 \documentclass[11pt]{article}
\usepackage{bbm}
 \usepackage{amssymb}
\usepackage{amssymb, amsthm, amsmath, amscd}
\usepackage[usenames,dvipsnames]{color}

\newtheorem{thm}{Theorem}[section]
\newtheorem{lem}[thm]{Lemma}
\newtheorem{prop}[thm]{Proposition}
\newtheorem{cor}[thm]{Corollary}
\theoremstyle{definition}
\newtheorem{NN}[thm]{}
\theoremstyle{definition}
\newtheorem{df}[thm]{Definition}
\theoremstyle{definition}
\newtheorem{rem}[thm]{Remark}
\newtheorem{nota}[thm]{Notation}
\theoremstyle{definition}

\renewcommand{\phi}{\varphi}

\newcommand{\N}{\mathbb{N}}
\newcommand{\Z}{\mathbb{Z}}
\newcommand{\Q}{\mathbb{Q}}
\newcommand{\R}{\mathbb{R}}
\newcommand{\C}{\mathbb{C}}
\newcommand{\T}{\mathbb{T}}

\numberwithin{equation}{section}

\newcommand{\Aff}{\operatorname{Aff}}
\newcommand{\Inf}{\operatorname{Inf}}

\newcommand{\id}{\operatorname{id}}

\newcommand{\aff}{\rm aff}

\newcommand{\morp}{contractive completely positive linear map}

\newcommand{\hm}{homomorphism}
\newcommand{\dt}{\delta}
\newcommand{\ep}{\varepsilon}

\newcommand{\sg}{\sigma}
\newcommand{\td}{\tilde}

\newcommand{\lr}{\longrightarrow}
\newcommand{\ld}{\lambda}
\newcommand{\cd}{\cdots}

\newcommand{\qq}{{\quad \quad}}
\newcommand{\sbs}{{\subset}}
\newcommand{\sps}{{\supset}}

\newcommand{\tht}{\theta}

\newcommand{\bb}{\mathfrak{b}}
\newcommand{\cc}{\mathfrak{c}}
\newcommand{\dd}{\mathfrak{d}}
\newcommand{\LD}{\Lambda}
\newcommand{\GM}{\Gamma}
\newcommand{\gm}{\gamma}
\newcommand{\sm}{\sigma}

\newcommand{\DT}{\Delta}

\newcommand{\spd}{(\diamondsuit)}
\newcommand{\spdd}{(\diamondsuit\diamondsuit)}
\newcommand{\spddd}{(\diamondsuit\diamondsuit\diamondsuit)}

\newcommand{\Tht}{\Theta}

\newcommand{\e}{\mbox{\large \bf 1}}
\newcommand{\0}{\mbox{\large \bf 0}}
\newcommand{\la}{\langle}
\newcommand{\ra}{\rangle}
\newcommand{\andeqn}{\,\,\,{\rm and}\,\,\,}
\newcommand{\rforal}{\,\,\,{\rm for\,\,\,all}\,\,\,}
\newcommand{\CA}{$C^*$-algebra}
\newcommand{\SCA}{$C^*$-subalgebra}

\newcommand{\af}{{\alpha}}
\newcommand{\bt}{{\beta}}
\newcommand{\dist}{{\rm dist}}

\newcommand{\diag}{{\rm diag}}

\newcommand{\wilog}{without loss of generality}
\newcommand{\Wlog}{Without loss of generality}

\newcommand{\D}{\mathbb D}
\newcommand{\beq}{\begin{eqnarray}}
\newcommand{\eneq}{\end{eqnarray}}
\newcommand{\tforal}{\,\,\,\text{for\,\,\,all}\,\,\,}
\newcommand{\tand}{\,\,\,\text{and}\,\,\,}

\usepackage{amsfonts}
\usepackage{mathrsfs}
\usepackage{textcomp}
\usepackage[all]{xy}

\title{A classification of finite simple amenable ${\cal Z}$-stable C*-algebras}
\author{Guihua Gong,  Huaxin Lin and Zhuang Niu
 }
\date{
}
\begin{document}

\maketitle

\begin{abstract}
%Let $\mathcal{S}$ be the class of all unital simple separable amenable stably finite Jiang-Su stable \CA s.
%We classify a class $\mathcal{A}$ of \CA s contained in $\mathcal{S}$ using the Elliott invariant. The class $\mathcal{A}$ is maximal in the sense that the Elliott invariant of every algebra in $\mathcal{S}$ is isomorphic to the Elliott invariant of some algebra in $\mathcal{A}$. The class $\mathcal{A}$ also contains the largest class of such algebras classified so far, namely those in $\mathcal{S}$ that have finite rational tracial rank and satisfy the Universal Coefficient Theorem (UCT).

We present a classification theorem for
a class of unital simple separable amenable ${\cal Z}$-stable  $C^*$-algebras by the Elliott invariant.
This class of simple \CA s exhausts all possible values of the Elliott invariant for unital stably finite simple separable amenable ${\cal Z}$-stable \CA s.  Moreover,  it contains all unital simple separable  amenable \CA s which satisfy the UCT and have finite  rational tracial rank.
% at most one.%generalized tracial rank one or zero.

\end{abstract}

\section{Introduction}
%%%%
%I will put my version of the introduction here -msun
The concept of a \CA\, exists harmoniously in many areas of mathematics. The original definition of \CA s axiomatized the norm closed self-adjoint subalgebras of $B(H),$ the algebra of all bounded linear operators on a Hilbert space $H.$ Thus \CA s are operator algebras. The study of \CA s may also be viewed as the study of a non-commutative analogue of topology. This is because every unital commutative \CA\, is isomorphic to $C(X),$ the algebra of continuous functions on a compact Hausdorff space $X$ (by the Gelfand transform).
%with values in the field of complex numbers $\mathbb{C}$. (
%The isomorphism is given by the Gelfand transform).
If we take our space $X$ and equip it with a group action via homeomorphisms, we enter the realm of topological dynamical systems, where remarkable progress has been made by considering the transformation \CA\, $C(X)\rtimes G$ arising from the crossed product construction. Analogously, we can consider the study of general crossed products as the study of non-commutative topological dynamical systems. Then there is also  naturally appeared the group \CA\, of $G$, which is fundamental in the study of abstract harmonic analysis.
\CA\, theory is also closely related to the non-commutative geometry.  Since von Neumann algebras are \CA s,
there are deep and inseparable interactions between \CA\, theory and von Neumann algebra theory.
The list goes on. Therefore, naturally,  it would be of great interest to classify \CA s. We have already seen some breathtaking advancements.

Early classification theorems include the work of Glimm in the 1950s who classified the Uniformly Hyperfinite algebras (UHF-algebras) by supernatural numbers.  We then saw Dixmier's and Brattelli's work on matroid \CA s and AF-algebras respectively. G. A. Elliott gave a  classification of AF-algebras by dimension groups in 1976 (\cite{Ell-76}) using what is  now known as Elliott's intertwining argument. By 1989,   Elliott had began his classification program by classifying simple A$\T$-algebras of real rank zero by scaled ordered K-theory.  Since then there has been rapid progress in the program to classify separable amenable simple \CA s now known as the Elliott program.  Elliott-Gong (\cite{EG-RR0AH}) and Elliott-Gong-Li (\cite{EGL-AH}) (together with a reduction theorem by Gong (\cite{Gong-AH})) classified the unital simple AH-algebras with no dimension growth by the Elliott invariant (see \ref{DEll} below).

The Elliott's  intertwining argument gave a framework to carry out further classification and focused one's efforts on studying the maps between algebras and invariants of certain building blocks (for AF algebras these building blocks would be finite dimensional algebras and for AH algebras these building blocks would be  certain homogeneous \CA s). In particular, one wants to know when maps between invariants are induced by maps between building blocks (sometimes referred to as an existence theorem) and to know when maps between the building blocks are approximately unitarily equivalent (often called an uniqueness theorem).   To classify \CA s without assuming some particular inductive limit structure, one would like to
establish abstract existence theorems
and uniqueness theorems.   These efforts became the engine for these
rapid
developments
(\cite{Lnjotuni}, \cite{LinTAF2}, \cite{LnAUCT} and  \cite{DE} for example). Both existence theorems and uniqueness theorems
used $KL$-theory  (\cite{Ror-KL-I}) and the total $K$-theory ($\underline{K}(A)$) developed by
Dadarlat and Loring (\cite{DL}).   These existence and uniqueness theorems provide not only the technical tools for the classification program but also provide the foundation to understand the morphisms in the category of \CA s.

 The rapid  developments mentioned above includes  the Kirchberg-Phillips (\cite{Kirch-Infty}, \cite{KP0}, \cite{Ph1})  classification of purely infinite simple separable amenable \CA s which satisfy the UCT by their $K$-theory.  There is also the classification of unital simple amenable C*-algebras in the UCT class which have finite tracial rank
%by Lin
(\cite{LinTAF2}, \cite{Lnduke} , \cite{LinTAI} and \cite{Lin-LAH}).

On the other hand,  it had been suggested in  \cite{DNNP-AH} and \cite{BDR} that unital simple AH-algebras without a dimension growth condition might behave differently.  It was Villadsen (\cite{Vill-perf} and \cite{Vill-sr}) who showed that unital simple AH-algebras may have perforated $K_0$-groups and may have stable rank
of any  positive integer values.  M. R{\o}rdam exhibited a unital separable simple C*-algebra which is finite but not stably finite
 (\cite{Ror-infproj}).
 %C*-algebra with an infinite projection and a finite projection? I have an impression that finite but not stably finite is something known earlier, but I need to have a look.)}.
It was Toms (\cite{Toms-Ann}) who showed that there are unital simple AH-algebras of stable rank one with the same
Elliott invariant that are not isomorphic.
Before that,  Jiang-Su (\cite{JS}) constructed a unital simple ASH-algebra ${\cal Z}$ of stable rank one which has the same Elliott invariant as that of $\C.$ In particular,  ${\cal Z}$ has no non-trivial projections. If $A$ is a unital separable
amenable simple \CA\, which belongs to a classifiable class, then one should expect $A\otimes {\cal Z}\cong A,$ since
they have the same $K$-theory.
Such \CA s are called ${\cal Z}$-stable \CA s. The existence of non-elementary simple $C^*$-algebras which were not ${\cal Z}$ stable was first proved  by
Gong-Jiang-Su (see \cite{GJS}).  Toms' counterexample is in particular not ${\cal Z}$-stable. Thus ${\cal Z}$-stability
should be added to
%was added  to our
into the
hypotheses if one uses the conventional Elliott invariant.  (The classes already classified up to that  time can easily be shown to be ${\cal Z}$-stable.)

% of Elliott's program.

%but had the shortcoming of not containing ${\cal Z}$ itself.
The next development in this direction came from a new approach of W. Winter who made use of the ${\cal Z}$-stability assumption in a  remarkably  innovative way (\cite{Winter-Z}). His idea was to view $A\otimes {\cal Z}$ as a (stationary) inductive limit of algebras of paths in $A\otimes Q$ with endpoints in $A\otimes M_{\bf p}$ and $A\otimes M_{\bf q}$ (where $Q$ is the UHF-algebra with its $K_0(Q)=\Q,$ ${\bf p}$ and ${\bf q}$ are coprime supernatural numbers and $M_{\bf p}$ and $M_{\bf q}$ their associated UHF algebras).
Suppose that the endpoint algebras are classifiable.
Winter  showed that  if somehow there is a continuous path of  isomorphisms
from one endpoint to the other,
%connected continuously and the class of all $A\otimes M_{\bf p}$ could be classified,
then   the class of all $A$ can  also be classified.
Winter's procedure  provided a new framework to carry out classification, however to actually execute the continuation from endpoint to endpoint alluded to above, one has to establish new types of uniqueness and existence theorems.
 In other words,  just like before, in Elliott's intertwining argument,  the new procedure  ultimately, but not surprisingly,
 depends on some existence and uniqueness theorems. However, this time we need the uniqueness and existence theorems
  with respect to asymptotic unitary equivalence of the maps involved rather approximate unitary equivalence, which is significantly more demanding. For example, in the case of the existence theorem, we need to construct maps which lift a prescribed $KK$-element
rather than a $KL$-element, which was once thought  to be out of reach for general stably finite algebras since the $KK$-functor does not preserve inductive limits. It was an unexpected usage of the Basic Homotopy Lemma that made this possible.  Moreover, the existence theorem also needs to respect the prescribed rotation related map.  The existence
 theorems
are very different from ones in the  early study. Inevitably,  the uniqueness theorem also
becomes more complicated (again the Basic Homotopy Lemma plays the key role).
Once  we overcame these new hurdles,
 the class of all ${\cal Z}$-stable \CA s $A$ whose tensor products with all UHF-algebras of infinite type are of tracial rank zero were shown to be classified by the Elliott invariant in \cite{Winter-Z}, \cite{Lin-App} and
 \cite{L-N}. This class was then expanded to the class ${\cal B}$ of unital simple separable amenable
 ${\cal Z}$-stable  \CA s which satisfy the UCT and
 have finite tracial rank after tensoring some
 %such that $A\otimes U$ have  finite tracial rank for all
infinite dimensional UHF-algebra.  It is shown in \cite{Lnclasn}
%(see also \cite{L-N} for the class
%of those $A$ such that $A\otimes U$ has tracial rank at most one)
that \CA s in the class ${\cal B}$ can be classified up to
isomorphism by the Elliott invariant.  These constitute a significant expansion of the class of \CA s which can be classified because there are a great amount of unital simple \CA s which do not have finite tracial rank whose tensored products with a UHF-algebra $U$ do. For example the Jiang-Su algebra ${\cal Z}$ is projectionless, but ${\cal Z}\otimes U\cong U$ for any infinite dimensional UHF-algebra $U.$ In fact, the class ${\cal B}$ exhausts all those Elliott invariant with weakly unperforated
simple rational Riesz groups as $K_0$-groups under the restriction that the maps from tracial state spaces to state spaces of $K_0$
map the extremal points to extremal points (\cite{LNjfa}).  The class ${\cal B}$ not only contains all unital simple separable amenable
\CA s with finite tracial rank in the UCT class  and  the Jiang-Su algebra but also contains many other
simple \CA s.   More importantly it unifies the previously classified classes such as the so-called dimension drop algebras as well as those dimension drop circle algebras which were known not to be AH-algebras (\cite{Lnjfa2010}). However the restriction on the parings of tracial state spaces and state spaces of $K_0$-groups prevents the class ${\cal B}$ from including those inductive limits of  so-called ``point-line" algebras, which we called the Elliott-Thomsen building blocks and leads us to the main goal of our paper.

The goal of this article is to give a classification of a certain class of unital simple separable amenable \CA s satisfying the UCT by the Elliott invariant. This class is significant because it exhausts all possible values of the Elliott invariant for all unital separable simple
amenable stably finite ${\cal Z}$-stable \CA s and it also contains the class ${\cal B}$ properly.
We introduce a class of separable simple \CA s which will be called simple \CA s of \emph{generalized} tracial rank at most one.
The definition follows the same spirit as that of tracial rank one  (or zero). But, instead of using only finite direct sums of
matrix algebras of continuous functions on a one-dimensional finite CW complex, we use certain C*-subalgebras of interval algebras. These were first introduced into the Elliott program by Elliott and Thomsen (\cite{ET-PL}), sometimes also called
one dimensional non-commutative  CW complexes (NCCW), as a model to approximate \CA s (tracially).  This class will be denoted by ${\cal B}_1$ (see \ref{DB1} below).
Denote by ${\cal N}_1$  the family of unital separable amenable simple \CA s  $A$ which satisfy the UCT such that $A\otimes Q\in {\cal B}_1,$  and by ${\cal N}_0$ the subclass of those \CA s $A$ such that
$A\otimes Q\in {\cal B}_0,$ where $Q$ is the UHF-algebra with $K_0(Q)=\Q$. %$(K_0(Q), K_0(Q)_+, [1_Q]=(\Q,\Q_+,1).$
The main theorem of this article has two parts. The first part states (see \ref{MFTh}) that
if $A$ and $B$ are two unital separable simple ${\cal Z}$-stable \CA s in ${\cal N}_1.$ Then
$A\cong B$ if and only if ${\rm Ell}(A)\cong {\rm Ell}(B)$ (see the definition \ref{DEll} below).  As a consequence,
two unital \CA s in ${\cal B}_1$  which satisfy the UCT are isomorphic if and only if they have the same
Elliott invariant. The second part of the main theorem states that, given any  countable weakly unperforated simple ordered
group $G_0$ with order unit $e,$ any countable abelian  group $G_1,$ and any metrizable Choquet simple $T$ and
any surjective continuous affine map $s: T\to S_e(G_0)$ (the state space of $G_0$), there exists a unital \CA\, $A\in {\cal N}_1$ such that
$$
{\rm Ell}(A)=(K_0(A), K_0(A)_+, [1_A], K_1(A), T(A), r_A)=(G_0, (G_0)_+, e, G_1, T, s).
$$

The article is organized as follows.
Section 2 serves  as  preliminaries for the article which contains a number of conventions.
In Section 3, we study the class of Elliott-Thomsen building blocks, denoted by $\mathcal C$, a class of unital C*-subalgebras of finite direct sums of interval algebras. We also use ${\cal C}_0$ for the subclass of \CA s in ${\cal C}$ with trivial $K_1.$
The Elliott-Thomsen building blocks are also called point-line algebras, or one dimensional non-commutative CW complexes as studied in \cite{ELP1} and \cite{ELP2}. %We  compute the exponential rank  and exponential length  of unitaries in
%the closure of commutator subgroups as well as ordered $K$-theory  of \CA s in ${\cal C}.$ Other properties are also
%presented.
Section 4 and 5 discuss the uniqueness theorem for maps from \CA s in ${\cal C}$ to finite dimensional
\CA s.  Section 8 presents a uniqueness theorem for maps from a C*-algebra in ${\cal C}$ to another C*-algebra in ${\cal C}.$ This is done by using
a homotopy lemma established in  Section 6 and existence theorems established in Section 7 to bridge uniqueness theorems in Section 4 and 5
to ones in Section 8.  In Section 9, the classes ${\cal B}_1$ and ${\cal B}_0$ are introduced.   Properties  of \CA s in
class ${\cal B}_1$
are discussed in Section 9, 10 and 11. \CA s in ${\cal B}_1$ are simple and are of  generalized tracial rank one (or zero).
These unital simple \CA s may also be characterized as tracially approximated by subhomogeneous \CA s
with one dimensional spectra.
%SCA s of $C([0,1],F)$ where $F$ is a finite
%dimensional \CA s.
%in ${\cal C}$  (and ${\cal C}_0,$
%a subclass of \CA s in ${\cal C}$ with trivial $K_1$-groups).
We show, for example, in Section 10, they are ${\cal Z}$-stable.
Section 12 is dedicated to the main uniqueness theorem for \CA s in ${\cal B}_0$ used in the isomorphism theorem in Section 21. Sections 13 and 14 are devoted to the range theorem. It includes one of the main results:
%Let  ${\cal N}_1$  be the class of  simple unital amenable \CA s  which satisfy the UCT and whose tonsored products with UHF-algebras of infinite type
%are in ${\cal B}.$
 We show in Section 13 that, given any six-tuples of possible Elliott invariant
for unital separable simple ${\cal Z}$-stable \CA s,
there is a unital simple ${\cal Z}$-stable \CA\, in ${\cal N}_1$ whose Elliott invariant is exactly as the given one.  Section 14
gives a similar construction for a set of restricted Elliott invariants. In this case, simple \CA s constructed
are inductive limits of \CA s which are finite direct sums of some homogeneous \CA s and \CA s in ${\cal C}_0.$ This subclass plays an important role in this article.
It should  be pointed out, however, that there are unital simple ${\cal Z}$-stable \CA s in ${\cal N}_1$ which can not be written as inductive limits
of finite direct sums of homogeneous \CA s and \CA s in ${\cal C}_0$ (or in ${\cal C}_1$). Section 15 to 19 could all be described as parts of existence theorems. These deal with
the issues of existence  for maps from \CA s in ${\cal C}$ to finite dimensional \CA s and then to
\CA s in ${\cal C}$ which
match the prescribed $K_0$-maps and tracial information. Ordered structure and combined simplex information
of these \CA s
become complicated.  We also need to consider maps
from homogeneous \CA s to \CA s in ${\cal C}.$ The mixture with higher dimensional CW complexes does not ease the difficulties. However, in Section 18, we show that,  at least under certain restrictions, any given compatible triple which consists of a strictly positive $KL$-element,
a map on the tracial state space and a \hm\, on a quotient of the unitary group,  it is possible to construct a
\hm\, between \CA s in ${\cal B}_0$ which matches the triple. Variations of this are also discussed.
In Section 19, we show that ${\cal N}_1={\cal N}_0.$

In Section 20, we continue to study the existence theorem.
In Section 21, we show that any  unital simple \CA s with the form $A\otimes U$ with $A\in {\cal B}_0$ is isomorphic
to a unital \CA s constructed in Section 14, and two such \CA s are isomorphic  if  (and only if)
they have the same Elliott invariant.  This isomorphism theorem is special but  is also    the foundation of our main
isomorphism theorem in Section 29.  In the next  seven sections, we show that Winter's strategy may be carried out
in the general case.    This  requires   a much more sensitive uniqueness and existence theorem.
The uniqueness theorem now requires an asymptotic unitary equivalence theorem which is proved in Section 27.
 To do this, we first need another Basic Homotopy Lemma. Sections 22, 24 and 25 establish the needed homotopy lemma while Section 22 serves as the existence
theorem for the homotopy lemma.
%The asymptotic unitary equivalence theorem is established in Section 26.
Sections 26 and 28 are for the rotation maps and existence theorem.
 In Section 29, we show that two  ${\cal Z}$-stable \CA s in ${\cal N}_0$ are isomorphic
if and only if they have the same Elliott invariant.
Since we have shown in section 19 that
${\cal N}_0={\cal N}_1.$ This completes the proof of the main isomorphism theorem for  all ${\cal Z}$-stable
\CA s in ${\cal N}_1.$

%{\done yet}

{\bf Acknowledgements}:
A large part of this article were written during the summers  of 2012, 2013 and 2014
when all three authors visited The Research Center for Operator Algebras in East China Normal University.
They were partially supported by the center. Both the first named author and  the second named author have been supported partially by NSF grants.
 The work of the third named author has been partially supported by a NSERC Discovery Grant, a Start-Up Grant from the University of Wyoming, and a Simons Foundation Collaboration Grant.
 Authors would like to thank Michael Yuan Sun for his efforts to read  part of the earlier version of this research and
for his comments.
% on introduction.

 \section{Preliminaries}

\begin{df}\label{Du}
{\rm
Let $A$ be a unital \CA. Denote by $U(A)$ the unitary group of $A$, and denote by  $U_0(A)$ the normal subgroup of $U(A)$ consisting of those unitaries which are in the  connected component  of $U(A)$ containing $1_A.$  Denote by
$DU(A)$ the commutator subgroup of $U_0(A)$ and $CU(A)$ the closure of $DU(A)$ in $U(A).$
}
\end{df}
\begin{df}\label{Aq}
{\rm Let $A$ be a unital \CA\, and let $T(A)$ denote the simplex of  tracial states of $A.$
Let $\tau\in T(A).$ We say that $\tau$ is faithful if $\tau(a)>0$ for all $a\in A_+\setminus\{0\}.$
Denote by $T_f(A)$ the set of all faithful tracial states.

Denote by $\Aff(T(A))$ the space of all real continuous affine functions on $T(A)$ and
denote by ${\rm LAff}_b(T(A))$  the set of all bounded lower-semi-continuous real affine functions on $T(A).$

Suppose that $T(A)\not=\emptyset.$ There is an affine  map
$r_{\aff}: A_{s.a.}\to \Aff(T(A))$ defined by
$$
r_{\aff}(a)(\tau)=\hat{a}(\tau)=\tau(a)\tforal \tau\in T(A)
$$
and for all $a\in A_{s.a.}.$ Denote by $A_{s.a.}^q$ the space  $r_{\aff}(A_{s.a.})$ and
$A_+^q=r_{\aff}(A_+).$

{For each integer $n\ge 1$ and $a\in M_n(A),$
write $\tau(a)=(\tau\otimes \mathrm{Tr})(a),$ where $\mathrm{Tr}$ is the (non-normalized) standard trace on $M_n.$}
}

\end{df}

\begin{df}\label{Drho}
Let $A$ be a unital stably finite \CA\, with $T(A)\not=\emptyset.$ Denote by $\rho_A: K_0(A)\to \Aff(T(A))$ the  order preserving \hm\,  defined by $\rho_A([p])=\tau(p)$ for any projection $p\in M_n(A),$
$n=1,2,...$ (see the above convention).
%{Let $A$ be a unital stably finite \CA.}
A map $s: K_0(A)\to \R$ is said to be a state if $s$ is an order preserving \hm\, such that
$s([1_A])=1.$ The set of states on $K_0(A)$ is denoted by $S_{[1_A]}(K_0(A)).$

{Denote by $r_A: T(A)\to S_{[1_A]}(K_0(A))$ the map defined by $r_A(\tau)([p])=\tau(p)$ for all projections
$p\in M_n(A)$ (for all integer $n$) and for all $\tau\in T(A).$ }
\end{df}

{
\begin{df}\label{DEll}
Let $A$ be a unital simple \CA. The Elliott invariant of $A$, denote by ${\rm Ell}(A),$ is the following
six tuple
$$
{\rm Ell}(A)=(K_0(A), K_0(A)_+, [1_A], K_1(A), T(A), r_A).
$$
Suppose that $B$ is another unital simple \CA. We write  ${\rm Ell}(A)\cong {\rm Ell}(B),$ if
there is an order isomorphism $\kappa_0: K_0(A)\to K_0(B)$
such that $\kappa_0([1_A])=[1_B],$ an isomorphism $\kappa_1:
K_1(A)\to K_1(B)$ and an affine homeomorphism
$\kappa_\rho: T(B)\to T(A)$ such that
$$
r_A(\kappa_\rho(t))(x)=r_B(t)(\kappa_0(x))\rforal x\in K_0(A)\andeqn \rforal t\in T(B).
$$
\end{df}
}

\begin{df}\label{Dball}
Let $X$ be a compact metric space, $x\in X$ be a point and let $r>0.$
Denote by $B(x, r)=\{y\in X: {\rm dist}(x, y)<r\}.$

Let $\ep>0.$ Define $f_\ep\in C_0((0,\infty))$ to be a function with
 $f_\ep(t)=0$ if $t\in [0,\ep/2],$ $f_\ep(t)=1$ if $t\in [\ep,\infty)$ and
$f_\ep(t)=2{t-\ep/2\over{\ep}}$  if $t\in [\ep/2, \ep].$ Note that$0\le f\le 1$ and $f_\ep f_{\ep/2}=f_\ep.$
\end{df}

\begin{df}\label{DW(A)}
Let $A$ be a \CA. Let $a, b\in M_n(A)_+$. Following Cuntz (\cite{Cu1}), we write
$a\lesssim b$ if there exists a sequence $(x_n)\subset M_n(A)$ such that
$\lim_{n\to\infty} x_n^*bx_n=a.$ If $a\lesssim b$ and $b\lesssim a,$ then
we write $a\sim b.$ The relation ``$\sim$" is an equivalence relation.
Denote by $W(A)$ the Cuntz semi-group of the equivalence classes of positive
elements in $\cup_{m=1}^{\infty} M_m(A)$ with orthogonal
addition (i.e., $[a+b]=[a\oplus b]$).

If $p, q\in M_n(A)$ are projections, then $p\lesssim q$ if and only if $p$ is Murray
%Munrray
-von Neumann equivalent to a subprojection of $q$. In particular,  when $A$ is stably finite, then $p\sim q$ if and only if they are von Neumann equivalent.
% {\color{Green} The last sentence is not true in general. See below.}

Denote by $QT(A)$ the set of normalized quasi-traces on $A.$
%Suppose that $\tau\in QT(A),$  then $\tau$ can be extended to a quasi-trace on $M_n(A)$ by
%define $\tau((a_{i,j})_{n\times n})=\sum_{i=1}^n \tau(a_{i,i})$ for all $(a_{i,j})_{n\times n}\in M_n(A).$
For $a\in A_+$ and $\tau\in QT(A),$ define
$$
d_{\tau}(a)=\lim_{\ep\to 0} \tau(f_{\ep}(a)).
$$
Suppose that $QT(A)\not=\emptyset.$
We say $A$ has strict comparison for positive elements, if for any $a, b\in M_n(A)$ (for all integer $n\ge 1$)
$d_\tau(a)<d_\tau(b)$ for $\tau\in QT(M_n(A))$  implies $a\lesssim b.$
\end{df}

%{\color{Green}
%\begin{rem}
%Let $A$ be a C*-algebra. If $p, q\in M_n(A)$ are two projections,
%then $p\lesssim q$ if and only if $p$ is Murray-von Neumann equivalent to a subprojection of $q$. But in general, $p\sim q$ does not imply $p$ is Muaary-von Neumann equivalent to $q$. For instance, all projections in $\mathcal O_\infty$ are Cuntz equivalent, but $K_0(\mathcal O_\infty)\cong\Z$.
%\end{rem}
%}

\begin{df}\label{Aq1}
Let $A$ be a \CA.  Denote by $A^{\bf 1}$ the unit ball of $A.$
$A_+^{q, {\bf 1}}$ is the image of the intersection of $A_+\cap A^{{\bf 1}}$ in $A_+^q.$
\end{df}

\begin{df}\label{dInn}
%{\rm
%{
Let $A$ be a unital C*-algebra and let $u\in U(A).$ We write
${\rm Ad}\, u$ for the automorphism $a\mapsto u^*au$ for all $a\in A.$
Suppose $B\subseteq A$ is a unital C*-subalgebra.
Denote by $\overline{{\rm Inn}}(B,A)$ the set of all those monomorphisms $\phi: B\to A$ such that
there exists a sequence of unitaries $\{u_n\}\subset {{A}}$ so that
$\phi(b)=\lim_{n\to\infty} u_n^*bu_n$ for all $b\in B.$
%}
%}
\end{df}

\begin{df}\label{dfcalN}
%{\rm {
Denote by ${\cal N}$ the class of separable amenable \CA s which satisfy the Universal Coefficient Theorem (UCT).

Denote by ${\cal Z}$ the Jiang-Su algebra (\cite{JS}). % of unital simple projectionless  \CA\,
Note that $\mathcal Z$ has a unique trace and $K_i({\cal Z})=K_i(\C)$ ($i=0,1$).
A C*-algebra $A$ is said to be ${\cal Z}$-stable if $A\cong A\otimes {\cal Z}.$
%{\color{green} (Should one insert simplicity?)}
%}
%}
\end{df}

\begin{df}\label{DKLtriple}
{\rm

Let $A$ be a unital C*-algebra. Recall that, following D\u{a}d\u{a}rlat and Loring (\cite{DL}), one defines
\begin{equation}\label{Dbeta-5}
\underline{K}(A)=\bigoplus_{i=0,1}K_i(A)\oplus\bigoplus_{i=0,1}
\bigoplus_{k\ge 2}K_i(A,\Z/k\Z).
\end{equation}
There is a commutative \CA\, $C_k$ such that one may identify
$K_i(A\otimes C_k)$ with $K_i(A, \Z/k\Z).$
{ Let $A$ be a unital separable amenable C*-algebra, and let $B$ be a $\sigma$-unital C*-algebra.
Following R\o rdam (\cite{Ror-KL-I}),  $KL(A,B)$ is the quotient of $KK(A,B)$ by those elements
represented by limits of trivial extensions (see \cite{LnAUCT}). In the case that $A$ satisfies the UCT,
R\o rdam defines $KL(A,B)=KK(A,B)/{\cal P},$ where
${\cal P}$ is the subgroup corresponding to the pure extensions of the $K_*(A)$ by $K_*(B).$
In \cite{DL}, D\u{a}d\u{a}rlat and Loring  proved that
\beq\label{DKL-2}
KL(A,B)={\rm Hom}_{\Lambda}(\underline{K}(A), \underline{K}(B)).
\eneq
}
{
Now suppose that $A$ is stably finite. Denote by $KK(A,B)^{++}$ the set of those elements
$\kappa\in KK(A,B)$
such that $\kappa(K_0(A)\setminus \{0\})\subset K_0(B)_+\setminus\{0\}.$
Suppose further that both $A$ and $B$ are unital. Denote by $KK_e(A,B)^{++}$ the subset of
those $\kappa\in KK(A,B)^{++}$ such that $\kappa([1_A])=[1_B].$
Denote by $KL_e(A,B)^{++}$ the image of $KK_e(A,B)^{++}$ in $KL(A,B).$ }

}
\end{df}

\begin{df}\label{Dmap}
Let $A$  and $B$ be two \CA s and let $\phi: A\to B$ be a map.
{{We will sometime, without warning, continue to use $\phi$ for the induced
map $\phi\otimes {\rm id}_{M_n}: A\otimes M_n\to B\otimes M_n.$
Also $\phi\otimes 1_{M_n}: A\to B\otimes M_n$ is  used for the amplification
which maps $a$ to $a\otimes 1_{M_n},$ the diagonal element with $a$ repeated $n$ times.}}
Throughout the paper,  if $\phi$ is a \hm, we will use $\phi_{*i}: K_i(A)\to K_i(B)$, $i=0, 1$, for the induced \hm.
We will use $[\phi]$ for the element in $KL(A,B)$ (or $KK(A,B)$ if there is no confusion) which is induced by $\phi.$
Suppose that $T(A)\not=\emptyset$ and $T(B)\not=\emptyset.$ Then $\phi$ induces
an affine map $\phi_T: T(B)\to T(A)$ defined by $\phi_T(\tau)(a)=\tau(\phi(a))$ for all $\tau\in T(B)$ and
$a\in A_{s.a.}.$ Denote by $\phi^{\sharp}: \Aff(T(A))\to \Aff(T(B))$ the affine continuous map defined by
$\phi^{\sharp}(f)(\tau)=f(\phi_T(\tau))$ for all $f\in \Aff(T(A))$ and $\tau\in T(B).$

\end{df}

%Let $\dt>0$ be a positive number, ${\cal F}\subset A$ be a finite subset and let
%${\cal G}\subset \underline{K}(A)$ be a finite.

\begin{df}\label{KLtriple}
Let $A$ be a unital separable amenable \CA\, and let  $x\in A.$ Suppose
that
$\|xx^*-1\|<1$ and $\|x^*x-1\|<1.$ Then $x|x|^{-1}$ is a unitary.
Let us use $\langle x \rangle $ to denote $x|x|^{-1}.$

Let ${\cal F}\subset A$ be a finite subset and $\ep>0$ be a positive number.
We say a map $L: A\to B$ is ${\cal F}$-$\ep$-multiplicative if
$$
\|L(xy)-L(x)L(y)\|<\ep\rforal x,\, y\in {\cal F}.
$$

Let ${\cal
P}\subset \underline{K}(A)$ be a finite subset. There is $\ep>0$ and a finite subset
${\cal F}$ satisfying the following: for any unital \CA\, $B$ and
any unital ${\cal F}$-$\ep$-multiplicative \morp\, $L : A\to B,$ $L$
induces a \hm\, $[L]$ defined on  $G({\cal P}),$ where $G({\cal
P})$ is the subgroup generated by ${\cal P},$ to $\underline{K}(B)$
such that
\beq\label{KLtriple-1}
\|L(p)-q\|<1\andeqn \|\langle
L(u)\rangle -v\|<1,
\eneq
where $q\in M_n(B)$ (for some $n\ge 1$) is a projection such that $[q]=[L]([p])$ in $K_0(B)$  and
$v\in M_n(B)$  is a unitary such that
$[v]=[L]([u])$ for all $[p]\in {\cal P}\cap K_0(A)$ and $[u]\in
{\cal P}\cap K_1(A).$ This also applies to ${\cal P}\cap
K_i(A, \Z/k\Z)$ with a necessary modification, by replacing $L,$ by
$L\otimes {\rm id}_{C_k},$ where $C_k$ is the commutative C*-algebra described in \ref{DKLtriple}, for
example. Such a triple $(\ep, {\cal F}, {\cal P})$ is called a
$KL$-triple for $A.$

Suppose that $K_i(A)$ is finitely generated.
Then, by \cite{DL}, there  is  $n_0\ge 1$ such that every element
$\kappa\in \mathrm{Hom}_{\Lambda}(\underline{K}(A), \underline{K}(B))$ is
determined by $\kappa|_{G_{n_0}},$ where
$G_{n_0}=\bigoplus_{i=0,1}K_i(A) \oplus \bigoplus_{i=0,1}\bigoplus_{k=2}^{n_0}K_i(A,\Z/k\Z).$
 %{\color{green} (``eqnarray" and ``equation" have different spacings. )}
 Therefore, for some large ${\cal P},$  if $(\ep, {\cal F}, {\cal P})$ is a $KL$-triple for $A,$
 then $[L]$ defines an element in $KL(A,B)=KK(A,B).$ In this case, we say $(\ep, {\cal F})$ is
 a $KK$-pair.
 \end{df}

\begin{df}\label{Dbeta}
{\rm Let $A$ be a unital \CA. Consider the tensor product
$A\otimes C(\T).$
By the K\"{u}nneth Theorem, the tensor product induces two canonical injective \hm s:
\begin{equation}\label{Dbeta-1}
\bt^{(0)}: K_0(A)\to K_1(A\otimes C(\T))\quad\mathrm{and}\quad
\bt^{(1)}: K_1(A)\to K_0(A\otimes C(\T)).
\end{equation}
In this way, one may write
\begin{equation}\label{Dbeta-2}
K_i(A\otimes C(\T))=K_i(A)\oplus \bt^{(i-1)}(K_{i-1}(A)),\,\,\,i=0,1.
\end{equation}
For each $i\ge 2,$ one also obtains the following injective \hm s:
\begin{equation}\label{Dbeta-3}
\bt_k^{(i)}: K_i(A, \Z/k\Z)\to K_{i-1}(A\otimes C(\T),\Z/k\Z),\,\,\,i=0, 1.
\end{equation}
Thus one may write
\begin{equation}\label{Dbeta-4}
K_i(A\otimes C(\T),\Z/k\Z)=K_i(A,\Z/k\Z)\oplus
\bt_k^{(i-1)}(K_{i-1}(A,\Z/k\Z)),\,\,\,i=0,1.
\end{equation}

%Recall, following Dadarlat and Loring, one write that
%\beq\label{Dbeta-5}
%\underline{K}(A)=\bigoplus_{i=0,1}K_i(A)\bigoplus\bigoplus_{i=0,1}
%\bigoplus_{k\ge 2}K_i(A,\Z/k\Z).
%\eneq

If $x\in \underline{K}(A),$ we use  ${\boldsymbol{\bt}}(x)$ for $\bt^{(i)}(x)$
if $x\in K_i(A)$ and
$\bt_k^{(i)}(x)$ if $x\in K_i(A,\Z/k\Z).$ So one has an injective \hm\,
\beq\label{Dbeta-6}
{\boldsymbol{\bt}}: \underline{K}(A)\to \underline{K}(A\otimes C(\T))
\eneq
and writes
\beq\label{Dbeta-7}
\underline{K}(A\otimes C(\T))=\underline{K}(A)\oplus {\boldsymbol{\bt}}(\underline{K}(A)).
\eneq

Let $h: A\otimes C(\T)\to B$ be a unital \hm. Then
$h$ induces a \hm\, $h_{*i,k}: K_i(A\otimes C(\T),\Z/k\Z)\to K_i(B,\Z/k\Z),$
$k=0,2,3,...$ and $i=0,1.$
Suppose that $\phi: A\to B$ is a unital \hm\, and $v\in U(B)$ is a unitary
such that $\phi(a)v=v\phi(a)$ for all $a\in A.$
Then $\phi$ and $v$ induce a unital \hm\,
$h: A\otimes C(\T)\to B$ by $h(a\otimes z)=\phi(a)v$ for all $a\in A,$ where
$z\in C(\T)$ is the identity function on the unit circle $\T.$
We use ${\rm Bott}(\phi,\, v)$ for all \hm s $h_{*i-1,k}\circ \bt_k^{(i)}$ and we write
\beq\label{Dbeta-8}
{\rm Bott}(\phi, \,v)=0
\eneq
if $h_{*i,k}\circ \bt_k^{(i)}=0$ for all $k$ and $i.$  %It seems helpful to point out that
In particular, since $A$ is unital, \eqref{Dbeta-8} implies that $[v]=0$ in $K_1(B)$.
We also use ${\rm bott}_i(\phi, \, v)$ for
$h_{*i-1}\circ \bt^{(i)},$ $i=0,1.$

Suppose that  $A$ is a unital separable amenable \CA, $\phi: A\to B$ is a \hm\, and $v\in B$ is a unitary.
For any $\ep>0$ and  any finite subset ${\cal F}\subset A,$ there is $\dt>0$ and a finite subset ${\cal G}\subset A$ such that if
\beq\label{Dbeta-9}
\|[\phi(g), \, v]\|<\dt\rforal f\in {\cal G},
\eneq
then, by 2.8 of \cite{LnHomtp}, there exists a unital ${\cal F}$-$\ep$-multiplicative \morp\, $L: A\otimes C(\T)\to B$ such that
\beq\label{Dbeta-10}
\|L(f)-\phi(f)\|<\ep\rforal f\in {\cal F}\andeqn
\|L(1\otimes z)-v\|<\ep,
\eneq
where $z\in C(\T)$ is the standard unitary generator of $C(\T).$
Therefore, for each finite subset ${\cal Q}\subset \underline{K}(A\otimes C(\T)),$ there is $\dt>0$ and a finite subset ${\cal G}$ such that,  when (\ref{Dbeta-9}) holds,
$[L]|_{{\cal Q}}$ is well defined.
Let ${\cal P}\subset \underline{K}(A)$ be a finite subset.
There are $\dt_{\cal P}>0$ and a finite subset ${\cal F}_{\cal P}$ satisfying the following: if (\ref{Dbeta-9}) holds for $\dt_{\cal P}$ (in place of $\dt$) and ${\cal F}_{\cal P}$ (in place of ${\cal F}$),
then
$[L]|_{{\boldsymbol{\bt}}({\cal P})}$ is well defined. In this case, we will write
\beq\label{Dbeta-11}
{\rm Bott}(\phi,\, v)|_{\cal P}(x)=[L]|_{{\boldsymbol{\bt}}({\cal P})}(x)
\eneq
for all $x\in {\cal P}.$  In particular,
when
$
[L]|_{{\boldsymbol{\bt}}({\cal P})}=0,
$
we will write
\beq\label{Dbeta-12}
{\rm Bott}(\phi,\, v)|_{\cal P}=0.
\eneq
When $K_*(A)$ is finitely generated, $KL(A, B)=\mathrm{Hom}_{\Lambda}(\underline{K}(A), \underline{K}(B))$ is determined
by a finitely generated subgroup of $\underline{K}(A)$ (see \cite{DL}).
Let ${\cal P}$ be a finite subset which generates this subgroup.
Then, in this case, instead of (\ref{Dbeta-13}), we may write
that
\beq\label{Dbeta-13}
{\rm Bott}(\phi, \, v)=0.
\eneq
In general, if ${\cal P}\subset K_0(A),$
we will write
\beq\label{Dbeta-14}
{\rm bott}_0(\phi, \,v)|_{\cal P}={\rm Bott}(\phi,\, v)|_{\cal P}
\eneq
and if ${\cal P}\subset K_1(A),$
we will write
\beq\label{Dbeta-15}
{\rm bott}_1(\phi,\, v)|_{\cal P}={\rm Bott}(\phi, \, v)|_{\cal P}.
\eneq
}
\end{df}

\begin{df}\label{expLR}
{\rm
Let $A$ be a unital \CA.  Each element $u\in U_0(A)$ can be written as  $u= e^{ih_1}e^{ih_2}\cdots e^{ih_k}$ for $h_1, h_2,..., h_k\in A_{s,a}$.
We write that ${\rm cer}(u)\le k$ if $u=e^{ih_1}e^{ih_2}\cdots e^{ih_k}$ for some  selfadjoint elements $h_1, h_2, ..., h_k.$
We  write that ${\rm cer}(u)=k$ if ${{\rm cer}}(u)\le k$ and
$u$ is {\it not} a norm limit of unitaries $\{u_n\}$ with $cer(u_n)\le k-1.$
We write ${\rm cer}(u)=k+\ep$ if ${\rm cer}(u)\not=k$ and there exists a sequence of unitaries
$\{u_n\}\subset A$ such that $u_n\in U_0(A)$ with ${\rm cer}(u_n)\le k.$

%Define the exponential rank of $u$ (denoted by $cer(u)$) to be $k$ if $u= e^{ih_1}e^{ih_2}\cdots e^{ih_k}$ for some $h_1, h_2, ..., h_k$, and there is a $\dt>0$ such that for any $v\in U_0(A)$ with $\|u-v\|<\dt$, there are no $h'_1, h'_2,..., h'_{k-1}\in A_{s,a}$ satisfying $v=e^{ih'_1}e^{ih'_2}\cdots e^{ih_{k-1}}$.  {\color{Green} Is it the same to say $\mathrm{cer}(u)$ is the smallest $k$ such that $u= e^{ih_1}e^{ih_2}\cdots e^{ih_k}$ for some $h_1, h_2, ..., h_k$?}

%The exponential rank of $u$ is said to be $k+\ep$ if for any $\dt>0$, there are $h_1, h_2,\cdots, h_k\in A_{s,a}$ such that
%$$\|u-e^{ih_1}e^{ih_2}\cdots e^{ih_k}\|<\dt$$
%but there are no $h'_1, h'_2,..., h'_{k}\in A_{s,a}$ satisfying $$u=e^{ih'_1}e^{ih'_2}\cdots e^{ih'_{k}}.$$
Define exponential length of $u$ by
$${\rm cel}(u)=\inf{\big{\{}}\mbox{length of } (u(t))_{0\leq t\leq 1}~|~~ u(t)\in U_0(A),~u_0=u,~ u_1=1{\big{\}}}.$$
Obviously if $u= e^{ih_1}e^{ih_2}\cdots e^{ih_k}$, then (see \cite{Ringrose-cel})
$${\rm cel}(u) \leq \|h_1\|+\|h_2\|+\cdots + \|h_k\|~.$$}
\end{df}

\begin{df}\label{Dcu}
%{\rm
Suppose that $A$ is a unital \CA\, with $T(A)\not=\emptyset.$
Recall that $CU(A)$ is the closure of commutator subgroup of $U_0(A).$
Let $u\in U(A).$  We use
${\bar u}$ for the image in $U(A)/CU(A).$  It was proved in
\cite{Thomsen-rims} that there is
a splitting short exact sequence:
\begin{equation}\label{Dcu-1}
0\to \Aff(T(A))/\overline{\rho_A(K_0(A))}\to \bigcup_{n=1}^{\infty}U(M_n(A))/CU(M_n(A))\to
K_1(A)\to 0.
\end{equation}
Let $J_c$ be a  fixed  splitting map.  Then one may write
\begin{equation}\label{Duc-2}
 \bigcup_{n=1}^{\infty}U(M_n(A))/CU(M_n(A))=\Aff(T(A))/\overline{\rho_A(K_0(A))}\oplus J_c(K_1(A)).
\end{equation}
If $A$ has stable rank $k,$ then $K_1(A)=U(M_k(A))/U_0(M_k(A)).$
Note
$${\rm csr}(C(\T, A))\le tsr(A)+1=k+1.$$
%{\color{green} (I changed ``str" to ``sr")}
It follows from Theorem 3.10  of \cite{GLX-ER} that
\begin{equation}\label{Dcu-3}
\bigcup_{n=1}^{\infty}U_0(M_n(A))/CU(M_n(A))=U_0(M_k(A))/CU(M_k(A)).
\end{equation}
Then one has the following splitting short exact sequence
\begin{equation}\label{Dcu-4}
0\to \Aff(T(A))/\overline{\rho_A(K_0(A))}\to U(M_k(A))/CU(M_k(A))\to U(M_k(A))/U_0(M_k(A))\to 0
\end{equation}
and one may write
\begin{eqnarray}\label{Dcu-5}
 && U(M_k(A))/CU(M_k(A))=
\Aff(T(A))/\overline{\rho_A(K_0(A))}\oplus J_c(K_1(A))\\
 &&=\Aff(T(A))/\overline{\rho_A(K_0(A))}\oplus J_c(U(M_k(A))/U_0(M_k(A))). \label{Dcu-6}
\end{eqnarray}

For each piecewise smooth and continuous path $\{u(t): t\in [0,1]\}\subset U(M_k(A)),$ define % for each $\tau\in T(A),$
$$
D_A(\{u(t)\})(\tau)={1\over{2\pi i}}\int_0^1 \tau({du(t)\over{dt}}u^*(t))dt,\quad \tau\in T(A).
$$
For each $\{u(t)\},$ the map $D_A(\{u\})$ is a real continuous affine function on $T(A).$
Let $$\overline{D}_A: U_0(M_k(A))/CU(M_k(A))\to \Aff(T(A))/\overline{\rho_A(K_0(A))}$$ denote the de la Harpe and
Skandalis (\cite{HS}) determinant given by
$$
\overline{D}_A(\bar u)=D_A(\{u\})+\overline{\rho_A(K_0(A))},\quad  u\in U_0(M_k(A)),
$$
where $\{u(t): t\in [0,1]\}\subset M_k(A)$ is a piecewise
smooth and continuous path of unitaries with $u(0)=1$ and  $u(1)=u.$  It is known that the de la Harp and Skandalis determinant is independent
of the choice of representatives for $\bar u$ and the choice of  path $\{u(t)\}.$
Define
\begin{equation}\label{July17-1}
\|\overline{D}_A({\bar u})\|=\inf\{\|D_A(\{v\})\|: v(0)=1,\,\,\, v(1)=v\andeqn {\bar v}={\bar u}\},
\end{equation}
where $\|D_A(\{v\})\|=\sup_{\tau\in T(A)}\|D_A(\{v\})(\tau)\|.$

We will fix a metric on $U(M_k(A))/CU(M_k(A)).$ Suppose that  $u, v\in U(M_k(A)).$
Define
\begin{equation}\label{July16-1}
{\rm dist}({\bar u}, {\bar v})=
\left\{
\begin{array}{cl}
2, & \textrm{if $uv^*\not\in U_0(M_k(A))$,}\\
 \|\overline{D}_A(uv^*)\|, & \textrm{otherwise}.
\end{array}
\right.
\end{equation}
%If $uv^*\not\in U_0(M_k(A)),$ then define ${\rm dist}({\bar u}, {\bar v})=2$;
%if $uv^*\in U_0(M_k(A)),$ define
%\beq\label{July16-1}
%{\rm dist}({\bar u}, {\bar v})=\|\overline{D}_A(uv^*)\|.
%\eneq
This is a metric. Note that, if $u, v\in U_0(M_k(A)),$ then
\begin{equation}\label{july19-1}
{\rm dist}(\overline{uv},\overline{1_k})\le {\rm dist}({\bar u}, \overline{1_k})+{\rm dist}({\bar v}, \overline{1_k}).
\end{equation}

%}
\end{df}

%\begin{df}\label{Dcu}
%{\rm
%Let $CU(A)$ denote the closure of commutator subgroup of the unitary group of $A$.
%}
%\end{df}

\begin{df}\label{DLddag}
Let $A$ be a unital separable amenable \CA.
For any finite subset ${\cal U}\subset U(A),$ there exists $\dt>0$ and  a finite subset ${\cal G}\subset  A$
satisfying the following:
If $B$ is another unital \CA\, and if $L: A\to B$ is an ${\cal F}$-$\ep$-multiplicative \morp,
then
$\overline{\langle L(u)\rangle}$ is a well defined element in $U(B)/CU(B)$  for all $u\in {\cal U}.$ We will write
 $L^{\ddag}({\bar u})=\overline{\langle L(u)\rangle}.$  Let $G({\cal U})$ be the subgroup generated by ${\cal U}.$
 We may assume that $L^{\ddag}$ is a well-defined \hm\, on $G({\cal U})$ so that
 $L^{\ddag}(u)=\overline{\langle L(u)\rangle}$ for all $u\in {\cal U}.$ In what follows, whenever we write
 $L^{\ddag},$ we mean that $\ep$ is small enough and ${\cal F}$ is large enough so that $L^{\ddag}$ is well defined
 (see Appendix in \cite{Lin-LAH}).
 Moreover, for an integer $k\ge 1,$ we will also use $L^{\ddag}$ for the map on $U(M_k(A))/CU(M_k(A))$ induced by $L\otimes {\rm id}_{M_k}.$  In particular, when $L$ is a unital \hm, the map $L^{\ddag}$ is well defined
 on $U(A)/CU(A).$

\end{df}

\begin{df}\label{Mappingtorus}
Let $C$ and $B$ be  unital \CA s and let $\phi_1,\phi_2: C\to B$ be two monomorphisms.
Define
\begin{equation}\label{Maptorus-1}
M_{\phi_1, \phi_2}=\{(f,c): C([0,1], B)\oplus C: f(0)=\phi_1(c)\andeqn f(1)=\phi_2(c)\}.
\end{equation}
Denote by $\pi_t: M_{\phi_1, \phi_2}\to B$  the point evaluation at $t\in [0,1].$
One has the following short exact sequence:
\begin{equation*}
0\to SB\stackrel{\imath}{\to}M_{\phi_1, \phi_2} \stackrel{\pi_e}{\to} C\to 0,
\end{equation*}
where $\imath: SB\to M_{\phi, \psi}$ is the embedding and $\pi_e$ is the
quotient map from $M_{\phi, \psi}$ to $C.$
Denote by $\pi_0, \pi_1: M_{\phi, \psi}\to C$ by
the point evaluations at $0$ and $1,$ respectively. Since  both $\phi_1$ and $\phi_2$ are injective, one may identify  $\pi_e$ by
the point-evaluation at $0$ for convenience.

{Suppose that $[\phi]=[\psi]$ in $KL(C,B).$ Then $M_{\phi, \psi}$ corresponds  a trivial element in $KL(A,B).$ In particular, the corresponding extensions
$$
0\to K_i(B) \stackrel{\imath_*}{\to} K_i(M_{\phi, \psi})\stackrel{\pi_e}{\to}K_i( C)\to 0\,\,\,\,\,\,\,{\rm (}i=0,1{\rm )}
$$
are pure.}
\end{df}

\begin{df}\label{DRphipsi}
{Suppose that $T(B)\not=\emptyset.$ Let $u\in M_l(M_{\phi, \psi})$ (for some integer $l\ge 1$) be a unitary
which is a piecewise smooth continuous  function on $[0,1].$
Then
$$
D_B(\{u(t)\})(\tau)={1\over{2\pi} i}\int_0^1 \tau({du(t)\over{dt}}u^*(t))dt\, \tforal \tau\in T(B).
$$
(see \ref{Aq} for the extension of $\tau$ on $M_l(B)$) as defined in \ref{Dcu}.
}
{Suppose that $\tau\circ \phi=\tau\circ \psi$ for all $\tau\in T(B).$
Then there exists a \hm\,
$$
R_{\phi, \psi}: K_1(M_{\phi, \psi})\to \Aff(T(B))
$$
defined by $R_{\phi, \psi}([u])(\tau)=D_B(\{u(t)\}){{(\tau)}}$ as above which is independent  of the choice of the piecewise smooth
paths $u$ in $[u].$
We have the following commutative diagram:}
$$
\begin{array}{ccc}
K_0(B) & \stackrel{\imath_*}{\longrightarrow} & K_1(M_{\phi, \psi})\\
 \rho_B\searrow && \swarrow R_{\phi, \psi}\\
  & \Aff(T(B))
  \end{array}
  $$
{Suppose, in addition,  that $[\phi_1]=[\phi_2]$ in $KK(C,B).$  Then the following exact sequence splits:
\begin{equation}\label{Aug2-2}
0\to \underline{K}(SB)\to \underline{K}(M_{\phi_1, \phi_2})\overset{[\pi_e]}{\rightleftharpoons}_{\theta} \underline{K}(C)\to 0.
\end{equation}
We may assume that $[\pi_0]\circ [\theta]=[\phi_1]$ and
$[\pi_1]\circ [\theta]=[\phi_2].$
In particular, one may write $K_1(M_{\phi, \psi})=K_0(B)\oplus K_1(C).$ Then we obtain a \hm\,
$$
R_{\phi, \psi}\circ \theta|_{K_1(C)}: K_1(C)\to \Aff(T(B)).
$$
We say the rotation map vanishes if there exists such  a map
$\theta$ such above that $R_{\phi, \psi}\circ \theta|_{K_1(C)}=0.$}

{Denote by ${\cal R}_0$ the set of those elements $\lambda\in {\rm Hom}(K_1(C), \Aff(T(B)))$ for which there is a \hm\,
$h: K_1(C)\to K_0(B)$ such that $\lambda=\rho_B\circ h.$ It is a subgroup of ${\rm Hom}(K_1(C), \Aff(T(B))).$
One has a well-defined element ${\overline{R_{\phi, \psi}}}\in {\rm Hom}(K_1(C), \Aff(T(B)))/{\cal R}_0$
(which is independent of the choice of $\theta$).}

In this case, there exists a \hm\, $\theta'_1: K_1(C)\to K_1(M_{\phi, \psi})$ such that
$(\pi_e)_{*1}\circ \theta_1'={\rm id}_{K_1(C)}$ and $R_{\phi, \psi}\circ \theta'_1\in {\cal R}_0$ if and only
if there is $\Theta\in {\rm Hom}_{\Lambda}(\underline{K}(C),\underline{K}(M_{\phi, \psi}))$ such that
$$
[\pi_e]\circ \Theta=[{\rm id}_{C}]\,\,\,\rm{in}\,\,\,KK(C,B)\andeqn R_{\phi, \psi}\circ \Theta|_{K_1(C)}=0.
$$
In other words, $\overline{R_{\phi, \psi}}=0$ if and only if there is $\Theta$ described above
such that $R_{\phi, \psi}\circ \Theta|_{K_1(C)}=0.$  When $\overline{R_{\phi, \psi}}=0,$ one has that
$\theta(K_1(C))\subset {\rm ker}R_{\phi, \psi}$ for some $\theta$ so that (\ref{Aug2-2}) holds. In this case
$\theta$ also gives the following:
$$
{\rm ker}R_{\phi, \psi}={\rm ker}\rho_B\oplus K_1(C).
$$

\end{df}

%\begin{df}\label{Dball}
%Let $X$ be a compact metric space, $x\in X$ be a point and let $r>0.$
%Denote by $B(x, r)=\{y\in X: {\rm dist}(x, y)<r\}.$

%Let $\ep>0.$ Define $f_\ep\in C_0((0,\infty))$ to be a function with
%$0\le f_\ep(t)\le 1,$ $f_\ep(t)=0$ if $t\in [0,\ep/2],$ $f_\ep(t)=1$ if $t\in [\ep,\infty)$ and
%$f_\ep(t)=2{t-\ep/2\over{\ep}}$  if $t\in [\ep/2, \ep].$ Note that $f_\ep f_{\ep/2}=f_\ep.$
%\end{df}

\begin{df}\label{appep}
Let $C$ be a  \CA, let $a, b\in C$ be two elements and let $\ep>0.$
We write $a\approx_\ep b$ if $\|a-b\|<\ep.$
Suppose that $A$ is another \CA,  $L_1, L_2: C\to A$ are two maps and ${\cal F}\subset C$ is a subset.
%and suppose
%that  $\ep>0.$
We write
\beq\label{appep-1}
L_1\approx_{\ep} L_2\,\,\,{\rm on}\,\,\, {\cal F},
\eneq
if $\|L_1(c)-L_2(c)\|<\ep$ for all $c\in {\cal F}.$

\end{df}

\begin{df}\label{Dfull}
Let $A$ and $B$ be \CA s, and assume that $B$ is unital. Let $\mathcal H\subseteq A_+\setminus \{0\}$ be a finite subset, let $T: A_+\setminus\{0\}\to \R_+\setminus\{0\}$ and let $N: A_+\setminus\{0\}\to \N$ be two maps. Then a map $L: A\to B$ is said to be $T\times N$-$\mathcal H$-full if for any $h\in \mathcal H$, if there are $b_1, b_2, ..., b_{N(h)}\in B$ such that $\|b_i\|\leq T(h)$ and
$$\sum_{i=1}^{N(h)} b_i^*L(h)b_i=1_B.$$
\end{df}

\section{The Elliott-Thomsen  building blocks}
%{Non-commutative finite CW compleces of dimension 1}

 To generalize the class of \CA s of tracial rank at most one, we naturally consider
 all  subhomogeneous \CA s with one dimensional spectrum  which, in particular,
 include circle algebras as well as dimension drop algebras. We begin, however, with the following
 special form:

\begin{df}[See \cite{ET-PL} and  \cite{point-line}]\label{DfC1}
{\rm
Let $F_1$ and $F_2$ be two finite dimensional \CA s.
Suppose that there are two unital \hm s
$\phi_0, \phi_1: F_1\to F_2.$
Denote the mapping torus $M_{\phi_1, \phi_2}$ by
$$
A=A(F_1, F_2,\phi_0, \phi_1)
=\{(f,g)\in  C([0,1], F_2)\oplus F_1: f(0)=\phi_0(g)\andeqn f(1)=\phi_1(g)\}.
$$

These \CA s have been introduced into the Elliott program by Elliott and Thomsen (\cite{ET-PL}), and in \cite{point-line}, Elliott used this class of C*-algebras {{and some other building blocks with 2-dimensional spectra}} to realize any weakly unperforated ordered group as the $K_0$-group of a simple ASH C*-algebra. Denote by ${\cal C}$ the class of all unital \CA s of the form $A=A(F_1, F_2, \phi_0, \phi_1)$ and  all finite dimensional \CA s.
These \CA s will be called Elliott-Thomsen building blocks.

A unital \CA\, $C\in {\cal C}$ is said to be {\it minimal} if
it is not a direct sum of more than one copy of \CA s in ${\cal C}.$
If $A\in {\cal C}$ is minimal and  is not finite dimensional, in what follows, we may  assume that ${\rm ker}\phi_0\cap {\rm ker}\phi_1=\{0\}.$
In general, if $A\in {\cal C},$ and ${\rm ker}\phi_0\cap {\rm ker}\phi_1\not=\{0\}.$ Then, we can write $A=A_1 \oplus ({\rm ker}\phi_0\cap {\rm ker}\phi_1)$, where $A_1=A(F'_1, F_2,\phi'_0, \phi'_1)$,  $F_1=F'_1\oplus ({\rm ker}\phi_0\cap {\rm ker}\phi_1)$ and $\phi'_i=\phi_i|_{F'_1}$ for $i=1,2$. (Note that $A_1=A(F'_1, F_2,\phi'_0, \phi'_1)$ satisfies the condition  ${\rm ker}\phi'_0\cap {\rm ker}\phi'_1 =\{0\},$ and that ${\rm ker}\phi_0\cap {\rm ker}\phi_1$ is a finite dimensional \CA.)

 %%%%%%%%%%%%%%

For $t\in (0,1),$ define $\pi_t: A\to F_2$ by $\pi_t((f,g))=f(t)$ for all $(f,g)\in A.$
For  $t=0,$ define $\pi_0: A\to \phi_0(F_1)\subset F_2$ by $\pi_0((f, g))=\phi_0(g)$ for all $(f,g)\in A.$
For $t=1,$ define $\pi_1: A\to \phi_1(F_1)\subset F_2$ by $\pi_1((f,g))=\phi_1(g))$ for all $(f,g)\in A.$
In what follows, we will call $\pi_t$  a point-evaluation of $A$ at $t.$
There is a canonical map $\pi_e: A \to F_1$ defined by $\pi_e(f,g)=g$ for all
pair $(f, g)\in A.$  It is a surjective map.
{\it The notation $\pi_e$ will be used for this map throughout  this paper.}

If $A\in {\cal C}$, then $A$ is the pull-back of
\begin{equation}\label{pull-back}
\xymatrix{
A \ar@{-->}[rr] \ar@{-->}[d]^-{\pi_e}  && C([0,1], F_2) \ar[d]^-{(\pi_0, \pi_1)} \\
F_1 \ar[rr]^-{(\phi_0, \phi_1)} & & F_2 \oplus F_2
}
\end{equation}
%
%\beq\label{pull-back}
%\begin{array}{ccc}
%\hspace{-0.15in}A & \cdots \to  &C([0,1], F_2)\\
%\downarrow_{\pi_e} && \hspace{0.4in}\downarrow_{(\pi_0, \pi_1)}\\
%\hspace{-0.15in}F_1 &\longrightarrow&  F_2\oplus F_2.
%\end{array}
%\eneq
Every such pull-back is an algebra in ${\cal C}.$
Infinite dimensional  C*-algebras in ${\cal C}$ are also called {\it one-dimensional non-commutative finite CW complexes} (NCCW)
(see \cite{ELP1} and \cite{ELP2}). %{\color{Green} (This diagram contains the case of finite dimensional C*-algebras---just let $F_2=\{0\}$.)}

Denote by ${\cal C}_0$ the sub-class of those \CA s $A$ in ${\cal C}$ such that $K_1(A)=\{0\}.$

{It follows from  Theorem 6.22 of \cite{ELP1} that \CA s in ${\cal C}$ are semiprojective, an important feature that we will use
later without warning.}
}

\end{df}

\begin{lem}\label{2Lg1}
Let $f\in C([0,1], M_k(\C))$ and let $a_0,a_1\in M_k(\C)$ be invertible elements with
$$\|a_0-f(0)\|<\ep\andeqn  \|a_1-f(1)\|<\ep.$$
Then there exists an invertible element  $g\in C([0,1], M_k(\C))$ such that $g(0)=a_0,~ g(1)=a_1$ and
$$\|f(t)-g(t)\|<\ep \qq \tforal t\in [0,1].$$
\end{lem}
\begin{proof}
Let $S\sbs M_k(\C)$ be the set consisting of all singular matrices.  Then $M_k(\C)$ is a $2k^2$-dimensional differential manifold (diffeomorphic to $\R^{2k^2}$) with $S$ being finite union of closed submanifold of codimension at least two. Since each continuous map between two differential manifolds (perhaps with boundary) can be approximated
arbitrarily well by smooth map{s}, we can find $f_1\in C^{\infty}([0,1],M_k(\C))$ with $f_1(0)=a_0$ and $f_1(1)=a_1$ and $\|f_1(t)-f(t)\|<\ep'<\ep$.  Apply  the relative version of the transversality theorem---the corollary on  page 73 of \cite{Guillemin-Pollack} and its proof (see page 70 and 68 of \cite{Guillemin-Pollack}), for example, with $Z=S, Y=M_k(\C), X=[0,1]$ with $\partial X=\{0,1\}$, one can find $g\in C^{\infty}([0,1],M_k(\C))$ with $g(0)=f_1(0)$,
 $g(1)=f_1(1)$, $g(t)\notin S$ and $\|g(t)-f_1(t)\|< \ep-\ep'$.
\end{proof}

\begin{prop}\label{2pg3}
If  $A\in {\cal C}$, then $A$ has stable rank one.
\end{prop}
\begin{proof}
Let $(f,a)$ be in $A$ with $f\in C([0,1], F_2)$ and $a\in F_1$  with $f(0)=\phi_0(a)$ and $f(1)=\phi_1(a)$.  For any $\ep>0$, since $F_1$ is a finite dimensional \CA, there is an invertible element $b\in F_1$ such that $\|b-a\|<\ep$.  Since $\phi_0, \phi_1$ are unital, $\phi_0(b)$ and $\phi_1(b)$ are invertible,
$$
\|\phi_0(b)-f(0)\|<\ep\andeqn \|\phi_1(b)-f(1)\|<\ep.
$$
By Lemma \ref{2Lg1},
%1.4,
there exists an invertible element $ g\in C([0,1], F_2)$ (applied to each direct summand of $F_2$) such that $g(0)=\phi_0(b),~~ g(1)=\phi_1(b)$ and
$$\|g-f\|<\ep.$$
This is what as desired.
\end{proof}

\begin{NN}\label{2Rg10}
Let $F_1=M_{R(1)}\oplus M_{R(2)}\oplus \cdots \oplus M_{R(l)}$ and $M_2=M_{r(1)}\oplus M_{r(2)}\oplus \cdots \oplus M_{r(k)}.$ $A=(F_1, F_2,\phi_0, \phi_1).$
Denote by $m$ the greatest common factor of
${\big \{}R(1),R(2), \cdots, R(l){\big \}}$.  Then each $r(j)$ is also a multiple of $m$.
 Let $\widetilde{F}_1=M_{\frac{R(1)}{m} }(\C)\oplus M_{\frac{R(2)}{m} }(\C)\oplus \cdots \oplus M_{\frac{R(l)}{m} }(\C)$
 and
 $\widetilde{F}_2=M_{\frac{r(1)}{m} }(\C)\oplus M_{\frac{r(2)}{m} }(\C)\oplus \cdots \oplus M_{\frac{r(k)}{m} }(\C)$.  Let ${\widetilde{\phi}_0},\widetilde{\phi}_1:~ \widetilde{F}_1\to \widetilde{F}_2$ be maps such that
\vspace{-0.1in} \beq\nonumber
 \widetilde{\phi}_{0*}, \widetilde{\phi}_{1*}: ~ K_0(\widetilde{F}_1)=\Z^l \longrightarrow \Z^k
 \eneq
 satisfying $\widetilde{\phi}_{0*}=(a_{ij})_{k\times l}$ and $\widetilde{\phi}_{1*}=(b_{ij})_{k\times l}$.  That is, {the maps which  are represented by the same matrices represent $\phi_{0*}$ and $ \phi_{1*}$.}

By \ref{2Pg12},
$$A(F_1, F_2, \phi_0, \phi_1)\cong M_m(A(\tilde{F}_1, \tilde{F}_2, \tilde{\phi}_0, \tilde{\phi}_1)).$$
\end{NN}

\begin{prop}\label{2Lg13}
Let $A=A(F_1, F_2, \phi_0, \phi_1)$.  Then $K_1(A)=\Z^l/\mathrm{Im}({\phi_0}_{*0}-{\phi_1}_{*0})$ and

\beq\label{2Rg12-n1}
K_0(A)\cong \left\{
\left(\!\!
         \begin{array}{c}
           v_1 \\
           v_2 \\
           \vdots \\
           v_l
         \end{array}\!\!
       \right) \in \Z^l~ , ~~  \phi_{~\!\!\!\!_{0*}}\!\!\left(\!\!
         \begin{array}{c}
           v_1 \\
           v_2 \\
           \vdots \\
           v_l
         \end{array}\!\!
       \right)=   \phi_{~\!\!\!\!_{1\!*}}\!\!\left(\!\!
         \begin{array}{c}
           v_1 \\
           v_2 \\
           \vdots \\
           v_l
         \end{array}\!\!
       \right)
\right\}
\eneq
with positive cone being $K_0(A) \cap \Z^l_+$, and {scale} $\left(\!\!
         \begin{array}{c}
           R_1 \\
          R_2\\
           \vdots \\
           R_l
         \end{array}\!\!
       \right)\in \Z^l$, where
$\Z_+^l=\left\{
\left(\!\!
         \begin{array}{c}
           v_1 \\
           v_2 \\
           \vdots \\
           v_l
         \end{array}\!\!
       \right)~;~~ v_i\geq 0
       \right\} \sbs \Z^l~.
$
Moreover, the map $\pi_e: A\to F_1$ induces the natural order embedding
$(\pi_e)_{*0}: K_0(A)\to K_0(F_1)=\Z^l,$ {in particular,  ${\rm ker}\rho_A=\{0\}.$}

Furthermore, if $K_1(A)=\{0\},$ then $\Z^l/K_0(A)\cong K_1(C_0((0,1), F_2))$ is torsion free, 
%(The quotient group
%$K_0(F_1)/{{(\pi_e)_{*0}(K_0(A))}}$ is always torsion free, since $K_0(F_1)/(\pi_e)_{*0}(K_0(A))$ is a subgroup of $K_0(F_2)$. It happens that $K_0(F_1)/{{(\pi_e)_{*0}}}(K_0(A))=K_0(F_2)$ holds if and only if $K_1(A)=\{0\}$),
and in this case, {{
$[\pi_e]|_{K_i(A, \Z/k\Z)}$ is injective for all $k\ge 2$ and $i=0,1.$}}

\end{prop}

\begin{proof}
Most of these assertions are known. We sketch the proof here. 
%Note that if ${\pi_e}_*(p)\in\Z_+^l$ for $p\in K_0(A)$, then $p$ must be positive. 
Consider the short exact sequence
\begin{displaymath}\label{sixterm-1}
\xymatrix{
0\ar[r] & \mathrm{C}_0((0,1), F_2) \ar[r] & A \ar[r]^{\pi_e} & F_1\ar[r] & 0.}
\end{displaymath}
We obtain 
\begin{equation}\label{sixterm}
\xymatrix{
0 \ar[r] & K_0(A) \ar[r]^{{\pi_e}_*} & K_0(F_1) \ar[r] & K_0(F_2) \ar[r] & K_1(A)\ar[r] & 0,
}
\end{equation}
where the map $K_0(F_1)\to K_0(F_2)$ is given by ${\phi_0}_{*0}-{\phi_1}_{*0}$.
In particular, $(\pi_e)_{*0}$ is injective.   
%If $x\in (\pi_e)_{*0}(K_0(A))_+,$ then $(\phi_0)_{*0}(x)=(\phi_1)_{*0}(x).$ 
If $p\in \pi_e(M_m(A))$ is a projection, then $(\phi_0)_{*0}(p)=(\phi_1)_{*0}(p).$
Therefore $\phi_0(p)$ and $\phi_1(p)$ have the same rank. It follows  that there is a projection 
$q\in M_m(A)$ such that $\pi_e(q)=p.$ This implies that $(\pi_e)_{*0}$ is an order embedding. 
This also implies that ${\rm ker}\rho_A=\{0\}.$  Other descriptions of $K_i(A)$ ($i=0,1$) also follow.
The quotient group
$K_0(F_1)/{{(\pi_e)_{*0}(K_0(A))}}$ is always torsion free, since $K_0(F_1)/(\pi_e)_{*0}(K_0(A))$ is a subgroup of $K_0(F_2)$. It happens that $K_0(F_1)/{{(\pi_e)_{*0}}}(K_0(A))=K_0(F_2)$ holds if and only if $K_1(A)=\{0\}.$

 In the case that $K_1(A)=\{0\},$
one also computes that $K_0(A, \Z/k\Z)$ may be identified with
$K_0(A)/kK_0(A)$ and $K_1(A,\Z/k\Z)=\{0\}$ for all $k\ge 2.$
%Moreover $K_0(F_1)/(\pi_e)_{*0}(K_0(A))\cong K_0(F_2).$

{{To see that $[\pi_e]|_{K_0(A,\Z/k\Z)}$ is injective for $k\ge 2,$  let
${\bar x}\in K_0(A,\Z/k\Z)=K_0(A)/kK_0(A)$ and let
$x\in K_0(A)$ be such that its image in $K_0(A)/kK_0(A)$ is ${\bar x}.$
If $[\pi_e]({\bar x})=0,$ then $(\pi_e)_{*0}(x)\in kK_0(F_1).$ Let
$y\in K_0(F_1)$ be such that $ky=(\pi_e)_{*0}(x).$ Then $k{\bar y}={\overline{(\pi_e)_{*0}(x)}}=0$ in
$K_0(F_1)/(\pi_e)_{*0}(K_0(A)).$ This implies that $K_0(F_2)$ has torsion. A contradiction.
Therefore $[\pi_e]|_{K_0(A/\Z/k\Z)}$ is injective for $k\ge 2.$ }}
  Then {{ the rest of}} proposition also  follows from \eqref{sixterm}.
\end{proof}

\begin{NN}\label{2Rg13}
{\rm
Let us describe how to identify a unitary $u\in A$ with $[u]\in K_1(A)=\Z^k/(\phi_{1*}-\phi_{0*})(\Z^l)$.
Let $u=(f,a)\in A$, where $a=(a_1,a_2,..., a_l)\in \bigoplus_{j=1}^l M_{R(j)}(\C)=F_1$.  Note that every unitary in $M_N(\C)$ can be written as $e^{ih}$ {{for some self-adjoint element $h\in M_N(\C)).$}}  Write $h=(h_1, h_2, ..., h_l)\in \bigoplus_{j=1}^l M_{\{ j\}}(\C)=F_1$ such that $a=e^{ih}$.  Let $(g, h)\in A$ be defined by
$$
g(t)=\left\{\begin{array}{lll}
              \phi_0(h)(1-2t), &  & 0\leq t\leq \frac{1}{2}, \\
                &   &  \\
                \phi_1(h)(2t-1), &  & \frac{1}{2}< t\leq 1.
            \end{array}
\right.
$$

Obviously, $(g, h)$ is well defined and inside $A$.  Let $v=(f_1, 1_{F_1})$, where $f_1(t)=f(t)e^{-ig(t)}$.
Let $U(s)=(f(t)e^{-ig(t)s}, e^{i h(1-s)})$ for $s\in [0,1].$ Then $\{U(s): s\in [0,1]\}$ is a continuous path
of unitaries in $A$ with $U(0)=u$ and $U(1)=(f_1, 1_{F_1}).$ Therefore $uv^*\in U_0(A).$
In particular, $[u]=[v] \in K_1(A)$, and $v=(v_1,v_2,...,v_k)$ with $v_j(0)=v_j(1)=1$.  Write
$f=(a_1,a_2,...,a_k)$ and
$fe^{-igs}=(b_1(s),b_2(s),...,b_k(s)),$ where $a_j,\,b_j(s)\in C([0,1], M_{r(j)})$ for all $s\in [0,1].$
Suppose that ${\rm det}(a_m(t))=1$ for all $t\in [0,1],$   $m=1,2,...,k,$
and ${\rm det}(e^{i h_j})=1$ for $j=1,2,...,l.$
Define $B_m(s)={\rm det}(b_m(s))^*,$  $B(s)=(B_1(s), B_2(s),...,B_k(s))$  and
$$d(s)=({\rm det}(e^{i h_1(1-s)}),{\rm det}(e^{i h_2(1-s)}),..., {\rm det}(e^{i h_k(1-s)})).$$
 Define
$W(s)=(B_j(s)U(s), d(s)^*e^{ih(1-s)}).$ Then $\{W(s): s\in [0,1]\}$ is a continuous path of unitaries in $A$
such that $W(0)=u$ and $U(1)=(B(1)f_1,1).$ Note that ${\rm det}(B_m(1) b_m(1))=1$ for
$1\le m\le k.$
 In other words, if ${\rm det}(u(\psi))=1$ for all irreducible representations $\psi,$ then we may also assume that
 ${\rm det}(v(\psi))=1$ for $\psi.$
Note that  $v\in C_0((0,1), F_2\widetilde{)}$, the unitilization of the ideal $C_0((0,1), F_2)\sbs A$.  Hence $u$ defines an element $(s_1,s_2,..., s_k)\in \Z^k=K_1(C_0((0,1), F_2))$, where $s_j$ is winding number of the map
$$t\in[0,1]\lr \mathrm{det}(v_j(t))\in \mathbb T\subseteq \C.$$ In particular, if ${\rm det}(u(\psi))=1$
for all irreducible representations $\psi,$ then $u\in U_0(A).$
Such a $k$-tuple gives an element
$$[(s_1,s_2,..., s_k)]\in \Z^k/(\phi_{1*}-\phi_{0*})(\Z^l)~.
$$
}
\end{NN}

\begin{lem}\label{2Lg8}
Let $A=A(F_1, F_2, \phi_0, \phi_1)$ be as in Definition \ref{DfC1}. A unitary $u\in U(A)$ is in $CU(A)$ if and only if for each irreducible representation $\psi$ of $A$, one has that $\mathrm{det}(\psi(u))=1$.
\end{lem}

\begin{proof}
Obviously, the condition is necessary. One only has to show that the condition is also sufficient.
From the discussion \ref{2Rg13} above, if $u\in U(A)$ with ${\rm det}(\psi(u))=1$ for all
irreducible representations $\psi,$ then  $u\in U_0(A).$

Write $F_1=M_{R(1)}\oplus R_{R(2)}\oplus \cdots \oplus M_{R(l)}$
and $F_2=M_{r(1)}\oplus M_{r(2)}\oplus \cdots \oplus M_{r(k)}.$
Then, since $u\in U_0(A),$ we may write that $u=\prod_{n=1}^m \exp(i 2\pi h_n)$ for some  $h_n\in A_{s.a.},$ $n=1,2,...,m.$
We may write $h_n=(h_{nI},h_{nq})\in A$ with $\phi_i(h_{nq})=h_{nI}(i),$ $i=0,1.$
For any irreducible representation $\psi$ {{of A}}, $\sum_{j=1}^n {\rm Tr}_\psi\psi(h_j)\in \Z,$ where
${\rm Tr}_\psi$ is the standard (unnormalized) trace on $\psi(A)\cong M_{n(\psi)}$ for some
integer $n(\psi).$
 Put $H_j(t)={\rm Tr}_j(\sum_{n=1}^m \pi_j(h_n(t)))$ and $H(t)=(H_1, H_2,...,H_k)$ for $t\in [0,1],$ where ${\rm Tr}_j$ is
 the standard trace on $M_{r(j)}$ and $\pi_j: F_2\to M_{r(j)}$ is the projection map,  $j=1,2,...,k.$ Put $a(\psi)=\sum_{j=1}^n {\rm Tr}_\psi(h_j)$ for
 those $\psi$ corresponding to  irreducible representations
of $F_1.$ {{Since $\sum_{j=1}^n {\rm Tr}_\psi\psi(h_j)\in \Z,$ there is a projection}}  $p\in M_N(F_1)$ such that $\mathrm{Tr}_{\psi}\psi(p)=a(\psi)$ for some $N\ge 1.$
Then, ${\rm Tr}_j'(\pi_j'\circ \phi_0(p))=H_j(0)=H_j(1)={\rm Tr}_j'({{\phi}}_1(p)),$ where
${\rm Tr}_j'$ is the standard trace on $M_N(M_{r(j)})$ and $\pi_j': M_N(F_2)\to M_N(M_{r(j)})$ is the projection map.
There is a projection $P_j\in M_N(C([0,1], M_{r(j)}))$ such that $P_j(0)=\pi_j\circ \phi_0(p)$ and
$P_j(1)=\pi_j\circ \phi_1(p).$ Put $P\in M_N(C([0,1], F_1))$ such that $\pi_j'(P)=P_j$ and
$e=(P, p).$ Then $e\in M_N(A).$
Consider
the continuous path $u(t)=\prod_{n=1}^m \exp(i 2\pi h_nt)$ for $t\in [0,1].$
Then  $u(0)=1,$  $u(1)=u$ and
$$
\tau({du(t)\over{d{{t}}}}u^*(t))=i 2\pi \sum_{n=1}^m\tau(h_n)\rforal \tau\in T(A).
$$
But,  for  all $a\in A$ and $\tau\in T(A),$
$$
\tau(a)=\sum_{s=1}^l \af_s {\rm tr}_s (a)+\sum_{j=1}^k \int_{(0,1)} {\rm tr}_j(\pi_n(a)) d\mu_j(t),
$$
where $\mu_j$ is a Borel measure on $(0,1),$ ${\rm tr}_s$ is the tracial state on $M_{R(s)}$ and
${\rm tr}_j$ is a tracial state on $M_{r(j)},$ $\af_i\ge 0$ and
$\sum_{s=1}^l\af_s+\sum_{j=1}^k \|\mu_j\|=1.$
We compute that
$$
{1\over{2\pi i}}\int_0^1 \tau({du(t)\over{dt}}u(t)^*)dt=\tau(e)\tforal \tau\in T(A).
$$
In other words, $D_A(\{u(t)\})\in \rho_A(K_0(A))$ (see \ref{Dcu}).
It follows from \cite{Thomsen-rims} and the fact that $A$ has stable rank one that $u\in CU(A)$ (see also \cite{GLX-ER}).
\end{proof}

The following is known (see \cite{Thomsen-Hom-AT}, \cite{Thomsen-UOrbit-IA}, and  \cite{Phil-cer}). % {\bf I am not sure the references--L}
\begin{lem}\label{2Lgsphillips}
Let $u$ be a unitary in $C([0,1], M_n).$
Then, for any $\ep>0,$ there exist continuous functions
$h_j\in C([0,1])_{s.a.}$ such that
$$
\|u-u_1\|<\ep,
$$
where $u_1=\exp(i \pi H),$ $H=\sum_{j=1}^n h_jp_j$ and $\{p_1, p_2,...,p_n\}$ is a set of mutually orthogonal rank one projections in $C([0,1], M_n),$
and $\exp(i \pi h_j(t))\not=\exp(i\pi h_k(t))$ if $j\not=k$ for all
{{$t\in (0,1),$ $u_1(0)=u(0)$  and $u_1(1)=u(1).$}}
%Moreover,  suppose that $u(0)=\sum_{j=1}^n \exp(i a_j)p_j(0)$ for some real number $a_j$ which are distinct,
%we may assume that $h_j(0)=a_j.$

 Furthermore,
if ${\rm det}(u(t))=1$ for all $t\in [0,1],$ then
we may also assume that ${\rm det}(u_1(t))=1$ for all $t\in [0,1].$

\end{lem}

%\begin{NN}
%{\rm
%It is well known (see [Thomsen] and [Ell 1-2]) that for any unitary$u\in U(M_k(C[0,1]))$, and $\ep>0$, there are real valued functions $h_1, h_2, \cdots, h_k:~ [0,1]\to \R$ and $v\in U(M_k(C[0,1]))$ such that
%$$\|u(t) -v(t)\left(
%\begin{array}{cccc}
 % e^{2\pi i h_1(t)} &   &   &   \\
 %   & e^{2\pi i h_2(t)}  &  &   \\
 %   &   & \ddots &   \\
 %   &   &   & e^{2\pi i h_k(t)}
%\end{array}
%\right)v^*(t)\|<\ep.$$
%In fact,, one can approximate $u$ by a unitary $u'$ within $\ep$ such that $u'(t)$ has distinct eigenvalues at each $t\in [0,1]$ (see [Choi-Ell]).  Then one can write the eigenvalues of $u'(t)$ in the form of $e^{2\pi i h_j(t)}$ for continuous functions $h_j: [0,1]\to \R$.

%Also if $u$ satisfies the condition that $det(u(t))=1$, then one can approximate $u$ by $u'$ such that $u'(t)$ has distinct eigenvalues at each $t\in [0,1]$ and $det(u'(t))=1$.

\begin{rem}\label{2Rg5}
It follows  from the above (also see Theorem 1-2 \cite{Thomsen-Hom-AT} and \cite{Thomsen-UOrbit-IA}% [Thomsen1-2,Thomsen, K.: Homomorphisms between finite direct sums of circle algebras. Aarhus: Matematisk Institut, Linear and multilinear algebra, and Inductive limits of interval algebras: unitary orbits of positive elements ]
) that  $u, v \in U(M_k(C[0,1]))$ are approximately unitarily equivalent (that is, for any $\ep$, there is $w\in U(M_k(C[0,1]))$ such that $\|u-wvw^*\|<\ep$) if and only if for each $t\in [0,1]$,  $u(t)$ and $v(t)$ have the same set of eigenvalues counting multiplicities (see also Lemma 3.1 of \cite{Lin-PNAS}).
%Spectrum$(u(t))=$Spectrum$(v(t))$.
\end{rem}

\begin{lem}\label{2Lg6}
Let $A=A(F_1, F_2, \phi_0, \phi_1)$ be as in Definition \ref{DfC1}. For any unitary $(f,a)\in U(A), ~\ep>0$, there is a unitary $(g,a)\in U(A)$ such that $\|g-f\|<\ep$ and, for each block $M_{r(j)} \sbs F_2= \bigoplus_{j=1}^k M_{r(j)}$, there are real valued functions $h^j_1,h^j_2,..., h^j_{k_j}: [0,1]\to \R$  {{ such that
$g_j=\sum_{j=1}^n h_jp_j$ and $\{p_1, p_2,...,p_n\}$ is a set of mutually orthogonal rank one projections in $C([0,1], M_{r(j)})$ and $\exp(i \pi h_j(t))\not=\exp(i\pi h_k(t))$ if $j\not=k$ for all
$t\in (0,1).$ Moreover, if $(f,a)\in CU(A),$  one can choose $(g,a)\in CU(A).$}}
%%
%such that
%
%and $v\in U(C[0,1], M_{r(j)})$ such that
%$g_j=
%
%{\red{, for all $t\in [0,1],$}}
%\begin{equation}\label{decomp-g}
%g^j(t)=v(t)\left(
%\begin{array}{cccc}
 % e^{2\pi i h^j_1(t)} &   &   &   \\
   % & e^{2\pi i h^j_2(t)}  &  &   \\
    %&   & \ddots &   \\
    %&   &   & e^{2\pi i h^j_{k_j}(t)}
%\end{array}
%\right)v^*(t),
%\end{equation}
%where $g^j=\pi_j(g)$ and where $\pi_j: C([0,1], F_2)\to  C([0,1], M_{r(j)})$ is the projection map,
%$j=1,2,..., k.$ {\red{Moreover, if $(f,a)\in CU(A),$ one can choose $(g,a)\in CU(A).$}}
\end{lem}

\begin{proof}
{{By \ref{2Lgsphillips}, f}}or each unitary $f^j\in C([0,1], M_{r(j)})$, one can approximates $f^j$ by $g_1^j$ to within
{{$\ep$}} such that $g^j(0)=f^j(0)$, $g^j(1)=f^j(1)$ and for each $t$ in the open interval $(0,1)$, $g^j(t)$ has distinct eigenvalues.   {{If $(f,a)\in U(A),$ then  $(g,a)\in U(A)$ too.}} {{Combining \ref{2Lgsphillips} with  \ref{2Lg8}, the last statement also follows.}}
% Then, by \ref{2Lgsphillips}, %1.10,
%the unitary  $g_1^j$ can be approximated by $g^j$ within $\ep/2$ which can be written as the desired form (where for each $t\in [0,1]$ $g_1^j(t)$ and $g^j(t)$ has same spectra $\{e^{2\pi i h^j_l(t)} \}_{l=1}^{k_j}$).
\end{proof}

%{\bf The above proof did not really say how to do this from \ref{2Lgsphillips}. Answer: I modified the proof and should be OK now.}

\begin{rem}\label{2Rg7}
In %\eqref{decomp-g}, %the form of $g^j(t)$ in \ref{2Lg6},
{{(\ref{2Lg6}),}}
one may assume that $h^j_1(0),h^j_2(0), ..., h^j_{k_j}(0)\in [0,1)$. For an arbitrarily small $t\in (0,\dt)$, one may assume
\begin{equation}\label{eqn-G-ast}
\max \{ h^j_i(t);~ 1\leq i\leq k_j \}-\min\{ h^j_i(t);~ 1\leq i\leq k_j \}<1
\end{equation}
and $h^j_{i_1}(t)\not= h^j_{i_2}(t)$ for $i_1\not= i_2$.  From the choice of $g^j$, we know that for any $t\in (0,1)$, $e^{2\pi i h^j_{i_1}(t)}\not= e^{2\pi i h^j_{i_2}(t)}$.  That is, $h^j_{i_1}(t) - h^j_{i_2}(t)\not\in \Z$.  This implies that \eqref{eqn-G-ast} holds for all $t\in (0,1)$. Hence, one may assume that
\begin{equation}\label{eqn-G-ast-01}
\max \{ h^j_i(1);~ 1\leq i\leq k_j \}-\min\{ h^j_i(1);~ 1\leq i\leq k_j \}\leq1.
\end{equation}
\end{rem}

%Let $CU(A)$ denote the closure of commutator subgroup of the unitary group of $A$.

\begin{lem}\label{2Lg9}
For any $u\in CU(A)$, one has that $\mathrm{cer}(u)\leq 2+\ep$ and $\mathrm{cel}(u)\leq 4\pi$.
{Moreover, there exists a continuous path of unitaries $\{u(t): t\in [0,1]\}\subset CU(A)$  with length
at most $4\pi$ such that $u(0)=1_A$ and $u(1)=u.$}
\end{lem}

\begin{proof}
%{\bf Step 1.}
{{ Case (i): the case that}}
%Assume
$u=(f,a)$ with $a=1$.  \, As in
%\ref{2Lgsphillips}, \ref{2Rg5} and
\ref{2Lg6},
 up to approximation within arbitrarily small pre-given number $\ep>0$, $u$ is unitarily equivalent to $v=(g,a)\in CU(A)$ with
$$g_j(t) =\diag{\big (}e^{2\pi i h^j_1(t)}, e^{2\pi i h^j_2(t)},..., e^{2\pi i h^j_{k_j}(t)}{\big )}$$
with distinct eigenvalues for each $t\in (0,1)$.  (Note that since $f(0)=f(1)=1=g(0)=g(1)$, the unitary to intertwining the approximation of $u$ and $v$ can be chosen to be $1$ at $t=0, 1$, and therefore, the unitary is in $A=A(F_1, F_2, \phi_0, \phi_1)$.)

Furthermore, one can assume {{that}}
$$h^j_1(0)=h^j_1(0)=\cdots = h^j_{k_j}(0)=0.$$
Since $\mathrm{det}(g_j(t))=1$ for all $t\in [0,1]$, one has that $h^j_1(t)+h^j_2(t)+ \cdots + h^j_{k_j}(t)\in \Z$.

By the continuity of each $h^j_s(t)$, we know that
\begin{equation}\label{eq-G-0819-001}
\sum_{s=1}^{k_j} h^j_s(t) = 0.
\end{equation}
Furthermore, by $h^j_s(1)\in \Z$ (since $g_j(1)=1$), we know that $h^j_s(1)=0$ for all $s\in \{1, 2, \cdots, k_j\}$.  Otherwise, ${\displaystyle \min_s } {\big \{ }h^j_s(1) {\big \} }\leq -1$ and ${\displaystyle \max_s } {\big \{ } h^j_s(1){\big \} } \geq 1$ which implies that
\begin{eqnarray*}
{\displaystyle \max_s } {\big \{ } h^j_s(1){\big \} }-
{\displaystyle \min_s } {\big \{ }h^j_s(1) {\big \} } \geq 2,
\end{eqnarray*}
and this contradicts to Remark \ref{2Rg7}.  That is, one has proved that $h={\big (} (h^1, h^2, \cdots, h^k), 0{\big )}$, where $h^j(t)=\diag {\big (} h^j_1(t), h^j_2(t), \cdots, h^j_{k_j}(t) {\big )}$ is an element in $A=A(F_1, F_2, \phi_0, \phi_1)$ with $h(0)=h(1)=0$.
 As $g=e^{2\pi ih}$, we have $\mathrm{cer}(u)\leq 1+\ep$.
{We also have ${\rm tr}(h(t))=0$ for all $t.$}

It follows from \eqref{eq-G-0819-001} above and ${\displaystyle \max_s } {\big \{ } h^j_s(t){\big \} }-
{\displaystyle \min_s } {\big \{ }h^j_s(t) {\big \} } \leq 1~$ (see \eqref{eqn-G-ast-01} in Remark \ref{2Rg7}) that
\begin{equation*}
h^j_s(t)~ \sbs ~(-1, 1),\quad t\in [0,1],\ s=1,2,\cdots, k_j~.
\end{equation*}
Hence $\|2\pi h\|\leq 2\pi$ which implies $\mathrm{cel}(u) \leq 2\pi$.
Moreover let $u(s)=\exp(i s2h).$ Then $u(0)=1$ and $u(1/2)=u.$  Since ${\rm tr}(s2h(t))=2s\cdot{\rm tr}(h(t))=0$
for all $t,$ one has that $u(s)\in CU(A)$ for all $s\in [0,1/2].$

%{\bf Step 2.}
{{Case (ii):}} The general case.  Since $a=(a^1, a^2, \cdots, a^l)$ with $\mathrm{det}(a^j)=1$ for $a^j\in F^j_1$.  So $a^j=\exp({2\pi i h^j})$ for $h^j\in F^j_1$ with $\mathrm{tr}(h^j)=0$ and $\|h^j\|<1$.  Define $H\in A(F_1, F_2, \phi_0, \phi_1)$ by
\begin{displaymath}
H(t)=
\left\{
\begin{array}{lcl}
  \phi_0(h^1, h^2, ... , h^l)\cdot (1-2t), & \textrm{if $0\leq t\leq \frac{1}{2}$}, \\
   & & \\
  \phi_1(h^1, h^2, ... , h^l)\cdot(2t-1),   & \textrm{if $\frac{1}{2}< t\leq 1$}.
\end{array}
\right.
\end{displaymath}
(Note $H(\frac12)=0,~ H(0)= \phi_0(h^1, h^2, \cdots, h^l)$, $ H(1)= \phi_1(h^1, h^2, \cdots, h^l)$ and therefore $H\in A(F_1, F_2, \phi_0, \phi_1)$.
{Moreover ${\rm tr}(H(t))=0$ for all $t.$} Then $u'=u \cdot \exp({-2\pi i H})\in A(F_1, F_2, \phi_0, \phi_1)$ with
$u'(0)=u'(1)=1$.  By {{Case (i),}} $\mathrm{cer}(u')\leq 1+\ep$ and $\mathrm{cel}(u')\leq 2\pi$,  we have $\mathrm{cer}(u)\leq 2+\ep$ and $\mathrm{cel}(u)\leq 2\pi+2\pi \|H\| \leq 4\pi.$ Furthermore, we
note that $\exp(-2\pi s H)\in CU(A)$ for all $s$ as in {{Case (i)}}.
%{\bf Step} 1. }
\end{proof}

\begin{NN}\label{2Rg11}
Let $F_1=M_{R_1 }(\C)\oplus M_{R_2 }(\C)\oplus \cdots \oplus M_{R_l}(\C),$  let
 $F_2=M_{r_1}(\C)\oplus M_{r_2 }(\C)\oplus \cdots \oplus M_{r_k }(\C)$  {{and}}  let $\phi_0,\,\,
 %let
  \phi_1:~ F_1\to F_2$ be unital homomorphisms, where $R_j$ and $r_i$ are positive integers.  Then $\phi_0$ and $\phi_1$ induce homomorphisms $$\phi_{0*}, \phi_{1*}: ~ K_0(F_1)=\Z^l \longrightarrow K_0(F_2)=\Z^k$$
 by matrices $(a_{ij})_{k\times l}$ and $(b_{ij})_{k\times l}$, respectively, where {$r_i=\sum_{j=1}^l a_{ij}R_j$ for $i=1,2,...,k.$}
\end{NN}

\begin{prop}\label{2Pg12}
For fixed finite dimensional C*-algebras $F_1, F_2$, the C*-algebra $A=A(F_1, F_2, \phi_0, \phi_1)$ is completely determined by $\phi_{0*}, \phi_{1*}: \Z^l \longrightarrow \Z^k$ {{(up to isomorphisms)}}.
\end{prop}

\begin{proof} Let $B=A(F_1, F_2, \phi'_0, \phi'_1)$ with $\phi'_{0*}=\phi_{0*}, \phi'_{1*}=\phi_{1*}$.  It is well known %(see, for instance, \cite{Bla-Ktheory}{---Should it refer to Elliott, if we should refer at all?})
that there exist two unitaries $u_0, u_1\in F_2$ such that
$$u_0 \phi_0(a) u_0^*=\phi'_0(a), \quad a\in F_1\andeqn
%and
%$$
u_1 \phi_1(a) u_1^*=\phi'_1(a), \quad a\in F_1.$$
Since $U(F_2)$ is path connected, there is a unitary path $u:~ [0,1]\to U(F_2)$ with $u(0)=u_0$ and $u(1)=u_1$.  Define $\phi:~ A \to B$ by
$$\phi(f,a)=(g,c),$$
where $g(t)=u(t) f(t) u(t)^*$.
 Then a straightforward calculation shows that the map $\phi$ is a *-isomorphism.
\end{proof}

\begin{thm}\label{FG-Ratn}
Let $A=A(F_1, F_2, \phi_1, \phi_2)$ be in ${\cal C}.$  Then $K_0(A)_+$ is finitely generated by
its minimal elements; in other words, there is an integer $m\ge 1$ and finitely many
minimal projections of $M_m(A)$ such that these minimal projections generate the positive cone
$K_0(A)_+.$
%Then there are $$g_1, g_2, ..., g_n\in K_0(A)^+$$ such that for any
%$g\in K_0(A)^+$, there are positive real numbers $r_1, r_2, ..., r_n$ such that
%$$g=r_1g_1+ r_2g_2+\cdots+r_ng_n.$$
\end{thm}
\begin{proof}
We first show that
  $K_0(A)_+\setminus\{0\}$ has only finitely many minimal elements.

  Suppose otherwise that $\{q_n\}$ is an infinite set of minimal elements  of $K_0(A)_+\setminus \{0\}.$
Write $q_n=(m(1,n), m(2,n), ..., m({{l}}, n))\in \Z^l_+,$ where $m(i,n)$ are non-negative integers,
$i=1,2,...,l$ and $n=1,2,....$ If there is an integer $M\ge 1$ such that
$m(i,n)\le M$ for all $i$ and $n,$ then $\{q_n\}$ is a finite set. So we may assume
that $\{m(i,n)\}$ is unbounded for some $1\le i\le l.$
There is a subsequence of $\{n_k\}$ such that
$\lim_{k\to\infty} m(i,n_k)=\infty.$ To simplify the notation, without loss of generality, we may assume
that  $\lim_{n \to\infty}m(i,n)=\infty.$ We may assume that, for some $j,$ $\{m(j,n)\}$ is bounded.
Otherwise, by passing to a subsequence,  we may assume that $\lim_{n\to\infty}m(i,n)=\infty$ for
all $i\in \{1,2,...,l\}.$
Therefore  $\lim_{n\to\infty}m(i,n)-m(i,1)=\infty.$
It follows that, for some $n\ge 1,$
$m(i,n)>m(i,1)$ for all $i\in\{1,2,...,l\}.$ Therefore $q_n\ge q_1$ which contradicts the fact
that $q_n$ is minimal.
By passing to a subsequence, we may write $\{1,2,...,l\}=N\sqcup B$ such that
$\lim_{n\to\infty} m(i,n)=\infty$ if $i\in N$ and
$\{m(i,n)\}$ is bounded if $i\in B.$ Therefore $\{m(j,n)\}$ has only finitely many different values if $j\in B.$
Thus, by passing to a subsequence  again, we may assume that
$m(j,n)=m(j,1)$ if $j\in B.$ Therefore, for some $n>1,$
$m(i,n)>m(i,1)$ for all $n$ if $i\in N$ and $m(j,n)=m(j,1)$ for all $n$ if $j\in B.$
It follows that $q_n\ge q_1.$ This is impossible since $q_n$ is minimal.
This shows that $K_0(A)_+$ has only finitely many minimal elements.

To show that $K_0(A)_+$ is generated by these minimal elements,
fix an element $q\in K_0(A)_+\setminus \{0\}.$   It suffices to show that $q$ is a finite sum
of minimal elements in $K_0(A)_+.$ If $q$ is not minimal, consider
the set of all elements in $K_0(A)_+\setminus\{0\}$ which are smaller than $q.$
This set is finite.  Choose one which is minimal among them, say $p_1.$ Then $p_1$ is minimal  element
in $K_0(A)_+\setminus \{0\},$ otherwise there is one smaller than $p_1.$  Since $q$ is not minimal, $q\not=p_1.$
Consider $q-p_1\in K_0(A)_+\setminus\{0\}.$ If $q-p_1$ is minimal, then
$q=p_1+(q-p_1).$ Otherwise, we repeat the same argument to obtain a minimal element
$p_2\le q-p_1.$ If $q-p_1-p_2$ is minimal, then $q=p_1+p_2+(q-p_1-p_2).$
Otherwise we repeat the same argument. This process is finite. Therefore $q$ is a finite sum of
minimal elements in $K_0(A)_+\setminus \{0\}.$
%It is enough to show that the cone
%$$V_\Q:=\{\xi \in (\Q^+)^l:  (\phi_{0})_*\xi=(\psi_{1})_*\xi\ \}$$
%is generated by finitely many elements of $K_0(A)^+$.
%
%
%Denote by $u=(u_1, u_2, ..., u_l)$ the standard order unit of $K_0(A)$, and consider
%\begin{displaymath}%
%S=\{x=(x_1, x_2, ..., x_l)\in\R^l: \left\{
%\begin{array}{l}
%\left<x-u, u\right>=0,\\
%((\phi_{0})_*-(\psi_{1})_*) (x_1, x_2, ..., x_l)^{\mathbf T}=0\\
%x_i\geq 0,\quad i=1, ..., l,
%\end{array}
%\right. \}.
%\end{displaymath}
%Then $S$ is a section of the cone
%$$V_\R:=\{\xi \in (\R^+)^l:  (\phi_{0})_*\xi=(\psi_{1})_*\xi\ \}.$$
%Indeed, for any $x=(x_1, x_2, ..., x_l)\in V_\R $, since $u_i>0$, $1\leq i\leq l$, one has that $\left<x, u\right>> 0.$ Put
%$$\lambda=\frac{ \left<u, u\right>}{\left<x, u\right>}>0,$$ and a straightforward calculation shows that
%$$\left<\lambda x-u, u\right>=0.$$
%Therefore $S$ is a section of $V_\R$.
%
%On the other hand, the set $S$ is a convex  set with finitely many extreme points. Denote the extreme points by $g'_1, g'_2, ..., g'_n$ for some $n\in\N$. Then the cone $V_\R$ is generated by $g'_1, g'_2, ..., g'_n$ in $\R^l$. Since the entries of $u$ and $(\phi_{0})_*-(\psi_{1})_*$ are integers, the entries of each vector $g'_i$, $1\leq i\leq n$, are rational. Therefore the cone $V_\Q$ is generated by $g'_1, g'_2, ..., g'_n$ in $\Q^l$.
%
%Since $g'_1, g'_2, ..., g'_n$ are rational, there is $M\in \N$ such that $$Mg'_1, Mg'_2, ..., Mg'_n\in (\Z^+)^l.$$ Denote by $$g_i=Mg'_i,\quad 1\leq i\leq n.$$ It is clear that $g_i\in K_0(A)^+$ and $\{g_1, g_2, ..., g_n\}$ still generates the cone $V_\Q$.
\end{proof}

%\begin{thm}\label{FG-PosCone}
%The positive cone of the ordered $K_0$-group of an Elliott-Thomsen algebra is finitely generated.
%\end{thm}
%\begin{proof}
%Let $A=A(F_1, F_2, \phi_1, \phi_2)$ be an Elliott-Thomsen algebra. By Lemma \ref{FG-Ratn}, there are $$g_1, g_2, ..., g_n\in K_0(A)^+$$ such that for any
%$g\in K_0(A)^+$, there are positive real numbers $r_1, r_2, ..., r_n$ such that
%$$g=r_1g_1+ r_2g_2+\cdots+r_ng_n.$$
%Fix $\{g_1, g_2, ..., g_n\}$. Note that the set
%$$G_0:=\{\sum_{i=1}^n r_ig_i: r_i\in[0, 1]\}\cap\Z^l$$
%is finite, so is the set $G_0\cap K_0(A)^+$. Then, one asserts that the finite set $G_0\cap K_0(A)^+$ generates the cone $K_0(A)^+$.

%Let $g\in K_0(A)^+$ be arbitrary. By Lemma \ref{FG-Ratn}, there are positive rational numbers $r_1, r_2, ..., r_n$ such that
%$$g=s_1g_1+s_2g_2+\cdots+s_ng_n$$
%and hence
%$$g=\lfloor s_1\rfloor g_1+\lfloor s_2\rfloor g_2+\cdots+\lfloor s_n\rfloor g_n+\{ s_1\} g_1+\{ s_2\} g_2+\cdots+\{ s_n\} g_n,$$
%where $s=\lfloor s\rfloor+\{s\}$ and $\lfloor s\rfloor$ is the integer part of $s$.
%Note that
%$$\{ s_1\} g_1+\{ s_2\} g_2+\cdots+\{ s_n\} g_n\in K_0(A)\cap (\Z^+)^l= K_0(A)^+,$$ and hence
%$$\{ s_1\} g_1+\{ s_2\} g_2+\cdots+\{ s_n\} g_n\in G_0\cap K_0(A)^+,$$
%Since $g_i\in G_0\cap K_0(A)^+$, $1\leq i\leq n$, this implies that $g$ is in the cone generated by $G_0\cap K_0(A)^+$. Since $g$ is arbitrary, the cone $K_0(A)^+$ is generated by $G_0\cap K_0(A)^+$, as desired.
%\end{proof}

\begin{thm}\label{2Tg14}
The exponential rank of $A=A(F_1, F_2, \phi_0, \phi_1)$ is at most $3+\ep$.
\end{thm}

\begin{proof}
For each unitary $u\in U_0(A)$, by \ref{2Rg13}, one can write $u=ve^{ih}$, where $v=(g,h)$ with $v(0)=v(1) =1\in F_2$.  So we only need to prove the exponential rank of $v$ is at most $2+\ep$.  Consider $v$ as an element in $C_0((0,1), F_2)\widetilde{)}$ which defines an element $(s_1,s_2,..., s_k)\in \Z^k=K_1(C_0((0,1), F_2))$.  Since $[v]=0$ in $K_1(A)$, there are $(m_1,m_2, ..., m_l)\in \Z^l$ such that
$$
(s_1, s_2,...,s_k)=((\phi_1)_{*0}-(\phi_0)_{*0})((m_1,m_2,...,m_l)).
$$
Note that
$$
(\phi_0)_{*0}((R_1, R_2,...,R_l))=(\phi_1)_{*0}((R_1,R_2,...,R_l))=(r_1,r_2,...,r_k)=[{\bf1}_{F_2}]\in K_0(F_2).
$$
%$\left(\!\!
 %        \begin{array}{c}
   %        s_1 \\
     %      s_2 \\
       %    \vdots \\
          % s_k
        % \end{array}\!\!
       %\right)
       %=(\phi_{1*}-\phi_{0*})
      % \left(\!\!
        % \begin{array}{c}
          % m_1 \\
          % m_2 \\
          % \vdots \\
          % m_l
         %\end{array}\!\!
     %  \right)$.
      % Note that $\phi_{0*}
       %\left(\!\!
        % \begin{array}{c}
          % R_1 \\
         % R_2\\
          % \vdots \\
           %R_l
        % \end{array}\!\!
       %\right)
       %=\phi_{1*}
       %\left(\!\!
         %\begin{array}{c}
         %  R_1\\
         % R_2\\
          % \vdots \\
         %  R_l
       %  \end{array}\!\!
      % \right)
      % =
      % \left(\!\!
      %   \begin{array}{c}
       %    {r_1} \\
        %  {r_2}\\
         %  \vdots \\
         % {r_k}
        % \end{array}\!\!
      % \right)=[\mbox{\large \bf 1}_{F_2}]\in K_0(F_2)$.

Increasing $(m_1,m_2,\cd, m_l)$ by adding a positive multiple of ${\big (}R_1,R_2,...,R_l{\big )}$, we can assume $m_j\geq 0$ for all $j\in \{ 1, 2, ..., l\}$.  Let $a=(m_1P_1, m_2P_2, ..., m_lP_l)$, where
$$P_j=\left(
        \begin{array}{cccc}
          1 &   &   &   \\
            & 0 &  &  \\
           &  & \ddots &  \\
           &  &  & 0
        \end{array}
      \right)
      \in M_{\{ j\} }(\C)\subset F_1~.$$
Let $h$ be defined by
$$
h(t)=\left\{\begin{array}{lll}
              \phi_0(a)(1-2t), &  & 0\leq t\leq \frac12, \\
                &   &  \\
                \phi_1(a)(2t-1), &  & \frac12< t\leq 1.
            \end{array}
\right.
$$
Then $(h,a)$ defines a selfadjoint element in $A$.  One also has $e^{2\pi ih}\in C_0((0,1), F_2\widetilde{)}$, since $e^{2\pi ih(0)}=e^{2\pi ih(1)}=1$.  Furthermore, $e^{2\pi ih}$ defines
$(\phi_{1*}-\phi_{0*})((m_1,m_2,...,m_l))\in \Z^k$
      % \left(\!\!
        % \begin{array}{c}
          % m_1 \\
          % m_2 \\
          % \vdots \\
          % m_l
        % \end{array}\!\!
      % \right)\in \Z^k$
as an element in $K_1(C_0((0,1), F_2))=\Z^k$.  Let $w=v e^{-2\pi ih}$.  Then $w$ satisfies $w(0)=w(1)=1$ and $w\in C_0((0,1), F_2\widetilde{)}$ defines
$$(0,0,..., 0)\in K_1(C_0((0,1), F_2))$$
up to an approximation within  a sufficiently small  $\ep$, one can assume that $w=(w_1,w_2,\cd, w_k)$ such that for all $j=1,2, \cd, k$,

\begin{enumerate}
\item $w_j(t)=\diag{\big(}e^{2\pi ih^j_1(t)}, e^{2\pi ih^j_2(t)},...,e^{2\pi ih^j_{k_j}(t)}{\big )}$;
\item the numbers $e^{2\pi ih^j_1(t)}, e^{2\pi ih^j_2(t)},... ,e^{2\pi ih^j_{k_j}(t)}$ are distinct for all $t\in (0,1)$;
%
%\vspace{.1in}
%
\item $h_1^j(0)=h_2^j(0)=\cd= h_{k_j}^j(0)=0$.
\end{enumerate}
Since $w_j(1)=1$, one has that  $h_i^j(1)\in \Z.$ %$h_i^j(t)\in \Z$.

On the other hand, the unitary $w$ defines
$$ h^j_1(t)+h^j_2(t)+\cd +h^j_{k_j}(t)\in \Z\cong K_1(C_0((0,1), F_2^j))$$
which is zero by the property of $w$.  From \eqref{eqn-G-ast-01}, %1.9,
one has $h^j_{i_1}(1)-h^j_{i_2}(1)\leq 1$. This implies
$$ h^j_1(1)=h^j_2(1)=\cd=h^j_{k_j}(1)=0.$$ Hence $h={\big(} (h^1, h^2, ..., h^k), 0 {\big )}$ defines a {{selfadjoint}} element in $A$ and $w=e^{2\pi ih}$.
\end{proof}

\begin{NN}\label{2Rg15}
Let $A=A(F_1, F_2, \phi_0, \phi_1)\in {\cal C},$
where $F_1=M_{r_1}\oplus M_{r_2}\oplus\cdots \oplus M_{r_l}$ and
$F_2=M_{R_1}\oplus M_{R_2}\oplus \cdots \oplus M_{R_k}.$ Let us calculate the Cuntz semigroup of $A$.  It is well known the extreme points of $\mathrm{T}(A)$  are canonically one-to-one corresponding to the  %{is as same as the set of}
irreducible representations of $A$, which are given by
$$\coprod_{j=1}^k (0,1)_j \cup \{ \rho_1, \rho_2, ..., \rho_l\} = \mathrm{Irr}(A),$$
where $(0,1)_j$ is the same open interval $(0,1)$.  We use subscript $j$ to indicate the $j$-th copy.

The affine space $\mathrm{Aff}(\mathrm{T}(A))$ can be identified with the subset of
$$\bigoplus_{j=1}^k C([0,1]_j,\R)\oplus \underbrace{(\R\oplus\R\oplus\cd\oplus \R)}_{l~~ copies}$$
consisting of $(f_1,f_2,..., f_k, g_1, g_2,...,g_l)$ satisfying the condition
$$f_i(0)=\frac1{R_i}\sum_{j=1}^l a_{ij}g_j\cdot r_j \qq\qq \mbox{and} \qq\qq
f_i(1)=\frac1{R_i}\sum_{j=1}^l b_{ij}g_j\cdot r_j,$$
where $(a_{ij})_{k\times l}=\phi_{0*}$ and $(b_{ij})_{k\times l}=\phi_{1*}$ as in \ref{2Rg11}. %1.13.

For any selfadjoint $h\in (A\otimes {\cal K})_+$, one can define a map $D_h:~ \mathrm{Irr}(A) \to \Z_+\cup\{ \infty\}$ as for any $\pi\in \mathrm{Irr}(A)$
$$D_h(\pi)=\lim\limits_{n\to\infty}\mathrm{Tr}(\pi\otimes\id_{\cal K}(h^{\frac1n}))~,$$
where Tr is the unnormalized trace.  Then  we write $D_h=(D_h^1,D_h^2,...,D_h^k,D_h(\rho_1),D_h(\rho_2),...,D_h(\rho_l)),$ where $D_h^i=D_h|_{(0,1)_i},$ which satisfies the following conditions
\begin{enumerate}
\item $D_h^j$ is lower semi-continuous on each $(0,1)_j$,
\item $\liminf_{t\to 0}D_h^i(t) \geq \sum_j a_{ij} D_h(\rho_j)$ and $\liminf_{t\to 1}D_h^i(t) \geq \sum_j b_{ij} D_h(\rho_l)$.
\end{enumerate}
It is straight forward to verify that the image of the map $h\in (A\otimes {\cal K})_+ \to D_h$ is the subset of $\mbox{Map} (\mathrm{Irr}(A), \Z_+\cup\{ \infty\})$ consisting elements satisfying the above two conditions.

Note that $D_h(\pi)=\mbox{rank}(\pi\otimes \id_{\cal K}(h))$ for each $\pi\in \mathrm{Irr}(A)$ and $h\in (A\otimes {\cal K})_+$.

The following result is well known to experts (for example, see \cite{CEI-CuntzSG}).
\end{NN}

\begin{thm}\label{2Tg16}
 Let $A=A(F_1, F_2, \phi_0, \phi_1)\in {\cal C}$ and let $n\ge 1$ be an integer.

\noindent {\rm (a)} The following are equivalent:
\begin{enumerate}
\item $h\in M_n(A)_+$ is Cuntz equivalent to a projection;
\item $0$ is an isolate point in the spectrum of $h$;
\item  $D_h^j$ continuous on each $(0,1)_j$, $\liminf_{t\to 0}D_h^i(t) = \sum_j a_{ij} r_h^j$ and $\liminf_{t\to 1}D_h^i(t) = \sum_j b_{ij} r_h^j$.
\end{enumerate}

\noindent {\rm (b)} For $h_1, h_2 \in M_n(A)_+$,  $h_1\lesssim h_2$
%$h_1$ is Cuntz sub-equivalent to $h_2$ (denote by $h_1\prec h_2$ )
if and only if $D_{h_1}(\pi) \leq  D_{h_2}(\pi)$ for each $\pi\in \mathrm{Irr}(A)$. In particular, $A$ has strictly comparison for positive elements.

\end{thm}

\begin{proof}
For part (b), obviously, $h_1 {\lesssim} h_2$ implies $D_{h_1}(\phi)\leq D_{h_2}(\phi)$ for each $\phi\in Irr(A)$. Conversely assume that $h_1=(f,a)$ and $h_2=(g,b)$ satisfy that $D_{h_1}(\phi)\leq D_{h_2}(\phi)$ for each $\phi\in Irr(A)$.

First, there are strictly positive  {{functions}}
%in
$s_1,s_2\in C_0((0,1])$
such that $s_1(a)$ and $s_2(b)$ are projections in $F_1.$
Note that $s_i(h_i)$ are Cuntz equivalent to $h_i$ ($i=1,2$).
By replacing $h_i$ by $s_i(h_i),$ \wilog\, we may assume that $a$ and $b$ are projections.

Fix $1/4>\ep>0.$ There exists $\dt_1>0$ such that
\beq\label{N318-1}
\|f(t)-\phi_0(a)\|<\ep/64\,\,{\rm for\,\,all}\,\, t\in [0,2\dt_1]\andeqn \hspace{-0.05in}\|f(t)-\phi_1(a)\|<\ep/16\,\,{\rm for\,\,all}\,\, \hspace{-0.02in}t\in [1-2\dt_1,1].
\eneq
It follows that $f_{\ep/8}(f(t))$ is a projection in $[0,2\dt_1]$ and $[1-2\dt_1,1].$
Put $h_0=f_{\ep/8}(h_1).$  Note that $h_0=(f_{\ep/8}(f), a).$
 Then $D^i_{h_0}(\pi)$ is constant in $(0, 2\dt_1]\subset (0,1)_i$ and
$[1-2\dt_1,1]\subset (0,1)_i,$ $i=1,2,...,k.$
Choose $\dt_2>0$ such that
\beq\label{N318-2}
D_{h_2}^i(t)\ge d_{Tr_i}(\phi_0(b))\rforal t\in (0,2\dt_2]_i\andeqn D_{h_2}(t)\ge d_{Tr}(\phi_1(b))\rforal t\in [2\dt_2,1)_i,
\eneq
where $Tr_i$ is the standard trace on $M_{R_i},$ $i=1,2,...,k.$
%There is  $\dt>0$ with $\dt<\min\{\dt_1, \dt_2\}$ and a homeomorphism $s:  [\dt, 1-\dt] \to [0,1]$
%such that
%\beq\label{Lg17-n3}
%s|_{[2\dt_1, 1-2\dt_1]}={\rm id}|_{[2\dt_1, 1-2\dt_1]}\andeqn \|c_i-h_i\|<\ep/64.
%\eneq
%, \dt_0/16\}\andeqn \|f_{\ep/2}(c_i)-f_{\ep/2}(h_i)\|<\min\{\ep/16, \dt_0/16\},\,\,\,i=1,2,
%\eneq
%where $c_i(t)=h_i(s(t))$ if $t\in [\dt, 1-\dt]$ and $c_i(t)=h_i(0)$ for all $t\in [0,\dt]$ and
%$c_i(t)=h_i(1)$ for all $t\in [1-\dt,1].$
%It follows that
%\beq\label{Lg17-n4}
%f_{\ep/8}(h_0)\lesssim c_1.
%\eneq
Choose $\dt=\min\{\dt_1, \dt_2\}.$ Since $a\lesssim b$ in $F_1,$  there is a  unitary $u_e\in F_1$ such that
$u_e^*au_e=q\le b,$ where $q\le b$ is a sub-projection.
Let $u_0=\phi_0(u_e)$ and $u_1=\phi_1(u_e)\in F_2.$ Then one can find  a unitary $u\in A$ such that
$u(0)=u_0$ and $u(1)=u_1.$
\Wlog, by replacing $h_0$ by $u^*h_0u,$ we may assume that $a=q.$

 Let $\bt: [0,1] \to [0,1]$ be a continuous function which is $1$ on the boundary and $0$ on $[\dt, 1-\dt]$. Let $h_3=(f_2, b-a)$ with $f_2(t)=\bt(t)\phi_0(b-a)$ for $t\in [0,\dt]$, $f_2(t)=\bt(t)\phi_1(b-a)$ for $t\in [1-\dt,\dt]$, and $f_2(t)=0$ for  $t \in[\dt, 1-\dt]$. Define $h_1'=h_0+h_3.$  Note that $h_1'$ has the form
 $(f',b)$ for some $f'\in C([0,1], F_2).$ Let $\pi_i: F_2\to M_{R_i}$ be the projection map.
  Then
 \beq\label{N318-3}
 h_0\le h_1',\,\,\, D_{h_1'}^i(\pi)=D_{h_0}^i(\pi)\rforal \pi\in (\dt, 1-\dt)\subset (0,1)_i\\
  {\rm rank}( h_1'(\pi))\le {\rm rank} (h_0(t))+{\rm rank}(h_3(t))\le {\rm rank}(\pi_i(\phi_0(a)))+{\rm rank}(\pi_i(\phi_0(b-a)))\\
  ={\rm rank}(\pi_i(\phi_0(b)))\le D_{h_2}^i(\pi)\rforal \pi\in (0,\dt)\subset (0,1)_i\andeqn\\
{\rm rank}( h_1'(\pi))\le {\rm rank} (h_0(t))+{\rm rank}(h_3(t))\le {\rm rank}(\pi_i(\phi_1(a)))+{\rm rank}(\pi_i(\phi_1(b-a)))\\
  ={\rm rank}(\pi_i(\phi_1(b)))\le D_{h_2}^i(\pi)\rforal \pi\in (1-\dt,1)\subset (0,1)_i.
 \eneq
It follows that
\beq\label{N318-4}
D_{h_1'}(\pi)\le D_{h_2}(\pi)\rforal \pi\in Irr(A).
\eneq
It follows from  (\ref{N318-4}) and  Theorem 1.1 of \cite{RL0} that
$h_1'\lesssim h_2$  in $C([0,1], F_2).$ Since $C([0,1], F_2)$ has stable rank one, there is a unitary
$w\in C([0,1], F_2)$ such that $w^*h_1'w=h_4\in \overline{h_2C([0,1], F_2)h_2}.$
Since $h_1'=(f', b),$ $h_2=(g,b)$ and  $b$ is a projection, $h_4(0)=\phi_0(b)$ and $h_4=\phi_1(b).$
In other words, $h_4\in A.$ In particular, $h_4\in \overline{h_2Ah_2}$ and $h_4\lesssim h_2.$
Note that $w^*\phi_i(b)w=\phi_i(b),$ $i=0,1.$ Using a continuous path of unitaries
which commutes $\phi_i(b)$ ($i=0,1$) and connects to the identity, it is easy to find a sequence of unitaries $u_n\in C([0,1], F_2)$
such that $u_n(0)=u_n(1)=1_{F_2}$ such that
\beq\label{N318-5}
\lim_{n\to\infty}u_n^*h_1'u_n=h_4.
\eneq
Then $u_n\in A.$  We also have that $u_nh_4u_n^*\to h_1'.$ Thus $h_1'\lesssim h_4.$ It follows that
\beq\label{N318-5-1}
f_{\ep}(h_1)\lesssim h_0\lesssim h_1'\lesssim h_4\lesssim h_2
\eneq
for all $\ep>0.$  This implies that  $h_1\lesssim h_2$ and part (b) follows.

To prove part (a), we note that (1) and (2) are obviously equivalent and both implies (3).
That of (3) implies (1) follows from the computation of $K_0(A)$ in \ref{2Lg13} and part (b).
%{\color{Green} Why do we have two proofs? The equivalent of 1-2-3 in Part a follows from a result of Alin-Brown. Also, the statement of Theorem \ref{2Tg17} is weaker than that of \ref{2Tg16}.}
\end{proof}

\begin{lem}\label{cut-full-pj}
Let $C\in \mathcal C$, and let $p\in C$ be a projection. Then $pCp\in \mathcal C.$
%{\bf Does this part need $p$ to be full?}
Moreover, if $p$ is full and $C\in\mathcal C_0$, then $pCp\in\mathcal C_0.$
\end{lem}
\begin{proof}
We may assume that $C$ is not of finite dimensional.
Write $C=C(F_1, F_2, \phi_0, \phi_1)$. Denote by $p_e=\pi_e(p)$, where $\pi_e: C\to F_1$ is the map defined in %\ref{Cendmap}.
\ref{DfC1}.

For each $t\in [0,1],$ write
$\pi_t(p)=p(t)$ and ${\tilde p}\in C([0,1], F_2)$ such that $\pi_t({\tilde p})=p(t)$ for all $t\in [0,1].$
Then $\phi_0(p_e)=p(0)$, $\phi_1(p_e)=p(1),$ and
\begin{equation}\label{cfull-1}
pCp=\{(f,g)\in C: f(t)\in p(t)F_2p(t),
\andeqn g\in p_1F_1p_1\}.
\end{equation}

%Since $p$ is full, the projection $p_1$ is full in $F_1$. For each $t\in [0,1],$ Write
%$\pi_t(p)=p(t)$ and ${\tilde p}\in C([0,1], F_2)$ such that $\pi_t({\tilde p})=p(t)$ for all $t\in [0,1].$
%Then $\phi_0(p_1)=p(0)$ and $\phi_1(p_1)=p(1).$
%Moreover
%\beq\label{cfull-1}
%pCp=\{(f,g)\in C: f(t)\in p(t)F_2p(t),
%\andeqn g\in p_1F_1p_1\}.
%\eneq
%The assumption that $p$ is full implies that
%$p(t)$ is full in $F_2$ for each $t\in [0,1]$ and $p_1$ is full in $F_1.$
Put $p_0=p(0).$
There is a unitary $W\in C([0,1], F_2)$ such that
$W^*{\tilde p}W=p_0.$
Define $\Phi: {\tilde p}C([0,1], F_2){\tilde p}\to
C([0,1], p_0F_2p_0)$ by $\Phi(f)=W^*fW$ for all $f\in  {\tilde p}C([0,1], F_2){\tilde p}.$
 Put
 %$\Phi_t=\pi_t\circ \Phi,$
$F_1'=p_1F_1p_1$ and $F_2'=p_0F_2p_0.$
Define $\psi_0={\rm Ad}\, W(0) \circ \phi_0|_{F_1'}$ and $\psi_1={\rm Ad}\, W(1) \circ\phi_1|_{F_1'}.$ Put
$$
C_1=\{(f,g)\in C([0,1], F_2')\oplus F_1': f(0)=\psi_0(g)\andeqn f(1)=\psi_1(g)\},
$$
and note that $C_1\in\mathcal C$.
Define $\Psi: pCp\to C_1$ by
\beq\label{cfull-2}
\Psi((f,g))=(\Phi(f),g)\rforal f\in {\tilde p}C([0,1], F_2){\tilde p}\andeqn g\in F_1'.
\eneq
It is ready to verify that $\Psi$ is an isomorphism.

If $p$ is full and $C\in {\cal C}_0,$  then, by a result of Brown (\cite{Brown-Hereditary}), the hereditary {{\SCA\,}} $pCp$ is stably isomorphic to $C,$ and hence $K_1(pCp)=K_1(C)=\{0\}$; that is, $pCp\in \mathcal C_0$.
%Also consider the restriction of $p$ to $\mathrm{C}([0, 1], F_2)$, and denote it by $\tilde{p}_2$. Then $pCp$ is isomorphic to the mapping torus
%$$\{(f, g)\in (\tilde{p}_2(\mathrm{C}([0, 1], F_2))\tilde{p}_2)\oplus(p_1F_1p_1):\ f(0)=\phi_0(g),\ f(1)=\phi_1(g) \}.$$
%Note that $\tilde{p}_2(\mathrm{C}([0, 1], F_2))\tilde{p}_2$ is isomorphic to
%$\mathrm{C}([0, 1], p_2F_2p_2))$, where $p_2$ is a full projection of $F_2$. Thus,
%the C*-algebra $pCp$ has the form $C(p_1F_1p_1, p_2F_2p_2, \psi_0, \psi_1)$, where $\psi_0, \psi_1: p_1F_1p_1\to p_2F_2p_2$ are the map induced by $\phi_1$ and $\phi_2$, respectively. Moreover, since $p$ is full, one has that $K_1(pCp)\cong K_1(C)=\{0\}$, and hence $pCp\in \mathcal C_0$. (In fact, since $p_1$ and $p_2$ are full projections, one has that $[\phi_0]_0=[\psi_0]_0$ and $[\phi_1]_0=[\psi_1]_0$. Therefore the exponential map associated to $C(p_1F_1p_1, p_2F_2p_2, \psi_0, \psi_1)$ is the same as the exponential map associated to $C(F_1, F_2, \phi_0, \phi_1)$, and hence is surjective.)
\end{proof}

Class ${\cal C}$ and ${\cal C}_0$ are not closed under quotient. However, we have the following:

\begin{lem}\label{subapprox}
Any quotient of a C*-algebra in ${\mathcal C}$ (or in ${\cal C}_0$)  can be locally approximated by C*-algebras in ${\cal C}$ (or in ${\cal C}_0).$
%; any quotient of a C*-algebra in ${\mathcal C}_0$ can be locally approximated by C*-algebras in ${\cal C}_0$
More precisely,
let  $A\in {\cal C}$ (or $A\in {\cal C}_0$), let $B$ be a quotient of $A,$ let ${\cal F}\subset B$ be a finite set and let
$\ep>0,$ there exists a unital \SCA\, $B_0\subset B$ with $B_0\in {\cal C}$ (or $B_0\in {\cal C}_0$)  such that
\beq\nonumber
{\rm dist}(x, B_0)<\ep \rforal x\in {\cal F}.
\eneq
\end{lem}

\begin{proof}
Let $A\in\mathcal C.$ We may consider only those $A$'s which are not finite dimensional. Let $I$ be an ideal of $A$. Write $A=A(E, F, \phi_0, \phi_1)$, where $E=E_1\oplus\cdots\oplus E_l$, where $E_i\cong \mathrm{M}_{k_i}$, and  $F=F_1\oplus\cdots\oplus F_s$ with $F_j\cong\mathrm{M}_{m_j}$
Let $J=\{f\in C([0,1], F): f(0)=f(1)=0\}\subset A.$
%We may write $J=\bigoplus_{j=1}^sC_0((0,1),F_j).$ 
As before,
we may write $[0,1]_j$ for the the spectrum of the $j$-th summand of $C([0,1],F_j),$ whenever it is convenient.
Put $\phi_{i,j}=\pi_j\circ \phi_i: E\to F_j,$ where $\pi_j: F\to F_j$ is the quotient map, 
$i=0,1.$ 
%Then  $A/I$ {\it may}  be written  (with a re-indexing) as 
%those $(f,a),$  where  $a\in {\tilde E}$ and 
%${\tilde E}=\bigoplus_{i=1}^{l'} E_i$ ($l'\le l$), $f\in \bigoplus_{j=1}^{s'} C({\tilde I}_j, F_j),$
%where ${\tilde I}_j\subset [0,1]_j$ is a compact subset, $f(0_j)=\phi_{0,j}(a),$ 
%if $0_j\in {\tilde I}_j$ and 
%$f(1_j)=\phi_{1,j}(a)$ if $1_j\in {\tilde I}_j,$ $j=1,2,...,s$ (see \ref{2Rg15} for notation). 
%$A_0\oplus A_1$ with $A_1$  being the following form 
Then  $A/I$ {\it may}  be written  (with a re-indexing) as
%the spectrum of the quotient $S/I$ is a closed subset of the spectrum of $S$, and the C*-algebra $S/I$ has the fo
$$\{(f,a): a\in {\tilde E}, f\in \bigoplus_{1\le j\le s'} C({\tilde I}_j, F_j), f(0_j)={\tilde \phi}_{0,j}(a), \,{\rm if}\,\, 0_j\in {\tilde I}_j,
f(1_j)={\tilde \phi}_{1,j}(a), {\rm if}\,\, 1_j\in {\tilde I}_j\},
$$
where $s'\le s,$ $l\le l',$ ${\tilde E}=\bigoplus_{i=1}^{l'} E_i$  and ${\tilde \phi}_{i,j}=\phi_{i,j}|_{\tilde{E}}$ ($l'\le l$)  
and ${\tilde I}_j\subset [0,1]_j$ is a compact subset.
%%
%$$\{(f_1\oplus\cdots\oplus f_l', a): f_i\in \mathrm{C}(\tilde{I}_i, F_i), a\in \tilde{E}, f_1(0)\oplus\cdots\oplus f_l(0)={\tilde \phi}_0(a),\ f_1(1)\oplus\cdots\oplus f_{l'}(1)={\tilde \phi}_1(a) \},$$
%where $\tilde{I}_i=[0, 1]_i\setminus\bigcup_{n=1}^\infty (a_{i, n}, b_{i, n})$, and $(a_{i, n}, b_{i, n})$, $n=1, 2, ...$ are
%mutually disjoint open subintervals of $(0, 1)$, $\tilde{E}=E_{1}\oplus\cdots\oplus E_{l'}$
%(with $l'\le l$ and with a re-indexing) 
% ${\tilde E}$ is a quotient of $E$ (or a \SCA) and ${\tilde \phi}_i=(\phi_i)|_{{\tilde E}}$ ($i=0,1$), 
% and $A_0\cong \bigoplus_{j=1}^{m'}C(J_{j}, F_j'),$ where 
% $J_{j}\subset (0,1)$ is a compact subset and $F_j'$ is a simple finite dimensional \CA. 

It follows from \cite{Freu} that there is a sequence of $X_{n,j}$ which is a finite disjoint union of closed intervals
(includes points) such that ${\tilde I}_j$ is an inverse limit of $X_{n,j}$ and each  map $s_{n,j}: {\tilde I}_j\to X_{n,j}$ is {\it surjective}. Moreover $X_{n,j}$ can be identified  with disjoint union of closed subintervals of $[0,1].$  With this identification, 
we may further assume that $s_{n,j}(0)=0,$ if $0\in {\tilde I}_j$ and $s_{n,j}(1)=1,$ if $1\in {\tilde I}_j.$
Indeed, if ${\tilde I}_j$ contains an interval $[0,d_{n,j}'],$ we may assume that one of subintervals of $X_{n,j}$
is an interval with the form $[0, d_{n,j}]$ and $s_{n,j}(0)=0.$ Otherwise, there exist two sequences 
$\dt_{n,j,2}>\dt_{n,j,1} >0$ with $\lim_{n\to\infty}\dt_{n,j,2}=0$
% for each $j$ 
such that 
${\tilde I}_j\cap (\dt_{n,j,1},\dt_{n,j,2})=\emptyset.$ By redefining $s_{n,j}$ so that 
$s_{n,j}({\tilde I}_j\cap [0, \dt_{n,j,1}])=\{0\}$ and keeping $s_{n,j}|_{{\tilde I}_j\cap [\dt_{n,j,2},1]},$ we can assume 
that $s_{n,j}(0)=0.$ This can also be done at $1.$
Write $C({\tilde I}_j, F_j)=\overline{\cup_{n=1}^{\infty}(C(X_{n,j}, F_j))}.$ 
% and 
%$s_{n,j}^*: C(X_{n,j}, F_j)$ for the embedding, $j=1,2,...,s'.$
Let $s_n: \bigoplus_{j=1}^{s'} C(X_{n,j}, F_j)\to \bigoplus_{j=1}^{s'}C({\tilde I}_j, F_j)$  be 
the map induced by $s_{n,j}.$  Put $C_j=C(X_{n,j}, F_j).$ 
Then, for all sufficiently large $n,$  for each $f\in {\cal F},$ there is $g\in \bigoplus_{j=1}^{s'} C_j$ such that
\beq\label{1511-n2}
\|f|_{{\tilde I}_j}-s_{n,j}^*(g|_{X_{n,j}})\|<\ep/4,\,\,  g(0_j)=f(0_j),\, {\rm if}\,\, 0_j\in {\tilde I}_j \andeqn  g(1_j)=f(1_j),\,
{\rm if}\,\,\, 1_j\in {\tilde I}_j.
\eneq
%for all $f\in {\cal F}.$
%\in_{\ep/2} C_j\subset C({\tilde I}_j, F_j)
%\rforal f\in {\cal F}.
%\eneq
Note that $C_j$ is a unital \SCA\, of $C({\tilde I}_j, F_j),$ $j=1,2,...,s.$
%Similarly, $A_0=\lim_{n\to\infty}(\bigoplus_{j=1}^{m'} C(Y_{n,j}, F_j'), \imath_{n,j})$ 
%where each $Y_{n,j}$ is a disjoint union of close subintervals (including points) of $(0,1)$ 
%and $\imath_{n,j}$ is a embedding. Let $\imath_n=\bigoplus_{j=1}^{m'} \imath_{n,j}.$ 
%There exists $B=\bigoplus_{j=1}^{m'} C(Y_{n',j}, F_j')$ for some integer $n'$ and 
%$h\in B$ such that
%\beq\label{1511-sec3-1}
%\|f|_{J_i}-\imath_{n,j}(h|_{Y_{n,j}})\|<\ep/4\tforal f\in {\cal F}.
%\eneq
Define
$$
C=\{(f, a):  f\in \bigoplus_{j=1}^{s'}C_j, f(s_{n,j}(0_j))={\tilde \phi}_{0,j}(a),\,
{\rm if}\,\, 0_j\in {\tilde I}_j \andeqn f(s_{n,j}(1_j))={\tilde \phi}_{1,j}(a)\}.
$$
%\{(f_1\oplus\cdots\oplus f_s, a): f_i\in C_i, a\in \tilde{E}, f_1(0)\oplus\cdots\oplus f_l(0)=\phi_0(a),\ f_1(1)\oplus\cdots\oplus f_s(1)=\phi_1(a) \}
%which is a \SCA\, of $A/I.$
%{{Let}}  $Y_1,Y_2,...,Y_m$  be the connected components of $X_{n_j,j}$ which are isomorphic to closed interval %$[0,1]$
%such that $s_{n_j,j}(0)$ and $s_{n_j,j}(1)$ (in case $0$ or $1$ is not in ${\tilde I}_j$ we do not need
%consider that point ) are both not in those intervals.
Then $g$ in  (\ref{1511-n2}) is in $C$ and ${\cal F}\subset_{\ep} C.$  Moreover, $C\in {\cal C}$  and 
%\cong C_1\oplus C_2,$ where $C_1$ is a finite  direct sum of \SCA s which have finite dimension or  have the form
%$C(Y_i, F')\cong C([0,1], F')$ and $C_2\in {\cal C}.$ 
%Moreover, 
$C$ is isomorphic to a quotient of $A.$
%By (\ref{1511-n2}), 
%\beq\label{150103-p2}
%{\cal F}\subset_{\ep} C.
%\eneq
%Note also $C$ is also isomorphic to a quotient of $A.$
This proves the lemma in the case that $A\in {\cal C}.$

Now suppose that $A\in {\cal C}_0.$  
Since $C$ is isomorphic to a quotient of $A,$ 
%One realizes that $C$  constructed above is isomorphic
%to a quotient of $A.$ Therefore 
it suffices to show that, for any ideal $I\subset A,$ if $K_1(A)=\{0\},$ then
$K_1(A/I)=\{0\}.$
To this end, we consider the following six-term exact sequence:
$$
\xymatrix{
K_0(I) \ar[r] & K_0(A) \ar[r] & K_0(A/I) \ar[d] \\
K_1(A/I) \ar[u]^{\mathrm{\dt_1}} & K_1(A) \ar[l]  & K_1(I) . \ar[l]
}
$$

By \ref{2pg3}, $A$ and $A/I$ has stable rank one,  it follows from Proposition 4 of \cite{LR} that
$\dt_1=0.$ Since $K_1(A)=\{0\},$ it follows that  $K_1(A/I)=\{0\}.$ Lemma follows.
\end{proof}

%Finally, 
As the end of this section, we would like to return to the beginning of this section by stating the following proposition which will not be used.
\begin{prop}[Theorem 2.15 of \cite{ENSTW}]\label{ASCAs}
Let $A$ be a unital \CA\, which is a subhomogeneous \CA\, with one dimensional spectrum.
Then, for any finite subset ${\cal F}\subset A$ and any $\ep>0,$ there exists  unital  \SCA\, $B$ of $A$ which
is in ${\cal C}$ such that
\begin{equation*}
{\rm dist}(x, B)<\ep\tforal x\in {\cal F}.
\end{equation*}
\end{prop}

\begin{proof}
%We may assume that $A$ is a unital \SCA\, of $A_0=C([0,1], F),$ where $F$ is a finite dimensional \CA.
We use the fact that $A$ is an inductive limit of \CA s in ${\cal C}$ (\cite{ENSTW}).
Therefore, there is \CA\, $C\in {\cal C}$ and a unital \hm\, $\phi: C\to A$ such that
\begin{equation*}
{\rm dist}(x, \phi(C))<\ep/2\rforal x\in {\cal F}.
\end{equation*}
Then we apply Lemma \ref{subapprox}.
\end{proof}

%%%%%%%%%%%%%%%%%%%%%%%%%%%%%%%%%%%%%%%%%%%%%%%%%%%%%%%%%%

\section{Maps to finite dimensional \CA s}

\begin{lem}\label{8-N-0}
Let $z_1, z_2,..., z_n$ {be} positive integers which may not be distinct. There is a positive integer
{{$T=n\cdot \max\{z_iz_j: 1\le i,j\le n\}$}}
%depending on $(z_1, z_2,..., z_n)$
such that for any two nonnegative integer linear {combinations}  $a=\sum_{i=1}^n a_i\cdot z_i$ and $b=\sum_{i=1}^n b_i\cdot z_i$, there are two combinations $a'=\sum_{i=1}^n a'_i\cdot z_i$ and $b'=\sum_{i=1}^n b'_i\cdot z_i$ with $a'=b'$, $0\leq a'_i\leq a_i$, $0\leq b'_i\leq b_i$, and ${\rm min}\{a-a', b-b'\}\leq T$.

{Consequently, if $\dt>0$ and $|a-b|<\dt,$ we also have  that $\max\{a-a', b-b'\}<\dt+T.$}
\end{lem}

\begin{proof}
{To prove the first part, let} $T=n\cdot{\rm max}_{i,j}\{z_iz_j\}$. It is enough to prove that if  $a,b >T,$ 
%and $b>T$, 
then there are nonzero $0<a'=\sum_{i=1}^n a'_i\cdot z_i=b'=\sum_{i=1}^n b'_i\cdot z_i$ with $0\leq a'_i\leq a_i$, $0\leq b'_i\leq b_i$. But if both $a>T$ and $b>T$, then there are two (not necessary distinct) index $i, j$, with $a_i\ge z_j$ and $b_j\ge z_i$.
{Then choose} $a^{(1)'}=a_i'z_i$ and $b^{(1)'}=b_j'z_j$ {with $a_i'=z_j$ and $b_j'=z_i.$}
{If ${\rm min}\{a-a^{(1)'},b-b^{(1)'}\}\le T.$ Then we are done. If not, we repeat this on $a-a^{(1)}$ and $b-b^{(1)}$ and obtain
$a^{(2)}\le a-a^{(1)}$ and $b^{(2)}\le b-b^{(1)}$ such that $a^{(2)}=b^{(2)}.$  Put $a^{(2)'}=a^{(1)'}+a^{(2)}$ and
$b^{(2)'}=b^{(1)'}+b^{(2)}.$ Note we have $a^{(2)'}=b^{(2)'}$ and
we can also have $a^{(2)'}=\sum_ia_i^{(2)}z_i$ and $b^{(2)'}=\sum_ib^{(2)}_iz_i$ with
$0\le a_i^{(2)}\le a_i$ and $0\le b_i^{(2)}\le b_i'$ for all $i.$ If $\min\{a-a^{(2)'}, b-b^{(2)'}\}\le T,$ then we are done.
Otherwise, we continue.  An inductive argument shows the first part of lemma follows.}

To see the second part, assume that $a-a'\le T.$ Then $b-b'<|a-b|+T.$
\end{proof}

%{\bf {{I have a reference for the following.  But what would be the original reference?--L
%---I will check this tomorrow afternoon, if I got the time }} }

%The following is convent in this section and follows from Lemma 2.15 of \cite{Li-interval} directly.{\bf Checked?
%Could we also remove "easily"---I think I do have reference for this--L}

\begin{thm}\label{uniCMn} {\rm (see 2.10 of \cite{Lin-AU11}, Theorem 4.6 of \cite{Lntams07}  and 2.15 of \cite{Li-interval})}
 Let $X$ be a connected compact metric space, and {let $C= C(X)$}. Let $\mathcal F\subseteq C$ be a finite set, and let $\epsilon>0$ be a constant. There is a finite set $\mathcal H_1\subseteq C^+$ such that for any $\sigma_1>0$ there is a finite subset $\mathcal H_2\subseteq C$ and $\sigma_2>0$ such that for any  {unital }homomorphisms $\phi, \psi: C\to M_n$ for a matrix algebra $M_n$ satisfying
\begin{enumerate}
\item $\phi(h)>\sigma_1$ and $\psi(h)>\sigma_1$ for any $h\in\mathcal H_1$, and
\item $|{\rm tr}\circ \phi(h)-{\rm tr}\circ \psi(h)|<\sigma_2$ for any $h\in \mathcal H_2$,
\end{enumerate}
then there is a unitary $u\in M_n$ such that
$$\|\phi(f)-u^*\psi(f)u\|<\epsilon\quad \textrm{for any $f\in\mathcal F$}.$$
\end{thm}

We would also like to state another version of the above theorem.

\begin{thm}\label{Oldthm201408}
Let $A=C(X),$ where $X$ is  a compact metric space and let $\Delta:  A_+^{q, \bf 1}\setminus \{0\}\to (0,1)$ be an order preserving map.

For any $\ep>0,$ any finite set ${\cal F}\subset A,$   there exist a finite set
${\cal P}\subset \underline{K}(A),$ a finite set
${\cal H}_1\subset A_+^{\bf 1}\setminus \{0\},$
a finite set ${\cal H}_2\subset A_{s.a.}$ and  $\dt>0$ satisfying the following:
If $\phi_1, \phi_2: A\to M_n$ (for some integer $n\ge 1$) are two unital \hm s such that
\begin{eqnarray*}
&&[\phi_1]|_{\cal P}=[\phi_2]|_{\cal P},\\
&&\tau\circ \phi_1 (h)\ge \Delta(\hat{h})\tforal h\in {\cal H}_1\andeqn\\
&&|\tau\circ \phi_1(g)-\tau\circ \phi_2(g)|<\dt\tforal g\in {\cal H}_2,
\end{eqnarray*}
then, there exist   a unitary $u\in M_n$ such that
\beq\label{8-NN-1}
\|{\rm Ad}\, u\circ \phi_1(f)-\phi_2(f)\|<\ep\tforal f\in {\cal F}.
\eneq
\end{thm}

\begin{rem}\label{remark831KL}
{\rm
Let $X$ be a compact metric space and let $A=C(X).$ Suppose that $\phi_1, \phi_2: A\to M_n$ are  two  \hm s.
Then $(\phi_i)_{*1}=0,$ $i=1,2.$  Let ${\cal P}\subset \underline{K}(A)$ be a finite subset and
let $G$ be the subgroup generated by ${\cal P}.$
%Let ${\cal F}\subset C(X)$ be a finite subset.
There exists a finite CW complex $Y$ and a unital \hm\, $h: C(Y)\to C(X)$ such that
$G\subset [h](\underline{K}(C(Y))).$
Write $K_0(C(Y))=\Z^k\oplus {\rm ker}\rho_{C(Y)},$ where $\Z^k$ is generated by
mutually orthogonal projections $\{p_1,p_2,...,p_k\}$ which correspond to $k$ different
path connected components $Y_1, Y_2,..., Y_k$ of $Y.$  Fix  $\xi_i\in Y_i,$
Let $C_i=C_0(Y_i\setminus \{\xi_i\}),$ $i=1,2,...,k.$   Since $Y_i$ is path connected, by
considering the point-evaluation at $\xi_i,$
it is easy to see that, for any $\phi: C(Y)\to M_n,$ $[\phi]|_{\underline{K}(C_i)}=0.$
Let ${\bar G}=[h](\underline{K}(C(Y))).$ Suppose that
$\tau\circ \phi_1(p_i)=\tau\circ \phi_2(p_i),$ $i=1,2,...,k.$
Then, from the above, one computes that
\beq\label{KL-1408}
[\phi_1]|_{\cal P}=[\phi_2]|_{\cal P}.
\eneq
We will use this fact in the next proof.
}
\end{rem}

\begin{lem}\label{Aug-N-1}
Let $X$ be a compact metric space, let $F$ be a finite dimensional \CA\, and let $A=PC(X,F)P,$
where $P\in C(X,F)$ is a projection.
Let $\Delta: A_+^{q,{\bf 1}}\setminus \{0\}\to (0,1)$ be an order preserving map.

For any $\ep>0,$ any finite subset ${\cal F}\subset A$  and any $\sigma>0,$  there exists a finite subset
${\cal H}_1\subset A_+^{\bf 1}\setminus \{0\},$
a finite subset ${\cal H}_2\subset A_{s.a.}$ and  $\dt>0$  satisfying the following:
If $\phi_1, \phi_2: A\to M_n$ (for some integer $n\ge 1$) are two unital \hm s such that
\begin{eqnarray*}%\label{8-N-1-1}
\tau\circ \phi_1 (h)\ge  \Delta(\hat{h})\tforal h\in {\cal H}_1,\andeqn\\
|\tau\circ \phi_1(g)-\tau\circ \phi_2(g)|<\dt\tforal g\in {\cal H}_2,
\end{eqnarray*}
then, there exist  a projection $p\in M_n,$  a unital \hm\, $H: A\to pM_np, $
unital \hm s $h_1, h_2: A\to (1-p)M_n(1-p)$ and a unitary $u\in M_n$ such that
\begin{eqnarray*}%\label{8-N-1-2}
&&\|{\rm Ad}\, u\circ \phi_1(f)-(h_1(f)+H(f))\|<\ep,\\
&& \|\phi_2(f)-(h_2(f)+H(f)\|<\ep\tforal f\in {\cal F},\\
&&\andeqn
\tau(1-p)<\sigma,
\end{eqnarray*}
where $\tau$ is the tracial state of $M_n.$

\end{lem}

\begin{proof}
%We will prove only the case that $A$ has infinite dimension.
%The case that $A$ is finite dimensional  is easy to prove and we will leave it to the reader.

We first prove the case that $A=C(X).$
%Note we now assume that $X$ has infinitely many points.

Let $\Delta_1=(1/2)\Delta.$ Let ${\cal P}\subset \underline{K}(A)$ be a finite set,
${\cal H}_1'\subset A_+^{\bf 1}\setminus \{0\}$
(in place of ${\cal H}_1$)
be a finite set, ${\cal H}_2'\subset  A_{s.a.}$
(in place of ${\cal H}_2$)  be a finite set
and $\dt_1>0$ (in place of $\dt$) required by \ref{Oldthm201408} for
$\ep/2$ (in place of $\ep$), ${\cal F}$ and  $\Delta_1.$

%Let ${\bar G}$ be the subgroup generated by ${\cal P}.$
%There is an integer $k_0\ge 1$ such that
%\beq\label{oldthm1408-p}
%{\bar G}\cap K_*(A, \Z/k\Z)=\{0\}\rforal k\ge k_0.
%\eneq
%Choose $K=k_0!.$
%Choose an integer $N_0\ge 1$ such that $1/N_0<\min\{\sigma\cdot \dt_1, 1\}/16.$
%
%Since $X$ has infinitely many points, by choosing larger ${\cal H}_1',$ we may assume
%that $n\ge N_0^2(K+1)^2.$

Without loss of generality, we may assume that $1_A\in {\cal F},$ $1_A\in {\cal H}_1'\subset {\cal H}_2'$ and
${\cal H}_2'\subset A_+^{\bf 1}\setminus \{0\}.$
So, in what follows, ${\cal H}_2'\subset A_+^{\bf 1}\setminus \{0\}.$
Put
\beq\label{8-31-1}
\sigma_0=\min\{\Delta_1(\hat{g}): g\in {\cal H}_2'\}.
\eneq

%Let ${\cal P}_0={\cal P}\cap K_0(A)$ and $G_0$ be the subgroup generated by ${\cal P}_0.$
Let $G$ be the subgroup generated by ${\cal P}$ and let ${\bar G}$ be defined in \ref{remark831KL}.
Let ${\cal P}_0$ be a set of generators of ${\bar G}\cap K_0(A).$
Without loss of generality,  we may assume that
${\cal P}_0=\{p_1,p_2,...,p_{k_1}\}\cup\{z_1,z_2,...,z_{k_2}\},$ where
$p_i\in C(X)$ are projections (corresponding to clopen subsets) and $z_j\in {\rm ker}\rho_A(K_0(A)).$

%There is $d_0>0$ such that
%\beq\label{oldthm1408-1}
%|f(x)-f(x')|<\min\{\ep/16, \sigma_0/4\} \rforal f\in {\cal F}\cup {\cal H}_2'.
%\eneq
%if $x, x'\in X$ and ${\rm dist}(x, x')<d_0.$

%A subset $Y$ of $X$ is said to be $d_0/2$-connected, if for any two points $x, x'\in Y,$ there are
%$x_1, x_2,...,x_m\in Y$ such that
%${\rm dist}(x, x_1)<d_0/2,$ ${\rm dist}(x_m, x')<d_0/2$ and ${\rm dist}(x_i, x_{i+1})<d_0/2,$
%$i=1,2,...,m-1.$
%Let $Y_i$ be the clopen subset corresponding to the projection $p_i,$ $i=1,2,...,k_1.$
Without loss of generality, we may assume that
%$Y_i$ are $d_0/2$-connected and
$\{p_i: 1\le i\le k_1\}$ is a set of mutually orthogonal
projections such that $1_A=\sum_{i=1}^{k_1}p_i.$

%There are disjoint  subset $Z_1, Z_2,...,Z_m$ of $X$ such that each $Z_j$ contains
%at least one non-empty open balls $O_j\subset Z_j$ such that the closure  of $O_j$ contained
%in an open ball $O_j'\subset Z_j,$
%${\rm diam}(Z_j)<d_0/2,$ $j=1,2,...,m,$ and $\sqcap_{j=1}^m Z_j=X.$
%There are
%Choose $h_j\in A_+^{\bf 1}\setminus \{0\}$ such that $h_j(x)=1$ if $x\in O_j,$
%and $h_j(x)=0$ if $x\not\in Z_j.$
%Put ${\cal H}_1''=\{h_j: 1\le j\le m\}.$

%Since $Y_i$'s are clopen subsets, we may assume that each $Z_j$ is a subset of $Y_i$ for some $i.$

Let ${\cal H}_1={\cal H}_1'\cup\{p_i: 1\le i\le k_1\}$
%\cup {\cal H}_1''$
and ${\cal H}_2={\cal H}_2'\cup {\cal H}_1.$
Let $\sigma_1=\min\{\Delta(\hat{g}): g\in {\cal H}_2\}.$
Choose $\dt=\min\{\sigma_0\cdot \sigma/4k_1, \sigma_0\cdot \dt_1/4k_1, \sigma_1/16k_1\}.$
%\sigma_0/4k_1, 1/N_0^2k_1, \sigma_1/16k_1\}.$

Suppose now that $\phi_1, \phi_2: A\to M_n$ are two unital \hm s described in the lemma
for the above ${\cal H}_1,$ ${\cal H}_2$ and $\Delta.$

We may write $\phi_j(f)=\sum_{k=1}^{n}f(x_{k,j})q_{k,j}$ for all $f\in C(X),$
where $\{q_{k,j}: 1\le k\le n\}$ ($j=1,2$) is a set of mutually orthogonal rank one projections
and $x_{k,j}\in X.$
We have
\beq\label{8-31-2}
|\tau\circ \phi_1(p_i)-\tau\circ \phi_2(p_i)|<\dt,\,\,\,i=1,2,...,k_1,
\eneq
where $\tau$ is the tracial state on $M_n.$
Therefore, there exists a projection $P_{0,j}\in M_n$ such that
\beq\label{8-31-3}
\tau(P_{0,j})<k_1\dt<\sigma_0\cdot \sigma,\,\,\, j=1,2,
\eneq
${\rm rank}(P_{0,1})={\rm rank}(P_{0,2}),$
unital \hm s $\phi_{1,0}: A\to  P_{0,1}M_nP_{0,1}, $ $\phi_{2,0}: A\to P_{0,2}M_nP_{0,2},$
$\phi_{1,1}: A\to (1-P_{0,1})M_n(1-P_{0,1})$ and $\phi_{1,2}: A\to (1-P_{0,2})M_n(1-P_{0,2})$
such that
\beq\label{8-31-4}
\phi_1=\phi_{1,0}\oplus \phi_{1,1},\,\,\, \phi_2=\phi_{2,0}\oplus \phi_{2,1},\\\label{n151110-1}
\tau\circ \phi_{1,1}(p_i)=\tau\circ \phi_{1,2}(p_i),\,\,\,i=1,2,...,k_1.
\eneq
By replacing $\phi_1$ by ${\rm Ad}\, v\circ \phi_1,$ simplifying the notation, without loss of generality,
we may assume that $P_{0,1}=P_{0,2}.$
It follows   from (\ref{n151110-1}) that (see \ref{remark831KL}) 
\beq\label{8-31-5}
[\phi_{1,1}]|_{\cal P}=[\phi_{2,1}]|_{\cal P}.
\eneq
By (\ref{8-31-3}) and the choice of $\sigma_0,$ we also have
\beq\label{8-31-6}
\tau \circ \phi_{1,1}(g) \ge \Delta_1(\hat{g})\rforal g\in {\cal H}_1'\andeqn\\
|\tau \circ \phi_{1,1}(g)-\tau\circ \phi_{1,2}(g)|<\sigma_0\cdot \dt_1\rforal g\in {\cal H}_2'.
\eneq
Therefore
\beq\label{8-31-6+}
t \circ \phi_{1,1}(g) \ge \Delta_1(\hat{g})\rforal g\in {\cal H}_1'\andeqn\\
|t \circ \phi_{1,1}(g)-t \circ \phi_{1,2}(g)|<\dt_1\rforal g\in {\cal H}_2',
\eneq
where $t$ is the tracial state on $(1-P_{1,0})M_n(1-P_{1,0}).$
By applying \ref{Oldthm201408},
there exists a unitary  $v_1\in (1-P_{1,0})M_n(1-P_{1,0})$ such that
\beq\label{8-31-7}
\|{\rm Ad}\, v_1\circ \phi_{1,1}(f)-\phi_{2,1}(f)\|<\ep/16\rforal f\in {\cal F}.
\eneq
Put $H=\phi_{2,1}$ and $p=P_{1,0}.$ The lemma for the case that $A=C(X)$ follows.

%We may write $\phi_{j,1}(f)=\sum_{k=1}^{N_1}f(x_{k,j})q_{k,j}$ for all $f\in C(X),$
%where $\{q_{k,j}: 1\le k\le n\}$ ($j=1,2$) is a set of mutually orthogonal rank one projections
%and $x_{k,j}\in X$ and $N_1\le n.$
%
%The condition that
%\beq\label{8-31-5}
%\tau\circ \phi_j(g)\ge \Delta(\hat{g})\rforal g\in {\cal H}_1 \,(j=1,2)
%\eneq
%implies that
%\beq\label{8-31-6}
%\tau\circ \phi_{1,1}(h_j),\,\,\tau\circ \phh_{2,1}(h_j)>\sigma_1/4
%\eneq
%Let $M_j$ be the number of points $x_{k,1}$ (counting multiplicity) in $Z_j,$ $j=1,2,...,m.$
%Then, by (\ref{8-31-6}),
%$M_j>N_0(K+1)^2.$
%Write $M_j=s_jK+r_j$ for some $0\le r_j<K.$ Then $s_j\ge N_0,$ $j=1,2,...,m.$
%Denote by
%Define $\psi_1: A\to (1-P_{0,1})M_n(1-P_{0,1})$ by
%$\psi_1

For the case that $A=M_r(C(X)),$ let $e_1\in M_r$ be a rank one projection.
Put $B=e_1(M_r(C(X)))e_1\cong C(X).$ Consider $\psi_1=\phi_1|_B$ and $\psi_2=\phi_2|_B.$
Since $\phi_1$ and $\phi_2$ are unital, ${\rm rank}(\psi_1(1_B))={\rm rank}(\psi_2(1_B)).$ By replacing $\psi_1$
by ${\rm Ad}\, v\circ \psi_1$ for some unitary $v\in M_n,$ we may assume that $\psi_1(1_B)=\psi_2(1_B).$
Then  what has been proved could be applied to $\psi_1$ and $\psi_2.$
The case  of $A=M_n(C(X))$ then  follows.

For the case $A=C(X,F),$ let $F=M_{r_1}\oplus M_{r_2}\oplus \cdots \oplus M_{r_k}.$
Denote by $E_i$ the projection of $A$ corresponding to the identity of $M_{r_i},$ $i=1,2,...,k.$
The same argument used the above shows that we can find a projection $Q_j\in M_n$ ($j=1,2$),
unital \hm s $\psi_{j,0}: A\to Q_jM_nQ_j$ and unital \hm s $\psi_{j,1}: A\to (1-Q_j)M_n(1-Q_j)$ such that
\beq\label{8-31-101}
{\rm rank}(Q_1)={\rm rank}(Q_2),\,\,\,
\tau(Q_j)<\dt,\\
\tau(\psi_{j,1}(E_i))=\tau(\psi_{j,2}(E_i)),\,\,\,i=1,2\andeqn
\phi_j=\psi_{j,0}\oplus \psi_{j,1},
\eneq
$j=1,2.$ By replacing $\phi_1$ by ${\rm Ad}\, v\circ \phi_1$ for a suitable unitary $v\in M_n,$ we may assume
$Q_1=Q_2$ and $\psi_{1,1}(E_i)=\psi_{2,1}(E_i),$ $i=1,2,...,k.$ Thus this case has been reduced to the case
that $A=M_r(C(X)).$
Therefore  the case that $A=C(X,F)$ also follows from the above proof.

Let us consider the case that $A=PC(X,F)P.$  Note that $\{{\rm rank}(P(x)): x\in X\}$ is a finite set of positive integers.
Therefore the set $Y=\{x\in X: {\rm rank}(P(x))>0\}$ is compact and open.  Then we may write  $A=PC(Y,F)P.$
Thus we may assume that $P(x)>0$ for all $x\in X.$ This implies easily that $P$ is a full projection of $C(Y,F).$
Then,  by \cite{Brown-Hereditary}, $A\otimes {\cal K}\cong C(Y,F)\otimes {\cal K}.$
It follows that there is an integer $m\ge 1$ and a projection $e\in M_m(A)$ such that
$eM_m(A)e\cong C(X,F).$ By extending $\phi_1$ and $\phi_2$ to maps from $eM_m(A)e,$ one easily sees
that the proof reduces to the case that $A=C(X,F).$
\end{proof}

\begin{cor}\label{Aug-N-2}
Let $X$ be a compact metric space, let $F$ be a finite dimensional \CA\, and let $A=PC(X,F)P,$
where $P\in C(X,F)$ is a projection.
Let $\Delta: A_+^{q,{\bf 1}}\setminus\{0\}\to (0,1)$ be an order preserving map.  Let $1>\af>1/2.$

For any $\ep>0,$ any finite subset ${\cal F}\subset A,$  any finite subset
${\cal H}_0\subset A_+^{\bf 1}\setminus \{0\}$ and any integer $K\ge 1.$   There is an integer $N\ge 1,$  a finite subset
${\cal H}_1\subset A_+^{\bf 1}\setminus \{0\},$
a finite subset ${\cal H}_2\subset A_{s.a.},$ $\dt>0$
%and  a finite subset ${\cal P}\subset A$ of projections
satisfying the following:
If $\phi_1, \phi_2: A\to M_n$ (for any integer $n\ge N$) are two unital \hm s such that
\begin{eqnarray*}%\label{8-N-1-1}
%[\phi_1(p)]=[\phi_2(p)]\tforal p\in {\cal P},\\
\tau\circ \phi_1 (h)\ge \Delta(\hat{h})\tforal h\in {\cal H}_1, {\tand}\\
|\tau\circ \phi_1(g)-\tau\circ \phi_2(g)|<\dt\tforal g\in {\cal H}_2,
\end{eqnarray*}
then, there exist  mutually orthogonal nonzero projections $e_0, e_1,e_2,...,e_K\in M_n$ such that
$e_1, e_2,...,e_K$ are equivalent, $e_0\lesssim e_1$   and $e_0+\sum_{i=1}^Ke_i=1_{M_n},$ and there are unital \hm s
$h_1, h_2: A\to e_0M_ne_0,$  $\psi: A\to e_1M_ne_1$
 and a unitary $u\in M_n$ such that
\begin{eqnarray*}%\label{8-N-1-2}
&&\|{\rm Ad}\, u\circ \phi_1(f)-(h_1(f)+{\rm diag}(\overbrace{\psi(f), \psi(f),...,\psi(f)}^K))\|<\ep,\\
 &&\|\phi_2(f)-(h_2(f)+{\rm diag}(\overbrace{\psi(f),\psi(f),...,\psi(f)}^K))\|<\ep\tforal f\in {\cal F},\\
&&\andeqn \tau\circ \psi(g)\ge \af{\Delta(\hat{g})\over{K}}\tforal g\in {\cal H}_0,
\end{eqnarray*}
where $\tau$ is the tracial state of $M_n.$

\end{cor}

\begin{proof}

By applying \ref{Aug-N-1}, it is easy to see that it suffices to prove the following statement:

Let $X,F,P$  $A$  and  $\af$ be as in the corollary.

Let $\ep>0,$  let ${\cal F}\subset A$ be a finite subset, let ${\cal H}_0\subset A_+^{\bf 1}\setminus \{0\}$ and let $K\ge 1.$  There is an integer $N\ge 1,$ a finite subset
${\cal H}_1\subset A_+^{\bf 1}\setminus \{0\}$ satisfying the following:
Suppose that $H: A\to M_n$ (for some $n\ge N$) is a unital \hm\, such that
\beq\label{8-N-1-3}
\tau\circ H(g)\ge \Delta(\hat{g})\rforal g\in {\cal H}_0.
\eneq
Then  there  are mutually orthogonal projections $e_0, e_1,e_2,...,e_K\in M_n,$ a unital \hm\, $\phi: A\to e_0M_ne_0$
and a unital \hm\, $\psi: A\to e_1M_ne_1$ such that
\beq\label{8-N-1-4}
\|H(f)-(\phi(f)\oplus {\rm diag}(\overbrace{\psi(f), \psi(f),...,\psi(f)}^K))\|<\ep\rforal f\in {\cal F},\\
\tau\circ \psi(g)\ge \af\Delta(\hat{g})/K\rforal g\in {\cal H}_0.
\eneq

We make one further reduction:  Using the argument at the end of the proof of \ref{Aug-N-1}, it suffices to prove
the above statement for $A=C(X).$

Put
\beq\label{8-N-1-5}
\sigma_0= ((1-\af)/4)\min\{\Delta(\hat{g}): g\in {\cal H}_0\}>0.
\eneq
Let $\ep_1=\min\{\ep/16, \sigma_0\}$ and let ${\cal F}_1={\cal F}\cup {\cal H}_0.$
Choose $d_0>0$ such that
\beq\label{8-N-1-6}
|f(x)-f(x')|<\ep_1\rforal f\in {\cal F}_1,
\eneq
provided that $x,x'\in X$ and ${\rm dist}(x,x')<d_0.$

Choose $\xi_1,\xi_2,...,\xi_m\in X$ such that $\cup_{j=1}^m B(\xi_j,d_0/2)\supset X,$
where $B(\xi, r)=\{x\in X: {\rm dist}(x, \xi)<r\}.$
There is $d_1>0$ such that $d_1<d_0/2$ and
\beq\label{8-N-1-7}
B(\xi_j, d_1)\cap B(\xi_i, d_1)=\emptyset,
\eneq
if $i\not=j.$
There is, for each $j,$  a function $h_j\in C(X)$ with $0\le h_j\le 1,$
$h_j(x)=1$ if $x\in B(\xi_j,d_1/2)$ and $h_j(x)=0$ if $x\not\in B(\xi_j,d_1).$
Define ${\cal H}_1={\cal H}_0\cup \{h_j: 1\le j\le m\}$ and put
\beq\label{8-N-1-8}
\sigma_1=\min\{\Delta(\hat{g}): g\in {\cal H}_1\}.
\eneq
Choose an integer $N_0\ge 1$ such that $1/N_0<\sigma_1 \cdot (1-\af)/4$ and $N=4m(N_0+1)^2(K+1)^2.$

Now let $H: C(X)\to M_n$ be a unital \hm\, with $n\ge  N$ satisfying the assumption (\ref{8-N-1-3}).
Let $Y_1=\overline{B(\xi_1, d_0/2)}\setminus \cup_{i=2}^m B(\xi_i, d_1),$
$Y_2=\overline{B(\xi_2,d_0/2)}\setminus (Y_1\cup\cup_{i=3}^m B(\xi_i, d_1),$
$Y_j=\overline{B(\xi_j, d_0/2)}\setminus (\cup_{i=1}^{j-1} Y_i\cup \cup_{i=j+1}^m B(\xi_i, d_1)),$
$j=1,2,...,m.$ Note that $Y_j\cap Y_i=\emptyset$ if $i\not=j$ and
$B(\xi_j,d_1)\subset Y_j.$
We write that
\beq\label{8-N-1-9}
H(f)=\sum_{i=1}^n f(x_i)p_i=\sum_{j=1}^m (\sum_{x_i\in Y_j} f(x_i)p_i)\rforal f\in C(X),
\eneq
where $\{p_1,p_2,...,p_n\}$ is a set of mutually orthogonal rank one projections in $M_n,$
$\{x_1,x_2,...,x_n\}\subset X.$
Let $R_j$ be  the cardinality of $\{x_i: x_i\in Y_j\}.$
Then, by (\ref{8-N-1-3}),
\beq\label{8-N-1-10}
R_j\ge N\tau\circ H(h_j)\ge N\Delta(\hat{h_j})\ge (N_0+1)^2 K\sigma_1\ge (N_0+1)K^2,\,\,\,j=1,2,...,m.
\eneq
Write $R_j=S_jK+r_j,$ where $S_j\ge N_0Km$ and $0\le r_j<K,$ $j=1,2,...,m.$
Choose $x_{j,1},x_{j,2},...,x_{j,r_j}\subset \{x_i\in Y_j\}$ and denote
$Z_j=\{x_{j,1},x_{j,2},...,x_{j,r_j}\},$ $j=1,2,...,m.$

Therefore we may write
\beq\label{8-N-1-11}
H(f)=\sum_{j=1}^m (\sum_{x_i\in Y_j\setminus Z_j} f(x_i)p_i)+\sum_{j=1}^m(\sum_{i=1}^{r_j}f(x_{j,i})p_{j,i})
\eneq
for $f\in C(X).$ Note that the cardinality of $\{x_i\in Y_j\setminus Z_j\}$ is $KS_j,$ $j=1,2,...,m.$
Define
\beq\label{8-N-1-12}
\Psi(f)=\sum_{j=1}^m f(\xi_j)P_j=\sum_{k=1}^K (\sum_{j=1}^m f(\xi_j)Q_{j,k})\rforal f\in C(X),
\eneq
where $P_j=\sum_{x_i\in Y_j\setminus Z_j}p_i=\sum_{k=1}^KQ_{j,k}$ and
${\rm rank}\,Q_{j,k}=S_j,$ $j=1,2,...,m.$
Put $e_0=\sum_{i=1}^m(\sum_{i=1}^{r_j}p_{j,i}),$
$e_k=\sum_{j=1}^m Q_{j,k},$ $k=1,2,...,K.$
Note that
\beq\label{8-N-1-13}
{\rm rank}(e_0)=\sum_{j=1}^mr_j<mK\andeqn {\rm rank}(e_k)=S_j\\
S_j\ge N_0mK>mK,\,\,\,j=1,2,..,K.
\eneq
It follows that $e_0\lesssim e_1$ and $e_i$ is equivalent to $e_1.$
Moreover, we may write
\beq\label{8-N-1-14}
\Psi(f)={\rm diag}(\overbrace{\psi(f), \psi(f),...,\psi(f)}^K)\rforal f\in A,
\eneq
where $\psi(f)=\sum_{j=1}^m f(\xi_j)Q_{j,1}$ for all $f\in A.$
We also estimate that
\beq\label{8-N-1-15}
\|H(f)-(\phi(f)\oplus {\rm diag}(\overbrace{\psi(f), \psi(f),...,\psi(f)}^K))\|<\ep_1\rforal f\in {\cal F}_1.
\eneq
We also compute that
\beq\label{8-N-1-16}
\tau\circ \psi(g)\ge (1/K)(\Delta(\hat{g})-\ep_1-{mK\over{N_0Km}} )\ge \af{\Delta(\hat{g})\over{K}}
%\tau\circ \psi(g)\ge (1/K)(\{\Delta(\hat{g}): g\in {\cal H}_0\}-\ep_1-{mK\over{N_0Km}} )\ge \af{\Delta(\hat{g})\over{K}}
\eneq
for all $g\in {\cal H}_0.$
\end{proof}

\begin{rem}\label{8-N-2-r}
{\rm
If we also assume that $X$ has infinitely many points, then Lemma \ref{Aug-N-2} holds without
mentioning the integer $N.$ This can be seen by taking larger ${\cal H}_1$ which  will force the integer
$n$ large.
}
\end{rem}

\begin{df}\label{8-N-3}
Denote by ${\cal {\bar D}}_0'$ the class of all
finite dimensional
 \CA s.
%with the form
%$PM_r(C(X))P,$ where $X$ is a compact metric space with finitely many path connected components, $r\ge 1$ is an integer and $P\in M_r(C(X))$ is a projection.
For $k\ge 1,$ denote by ${\cal {\bar D}}_k'$ the class of all \CA s with the form:
$$
A=\{(f,a)\in PC(X, F)P\oplus B: f{|_{Z}}=\Gamma(a)\},
$$
where $X$ is a compact  metric space,  $F$ is a finite dimensional \CA, $P\in C(X,F)$ is a projection, { $Z\subset X$ is a nonempty proper subset of $X,$
%pace with finitely many path connected components such that $X\setminus Z$ is path connected,
$B\in {\cal {\bar D}}_{k-1}'$ }
and $\Gamma: B\to PC({ Z},F)P$ is a unital \hm,
where we assume that there is $d_{X,{ Z}}>0$ { (we  will denote it by $d_X$ if there is no danger of confusion)} such that, for any $0<d\le d_X,$
there exists $s_*^d: \overline{X^d}\to { Z}$ such that
\beq\label{Amnotation-1}
s_*^d(x)=x\rforal x\in Z \andeqn \lim_{d\to 0} \|f|_{ Z}\circ s_*^d-f|_{\overline{X^d}}\|=0\rforal f\in C(X,F),
\eneq
where $X^d=\{x\in X: {\rm dist}(x, { Z})< d\}.$
We also assume that,  for any $0<d<d_X/2$ and for any $d>\dt>0,$
there is { a homeomorphism}  $r: X\setminus X^{d-\dt}\to X\setminus X^{d}$ such that
\beq\label{Amnotation-2+}
{\rm dist}(r(x), x)<\dt\rforal x\in X\setminus X^{d-\dt}.
\eneq
For the convenience, we will also assume, in addition, that
$F=M_{s(1)}\bigoplus M_{s(2)}\cdots \bigoplus M_{s(k_1)},$
$$
PC(X,F)P =\bigoplus_{j=1}^{k_1}P_jC(X,M_{s(j)})P_j,
$$
where each $P_j$ has a constant rank $r(j)>0.$
%$P(x)\not=0$ for all $x\in X,$

%{ (For example one can choose  $X$ to be any connected simplicial complex of dimension at least one and $Z$ a nonempty proper  sub complex of $X$ with $X\setminus Z$ connected.)}

%{ {Huaxin: I think in the above definition, for ${\cal D}_0$, we assume it is a general finite dimensional $C^*$ algebras, but when we define ${\cal D}_k$ from ${\cal D}_{k-1}$, we can assume the above $F$ is just one matrix algebra--not necessary direct sum of matrix algebras. This definition has same generality and will make induction argument simple in notation because we don't need to deal with several blocks of $F$. I will let you to decide ---Guihua}}

Let $A_m$ be a unital \CA\, in ${\cal {\bar D}}_m'.$
For $0\le k<m,$ let  $A_k\in {\cal {\bar D}}_k'$ such that
$A_{k+1}=\{(f, a)\in P_{k+1}C(X_{k+1}, F_{k+1})P_{k+1}\oplus A_k: f|_{ Z_{k+1}}=\Gamma_{k+1}(a)\},$
where $F_{k+1}$ is a finite dimensional \CA, $P_{k+1}\in C(X_{k+1}, F_{k+1})$ is a projection,
$\Gamma_{k+1}: A_k\to P_{k+1}C({ Z_{k+1}}, F_{k+1})P_{k+1}$ is a unital \hm,
$k=0,1,2,...,m-1.$ Denote by $\partial_{k+1}: P_{k+1}C(X_{k+1}, F_{k+1})P_{k+1}\to
Q_{k+1} C(Z_{k+1}, F_{k+1})Q_{k+1}$  the map defined by
$f\mapsto f|_{Z_{k+1}}$ ($Q_{k+1}=P_{k+1}|_{Z_{k+1}}$).
We use $\pi_e^{(k+1)}: A_{k+1}\to A_k$ for the quotient map and
$\lambda_{k+1}: A_{k+1}\to P_{k+1}C(X_{k+1}, F_{k+1})P_{k+1}$ for the map  defined by
$(f,a)\mapsto f.$

For each $k,$ one has the following commutative diagram:
\begin{equation}\label{pull-back-k}
\xymatrix{
A_{k} \ar@{-->}[rr]^{\lambda_k} \ar@{-->}[d]^-{\pi_e^{(k+1)}}  && P_kC(X_k, F_k)P_k \ar[d]^-{\partial_k} \\
A_{k-1}\ar[rr]^-{\Gamma_k} & & Q_kC({ Z_k}, F_k)Q_k,
}
\end{equation}

In general, suppose that $A=A_m\in{\cal {\bar D}}_m'$ is constructed in the following sequence
$$A_0=F_0,
%P_0M_r(X_0)P_0,
~~ A_1=A_0\oplus_{Q_1C(Z_1, F_1)Q_1}P_1C(X_1, F_1)P_1,~~A_2=A_1\oplus_{Q_2C(Z_2, F_2)Q_2}P_2C(X_2, F_2)P_2, \cdots,~~~~~~~~~~~~$$
$$ ~~~~~~~~~~~~~~~~~~~~~~ A_m=A_{m-1}\oplus_{Q_mC(Z_m, F_m)Q_m}P_mC(X_m, F_m)P_m,$$
where $Q_i=P_i|_{Z_i},$ $i=1,2,...,m.$
%Set $F_0=M_r.$ Then $A_0=P_0C(X_0, F_0)P_0.$
%Note that $X_0$ may {\it not} be connected.
With $\pi_e^{(k+1)}$ and $\lambda_{k}$ above we can define quotient map $\Pi_k: A=A_m \to A_k$ and homomorphism $\LD_k: A =A_m \to P_kC(X_k, F_k)P_k$ as follows:
$$\Pi_k=\pi_e^{(k+1)}\circ \pi_e^{(k+2)}\circ \cdots \circ \pi_e^{(m-1)}\circ \pi_e^{(m)}~~~~~~~~~~{\mbox {and}}~~~~~~~~~ \LD_k=\lambda_k\circ \Pi_k. $$
Combining all $\LD_k$ we get the inclusion homomorphism
$$\LD: A \to \bigoplus_{k=0}^mP_kC(X_k, F_k)P_k$$
with $X_0$ being single point set.
For each $k\ge 1,$ we may write
$$
P_k(C(X_k, F_k))P_k=P_{k,1}C(X_k, M_{s(k,1)})P_{k,1}\oplus P_{k,2}C(X_k, M_{s(k,2)})P_{k,2}\oplus\cdots P_{k,t_k}C(X_k, M_{s(k,t_k)})P_{t_k},
$$
where $P_{k,j}\in C(X_k, M_{s(k,j)})$ is a projection of rank $r(k,j)$ at each $x\in X_k.$
%For $k=0,$ we may write
%$$
%A_0=\bigoplus_{j=1}^{l(0)}P_{0,j}C(X_{0,j}, M_r)P_{0,j},
%$$
%where $X_{0,j}$ is connected, $X_0=\sqcup_{j=1}^{l(0)}X_{0,j}$ and
%$P_{0,j}\in C(X_{0,j}, M_r)$ is a projection of rank $r(0,j).$
%Let $F_k=\bigoplus_{j=1}^{t_k}M_{s(j)}$.
For each $x\in X_k$ and $j\leq t_k$ and $f\in A$, denoted by $\pi_{(x,j)}(f)\in M_{r(k,j)}$, the evaluation of
$j^{th}$ component of $\LD_k(f)$ at point $x.$
%$$\LD_k(f)\in C(X_k, F_k)=C(X_k, M_{s(1)})\oplus C(X_k, M_{s(2)})\oplus\cdots C(X_k, M_{s(t_k)})$$ at point $x\in X_k$.}~~
 Then for each pair $(x,j)$, $\pi_{(x,j)}$ is a finite dimensional representation of $A$, and furthermore if we assume $x\in X_k\setminus Z_k$, then $\pi_{(x,j)}$ is an irreducible representation.

 %One notices that ${\cal C}\subset {\cal D}_1'.$

%Put $F_A=A_0.$ Let $\pi_F=\pi_e^{(1)}: A_1\to A_0$ be the surjective \hm\, in the pull-back:
%\begin{equation}\label{pull-back-k}
%\xymatrix{
%A_1 \ar@{-->}[rr]^{\lambda_1} \ar@{-->}[d]^-{\pi_e^{(1)}}  && C(X_1, F_1) \ar[d]^-{\partial_1} \\
%A_0\ar[rr]^-{\Gamma_1} & & C({ Z_1}, F_1),
%}
%\end{equation}
%Note that if $A_1, A_2\in {\bar D}_k,$ then $A_1\oplus A_2\in {\bar D}_k.$
For each integer $m\ge 0,$ denote by ${\cal {\bar D}}_m$ the class of \CA s which are finite direct sums
of those \CA s in ${\cal {\bar D}}_m'.$
%We use ${\cal D}_0'={\cal {\bar D}}_0.$
 If in the above definition
of ${\cal {\bar D}}_k',$ in addition, we assume that $X$ is path connected,
%=\sqcup_{i=1}^m X_i,$ where each $X_i$ is path connected, $Z=\sqcup_{i=1}^m Z_i,$ where $Z_i\subset X_i$
$Z$  has finitely
many path connected components and $X_i\setminus Z_i$ is path connected,
then we will use ${\cal D}_k'$ for the resulting class of \CA s
and we use
${\cal D}_k$ for  the class of \CA s which are finite direct sums
of those \CA s in ${\cal D}_m'.$  Note that ${\cal D}_k\subset {\cal {\bar D}}_k$ and ${\cal C}\subset {\cal D}_1.$

\end{df}

%Let $A$ be a $m$-dimensional NCCW complex as above.  In what follows we denote $X=\I^m,$
%$Y=X\setminus \partial(X)\,(\cong S^{m})$ and $I=C_0(Y, F)\subset A.$
%For $d>0,$ define
%$X^d=\{x\in X: {\rm dist}(x, \partial(X))< d\}$ and $Y^d=X\setminus X^d.$
%Let
%\beq\label{dAd}
%A^d=\{(f,b)\in C(X^d, F)\oplus B: \pi_I(f)=\Gamma(b)\},
%\eneq
%where $\pi_I=\partial_m: C(\overline{X^d}, F)\to C(\partial(X),F)$ is  the quotient map induced by the same restriction on $\partial(X).$ There is a quotient map $\pi^d: A\to A^d$ by $\pi^d((f, b))=(f_{\overline{X^d}},b)$ for all $(f, b)\in A.$
%Note that $A/I\cong B.$
%We will also use $\pi_I$ for the quotient map.

%There is $d^Y>0$ such that, for any $0<d\le d^Y,$ there is a continuous map  $s^d_*: X^d\to \partial(X)$ such
%that $s^d_*(x)=x$ for all $x\in \partial(X)$ and ${\rm dist}(s^d_*(x),x)\le d$ for all $x\in X^d.$
%Define $s^d: A/I\to A^d$ by $s^d(b)=(\Gamma(f)\circ s^d_*, b)$ for all $b\in B.$
%For any $g\in A,$\beq\label{d-nccw-c}
%\lim_{d\to 0} \|s^d\circ \pi_I(g)-\pi^d(g)\|=0\rforal g\in A.
%\eneq
%
%Write $F=M_{r(1)}\oplus M_{r(2)}\oplus \cdots \oplus M_{r(k_2)}.$ Denote by $\pi_{x, j}: A\to M_{r_j}$ the point-evaluation
%defined by $p_j(\lambda(a)(x)),$ where $\lambda: A\to C(X, F)$ be give by the pull-bak map (as $\lambda_k$ in (\ref{pull-back-k})), $x\in X$ and $p_j: F\to M_{r(j)}$ is the projection map.

\begin{rem}\label{8-n-rn}
Let $A\in {\cal {\bar D}}_k$ (or $A\in {\cal D}_k$). It is easy to check that $C(\T)\otimes A\in {\cal {\bar D}}_{ k+1}$
(or $A\in {\cal D}_{k+1}$). First of all it is easy to check that if $F_0$ is a finite dimensional $C^*$ algebra, then $C(\T)\otimes F_0\in {\cal {\bar D}}_1$ by putting $F_1=F_0$ and $X_1=\T$ with $Z_1=\{1\}\subset \T$ and $\GM_1: F_0\to C(Z_1, F_1)\cong F_0$ to be the identity map. And if a  pair of spaces $(X_k, Z_k)$ satisfies   the condition described as pair $(X, Z)$ in the above definition, in particular, the existence of the retraction $s_*^d$ and homeomorphism $r$ as in \ref{Amnotation-1} and \ref{Amnotation-2+}, then the  pair  $(X_k \times \T , Z_k\times \T)$ also satisfies the same condition.

\end{rem}

\begin{rem}\label{8-n-K0nccw}
Let $A\in {\cal D}_k$ be a unital \CA.
There is an integer $l\ge 1$ such that $(K_0(A))/\ker(\rho_A)\subset \Z^l$ as an order group, where
$\Z^l$ is equipped with the usual order.
This has been proved for $k=0,1.$ Now assume that $k>1$ and $A\in {\cal  D}_k.$
Then the  above statement holds for any unital \CA\, in ${\cal D}_d$ with $0\le d<k.$ Consider the short exact sequence
\beq\label{K0nccw-1}
0\to  I
%C_0(\I^{k}\setminus \partial(I^k), F))
\to_j A\to C\to_{\pi_e} 0,
\eneq
where $C\in {\cal  D}_d$ for some  $d<k,$  $I=PC_0(X\setminus { Z}, F))P$
and $j: I\to A$ is the embedding.
Note $\rho_A\circ j_{*0}=0,$ { since $X$ is path connected and $Z\subset X$ is non empty}.  It follows that $\rho_A$ factors through $(\pi_e)_{*0}.$  Hence the map $K_0(A)/\ker (\rho_A) \to K_0(C)/\ker (\rho_C)$ is injective and consequently the statement holds for $A$.
%By the induction, it is easy
%to see that
%$\rho_A(K_0(A))\subset \rho_C(K_0(C))\subset \Z^l$ for some $l\ge 1.$
In particular for $A\in {\cal D}_k$ as notation in \ref{8-N-3}, we have
\beq\label{K0nccw-2}
\hspace{-0.2in}K_0(A) /\ker(\rho_A)\subset  (K_0(A_0))/\ker(\rho_{A_0}),\,K_0(A) /\ker(\rho_A)\cong (\Pi_0)_{*0}(K_0(A)) \subset K_0(A_0)\cong \Z^l.
\eneq
\end{rem}

\begin{lem}\label{8-N-4}
%Let $X$ be a compact metric space, let $Y\subset X$ be an open subset, $F_1$ and $F_2$ be two finite dimensional
%\CA s with a unital embedding $\imath: F_1\to F_2.$
Let $A\in {\cal {\bar D}}_k$ be a unital \CA.
%finite dimensional NCCW complex.
%sub-homogeneous \CA\, with the
 %property ($S_m$).
 Let $\Delta: A_+^{q, {\bf 1}}\setminus \{0\}\to (0,1)$ be an order preserving map.  Let $1>\af>1/2.$

Let $\ep>0,$ ${\cal F}\subset A$ be a finite subset, ${\cal H}_0\subset A_+^{\bf 1}\setminus \{0\}$ be a finite subset
and $K\ge 1$ be an integer. There exist  an integer $N\ge 1,$ $\dt>0,$  a finite subset ${\cal H}_1\subset A_+^{\bf 1}\setminus \{0\},$ a finite
subset ${\cal H}_2\subset A_{s.a.}$  satisfying the following:

If $\phi_1, \phi_2: A\to M_n$(for some integer $n\ge N$) are two unital \hm s such that
\beq\label{8-N-4-1}
\tau\circ \phi_1(g)\ge \Delta(\hat{g})\tforal g\in {\cal H}_1\tand\\
|\tau\circ \phi_1(g)-\tau\circ \phi_2(g)|<\dt\tforal g\in {\cal H}_2,
\eneq
where $\tau$ is the tracial state on $M_n,$ then
there exist mutually orthogonal projections $e_0, e_1, e_2,...,e_K\in M_n$ such that
 $e_1,e_2,...,e_K$ are mutually equivalent, $e_0\lesssim e_1,$ and $e_0+\sum_{i=1}^K e_i=1_{M_n},$
unital \hm s $h_1, h_2: A\to e_0M_ne_0,$ a unital \hm\, $\psi: A\to e_1M_ne_1$  and a unitary
$u\in M_n$ such that
\beq\label{8-N-4-3}
&&\|{\rm Ad}\, u \circ \phi_1(f)-(h_1(f)\oplus {\rm diag}(\overbrace{\psi(f), \psi(f),...,\psi(f)}^K))\|<\ep,\\
&&\|\phi_2(f)-(h_2(f)\oplus {\rm diag}(\overbrace{\psi(f), \psi(f),...,\psi(f)}^K))\|<\ep  \tforal f\in {\cal F}\\
&&\tand \tau\circ \psi(g)\ge \af{\Delta(\hat{g})\over{K}}\tforal g\in {\cal H}_0,
\eneq
where $\tau$ is the tracial state on $M_n.$
\end{lem}

\begin{proof}
%Let $A$ be an $k$-dimensional NCCW complex.
We will use induction on integer $k\ge 0.$ The case $k=0$ follows from \ref{Aug-N-2}. Now assume
that the lemma holds for integers $0\le k\le  m.$

  We assume that $A\in {\cal {\bar D}}_{m+1}'.$
  %is an $m+1$-dimensional NCCW complex as in \ref{8-N-3}
  %and we
  We will retain the notation for $A$ as an algebra in ${\cal {\bar D}}_{m+1}'$  in the  later part of  \ref{8-N-3}.
  % that
%$A\subset PC(X,F_2)P$ as a unital \SCA, where $X$ is a compact metric space and
In particular, $X_{m+1}=X,$ ${ Z_{m+1}= Z}$, $Y=X\setminus { {Z}},$
$X^0={ Z}=X\setminus Y,$ $I= PC_0(Y,F)P\subset A.$
We will write
\beq\label{adf}
A=\{(f,b)\in PC(X, F)P\oplus B: f|_{X^0}=\Gamma(b)\},
\eneq
where $B\in {\cal {\bar D}}_m'$  is a unital \CA\,  and will be identified with
$A/I.$
We also keep the notation $\lambda: A\to PC(X, F)P$ in the pull-back in the last part of \ref{8-N-3}.
We will use $f|_S$ for $\lambda(f)|_S$ for $f\in A$  and $S\subset X$ in the proof when there is no confusion.
%We assume that $A\cap I=I=\{f\in A: f|_{X^0}=0\}.$
%There exists  $d_0>0$ such that there exists a unital injective \hm\, $s: A/I\to A^{d_0}$ such that
%$\pi_I'\circ s={\rm id}_{A/I}$ (see \ref{8-N-3}).
Let ${ d_X}> 0$ be given in \ref{8-N-3}.
Denote by $\pi_I: A\to A/I$ the quotient map.
%We write $F=M_{r_1}\oplus M_{l_2}\oplus \cdots \oplus M_{l_{k_1}}.$

%By considering each summand of $A,$ to simplify notation, \wilog, we may assume
%that  $F=M_{r(k)'}$ and
We may write that
\beq\label{8-N-40-150825-1}
PC(X,F)P=\bigoplus_{j=1}^{k_2}P_{j}C(X, M_{s(k,j)})P_{j},
\eneq
where $P_{j}\in C(X, M_{s(k,j)})$ is a projection of rank $r(j)$ at each $x\in X.$
We may assume that $A/I$ has irreducible representations with rank $l_1, l_2,...,l_{k_1}.$
%and $P$ has rank $r(k)\le r(k)'.$
%$I\cong PC_0(Y, F_2)P$ such that
%$P_xF_2P_x$  is a
%\SCA\,  as well as a
%quotient of
%$F=M_{r(1)}\oplus M_{r(2)}\oplus \cdots \oplus M_{r(k_2)}$ (where
%$r(1),r(2),...,r(k_2)$ may not be distinct.)
Choose
%$T=k_1\cdot \max_{i,j}\{z_iz_j: z_i, z_j\in \{l_1,l_2,...,l_{k_1}, r(k)\}\}.$
$T=(k_1k_2)\cdot \max_{i,j}\{z_iz_j: z_i, z_j\in \{l_1,l_2,...,l_{k_1}, r(1),r(2),...,r(k_2)\}\}.$

 Choose $\bt=\sqrt{1-(1-\af)/8}=\sqrt{(7+\af)/8}.$  Note  that $1>\bt^2>\af.$
Fix $N_{00}\ge 4$ such that $1/N_{00}<{(1-\bt)\over{64}}.$

Fix $\ep>0,$ a finite subset ${\cal F}\subset A$ and a finite subset ${\cal H}_0\subset A_+^{\bf 1}\setminus \{0\}.$
We may assume that $1_A\in {\cal H}_0\subset {\cal F}.$ Without loss of generality, we may also assume that
${\cal F}\subset A_{s.a.}$ and $\|f\|\le 1$ for all $f\in {\cal F}.$
Write $I=\{f\in PC(X, F)P: f|_{X^0}=0\}.$
Put
\beq\label{8-N-4-3+}
\dt_{00}=\min\{\Delta(\hat{g})/2: g\in {\cal H}_0\}>0.
\eneq

There is $d>0$ such that
\beq\label{8-N-4-4}
%\|\pi_{x}(f)-\pi_{x'}(f)\|<\min\{\ep, \dt_{00}\}/256KN_{00}\rforal f\in {\cal F},
\|\pi_{x,j}(f)-\pi_{x',j}(f)\|<\min\{\ep, \dt_{00}\}/256KN_{00}\rforal f\in {\cal F},
\eneq
provided that ${\rm dist}(x, x')<d$ for any pair $x, x'\in X$
%$j=1,2,...,k_2$
(here we  identify $\pi_{x,j}(f)$ with $\pi_{x,j}({ \LD}(f))$
---see the { \ref{8-N-3}}).
Put $\ep_0=\min\{\ep, \dt_{00}\}/16KN_{00}.$

We also assume that, for any $x\in X^d=\{x\in X: {\rm dist}(x, X^0)<d\}$ {{choosing smaller $d$ if necessary,}}
\beq\label{8-N-4-4+}
%\|\pi_{x}\circ s^d\circ { (\LD(f)|_Z)}-\pi_{x}(f)\|<\ep_0/16,
\|\pi_{x,j}\circ s^d\circ { (\LD(f)|_Z)}-\pi_{x,j}(f)\|<\ep_0/16,
\eneq
{ where $s^d: QC(Z, F)Q\to P|_{X^d}C(\overline{X^d}, F)P|_{X^d}$ is induced by $s_*^d: \overline{X^d} \to Z$ (see \ref{8-N-3}). Note that $s^d$ also induces a map (still denoted by $s^d$)
$$s^d: B \to B\oplus_{QC(Z,F)Q}P|_{\overline{X^d}}C(\overline{X^d}, F)P|_{\overline{X^d}},$$}
where $Q=P|_Z.$
To simplify the notation, we may assume that $d<d_{ X}/2.$
%{ {what is $d_0$ here. I suggest to change to $d_X/2$.}}

%If $X^0$ is a clopen set of $X,$ Then $A=A/I\oplus I,$ where
%both of $A/I$ and $I$ have the property ($S_{m}$).  No induction is needed.
%Therefore we may assume that $X^{d/256}\setminus X^0$ has infinitely many points.

For any $b>0,$ as in \ref{8-N-3}, we will continue to use $X^b$ for
$\{x\in X: {\rm dist}(x, X^0)<b\}.$
%Choose $d/4<d_0 <d/2$ such that

Let  $Y_0=X\setminus X^{d/2}.$
%$ and $Y_0=\overline{X^{d/4}}.$
%be the closure of $\{x\in Y: {\rm dist}(x, X^0)<d/4\}.$
%Put $Y_1=\{x\in X: {\rm dist}(x, Y_0)\le d/16\}.$ Then $Y_1$  is compact.
%By the assumption, there is $\eta_0>0$ such that $\eta_0<d/32$ such
%that there is a homeomorphism $\xi_0: Y_{00}\to Y_0,$ where $Y_{00}=X\setminus X^{d/8+\eta_0/2}.$

%In what follows we assume that $Y_{00}\not=\emptyset.$  In the end of proof, we will take care of the case
%that $Y_{00}=\emptyset.$
%Note, by \ref{8-N-3},
%%We may assume that $Y_{00}$ has infinitely many points.  The case that $Y_{00}$ has finitely many points
%%will be dealt with at the end of this proof.

Put $B_{0,0}=PC(Y_0, { F})P.$
% where $P_0=P|_{Y_{00}}.$
Let ${\cal F}_{I,0}=\{f|_{Y_0}: f\in {\cal F}\}$ and let ${\cal H}_{0,I,0}=\{h|_{Y_0}: h\in {\cal H}_0\}.$

%Let $f_{0,0}\in C_0(Y)_+$ be such that $0\le  f_{0,0}\le 1,$ $f_{0,0}(x)=1$ if $x\in Y_0,$
%$f_{0,0}(x)=0$ if $x\in Y_0$ and $f_{0,0}(x)>0$ if ${\rm dist}(x, X^0)>d/4.$
%
%{ {I suggested the above paragraph to be changed to the following, otherwise we will have problem for the definition of $g'$ in next paragraph:}}

 {Let $f_{0,0}\in C_0(Y)_+$ be such that $0\le  f_{0,0}\le 1,$ $f_{0,0}(x)=1$ if $x\in X\setminus X^d,$
$f_{0,0}(x)=0$ if $x\notin Y_0$ and $f_{0,0}(x)>0$ if ${\rm dist}(x, X^0)>d/2.$}

Let $\Delta_{I,0}: (B_{0,0})_+^{q, {\bf 1}}\setminus \{0\}\to (0,1)$ be defined by
\beq\label{8-N-4-5}
\Delta_{I,0}(\hat{g})=\bt \Delta(\widehat{g'})\tforal g\in (B_{0,0})_+^{\bf 1}\setminus\{0\},
\eneq
where $g'=f_{0,0}\cdot P\cdot g$
%{ {this is not well define if not make the above change since for $x\in X^{d/2}\setminus \overline{X^{d/4}}$, $g(x)$ has %no defintion and $f_{00}(x)\not= 0$}}
which  is viewed as an element in $I_+^{\bf 1}.$
%Also  \andeqn
%D_{0,I}=\Delta_I(\hat{1_{B_0}}).
%\eneq
Note that if $g\not=0,$ then $f_{0,0}\cdot P\cdot g\not=0.$ So $\Delta_{I,0}: (B_{0,0})_+^{q, {\bf 1}}\to (0,1)$ is an order preserving
map.  Let $N^I\ge 1$ be an integer (in place of $N$) as required  by \ref{Aug-N-2} for $B_{0,0}$ (in place of $A$),
$\Delta_{I,0}$ (in place of $\Delta$), $\ep_0/16$ (in place of $\ep$), ${\cal F}_{I,0}$ (in place of ${\cal F}$) and
${\cal H}_{0,I,0}$ (in place of ${\cal H}_0$) { and $2K$ (in place of $K$)}.

Let ${\cal H}_{1, I,0}\subset (B_{0,0})_+^{\bf 1}\setminus \{0\}$ be a finite subset (in place of ${\cal H}_1$),
${\cal H}_{2,I,0}\subset (B_{0,0})_{s.a.}$ (in place of ${\cal H}_2$) and $\dt_{1}>0$ (in place of $\dt$)
{{be}} as required
by \ref{Aug-N-2} for $\ep_0/16$ (in place of $\ep$), ${\cal F}_{I,0}$ (in place of ${\cal F}$), $2K$ { (in place of $K$)}  and ${\cal H}_{0,I,0}$ associated
with $B_{0,0}$ (in place of $A$) and $\Delta_{I,0}$ (in place of $\Delta$) { and $\beta$ (in place of $\alpha$)}.
Without loss of generality, we may assume that
$\|g\|\le 1$ for all $g\in {\cal H}_{2,I,0}.$
%Since $Y^{d^Y}$ has infinitely many points, we may assume that

Let ${\cal F}_{\pi}=\pi_I({\cal F}).$
%\{f|_{X^0}: f\in {\cal F}\}.$
Let $g_0'\in C(X)_+$ with $0\le g_0'\le 1$ such that
$g_0'(x)=0$ if ${\rm dist}(x, X^0)<d/256$ and $g_0'(x)=1$ if ${{{\rm dist}(x,Y_0)\le d/16}}.$
Define $g_0=1_A-g_0'\cdot P.$ Since $g_0'\in I,$ we view $g_0$ as an element in $A.$
%Let $g_0\in C(X)_+$ with $0\le g_0\le 1$ such that
%$g_0(x)=1$ if ${\rm dist}(x, X^0)<d/256$ and $g_0(x)=0$ if $x\in Y_{00}.$
%Since $1-g_0\in I,$ we may view $g_0$ as an element of $A.$
Note that  $(1-g_0')P\cdot s^d(g)$ may be viewed as an element in
$P(C(X,F)P$ as $(1-g_0')(x)=0$ for $x\notin X^d.$ (Recall that $s^d:  QC(Z,F)Q \to PC(\overline{X^d}, F)P$ is also regarded as map  $$s^d: B=A/I \to B\oplus_{QC(Z,F)Q}P_{\overline{X^d}}C(\overline{X^d}, F)P|_{\overline{X^d}}.)$$
Hence
$g_0\cdot s^d(g)$ may be identified with
$((1-g_0')P\cdot s^d(g), g)\in A.$
Define
\beq\label{8-N-4-9}
\Delta_\pi(\hat{g})=\bt\Delta(\widehat{g_0\cdot s^d(g)})\rforal g\in (A/I)_+^{\bf 1}.
\eneq
{{We will late used the fact that $g_0(x)=0$ if ${\rm dist}(x, Y_0)\le d/16.$}}
%Furthermore for any $g\in B$, $g_0\cdot s^d(g) \in B\oplus_{C(Z,F)}C(\overline{X^d}, F)$ could be regarded as element in $A=B\oplus_{C(Z,F)}C(X, F)$, as $g_0(x)=0$ for $x\notin X^d$.)

 Note if $g$ is nonzero, so is $s^d(g).$  {{Since $g_0|_{X^0}=1,$}} we have that  $g_0 \cdot s^d(g)\not=0.$
It follows that $\Delta_\pi: (A/I)_+^{q, {\bf 1}}\setminus \{0\}\to (0,1)$ is an order preserving map.

Put ${\cal H}_{0,\pi}=\pi_I({\cal H}_0).$
%\{f|_{X^0}: f\in {\cal H}_0\}.$
Let $N^{{{\pi}}}\ge 1$ be the integer associated  with $A/I{ (=B)},$ $\Delta_\pi,$ $\ep_0/16,$ ${\cal F}_{\pi}$ and
${\cal H}_{0,\pi}$ (as required by the inductive assumption that the lemma holds for integer $m$).

 Let ${\cal H}_{1,\pi}\subset (A/I)_+^{q, {\bf 1}}$ be a finite subset
(in place of ${\cal H}_1$),
${\cal H}_{2,\pi}\subset A/I_{s.a.}$ be a finite subset (in place of ${\cal H}_2$) and let $\dt_2>0$ (in place of $\dt$)
be as required
{{by the inductive assumption that this  lemma for the case  that $k=m$}}  for $\ep_0/16$ (in place of $\ep$), ${\cal F}_{\pi}$ (in place of ${\cal F}$), ${\cal H}_{0, \pi}$ (in place of ${\cal H}_0$)  and $2K$ associated with $A/I$ (in place of $A$), $\Delta_{\pi}$ (in place of $\Delta$) and $\bt$ (in place of $\af$).
{{\Wlog,  we may assume that $\|h\|\le 1$ for all $h\in {\cal H}_{2,\pi}.$}}

Set ${{{\dt_{000}}}}=\min\{{{ \dt_1,\dt_2, \ep_0}}\}.$ There is an integer $N_0\ge 256$ such that
\beq\label{8-N-4-5+}
1/N_0<\Delta(\widehat{f_{0,0}\cdot P})\cdot{{\dt_{000}}}^2\cdot \min\{\Delta_{I,0}(\hat{g}): g\in {\cal H}_{1,I,0}\}\cdot
\min\{\Delta_{ \pi}(\hat{g}): g\in {\cal H}_{1,\pi}\}/64KN_{00}.
\eneq
Define $Y_k$ to be the closure of $\{y\in Y: {\rm dist}(y, Y_0)<{kd\over{64N_{ 0}^2}}\},$
%$Y_k=\{y\in Y: {\rm dist}(y, Y_0)\le {kd\over{64N_0^2}}\},$
$k=1,2,..., 4N_0^2.$

Let ${\cal F}_{I,k}=\{f|_{Y_k}: f\in {\cal F}\}$ and let ${\cal H}_{0,I,k}=\{h|_{Y_k}: h\in {\cal H}_0\}.$
Put $B_{0,k}=P|_{Y_k}C(Y_k, F)P|_{Y_k}.$
%where $P_k=P|_{Y_k}.$
Let $f_{0,k}\in C_0(Y)_+$ be such that $0\le f_{0,{ k}}\le 1,$ $f_{0,k}(x)=1$ if $x\in Y_{k-1},$
$f_{0,k}(x)=0$ if $x\not\in Y_k$ and $f_{0,k}(x)>0$ if ${\rm dist}(x, Y_0)<{kd\over{64N_0^2}},$ $k=1,2,...,4N_0^2.$

%We may assume that $X^{d/256}$ is an infinite set. Otherwise, $A\cong F_3\oplus C(Y', F_2)$
%for some finite dimensional \CA\, $F_3$ and for some compact subset $Y'\subset Y.$
%Therefore, the lemma follows from \ref{Aug-2}.

Let $r_k: Y_k\to Y_0$ be a homeomorphism such that
%$r_k(x)=x$ for $x\in Y_{00}$ and
\beq\label{8-N-4-6}
{\rm dist}(r_k(x),x)<d/4\tforal x\in Y_k,\,\,\,k=1,2,..., 4N_0^2
\eneq
(see \ref{8-N-3}).

Let $\Delta_{I,k}: (B_{0,k})_+^{q, {\bf 1}}\setminus\{0\} \to (0,1)$ be defined by
\beq\label{8-N-4-7}
\Delta_{I,k}(\hat{g})=\bt\Delta(\widehat{g'})\tforal g\in (B_{0,k})_+^{\bf 1},
\eneq
where $g'=(f_{0,k}\cdot P \cdot g)\circ r_k^{-1}$ which  is viewed as an element in $I_+^{\bf 1}.$
Note that if $g\not=0,$ then $f_{0,k}\cdot 1_A\cdot g\not=0.$ So $\Delta_{I,k}: (B_{0,k})_+^{q, {\bf 1}}\setminus \{0\}\to (0,1)$ is an order preserving
map.

Let ${\cal F}_{I,k}'=\{f\circ r_k: f\in {\cal F}_{I,0}\}$ and ${\cal H}_{0,I,k}'=\{g\circ r_k: g\in {\cal H}_{0,I,0}\},$
$k=1,2,...,4N_0^2.$

Any unital \hm\, $\Phi: B_{0,k}\to C$ (for any unital \CA\, $C$)
 induces  a unital \hm\, $\Psi: B_{0,0}\to C$ by $\Psi(f)=\Phi(f\circ r_k)$ for all $f\in B_{0,0}.$ Note also
that $f\mapsto f\circ r_k$ is an isomorphism from $B_{0,0}$ onto $B_{0,k}.$

Then, for $\ep_0/16$ (in place of $\ep$), ${\cal F}_{I,k}'$ (in place of ${\cal F}$), $K$  and ${\cal H}_{0,I,k}'$ associated
with $B_{0,k}$ (in place of $A$), $\Delta_{I,k}$ (in place of $\Delta$) and $\bt$(in place of $\af$),  to apply \ref{Aug-N-2},
one can choose ${\cal H}_{1,l,k}$ (in place of ${\cal H}_1$) to be $H_{1,{ I},0}\circ r_k,$
${\cal H}_{2,I,k}$ (in place of ${\cal H}_2$) to be ${\cal H}_{2,{ I},0}\circ r_k$ and $\dt_1$ (in place of $\dt$).

We also note that
\beq\label{8-N-4-8}
\|f-f|_{{{Y_0}}}\circ r_k\|<\min\{\ep, \dt_{00}\}/64KN_{00}\tforal f\in {\cal F}_{I,k}.
\eneq

{{Define
\begin{equation*}
\sigma_0=\min\{\min\{\Delta_{I,0}(\hat{g}): g\in {\cal H}_{1, I,0}\}, \min\{\Delta_{\pi}(\hat{g}): g\in {\cal H}_{1,\pi}\}\}.
\end{equation*}
}}
{{Choose an integer $N\ge N^\pi+N^I$ such that
\beq\label{8-N-4-13-}
{T\over{N}}< \sigma_0\cdot \min\{\dt_1/64, \dt_2/64, \ep_0/64K\}/N_{00}.
\eneq
}}

Put
\beq\label{8-N-4-10}
{\cal H}_1=\cup_{k=0}^{4N_0^2}\{f_{0,k}\cdot 1_A\cdot g\circ { r_k}: g\in {\cal H}_{1,I ,0}\}\cup
\{(g_0\cdot 1_A\cdot s^{ d}(g),g): g\in {\cal H}_{1,\pi}\}\cup {\cal G}_I\cup {\cal G}_\pi,
\eneq
where ${\cal G}_I$ is a finite subset in $\{g\cdot P: g\in C_0(Y): g|_{Y\setminus Y_0}=0\}$
{{which contains at least $N$ mutually orthogonal non-zero positive  scalar elements}} and  ${\cal G}_\pi$ is  a finite subset of $J$ {{which contains at least $N$ mutually orthogonal non-zero scalar elements,}}
where $J=\{(f\cdot P, b)\in A: f\in C(X): f|_{X\setminus X^{d/256}}=0\}$ such that $0\le f\le 1$ for all $f\in {\cal G}_I\cup {\cal G}_\pi.$
%to be determined later.
%Note ${\cal G}_I\cup {\cal G}_\pi\subset A.$

%For each $g\in {\cal H}_{2,I,k},$ there is $g^I\in I_{s.a.}$ such that $g^I|_{Y_k}=g$ and $\|g\|\le 1.$
With  the convention that $r_0: Y_0\to Y_0$ is the identity map, put
\beq\label{8-N-11}
{\cal H}_2'=\cup_{k=0}^{4N_0^2}\{f_{0,k}\cdot g\circ r_k: g\in {\cal H}_{2, I,0}\}\cup\{f_{0,k}, 0\le k\le 4N_0^2\}.
\eneq
Define $g_{0,k}'\in C(X)_+$ so that $0\le g_{0,k}'\le 1,$
$g_{0,k}'(x)=0$ if $x\not\in Y_k,$  $g_{0,k}'(x)=1$ if $x\in Y_{k-1}$
and $g_{0,k}'(x)>0$ if $x\in Y_k\setminus Y_{k-1},$  $k=1,2,....$ \\

%{
%{ It seems the $g_{0,k}'$ is the $f_{0,k}$ before. Why define it again}
%} \\

 Put $g_{0.k}=1_A-g_{0,k}'\cdot P.$ Note, since $g_{0,k}' \cdot P\in I,$
$g_{0,k}\in A.$
Define
\beq\label{8-N-12}
{\cal H}_2''=\cup_{k=1}^{4N_0^2}\{(g_{0,k}\cdot  s^{ d}(g),g)\in A: g\in {\cal H}_{2, \pi}\}\cup {\cal F}.
\eneq
%where ${\cal G}$ be a finite subset of $J,$
%where $J=\{f\cdot 1_A: f\in C(X): f|_{X\setminus X^{d/256}}=0\}$ such that $0\le f\le 1$ for all $f\in {\cal G}.$
Put
\beq\label{8-N-13}
{\cal H}_2={\cal H}_2'\cup{\cal H}_2''.
\eneq
%Define
%\begin{equation*}
%\sigma_0=\min\{\min\{\Delta_{I,0}(\hat{g}): g\in {\cal H}_{1, I,0}\}, \min\{\Delta_{\pi}(\hat{g}): g\in {\cal H}_{1,\pi}\}\}.
%\end{equation*}
%Choose an integer $N_1\ge N^\pi+N^I$ such that
%\beq\label{8-N-4-13-}
%{T\over{N_1}}< \sigma_0\cdot \min\{\dt_1/64, \dt_2/64, \ep_0/64K\}/N_{00}.
%\eneq
%Since $X^{d/256}$ has infinitely  many points, we may assume that
%Note that  we may assume that ${\cal G}_\pi$ contains  at least
%$N_1$ many positive and mutually orthogonal scalar functions in $J;$ and
%since we assume that $Y_{00}$ contains infinitely many points,
%we may also assume that ${\cal G}_I$ contains at least $N_1$ many positive and mutually
%orthogonal scalar functions.

Let
\begin{equation*}
\dt={\sigma_0\cdot \min\{\dt_1/64, \dt_2/64, \ep_0/64K\}\over{4KN_1N_{00}}}.
\end{equation*}
%Since {\color{Green} (Anything missing? Or just delete ``Since"?)}\\

%{ {What is $\sigma$, is $sigma=\epsilon_0$ ?}}\\

Now let $\phi_1,\phi_2: A\to M_n$ (for some integer $n\ge N$) be two unital \hm s such that
they satisfy the assumption for the above ${\cal H}_1,$ ${\cal H}_2$ and $\dt.$

Consider two finite  Borel  measures on $Y$ defined by
\beq\label{8-N-4-15}
\int_Y f\mu_i=\tau\circ \phi_i(f\cdot P)\rforal f\in C_0(Y),\,\,\,i=1,2,
\eneq
where $\tau$ is the tracial state on $M_n.$

Note that $\{Y_k\setminus Y_{k-1}: k=1,2,....,4N_0^2\}$ is a family of $4N_0^2$ disjoint Borel sets.
There are at least $2N_0$ of $k's$ such that
\beq\label{8-N-4-16}
\mu_1(Y_k\setminus Y_{k-1})<1/N_0.
\eneq
It is clear then there is at least one of them satisfies
\beq\label{8-N-4-17}
\mu_i(Y_k\setminus Y_{k-1})<1/N_0,\,\,\, i=1,2.
\eneq
Fix one of such $k.$

We may write that
\beq\label{8-N-4-18}
\phi_1=\Sigma^1_\pi\oplus \Sigma_b^1\oplus \Sigma^1_s\oplus \Sigma^1_I\andeqn \phi_2=\Sigma^2_\pi\oplus
\Sigma_b^2\oplus \Sigma^2_s\oplus \Sigma^2_I,
\eneq
where $\Sigma^1_I$  and $\Sigma^2_I$ are finite direct sums of the form $\pi_{x,j}$ for $x\in Y_{k-1},$
$\Sigma_s^1$ and $\Sigma_s^2$ are finite direct sums of the from $\pi_{x,j}$ for $x\in Y_k\setminus Y_{k-1},$
$\Sigma_b^1$ and $\Sigma_b^2$ are    finite direct sum of the form
$\pi_{x,j}$ for $x\in Y\setminus Y_k$
and $\Sigma^1_{\pi}$ and $\Sigma^2_{\pi}$ are finite direct sum of the form
 ${\bar\pi}_{x,i}$
 %for $x\in X^0$
 given by irreducible representations of $A/I$ {{(note that these $\pi_{x,j}$ or ${\bar \pi}_{x,i}$ can be repeated)}}.

Define $\psi_I^{1,0}, \psi_I^{2,0}: B_{0,k-1}\to M_n$ by
\beq\label{8-N-4-19}
\psi_I^{i,0}(f)=\Sigma_I^i(f)\tforal f\in B_{0,k-1},\,\,\,i=1,2.
\eneq
By the choice of ${\cal H}_2,$
we estimate that
\beq\label{8-N-4-20}
&&|\tau\circ \psi_I^{1,0}(1_{B_{0,k-1}})-\tau\circ \psi_I^{2,0}(1_{B_{0,k-1}})|\\
&\le& |\tau\circ \Sigma_I^1(f_{0,k}\cdot P)-\tau\circ \phi_1(f_{0,k}\cdot P)|+|\tau\circ \phi_1(f_{0,k}\cdot P)-\tau\circ \phi_2(f_{0,k}\cdot P)|\\
&+&|\tau\circ \phi_2(f_{0,k}\cdot P)-\Sigma_I^2(f_{0,k}\cdot P)|
<1/N_0+\dt+1/N_0\\
&\le & \dt+\Delta(\widehat{f_{00}}\cdot P)\dt_1\min\{\Delta_{I,0}(\hat{g}): g\in {\cal H}_{1,I,0}\}/32KN_{00}.
\eneq
It follows from \ref{8-N-0} that there are two mutually equivalent projections $p_{1,0}$ and $p_{2,0}\in M_n$ such that
$p_{i,0}$ commutes with $\psi_I^{i,0}(f)$ for all $f\in B_{0,k-1}$ and $p_{i,0}\psi_I^{i,0}(1_{B_{0,k-1}})=p_{i,0}.$ $i=1,2,$
and
\beq\label{8-N-4-21}
\hspace{-0.2in}0\le \tau\circ \psi_I^{i,0}(1_{B_{0,k-1}})-\tau(p_{i,0})<\dt+\Delta(\widehat{f_{00}\cdot P})\dt_1\min\{\Delta_{I,0}(\hat{g}): g\in {\cal H}_{1,I,0}\}/32{{K}}N_{00}+T/n,
\eneq
$i=1,2.$
Since $Y_0\subset Y_{k-1},$  ${\rm supp}(f_{00}){ =Y_0}\subset Y_{k-1}.$
Therefore, using  (\ref{8-N-4-5+}),
\beq\label{10-5-N1}
\tau\circ \psi_I^{1,0}(1_{B_{0,k-1}})\ge \tau\circ \psi_I^{1,0}(f_{00}\cdot P)\ge \Delta(\widehat{f_{00}\cdot P})>
\max\{\sigma_0, 8/(N_0\dt_1)\}.
\eneq
Hence
\beq\label{10-5-N2}
\tau(p_{2,0})>\max\{\sigma_0, 8/(N_0\dt_1)\}/2.
\eneq
Put $q_{i,0}=\psi_I^{i,0}(1_{B_{0,k-1}})-p_{i,0},$ $i=1,2.$
There is a unitary $U_1\in M_n$ such that $U_1^*p_{1,0}U_1=p_{2,0}.$
Define $\psi_I^1: B_{0,k-1}\to p_{2,0}M_np_{2,0}$ by
$\psi_I^1(f)=U_1^*p_{1,0}\psi_I^{1,0}(f)U_1$ for all $f\in B_{0,k-1}$ and define
$\psi_I^2: B_{0,k-1}\to p_{2,0}M_np_{2,0}$ by $\psi_I^2(f)=p_{2,0}\psi_I^{2,0}(f)$ for all $f\in B_{0,k-1}.$
We compute
that (using (\ref{8-N-4-21}) among other items)
\beq\label{8-N-4-22}
\tau\circ \psi_I^1(g\circ r_{k-1}^{-1}) &\ge & \Delta(\widehat{f_{0,k-1}}\cdot 1_A\cdot g)-\min\{\Delta(\hat{g}: g\in {\cal H}_{1,I,0}\}/2N_{00}\\
&>& \bt \Delta(\widehat{f_{0,k-1}}\cdot 1_A\cdot g)=\Delta_{I,0}(g\circ r_{k-1}^{-1})
\eneq
for all $g\in {\cal H}_{1,I,0}.$
Therefore
\beq\label{8-N-4-23}
t\circ \psi_I^1(g)\ge \Delta_{I,k}(g)\tforal g\in {\cal H}_{1,I,k},
\eneq
where $t$ is the tracial state on $p_{2,0}M_np_{2,0}.$
We also estimate that
%, with $g'=g\circ r_k^{-1},$
\beq\label{8-N-4-24}
\hspace{-0.5in}|t\circ \psi_I^1(g)-t\circ \psi_I^2(g)| &=& (1/\tau(|\tau\circ \psi_I^1(g)-\tau\circ \psi_I^2(g)|\\
&\le&  (1/\tau(p_{1,0}))|\tau\circ \psi_I^1(g)-\tau\circ \psi_I^{1,0}(g)|\\
&+& (1/\tau(p_{1,0}))|\tau\circ \psi_I^{1,0}(g)-
\tau(\phi_1(f_{0,k}\cdot 1_A\cdot  g))|\\
&+& (1/\tau(p_{2,0}))|\tau(\phi_1(f_{0,k}\cdot 1_A\cdot g))-\tau(\phi_2(f_{0,k}\cdot 1_A\cdot g))|\\
&+& (1/\tau(p_{2,0}))|\tau(\phi_2(f_{0,k}\cdot 1_A\cdot g)-\tau\circ \psi_I^{2,0}(g)|\\
&+&|\tau\circ \psi_I^{2,0}(g)-\tau\circ \psi_I^2(g)|\\\label{8-N-4-24+}
&<& (1/\tau(p_{2,0}))(\tau(q_{1,0})+1/N_0+ \dt +1/N_0+\tau(q_{2,0}))<\dt_1
%\\&<& (1/\sigma_0)(1/N_0+\dt+1/N_0)<\dt_1
\eneq
for all $g\in {\cal H}_{2,I,k}$ (the last step uses (\ref{10-5-N2})).  Note also that $p_{2,0}$ has rank at least $N^I.$
It follows (by applying \ref{Aug-N-2}) that there are mutually orthogonal projections $e_0^I,e_1^I,e_2^I,...,e_{2K}^l\in
p_{2,0}M_np_{2,0}$ such that $e_0^I+\sum_{i=1}^{2K}e_i^I=p_{2,0},$ $e_0^I\lesssim e_1^I$ and
$e_j^I$ are equivalent to $e_1^I,$ two  unital \hm s $\psi_{1,I,0}, \psi_{2,I,0}: {{B_{0,k-1}}}\to e_0^IM_ne_0^I,$ a
%B_{0,k}\to e_0^IM_ne_0^I,$ a
unital \hm\, $\psi_{I}: {{B_{0,k-1}}}\to e_1^IM_ne_1^I$ and a
%B_{0,k}\to e_1^IM_ne_1^I$ and a
unitary $u_1\in p_{2,0}M_np_{2,0},$
such that
\beq\label{8-N-4-25}
\|{\rm Ad}\, u_1\circ \psi_I^1(f)-(\psi_{1,I,0}(f)\oplus {\rm diag}(\overbrace{\psi_{I}(f), \psi_I(f),...,\psi_I(f)}^{2K}))\|<\ep_0/16\\
\andeqn\,\,\,
\|\psi_I^2(f)-(\psi_{2,I,0}(f)\oplus {\rm diag}(\overbrace{\psi_{I}(f), \psi_I(f),...,\psi_I(f)}^{2K}))\|<\ep_0/16
\eneq
for all $ f\in {{{\cal F}_{I,k-1}'}}.$
%{\cal F}_{I,k}'.$

 By (\ref{8-N-4-8}), the above implies that
\beq\label{8-N-4-26}
\|{\rm Ad}\, u_1\circ \psi_I^1(f)-(\psi_{1,I,0}(f)\oplus {\rm diag}(\overbrace{\psi_{I}(f), \psi_I(f),...,\psi_I(f)}^{2K}))\|<\ep_0/8\\
\label{8-N-4-26+}
\andeqn\,\,\,
\|\psi_I^2(f)-(\psi_{2,I,0}(f)\oplus {\rm diag}(\overbrace{\psi_{I}(f), \psi_I(f),...,\psi_I(f)}^{2K}))\|<\ep_0/8
\eneq
for all $ f\in {{{\cal F}_{I,k-1}}}.$
%{\cal F}_{I,0}.$

For each $x\in X\setminus Y_k$ such that $\pi_{x,j}$ appeared in ${{\Sigma_b^1}},$ or ${{\Sigma_b^2}},$
%$\Sigma_s^1,$ or $\Sigma_s^2,$
%there is ${\bar x}\in {\bar X}$ such that,
by (\ref{8-N-4-4}),
\beq\label{8-N-4-27}
\|\pi_{x,j}(f)-\pi_{x,j}\circ s^d\circ \pi_I(f))\|<\ep_0/16\tforal f\in {\cal F}.
\eneq
%where $\sum_{x,i'}$ is a finite direct sum of some ${\bar \pi}_{{\bar x},i}.$
%Denote by
Define $\Sigma_{\pi,b,i}=\Sigma_b^i\circ s^d: A/I\to M_n,$ $i=1,2.$
%the finite direct sum of $\sum_{i=1}^{k_2}\sum_{x,j,i}{\bar \pi}_{{\bar x}, j,i}\circ \pi_I$ for
%which $\pi_{x,j}$ appeared in $\Sigma_\pi^1,$ and
%$\Sigma_{\pi, b,2}$ the finite direct sum of $\sum_{i=1}^{k_2}\sum_{x,j,i}{\bar \pi}_{{\bar x}, j,i}\circ \pi_I$ for
%which $\pi_{x,j}$ appeared in $\Sigma_b^2.$

Define $\Phi_1: A/I\to (1-p_{2,0})M_n(1-p_{2,0})$ by
\beq\label{8-N-4-28}
\Phi_1(f)={\rm Ad}\, U_1\circ (\Sigma_\pi^1\oplus \Sigma_{\pi, b,1})(f)\rforal f\in A/I.
\eneq
Define $\Phi_2: A/I\to (1-p_{2,0})M_n(1-p_{2,0})$ by
\beq\label{8-N-4-29}
\Phi_2(f)=(\Sigma_\pi^1\oplus \Sigma_{\pi, b,2})(f)\rforal f\in A/I.
\eneq
Note that
\beq\label{8-N-4-29+}
\Phi_1(1_{A/I})=\Sigma_\pi^1(g_{0,k})\oplus \Sigma_b^1(g_{0,k})\andeqn
\Phi_2(1_{A/I})=\Sigma_\pi^1(g_{0,k})\oplus \Sigma_b^2(g_{0,k}).
\eneq
We estimate that
\beq\label{8-N-4-30}
\hspace{-0.8in}|\tau\circ \Phi_1(1_{A/I})-\tau\circ \Phi_2(1_{A/I})|
&\le & |\tau\circ \Phi_1(1_{A/I})-\tau\circ \phi_1(g_{0,k})|\\
&&\hspace{-1.3in}+|\tau\circ \phi_1(g_{0,k})-\tau\circ \phi_2(g_{0,k})|+
|\tau\circ \phi_2(g_{0,k})-\tau\circ \Phi_2(g_{0,k})|\\
&<& 1/N_0+\dt+1/N_0\\
&<&\dt+\Delta(\widehat{f_{00}\cdot P})\dt_{ 2}\min\{\Delta_{ \pi}(\hat{g}): g\in {\cal H}_{1,{ \pi}}\}/32N_{00}.\hspace{-0.4in}.
\eneq
It follows from \ref{8-N-0} that there are two mutually equivalent projections $p_{1,1}$ and
$p_{2,1}\in (1-p_{2,0})M_n(1-p_{2,0})$ such that
$p_{i,1}$ commutes with $\Phi_i(f)$ for all $f\in A/I$ and $p_{i,1}\Phi_i(1_{A/I})=p_{i,1}.$ $i=1,2,$
and
\beq\label{8-N-4-32}
0\le \tau\circ \Phi_i(1_{A/I})-\tau(p_{i,1})<\dt+\Delta(\widehat{f_{00}})\dt_{ 2}\min\{\Delta_{\pi}(\hat{g}): g\in {\cal H}_{1,\pi}\}/{ (32N_{00})}+T/n
\eneq
$i=1,2.$

%{ Since ${\rm supp} (g_0)\subset \overline{X^{d/256}}$,
{{Since $g_0(x)=0$ if ${\rm dist}(x, Y_0)\le d/16,$}} we have
$$\tau\circ \Phi_1(1_{A/I}) > \DT(\widehat{g_0\cdot s^d(1)})> \DT_{\pi}(\hat{1})\geq {{\max}}\{\sigma_0, 8/(N_0\dt_2)\},$$
and $\tau(p_{2,1})\geq \frac{\max(\sigma_0, 8/(N_0\dt_2))}{2}$.

Put $q_{i,1}=\Phi_i(1_{A/I})-p_{i,1},$ $i=1,2.$
There is a unitary $U_2\in (1-p_{2,0})M_n(1-p_{2,0})$ such that $U_2^*p_{1,1}U_2=p_{2,1}.$
Define $\Phi_\pi^1: A/I\to p_{2,1}M_np_{2,1}$ by
$\Phi_\pi^1(f)=U_2^*p_{1,1}\Phi_1(f)U_2$ for all $f\in B_{0,k-1}$ and define
$\Phi_\pi^2: A/I\to p_{2,1}M_np_{2,1}$ by $\Phi_\pi ^2(f)=p_{2,1}\Phi_2(f)$ for all $f\in A/I.$

We compute
that
%(using (\ref{8-N-4-21}) among other items)
%{ {I did not see the below has anything to do with (\ref{8-N-4-21})}}
\beq\label{8-N-4-33}
\tau\circ \Phi_\pi^1(g) &\ge & \Delta(\widehat{g_0\cdot 1_A\cdot s^{ d}(g)})-\sigma_0/N_{00}\\
&>& \bt \Delta(\widehat{g_0\cdot 1_A\cdot s^{ d}(g)})=\Delta_{\pi}(g)
\eneq
for all $g\in {\cal H}_{1,\pi}.$
Therefore
\beq\label{8-N-4-34}
t_1\circ \Phi_\pi^1(g)\ge \Delta_{\pi}(g)\tforal g\in {\cal H}_{1,\pi},
\eneq
where $t_1$ is the tracial state on $p_{2,1}M_np_{2,1}.$

We also estimate  (similar to the estimate of (\ref{8-N-4-24+})) that
\beq\label{8-N-4-35}
\hspace{-0.6in}|t_1\circ \Phi_\pi^1(g)-t_1\circ \Phi_\pi^2(g)| &=& (1/\tau(p_{2,1}))|\tau\circ \Phi_\pi^1(g)-\tau\circ \Phi_\pi^2(g)|\\
&\le&  (1/\tau(p_{2,1}))|\tau\circ \Phi_\pi^1(g)-\tau\circ \Phi_1(g)|\\
&+& (1/\tau(p_{2,1}))|\tau\circ \Phi_1(g)-\tau\circ \phi_1(g_{0,k}\cdot 1_A\cdot s^d(g))|\\
&+& (1/\tau(p_{2,1}))|\tau\circ \phi_1(g_{0,k}\cdot 1_A\cdot s^d(g))-\tau\circ \phi_2(g_{0,k}\cdot 1_A\cdot s^d(g))|\\
&+& (1/\tau(p_{2,1}))|\tau\circ \phi_2(g_{0,k}\cdot 1_A\cdot s^d(g))-\tau\circ \Phi_2(g)|\\
&+&(1/\tau(p_{2,1}))|\tau\circ \Phi_2(g)-\tau\circ \Phi_\pi^2(g)|\\
&<&  (1/\tau(p_{2,1}))(1/N_0+ \dt +1/N_0)<\dt_{ 2}
%\\&<& (1/\sigma_0)(1/N_0+\dt+1/N_0)<\dt_1
\eneq
for all $g\in {\cal H}_{2,\pi}.$
It follows from \ref{Aug-N-2} that there are mutually orthogonal projections $e_0^\pi,e_1^\pi,e_2^\pi,...,e_{2K}^\pi\in
p_{2,1}M_np_{2,1}$ such that $e_0^\pi\lesssim e_1^\pi$ and
$e_j^\pi$ are equivalent to $e_1^\pi,$ two  unital \hm s $\psi_{1,\pi,0}, \psi_{2,\pi,0}: A/I\to e_0^\pi M_ne_0^\pi,$ a
unital \hm\, $\psi_{\pi}: A/I\to e_1^\pi M_ne_1^\pi$ and a
unitary $u_2\in p_{2,1}M_np_{2,1},$
such that
\beq\label{8-N-4-36}
\|{\rm Ad}\, u_2\circ \Phi_\pi^1(f)-(\psi_{1,\pi,0}(f)\oplus {\rm diag}(\overbrace{\psi_{\pi}(f), \psi_\pi(f),...,\psi_\pi(f)}^{2K}))\|<\ep_0/16\\
\andeqn\,\,\,
\|\Phi_\pi^2(f)-(\psi_{2,\pi,0}(f)\oplus {\rm diag}(\overbrace{\psi_{\pi}(f), \psi_\pi(f),...,\psi_\pi(f)}^{2K}))\|<\ep_0/16
\eneq
for all $ f\in {\cal F}_{\pi}.$
% Moreover,
%\beq\label{8-N-4-25+}
%\tau\circ\psi_I(g)\ge \bt  {\Delta_{I,0}(\hat{g})\over{K}}\tforal g\in {\cal H}_{0,I,k}'.
%\eneq
Let ${{\psi_\pi^{1}}}: A\to  p_{2,1}M_np_{2,1}$ by
$\psi_\pi^{1}(f)={\rm Ad}\, u_2\circ{\rm Ad}\, U_2 (p_{2,1}(\Sigma_\pi^1\oplus \Sigma_b^1)(f))$ and
define $\psi_\pi^{2}: A\to p_{2,1}M_np_{2,1}$ by
$\psi_\pi^{2}(f)=p_{2,1}(\Sigma_\pi^2\oplus \Sigma_b^2)(f)$
 for all $f\in A.$
Then, by (\ref{8-N-4-27}),
\beq\label{8-N-4-37}
\|\psi_\pi^{i}(f)-(\psi_{i,\pi,0}\circ\pi_I(f)\oplus {\rm diag}(\overbrace{\psi_{\pi}(\pi_I(f)), \psi_\pi(\pi_I(f)),...,\psi_\pi(\pi_I(f))}^{2K}))\|<\ep_0/8
\eneq
for all $f\in {\cal F},$ $i=1,2.$

Put $e_i=e_{2i-1}^{ I}\oplus e_{2i}^I\oplus e_{2i-1}^\pi\oplus e_{2i}^\pi,$ $i=1,2,...,K.$
Define $\psi: A\to  e_1M_ne_1$ by
$$
\psi(f)={\rm diag}(\psi_I(f|_{{{Y_{k-1}}}}),\psi_I(f|_{{{Y_{k-1}}}}), \psi_\pi\circ \pi_I(f), \psi_\pi\circ \pi_I(f))
$$ for all $f\in A.$
%\psi_I(f|_{Y_k})+\psi_\pi\circ \pi_I(f)$ for all $f\in A.$
 By (\ref{8-N-4-16}), (\ref{8-N-4-21}) and (\ref{8-N-4-32}),
 \beq\label{8-N-4-38}
 \tau(q_{i,0})+\tau(q_{i,1})+\tau(\Sigma_s^i(1_A))<  1/64K+ 1/N_0+1/64K<1/16K.
 \eneq
We have, for $f\in A,$
\beq\label{8-N-4-39}
\phi_2(f)={{\psi_\pi^{2}}}(f)\oplus q_{2,1}(\Sigma_\pi^2+\Sigma_b^2)(f)\oplus \Sigma_s^2(f)\oplus \psi_I^2(f|_{{{Y_{k-1}}}})\oplus q_{2,0}{{\psi_I^{2,0}}}(f{{|_{Y_{k-1}}}}).
%\Phi_J^{2,0}(f{\red{|_{Y_{k-1}}}}).
\eneq
Put $e_0=e_0^I\oplus e_0^\pi+q_{2,1}+\Sigma_s^2(1_A)+q_{2,0}.$
Then
\beq\label{8-N-4-50}
\tau(e_0)<\tau(e_0^I)+\tau(e_0^\pi)+1/16K\le \tau(e_1).
\eneq
In other words, $e_0\lesssim e_1.$ Moreover $e_1$ is equivalent to each $e_i,$ $i=1,2,...,K.$
{{Define $h_2: A\to e_0M_ne_0$ by,  for each $f\in A,$
$$
h_2(f)=\psi_{2,I,0}(f|_{Y_{k-1}})\oplus \psi_{2,\pi,0}(\pi_I(f)) \oplus q_{2,1}(\Sigma_\pi^2+\Sigma_b^2)(f)\oplus \Sigma_s^2(f)\oplus q_{2,0}\psi_I^{2,0}(f|_{Y_{k-1}}).
$$
}}
It follows from {{(\ref{8-N-4-26+}), (\ref{8-N-4-37}) and}} above that
%there is a unital \hm\, $h_1: A\to e_0M_ne_0$
%such that
\beq\label{8-N-4-51}
\|\phi_2(f)-(h_{{2}}(f)\oplus {\rm diag}(\overbrace{\psi(f), \psi(f),...,\psi(f)}^K)\|<\ep_0/8\rforal f\in {\cal F}.
\eneq
Similarly, there exists a unitary $U\in M_n$ and a unital \hm\, $h_2: {{A\to}} e_0M_ne_0$ such that
\beq\label{8-N-4-52}
\|{\rm Ad}\, U\circ \phi_1(f)-(h_{{1}}(f)\oplus{\rm diag}(\overbrace{\psi(f), \psi(f),...,\psi(f)}^K)\|<\ep_0/8\rforal f\in {\cal F}.
\eneq
Since we also assume that ${\cal H}_0\subset {\cal F}$ in the above proof, it is easy to check, by the choice of $\ep_0$
and $\bt,$
that
\beq\label{8-N-4-53}
\tau\circ \psi(g)\ge \af{\Delta(\hat(g))\over{K}}\tforal g\in {\cal H}_0.
\eneq

\end{proof}

%xxxxxxxxx

\begin{rem}\label{8-N-4-r}
{\rm
If we assume that $A$ is infinite dimensional, then Lemma \ref{8-N-4} still holds  without the assumption about the integer $N.$  This could be easily seen by taking a larger ${\cal H}_1.$
%It is perhaps important to note that, in the proof of Lemma \ref{8-N-4} , we did not use the assumption that $X\setminus Z$ is connected.  We did use (\ref{Amnotation-1}) and (\ref{Amnotation-2+}).
%Moreover we do not need to assume that $X\setminus Z$ is path connected for Lemma \ref{8-N-4}.
%The condition that $X\setminus Z$ will used in section 16 though.

}
\end{rem}

The following is known and is taken from Theorem 3.9 of \cite{LnAUCT}

\begin{thm}\label{Lnuct}
Let $A$ be a unital separable amenable \CA\,  which satisfies the UCT and let $B$ be a unital \CA.  Suppose that
$h_1,\, h_2: A\to B$ are two \hm s such that
$$
[h_1]=[h_2]\,\,\,{\rm in}\,\,\, KL(A,B).
$$
Suppose that $h_0: A\to B$ is a unital full monomorphism. Then, for any $\ep>0$ and   any
finite subset ${\cal F}\subset A,$ there exits an integer $n\ge 1$ and a unitary $W\in M_{n+1}(B)$
such that
$$
\|W^*{\rm diag}(h_1(a),h_0(a),...,h_0(a))W-{\rm diag}(h_2(a),h_0(a),...,h_0(a))\|<\ep
$$
for all $a\in {\cal F},$  {{and
$W^*pW=q,$}} where
\beq\nonumber
p={\rm diag}(h_1(1_A),h_0(1_A),...,h_0(1_A))\andeqn q={\rm diag}(h_2(1_A),h_0(1_A),...,h_0(1_A)).
\eneq
In particular,  if $h_1(1_A)=h_2(1_A),$ $W\in U(pM_{n+1}(B)p).$
\end{thm}

\begin{proof}
This is a slight variation of Theorem 3.9 of \cite{LnAUCT}.  If $h_1$ and $h_2$ are both unital, then
it is exactly the same as Theorem 3.9 of \cite{LnAUCT}.  So suppose that $h_1$ is not unital. Let $A'=\C\oplus A.$
Choose $p_0=1_B-h_1(1_A)$ and $p_1={\rm diag}(p_0, 1_B).$ Put $B'=p_1M_2(B)p_1.$
Define $h_1': A'\to B'$ by $h_1'(\lambda \oplus a)=\lambda\cdot {\rm diag}(p_0,p_0)\oplus h_1(a)$ for all
$\lambda\in \C$ and $a\in A,$  and define $h_2': A'\to B'$ by
$h_2'(\lambda \oplus a)=\lambda \cdot {\rm diag}(p_0, 1_B-h_2(1_A))\oplus h_2(a)$ for all $\lambda\in \C$ and $a\in A.$
Then $h_1'$ and $h_2'$ are unital and $[h_1']=[h_2']$ in $KL(A', B').$ Define
$h_0': A'\to B'$ by $h_0'(\lambda\oplus a)=\lambda\cdot p_0\oplus h_0(a)$ for all $\lambda\in \C$ and $a\in A.$
Note that $h_0'$ is full in $B'.$
So, Theorem 3.9 of \cite{LnAUCT} applies. It follows that there is an integer $n\ge 1$ and a unitary
$W'\in M_{n+1}(B')$ such that
\beq\nonumber
\|(W')^*{\rm diag}(h_1'(a),h_0'(a),...,h_0'(a))W'-{\rm diag}(h_2'(a),h_0'(a),...,h_0'(a))\|<\min\{1/2, \ep/2\}
\eneq
for all $a\in {\cal F}\cup\{1_A\}.$  In particular,
\beq\label{lnauct-3}
\|(W')^*pW'-q\|<\min\{1/2,\ep/2\}.
\eneq
There is a unitary $W_1\in M_{n+1}(B')$ such that
\beq\label{linauct-3}
\|W_1-1_{M_{n+1}(B')}\|<\ep/2\andeqn W_1^*(W')^*pW'W_1=q.
\eneq
Put $W=W'W_1.$ Then
\beq\label{linatct-4}
\|W^*{\rm diag}(h_1(a),h_0(a),...,h_0(a))W-{\rm diag}(h_2(a),h_0(a),...,h_0(a))\|<\ep
\eneq
 for all $a\in {\cal F}.$ Lemma follows.
\end{proof}

\begin{lem}\label{Lauct2} {\rm (cf.{{ 5.3 of \cite{Lnjotuni}, Theorem 3.1 of}} \cite{GL-almost-map}, \cite{Da1}, 5.9 of \cite{LnAUCT} and
{{Theorem 7.1 of \cite{Lin-hmtp}}})}
  %\\
Let $A$ be a unital separable amenable  \CA\, which satisfies the UCT and let $\Delta: A_+^{q, {\bf 1}}\setminus\{0\} \to (0,1)$ be an order preserving
map.  For any $\ep>0$ and any finite subset ${\cal F}\subset A,$ there exists
$\dt>0,$ a finite subset ${\cal G}\subset A,$
a finite subset ${\cal P}\subset \underline{K}(A),$ a finite subset ${\cal H}\subset A_+^{\bf 1}\setminus \{0\}$ and an integer $K\ge 1$ satisfying the following:
For any two unital ${\cal G}$-$\dt$-multiplicative \morp s $\phi_1, \phi_2: A\to M_n$ (for some integer $n$) and any
unital ${\cal G}$-$\dt$-multiplicative \morp\, $\psi: A\to M_m$ with $m\ge n$ such that
\beq\label{Lauct-1}
\tau\circ \psi(g)\ge \Delta(\hat{g})\tforal g\in {\cal H}\tand\,
[\phi_1]|_{\cal P}=[\phi_2]|_{\cal P},
\eneq
there exists a unitary $U\in M_{Km+n}$ such that
\beq\label{Lauct-2}
\|{\rm Ad}\, U\circ (\phi_1\oplus \Psi)(f)-(\phi_2\oplus \Psi)(f)\|<\ep\tforal f\in A,
\eneq
where
$$
\Psi(f)={\rm diag}(\overbrace{\psi(f), \psi(f),...,\psi(f)}^K)\tforal f\in A.
$$
\end{lem}

\begin{proof}
This follows from \ref{Lnuct}.
%Theorem 3.9 of \cite{Lauct}. Note that Theorem 3.9 holds for $A$ without assuming $h_1$ and $h_2$
%are unital by choosing a larger $K$ if necessary (but integer $K$ depends only on the smallest rank of irreducible representations).

Fix $\Delta$ as given.
Suppose the lemma is false. Then there exist  $\ep_0>0$ and a finite subset
${\cal F}_0\subset A,$  an increasing sequence of finite subsets
$\{{\cal P}_n\}$ of $\underline{K}(A)$ whose union is $\underline{K}(A),$  an increasing sequence
of finite subsets $\{{\cal H}_{n}\}\subset A_+^{\bf 1}\setminus \{0\}$ whose union is dense
in $A_+^{\bf 1}$ {{and if $a\in {\cal H}_n$ and $f_{1/2}(a)\not=0,$ then
$f_{1/2}(a)\in {\cal H}_{n+1},$}}
three sequences of increasing   integers  $\{R(n)\}$ and  $\{r(n)\},\,\{s(n)\}$ (with $s(n)\ge r(n)$),
 two sequences of \morp s $\phi_{1,n}, \phi_{2,n}, : A\to M_{r(n)}$
 with the properties that
 \beq\label{Lauct-3-}
 [\phi_{1,n}]|_{{\cal P}_n}&=&[\phi_{2,n}]|_{{\cal P}_n}\andeqn\\
 \lim_{n\to\infty}\|\phi_{i,n}(ab)-\phi_{i,n}(a)\phi_{i,n}(b)\|&=&0\rforal a,b \in A,\,\,\,i=1,2,
 \eneq
 a sequence of unital \morp s
 $\psi_n: A\to M_{s(n)}$ with the properties
 that
 \beq\label{Lauct-3}
&& \tau_n\circ \psi_n(g)\ge \Delta(\hat{g})\rforal g\in {\cal H}_{n}\andeqn\\
 && \lim_{n\to\infty}\|\psi_{n}(ab)-\psi_{n}(a)\psi_{n}(b)\|=0\rforal a,b \in A
 \eneq
such that
\beq\label{Lauct-4}
\inf\{ \sup\{\|{\rm Ad}\, U_n\circ (\phi_{1,n}(f)\oplus {\tilde \psi}^{R(n)}_n(f))-(\phi_{2,n}(f)\oplus {\tilde \psi}^{R(n)}_n(f)\|: f\in {\cal F}\}\}\ge \ep_0,
\eneq
where $\tau_n$ is the normalized trace  on $M_{r(n)},$
${\tilde \psi}^{(R(n))}(f)={\rm diag}(\overbrace{\psi_n(f),\psi_n(f),...,\psi_n(f)}^{R(n)})$ for all $f\in A,$
and the infimum  is taken among all unitaries $U_n\in M_{r(n)}.$

Note that, by (\ref{Lauct-3}), since $\{{\cal H}_n\}$ is increasing, for any $g\in {\cal H}_{n}\subset A_+^{\bf 1},$  we compute that
\beq\label{Lauct-5}
\tau_m(p_m)\ge \Delta(\hat{g})/2,
\eneq
where $p_m$ is the spectral projection of $\psi_m(g)$ corresponding to the subset $\{\lambda>\Delta(\hat{g})/2\}$
for all $m\ge n.$
It follows that (for all sufficiently large $m$)  there are  element $x_{g,i,m}\in M_{s(m)}$ with $\|x_{g, i,m}\|\le 1/\Delta(\hat{g}),$
$i=1,2,...,N(g)$ such that
\beq\label{Lauct-6}
\sum_{i=1}^{N(g)} x_{g,i,m}^*\psi_m(g)x_{g,i,m}=1_{s(m)},
\eneq
where $1\le N(g)\le 1/\Delta(\hat{g})+1.$
Put $X_{g,i}=\{x_{g,i,m}\},$ $i=1,2,...,N(g).$ Then $X_{g,i}\in \prod_{n=1}^{\infty}M_{r(n)}.$
Let $Q(\{M_{r(n)}\})=\prod_{n=1}^{\infty}M_{r(n)}/\bigoplus_{n=1}^{\infty}M_{r(n)},$
$Q(\{M_{s(n)}\})=\prod_{n=1}^{\infty}M_{s(n)}/\bigoplus_{n=1}^{\infty}M_{s(n)},$ and let
$\Pi_1: \prod_{n=1}^{\infty}M_{r(n)}\to \prod_{n=1}^{\infty}M_{r(n)}/\bigoplus_{n=1}^{\infty}M_{r(n)},$
$\Pi_2:\prod_{n=1}^{\infty}M_{s(n)}\to \prod_{n=1}^{\infty}M_{s(n)}/\bigoplus_{n=1}^{\infty}M_{s(n)}$ be  quotient maps.
Denote by $\Phi_i: A\to Q(\{M_{r(n)}\})$ the \hm s
$\Pi_1\circ \{\phi_{i,n}\}$ and denote by ${\bar \psi}: A\to Q(\{M_{s(n)}\})$ the \hm\, $\Pi_2\circ \{\psi_n\}.$
For each $g\in \cup_{n=1}^{\infty} {\cal H}_n,$
\beq\label{Lauct-7}
\sum_{i=1}^{N(g)}\Pi_2(X_{i,g})^*{\bar \psi}(g)\Pi_2(X_{i,g})=1_{Q(\{M_{s(n)}\})}.
\eneq
{{Note that if $g\in \cup_{n=1}{\cal H}_n$ and $f_{1/2}(g)\not=0,$ then $g_{1/2}(g)\in \cup_{n=1}{\cal H}_n,$ and
$\cup_{n=1}{\cal H}_n$ is dense in $A_+^{\bf 1}.$}} This implies  that ${\bar \psi}$ is full.
Note that both $\prod_{n=1}^{\infty} M_{r(n)}$ and $Q(\{M_{r(n)}\}$ have stable rank one and real rank zero.
One computes that
\beq\label{Lauct-10}
[\Phi_1]=[\Phi_2]\,\,\,{\rm in}\,\,\, KL(A, Q(\{M_{r(n)}\}).
\eneq
By applying Theorem \ref{Lnuct} (Theorem 3.9 of \cite{LnAUCT}), there exists an integer $K\ge 1$ and a unitary $U\in PM_{K+1}(Q(\{M_{s(n)}\})P,$ where $P={\rm diag}(1_{Q(\{M_{r(n)}\}}, 1_{M_K(Q(\{M_{s(n)}\}))}),$
such that
\beq\label{Lauct-11}
\|{\rm Ad}\, U\circ (\Phi_1(f)\oplus {\rm diag}(\overbrace{{\bar \psi}(f),...,{\bar \psi}(f)}^K)
-(\Phi_2(f)\oplus {\rm diag}(\overbrace{{\bar \psi}(f),...,{\bar \psi}(f)}^K)\|<\ep_0/2
\eneq
for all $f\in {\cal F}_0.$
It follows that there are unitaries
$$
\{U_n\}\in \prod_{n=1}^{\infty}M_{r(n) +Ks(n)}
$$
such that, for all
large $n,$
\beq\label{Lauct-12}
\hspace{-0.4in}\|{\rm Ad}\, U_n\circ \phi_{1,n}(f)\oplus {\rm diag}(\overbrace{\psi_n(f),...,\psi_n(f)}^K)-
\phi_{2,n}(f)\oplus {\rm diag}(\overbrace{\psi_n(f),...,\psi_n(f)}^K)\|<\ep_0/2
\eneq
for all $f\in {\cal F}_0.$  This leads to a contradiction with (\ref{Lauct-4}) when we choose $n$ with $R(n)\ge K.$
\end{proof}

\begin{rem}
This lemma holds in a much more general  setting and variations of it has appeared. We state this version here for  our immediate purpose (see \ref{Suni} and part (1) of \ref{Rsuni} for more comments).
%In an earlier draft, it was stated for $A$ being residually finite dimensional. But the same proof
%works for general cases. The argument to obtain an independent integer $K$ as stated here has been appeared
%before (see the proof of Theorem 4.8 of \cite{Lnjotuni} and \cite{Da1}). {\red{ That part of the  proof presented here  is a repetition of some previous
%proof. }}
%The condition on trace (the first part of (\ref{Lauct-1})) is used for fullness of the map ${\bar \psi}.$
%It can be replaced by the one stated in Theorem 5.9 of \cite{LnAUCT} when the target algebras
%have no traces. In fact, the statement of Theorem 5.9 of \cite{LnAUCT} holds with the additional
%feature that the integer $k$ there does not depend on maps  $\phi,$ $\psi$ and $\sigma$ as long as
%they satisfy other conditions stated in Theorem 5.9 of \cite{LnAUCT} as that proof referred to that
%of Theorem 4.8 of \cite{Lnjotuni} which had dealt with the independence of integer $k.$
%Note that the target algebras have to be in ${\bf C}_{l,b,M}$ for given $l, b$ and $M$ as stated in
%Theorem 5.9 of \cite{LnAUCT}.
\end{rem}

It should be noted that, in the following statement the integer $L$ and $\Psi$ depend not only on $\ep,$ ${\cal F},$
${\cal G},$  but also depend on $B$ as well as $\phi_1$ and $\phi_2.$

\begin{lem}\label{Newstableuniq}
Let $C$ be a unital amenable separable residually finite dimensional \CA\, which satisfies the UCT. For any $\ep>0,$ any finite subset
${\cal F}\subset C,$ there exists a finite subset ${\cal G}\subset C,$
$\dt>0,$ a finite subset ${\cal P}\subset \underline{K}(C)$
%a finite subset ${\cal H}\subset A_+\setminus\{0\}$
%and an integer $L\ge 1$
satisfying the following:
For any unital ${\cal G}$-$\dt$-multiplicative \morp s $\phi_1, \phi_2:
C\to A$ (for any unital \CA\, $A$)
%and any unital \hm\,
%$\psi: A\to M_k\subset M_k(B)$
such that
\beq\label{New-1}
[\phi_1]|_{\cal P}=[\phi_2]|_{\cal P},
%tr(\psi(h))>0,
\eneq
there exist an integer $L\ge 1,$ a  unital \hm~ $\Psi: C\to M_L(\C)\subset M_L(A)$, and  a unitary $U\in U(M_{L+1}(A))$   such that
\beq\label{New-2}
\|{\rm Ad}\, U\circ {\rm diag}(\phi_1(f), \Psi(f))-{\rm diag}(\phi_2(f),
\Psi(f)\|<\ep
\eneq
for all $f\in {\cal F}.$
%where
%$\Psi: C\to M_L(\C)\subset M_L(A)$ is a unital \hm.
%Moreover one can choose $U\in U_0(M_{L+1}(B)).$
\end{lem}

\begin{proof}
The proof is almost the same as that of Theorem 9.2 of \cite{Linajm}.
%Let $\gamma_n$ be a separating sequence of finite dimensional representations of $C.$
%Let $D_n=\gamma_n(C).$ There is, for each $n,$ an integer $L_n\ge 1$ such that
%$D_n\subset M_{L_n}(\C)$ with $1_{D_n}=1_{M_{L_n}}.$ So we may view $\gamma_n: C\to M_{L_n}(\C)$
%is a unital \hm.
%We first proof the first part of the statement (we will show that we can choose $U\in U_0(M_{L+1}(B))$
%at the end of the proof).
Suppose that the lemma is false.
We then obtain  positive number $\ep_0>0$ and a finite ${\cal F}_0\subset C,$ a sequence  of finite
subsets ${\cal P}_n\subset \underline{K}(C)$ with ${\cal P}_n\subset {\cal P}_{n+1}$ and
$\cup_n {\cal P}_{n+1}=\underline{K}(C),$ a sequence of unital \CA s $\{A_n\},$ a sequence
of unital \morp s $\{L_n^{(1)}\}$ and $\{L_n^{(2)}\}$ (from $C$ to $A_n$)
such that
\beq\label{12/31/24-1}
&&\hspace{-0.3in}\lim_{n\to\infty}\|L_n^{(i)}(ab)-L_n^{(i)}(a)L_n^{(i)}(b)\|=0\rforal a, b\in C,\\
&&\hspace{-0.3in}{[}L_n^{(1)}{]}|_{{\cal P}_n}={[}L_n^{(2)}{]}|_{{\cal P}_n},\\
&&\hspace{-0.3in}\inf\{\sup\{\|u_n^*{\rm diag}(L_n^{(1)}(a), \Psi_n(a)){{u_n}}-{\rm diag}(L_n^{(2)}(a), \Psi_n(a)\|: a\in {\cal F}_0\}\ge \ep_0,
\eneq
where the infimum
is taken among all integers $k >1,$  all possible unital homomorphisms
$\Psi_n: C\to  M_k(\C)$ and all possible unitaries $U \in M_{k+1}(A_n).$  We may assume
that $1_C\in {\cal F}.$
Define $B_n=A_n\otimes {\cal K},$ $B=\prod_{n=1}^{\infty} B_n$ and $Q_1=B/\bigoplus_{n=1}^{\infty} B_n.$
Let $\pi: B\to Q_1$ be the quotient map. Define $\phi_j:C\to B$ by $\phi_j(a)=\{L_n^{(j)}(a)\}$ and
define  $\bar{\phi_j}=\pi\circ \phi_j,$ $j=0,1.$  Note that $\bar{\phi_j}: C\to Q_1$ is a \hm\,
%by  2.9 of \cite{GLjot} (see also p.990 of \cite{GLitj}),  since $B_n$ is stable,
as in the proof of 9.2 of \cite{Linajm},
we have
$$
[{\bar{ \phi_1}}]=[\bar{\phi_2}]\,\,\, KL(C, Q).
$$
%

%We then proceed the proof of
%Theorem 9.2 of \cite{Linajm} to (e 9.295) in the proof of Theorem 9.2 of \cite{Linajm}.
%We then obtain $\ep_0>0$ and a finite finite
Fix an irreducible representation  $\phi_0': C\to M_{r}.$
Denote by $p_n$ the unit of  the unitization ${\tilde B_n}$ of $B_n,$ $n=1,2,....$
%and $q_n=p_n-1_{A_n}.$
Define a  \hm\,
$\phi_0^{(n)}: C\to M_r({\tilde B}_n)=M_r\otimes {\tilde B}_n$ by
$\phi_0^{(n)}(c)=\phi_0'(c)\otimes 1_{{\tilde B}}$ for all $c\in C.$ Put
$$
e_A=\{1_{A_n}\},\,\,  P=\{1_{M_r({\tilde B}_n)}\}+e_A.
$$
Put also  $Q_2=\pi(P)M_{r+1}({\tilde Q_1})\pi(P)$ and
define
${\bar \phi}_j={\bar \phi}'_j\oplus \pi\circ \{\phi_0^{(n)}\},$ $j=1,2.$
Then
\beq\label{Jukly31-1}
[{\bar \phi}_1]=[{\bar \phi}_2]\,\,\,{\rm in}\,\,\, KK(C, Q_2).
\eneq
The point to add $\pi\circ \{\phi_0^{(n)}\}$ is that, now, ${\bar \phi}_1$ and ${\bar \phi}_2$ are unital.
%Let $\gamma_n': C\to M_{L_n}(\C)\subset M_{L_n}(A_n)\subset A_n\otimes {\cal K}=B_n$ be defined by $\gamma_n.$
%Define $\Psi: C\to Q_1
It follows from Theorem 4.3 of \cite{Da1} that there exists an integer $K>0,$ a unitary $u\in M_{1+K}(Q_1)$ and
a unital \hm\, $\psi: C\to M_{K}(\C)\subset M_{K}(Q_2)$ (by identifying $M_{K}(\C)$ with the unital subalgebra of
$M_K(Q_2)$) such that
\begin{equation*}
{\rm Ad}\, u\circ {\rm diag}({\bar \phi}_1, \psi)\approx_{\ep_0/4} {\rm diag}({\bar \phi}_2, \psi) \,\,\,{\rm on}\,\,\, {\cal F}_0.
\end{equation*}
There exists a unitary $V=\{V_n\}\in M_{1+K}(PM_{r+1}({\tilde B})P)$ such that
$\pi(V)=u.$ It follows (by identifying $M_K(\C)$ with $M_K(\C)\otimes 1_{Q_2}$) that for all sufficiently large $n,$
\begin{equation*}
{\rm Ad}\, V_n\circ {\rm diag}(L_1^{(n)}\oplus \phi_n^{(n)}, \psi)\approx_{\ep_0/3} {\rm diag}(L_2^{(n)}\oplus \phi_n^{(n)}, \psi)
\,\,\,{\rm on}\,\,\, {\cal F}_0.
\end{equation*}
Denote by, for each integer $k\ge 1,$  $e_{n,k,0}={\rm diag}(\overbrace{1_{A_n}, 1_{A_n},...,1_{A_n}}^k)\in A_n\otimes {\cal K}=B_n,$
\beq\label{July31-3}
e_{n,k}'&=&{\rm diag}(1_{A_n},\overbrace{e_{n,k,0},e_{n,k,0},...,e_{n,k,0}}^r)\in PM_{1+r}(B_n)P\andeqn\\
e_{n,k}''&=& {\rm diag}(\overbrace{e_{n,k}',e_{n,k}',...,e_{n,k}'}^K)\in M_{K}(PM_{1+r}(B_n)P).
\eneq
It should be noted that
$e_{n,k}''$ commutes with $\psi$ and $e_{n,k}'$ commutes with $\phi_n^{(0)}.$ Put $e_{n,k}=e_{n,k}'\oplus e_{n,k}''$
in $M_{1+K}(PM_{1+r}(B_n)P).$ Then $\{e_{n,k}\}$ forms an approximate identity for
$M_{1+K}(PM_{1+r}(B_n)P).$ Note that $V_n\in M_{1+K}(PM_{r+1}({\tilde B})P).$ It is easy to check that
\beq\label{Aug1st-4}
\lim_{k\to\infty}\|[V_n, e_{n,k}]\|=0.
\eneq
It follows that there exists a unitary $U_{n,k}\in e_{n,k}M_{1+K}(PM_{1+r}(B_n)P)e_{n,k}$  for each $n$ and $k$ such that
\beq\label{Aug1st-5}
\lim_{k\to\infty}\|e_{n,k}V_ne_{n,k}-U_{n,k}\|=0.
\eneq
For each $k,$ there is $N(k)=rk+K(rk+1)$ such that
\beq\label{Aug1st-6}
M_{N(k)}(A_n)=((e_{n,k}'-1_{A_n})\oplus e_{n,k}'')M_{1+K}(PM_{r+1}(B_n)P)((e_{n,k}'-1_{A_n})\oplus e_{n,k}'').
\eneq
Moreover $e_{n,k}M_{1+K}(PM_{1+r}(B_n)P)e_{n,k}=M_{N(k)+1}(A_n).$
Define
$\Psi_n(c)=(e_{n,k}'-1_{A_n})\phi_0^{(n)}(c)(e_{n,k}'-1_{A_n})\oplus e_{n,k}''\psi(c)e_{n,k}$ for $c\in C.$ Then, for all large
$k$ and large $n,$
\beq\label{Aug1st-7}
{\rm Ad}\, U_n\circ {\rm diag}(L_1^{(n)}, \Psi_n)\approx_{\ep_0/2}{\rm diag}(L_2^{(n)}, \Psi_n)\,\,\,{\rm on}\,\,\, {\cal F}_0.
\eneq
This gives a contradiction.
% (to the same assumption as (e.9.291) in the proof of Theorem 9.2 of \cite{Linajm}).}
%which proves the first part of the proof.
%To see that one can choose $U\in U_0(M_{L+1}(A))$ in the statement,
%suppose that (\ref{New-2}) holds.
%Choose that
%
%
\end{proof}

\begin{thm}\label{UniqAtoM}
Let $A\in {\cal {\bar D}}_s$ be a unital subhomogeneous \CA\, {{and}}
let $\Delta: A_+^{q,{\bf 1}}\setminus \{0\}\to (0,1)$ be an order preserving map.

For any $\ep>0$ and finite subset ${\cal F}\subset A,$ there exists $\dt>0,$ a finite subset
${\cal P}\subset \underline{K}(A),$ a finite subset ${\cal H}_1\subset A_+^{\bf 1}\setminus \{0\}$ and
a finite subset ${\cal H}_2\subset A_{s.a.}$ satisfying the following:

If $\phi_1, \phi_2: A\to M_n$ are two unital \hm s such that
\beq\label{UnAM-1}
&&[\phi_1]|_{\cal P}=[\phi_2]|_{\cal P},\\\label{UNAM-1+}
&&\tau\circ \phi_1(g)\ge \Delta(\hat{g})\tforal g\in {\cal H}_1\tand\\
&&|\tau\circ \phi_1(h)-\tau\circ \phi_2(h)|<\dt\tforal h\in {\cal H}_2,
\eneq
then there exists a unitary $u\in M_n$ such that
\beq\label{UnAM-2}
\|{\rm Ad}\, u\circ \phi_1(f)-\phi_2(f)\|<\ep\tforal f\in {\cal F}.
\eneq
\end{thm}

\begin{proof}
If $A$ has finite dimensional, the lemma is known. So, in what follows, we will assume
that $A$ is infinite dimensional.

Define $\Delta_0: A_+^{q, {\bf 1}}\setminus \{0\}\to (0,1)$ by
$\Delta_0=(3/4)\Delta.$
Fix $\ep>0$ and a finite subset ${\cal F}\subset A.$
Let ${\cal P}\subset \underline{K}(A)$ be a finite subset, ${\cal H}_0\subset A_+^{q, {\bf 1}}\setminus \{0\}$
(in place of ${\cal H}$) be a finite subset and an integer $K\ge 1$ be {{as}} required by \ref{Lauct2} for $\ep/2$ (in place of $\ep$),
${\cal F}$ and $\Delta_0.$

Choose $\ep_0>0$  and a finite subset ${\cal G}\subset A$ such that $\ep_0<\ep$ and
\beq\label{UnAM-5}
[\Phi_1']|_{\cal P}=[\Phi_2']|_{\cal P}
\eneq
for any pair of unital \hm s from $A,$  provided that
\beq\label{UnAM-6}
\|\Phi_1'(g)-\Phi_2'(g)\|<\ep_0\rforal g\in {\cal G}.
\eneq
We may assume that ${\cal F}\subset {\cal G}$ and $\ep_0<\ep/2.$

Let $\af=3/4.$  Let $N\ge 1$ be an integer, $\dt_1>0$ (in place of $\dt$), ${\cal H}_1\subset A_+^{\bf 1}\setminus \{0\}$  be a finite subset and ${\cal H}_2\subset A_{s.a.}$ be a finite subset {{as}} required by \ref{8-N-4}
for $\ep_0/2$ (in place $\ep$), ${\cal G}$ (in place of ${\cal F}$), ${\cal H}_0,$ $K$ and $\Delta_0$ (in place of $\Delta$).
By choosing larger ${\cal H}_1,$ since $A$ has infinite dimensional, we may assume
that ${\cal H}_1$ contains at least $N$ many mutually orthogonal non-zero positive elements.

Now suppose that $\phi_1, \phi_2$ are two unital \hm s satisfying the assumption
for the above ${\cal P},$ ${\cal H}_1$ and ${\cal H}_2.$
The assumption (\ref{UNAM-1+}) implies that $n\ge N.$
By applying \ref{8-N-4}, we obtain a unitary $u_1\in M_n,$ mutually orthogonal non-zero projections
$e_0, e_1, e_2,...,e_K\in M_n$ with $\sum_{i=0}^K e_i=1_{M_n},$ $e_0\lesssim e_1,$
$e_1$ are equivalent to $e_i,$ $i=1,2,...,K,$ unital \hm s $\Phi_1,\Phi_2: A\to e_0M_ne_0$ and
a unital \hm\, $\psi: A\to e_1M_ne_1$ such that
\beq\label{UnAM-7}
\|{\rm Ad}\, u_1\circ \phi_1(f)-(\Phi_1(f)\oplus \Psi(f))\|<\ep_0/2\rforal f\in {\cal G},\\\label{UnAM-7+}
\|\phi_2(f)-(\Phi_2(f)\oplus \Psi(f))\|<\ep_0/2\rforal f\in {\cal G}\andeqn\\
\tau\circ \psi(g)\ge (3/4)\Delta(\hat{g)}/K\rforal g\in {\cal H}_0,
\eneq
where $\Psi(a)={\rm diag}(\overbrace{\psi(a), \psi(a),...,\psi(a)}^K)$ for all $a\in A$ and
$\tau$ is the tracial state on $M_n.$

Since $[\phi_1]|_{\cal P}=[\phi_2]|_{\cal P},$ by the choice of $\ep_0$ and ${\cal G},$ we compute that
\beq\label{UNAM-8}
[\Phi_1]|_{\cal P}=[\Phi_2]|_{\cal P}.
\eneq
Moreover,
\beq\label{UnAM-9}
t\circ \psi(g)\ge (3/4)\Delta(\hat{g})\rforal g\in {\cal H}_0,
\eneq
if $t$ is the tracial state of $e_1M_ne_1.$
 By \ref{Lauct2}, there is a unitary $u_2\in M_n$ such that
 \beq\label{UnAM-10}
 \|{\rm Ad}\, u_2\circ (\Phi_1\oplus \Psi)(f)-(\Phi_1\oplus \Psi(f))\|<\ep/2\tforal f\in {\cal F}.
 \eneq
 Put $U=u_2u_1.$ Then, by (\ref{UnAM-7}), (\ref{UnAM-7+}) and (\ref{UnAM-10}),
 \beq\label{UnAM-11}
 \|{\rm Ad}\, U\circ \phi_1(f)-\phi_2(f)\|<\ep_0/2+\ep/2+\ep_0/2<\ep\rforal f\in {\cal F}.
 \eneq
\end{proof}

\begin{lem}\label{Lfullab}
Let $A\in {\cal D}_s$ be a unital \CA\, and let $\Delta: A_+^{q, \bf 1}\setminus \{0\}\to (0,1)$
{{be an order preserving map.}}
Let ${\cal P}_0\subset K_0(A)$ be a finite subset.
Then there exists an integer $N({\cal P}_0){{\ge 1}}$ and a finite subset ${\cal H}\subset A_+^{\bf 1}\setminus \{0\}$ satisfying the following:
For any unital \hm\, $\phi: A\to M_k$ (for some $k\ge 1$) and any unital \hm\,
$\psi: A\to M_R$ for some integer $R\ge N({\cal P}_0)k$
such that
\beq\label{Lfullab-0}
\tau\circ \psi(g)\ge \Delta(\hat{g})\tforal g\in {\cal H},
\eneq
there exists a unital \hm\,
$h_0: A\to M_{R-k}$ such that
\beq\label{Lfullab-1}
[\phi\oplus h_0]|_{{\cal P}_0}=[\psi]|_{{\cal P}_0}.
\eneq
\end{lem}

\begin{proof}
Let $G_0$ be a subgroup of $K_0(A)$ generated by ${\cal P}_0.$   We may also assume, without loss of generality, that
${\cal P}_0=\{[p_1],[p_2],...,[p_{m_1}]\}\cup \{z_1,z_2,...,z_{m_2}\},$
where $p_1, p_2,...,p_{m_1}\in M_l(A)$ are  projections
 and $z_j\in {\rm ker}\rho_A,$ $j=1,2,...,m_2.$
%Let $G_{1,t}={\rm Tor}(K_1(A))\cap G.$ Let $j\ge 1$ be an integer such that
%$jx=0$ for all $x\in G_1.$ Put $J=j!.$

%Suppose that $A\in {\cal D}_m'.$
%has property ($S_m$).
We prove the lemma  by induction.
Assume first
%that $A\in {\cal D}_0.$  We may  write
%{\red{Case (1)}}:
$A=PC(X,F)P{,}$ {{w}}here $X$ is path connected. This, of course, includes the case that $X$ is a single point.
There is $d>0$ such that
\beq\label{Lfullab-5}
\|\pi_{x,j}\circ p_i-\pi_{x',j}\circ p_i\|<1/2,\,\,\, i=1,2,...,m_1,
\eneq
provided that ${\rm dist}(x,x')<d,$ where $\pi_{x,j}$ is identified with $\pi_{x,j}\otimes {\rm id}_{M_l}.$
Since $X$ is compact, we may assume that $\{x_1,x_2,...,x_{m_3}\}$ is a $d/2$-dense set.
Write $P_{x_i}FP_{x_i}=M_{r(i,1)}\oplus M_{r(i,2)}\oplus \cdots \oplus M_{r(i,k(x_i))},$ $i=1,2,...,m_3.$

There are $h_{i,j}\in C(X)$ with $0\le h_{i,j}\le 1,$ $h_{i,j}(x_i)=1_{M_{r(i,j)}}$ and
$h_{i,j}h_{i',j'}=0$ if $(i,j)\not=(i',j').$  Moreover we assume
that $h_{i,j}(x)=0$ if ${\rm dist}(x, x_i)\ge d.$

Put $g_{i,j}=h_{i,j}\cdot P\in A,$ $j=1,2,...,k(x_i), i=1,2,...,m_3,$
Let
\beq\label{Lfullab-6}
\sigma_0=\min\{\Delta(\hat{h_{i,j}}): 1\le j\le k(x_i),\,1\le i\le m_3\}\andeqn N({\cal P}_0)\ge 2/\sigma_0.
\eneq
Put ${\cal H}=\{h_{i,j}: 1\le j\le k(x_i),\, 1\le i \le m_3\}.$

Now suppose that $\phi: A\to M_k$ and $\psi: A\to M_R$ with $R\ge N({\cal P}_0)k$ and
\beq\label{Lfullab-7}
\tau\circ \psi(g)\ge \Delta(\hat{g})\rforal g\in {\cal H}.
\eneq
Write $\phi=\bigoplus_{i,j}^{m_3}\Pi_{y_i,j},$ where $\Pi_{y_i,j(i)}$ is $T_{i,j}$ copies of $\pi_{y_i,j}.$
 Note $k-T_i>0$ for all $i.$
Since $R\ge N({\cal P}_0)k,$ (\ref{Lfullab-7}) implies
that $\psi$ may be viewed as direct sum of at least
\beq\label{Lfullab-8}
\Delta(\widehat{h_{j,i}})\cdot (2k/\sigma_0)>2k
\eneq
copies of $\pi_{x,j}$ with ${\rm dist}(x, x_i)<d,$ $i=1,2,...,m_3.$
Rewrite
$\psi=\Sigma_1\oplus \Sigma_2,$ where
$\Sigma_1$ contains exactly $T_{i,j}$ copies of $\pi_{x,j}$ with ${\rm dist}(x, x_i)<d$ for each $i$ and $j.$
Then
\beq\label{Lfullab-9}
{\rm rank} \, \Sigma_1(p_i)={\rm rank} \phi(p_i),\,\,\,i=1,2,...,m_1.
\eneq
Put $h_0=\Sigma_2.$
Note for any unital \hm\, $h: A\to M_n,$ $[h(z)]=0$ for all $z\in {\rm ker}\rho_A.$
%{\red{This prove Case (i).}}

%{\red{By choosing a large ${\cal H}$ and  consider each component, it is also clear that
%the case $A=PC(X,F)P,$   where $X$ is a disjoint union of finitely many path component,
This proves the case that $A=PC(X,F)P$ as above, in particular,  the case that $A\in {\cal D}_0.$
%the case that $A\in {\cal D}_0$ follows.}}
%has property ($S_0$).

Now assume the lemma holds for any \CA\, $A\in {\cal D}_m.$
%with the property ($S_m$).

Let $A$ be a C*-algebra in ${\cal D}_{m+1}'.$
% having property ($S_{m+1}$).
We assume that $A\subset PC(X,F)P\oplus B$ is a unital \SCA\,
and $I=\{f\in PC(X,F)P: f|_{X^0}=0\},$ where $X^0=X\setminus Y$ and $Y$ is an open subset of $X$
and $B\in {\cal D}_m'$ and $A/I\cong B.$
We assume that, {{if ${\rm dist}(x,x')<2d,$}}
\beq\label{Lfullab-10}
\|\pi_{x, j}(p_i)-\pi_{x',j}(p_i)\|<1/2\andeqn
\|\pi_{x,j}\circ s\circ  \pi_I(p_i)-\pi_{x,j}(p_i)\|<1/2,
\eneq
%provided that ${\rm dist}(x,x')<2d,$
where $s: A/I\to A^{d}=\{f|_{X^d}: f\in A\}$ is an injective \hm\, given by \ref{8-N-3}. We also assume
that $2d<d^Y.$
Define $\Delta_\pi: (A/I)_+^{q, {\bf 1}}\setminus \{0\}\to (0,1)$ by
\beq\label{Lfullab-11}
\Delta_\pi(\hat{g})=\Delta(\widehat{g_0\cdot P\cdot s(g)})\rforal g\in (A/I)_+^{\bf 1}\setminus \{0\},
\eneq
where $g_0\in C(X^d)_+$ with $0\le g\le 1,$ $g_0(x)=1$ if $x\in X^0,$
$g_0(x)>0$ if ${\rm dist}(x, X^0)<d/2$ and
$g_0(x)=0$ if ${\rm dist}(x, X^d)\ge d/2\}.$

Note that $g_0\cdot s(g)>0$ if $g\in (A/I)_+\setminus \{0\}.$ Therefore $\Delta_\pi$ is indeed an order preserving map
from $(A/I)_+^{q, {\bf 1}}\setminus \{0\}$ into $(0,1).$

Note that $A/I\in {\cal D}_m'.$
%has property ($S_m$).
By the inductive assumption,
there is an integer $N_\pi({\cal P}_0)\ge 1,$ a finite subset ${\cal H}_\pi\subset (A/I)_+^{\bf 1}\setminus \{0\}$ satisfying the following:
if $\phi': A/I\to M_{k'}$ is a unital \hm\, and $\psi': A/I\to M_{R'}$ is a unital \hm\, for some $R'>N_\pi({\cal P}_0)k'$ such that
\begin{equation*}
t\circ \psi'(\hat{g})\ge \Delta_\pi(\hat{g})\rforal g\in {\cal H}_\pi,
\end{equation*}
where $t$ is the tracial state of $N_{R'},$ then there exists a unital \hm\, $h_\pi: A/I\to M_{R'-k'}$ such that
\begin{equation*}
(\phi'\oplus h_\pi)_{*0}|_{\bar{\cal P}_0}=(\psi_\pi)_{*0}|_{\bar{\cal P}_0},
\end{equation*}
where $\bar{\cal P}_0=\{(\pi_I)_{*0}(p): p\in {\cal P}_0\}.$

Let $r: Y^{d/2}\to Y^d$ be a homeomorphism such that ${\rm dist}(r(x), x)<d$ for all $x\in Y^{d/2}.$

Let $C=\{f|_{Y^d}: f\in I\}.$
Define $\Delta_I: C_+^{q, {\bf 1}}\setminus \{0\}\to (0,1)$ by
\begin{equation}
\Delta_I(\hat{g})=\Delta(\widehat{f_0\cdot g\circ r})\rforal g\in  C_+^{q,{\bf 1}}\setminus \{0\},
\end{equation}
where $f_0\in C_0(Y)_+$ with $0\le f_0\le 1,$ $f_0(x)=1$ if $x\in Y ^d,$ $f_0(x)=0$ if ${\rm dist}(x, X^0)\le d/2$
and $f_0(x)>0$ if ${\rm dist}(x, X^0)>d/2\}.$

Note that $C=P|_{Y^d} C(Y_d, F)P_{Y^d}.$  By what has been proved, there is an integer $N_I({\cal P}_0)\ge 1$ and a finite subset ${\cal H}_I\subset C_+^{\bf 1}\setminus \{0\}$ satisfying the following:
if $\phi'': C\to M_{k''}$ is a unital \hm\, and $\psi'': C\to M_{R''}$ (for some $R''\ge N_I({\cal P}_0)k''$) is another  unital \hm\, such that
\begin{equation*}
t\circ \psi''(\hat{g})\ge \Delta_I(\hat{g})\tforal g\in {\cal H}_I,
\end{equation*}
where $t$ is the tracial state on $M_{R''},$ then there exists a unital \hm\, $h'': C\to M_{R''-k''}$ such that
\beq\label{Lfullab-14}
(\phi''\oplus h'')_{*0}|_{{\cal P}_0}=(\psi''|_{*0})|_{{\cal P}_0}.
\eneq

Put
\beq\label{Lfullab-13+}
\sigma=\min\{\min\{\Delta_\pi(\hat{g}): g\in {\cal H}_\pi\},\min\{\Delta_I(\hat{g}): g\in {\cal H}_\pi\}\}>0.
\eneq

Let $N=(N_\pi({\cal P}_0)+N_I({\cal P}_0))/\sigma$ and let
\beq\label{Lfullab-15}
{\cal H}=\{g_0\circ s(g): g\in {\cal H}_\pi\}\cup \{f_0\cdot g\circ r\}.
\eneq

Now suppose that $\phi$ and $\psi$ satisfy the assumption for $N=N({\cal P}_0)$ and
${\cal H}$ as above.
We may write $\phi=\Sigma_{\phi, \pi}\oplus \Sigma_{\phi, I},$
where $\Sigma_{\phi, \pi}$ corresponds to a finite direct sum of irreducible representations
of $A$ which factors through $A/I$ and $\Sigma_{\phi, I}$ corresponds to the finite direct sum
of irreducible representations of $I.$
We also write
\beq\label{Lfullab-16}
\psi=\Sigma_{\psi, \pi}\oplus \Sigma_{\psi, b}\oplus \Sigma_{\psi, I'},
\eneq
where $\Sigma_{\psi, \pi}$ corresponds to the finite direct sum of irreducible representations
of $A$ which factors through $A/I,$ $\Sigma_{\psi, b}$ which corresponds to the finite direct sum
of irreducible representations  which factors through point-evaluations at $x\in Y$ with ${\rm dist}(x, X^0)<d/2$
and $\Sigma_{\psi, I'}$ corresponds the rest of irreducible representations (which
can be factors through point-evaluations at $x\in Y$ with ${\rm dist}(x, X^0)\ge d/2$).

Put $q_\pi=(\Sigma_{\psi, \pi}\oplus \Sigma_{\psi,b})(1_A)$ and
$k'={\rm rank} \Sigma_{\phi,\pi}(1_A).$
Define $\psi_\pi: A/I\to M_{{\rm rank}(q_\pi)}$ by
$\psi_\pi(a)=(\Sigma_{\psi,\pi}\oplus \Sigma_{\psi,b})\circ \pi_I\circ s(a)\rforal a\in A/I.$
Then
\beq\label{Lfullab-17-}
t\circ \psi_\pi(g)\ge \tau\circ \psi_\pi(g)\ge \Delta(\widehat{g_0\cdot s(\pi_I(g))})=\Delta_\pi(\hat{g})\tforal g\in {\cal H}_\pi.
\eneq
where $t$ is the tracial state on $M_{{\rm rank}(q_\pi)}.$
Note that
\beq\label{Lfullab-17}
\tau\circ \psi(g_0)\ge \Delta(\hat{g_0})=\Delta_\pi(\widehat{1_{A/I}}),
\eneq
Therefore
\beq\label{Lfullab-18}
{\rm rank} (q_\pi)\ge R\Delta_\pi(1_{A/I})\ge N_\pi({\cal P}_0)k'.
\eneq
By the inductive assumption, there is a unital \hm\, $h_\pi: A/I\to M_{{\rm rank}q_\pi-k'}$ such that
\beq\label{Lfullab-19}
(\Sigma_{\phi,\pi}\oplus h_\pi)_{*0}|_{{\cal P}_0}=(\psi_\pi)_{*0}|_{{\cal P}_0}.
\eneq

Put $q_I=\Sigma_{\psi, I'}(1_A)$ and $k''={\rm rank}(\Sigma_{\phi, I}(1_A)).$
Define $\psi_I: C\to M_{{\rm rank}q_I}$ by
\beq\label{Lfullab-20}
\psi_I(a)=\Sigma_{\psi,I'}(a\circ r)\rforal a\in C.
\eneq
Then
\beq\label{Lfullab-21}
t\circ \psi_I(g)\ge \tau\circ \Sigma_{\psi, I'}(a\circ r)\ge \psi(f_0\cdot a\circ r)\ge
\Delta(\widehat{f_0\cdot g\circ r})=\Delta_I(\hat{g})\rforal g\in {\cal H}_I,
\eneq
where $t$ is the tracial state on $M_{{\rm rank}(q_I)}.$
Note that
\beq\label{Lfullab-22}
\tau\circ \psi(f_0)\ge \Delta(\hat{f_0})=\Delta_I(\hat{1_A}).
\eneq
Therefore
\beq\label{Lfullab-23}
{\rm rank}(q_I)\ge R\Delta_I(1_A)\ge N_I({\cal P}_0)k''.
\eneq
There is $0<d_1<d<d^Y$ such that all irreducible representations appeared in $\Sigma_{\phi,I}$
factor through point-evaluations at $x$ with ${\rm dist}(x, X^d)\ge d_1.$
Put $r': Y^{d_1}\to Y^d.$
Define $\phi_I: C\to M_{k''}$ by $\phi_I(f)=\Sigma_{\phi, I}(f\circ r').$

By the inductive assumption, there is a unital \hm\, $h_I: C\to M_{{\rm rank}(q_I)-k''}$ such that
\beq\label{Lfullab-24}
(\Sigma_{\phi, I}\oplus h_I)_{*0}|_{{\cal P}_0}=(\psi_I)_{*0}|_{{\cal P}_0}.
\eneq
Define $h: A\to M_{R-k}$ by $h(a)=h_\pi(\pi_I(a))\oplus h_I(a|_{Y^d})\oplus \Sigma_{\psi,b}(a) $ for all $a\in A.$
Then, for each $i,$
\beq\label{Lfullab-25}
\hspace{-0.7in}{\rm rank}\phi(p_i)+{\rm rank}h(p_i)&=&{\rm rank}(\Sigma_{\phi, \pi}(p_i))+
{\rm rank}(\Sigma_{\phi, I}(p_i))\\
&+& {\rm rank}(\Sigma_{\psi, b}(p_i))+{\rm rank}h_\pi(p_i)+{\rm rank}h_I(p_i)\\
&=& {\rm rank}\psi_\pi(p_i)+{\rm rank}(\Sigma_{\psi, b}(p_i))+{\rm rank}\psi_I(p_i)\\
&=&{\rm rank} \psi(p_i),\,\,\, i=1,2,...,m_1.
\eneq
Since $(\phi)_{*0}(z_j)=h_{*0}(z_j)=\psi_{*0}(z_j),$ $j=1,2,...,m_2,$ we conclude that
\beq\label{Lfullab-26}
(\phi\oplus h)_{*0}|_{{\cal P}_0}=\psi_{*0}|_{{\cal P}_0}.
\eneq
This completes the induction process.
\end{proof}

\begin{lem}\label{fullabs}
%Let $A=A_0$ or $A=A_0\otimes C(\T),$ where $A_0\in {\cal C}_0$ be a unital separable \CA\,
Let $A\in {\cal D}_s$ be a unital \CA\, and let
$\Delta: A^{q,{\bf 1}}_+\setminus \{0\}\to (0,1)$ be
%a positive
an order preserving map.
For any $\ep>0$ and any finite subset ${\cal F},$ there exist a finite
subset ${\cal H}\subset A_+^{\bf 1}\setminus \{0\}$ and an integer
$L\ge 1$ satisfying the following:
For any unital \hm\, $\phi: A\to M_{k}$ and
any unital \hm\, $\psi: A\to M_{R}$ for some $R\ge Lk$ such that
\beq\label{fullabs-1}
{\rm tr}\circ \psi(h)\ge \Delta(\hat{h})\tforal h\in {\cal H},
\eneq
there exist a unital \hm\, $\phi_0: A\to M_{R-k}$ and a unitary $u\in M_{R}$ such that
\beq\label{fullabs-2}
\|{\rm Ad}\, u\circ {\rm diag}(\phi(f), \phi_0(f))-
\psi(f)\|<\ep
\eneq
for all $f\in {\cal F}.$

\end{lem}

\begin{proof}
Let $\dt>0,$ ${\cal P}\subset \underline{K}(A)$ be a finite subset, ${\cal H}_1\subset A_+^{\bf 1}\setminus \{0\}$
be a finite
subset, ${\cal H}_2\subset A_{s.a.}$ be a finite subset  and $N_0$ be an integer {{as}} required by \ref{UniqAtoM} for $\ep/4$ (in place of $\ep$),
${\cal F},$ $(1/2)\Delta$ and $A.$
Without loss of generality, we may assume that ${\cal H}_2\subset A_+^{\bf 1}\setminus \{0\}.$
%By choosing larger ${\cal H}_1,$ we may assume that
Let $\sigma_0=\min\{\min\{\Delta(\hat{g}): g\in {\cal H}_1\}, \min\{\Delta(\hat{g}): g\in {\cal H}_2\}\}$

Let $G$ be a subgroup of $\underline{K}(A)$ generated by ${\cal P}.$ Put
${\cal P}_0={\cal P}\cap K_0(A).$  We may also assume, without loss of generality, that
${\cal P}_0=\{[p_1],[p_2],...,[p_{m_1}]\}\cup \{z_1,z_2,...,z_{m_2}\},$
where $p_1, p_2,...,p_{m_1}$ are projections in $M_l(A)$
 and $z_j\in {\rm ker}\rho_A,$ $j=1,2,...,m_2.$
Let $j\ge 1$ be an integer such that $K_0(A, \Z/j'\Z)\cap G=\emptyset$ for all $j'\ge j.$
Put $J=j!.$

Let $N({\cal P}_0)\ge 1$ be an integer and ${\cal H}_3\subset A_+^{\bf 1}\setminus \{0\}$ be a  finite subset
{{as}} required by \ref{Lfullab} for ${\cal P}_0.$

Let $p_s=(a^{(s)}_{i,j})_{l\times l},$ $s=1,2,...,m_1$, and
choose $\ep_0>0$ and a finite subset ${\cal F}_0$ such that
\beq\label{fullabs-9}
[\psi']|_{\cal P}=[\psi'']|_{\cal P}{{,}}
\eneq
provided that $\|\psi'(a)-\psi''(a)\|<\ep_0$ for all $a\in {\cal F}_0.$

Put ${\cal F}_2={\cal F}\cup {\cal F}_1\cup {\cal H}_2$ and put $\ep_1=\min\{\ep/16, \ep_1/2\}.$
Let $K>8((N({\cal P}_0)+1)(J+1)/\dt\sigma)$ be an integer.
Let ${\cal H}_0={\cal H}_1\cup {\cal H}_3.$
Let $N_1\ge 1$ (in place of $N$) be an integer, $\dt_1>0$ (in place of $\dt$), ${\cal H}_4\subset
A_+^{\bf 1}\setminus \{0\}$ (in place of ${\cal H}_1$) be a finite subset and ${\cal H}_5\subset A_{s.a.}$
(in place of ${\cal H}_2$) be a finite subset {{as}} required by \ref{8-N-4} for $\ep_1$ (in place of $\ep$)
${\cal F}_2$ (in place of ${\cal F}$),  ${\cal H}_0$ and $K.$
Let $L=K(K+1)$ and let ${\cal H}={\cal H}_4\cup {\cal H}_0$ as well as $\af=15/16.$
Suppose that $\phi$ and $\psi$ {{satisfy}}  the assumption for the above $L$ and ${\cal H}.$

Then, by applying \ref{8-N-4}, there are mutually orthogonal projections
$e_0, e_1, e_2,...,e_K\in M_R$ such that $e_0\lesssim e_1, $ $e_i$ is equivalent to $e_1,$ $i=1,2,...,K,$
a unital \hm\, $\psi_0: A\to e_0M_Re_0$ and a unital \hm\, $\psi_1: A\to e_1M_Re_1$ such that
\beq\label{fullabs-10}
\|\psi(a)-(\psi_0(a)\oplus {\rm diag}(\overbrace{\psi_1(a),\psi_1(a),...,\psi_1(a)}^{K}))\|<\ep_1\tforal a\in {\cal F}_2\\
\andeqn
\tau\circ \psi_1(g)\ge (15/16){\Delta(\hat{g})\over{K}}\rforal g\in {\cal H}_0.
\eneq
Put $\Psi=\psi_0\oplus {\rm diag}(\overbrace{\psi_1(a),\psi_1(a),...,\psi_1(a)}^{K})).$
We compute that
\beq\label{fullabs-10+}
[\Psi]|_{\cal P}=[\psi]|_{\cal P}.
\eneq
Let $R_0={\rm rank}(e_1).$ Then
\beq\label{fullabs-11}
R_0=R\tau\circ \psi_1(1_A)&\ge& Lk(15/16) {\Delta(\widehat{1_A})\over{K}}\ge k(K+1)(15/16)\Delta(\widehat{1_A})\\
&\ge & k(15/16)8N({\cal P}_0)(J+1)/\dt.
\eneq
Moreover,
\beq\label{fullabs-12}
t\circ \psi_1(\hat{g})\ge (15/16)\Delta(\hat{g})\tforal g\in {\cal H}_3,
\eneq
where $t$ is the tracial state on $M_{R_1}.$
It follows from \ref{Lfullab} that there exists a unital \hm\, $h_0: A\to M_{R_0-k}$ such that
\beq\label{fullabs-13}
(\phi\oplus h_0)_{*0}|_{{\cal P}_0}=(\psi_1)_{{\cal P}_0}.
\eneq
Put
\beq\label{fullabs-14}
h_1=h_0\oplus {\rm diag}(\overbrace{\phi\oplus h_0,\phi\oplus h_0,...,\phi\oplus h_0}^{J-1})\andeqn
\psi_2={\rm diag}(\overbrace{\psi_1, \psi_1,...,\psi_1}^J)
\eneq
Then
\beq\label{fullabs-15}
[\phi\oplus h_1]|_{\cal P}=[\psi_2]|_{\cal P}.
\eneq
Put $\Psi'={\rm diag}(\overbrace{\psi_1,\psi_1,...,\psi_1}^{K-J}).$
Let $\phi_0=h_1\oplus \psi_0\oplus \Psi'.$ Then
\beq\label{fullabs-16}
[\phi\oplus \phi_0]|_{\cal P}=[\psi_0\oplus \psi_2\oplus \Psi']|_{\cal P}=[\psi]|_{\cal P}
\eneq
Since $J/(K-J)<\dt$ and by (\ref{fullabs-1}), by applying \ref{UniqAtoM}, there is a unitary $u\in M_R$ such that
\beq\label{fullabs-17}
\|{\rm Ad}\, u\circ (\phi(f)\oplus \phi_0(f))-\psi(f)\|<\ep\tforal f\in {\cal F}.
\eneq
\end{proof}

\section{Almost multiplicative maps to finite dimensional \CA s}

%\begin{df}
%A separable \CA\, $A$ is said to have bounded finite dimensional representations,
%if every irreducible representation $\pi$ of $A$ is finite dimensional and there is
%an integer $d\ge 1$ such that ${\rm dim} \pi(A)\le d.$
%Denote by ${\cal A}_1$ the class of all \CA s with bounded finite dimensional representations.
%\end{df}

%\begin{lem}\label{htoMn}
%Let $A\in {\cal A}_1$ be a unital \CA.
%For any $\ep>0$ and any finite subset ${\cal F}\subset A,$ there exists
%a finite subset ${\cal H}_1\subset A_+$ satisfying the following:
 % For any $\sigma>0,$ there exists a finite subset ${\cal P}\subset \underline{K}(A),$ a finite subset ${\cal H}_2\subset A_{s.a.}$ and $\dt>0$ satisfying the following:

%For any integer $n\ge 1$ and any two unital \hm s $\phi, \psi:
%A\to M_n,$ if
%\beq\label{htoMn-1}
%[\phi]|_{\cal P}&=&[\psi]|_{\cal P}\\
%|\tau\circ \phi(a)-\tau\circ \psi(a)|&<&\dt\tforal a\in {\cal H}_2\\
%\tau\circ \phi(b)\ge \sigma&&\andeqn
%\tau\circ \psi(b)\ge \sigma
%\eneq
%for all $b\in {\cal H}_1.$
%Then, there exists a unitary
%$u\in M_n$ such that
%\beq\label{htoMn-2}
%\|{\rm Ad}\, u\circ \phi(f)-\psi(f)\|<\ep\tforal f\in {\cal F}.
%\eneq

%\end{lem}

Note that in the following statement,  $n$ is given and $(\mathcal G, \delta)$ depends on $n$.

\begin{lem}\label{LtoMN}
Let $n\ge 1$ be an integer and let $A$ be a unital separable  \CA.
For any $\ep>0$ and any finite subset ${\cal F}\subset A,$ there exist
$\dt>0$ and a finite subset ${\cal G}\subset A$ such that, for any unital ${\cal G}$-$\dt$-multiplicative \morp\, $\phi: A\to M_n,$ there exists a
unital \hm\, $\psi: A\to M_n$ such that
$$
\|\phi(a)-\psi(a)\|<\ep\rforal a\in {\cal F}.
$$
\end{lem}

\begin{proof}

Suppose that the lemma is not true for  certain finite set ${\cal F} \sbs A$ and $\ep_0{{>0}}$. Let $\{{\cal G}_k\}_{k=1}^{\infty}$ be a sequence of finite subsets of $A$ with ${\cal G}_k\subset  {\cal G}_{k+1}$ and $\overline{\cup_k {\cal G}_k}=A$
and let $\{\dt_k\}$ be a monotone decreasing sequence of positive numbers with $\dt_k \rightarrow 0$.
Since the lemma is assumed not to be true, there are unital ${\cal G}_k$-$\dt_k$-multiplicative \morp s $\phi_k: A\to M_n$ such that
\beq\label{150103-sec5-1}
 \mbox{inf}\{\mbox{max}_{a\in {\cal F}}\|\phi_k(a)-\psi(a)\|:~~~\psi: A\to M_n ~\mbox{homomorphisms} \} \geq \ep_0.
\eneq
For each pair $(i,j)$ with $1\leq ~i,j~\leq n$, let $l^{i,j}: M_n \to \C$ be the map defined by taking matrix $a\in M_n$ to the entry of $i^{th}$ row and $j^{th}$ column of $a$. Let $\phi_k^{i,j} =l^{i,j}\circ \phi_k: A\to \C$. Note that the unit ball of the dual space of $A$ (as Banach space) is weak$*$-compact.  Since $A$ is separable, there is a subsequence (instead of subnet) of $\{\phi_k\}$ (still denote by $\phi_k$) such that $\{\phi_k^{i,j}\}$ is weak$*$ convergent for all $i,j$. In other words,  $\{\phi_k\}$ convergent pointwise. Let $\psi_0$ be the the limit. Then $\psi_0$ is a homomorphism and for $k$ large enough, we have
$$\|\phi_k(a)-\psi_0(a)\|< \ep_0, ~~\rforal a\in F.$$
This is a contradiction to (\ref{150103-sec5-1}) above.
\end{proof}

%{\bf I have a proof and can be easily copied to here-----do not know if it will be absolutely necessary here. Answer: I add the proof for the general case of $A$. I think it is necessary.}

\begin{lem} {\rm (cf. Lemma 4.5 of \cite{Lncrell})}\label{0dig}
Let $A$ be a unital \CA\, arising from a locally trivial continuous field of $M_n$ over a compact metric space $X.$ Let $T$ be a finite subset of tracial states on $A.$ For any finite subset ${\cal F}\subset A$ and for any $\ep>0$ and any $\sigma>0,$ there is an ideal $J\subset A$ such that
$\|\tau|_J\|<\sigma$ for all $\tau\in T,$ a finite dimensional \SCA\, $C\subset A/J$ and a unital \hm\, $\pi_0$
{{from $A/J$}}
%with domain $A$
such that
\beq\label{0dig-1}
{\rm dist}(\pi(x), C)<\ep\tforal x\in {\cal F} {{\andeqn}}
%\\\label{0dig-2}
\pi_0(A/J)=\pi_0(C)\cong C,
%\andeqn {\rm ker}\pi_0\supset J,
\eneq
where $\pi: A\to A/J$ is the quotient map.
\end{lem}
\begin{proof}
This  follows from
Lemma 4.5 of \cite{Lncrell}.
Choose $x_i\in F_i$
%The only  item In that proof, one can choose
%$F_1'=\cup_{F_j\in B(\xi_1, \dt_i)} F_j,$ $F_i'=\cup_{F_j\in B(\xi_i, \dt_i)\setminus \cup_{l\le i-1} F_l'}F_l' $ and use $F_i'$ instead of $F_j,$
%
% only the second part of (\ref{0dig-1})
%{\color{Green} (0dig-2 is changed to 0dig-2)}
%is required to be proved.
%However, this is also included in the proof of Lemma 4.5
%of \cite{Lncrell}. Note that the point $\xi_j\in F_j,$
$j=1,2,...,k.$ One can then choose
$\pi_0=\bigoplus_{j=1}^k \pi_{x_j}.$
\end{proof}

\begin{lem} {\rm (cf. Lemma 4.7 of \cite{Lncrell})}\label{abdig}
Let $A$ be a unital  separable subhomogeneous C*-algebra. % \CA\, whose irreducible representations have bounded dimensions.
Let $T\subset T(A)$ be a finite subset. For any finite subset ${\cal F}\subset A,$ $\ep>0$ and $\sigma>0,$ there is an ideal $J\subset A$ such that
$\|\tau|_J\|<\sigma$ for all $\tau\in T,$ a finite dimensional \SCA\, $C\subset A/J$ and  a unital \hm\, $\pi_0$
{{from $A/J$}}
%with domain $A$
such that
\beq\label{abdig-1}
{\rm dist}(\pi(x), C)<\ep\tforal x\in {\cal F} {{\andeqn}}
%\\\label{abdig-2}
\pi_0(A/J)=\pi_0(C)\cong C,
%\andeqn {\rm ker}\pi_0\supset J,
\eneq
where $\pi: A\to A/J$ is the quotient map.
\end{lem}

\begin{proof}
The proof is in fact contained in that of Lemma 4.7 of \cite{Lncrell}. Each time Lemma 4.5 of \cite{Lncrell} applied, one can apply
 \ref{0dig} instead. The complete proof is omitted here.
\end{proof}

%The independent of $\pi_0$ from $T$ and $\sigma$ is not used in the next proof.

\begin{lem}\label{LtoMn}
Let $A$ be a unital subhomogeneous C*-algebra. %\CA\, whose irreducible representations have bounded dimensions.
Let $\ep>0,$ let ${\cal F}\subset A$ be a finite subset and let $\sigma_0>0.$
There exist
%a finite subset ${\cal H}$ satisfies the following:
 %For any $\sigma>0,$ there exists
 $\dt>0$ and a finite subset ${\cal G}\subset A$ satisfying the following:
Suppose that $\phi: A\to M_n$ (for some integer $n\ge 1$) is a
${\cal G}$-$\dt$-multiplicative \morp.
%such that
%\beq\label{LtoMn-1}
%\tau\circ \phi(b)\ge \sigma
%\eneq
%for all $b\in {\cal H}.$
Then,
there exists a projection $p\in M_n$  and a unital \hm\,
$\phi_0: A\to pM_np$ such that
\beq\label{LtoMn-2}
\|p\phi(a)-\phi(a)p\|&<&\ep\tforal a\in {\cal F},\\
\|\phi(a)-[(1-p)\phi(a)(1-p)+\phi_0(a)]\|&<&\ep\tforal  a\in {\cal F}\tand\\
tr(1-p)&<&\sigma_0,
\eneq
where $tr$ is the normalized trace on $M_n.$
\end{lem}

\begin{proof}
We assume that the lemma is false.
Then there exists $\ep_0>0,$ a finite subset ${\cal F}_0,$ a positive number $\sigma_0>0,$
an increasing sequence
of finite subsets ${\cal G}_n\subset A$ such that ${\cal G}_n\subset  {\cal G}_{n+1}$ and such that
$\cup_{n=1}{\cal G}_n$ is dense in $A,$ a sequence of decreasing positive
numbers $\{\dt_n\}$ with $\sum_{n=1}^{\infty}\dt_n<\infty,$
a sequence of integers $\{m(n)\}$ and a sequence
of unital ${\cal G}_n$-$\dt_n$-multiplicative \morp s $\phi_n: A\to M_{m(n)}$ satisfying the following:
\beq\label{LtoMn-6}
\inf \{\max\{\|\phi_n(a)-[(1-p)\phi_n(a)(1-p)+\phi_0(a)\|: a\in {\cal F}_0\}\}\ge \ep_0
\eneq
where infimum is take among all projections $p\in M_{m(n)}$
with $tr_n(1-p)<\sigma_0,$ where $tr_n$ is the normalized
trace on $M_{m(n)}$ and all possible \hm s $\phi_0: A\to pM_{m(n)}p.$
By the virtue of \ref{LtoMN},
one may also assume that $m(n)\to \infty$ as $n\to\infty.$

%Let $tr_n$ be the normalized trace on $M_{m(n)}.$
Note  that  $\{tr_n\circ \phi_n\}$ is a sequence of (not necessary tracial) states of $A.$ Let
$t_0$ be a weak limit of $\{tr_n\circ \phi_n\}.$ Since $A$ is separable, there is a subsequence (instead of subnet) of
 $\{tr_n\circ \phi_n\}$ converging to $t_0$.
%By passing to a subsequence, to simplify the notation,
Without loss of generality, we may assume that
$tr_n\circ \phi_n$ converges to $t_0.$ By the {{${\cal G}_n$-$\dt_n$}}-multiplicativity of $\phi_n$, we know that $t_0$ is a tracial state on $A$.

Consider the ideal $\bigoplus_{n=1}^{\infty}(\{M_{m(n)}\}),$ where
$$
\bigoplus_{n=1}^{\infty}(\{M_{m(n)}\})=\{\{a_n\}: a_n\in M_{m(n)}\andeqn \lim_{n\to \infty}\|a_n\|=0\}.
$$
Denote by $Q$ the quotient
$\prod_{n=1}^{\infty}(\{M_{m(n)}\}/\bigoplus_{n=1}^{\infty}(\{M_{m(n)}\}).$
Let
$\pi_\omega: \prod_{n=1}^{\infty}(\{M_{m(n)}\})\to  Q$ be the quotient map.
%\prod_{n=1}^{\infty}(M_{m(n)})/\bigoplus_{n=1}^{\infty}(\{ M_{m(n)}\})$ be the quotient map.
Let $A_0=\{\pi_{\omega}(\{\phi_n(f)\}): f\in A\}$ which is a  subalgebra of
$Q.$
%$l^{\infty}(M_{m(n)})/c_0(\{ M_{m(n)}\})$. Denote $\Psi=\pi_\omega\circ \{\phi_n\}.$
Then $\Psi$ is a unital \hm~ from $A$ to $\prod_{n=1}^{\infty}(M_{m(n)})/\bigoplus(\{ M_{m(n)}\})$ with $\Psi(A)=\pi_\omega(A_0)$. If $a\in A$ has  zero image in $\pi_\omega(A_0)$, that is, $\phi_n(a) \to 0$, then $t_0(a)=\lim_{n\to \infty} tr_n(\phi_n(a))=0$.
  So we  may view $t_0$ as a state on $\pi_\omega(A_0)=\Psi(A).$
%Note that $t_0$ vanishes on $c_0(\{M_{m(n)}\}).$
%Let $\gamma_n: l^{\infty}(\{M_{m(n)}\}\to M_{m(n)}$ be the projection
%on the $n$-coordinate. Define $\tau_n=tr_n\circ \gamma_n.$ Then
%each $\tau_n$ is a tracial state of $l^{\infty}(\{M_{m(n)}\}).$ Let
%$\tau_0$ be a weak limit of $\{tau_n\}.$ One verifies that $\tau_0$ is a tracial state of $l^{\infty}(\{M_{m(n)}\}).$ By passing another subsequence, to simplify the notation, without loss of generality, one may further assume that
%$$
%\lim_{n\to\infty}\tau_n(x)=\tau_0(x)\tforal x\in l^{\infty}(\{M_{m(n)}\}).
%$$

%Let $\omega\in \bt N\setminus \N$ be a free ultrafilter.
%It should be noted that $\tau_0\circ \Psi$ may be identified with
%$t_0$ and it is a tracial state on $\Psi(A).$

It follows from Lemma \ref{abdig} that there is an
% (two-sided closed) 
ideal $I\subset \Psi(A)$ and a finite dimensional \SCA\, $B\subset \Psi(A)/I$
and a unital \hm\, $\pi_{00}:\Psi(A)/I \to B$ such that
\beq\label{LtoMn-10}
{\rm dist}(\pi_I\circ \Psi(f), B)&<&\ep_0/16\tforal f\in {\cal F}_0,\\\label{LtoMn-11}
\|(t_0)|_I\|&<&\sigma_0/2\andeqn
%\\\label{LtoMn-11+}
\pi_{00}|_B=\mbox{id}_B.
\eneq
Note that $\pi_{00}$ can be regarded as map from $A$ to $B$ {{with}} ${\rm ker} \pi_{00}\supset I$.
There is, for each $f\in {\cal F}_0,$ an element $b_f\in B$ such that
\beq\label{LtoMn-12}
\|\pi_I\circ \Psi(f)-b_f\|<\ep_0/16.
\eneq
Put $C'=B+I$ and $I_0=\Psi^{-1}(I)$ and $C_1=\Psi^{-1}(C').$ For each $f\in {\cal F}_0,$ there exists
$a_f\in C_1\subset A$ such that
\beq\label{LtoMn13}
\|f-a_f\|<\ep_0/16\andeqn \pi_I\circ \Psi(a_f)=b_f.
\eneq
Let $a\in (I_0)_+$ be a strictly positive element and let
$J=\overline{\Psi(a)Q\Psi(a)}$ be the hereditary \SCA\, of $Q$ generated by $\Psi(a).$
%By (\ref{LtoMn-11}),
%\beq\label{LtoMn-14}
%\|(\tau_0)|_J\|<\sigma_0/2.
%\eneq
Put $C_2=\Psi(C_1)+J.$ Then $J$ is an ideal of $C_2.$
Denote by $\pi_J: C_2\to B$ the quotient map.
Since $Q$ and $J$ have real rank zero and
$C_2/J$ has finite {{dimension}}, by Lemma 5.2  of \cite{LnTAF}, $C_2$ has real rank zero. It follows that
$$
0\to J\to C_2\to B\to 0
$$
is a quasidiagonal extension. As in Lemma 4.9 of \cite{Lncrell}, there is a projection
$P\in J$ and a unital \hm\, $\psi_0: B\to (1-P)C_2(1-P)$ such that
\beq\label{LtoMn-15}
\|P\Psi(a_f)-\Psi(a_f)P\|<\ep_0/8\andeqn \|\Psi(a_f)-[P\Psi(a_f)P+\psi_0\circ \pi_J\circ\Psi(a_f)]\|<\ep_0/8
\eneq
for all $f\in {\cal F}_0.$ Let $H: A\to \psi_0(B)$ be defined by
$H=\psi_0\circ \pi_{00}\circ \pi_I\circ \Psi$.
%Since ${\rm ker}\pi_{00}\supset I$ and $\pi_{00}(\Psi(A))=\pi_{00}(B)\cong B,$
%there is a unital \hm\, $H: A\to \psi_0(B)$ such that $H|_B=\psi_0.$
One estimates that
\beq\label{LtoMn-16}
\|P\Psi(f)-\Psi(f)P\|<\ep_0/2\andeqn\\\label{LtoMn-17}
\|\Psi(f)-[P\Psi(f)P+H(f)]\|<\ep_0/2
\eneq
for all $f\in {\cal F}_0.$ Note that ${\rm dim}H(A)<\infty,$ and that $H(A)\sbs Q.$
%=l^{\infty}(\{M_{m(n)}\}/c_0(\{M_{m(n)}\}).$
There is a
\hm\, $H_1: H(A)\to \prod_{n=1}^{\infty}(\{M_{m(n)}\})$ such that
$\pi\circ H_1\circ H=H.$ One may write $H_1=\{h_n\},$ where
each $h_n: H(A)\to M_{m(n)}$ is a (not necessary unital) \hm, $n=1,2,....$
There is also a sequence of projections $q_n\in M_{m(n)}$ such that
$\pi(\{q_n\})=P.$ Let $p_n=1-q_n,$ $n=1,2,....$ Then, for sufficiently large $n,$
by (\ref{LtoMn-16}) and (\ref{LtoMn-17}),
\beq\label{LtoMn-18}
\|(1-p_n)\phi_n(f)-\phi_n(f)(1-p_n)\|<\ep_0\andeqn\\
\|\phi_n(f)-[(1-p_n)\phi_n(f)(1-p_n)+h_n\circ H(f)]\|<\ep_0
\eneq
for all $f\in {\cal F}_0.$ Moreover, since $P\in J,$ for any $\eta>0,$ there is $b\in I_0$ with
$0\le b\le 1$ such that
\begin{equation*}
\|\Psi(b)P-P\|<\eta.
\end{equation*}
However, by (\ref{LtoMn-11}),
\beq\label{LtoMn-20}
0<t_0(\Psi(b))<\sigma_0/2\tforal b\in I_0\,\,\, {\rm with}\,\,\, 0\le n\le 1.
\eneq
By choosing sufficiently small $\eta,$
for all sufficiently large $n,$
\begin{equation*}
tr_n(1-p_n)<\sigma_0.
\end{equation*}
This contradicts with (\ref{LtoMn-6}).
\end{proof}

%\begin{lem}\label{numbercomp}
%%Let $z_1, z_2,\cdots, z_n$ are positive integers. There is a positive integers $T$ depending on $(z_1, z_2,\cdots, z_n)$ such that for any two nonnegative integer linear combination  $a=\sum_{i=1}^n a_i\cdot z_i$ and $b=\sum_{i=1}^n b_i\cdot z_i$, there are two combinations $a'=\sum_{i=1}^n a'_i\cdot z_i$ and $b'=\sum_{i=1}^n b'_i\cdot z_i$ with $a'=b'$, $0\leq a'_i\leq a_i$, $0\leq b'_i\leq b_i$, and ${\rm min}\{a-a', b-b'\}\leq T$.

%\end{lem}

%\begin{proof}
%
%Let $T=n\cdot{\rm max}_{i,j}\{z_iz_j\}$. It is enough to prove that if both $a>T$ and $b>T$, then there are nonzero $0<a'=\sum_{i=1}^n a'_i\cdot z_i=b'=\sum_{i=1}^n b'_i\cdot z_i$ with $0\leq a'_i\leq a_i$, $0\leq b'_i\leq b_i$. But if both $a>T$ and $b>T$, then there are two (not necessary distinct) index $i, j$, with $a_i\ge z_j$ and $b_j\ge z_i$. And therefore we can define $a'=a'_iz_i=b'=b'_jz_j$, by chose $a'_i=z_j$ and $b'_j=z_i$.
%
%\end{proof}

\begin{cor}\label{LtoMncor}
Let $A$ be a unital subhomogeneous C*-algebra.  %\CA\, whose irreducible representations have bounded dimensions.
Let $\eta>0,$ let ${\cal E}\subset A$ be a finite subset and let $\eta_0>0.$
There exist
 $\dt>0$ and a finite subset ${\cal G}\subset A$ satisfying the following:
Suppose that $\phi,~ \psi: A\to M_n$ (for some integer $n\ge 1$) are two
${\cal G}$-$\dt$-multiplicative \morp s.
Then,
there exist  projections $p, q\in M_n$ with ${\rm rank}(p)={\rm rank}(q)$ and  unital \hm s\,
$\phi_0: A\to pM_np$ and $\psi_0: A\to qM_nq$ such that
$$
\|p\phi(a)-\phi(a)p\|<\eta,\quad \|q\psi(a)-\psi(a)q\|<\eta,\quad a\in {\cal E},$$
$$\|\phi(a)-[(1-p)\phi(a)(1-p)+\phi_0(a)]\|<\eta,\quad \|\psi(a)-[(1-q)\psi(a)(1-q)+\psi_0(a)]\|<\eta,\quad  a\in {\cal E}$$
$$\andeqn tr(1-p)=tr(1-q)<\eta_0,$$
where $tr$ is the normalized trace on $M_n.$
\end{cor}

For convenience of future use, we used $\eta$, $\eta_0$ and ${\cal E}$ to replace $\ep$, $\sg_0$ and ${\cal F}$ in \ref{LtoMn}.

\begin{proof}

By \ref{LtoMn}, we can get such decomposition for $\phi$ and $\psi$ separately, then the only missing part  is that ${\rm rank}(p)={\rm rank}(q)$.
Let $\{z_1,z_2,...,z_m\}$ be the set of ranks of irreducible representations of $A$ and let
$T$ be the number given by \ref{8-N-0} corresponding to  $\{z_1,z_2,...,z_m\}.$
%) (in the case $A=A(F_1,F_2,\phi_0,\phi_1)$ or  $A=A(F_1,F_2,\phi_0,\phi_1)\otimes C(S^1)$ the set is $\{k_1,k_2,\cdots,k_p,l_1,l_2,\cdots, l_l\}$).
We apply \ref{LtoMn} to $\eta_0/2$ instead of $\sg_0$ (and, $\eta$ and ${\cal E}$ in places of  $\ep$ and ${\cal F}$). By \ref{LtoMN}, we can assume the size $n$ of matrix $M_n$ is so large that $\frac{T}n<\eta_0/2$. By \ref{8-N-0}, we can take sub-representations out of $\phi_0$ and $\psi_0$ (one of them has size at most $T$) so that the remainder of $\phi_0$ and $\psi_0$ have same size---that is for ${\rm rank}({\rm new}~ p)={\rm rank}({\rm new}~ q)$, and $\mathrm{tr}(1-({\rm new}~p))=\mathrm{tr}(1-({\rm new}~q))<\frac{\eta_0}2+\frac{T}n<\eta_0$.
\end{proof}

\begin{lem}\label{Combinerep}
%Let $A=A_0$ or $A=A_0\otimes C(\T),$ where $A_0\in {\cal C}_0,$ be a unital \CA.
Let $A\in \overline{ {\cal D}}_s$ be an infinite dimensional unital \CA, let $\ep>0$   and let
${\cal F}\subset A$ be  a finite subset.
 let $\ep_0>0$ and let ${\cal G}_0\subset A$ be a finite subset.
Let $\Delta: A_+^{q, {\bf 1}}\setminus\{0\}\to (0,1)$ be
%a positive
{{an order preserving}} map.

Suppose that
${\cal H}_1\subset A_+^{\bf 1}\setminus \{0\}$ is a finite subset,
$\ep_1>0$ is a
positive number and
$K\ge 1$ is an integer.
There exists $\dt>0,$  $\sg>0$ and  a finite subset
${\cal G}\subset A$ and a finite subset ${\cal H}_2\subset A_+^{\bf 1}\setminus\{0\}$ satisfying the following:
Suppose that $L_1, L_2: A\to M_n$ (for some integer $n\ge 1$) are unital ${\cal G}$-$\dt$-multiplicative \morp s
\beq\label{Comb-1}
tr\circ L_1(h)\ge \Delta({\hat h})~\mbox{and}~tr\circ L_2(h)\ge \Delta({\hat h})~~\tforal h\in {\cal H}_2,\ and
\eneq
\beq\label{Comb-2}
|tr\circ L_1(h)-tr\circ L_2(h)|<\sg \rforal ~h\in {\cal H}_2.
\eneq
Then there exist mutually orthogonal  projections $e_0, e_1,e_2,...,e_K\in
M_n$ such that
$e_1, e_2,...,e_K$ are pairwise equivalent, $e_0\lesssim e_1,$
$tr(e_0)<\ep_1$ and $e_0+\sum_{i=1}^K e_i=1,$ and there exist  unital
${\cal G}_0$-$\ep_0$-multiplicative \morp s
 $\psi_1,\psi_2: A\to e_0M_ke_0$,  a unital \hm\, $\psi: A\to e_1M_ke_1$, and unitary $u\in M_n$ such that
\beq\label{Com-3}
\|L_1(f)-{\rm diag}(\psi_1(f), \overbrace{\psi(f),\psi(f),...,\psi(f)}^K)\|<\ep \tand
\eneq
\beq\label{Com-4}
\|uL_2(f)u^*-{\rm diag}(\psi_2(f), \overbrace{\psi(f),\psi(f),...,\psi(f)}^K)\|<\ep
\eneq
for all $f\in {\cal F},$ where $tr$ is the tracial state of $M_n.$
Moreover,
\beq\label{Com-5}
tr(\psi(g))\ge {\Delta({\hat g})\over{3K}}\tforal g\in {\cal H}_1.
\eneq
\end{lem}

\begin{proof}
First note that the following statement is evident. For any $C^*$-algebra $A$, a finite subset ${\cal G}_0\subset A$  and  $\ep_0>0$, there are  finite subset ${\cal G}'\subset A$ with {{$\delta'>0,$   ${\cal F}'\subset A$  and}} with $\ep'>0$ satisfying the following condition.  If $L: A\to B$ is a unital ${\cal G}'$-$\dt'$-multiplicative \morp\,,  $p_0,p_1\in B$ are projections with $p_0+p_1= 1_B$, and $L'_0: A\to p_0Bp_0$, $L'_1: A\to p_1Bp_1$ are linear maps with
$$\|L(f)-{\rm diag}(L'_0(f),L'_1(f))\|< \ep'\rforal f\in {\cal F}',$$
then both $L'_0$ and $L'_1$ are $\ep_0$-${\cal G}_0$-multiplicative. Therefore if $\ep$ is sufficiently small and ${\cal F}$ is sufficiently large relative to $(\ep_0, {\cal G}_0)$, {{t}}hen (\ref{Com-3}) and (\ref{Com-3}) imply $\psi_1$ and $\psi_2$
are {{${\cal G}_0$-$\ep_0$}}-multiplicative.

Put
\beq\label{8-Nsec4-1-1}
\ep_1=\min\{\ep/16, \ep'/16, \frac{1}{64K}{\rm min}\{\Delta ({\hat h}), h\in {\cal H}_1\}\}.
\eneq
Let $\Delta_1=(3/4)\Delta.$  Let $\dt_1>0$ (in place of $\dt$), ${\cal H}_{1,0}\subset A_+^{\bf 1}\setminus \{0\}$ (in place of ${\cal H}_1$) be a finite subset, ${\cal H}_{2,0}\subset A_{s.a.}$ (in place of ${\cal H}_2$) be a finite  subset {{as}}
required by \ref{8-N-4} (see also its remark \ref{8-N-4-r}) for $\ep_1$ (in place of $\ep$),  ${\cal F}\cup {\cal F}'$
{{(in place of ${\cal F}$),}}
$2K$ (in place of $K$), $\Delta_1$ and $A$ as well as $\af=(3/4).$  We may assume that
${\cal H}_{0,2}\subset A_+^{\bf 1}\setminus \{0\}.$

Let $\eta_0=\min\{\dt_1/16, \ep_1/16, \min\{\tau(h): h\in {\cal H}_{1,0}\cup {\cal H}_{2,0}\}.$
{{Let}} $\dt_2>0$ (in place of $\dt$)  and let  ${\cal G}_1\subset A$ (in place of ${\cal G}$) be the finite subset {{as}} required by \ref{LtoMncor} for
$\eta=\eta_0\cdot \min\{\ep_1,\dt_1/4\},$  $\eta_0$ and  ${\cal E}={\cal F}\cup {\cal H}_{1,0}\cup {\cal H}_{2,0}\cup {\cal H}_1.$
Let  $\dt=\eta_0\cdot \min\{\dt_2/2, \dt_1/2,\},$ $\sigma=\min\{\eta_0/2, \eta_1/2\},$
let ${\cal G}={\cal G}'\cup {\cal G}_1\cup {\cal F}\cup {\cal F}'\cup {\cal E}$ and
let ${\cal H}_2={\cal H}_{1,0}\cup {\cal H}_{2,0}\cup {\cal H}_1.$

Now suppose that $L_1$ and $L_2$ satisfy the assumption of the lemma with respect to the above $\dt,$ $\sigma$ and
${\cal G}$  and ${\cal H}_2.$
It follows from \ref{LtoMncor} that there exists a projection $p\in M_n,$
two unital \hm s $\phi_1, \phi_2: A\to pM_np$
and a unitary $u_1\in M_n$ such that
\beq\label{8-Nsec4-1-2}
\|{\rm Ad}\, u_1\circ L_1(a)-((1-p)u_1^*L_1(a)u_1(1-p)+\phi_1(a))\|<\eta,\\
\|L_2(a)-((1-p)L_2(a)(1-p)+\phi_2(a)\|<\eta
\eneq
for all $a\in {\cal E}$ and
\beq\label{8-Nsec4-1-3}
\tau(1-p)<\eta_0,
\eneq
where $\tau$ is the tracial state on $M_n.$

We compute that
\beq\label{8-Nsec4-1-4}
\tau\circ \phi_1(g)\ge \Delta(\hat{g})-\eta-\eta_0\ge (3/4)\Delta(\hat{g}) \tforal g\in {\cal H}_{1,0}\andeqn\\
|\tau\circ \phi_1(g)-\tau\circ \phi_2(g)|<2\eta+2\eta_0+\dt<\eta_0 \dt_1\tforal g\in {{{\cal H}_{2,0}.}}
\eneq

It follows from \ref{8-N-4} (and its remark (\ref{8-N-4-r})) that  there {{exist}}
mutually orthogonal projections $q_0, q_1,...,q_{2K}\in pM_np$ such that
$q_0\lesssim q_1$ and $q_i$ is equivalent to $q_1$ for all $i=1,2,...,2K,$ two unital \hm s
$\phi_{1,0}, \phi_{2,0}: A\to q_0M_nq_0,$ a unital \hm\, $\psi': A\to e_1M_ne_1$ and a unitary
$u_2\in pM_np$ such that
\beq\label{8-Nsec4-1-5}
\|{\rm Ad}\, u_2\circ \phi_1(a)-(\phi_{1,0}\oplus {\rm diag}(\overbrace{\psi'(a),\psi'(a), ...,\psi'(a)}^{2K}))\|<\ep_1\\
\andeqn\|\phi_2(a)-(\phi_{2,0}\oplus {\rm diag}(\overbrace{\psi'(a),\psi'(a), ...,\psi'(a)}^{2K}))\|<\ep_1
\eneq
 for all $a\in {\cal F}\cup {\cal F}'.$ Moreover,
 \beq\label{8-Nsec4-1-6}
 \tau\circ\psi(g)\ge (3/4)^2{\Delta(\hat{g}){2K}}\rforal g\in {\cal H}_1.
 \eneq
Let $u=((1-p)+u_2)u_1,$ $e_0=(1-p)\oplus q_0,$
$e_i=q_{2i-1}\oplus q_{2i},$ $i=1,2,...,K,$ let $\psi_1(a)=(1-p)u_1^*L_1(a)u_1(1-p)\oplus \psi_{1,0},$
$\psi_2(a)=(1-p)L_2(a)(1-p)\oplus \psi_{2,0}$ and $\psi(a)={\rm diag}(\psi'(a)\oplus \psi'(a))$ for $a\in A.$
{{Then}}
\beq\label{8-Nsec4-1-7}
\|{\rm Ad}\, u\circ L_1(f)-(\psi_1(f)\oplus {\rm diag}(\overbrace{\psi(f), \psi(f),...,\psi(f)}^K))\|<\ep\\
\andeqn \|L_2(f)-(\psi_2(f)\oplus {\rm diag}(\overbrace{\psi(f), \psi(f),...,\psi(f)}^K))\|<\ep
\eneq
for all $f\in {\cal F},$
\beq\label{8-Nsec4-1-8}
\tau\circ \psi(g)\ge {\Delta(\hat{g})\over{3K}}\rforal g\in {\cal H}_1.
\eneq
Moreover $\psi_1$ and $\psi_1$ are ${\cal G}_0$-$\ep_0$-multiplicative.
\end{proof}

%{\bf These may not be the best way to state them, one can change the proofs late to fit different way to state them. As this moment, I thought that I should state this way. Answer: Guihua changed the statement to get simultaneous decomposition. }

%\begin{df}\label{Dfull}
%Let $A$ be a unital separable \CA. Let $T_1: A_+\setminus \{0\}\to R_+\setminus \{0\}$ and $T_2: A_+\setminus \{0\}\to \N$ be two maps.
%Let $B$ be a unital \CA\, and $\phi: A\to B$ be a map.
%We say $\phi$ is $T_1\times T_2$-full if, for each $a\in A_+\setminus \{0\},$
%there exist $x_1, x_2,...,x_{T_2(a)}\in A$ with $\|x_i\|\le T_1(a),$
%$i=1,2,..., T_2(a)$ such that
%\beq\label{Dfull-1}
%\sum_{i=1}^{T_2(a)} x_i^*\phi(a)x_i=1_B.
%\eneq
%\end{df}
\begin{cor}\label{repcor}
Let  $A_0\in {\overline {\cal D}}_s$ be a unital \CA,
%$A=A_0\otimes C(\T),$ where $A_0\in {\cal C}_0,$ be a unital \CA,
let $\ep>0$ and
let ${\cal F}\subset A$ be a finite subset.
Let $\Delta: (A_0)_+^{q, {\bf 1}}\setminus \{0\}\to (0,1)$ be an order preserving  map.

Suppose that ${\cal H}_1\subset (A_0)_+^{\bf 1}\setminus \{0\}$ is a finite subset,  $\sigma>0$
is  positive number and $n\ge 1$ is an integer.
There exists a finite subset ${\cal H}_2\subset (A_0)_+^{\bf 1}\setminus\{0\}$ satisfying the following:
Suppose that $\phi: A=A_0\otimes C(\T)\to M_k$ (for some integer $k\ge 1$) is a unital \hm\, and
\beq\label{repcor-1}
tr\circ \phi(h\otimes 1)\ge \Delta(\hat{h})\tforal h\in {\cal H}_2.
\eneq
Then there exist mutually orthogonal  projections $e_0,e_1,e_2,...,e_n\in
M_k$ such that
$e_1, e_2,...,e_n$ are equivalent and $\sum_{i=0}^n e_i=1,$ and there {{exists}}
%a
unital \hm s $\psi_0: A= A_0\otimes C(\T)\to e_0M_ke_0$ and $\psi: A=A_0\otimes C(\T)\to e_1M_ke_1$ such that
one may write that
\beq\label{repcor-2}
&&\|\phi(f)-{\rm diag}( \psi_0(f),\overbrace{\psi(f),\psi(f),...,\psi(f)}^n)\|<\ep\\
&&\tand {\rm tr}(e_0)<\sigma
\eneq
for all $f\in {\cal F},$ where $tr$ is the tracial state on $M_k.$
Moreover,
\beq\label{repcor-3}
tr(\psi(g\otimes 1))\ge {\Delta(\hat{g})\over{2n}}\tforal g\in {\cal H}_1.
\eneq

\end{cor}

The following is  {{taken from 9.4 of \cite{LinTAI}}}.
%well-known fact which can be easily proved directly and which had been embedded into different proofs.

\begin{lem}\label{measureexistence}
Let $A$ be a unital separable \CA. For any $\ep>0,$  any  finite subset ${\cal H}\subset A_{s.a.},$ there exists a finite subset ${\cal G}\subset A$ and $\dt>0$  satisfying the following:
Suppose that $\phi: A\to B$ {{(for some unital  \CA\, $B$)}} is a ${\cal G}$-$\dt$-multiplicative
\morp\, and $t\in T(B)$ is a tracial state of $B.$ Then, there exists
a tracial state $\tau\in T(A)$ such that
\beq\label{mext-1}
|t\circ\phi(h)-\tau(h)|<\ep\tforal h\in {\cal H}.
\eneq
\end{lem}

%\begin{proof}
%{\color{Green} It follows from the same argument as that of Lemma \ref{LtoMN}.}
%\end{proof}
%{\bf The above  was proved somewhere. Need a reference.}
% We might need it in the next proof. Answer: We don't need it for next proof, may be you need it somewhere else.}

%{\bf To prove the following, we might need a simultaneously dig as in \ref{LtoMn} which can be done. We probably need \ref{measureexistence} here too.  The question is: do we really need that? }

\begin{thm}\label{Newunique1}
%Let $A=A_0$ or  $A=A_0\otimes C(\T),$ where $A_0\in {\cal C}_0,$ be a unital \CA\, or {\color{green} let $A=C(X)$ for a compact metrizable space $X$}.
Let $A\in {\overline{\cal D}}_s$ be a unital \CA.
Let $ \Delta: A_+^{ q, {\bf 1}}\setminus \{0\}\to (0,1)$ be an order preserving map.

Let  $\ep>0$ and let ${\cal F}\subset A$ be a finite subset. There exists a finite subset
${\cal H}_1\subset A_+^{{ {\bf 1}}}\setminus \{0\},$
%satisfying the following:
%For any $\sigma_1>0,$ there exists a finite subset
a finite subset ${\cal G}\subset A,$
$\dt>0,$ a finite subset ${\cal P}\subset \underline{K}(A),$
a finite subset ${\cal H}_2\subset A_{s.a.}$ and $\sigma>0$ satisfying the following:
Suppose that $L_1, L_2: A\to M_k$ (for some integer $k\ge 1$) are two unital
${\cal G}$-$\dt$-multiplicative \morp s such that
\beq\label{Newu-1}
[L_1]|_{\cal P}&=&[L_2]|_{\cal P},\\\label{Newu-1+}
{\rm tr}\circ L_1(h)&\ge& \Delta({\hat{h}})\,\,\,
%tr\circ L_2(h)\ge  \Delta(\hat{h})
 \tforal h\in {\cal H}_1\tand\\
|{\rm tr}\circ L_1(h)-\mathrm{tr}\circ L_2(h)|&<&\sigma\tforal h\in {\cal H}_2,
\eneq
then there exists a unitary $u\in M_k$ such that
\beq\label{Newu-4}
\|{\rm Ad}\, u\circ L_1(f)-L_2(f)|<\ep\tforal f\in {\cal F}.
\eneq
\end{thm}

\begin{proof}
The proof is exactly the same as that of \ref{UniqAtoM}.  As in the proof of \ref{UniqAtoM}, we will use \ref{Lauct2}. However, here we will use \ref{Combinerep}  instead of \ref{8-N-4}.
%{\color{green} If $A=C(X)$, then it follows from Theorem 2.10 of \cite{Lin-AU11}.} If  $A=A_0$ or  $A=A_0\otimes C(\T),$ where $A_0\in {\cal C}_0,$, then this follows from  \ref{Oldstableuniq} and \ref{Combinerep}.
%{\bf  This is the Main Goal of this section. Please write a proof. Answer: Since I modified \ref{Combinerep}, now should be OK}
\end{proof}

%\begin{cor}\label{CNewuni1}
%The statement of Theorem {\rm \ref{Newunique1}} holds if {\rm (\ref{Newu-1+})} is replaced by
%\beq\label{Newu-5}
%{\rm tr}\circ L_1(h)\ge \Delta(\hat{h}) \tforal h\in {\cal H}_1.
%\eneq
%\end{cor}
%\begin{lem}\label{Newuni2}
%Let $A\in {\cal A}_1$ be a unital \CA. let $\ep>0$ and let ${\cal F}\subset A$ be a finite subset. There exists a finite subset
%${\cal H}_1\subset A_+\setminus \{0\}$ satisfying the following:
%For any $\sigma_1>0,$ there exists a finite subset
%${\cal H}_2\subset A_{s.a.}$ and $\sigma_2>0$ satisfying the following:
%Suppose that $L_1, L_2: A\to M_k$ (for some integer $k\ge 1$) are two unital
%\hm s such that
%\beq\label{New2-1}
%tr\circ L_1(h)\ge \sigma_1,\,\,\,tr\circ L_2(h)\ge \sigma_1\tforal h\in {\cal H}_1\andeqn\\
%|tr\circ L_1(h)-tr\circ L_2(h)|<\sigma_2\tforal h\in {\cal H}_2,
%\eneq
%then there exists a unitary $u\in M_k$ such that
%\beq\label{New2-2}
%\|{\rm Ad}\, u\circ L_1(f)-L_2(f)|<\ep\tforal f\in {\cal F}.
%\eneq
%
%\end{lem}

%%%%%%%%%%%%%%%%%%%%%%%%%%
%%%%%%%%%%%%%
%%%%%%%%%
%%%%%%%%%%%%%%%%%%%%%%%%%%%%%%%%%%%%%%%%%%%%%%%%%%%%%%%%%%%%%

%%%%%%%%%%%%%%

%%%%%%%%%%%%%%%

%%%%%%%%%%

%%%%%%%%%%%%%

%%%%%%%%%%%

\section{Homotopy Lemma in finite dimensional \CA s}

\begin{lem}\label{homhom}
Let $A$ be a unital separable \CA\, and let $\phi: A\to M_k$ (for some integer $k\ge 1$) be a unital linear map. Suppose that $u\in M_k$ is a unitary
such that
\beq\label{homhom-1}
u\phi(a)=\phi(a)u\tforal a\in A.
\eneq
Then there exists a continuous path of unitaries
$\{u_t: t\in [0,1]\}\subset M_k$ such that
\beq\label{homhom-2}
u_0=u,\,\,\, u_1=1,\,\,\,
u_t\phi(a)=\phi(a)u_t\tforal a\in A
\eneq
and for all $t\in [0,1].$ Moreover,
\beq\label{homhom-3}
{\rm length}(\{u_t\})\le \pi.
\eneq
\end{lem}

\begin{proof}
There is $d>0$ such that spectrum of $u$ has a gap containing an arc with length at least $d.$ There is a continuous function $h$  from $sp(u)$ to
$[-\pi, \pi]$ such that
\beq\label{homhom-4}
\exp(i h(u))=u.
\eneq
Therefore
\beq\label{homhom-5}
\phi(a)h(u)=h(u)\phi(a)\tforal a\in A.
\eneq
Note that $h(u)\in (M_k)_{s.a.}$ and $\|h(u)\|\le \pi.$
Define $u_t=\exp(i(1-t)h(u))$ ($t\in [0,1]$). Then
$
u_0=u\andeqn u_1=1.
$
Also
$$
u_t\phi(a)=\phi(a)u_t
$$
for all $a\in A$ and $t\in [0,1].$
Moreover, one has
$$
{\rm lenghth}(\{u_t\})\le \pi,
$$
as desired.
\end{proof}

%\begin{lem}\label{MExt}
%Let $A\in {\cal A}_1$ be a unital \CA, let $\sigma>0,$ let ${\cal H}\subset A_{s.a.}$ be a finite subset. There exists a finite subset ${\cal G}\subset A$ and $\dt>0$ satisfying the following:
%Suppose that $L: A\to M_k$ (for some integer $k\ge 1$) is a unital $\dt$-${\cal G}$-multiplicative \morp. Then there is a unital \hm\,
%$\psi: A\to M_k$ such that
%\beq\label{MExt-1}
%|tr\circ L(f)-tr\circ \psi(f)|<\sigma\tforal f\in {\cal H}.
%\eneq
%\end{lem}
%
%%\begin{proof}
%It follows from \ref{measureexistence}
%\end{proof}

\begin{lem}\label{EXTMM}
Let $A\in {\overline{\cal D}}_s$ be a unital C*-algebra, let ${\cal H}\subset (A\otimes C(\T))_{s.a.}$ be a finite subset,
  let $1> \sigma>0$ be a positive number and let $\Delta: A_+^{q,\bf 1}\setminus \{0\}\to (0,1)$ be an o{{r}}der preserving map.  Let  $\ep>0,$ ${\cal G}_0\subset A\otimes C(\T)$ be a finite subset, ${\cal P}_0, {\cal P}_1\subset \underline{K}(A)$ be finite subsets and
  ${\cal P}={\cal P}_0\cup \bt({\cal P}_1)\subset \underline{K}(A\otimes C(\T)).$
There exist $\dt>0,$ a finite subset ${\cal G}\subset A\otimes C(\T)$
%a finite subset ${\cal G}_0\subset A$
and a finite subset ${\cal H}_1\subset (A\otimes C(\T))_+^{\bf 1}\setminus \{0\}$ satisfying the following:
Suppose that $L: A\otimes C(\T)\to M_k$ (for some integer $n\ge 1$) is a
${\cal G}$-$\dt$-multiplicative \morp\,
%and $\phi: A\to M_k$ is a unital \hm\,
such that
\beq\label{EXTMM-1}
%\|L(g\otimes 1)-\phi(g)\|&<&\dt\tforal g\in {\cal G}_0,\\\label{EEXTMM-1+}
{\rm tr}\circ L(h) &\ge& \Delta(\hat{h})\tforal h\in {\cal H}_1,
\tand\\
{[L]}|_{\bt ({\cal  P}_1)}&=&0.
\eneq
Then there exists a unital ${\cal G}_0$-$\dt$-multiplicative \morp\, $\psi: A\otimes C(\T)\to M_k$ such that
$u=\psi(1\otimes z)$ is a unitary,
\beq\label{EXTMM-2}
u\psi(a\otimes 1)&=&\psi(a\otimes 1)u\tforal a\in A\\
{[L]}|_{\cal P}&=&[\psi]|_{\cal P}\tand\\
|{\rm tr}\circ L(h)-{\rm tr}\circ \psi(h)|&<&\sigma\tforal h\in {\cal H}.
\eneq

\end{lem}

\begin{proof}
Let  ${\cal H}$ and $\sigma_0,$ $\ep$ and ${\cal G}_0$ are given.
Without loss of generality, we may assume that
${\cal H}\subset {\cal G}_0$ which is in the unit ball of $A$ and
$\sigma<\ep/4.$ %{\color{Green} (Either something here or something wrong below)}
We may also assume that
$$
{\cal G}_0=\{g\otimes f: g\in {\cal G}_{0A}\andeqn f\in {\cal G}_{1T}\},
$$
where ${\cal G}_{0A}\subset A$ and ${\cal G}_{1T}\subset C(\T)$ are finite
subsets. To simplify matter further, we may assume, without loss of generality, that
${\cal G}_{1T}=\{1_{C(\T)}, z\},$ where $z\in C(\T)$ is the standard unitary generator.

We may assume that ${\cal G}_{0A}$ is sufficiently large and $\ep$ is sufficiently small such that
$[L_1]|_{\cal P}$ is well defined for any unital ${\cal G}_0$-$\ep$-multiplicative \morp\, from $A\otimes C(\T)$ and
\beq\label{Extmm-3}
[L_1]|_{{\cal P}_0}=[L_2]|_{{\cal P}_0}
\eneq
for any unital ${\cal G}_{0A}$-$\ep$-multiplicative \morp\, $L_2$
from $A\otimes C(\T)$ such that
\beq\label{Extmm-4}
L_1\approx_{\ep} L_2\,\,\,{\rm on}\,\,\, {\cal G}_{0A}.
\eneq
%We may also assume that $\ep<\sigma.$ {\color{Green} (Either something here or something wrong above.)}

Let $n$ be an integer such that $1/n<\sigma/2.$
Note that $A\otimes C(\T)\in {\overline {\cal D}}_s.$

Let $\dt>0,$ ${\cal G}\subset A\otimes C(\T)$ and ${\cal H}_1\subset
A\otimes C(\T)^{\bf 1}_+\setminus \{0\}$ (in place of ${\cal H}_2$) be finite subsets required by \ref{Combinerep} for $A\otimes C(\T)$ (in place of $A$), $\ep/2$ (in place of $\ep$), ${\cal G}_0$ (in place of ${\cal F}$),
${\cal H}$ (in place of ${\cal H}_1$) and $\Delta.$
Now suppose that $L: A\otimes C(\T)\to M_k$ satisfies the assumption
for the above $\dt,$ ${\cal G}$ and ${\cal H}_1.$
It follows from \ref{Combinerep} that there is a projection $e_0\in M_k$ and a ${\cal G}_0$-$\ep/2$-multiplicative \morp $\psi_0: A\otimes C(\T)\to e_0M_ke_0$ and a unital \hm\, $\psi_1: A\otimes C(\T)\to (1-e_0)M_k(1-e_0)$
such that
\beq\label{Exmm-4}
{\rm tr}(e_0) &<&1/n<\sigma,\\
\|L(a)-\psi_0(a)\oplus \psi_1(a)\|&<&\ep\tforal a\in {\cal G}_0.
\eneq
 Define $\psi: A\otimes C(\T)\to M_k$ by
 $\psi(a)=\psi_0(a)\oplus \psi_1(a)$ for all $a\in A$ and
 $\psi(1\otimes z)=e_0\oplus \psi_1(1\otimes z).$

Put $u=\psi(1\otimes z).$ One verifies that this $\psi$ and $u$ satisfy all requirements.
\end{proof}

%{\bf It is easier than the first look-----applying \ref{Combinerep}.}

%\begin{lem}\label{MExtmeasure}
%Let $A\in {\cal A}_1$ be a unital \CA\, let ${\cal H}\subset A_+\setminus \{0\}$ be a finite subset and let $1>r>0.$ For any $1>\sigma>0,$ there exists  $\dt>0,$ a finite subset ${\cal G}\subset A$  satisfying the following:
%Suppose that $L: A\to M_k$ (for some integer $n\ge 1$) is a $\dt$-${\cal G}$-multiplicative \morp\, such that
%\beq\label{Mextm-1}
%{\rm tr}\circ L(h)>\sigma \tforal h\in {\cal H}.
%\eneq
%Then there exists a faithful tracial state $t$ such that
%\beq\label{Mextm-2}
%{\rm tr}\circ L(h)\ge r t(h)\tforal h\in {\cal H}.
%\eneq
%\end{lem}

\begin{lem}\label{homfull}
Let $A\in {\overline{\cal D}}_s$ be a unital \CA\, and
%Let $A=A_0,$ or $A=A_0\otimes C(\T)$ for some $A_0\in {\cal C}_0$, {\color{green} or let $A=C(X)$ for a metric space $X$.} % be a unital separable \CA\, and
let $\Delta: (A\otimes C(\T))_+^{q, {\bf 1}}\setminus\{0\}\to (0,1)$ be an order  preserving map.
Let $\ep>0$ and let
${\cal F}\subset A$ be a finite subset. There exists a finite subset
${\cal H}_1\subset A_+^{\bf 1}\setminus \{0\},$ a finite subset
${\cal H}_2\subset C(\T)_+^{\bf 1}\setminus \{0\},$
a finite subset ${\cal G}\subset A,$ $\dt>0$
%a finite
%subset ${\cal H}_2\subset C(\T)_+\setminus\{0\}$
and a finite subset ${\cal P}\subset \underline{K}(A)$ such that,
if $L: A\otimes C(\T)\to M_k$ (for some integer $k\ge 1$) is ${\cal G}'$-$\dt$-multiplicative \morp,
where ${\cal G}'=\{g\otimes f: g\in {\cal G}, f=\{1, z,z^*\}\}$ and $u\in M_k$ is a unitary
such that
\beq\label{homfull-1}
%\tand
\|L(1\otimes z)-u\|&<&\dt,
\\\label{homfull-2}
[L]|_{\bt(\cal P)}&=&0\tand\\\label{homfull-3}
{\rm tr}\circ L(h_1\otimes h_2)&\ge &\Delta(\widehat{h_1\otimes h_2})
\eneq
for all $h_1\in {\cal H}_1$ and $h_2\in {\cal H}_2,$ then
there exists a continuous path of unitaries $\{u_t: t\in [0,1]\}\subset M_k$
with $u_0=u$ and $u_1=1$ such that
\beq\label{homfull-3+}
\|L(f\otimes 1)u_t-u_tL(f\otimes 1)\|<\ep\tforal f\in {\cal F}
\eneq
and $t\in [0,1].$
Moreover,
\beq\label{homfull-4}
{\rm length}(\{u_t\})\le \pi+\ep.
\eneq
\end{lem}

\begin{proof}
%Let $z\in C(\T)$ be the identity function on the unit circle.
%{\color{green} Assume that $A=A_0,$ or $A=A_0\otimes C(\T)$ for some $A_0\in {\cal C}_0$.}
Let $\Delta_1=(1/2)\Delta,$
${\cal F}_0=\{f\otimes 1: 1\otimes z: f\in {\cal F}\}$ and
let $B=A\otimes C(\T).$  Then $B\in {\bar{\cal D}}_s.$ Let ${\cal H}'\subset B_+^{\bf 1}\setminus\{0\}$ (in place of ${\cal H}_1$),  {{ ${\cal H}_0\subset A_{s.a.}$}} {{(in place of ${\cal H}_2$) be a finite subset,}} ${\cal G}_1\subset A\otimes C(\T)$ (in place of ${\cal G}$), $\dt_1>0$ (in place of $\dt$), ${\cal P}'\subset \underline{K}(B)$ (in place of ${\cal P}$) be the finite sets and constants {{as}} required by \ref{Newunique1} (for $B$ instead of $A$) for $\ep/16$ (in place of
$\ep$), ${\cal F}_0$ (in place of ${\cal F}$) and $\Delta.$
 Without loss of generality,
 {{we may assume that ${\cal H}_0\subset B_+^{\bf 1}\setminus \{0\}$ and further, to simplify notation,
 we may assume that ${\cal H}_0={\cal H}'.$}}  We may assume
that there are finite sets ${\cal  H}_1'\subset A_+^{\bf 1}\setminus \{0\}$,
${\cal  H}_2'\subset C(\T)_+^{\bf 1}\setminus \{0\}$,
and
${\cal G}_1'\subset A$
such that
$$
{\cal  H}'=\{h_1\otimes h_2: h_1\in {\cal H}_1'\andeqn h_2\in {\cal H}_2'\}
$$
%Let $\tau\in T(A\otimes C(\T))$ be a faithful tracial state and let $1>r>0.$
%Set
%\beq\label{homfull-5}
%\sigma_1=\inf\{\Delta(\hat{b}): b\in {\cal H}'\}.
%\eneq
and
${\cal G}_1=\{g\otimes f: g\in {\cal G}_1'\andeqn f\in \{1,z,z^*\}\}$. We may also assume
that $1_A\in {\cal H}_1'$ and $1_{C(\T)}\in {\cal H}_2'.$

%Let ${\cal H}_3={\cal H}_1'\cup {\cal H}_2'.$

%Let ${\cal G}''\subset A\otimes C(\T)$ (in place of ${\cal G}$) be a finite subset, $\dt_2>0$ (in place of $\dt$) required by \ref{MExt}
%for ${\cal H}_3$ (in place of ${\cal H}$) and $\sigma_3$ (in place of $\sigma$),
%where $\sigma_3=\min\{\sigma_1/3, \sigma_2/3\}.$
%Choose ${\cal G}_1'={\cal G}'\cup {\cal G}_2$ and
%$\dt_3=\min\{\dt_1/2, \dt_2/2\}.$
Without loss of generality, one may assume that
\beq\label{homfull-8}
{\cal P}'={\cal P}_0\sqcup {\cal P}_1,
\eneq
where ${\cal P}_0\subset \underline{K}(A)$ and
${\cal P}_1\subset {\boldsymbol{\bt}}(\underline{K}(A))$ are finite subsets.
Let ${\cal P}\subset \underline{K}(A)$ be a finite subset such that
${\boldsymbol{\bt}}({\cal P})={\cal P}_1.$
Let
\beq\label{homfull-8+1}
\sigma=\min\{\Delta_1(\hat{h}): h\in {\cal H}'\}.
\eneq
%Let $\ep_1=\min\{\ep/16, \dt_1/4, \sigma_1/16\}.$

There is $\dt_2>0$ (in place of $\dt$) with $\dt_2<\ep/16,$ a finite subset ${\cal G}_2\subset A\otimes C(\T)$(in place of ${\cal G}$) and
a finite subset ${\cal H}_3\subset (A\otimes C(\T))_+^{\bf 1}\setminus \{0\}$ (in place of ${\cal H}_1$)
required by \ref{EXTMM} for $\sigma,$ $\Delta,$  ${\cal H}'$(in place of
${\cal H}$), $\min\{\ep/16,\dt_1/2\}$ (in place of $\ep$),
${\cal G}_1$ (in place of ${\cal G}_0$),  ${\cal P}_0$ and ${\cal P}$
(in place of ${\cal P}_1$).
We may also assume that
$$
{\cal G}_2=\{g\otimes f: g\in {\cal G}_2'\andeqn f\in \{1, z,z^*\}{{\}}}
$$
for a finite set ${\cal G}_2'\subset A$.
%Let $\dt_1>\dt_2>0$  and ${\cal G}_0'\subset A$ satisfying the following:
%For any $\dt_1$-${\cal G}_0''$-multiplicative \morp\, $L': A\otimes C(\T)\to C$ (for any unital \CA\, $C$),
%where ${\cal G}_0''=\{g\otimes 1: g\in {\cal G}_0'\},$ such that
%$$
%\|L'(g\otimes 1)w-L'(g\otimes 1)w\|<\dt_2\tforal g\in {\cal G}_0'
%$$
%and for any unitary $w\in C$ with
%$$
%\|L(1\otimes z)-w\|<\dt_2,
%$$
%then $L$ is $\dt_1$-${\cal G}_1$-multiplicative.
We may further assume that
$$
{\cal H}_3=\{h_1\otimes h_2: h_1\in {\cal H}_4\andeqn h_2\in {\cal H}_5\}
$$
{{for}} finite sets ${\cal H}_4\subset A_+^{\bf 1}\setminus \{0\}$ and
${\cal H}_5\subset C(\T)_+^{\bf 1}\setminus \{0\}$.
Let ${\cal G}={\cal F}\cup {\cal G}_1'\cup {\cal G}_2',$ $\dt=\min\{\dt_1/2,\dt_2/2, \ep/16\},$
${\cal H}_1={\cal H}_1'\cup {\cal H}_4$ and
${\cal H}_2={\cal H}_2'\cup {\cal H}_5.$

Now suppose that one has a linear map $L: A\otimes C(\T)\to M_k$ and a unitary $u\in M_k$ satisfy the assumption with
the above ${\cal H}_1,$ ${\cal H}_2,$ ${\cal G},$ ${\cal P},$  $\dt$ and $\sigma.$
%As mentioned above, $L$ is $\dt_1$-${\cal G}_1$-multiplicative.
It follows from \ref{EXTMM} that there is a unital  {{${\cal G}_1$-$\min\{\ep/16, \dt_1/2\}$}}-multiplicative
\morp\,
$\psi: A\otimes C(\T)\to M_k$ such that $w=\psi(1\otimes z)$ is a unitary,
\beq\label{homfull-9}
w\psi(g\otimes 1)&=&\psi(g\otimes 1)w\tforal g\in A,\\\label{homfull-9n}
[\psi]|_{\cal P'}&=&[L]|_{\cal P'},\\\label{homfull-9+}
|{\rm tr}\circ L(g)-{\rm tr}\circ \psi(g)|&<&\sigma\tforal g\in {\cal H}_3.
\eneq
It follows that
\beq\label{homfull-10}
\mathrm{tr}\circ \psi(h)\ge \mathrm{tr}\circ L(h)-\sigma\ge \Delta_1(\hat{h})
\eneq
for all $h\in {\cal H}'.$
Combining (\ref{homfull-9}), (\ref{homfull-2}), (\ref{homfull-3}), (\ref{homfull-9+}) and (\ref{homfull-3+}), by applying \ref{Newunique1},
%and \ref{CNewuni1},
one obtains a unitary $U\in M_k$ such that
\beq\label{homfull-12}
\|{\rm Ad}\,U\circ \psi(f)-L(f)\|<\ep/16\tforal f\in {\cal F}_0.
\eneq
Let $w_1={\rm Ad}\,U\circ \phi(1\otimes z).$ Then
\beq\label{homfull-13}
\|u-w\| &\le & \|u-L(1\otimes z)\|+\|L(1\otimes z)-{\rm Ad}\, U\circ \psi(1\otimes z)\|\\
&<&\dt+\ep/16<\ep/8.
\eneq
Then there is a continuous path of unitaries $\{u_t\in [0,1/2]\}\subset M_k$
such that
\beq\label{homfull-13+}
\|u_t-u\|<\ep/8,\,\,\,\|u_t-w\|<\ep/8,\,\,\, u_0=u,\,\,\, u_{1/2}=w\\
\andeqn
{\rm length}(\{u_t: t\in [0,1/2]\})<\ep\pi/8.
\eneq
It follows from \ref{homhom} that there exists a continuous path of unitaries $\{u_t: t\in [1/2,1]\}\subset M_k$ such that
\beq\label{homfull-14}
u_{1/2}=w,\,\,\, u_1=1\andeqn u_t({\rm Ad}\, U\circ \phi(f\otimes 1))=({\rm Ad}\, U\circ \phi(f\otimes 1))u_t
\eneq
for all $t\in [1/2,1]$ and $f\in A\otimes 1.$
Moreover,
\beq\label{homfull-15}
{\rm length}(\{u_t: t\in [1/2,1]\})\le \pi.
\eneq
It follows that
\beq\label{homfull-16}
{\rm length}(\{u_t: t\in [0,1]\}) \le \pi+\ep\pi/6.
\eneq
Furthermore,
\beq\label{homfull-17}
\|u_tL(f\otimes 1)-L(f\otimes 1)u_t\|<\ep\tforal f\in {\cal F}
\eneq
and $t\in [0,1].$
%
%{\color{green}Assume that $A=C(X)$. Then the Lemma follows from exactly the same argument as above, with \ref{Newunique1} being replaced by Theorem 2.10 of \cite{Lin-AU11}.}
\end{proof}

\begin{lem}\label{oldnuclearity} {\rm (Lemma 2.8 of \cite{LnHomtp} )}
Let $A$ be a unital  {{ separable amenable}} \CA. Let $\ep>0,$ let
${\cal F}_0\subset A$ be a finite subset and let ${\cal F}\subset A\otimes C(\T)$ be a finite subset. There exists a finite subset ${\cal G}\subset A$ and $\dt>0$ satisfying the following:
For any ${\cal G}$-$\dt$-multiplicative \morp\, $\phi: A\to B$ (for some unital \CA\, $B$) and any unitary $u\in  B$ such that
\beq\label{oldnuc-1}
\|\phi(g)u-u\phi(g)\|<\dt\tforal g\in {\cal G},
\eneq
there exists a unital ${\cal F}$-$\ep$-multiplicative \morp\, $L: A\otimes C(\T)\to B$ such that
\beq\label{oldnuc-2}
\|\phi(f)-L(f\otimes 1)\|<\ep\andeqn
\|L(1\otimes z)-u\|<\ep
\eneq
for all $f\in {\cal F}_0,$ where $z\in C(\T)$ is the identity function
on the unit circle.
\end{lem}

\begin{lem}\label{changespectrum}
%{\color{green} Let $A=A_0,$ or $A=A_0\otimes C(\T)$ for some $A_0\in {\cal C}_0$, or let $A=C(X)$ for a metric space $X$.}  %
Let $A\in {\bar {\cal D}}_s$ be a unital \CA,
Let  $\ep>0$  and let ${\cal F}\subset A$ be a finite subset.
 Let ${\cal H}_1\subset A_+^{\bf 1}\setminus\{0\}$ and let  ${\cal H}_2\subset C(\T)_+^{\bf 1}\setminus \{0\}$  be  finite subsets.
For any order preserving map $\Delta: A_+^{q, {\bf 1}}\setminus \{0\}\to (0,1),$ there exists a finite subset ${\cal G}\subset A,$
a finite
subset ${\cal H}_1'\subset A_+\setminus \{0\}$
and $\dt>0$ satisfy the following: for any unital ${\cal G}$-$\dt$-multiplicative \morp\, $\phi: A\to M_k$ (for some integer $k\ge 1$) and any unitary $u\in M_k$ such that
\beq\label{change-1}
\|u\phi(g)-\phi(g)u\|<\dt\tforal g\in {\cal G} \tand \\ \label{change-1+}
tr\circ \phi(h)\ge \Delta(\hat{h})\tforal h\in {\cal H}_1',
\eneq
there exists a continuous path of unitaries
$\{u_t: t\in [0,1]\}\subset M_k$  such that
\beq\label{change-2}
&&u_0=u,\,\,\, u_1=w,\,\,\,\|u_t\phi(f)-\phi(f)u_t\|<\ep\tforal f\in {\cal G}\tand t\in [0,1],\\
%\eneq
%for all $f\in {\cal G}$ and $t\in [0,1],$
%\beq
%\\\label{change-3}
&&\hspace{0.4in}\mathrm{tr}\circ L(h_1\otimes h_2)\ge  \Delta(\widehat{h_1}) \tau_m(h_2)/4
\eneq
for all $h_1\in {\cal H}_1$ and $h_2\in {\cal H}_2,$
where $L: A\otimes C(\T)\to M_k$ is a \morp\,
such that
\beq\label{changes-4}
%&&
\|L(f\otimes 1)-\phi(f)\|<\ep \tforal f\in {\cal F}, \tand
%&&
\|L(1\otimes z)-w\|<\ep,
\eneq
and
$\tau_m$ is the tracial state on $C(\T)$ induced by the Lebesgue
%%%Lesbegue 
measure
on the circle.
Moreover,
\beq\label{changes-5}
{\rm length}(\{u_t\})\le \pi+\ep.
\eneq

\end{lem}

\begin{proof}
%{\color{green} Assume that $A=A_0,$ or $A=A_0\otimes C(\T)$ for some $A_0\in {\cal C}_0$}
There exists an integer $n\ge 1$ such that
\beq\label{changes-6}
 (1/n)\sum_{j=1}^n f(e^{\theta+j 2\pi i/n})\ge (63/64)\tau_m(f)
\eneq
for all $f\in {\cal H}_2$ and for any $\theta\in [-\pi, \pi].$
We may also assume that $16\pi/n<\ep.$

%Let $1>\sigma>0.$
%Let ${\cal H}_1'\subset A_+\setminus \{0\}$ be a finite subset
%required by \ref{MExtmeasure}  for $r=3/4,$ $\sigma/4$ and ${\cal H}_1$ (in place of ${\cal H}$).  Without loss of generality, one may assume that
%${\cal H}_1\subset {\cal H}_1'.$
%Let $\dt_0>0$ and let ${\cal G}_0$ be a finite required by \ref{MExtmeasure} for ${\cal H}_1$ (in place of ${\cal H}$), $r=3/4$ and $\sigma.$
%Let $t$ be a faithful tracial state on $A$ given by
%\ref{MExtmeasure} corresponding  ${\cal H}_1$ (in place of ${\cal H}$) and
%$r=3/4.$

Let
$$
\sigma_1=(1/2^{10})\inf\{t(h): h\in {\cal H}_1\}\cdot  \inf\{\tau_m(g): g\in {\cal H}_2\}.
$$

%Define
%\beq\label{changes-7}
%\tau_0=t\otimes \tau_m\,\,\,{\rm on}\,\,\, A\otimes C(\T).
%\eneq
Let ${\cal F}'=\{f\otimes 1, f\otimes z: f\in {\cal F}\cup {\cal H}_1\}.$
Let $\dt_1>0$ (in place of $\dt$) and ${\cal G}_1\subset A\otimes C(\T)$ (in place of ${\cal G}$) be a finite subset
{{as}} required by \ref{LtoMn} for $\ep/32$
(in place of $\ep$), ${\cal F}'$ (in place of ${\cal F}$) and $\sigma_1/16$ (in place of $\sigma_0$).
Without loss of generality, one may assume that, for a finite set ${\cal G}_2\subset A,$
$$
{\cal G}_1=\{g\otimes 1, 1\otimes z: g\in {\cal G}_2\}.
$$
%for a finite set ${\cal G}_2\subset A$.

Let ${\cal H}_1'\subset A_+\setminus \{0\}$  (in place of ${\cal H}_2$) be a finite subset {{as}} required by
\ref{repcor} for $\min\{\ep/32,\sigma_1/16\}$ (in place of $\ep$), ${\cal F}\cup {\cal H}_1$ (in place of ${\cal F}$), ${\cal H}_1$ (in place of ${\cal H}$), $(190/258)\Delta$ (in place of $\Delta$) and  $\sigma_1/16$ (in place of $\sigma$) and integer
$n.$

Put
$$
{\cal H}'=\{h_1\otimes h_2, h_1\otimes 1, 1\otimes h_2: h_1\in {\cal H}_1\andeqn h_2\in {\cal H}_2\}.
$$

Let ${\cal G}_3={\cal G}_2\cup {\cal H}_1\cup {\cal H}_1'.$ To simplify the notation, without loss of generality, one may assume that
${\cal G}_3$ and ${\cal F}'$ are all in the unit ball of $A\otimes C(\T).$ Let
$\dt_2=\min\{\ep/64,\dt_1/2, \sigma_1/16\}.$

Let ${\cal G}_4\subset A$ be a finite subset (in place of ${\cal G}$) and let
$\dt_3$ (in place of $\dt$) be positive as required by \ref{oldnuclearity} for ${\cal G}_3$ (in place of ${\cal F}_0$), ${\cal F}'$ (in place of ${\cal F}$),  and $\dt_2$ (in place of $\ep$).

Let ${\cal G}={\cal G}_4\cup {\cal G}_3$ and
$\dt=\min\{\dt_1/4, \dt_2/2, \dt_3/2\}.$
Now let $\phi: A\to M_k$ be a unital ${\cal G}$-$\dt$-multiplicative \morp\,,  and let $u\in M_k$ be a unitary
such that (\ref{change-1}) and (\ref{change-1+}) hold for the above $\dt,$ $\sigma,$ ${\cal G}$ and
${\cal H}_1'.$
%By \ref{MExtmeasure},
%\beq\label{changes-8-1}
%{\rm tr}\circ \phi(h)\ge (3/4)t(h)\rforal h\in {\cal H}_1.
%\eneq

It follows from \ref{oldnuclearity} that there
exists a ${\cal G}_3$-$\dt_2$-multiplicative \morp\,
$L_1: A\otimes C(\T)\to M_k$ such that
\beq\label{changes-8}
\|L_1(g\otimes 1)-\phi(g)\|<\dt_2\tforal g\in {\cal G}_2\andeqn
\|L_1(1\otimes z)-u\|<\dt_2.
\eneq
%It follows from \ref{MExtmeasure} that
We then have that
\beq\label{changes-8+1}
{\rm tr}\circ L_1(h\otimes 1)&\ge &{\rm tr}\circ \phi(h)-\dt_2\\\label{changes-8+2}
&\ge & \Delta(\hat{h})-\sigma_1/16\ge (191/256)\Delta(\hat{h})
\eneq
for all $h\in {\cal H}_1.$
%One also has that
%\beq\label{changes-8+2}
%{\rm tr}\circ L_1(h\otimes 1) &\ge & {\rm tr}\circ \phi(h)-\sigma_1/4\\
%&\ge &  \sigma-\sigma_1/4>\sigma/2
%\eneq
%for all $h\in {\cal H}_1.$
It follows from \ref{LtoMn} that there exist a projection
$p\in M_k$ and a unital \hm\, $\psi: A\otimes C(\T)\to pM_kp$ such that
\beq\label{changs-9}
&&\|pL_1(f)-L_1(f)p\|<\min\{\ep/32,\sigma_1/16\}\tforal f\in {\cal F}',\\\label{changes-9+}
&&\|L_1(f)-((1-p)L_1(f)(1-p)+\psi(f))\|<\min\{\ep/32,\sigma_1/16\}\tforal f\in {\cal F}' \\\label{changes-10}
&&\andeqn {\rm tr}(1-p)<\sigma_1/16.
\eneq
Note that $pM_kp\cong M_m$ for some $m\le k.$ It follows from (\ref{changes-8+2}), (\ref{change-1+}), (\ref{changes-9+})  and (\ref{changes-10}) that
\beq\label{changes-11}
{\rm tr}\circ \psi(h)\ge (191/256)\Delta(\hat{h})-\sigma_1/16-\sigma_1/16 \ge (190/256)\Delta(\hat{h})
%\eneq
\rforal h\in {\cal H}_1.
\eneq

By \ref{repcor}, there are mutually orthogonal
projections $e_0, e_1, e_2,...,e_n\in pM_kp$  such that  $e_1, e_2,...,e_n$ are equivalent and there are unital \hm s $\psi_0: A\otimes C(\T)\to e_0M_ke_0$ and
$\psi_1: A\otimes C(\T)\to e_1M_ke_1$ such that
\beq\label{changes-12}
\|\psi(f)-{\rm diag}(\psi_0(f),\overbrace{\psi_1(f),...,\psi_1(f)}^n)\|&<&\min\{\ep/32, \sigma_1/6\}\tforal
f\in {\cal F}_1\\
\andeqn{\rm tr}(e_0)&<&\sigma_1/16.
\eneq
%Moreover,
%\beq\label{changes-13}
%{\rm tr}\psi_1(h\otimes 1)\ge (3/4)t(h)/2n\tforal h\in {\cal H}_1.
%\eneq
Let $w_0'=\psi_1(1\otimes z).$ One may write
$$
w_0'={\rm diag}(\exp(ia_1),\exp(ia_2),...,\exp(ia_n)),
$$
where $a_j\in e_jM_ke_j$ is
 a selfadjoint element
with $\|a_j\|\le \pi.$
%By linear algebra, it is easy to find
Choose a continuous path of unitaries
$\{w_{t,j}': t\in [0,1]\}\subset  e_jM_ke_j$ such that
\beq\label{changes-14}
w_{0,j}'=\exp(ia_j),\,\,\, w_{1,j}'=\exp(i(2\pi j/n)),
\andeqn {\rm length}(\{w_{t,j}'\})\le \pi+\ep/4.
\eneq
Moreover, one can choose such $w_{t,j}'$
%is
{{in}} the commutant of $\psi_1(A)$. % that it commutes with every element in $\psi_1(f),$ $f\in A.$
There is a unitary $w_0''\in (1-p)M_k(1-p)$ such that
\beq\label{changes-15}
\|w_0''-(1-p)L_1(1\otimes z)(1-p)\|<\ep/16.
\eneq
Put
\beq\label{changes-16}
u_0'=w_0''\oplus \psi_0(1\otimes z)\oplus w_0'.
\eneq
Then $u_0'$ is a unitary and
\beq\label{changes-17}
\|u-u_0'\| &\le & \|u-L_1(1\otimes z)\|+\|L_1(1\otimes z)-u_0'\|\\
&\le & \dt_2+\ep/16<\ep/8.
\eneq
One obtains a continuous path of unitaries $\{w_t\in [0,1]\}\subset M_k$ such that
\beq\label{changes-18}
w_0=u,\,\,\,w_1=w_0''\oplus \psi_0(1\otimes z)\oplus {\rm diag}(w_{1,1}', w_{1,2}',...,w_{1,n}')\\
\|w_t\phi(f)-\phi(f)w_t\|<\ep\tforal f\in {\cal F}, \andeqn
 {\rm length}(\{w_t\})\le \pi+\ep.
\eneq
Define $L: A\otimes C(\T)\to M_k$ by
\vspace{-0.1in}\beq
L(a\otimes f)=(1-p)L_1(a\otimes f)(1-p)\oplus
({\rm diag}(\psi_0(a), \overbrace{\psi_1(a),...,\psi_1(a)}^n)f(w_1))
\eneq
for all $a\in A$ and $f\in C(\T).$
It follows that
\beq\label{changes-19}
\|L(f\otimes 1)-\phi(f)\|<\ep\tforal f\in {\cal F}\andeqn
\|L(1\otimes z)-w_1\|<\ep.
\eneq
One also has
that (note that $w_{1,j}'$ is scalar)
\beq\label{changes-20}
\hspace{-0.4in}{\rm tr}\circ L(h_1\otimes h_2) &\ge & {\rm tr}(\psi_0(h_1)+{\rm  tr}(\diag(\overbrace{\psi_1(a),...,\psi_1(a)}^n)f(w_1))\\
%n{\rm tr}(\psi_1(h_1\otimes 1))){\rm tr}(h_2(w_1))\\
&\ge &{\rm tr}\circ \psi(h_1)\sum_{j=1}^n h_2(e^{i2\pi j/n})
\ge  (190/256)\Delta(\widehat{h_1})\sum_{j=1}^n h_2(e^{i2\pi j/n})\\
&\ge & (190/256)\Delta(\widehat{h_1})(63/64)\tau_m(h_2)
%&\ge & (190/256)\Delta(\widehat{h_1})({1-\sigma_1/16\over{n}})\sum_{j=1}^n h_2(e^{i2\pi j/n})-\sigma_1/6\\
%&\ge & {\rm tr}\circ \psi(h_1)({1-\sigma_1/16\over{n}})\sum_{j=1}^n h_2(e^{i2\pi j/n})-\sigma_1/6\\
%&\ge & (190/256)\Delta(\widehat{h_1})({1-\sigma_1/16\over{n}})\sum_{j=1}^n h_2(e^{i2\pi j/n})-\sigma_1/6\\
%&\ge & (190/256)\Delta(\widehat{h_1})(63/64)(1-\sigma_1/16)\tau_m(h_2)-\sigma_1/6\\
%&\ge & (190/256)\Delta(\widehat{h_1})({(63/64)(1-1/2^{14}})\tau_m(h_2)-(1/2^{12})t(h_1)\tau_m(h_2)\\
%&\ge & 
\ge \Delta(\widehat{h_1})\cdot \tau_m(h_2)/4
\eneq
for all $h_1\in {\cal H}_1$ and $h_2\in {\cal H}_2.$
%{\color{green} If $A=C(X)$, the Lemma follows from the exact same argument with \ref{repcor} being replaced by ??.}
\end{proof}

\begin{df}\label{Ddel}
Let $A$ be a unital \CA\, with $T(A)\not=\emptyset$ and let $\Delta: A_+^{q,{\bf 1}}\setminus \{0\}\to (0,1)$ be an order preserving
map. Suppose that $\tau_m: C(\T)\to \C$ is the tracial state given by the normalized Lebesgue
%%%Lesbegue
 measure.
Define
$\Delta_1: (A\otimes C(\T))_+^{q, {\bf 1}}\setminus\{0\}\to (0,1)$ by
\beq\label{Ddel-1}
\hspace{-0.3in}\Delta_1(\hat{h})=\sup\{ {\Delta(h_1)\tau_m(h_2)\over{4}}:\,
\hat{h}\ge \widehat{h_1\otimes h_2}\andeqn h_1\in A_+\setminus \{0\},\,\,\, h_2\in C(\T)_+\setminus\{0\}\}.
\eneq
\end{df}

\begin{lem}\label{homotopy1}
%{\color{green} Let $A=A_0,$ or $A=A_0\otimes C(\T)$ for some $A_0\in {\cal C}_0$, or let $A=C(X)$ for a metric space $X$.}
%Let $A=A_0,$ or $A=A_0\otimes C(\T)$ for some $A_0\in {\cal A}_1$ be a unital \CA\, and
Let $A\in  {\overline {\cal D}}_s$ be a unital \CA.
Let $\Delta: A_+^{q, {\bf 1}}\setminus \{0\}\to (0,1)$ be an order preserving  map.
For any $\ep>0$ and any finite subset ${\cal F}\subset A,$ there exists a finite subset ${\cal H}\subset A_+^{\bf 1}\setminus\{0\},$
 $\dt>0,$  a finite subset
${\cal G}\subset A$  and a finite subset ${\cal P}\subset \underline{K}(A)$ satisfying the following:
For any unital  ${\cal G}$-$\dt$-multiplicative \morp\, $\phi: A\to M_k$ (for some integer $k\ge 1$) and any unitary $v\in M_k$ such that
\beq\label{homot-1}
&&tr\circ \phi(h)\ge \Delta(\hat{h}) \tforal h\in {\cal H},\\\label{homot-2}
&&\|\phi(g)v-v\phi(g)\|<\dt\tforal g\in {\cal G}\tand\\\label{homot-2+}
&&{\rm  Bott}(\phi,v)|_{\cal P}=\{0\},
\eneq
then there exists a continuous path of unitaries
$\{u_t: t\in [0,1]\}\subset M_k$ such that
\beq\label{homot-3}
u_0=v,\,\,u_1=1,\tand \|\phi(f)u_t-u_t\phi(f)\|<\ep
\eneq
for all $t\in [0,1]$ and $f\in {\cal F}.$ Moreover,
\beq\label{homot-4}
{\rm length}(\{u_t\})\le 2\pi+\ep.
\eneq

\end{lem}

\begin{proof}
Let $\Delta_1$ be as in \ref{Ddel}.
Let ${\cal H}_1\subset A_+^{\bf 1}\setminus \{0\}$ and ${\cal H}_2\subset C(\T)_+^{\bf 1}\setminus \{0\}$ be finite subsets,
${\cal G}_1\subset A$ (in place of ${\cal G}$) be a finite subset, $\dt_1>0$ (in place of $\dt$) and  ${\cal P}\subset \underline{K}(A)$ be a finite subset required by \ref{homfull} for $\ep/4$ (in place of $\ep$),
 ${\cal F}$ and $\Delta_1.$

%Let $1>\sigma>0$ and let $\tau$ be a faithful tracial state given by \ref{MExtmeasure} for $r=3/4.$
%There exists a finite subset ${\cal G}_1\subset A$ (in place of ${\cal G}$),
%$\dt_1>0$ (in place of $\dt$)
 %a finite subset ${\cal H}_3\subset C(\T)_+\setminus \{0\}$
% and a finite subset ${\cal P}\subset \underline{K}(A)$ required by \ref{homfull} for $\ep/4$ (in place of $\ep$), %${\cal F},$  above ${\cal H}_1,$ ${\cal H}_2$ and $\Delta_1.$

  Let ${\cal G}_2\subset A$ (in place of ${\cal G}$) be a finite subset, ${\cal H}_1'\subset A_+\setminus \{0\}$ be a finite subset, $\dt_2>0$ (in place of $\dt$) be  as required by
\ref{changespectrum} for $\min\{\ep/16, \dt_1/2\}$ (in place of $\ep$),
${\cal G}_1\cup {\cal F}$ (in place of ${\cal F}$) and ${\cal H}_1$ and ${\cal H}_2.$

Let ${\cal G}={\cal G}_2\cup {\cal G}_1\subset {\cal F}$ and
let $\dt=\min\{\dt_2, \ep/16\}.$  Let ${\cal H}={\cal H}_1.$

Now suppose that $\phi: A\to M_k$ is a unital ${\cal G}$-$\dt$-multiplicative \morp\, and $u\in M_k$ is a unitary which satisfy the assumption for the above ${\cal H},$ $\dt,$
${\cal G}$ and  ${\cal P}.$

By applying \ref{changespectrum}, one obtains a continuous path of unitaries $\{u_t: t\in [0,1/2]\}\subset M_k$ such that
\beq\label{homot-5}
u_0=u,\,\,\, u_1=w,\,\,\,\|u_t\phi(g)-\phi(g)u_t\|<\min\{\dt_1, \ep/4\}
\eneq
for all $g\in {\cal G}_1\cup {\cal F}$ and $t\in [0,1/2].$
Moreover, there is a unital \morp\, $L: A\otimes C(\T)\to M_k$ such that
\beq\label{homot-6}
&&\|L(g\otimes 1)-\phi(g)\|<\min\{\dt_1, \ep/4\}\tforal g\in {\cal G}_1\cup {\cal F},\\\label{homot-7}
&&\|L(1\otimes z)-w\|<\min\{\dt_1, \ep/4\}\\\label{homot-8}
&&\andeqn {\rm tr}\circ L(h_1\otimes h_2)\ge \Delta(h_1)\tau_m(h_2)/4
\eneq
for all $h_1\in {\cal H}_1$ and $h_2\in {\cal H}_2.$ Furthermore,
\beq\label{homot-12}
{\rm length}(\{u_t:t\in [0,1/2]\})\le \pi+\ep/4.
\eneq

Note that
\beq\label{homot-13}
[L]|_{\bt({\cal P})}={\rm Bott}(\phi, w)|_{\cal P}=
{\rm Bott}(\phi, u)|_{\cal P}=0.
\eneq

By (\ref{homot-6}), (\ref{homot-7}), (\ref{homot-13}) and (\ref{homot-8}),
applying \ref{homfull},
one obtains a continuous path of unitaries
$\{u_t\in [1/2, 1]\}\subset M_k$ such that
\beq\label{homot-14}
u_{1/2}=w,\,\,\,u_1=1,\,\,\,
\|u_t\phi(f)-\phi(f)u_t\|<\ep/4\tforal f\in {\cal F}\\
\andeqn {\rm length}(\{u_t: t\in [1/2, 1]\})\le \pi+\ep/4.
\eneq
Therefore $\{u_t: t\in [0,1]\}\subset M_k$ is a continuous path of unitaries in $M_k$ with $u_0=u$ and $u_1=1$ such that
\beq\label{homot-15}
\|u_t\phi(f)-\phi(f)u_t\|<\ep\tforal f\in {\cal F}\\
\andeqn
{\rm length}(\{u_t: t\in [0,1]\})\le 2\pi+\ep.
\eneq
\end{proof}

\section{An Existence Theorem for  Bott maps }

\begin{lem}\label{Ext1}
Let $A$ be a unital {{amenable}} residually finite dimensional  C*-algebra which satisfies the UCT, % \CA\, whose ranks of irreducible representations are bounded,
let
$G=\Z^r\oplus {\rm Tor}(G)\subset K_0(A)$ be a finitely generated subgroup with $[1_A]\in G$ and
let
$J_0, J_1\ge 0$ be integers.

For any $\dt>0,$ any finite subset ${\cal G}\subset A$ and any finite subset ${\cal P}\subset \underline{K}(A)$ with ${\cal P}\cap K_0(A)\subset G,$
there exist integers $N_0, N_1,...,N_k$ and
unital \hm s $h_j: A\to M_{N_j},$ $j=1,2,...,k,$ satisfying the following:

For any $\kappa\in Hom_{\Lambda}(\underline{K}(A), \underline{K}({\cal K})),$
with $|\kappa([1_A])|=J_1$ and
\beq\label{Ext1-1}
J_0=\max\{|\kappa(g_i)|: g_i=(\overbrace{0,...,0}^{i-1}, 1, 0,...,0)\in \Z^r: 1\le i\le r\},
\eneq
there exists a ${\cal G}$-$\dt$-multiplicative \morp\, $\Phi: A\to M_{N_0+\kappa([1_A])}$
such that
\beq\label{Ext1-2}
[\Phi]|_{\cal P}=(\kappa+[h_1]+[h_2]+\cdots +[h_k])|_{\cal P}.
\eneq

(Note that $N_0=\sum_{i=1}^kN_i$).

\end{lem}

\begin{proof}
It follows from 6.1.11 of \cite{Lnbok} (see also \cite{LinTAF2} and \cite{DE}) that, for each such $\kappa,$
there is a unital ${\cal G}$-$\dt$-multiplicative \morp\, $L_{\kappa}: A\to M_{n(\kappa)}$ (for integer $n(\kappa)\ge 1$) such that
\beq\label{Ext1-3}
[L_{\kappa}]|_{\cal P}=(\kappa+[h_{\kappa}])|_{\cal P},
\eneq
where $h_{\kappa}: A\to M_{N_{\kappa}}$ is a unital \hm.
There are only finitely many different $\kappa|_{\cal P}$ so that
(\ref{Ext1-1}) holds, say $\kappa_1, \kappa_2,...,\kappa_k.$
Set $h_i=h_{\kappa_i},$ $i=1,2,...,k.$ Let $N_i=N_{\kappa_i},$ $i=1,2,...k.$
Note that $N_i=J_1+n(\kappa_i),$ if $\kappa([1_A])=J_1,$
and $N_i=-J_1+n(\kappa_i),$ if $\kappa([1_A])=-J_1.$
Define
$$
N_0=\sum_{i=1}^k N_i.
$$
If $\kappa=\kappa_i,$ define $\Phi: A\to M_{N_0+\kappa([1_A])}$ by
$$
\Phi=L_{\kappa_i}+\sum_{j\not=i}h_j.
$$
The lemma follows.
\end{proof}

\begin{lem}\label{Ext2}
Let $A$ be a unital \CA\, as in \ref{Ext1} and let $[1_A]\in G=\Z^r\oplus {\rm Tor}(G)\subset K_0(A)$ be a finitely generated subgroup. There exists $\Lambda_i\ge 0,$ $i=1,2,...,r,$ satisfying the following:
For any $\dt>0,$ any finite subset ${\cal G}\subset A$ and any finite subset ${\cal P}\subset \underline{K}(A)$ with
${\cal P}\cap K_0(A)\subset G,$  there exist integers
$N(\dt, {\cal G}, {\cal P}, i)\ge 1,$  $i=1,2,...,r,$ satisfying the following:

Let $\kappa\in Hom_{\Lambda}(\underline{K}(A), \underline{K}({\cal K}))$
and
$
S_i=\kappa(g_i),
$
where $g_i=(\overbrace{0,...,0}^{i-1},1,0,...,0)\in \Z^r,$ there exists
a unital ${\cal G}$-$\dt$-multiplicative \morp\, $L: A\to M_{N_1}$ and
a \hm\, $h: A\to M_{N_1}$ such that
\beq\label{EXt2-1}
[L]|_{\cal P}=(\kappa+[h])|_{\cal P},
\eneq
where $N_1=\sum_{i=1}^r (N(\dt, {\cal G},{\cal P}, i)\pm\Lambda_i)\cdot |S_i.|$
\end{lem}

\begin{proof}
Let $\psi_i^{+}: G\to \Z$ be a \hm\, defined by
$\psi_i^{+}(g_i)=1,$ $\psi_i^{+}(g_j)=0,$ if $j\not=i,$ and
 $\psi_i^{+}|_{{\rm Tor}(G)}=0,$ and let $\psi_i^{-}(g_i)=-1$ and $\psi_i^{-}(g_j)=0,$ if $j\not=i,$  and $\psi_i^{-}|_{{\rm Tor}(G)}=0,$  $i=1,2,...,r.$
 Note that $\psi_i^-=-\psi_i^+,$ $i=1,2,...,r.$
 Let $\Lambda_i=|\psi_i^+([1_A])|,$ $i=1,2,...,r.$

Let $\kappa_i^{+}, \kappa_i^-\in \mathrm{Hom}_{\Lambda}(\underline{K}(A), \underline{K}({\cal K}){{)}}$ be such that
$\kappa_i^+|_G=\psi_i^+$ and $\kappa_i^-=\psi_i^-,$ $i=1,2,...,r.$
Let $N_0(i)\ge 1$ (in place of $N_0$) be required by \ref{Ext1} for
$\dt,$ ${\cal G},$ $J_0=1$ and $J_1=M_i.$
Define
$N(\dt, {\cal G},{\cal P}, i)=N_0(i),$ $i=1,2,...,r.$

Let $\kappa\in \mathrm{Hom}_{\Lambda}(\underline{K}(A), \underline{K}({\cal K})).$
Then
$\kappa|_G=\sum_{i=1}^rS_i\psi_i^+,$
where $S_i=\kappa(g_i),$ $i=1,2,...,r.$
%Therefore, without loss of generality, we may assume
%that
%$\kappa=S_i\kappa_i^+$ for some $\kappa_i^+$ $(i=1,2,...,r$) as denoted above.
By applying \ref{Ext1}, one obtains ${\cal G}$-$\dt$-multiplicative
\morp s $L_i^{\pm}: A\to M_{N_0(i)+\kappa_i^{\pm}([1_A])}$  and a \hm\, $h_i^{\pm}: A\to M_{N_0(i)}$ such that
\beq\label{Ext2-3}
[L_i^{\pm}]|_{\cal P}=(\kappa_i^{\pm}+[h_i^{\pm}])|_{\cal P},\,\,\,i=1,2,...,r.
\eneq
Define
$L=\sum_{i=1}^r L_i^{\pm,|S_i|},$
where $L^{\pm, |S_i|}: A\to M_{|S_i|N_0(i)}$ is defined by
$$L^{\pm, |S_i|}(a)={\rm diag}(\overbrace{L_i^{\pm}(a), ...,L_i^{\pm}(a)}^{|S_i|})$$ for
all $a\in A.$
One checks that
$L: A\to M_{N_1},$ where $N_1=\sum_{i=1}^r|S_i|(\Lambda_i'+N(\dt, {\cal G}, {\cal P}, i))$ with $\Lambda_i'=\psi^+_i([1_A])$ if $S_i>0,$ or $\Lambda_i'=-\psi_+^+([1_A])$ if $S_i<0,$ is a
unital ${\cal G}$-$\dt$-multiplicative \morp\, and
$$
[L]|_{\cal P}=(\kappa+[h])|_{\cal P}
$$
for some \hm\, $h: A\to M_{N_1}.$

\end{proof}

\begin{rem}
Note that in the above proof, $h$ may be written as $\sum_{i=1}^r h_i^{\pm, |S_i|}.$
Let us agree that $h_i^{\pm, 0}$ is zero.
Let $S\ge 1$ be  a fixe integer.
%Choose new $N_1$ to be $\sum_{i=1}^r (N(\dt, {\cal G},{\cal P}, i)\pm\Lambda_i)\cdot S^
Let us assume that $\max\{|S_i|\}\le  S$ as a restriction on $\kappa.$
Choose new $N_1$ to be $\sum_{i=1}^r (N(\dt, {\cal G},{\cal P}, i)\pm\Lambda_i)\cdot S^r.$
Let $h_0=\sum_{0\le s_{i}\le S}\sum_{i=1}^r h_i^{\pm, s_{i}}: A\to M_{N_1}$ Then $h_0$ is independent
of $\kappa.$
Choose $\Phi=L\oplus \sum_{s_i\not |S_i|} h_i^{\pm, s_i}.$ Then $\Phi: A\to M_{N_1}$ and
 $[\Phi]|_{\cal P}=(\kappa+ [h_0])_{\cal P}.$ In other words, $h_0$ can be chosen before $\kappa$ is given
 as long as $\kappa$ has the restriction as mentioned above (as in  \ref{Ext1}).

\end{rem}

\begin{lem}\label{Ext3}
%{\color{green} Let $A=A_0,$ or $A=A_0\otimes C(\T)$ for some $A_0\in {\cal C}_0$, or let $A=C(X)$ for a compact metric space $X$}  such that $K_i(A)$
Let $A\in {\overline{{\cal D}}}_s$ be a unital \CA\, and let ${\cal P}\subset \underline{K}(A)$ be a finite subset.
Suppose that $G\subset \underline{K}(A)$ is the group generated by ${\cal P}$, and
$G_1=G\cap K_1(A)=
%is finitely generated ($i=0,1$) so that
%$K_1(A)=\
Z^r\oplus {\rm Tor}(K_1(A)).$
%Suppose that $A$ has the property (P).
Let ${\cal F}\subset A,$ let
$\ep>0$, and let $\Delta: A_+^{q, {\bf 1}}\setminus \{0\}\to (0,1)$ be an order preserving map.

%For any finite subset ${\cal P}\subset \underline{K}(A),$ t
There exist
$\dt>0,$  a finite subset ${\cal G}\subset A,$ a
 finite subset ${\cal H}\subset A_+^{{{\bf 1}}}\setminus \{0\}$, and an integer $N\ge 1$ satisfying the following: Let $\kappa\in KK(A\otimes C(\T), \mathbb C)$ and put
\begin{equation}\label{Ext3-2}
K=\max\{|\kappa(\boldsymbol{\bt}(g_i))|:1\le i\le r\},
\end{equation}
where $g_i=(\overbrace{0,...,0}^{i-1},1,0,...,0)\in \Z^r.$ Then for any unital ${\cal G}$-$\dt$-multiplicative
\morp\, $\phi: A\to M_R$ such that  $R\ge N(K+1)$ and
\begin{equation}\label{Ext3-1}
{\rm tr}\circ \phi(h)\ge \Delta(\hat{h})\tforal h\in {\cal H},
\end{equation}
% and for any
%$\kappa\in KK(A\otimes C(\T), M_R)$ with
%\beq\label{Ext3-2}
%K=\max\{|\kappa(\boldsymbol{\bt}(g_i))|:1\le i\le r\},
%\eneq
%where $g_i=(\overbrace{0,...,0}^{i-1},1,0,...,0)\in \Z^r,$
there exists a
unitary $u\in M_R$ such that
\beq\label{EXt3-3}
&&\|[\phi(f),\, u]\|<\ep\tforal f\in {\cal F}
\andeqn\\
&&{\rm Bott}(\phi,\, u)|_{\cal P}=\kappa\circ {\boldsymbol{\bt}}|_{\cal P}.
\eneq
\end{lem}

\begin{proof}

%{\color{green} Let $A=A_0,$ or $A=A_0\otimes C(\T)$ for some $A_0\in {\cal C}_0$.}

%Let $N(p)\ge 1$ be the integer given by the property (P) (see \ref{prime}).
To simplify notation, without loss of generality,  we may assume that ${\cal F}$ is a subset of the unit ball.
Let $\Delta_1=(1/8)\Delta$ and $\Delta_2=(1/16)\Delta.$

Let $\ep_0>0$ and let ${\cal G}_0\subset A$ be a finite subset satisfy the following:
If $\phi': A\to B$ (for any unital \CA\, $B$) is a unital ${\cal G}_0$-$\ep_0$-multiplicative \morp\, and $u'\in B$ is a unitary such that
\beq\label{Ext3-n+1}
\|\phi'(g)u'-u'\phi'(g)\|<4\ep_0\tforal g\in {\cal G}_0,
\eneq
then ${\rm Bott}(\phi',\, u')|_{\cal P}$ is well defined. Moreover,
if $\phi': A\to B$ is another unital ${\cal G}_0$-$\ep_0$-multiplicative \morp\,, then
\beq\label{Ext3-n-2}
{\rm Bott}(\phi', \, u')|_{\cal P}={\rm Bott}(\phi'', \, u'')|_{\cal P},
\eneq
provided that
\beq\label{Ext3-n-3}
\|\phi'(g)-\phi''(g)\|<4\ep_0\andeqn \|u'-u''\|<4\ep_0\tforal g\in {\cal G}_0.
\eneq
We may assume that $1_A\in {\cal G}_0.$
Let
$$
{\cal G}_0'=\{g\otimes f: g\in {\cal G}_0\} \andeqn f=\{1_{C(\T)}, z, z^*\},
$$
where $z$ is the identity function on the unit circle $\T.$
We also assume that if $\Psi': A\otimes C(\T)\to C$ for a unital \CA\, $C$
is a ${\cal G}_0'$-$\ep_0$-multiplicative \morp, then there exists a unitary $u'\in C$ such that
\beq\label{Ext3-n-4}
\|\Psi'(1\otimes z)-u'\|<4\ep_0.
\eneq

Without loss of generality,  we may assume that
${\cal G}_0$  is in the unital ball of $A.$
Let $\ep_1=\min\{\ep/64, \ep_0/512\}$ and ${\cal F}_1={\cal F}\cup {\cal G}_0.$

Let ${\cal H}_0\subset A_+^{\bf 1}\setminus \{0\}$ (in place of ${\cal H}$) be a finite subset and  $L\ge 1$ be an integer {{as}} required by \ref{fullabs} for $\ep_1$ (in place of $\ep$) and ${\cal F}_1$ (in place of ${\cal F}$) as well as $\Delta_2$ (in place of $\Delta$).

Let ${\cal H}_1\subset A_+^{\bf 1}\setminus\{0\}$,
${\cal G}_1\subset A$ (in place of ${\cal G}$),
$\dt_1>0$ (in place of $\dt$), ${\cal P}_1\subset  \underline{K}(A)$
(in place of ${\cal P}$), ${\cal H}_2\subset A_{s.a.}$, and $1>\sigma>0$ be required by \ref{Newunique1}
for $\ep_1$ (in place of $\ep$),  ${\cal F}_1$ (in place of ${\cal F}$)
and $\Delta_1.$  We may assume that $[1_A]\in {\cal P}_2,$ ${\cal H}_2$ is in the unit ball of $A$
and ${\cal H}_0\subset {\cal H}_1.$

%We may also assume that ${\cal H}_1$ is required by \ref{Newuni2} for $\ep/16$ and ${\cal F},$
%${\cal G}_1,$ $\dt_1$ and ${\cal H}_2$ are required by \ref{Newunique1}
%for $\ep/16,$ ${\cal F},$ ${\cal H}_1$ and $\sigma_1.$

%Let ${\cal H}_3=\{h_+, h_-: h\in {\cal H}_2.$

Without loss of generality, we may assume that
$\dt_1<\ep_1/16, \sigma<\ep_1/16$, and ${\cal F}_1\subset {\cal G}_1.$
Put ${\cal P}_2={\cal P}\cup {\cal P}_1$ and ${\cal H}_1'={\cal H}_1\cup {\cal H}_0.$

Suppose that $A$ has irreducible representations of rank
$r_1,r_2,...,r_k.$
%Put $d=\sum_{i=1}^k r_i.$
Fix an irreducible representation
$\pi_0: A\to M_{r_1}.$
Let $N(p)\ge 1$ (in place of $N({\cal P}_0))$ and ${\cal H}_0'\subset A_+^{\bf 1}\setminus \{0\}$ (in place of ${\cal H}$) be a finite subset
required by \ref{Lfullab} for  $\{1_A\}$ (in place of ${\cal P}_0$) and $(1/16)\Delta.$ Let ${\cal H}_1'={\cal H}_1\cup{\cal H}_0'.$

 Let $G_0=G\cap K_0(A)$ and
write $G_0=\Z^{s_1}\oplus \Z^{s_2}\oplus {\rm Tor}(G_0),$ where $\Z^{s_2}\oplus {\rm Tor}(G_0)
\subset {\rm ker}\rho_A.$
Let $x_j=(\overbrace{0,...,0}^{j-1},1,0,...,0)\in \Z^{s_1}\oplus \Z^{s_2},$ $j=1,2,...,s_2.$ Note that
$A\otimes C(\T)\in {\overline{\cal D}}_s$ and
$A\otimes C(\T)$ has irreducible representations of rank
$r_1, r_2,...,r_k.$
Let
$$
{\bar r}=\max \{ |(\pi_0)_{*0}(x_j)|: 0\le j\le s_1+s_2\}.
$$

Let ${\cal P}_3\subset \underline{K}(A\otimes C(\T))$ be a finite subset
set containing ${\cal P}_2,$ $\{{\boldsymbol{\bt}}(g_j): 1\le j\le r\}$ and a finite subset which generates
${\boldsymbol{\bt}}({\rm Tor}(G_1)).$  Choose $\dt_2>0$ and finite subset
$$
\overline{{\cal G}}=\{g\otimes f: g\in {\cal G}_2,\,\,\, f\in \{1, z, z^*\}\}
$$
in $A\otimes C(\T),$ where
${\cal G}_2\subset A$ is a finite subset
% and $z\in C(\T)$ is the identity function on the unit circle,
such that, for any unital $\overline{{\cal G}}$-$\dt_2$-multiplicative \morp\,
$\Phi': A\otimes C(\T)\to C$ (for any unital \CA\, $C$ with $T(C)\not=\emptyset$),
$[\Phi']|_{{\cal P}_3}$ is well defined and
%\beq\label{Ext3-5}
%|\tau([\Phi](x_j))|<\sigma_2/16\andeqn |\tau([\Phi](\boldsymbol{\bt}(g_i)))|<\sigma_2/16,
%\eneq
%$j=1,2,...,s_2$ and $i=1,2,...,r.$ Moreover, by choosing even smaller
%$\dt_2$ and larger ${\cal G}_2,$ we may further assume that
\beq\label{Ext3-6}
[\Phi']|_{{\rm Tor}(G_0)\oplus \boldsymbol{\bt}({\rm Tor}(G_1)}=0.
\eneq
We may assume ${\cal G}_2\supset {\cal G}_1\cup {\cal F}_1.$

Let $\sigma_1=\min\{\Delta_2(\hat{h}): h\in {\cal H}_1'\}.$
Note $K_0(A\otimes C(\T))=K_0(A)\oplus \boldsymbol{\bt}(K_1(A))$ and
$\underline{K}(A\otimes C(\T))=\underline{K}(A)\oplus \boldsymbol{\bt}(\underline{K}(A)).$
Consider the subgroup of
$K_0(A\otimes C(\T))$ given by
$$
\Z^{s_1}\oplus \Z^{s_2}\oplus \Z^r\oplus {\rm Tor}(K_0(A)\oplus {\boldsymbol{\bt}}({\rm Tor}(K_1(A)).
$$

Let $\dt_3=\min\{\dt_1, \dt_2\}.$ Let $N(\dt_3, \overline{{\cal G}}, {\cal P}_3,i)$ and $\Lambda_i,$
$i=1,2,...,s_1+s_2+r,$ be required by \ref{Ext2} for $A\otimes C(\T)$.
Choose an integer $n_1\ge N(p)$ such that
\beq\label{Ext3-7}
{(\sum_{i=1}^{s_1+s_2+r} N(\dt_3, \overline{{\cal G}}, {\cal P}_3,i)+1+\Lambda_i)N(p)\over{n_1-1}}<\min\{\sigma/16, \sigma_1/2\}.
\eneq
Choose $n>n_1$ such that
\beq\label{Ext3-7+1}
{n_1+2\over{n}}<\min\{\sigma/16, \sigma_1/2, 1/(L+1)\}.
\eneq

Let $\ep_2>0$ and let ${\cal F}_2\subset A$ be a finite subset
such that
$[\Psi]|_{{\cal P}_2}$ is well defined.

%$\ep_2<\ep/16,$ ${\cal F}\subset {\cal F}_1$ and
%$(\ep_2, {\cal F}_1)$ is a KK-triple.

Let $\ep_3=\min\{\ep_2/2, \ep_1\}$ and ${\cal F}_3={\cal F}_1\cup {\cal F}_2.$

Denote by $\dt_4>0$ (in place of $\dt$), ${\cal G}_3\subset A$
(in place of ${\cal G}$), ${\cal H}_3\subset A_+^{\bf 1}\setminus\{0\}$ (in place of ${\cal H}_2$) the finite sets
{{as}} required by \ref{Combinerep} for $\ep_3$ (in place of $\ep$),
${\cal F}_3\cup {\cal H}_1'$ (in place of ${\cal F}$), $\dt_3/2$ (in place of $\ep_0$),
${\cal G}_2$ (in place of ${\cal G}_0$),
%(in place of ${\cal F}$),
$\Delta,$ ${\cal H}_1'$ (in place of ${\cal H}$),
$\min\{\sigma/16, \sigma_1/2\}$  (in place of $\sigma$) and $n^2$ (in place of $K$) required by
\ref{Combinerep} (with $L_1=L_2$).
%Without loss of generality, we may assume that
%$$
%{\cal G}_3=\{g\otimes f: g\in {\cal G}_4\andeqn  f\in \{1, z, z^*\}\},
%$$
%where ${\cal G}_4\subset A$ is a finite subset.

Set ${\cal G}={\cal F}_3\cup {\cal G}_1\cup {\cal G}_2\cup {\cal G}_3$, set  $\dt=\min\{\ep_3/16, \dt_4, \dt_3/16\}{{,}}$
and
set ${\cal G}_5=\{g\otimes f: g\in {\cal G}_4,\,\,\,f\in\{1, z, z^*\}\}.$
%Without loss of generality, we may assume
%that
%if $L': A\otimes C(\T)\to C$ is a unital $\dt_4$-${\cal G}_4$-multiplicative \morp\, (for  any unital \CA\, $C$),
%then there exists a unitary $u'\in C$ such that
%$$
%\|L'(1\otimes z)-u'\|<\min\{\ep/16, 1/16\}.
%$$

%Choose $\dt>0$ and a finite subset ${\cal G}\subset A$ satisfying the following: For any unital $\dt$-${\cal G}$-multiplicative \morp\,
%$L: A\to C$ (for any unital \CA\, $C$) and any unitary $u'\in C$ with
%$\|L(g),\, u']\|<\dt$ for all $g\in {\cal G},$ there is a unital $\dt_4$-${\cal G}_5$-multiplicative \morp\,
%$L': A\otimes C(\T)\to C$ such that
%$$
%\|L'(g\otimes 1)-L(g)\|<\dt_5/2\tforal g\in {\cal G}_4\andeqn
%\|L'(1\otimes z)-u'\|<\dt_5/2.
%$$
Let ${\cal H}={\cal H}_1'\cup {\cal H}_3.$  Define $N_0=(n+1)N(p)(\sum_{i=1}^{s_1+s_2+r}N(\dt_3, {\cal G}_0, {\cal P}_3,i)+\Lambda_i+1)$ and define $N=N_0+N_0{\bar r}.$
Fix any $\kappa\in KK(A\otimes C(\T), \C)$
with
$$
K=\max\{|\kappa({\boldsymbol{\bt}}(g_i)|: 1\le j\le r\}.
$$
Let $R>N(K+1).$

Suppose
that $\phi:A\to M_R$ is a unital ${\cal G}$-$\dt$-multiplicative
\morp\, such that
\beq\label{Ext3-8}
{\rm tr}\circ \phi(h)\ge \Delta(\hat{h})\tforal h\in {\cal H}.
\eneq
Then, by \ref{Combinerep}, there exist mutually orthogonal projections
$e_0, e_1,e_2,...,e_{n^2}\in M_R$ such that $e_1, e_2,...,e_{n^2}$ are equivalent,
${\rm tr}(e_0)<\min\{\sigma/64, \sigma_1/4\}$  and $e_0+\sum_{i=1}^{n^2}e_i=1_{M_R},$ and there exists a unital ${\cal G}_2$-$\dt_3/2$-multiplicative \morp\, $\psi_0: A\to e_0M_Re_0$ and a unital \hm\, $\psi: A\to e_1M_Re_1$ such that
\beq\label{Ext3-9}
\|\phi(f)-(\psi_0(f)\oplus \overbrace{\psi(f),\psi(f),...,\psi(f)}^{n^2})\|<\ep_3\tforal f\in {\cal F}_3\andeqn\\
{\rm tr}\circ \psi(h)\ge \Delta(\hat{h})/3n^2\tforal h\in {\cal H}_1.
\eneq
Let $\alpha\in \mathrm{Hom}_{\Lambda}(\underline{K}(A\otimes C(\T)), \underline{K}(M_r))$ be defined as follows:
$\alpha|_{\underline{K}(A)}=[\pi_0]$ and
$\alpha|_{\boldsymbol{\bt}(\underline{K}(A))}
=\kappa|_{\boldsymbol{\bt}(\underline{K}(A))}.$
Let
$$
\max\{|\kappa\circ \boldsymbol{\bt}(g_i)|: i=1,2,...,r, |\pi_0(x_j)|:1\le j\le s_1+s_2\}\le \max\{K, {\bar r}\}.
$$

Applying \ref{Ext2}, we obtain a unital ${\cal G}$-$\dt_3$-multiplicative \morp\,
$\Psi: A\otimes C(\T)\to M_{N_1'},$ where $N_1'\le N_1=\sum_{j=1}^{s_1+s_2+r} N(\dt_3, {\cal G}_0, {\cal P}_3,j)+\Lambda_i)\max\{K, {\bar r}\},$
and a \hm\, $H_0: A\otimes C(\T)\to H_0(1_A)M_{N_1'}H_0(1_A)$ such that
such that
\beq\label{Ext3-10-}
[\Psi]|_{{\cal P}_3}=(\af+[H_0])|_{{\cal P}_3}.
\eneq
In particular, since $[1_A]\in {\cal P}_2\subset {\cal P}_3,$
%\begin{equation*}
${\rm rank} \Psi(1_A)=r_1+{\rm rank} (H_0).$
%\end{equation*}
Note that
\beq\label{Ext3-10}
{N_1'+N(p) \over{R}}\le {N_1+N(p)\over{N(K+1)}}<1/(n+1).
\eneq
Let $R_1={\rm rank}\,e_1.$
Then $R_1\ge R/(n+1).$  So, from (\ref{Ext3-10}), one has that $R_1\ge N_1+N(p).$
In other words, $R_1-N_1'\ge N(p).$
%Since $A$ has property (P), there are non-negative integers $m_1,...,m_k$ such that
%$\sum_{i=1}^k m_ir_i=R_1-N_1.$
Note that
\begin{equation*}
t\circ  \psi(\hat{g})\ge (1/3)\Delta(\hat{g})\ge \Delta_2(\hat{g})\rforal g\in {\cal H}_0',
\end{equation*}
%where ${\tilde \psi}(a)=\diag(\phi(a),\phi(a),...,\phi(a)),$ where $\phi(a)$ repeats $n^2$-times and
where $t$ is the tracial state on $M_{R_1}.$
By applying \ref{Lfullab} to the case that $\pi_0\oplus H_0$ (in place $\phi$), ${\tilde \psi},$ 
where ${\tilde \psi}$ is an amplification of $\psi$ with $\psi$ repeated $n$ times, and ${\cal P}_0=\{[1_A]\},$
%$R_1=\sum_{i=1}^kT_ir_i$ with $T_i\ge 1.$
%Therefore $nT_i\ge N_1'.$
%It follows that
we obtain  a unital \hm\, $h_0: A\otimes C(\T)\to M_{nR_1-N_1'} .$
Define $\psi_0': A\otimes C(\T)\to e_0M_Re_0$ by
$\psi_0'(a\otimes f)=\psi_0(a)\cdot f(1)\cdot e_0$ for all $a\in A$ and $f\in C(\T),$ where $1\in  \T.$
Define $\psi': A\otimes C(\T)\to e_1M_Re_1$ by
$\psi'(a\otimes f)=\psi(a)\cdot f(1)\cdot e_0$ for all $a\in A$ and $f\in C(\T).$
Put $E_1={\rm diag}(e_1,e_2,...,e_{nn_1}).$

Define $L_1: A\to E_1M_RE_1$ by
$L_1(a)=\pi_0(a)\oplus H_0|_{A}(a)\oplus h_0(a\otimes 1)\oplus (\overbrace{\psi(f),...,\psi(f)}^{n(n_1-1)})$ for $a\in A$, and define
$L_2: A\to E_1M_RE_1$ by
$L_2(a) =\Psi(a\otimes 1)\oplus h_0(a\otimes 1)\oplus (\overbrace{\psi(f),...,\psi(f)}^{n(n_1-1)})$ for $a\in A.$
Note that
\beq\label{Ext3-12}
&&[L_1]|_{{\cal P}_1}=[L_2]|_{{\cal P}_1},\\
&&{\rm tr}\circ L_1(h)\ge \Delta_1(\hat{h}),\,\,\,{\rm tr}\circ L_2(h)\ge \Delta_1(\hat{h})\tforal h\in {\cal H}_1, \andeqn\\
&&|{\rm tr}\circ L_1(g)-{\rm tr}\circ L_2(g)|<\sigma \tforal g\in {\cal H}_2.
\eneq
It follows from \ref{Newunique1} that
there exists a unitary $w_1\in E_1M_RE_1$ such that
\beq\label{Ext3-13}
\|{\rm ad}\, w_1\circ L_2(a)-L_1(a)\|<\ep_1\tforal a\in {\cal F}_1.
\eneq
Define $E_2=(e_1+e_2+\cdots +e_{n^2})$ and define
$\Phi: A\to E_2M_RE_2$ by
\beq\label{Ext3-12+n-1}
\Phi(f)(a)={\rm diag}(\overbrace{\psi(a),\psi(a),...,\psi(a)}^{n^2})\tforal a\in A.
\eneq
Then
\beq\label{EXt3-12+n-2}
{\rm tr}\circ \Phi(h)\ge \Delta_2(\hat{h})\tforal h\in {\cal H}_0.
\eneq
By (\ref{Ext3-7+1}), one has that ${n\over{n_1+2}}>L+1.$ By applying \ref{fullabs}, we obtain a unitary $w_2\in E_2M_RE_2$ and
a unital \hm\, $H_1: A\to (E_2-E_1)M_R(E_2-E_1)$ such that
\beq\label{Ext3-14}
\|{\rm ad}\, w_2\circ {\rm diag}(L_1(a), H_1(a))-\Phi(a)\|<\ep_1\tforal a\in {\cal F}_1.
\eneq
Put
$$
w=(e_0\oplus w_1\oplus (E_2-E_1))(e_0\oplus w_2)\in M_R.
$$
Define $H_1': A\otimes C(\T)\to (E_2-E_1)M_R(E_2-E_1)$ by
$H_1'(a\otimes f)=H_1(a)\cdot f(1)\cdot (E_2-E_1)$ for all $a\in A$ and
$f\in C(\T).$
Define $\Psi_1: A\to M_R$ by
\beq\label{Ext3-11}
\Psi_1(f)=\psi_0'(f)\oplus \Psi(f)\oplus h_0\oplus \overbrace{\psi'(f),...,\psi'(f)}^{n_1-1})\oplus H_1'(f)\tforal f\in A\otimes C(\T).
\eneq
It follows from (\ref{Ext3-13}), (\ref{Ext3-14}) and (\ref{Ext3-9}) that
\beq\label{Ext3-15}
\|\phi(a)-w^*\Psi_1(a\otimes 1)w\|<\ep_1+\ep_1+\ep_3\tforal a\in {\cal F}.
\eneq
Now pick a unitary $v\in M_R$ such that
\beq\label{Ext3-15+}
\|\Psi_1(1\otimes z)-v\|<4\ep_1.
\eneq
Put $u=w^*vw.$
Then, we estimate  that
\beq\label{Ext3-16}
\|[\phi(a), \, u]\|<\min\{\ep, \ep_0)\tforal a\in {\cal F}_1.
\eneq
Moreover, by (\ref{Ext3-13}),(\ref{Ext3-10-}) and (\ref{Ext3-n-2}), one has
\beq\label{Ext3-17}
{\rm Bott}(\phi,\, u)|_{{\cal P}}=\kappa\circ \boldsymbol{\bt }|_{\cal P}.
\eneq

%{\color{green} If $A=C(X)$, then the Lemma follows from exact the same argument with \ref{fullabs} being replaced by ?? and \ref{Newunique1} replaced by Theorem 2.10 of \cite{Lin-AU11}.}
\end{proof}

\section{A Uniqueness Theorem for \CA s
%whose irreducible representation have bounded
%dimensions {\color{Green} (
 in $\mathcal D_s$}

The main goal of this section is to prove Theorem \ref{UniqN1}.
\begin{df}\label{1endpoints}
Let $A$ be a unital \CA\, and let $C\in  {\cal C},$
where
$C=C(F_1, F_2, \phi_0, \phi_1)$ is as in \ref{DfC1}.
Suppose that $L: A\to C$ is a \morp.
Define $L_e=\pi_e\circ L.$
Then
 $L_e: A\to F_1$ is a \morp\, such that
\beq\label{1end-1}
\phi_0\circ L_e=\pi_0\circ L\andeqn
\phi_1\circ L_e=\pi_1\circ L.
\eneq
Moreover, if $\dt>0$ and ${\cal G}\subset A$ and $L$ is ${\cal G}$-$\dt$-multiplicative, then $L_e$ is also ${\cal G}$-$\dt$-multiplicative.
\end{df}

\begin{lem}\label{endpoints}
Let $A$ be a unital \CA\, and let $C\in  {\cal C},$
where
$C=C(F_1, F_2, \phi_0, \phi_1)$ is as in  {\rm \ref{DfC1}}.
Let  $L_1, L_2: A\to C$ be  two unital completely positive linear maps, let  $\ep>0$ and let ${\cal F}\subset A$ be a subset.
 Suppose that there are unitaries $w_0\in \pi_0(C)\subset F_2$  and $w_1\in \pi_1(C)\subset F_2$ such that
 \beq\label{endpoints-1}
 &&\|w_0^*\pi_0\circ L_1(a)w_0-\pi_0\circ L_2(a)\|<\ep\andeqn\\
 % \eneq
 % \beq
 \label{endpoints-1+1}
 && \|w_1^*\pi_1\circ L_1(a)w_1-\pi_1\circ L_2(a)\|<\ep \tforal a\in {\cal F}.
 \eneq
 Then there exists a unitary $u\in F_1$ such that
 \beq\label{endpoints-2}
 &&\|\phi_0(u)^*\pi_0\circ L_1(a)\phi_0(u)-\pi_0\circ L_2(a)\|<\ep, \andeqn
 %\eneq
% \beq
\\\label{endpoints-2+1}
&&\|\phi_1(u)^*\pi_1\circ L_1(a)\phi_1(u)-\pi_1\circ L_2(a)\|<\ep\tforal a\in {\cal F}.
\eneq
\end{lem}

\begin{proof}
Write $F_1%=\bigoplus_{i=1}^k F_1^i
=M_{n_1}\oplus M_{n_2}\oplus\cdots \oplus M_{n_k}$ and
$F_2=M_{r_1}\oplus M_{r_2}\oplus\cdots \oplus M_{r_l}.$
We may assume that, ${\rm ker}\phi_0\cap {\rm ker}\phi_1=\{0\}.$ (See \ref{DfC1}.)

We may assume that there are $k(0)$ and $k(1)$ such that $\phi_0|_{M_{n_i}}$ is injective, $i=1,2,...,k(0)$ with $k(0)\le k,$
$\phi_0|_{M_{n_i}}=0$ if $i> k(0),$ and
$\phi_1|_{M_{n_i}}$ is injective, $i=k(1), k(1)+1,...,k$ with $k(1)\le k,$
$\phi_1|_{M_{n_i}}=0,$ if $i<k(1).$
Write $F_{1,0}=\bigoplus_{i=1}^{k(0)}M_{n_i}$ and $F_{1,1}=\bigoplus_{j=k(1)}^{k}M_{n_{j}}.$
Note that $k(1)\le k(0)+1,$ $\phi_0|_{F_{1,0}}$ and $\phi_1|_{F_{1,1}}$ are injective.
Note $\phi_0(F_{1,0})=\phi_0(F_1)=\pi_0(C)$ and $\phi_1(F_{1,1})=\phi_1(F_1)=\pi_1(C).$
Let $\psi_0=(\phi_0|_{F_{1,0}})^{-1}$ and $\psi_1=(\phi_1|_{F_{1,1}})^{-1}.$

%Fix $a\in A.$ Let $g_{a,i}\in F_1$ such that
%$\phi_0(g_{a,i})=\pi_0\circ L_i(a)$ and $\phi_1(g_{a,i})=\pi_1\circ L_i(a),$ $i=1,2.$
%The existence of such $g_{a,i}$ follows from the fact that $L_i(a)\in C.$
%Write
For each fixed $a\in A$, since $L_i(a)\in C~(i=0,1)$, there are elements
$$
g_{a,i}=g_{a,i,1}\oplus g_{a,i,2}\oplus \cd \oplus g_{a,i,k(0)}\oplus\cdots \oplus g_{a,i,k}\in F_1,
$$
such that $\phi_0(g_{a,i})=\pi_0\circ L_i(a)$ and $\phi_1(g_{a,i})=\pi_1\circ L_i(a),$ $i=1,2$,
 where $g_{a,i,j}\in M_{n_j},$ $j=1,2,...,k$ and $i=1,2.$ Note that such $g_{a,i}$ is unique since ${\rm ker}\phi_0\cap {\rm ker}\phi_1=\{0\}.$
 Since $w_0\in \pi_0(C)=\phi_0(F_1)$, there is a unitary
 $$u_0=u_{0,1}\oplus u_{0,2} \oplus \cd \oplus u_{0, k(0)}\oplus \cd \oplus u_{0,k}$$
 such that $\phi_0(u_0)=w_0$. Note that the first $k(0)$ components of $u_0$ is uniquely determined by $w_0$ (since $\phi_0$ is injective on this part) and the components after $k(0)$'th component can be {chosen} arbitrarily (since $\phi_0=0$ on this part). Similarly there exists
 $$u_1=u_{1,1}\oplus u_{1,2} \oplus \cd \oplus u_{1, k(1)}\oplus \cd \oplus u_{1,k}$$
 such that  $\phi_1(u_1)=w_1$

 Now by \ref{endpoints-1} and \ref{endpoints-1+1}, we have
 \beq\label{endpoints-3}
 \|\phi_0(u_0)^*\phi_0(g_{a,1})\phi_0(u_0)-\phi_0(g_{a,2})\|<\ep\andeqn\\
\|\phi_1(u_1)^*\phi_1(g_{a,1})\phi_1(u_1)-\phi_1(g_{a,2}))\|<\ep\tforal a\in {\cal F}.
\eneq

Since $\phi_0$ is injective on $F_1^i$ for  $i\leq k(0)$   and $\phi_1$ is injective on $F_1^i$ for $i>k(0)$ (note that we use $k(1)\leq k(0)+1$),
we have
\beq\label{endpoints-3+1}
 \|(u_{0,i})^*(g_{a,1,i})u_{0,i}-(g_{a,2,i})\|<\ep\quad \tforal i\leq k(0), \andeqn\\
 %\eneq
% \beq
\label{endpoints-3+2}
\|(u_{1,i})^*(g_{a,1,i})u_{1,i}-(g_{a,2,i})\|<\ep\quad \tforal i>k(0), \tforal a\in {\cal F}.
\eneq
Let $u=u_{0,1}\oplus \cd \oplus u_{0,k(0)}\oplus u_{1,k(0)+1} \oplus \cd \oplus u_{1,k} \in F_1$---that is, for the first $k(0)$'s components of $u$, we  use $u_0$'s corresponding components,  and for the last $k-k(0)$ components of $u$, we use $u_1$'s. From \ref{endpoints-3+1} and \ref{endpoints-3+1}.
we have
$$\|u^*g_{a,1}u-g_{a,2}\|<\ep~~\tforal a\in {\cal F}.$$
Apply $\phi_0$ and $\phi_1$ to the above inequality, we get \ref{endpoints-2} and \ref{endpoints-2+1} as desired.
\end{proof}

\begin{lem}\label{CUend}
Let $A$ be a unital \CA\,
%in ${\cal A}_s.$
and
let $C\in  {\cal C},$
where
$C=C(F_1, F_2, \phi_0, \phi_1)$ is as in {\rm  \ref{DfC1}}.
Suppose $L, \Psi: A\to C$ are  two ${\cal G}$-$\dt$-multiplicative
\morp s for some $1>\dt>0$ and a subset ${\cal G}\subset A.$
Suppose that $g\in U(A),$ $1/2>\gamma>0$ and there is $v\in CU(C)$ such that
\beq\label{CUend-1}
\|\langle L (g^*)\rangle \langle \Psi(g)\rangle -v\|<\gamma.
\eneq
Then, there is $v_e\in CU(F_1)$ such that
\beq\label{Cuend-2}
\|\langle L_e(g^*)\rangle \langle \Psi_e(g) \rangle -v_e\|<\gamma,\,\, \pi_0(v_e)=v\andeqn
\phi_1(v_e)=\pi_1(v).
\eneq
\end{lem}

\begin{proof}
The element $v_e=\pi_e(v)$ satisfy the condition of the lemma.
\end{proof}

\begin{thm}\label{UniqN1}
Let $A\in {{{\overline{{\cal D}}}}}_s$ be a unital \CA\, with finitely generated $K_i(A)$ ($i=0,1$).
%$A=A_0,$ or $A=A_0\otimes C(\T)$ for some $A_0\in {\cal C}_0$, or let $A=C(X)$ for a metric space $X$.
%Suppose that $A$ either satisfies the property (P), or $K_1(A)=\{0\}.$
Let ${\cal F}\subset A,$ let
$\ep>0$ be a positive number and let $\Delta: A_+^{q, {\bf 1}}\setminus \{0\}\to (0,1)$  be an order preserving map. There exists  a finite subset ${\cal H}_1\subset A_+^{\bf 1}\setminus \{0\},$
$\gamma_1>0,$ $\gamma_2>0,$ $\dt>0,$ a finite subset
${\cal G}\subset A$ and a finite subset ${\cal P}\subset \underline{K}(A),$ a finite subset ${\cal H}_2\subset A$, a finite subset ${\cal U}\subset J_c(K_1(A))$ {\rm (see \ref{Dcu})} for which $[{\cal U}]\subset {\cal P}$, and $N\in \mathbb N$ satisfying the following:
For any unital ${\cal G}$-$\dt$-multiplicative \morp s $\phi, \psi: A\to C$
for some $C\in {\cal C}$ such that
\begin{equation}\label{Uni1-1}
[\phi]|_{\cal P}=[\psi]|_{\cal P},
\end{equation}
\begin{equation}\label{Uni1-2}
\tau(\phi(a))\ge \Delta(a),\,\,\, \tau(\psi(a))\ge \Delta(a),\quad \textrm{for all $\tau\in T(C)\ a\in {\cal H}_1$},
\end{equation}
\begin{equation}\label{Uni1-3}
|\tau\circ \phi(a)-\tau\circ \psi(a)|<\gamma_1 \tforal a\in {\cal H}_2,
\end{equation}
and
\begin{equation}\label{Uni1-3+1}
{\rm dist}(\phi^{\ddag}(u), \psi^{\ddag}(u))<\gamma_2 \tforal u\in {\cal U},
\end{equation}
there exists a unitary $W\in C\otimes M_{N}$ such that
\begin{equation}\label{Uni1-4}
\|W(\phi(f)\otimes 1_{M_{{N}}})W^*-(\psi(f)\otimes 1_{M_{N}})\|<\ep,\tforal f\in {\cal F}.
\end{equation}
\end{thm}

\begin{proof}
%{\bf To simplify the proof, in what follows, I assume that $K_i(A)$ is finitely generated---we worry the general case later.}
%Without loss of generality, one may assume that $\underline{K}(A)$ is finitely generated (if $A=A_0$ or $A=A_0\otimes\mathrm{C}(\mathbb T)$ for some $A_0\in\mathcal C_0$, then $K_i(A)$ is finitely generated; if $A=C(X)$, then it is an inductive limit of commutative C*-algebras with spectra finite CW-complex).

%There is $n_0$ such that
%$n_0x=0$ for all $x\in K_i(A\otimes C(\T)),$ $i=0,1.$
%%$$\underline{K}(A\otimes C(\mathbb T))\cap K(A\otimes C(\mathbb T), \mathbb Z/n\mathbb Z)=\{0\},\quad n\geq n_0.$$

Since $K_*(A)$ is finitely generated,  there is $n_0$ such that $\kappa\in\mathrm{Hom}_\Lambda(\underline{K}(A), \underline{K}(C))$ is determined by its restriction to $K(A, \Z/n\Z)$, $n=0,..., n_0$.
Set $N=n_0!$.
Let  $\Delta_1$ be defined in \ref{Ddel} for the given $\Delta.$

Let ${\cal H}_1'\subset A_+\setminus \{0\}$ (in place of ${\cal H}_1$)
for $\ep/32$  (in place of $\ep$) and ${\cal F}$ required by \ref{homotopy1}.

%Let $\eta_1=\eta_1'/3.$
 Let
 $\dt_1>0$ (in place of $\dt$), ${\cal G}_1\subset A$ (in place of ${\cal G}$) be a finite subset and let ${\cal P}_0\subset \underline{K}(A)$ (in place of ${\cal P}$) be a finite subset required by \ref{homotopy1}
for $\ep/32$ (in place of $\ep$), ${\cal F}$ and  $\Delta_1.$
We may assume that $\dt_1<\ep/32$ and $(2\dt_1, {\cal G}_1)$
is a $KK$-pair (see the end of \ref{KLtriple}).
%Let ${\bar G}$ be the group generated by ${\cal P}_0.$
%Put $G_1={\bar G}\cap K_1(A).$

Moreover, we may assume that $\dt_1$ is sufficiently small  that
if $\|uv-vu\|<3\dt_1,$ then the Exel formula
$$
\tau({\rm bott}_1(u,v))={1\over{2\pi\sqrt{-1}}}(\tau(\log(u^*vuv^*))
$$
holds for any pair of unitaries $u$ and $v$ in any unital \CA\, $C$ with tracial rank zero and any $\tau\in T(C)$ (see Theorem 3.6 of \cite{Linajm}). Moreover
if $\|v_1-v_2\|<3\dt_1,$ then
$$
{\rm bott}_1(u,v_1)={\rm bott}_1(u,v_2).
$$

Let  $g_1, g_2,...,g_{k(A)}\in U(M_{m(A)}(A))$ $(m(A)\ge 1$ is an integer) be a finite subset such that
$\{\bar{g_1}, \bar{g_2},..., \bar{g}_{k(A)}\}\subset J_c(K_1(A))$ and such that $\{[g_1], [g_2],...,[g_{k(A)}]\}$ forms  a set of generators for $K_1(A).$
Let ${\cal U}=\{\bar{g_1}, \bar{g_2},..., \bar{g}_{k(A)}\}\subset J_c(K_1(A))$ be a finite subset.

Let ${\cal U}_0\subset A$ be a finite subset
such that
$$
\{g_1, g_2,...,g_{k(A)}\}\subseteq\{(a_{i,j}): a_{i,j}\in {\cal U}_0\}.
$$

Let $\dt_u=\min\{1/256m(A)^2, \dt_1/16m(A)^2\},$
${\cal G}_u={\cal F}\cup{\cal G}_1\cup {\cal U}_0$ and
let ${\cal P}_u={\cal P}_0.$
%\cup {\cal P}_1.$

%We now assume that $A$ has property (P) and $A$ is infinite dimensional.

Let $\dt_2>0$ (in place of $\dt$), ${\cal G}_2\subset A$ (in place of ${\cal G}$), ${\cal H}_2'\subset A_+\setminus \{0\}$ (in place of ${\cal H}$), $N_1\ge 1$ (in place of $N$) be the finite subsets and the constants as required by \ref{Ext3} %{\color{green} (the statement of \ref{Ext3} requires that $A=A_0$ or $A=A_0\otimes\mathrm{C}(\mathbb T)$ for some $A \in\mathcal C_0$)}
for the data
$\dt_u$ (in place of $\ep$), ${\cal G}_u$ (in place of ${\cal F}$), ${\cal P}_u$ (in place of ${\cal P}$)  and
$\Delta$ and with $\bar{g}_j$ (in place of $g_j$), $j=1,2,...,k(A)$ (with $k(A)=r$).

Let $d=\min\{\Delta(\hat{h}): h\in {\cal H}_2'\}.$
 Let $\dt_3>0$  and let ${\cal G}_3\subset A\otimes C(\T)$
 be finite subset satisfying the following:
 For any ${\cal G}_3$-$\dt_3$-multiplicative \morp\, $L': A\otimes C(\T)\to C'$ (for any unital \CA\, $C'$ with $T(C')\not=\emptyset$),
 \beq\label{Uni-10}
 |\tau([L]({\boldsymbol{\bt}}(\bar{g}_j))|<d/8,\,\,\,j=1,2,...,k(A).
 \eneq

Without loss of generality, we may assume
that
$$
{\cal G}_3=\{g\otimes z: g\in {\cal G}_3'\andeqn z\in \{1,z, z^*\}\},
$$
where
${\cal G}_3'\subset A$ is a finite subset (by  choosing a smaller $\dt_3$ and large ${\cal G}_3'$).

  Let $\ep_1''=\min\{d/27m(A)^2, \dt_u/2,  \dt_2/2m(A)^2, \dt_3/2m(A)^2\}$ and let ${\bar \ep}_1>0$ (in place of $\dt$) and ${\cal G}_4\subset A$ (in place of ${\cal G}$) be a finite subset as
required by \ref{oldnuclearity} for $\ep_1''$ (in place of $\ep$) and ${\cal G}_u\cup {\cal G}_3'.$
Put
$
\ep_1=\min\{\ep_1', \ep_1'',{\bar \ep_1}\}.
$
Let ${\cal G}_5={\cal G}_u\cup {\cal G}_3'\cup {\cal G}_4.$

%Let ${\cal H}_3'\subset A_+\setminus \{0\}$  be required by  Theorem \ref{Newunique1}  %{\color{green}(the statement of \ref{Newunique1} requires that $A=A_0$ or $A=A_0\otimes C(\mathbb T)$ for some $A_0\in\mathcal C_0$)}
%for $\ep_1/4$ (in place $\ep$) and ${\cal G}_5$ (in place of ${\cal F}$).

%Let $\eta_3''>0$ (in place of $\eta_1$) be required by 2.12 of \cite{ChomC1}
%for $\ep_1/4$ (in place of $\ep$) and ${\cal G}_5$
%(in place of ${\cal F}$).

%Let $\sigma_3=\min\{\Delta(\hat{h}): h\in {\cal H}_3'\}.$
% Let $\dt_4>0$ (in place of $\dt$),  ${\cal G}_6\subset A$  be a finite subset (in place of ${\cal G}$), ${\cal H}_4'\subset A_{s.a.}$ be a finite subset, let ${\cal P}_3\subset \underline{K}(A)$ (in place of ${\cal P}$) and $\sigma_4>0$ (in place of $\sigma_2$) be required by \ref{Newunique1}  for $\ep_1/16$ (in place of $\ep$), ${\cal G}_5$ (in place of ${\cal F}$) and $\sigma_3$ (in place of $\sigma_1$).

Let $\mathcal H_3'\subseteq A_+^{\bf 1}\setminus \{0\}$ (in place of $\mathcal H_1$), $\dt_4>0$ (in place of $\dt$),  ${\cal G}_6\subset A$  (in place of ${\cal G}$), ${\cal H}_4'\subset A_{s.a.}$ (in place of $\mathcal H_2$), ${\cal P}_1\subset \underline{K}(A)$ (in place of ${\cal P}$) and $\sigma_4>0$ (in place of $\sigma_2$) be the finite subsets and constants  as required by Theorem \ref{Newunique1} with respect to $\ep_1/4$ (in place $\ep$) and ${\cal G}_5$ (in place of ${\cal F}$) and $\Delta$.

Choose $N_2\ge N_1$ such that $(k(A)+1)/N_2<d/8.$
 Choose ${\cal H}_5'\subset A_+^{\bf 1} \setminus \{0\}$ and
 $\dt_5>0$ and a finite subset ${\cal G}_7\subset A$ such that, for any $M_m$ and unital ${\cal G}_7$-$\dt_5$-multiplicative \morp\, $L': A\to M_m,$
 if
 ${\rm tr}\circ L'(h)>0\tforal h\in {\cal H}_5',$
 then
 $m\ge N_2((8/d)+1).$

Put $\dt=\min\{\ep_1/16,  \dt_4/4m(A)^2, \dt_5/4m(A)^2\},$
${\cal G}={\cal G}_5\cup {\cal G}_6\cup {\cal G}_7$, and  ${\cal P}={\cal P}_u\cup {\cal P}_1.$
%Let ${\bar G}_1$ be the subgroup generated by ${\cal P}.$
%We then view ${\bar G}_1\cup {\boldsymbl{\bt}}({\bar G}_1)$ as a subgroup of $\underline{K}(A\otimes C(\T)).$
%There is an integer $n_0$ such that
%\beq\label{UniqN-int}
%
%\eneq
Put
$$
{\cal H}_1={\cal H}_1'\cup {\cal H}_2'\cup {\cal H}_3'\cup {\cal H}_4'\cup {\cal H}_{{5}}'
$$
and let ${\cal H}_2={\cal H}_4'.$
Let $\gamma_1=\sigma_4$ and let
$0<\gamma_2<\min\{d/16m(A)^2, \dt_u/9m(A)^2, 1/256m(A)^2\}.$

Now suppose that $C\in {\cal C}$ and $\phi, \psi: A\to C$ are two unital
${\cal G}$-$\dt$-multiplicative \morp s satisfying the condition of the theorem for the given $\Delta,$ ${\cal H}_1,$ $\dt,$ ${\cal G},$ ${\cal P},$ ${\cal H}_2,$ $\gamma_1,$ $\gamma_2$ and ${\cal U}.$

Let
$$
0=t_0<t_1<\cdots <t_n=1
$$
be a partition of $[0,1]$ so that
\beq\label{Uni-11}
\|\pi_{t}\circ \phi(g)-\pi_{t'}\circ \phi(g)\|<\ep_1/16\andeqn
\|\pi_{t}\circ \psi(g)-\pi_{t'}\circ \psi(g)\|<\ep_1/16
\eneq
for all $g\in {\cal G},$ provided $t, t'\in [t_{i-1}, t_i],$ $i=1,2,...,n.$

We write $C=A(F_1, F_2, h_0, h_1),$
$F_1=M_{m_1}\oplus M_{m_2}\oplus \cdots \oplus M_{m_{F(1)}}$ and
$F_2=M_{n_1}\oplus M_{n_2}\oplus \cdots \oplus M_{n_{F(2)}}.$ By the choice of ${\cal H}_5',$ one has that
\beq\label{Uni11+1}
n_j\ge N_2(8/d+1)\andeqn m_s\ge  N_2(8/d+1),\,\,\, 1\le j \le  F(2), \,\,\, 1\le s\le F(1).
\eneq
Applying Theorem \ref{Newunique1}, one obtains a unitary
$w_i\in F_2,$ if $0<i<n,$ $w_0\in h_0(F_1),$ if $i=0,$
and $w_1\in h_1(F_1),$ if $i=1,$
such that
\beq\label{Uni-12}
\|w_i\pi_{t_i}\circ \phi(g)w_i^*-\pi_{t_i}\circ \psi(g)\|<\ep_1/16\tforal g\in {\cal G}_5.
\eneq

It follows from \ref{endpoints} that  we may assume that there is a unitary  $w_e\in F_1$ such that $h_0(w_e)=w_0$ and $h_1(w_e)=w_n.$

By \eqref{Uni1-3+1}, there is a unitary  $\omega_j\in M_{m(A)}(C)$
such that $\omega_j\in CU(M_{m(A)}(C))$ and
\beq\label{Uni-13}
\|\langle (\phi\otimes {\rm id}_{M_{m(A)}}(g_j^*)\rangle \langle (\psi\otimes {\rm id}_{M_{m(A)}})(g_j)\rangle  -\omega_j\|<\gamma_2,\,\,\,j=1,2,...,k(A).
\eneq
({\it note that we now have $w_i$ as well as $\omega_i$ in the proof}.) Write
$$
\omega_j=\prod_{l=1}^{e(j)}\exp(\sqrt{-1}a_j^{(l)})
$$
for some selfadjoint element $a_j^{(l)}\in M_{m(A)}(C),$
$l=1,2,...,e(j),$ $j=1,2,...,k(A).$
Write
\beq\nonumber
a_j^{(l)}=(a_j^{(l,1)}, a_j^{(l,2)},...,a_j^{(l,n_{F(2)})})\andeqn \omega_j=(\omega_{j,1},\omega_{j,2},...,\omega_{j,F(2)})
\eneq
in $C([0,1], F_2)=C([0,1],M_{n_1})\oplus \cdots \oplus C([0,1],M_{n_{{F(2)}}}),$
where $\omega_{j,s}=\exp(\sqrt{-1}a_j^{(l,s)}),$ $s=1,2,...,F(2).$
Then
$$
\sum_{l=1}^{e(j)}{n_s(t_s\otimes {\rm Tr}_{m(A)})(a_j^{(l,s)}(t))\over{2\pi}}\in \Z,\quad t\in [0,1],
$$
where $t_s$ is the normalized trace on $M_{n_s},$
$s=1,2,...,F(2).$
In particular,
\beq\label{Uni-15}
\sum_{l=1}^{e(j)}n_s(t\otimes {\rm Tr}_{m(A)})(a_j^{(l,s)}(t))=\sum_{l=1}^{e(j)}n_s(t\otimes {\rm Tr}_{m(A)})(a_j^{(l,s)}(t'))
\rforal t, t''\in [0,1].
\eneq
%for all $t, t'\in [0,1].$

Let $W_i=w_i\otimes {\rm id}_{M_{m(A)}},$ $i=0,1,....,n$
and $W_e=w_e\otimes {\rm id}_{M_{m(A)}(F_1)}.$
Then
\begin{eqnarray}\label{Ui-16}
&&\|\pi_i(\langle \phi\otimes {\rm id}_{M_{m(A)}})(g_j^*)\rangle) W_i(\pi_i(\langle \phi\otimes {\rm id}_{M_{m(A)}})(g_j)\rangle)W_i^*-\omega_j(t_i)\|\\
&&\hspace{0.2in}<3m(A)^2\ep_1+2\gamma_2<1/32.
\end{eqnarray}
We also have (with $\phi_e=\pi_e\circ \phi$)
\beq\label{Uni-17}
\hspace{-0.4in}\|\langle \phi_e\otimes {\rm id}_{M_{m(A)}})(g_j^*)\rangle W_e(\langle \phi_e\otimes {\rm id}_{M_{m(A)}})(g_j)\rangle)W_e^*-\pi_e(\omega_j)\|<3m(A)^2\ep_1+2\gamma_2<1/32.
\eneq

It follows from \eqref{Ui-16} that there exists selfadjoint elements $b_{i,j}\in M_{m(A)}(F_2)$ such that
\beq\label{Uni-18}
\hspace{-0.2in}\exp(\sqrt{-1}b_{i,j})=\omega_j(t_i)^*(\pi_i\langle \phi\otimes {\rm id}_{M_{m(A)}})(g_j^*)\rangle)W_i(\pi_i(\langle \phi\otimes {\rm id}_{M_{m(A)}})(g_j)\rangle) W_i^*,
\eneq
and $b_{e,j}\in M_{m(A)}(F_1)$ such that
\beq\label{Uni-19}
\hspace{-0.2in}\exp(\sqrt{-1}b_{e,j})=\pi_e(\omega_j)^*(\pi_e\langle \phi\otimes {\rm id}_{M_{m(A)}})(g_j^*)\rangle)W_e(\pi_e(\langle \phi\otimes {\rm id}_{M_{m(A)}})(g_j)\rangle) W_e^*,
\eneq
and
\beq\label{Uni-20}
\|b_{i,j}\|<2\arcsin (3m(A)^2\ep_1/4+2\gamma_2),\,\,\,j=1,2,...,k(A),\,i=0,1,...,n,e.
\eneq
Write
\beq\nonumber
b_{i,j}=(b_{i,j}^{(1)},b_{i,j}^{(2)},...,b_{i,j}^{F(2)})\in F_2\andeqn
b_{e,j}=(b_{e,j}^{(1)}, b_{e,j}^{(2)},...,b_{e,j}^{(F(1))})\in F_1.
\eneq
We have that
\beq\label{Uni-20+}
h_0(b_{e,j})=b_{0,j}\andeqn h_1(b_{e,j})=b_{n,j}.
\eneq
Note that
\beq\label{Uni-21}
(\pi_i(\langle \phi\otimes {\rm id}_{M_{m(A)}}(g_j^*)\rangle))W_i(\pi_i(\langle \phi\otimes {\rm id}_{M_{m(A)}})(g_j)\rangle) W_i^*
=\pi_i(\omega_j)\exp(\sqrt{-1}b_{i,j}),
\eneq
$j=1,2,...,k(A)$ and $i=0,1,...,n,e.$

Then,
\beq\label{Uni-22}
{n_s\over{2\pi}}(t_s\otimes {\rm Tr}_{M_{m(A)}})(b_{i,j}^{(s)})\in \Z,
\eneq
where $t_s$ is the normalized trace on $M_{n_s},$
$s=1,2,...,F(2),$ $j=1,2,...,k(A),$ and $i=0,1,...,n.$
We also have
\beq\label{Uni-23}
{m_s\over{2\pi}}(t_s\otimes {\rm Tr}_{M_{m(A)}})(b_{e,j}^{(s)})\in \Z,
\eneq
where $t_s$ is the normalized trace on $M_{m_s},$ $s=1,2,...,F(1),$
$j=1,2,...,k(A).$
Put
$$
\lambda_{i,j}^{(s)}={n_s\over{2\pi}}(t_s\otimes {\rm Tr}_{M_{m(A)}})(b_{i,j}^{(s)})\in \Z,
$$
where $t_s$ is the normalized trace on $M_{n_s},$
$s=1,2,...,n,$ $j=1,2,...,k(A)$ and  $i=0,1,2,...,n.$

Put
$$
\lambda_{e,j}^{(s)}={m_s\over{2\pi}}(t_s\otimes {\rm Tr}_{M_{m(A)}})(b_{e,j}^{(s)})\in \Z
$$
where $t_s$ is the normalized trace on $M_{m_s},$
$s=1,2,...,F(1)$ and $j=1,2,...,k(A).$
Denote
\beq\nonumber
\lambda_{i,j}=(\lambda_{i,j}^{(1)}, \lambda_{i,j}^{(2)},...,\lambda_{i,j}^{(F(2)})\in \Z^{F(2)},\andeqn
\lambda_{e,j}=(\lambda_{e,j}^{(1)},\lambda_{e,j}^{(2)},...,\lambda_{e,j}^{(F(1)})\in \Z^{F(1)}.
\eneq
 %Moreover
 We have, by (\ref{Uni-20}),
\beq\label{Uni-24}
|{\lambda_{i,j}^{(s)}\over{n_s}}|&<&d/4, \,\,\,s=1,2,...,F(2),\andeqn\\\label{Uni-24+}
|{\lambda_{e,j}^{(s)}\over{m_s}}|&<&d/4,\,\,\,s=1,2,...,F(1),
\eneq
 $j=1,2,...,k(A),$ $i=0,1,2,...,n.$

Define $\af_i^{(0,1)}: K_1(A)\to \Z^{F(2)}$ by mapping $[g_j]$ to $\lambda_{i,j},$ $j=1,2,...,k(A)$, $i=0,1,2,...,n,$ and
define
$\af_e^{(0,1)}: K_1(A)\to \Z^{F(1)}$ by
mapping $[g_j]$ to $\lambda_{e,j},$ $j=1,2,...,k(A).$
We write $K_0(A\otimes C(\T))=K_0(A)\oplus {\boldsymbol{\bt}}(K_1(A)))$
(see   \ref{Dbeta}
%2.10 of \cite{LnHomtp}
for the definition
of ${\boldsymbol{\bt}}$).
Define $\af_i: K_*(A\otimes C(\T))\to K_*(F_2)$ as follows:
On $K_0(A\otimes C(\T)),$ define
\beq\label{Uni-25}
\af_i|_{K_0(A)}=[\pi_i\circ \phi]|_{K_0(A)},\,\,\,
\af_i|_{{\boldsymbol{\bt}}(K_1(A))}=\af_i\circ {\boldsymbol{\bt}}|_{K_1(A)}=\af_i^{(0,1)}
\eneq
and on $K_1(A\otimes C(\T)),$ define
\beq\label{Uni-25+}
\af_i|_{K_1(A\otimes C(\T))}=0,\,\,\,i=0,1,2,...,n.
%\af_i|_{K_1(C(\T)\otimes C(X_0))}=0,\\
%\af_i|_{K_*(C(\T)\otimes C(X_0),\Z/m\Z}=0,\,\,\,m=2,3,...,K.
\eneq
%$i=0,1,2,...,n.$

Also define  $\af_e\in {\rm Hom}(K_*(A\otimes C(\T)), K_*(F_1)),$ by
\beq\label{Uni-26}
\af_e|_{K_0(A)}=[\pi_e\circ \phi]|_{K_0(A)},\,\,\,
\af_e|_{{\boldsymbol{\bt}}(K_1(A))}=\af_i\circ {\boldsymbol{\bt}}|_{K_1(A)}=\af_e^{(0,1)}
\eneq
on $K_0(A\otimes C(\T))$ and  $(\af_e)|_{K_1(A\otimes C(\T))}=0.$
%{\color{green}
Note that
\beq\label{Uni-26+0}
(h_0)_{*}\circ \af_e=\af_0\andeqn (h_1)_{*}\circ\af_e=\af_n.
\eneq
Since $A\otimes C(\mathbb{T})$ satisfies the UCT,  the map $\alpha_e$ can be lifted to an element of $KK(A\otimes C(\mathbb T), F_1)$ which is still denoted by $\alpha_e$. Then define
\beq\label{Uni-26+}
\af_0=\af_e\times [h_0] \andeqn \af_n=\af_e\times [h_1]
\eneq
in $KK(A\otimes C(\mathbb T), F_2)$.
%Note that these $KK$-elements are the lifting of the original maps $\alpha_0$ and $\alpha_n$ respectively.
For $i=1, ..., n-1$, also pick a lifting of $\alpha_i$ in $KK(A\otimes C(\mathbb T), F_2)$, and still denote it by $\alpha_i$.
%}
%Since we assume that $K_*(A)$ is torsion free, $\af_i$ defines an element in $KK(A\otimes C(\T), F_2),$ if
%$i=1,2,...,n,$ and $\af_e$ defines an element in $KK(A\otimes C(\T), F_1).$
%Note that
%\beq\label{Uni-26+}
%[\af_e]\times [h_0]=\af_0 \andeqn [\af_e]\times [h_1]=\af_n.
%\eneq
We estimate that
$$
\|(w_{i}^*w_{i+1})\pi_{t_i}\circ \phi(g)-\pi_{t_i}\circ \phi(g)(w_{i}^*w_{i+1})\|<\ep_1/4\tforal g\in {\cal G}_5,
$$
$i=0,1,...,n-1.$
Let $\Lambda_{i,i+1}: C(\T)\otimes A\to F_2$ be a unital \morp\, given
by the pair $w_{i}^*w_{i+1}$ and $\pi_{t_i}\circ \phi$  (by \ref{oldnuclearity}, see 2.8 of \cite{LnHomtp}).
Denote $V_{i,j}=\langle \pi_{t_i}\circ \phi\otimes {\rm id}_{M_{m(A)}}(g_j) \rangle,$ $j=1,2,...,k(A)$ and $i=0,1,2,...,n-1.$

Write
$$
V_{i,j}=(V_{i,j, 1}, V_{i,j,2},..., V_{i,j,F(2)})\in F_2,\,\,\, j=1,2,...,k(A),\,\,\, i=0,1,2,...,n.
$$
Similarly, write
$$
W_i=(W_{i,1},W_{i,2},...,W_{i,F(2)})\in F_2,\,\,\, i=0,1,2,...,n.
$$
We have
\beq\label{Uni-27}
\|W_{i}V_{i,j}^*W_{i}^* V_{i,j}V_{i,j}^*W_{i+1}V_{i,j}W_{i+1}^*-1\|<1/16\\
\|W_{i}V_{i,j}^*W_{i}^*V_{i,j}V_{i+1,j}^*
W_{i+1} V_{i+1,j}W_{i+1}^*-1\|<1/16
\eneq
and there is a continuous path $Z(t)$
of unitaries  such that $Z(0)=V_{i,j}$ and $Z(1)=V_{i+1,j}.$ Since
$$
\|V_{i,j}-V_{i+1,j}\|<\dt_1/12,\,\,\,j=1,2,...,k(A),
$$
we may assume that $\|Z(t)-Z(1)\|<\dt_1/6$ for all $t\in [0,1].$
We also write
$$
Z(t)=(Z_1(t), Z_2(t), ...,Z_{F(2)}(t))\in F_2\andeqn t\in [0,1].
$$
We obtain a continuous path
$$
W_{i}V_{i,j}^*W_{i}^*V_{i,j}Z(t)^*W_{i+1} Z(t)W_{i+1}^*
$$
which is in $CU(M_{nm(A)})$
for all $t\in  [0,1]$ and
$$
\|W_{i}V_{i,j}^*W_{i}^*V_{i,j}Z(t)^*W_{i+1} Z(t)W_{i+1}^*-1\|<1/8\tforal t\in [0,1].
$$
It follows that
$$
(1/2\pi\sqrt{-1})(t_s\otimes {\rm Tr}_{M_{m(A)}})[\log(W_{i,s}V_{i,j,s}^*W_{i,s}^*V_{i,j,s}Z_s(t)^*W_{i+1,s}Z_s(t)W_{i+1,s}^*)]
$$
is a constant integer,  where $t_s$ is the normalized trace on $M_{n_s}.$
In particular,
\beq\label{Uni-28}
&&(1/2\pi\sqrt{-1})(t_s\otimes {\rm Tr}_{M_{m(A)}})(\log(W_{i,s}V_{i,j,s}^*W_{i,s}^*W_{i+1,s} V_{i,j,s}W_{i+1,s}^*))\\\label{Uni-28+}
&&=(1/2\pi\sqrt{-1})(t_s\otimes {\rm Tr}_{M_{m(A)}})(\log(W_{i,s}V_{i,j,s}^*W_{i,s}^*V_{i, j}V_{i+1,j,s}^*W_{i +1}V_{i,j,s}W_{i+1,s}^*)).
\eneq

One also has
\begin{eqnarray}\label{Uni-29}
&&W_{i}V_{i,j}^*W_{i}^*V_{i, j}V_{i+1,j}^*W_{i+1} V_{i+1,j}W_{i+1}^*\\
&&=(\omega_j(t_{i})\exp(\sqrt{-1}b_{i,j}))^*\omega_j(t_i)
\exp(\sqrt{-1}b_{i+1,j})\\\label{Uni-30}
&&=\exp(-\sqrt{-1}b_{i,j})\omega_j(t_{i})^*\omega_j(t_{i+1})
\exp(\sqrt{-1}b_{i+1,j}).
\end{eqnarray}
Note that, by (\ref{Uni-13}) and (\ref{Uni-11}), for $t\in [t_i, t_{i+1}],$
\begin{equation}
\|\omega_j(t_{i})^*\omega_j(t)-1\|<3(3\ep_1'+2\gamma_2)<3/32,
\end{equation}
$j=1,2,...,k(A),$ $i=0,1,..., n-1.$

By Lemma 3.5 of \cite{Lin-AU11},
\begin{equation}\label{LeasyApp}
(t_s\otimes {\rm Tr}_{m(A)})(\log(\omega_{j,s}(t_{i})^*\omega_{j,s}(t_{i+1})))=0.
\end{equation}
It follows that (by the Exel formula (see \cite{HL}), using (\ref{Uni-28+}), (\ref{Uni-30}) and (\ref{LeasyApp}))
\beq\label{Uni-31}
&&\hspace{-0.6in}(t\otimes {\rm Tr}_{m(A)})({\rm bott}_1(V_{i,j}, W_{i}^*W_{i+1}))\\
 \hspace{-0.2in}&=&
({1\over{2\pi \sqrt{-1}}})(t\otimes {\rm Tr}_{m(A)})(\log(V_{i,j}^*W_{i}^*W_{i+1}V_{i,j}W_{i+1}^*W_{i}))\\
 \hspace{-0.2in}&=&({1\over{2\pi \sqrt{-1}}})(t\otimes {\rm Tr}_{m(A)})(\log(W_{i}V_{i,j}^*W_{i}^*W_{i+1}V_{i,j}W_{i+1}^*))\\
&=&({1\over{2\pi \sqrt{-1}}})(t\otimes {\rm Tr}_{m(A)})(\log(W_{i}V_{i,j}^*W_{i}^*V_{i,j}V_{i+1,j}^*
W_{i+1}V_{i+1,j}W_{i+1}^*))\\
&=& ({1\over{2\pi \sqrt{-1}}})(t\otimes {\rm Tr}_{m(A)})(\log(\exp(-\sqrt{-1}b_{i,j})\omega_j(t_{i})^*
\omega_j(t_{i+1})\exp(\sqrt{-1}b_{i+1,j}))\\
&=& ({1\over{2\pi \sqrt{-1}}})[(t\otimes {\rm Tr}_{k(n)})(-\sqrt{-1}b_{i,j})+(t\otimes {\rm Tr}_{k(n)})(\log(\omega_j(t_{i})^*\omega_j(t_{i+1}))\\
&&\hspace{1.4in}+(t\otimes {\rm Tr}_{k(n)})(\sqrt{-1}b_{i,j})]\\
%&=& ({1\over{2\pi \sqrt{-1}}})
&=&{1\over{2\pi}}(t\otimes {\rm Tr}_{k(n)})(-b_{i,j}+b_{i+1,j})
\eneq
for all $t\in T(F_2).$
In other words,
\beq\label{Uni-32}
{\rm bott}_1(V_{i,j}, W_{i}^*W_{i+1}))=-\lambda_{i,j}+\lambda_{i+1,j}
\eneq
$j=1,2,...,m(A),$  $i=0,1,...,  n-1.$

%{\color{green}
Consider $\alpha_0, ..., \alpha_n\in KK(A\otimes C(\mathbb T), F_2)$ and $\alpha_e\in KK(A\otimes C(\mathbb T), F_1)$. Note that
$$ \left|\alpha_i(g_j)\right|=\left|\lambda_{i, j}\right|,$$
and by \eqref{Uni-24}, one has
$$m_s, n_j\ge N_2(8/d+1).$$
Applying \ref{Ext3}
(using \eqref{Uni-24+},  among other items), there are unitaries
%{\color{green}
$z_i\in F_2$, $i=1,2,...,n-1$, and $z_e\in F_1$ such that
\begin{eqnarray}\label{Uni-33}
&&\|[z_i,\,\pi_{t_i}\circ\phi(g)]\|<\dt_u\tforal g\in {\cal G}_u\\
%\Andean
%{\rm Bott}(z_i, \pi_{t_i}\circ \phi)=\af_i.
&& {\rm Bott}(z_i, \pi_{t_i}\circ \phi)=\alpha_i\andeqn
 {\rm Bott}(z_e, \pi_e\circ \phi)=\alpha_e.
\end{eqnarray}
Put $$z_0=h_0(z_e)\quad \mathrm{and} \quad z_n=h_1(z_e).$$
One verifies (by \eqref{Uni-26+}) that
\begin{eqnarray}\label{Uni-34}
%{\rm Bott}(h_j(z_e), \pi_{t_j}\circ \phi)=\af_i\times [h_j],\,\,\,j=0,1.
{\rm Bott}(z_0, \pi_{t_0}\circ \phi)=\alpha_0\andeqn
{\rm Bott}(z_n, \pi_{t_n}\circ \phi)=\alpha_n.
\end{eqnarray}

%Replacing $z_0$ and $z_n,$ if necessarily, we may assume that $z_0=h_0(z_e)$ and $z_n=h_1(z_e).$
Let $U_{i,i+1}=z_{i}(w_{i})^*w_{i+1}(z_{i+1})^*,$
$i=0,1,2,...,n-1.$
Then
\beq\label{NT-8}
\|[U_{i,i+1}, \pi_{t_i}\circ \phi(g)]\|<\min\{\dt_1, \dt_2\},\quad g\in {\cal G}_u,\ i=0,1,2,...,n-1.
\eneq
Moreover, for $i=0,1, 2, ..., n-1$,
%{\color{green}
\begin{eqnarray*}
{\rm bott}_1(U_{i,i+1}, \pi_{t_i}\circ \phi)&=&
{\rm bott}_1(z_{i},\, \pi_{t_i}\circ \phi))+{\rm bott}_1((w_{i}^*w_{i+1}, \pi_{t_i}\circ \phi))\\
&&+{\rm bott}_1((z_{i+1}^*,\, \pi_{t_i}\circ\phi)\\
&=& (\lambda_{i,j})+(-\lambda_{i,j}+\lambda_{i+1,j})+(-\lambda_{i+1,j})\\
&=&0.
\end{eqnarray*}
Note that for any $x\in\bigoplus_{*=0, 1}\bigoplus_{k=1}^{n_0} K_*(A\otimes C(\mathbb T),\mathbb Z/k\mathbb Z)$, one has $N x=0$.

Therefore
\begin{equation}\label{Uni-37}
{\rm Bott}((\underbrace{U_{i,i+1}, ..., U_{i, i+1}}_N), (\underbrace{\pi_{t_i}\circ \phi, ..., \pi_{t_i}\circ\phi}_N))|_{\cal P}=N{\rm Bott}(U_{i,i+1}, \pi_{t_i}\circ \phi)|_{\cal P} =0,
\end{equation}
$i=0, 1, 2,...,n-1.$
%}
%{and
%\begin{eqnarray}
%\textrm{Bott}(U_{n-1, n}, \pi_{t_{n-1}}\circ\phi) &=& \tilde{\kappa}_{n-1} +\Lambda_{n-1, n}- {\rm Bott}((z_{n},\, \pi_{1}\circ\phi) \\
% &=& \tilde{\alpha}_{n-1}+\beta_{n-1}+\Lambda_{n-1, n}-{\rm Bott}((z_{n},\, \pi_{1}\circ\phi)\\
% &=&0.
%\end{eqnarray}
%}
%\beq
%\hspace{-0.4in}{\rm bott}_1(U_{i,i+1}, \pi_{t_i}\circ \phi)(g_j)&=&
%{\rm bott}_1(z_{i},\, \pi_{t_i}\circ \phi))(g_j)\\
%&&+{\rm bott}_1((w_{i}^*w_{i+1}, \pi_{t_i}\circ \phi))(g_j)\\
%&&\hspace{0.2in}+{\rm bott}_1((z_{i+1}^*,\, \pi_{t_i}\circ\phi)(g_j)\\
%&=& (\lambda_{i,j})+(-\lambda_{i,j}+\lambda_{i+1,j})+(-\lambda_{i+1,j})\\
%&=&0
%\eneq
%$j=1,2,...,k(A).$
%Note also that $K_1(F_2)=\{0\}.$ It follows that, since $K_*(A)$ is torsion free,
%That is,
%\beq\label{Uni-37}
%{\rm Bott}(U_{i,i+1}, \pi_i\circ \phi)|_{\cal P}=0,\,\,\, i=0, 1, 2,...,n-1.
%\eneq
Note that, by the assumption \eqref{Uni1-2},
\beq\label{Uni-38}
t_s\circ \pi_t\circ \phi(h)\ge \Delta(\hat{h})\tforal h\in {\cal H}_1',
\eneq
where $t_s$ is the normalized trace on $M_{n_s},$ $1\le s\le F(2).$

By applying \ref{homotopy1}, using \eqref{Uni-38}, \eqref{NT-8} and \eqref{Uni-37}, there exists
a continuous path of unitaries, $\{\tilde{U}_{i,i+1}(t): t\in [t_i, t_{i+1}]\}\subset F_2\otimes M_N(\mathbb C)$ such that
\begin{equation}\label{Uni-39}
\tilde{U}_{i,i+1}(t_i)={\rm id}_{F_2\otimes M_N(\mathbb C)},\,\,\, \tilde{U}_{i, i+1}(t_{i+1})=(z_iw_i^*w_{i+1}z^*_{i+1})\otimes 1_{M_N(\mathbb C)},
\end{equation}
and
\begin{equation}\label{Uni-39+}
\|\tilde{U}_{i, i+1}(t)(\underbrace{\pi_{t_i}\circ \phi(f), ..., {{\pi}}_{t_i}\circ\phi(f)}_N)\tilde{U}_{i, i+1}(t)^*-(\underbrace{\pi_{t_i}\circ \phi(f), ..., {{\pi}}_{t_i}\circ\phi(f)}_N)\|<\ep/32
\end{equation}
for all $f\in {\cal F}$ and for all $t\in [t_i, t_{i+1}].$ Define
$W\in C\otimes M_N$ by
\beq\label{Uni-40}
W(t)=(w_iz_i^*\otimes 1_{M_N})\tilde{U}_{i, i+1}(t)\tforal t\in [t_i, t_{i+1}],
\eneq
$i=0,1,...,n-1.$ Note that $W(t_i)=w_iz_i^*\otimes 1_{M_N},$ $i=0,1,...,n.$
Note also that
$$W(0)=w_0z_0^*\otimes 1_{M_N} =h_0(w_ez_e^*)\otimes 1_{M_N}$$ and
$$W(1)=w_nz_n^*\otimes 1_{M_N}=h_1(w_ez_e^*)\otimes 1_{M_N}.$$
So $W\in C\otimes M_N.$
One then checks that, by  (\ref{Uni-11}), (\ref{Uni-39+}) , (\ref{Uni-33}) and (\ref{Uni-12}),
\begin{eqnarray}
&&\|W(t)((\pi_t\circ \phi)(f)\otimes 1_{M_N})W(t)^*-(\pi_t\circ \psi)(f)\otimes 1_{M_N}\|\\
&<&\|W(t)((\pi_t\circ \phi)(f)\otimes 1_{M_N})W(t)^*-W(t)((\pi_{t_i}\circ \phi)(f)\otimes 1_{M_N})W^*(t)\|\\
 &&+\|W(t)(\pi_{t_i}\circ \phi)(f)W(t)^*-W(t_i)(\pi_{t_i}\circ \phi)(f) W(t_i)^*\|\\
 &&+\|W(t_i)((\pi_{t_i}\circ \phi)(f)\otimes 1_{M_N}) W(t_i)^*-(w_i(\pi_{t_i}\circ \phi)(f)w_i^*)\otimes 1_{M_N}\|\\
 &&+\|w_i(\pi_{t_i}\circ \phi)(f)w_i^*-\pi_{t_i}\circ \psi(f)\|\\
 && +\|\pi_{t_i}\circ \psi(f)-\pi_t\circ \phi(f)\|\\
 &<&\ep_1/16+\ep/32+\dt_u+\ep_1/16+\ep_1/16<\ep
\end{eqnarray}
 for all $f\in {\cal F}$ and for $t\in [t_i, t_{i+1}]$.
\end{proof}

\begin{rem}\label{Rem-on-N}
With a minor modification, the proof also works without assuming that $K_*(A)$ is  finitely generated.
In Theorem \ref{UniqN1}, the multiplicity $N$ only depends on $\underline{K}(A)$ as $\underline{K}(A)$ is finitely generated. However, if $K_*(A)$ is not finitely generated, it also depends on $\mathcal F$ and $\epsilon$. Moreover, if $K_*(A)$ is torsion free, or if $K_1(C)=0$, then the multiplicity $N$ can be chosen to be $1$.
\end{rem}

\begin{cor}\label{CUniN1}
Theorem \ref{UniqN1} holds if $A$ is replaced by $M_m(A)$ for any
integer $m\ge 1.$

\end{cor}

\section{\CA s in ${\cal B}_1$}

\begin{df}\label{DB1}
Let $A$ be a unital infinite dimensional simple \CA. We say $A\in {\cal B}_1$ if the following hold:
Let $\ep>0,$ let $a\in A_+\setminus \{0\}$ and let ${\cal F}\subset A$ be a
finite subset.
There exists a nonzero projection $p\in A$ and a \SCA\, $C\in {\cal C}$
%with $K_1(C)=\{0\}$
%(and $K_0(C)$ is torsion free)
with $1_C=p$ such that
\begin{eqnarray}
&&\|xp-px\|<\ep\tforal x\in {\cal F},\\
&&{\rm dist}(pxp, C)<\ep\tforal x\in {\cal F}\andeqn \\
&&1-p\lesssim a.
\end{eqnarray}

In the above, if $C$ can always be chosen in $\mathcal C_0$, that is, $K_0(C)=\{0\},$ then we say that $A\in {\cal B}_0.$
\end{df}

\begin{df}\label{DgTR}
Let $A$ be a unital simple \CA. We say $A$ has the generalized  tracial rank at most one,  if the following hold:

Let $\ep>0,$ let $a\in A_+\setminus \{0\}$ and let ${\cal F}\subset A$ be a
finite subset.
There exists a nonzero projection $p\in A$ and a \SCA\, $C$ which is a subhomogeneous \CA\, with one dimensional
spectrum, or $C$ is finite dimensional \CA\,
%and
%with $K_1(C)=\{0\}$
%(and $K_0(C)$ is torsion free)
with $1_C=p$ such that
\begin{eqnarray}\label{Dgtr-1}
&&\|xp-px\|<\ep\tforal x\in {\cal F},\\\label{Dgtr-2}
&&{\rm dist}(pxp, C)<\ep\tforal x\in {\cal F}\andeqn \\\label{Dgtr-3}
&&1-p\lesssim a.
\end{eqnarray}

In this case, we write $gTR(A)\le 1.$
\end{df}

\begin{rem} It follows from  {\rm \ref{ASCAs}} that $gTR(A)\le 1$ if and only if $A\in {\cal B}_1.$
\end{rem}

%If $A\in {\cal B}_1,$ we  say $A$ has generalized tracial rank at most one.  In this case, we will write $gTR(A)\le 1.$

%If in the above definition, only (\ref{Dgtr-1}) and (\ref{Dgtr-2}) hold, then we say $A$ has the property
%($L_b$).

Let ${\cal D}$ be a class of unital \CA s.
%We will use the following  general definition (see Definition 2.2 of \cite{EN-Tapprox}).
\begin{df}\label{DATD}
Let $A$ be a unital simple \CA. We say $A$ is tracially approximately in ${\cal D},$ denoted by $A\in\mathrm{TA}{\mathcal D}$, if
 the following hold:

 For any $\ep>0,$ any $a\in A_+\setminus \{0\}$ and any finite subset ${\cal F}\subset A$,
there exist a nonzero projection $p\in A$ and a \SCA\, $C\in {\cal D}$
%with $K_1(C)=\{0\}$
%(and $K_0(C)$ is torsion free)
%and
with $1_C=p$ such that
\begin{eqnarray}\label{DtrD-1}
&&\|xp-px\|<\ep\tforal x\in {\cal F}, \\\label{DtrD-2}
&&{\rm dist}(pxp, C)<\ep\tforal x\in {\cal F},\  \textrm{and} \\
&&1-p\lesssim a.
\end{eqnarray}
(see Definition 2.2 of \cite{EN-Tapprox}).
Note that $\mathcal B_0=\mathrm{TA}\mathcal C_0$ and $\mathcal B_1=\mathrm{TA}\mathcal C$.
%\end{df}
If in the above definition, only (\ref{DtrD-1}) and (\ref{DtrD-2}) hold, then we say $A$ has the property
($L_{\cal D}$).

\end{df}

The following proposition was first appear in an unpublished paper of the second named author
distributed in 1998.
%Note that the following holds for any unital simple \CA\, $A$ which is tracially in
%${\cal D},$ where ${\cal D}$ is a class of unital \CA s, if we replace $C$ by a \CA\, in ${\cal D}.$ A proof using asymptotical sequence argument was given in \cite{EN-Tapprox}.

\begin{prop}\label{Adset9201411}
Let $A$ be a unital simple \CA\, which has the property  {\rm ($L_{\cal D}$)}.
Then, for any $\ep>0$ and any finite subset ${\cal F}\subset A,$ there exists a projection
$p\in A$ and a \SCA\, $C\in {\cal D}$ with $1_C=p$ such that
\beq\label{Addb-1}
&& \|[x,\, p]\|<\ep\tforal f\in {\cal F}, \\
&& {\rm dist}(pxp, C)<\ep,\ \textrm{and}, \\
&& \|pxp\|\ge \|x\|-\ep\tforal x\in {\cal F}.
\eneq
\end{prop}

\begin{proof}
Fix $\ep>0$ and a finite subset ${\cal F}\subset A.$
It is clear that, \wilog, we may assume that $x\in A_+$ and $\|x\|=1$ for all $x\in {\cal F}.$
Let $f\in C(([0,1])_+$ be such that $f(t)=0$ if $t\in [0,1-\ep/2],$ $f(t)=1$ if $t\in [1-\ep/4, 1]$
and $f(t)$ is linear in $(1-\ep/2,1-\ep/4).$
For each such $x\in {\cal F},$ there exist $y_1(x), y_2(x),...,y_{k(x),x}\in A$ such that
\beq\label{Addb-2}
\sum_{i=1}^{k(x)}y_i(x)^*f(x)y_i(x)=1_A.
\eneq
Put
$$
\sigma(x)={1\over{k(x)(\max\{1,\|y_i(x)\|^2\}}}\andeqn
\sigma=\min\{\sigma(x): x\in {\cal F}\}.
$$
%Let $f\in C(([0,1])_+$ be such that $f(t)=0$ if $t\in [0,1-\ep/2],$ $f(t)=1$ if $t\in [1-\ep/4, 1]$
%and $f(t)$ is linear in $(1-\ep/2,1-\ep/4).$
Choose $\dt>0$ and a large finite subset ${\cal F}_1\supset {\cal F}$ such that, for any projection $q,$
\beq\label{Addb-3}
\|qy-yq\|<\dt\rforal y\in {\cal F}_1
\eneq
implies that, for all $x\in {\cal F},$
\beq\label{Addb-4}
\|\sum_{i=1}^{k(x)}qy_i(x)^*qf(x)qy_i(x)q-q\|<\sigma/16 \andeqn \|f(qxq)-qf(x)q\|<\sigma/16.
\eneq
Now, since $A$ has property ($L_{\cal D}$), there is a projection $p\in A$ and a \SCA\, $C\in {\cal D}$ with $1_C=p$
such that
\beq\label{Addb-5}
\|py-yp\|<\min\{\ep/2, \dt\} \andeqn
{\rm dist}(pyp, C)<\ep\rforal y\in {\cal F}_1.
\eneq
Therefore, by the choice of $\dt$ and ${\cal F}_1,$
\beq\label{Addb-6}
\|\sum_{i=1}^{k(x)}py_i(x)pf(x)py_i(x)p\|\ge 1-\sigma/16\andeqn\\
 \|f(pxp)-pf(x)p\|<\sigma/16 \rforal x\in {\cal F}.
\eneq
Therefore
\beq\label{Addb-7}
\|pf(x)p\|\ge 15\sigma/16\rforal x\in {\cal F}.
\eneq
It follows that
\beq\label{Addb-8}
\|f(pxp)\|\ge 14\sigma/16.
\eneq
Therefore
$$
\|pxp\|>1-\ep\rforal x\in {\cal F}.
$$
\end{proof}

\begin{thm}\label{B1sp}
Let $A$ be a unital separable simple \CA\, in ${\cal B}_1$ (or in ${\cal B}_0$).  Then either $A$ is an inductive limit of
unital \CA s  $A_n\in {\cal C}$
(or $A_n\in {\cal C}_0$)
with $A=\overline{\bigcup_{n=1}^{\infty} A_n}$, or $A$  has the property (SP).
%If $A\in {\cal B}_0,$ then either $A$ is an inductive limit of unital \CA s $A_n\in {\cal C}_0$ with $K_1(A_n)=\{0\}.$
\end{thm}

\begin{proof}
This follows from Definition \ref{DB1} immediately.  %We may assume that $A$ is not one dimensional.
Let ${\cal F}_1, {\cal F}_2,...,{\cal F}_n,...$ be a sequence of increasing
finite subsets of the unit ball of $A$ whose union is dense in the unit ball.   If $A$ does not have property (SP), then there is a non-zero positive element $a\in A$ such that $\overline{aAa}\not=A$ and
$\overline{aAa}$ has no non-zero projection.  Then, for each $n\ge 1,$  there is a projection
$1_A-p_n\lesssim a$ and a \SCA\, $C_n\in {\cal C}$ (or ${\cal C}_0$) such that
$1_{C_n}=p_n$ and
\beq\label{B1sp-1}
\|p_nx-xp_n\|<1/2^n \andeqn {\rm dist}(p_nxp_n, C_n)<1/2^n \rforal x\in {\cal F}_n.
\eneq
Since $\overline{aAa}$ does not have any non-zero projection, one has $1_A-p_n=0.$
In other words,
$1_A=p_n$ and
\beq\label{B1sp-2}
{\rm dist}(x, C_n)<1/2^n\rforal x\in {\cal F}_n, n=1,2,....
\eneq
It follows that $\overline{\bigcup_{n=1}^{\infty} C_n}=A.$
Since each \CA\, $C_n$ is semiprojective (see  \ref{DfC1}, also Theorem 6.22 of \cite{ELP1}),
$A$ is in fact an inductive limit of \CA s in $C_0$ (or $\mathcal C_0$).
\end{proof}

\begin{thm}\label{B1stablerk}
Let $A\in {\cal B}_1.$ Then $A$ has stable rank one.
\end{thm}
\begin{proof}
This follows from (\ref{2pg3}) and 4.3 of \cite{EN-Tapprox} (or a similar result in \cite{Fan-sr1}).
\end{proof}

%{\bf I inserted the following two lemmas---N.----I made some changes--L}

\begin{lem}\label{MF}
Let ${\cal D}$ be a family of unital separable  \CA s which are  residually finite dimensional.
Any unital separable simple C*-algebra
% in $TA{\cal C}'$
with property $L_{\cal D}$ can be embedded in $\prod M_{r(n)}/\bigoplus M_{r(n)}$ for some sequence of integers $\{r(n)\}.$
\end{lem}
\begin{proof}
Let $A$ be a unital separable simple C*-algebra  with property ($L_{\cal D}$).
%in
%{$TA{\cal C}'.$}
%$\mathcal B_1$.
Let $\mathcal F_1\subseteq\mathcal F_2\subseteq\cdots\subseteq \mathcal F_i\subseteq\cdots$ be an increasing sequence of finite subsets of $A$ with the union being dense in  $A$.
Since
%$A\in \mathcal B_1$,
$A$ has property ($L_{\cal D}$),
%{is in $TA{\cal C}',$}
for each $n,$ there is a projection $p_n\in A$ and $C_n\subseteq A$ with $1_{C_n}=p_n$ and $C_n\in {\mathcal D}$ such that
\beq\label{MF-1}
\hspace{-0.1in}\| p_nf-fp_n \|<1/2^{n+2},\,\, \|p_nfp_n\|\ge \|f\|-1/2^{n+2} \andeqn  p_nfp_n\in_{1/2^{n+2}} C_n\,\, \rforal f\in\mathcal F_n.
\eneq
For each $a\in {\cal F}_n,$ there exists $c(a)\in C_n$ such that
$\|p_nap_n-c(a)\|<1/2^{n+1}.$
  There are unital \hm s $\pi_n': C_n\to B_n,$ where $B_n$ is a finite dimensional \SCA\, such that
\beq\label{MF-4}
\|\pi_n'(c(a))\|=\|c(a))\|\ge \|a\|-1/2^{n+1} \rforal \in {\cal F}_n,\,\,\,\,n=1,2,....
\eneq
There is an integer $r(n)\ge 1$ such that $B_n$ is unitally embedded into
$M_{r(n)}.$ Denote by $\pi_n: C_n\to M_{r(n)}$ the composition of $\pi_n'$ and the embedding.
Note $C_n\subset p_nAp_n.$  Then there is \morp\, $\Phi_n': p_nAp_n\to M_{r(n)}$ such that
\beq\label{MF-4+}
\Phi_n'|_{C_n}=\pi_n.
%\|\Phi_n(c(a))-\pi_n(c(a))\|<1/2^{n+2}\rforal a\in {\cal F}_n.
\eneq
Define $\Phi_n: A\to M_{r(n)}$ by $\Phi_n(a)=\Phi_n'(p_nap_n)$ for all $a\in A.$  It is a unital \morp.
Moreover,
\beq\label{MF-2}
\|\Phi_n(p_nap_n)-\Phi_n(c(a))\|<1/2^{n+1}\rforal a\in {\cal F}_n,
\eneq
$n=1,2,...,\infty.$
Combining with (\ref{MF-1}), we obtain that
\beq\label{MF-2ad1411}
 \|\Phi_n(f)\|\ge \|f\|-1/2^{n}\tforal f\in {\cal F}_n,\,\,\,\, n=1,2,....
\eneq
%Note that the map $f\mapsto p_nfp_n$ is a \morp.}
%For  any $a\not\in {\cal F},$ choose an element $\Phi_n(a)\in C_n$
%such that $\|\Phi_n(a)\|=\|a\|.$
 Define $\Phi: A\to \prod_{n=1}^{\infty} M_{r(n)}$ by
$\Phi(a)=\{\pi_n'\circ \Phi_n(a)\}$ for all $a\in A.$  Let
$$
\Pi: \prod_{n=1}^{\infty} M_{r(n)}\to \prod_{n=1}^{\infty} M_{r(n)}/\bigoplus_{n=1}^{\infty}M_{r(n)}
$$
be the quotient map.
Put $\Psi=\Pi\circ \Phi.$ One easily checks that $\Psi$ is in fact a unital \hm.
Moreover, by (\ref{MF-4}), $\Psi$ is not zero. Since $A$ is simple,
it is a unital monomorphism.
%
%Since $A$ is unital and simple, it is easy to see that
%$\liminf_{n\to\infty}\|p_nfp_n\|=\|f\|$ for all $f\in A.$
%It follows (\ref{MF-1}) that there is a unital  monomorphism
%$\phi: A\to \prod_{n=1}^{\infty} C_n/\bigoplus_{n=1}^{\infty} C_n.$
%It follows
%There then exists $\mathcal F$-$1/2^{i-1}$-multiplicative linear map
%$\Phi_i:A \to C$. Pick an arbitrary unital homomorphism $\pi_i: C\to M_{n_i}$ for some $n_i$.
%Then the map $\pi_i\circ\Phi_i: A\to M_{n_i}$ is unital and
%$\mathcal F_i$-$1/2^{i-1}$-multiplicative. Then the map
%$(\pi_1\circ\Phi_1, \pi_2\circ\Phi_2, ..., \pi_i\circ\Phi_i, ...)$ induces a
%unital homomorphism from $A$ to $\prod M_{n_i}/\bigoplus M_{n_i}$. Since $A$ is
%simple, the map is an embedding.
\end{proof}

\begin{lem}\label{Popa}
Let ${\cal D}$ be the family of unital separable  residually finite dimensional \CA s  and let
$A$ be a unital simple separable \CA\, which has  property ($L_{\cal D}$) and  property  (SP). Then
$A$ satisfies  the   following Popa condition:
Let $\ep>0$ and let ${\cal F}\subset A$ be a finite subset.
There exists a finite dimensional \SCA\, $F\subset A$  with $P=1_F$ such that
\beq\label{Popa-1}
\|[P, x]\|<\ep,\,\,\, PxP\in_{\ep}F\andeqn \|PxP\|\ge \|x\|-\ep
\eneq
for all $x\in {\cal F}.$  In particular, if $A\in {\cal B}_1,$ then $A$ satisfies the Popa condition.
\end{lem}

\begin{proof}
We may assume that ${\cal F}\subset A^{\bf 1}$ and $0<\ep<1/2.$
Without loss of generality, we may assume that
$$
d=\min\{\|x\|: x\in {\cal F}\}>0.
$$
%Since $A\in \mathcal B_1$, there is a projection $p\in A$ and a \SCA\, $C\subseteq A$ with $C\in\mathcal C$ and $p=1_C$ such that
Since $A$ has property ($L_{\cal D}$), there is a projection $p\in A$ and a \SCA\, $D\subseteq A$ with $D\in\mathcal {\cal D}$ and $p=1_D$ such that

\beq\label{popa-0}
\|px-xp\|<d\epsilon/16,
\,\, pxp\in_{d\epsilon/16} D\andeqn
\|pxp\|\ge (1-\ep/16)\|x\|
\eneq
for all $x\in {\cal F}$ (see  \ref{Adset9201411}).

Let ${\cal F}'\subset D$ be a finite subset such that, for each $x\in {\cal F},$ there exists $x'\in {\cal F}'$ such
that $\|pxp-x'\|<d\ep/16.$
Since $D\in {\cal D},$ there is a unital surjective \hm\,
$\pi: D\to D/{\rm ker}\pi$ such that $F_1:=D/{\rm ker}\pi$ is a finite dimensional \CA\, and
\beq\label{popa-2}
 \|\pi(x')\|\ge (1-\ep/16)\|x'\|\rforal x'\in {\cal F}'.
\eneq
%Without loss of generality, we may assume that
%$$
%d=\min\{\|x'\|: x'\in {\cal F}'\}>0.
%$$
Let $B=\overline{{\rm ker}\pi A{\rm ker}\pi}.$ $B$ is a hereditary \SCA\, of $A.$
Let $C$ be the closure of $D+B.$
Note that $1_C=1_D=p.$ As in the proof of 5.2 of \cite{LnTAF},   $B$ is an ideal of $C$ and
$C/B\cong D/{\rm ker}\pi=F_1.$
The lemma then follows from Lemma 2.1 of \cite{Niu}. In fact, since $pAp$ has property (SP), by
Lemma 2.1 of \cite{Niu}, there a projection $P\in pAp$  and a monomorphism $h: F_1\to PAP$ such that
\beq\label{popa-n3}
h(1_{F})=P,\,\,\, \|Px'-x'P\|<\ep/16\\
\andeqn \|h\circ \pi(x')-Px'P\|<\ep\cdot d/16
%<(\ep/16)\|x\|
\eneq
for all  $x'\in {\cal F}'.$

Put $F=h(F_1).$ Then, one estimates that, for all $x\in {\cal F},$
\beq\label{popa-n4}
&&\|Px-xP\|\le \|Ppx-Px'\|+\|Px'+x'P\|+\|x'P-xpP\|\\
&&<\ep/16+\ep/16 +\ep/16<\ep,\\
%&&\andeqn
&&PxP\approx_{\ep/16} Px'P\in_{\ep/16} {{F_1}}\andeqn\\
&& \|PxP\|=\|PpxpP\|\ge \|Px'P\|-d\ep/16\ge \|h\circ \pi(x')\|-d\ep/8\\
&&=\|\pi(x')\|-d\ep/8
\ge \|x'\|-d\ep/16-d\ep/8\ge \|pxp\|-d\ep/4\\
&&\ge (1-\ep/16)\|x\|-d\ep/4
\ge \|x\|-\ep.
\eneq
\end{proof}

\begin{thm}\label{B1hered}
Let $A\in {\cal B}_1$ (or $A\in {\cal B}_0$). Then, for any projection $p\in A,$ one has that $pAp\in {\cal B}_1$ (or $pAp\in {\cal B}_0$).
\end{thm}

\begin{proof}
Let $1/4>\ep>0,$ let $a\in (pAp)_+\setminus \{0\}$ and let ${\cal F}\subset pAp$ be a finite subset.
Since $A$ is unital and simple, there are $x_1,x_2, ...,x_m\in A$ such that
\beq\label{B1hered-1-1}
\sum_{i=1}^m x_i^*px_i=1_A.
\eneq
Put ${\cal F}_1=\{p, x_1,x_2,...,x_m, x_1^*,x_2^*,...,x_m^*\}\cup {\cal F}.$
Let $K=m^2\max\{\|x\|: x\in {\cal F}_1\}.$
  Since $A\in {\cal B}_1,$ there is a projection $e\in A$ and a unital \SCA\, $C_1\in {\cal C}$ (or $C_1\in {\cal C}_0$)
with $1_{C_1}=e$ such that
\beq\label{B1hered-1}
\|xe-ex\|&<&\ep/64(K+1)\tforal x\in {\cal F}_1,\\
{\rm dist}(exe, C_1)&<&\ep/64(K+1)\tforal x\in  {\cal F}_1\andeqn\\
1-e &\lesssim & a.
\eneq
Since $p\in {\cal F}_1,$ there is a projection $q\in C_1$ such that
\beq\label{B1hered-2}
\|epe-q\|<\ep/64(K+1).
\eneq
It follows that
$$
\|pep-q\|<\ep/32(K+1).
$$
Moreover, there are $y_1,y_2,...,y_m\in C_1$ such that
\beq\label{B1hered-2+1}
\|\sum_{i=1}^m y_i^*qy_i-e\|<\ep.
\eneq
It follows that $q$ is full in $C_1.$ It follows from \ref{cut-full-pj} that
$qC_1q\in {\cal C}$(or $qC_1q\in {\cal C}_0$).
%and $K_1(qC_1q)=\{0\}.$
There is a unitary $u\in A$ such that
$$
\|u-1\|<\ep/16(K+1)\andeqn u^*qu\le p.
$$
Put $q_1=u^*qu$ and
$C=u^*qC_1qu.$ Then $C\in {\cal C}$ (or $C\in {\cal C}_0$) and $1_C=q_1.$
We also have
\beq\label{B1hered-3=}
\|epe-q_1\|<\ep/64(K+1)+\ep/8(K+1)=9\ep/64(K+1).
\eneq
If $x\in {\cal F},$
then
\beq\label{B1hered-3}
\|q_1x-xq_1\| &\le& 2\|(q_1-epe)x\|+\|epex-xepe\|\\
&<& 18\ep /64+ \ep/16(K+1)<\ep\tforal x\in {\cal F}.
\eneq
Similarly, we estimate that
\beq\label{B1hered-4}
{\rm dist}(qxq, C)<\ep\tforal x\in {\cal F}.
\eneq
We also have
\beq\label{B1hered-5}
\|(p-q_1)-(p-pep)\|=\|q_1-pep\|<\ep/32(K+1)+\ep/8(K+1)=5\ep/32(K+1).
\eneq
Put $\eta=5\ep/32(K+1)<1/16.$
Let $f_{\eta/2}(t)\in C_0((0,\infty)$ be as in \ref{Dball}.
Then, by 2.2 of \cite{RorUHF2},
\beq\label{B1hered-6}
p-q_1=f_{\eta}(p-q_1)\lesssim p-pep\lesssim 1-e\lesssim a.
\eneq
This shows that $pAp\in {\cal B}_1.$
\end{proof}

%{\bf I inserted a proof for the following theorem. N-----I made significant changes here. But have not delete old stuff. {Now it is deleted---L.}}

\begin{thm}\label{Comparison}
Let ${\cal D}$ be a class of until \CA s  which is closed under tensor products with a finite dimensional C*-algebra and which satisfies the strict comparison property for positive elements (see \ref{DW(A)}).
 Let $A$ be a unital simple C*-algebra in the class $\mathrm{TA}{\cal D},$  then $A$ has strict comparison for positive elements. In particular, if $A\in\mathcal B_1$, then $A$ has strictly comparison for positive elements and
$K_0(A)$ is weakly unperforated.
\end{thm}

%\begin{proof}
%This follows from the proof of Theorem 5.5 of \cite{Lin-LAH}.
%Note the proof is easier than that of Theorem 5.5 of \cite{Lin-LAH}.
%In fact, the \SCA\, $C$ appear in that proof will be in ${\cal D}$ and have strict comparison. In particular,
%we do not need to worry about the value of $9d(X)/m$ here.
%\end{proof}
%
%{{\bf Zhuang: I suggest that the following is reduced to a corollary in the light of the above.}}
%{\color{Green} In fact, the proof of \ref{Comparison} is in fact a proof of \ref{ComparisonR}}
%\begin{thm}\label{Comparison}
%Let $A\in {\cal B}_1.$ Then $A$ has strictly comparison for positive elements and
%$K_0(A)$ is weakly unperforated.
%\end{thm}

\begin{proof}
%Let $a, b\in A_+\setminus\{0\}.$ Suppose that
%$$
%d_{\tau}(a)<d_{\tau}(b)\tforal \tau\in T(C).
%$$
%Since $A$ is simple.

%Note that although $A$ may not be exact, any quasi-trace of $A$ is indeed a trace. (This follows from the fact that any C*-algebras in $\mathcal C_0$ has this property. We leave the proof to the readers).
By a result of R\o dram (see, for example,  Corollary 4.6 of \cite{Ror-srZ}), to show that $A$ has strict comparison for positive elements, it is enough to show that $W(A)$ is almost unperforated, i.e., for any positive elements $a, b$ in a matrix algebra over $A$, if $(n+1) [a] \leq n [b]$ for some $n\in\mathbb N$, then $[a]\le [b]$.

Let $a, b$ be such positive elements. Since any matrix algebra over $A$ is still in $\mathrm{TA}{\mathcal D}$, let us assume that $a, b\in A$.

First we consider the case that $A$ does not have (SP) property. In this case, by the proof of  \ref{B1sp},
$A=\overline{\cup_{n=1}^{\infty} A_n},$ where $A_n\in {\cal D}$ ($\{A_n\}$ be not be increasing).

Without loss of generality, we may assume that $0\le a, \, b\le 1.$
 Let $\ep>0.$ It follows from an argument of Rordam (see Lemma 5.6 of \cite{Niu-MD}) that there exists an integer $m\ge 1,$
 $a', b'\in A_m$ such that
 \beq\label{Compr-n-1}
 \|a'-a\|<\ep/2,\,\,\, \|b'-b\|<\ep/2,\,\,\, b'\lesssim b \andeqn\\
 {\rm diag}(\overbrace{f_{\ep/2}(a'), f_{\ep/2}(a'),...,f_{\ep/2}(a')}^{n+1})\lesssim
 {\rm diag}(\overbrace{b',b',...,b'}^n)\,\,\,{\rm in}\,\,\, A_m.
 \eneq
Since $A_m$ has strict comparison (see part (b) of \ref{2Tg16}), one has
\beq\label{Compr-n-2}
f_{\ep/2}(a')\lesssim b'\,\,\,{\rm in}\,\,\, A_m.
\eneq
It follows, using 2.1 of \cite{RorUHF2}, that
\beq\label{Compri-n-3}
f_\ep(a)\lesssim f_{\ep/2}(f_{\ep/2}((a))\lesssim f_{\ep/2}(a')\lesssim b'\lesssim b
\eneq
for every $\ep>0.$
It follows that $a\lesssim b.$

Now we assume that $A$ has (SP).
Let $1/4>\ep>0.$  We may further assume that $\|b\|=1.$
Since $A$ has (SP) and simple, there are mutually orthogonal and mutually equivalent non-zero projections
$e_1,e_2,...,e_{n+1} \in \overline{f_{3/4}(b)Af_{3/4}(b)}.$
Put $E=e_1+e_2+\cdots +e_{n+1}.$
By 2.4 of \cite{RorUHF2}, we also have that
\beq\label{Compr-n2}
(n+1)[f_{\ep/2}(a)]\le n[f_{\dt}(b)]
\eneq
for some $\ep>\dt>0.$  Put $\eta=\min\{\ep/4, \dt/4, 1/8\}.$
It follows from  the definition \ref{DB1} that there is a \SCA\, $C=pAp\oplus S$ with $S\in\mathcal D$ and
$a', b', E', e_i'\in C$ ($i=1,2,...,n+1$) such that
$0\le a', b'\le 1$ and $E', e_i'$ are projections in $C$,
\beq\label{Compri-n3}
||a-a'||<\eta,\  b'\lesssim f_{\dt}(b), \|f_{1/2}(b') E'-E'\|<\eta,\\
E'=\sum_{i=1}^{n+1}e_i',\,\,\,\|e_i-e_i'\|<\eta \andeqn  \|E-E'\|<\eta<1,
\eneq
and
\begin{equation}\label{Compri-n4}
{\rm diag}(\overbrace{f_{\ep/2}(a'),f_{\ep/2}(a'),...,f_{\ep/2}(a')}^{n+1}) \lesssim
{\rm diag}(\overbrace{b',b',...,b'}^n)\andeqn
(n+1)[E']< [b'] \quad \textrm{in $C$}
\end{equation}
(see Lemma 5.6 of \cite{Niu-MD}).
Moreover, the projection $p$ can be chosen so that $p\lesssim e_1$.
From (\ref{Compri-n3}), there is a projection ${e_i''}$, ${E}''\in \overline{f_{1/2}(b')Cf_{1/2}(b')}$
($i=1,2,...,n+1$) such that
$\|E'-E''\|<2\eta,$  $\|e_i''-e_i'\|<2\eta,$ $i=1,2,...,n+1$ and
$E''=\sum_{i=1}^{n+1}e_i''$ (we also assume that $e_1'',e_2'',...,e_{n+1}''$ are mutually orthogonal).
Note that $e_i'$ and
$e_i''$ are equivalent.
Choose a function $g\in C_0((0,1])_+$ with $g\le 1$ such that $g(b')f_{1/2}(b')=f_{1/2}(b')$
and $[g(b')]=[b']$ in $W(C).$ In particular,
$g(b')E''=E''.$

Write
$$a'=a'_0 \oplus a'_1,\quad g(b')=b'_0 \oplus b'_1, e_i''=e_{i,0}\oplus e_{i,1}\quad\textrm{and}\quad E''=E'_0 \oplus E'_1$$
with $a'_0, b'_0, E_0', e_{i,0}\in pAp$ and $a'_1, b'_1, e_{i,1},\, E_1' \in S,\, i=1,2,...,n+1$.
Note that $E_1'b_1'=E_1'b_1'=E_1'.$
This, in particular, implies that
\begin{equation}\label{12/27/04-comp-1}
\tau(b_1')\ge (n+1)\tau(e_{1,1})\rforal \tau\in T(S).
\end{equation}
It follows from \eqref{Compri-n4} that
$$
d_\tau(f_{\ep/2}(a_1')) \leq \frac{n}{n+1} d_\tau (b_1'),
%=\frac{n}{n+1} \tau ({g(b')}) ,
\rforal \tau\in\mathrm{T}(S).
$$
%However, since $(n+1)[e'']< [g(b')]$ in $C$ (hence in $S$), one has that for any $\tau\in \mathrm T(S)$,
Since $(b_1'-e_{1,1})e_{1,1}=0$ and $b_1'=(b_1'-e_{1,1})+e_{1,1},$ for all $\tau\in T(S),$
$$
d_\tau((b'_1-e_{1,1}))=d_\tau(b'_1)-\tau(e_{1,1})>d_\tau(b'_1)-\frac{1}{n+1}d_\tau(b'_1)\geq d_\tau(f_{\ep/2}(a_1'))).
$$
Since $S$ has the strict comparison, one has
$$
f_{\ep/2}(a_1')_+\lesssim  (b_1'-e_1'),
$$
and therefore
\beq\label{Compari-n6}
&&f_{\ep}(a)\lesssim f_{\ep/2}(a')\lesssim  p\oplus f_{\ep/2}(a_1')\lesssim p\oplus(b'_1-e_{1,1})\\
&& \lesssim  e_1 \oplus(b_1'-e_{1,1})\lesssim e_1\oplus (b_1'-e_{1,1})+(b_0'-e_{1,0})\\
&&\sim e_1''\oplus (g(b')-e_1'')\sim g(b')\sim b'\lesssim b.
\eneq
Since $\epsilon$ is arbitrary, one has that $a\lesssim b.$

Hence one always has that $a\lesssim b$, and therefore $W(A)$ is almost unperforated.
%Since each quasi-trace of $A$ is a trace, one concludes that $A$ has strict comparison on positive elements.
\end{proof}

\begin{lem}\label{Tapprox}
Let $\mathcal D$ be a class of unital amenable C*-algebras, let $A$ be a separable unital C*-algebra which
is $TA{\cal D}$ and
% can be tracially approximated by the C*-algebras in $\mathcal D$.
let $C$ be a unital  ( amenable)  C*-algebra.

Let $\mathcal F, \mathcal G\subseteq C$ be finite subsets, let $\epsilon>0$ and $\delta>0$ be positive. Let $\mathcal H\subseteq C_+^{\bf 1}$ be a finite subset, and let $T: C_+\setminus\{0\}\to \mathbb R_+\setminus\{0\}$ and $N: C_+\setminus\{0\}\to\mathbb N$ be maps. Let $\Delta: C_+^{q, 1}\setminus\{0\}\to (0, 1)$ be an order preserving map. Let $\mathcal H_1\subseteq C_{+}^{\bf 1}\setminus \{0\}$, $\mathcal H_2\subseteq C_{s.a.}$ and $\mathcal U\subseteq U(C)/CU(C)$ be finite subsets. Let $\sigma_1>0$ and $\sigma_2>0$ be constants.   Let $\phi, \psi: C\to A$ be two unital $\mathcal G$-$\dt$-multiplicative linear maps such that
\begin{enumerate}
\item\label{cond-708-001} $\phi$ and $\psi$ are $T\times N$-$\mathcal H$-full  (see the definition \ref{Dfull}),
\item\label{cond-708-002} $\tau\circ\phi(c)>\Delta(\hat{c})$ and $\tau\circ\psi(c)>\Delta(\hat{c})$ for any $c\in\mathcal H_1$,
\item\label{cond-708-003} $|\tau\circ\phi(c)-\tau\circ\psi(c)|<\sigma_1$ for any $\tau\in \mathrm{T}(A)$ and any $c\in \mathcal H_2$,
\item\label{cond-708-004} $\mathrm{dist}(\phi^{\ddagger}(u), \psi^{\ddagger}(u))<\sigma_2$ for any $u\in\mathcal U$.
\end{enumerate}

Then, for any finite subset $\mathcal F'\subseteq A$ and $\epsilon'>0$, there exists a C*-subalgebra $D\subseteq A$ with $D\in\mathcal D$ such that if $p=1_D$, then, for any $a\in\mathcal F'$,
%\begin{enumerate}
%\item

\indent
{\rm (a)} $\|pa-ap\|<\epsilon'$,
%\item

\indent
{\rm (b)} $pap\in_{\epsilon'} D$,
%\item

\indent
{\rm (c)} $\tau(1-p)<\epsilon'$, for any $\tau\in T(A)$.
%\end{enumerate}

There are also (completely positive) linear map $j_0: A \to (1-p)A(1-p)$ and
a unital \morp\, $j_1: A\to D$ such that 
\beq\nonumber
j_0(a)=(1-p)a(1-p) \rforal a\in A\andeqn\\\nonumber
\|j_1(a)-pap\|<3\epsilon' \rforal a\in {\mathcal F}.
\eneq

Moreover, define
$$
\phi_0=j_0\circ\phi,\,\, \psi_0=j_0\circ\psi,
\phi_1=j_1\circ\phi\quad\textrm{and}\quad \psi_1=j_1\circ\psi.
$$
With a sufficiently large $\mathcal F'$ and small enough $\epsilon'$, one has that $\phi_0$, $\psi_0$, $\phi_1$ and $\psi_1$ are $\mathcal G$-$2\dt$-multiplicative and
\begin{enumerate}\setcounter{enumi}{4}
\item\label{concl-708-001} $\|\phi(c)-(\phi_0(c)\oplus\phi_1(c))\|<\epsilon$ and $\|\psi(c)-(\psi_0(c)\oplus\psi_1(c))\|<\epsilon$, for any $c\in\mathcal F$,
\item\label{concl-708-002} $\phi_0, \psi_0$ and $\phi_1, \psi_1$ are $2T\times N$-$\mathcal H$-full,
\item\label{concl-708-003} $\tau\circ\phi_1(c)>\Delta(\hat{c})/2$ and $\tau\circ\psi_1(c)>\Delta(\hat{c})/2$ for any $c\in\mathcal H_1$,
\item\label{concl-708-004} $|\tau\circ\phi_1(c)-\tau\circ\psi_1(c)|<2\sigma_1$ for any $\tau\in \mathrm{T}(D)$ and any $c\in \mathcal H_2$,
\item\label{concl-708-005} $\mathrm{dist}(\phi_{{i}}^{\ddagger}(u), \psi_{{i}}^{\ddagger}(u))<2\sigma_2$ for any $u\in\mathcal U,$
{$i=0,1.$}
\end{enumerate}
\end{lem}

\begin{proof}
Without loss of generality, one may assume that each element of $\mathcal F$, $\mathcal G,$
${\cal H}_2,$ or $\mathcal F'$ has norm at most one and $1_A\in \mathcal F'$.

Since $\phi$ and $\psi$ are $T\times N$-$\mathcal H$-full, for each $h\in\mathcal H$, there are $a_{1,h}, ..., a_{N(h),h}$ and $b_{1,h}, ..., b_{N(h),h}$ in $A$ with $\|a_{i,h}\|, \|b_{i,h}\|\le T(h)$ such that
\beq\label{911-n1}
\sum_{i=1}^{N(h)} a^*_{i,h}\phi(h)a_{i,h}=1_A \quad\textrm{and}\quad \sum_{i=1}^{N(h)} b^*_{i,h}\psi(h)b_{i,h}=1_A.
\eneq
Put
$
d_0=\min\{\Delta(\hat{h}): h\in {\cal H}_1\}.
$
By  (\ref{cond-708-003}), there are, for each $c\in {\cal H}_2,$
$x_{1,c},x_{2,c},...,x_{t(c), c}\in A$ such that
\beq\label{911-n2}
\|\sum_{i=1}^{t(c)} x_{i, c}^*x_{i,c}-\phi(c)\|<\sigma_1\andeqn
\|\sum_{i=1}^{t(c)} x_{i,c}x_{i,c}^*-\psi(c)\|<\sigma_1.
\eneq
Let
$$
t({\cal H}_2)=\max \{\|x_{i,c}\|: 1\le i\le t(c): c\in {\cal H}_2\}.
$$
For the given finite subset $\mathcal F'\subseteq A$ and given $\epsilon'>0$, since $A$ can be tracially approximated by the \CA s
%\, C*-algebras
in the class $\mathcal D$, there exists a C*-subalgebra $D\subseteq A$ with $D\in\mathcal D$ such that if $p=1_D$, then, for any $a\in\mathcal F'$,
%\begin{enumerate}
%\item

(i) $\|pa-ap\|<\epsilon'$,
%\item

(ii) $pap\in_{\epsilon'} D$,

(iii)  $\tau(1-p)<\epsilon'$, for any $\tau\in T(A)$.
%\end{enumerate}

To make ${\cal F}'$ large and $\ep'$ small,
we may assume that ${\cal F}'$ also contains ${\cal F},$ ${\cal G},$ ${\cal H},$
$\phi({\cal G}\cup {\cal F}),$ $\psi({\cal G}\cup {\cal F}),$
$\phi({\cal H}),$ $\psi({\cal H}),$ ${\cal H}_1,$ ${\cal H}_2,$
$x_{i,c}, x_{i,c}^*,$ $i=1,2,...,t(c)$ and $c\in {\cal H}_2,$ as well as
$a_{i,h}, a_{i,h}^*, b_{i,h}, b_{i,h}^*,$
$i=1,2,...,N(h)$ and $h\in {\cal H},$ and
$$\ep'<\min\{\min\{1/64(T(h)+1)(N(h)+1): h\in {\cal H}\}, \ep, \dt,d_0, \sigma_1, \sigma_2\}/64(t({\cal H}_2)+1)^2.$$

For each $a\in\mathcal F'$, choose $d_a\in D$ such that $\|pap-d_a\|<\epsilon'$ (choose $d_{1_A}=1_D$).  Consider the finite subset $\{d_ad_b:\ a, b\in\mathcal F'\}\subseteq D$. Since $D$ is an amenable \SCA\,  of $pAp$, there is unital completely positive linear map $L: pAp\to D$ such that
$$\|L(d_ad_b)-d_ad_b\|<\epsilon',\quad a, b\in\mathcal F'.$$
Define $j_1: A\to D$ by $j_1(a)=L(pap).$
Then, for any $a\in\mathcal F'$, one has
\begin{eqnarray}
\|j_1(a)-pap\| & = &\| L(pap)-pap\|=
%& = &
\|L(d_a)-d_a\| + 2\epsilon'={{2\ep'}}.
%& = & 3\epsilon'.
\end{eqnarray}
Note that $j_0$ and $j_1$ are  {{${\cal F'}$-$7\epsilon'$}}
%-$\mathcal F'$
-multiplicative, and
$$\|a-j_0(a)\oplus j_1(a)\|<4\epsilon' \rforal  a\in\mathcal F'.$$

%Consider the maps
Define
$$
\phi_0=j_0\circ\phi,\,\, \psi_0=j_0\circ\psi,
\phi_1=j_1\circ\phi\quad\textrm{and}\quad \psi_1=j_1\circ\psi.
$$
Then (by the choices of ${\cal F}'$ and $\ep'$),
% by choosing $\mathcal F'$ sufficiently large (containing $\phi(\mathcal G\cup\mathcal F)\cup\psi(\mathcal G\cup\mathcal F)$) and $\epsilon'$ sufficiently small (less than $\min\{\delta/7, \epsilon\}$),
maps $\phi_0$, $\psi_0$, $\phi_1$, and $\phi_1$ are $\mathcal G$-$2\dt$-multiplicative, and for any $c\in\mathcal F$,
\begin{equation}\label{app-709-0002}
\|\phi(c)-(\phi_0(c)\oplus\phi_1(c))\|<\epsilon\quad\textrm{and}\quad \|\psi(c)-(\psi_0(c)\oplus\psi_1(c))\|<\epsilon.
\end{equation}

%Since $\phi$ and $\psi$ are $T\times N$-$\mathcal H$-full, for each $h\in\mathcal H$, there are $a_1, ..., a_{N(h)}$ and $b_1, ..., b_{N(h)}$ in $A$ with $\|a_i\|, \|b_i\|<T(h)$ such that
%$$\sum_{i=1}^{N(h)} a^*_i\phi(h)a_i=1_A \quad\textrm{and}\quad \sum_{i=1}^{N(h)} b^*_i\psi(h)b_i=1_A.$$
Apply $j_1$ on both sides of both equations in (\ref{911-n1}).
%By increase $\mathcal F'$ and decrease $\epsilon'$ (recall that $j_1$ is $7\epsilon'$-$\mathcal F'$-multiplicative),
One obtains two invertible elements  $e_h:=\sum_{i=1}^{N(h)} j_1(a_i^*)\phi_1(h)j_1(a_i)$ and $f_h:=\sum_{i=1}^{N(h)} j_1(b^*_i)\psi_1(h)j_1(b_i)$  such that
$$ |\|e_h^{-\frac{1}{2}}\|-1|<1 \quad\textrm{and}\quad  |\|f_h^{-\frac{1}{2}}\|-1|<1.$$
Note that
$$\sum_{i=1}^{N(h)} e_h^{-\frac{1}{2}}j_1(a_i^*)\phi_1(h)j_1(a_i)e_h^{-\frac{1}{2}}=1_D,
%\quad\textrm{and}\quad
\sum_{i=1}^{N(h)} f_h^{-\frac{1}{2}}j_1(b_i^*)\psi_1(h)j_1(b_i)f_h^{-\frac{1}{2}}=1_D,$$
%and
$$\|j_1(a_i)e_h^{-\frac{1}{2}}\|<2 T(h) \quad\mathrm{and}\quad \|j_1(b_i)f_h^{-\frac{1}{2}}\|<2 T(h).$$
Therefore, $\phi_1$ and $\psi_1$ are $2T\times N$-$\mathcal H$-full, and this proves \eqref{concl-708-002}. The same calculation also shows that $\phi_0$ and $\psi_0$ are $2T\times N$-$\mathcal H$-full.
Note that we have shown (\ref{concl-708-001}) and (\ref{concl-708-001}) hold.
%Note that increasing $\mathcal F'$ and decreasing $\epsilon'$ preserve \eqref{app-709-0002}.
Since $\ep'<d_0/16,$ it is also easy to check that  \eqref{concl-708-003} holds.

To see \eqref{concl-708-004}, one notes that
\beq\label{911-n4}
\|\sum_{i=1}^{t(c)} d_{x_{i,c}}^*d_{x_{i,d}}-\phi_1(c)\|<2\sigma_1\andeqn
\|\sum_{i=1}^{t(c)}d_{x_{i,c}}d_{x_{i,d}}^*-\psi_1(c)\|<2\sigma_1
\eneq
for all $c\in {\cal H}_2.$ Then \eqref{concl-708-004} also holds.

Let us show that \eqref{concl-708-005} holds with sufficiently large $\mathcal F'$ and sufficiently small $\epsilon'$.

Choose unitaries $u_1, u_2, ..., u_n\in C$ such that $\mathcal U=\{\overline{u_1}, \overline{u_2}, ..., \overline{u_n}\}$. Pick unitaries $w_1, w_2, ..., w_n\in A$ such that
each $w_i$ is a commutator and
$$\mathrm{dist}(\left<\phi(u_i)\right>\left<\psi(u^*_i)\right>, w_i)<\sigma_2.$$
%where $\left<z\right>=z(z^*z)^{-\frac{1}{2}}$ for any invertible $z$.
Choose $\mathcal F'$ sufficiently large and $\epsilon'$ sufficiently small such that there are commutators $w_1', w_2, ..., w_n'\in CU(D)$  and commutators $w_1'', w_2'',...,w_n''\in (1-p)A(1-p)$ satisfying
$$\|j_1(w_i)-w_i'\|<\sigma_2/2 \andeqn \|j_0(w_i)-w_i''\|<\sigma_2/2, \quad 1\leq i\leq n,$$
(see Appendix of \cite{Lin-LAH}) and
$$\| \left<\phi_{k}(u_i)\right>\left<\psi_{k}(u^*_i)\right> - j_{k}(\left<\phi(u_i)\right>\left<\psi(u^*_i)\right>)\|<\sigma_2/2, \quad 1\leq i\leq n \andeqn k=0,1.$$
Then
\begin{eqnarray}
&&\|  \left<\phi_{k}(u_i)\right>\left<\psi_{k}(u^*_i)\right> - w_i'  \| \\
& \leq & \|  \left< \phi_{k}(u_i)\right>\left<\psi_{k}(u^*_i)\right> - j_{k}(w_i) \| + \| j_{k}(w_i)  -w_i'\| \\
&\leq & \|  \left< \phi_{k}(u_i)\right>\left<\psi_{k}(u^*_i)\right> - j_{k}(\left<\phi(u_i)\right>\left<\psi(u^*_i)\right>)\| + \\
&& \| j_{k}(\left<\phi(u_i)\right>\left<\psi(u^*_i)\right>)-j_{k}(w_i)\| + \sigma_2/2\\
&\leq& 2\sigma_2,\,\,\, k=0,1.
\end{eqnarray}
This proves \eqref{concl-708-005}.
\end{proof}

\section{${\cal Z}$-stability}

\begin{lem}\label{Affon1}
Let $A\in {\cal B}_1$ {\rm (or ${\cal B}_0$)} be a unital infinite dimensional simple \CA. Then,  for any $\ep>0,$ any $a\in A_+\setminus\{0\},$ any finite subset ${\cal F}\subset A$ and any
integer $N\ge 1,$ there exists a projection $p\in A$ and a \SCA\, $C\in {\cal C}$  (or ${\cal C}_0$)  with $1_C=p$ satisfies the following:

{\rm (1)} ${\rm dim}(\pi(C))\ge N^2$ for every irreducible representation $\pi$ of $C;$
%{\color{Green} (why use $N^2$ rather than $N$?)}

{\rm (2)} $\|px-xp\|<\ep$ for all $x\in {\cal F};$

{\rm (3)} ${\rm dist}(pxp, C)<\ep$ for all $x\in {\cal F}$ and

{\rm (4)} $1-p\lesssim a.$

\end{lem}

\begin{proof}
Since $A$ is an  infinite dimensional simple \CA, there are
$N+1$ mutually orthogonal non-zero positive elements
$a_1, a_2,...,a_{N+1}$ in $A.$  Since $A$ is simple,
there are $x_{i,j}\in A,$ $j=1,2,...,k(i),$ $i=1,2,...,N+1,$ such that
$$
\sum_{j=1}^{k(i)}x_{i,j}^*a_ix_{i,j}=1_A.
$$
Let
$$
K=(N+1)\max\{\|x_{i,j}\|+1: 1\le j\le k(i),\,\,\, 1\le i\le N+1\}.
$$
Let $a_0\in A_+\setminus \{0\}$ be such that
$a_0\lesssim a_i$ for all $1\le i\le N+1.$ Since $\overline{a_0Aa_0}$ is also an infinite dimensional simple \CA, one obtains $a_{01}, a_{02}\in \overline{a_0Aa_0}$ which are mutually orthogonal and nonzero. One then
obtains a non-zero element $a\in \overline{a_{01}Aa_{01}}$ such that $a\lesssim a_{02}.$

Let
$$
{\cal F}=\{a_i: 1\le i\le N+1\}\cup\{x_{i,j}: 1\le j\le k(i),\,1\le i\le N+1\}\cup \{a\}.
$$
Now since $A\in {\cal B}_1,$ there is a projection $p\in A$ and $C\in {\cal C}$ with $1_C=p$ such that

(1) $\|xp-px\|<\min\{1/2,\ep\}/2K$ for all $x\in {\cal F},$

(2) ${\rm dist}(pxp, C)<\min\{1/2, \ep\}/2K$ for all $x\in {\cal F}$, and

(3) $1-p\lesssim a.$

Thus, with a standard computation, we obtain mutually orthogonal
 non-zero positive elements $b_1,b_2,...,b_{N+1}\in C$ and
 $y_{i,j}, \in C$ ($1\le j\le k(i)),$ $i=1,2,...,N+1,$ such that
\beq\label{Affon1-2}
\|\sum_{j=1}^{k(i)}y_{i,j}^*b_iy_{i,j}-p\|<\min\{1/2, \ep/2\}.
\eneq
For each $i,$ we find another element $z_i\in C$ such that
\beq\label{Affon1-3}
\sum_{j=1}^{k(i)}z_i^*y_{i,j}b_iy_{i,j}z_i=p.
\eneq
Let $\pi$ be an irreducible representation of $C.$ Then
by (\ref{Affon1-3}),
\beq\label{Affon1-4}
\sum_{j=1}^{k(i)}\pi(z_i^*y_{i,j})\pi(b_i)\pi(y_{i,j}z_i)=\pi(p).
\eneq
Therefore $\pi(b_1), \pi(b_2),...,\pi(b_{N+1})$ are mutually orthogonal
non-zero positive elements in $\pi(A).$ Then (\ref{Affon1-4}) implies
that $\pi(C)\cong M_n$ with $n\ge N+1.$ This proves the lemma.

\end{proof}

\begin{cor}\label{CTdense}
Let $A\in {\cal B}_1$ be a unital simple \CA. Then, for any $\ep>0$ and $f\in \Aff(T(A))_{++},$ there exists
a \SCA\, $C\in {\cal C}_0$ in $A,$ an element $c\in C_+$ such that
\beq\label{CTdense-1}
{\rm dim}\pi(C)\ge  (4/\ep)^2\,\,\,\text{for \,each \,irreducible\,representation} \,\,\pi\,\,\,of \,\,C,\\
0<\tau(f)-\tau(c)<\ep/2\rforal \tau\in T(A).
\eneq
\end{cor}

The following is  known.

\begin{lem}\label{Affon1+}
Let $C=M_n([0,1])$ and $g\in {\rm LAff}_b(T(C))_+.$ Then there exists $a\in
C_+$ with $0\le a\le 1$ such that
$$
0\le g(t)-d_t(a)\le 1/n \tforal t\in [0,1],
$$
where $d_t(a)=\lim_{k\to\infty} a^{1/k}(t)$ for all $t\in [0,1].$
\end{lem}

\begin{proof}
We will use the proof of Lemma 5.2 of \cite{BPT}.
For each $0\le i\le  n-1,$ define
$$
X_i=\{t\in [0,1]: g(t)>i/n\}.
$$
Since $g$ is lower semi-continuous, $X_i$ is open in $[0,1].$ There
is a continuous function $g_i\in C([0,1])_+$ with $0\le g_i\le 1$
such that
$$
\{t\in [0,1]: g_i(t)\not=0\}=X_i,\,\,\,i=0,1,...,n-1.
$$
Let $e_1,e_2,...,e_n$ be $n$ mutually orthogonal rank one
projections in $C=M_n(C([0,1]).$ Define \beq\label{Affon1++1}
a=\sum_{i=1}^{n-1} g_i e_i\in C. \eneq Then $0\le a\le 1.$ Put
$Y_i=\{t\in [0,1]: (i+1)/n\ge g_i(t)>i/n\}=X_i\setminus
\bigcup_{j>i} X_j,$ $i=0,1,2,...,n-1.$ These are mutually disjoint
sets. Note that
$$
[0,1]=([0,1]\setminus X_0)\cup \bigcup_{i=0}^{n-1}Y_i.
$$
If $x\in ([0,1]\setminus X_0)\cup Y_0,$ $d_t(a)=0.$ So $0\le g(t)-d_t(a)(t)\le 1/n$ for all such $t.$
If $t\in Y_j,$
\beq\label{Affon1++2}
d_t(a)=j/n.
\eneq
Then
\beq\label{Affon1++3}
0\le g(t)-d_t(a)\le 1/n\tforal t\in Y_j.
\eneq
It follows that
\beq\label{Affon1++4}
0\le g(t)-d_t(a)\le 1/n\tforal t\in [0,1].
\eneq
\end{proof}

\begin{lem}\label{Affon2}
Let $F_1$ and $F_2$ be two finite dimensional \CA s such that each
simple summand of $F_1$ and $F_2$ has rank at least $k$, where $k\ge
1$ is an integer.  Let $\phi_0, \phi_1: F_1\to F_2$ be unital \hm s.
Let $C=A(\phi_0, \phi_1, F_1, F_2).$ Then, for any $f\in {\rm
LAff}_b(T(C))_+$ with $0\le f \le 1,$ there exists a positive
element $a\in M_2(C)$ such that
$$
\max_{\tau\in T(C)}|d_{\tau}(a)-f(\tau)|\le 2/k.
$$
\end{lem}

\begin{proof}
Let $I=\{g\in C: g(0)=g(1)=0\}.$ Note that $C/I$ is a finite
dimensional \CA. Let
$$
T=\{\tau\in T(C):\, {\rm ker}\tau\supset I\}.
$$
Then $T$ may be identified with $T(C/I).$ Let $f\in  {\rm
LAff}_b(T(C))_+$ with $0\le f \le 1.$ It is easy to see that
there exists $b\in (C/I)_+$ for some integer $m_1\ge 1$ such that
\beq\label{Affon2-1} 0\le f(\tau)-d_{\tau}(b) \rforal \tau\in T\andeqn
\max\{f(\tau)-d_\tau(b): \tau\in T\}\le 1/k, \eneq {and furthermore,
if $f(\tau)>0$, then $f(\tau)-d_\tau (b)>0.$}  { As matter of fact,
%in above
we can choose $b$ to be rank $r$ projection in $i^{th}$
block $M_{R(s)}$ of $F_1$ and $\frac
rm<f(\tau_{q(s)})\leq\frac {r+1}m$, where $\tau_{q(s)}$ is the normalized
trace of $M_{R(s)} $ regarded as an element in $T(C/I)$ (since $m\ge k$). For such a
choice, we have $d_\tau (b)= \tau (b)$ for all $\tau$.}
%We may write
%$C/I=\phi_0(F_1)\oplus \phi_1(F_2).$ Accordingly, $b=b_0\oplus b_1.$
{Note that $b=(b_{q(1)}, b_{q(2)},...,b_{q(l)})\in C/I=F_1=M_{R(1)}\oplus M_{R(2)}\oplus\cdots \oplus M_{R(l)}$. Let $b_0=\phi_0(b)\in F_2$ and
$b_1=\phi_1(b)\in F_2.$}
% We will view $b_0, b_1\in M_{m_1}(F_2).$
Write $F_2=M_{r_1}\oplus M_{r_2}\oplus\cdots\oplus M_{r_m}.$ Write
$b_0=b_{0,1}\oplus b_{0,2}\oplus \dots \oplus b_{0,r_m}$ and
$b_1=b_{1,1}\oplus b_{1,2}\oplus\cdots \oplus b_{1,r_m},$ where
$b_{0,j},\,b_{1,j}\in M_{r_j},$ $j=1,2,...,m.$
 Let $\tau_{t,j}=tr_j\circ \Psi_j\circ \pi_t,$ where $tr_j$ is the normalized trace on $M_{r_j},$
$\Psi_j: F_2\to M_{r_j}$ is the quotient map and
$\pi_t: A\to F_2$ is the evaluation at $t\in (0,1).$  Since $f$ is lower semi-continuous on $T(C),$
$$
\liminf_{t\to 0} f(\tau_{t,j})\ge tr_j(b_{0,j})\andeqn \liminf_{t\to
1} f(\tau_{t,j})\ge tr_j(b_{1,j}).
$$
Note that $tr_j(b_{0,j})=\sum_{s=1}^l \af_s tr_{q(s)}(b_{q(s)})$ for some $\af_s\ge 0$ with
$\sum_{s=1}^l\af_s=1.$ Therefore, if $tr_j(b_{0,j})>0,$ then
$$
f(\sum_{s=1}^l \af_s tr_{q(s)})>\sum_{s=1}^l \af_s tr_{q(s)}(b_{q(s)})=tr_j(b_{0,j}).
$$
{Hence, if $tr_j(b_{0,j})>0$ (or $tr_j(b_{1,j})>0$),
then $\liminf_{t\to 0} f(\tau_{t,j})> tr_j(b_{0,j})$ (or
$\liminf_{t\to 1} f(\tau_{t,j})> tr_j(b_{1,j})$).}
 Therefore, there exists $1/8>\dt>0,$ such that
\beq\label{Affon2-2}
f(\tau_{t,j}) &\ge& tr_j(b_{0,j})\tforal t\in (0,2\dt)\andeqn\\
f(\tau_{t,j}) &\ge&  tr_j(b_{1,j})\tforal t\in (1-2\dt,1),
\,\,\,j=1,2,...,l. \eneq Let \beq\label{Affon2-3}
c(t)&=&({\dt-t\over{\dt}})b_0\,\,\,  \text{if}\,\,\,t\in [0,
\dt)\\\label{Affon2-3+1} c(t)&=& 0 \,\,\,\text{if}\,\,\, t\in [\dt,
1-\dt]\andeqn\\\label{Affon2-3+2}
c(t)&=&({t-1+\dt\over{\dt}})b_1\tforal t\in (1-\dt, 1]
\eneq
Note
that $c\in A.$ Define
\beq\label{Affon2-4} g_j(0)&=&
0\\\label{Affon2-4+1} g_j(t)&=&f(\tau_{t,j})-tr_j(b_{0,j})\tforal
t\in (0, \dt]\\\label{Affon2-4+2} g_j(t)&=& f(\tau_{t,j}) \tforal
t\in (\dt, 1-\dt)\\\label{Affon2-4+3}
g_j(t)&=& f(\tau_{t,j})-tr_j(b_{1,j})\tforal  t\in [1-\dt, 1)\andeqn\\
g_j(1)&=&0.
\eneq
One verifies that $g_j$ is lower semi-continuous on
$[0,1].$ It follows \ref{Affon1+} that there exists $a_1\in C([0,1],
F_2)_+$
%$a_1\in M_{m_2}(C([0,1])_+$ for some $m_2\ge 1$
such that
\begin{equation}\label{Affon2-5}
0\le g_j(t)-d_{tr_{t,j}}(a_1)\le 1/r_j\leq1/k \tforal t\in [0,1].
\end{equation}
Note that $a_1(0)=0$ and $a_1(1)=0.$
Therefore  $a_1\in I\subset C.$  Now let $a=c\oplus a_1\in M_2(C).$
 Note that
\begin{eqnarray*}%\label{Affon2-5}
d_{\tau}(a)&=& d_{\tau}(c)+d_{\tau}(a_1)=d_{\tau}(b)\,\,\, \text{if}\,\,\, t\in T,\\
d_{tr_{t,j}}(a)&=& d_{tr_{t,j}}(c)+d_{tr_{t,j}}(a_1)=d_{t,j}(b_0)+d_{tr_{t,j}}(a_1)\tforal t\in (0, \dt),\\
d_{tr_{t,j}}(a)&=& d_{tr_{t,j}}(a_1)\tforal t\in [\dt, 1-\dt]\andeqn\\
d_{tr_{t,j}}(a)&=&d_{tr_{t,j}}(c)+d_{tr_{t,j}}(a_1)=d_{t,j}(b_1)+d_{tr_{t,j}}(a_1) \tforal t\in (1-\dt, 1).
\end{eqnarray*}
Then, combining with \eqref{Affon2-1}, \eqref{Affon2-3}, \eqref{Affon2-3+1}, \eqref{Affon2-3+2}, \eqref{Affon2-4},
\eqref{Affon2-4+1}, \eqref{Affon2-4+2}, \eqref{Affon2-4+3} and \eqref{Affon2-5}, we have
\beq\label{Affon2-6}
0\le f(\tau)-d_{\tau}(a)\le 2/k\tforal \tau\in T \andeqn \tau=tr_{t,j}, j=1,2,...,l, \,\,\,t\in (0,1).
\eneq
Since $T\cup \{tr_{t,j}: 1\le j\le l, \andeqn t\in (0,1)\}$ contains all extremal points of $T(C),$ we conclude that
\beq\label{Affon2-7}
0\le f(\tau)-d_{\tau}(a)\le 2/k\tforal \tau\in T(C).
\eneq
\end{proof}

{
\begin{thm}\label{Aff=W-L}
Let $A\in {\cal B}_1$ be a unital simple \CA. Then
the map $W(A)\to V(A)\sqcup {\rm LAff}_b(A)_{++}$ is surjective.
\end{thm}
 }
 \begin{proof}
The proof follows the same lines of Theorem 5.2 of \cite{BPT}.
 It suffices to show that the map $a\mapsto d_\tau(a)$ is surjective from $W(A)$ onto ${\rm LAff}_b(T(A)).$
Let $f\in {\rm LAff}_b(A)_+$ with $f(\tau)>0$ for all  $\tau\in T(A).$
We may assume that $f(\tau)\le 1$ for all $\tau\in T(A).$
 As in the proof of 5.2 of \cite{BPT}, it suffices to find a sequence of $a_i\in M_2(A)_+$ such that
 $a_i\lesssim a_{i+1},$
 $[a_n]\not=[a_{n+1}]$ (in $W(A)$) and
 $$
 {\lim_{n\to\infty}d_{\tau}(a_n)=f(\tau)\rforal \tau\in T(A).}
 $$
Using the semi-continuity of $f,$  we find a sequence $f_n\in
\Aff(T(A))_{++}$ such that \beq\label{Aff=W-L-1}
f_n(\tau)< f_{n+1}(\tau)\rforal \tau\in T(A),\,\,\,n=1,2,....\\
\lim_{n\to\infty} f_n(\tau)=f(\tau)\tforal \tau\in T(A).
 \eneq
Since $f_{n+1}-f_n$ is continuous and strictly positive on the
compact set $T,$ there is $\ep_n>0$ such that
$(f_n-f_{n+1})(\tau)>\ep_n$ for all $\tau\in T(A),$ $n=1,2,....$
It follows  from \ref{CTdense}, for each $n,$ there is a \SCA\,
$C_n$ of $A$ with $C_n\in {\cal C}$ and an element $b_n\in (C_n)_+$
such that \beq\label{Aff=W-L-2}
{\rm dim}\pi(C_n)\ge  (16/\ep_n)^2\,\,\,\text{for \,each \,irreducible\,representation} \,\,\pi\,\,\,of \,\,C_n,\\
0<\tau(f_n)-\tau(b_n)<\ep_n/4\rforal \tau\in T(A).
\eneq
Applying \ref{Affon2}, one obtains an element $a_n\in M_2(C_n)_+$
such that \beq\label{Aff=W-L-3} 0<t(b_n)-d_t(a_n)<\ep_n/4\tforal
t\in T(C_n). \eneq It follows that \beq\label{Aff=W-L-4}
0<\tau(f_n)-d_\tau(a_n)<\ep_n/2\tforal \tau\in T(A). \eneq One then
checks that $\lim_{n\to\infty}d_\tau(a_n)=f(\tau)$ for all $\tau\in
T(A).$ Moreover, $d_\tau(a_n)<d_\tau(a_{n+1})$ for all $\tau\in
T(A),$ $n=1,2,....$ It follows  from \ref{Comparison} that
$a_n\lesssim a_{n+1},$ $[a_n]\not=[a_{n+1}],$ $n=1,2,....$  This
ends the proof.
\end{proof}

\begin{thm}\label{Aff=W}
Let $A\in {\cal B}_1$ be a unital simple \CA. Then $W(A)$ has
$0$-almost divisible property.
%Then $W(A)=V(A)\sqcup LAff_b(T(C)).$
\end{thm}

\begin{proof}
Let $a\in M_n(A)_+\setminus \{0\}$ and $k\ge 1$ be an integer. We
need to show that there exists an element $x\in M_{m'}(A)_+$ for
some $m'\ge 1$ such that \beq\label{Aff=W-1} k [x] \le [a] \le (k+1)
[ x ] \eneq in $W(A).$ {
%There is $\sigma>0$ such that
%$d_\tau(a)>\sigma$ for all $\tau\in T(A).$
It follows from \ref{Aff=W-L} that, since $kd_\tau(a)/(k^2+1)\in {\rm
LAff}_b(T(A)),$ there is $x\in M_{2n}(A)_+$ such that
\beq\label{Aff=W-2} d_\tau(x)=kd_\tau(a)/(k^2+1)\rforal \tau\in
T(A). \eneq Then, \beq\label{Aff=W-3}
kd_\tau(x)<d_\tau(a)<(k+1)d_\tau(x)\tforal \tau\in T(A). \eneq It
follows from \ref{Comparison} that \beq\label{Aff=W-4} k[x]\le [a]\le
(k+1) [x]. \eneq }
\end{proof}

\begin{thm}\label{Zstable}
Let $A\in {\cal B}_1$ be a unital separable simple amenable \CA. Then
$A\otimes {\cal Z}\cong A.$
\end{thm}

\begin{proof}
Since $A\in {\cal B}_1,$ $A$ has finite weak tracial nuclear dimension (see 8.1 of \cite{Lin-LAH}).
By \ref{Comparison}, $A$ has the strict comparison property for positive elements.
 Note, by \ref{B1hered}, every unital hereditary \SCA\, of $A$ is in ${\cal B}_1.$
 Thus, by \ref{Aff=W},  its Cuntz semigroup also has $0$-almost divisibility. It follows from 8.3 of \cite{Lin-LAH} that $A$ is ${\cal Z}$-stable.
\end{proof}

\section{The unitary groups}

\begin{thm} {\rm (cf. Theorem 6.5 of \cite{LinTAI})}\label{UL1}
Let $K\in \N$ be an integer and let ${\cal B}$ be a class of unital \CA s which has the property
that ${\rm cer}(B)\le K$
%{\color{Green} ($K+\epsilon?)$}
for all $B\in {\cal  B}.$
Let $A$ be a unital simple \CA\, which is tracially in ${\cal B}$ and let $u\in U_0(A).$
Then, for any $\ep>0,$  there exists a unitary
$u_1, u_2\in A$ such that $u_1$ has exponential length no more than $2\pi,$ $u_2$ has exponential rank $K$ and
$$
\|u-u_1u_2\|<\ep.
$$
Moreover, ${\rm cer}(A)\le K+2+\ep.$
\end{thm}

\begin{proof}
The proof is exactly the same as that of Theorem 6.5 of \cite{LinTAI}.
\end{proof}

\begin{cor}\label{cerB1}
Any C*-algebra in the class ${\cal B}_1$ has exponential rank at most $5+\epsilon$.
\end{cor}
\begin{proof}
By Theorem \ref{2Tg14},  \CA s in $\mathcal C_0$ has exponential rank at most $3+\epsilon$. Therefore, by Theorem \ref{UL1}, any \CA\,  in $\mathcal B_1$ has exponential rank at most $5+\epsilon$.
\end{proof}

\begin{thm}\label{Ulength}
Let $L>0$ be a positive number and let ${\cal B}$ be a class of unital \CA s such that ${\rm cel}(v)\le L$ for
every  unitary $v$ {in their closure of commutator subgroups}.
Let $A$ be a unital simple \CA\,  which is tracially in ${\cal B}$  and let $u\in CU(A).$
Then $u\in U_0(A)$ and ${\rm cel}(u)\le 3\pi+L.$
\end{thm}

 \begin{proof}
 Let $1>\ep>0.$
 There are $v_1, v_2,....,v_k\in U(A)$ such that
 \beq\label{Ulength-1}
 \|u-v_1v_2\cdots v_k\|<\ep/16
 \eneq
 and $v_i=a_ib_ia_i^*b_i^*,$ where $a_i, b_i\in U(A).$
 Let $N$ be an integer in Lemma 6.4 of \cite{LinTAI} (for $L=8\pi+\ep$).
 Since $A$ is tracially in ${\cal B},$ there is a projection $p\in A$ a unital \SCA\, in ${\cal B}$ with $1_B=p$
 such that
 \beq\label{Ulnegth-2}
 \|a_i-(a_i'\oplus a_i'')\|<\ep/{32}k,\,\,\,\|b_i-(b_i'\oplus b_i'')\|<\ep/{32}k,\,\,\,i=1,2,...,k\\
 \|u-\prod_{i=1}^k (a_i'b_i'(a_i')^*(b_i')^*\oplus a_i''b_i''(a_i'')^*(b_i'')^*\|<\ep/8,
 \eneq
 where $a_i', b_i'\in U((1-p)A(1-p)),$ $a_i'', b_i''\in U_0(B)$ and
 $6N[1-p]\le [p].$
 Put
 \beq\label{Ulength-3}
 w=\prod_{i=1}^k a_i'{{b_i'}}(a_i')^*(b_i')^*
 \andeqn z=\prod_{i=1}^k a_i''b_i''(a_i''){{(b_i'')^*.}}
 \eneq
 Then $z\in  CU(B).$ Therefore ${\rm cel}_{B}(z)\le L$  in $B\subset pAp.$
 It is standard to show that
 $$
 a_i'b_i'(a_i')^*(b_i)^*\oplus (1-p)\oplus (1-p)
 $$
 is in $U_0(M_4((1-p)A(-1p)))$ and it has exponential length no more than $4(2 \pi)+2\ep/16k.$
 This implies
 $$
 {\rm cel}(w\oplus (1-p)\oplus (1-p))\le 8\pi k+\ep/4
 $$
 in $U(M_{{4}}((1-p)A(1-p))).$ It follows from Lemma 6.4 of \cite{LinTAI} that
 \beq\label{Ulength-4}
 {\rm cel}(w\oplus p)\le 2\pi+\ep/4.
 \eneq
 It follows that
 $$
 {\rm cel}((w\oplus p)((1-p)\oplus z))<2\pi+\ep/4+L+\ep/16.
 $$
 {{Therefore}} ${\rm cel}(u)\le 2\pi+L+\ep.$
  \end{proof}

%{\color{green}
\begin{cor}\label{celB1}
Let $A\in\mathcal B_1$, and let $u\in CU(A)$. Then $u\in U_0$ and $\mathrm{cel}(u)\leq 7\pi$.
\end{cor}
\begin{proof}
{{This}}  follows from Lemma \ref{2Lg9} and Theorem \ref{Ulength}.
\end{proof}
%}

\begin{lem}\label{UCUdiv}
Let $A$ be a unital \CA,\,let $U$ be an infinite dimensional UHF-algebra  and let $B=A\otimes U.$
%Suppose that there is an integer $K>0$ such that
%${\rm cel }(u)\le K$ for all $u\in CU(A).$
Then $U_0(B)/CU(B)$ is torsion free  and divisible.

\end{lem}

\begin{proof}
%{\color{green}
Since $B=A\otimes U$, and $U$ is a UHF-algebra with infinite type,
%one has that
{{for}} any projection $p$ in a matrix algebra over $B$, there are  projections $e_1, e_2, ..., e_m\in B$ for some $m$ such that
$$p=e_1\oplus e_2 \oplus\cdots\oplus e_m.$$
Therefore $\rho_B(K_0(B))$ is spanned by the image of the projections in $B$, and by Theorem 3.2 of \cite{Thomsen-rims}, the group $U_0(B)/CU(B)$ is isomorphic to $\Aff(T(B))/\overline{\rho_B(K_0(B))}$. Since $\Aff(T(B))$ is divisible (it is a vector space), so is its quotient.
%}

We will show that $\Aff(T(B))/\overline{\rho_B(K_0(B))}$ is torsion
free. Suppose that $a\in \Aff(T(B))$ so that ${n}a\in
\overline{\rho_B(K_0(B))}$ for some integer ${n}>1.$
Let $\ep>0.$ There exists a pair of projections  $p, q\in M_n(B)$ such that
$$
|na(\tau)-\tau(p)-\tau(q)|<\ep/2
$$
for all $\tau\in T({B}),$  regarded as unnormalized trace
on $M_n(B)$. Since $B=A\otimes U,$ there are mutually orthogonal
projections $p_1, p_2,...,p_n, p_{n+1}, q_1, q_2, ...,q_n, q_{n+1}$ such that
$$ p_1+p_2+\cdots+p_n+p_{n+1}=p\andeqn q_1+q_2+\cdots +q_n+q_{n+1}=q$$
\beq
[p_1]=[p_j],\,\,\,[q_1]=[q_j],\,\,\, j=1,2,...,n\andeqn \tau(p_{n+1}),\, \tau(q_{n+1})<\ep/2\tforal \tau\in T(B).
\eneq
It follows that
$$
|a(\tau)-(\tau(p_1)-\tau(q_1))|<\ep\tforal \tau\in T(B).
$$
This implies that $a\in \overline{\rho_B(K_0(B))}.$ Therefore $\Aff(T(B))/\overline{\rho_B(K_0(B))}$ is torsion
free and the lemma follows.
\end{proof}

\begin{thm}\label{Utorlength}
Let $A$ be a unital \CA\, such that
there is an integer $K>0$ such that
${\rm cel}(u)\le K$ for all $u\in CU(A).$
Suppose that  $U_0(A)/ CU(A)$ is torsion free and suppose that
 $u,\, v\in U(A)$ such that $u^*v\in U_0(A).$
Suppose also that there is $k\in \N$ such that ${\rm cel}((u^k)^*v^k)<L$ for some $L>0.$
Then
\beq\label{Utorl-1}
{\rm cel}(u^*v)\le K+L/k.
\eneq
%Moreover, there is $y\in U_0(A)$ with ${\rm cel}(y)\le L/k$ such that $\overline{u^*y}={\bar y}$ in $U(A)/CU(A).$
\end{thm}

\begin{proof}
It follows from \cite{Ringrose-cel} that, for any $\ep>0,$   there are $a_1, a_2,...,a_N\in A_{s.a.}$ such that
\beq\label{Utorl-2}
(u^k)^*v^k=\prod_{j=1}^N \exp(\sqrt{-1}a_j)\andeqn \sum_{j=1}^N \|a_j\|\le L+\ep/2.
\eneq
Choose
$
w=\prod_{j=1}^N\exp(-\sqrt{-1}a_j/k).
$
Then
$
(u^*vw)^k\in CU(A).
$
 Since $U_0(A)/CU(A)$ is assumed to be torsion free, it follows  that
\beq\label{Utorl-3}
u^*vw\in CU(A).
\eneq
Thus,
$
{\rm cel}(u^*vw)\le K.
$
Note that
$
{\rm cel}(w)\le L/k +\ep/2k,
$
It follows that
$$
{\rm cel}(u^*v)\le K+L/k+\ep/2k.
$$
\end{proof}

\begin{cor}\label{Unotrosion}
Let $A$ be a unital simple \CA\, in ${\cal B}_1,$ $B=A\otimes U,$ where
$U$ is a UHF-algebra
%with
{{of}}  infinite type.
% $B.$
Then
\begin{enumerate}
\item $U_0(B)/CU(B)$ is torsion free and divisible; and
\item if $u, v\in U(B)$ with ${\rm cel}((u^*)^kv^k)\le L$ for some integer $k>0,$  then
$$
{\rm cel}(u^*v)\le {7}\pi+L/k.
$$
\end{enumerate}
\end{cor}
%{\color{green}
\begin{proof}
The lemma follows from Lemma \ref{UCUdiv}, Corollary \ref{celB1}, and
Theorem \ref{Utorlength}.
\end{proof}
%}

\begin{cor}\label{Ukerdiv}
Let $A_n$ be a sequence of unital separable simple \CA s in ${\cal B}_1$ and let $B_n=A_n\otimes U.$
Then the kernel of the map
\beq\label{Ukerdiv-1}
K_1(\prod_n B_n)\to \prod^b_n K_1(B_n)\to 0
\eneq
is a divisible and torsion free.
\end{cor}

\begin{proof}
By  Corollary \ref{cerB1}, the exponential rank of each $B_n$ is
bounded by  $6.$  Since each $B_n$ has stable rank one, by (2) of
Proposition 2.1 of \cite{GL-almost-map}, the kernel of the map in (\ref{Ukerdiv-1}) is divisible.
Cor \ref{Ukerdiv} implies that $\{A_n\}$ has exponential length divisible rank
$7\pi+L/k$ (see Definition 2.1 of \cite{GL-almost-map}).
It follows Lemma 2.2 of \cite{GL-almost-map} that  the kernel is also torsion free.

%Suppose that $\{u_n\}\in U(M_{{ K}}(\prod_n B_n))$ such that
%$[\{u_n\}]$ is in the kernel  and $k[\{u_n\}]=0$ for some integer
%$k>0.$ By changing notation, without loss of generality, we may
%assume that $\{u^k_n\}\in U_0(M_K(\prod_n B_n))$   for some integer $K\ge 1.$  Also $u_n\in
%U_0(B_n)$  for each $n$. However, the fact that $\{u_n^k\}\in
%U_0(M_K(\prod_n B_n))$ implies that there is $L>0$ such that
%${\rm cel }(u_n^k)\le L$ for all $n.$  It follows from
%\ref{Unotrosion} that
%$$
%{\rm cel}(u_n)\le {7}\pi+L/k+\pi/4\tforal n.
%$$
%Consequently, $\{u_n\}\in U_0(M_K(\prod_n B_n)).$ Therefore $[\{u_n\}]=0$ in $K_1(\prod_n B_n).$ So the kernel
%is torsion free.
\end{proof}

\begin{lem}\label{ph}
Let $K\ge 1$ be an integer.  Let $A$ be a  unital  simple \CA\, in
${\cal B}_1.$  Let  $e\in A$ be a projection and let $u\in
U_0(eAe).$ Suppose that $w=u+(1-e)$ and suppose   $\eta\in (0,2].$ Suppose
also that \beq\label{PH-1} [1-e]\le  K[e]\,\,\,{\rm in}\,\,\, K_0(A)
\andeqn {\rm dist}({\bar w},{\bar 1})\leq \eta. \eneq

Then, if $\eta<2,$ one has
$$
{\rm cel}_{eAe}(u)<({K\pi\over{2}}+1/16)\eta+6\pi \andeqn {\rm
dist}({\bar u}, {\bar e})<(K+1/8)\eta,
$$
and if $\eta=2,$ one has
$$
{\rm cel}_{eAe}(u)<{K\pi\over{2}}{\rm cel}(w)+1/16+6\pi.
$$
\end{lem}

\begin{proof}
We  assume that (\ref{PH-1}) holds. Note that $\eta\le 2.$ Put
$L={\rm cel}(w).$
We first consider the case that $\eta<2.$ There is a projection
$e'\in M_2(A)$ such that
$$
[(1-e)+e']=K[e].
$$
To simplify notation,  by \ref{B1hered} and by replacing $A$ by
$(1_A+e')M_2(A)(1_A+e')$ and $w$ by $w+e',$ without loss of
generality, we may now assume that
\beq\label{PH-10}
[1-e]=K[e]\andeqn {\rm dist}({\bar w},{\bar 1})<\eta.
\eneq
There is $R_1>1$ such that
$\max\{L/R_1,2/R_1,\eta\pi/R_1\}<\min\{\eta/64, 1/16\pi\}.$

For
any ${\eta\over{32K(K+1)\pi}}>\ep>0$ with $\ep+\eta<2,$  since
$TR(A)\le 1,$ there exists a projection $p\in A$ and a \SCA\, $D\in
{\cal C}$ with $1_D=p$ such that
\begin{enumerate}

\item $\|[p,\,x]\|<\ep$ for $x\in \{u, w, e, (1-e)\},$

\item $pwp, pup, pep, p(1-e)p\in_{\ep} D,$

\item there is a  projection $q\in D$ and a unitary $z_1\in qDq$
such that $\|q-pep\|<\ep,$ $\|z_1-quq\|<\ep,$ $\|z_1\oplus (p-
q)-pwp\|<\ep$ and $\|z_1\oplus (p-q)-c_1\|<\ep+{\eta};$

\item  there is a  projection $q_0\in (1-p)A(1-p)$ and a unitary
$z_0\in q_0Aq_0$ such that\\ $\|q_0-(1-p)e(1-p)\|<\ep,$
$\|z_0-(1-p)u(1-p)\|<\ep,$ $\|z_0\oplus
(1-p-q_0)-(1-p)w(1-p)\|<\ep,$ $\|z_0\oplus
(1-p-q_0)-c_0\|<\ep+{\eta},$

\item $[p-q]=K [q]$ in $K_0(D),$ $[(1-p)-q_0]=K[q_0]$ in $K_0(A);$

\item $2(K+1)R_1[1-p]<  [p]$ in $K_0(A);$

\item $ {\rm cel}_{(1-p)A(1-p)}(z_0\oplus (1-p-q_0))\le L+\ep,$

\end{enumerate}
where
 $c_1\in CU(D)$ and
$c_0\in CU((1-p)A(1-p)).$

By  Lemma \ref{2Lg8}, one has that ${\rm det}_D(c_1)=1.$ Since $\ep+\eta<2,$
there is $h\in D_{s.a}$ with $\|h\|\le
2\arcsin({\ep+\eta\over{2}})$ such that (by (3) above)
\beq\label{ph-5}
(z_1\oplus (p-q))\exp(ih)=c_1.
\eneq
It follows that
\beq\label{ph-6}
{\rm det}_D((z_1\oplus(p-q))\exp(ih))=1,
\eneq
or
\beq\label{ph-6+1}
D_D(z_1\oplus (p-q)\exp(ih))(t)=0\tforal t\in T(D).
\eneq
It follows that
\beq\label{ph-6-2}
|D_D(z_1\oplus (p-q))(t)|\le 2\arcsin({\ep+\eta\over{2}})\tforal \tau\in T(D).
\eneq
By (5) above, one obtains that
\beq\label{ph-7}
|D_{qDq}(z_1)(t)|\le K2\arcsin({\ep+\eta\over{2}})\tforal t\in T(qDq).
\eneq
If $2K\arcsin({\ep+\eta\over{2}})\ge \pi,$ then
$
2K({\ep+\eta\over{2}}){\pi\over{2}}\ge \pi.
$
It follows that
\beq\label{ph-8-}
K(\ep+\eta)\ge 2\ge {\rm dist}(\overline{z_1}, {\overline{q}}).
\eneq
 Since those unitaries in $D$ with ${\rm det}(u)=1$ (for all
points) are in $CU(D)$ (see  Lemma \ref{2Lg8}), from
(\ref{ph-7}), one computes that, when
$2K\arcsin({\ep+\eta\over{2}})< \pi,$
\beq\label{ph-8}
{\rm dist}(\overline{z_1},
{\overline{q}})<2\sin(K\arcsin({\ep+\eta\over{2}}))\le
K(\ep+\eta).
\eneq
By combining both (\ref{ph-8-}) and (\ref{ph-8}), one obtains that
\beq\label{ph-8+}
{\rm dist}(\overline{z_1}, {\overline{q}})\le K(\ep+\eta) \le
K\eta+{\eta\over{32(K+1)\pi}}.
\eneq
By (\ref{ph-7}), it follows from  Lemma \ref{2Lg9} that
\beq\label{ph-8+1}
\hspace{-0.3in}{\rm cel}_{qDq}(z_1)\le 2K\arcsin{\ep+\eta\over{2}} +4\pi \le
K(\ep+\eta){\pi\over{2}}+4\pi\le
(K{\pi\over{2}}+{1\over{64(K+1)}})\eta +4\pi .
\eneq
By (5) and (6) above,
$$
(K+1)[q]=[p-q]+[q]=[p]>2(K+1)R_1[1-p].
$$
Since $K_0(A)$ is weakly unperforated, one has
\beq\label{ph-9}
2R_1[1-p]<[q].
\eneq

There is a unitary $v\in A$ such that
%\beq\label{ph-10}
$v^*(1-p-q_0)v\le q.$
%\eneq
Put $v_1=q_0\oplus (1-p-q_0)v.$ Then
\beq\label{ph-11}
v_1^*(z_0\oplus (1-p-q_0))v_1=z_0\oplus v^*(1-p-q_0)v.
\eneq
Note that
\beq\label{ph-10+n}
\|(z_0\oplus v^*(1-p-q_0)v)v_1^*c_0^*v_1-q_0\oplus
v^*(1-p-q_0)v\|<\ep+\eta.
\eneq
Moreover, by (7) above,
\beq\label{ph-10+}
{\rm cel}(z_0\oplus v^*(1-p-q_0)v)\le L+\ep,
\eneq
It follows from (\ref{ph-9}) and Lemma 6.4 of \cite{LinTAI} that
\beq\label{ph-10+1}
{\rm cel}_{(q_0+q)A(q_0+q)}(z_0\oplus q)\le 2\pi+(L+\ep)/R_1.
\eneq
Therefore, combining (\ref{ph-8+1}),
\beq\label{ph-10+2}
{\rm cel}_{(q_0+q)A(q_0+q)}(z_0+z_1)\le
2\pi+(L+\ep)/R_1+(K{\pi\over{2}}+{1\over{64(K+1)}})\eta +6\pi.
\eneq
By (\ref{ph-10+}), (\ref{ph-9}) and Lemma 3.1 of \cite{Lin-hmtp},
in $U_0((q_0+q)A(q_0+q))/CU((q_0+q)A(q_0+q)),$
\beq\label{ph-13}
{\rm dist}(\overline{z_0+q},
{\overline{q_0+q}})<{(L+\ep)\over{R_1}}.
\eneq
Therefore, by (\ref{ph-8}) and (\ref{ph-13}),
\beq\label{ph-14}
{\rm dist}({\overline{z_0\oplus z_1}},
{\overline{q_0+q}})<{(L+\ep)\over{R_1}}+K\eta+{\eta\over{32(K+1)\pi}}<(K+1/16)\eta.
\eneq
We note that
\beq
\|e-(q_0+q)\|<2\ep\andeqn \|u-(z_0+z_1)\|<2\ep.
\eneq
It follows that
\beq\label{ph-15}
{\rm dist}({\bar u}, {\bar e})<4\ep +(K+1/16)\eta<(K+1/8)\eta.
\eneq

Similarly, by (\ref{ph-10+2}),
\beq\label{ph-15+1}
{\rm cel}_{eAe}(u)&\le& 4\ep\pi+2\pi+(L+\ep)/R_1+(K{\pi\over{2}}+{1\over{64(K+1)}})\eta
+4\pi\\
&<&(K{\pi\over{2}}+1/16)\eta+6\pi.
\eneq

This proves the case that $\eta<2.$

Now suppose that $\eta=2.$
%Let $\Lambda={{\rm cel}(w)\over{K}}.$
%Let $L={\rm cel}(w).$
Define $R=[{\rm cel}(w)+1].$ Note that ${{\rm cel}(w)\over{R}}<1.$
%Choose an integer $R\ge 4K$ such that $L/R<1/2.$
There is a projection $e'\in M_{R+1}(A)$ such that
$$
[(1-e)+e']=(K+RK)[e].
$$
It follows from Lemma 3.1 of \cite{Lin-hmtp} that
\beq\label{PH-111}
{\rm dist}(\overline{w\oplus e'},{\overline{ 1_A+e'}})<{{\rm
cel}(w)\over{R+1}}.
\eneq
Put $K_1=K(R+1).$ To simplify notation,  without loss of
generality, we may now assume that
\beq\label{PH-102}
[1-e]=K_1[e]\andeqn {\rm dist}({\bar w},{\bar 1})<{{\rm
cel}(w)\over{R+1}}.
%\andeqn {\rm
%cel}(w)\le 2\arcsin(L/2R)+2\pi.
\eneq

It follows from the first part of the lemma that
\beq\label{PH-103}
\hspace{-0.3in}{\rm cel}_{eAe}(u)<
({K_1\pi\over{2}}+{1\over{16}}){{\rm cel}(w)\over{R+1}}+6\pi
\le {K\pi{\rm cel}(w)\over{2}}+{1\over{16}}+6\pi.
\eneq

\end{proof}

%\begin{cor}\label{UCUinj}
%Let $A\in {\cal B}_1$ be a unital  simple \CA\, and let
%$e\in A$ be a non-zero projection. Then the map $u\mapsto u+(1-e)$
%induces an injective contractible \hm\,  from $U(eAe)/CU(eAe)$ onto $U(A)/CU(A).$
%\end{cor}
%
%\begin{proof}
%This was originally proved here following \ref{ph}. However, by \ref{???},
%$A$ has stable rank one. Thus this follows from Theorem 4.6 of \cite{GLX-ER}.
%\end{proof}

%\begin{proof}
%Let $j$ be the map.
%Fix a unitary $u\in eAe$ so that ${\bar u}\in {\rm ker}\,
%j.$ We will show that $u\in CU(eAe).$%
%
%There is an integer $K\ge 1$ such that
%$$
%K[e]\ge [1-e]\,\,\,{\rm in}\,\,\, K_0(A).
%$$
%Put $v=u+(1-e).$ Let $1>\ep>0.$ Since ${\bar u}\in {\rm ker}j,$
%$v\in CU(A).$ In particular,
%$$
%{\rm dist}({\bar v}, \bar{1})<\ep/(K\pi/2+1).
%$$
%It follows from Lemma \ref{ph} that
%$$
%{\rm dist}({\bar u}, {\bar
%e})<({K\pi\over{2}}+1/16)(\ep/(K\pi/2+1))<\ep.
%$$
%It then follows that
%$$
%u\in CU(eAe).
%$$
%\end{proof}

\begin{thm}\label{UCUiso}
%Let $A\in {\cal B}_1$ be a unital  simple \CA\,
{Let $A$ be a unital  simple \CA\, of stable rank one} and let
$e\in A$ be a non-zero projection. Then the map $u\mapsto u+(1-e)$
induces an isomorphism $j$ from $U(eAe)/CU(eAe)$ onto $U(A)/CU(A).$
\end{thm}

\begin{proof}
This was originally proved here following \ref{ph} (with additional assumption that $A\in {\cal B}_1${{).}}  However, by \ref{2pg3},
$A$ has stable rank one. Thus this follows from Theorem 4.6 of \cite{GLX-ER}.
\end{proof}

\begin{cor}\label{c1}
Let {$A$ be a unital  simple \CA\, of stable rank one.} Then
the map $j: a\to {\rm diag}(a, \overbrace{1,1,...,1}^m)$ from $A$ to
$M_n(A)$ induces an isomorphism from $U(A)/CU(A)$ onto
$U(M_n(A))/CU(M_n(A))$ for any integer $n\ge 1.$
\end{cor}

\begin{proof}
This follows from \ref{UCUiso} but also follows from 3.11 of \cite{GLX-ER}.
\end{proof}

\section{A Uniqueness Theorem for \CA s in ${\cal B}_0$ }

%The following lemma follows directly from the definition of ${\cal B}_1$ and the approximate divisibility of UHF-algebras.

%\begin{prop}\label{decomp}
%Let $A_1$ be a unital separable simple \CA\, in ${\cal B}_1,$
%let $U$ be a UHF-algebra of infinite type and $A=A_1\otimes U.$
%Then, for any $\ep>0,$ any  finite subset ${\cal F}\subset A,$  and
%any integer $n\ge 1,$
%there exist mutually orthogonal  projections $p_0,p_1,...,p_n\in  A,$  unitaries $u_j\in A,$ $j=1,2,...,n-1,$
%a \SCA\, $C\in {\cal C}$ with $1_C=p_1,$  unital ${\cal F}$-$\ep$-multiplicative \morp s $\phi_0: A\to p_0Ap_0$
%and $\phi_1:  A\to C$ such that
%\beq\label{decomp-1}
%&&p_0\lesssim p_1, u_j^*p_1u_j=p_{j+1},,\,\,\, j=1,2,...,n-1;\\
%&&\|x-\phi_0(x)\oplus  \phi_1(x)\oplus ({\rm Ad}\, u_1\circ \phi_1(x))\oplus \cdots \oplus ({\rm Ad}\, u_{n-1}\circ \phi_1(x))\|<\ep
%\eneq
%for all $ x\in {\cal F}.$
%\end{prop}

%\begin{df}
%Let $A$ and $B$ be C*-algebras, and assume that $B$ is unital. Let $\mathcal H\subseteq A_+$ be a finite subset, let $T: A_+\setminus\{0\}\to \R_+\setminus\{0\}$ and let $N: A_+\setminus\{0\}\to \N$ be two maps. Then a map $L: A\to B$ is said to be $T\times N$-$\mathcal H$-full if for any $h\in \mathcal H$, if there are $b_1, b_2, ..., b_{N(h)}\in B$ such that $\|b_i\|\leq T(h)$ and
%$$\sum_{i=1}^{N(h)} b_i^*L(h)b_i=1_B.$$
%\end{df}

{{The following follows from Theorem 7.1 of \cite{Lin-hmtp}.}}

\begin{thm}\label{Suni} {{\rm  (cf.  Theorem Theorem 7.1 of \cite{Lin-hmtp})}}
Let  $A$ be a unital separable amenable \CA\, which satisfies the UCT, {{let}} $T\times N: A_+\setminus \{0\}
\to \R_+\setminus \{0\}\times \N$ be a map  and let ${\bf L}: U(M_{\infty}(A))\to
\R_+$ be another map. For any $\ep>0$ and any finite subset ${\cal F}\subset A,$ there exists  $\dt>0,$ a finite subset
${\cal G}\subset A,$ a finite subset ${\cal H}\subset A_+\setminus \{0\},$ a finite subset ${\cal U}\subset
{\cup_{m=1}U(M_m(A))},$ a finite subset ${\cal P}\subset \underline{K}(A)$ and an integer $n>0$ satisfying the following:
for any unital separable simple \CA\, $C$ in ${\cal B}_1,$ if $\phi, \psi, \sigma: A\to B
=C\otimes U,$ where
$U$ is a UHF-algebra of infinite type,
are three ${\cal G}$-$\dt$-multiplicative \morp s { and  $\sigma$ is}  unital and $T\times N$-${\cal H}$-full ({{see}} \ref{Dfull}) with
{the properties that}
\beq\label{Suni-1}
[\phi]|_{\cal P}=[\psi]|_{\cal P}\tand
{\rm cel}(\langle \psi (u)\rangle^*\langle \phi(v)\rangle)\le {\bf L}(v)
\eneq
for all $v\in {\cal U}$ and $\sigma$ is unital,  there exists a unitary $u\in M_{n+1}(B)$ such that
$$
\|u^*{\rm diag}(\phi(a), \overline{\sigma}(a))u-{\rm diag}(\psi(a), \overline{\sigma}(a))\|<\ep
$$
for all $a\in {\cal F},$ where
$$
\overline{\sigma}(a)={\rm diag}(\overbrace{\sigma(a),\sigma(a),...,\sigma(a)}^n)\tforal a\in A.
$$
\end{thm}

\begin{proof}
Note that $B$  has stable rank one,  ${\rm cer}(M_m(B))\le  6$ (see \ref{cerB1}), exponential length divisible rank
$E(L,k)=7\pi +L/k$ ((2) of \ref{Unotrosion}) and $K_0(B)$ are weakly unperforated (in particular,
it has $K_0$-divisible rank $T(L,k)=L+1$).
The class of unital \CA s have just mentioned properties will be denoted by ${\bf C}.$
We sketch the proof below.
Suppose that the theorem is false.  Then there exists  $\ep_0>0$ and a finite subset ${\cal F}\subset A$
such that there are a sequence of positive numbers $\{\dt_n\}$ with $\dt_n\searrow 0,$ an increasing sequence
$\{{\cal G}_n\} \subset A$ of finite subsets such that $\cup_n {\cal G}_n$ is dense in $A,$
an increasing sequence $\{{\cal P}_n\}\subset \underline{K}(A)$ of finite subsets  with
$\cup_{n=1}{\cal P}_n=\underline{K}(A),$  an increasing sequence of finite subsets
$\{{\cal U}_n\}\subset U(M_{\infty}(A))$ such that $\cup_{n=1}{\cal U}_n\cap U(M_m(A))$ is dense
in $U(M_m(A))$ for each integer $m\ge1,$  an increasing sequence of finite subsets
$\{{\cal H}_n\}\subset  A_+^{\bf 1}\setminus \{0\}$ such that if $a\in {\cal H}_n$
and $f_{1/2}(a)\not=0,$ then
$f_{1/2}(a)\in {\cal H}_{n+1}$ and $\cup_{n=1}{\cal H}_n$ is dense
in $A_+^{\bf 1},$  a sequences of integers $\{k(n)\}$ with
$\lim_{n\to\infty} k(n)=\infty$,  a sequence of unital \CA s $B_n\in {\bf C},$
%{\bf C}_{r_0, r_1, T, s, R},$
two sequences of ${\cal G}_n$-$\dt_n$-multiplicative \morp s $\phi_{n}, \psi_{n}: A\to B_n,$
%where $D_n$ is a UHF-algebra of infinite type,
such that
\beq\label{stableun2-1}
[\phi_{n}]|_{{\cal P}_n}=[\psi_{n}]|_{{\cal P}_n}\andeqn
{\rm cel}(\langle \phi_{n}(u)\rangle \langle \psi_{n}(u^*)\rangle)\le {\bf L}(u)
\eneq
for all $u\in {\cal U}_n$
and a sequence of unital ${\cal G}_n$-$\dt_n$-multiplicative \morp\, $\sigma_n: A\to
B_n$ which is also $F$-${\cal H}_n$-full satisfying
\beq\label{stableun2-2}
&&\hspace{0.3in} \inf\{\sup\|v_n^*{\diag}( \phi_{n}(a), S_n(a))v_n-{\diag}(\psi_n(a), S_n(a))\|: a\in {\cal F}\|\}\ge \ep_0,
\eneq
where the infimum is taken among all unitaries $v_n\in M_{k(n)l(n)+1}(B_n)$
and where 
$$
S_n(a)={\diag}(\overbrace{\sigma_n(a),\sigma_n(a),...,\sigma_n(a)}^{k(n)})\rforal a\in A.
$$

Let $C_0=\bigoplus_{n=1}^{\infty}B_n,$ $C=\prod_{n=1}^{\infty}B_n,$ $Q(C)=C/C_0$ and
$\pi: C\to Q(C)$ is the quotient map.  Define $\Phi, \Psi, S: A\to C$ by
$\Phi(a)=\{\phi_n(a)\},$ $\Psi(a)=\{\psi_n(a)\}$  and  $S(a)=\{\sigma_n(a)\}$
all $a\in A.$ Note that $\pi\circ \Phi,$ $\pi\circ \Psi$ and $\pi\circ S$ are \hm s.

As in the proof of 7.1 of \cite{Lin-hmtp}, since $B_n\in {\bf C},$ one computes that
\beq\label{stableun2-3}
[\pi\circ \Phi]=[\pi\circ \Psi]\,\,\,{\rm in}\,\,\, KL(A, Q(C)).
\eneq
One can also check that $\pi\circ S$ is a full \hm.
It follows from \ref{Lnuct} that there exists an integer $K\ge 1$ and a unitary $U\in M_{K+1}(Q(C))$
such that
\beq\label{stableun2-4}
\|U^*{\diag}(\pi\circ \Phi(a), \Sigma(a))U-{\diag}(\pi\circ \Psi(a))\|<\ep_0/4\rforal a\in {\cal F}.
\eneq
It follows that there exists a  $V=\{v_n\}\in C$  and an integer $N\ge 1$ such that, for any $n\ge N,$
$k(n)\ge K$ such that
\beq\label{stableun2-5}
\|v_n^*{\diag}(\phi_n(a), {\overline{\sigma}_n}(a))v_n-{\diag}(\psi_n(a), {\overline{\sigma}_n}(a))\|<\ep_0/2
\eneq
for all $a\in {\cal F},$ where
$$
{\overline{\sigma}_n}(a)={\diag}(\overbrace{\sigma_n(a), \sigma_n(a),...,\sigma_n(a)}^K)\rforal a\in A.
$$
This contradicts with (\ref{stableun2-2}).
%Therefore Theorem 7.1 of \cite{Lin-hmtp} applies.
%Note that we may choose  $n=lk$ in 7.1 of \cite{Lin-hmtp}.
\end{proof}

\begin{rem}\label{Rsuni}

(1) The assumption that $B=C\otimes U$ where $C\in {\cal B}_1$ and $U$ is an infinite  UHF-algebra
is just for our purpose in this paper. The same theorem holds for a much more general setting.
%What are needed in the proof is that (\ref{Suni-8}) and  (\ref{Suni-9}) hold.
Let ${\bf R}, {\bf r}: \N\to \N,$ ${\bf T}: \N^2\to \N,$ ${\bf E}: R_+\times \N\to \R_+$
be fixed maps.  Let ${\bf C}_{{\bf R, r, T, E}}$ be  the class of
unital \CA s which have  $K_i$-${\bf r}$-cancellation, $K_i$-divisible rank ${\bf T}$  and
exponential length divisible length ${\bf E},$ and ${\rm cer}(M_m(B))\le {\bf R}$
for all $B\in {\bf C}_{{\bf R, r, T, E}}.$ Then Theorem \ref{Suni} holds
for $B\in {\bf C}_{{\bf R, r, T, E}}$
% Theorem 2.1 of \cite{GL-almost-map} implies
%that  (\ref{Suni-8}) and  (\ref{Suni-9}) hold.
just as  7.1 of \cite{Lin-hmtp}.
One should also note that the proof  is  not new. It is just the same as the combination of
those of Theorem 5.9 of \cite{LnAUCT} and Theorem 4.8 of \cite{Lnjotuni}. One should also note
that $\phi$ and $\psi$ are not assumed to be unital. Thus Theorem \ref{Lauct2} can also  viewed
as a special case of Theorem \ref{Suni}.
%integer $l$ in the statement of 7.1 of \cite{Lin-hmtp} can be replaced by an arbitrary integer,
%in particular, Theorem \ref{Lauct2} also follows from 7.1 of \cite{Lin-hmtp} given \ref{Fullmeasure} below.

(2) \, \,Suppose that there exists an integer $n_0\ge 1$ such that
$U(M_{n_0}(A))/{U_0(M_{n_0}(A))}\to U(M_{n_0+k}(A))/U_0(M_{n_0+k}(A))$ is an isomorphism
for all $k\ge 1.$ Then ${\bf L}$ may be replaced by a map from $U(M_{n_0}(A))$ to $\R_+$,  and $\mathcal U$ can be chosen in {$U(M_{n_0}(A))$}.

{
 Moreover,
 the condition that
${\rm cel}(\la \phi(u)\ra \la \psi(u)^*)\le {\bf L}(u)$ may,  in practice,  be replaced
by a stronger condition that,  for all $u\in \overline{{\cal U}},$
\beq\label{Rsuni-1}
{\rm dist}(\phi^{\ddag}(u), \psi^{\ddag}(u))<L
\eneq
where $\overline{{\cal U}}\subset \cup_{m=1}^{\infty}U(M_m(A))/CU(M_m(A))$ is a finite subset
and where $L<2$ is a given constant. }

To see this, let ${\cal U}$ be a finite subset of $U(M_m(A))$  for
some large $m$ whose image in the group $U(M_m(A))/CU(M_m(A))$ is
$\overline{{\cal U}}.$ Then (\ref{Rsuni-1}) implies that
\beq\label{Rsuni-2} \|\la \phi(u)\ra\la \psi(u^*)\ra -v\|<{{2}}
\eneq
for some $v\in CU(M_m(A)),$ provided that $\dt$ is
sufficiently small and ${\cal G}$ is sufficiently large. Since ${\rm
cel}(v)\le {7}\pi,$ we conclude that
\beq\label{Rsuni-3}
{\rm cel}(\la \phi(u)\ra \la \psi(u^*)\ra)\le \pi+{7}\pi +1
\eneq for
all $u\in {\cal U},${and take ${\bf L}: U_\infty(A) \to \R_+$ to
be constant ${{8\pi}} +1$.}
Furthermore,  we may assume $$\overline{{\cal U}}\subset U(M_{n_0}(A))/CU(M_{n_0}(C)),$$ if
$K_1(C)=U(M_{n_0}(C))/U_0(M_{n_0}(C)).$

\end{rem}

\begin{lem}\label{Ldet}
Let $A$ be a unital separable simple \CA\, with $T(A)\not=\emptyset.$  There exists a map
$\Delta_0: A_+^{q, {\bf 1}}\setminus \{0\}\to (0,1)$ satisfying the following:
For any finite subset ${\cal H}\subset A_+^{\bf 1}\setminus \{0\},$ there exits $\dt>0$ and a finite subset ${\cal G}\subset
A$ such that, for any unital \CA\, $B$ with $T(B)\not=\emptyset$ and any unital ${\cal G}$-$\dt$-multiplicative
\morp\, $\phi: A\to B,$ one has
\beq\label{Ldet-1}
\tau\circ \phi(h)\ge \Delta_0(\hat{h})/2\tforal h\in {\cal H}
\eneq
for all $\tau\in T(B).$ {Moreover, one may assume that $\Delta_0(\widehat{1_A})=3/4$.}
\end{lem}

\begin{proof}

Define, for each $h\in A_+^{\bf 1}\setminus \{0\},$
\beq\label{Ldet-2}
\Delta_0(\hat{h})=\min\{3/4, \inf\{\tau(h): \tau\in T(A)\}\}.
\eneq

Let ${\cal H}\subset A_+^{\bf 1}\setminus\{0\}$ be a finite subset.
Define
\beq\label{Ldet-3}
d=\min\{\Delta_0(\hat{h})/4: h\in {\cal H}\}>0{.}
\eneq

Let $\dt>0$ and let ${\cal G}\subset A$ be a finite subset {{as}}
required by \ref{measureexistence} for $\ep=d$ and ${\cal F}={\cal H}.$

Suppose that $\phi: A\to B$ is a unital ${\cal G}$-$\dt$-multiplicative \morp. Then, 
by \ref{measureexistence}, for each $t\in T(B),$ there exists $\tau\in T(A)$ such that
\beq\label{Ldet-4}
|t\circ \phi(h)-\tau(h)|<d\tforal h\in {\cal H}.
\eneq
It follows that
%\beq\label{Ldet-5}
$
t\circ \phi(h)>\tau(h)-d\tforal h\in {\cal H}.
$
%\eneq
Thus
\beq\label{Ldet-6}
t\circ \phi(h)>\Delta_0(\hat{h})-d>\Delta_0(\hat{h})/2\tforal h\in {\cal H} {{\andeqn \rforal t\in T(B)}}.
\eneq
%for all $t\in T(B).$
\end{proof}

%{ It also seems that one does not need to assume that $C$ is a commutative C*-algebra. So, I proposed a new statement, but with the old one remain---N}
%{However, I think that I proved the following somewhere, or so I think---L.})

\begin{lem}\label{Fullmeasure}
Let $C$ be a unital C*-algebra, and let $\Delta: C_+^{q, {\bf 1}}\setminus \{0\}\to (0,1)$ be an order preserving
map.  There exists  a map  $T\times N: C_+\setminus \{0\}\to \R_+\setminus \{0\}\times \N$
satisfying the following:
For any finite subset ${\cal H} \subset C_+^{\bf 1}\setminus \{0\}$ and any unital C*-algebra $A$ 
 with the strict comparison of positive elements, if $\phi: C \to A$ is a unital \morp\ satisfying
\beq\label{Fullm-1}
\tau\circ \phi(h)\ge \Delta(\hat{h})\tforal h\in {\cal H} \tforal \tau\in\mathrm T(A),
\eneq
then $\phi$ is $(T\times N)$-${\cal H}$-full.
\end{lem}

%\begin{lem}\label{Fullmeasure}
%Let $X$ be a compact metric space and $\Delta: C(X)_+^{{\bf 1}}\setminus \{0\}\to (0,1)$ be an order preserving
%map.  There exists  a map  $T\times N: C(X)_+\setminus \{0\}\to \R_+\setminus \{0\}\times \N$
%satisfying the following:
%For any finite subset ${\cal H}_1\subset C(X)_+^{\bf 1}\setminus \{0\},$ there exists a finite
%subset ${\cal H}_2\subset C(X)_+^{\bf 1}\setminus \{0\}$ such that if $\phi: C(X)\to A$ (for any unital infinite dimensional simple \CA\, $A$ with the strict comparison for positive elements) is a unital \morp\, with
%\beq\label{Fullm-1}
%\tau\circ \phi(h)\ge \Delta(\hat{h})\tforal h\in {\cal H}_2
%\eneq
%and for all $\tau\in T(A),$ then $\phi$ is $(T\times N)$-${\cal H}_1$-full.
%\end{lem}

\begin{proof}
For each $\delta\in(0, 1)$, let $g_\delta:[0, 1] \to [0, +\infty)$ be the continuous function defined by
$$g_{\delta}(t)=\left\{
\begin{array}{ll}
0 & \textrm{if $t\in[0, \frac{\delta}{4}]$},\\
\frac{f_{\dt/2}(t)}{t} & \textrm{otherwise},
\end{array}
\right.$$
where $f_{\dt/2}$ is as defined in \ref{Dball}. Note that
\begin{equation}\label{factor}
g_\delta(t) t  =  f_{\dt/2}(t)\rforal t\in[0, 1].
\end{equation}
Let $h\in C_+^{\bf 1}\setminus\{0\}$. Then define
$$
T(h)= \|(g_{\Delta(\widehat{h})})^{\frac{1}{2}}\|{=\frac 2\dt}\andeqn N(h)=\lceil\frac{2}{\Delta(\widehat{h})}\rceil.
$$
Then the function $T\times N$ satisfies the lemma.

Indeed, let $\mathcal H\subseteq C^1_+\setminus \{0\}$ be a finite subset. Let $A$ be a unital 
\CA\,  with
the strict comparison for positive elements, and let $\phi: C\to A$ be a unital positive linear map satisfying
\begin{equation}\label{eq-lb-tr}
\tau\circ\phi(h)\geq\Delta(\hat{h})\rforal h\in \mathcal H \rforal \tau\in \mathrm{T}(A).
\end{equation}
%Consider the positive element
%$(\phi(h)-\frac{\Delta(\hat{h})}{2})_+.$
Put $\phi(h)^{-}=(\phi(h)-\frac{\Delta(\hat{h})}{2})_+.$
Then, by  \eqref{eq-lb-tr}, one has that, since $0\le h\le 1,$
$$d_\tau(\phi(h)^{-})=d_\tau((\phi(h)-\frac{\Delta(\hat{h})}{2})_+)\ge \tau((\phi(h)-\frac{\Delta(\hat{h})}{2})_+)\geq \frac{\Delta(\hat{h})}{2}\rforal \tau\in\mathrm{T}(A).$$
 The above also shows that $\phi(h)^{-}\not=0.$
Since $A$ has the strict comparison for positive elements, one has
$K\left<(\phi(h)^{-}\right> > \left< 1_A \right>,$ where $K=\lceil \frac{2}{\Delta(\hat{h})}\rceil$
and where $\la x\ra$ denotes the class of $x$ in  $W(A).$

Therefore
there is a partial isometry $v=(v_{ij})_{K\times K}\in \mathrm{M}_{K}(A)$ such that
$$
vv^*=1_A\quad\textrm{and}\quad v^*v\in \overline{(\phi(h)^{-}\otimes 1_K) M_K(A)(\phi(h)^{-}\otimes 1_K)}.
$$
%Put $d_K(h,f)=f_{\dt/2}(\phi(h))\otimes 1_K.$ 
Note  that
$
c(f_{\dt/2}(\phi(h))\otimes 1_K)=c(f_{\dt/2}(\phi(h))\otimes 1_K)=c\rforal c\in\overline{(\phi(h)^{-}\otimes 1_K)M_K(A)(\phi(h)^{-}\otimes 1_K)},
$
where $\dt=\Delta(\hat{h}),$
%$$(\bigoplus_K h_{\Delta(\hat{f})}(\phi(f))) c =c (\bigoplus_K h_{\Delta(\hat{f})}(\phi(f))) = c\rforal c\in\textrm{Her}(\bigoplus_K (\phi(f)-\frac{\Delta(\hat{f})}{2})_+).$$ In particular,
%$$(\bigoplus_K h_{\Delta(\hat{f})}(\phi(f))) (v^*v) =(v^*v) (\bigoplus_K h_{\Delta(\hat{f})}(\phi(f))) = v^*v,$$
and therefore
$$
v(f_{\dt/2}(\phi(h))\otimes 1_K)v^*=vv^*=1_A.
$$
Consider the upper-left corner of $M_K(A),$ one has that
$$\sum_{i=1}^K v_{1,i}f_{\dt/2}(\phi(h))v_{1,i}^*=1_A,$$
and therefore, by \eqref{factor}, one has
$$
\sum_{i=1}^K v_{1,i} (g_{\Delta(\hat{f})}(\phi(h)))^{\frac{1}{2}}\phi(h)(g_{\Delta(\hat{f})}(\phi(h)))^{\frac{1}{2}}v_{1,i}^*=1_A.
$$
Since $v$ is a partial isometry, one has that $\|v_{i,j}\|\leq 1$, $i,j=1, ..., K$, and therefore
$$\|v_{1,i} (g_{\Delta(\hat{f})}(\phi(f)))^{\frac{1}{2}}\|\leq \|(g_{\Delta(\hat{f})}(\phi(f)))^{\frac{1}{2}}\|\leq \|(g_{\Delta(\hat{f})})^{\frac{1}{2}}\|=T(f).$$ Hence the map $\phi$ is $T\times N$-$\mathcal H$-full, as desired.

\end{proof}

%{\color{green} The following seems to be a corollary of the proof above, rather than the statement. So is the next corollary. (Maybe one need to be more careful: in general, C*-algebras in $\mathcal B_1$ might not be tracially approximately divisible. So, maybe one should assume that $A$ is UHF-stable.

%The torsion-free condition should be removed once such condition in \ref{UniqN1} is removed.

\begin{thm}\label{UniCtoA}
Let $C$ be a unital \CA\, in ${\overline{\cal D}}_s$ {\rm (see \ref{8-N-3})}.
% with finitely generated $K_*(C).$}
%with the form
%Let $C'\in \mathcal C_0$ be a unital \CA\,
%with torsion free $K_i(C')$
%and let
%$C=C'\otimes C(X),$ where $C'\in {\cal C}_0$ and
%where $X$ is a point, or $X=\T, $ or $C=PM_n(C(Y))P,$
%where $Y$ is a finite dimensional metric space and $P\in M_n(C(Y))$ be a projection.}
Let ${\cal F}\subset C$ be a finite subset, let $\ep>0$ be a
positive number and let $\Delta: C_+^{q, {\bf 1}}\setminus \{0\}\to
(0,1)$ be an {order preserving } map. There exists  a finite
subset ${\cal H}_1\subset C_+^{\bf 1}\setminus \{0\},$ there exists
$\gamma_1>0,$ $1>\gamma_2>0,$ $\dt>0,$ a finite subset ${\cal
G}\subset C$ and a finite subset ${\cal P}\subset \underline{K}(C),$
a finite subset ${\cal H}_2\subset C_{s.a.}$ and a finite subset
${\cal U}\subset \cup_{m=1}^{\infty}U(M_m(C))/CU(M_m(C))$ for which
$[{\cal U}]\subset {\cal P}$ satisfying the following: For any
unital ${\cal G}$-$\dt$-multiplicative \morp s $\phi, \psi: C\to A$,
where $A=A_1\otimes U$ for some
%unital separable simple C*-algebra
$A_1\in {\mathcal B_1}$ and a UHF-algebra $U$ of infinite type
satisfying
\beq\label{CUni-1}
&& [\phi]|_{\cal P}=[\psi]|_{\cal P},\\
\label{CUni-2}
&&\tau(\phi(a))\ge \Delta(\hat{a}),\,\,\,
\tau(\psi(a))\ge \Delta(\hat{a})
 \eneq
 for all $\tau\in T(A)$ and
for all $a\in {\cal H}_1,$
\beq\label{CUni-3}
 |\tau\circ\phi(a)-\tau\circ \psi(a)|&<&\gamma_1\tforal a\in {\cal
H}_2\tand\\\label{CUni1-3+1}
{\rm dist}(\phi^{\ddag}(u),
\psi^{\ddag}(u))&<&\gamma_2\tforal u\in {\cal U},
\eneq
there exists a
unitary $W\in A$ such that \beq\label{CUni-4}
\|W^*\phi(f)W-\psi(f)\|<\ep\tforal f\in {\cal F}. \eneq
\end{thm}

\begin{proof}
Let $T'\times N: C_+\setminus\{0\}\to\mathbb R_+\{0\}\times\mathbb N$ be the map of Lemma \ref{Fullmeasure} with respect to $C$ and $\Delta/4$. Let $T=2T'.$
%We may also assume that there is $n_0\ge 1$ such that the map $U(M_{n_0}(C))/U_0(M_{n_0}(C))$ to
%$U(M_{n_0+k}(A))/U_0(M_{n_0+k}(C))$ is isomorphism for all $k\ge 1.$

%For each $u\in U(M_2(A))$,
Define
${\bf L}=1.$
Let $\delta_0>0$ (in place of $\delta$), $\mathcal G_0\subseteq C$ (in place of $\mathcal G$), $\mathcal H_0\subseteq C_+\setminus\{0\}$ (in place of $\mathcal H$), $\mathcal U_0\subseteq {U(M_{n_0}(C))/CU(M_{n_0}(C))}$ (in place of $\mathcal U$), $\mathcal P_0\subseteq \underline{K}(C)$ (in place of $\mathcal P$) {be finite subsets} and {$n_1$ (in place of $n$) be an integer as required by} Theorem \ref{Suni} with respect to $C$ (in place of $A$), $T\times N$, ${\bf L}$, $\mathcal F$ and $\epsilon/2$ {----(see (2) {{of}} \ref{Rsuni})}.

Let $\mathcal H_{1, 1}\subseteq C_+^1\setminus\{0\}$ (in place of $\mathcal H_1$), $\mathcal H_{1, 2}\subseteq A$ (in place of $\mathcal H_2$), $\gamma_{1, 1}>0$ (in place of $\gamma_1$), $\gamma_{1, 2}>0$ (in place of $\gamma_2$), $\delta_{1}>0$ (in place of $\delta$), $\mathcal G_1\subseteq C$ (in place of $\mathcal G$), $\mathcal P_1\subseteq\underline{K}(C)$ (in place of $\mathcal P$), $\mathcal U_1\subseteq J_c(K_1(C))$ (in place of $\mathcal U$)
{and  $n_2$ (in place of $N$)} be the finite subsets and constants of Theorem \ref{UniqN1} with respect to $C$ (in place of $A$), $\Delta/4$, $\mathcal F$ and $\epsilon/4$.

Put $\mathcal G=\mathcal G_0\cup\mathcal G_1$, $\delta=\min\{\delta_0/4, \delta_1/4\}$, $\mathcal P=\mathcal P_0\cup\mathcal P_1$, $\mathcal H_1=\mathcal H_{1, 1}$, $\mathcal H_2=\mathcal H_{1, 2}$, $\mathcal U=\mathcal U_0\cup\mathcal U_1$, $\gamma_1=\gamma_{1, 1}{/2}$, $\gamma_2=\gamma_{1, 2}{/2}$. One asserts these are desired finite subsets and constants (for $\mathcal F$ and $\epsilon$).  {We may assume that
$\gamma_2<1/4.$}

In fact, let $A=A_1\otimes U$, where $A\in\mathcal B_1$ and $U$ is a UHF-algebra
{of infinite type}. Let $\phi, \psi: C\to A$ be $\mathcal G$-$\dt$-multiplicative maps satisfying {{(\ref{CUni-1})
to (\ref{CUni1-3+1}) for the above chosen ${\cal G},$ ${\cal H}_1,$ ${\cal P},$ ${\cal H}_2,$ ${\cal U},$
$\gamma_1$ and $\gamma_2.$}}
%
%the assumption
%of the l
%\begin{eqnarray}
%&&[\phi]|_{\cal P}=[\psi]|_{\cal P} \label{CUni-1-1}\\
%&&\tau(\phi(a))\geq \Delta(a)\ \textrm{and}\  \tau(\psi(a))\geq \Delta(a) \rforal \tau\in T(A)\rforal a\in {\cal H}_1, \label{CUni-2-2}\\
%&&|\tau\circ \phi(a)-\tau\circ \psi(a)|<\gamma_1\rforal a\in {\cal H}_2 \label{CUni-3-3}\\
%&&{\rm dist}(\phi^{\ddag}(u), \psi^{\ddag}(u))<\gamma_2\rforal u\in {\cal U} \label{CUni-4-4}.
%\end{eqnarray}

Since $A=A_1\otimes U$, $A\cong A\otimes U.$ Moreover,
$j\circ \imath: A\to A$ is approximately inner, where
$\imath: A\to A\otimes U$ is defined by
$a\mapsto  a\otimes 1_U$  and $j: A\otimes U\to A$ is an isomorphism.
Thus, we may assume that $A=A_1\otimes U\otimes U=A_2\otimes U,$ where $A_2=A_1\otimes U.$
Moreover, without loss of generality, we may assume that the images of $\phi$ and $\psi$
are in $A_2.$
%Put
%\beq\label{nn107-d=}
%d=\inf\{\Delta(\hat{h})/4: h\in {\cal H}_1\setminus \{0\}\}.
%\eneq
%{Therefore, for any $\eta_0>0$ with
%$\eta_0<\ep/8$ and finite subset ${\cal G}'\subset C$ with ${\cal G}'\supset {\cal G}\cup {\cal F},$
%there exists a unitary $W\in A\otimes U$ such that
%\beq\label{n107-2}
%\|\phi(c)-{\rm Ad}\, W\circ j\circ  (\phi_1'(c)+\phi_2'(c)+\cdots +\phi_m'(c))\|<\eta_0\andeqn\\
%\|\psi(c)-{\rm Ad}\, W\circ j\circ (\psi_1'(c)+\psi_2'(c)+\cdots +\psi_m'(c))\|<\eta_0
%\eneq
%Put $e_i={\rm Ad}\, W\circ j(1_A\otimes e_i'),$
%$\phi_i={\rm Ad}\, W\circ j\circ \phi_i'$ and $\psi_i={\rm Ad}\, W\circ j\circ \psi_i',$ $i=1,2,...,m.$
%there are projections $e_1, e_2, ..., e_{m}\in A$ and unitaries $u_1=1, u_2, ..., u_m\in A$ for some $m\geq n+1$ such that
%{enumerate}
%\beq\nonumber
%{e_1+\cdots+e_m=1_A\andeqn
%e_i=u_i^*e_1u_i \rforal 1\leq i\leq m,}
%\eneq
%\end{enumerate}
%and  moreover,
%there are unital $\delta$-$\mathcal G$-multiplicative maps
%$$\phi_i, \psi_i: C\to e_iAe_i,\quad 0\leq i\leq n$$ such that
%\begin{enumerate}
%{It is also clear that there are unitaries $u_i\in A$ such that $e_i=u_i^*e_1u_i$ and
%\beq\nonumber
%\phi_i=u_i^*\phi_1 u_i,\,\,\psi_i=u_i^*\psi_1 u_i,\,\,\,1\leq i\leq m,
%\eneq
%Moreover, for any $a\in\mathcal {\cal G}'$, one has
  %       \begin{eqnarray}\label{decp-uhf-CA}
    %     \|\phi(a)-\phi_1(a)\oplus\cdots\oplus\phi_m(a)\|&<&\eta_0,\\
     %    \|\psi(a)-\psi_1(a)\oplus\cdots\oplus\psi_m(a)\|&<&\eta_0.
      %   \end{eqnarray}
%By choosing sufficiently small $\eta_0,$  with the image of $\phi$ and $\psi$, one may also assume that
Since $A_2\in\mathcal B_1$,   for a finite subset ${\cal G}''\subseteq A_2$
and $\delta'>0$, there is a projection $p\in A_2$, a C*-subalgebra
$D\in \mathcal C_1$
%{ You want $D\in {\mathcal C}$ not
%$\mathcal C_0$, otherwise you need $A_2\in {\cal B}_0$}
with $p=1_D$
such that
\begin{enumerate}
\item $\|pg-gp\|<\delta'$ for any $g\in\mathcal G'',$
\item $pgp\in_{\delta'} D$,
\item $\tau(1-p)<\min\{\delta', \gamma_1/4, 1/8m\}$ for any $\tau\in\mathrm{T}(A)$.
\end{enumerate}
{Define $j_0: A_2\to (1-p)A_2(1-p)$ by $j_0(a)=(1-p)a(1-p)$ for all $a\in A_2.$
For any $\ep''>0$ and any finite subset ${\cal F}''\subset A_2,$
there is also a unital \morp\, $j_1: A\to D$ such that
$\|j_1(a)-pap\|<\ep''$  for all $a\in {\cal F}'',$ provided that ${\cal G}''$ is sufficiently large and
$\dt'$ is sufficiently small. }
Therefore, in particular, we may assume that
\beq\label{nn107-1-}
&&\|\phi(c)-(j_0\circ \phi(c)\oplus j_1\circ \phi(c))\|<\ep/16\andeqn\\
&&\|\psi(c)-(j_0\circ \psi(c)\oplus j_1\circ \psi(c))\|<\ep/16
\eneq
for all $c\in {\cal F}.$

{Choose an integer $m\ge 2(n_1+1)n_2$ and mutually orthogonal and
mutually equivalent projections $e_1,e_2...,e_m\in U$ with
$\sum_{i=1}^me_i=1_U.$} {Define $\phi_i', \psi_i': C\to A\otimes U$
by $\phi_i'(c)=\phi(c)\otimes e_i$ and $\psi_i'(c)=\psi(c)\otimes
e_i$ for all $c\in C,$ $i=1,2,...,m.$ Note that \beq\label{n107-1}
[\phi_1']|_{{\cal P}}=[\phi_i']|_{\cal P}=[\psi'_1]|_{\cal
P}=[\psi_i']|_{\cal P}, \eneq $i=1,2,...,m.$ } {Note also that $\phi_i',
\psi_i': C\to e_iAe_i$ are ${\cal G}$-$\dt$-multiplicative.}

Write $m=kn_2+r,$ where $k\ge n_1+1$ and $r<n_2$ are integers.
%{ I think that $r$ should not be so much small for using
%\ref{ph}. I think it should be like $\frac 1{2(n_1+1)} <\frac rm
%<\frac 1{n_1+1}$ or something similar, see comment below.}
Define\\
$
{\tilde \phi}, {\tilde \psi}: C\to (1-p)A_2(1-p)\oplus \bigoplus_{i=kn_2+1}^m A_2\otimes e_i
$
by
\beq\label{nn107-1}
&&{\tilde\phi}(c)=j_0\circ \phi(c)\oplus \sum_{i=kn_2+1}^m  j_1\circ \phi(c)\otimes e_i\andeqn\\
&&{\tilde \psi}(c)=j_0\circ \psi(c)\oplus \sum_{i=kn_2+1}^m  j_1\circ
\psi(c)\otimes e_i \eneq for all $c\in C.$  With sufficiently large
${\cal G}''$ and small $\dt',$ we may assume that ${\tilde \phi}$
and ${\tilde \psi}$ are ${\cal G}$-$2\dt$-multiplicative
%%%multiplictative
 and, by
(\ref{n107-1}), \beq\label{nn107-2} [{\tilde \phi}]|_{{\cal
P}}=[{\tilde \psi}]|_{{\cal P}}. \eneq Moreover,  by \ref{Tapprox},
we may further assume that
%{ \ref{Tapprox} did not talk about
%this part, it only says about the part of $j_1\circ \phi$ and
%$j_1\circ \psi$. I think we need to to use \ref{ph}, knowing from
%\ref{nn107-10++} that ${\rm dist}({\tilde \phi}^{\ddag}({\bar
%v}\oplus 1), {\tilde \psi}^{\ddag}({\bar v}\oplus 1))$ are very
%small, but need to do more detailed estimation}
\beq\label{nn107-2+}
{\rm dist}({\tilde \phi}^{\ddag}({\bar v}), {\tilde
\psi}^{\ddag}({\bar v}))<\gamma_2\le {\bf L} \eneq for all ${\bar
v}\in {\cal U}.$ Define  $\phi_i^1,\psi_i^1: C\to D\otimes e_i$ by
$\phi_i^1(c)=j_1\circ \phi(c)\otimes e_i$ and $\psi_i^1(c)=j_1\circ
\psi(c)\otimes e_i.$
%By (\ref{nn107-d=}) ,  (3) above,  and \
By \ref{Tapprox} and by choosing even larger ${\cal G}''$ and smaller $\dt',$ we may assume
that
\beq\label{nn107-10}
\tau\circ \phi_i^1(h)\ge \Delta(\hat{h})/2\andeqn \tau(\psi_i^1(h))\ge \Delta(\hat{h})/2\rforal h\in {\cal H}_1
\eneq
and for all $\tau\in T(pAp\otimes e_i),$
\beq\label{nn107-10+}
|t\circ \phi_i^1(c)-t\circ \psi_i^1(c)|<\gamma_{1,1}\rforal c\in {\cal H}_2\andeqn
 \eneq
\beq\label{nn107-10++}
{\rm dist}((\phi_i^1)^{\ddag}({\bar v}),(\psi_i^1)^{\ddag}({\bar v})) <\gamma_{1,2} \rforal {\bar v}\in {\cal U}.
 \eneq
By applying \ref{Fullmeasure},  $\phi_i^1$ and $\psi_i^1$ are $T\times N$-${\cal H}_1$-full.
Moreover, we may also assume that
$\phi_i^1$ and $\psi_i^1$ are ${\cal G}$-$2\dt$-multiplicative \morp s and
\beq\label{nn107-10+++}
[\phi_i^1]|_{\cal P}=[\psi_i^1]|_{\cal P}.
\eneq

Define
$\Phi, \Psi: C\to  \bigoplus_{i=1}^{kn_2} D\otimes e_i$
by
\beq\label{nn107-3}
\Phi(c)=\bigoplus_{i=1}^{kn_2} \phi_i^1(c)\andeqn
\Psi(c)=\bigoplus_{i=1}^{kn_2} \psi_i^1(c)
 \eneq
for all $c\in C.$
By (\ref{nn107-10+++}), (\ref{nn107-10}), (\ref{nn107-10+}), (\ref{nn107-10++}) and
by \ref{UniqN1}, there exits a unitary $W_1\in (\sum_{i=1}^{kn_2}p\otimes e_i)(A_2\otimes U)(\sum_{i=1}^{kn_2}p\otimes e_i)$ such that
\beq\label{nn107-15}
\|W_1^*\Phi(c)W_1-\Psi(c)\|<\ep/4\rforal c\in {\cal F}.
\eneq
 Note  that
 \beq\label{nn107-16}
  \tau(1-p)+\sum_{kn_2+1}^m\tau(e_i)<(1/m)+(r/m)\le n_2/m
 \eneq
  for all $\tau\in T(A).$  Note also that $k\ge n_1.$
  By (\ref{nn107-2}), (\ref{nn107-2+}), since $\psi_i^1$ is $T\times N$-${\cal H}_1$-multiplicative, by applying \ref{Suni}, there exists a unitary $W_2\in A$ such that
  \beq\label{nn107-17}
  \|W_2^*({\tilde \phi}(c)\oplus \Psi(c))W_1-({\tilde \psi}(c)\oplus \Psi(c))\|<\ep/2
  \rforal c\in {\cal F}.
  \eneq
  Set
  $$
  W=({\rm diag}(1-p, e_{kn_2+1}, e_{kn_2+2},...,e_m)\oplus W_1)W_2.
  $$
  Then
  we compute that
  \beq\label{nn107-18}
  \|W^*({\tilde \phi}(c)\oplus \Phi(c))W-({\tilde \psi}(c)\oplus \Psi(c))\|<\ep/2+\ep/4
  \eneq
  for all $c\in {\cal F}.$  By (\ref{nn107-1-}),
  we have
  \beq\label{nn107-19}
  \|W^*\phi(c)W-\psi(c)\|<\ep\rforal c\in {\cal F},
  \eneq
  as desired.

\end{proof}

\begin{rem}\label{RemUniCtoA}

{\rm
First we note that,  given (\ref{CUni-3}) and with an arrangement at the beginning of the proof, 
we only need   one inequality  in  (\ref{CUni-2}).  %If  $K_1(C)=\{0\}$ or $K_1(C)$ is torsion, then one does not need ${\cal U}$ and $\gamma_2$ in the statement.
%{Moreover,
Note also that the condition that ${\cal U}\subset \cup_{m=1}^{\infty}U(M_m(C))/CU(M_m(C))$ can be replaced by
${\cal U}\subset J_c(K_1(C)).$ To see this, we note that we may write, by \ref{Dcu},
\beq\label{RemUniC-1}
\bigcup_{m=1}^{\infty}U(M_m(C))/CU(M_m(C))=\Aff(T(C))/\overline{\rho_C(K_0(C))}\oplus J_c(K_1(C)).
\eneq
So, without loss of generality, we may assume that
${\cal U}={\cal U}_0\cup {\cal U}_1,$ where ${\cal U}_0\subset \Aff(T(C))/\overline{\rho_C(K_0(C))}$
and ${\cal U}_1\subset J_c(K_1(C)).$
Using the de La Harpe and Skandalis determinant (see again \ref{Dcu}),
with sufficiently small $\dt,$ $\gamma_1$ and sufficiently large ${\cal G}$ and ${\cal H}_1,$  the condition (\ref{CUni-3})
implies that
\beq\label{RemUniC-2}
{\rm dist}(\phi^{\ddag}(u),\psi^{\ddag}(u))<\gamma_2 \rforal u\in {\cal U}_0.
\eneq
This  particularly implies that, when $K_1(C)=0,$  ${\cal U}$ and $\gamma_2$ are not needed in the
statement of \ref{UniCtoA}. When $K_1(C)$ is torsion, $\phi^{\ddag}|_{J_c(K_1(C))}=0,$
since $\Aff(T(A))/\overline{\rho_A(\rho_A(K_0(A))}$ is torsion free (see \ref{UCUdiv}).

In general, we can further assume that ${\cal U}\subset J_c(K_1(C))$ generates a free group.
Let $G({\cal U})$ be the subgroup generated by ${\cal U}.$
Write $G({\cal U})=G_1\oplus G_t,$ where $G_0$ is free and $G_t$ is torsion.
Let $g_1, g_2,...,g_{k_1}$ be the generators of $G_1$ and $f_1, f_2,..., f_{k_2}$ be the generators of $G_t.$ As in the proof, we may assume that $[x]\in {\cal P}$ for all
$x\in {\cal U}',$ where ${\cal U}'$ is a finite subset of unitaries
such that $\overline{[{\cal U'}]}={\cal U}.$  By choosing even larger ${\cal P}$ and
${\cal U},$ we may assume that $g_i, f_j\in {\cal U}.$  There is
an integer $m(j)\ge 1$ such that $m(j)f_j=0.$ It follows
that $m(j)\phi^{\ddag}(f_j)=m(j)\psi^{\ddag}(f_j)=0.$ With  sufficiently small $\dt,$ $\gamma_1$ and sufficiently large ${\cal G}, $ condition (\ref{CUni-1}) implies
that
$$
\pi\circ \psi^{\ddag}(f_j)=\pi\circ \phi^{\ddag}(f_j),
$$
where $\pi: \bigcup_{m=1}^{\infty} U(M_m(A))/CU(M_m(A))\to K_1(A)$ is the quotient map.
Therefore
\beq\label{141105-1}
J_c\circ \pi\circ \psi^{\ddag}(f_j)=J_c\circ \pi\circ \phi^{\ddag}(f_j).
\eneq
Note $m(j)J_c\circ \pi\circ \phi^{\ddag}(f_j)=0.$
It follows that
$$
m(j)(\phi^{\ddag}(f_j)-J_c\circ \pi\circ \phi^{\ddag}(f_j))=0.
$$
However,
$$
\phi^{\ddag}(f_j)-J_c\circ \pi\circ \psi^{\ddag}(f_j)\in {\rm Aff}(T(A))/\overline{\rho_A(K_0(A))}.
$$
Since $A_1\in {\cal B}_1,$ by \ref{Unotrosion}, ${\rm Aff}(T(A))/\overline{\rho_A(K_0(A))}$ is torsion free. Therefore
$$
\phi^{\ddag}(f_j)-J_c\circ \pi\circ \phi^{\ddag}(f_j)=0
$$
Similarly,
$
\psi^{\ddag}(f_j)-J_c\circ \pi\circ \psi^{\ddag}(f_j)=0.
$
Thus, by (\ref{141105-1}),
\beq\label{141105-2}
\phi^{\ddag}(f_j)=\psi^{\ddag}(f_j).
\eneq
In other words, with sufficiently small $\dt,$ $\gamma_1$ and sufficiently large ${\cal G}, $ (\ref{CUni-1}) implies (\ref{141105-2}).  So $G_t$ can be dropped. Therefore
 in the statement in \ref{UniCtoA}, we may assume
that ${\cal U}$ generates a free subgroup of $J_c(K_1(C)).$

Furthermore, if $C$ has stable rank $k,$ then $J_c(K_1(C))$ may be replaced by $J_c(U(M_k(C))/U_0(M_k(C))),$ see
(\ref{Dcu-6}).
In the case that  $C$ has stable rank one,
then ${\cal U}$ may be assumed to be a subset of $J_c(U(C)/U_0(C)).$ In the case that $C=C'\otimes C(\T)$ for
some $C'$ with stable rank one,  then stable rank of $C$ is no more than $2.$ Therefore, in this case,
${\cal U}$ may be assumed to be in $J_c(U(M_2(C))/U_0(M_2(C)),$ or  $U(M_2(C))/CU(M_2(C)).$

}

%\end{rem}

%\begin{rem}\label{Re2MUN}
{\rm
The introduction of $\Delta$ and  condition (\ref{CUni-2}) are for convenience which can be replaced by original
fullness condition. $\Delta$ can be replaced by a map $T\times N: C_+\setminus \{0\}$ and condition
(\ref{CUni-2}) is replaced by
$\phi$ and $\psi$ are $T\times N$-${\cal H}_1$-full as indicated in the proof.
}
\end{rem}

%{\color{green} With the new forms of Theorem \ref{UniqN1}, the next Corollary can be merged with \ref{UniCtoA}.}

%\begin{cor}\label{CCtoAuni}
%Let $C=PM_m(C(X))P,$ where $X$ is a finite CW complex and $P\in M_m(C(X))$ is a projection.
%Let ${\cal F}\subset C,$ let
%$\ep>0$ be a positive number and let $\Delta: C_+^q\setminus \{0\}\to (0,1).$  There exists  a finite subset
%${\cal H}_1\subset C_+^{\bf 1}\setminus \{0\},$
%there exists $\gamma_1>0,$ $\gamma_2>,$ $\dt>0,$ a finite subset
%${\cal G}\subset C$ and a finite subset ${\cal P}\subset \underline{K}(C),$ a finite subset ${\cal H}_2\subset C_{s.a.}$ and a finite subset ${\cal U}\subset J_c(K_1(C))$ for which $[{\cal U}]\subset {\cal P}$ satisfying the following:
%For any unital $\dt$-${\cal G}$-multiplicative \morp s $\phi, \psi: C\to A$, where {\color{green} $A=A_1\otimes U$}
%for some unital separable simple \CA\, $A_1\in {\mathcal B_1}$ such that
%\beq\label{CUni-1}
%[\phi]|_{\cal P}=[\psi]|_{\cal P},
%\eneq
%\beq\label{CUni-2}
%\tau(\phi(a))\ge \Delta(a),\,\,\, \tau(\psi(a))\ge \Delta(a)
%\eneq
%for all $\tau\in T(A)$ and for all $a\in {\cal H}_1,$
%\beq\label{CUni-3}
%|\tau\circ \phi(a)-\tau\circ \psi(a)|<\gamma_1\tforal a\in {\cal H}_2\tand\\\label{Uni1-3+1}
%{\rm dist}(\phi^{\ddag}(u), \psi^{\ddag}(u))<\gamma_2\tforal u\in {\cal U},
%\eneq
%there exists a unitary $W\in A$ such that
%\beq\label{CUni-4}
%\|W^*\phi(f)W-\psi(f)\|<\ep\tforal f\in {\cal F}.
%\eneq
%\end{cor}

%\begin{rem}
%{\rm If $K_1(C)$ is torsion, then one does not need ${\cal U}$ and $\gamma_2$ in the statement.}
%
%\end{rem}

\begin{cor}\label{Unitaryuni}
Let
$\ep>0$ be a positive number and let $\Delta: C(\T)_+^{q, {\bf 1}}\setminus \{0\}\to (0,1)$ be {{an order
preserving}} map.  There exists  a finite subset
${\cal H}_1\subset C(\T)_+^{\bf 1}\setminus \{0\},$
there exists $\gamma_1>0,$ $1>\gamma_2>0$ and any finite subset ${\cal H}_2\subset C(\T)_{s.a.}$
satisfying the following:
For any two unitaries $u_1$ and $u_2$ in a unital separable simple
\CA\, $A\in {\cal B}_1$
such that
\beq\label{UUni-1}
[u_1]=[u_2]\in K_1(A) {{\tand}}
%\eneq
%\beq\label{UUni-2}
\tau(f(u_1)),\,\tau(f(u_2))\ge \Delta(\hat{f})
\eneq
for all $\tau\in T(C)$ and for all $f\in {\cal H}_1,$
\beq\label{UUni-3}
|\tau(g(u_1))-\tau(g)(u_2)|<\gamma_1\tforal g\in {\cal H}_2\tand
%\\\label{CCUni1-3+1}
{\rm dist}(\bar{u_1}, \bar{u_2})<\gamma_2,
\eneq
there exists a unitary $W\in C$ such that
\beq\label{UUni-4}
\|W^*u_1W-u_2\|<\ep.
\eneq
\end{cor}

\begin{lem}\label{tensorprod}
Let $A$ be a \CA\, and $X$ be a compact metric space.
Suppose that $y\in (A\otimes C(X))_+\setminus \{0\}.$
Then exists $a(y)\in A_+\setminus \{0\},$ $f(y)\in C(X)_+\setminus \{0\}$ and
$r_y\in A\otimes C(X)$ such that $\|a(y)\|\le \|y\|,$ $\|f(y)\|\le 1,$ $\|r_y\|\le \|a\|$ and
$r_y^*yr_y=a(y)\otimes f(y).$
\end{lem}

\begin{proof}
Identify $A\otimes C(X)$ with $C(X,A).$
Let $x_0\in X$ such that $\|y(x_0)\|=\|y\|.$
There is $\dt>0$  such that $\|y(x)-y(x_0)\|<\|y\|/16$ for   $x\in B(x_0, 2\dt).$
Let $Y=\overline{B(x_0, \dt)}.$  Let $z(x)=f_{\|y\|/4}(y(x_0))$ for all $x\in Y.$
Note $z\not=0.$
By 2.2 of \cite{RorUHF2}, there exist $r\in C(Y, A)$ with $\|r\|\le \|y\|$ such that
$r^*(y|_Y)r=z.$ Choose $g\in C(X)_+\setminus \{0\}$ such that
$0\le g\le 1,$ $g(x)=0$ if ${\rm dist}(x, x_0)\ge \dt$ and $g(x)=1$ if ${\rm dist}(x, x_0)\le \dt/2.$
Then one may view $rg^{1/2}, zg\in C(X, A).$ Put $r_y=rg^{1/2},$ $a(y)=f_{\|y\|/4}(y(x_0))$ and
$f(y)=g.$ Then
$$
r_y^*yr_y=zg=a(y)\otimes f.
$$
\end{proof}

\begin{thm}\label{MUN1}
Let $A_1\in {\cal B}_1$
% {{{\cal B}_0}}$
be a unital simple \CA\, which satisfies the UCT,
 %$A_1=A_0\otimes U_0,$ where $U_0$ is a UHF-algebra of infinite type and
 $A=A_1\otimes C(X),$ where $X$ be a point or $X=\T.$
For any $\ep>0,$ any finite subset ${\cal F}\subset A$  and any order preserving map $\Delta:
C(X)_+^{{\bf 1}}\setminus \{0\}\to (0,1),$
there exists
$\dt>0,$ a finite subset ${\cal G}\subset A,$  $\sigma_1, \sigma_2>0,$ a finite subset ${\cal P}\subset \underline{K}(A),$ a finite subset ${\cal H}_1\subset C(X)_+^{\bf 1}\setminus\{0\},$ a finite subset ${\cal U}\subset {U(M_2(A))/CU(M_2(A))}$
{(---see \ref{ReMUN1})}
and a finite subset ${\cal H}_2\in A_{s.a}$ satisfying the following:

Let $B'\in {\cal B}_1,$ let $B=B'\otimes U$ for some  UHF-algebra $U$ of infinite type and let $\phi, \psi: A \to B$ be two unital ${\cal G}$-$\dt$-multiplicative \morp s such that
\beq\label{MUN1-1}
[\phi]|_{\cal P}&=&[\psi]|_{\cal P},\\\label{MUN1-1+}
\tau\circ \phi(1\otimes h) &\ge & \Delta(h)\tforal  h\in {\cal H}_1\tand \tau\in T(B),\\
|\tau\circ \phi(a)-\tau\circ \psi(a)|&<&\sigma_1\tforal a\in {\cal H}_2\tand\\\label{MUN1-1+2}
{\rm dist}(\phi^{\ddag}({\bar{u}}), \psi^{\ddag}({\bar u}))&<&\sigma_2\tforal {\bar u}\in {\cal U}.
%{\rm dist}(\overline{\langle \phi(u)\rangle}, \overline{ \langle \psi(u)\rangle})&<&\sigma_2\tforal u\in {\cal U}.
\eneq
Then there exists a unitary $u\in U(B)$ such that
\beq\label{ }
\|{\rm Ad}\, u\circ \phi(f)-\psi(f)\|<\ep\tforal f\in {\cal F}.
\eneq

\end{thm}

\begin{proof}
Let $\ep>0$ and let ${\cal F}\subset A$ be a finite subset.
Without loss of generality, we  may assume that
$$
{\cal F}=\{ a\otimes f: a\in {\cal F}_1\andeqn f\in {\cal F}_2\},
$$
where ${\cal F}_1\subset A$ is a finite subset and ${\cal F}_2\subset C(X)$
is also a finite subset. We further assume that
${\cal F}_1$ and ${\cal F}_2$ are  in the unit ball  of $A$ and $C(X),$ respectively.

{Let ${\bf L}=1.$}
%${\bf L}: U(M_2(A))\to \R_+$ by ${\bf L}(v)=8\pi+1$ for all $v\in U(M_2(A)).$
%be defined as follows.
%If $u\in U_0(A),$ define ${\bf L}(u)=4\pi({\rm cel}(u)+1)+1;$ if
%$u\in U(A)\setminus U_0(A),$ define
%${\bf L}(u)=8\pi+1.$
Let $\Delta: C(X)_+^{\bf 1}\setminus \{0\}\to (0,1)$ be an order preserving map.
Let $T'\times N': C(X)_+\setminus \{0\}\to \R_+\setminus \{0\}\times \N$ be a map  given by \ref{Fullmeasure}  with respect to  $3\Delta/16.$
Since $A_1$ is a unital separable simple \CA, the identity map on $A_1$  is
$T''\times N''$-full for some $T''\times N'': (A_1)_+\setminus \{0\}\to \R_+\setminus \{0\}\times \N.$

Define  a map $T\times N: A_+\setminus \{0\}\to \R_+\setminus \{0\}\times \N$ as follows:
For any $y\in A_+\setminus \{0\},$ {{by \ref{tensorprod}}},
there exist $a(y)\in (A_1)_+\setminus \{0\},$  $f(y)\in C(X)_+\setminus \{0\}$
{{and $r_y\in A\otimes C(X)$  with $\|r_y\|\le \|y\|$
 such that $r_y^*yr_y=a(y)\otimes f(y),$  $\|a(y)\|\le \|y\|$ and $\|f(y)\|\le 1.$ }}
%with $(a(y)\otimes f(y))^{1/2}\le y.$

There are $x_{a(y),1},x_{a(y),2},...,x_{a(y), N''(a(y))}\in A_1$ with
$\max\{\|x_{a(y),i}\|: 1\le i\le N''(a(y))\}=T''(a(y))$ such that
$$
\sum_{i=1}^{N''({a(y)})}x_{a(y),i}^*a(y)x_{a(y),i}=1_{A_1}.
$$
Then define
\beq\label{MUN1-6-3}
\hspace{-0.2in}(T\times N)(y)=(1+\max\{T''(a(y)), T'(f(y))\} \cdot \max\{1, \|y\|\}, N''(a(y)) \cdot N'(f(y))).
\eneq
%There are $x_{a,1},x_{a,2},...,x_{a, N''(a)}\in A_1$ with $\max\{\|x_{a,i}\|: 1\le i\le N''(a)\}=T''(a)$ such that
%$$
%\sum_{i=1}^{N''(a)}{x_{a,i}^*ax_{a,i}=1_{A_1}
%$$
%and there are $y_{f,1},y_{f,2},..., y_{f, N'(b)}\in C(X)$ with $\max\{\|y_{f,j}\|:1\le j\le N'(f)\}=T'(f)$
%such that
%$$
%\sum_{j=1}  ^{N'(f)} y_{f,j}^*fy_{f,j}=1_{C(X)}.
%$$
%Let $z_{i,j}=x_{a,i}\otimes y_{f,j},$ $1\le i\le N''(a),$ $1\le j\le N'(f).$ Then
%\beq\label{MUN1-6-2}
%\sum_{i,j}z_{i,j}^*(a\otimes f) z_{i,j}=1_A
%\eneq
%Define
%\beq\label{MUN1-6-3}
%(T\times N)(x)=(\max\{T''(a), T'(f)\}\,,N''(a)N'(f)).
%\eneq
%We have
%\beq\label{MUN1-6-4}
%(a\otimes b)^{1/4} x(a\otimes b)^{1/4}\ge a\otimes b.
%\eneq

 Let $\ep/16>\dt_1>0$ (in place of $\dt$), ${\cal G}_1\subset A$ (in place of ${\cal G})$,
 ${\cal H}_0\subset A_+\setminus \{0\}$ (in place of $\mathcal H$),
${\cal U}_1\subset {U(M_2(A))/CU(M_2(A))}$ (in place of ${\cal U}$---see  (2) of \ref{Rsuni}), ${\cal P}_1\subset \underline{K}(A)$ (in place of ${\cal P}$) and $n\ge 1$ be the finite subsets and constants as required by \ref{Suni} for $A,$ ${\bf L},$ $\ep/16$ (in place of $\ep$), ${\cal F}$ and $T\times N.$
Without loss of generality, we may assume that $\dt_1<\ep$ and
\beq\label{MUN1-6}
{\cal G}_1=\{a\otimes g: a\in {\cal G}_1'\andeqn g\in {\cal G}_1''\},
\eneq
where ${\cal G}_1'\subset A_1$  and
${\cal G}_1''\subset C(X)$ are finite subsets. We may further assume
that ${\cal F}_1\subset {\cal G}_1'$ and ${\cal F}_2'\subset {\cal G}_1'',$
and both are in the unit ball.  In particular,
${\cal F}\subset {\cal G}_1.$
We may also assume that
\beq\label{MUN1-6+}
{\cal H}_0=\{a\otimes f: a\in {\cal H}_0'\andeqn f\in {\cal H}_0''\},
\eneq
where ${\cal H}_0'\subset (A_1)_+\setminus \{0\}$ and ${\cal H}_0''\subset C(X)_+\setminus \{0\}$ be finite
subsets. {For convenience, we may further assume that,  {{for the above  integer $n\ge 1,$}}
\beq\label{n107-n1}
1/n<\inf\{\Delta(h): h\in {\cal H}_0''\}/16.
\eneq
}
{Let ${\cal U}_1=\{\bar{v_1},\bar{v_2},...,\bar{v_K}\},$ where
$v_1,v_2,...,v_K\in U(M_2(A)).$ Put
${\cal U}_0=\{v_1,v_2,...,v_K\}.$}
%For each $v\in {\cal U}_0,$ there are $u_{v,1},u_{v,2},..., u_{v, k(v)}\in U_0(M_2(A))$
%with $k(v)\le {\rm cel}(v)+1$ such that
%\beq\label{MUN1-6++}
%u_{v,1}=v,\,\,\, \|u_{v,j}-u_{v,j+1}\|<1,\,\,\, j=1,2,...,k(v)\andeqn u_{v,k(v)}=1_{M_2(A)}.
%\eneq
Choose a finite subset ${\cal G}_u'\subset A$ such that
\beq\label{MUN1-6+++}
{v_j\in \{(a_{i,j})_{1\le i, j\le 2}: a_{i,j}\in {\cal G}_u'\}\tforal v_j\in {\cal U}_0.}
\eneq
Choose $\dt_1'>0$ and a sufficiently large finite subset ${\cal G}_v\subset A_1$ which satisfying the following:
If $p\in A_1$ is a projection such that
$$
\|px-xp\|<\dt_1'\tforal x\in {\cal G}_v,
$$
then {there are unitaries $w_j\in ({\rm diag}(p,p)\otimes
1_{C(X)})M_2(A)({\rm diag}(p,p)\otimes 1_{C(X)})$ such that
\beq\label{n107-n2} \|{\rm diag}(p,p)v_j{\rm
diag}(p,p)-w_j\|<\dt_1/16n\tforal v_j\in {\cal U}_0\andeqn
j=1,2,...,K. \eneq }
Let
\beq\nonumber
{\cal G}_2'={\cal F}_1\cup {\cal G}_1'\cup {\cal H}_0'\cup {\cal G}_u\cup \{a(y),x_{a(y),j}, x_{a(y),j}^*: y\in {\cal H}_0'\},
\andeqn\\
M_1=64(\max\{\|x_{a(y), j}\|: y\in {\cal H}_0'\}+1) \cdot \max\{N(y): y\in {\cal H}_0'\}.
\eneq
Put $\dt_1''=\min\{\dt_1', \dt_1\}/(64(n+1)M_1).$
%It follows from \ref{decomp}
{Since $A_1\in  {\cal B}_1,$
%{\red{{\cal B}_0}},$
there exists mutually orthogonal
projections $p_0', p_1'\in A_1,$
%a unitary $u_j'\in A_1,$
%$j=1,2,...,2n-1,$
a \SCA\, $C\in {{{\cal C}}}$
%$K_1(C)=\{0\}$
and  $1_{C}=p_1$, unital {$\dt_1''/16$-${\cal G}_2'$}-multiplicative
\morp s $\imath_{00}': A_1\to p_0'A_1p_0'$ and $\imath_{01}': A_1\to
C$ such that
\beq\label{MUN1-7} {{\rm diag}(\overbrace{p_0',
p_0',...,p_0'}^{n+1})\lesssim p_1'\andeqn \|x-\imath_{00}'(x)\oplus
\imath_{01}'(x)\|<\dt_1''}
\eneq
% { the above is not enough,
%$p_0'$ should be small but not to small to apply \ref{ph}. I suggest
%%to change to ${{\rm diag}(\overbrace{p_0',
%p_0',...,p_0'}^{n+1})\lesssim p_1'\lesssim {\rm
%diag}(\overbrace{p_0', p_0',...,p_0'}^{2(n+1)})} $}
for all $x\in
{\cal G}_1'\cup {\cal G}_2',$ where $\imath_{00}'(a)=p_0ap_0$ for
all $a\in A_1.$} Define $p_0=p_0'\otimes 1_{C(X)},$ $p_1=p_1'\otimes
1_{C(X)}.$

\Wlog, we may assume that $p_0'\not=0.$  Since $A_1$  is simple, there is an integer $N_0>1$ such that
\beq\label{MUN1-n7}
N_0[p_0']\ge [p_1']\,\,\,{\rm in}\,\,\, W(A_1).
\eneq
This also implies that
\beq\label{MUN1-n8}
N_0[p_0]\ge [p_1].
\eneq

 Define
$\imath_{00}: A\to p_0Ap_0$ by $\imath_{00}(a\otimes f)=\imath_{00}'(a)\otimes f$ and $\imath_{01}: A\to C\otimes C(X)$ by
$\imath_{01}(a\otimes f)=\imath_{01}'(a)\otimes f$ for all $a\in A_1$ and
$f\in C(X).$
Define $L_0: A\to A$ by
$$
L_0(a)=\imath_{00}(a)\oplus \imath_{01}(a) \rforal a\in A.
$$
{For each $v_j\in {\cal U}_0,$ there exists a unitary $w_j\in M_2(p_0Ap_0)$ such that
\beq\label{n107-nn1}
\|{\rm diag}(p_0,p_0)v_j{\rm diag}(p_0,p_0)-w_j\|<\dt_1/16n,\,\,\,j=1,2,....
\eneq
}
Note that $C_1=C\otimes C(X){\subset  A}.$
%{\color{green} Since the $K_0(C)$ is torsion free and $K_1(C)=\{0\}$, one has that $K_*(C_1)$ is torsion free.} We may also assume that $C=M_m(C_0'),$
%where $m\ge 1$ is an integer and $C_0'\in {\cal C}_0$ so that $C_0'$ has the property (P).{\color{green} (Lemma \ref{prime-C0}.)}

Let $\imath_0: C\to A$ be the natural embedding as $C$ is a unital \SCA\, of $p_1A_1p_1.$ Let $\imath_0^{\sharp}: C^q\to A_1^q$ be
defined by $\imath_0^{\sharp}(\hat{c})=\hat{c}$ for $c\in C.$
Let $\Delta_0: A_1^{q,{\bf 1}}{{\setminus \{0\}}}\to (0,1)$ be the map given by \ref{Ldet}
%{\red{corresponding to the map $T''\times N''$}}
and
define
$$
\Delta_1(\hat{h})
=\sup\{\Delta_0(\imath_0^{\sharp}(\hat{h_1}))\Delta(h_2)/4: h\ge h_1\otimes h_2,\,\,\,
h_1\in C\setminus \{0\}\andeqn h_2\in C(X)_+\setminus \{0\}\}
$$
for all $h\in (C_1)_+\setminus \{0\}.$

Let ${\cal G}_3'=\imath_{0,1}'({\cal G}_1'\cup {\cal G}_2').$
Let ${\cal G}_3=\{a\otimes f: a\in {\cal G}_3'\andeqn f\in {\cal G}_1''\}.$
Let ${\cal H}_3\subset (C_1)_+\setminus \{0\}$  (in place of ${\cal H}_1$), $\gamma_1'>0$ (in place of $\gamma_1$), $\gamma_2'>0$ (in place of $\gamma_2$), $\dt_2>0$ (in place of $\dt$), ${\cal G}_4\subset C_1$ (in place of ${\cal G}$), ${\cal P}_2\subset \underline{K}(C_1)$ (in place of ${\cal P}$), ${\cal H}_4'\subset (C_1)_{s.a.}$ (in place of ${\cal H}_2$),
$\overline{{\cal U}_2}\subset {U(M_2(C_1))/U_0(M_2(C_1))}$  {(in place of ${\cal U}$---see \ref{RemUniCtoA})} be the finite subsets and constants as  required by {\ref{UniCtoA}} for $\dt_1/16$ (in place of $\ep$)
and ${\cal G}_3$ (in place of ${\cal F}$),
$\Delta_1/2$ (in place of $\Delta$)  and for $C_1$ (in place of $A$).

Let ${\cal U}_2\subset U(M_2(C_1))$ be a finite subset {which  has an one-to-one correspondence to its image
in $U(M_2(C_1))/CU(M_2(C_1))$ which is exactly  $\overline{{\cal U}_2}.$}  We also assume
that $\{[u]: u\in {\cal U}_2\}\subset {\cal P}_2.$

Without loss of generality, we may assume that
\beq\label{MUN-8}
{\cal H}_3=\{h_1\otimes h_2: h_1\in {\cal H}_3'\andeqn h_2\in {\cal H}_3''\},
\eneq
where ${\cal H}_3'\subset C_+\setminus \{0\}$ and ${\cal H}_3''\subset C(X)_+\setminus \{0\}$ are
finite subsets, and
\beq\label{MUN-8+1}
{\cal G}_4=\{a\otimes f: a\in {\cal G}_4'\andeqn f\in {\cal G}_4''\},
\eneq
where ${\cal G}_4'\subset C$ and ${\cal G}_4''\subset C(X)$ are finite subsets.

%For each $x\in {\cal H}_3'\subset C\subset A_1,$ there are

Let $\dt_3>0$ (in place of $\dt$) and let ${\cal G}_5\subset A_1$ (in place of ${\cal G}$) be the finite subset
as required by \ref{Ldet} for $\Delta_0$ and ${\cal H}_0'\cup {\cal H}_3'$.

%Let ${\cal H}_4={\cal H}_0\cup {\cal H}_3.$

Set
\beq\label{UNN1-9}
\dt={\min\{1/16, \ep/16, \dt_1, \dt_2, \dt_3\}\over{128N_0(n+1)(T({\cal H}_4)+N({\cal H}_4))}},
\eneq
and
set
$$
{\cal G}_6=\{{\cal G}_2'\cup L_0({\cal G}_2')\cup \imath_{01}'({\cal G}_2')\cup\bigcup_{j=1}^{n-1}\{{\rm Ad}\, u_j\circ \imath_{01}'({\cal G}_2')\}\cup {\cal G}_4'\cup {\cal G}_5\andeqn
$$
$$
{\cal G}= \{a\otimes f: a\in {\cal G}_6  \andeqn f\in {\cal G}_1''\cup {\cal H}_0''\cup {\cal G}_4''\}\cup\{p_j: 0\le j\le 1\}\cup\{v_j,w_j: 1\le j\le K\} .
$$
To simplify notation, without loss of generality, we may assume that ${\cal G}\subset A^{\bf 1}.$
% {\color{green} We may assume assume that $\delta$ sufficiently small and $\mathcal G$ sufficiently large so that the restriction of any $2\delta$-$\mathcal G$-multiplicative map $L: A\to B$ to $C\otimes 1_{C(X)}$  is $T''\times N''$-$\imath_{01}(\mathcal H_0')$-full.}
Let ${\cal P}={\cal P}_1\cup \{p_j: 0\le j\le 1\}\cup [\imath]({\cal P}_2),$ where $\imath: C_1\to A$ is the embedding.
Let ${\cal H}_1={\cal H}_0''\cup {\cal H}_3''.$

{Let ${\cal U}_2'=\{{\rm diag}(1-p_1,1-p_1)+w: w\in {\cal U}_2\}$ and let
${\cal U}_0''=\{ w_j+{\rm diag}(p_1, p_1): 1\le j\le K\}.$
Let ${\cal U}=\{\bar{v}: v\in {\cal U}_1\cup {\cal U}_2'\}$
and let } ${\cal H}_2={\cal H}_4''.$
Let $\sigma_1=\min\{{1\over{4n}},{\gamma_1'\over{16nN_0}}\}$ and $\sigma_2=\min\{{1\over{16nN_0}},{\gamma_2'\over{16nN_0}}\}.$

Now we assume that $B$ is as in the statement, $\phi,\,\,\psi: A\to B$ are two unital ${\cal G}$-$\dt$-multiplicative
\morp s satisfying the assumption for  the above defined $\dt,$ ${\cal G},$
${\cal P},$ ${\cal H}_1,$ ${\cal U},$ ${\cal H}_2,$ $\sigma_1$ and $\sigma_2.$

Note that $B'$ is in  ${\mathcal B_1}$  {and $B=B'\otimes U.$
We may also write $B=B_1\otimes U,$ where $B_1=B'\otimes U,$ since $U$ is strongly self absorbing.
Without loss of generality, by the fact that $U$ is strongly self absorbing, we may assume
that the image of both $\phi$ and $\psi$ are in $B_1.$}
{By replacing $\psi$ by ${\rm Ad}\, u_0$ for some unitary $u_0\in U(B_1)$ if necessary, we may assume
that
\beq\label{n107-11}
\psi\circ \imath_{00}'(1_A)=\phi\circ \imath_{00}'(1_A)=q.
\eneq
There is an integer  $m\ge n$ and mutually orthogonal and mutually unitarily equivalent
projections $e_1,...,e_{m}\in U$ such that $\sum_{j=1}^{m}e_i=1_U.$}
{Define $\phi_0', \psi_0': A\to qB_1q\otimes 1_U$ (see (\ref{n107-11})) by
\beq\label{n107-12}
\phi_0'(a)=\phi\circ \imath_{00}(a)\otimes 1_U\andeqn
\psi_0'(a)=\psi\circ \imath_{00}(a)\otimes 1_U
\eneq
for all $a\in A.$ Define $\Phi', \Psi: A\to (1-q)B_1(1-q)\otimes 1_U$ by
\beq\label{n107-13}
\Phi'(a)=\phi\circ \imath_{01}(a)\otimes 1_U\andeqn
\Psi'(a)=\phi\circ \imath_{01}(a)\otimes 1_U
\eneq
 for all $a\in A.$
 Define $\Phi, \Psi: C_1\to (1-q)B_1(1-q)$ by
 \beq\label{n107-14}
 \Phi=\phi\circ \imath\andeqn \Psi=\psi\circ \imath.
 \eneq
Define $\psi_i': A\to (1-q)B_1(1-q)\otimes e_i$ by
\beq\label{n107-14+}
\psi_i'(a)=\Psi'(a)e_i\rforal a\in A.
\eneq
}
{Note, by the choice of $\dt$ and ${\cal G},$ $\Phi$ and $\Psi$ are
${\cal G}$-$\dt$-multiplicative. By (the proof of) \ref{Ldet},
\beq\label{n107-15} \tau(\Phi(h))\ge \Delta_1(\hat{h})/2\rforal h\in
{\cal H}_3. \eneq By assumptions, for all $\tau\in T(B_1),$
\begin{equation}\label{n107-16}
|\tau(\Phi(c))-\tau(\Psi(c))|<\sigma_1\rforal c\in {\cal H}_4''.
\end{equation}
Therefore, for all $t\in T((1-q)B_1(1-q)),$
\begin{equation}\label{n107-17}
|t(\Phi(c))-t(\Psi(c))|<\gamma_1'\rforal c\in{\cal H}_4''.
\end{equation}
Since $[\imath]({\cal P}_2)\subset {\cal P},$ by
assumptions, one  also has \beq\label{n107-18} [\Phi]|_{{\cal
P}_2}=[\Psi]|_{{\cal P}_2}. \eneq One also computes that (as
elements in $M_2((1-q)B_1)(1-q))$)
\beq\label{n107-19}
{\rm
dist}(\Phi^{\ddag}({\bar v})\Psi^{\ddag}({\bar
v^*}))<\gamma_2'\rforal v\in {\cal U}_2.
\eneq
By the choices of
$\dt,$ ${\cal G}_4,$ $\gamma_1',$ $\gamma_2',$ ${\cal P}_2,$ ${\cal
H}_3,$ ${\cal H}_4''$ and $\overline{{\cal U}_2},$ and by applying
\ref{UniCtoA}, there exists $u_1\in (1-q)B_1(1-q)$ such that
\beq\label{n107-20}
\|u_1^*\Phi(c)u_1-\Psi(c)\|<\dt_1/16\rforal c\in
{\cal G}_3. \eneq Thus, by  (\ref{MUN1-7}), \beq\label{n107-21}
\|u_1^*\Phi'(a)u_1-\Psi'(a)\|<\dt_1/16 +\dt'' \rforal a\in {\cal
G}_1.
\eneq }
We check that, by  assumption  (\ref{MUN1-1+}) and
(\ref{n107-n1}),
\beq\label{n107-22}
\tau(\Psi'(h))\ge
15\Delta(\hat{h})/16 \rforal h\in {\cal H}_0'' \andeqn\rforal
\tau\in T(B_1). \eneq
By \ref{Fullmeasure}, it follows that  $\Psi'$
is $T'\times N'$-${\cal H}_0''$-full. Note that \beq\label{n107-23}
\|\sum_{i=1}^{N''(a(y))}\Psi'(x_{a(y),i}^*)\Psi'(a(y))\Psi'(x_{a(y),i})-(1-q)\|<
4N''(a(y))T''(a(y))\dt \eneq for all $a\in {\cal H}_0'.$ We conclude
that $\Psi'$ is $T\times N$-${\cal H}_0$-full. It follows that
$\psi_i'$ is $T\times N$-${\cal H}_0$-full, $i=1,2,...,m.$ It is
easy to see that
\beq\label{n107-24}
[\psi_0']|_{{\cal
P}_1}=[\phi_0']|_{{\cal P}_1}.
\eneq
Assumptions also imply that
\beq\label{n107-25}
{\rm dist}(\phi^{\ddag}(\overline{w_j+{\rm
diag}(p_1, p_1)}), \psi^{\ddag}(\overline{w_j+{\rm diag}(p_1,
p_1)}))<\sigma_2
\eneq $j=1,2,...,K.$ By \ref{MUN1-n8}, applying (\ref{ph})
%{
%for apply \ref{ph}, we need the supporting projection of $w_j$ not
%too much smaller that $p_1$, please revise}
and then
(\ref{n107-nn1}), we compute that \beq\label{n107-26} {\rm
dist}((\phi_0')^{\ddag}(v_j),(\psi_0')^{\ddag}(v_j))<\gamma_2',\,\,\,
j=1,2,...,K. \eneq It follows from \ref{Suni} and its remark
\ref{Rsuni} that there exists a unitary $u_2\in B$ such that
\beq\label{n107-27} \|u_2^*(\phi'(a)\oplus \psi_1'(a)\oplus\cdots
\oplus \psi_n'(a))u_2-(\psi_0'(a)\oplus \psi_1'(a)\oplus \cdots
\oplus \psi_n'(a))\| <\ep/16 \eneq for all $a\in {\cal F}.$ In other
words, \beq\label{n107-28} \|u_2^*(\phi'_0(a)\oplus
\Psi'(a))u_2-\psi_0'(a)\oplus \Psi'(a)\|<\ep/16\rforal a\in {\cal
F}. \eneq Thus, by (\ref{n107-21}), \beq\label{n107-29}
\|u_2^*(\phi_0'(a)\oplus u_1^*\Phi'(a)u_1)u_2-\psi_0(a)\oplus
\Psi(a)\|<\ep/16+\dt_1/16+\dt''\rforal a\in {\cal F}. \eneq Let
$u=(q+u_1)u_2\in U(B).$ Then,  by \eqref{MUN1-7}, one has
\begin{equation}\label{n107-30}
\|u^*\phi(a)u-\psi(a)\|<\ep\rforal a\in {\cal F},
\end{equation}
as desired.
\end{proof}

\begin{rem}\label{ReMUN1}
As in  remark \ref{RemUniCtoA}, the condition that ${\cal U}\subset U(M_2(A))/CU(M_2(A))$ can be replaced
by ${\cal U}\subset J_c(U(M_2(A))/U_0(M_2(A))).$ Moreover, if $X$ is a point, or  equivalently,
$A\in {\cal B}_1,$ we may take ${\cal U}\subset J_c(U(A)/U_0(A)),$ since $A$ has stable rank one. Furthermore,
in this case , we do not need $\Delta$ and (\ref{MUN1-1+}).
%It should also be that,  in the proof of \ref{MUN1},  we do not assume that the algebra
%$D$ has the property that $K_1(D)=\{0\},$ nor we assume that $K_1(C)=\{0\}.$
%Therefore Theorem \ref{MUN1} holds for
%the case that $A$ is a unital simple \CA\, which is ``tracially ${\cal C}_0$''.
Let $G$ be the subgroup generated by the finite subset ${\cal U}\subset J_c(K_1(A))$ and
$G=G_1\oplus Tor(G_1),$ where $G_1$ is free.
 Since $\Aff(T(B))/{\overline{\rho_B(K_0(B))}}$ is torsion free (by \ref{UCUdiv}),
we may further assume that $G$ is free.

In the case that $A=\overline{\bigcup_{n=1}^{\infty} A_n},$ in the theorem above, one may choose ${\cal U}$ to be
in $A_n$ for some sufficiently large $n.$
\end{rem}

\section{The range of invariant}
% and models for the C*-algebras in the class $\mathcal B_0$ with (SP)}

\begin{nota}\label{homrestr}
Let $A$ be  a unital subhomogeneous C*-algebra, that is,  the maximal dimension of irreducible representations of $A$  is finite. Let us use $RF(A)$ to denote the set of equivalence classes of all (not necessarily irreducible) finite dimensional representations. In this section,  by $\{y\}=\{x_1,x_2, ...,  x_k\}$, we mean that the representation corresponding to $y$ is direct sum of the representations $x_1,x_2, ..., x_k$.  If some of  $x$ repeats $k$ times, then  we use $x^{\sim k}$ to denote this.
%In other words, $\{x_1^{\sim k_1},x_2^{\sim k_2}, \cd, x_n^{\sim k_n}\} $ is a finite dimensional
%representation equivalent to the direct sum over all of $k_i$ copies of the representation corresponding  to $x_i$ for  $1\le i\le n.$
That is, $\{y\}$ may be written as $\{z_1^{\sim k_1}, z_2^{\sim k_2},...,z_m^{\sim k_m}\},$ where
for each $j,$ $z_j=x_i$ for some $i.$
Note that we do not insist that each $z_i$ should be irreducible.

%n Let $A$ be  a subhomogeneous algebras whose maximal dimension of irreducible representations  is a finite number.
Let us use $Sp(A)$ to denote the set of {equivalence classes of} all irreducible representations  of $A$, which may be viewed as a subset of $RF(A)$. Since $A$ is of type I,  the set $Sp(A)$ has an one-to-one correspondence  to the set of primitive  ideals of $A.$
%n, which has a nature structure of topology (though this topology is in general not Hausdorff).
Let $X\subset Sp(A)$ be a closed subset, then $X$ corresponds to the idea $I_X=\bigcap_{\psi\in X} \ker \psi$. In this section, let us use  $A|_X$ to denote the quotient algebra $A/I_X$.
%It is  convenient to use  $A|_X$ to denote the quotient algebra $A/I_X$, and we will do so in this section.
If $\phi: B \to A$ is a homomorphism, then we will use $\phi|_X: B \to A|_X$ to denote the composition $\pi\circ \phi$, where $\pi: A \to A|_X$ is the quotient map. As usual, if $B_1$ is a subset of $B$, we will also use $\phi|_{B_1}$ to denote the restriction of $\phi$ on $B_1$. These two notation will not be confused, since it will be clear from content which notation we refer to.

If $ \phi: A \to B$ is a homomorphism, then we write  $Sp(\phi)= \{x\in Sp(A): \ker \phi \subseteq \ker x\}$.

\end{nota}

\begin{nota}\label{ktimes}
 In this section we will use the concept   of sets with multiplicity. Therefore $X_1=\{x,x,x,y\}$ is different from $X_2=\{x,y\}.$
 %, since in the first set, $x$ appears  three times and in the second set it appears only once.
 We will also  use $x^{\sim k}$ for a simplified notation for $\{\underbrace{x,x,...,x}_k\}$ (see 1.1.7 of \cite{Gong-AH}). For example, $\{x^{\sim 2}, y^{\sim 3}\}=\{x,x,y,y,y\}$. Let $X=\{x_1^{\sim i_1}, x_2^{\sim i_2}, ... , x_n^{\sim i_n}\}$ and $Y= \{x_1^{\sim j_1}, x_2^{\sim j_2}, ... , x_n^{\sim j_n}\}$ (some of the $i_k$'s (or $j_k$'s) may be zero which means the element $x_k$  does not appear in the set $X$ (or $Y$)). If $i_k\leq j_k$ for all $k=1,2, ... , n$, then we say that $X\subset Y$ (see 3.21 of \cite{Gong-AH}). We  define
 $$
 X\cup Y=\{x_1^{\sim \max(i_1,j_1)},x_2^{\sim \max(i_2,j_2)}, ... , x_n^{\sim \max(i_n,j_n)}\}.
 $$
   By $X^{\sim k}$ we mean the set $\{x_1^{\sim ki_1},x_2^{\sim ki_2}, ... , x_n^{\sim ki_n}\}$.

% Let $A$ be  a subhomogeneous algebras whose maximal dimension of irreducible representations  is a finite number. We use $RF(A)$ to denote the set of equivalence classes of all (not necessarily irreducible) finite dimensional representations. In this section,  by $\{y\}=\{x_1,x_2, \cd, x_k\}$, we mean the representation corresponding to $y$ is direct sum of the representations $x_1,x_2,\cd, x_k$.  If some of  $x$ repeats $k$ times, then  we use $x^{\sim k}$ to denote it.
%In other words, $\{x_1^{\sim k_1},x_2^{\sim k_2}, \cd, x_n^{\sim k_n}\} $ is a finite dimensional
%representation equivalent to the direct sum over all of $k_i$ copies of the representation corresponding  to $x_i$ for  $1\le i\le n.$
%That is, $\{y\}$ may be written as $\{z_1^{\sim k_1}, z_2^{\sim k_2},...,z_m^{\sim k_m}\},$ where
%for each $j,$ $z_j=x_i$ for some $i.$
%It should be noted that we do not insist any $z_i$ should be irreducible.

 If $\phi:A \to M_m(\C)$ is a homomorphism, let us use $SP(\phi)$ to denote the corresponding equivalent class of $\phi$ in $RF(A)$. Any finite subset of $RF(A)$ also defines an element in $RF(A)$ which is the equivalent class of  the direct sum of all corresponding representations in the set with correct multiplicities. If both $X$ and $Y$  are finite subsets of $RF(A)$ with multiplicities, we say $X\subset Y$ if  the representation corresponding to  $X$ is equivalent to a sub-representation of that corresponding to   $Y$. That is,  if  we rewrite $X$ and $Y$ as $X=\{x_1^{\sim i_1}, x_2^{\sim i_2},\cd, x_n^{\sim i_n}\}$ and $Y= \{x_1^{\sim j_1}, x_2^{\sim j_2},\cd, x_n^{\sim j_n}\}$, with $x_i$ being irreducible representation, then $i_k\leq j_k$ for each $k\in \{1,2,\cd,n\}$. It is clear that for two homomorphisms $\phi_1: A \to M_{m_1}(\C)$ and $\phi_2: A \to M_{m_2} (\C)$,  we have $SP(\phi_1)\subset SP(\phi_2)$ if and only if  $\phi_1$ is equivalent to sub-representation of $\phi_2$. Strictly speaking, an element in $RF(A)$ is regarded as a set with multiplicity whose elements are in $Sp(A)$. But when we write $X=\{z_1^{\sim k_1}, z_2^{\sim k_2},\cd, z_m^{\sim k_m}\},$  we do not insist that $z_i$ itself in $Sp(A)$. It may be a list of several elements in $Sp(A)$---that is, we do not insists $z_i$ to be irreducible (but as we know, it can always be decomposed into irreducible ones). So in this notation, we do not differentiate $\{x\}$ and $x$, both give same element in $RF(A)$ and same set with multiplicity whose elements in $Sp(A)$.

 Comparing  with notation in \ref{homrestr}, if $\phi: A\to M_m(\C)$ is a homomorphism and if $SP(\phi)=\{x_1^{\sim k_1},x_2^{\sim k_2}, ... ,x_i^{\sim k_i}\}$, with $x_1, x_2, ... ,x_i$ being irreducible representations, then \\
 $Sp(\phi)=\{x_1,x_2,...,x_i\}\subset Sp(A).$ So $Sp(\phi)$ is an ordinary set which is a subset of $Sp(A)$, while $SP(\phi)$ is a set with multiplicity, whose elements are also elements in ${Sp(A)}$.
%{\bf {Suggestion:   Switch $sp(\phi)$ and $SP(\phi).$  So $sp(A)$ and $sp(\phi)$ are ordinary sets and
%$SP(\phi)$ for more complicated sets.---Done}}
\end{nota}

\begin{nota}\label{density}
Let us recall some notations from \ref{8-N-3}. Suppose that $A=A_m\in{\cal D}_m$  constructed in the following sequence
\beq\nonumber
A_0=F_0,~~ A_1=A_0\oplus_{Q_1C(Z_1, F_1)Q_1}P_1C(X_1, F_1)P_1,\\
~~A_2=A_1\oplus_{Q_2C(Z_2, F_2)Q_2}P_2C(X_2, F_2)P_2, \quad \cdots ~~~~~~~~~~~~
\eneq
is as in \ref{8-N-3}.
Let $\LD: A \to \bigoplus_{k=0}^mP_kC(X_k, F_k)P_k,$ be the inclusion homomorphism as \ref{8-N-3}.  Let  $\pi_{(x,j)}$, where $x\in X_k$ and $j$ the positive integer represented $j$-th block of
$A_k=\bigoplus_{j=1}^{l_k}P_{k,j}C(X_k, F_k)P_{k,j}$,   be the  finite dimensional representation of $A$ as in \ref{8-N-3}. According to \ref{homrestr}, one has that $\pi_{(x,j)}\in RF(A)$.

 Suppose all $X_k$ are metric spaces. Let $\phi: A \to M_{\bullet}(\C)$ be a homomorphism. We say that $SP(\phi)$ is $\dt$-dense in $Sp(A)$ if for each irreducible representation $\theta$ of $A$, there are two points $x, y, \in X_k$ and $j$ (same $k$ and $j$) such that $dist (x,y) <\dt$, and such that
$\tht \subset \pi_{(y,j)}$ and $\pi_{(x,j)}\subset SP(\phi)$. This will be used in this section and next. 

If $X_k=[0,1]$, then we use the ordinary metric of interval $[0,1]$.
\end{nota}

\begin{NN}\label{simplelimit}

Let $A=\varinjlim (A_n,\phi_{n,m})$ be an inductive limit with all $A_n\in D_{l(n)}$ for some $l(n)$. The limit $C^*$ algebra is simple if the inductive system satisfy the following condition: for any $n$ and $\dt>0$, there is an integer $m>n$ such that for any $y\in Sp(A_m)$, $SP(\phi_{n,m}|_y)$ is $\dt$-dense---consequently, for any $m'>m$ and any $y\in Sp(A_{m'})$, $SP(\phi_{n,m'}|_y)$ is $\dt$-dense. The proof of this claim is standard. Indeed, for any nonzero element $f\in A_n$, it must be nonzero on a open ball of radius $\dt$ (for some $\dt>0$) of a certain $X_k$ (a space appears in the construction of $A_n$). Then, by the given condition, there is an $m$, such that for any $m'>m$ and any  irreducible representation $\tht$ of $A_{m'}$, one has $\tht(\phi_{n,m'}(f))\not=0$. {It follows that the ideal $I$ generated by $\phi_{n,m'}(f))$ in $A_{m'}$ equals $A_{m'}$---otherwise, any irreducible representation $\tht$ of $A_{m'}/I\not=0$ (which is also a irreducible representation of $A_{m'}$) satisfies $\tht(\phi_{n,m'}(f))=0$, and this contradicts with the fact that $\tht(\phi_{n,m'}(f))\not=0$ for all irreducible representation $\tht$.} It is well known that an ideal $I$ of the limit C*-algebra $\lim(A_n,\phi_{n,m})$ is always the limit of ideals $\lim(I_n, \phi_{n,m}|_{I_n})$, where $I_n=\phi_{n,\infty}^{-1}(I)$. Consequently, the  $C^*$-algebra of the given inductive system is  simple.
%{\color{Green} (One should be careful for the last sentence. It is only shown that the algebraic limit is topologically simple inside itself. To see the C*-algebra limit is simple, it seems that some argument should be inserted here. But why not we just use the well-known fact that the ideal of the limit C*-algebra is always a limit of ideals.)}
%{\red{The original statement is same thing as Zhuang suggested, anyway I changed it Zhuang's suggestion--Guihua

\end{NN}

\begin{df}\label{Class0}
Denote by ${\cal N}_0$ the class of those unital  simple \CA s
$A$  in ${\cal N}$ for which $A\otimes U\in {\cal
N}\cap {\cal B}_0$, for any  UHF-algebra $U$ of
infinite type (see \ref{dfcalN} for definition of class ${\cal N}$).

Denote by ${\cal N}_1$ the class of those unital \CA s $A$ in ${\cal N}$ for which $A\otimes U\in {\cal N}\cap {\cal B}_1$, for any  UHF-algebra of
infinite type. In section 19, we will show that ${\cal N}_1={\cal N}_0.$

Also denoted by ${\cal N}_0^{\cal Z}$ (respectively ${\cal N}_1^{\cal Z}$) the class of all ${\cal Z}$-stable $C^*$-algebras in ${\cal N}_0$ (respectively ${\cal N}_1$).

\end{df}

%xxxxxxxabc

\begin{NN}\label{range 0.1}
Let $(G,G_+,u)$ be a scaled ordered abelian group $(G,G_+)$ with order
unit $u\in G_+\setminus \{0\}$.
{ The scale is given by
$\{g\in G_+: g\le u\}$. Sometimes we will also call $u$ the scale of the group.}
Let $S(G):=S_u(G)$ be the state space of $G.$
% Denote by $S(G)$ the state space
%of $G$ -- that is, $f\in S(G)$ is given by group homomorphism $f:~
%G\to \R$ with $f(G^+)\sbseq \R_+$ and $f(u)=1$.
Suppose that
$((G,G_+,u),K,\Delta,r)$ is a weakly  unperforated Elliott invariant---that is, $(G,G_+,u)$ is a simple scaled ordered countable group,
%{\color{Green}{which is weakly unpeforated}},
$K$ is a countable
abelian group, $\DT$ is a metrizable Choquet simplex, and $r: \DT \to S(G)$
is a surjective affine map such that for any $x\in G$,
\beq\label{1508/star13-1}
 x\in G_+\setminus \{0\} \qq \mbox{ if and only if } \qq r(\tau)(x)>0
\qq \mbox{ for all } \qq\tau\in \DT.
\eneq
The above condition (\ref{1508/star13-1}) is also called  weakly  unperforated for the
simple ordered group. (Note that this condition is equivalent to that
$x\in G_+\setminus \{0\}$ if and only if for any $f\in S(G)$,
$f(x)>0$.  The latter condition does not mention Choquet simplex
$\DT$.) {In this paper, we only consider the Elliott invariant for stably finite simple $C^*$-algebras and therefore $\DT$ is not empty. The above weakly unperforated condition is also equivalent to that $x>0$ if $nx>0$ for some positive integer $n$.}

In this section, we will show that for any weakly unperforated Elliott invariant $((G,G_+,u),K,
\Delta,r)$, there is a unital simple C*-algebra $A$ in the class ${\cal N}_0^{\cal Z}$ such that $$ ((K_0(A),K_0(A)_+,
[{\mbox{\large \bf 1}}_A]),K_1(A),T(A),r_A)\cong ((G,G_+,u),K,
\Delta,r).
$$
In \cite{point-line}, Elliott constructed an inductive limit $A$ with
${\rm Ell}(A)=((G,G_+,u),K,\Delta, r)$.  In this section, we will use
modified building blocks to construct a simple \CA\,  $A$ such that
${\rm Ell}(A)$ is as described, and, in addition, as a direct consequence of the construction, one has that $A\in {\cal N}_0^{\cal Z}$. %, which
%-- that is, $A\otimes M_{\mathfrak{p}}\in {\cal N}\cap{\cal B}_0$
%for every supernatural number ${\mathfrak{p}}$ which
%, tracially
%approximately Elliott-Thomsen trivial $K_1$ algebra. (In fact, this
%will be a direct consequence of the construction.

\end{NN}

\begin{NN}\label{range 0.2}  Our construction will be a modification of the Elliott
construction mentioned above.  As matter
of fact, for the case that $K=\{0\}$ and $G$ torsion free, our
construction uses the same building blocks in ${\cal C}_0,$ as in \cite{point-line}.  We will repeat a part of the
construction of Elliott for this case.  There are two steps in
Elliott construction:

Step 1.  Construct an inductive limit
$$ C_1\lr C_2 \lr \cd\lr C $$
with inductive limit of ideals
$$ I_1\lr I_2 \lr \cd\lr I $$
such that the non-simple limit $C$ has the described Elliott invariant
and the quotient $C/I$ is a simple AF algebra. For the case of $K=\{0\}$ and $G$ torsion free, we will use the notation $C_n$ and $C$ for the construction, and reserve $A_n$ and $A$ for the general case.

\vspace{0.1in}

Step 2.  Modify the above inductive limit to make $C$ (or $A$ in the general case) simple without
changing the Elliott invariant of $C$ (or $A$).

For reader's convenience, we will repeat Step 1 of Elliott construction with minimum modification. For Step 2, we will use a slightly different way to modify the inductive limit which will be more suitable for our purpose---that is, to construct an inductive limit $A\in {\cal
N}_0^{\cal Z}$ with  possible nontrivial $K_1$ and nontrivial
${\rm Tor}(K_0(A))$.

\end{NN}

\begin{NN}\label{range 0.3}

Let $\rho:~ G\to \Aff(\DT)$ be the
dual map of $r:\DT \to S(G)$.  That is, for every $g\in G, ~
\tau\in \DT,$
$$
\rho (g)(\tau)= r(\tau)(g)\in \R.
$$
%Then by the condition of weakly unperforated $(*)$ of \ref{range 0.1},
Since $G$ is weakly unperforated, one has that $g\in G_+\setminus \{0\}$ if and only if $\rho (g)(\tau)>0$
for all $\tau\in \DT$.  Note that $\Aff(\DT)$ is an ordered vector
space with {the strict order}, i.e., $f\in \Aff(\DT)_+\setminus\{0\}$ if and only if $f(\tau)
>0$ for all $\tau\in\DT$.  Since $\DT$ is a metrizable compact convex set,
there is a countable dense subgroup $G^1\subset \Aff(\DT).$   Put
$H=G\oplus G^1$ and define  { ${\tilde \rho}: H\to \Aff(\DT)$ by
${\tilde \rho}((g,f))(\tau)=\rho(g)(\tau)+f(\tau)$ for all $(g,f)\in G\oplus G^1$ and
$\tau\in \DT.$ Define
$H_+\setminus\{0\}$ to be the set of  elements
$(g,f)\in G\oplus G^1$ with ${\tilde \rho}((g,f))(\tau)>0$ for all $\tau\in \DT.$}
%$$
%\qq\qq\qq\qq\rho (g)(\tau)+f(\tau)>0\qq\qq \mbox{for all } \qq\qq
%\tau\in \DT.
%$$
The order unit (or scale) $u\in G_+$ could be regarded as $(u,0)\in G\oplus G^1=H$
as the order unit of $H_+$ (still denote it by $u$).  Then $(H, H_+,u)$
is a simple  ordered group.   Using the fact that ${\tilde \rho}(H)$ is dense and therefore has the Riesz
interpolation property,
%Using the fact I
it is straight forward to prove that
$(H, H_+,u)$ is a  Riesz group. We give a brief proof of this fact as below. Let $a_1, a_2, b_1, b_2\in H$ with $a_i<b_j$ for $i,j\in\{1,2\}$, then for each $\tau\in\DT$, $\max\{\tilde\rho(a_1)(\tau),\tilde\rho(a_2)(\tau)\}< \min\{\tilde\rho(b_1)(\tau),\tilde\rho(b_2)(\tau)\}$. Let $d=\min_{\tau\in \DT}\left( {\min}\{\tilde\rho(b_1)(\tau),\tilde\rho(b_2)(\tau)\}- \max\{\tilde\rho(a_1)(\tau),\tilde\rho(a_2)(\tau)\}\right) >0$. Since ${\tilde \rho}(H)$ is dense in $\DT$, there is a $c\in H$ such that $\max\{\tilde\rho(a_1)(\tau),\tilde\rho(a_2)(\tau)\}< \tilde\rho(c)(\tau) < \min\{\tilde\rho(b_1)(\tau),\tilde\rho(b_2)(\tau)\}$ for all $\tau\in \DT$. Consequently, $a_i<c<b_j$ for $i,j\in\{1,2\}$.  As a
direct summand of $H$, the subgroup $G$ is {\it relatively divisible}
subgroup of $H,$ i.e., if $g\in G,$  $m\in \N\setminus \{0\},$ and $h\in H$
such that $g=mh,$ then there is $g'\in G$ such that $g=mg'.$
Note that $G\subsetneqq H$.

Now we assume, until \ref{range 0.14}, that $G$ is torsion free and $K=0.$
Then  $H$ is also torsion free. Therefore $H$ is a dimension group.

%\vspace  {0.2in}
\end{NN}

 % abcxxxxxx

\begin{NN}\label{range 0.4}  In \ref{range 0.3} we can choose the dense
subgroup $G^1\subset  \Aff(\DT)$ to contain at least three elements
$x,y,z\in \Aff(\DT) $ such that $x,y$ and $z$ are $\Q$-linearly
independent.
% -- that is, if $\frac{n_1}{q_1} x +\frac{n_2}{q_2}y
%+\frac{n_3}{q_3}z=0$ for $n_1,n_2, n_3\in \Z;
%q_1,q_2,q_3\in\Z\setminus \{0\}$, then $n_1=n_2=n_3=0$.
With this
choice, when we write $H$ as inductive limit
$$
H_1\lr H_2\lr \cd
$$
of finite  direct sums of copies of ordered group $(\Z,)_+)$ as in Theorem 2.2
of \cite{EHS} , we can assume all $H_n$ have
at lease three copies of $\Z$.

Note that  the homomorphism
$$
\gm_{n,n+1}:~ H_n=\Z^{p_n}\lr H_{n+1}=\Z^{p_{n+1}}
$$ is given by a $p_{n+1}\times p_n$ matrix $\cc=(c_{ij})$ of nonnegative integers,
where $i=1,2,...,p_{n+1};~ j=1,2,...,p_n$ and $c_{ij}\in
)_+:=\{0,1,2,...\}.$ For $M>0$, if all $c_{ij}\geq M$, then we will
say $\gm_{n,n+1}$ is at least $M$-large or has multiplicity at least
$M$.  Note that since $H$ is a simple group, passing to {subsequences}, we
{may assume  that} at each step $\gm_{n,n+1}$ is at least $M_n$-large for
arbitrary choice of $M_n$ depending on our construction up to step
$n$.

\end{NN}

\begin{NN}\label{range 0.5} As in \ref{range  0.3}  and \ref{range  0.4}, we have $G\subseteq H$ with
$G_+=H_+\cap G$ and both $G$ and $H$ share the same order unit $u\in G\subseteq H$.  As in \ref{range 0.4}, write $H$ as inductive limit of
$H_n$---finite direct sum of ordered groups $(\Z, \Z_+)$.  Let $G_n={\gm_{n, \infty}^{-1}(\gm_{n, \infty}(H_n)\cap G)},$
%H_n\cap\gm_{n,\infty}^{-1}(G)$,
 where $\gm_{n,\infty}:~ H_n\to H$ is induced map to the limit.  We may
assume $u\in G_n\subseteq  H_n$ for each $n$ and $(G_n)_+=(H_n)_+\cap G_n$.
Since $G$ is a relatively divisible subgroup of $H$ {and $H$ is torsion free}, the quotient
$H_n/G_n$ is torsion free group and therefor a direct sum of copies
of $\Z$, denoted by $\Z^{l_n}$. Then we have the following commutative diagram
\begin{displaymath}
    \xymatrix{G_1 \ar[r]^{\gm_{_{12}}|_{_{G_1}}} \ar@{^{(}->}[d] & G_2 \ar[r]\ar@{^{(}->}[d]&\cd \ar[r]&G \ar@{^{(}->}[d]\\
         H_1 \ar[r]^{\gm_{_{12}}} \ar[d] & H_2 \ar[r]\ar[d]&\cd \ar[r]&H \ar[d]\\
         H_1/G_1 \ar[r]^{\td\gm_{_{12}}}  & H_2/G_2 \ar[r]&\cd \ar[r]&H/G }
\end{displaymath}
Let $H_n=\left(\Z^{p_n}, \Z^{p_n}_+, u_n\right)$, where
$u_n=([n,1],[n,2], ..., [n,p_n])\in (\Z_+\setminus \{0\})^{p_n}$.
Then $H_n$ can be realized as the $K_0$ group of
$F_n=\bigoplus_{i=1}^{p_n}M_{[n,i]}(\C)$---that is,
$$
(K_0(F_n),K_0(F_n)_+,{\e}_{F_n})=(H_n,(H_n)_+,u_n).
$$

If there are infinitely many $n$ such that the inclusions
%$G_n\hookrightarrow H_n$
$G_n\rightarrow H_n$
are isomorphisms, then, by passing to subsequence, we have that
%$G\hookrightarrow H$
$G_n\rightarrow H_n$
is also an isomorphism which contradict with $G\subsetneqq H$ in\ref{range 0.3}. Hence, without lose of generality, we can assume that for all $n$, $G_n\subsetneqq H_n$, and therefore $H_n/G_n\not= 0$.

To construct the C*-algebra with $K_0$ being $(G_n,(G_n)_+,u_n)$, we
consider the map $\pi:~ H_n\to H_n/G_n$ being identified as a map
(still denoted by)
$$\pi:~ \Z^{p_n}\lr \Z^{l_n}.$$
as in \cite{point-line}. We emphasis that $l_n>0$ for all $n$.
Such a map can be realized as difference of two
maps
$$\bb_0,~ \bb_1:~ \Z^{p_n}\lr \Z^{l_n}$$
corresponding to two $l_n\times p_n$ matrices of strictly positive
integers $\bb_0=(b_{0,ij})$ and $\bb_1=(b_{1,ij})$. That is,
$$
\qq\qq\pi\left(
     \begin{array}{c}
       t_1 \\
       t_2 \\
       \vdots \\
       t_p \\
     \end{array}
   \right)
   =
   \big(\bb_1-\bb_0\big)
   \left(
     \begin{array}{c}
       t_1 \\
       t_2 \\
       \vdots \\
       t_p \\
     \end{array}
   \right)\in \Z^{l_n},\qq
\mbox{for any} \quad \left(
     \begin{array}{c}
       t_1 \\
       t_2 \\
       \vdots \\
       t_p \\
     \end{array}
   \right)\in \Z^{p_n}.
$$
Note that $u_n=([n,1],[n,2], \cd, [n,p_n])\in G_n$ and hence
$\pi(u_n)=0$.  Consequently,
$$
 \bb_1
\left(
     \begin{array}{c}
       {[n,1]} \\
      {[n,2]} \\
       \vdots \\
       {[n,p_n]} \\
     \end{array}
   \right)
   =
\bb_0 \left(
     \begin{array}{c}
       [n,1] \\
       {[n,2]} \\
       \vdots \\
       {[n,p_n]} \\
     \end{array}
   \right)
\triangleq \left(
     \begin{array}{c}
       \{n,1\} \\
       \{n,2\}\\
       \vdots \\
       \{n,p_n\} \\
     \end{array}
   \right),
$$
i.e., denote that
$$
\{n,i\}\triangleq\sum_{j=1}^{p_n} b_{1,ij}[n,j]=\sum_{j=1}^{p_n}
b_{0,ij}[n,j].
$$

Let $E_n=\oplus_{i=1}^{l_n} M_{\{n,i\}}(\C)$.  We can choose any two
homomorphisms $\bt_0,\bt_1:~ F_n\to E_n$ such that $(\bt_0)_{*0}=\bb_0$
and $(\bt_1)_{*0}=\bb_1$.  Then define
$$C_n=\Big\{(f,a)\in C([0,1],E_n)\oplus F_n;~ f(0)=\bt_0(a), f(1)=\bt_1(a)\Big\},$$
which is $A(F_n, E_n, \bt_0,\bt_1)$ as in the definition of \ref{DfC1}.
Using the following six-term exact sequence
\begin{displaymath}
\xymatrix{
K_0(C_0\big((0,1),E_n\big) \ar[r]&K_0(C_n) \ar[r] &K_0(F_n)\ar[d]\\
K_1(F_n) \ar[u] &K_1(C_n)\ar[l]&K_1(C_0\big((0,1),E_n\big),\ar[l]}~
\end{displaymath}
and the fact that $\pi$ is surjective,
we have $K_1(C_n)=0$ and
$$
\Big(K_0(C_n),K_0(C_n)_+, \e_{C_n} \Big)=\Big({G_n},(G_n)_+,u_n\Big).
$$
Note that in the above diagram the map $K_0(F_n)=\Z^{p_n} \to
K_1(C_0\big((0,1),E_n\big))=\Z^{l_n}$ is given by $( \bb_1-\bb_0)\in
M_{l_n\times p_n}(\Z)$, which is surjective, as quotient map $H_n
\big(=\Z^{p_n}\big)\to H_n/G_n \big(=\Z^{l_n}\big)$.

As observed in \cite{point-line}, in the construction of $C_n$, we have the
freedom to choose the pair of the $K_0$ map $(\bt_0)_{*0}=\bb_0$ and
$(\bt_1)_{*0}=\bb_1$ as long as the difference is the same  map
$\pi:~ H_n(=\Z^{p_n})\lr H_n/G_n(=\Z^{l_n})$. For example, if
$(m_{ij})\in M_{l_n\times p_n}\big(\Z_+\setminus \{0\}\big)$ is any
$l_n\times p_n$ matrix of positive integer, then we can replace
$b_{0,ij}$ by $b_{0,ij}+m_{ij}$ and, at the same time, replace
$b_{1,ij}$ by $b_{1,ij}+m_{ij}$. That is, we can assume that each
entry of $\bb_0$ (and of $\bb_1$) is larger than any fixed integer
$M$ which depends on $C_{n-1}$ and $\psi_{n-1,n}: F_{n-1}\to F_n$.
Also, we can make all the entries of one column (say, the third
column) of both $\bb_0$ and $\bb_1$  much larger than all the
entries of another column
(say, the second), by choosing  $m_{i3}>>m_{j2}$ for all $i,~ j$.  \\

\end{NN}

\begin{NN}\label{range 0.5a}
For later use, we will also deal with the case that
$G_n=\Z^{\bullet}\oplus G_n'$ and $H_n=\Z^{\bullet}\oplus H_n'$,
with the inclusion map being identity for the first $\bullet$ copies
of $\Z$. In this case, the quotient map $H_n(=\Z^{p_n}) \to
H_n/G_n(=\Z^{l_n})$ given by the matrix $\bb_1-\bb_0$ maps the first
$\bullet$ copies of $\Z$ to zero. For this case, it will be much
more convenient to assume that the first $\bullet$-columns of  both
matrices $\bb_0$ and $\bb_1$ are zero and each entry of the last
$(p_n-\bullet)$-column of them are larger than any previously  given integer
$M$. Now we have that the entries of the matrices $\bb_0,\bb_1$ are strictly positive integers except the ones in the first $\bullet$-columns which are zero.

Consider the following diagram

\begin{displaymath}
    \xymatrix{\Z^{p_1^0}\oplus G'_1 \ar[r]^{\gm_{_{12}}|_{_{G_1}}} \ar@{^{(}->}[d] & \Z^{p_2^0}\oplus G'_2 \ar[r]\ar@{^{(}->}[d]&\cd \ar[r]&G \ar@{^{(}->}[d]\\
       \Z^{p_1^0}\oplus  H'_1 \ar[r]^{\gm_{_{12}}} \ar[d] & \Z^{p_2^0}\oplus H'_2 \ar[r]\ar[d]&\cd \ar[r]&H \ar[d]\\
         H_1/G_1 \ar[r]^{\td\gm_{_{12}}}  & H_2/G_2 \ar[r]&\cd \ar[r]&H/G. }
\end{displaymath}
That is $G_n=\Z^{p_n^0}\oplus G_n'$ and $H_n=\Z^{p_n^0}\oplus H_n'$,
with the inclusion map being identity for the first $p_n^0$ copies
of $\Z$. Write $p_n^1=p_n-p_n^0$. Then we assume that the entries of the matrices $\bb_0,\bb_1$ are strictly positive integers for last $p_n^1$-columns and are zeros for  the first $p_n^0$-columns. Note that since $l_n>0$ (see \ref{range 0.5}), we have $p_n^1>0$.

%The case {discussed in}  \ref{range 0.5} is a special case {of the current case with}   $p_n^0=0$ for all $n$.
{ {{The inductive limits $H=\varinjlim(H_n,\gm_{n,n+1})$ and $G=\varinjlim(G_n,\gm_{n,n+1}|_{G_n})$ (with $G_n\subset H_n$)  constructed in \ref{range 0.5} are in fact special cases of the current construction
when we assume that $p_n^0=0.$} }} One {notices} that, for the case $G_n=\Z^{p_n^0}\oplus G_n'$ and $H_n=\Z^{p_n^0}\oplus H_n'$,
with the inclusion map being identity for the first $p_n^0$ copies
of $\Z$, one can still use the construction of \ref{range 0.5} to make all the entries of $\bb_0$ and $\bb_1$ (not only the entries of last $p_n^1$) to be strictly positive (of course with first $p_n^0$ column of matrices $\bb_0$ and $\bb_1$ equal to each others).
%{\bf {Omit this: That is, it seems one only needs the construction for the special case of \ref{range 0.5}.} }
{However,}
%But
for the case that $p_n^0\not= 0$ for all $n$, we will get  better decomposition as we will discuss in the next section. In fact, for the case that $p_n^0\not= 0$ for all $n$, the construction is much simpler than the case that $p_n^0=0$ for all $n$.

If we write $F_n$ (in \ref{range 0.5})) as
$\bigoplus_{i=1}^{p_n^0}M_{[n,i]}(\C)\oplus F_n'$, where
$F_n'=\bigoplus_{i=p_n^0 +1}^{p_n}M_{[n,i]}(\C),$ then the maps
$\bt_0$ and $\bt_1$ are zero on the part
$\bigoplus_{i=1}^{p_n^0}M_{[n,i]}(\C)$. Moreover, the algebra
$C_n=\Big\{(f,a)\in C([0,1],E_n)\oplus F_n;~ f(0)=\bt_0(a),
f(1)=\bt_1(a)\Big\}$ in \ref{range 0.5} can be written as
$\bigoplus_{i=1}^{p_n^0}M_{[n,i]}(\C) \oplus C'$, where
$C'=\Big\{(f,a)\in C([0,1],E_n)\oplus F_n';~ f(0)=\bt_0(a),
f(1)=\bt_1(a)\Big\}$, which is $A(F_n', E_n,
\bt_0|_{F_n'},\bt_1|_{F_n'})$ as in the notation of \ref{DfC1}.

\end{NN}

\begin{NN}\label{conditions} Let us emphasis that once
$$(H_n,(H_n)_+,u_n)=\big(\Z^{p_n}, (\Z_+)^{p_n},
([n,1],[n,2], ..., [n,p_n])\big),~~~~~\bb_0,\bb_1: \Z^{p_n} \to
\Z^{l_n},$$
are fixed, then the algebras $F_n=\bigoplus_{i=1}^{p_n} M_{[n,i]}(\C)$, $E_n=\bigoplus_{i=1}^{l_n} M_{\{n,i\}}(\C),$ where \\
$\{n,i\}\triangleq\sum_{j=1}^{p_n} b_{1,ij}[n,j]=\sum_{j=1}^{p_n}
b_{0,ij}[n,j]$, and $C_n=A(F_n,E_n,\bt_0,\bt_1)$, with $K_0\bt_i=\bb_i$ ($i=0,1$), are determined  up to isomorphism.

To construct the inductive limit, we not only need to construct $C_n$'s (later on $A_n$'s for the  general case) but also
need to construct $\phi_{n,n+1}$   which realizes the corresponding K-theory map---that is $(\phi_{n,n+1})_{*0}=\gm_{_{n,n+1}}|_{G_n}$. In additional, we need to make the limit algebras to have the desired tracial state space.
%For this to be done,
{In order to do all these,} we need some extra conditions for $\bb_0, \bb_1$ for $C_{n+1}$ (or for $A_{n+1}$)  depending on $C_n$ (or $A_n$).
We will divide the construction into several steps with gradually stronger conditions for $\bb_0, \bb_1$ (for $C_{n+1}$)---of course depending on $C_n$ and the map $\gm_{n,n+1}: H_n\to H_{n+1}$,  to guarantee the construction can go through.

Let $G_n \subset H_n=\Z^{p_n}$,  $G_{n+1} \subset H_{n+1}=\Z^{p_{n+1}}$, and $\gm_{_{n,n+1}}:~ H_n\to H_{n+1}$, with $\gm_{_{n,n+1}} (G_n) \subset G_{n+1}$ be as in \ref{range 0.4} and \ref{range 0.5} (also see \ref{range 0.5a}). Then $\gm_{_{n,n+1}}$ induces a map $\td\gm_{_{n,n+1}}: H_n/G_n(=\Z^{l_n}) \to H_{n+1}/G_{n+1}= (\Z^{l_{n+1}})$.
Let $\gm_{_{n,n+1}}:~H_n\left(=\Z^{p_n}\right)\to
H_{n+1}\left(=\Z^{p_{n+1}}\right)$ be given by a matrix
$\cc=(c_{ij})\in M_{_{p_{n+1}\times p_n}}(\Z_+\setminus \{0\})$ and
$\td\gm_{_{n,n+1}}:~\Z^{l_n} \to \Z^{l_{n+1}}$ (as a map from
$H_n/G_n \to H_{n+1}/G_{n+1}$) be given by matrix $\dd=(d_{ij})$. Let $\pi_{n+1}: H_{n+1} (=\Z^{p_{n+1}}) \to H_{n+1}/G_{n+1}(=\Z^{l_{n+1}})$ be the quotient map.

Let us use $\bb'_0, \bb'_1:\Z^{p_{n+1}} \to \Z^{l_{n+1}}$ to denote the maps  
%defined during 
%
required for 
the construction of $C_{n+1}$, and reserve $\bb_0, \bb_1$ for $C_n$.  Of course $\pi_{n+1}=\bb'_1-\bb'_0$.

We will prove that if $\bb'_0, \bb'_1$ satisfy:
\beq\label{13spd-1}
\spd \qq\qq\qq\qq \td b_{0,ji},~\td b_{1,ji}~ > ~
\sum_{k=1}^{l_n}\big(|d_{jk}|+2\big)\max(b_{0,ki},~b_{1,ki})
\eneq
for all
$i\in\{1,2,...,p_n\}$ and for all  $j\in\{1,2,...,l_{n+1}\},$
where $\td \bb_0=\bb'_0\cdot \cc=(\td b_{0,ji})$ and $\td
\bb_1=\bb'_1\cdot \cc=(\td b_{1,ji})$, then one can construct the homomorphism $\phi_{n,n+1}: C_n \to C_{n+1}$ to realized the desired K-theory map (see \ref{range 0.6} below). If $\bb'_0, \bb'_1$ satisfy stronger condition:
 \beq
\spdd\qq\qq\qq\qq \td b_{0,ji}, \td
 b_{1,ji}
 > 2^{2n}\!\left(\sum_{k=1}^{l_n}
 (|d_{jk}|+2)\{n,k\}\!\!\right)
 \eneq
 for all
$i\in\{1,2,...,p_n\}$ and for all  $j\in\{1,2,...,l_{n+1}\}$, then we can prove the limit algebra constructed has the desired tracial state space (see \ref{range 0.8}---\ref{range 0.12}). (It follows from that $\{n,k\}\triangleq\sum_{i=1}^{p_n} b_{1,ki}[n,i]=\sum_{i=1}^{p_n}b_{0,ki}[n,i],$ we have that $\spdd$ is stronger than $\spd$.)

If we want to make the limit algebra to be simple, then we need to modify the connecting homomorphisms $\phi_{n,n+1}$ (but not the algebras $C_n$ and $C_{n+1}$), and we need even stronger condition for $\bb_0'=(b_{0,ij}')$ and $\bb_1'=(b_{1,ij}')$ below:
%{{as:}}
\beq
\td b_{0,il}> 2^{2n}\left(\sum_{k=1}^{l_{n}}(|d_{ik}|+2)\cdot \{n,k\}+
((n-1)\sum_{k=1}^{l_n}b_{0,kl})\cdot b'_{0,i1}
\right)\\
\hspace{-2.6in}\andeqn
%\qq\qq\qq
\td b_{1,il}>
2^{2n}\left(\sum_{k=1}^{l_{n}}(|d_{ik}|+2)\cdot \{n,k\}+
((n-1)\sum_{k=1}^{l_n}b_{0,kl})\cdot b'_{1,i1} \right),
\eneq
 for
each $l\in \{1,2,..., p_n\}$ and $i\in \{1,2,..., l_{n+1}\}$. But we will not do the modification for this special case since it only covers the case that the algebras have trivial $K_1$ group and torsion free $K_0$ group. Later on, we will deal with general case involving general $K_1$ group and $K_0$-group; % and general (not torsion free) $K_0$ group,
the conditions will be $\spddd_1$ and $\spddd_2$ specified in \ref{condition2} below.

Let us emphasis that when we modify inductive limit system to make it simple, we never change the algebras $C_n$ (or $A_n$ in the general case), what will be changed are the connecting homomorphisms.

%{\bf{Is it obvious that the latter conditions are stronger than previous ones?---Yes}}

\end{NN}

{\it Let $A$ and $B$ be two \CA s,  $\phi: A\to B$ be a \hm\, and $\pi\in RF(B).$
{For the rest of} this  section, we will use $\phi|_\pi$ for the composition $\pi\circ \phi,$ in particular,
in the following statement and its proof. This notation is consistent with  {\rm \ref{homrestr}}.}

{In the following lemma we will give the construction of  $\phi_{n,n+1}$  if $\bb'_0$ and $\bb'_1$ satisfy the condition $\spd$ in \ref{conditions}---of course, the condition depends on the previous  step.
So this lemma {provides the}
%that  how the
$n+1$-step of the  construction. Again, we first have $G,$ then obtains $H,$
$H_n$ and $G_n$ as constructed in \ref{range 0.3} and \ref{range 0.4}. }

\begin{lem}\label{range 0.6}
Let
\beq\nonumber
(H_n,(H_n)_+,u_n)=\big(\Z^{p_n}, (\Z_+)^{p_n},
([n,1],[n,2], ... , [n,p_n])\big),\\
F_n=\bigoplus_{i=1}^{p_n} M_{[n,i]}(\C),\,\,\, \bb_0,\bb_1: \Z^{p_n} \to
\Z^{l_n},\,\,\,
E_n=\bigoplus_{i=1}^{l_n} M_{\{n,i\}}(\C),~ \bt_0,\bt_1:
F_n \to E_n
\eneq
 with $(\bt_0)_{*0}=\bb_0, (\bt_1)_{*0}=\bb_1$, and
$C_n=A(F_n,E_n,\bt_0,\bt_1)$ with $K_0(C_n)=G_n$ be as in
{\rm \ref{range 0.5}} or more generally as in {\rm \ref{range 0.5a}}.
  Let
$$
(H_{n+1},H_{n+1}^+,u_{n+1})=\big(\Z^{p_{n+1}}, (\Z_+)^{p_{n+1}},
([{n+1},1],[{n+1},2], ..., [{n+1},p_{n+1}])\big),
$$
let  $\gm_{n,n+1}:~
H_n\to H_{n+1}$ be { an} ordered homomorphism with
$\phi_{n,n+1}(u_n)=u_{n+1}$ {\rm (}as in {\rm \ref{range  0.5}} or {\rm \ref{range  0.5a}} {\rm )}, let
$G_{n+1}\sbs H_{n+1}$ be a subgroup containing $u_{n+1}$ {\rm (}as in
{\rm \ref{range 0.4}}{\rm )} and satisfying $\gm_{n,n+1}(G_n)\sbs G_{n+1}$.  Let
$\pi_{n+1}: H_{n+1}(=\Z^{p_{n+1}}) \to H_{n+1}/G_{n+1}(=\Z^{l_{n+1}})$ denote the quotient map.

Suppose that the matrices $\bb'_0=(b'_{0,ij}),~\bb'_1=(b'_{1,ij}): \Z^{p_{n+1}} \to \Z^{l_{n+1}}$  satisfy $\bb'_1-\bb'_0=\pi_{n+1}$ and satisfy the condition $\spd$ in {\rm \ref{conditions}}. (As a convention, we assume that the entries of first $p_{n+1}^0$ columns of the matrices $\bb'_0$ and $\bb'_1$ are zeros and the entries of the last $p_{n+1}^1=p_{n+1}-p_{n+1}^0$ columns are strictly positive. Note that $p_{n+1}^0$ might be zero as in the special case {\rm \ref{range  0.5}}.) %which is the case  {\rm \ref{range  0.5}}---a  special case, )
%{\bf {This sentence is confusion: In section 3.10, all entries are strictly positive, in section 13.11, on the other, it is the situation described here. However,  this is contradicted to the first part of the statement, which seems to suggest
%that either way is fine.}}

Put
$F_{n+1}=\bigoplus_{i=1}^{p_{n+1}} M_{[{n+1},i]}(\C)$ and
$E_{n+1}=\bigoplus_{i=1}^{l_{n+1}} M_{\{{n+1},i\}}(\C)$ with $$\{n+1,i\}=\sum_{j=1}^{p_{n+1}} b'_{1,ij}[n+1,j]=\sum_{j=1}^{p_{n+1}}
b'_{0,ij}[n+1,j],$$
%unital
%homomorphisms
and pick unital homomorphisms $\bt_0',\bt_1': F_{n+1}\to E_{n+1}$ with $$(\bt_0')_{*0}=\bb'_0\quad \textrm{and}\quad (\bt_1')_{*0}=\bb'_1.$$
Set
$$C_{n+1}=A(F_{n+1},E_{n+1},\bt_0',\bt_1').$$

Then there is a  \hm~ $\phi_{n, n+1}:~
C_n\to C_{n+1}$ satisfying the following conditions:
\begin{enumerate}
\item[\rm{(1)}] $K_0(C_{n+1})=G_{n+1}$ as a scaled ordered group (as already verified in \ref{range 0.5}).
\item[\rm{(2)}] $(\phi_{n, n+1})_{*0}:~ K_0(C_n)=G_n \to K_0(C_{n+1})= G_{n+1}$
satisfies $(\phi_{n, n+1})_{*0}=\gm_{n, n+1}|_{G_n}$.
\item[\rm{(3)}] $\phi_{n, n+1}(C_0\big((0,1), E_n\big))\sbs C_0\big((0,1),
E_{n+1}\big)$.
\item[\rm{(4)}] Let $\td\phi_{n, n+1}:~ F_n\to F_{n+1}$ be the quotient map
induced by $\phi_{n, n+1}$ (note from (3), we know that this
quotient map exists), then
 $(\td\phi_{n, n+1})_{*0}=\gm_{n, n+1}:~ K_0(F_n)=H_n \to K_0(F_{n+1})=H_{n+1}$.
\item[\rm{(5)}] For each $y\in Sp(C_{n+1})$, $Sp(F_{n})\subset Sp(\phi_{n,n+1}|_y)$.
\item[\rm{(6)}] For each $j_0\in \{ 1, 2, ..., l_{n+1}\}, ~i_0\in \{ 1, 2, ... ,
l_{n}\}$, one of the following holds:
\begin{enumerate}
\item[{\rm (i)}]  for each $t\in (0,1)_{j_0}\sbs
Sp(I_{n+1})=\bigcup_{j=1}^{l_{n+1}}(0,1)_j\sbs Sp(C_{n+1})$,
$Sp(\phi_{n, n+1}|_t)\cap (0,1)_{i_0}$ contains $t\in (0,1)_{i_0}\sbs
Sp(C_n)$; or
\item[{\rm (ii)}] for each $t\in (0,1)_{j_0}\sbs
Sp(I_{n+1})=\bigcup_{j=1}^{l_{n+1}}(0,1)_j\sbs Sp(C_{n+1})$,
$Sp(\phi_{n, n+1}|_t)\cap (0,1)_{i_0}$ contains $1-t\in
(0,1)_{i_0}\sbs Sp(C_n)$.
\end{enumerate}

\end{enumerate}

 Consequently (as  $l_{n+1}>0$, that is $E_{n+1}\not=0$ or
 there is at least one interval in $Sp(C_{n+1})$),  the following is true:
 if $X\sbs Sp(C_{n+1})$ is $\dt$-dense,
then $\bigcup_{x\in X}Sp(\phi_{n, n+1}|_x)$ is $\dt$-dense in $Sp(C_n)$ (see \ref{density}).
%{\bf {{ You need to specify  $\dt$-density here too. One suggestion:  discuss this at earlier stage.
%I  could try. May be you want to do your way.-----Be aware:  $SP(\phi)$ is Not really %even a subset of
%$Sp(A_n)$}}}
(In general,    $\phi_{n, n+1}$ is
injective.)

%{Moreover, $\bt_0'$ and $\bt_1'$ have the property described in \ref{range 0.5} or
%\ref{range 0.5a} (but one of them throughout the construction).}
\end{lem}

%{\bf {Do we need to state something about $M$ large?, I don't think so}}

\begin{rem}\label{range 0.7} Let $I_n=C_0\big((0,1),E_n\big)$  and
$I_{n+1}=C_0\big((0,1),E_{n+1}\big)$.  If $\phi_{n, n+1}:~ C_n\to
C_{n+1}$
 is as described in \ref{range  0.6}, then we have the following exact sequences:
\begin{displaymath}
\xymatrix{0\ar[r] &
K_0(C_n)\ar[r]^{~}\ar@{->}[d]_{\gm_{_{n,n+1}}|_{_{_{G_n}}}} &
K_0(C_n/I_n) \ar@{->}[d]_{\gm_{_{n,n+1}}} \ar[r]^{~} & K_1(I_n)
\ar@{->}[d]_{\td\gm_{_{n,n+1}}} \ar[r]&0
\\
0\ar[r] & K_0(C_{n+1})\ar[r]^{~}& K_0(C_{n+1}/I_{n+1})
 \ar[r]^{~} & K_1(I_{n+1}) \ar[r] & 0}
\end{displaymath}
where $K_0(C_n)$ and $K_0(C_{n+1})$ are identified with $G_n$ and
$G_{n+1}$, $K_0(C_n/I_n)(=K_0(F_n))$ and
$K_0(C_{n+1}/I_{n+1})(=K_0(F_{n+1}))$ are identified with $H_n$ and
$H_{n+1}$, $K_1(I_n)$ is identified with $H_n/G_n$, and
$K_1(I_{n+1})$ is identified with $H_{n+1}/G_{n+1}$. Moreover,
$\td\gm_{_{n,n+1}}$ is induced by $\gm_{_{n,n+1}}:~ H_n\to H_{n+1}$.

Consider the matrix
$\cc=(c_{ij})\in M_{_{p_{n+1}\times p_n}}\Z_+\setminus \{0\})$ which is induced by the map
$\gm_{_{n,n+1}}: H_n\left(=\Z^{p_n}\right) \to H_{n+1}\left(=\Z^{p_{n+1}}\right)$
and consider the matrix $\dd=(d_{ij})$ which is induced by the map
$\td\gm_{_{n,n+1}}:~\Z^{l_n} \to \Z^{l_{n+1}}$ (as a map from
$H_n/G_n \to H_{n+1}/G_{n+1}$).
%{\color{Green} (It seems that the matrix $\dd$ is not used in this remark)} {\red{Answer :Yes,its entries $d_{jk}$ used below}}
Note that $\pi_n:~ H_n(=\Z^{p_n})\to  H_n/G_n(=\Z^{l_n})$ is given
by $\bb_1-\bb_0\in M_{l_n\times p_n}(\Z).$
{Here,  in the situation of \ref{range 0.5a},
 we assume the first
$p_n^0$ columns of both $\bb_0$ and $\bb_1$ are zero
and last $p_n^1$ columns are  strictly positive.}
 %%as in the
%case of \ref{range 0.5a}.
Let $\pi_{n+1}:~ H_{n+1}(=\Z^{p_{n+1}})
\to H_{n+1}/G_{n+1}(=\Z^{l_{n+1}})$ be the quotient map.
%If the
%algebra $A_{n+1}=A(F_{n+1},E_{n+1},\bt'_0,\bt'_1)$ and
%$\phi_{n,n+1}:~ A_n\to A_{n+1}$ are as in the lemma, and if
%$$
%(\bt'_0)_{*0}=\bb'_0\in M_{l_{n+1}\times p_{n+1}}(\Z^+)\andeqn
%(\bt'_1)_{*0}=\bb'_1\in M_{l_{n+1}\times p_{n+1}}(\Z^+),
%$$
%then $\pi_{n+1}=\bb'_1-\bb'_0$.
%As
%{In the case of
% \ref{range 0.5a},  we can write}
Write $G_{n+1}=\Z^{p_{n+1}^0}\oplus G_{n+1}'$,
$H_{n+1}=\Z^{p_{n+1}^0}\oplus H_{n+1}'$, then we can choose  both $\bb_0'$ and $\bb_1'$, with the first $p_{n+1}^0$ columns being zero and the last $p_{n+1}^1=p_{n+1}-p_{n+1}^0$ columns being strictly positive, such that $\pi_{n+1}=\bb'_1-\bb'_0$ and
 %{and insist that} first
%$\bullet\bullet$-columns of both $\bb_0'$ and $\bb_1'$ are zero.
%Denote that $\td \bb_1=\bb'_1\cdot \cc$ and $\td \bb_0=\bb'_0\cdot
%\cc$. In the proof of \ref{range 0.6}, we will  define
%$$
%E_{n+1}=\bigoplus _{i=1}^{l_n} M_{\{n+1,i\}}(\C),\qq \mbox{where
%integer }~\{n+1,i\}=\sum_{j=1}^{p_{n+1}}b'_{0,ij}[n+1,j]=
%\sum_{j=1}^{p_{n+1}}b'_{1,ij}[n+1,j],
%$$
% (this is true since $\bb'_1-\bb'_0=\pi_{n+1}$ which maps
% $([n+1,1],[n+1,2],\cd,[n+1,p_{n+1}])\in G_{n+1}$ into zero.)
%and  define homomorphisms $\bt'_0,\bt'_1:~ F_{n+1}\to E_{n+1}$ with
% $(\bt'_0)_{*0}=\bb'_0,~ (\bt'_1)_{*0}=\bb'_1$ and homomorphism
% $\phi_{_{n,n+1}}:~ A_n\to A_{n+1}=A(F_{n+1},E_{n+1},\bt'_0,\bt'_1)$
%as desired in the theorem. We will prove that such construction can
% be carried  out provided that
%two matrices of non negative integers $\bb'_1=(b'_{1,ij})$ and
%$\bb'_0=(b'_{0,ij})$ satisfy that $\pi_{n+1}=\bb'_1-\bb'_0$ and the
%following condition:
\beq
 \spd \qq\qq \td b_{0,ji},~\td b_{1,ji}~ > ~
\sum_{k=1}^{l_n}\big(|d_{jk}|+2\big)\max(b_{0,ki},~b_{1,ki})
\eneq
for all
$i\in\{1,2,...,p_n\}$ and for all  $j\in\{1,2,...,l_{n+1}\},$
where $\td \bb_0=\bb'_0\cdot \cc=(\td b_{0,ji})$ and $\td
\bb_1=\bb'_1\cdot \cc=(\td b_{1,ji})$.

%  If $l_{n+1}=0$
%(that is $G_{n+1}=\Z^{p_{n+1}^0}$ and
%$H_{n+1}=\Z^{p_{n+1}^0}$), then  both matrices $\td \bb_0$ and
%$\td \bb_1$,  as $l_{n+1}\times p_n$ matrix,  are  empty matrices---the matrices with no entry,  and hence the condition $\spd$ holds tautologically.

Note that  $l_{n+1}>0$ and  $p_{n+1}^1>0$. So we can make $\spd$ hold by only increasing the last
$p_{n+1}^1$ columns of the the matrices $\bb_0'$ and
$\bb_1'$---that is the first $p_{n+1}^0$ column of the matrices
are still kept to be  zero, since all the entries in $\cc$ are strictly positive.  Note that, even though the first $p_{n+1}^0$ column of $\bb_0'$ and
$\bb_1'$ (as $l_{n+1}\times p_{n+1}$ matrices) are zero, but all entries of $\td \bb_0$ and $\td \bb_1$ (as $l_{n+1}\times p_n$ matrix) has been
made strictly positive.
%Also note that if  $l_n=0$
%(that is $G_n=\Z^{p_n^0}$ and $H_n=\Z^{p_n^0}$), then  both
%matrices $\bb_0$ and $\bb_1$ are  empty matrices--the matrix with no entry, and the right %hand side is then the sum over empty set which is considered to be zero.
Again note that  the case of \ref{range 0.5} is a special case of \ref{range 0.5a} for $p_{n+1}^0=0$, so one does not need to deal with this case separately.
%{\bf {While I understand what this means, it certainly could cause some confusion.}}
% then we can define
%$$
%E_{n+1}=\bigoplus _{i=1}^{l_n} M_{\{n+1,i\}}(\C),\qq \mbox{where
%integer }~\{n+1,i\}=\sum_{j=1}^{p_{n+1}}b'_{0,ij}[n+1,j]=
%\sum_{j=1}^{p_{n+1}}b'_{1,ij}[n+1,j],
%$$
% (this is true since $\bb'_1-\bb'_0=\pi_{n+1}$ which takes
%$([n+1,1],[n+1,2],\cd,[n+1,p_{n+1}])\in G_{n+1}$ to zero.)
%and we can define homomorphisms $\bt'_0,\bt'_1:~ F_{n+1}\to E_{n+1}$
%with
%$K_0\bt'_0=\bb'_0,~ K_0\bt'_1=\bb'_1$ and homomorphism
%$\phi_{_{n,n+1}}:~ A_n\to A_{n+1}=A(F_{n+1},E_{n+1},\bt'_0,\bt'_1)$
%as desired in the theorem.

 \end{rem}

 %{\bf Proof of \ref{range 0.6}}
  \begin{proof}[Proof of \ref{range 0.6}]
 Suppose that
 $\bb_0',\bb_1'$  satisfy $\bb'_1-\bb'_0=\pi_{n+1}$ and
 (\ref{13spd-1}) (and the first $p_{n+1}^0$ columns to be zero).  Now let $E_{n+1},~\bt'_0,\bt'_1:~ F_{n+1}\to
 E_{n+1}$, and $C_{n+1}=A(F_{n+1},E_{n+1},\bt'_0,\bt'_1)$ be as
 constructed in \ref{range 0.5}.
 %{\bf {This is confusing since it suggests
 %that \ref{range 0.7} had constructed $E_{n+1}, \bt_0' $ and so on, and only $\phi_{n, %n+1}$
 %are the ones in questions. However, it seems that there was $\phi_{n, n+1}$ in \ref{range 0.7} too.
 %You also told me that $A_{n+1}$ needs to be constructed.   It is also unclear from the second sentence of \ref{range 0.7}
 %that how $A_{n+1}$ was constructed from the $n$-step, or the issue has been discussed in \ref{range 0.7}.
 %Or perhaps, the construction of $A_{n+1}$ is not issue--just take one from \ref{range 0.5}?}}
  We will define $\phi_{n,n+1}:~C_n\to C_{n+1}$
to satisfy Conditions (2)--(6)  of \ref{range  0.6} (Condition (1) is a property for $C_{n+1}$ which is verified in \ref{range 0.5}).

As usual, let us use $F_n^i$ (or $E_n^i$) to denote the $i$-th block of $F_n$ (or $E_n$).

There exists a unital \hm\,
$\td\phi_{n, n+1}:~ F_n\to F_{n+1}$  such that
$$(\td\phi_{n,n+1})_{*0}=
\gm_{n, n+1}:~ K_0(F_n)=H_n \to K_0(F_{n+1})=H_{n+1}, $$where
$\gm_{n,n+1}$ is defined in \ref{range 0.4}.
%(Such
%homomorphism exists  by [Elliott AF Classification].)
Note that
$Sp(C_{n+1})=\coprod_{j=1}^{l_{n+1}}(0,1)_j\bigcup Sp(F_{n+1})$ (see
\S 3). Write $Sp(F_{n+1})=(\tht_1',\tht_2',...,\tht'_{p_{n+1}})$ {and}
 $Sp(F_{n})=(\tht_1,\tht_2, ... ,\tht_{p_{n}})$. To define
$\phi_{n,n+1}:~C_n\to C_{n+1}$,  we need to specify {$\phi_{_{n,n+1}}|_y=e_y\circ \phi_{_{n,n+1}}$ below,} where for each
$y\in Sp(C_{n+1})$,
\begin{displaymath}
\xymatrix{\phi_{_{n,n+1}}|_y:~C_n\ar[r]&  C_{n+1}~ \ar[r]^{e_{_y}} &
C_{n+1}|_y}~,
\end{displaymath}
{and $e_y$ is the point-evaluation at $y.$}
%{{\bf  I have deleted (in \%) several lines}}
%That is, we need to specify $\phi_{_{n,n+1}}(f)(y)$ for each $f\in
%C_n$ and each $y\in Sp(C_{n+1})$.
%As we know,
%$Sp(C_{n+1})=\coprod_{j=1}^{l_{n+1}}(0,1)_j\bigcup Sp(F_{n+1})$.}}
%As
%condition (3)  $\phi_{_{n, n+1}}{(C_0\big((0,1), E_n\big))}\sbs
%C_0\big((0,1), E_{n+1}\big)$ required, homomorphism $\pi_{n+1} \circ \phi_{n,n+1}: C_n \to F_{n+1}$
%{has to take} the ideal $I_n=C_0\big((0,1), E_n\big)$ to zero, where
% $\pi_{n+1}: C_{n+1}  \to F_{n+1}=C_{n+1}/C_0\big((0,1), E_{n+1}\big)$ is the canonical quotient map, and therefore $\pi_{n+1} \circ \phi_{n,n+1}$ factors as
%\begin{displaymath}
%\xymatrix{ C_n\ar[r]&  F_n~ \ar[r]^{~~~\tiny \td\phi_{_{n,
%n+1}~~~~~~~~~}} & F_{n+1}.}
%\end{displaymath}
%{\bf{Are you verifying condition (3) above? }}
%
%Consequently,  for  for $\tht_j'\in Sp(F_{n+1})\sbs
%Sp(C_{n+1})$, map $\phi_{_{n, n+1}}|_{_{\tht'_j}}$ factors
%through $F_n$ {as below:}
%\begin{displaymath}
%\xymatrix{ C_n\ar[r]&  F_n~~ \ar[r]^{\tiny{\pi_{\tht_j'}\circ \td\phi_{_{n,
%n+1}}}} &~~~~ F_{n+1}|_{_{\tht'_j}}~.}
%\end{displaymath}
%where $\pi_{\tht_j'}: F_{n+1} \to F_{n+1}|_{\tht_j'}$ is the corresponding quotient map, see \ref{homrestr} for notation.
%
%If $l_{n+1}=0$, then $A_{n+1}=F_{n+1}$ (no $E_{n+1}$ appearing) and the definition of homomorphism $\phi_{n,n+1}$ is done.

%We need to specify for each $t\in (0,1)_{j_0}$,
%$$
%\phi_{_{n, n+1}}|_{_t} (f)\in M_{\{n+1,j\}}(\C)=E_{n+1}^{j_0}\sbs
%\bigoplus_{j=1}^{l_{n+1}}E_{n+1}^j=E_{n+1}, ~~~\forall f\in C_n.
%%$$
%We
%will {now} define
%another
{To actually construct $\phi_{n, n+1},$ we  first construct
a} homomorphism
$\psi:~ C([0,1],E_n)\to
C([0,1],E_{n+1})$.
{This can be done by define, for each $j,$
%We only need to define
the map}
$$
\qq\qq\qq\psi^j:~C([0,1],E_n)\to C([0,1],E^j_{n+1})\qq\qq \mbox{for
each }~ j=1,2,..., l_{n+1}.
$$

Let $(f_1,f_2,..., f_{l_n})\in C([0,1],E_n)$.  For any $k\in
\{1,2,..., l_{n}\}$, if $d_{jk}>0$, then let
$$F_k(t)=\diag\big(\underbrace{f_k(t),f_k(t),...,
f_{k}(t)}_{d_{jk}}\big) \in C\left([0,1],M_{d_{jk}\cdot
\{n,k\}}(\C)\right);$$ if $d_{jk}<0$, then let
$$F_k(t)=\diag\big(\underbrace{f_k(1-t),f_k(1-t),...,
f_{k}(1-t)}_{|d_{jk}|}\big) \in C\left([0,1],M_{|d_{jk}|\cdot
\{n,k\}}(\C)\right);$$ if $d_{jk}=0$, then let
$$F_k(t)=\diag\big(f_k(t),f_k(1-t)\big) \in C\left([0,1],M_{2\cdot
\{n,k\}}(\C)\right).$$

With the above notation, define
\beq\label{defpsi}
\psi^j(f_1,f_2,..., f_{l_n})(t)=\diag(F_1(t),F_2(t),..., F_{l_n}(t))
\in C\left([0,1],M_{\left({\sum_{k=1}^{l_n}d_k'}\cdot
\{n,k\}\right)}(\C)\right),
\eneq
where
$$
d_k'=\left\{\begin{array}{cc}
                |d_{jk}| & \mbox{ if } d_{jk}\not= 0, \\
                ~&\\
                2 & \mbox{ if } d_{jk} = 0.
              \end{array}\right.
$$
Note that
\beq\label{sizeE}
\{n+1,j\}=\sum_{l=1}^{p_{n+1}}b_{0,jl}'[n+1,l]
=\sum_{l=1}^{p_{n+1}}\sum_{i=1}^{p_{n}}b_{0,jl}'c_{li}[n,i]
=\sum_{i=1}^{p_{n}}\td b_{0,ji}[n,i].
\eneq
Recall that $\td \bb_0=\bb_0'\cdot \cc=\big(\td b_{0,ji}\big)$. From
 (\ref{13spd-1}),  \ref{sizeE}, $d_k'\leq |d_{kj}|+2$ and $\{n,k\}=\sum
b_{0,ki}[n,i]$, we have
$$
\{n+1,j\}=\sum_{i=1}^{p_{n}}\td b_{0,ji}[n,i]> \sum_{i=1}^{p_{n}} \big( \sum_{k=1}^{l_n}\big(|d_{jk}|+2\big)b_{0,ki}\big)[n,i]\geq \sum_{k=1}^{l_n}d_k'\{n,k\}.
$$
Hence the C*-algebra
$C\left([0,1],M_{\left({\sum_{k=1}^{l_n}d_k'}\cdot
\{n,k\}\right)}(\C)\right)$ can be regarded as a corner of the C*-algebra
$C([0,1],E_{n+1}^j)=C\left([0,1],M_{ \{n+1,j\}}(\C)\right)$, and
consequently, $\psi^j$ can be regarded as mapping into
$C([0,1],E_{n+1}^j)$.  Put all $\psi^j$ together we get $\psi:~
C([0,1],E_n)\to C([0,1],E_{n+1})$.

Define $\psi_0,\psi_1:~ C_n\to E_{n+1}$ to be
$$
\psi_0(f)=\psi(f)(0) \qq\qq \mbox{and}\qq\qq \psi_1(f)=\psi(f)(1)
$$
for any $f\in C_n\sbs C([0,1],E_n)$.  Since
$\psi_0(C_0((0,1),E_n))=0$ and $\psi_1(C_0((0,1),E_n))=0$, this
defines maps $\af_0,~\af_1:~ F_n\to E_{n+1}$.
 Note that for each $j\in\{1,2, ... , l_{n+1}\}$, the maps
 $\af_0^j,~ \af_1^j :~F_n\to E_{n+1} \to E_{n+1}^j$ have spectra
$$
SP(\af_0^j) =\left\{\tht_1^{\sim i_1},\tht_2^{\sim i_2},...,
\tht_{p_n}^{\sim i_{p_n}}\right\} \qq\mbox{ and }\qq SP(\af_1^j)
=\left\{\tht_1^{\sim i_1'},\tht_2^{\sim i_2'},..., \tht_{p_n}^{\sim
i_{p_n}'}\right\}
$$
respectively
(see \ref{ktimes} for the notation used here), where
$$i_l=\sum_{d_{jk}<0}|d_{jk}|b_{1,kl}+\sum_{d_{jk}>0}|d_{jk}|b_{0,kl}
+\sum_{d_{jk}=0}(b_{0,kl}+b_{1,kl})
$$
$$\qq\mbox{ and }\qq
i_l'=\sum_{d_{jk}>0}|d_{jk}|b_{1,kl}+\sum_{d_{jk}<0}|d_{jk}|b_{0,kl}
+\sum_{d_{jk}=0}(b_{0,kl}+b_{1,kl}).\qq\qq
$$
Note that $i_l'-i_l=\sum_{k=1}^{l_n}d_{jk}(b_{1,kl}-b_{0,kl})$, and note
that if $l\leq p_n^0$ in the case \ref{range 0.5a}, then
$b_{0,kl}=b_{1,kl}=0$ and consequently, $i_l=i_l'=0$. Put
\begin{displaymath}
\xymatrix{ \td\af_0=\bt_0'\circ \td\phi_{_{n,n+1}}:~F_n~
\ar[r]^{~~~~~~~~~~~~ \td\phi_{_{n, n+1}}} & F_{n+1} \ar[r]^{\tiny
\bt'_{0}}& E_{n+1}~}\andeqn\\\\
%\end{displaymath}
%and
%\begin{displaymath}
\xymatrix{ \td\af_1=\bt_1'\circ \td\phi_{_{n,n+1}}:~F_n~
\ar[r]^{~~~~~~~~~~~~ \td\phi_{_{n, n+1}}} & F_{n+1} \ar[r]^{\tiny
\bt'_{1}}& E_{n+1}~}.
\end{displaymath}
Then for each $j\in\{1,2,...,l_{n+1}\} $, the maps
$\td\af_0^j,~\td\af_1^j:~  F_n\to E_{n+1}\to E_{n+1}^j$ have spectra
$$%%%%%
SP(\td\af_0^j) =\left\{\tht_1^{\sim \bar i_1},\tht_2^{\sim \bar
i_2},..., \tht_{p_n}^{\sim \bar i_{p_n}}\right\} \qq\mbox{ and }\qq
SP(\td\af_1^j) =\left\{\tht_1^{\sim \bar i_1'},\tht_2^{\sim \bar
i_2'},..., \tht_{p_n}^{\sim \bar i_{p_n}'}\right\},
$$
where
$$\bar i_l=\sum_{k=1}^{p_{n+1}}b'_{0,jk}c_{kl}=\td b_{0,jl}
\qq\mbox{ and }\qq \bar i'_l=\sum_{k=1}^{p_{n+1}}b'_{1,jk}c_{kl}=\td
b_{1,jl}.
$$
From  (\ref{13spd-1}), we have that $\bar i_l> i_l$ and $\bar
i_l'> i_l'$. Furthermore $\bar i_l'-\bar
i_l=\sum_{k=1}^{p_{n+1}}(b'_{1,jk}-b'_{0,jk})c_{kl}$. Since
$(\bb_1'-\bb_0')\cc=\dd(\bb_1-\bb_0)$, we have that $\bar i_l'-\bar i_l=
i_l'- i_l$, and hence $\bar i_l'-i_l'= \bar i_l- i_l\triangleq r_l>0$.
Note that these numbers are defined for the homomorphisms $\af_0^j,~
\af_1^j, \td\af_0^j, \td\af_1^j :~F_n\to  E_{n+1}^j$. So strictly
speaking, $r_l$ should be denoted as $r_l^j> 0$.

Define a unital homomorphism $\Phi:~ C_n\to
C([0,1],E_{n+1})=\oplus_{j=1}^{l_{n+1}}C([0,1],E_{n+1}^j)$ by
$$\Phi^j(f_1,f_2,\cd,f_{l_n},a_1,a_2,... ,a_{p_n})=
\diag\big(\psi^j(f_1,f_2,\cd,f_{l_n}), a_1^{\sim r_1^j},a_2^{\sim
r_2^j},... ,a_{p_n}^{\sim r_{p_n}^j} \big).
$$
Again, define the maps $\Phi_0,~ \Phi_1:~ C_n\to E_{n+1}$ by
$$\Phi_0(F)=\Phi(F)(0)\andeqn\Phi_1(F)=\Phi(F)(1),$$
for $F=(f_1,f_2,...,f_{l_n},a_1,a_2,...,a_{p_n})\in C_n$.  These two
maps induce two quotient maps
$$ \td{\td\af}_0,~\td{\td\af}_1,~:F_n\to E_{n+1}.$$

From our calculation, for each $j\in\{1,2,...,l_{n+1}\} $, the map $
\td{\td\af}_0^j$ ($\td{\td\af}_1^j$ resp.) has same
spectrum as ${\td\af}_0^j$ (${\td\af}_1^j$ resp.)
does.  That is,
 $(\td{\td\af}_0^j)_{*0}=({\td\af}_0^j)_{*0}$ and $(\td{\td\af}_1^j)_{*0}=({\td\af}_1^j)_{*0}$.
 There are unitaries $U_0,~ U_1\in E_{n+1}$ such that $\mathrm{Ad}\,U_0\circ \td{\td\af}_0^j
 ={\td\af}_0^j$ and $\mathrm{Ad}\,U_1\circ \td{\td\af}_1^j
 ={\td\af}_1^j$.  Choose a unitary path $U\in C([0,1],E_{n+1})$ such that $U(0)=U_0$
 and $U(1)=U_1$.  Finally, set $\phi_{_{n, n+1}}:~ C_n\to C([0,1],E_{n+1})$ to
 be defined as
$\phi_{_{n, n+1}}=\mathrm{Ad}\,U\circ \Phi$.
From the construction, we have that $\psi(C_0\big((0,1), E_n\big))\sbs C_0\big((0,1),
E_{n+1}\big)$ and consequently, $\Phi(C_0\big((0,1), E_n\big))\sbs C_0\big((0,1),
E_{n+1}\big)$.

We conclude that $\phi_{_{n, n+1}}(I_n)\sbs I_{n+1}$.
%{Since $\psi$ is a $C([0,1])$-map and $U\in C([0,1], E_{n+1}),$ $\Phi$ is also a $C([0,1])$-map. }
Since $\mathrm{Ad}U(0)\circ
\td{\td\af}_0^j
 ={\td\af}_0^j=\bt_0'\td\phi_{_{n, n+1}}$ and $\mathrm{Ad}\,U(1)\circ \td{\td\af}_1^j
 ={\td\af}_1^j=\bt_1'\td\phi_{_{n, n+1}}$ we have that $\phi_{_{n, n+1}} (C_n)\subset C_{n+1}$ and furthermore,  the quotient map from
$C_n/I_n\to C_{n+1}/I_{n+1}$ induced by $\phi_{_{n, n+1}}$ is the
same as $\td\phi_{_{n, n+1}}$ (see definition of ${\td\af}_0^j$ and
${\td\af}_1^j$).
%Thus finish
Note that Condition (2)--(4) hold evidently.  {Condition (6) follows from (\ref{defpsi}) and the the definition of $F_k(t)$.}   If $x\in Sp(F_{n+1})\subset Sp(C_{n+1})$, then $Sp(F_n)\subset Sp(\phi_{n,n+1}|_x)(=Sp(\td \phi_{n,n+1}|_x)),$
%follows from
{by}  the fact that all entries of $\cc$ are strictly positive. If $x\in (0,1)_j=Sp(C_0((0,1), E_{n+1}^j)$, then each $\tht_i$, as the only element in $Sp(F_n^i)(\subset Sp(F_n))$,
%{{\bf This is the first time $F_n^i$ appears.
%Could I assume that it is a Simple direct summand of $F_n$?}})
 appears $r_i^j>0$ times in
$Sp(\phi_{n,n+1}|_x)$. Consequently, we also have $Sp(F_n)\subset Sp(\phi_{n,n+1}|_x)$. Hence Condition (5) holds.
This finishes the proof.
\end{proof}

\begin{NN}\label{range 0.8}  Let $\phi:~ C_n\to C_{n+1}$ be as in the proof above.
We will calculate the
map
$$\phi_{_{n,n+1}}^{\sharp}:~ \Aff(T(C_n))\to \Aff(T(C_{n+1})).$$
Recall from \ref{2Rg15} that
%3.20,
$$\Aff(T(C_n)) \sbs \bigoplus_{i=1}^{l_n}C([0,1]_i,\R)\oplus \R^{p_n}$$
consisting of $(f_1,f_2,...,f_{l_n},h_1,h_2,...,h_{p_n})$ which  satisfies
the conditions
\beq\label{1508/13star-2}
(*)\qq\qq\qq f_i(0)&=&\frac1{\{n,i\}}\sum b_{0,ij}h_j\cdot [n,j]\\\label{13stars-3}
(**)\qq\qq\qq f_i(1)&=&\frac1{\{n,i\}}\sum b_{1,ij}h_j\cdot [n,j],
\eneq
and $\Aff(T(C_{n+1})) \sbs \bigoplus_{i=1}^{l_{n+1}}C([0,1]_i,\R)\oplus
\R^{p_{n+1}}$ consisting of
$(f'_1,f'_2,...,f'_{l_{n+1}},h'_1,h'_2,...,h'_{p_{n+1}})$ which satisfies
$$f'_i(0)=\frac1{\{n+1,i\}}\sum b'_{0,ij}h'_j\cdot [n+1,j]
\qq\mbox{ and } \qq f'_i(1)=\frac1{\{n+1,i\}}\sum b'_{1,ij}h'_j\cdot
[n+1,j].
$$
Let
$$\phi_{_{n,n+1}}^{\sharp}(f_1,f_2,...,f_{l_n},h_1,h_2,...,h_{p_n})
=(f'_1,f'_2,...,f'_{l_{n+1}},h'_1,h'_2,...,h'_{p_{n+1}}).$$ Then
$$h_i'=\frac1{[n+1,i]}\sum_{j=1}^{p_n}c_{ij}h_j[n,j].$$
Recall that $(c_{ij})_{p_{n+1}\times p_n}$ is the matrix
corresponding to $(\td\phi_{_{n,n+1}})_{*0}=\gm_{_{n,n+1}}$ for
$\td\phi_{_{n,n+1}}: F_n\to F_{n+1}$, and note that since  $\td\phi_{_{n,n+1}}$ is unital, one has
$$\sum_{j=1}^{p_n}c_{ij}[n,j]=[n+1,j].$$  Also note that
\beq\nonumber
f'_i(t)=\frac1{\{n+1,i\}}\left\{\sum_{d_{ik}>0}d_{ik}f_k(t)\{n,k\}+
\sum_{d_{ik}<0}|d_{ik}|f_k(1-t)\{n,k\}+\right.\qq\qq\qq\qq\\\label{13starss}
\qq\qq\qq\qq +\left.\sum_{d_{ik}=0}(f_k(t)+f_k(1-t))\{n,k\}+
\sum_{l=1}^{p_n}r_l^ih_l[n,l]\right\}, \qq\qq
\eneq
where
 \begin{eqnarray}\label{rli}
r_l^i & = & \sum_{k=1}^{p_{n+1}}b_{0,ik}'c_{kl}-
\left(\sum_{d_{ik}<0}|d_{ik}|b_{1,kl}+
\sum_{d_{ik}>0}|d_{ik}|b_{0,kl}+ \sum_{d_{ik}=0}(b_{0,kl}+b_{1,kl})
\right)\nonumber \\
& = & \sum_{k=1}^{p_{n+1}}b_{1,ik}'c_{kl}-
\left(\sum_{d_{ik}>0}|d_{ik}|b_{1,kl}+
\sum_{d_{ik}<0}|d_{ik}|b_{0,kl}+ \sum_{d_{ik}=0}(b_{0,kl}+b_{1,kl})
\right).
\end{eqnarray}
It follows from the last paragraph of \ref{range  0.5} and \ref{range
0.5a} that, when we define $C_{n+1}$, we can always increase the entries
of the last $p_{n+1}^1=(p_{n+1}- p_{n+1}^0)$  columns of the matrices
$\bb_0'=(b_{0,ik}')$ and $\bb_1'=(b_{1,ik}')$ by adding an
arbitrarily (but same for $\bb_0'$ and $\bb_1'$)
 matrix $(m_{ik})_{l_{n+1}\times (p_{n+1}^1}$, with each $m_{ik}>0$
 sufficiently large, to the last $p_{n+1}^1$ columns of the matrices.
 In particular we can strengthen the requirement $\spd$ (see (\ref{13spd-1})) to
 \beq\label{13spdd}
\spdd \qq\qq \td b_{0,il}=\!\sum_{k=1}^{p_{_{n+1}}}b_{0,ik}'c_{kl}
 \!\!\left(\!\!\mbox{and } \td
 b_{1,il}=\sum_{k=1}^{p_{n+1}}b_{1,ik}'c_{kl}\!\!\right)
 > 2^{2n}\!\left(\sum_{k=1}^{l_n}
 (|d_{ik}|+2)\{n,k\}\!\!\right)
 \eneq
 for all
 $i\!\in\!\{1,...,l_{n+1}\}.$
 %(Again as the discussion about (\ref{13spd-1}) if $l_{n+1}=0$, then the condition holds tautologically, and if $l_n=0$, then the right hand side is a sum over empty set which is regarded to be zero.)
 %{\bf {Need to consider two cases. It should be made clearer, Answer: this is already discussed when we talked about (\ref{13spd-1})}}
This condition and \ref{rli} (and note that $b_{0,kl}\leq \{n,k\}$, $b_{1,kl}\leq \{n,k\}$ for any $k\leq l_n$) imply
\beq\label{13spdd'}
r_l^i \geq \frac{2^{2n}-1}{2^{2n}}\td b_{0,il} \qq \mbox{and
equivalently} \qq 0\leq \td b_{0,il}-r_l^i< \frac{1}{2^{2n}}\td
b_{0,il}.
\eneq
%%%%%%%%%%%%%%%
Recall that $\phi_{_{n,n+1}}^{\sharp}:~ \Aff(T(C_{n}))\to \Aff(T(C_{n+1}))$, and
$\td\phi_{_{n,n+1}}^{\sharp}:~ \Aff(T(F_{n}))\to \Aff(T(F_{n+1}))$ are the maps induced by corresponding homomorphisms (i.e., by the homomorphisms $\phi_{n,n+1}$ and $\td\phi_{_{n,n+1}}$; and same for $\pi_n^{\sharp}: ~ \Aff(T(C_{n}))\to \Aff(T(F_{n}))$ and
$\pi_{n+1}^{\sharp}: ~ \Aff(T(C_{n+1}))\to \Aff(T(F_{n+1}))$, which are the maps
induced by quotient
 maps
  $\pi_n: C_n\to F_n$ and $\pi_{n+1}: C_{n+1}\to F_{n+1},$ respectively.
    Since $\td\phi_{_{n,n+1}}\circ \pi_n=\pi_{n+1}\circ \phi_{_{n,n+1}}$, we have
$$
\pi_{n+1}^{\sharp}\circ\phi_{_{n,n+1}}^{\sharp}=\td\phi_{_{n,n+1}}^{\sharp}\circ\pi_n^{\sharp}:~ \Aff(T(C_{n+1}))\to \Aff(T(F_{n+1})).
$$

\end{NN}

\begin{NN}\label{range 0.9}
For each $n,$
%and the quotient map $\pi_n: A_n\to F_n$,
we will define
a map $\GM_n:~ \Aff(T(F_{n}))\to \Aff(T(C_{n}))$ which is a right inverse of
$\pi_n^{\sharp}:~ \Aff(T(C_{n}))\to \Aff(T(F_{n}))$---that is, $\pi_n^{\sharp}\circ
\GM_n=\id|_{_{\Aff(T(F_n))}}$ as below.

Recall that $C_n=A(F_n,E_n,\bt_0,\bt_1)$ with unital homomorphisms
$\bt_0,\bt_1: F_n\to E_n$ whose K-theory maps satisfy
$(\bt_0)_{*0}=\bb_0=(b_{0,ij})$ and  $(\bt_1)_{*0}=\bb_1=(b_{1,ij})$. For
$(h_1,h_2,...,h_{p_n})\in \R^{p_n}$, the linear maps $\bt_0^{\sharp},
\bt_1^{\sharp}:~ \Aff(T(F_n))=\R^{p_n} \to \Aff(T(E_n))=\R^{l_n}$ are given by the
matrices $\left(
            \begin{array}{c}
              \frac{b_{0,ij}[n,j]}{\{n,i\}} \\
            \end{array}
          \right)$
and $\left(
            \begin{array}{c}
              \frac{b_{1,ij}[n,j]}{\{n,i\}} \\
            \end{array}
          \right),$
where $i\in \{1,2,..., l_n\},~ j\in \{1,2,..., p_n\}$---that is,
$$
\bt_0^{\sharp}\left(
           \begin{array}{c}
             h_1 \\
             h_2 \\
             \vdots \\
             h_{p_n} \\
           \end{array}
         \right) =
\left(
           \begin{array}{c}
             \frac1{\{n,1\}}\sum_{j=1}^{p_n} b_{0,1j}[n,j]\cdot h_j \\
             \frac1{\{n,2\}}\sum_{j=1}^{p_n} b_{0,2j}[n,j]\cdot h_j \\
             \vdots \\
             \frac1{\{n,l_n\}}\sum_{j=1}^{p_n} b_{0,l_nj}[n,j]\cdot h_j \\
           \end{array}
         \right).
$$
(For  $\bt_1^{\sharp}$, one replaces $b_{0,ij}$ by $b_{1,ij}$ in the construction above).
For $h=(h_1,h_2,...,h_{p_n})$, set
$$\GM_n'(h)(t)= t\cdot \bt_1^{\sharp}(h)+(1-t)\cdot \bt_0^{\sharp}(h),$$ which
is an element in $C([0,1],\R^{l_n})=\bigoplus_{i=1}^{l_n}C([0,1]_i,\R)$.
Finally, define the map
%$$
\beq\nonumber
&&\hspace{-0.2in}\GM_n:~ \Aff(T(F_{n}))=\R^{p_n}\to \Aff(T(C_{n})) \sbs
\bigoplus_{j=1}^{l_n}C([0,1]_j,\R)\oplus \R^{p_n} \,\,\,{\rm by}\\
&&\GM_n(h)=(\GM_n'(h), h)\in \bigoplus_{j=1}^{l_n}C([0,1]_j,\R)\oplus
\R^{p_n}.
\eneq
%$$
Note that $\GM_n(h)\in \Aff(T(C_{n}))$ (see (\ref{1508/13star-2}) and (\ref{13stars-3}) in \ref{range 0.8}).
 Evidently, $\pi_n^{\sharp}\circ \GM_n=\id|_{_{\Aff(T(F_n))}}$.\\

\end{NN}

\begin{lem}\label{range 0.10}
 If condition $\spdd$ (see $(\ref{13spdd})$) holds, then for any $f\in
\Aff(T(C_{n}))$ with $\|f\|\leq 1$, and $f'=\phi_{_{n,n+1}}^{\sharp}(f)\in \Aff(T(C_{n+1}))$, we
have
$$\| \GM_{n+1}\circ\pi_{n+1}^{\sharp}(f')-f'\|<\frac2{2^{2n}}.$$
\end{lem}

\begin{proof}
Write 
\beq\nonumber
f=(f_1,f_2,...,f_{l_n},h_1,h_2,...,h_{p_n})\in \Aff(T(C_{n}))
%$$ and
%$$
\andeqn\\
f'=(f'_1,f'_2,...,f'_{l_{n+1}},h'_1,h'_2,...,h'_{p_{n+1}})\in
\Aff(T(C_{n+1})).
\eneq
%$$ 
Since
$\pi_{n+1}^{\sharp}\circ\GM_{n+1}=\id|_{_{\Aff(T(F_{n+1}))}}$, one has
$$\GM_{n+1}\circ\pi_{n+1}^{\sharp}(f'):=g'
:=(g'_1,g'_2,... ,g'_{l_{n+1}},h'_1,h'_2,... ,h'_{p_{n+1}});$$
%{\color{Green} (Do you mean $\GM_{n+1}(f')$ instead of $\GM_{n+1}\circ\pi_{n+1}^{\sharp}(f')$?)} {\red {Answer: I think it is OK since $\GM_{n+1}$ is from $\Aff(T(F_{n+1}))$ and $f'\in \Aff(T(C_{n+1}))$ only  $\pi_{n+1}^{\sharp}(f') \in \Aff(T(F_{n+1}))$---Gong}}
that is, $f'$ and $g'$ have the same boundary value
$(h'_1,h'_2,\cd,h'_{p_{n+1}})$.
%{\color{Green} (Don't understand this part, do you mean ``Since
%$$\Pi_{n+1}\circ\GM_{n+1}=\id|_{_{\Aff(T(F_{n+1}))}},$$-
%then
%$$g':=\GM_{n+1}\circ\Pi_{n+1}(f')=(g'_1,g'_2,\cd,g'_{l_{n+1}},h'_1,h'_2,\cd,h'_{p_{n+1}})$$---that %is, ..."?) Answer: Yes, this is what I mean, that $|Pi_{n+1}(g')=\Pi_{n+1}(f')$}

Note that the evaluations of $f'$ at zero, $(f'_1(0),f'_2(0),...,f'_{l_{n+1}}(0))$, and the evaluation at at one, $(f'_1(0),f'_2(0),...,f'_{l_{n+1}}(0))$
are
%(Note that the evaluation
%at $0$
%$(f'_1(0),f'_2(0),\cd,f'_{l_{n+1}}(0))$  at and also the evaluation at
%$1$) is
completely determined by $h'_1,h'_2,...,h'_{p_{n+1}}$.
From (\ref{13starss}) of \ref{range 0.8}, we also have
$$
f'_i(t)-f'_i(0)=\frac1{\{n+1,i\}}\left(\sum_{d_{ik}>0}d_{ik}f_k(t)\{n,k\}+
\sum_{d_{ik}<0}|d_{ik}|f_k(1-t)\{n,k\}\right.+\qq\qq\qq\qq\qq$$
$$\qq\qq\qq\qq\qq\qq\qq\left.+\sum_{d_{ik}=0}(f_k(t)+f_k(1-t))\{n,k\}-
\sum_{l=1}^{p_n}(\td b_{0,il}-r_l^i)h_l[n,l]\right).
$$
And from $(\ref{13spdd})$, one has
$$
\left|\sum_{d_{ik}>0}d_{ik}f_k(t)\{n,k\}+
\sum_{d_{ik}<0}|d_{ik}|f_k(1-t)\{n,k\}
+\sum_{d_{ik}=0}(f_k(t)+f_k(1-t))\{n,k\}\right|
$$
$$
\leq \sum_{k=1}^{l_n}(|d_{ik}|+2)\{n,k\} \leq \frac1{2^{2n}} \td
b_{0,il} < \frac1{2^{2n}}\cdot\{n+1,i\}.
$$
Note that $ \td b_{0,il}<\{n+1,i\},$  since $\sum_{l=1}^{p_n}\td
b_{0,il}[n,l]=\{n+1,i\}$ and $[n,l]{\ge 1}$.  Combining this with $(\ref{13spdd'})$, we
have
$$|f_i'(t)-f_i'(0)|<\frac2{2^{2n}}.$$
Similarly,  we have
$$|f_i'(t)-f_i'(1)|<\frac2{2^{2n}}.$$
But, by the definition of $\GM_{n+1},$ we have
$$g_i'(t)=tg_i'(1)+(1-t)g_i'(0).$$
Combining with $g_i'(0)=f_i'(0)$ and $g_i'(1)=f_i'(1)$, we have
$$|g_i'(t)-f_i'(t)|<\frac2{2^{2n}}\rforal i,$$
as desired.
\end{proof}

\begin{rem}\label{range 0.11} In the proof of Lemma \ref{range  0.10}, the key
point is that for any
$f'=\phi_{_{n,n+1}}^{\sharp}(f)\in \Aff(T(C_{n+1}))$, we have that all components of $f'$, as $\mathbb R^l_{n+1}$-valued functions on $(0,1)$, are
close to a constant function.
In the definition of $\phi_{_{n,n+1}}$ of the proof of \ref{range  0.6}, we
know that each component
$f_j'=\phi^j_{_{n,n+1}}(f_1,f_2,...,f_{l_n},a_1,a_2,...,a_{p_n})$ is
written as
$$
\diag\big(\psi^j(f_1,f_2,...,f_{l_n}), a_1^{\sim r_1^j},a_2^{\sim
r_2^j},...,a_{p_n}^{\sim r_{p_n}^j} \big)
$$
up to unitary equivalence.  The major part {of} $\diag\big( a_1^{\sim
r_1^j},a_2^{\sim r_2^j},...,a_{p_n}^{\sim r_{p_n}^j} \big)$
is constant (of course up to unitary equivalence)---that is,  each $a_i$ appeared above  is a fixed matrix as in $F_n$ and does not depend on $t\in (0,1)_j\subset Sp C_0((0,1), E_{n+1})$ as a function.
%{\bf {This is unclear as what it means}}
In fact, Condition $(\ref{13spdd})$
 implies that this part occupies more than $\frac{2^{2n}-1}{2^{2n}}$ of the
 whole size $\{n+1,j\}$ of $M_{\{n+1,j\}}(C[0,1])=C([0,1],E_{n+1}^j)$.
 %{\bf {If this is important, which  I guess
 %it is, then it is better spelled out here.}}
  In the rest of this section, we will use this argument several times. That is, some condition similar to $\spdd$ (see (\ref{13spdd})) will imply that the map $\phi_{n,n+1}^\#$ has the property that for any $f$, the image $f'=\phi_{n,n+1}^\#(f)$ is close to constant (on certain connected subsets of spectrum), and therefore it ``almost only" depends on its value at finitely many points (in the above case, it depend on its values on the spectrum of $F_n$). Then pass to the inductive limit and use the approximately intertwining argument, one can prove that the tracial state space of the limit algebra is as desired---in the case above, $\lim (AffTC_n, \phi_{n,n+1}^\#)=\lim (AffTF_n, \tilde\phi_{n,n+1}^\#)$ (see \ref{range 0.12} below).

%In this case we will not write the detailed calculation as above,  but only
%  specify the condition corresponding to $(\ref{13spdd})$.\\
%{\bf Almost constant=? major part constant}

%{\bf {Perhaps, a stupid question could have your attention:
%Do you mean you will apply Lemma \ref{range 0.10} without warning? If so, one can say things like:
%Since condition (\ref{13spdd}) holds, by Lemma \ref{range 0.10}, we have....
%}}
%{\bf { It is not clear to me that the above paragraph means that.  You probably means %somethings else.}}

\end{rem}

\begin{thm}\label{range 0.12} If $C_n$ and $\phi_{_{n,n+1}}:~ C_n\to C_{n+1}$ are
as in {\rm \ref{range  0.6}} and  {\rm \ref{range 0.7}} with condition $\spd$ (see $(\ref{13spd-1})$)
being replaced by the stronger condition $\spdd$ (see $(\ref{13spdd})$), then the inductive limit $C=\lim
(C_n,\phi_n)$ has the  Elliott invariants $(K_0(C),
K_0(C)_+,\e_C, T(C), r_{\!\!_C})=(G,G_+,\DT,r)$ (here we assume $G$ is
torsion free and $K_1(C)=\{0\}$).

\end{thm}

%{\bf {Does condition (\ref{13spdd}) stronger than that of (\ref{13spd-1})?
%If so, we should inform the reader. If not, a question may be raised: Could the same
%construction be carried out------It was stated before (\ref{13spd-1}), that construction %requires
%(\ref{13spd-1})?}}

\begin{proof}  From the construction, we have the following infinite  commutative diagram:
\begin{displaymath}
    \xymatrix{ I_1 \ar[r]\ar[d] & I_2 \ar[r]\ar[d]& I_3 \ar[d]\ar[r]&\cd I \\
               C_1 \ar[d]\ar[r] & C_2 \ar[r]\ar[d] &C_3 \ar[d]\ar[r]&\cd C \\
               C_1/I_1 \ar[r] & C_2/I_2 \ar[r] &C_3/I_3 \ar[r]&\cd C/I~,}
\end{displaymath}
where $C_n/I_n=F_n$.  Also from the construction, we have the
following diagram:
\begin{displaymath}
    \xymatrix{0 \ar[r]\ar[d] & 0 \ar[r]\ar[d]& 0 \ar[d]\ar[r]&\cd  \\ K_0(C_1)=G_1 \ar[r]\ar[d] & K_0(C_2)=G_2 \ar[r]\ar[d]& K_0(C_3)=G_3 \ar[d]\ar[r]&\cd  \\
               K_0(C_1/I_1)=H_1 \ar[d]\ar[r] & K_0(C_2/I_2)=H_2 \ar[r]\ar[d] &K_0(C_3/I_3)=H_3 \ar[d]\ar[r]&\cd \\
               K_1(I_1)\ar[d] \ar[r] & K_1(I_2)\ar[d] \ar[r] &K_1(I_3)\ar[d] \ar[r]&\cd \\
               0 \ar[r] & 0 \ar[r]& 0 \ar[r]&\cd ~.}
\end{displaymath}
So $K_0(C)=\lim(G_n,\gm_{_{n,n+1}}|_{_{G_n}})=G$.

The inductive limits also induce the following diagram:
\begin{displaymath}
    \xymatrix{
        \Aff(T(C_1)) \ar[r] \ar[d]_{\pi_1^\#} & \Aff(T(C_2)) \ar[r] \ar[d]_{\pi_2^\#}& \Aff(T(C_3)) \ar[r] \ar[d]_{\pi_3^\#}& \cd \Aff(T(C)))) \\
    \Aff(T(F_1)) \ar[r] & \Aff(T(F_2)) \ar[r]& \Aff(T(F_3)) \ar[r]& \cd \Aff(T((C/I)))~. }
\end{displaymath}
By the construction of $\GM_n$ and \ref{range  0.10} we have the following diagram
approximately intertwining:
\begin{displaymath}
    \xymatrix{
        \Aff(T(C_1)) \ar@/_/[d]_{\pi_1^\#}\ar[r]^{\phi_{1,2}^\#}  & \Aff(T(C_2)) \ar[r]^{\phi_{2,3}^\#} \ar@/_/[d]_{\pi_2^\#}& \Aff(T(C_3)) \ar[r] \ar@/_/[d]_{\pi_3^\#}& \cd \Aff(T(C ))\\
        \Aff(T(F_1)) \ar[r]^{\widetilde{\phi}_{1,2}^\#}\ar@/_/[u]_{\GM_1} & \Aff(T(F_2)) \ar[r]^{\tilde{\phi}_{2,3}^\#}\ar@/_/[u]_{\GM_2}& \Aff(T(F_3)) \ar[r]\ar@/_/[u]_{\GM_3}& \cd \Aff(T(C/I))~. }
\end{displaymath}
%{\bf {It would perhaps helpful to have  $\Lambda$ appears in the above digram, or mentioned.}}
Hence by (2.2)--(2.4) of \cite{Ell-RR0}, it induces unital ordered isomorphism between unital ordered real Banach spaces $\Aff(T(C))$ and
$\Aff(T((C/I)))=\Aff(\DT)$. This isomorphism $\pi_{\infty}^\#$ and its inverse $\GM_{\infty}$ are continuous with respect to the weak $*$ topology of $\Aff(T(C))$ and
$\Aff(T((C/I)))$---that is, for any net $l_\af\in \Aff(T(C))$
% {\color{Green} (why sequences are not enough?)} {\red{I think it is ok to use sequence,but I don't want to verify it}}
and $l\in \Aff(T(C))$, if $l_\af (\tau)\to l(\tau)$ for all $\tau \in T(C)$, then $(\pi_{\infty}^\# (l_\af))(\tau_1) \to (\pi_{\infty}^\# (l))(\tau_1)$ for all $\tau_1 \in T(C/I)$. Furthermore, $T(C)$ (or $T(C/I)$) is the set of all positive linear maps $\tau$ from $\Aff(T(C))$ (or $\Aff(T(C/I))$) to $\R$ with $\tau ({\bf 1})=1$, which is continuous with respect to weak $*$ topology on $\Aff(TC)$ (or $\Aff(T(C/I))$).        Consequently $T(C)=T(C/I)=\DT$, as they can be recovered from ${\Aff}(T(C))={\Aff(T(C/I))}={\Aff}(\DT)$ and as they are the sets of all positive unital linear weak $*$ continuous map from the same set ${\Aff(T(C))}=\Aff(T(C/I))={\Aff}(\DT)$ to $\R$.

Furthermore, the homomorphisms $\pi_i: C_i\to C_i/I_i, i=1,2,...,$ induce the inclusion map $\iota: K_0(C)= G\to K_0(C/I)=H$ and the isomorphism ${\Aff}(T(C)) \to {\Aff(T(C/I))}={\Aff}(\DT)$ by the above two approximately intertwining diagram, and therefore it is compatible with the map
$r_{\!\!_C}$ and $r$---that is $\pi_{\infty}^\#(r_{\!\!_C}(x))=r(K_0(\pi)(x))$ for any $x\in K_0(C)$. (Here, we use the facts that $r: G\to {\Aff}(\DT)$ factors through $H$ and $K_0(\pi)(x)\in G\subset H$ for $x\in K_0(C)$.)
%{\bf {  "Consequence" must be from a theorem---which is not \cite{Ell-RR0}, I %assume----therefore it requires a reference.
%I could not figure out from, say Afre..."}}
\end{proof}

\begin{NN}\label{range 0.13} From the construction, the algebra $C$ in \ref{range 0.12} has an ideal $I=\lim(I_n, \phi_{n,m}|_{I_n})$   and therefore  is
not a simple $C^*$-algebra.  However, one can modify the homomorphism $\phi_{_{n,n+1}}$ to a
map $\psi_{_{n,n+1}}$ to have the following properties:
\begin{enumerate}
\item[(1)]  $\psi_{_{n,n+1}}\sim_h \phi_{_{n,n+1}}$.
\item[(2)]  $\|\psi_{_{n,n+1}}^{\sharp}-\phi_{_{n,n+1}}^{\sharp}\|\leq \frac1{2^n}$.
\item[(3)]  $Sp(\psi_{_{n,n+1}}|_{\tht_1'})$ is $\frac1n$-dense in
$Sp(C_n)$,
%{\bf {Do you really mean this? I remember same issues appeared elsewhere.
%The point: we are dealing with non-Huasdorff spaces. As it stands, there is no true %metric on
%$Sp(A_n)$ unless $A_n$ degrades itself to some circles or interval algebras}}
where $\tht_1'$ is the first point in
$$Sp(F_{n+1})=\{\tht_1',\tht_2',\cd,\tht'_{_{p_{n+1}}}\}\sbs
Sp(C_{n+1})$$ (see \ref{density} for concept of $\dt$-dense).
%(by $\delta$-density
%(here $\delta=\frac1n$), we mean that any isolate points (only in the case of \ref{range %0.5a}) of $Sp(F_n)$ are in $SP(\psi_{_{n,n+1}}|_{\tht_1'})$, and %$SP(\psi_{_{n,n+1}}|_{\tht_1'})\cap Sp(C_0((0,1),E_n))$ is $\delta$-dense in %$Sp(C_0((0,1),E_n))$, which is disjoint union of several open intervals $(0,1)$ with %standard metric);
\item[(4)] Like $\phi_{n,n+1}$, for any $x\in Sp(C_{n+1})$, we have $Sp(F_n)\subset Sp(\psi_{n,n+1}|_x)$.
% {\bf {Something wrong with the sentence. I could not even guess}}
\item[(5)] Again like $\phi_{_{n,n+1}}$, the homomorphism
$\psi_{_{n,n+1}}$ is injective and furthermore, if $X\sbs
Sp(C_{n+1})$ is $\dt$-dense in $Sp(C_{n+1})$, then $\cup_{x\in
X}Sp(\psi_{_{n,n+1}}|_x)$ is $\dt$-dense in $Sp(C_{n})$.
\end{enumerate}
(See the end of \ref{homrestr}, \ref{ktimes} and \ref{density} for the notation used here.)
%and later on in \ref{range 0.29}
%we will discuss $\dt$-density at {different} but similar case.)
%----{\bf {It looks that this is the first time it appeared.
%Perhaps, definition should be given here instead}})

 For each $n$ and $x\in Sp(C_{n+2})$, combining (3) (for $\psi_{n,n+1}$) and (4) (for $\psi_{n+1,n+2}$), we have that
\beq\label{dense}
Sp(\psi_{n,n+2}|_x)=\bigcup_{y\in Sp(\psi_{n+1,n+2}|_x)} Sp(\psi_{n,n+1}|_y) \supset Sp(\psi_{n,n+1}|_{\tht_1'})
\eneq
is $1/n$-dense in $Sp(C_n)$

Combining with (5), for any $n$ and $m\geq n+2$, and any
$x\in Sp(C_m)$, $Sp(\psi_{_{n,m}}|_x)$ is $\frac1{m-2}$ dense in
$Sp(C_{n})$.  Consequently, the new inductive limit
$B=\lim(C_n,\psi_{_{n,n+1}})$ is simple by \ref{simplelimit}.  By (1) and (2),
$$(K_0(C), K_0(C)_+,\e_C, T(C), r_{\!\!_C})\cong (K_0(B), K_0(B)_+,\e_B, T(B), r_{\!\!_B}).$$
Since we will construct the algebra $C\in {\cal N}_0^{\cal Z}$ with general
weakly unperforated $K_0$ group (without the torsion free condition)
and general $K_1$ group, we will not give detailed construction of
the above special case.

The following theorem is in \S 3 of \cite{EG-RR0AH} (see Lemma 3.21 and
Corollary 3.22 there).

\end{NN}

\begin{prop}\label{range 0.14}
Let $X$ and $Y$ be path connected finite CW complexes of dimension
at most three, with base point $x_0\in X$ and $y_0\in Y$, such that
the cohomology groups $H^3(X)$ and $H^3(Y)$ are finite. Let $\af_0:~
K_0(C(X))\to K_0(C(Y))$ be a homomorphism satisfying that $\af_0$ is at least 12-large
and that $\af_0(K_0(C(X))_+\setminus \{0\})\sbs K_0(C(Y))_+\setminus \{0\}$, and let $\af_1:~
K_1(C(X))\to K_1(C(Y))$ be any homomorphism.  Let $P\in M_\infty
(C(X))$ be any {non-zero} projection and $Q\in M_\infty (C(Y))$ be a
projection with $\af_0([P])=[Q]$ (such projections always exist.)
Then there exists a unital homomorphism $\phi:~ PM_\infty (C(X))P\to
QM_\infty (C(Y))Q$ such that $\phi_{*0}=\af_0$ and
$\phi_{*1}=\af_1$, and such that
$$
\phi(PM_\infty (C_0(X\setminus \{x_0\}))P)\sbs QM_\infty
(C_0(Y\setminus \{y_0\}))Q$$
That is, if $f\in PM_\infty (C(X))P$ satisfies $f(x_0)=0$, then $\phi(f)(y_0)=0$.\\
\end{prop}

\begin{rem}\label{range 0.15}
 In the proof of the above proposition in \cite{EG-RR0AH},
 %[EG-Annal],
one reduced the general case to the case that $P$ is rank one trivial
projection ($PM_\infty (C(X))P=C(X)$). If one further assumes that
$\af_0$ is at least 13-large and $Y\not= \{pt\}$, then the
homomorphism $\phi$ in the proposition can be chosen to be
injective.  To prove this we choose a surjective homotopy trivial
continuous map $g:~ Y\to X$, which induces an injective homomorphism
$g^*:~ C(X)\to C(Y)$. Then apply the theorem to the new map (still denoted by $\alpha_0$) $
\af_0:=\af_0-(g^*)_{*0}:~ K_0(C(X))\to K_0(C(Y))$ which is at least
12-large and $\af_1$ to obtain $\phi_1:~ C(X)\to QM_\infty C(Y)Q$.
Then
 $\phi=\diag(\phi_1, g^*):~ C(X)\to (Q\oplus 1)M_\infty C(Y)(Q\oplus 1)$
 has the desired property.
\end{rem}

The following is perhaps known.
%{\bf Why italian?}

\begin{lem}\label{range 0.15a} Let $0\to E\to H \to H/E \to 0$ be a
short exact sequence of countable abelian groups with $H/E$ torsion
free. And let
 \begin{displaymath}
    \xymatrix{
        H_1 \ar[r]^{\gm'_{_{1,2}}} & H_2 \ar[r]^{\gm'_{_{2,3}}}&H_3
         \ar[r]^{\gm'_{_{3,4}}} &\cdots \ar[r]& H/E }
\end{displaymath}
be an inductive system with limit $H/E$ such that each $H_i$ is
finitely generated free abelian group.  Then there are an increasing
sequence of finitely generated subgroups $E_1\subset
E_2\subset\cd\subset E_n \subset \cd \subset E$ with
$E=\bigcup_{i=1}^{\infty}E_i$, and an inductive system
\begin{displaymath}
    \xymatrix{
       E_1 \oplus H_{1}\ar[r]^{\gm_{_{1,2}}} & E_2 \oplus H_{2} \ar[r]^{\gm_{_{2,3}}}&E_3 \oplus H_{3}
         \ar[r]^{\gm_{_{3,4}}} &\cdots \ar[r]& H }
\end{displaymath}
with limit $H$ with the following properties:

{\rm (i)} $\gm_{n,n+1} (E_n) \subset E_{n+1}$ and $\gm_{n,n+1}|_{E_n}$ is
the inclusion from $E_n$ to $E_{n+1}$.

{\rm (ii)} Let $\pi_{n+1}: E_{n+1}\oplus H_{n+1}\to H_{n+1}$ be the
canonical projection, then $\pi_{n+1}\circ
\gm_{n,n+1}|_{H_{n}}=\gm'_{n,n+1}.$

(Here, we do not assume $\gm'_{n,n+1}$ to be injective.)

\end{lem}

\begin{proof} Let $E=\{e_i\}_{i=1}^{\infty}$. We will
construct the system $(E_n\oplus H_{n},\gm_{n, n+1})$, inductively.
Let us assume we already have $E_n \subset E$ with
$\{e_1,e_2,\cd,e_n\}\subset E_n$ and the map $\gm_{n,\infty}:
E_n\oplus H_n \to H$ such that $\gm_{n,\infty}|_{E_n}$ is the
inclusion and $\pi\circ \gm_{n, \infty}|_{H_n} =\gm'_{n, \infty}$,
where $\pi: H \to H/E$ is the quotient map. Note that
$\gm'_{n+1,\infty}(H_{n+1})$ is a finitely generated free abelian
subgroup of $H/E$, one has a lifting map {$\gm_{n+1, \infty}: H_{n+1} \to H$ such
that $\pi\circ \gm_{n+1, \infty} =\gm'_{n+1, \infty}$.}
%from
%$\gm'_{n+1,\infty}(H_{n+1}) (\subset H/E)$ to $H$. Using this
%lifting, we obtain a map $\gm_{n+1, \infty}: H_{n+1} \to H$ such
%that $\pi\circ \gm_{n+1, \infty} =\gm'_{n+1, \infty}$.
For each
$h\in H_n$, we have  $\gm (h)
:=\gm_{n,\infty}(h)-\gm_{n+1,\infty}(\gm'_{n,n+1}(h))\in E$. Let
$E_{n+1}\subset E$ be
{the}
%a
finitely generated subgroup generated {by}
$E_n\cup\{e_{n+1}\}\cup \gm(H_n)$ and extend the map $\gm_{n+1,
\infty} $ on $E_{n+1}\oplus H_{n+1}$ by defining it to be inclusion
on $E_{n+1}$. And finally let $\gm_{n,n+1}: E_n\oplus H_n \to
E_{n+1}\oplus H_{n+1}$ be defined by $\gm_{n,n+1}(e,
h)=(e+\gm(h),\gm'_{n,n+1}(h))\in E_{n+1}\oplus H_{n+1}$ for each
$(e,h) \in E_n\oplus H_n $. Evidently, $\gm_{n, \infty} =\gm_{n+1,
\infty} \circ \gm_{n,n+1}$.
\end{proof}

\begin{NN}\label{range 0.16}
%Now let $((G,G_+,u),K,\Delta,r)$ be the one given in
%\ref{range  0.1}.  As in \ref{range  0.3}, let $\rho:~ G\to \Aff\DT$
%be dual to the map $r$. Let $G^1\subset  \Aff \DT$ be a dense subgroup
%with at least three $\Q$-linearly independent elements.  Again as in
%\ref{range  0.3}, let $H=G\oplus G^1$ {and let ${\tilde \rho}: H\to
%\Aff(T(\DT))$ be defined by ${\tilde \rho}((g,f))(\tau)=\rho(g)(\tau)+f(\tau)$ for all %$(g,f)\in H$ and
%$\tau\in \DT.$  Let
%$H^+\setminus \{0\}=\{(g,f): \tilde{\rho}((g,f))>0\}.$}
%be the  {set of
%%collection of $(g,f)\in G\oplus G^1$ with
%$$\qq\qq\qq\qq\qq \rho(g)(\tau)+f(\tau)>0 \qq\qq\qq\qq \mbox{for all } \tau\in \DT.$$
%The order unit $u\in G^+$ could be regarded as $(u,0)\in G\oplus G^1=H$
%as the order
%unit of $H^+$.  Now $(H,H^+,u)$ is a simple weakly unperforated group with
%Riesz decomposition property.  Note that  $H$ has torsion ${\rm Tor}(G)={\rm Tor}(H)$ here.
 % {\bf {While this could be easily done, originally in \ref{range 0.3}, $G$ is assumed to be torsion free and here
 % it is not----I have suggested make some relevant changes in \ref{range 0.3}. Suggest: remove the above and replaced by the foliowing.}}

 { Now   let $((G,G_+,u),K,\Delta,r)$ be the one given in
\ref{range  0.1} in general, i.e., $G$ may have torsion and $K$ may not be zero. Let $G^1\subset \Aff(\DT)$ be a dense subgroup
with at least three $\Q$-linearly independent elements and
$H=G\oplus G^1$ be as in \ref{range 0.3}. The order unit $u\in G_+$ could be regarded as $(u,0)\in G\oplus G^1=H$
as the order
 unit of $H_+$. Note that  $H$ has torsion ${\rm Tor}(G)={\rm Tor}(H)$ here.}

Then we have splitting short exact sequence
$$0\lr G\lr H\lr H/G(=G')\lr 0$$
with $H/{\rm Tor}(H)$ dimension group.
%By \ref{range 0.15a} applied
By applying \ref{range 0.15a}
%{\color{Green} (maybe it is helpful to point out that in \ref{range 0.15a}, if $H$ and $H/E$ are ordered group, and the order of $H$ is determined by the quotient map, then \ref{range 0.15a} actually gives an decomposition of $H$ as a ordered group)}  {\red{I did not mention this for the following reasons, 1)  when I apply it in below, I am not quite using the order that $E_n\oplus H_n$ is completely decided by the summand $H_n$, 2) not to important reason is that for such to hold, one needs the limit order group to be simple}}
to the short exact sequence $0 \to {\rm Tor}(H) \to H \to H/{\rm Tor}(H) \to 0$, we can write $H$ as inductive
limit of finitely generated abelian groups
\begin{displaymath}
    \xymatrix{
        H_1 \ar[r]^{\gm_{_{1,2}}} & H_2 \ar[r]^{\gm_{_{2,3}}}&H_3 \ar[r]^{\gm_{_{3,4}}} &\cd \ar[r] & H,}
\end{displaymath}
where $H_n=\oplus_{j=1}^{p_n} H_n^j$ with $H_n^1=\Z\oplus {\rm Tor}(H_n)$
and $H_n^i=\Z$ for all $i\geq 2.$
% (here we put the torsion part into
%first summand {\bf What do you mean by that? It looks the torsion in the second,}{ Answer: I mean $H_n^i$ is torsion free for $i\geq 2$ only $H_n^1$ has torsion--Gong}).
Presumingly the positive cone should be given by
$(H_n)_+\setminus \{0\}=
(\Z_+^{p_n}\setminus\{\underbrace{0,...,0}_{p_n}\}) \oplus
{\rm Tor}(H_n)$.
%But we change it a to a smaller positive cone with
But we change it to a smaller cone with
$$(H_n)_+=\big((\Z_+\setminus \{0\}\oplus {\rm Tor}(H_n)\cup \{0,0\}\big)\oplus \Z_+^{p_n-1}$$
with the order unit $u_n=\big(([n,1],\tau_n),[n,2],...,[n,p_n]\big)$,
where $\tau_n\in {\rm Tor}(H_n)$ and $[n,i]$ are positive integers.
Evidently, by simplicity of ordered group $H$, this modification
will not change the positive cone of the limit---in fact, since all the entries of  the map $\cc: \Z^{p_n} \to \Z^{p_{n+1}}$ are strictly positive,
any element in $(\Z_+^{p_n}\setminus\{\underbrace{0,...,0}_{p_n}\}) \oplus
{\rm Tor}(H_n)$ will be sent into $(\Z_+\setminus \{0\}\oplus {\rm Tor}(H_{n+1})\oplus (\Z_+\setminus \{0\})^{p_{n+1}-1}$ which is a subset of $H_{n+1}^+$, by $\gm_{n,n+1}$.  The map
$\gm_{_{n,n+1}}:~ H_n\to H_{n+1}$ could be represented by the
$p_{n+1}\times p_n$ matrix of homomorphism $\td \cc=(\td c_{ij})$,
where $\td c_{ij}:~ H_n^j\to H_{n+1}^i$.  If $i>1,~ j\geq 1$, then
each $\td c_{ij}=c_{ij}$ is a positive integer which defines a map
$\Z(=H_n^j)\to \Z(=H_{n+1}^j)$ by sending $m$ to $c_{ij}m$ (for
$j>1$) or a map $\Z\oplus {\rm Tor}(H_{n})\to \Z$ by sending
$(m,t)\in\Z\oplus {\rm Tor}(H_{n})$ to $c_{ij}m$ (note that the $i$-th
component of $\gm_{_{n,n+1}}({\rm Tor}(H_n))$ is $0$ if $i>1$). If $i=1$,
then $\td c_{ij}=c_{1j}+T_j$, where $T_j:~H_n^j \to {\rm Tor}(H_{n+1})\sbs
H_{n+1}^1$ and $c_{1j}$ is a positive integer. Since
$\gm_{_{n,n+1}}$ satisfies $\gm_{_{n,n+1}}({\rm Tor}(H_n))\sbs
{\rm Tor}(H_{n+1})$, it induces the  map
$$\gm'_{_{n,n+1}}:~ H_n/{\rm Tor}(H_n)=\Z^{p_n}\lr H_{n+1}/{\rm Tor}(H_{n+1})=\Z^{p_{n+1}}.$$
Then $\gm'_{_{n,n+1}}$ is given by the matrix $\cc=(c_{ij})$ of all
entries positive integers.  By passing to a subsequence, we can
require $c_{ij}$ to be larger than any previously given {number}
which
only depends on the construction up to $n^{th}$
step---we will specify how large $c_{ij}$ should be in \ref{condition2} below (see $\spddd_1$ there)\\

\end{NN}

\begin{NN}\label{range 0.17}

Let $G_n=H_n\cap \gm^{-1}_{_{n,\infty}}(G)$ with $(G_n)_+=(H_n)_+\cap
G_n$.  Then the unit  $u_n\in (H_n)_+$ is also the unit for $(G_n)_+$,
since $u\in G\subset H$.

From ${\rm Tor}(G)={\rm Tor}(H)$ we have
${\rm Tor}(G_n)={\rm Tor}(H_n)$. Furthermore we have the following commutative diagram:
\begin{displaymath}
    \xymatrix{0\ar[d] &0\ar[d]&&0\ar[d]\\
        G_1 \ar[r]^{\gm_{_{12}}|_{_{G_1}}} \ar[d] & G_2 \ar[r]\ar[d]&\cd \ar[r]&G \ar[d]\\
         H_1 \ar[r]^{\gm_{_{12}}} \ar[d] & H_2 \ar[r]\ar[d]&\cd \ar[r]&H \ar[d]\\
         H_1/G_1 \ar[r]^{\td\gm_{_{12}}} \ar[d] & H_2/G_2 \ar[r]\ar[d]&\cd \ar[r]&H/G \ar[d]\\
        0 &0&&0 }
\end{displaymath}
where $\td\gm_{_{n,n+1}}$ is induced by $\gm_{_{n,n+1}}$.  Note that
the inductive limit of quotient groups $H_1/G_1\to H_2/G_2\to \cd
\to H/G$ has no apparent order structure.
%{\bf Is this important?} { not important, but clear some possible confusion---Gong}

 Also write the group $K$ (in
\ref{range 0.16}) as inductive limit
\begin{displaymath}
    \xymatrix{
        K_1 \ar[r]^{\chi_{_{12}}}  & K_2 \ar[r]^{\chi_{_{23}}}&K_3 \ar[r]^{\chi_{_{34}}}&\cd \ar[r]& K ~,}
\end{displaymath}
where each $K_n$ is finitely generated.\\

\end{NN}

\begin{NN}\label{range 0.18}
~~Recall from \cite{EG-RR0AH}  that the finite CW complexes $T_{2,k}$  (or
$T_{3,k}$) is defined to be a 2-dimensional connected finite CW
complex with $H^2(T_{2,k})=\Z/k$ and
 $H^1(T_{2,k})=0$ (or 3-dimensional finite CW complex with $H^3(T_{3,k})=\Z/k$
 and $H^1(T_{3,k})=0=H^2(T_{3,k})$). (In {\cite{EG-RR0AH}} the spaces are denoted
 by $T_{II,k}$ and $T_{III,k}$.)  For each $n$, there is a space $X_n'$ which is
 of form
$$X_n'=S^1\vee S^1\vee\cd \vee S^1\vee T_{2,k_1} \vee T_{2,k_2} \vee \cd \vee T_{2,k_i} \vee T_{3,m_1} \vee T_{3,m_2}\vee \cd \vee T_{3,m_j}$$
with $K_0(C(X_1'))=H^1_n=\Z\oplus {\rm Tor}(H_n)$ and $K_1(C(X_1'))=K_n$. Let $x_n$ be the base point of $X_n'$ which is a common point of all spaces $S^1$, $T_{2,k}$, $T_{3,k}$ appeared above in the wedge $\vee$.
And there is a projection $P_n\in M_\infty (C(X_n))$ such that
$$[P_n]=([n,1],\tau_n)\in K_0(C(X_n))=\Z\oplus {\rm Tor}(K_0(C(X_n)))$$
where $([n,1],\tau_n)$ is the first component of unit $u_n\in H_n.$
{ We may assume that $P_n(x_n)=\e_{M_{[n, 1]}}$.}
%$ {\bf Is $\tau_n.$
%specified?} { No, just some element in the torsion part---Gong }.

Note that ${\rm rank}(P_n)=[n,1]$.
Assume  $P_n(x_n)=\e_{M_{[n,1]}(\C)}$,
where $M_{[n,1]}(\C)$ is identified with upper left corner of
$M_{\infty}(\C)$. Define $X_n=[0,1]\vee X_n'$ with $1\in [0,1]$
identified with the base point $x_n\in X_n'$.  We name $0\in [0,1],$
by the symbol $ \tht_1$. So $[0,1]$ is identified with $[\tht_1,\tht_1+1]$. Under this identification,   we have
$X_n=[\tht_1,\tht_1+1]\vee X_n'$.  $P_n\in M_{\infty}(C(X_n'))$ can
be extended to a projection, still called $P_n\in
M_{\infty}(C(X_n))$ by $P_n(\tht_1+t)=\e_{[n,1]}(\C)$ for each
$t\in(0,1)$. Now we will also call $\tht_1$ the base point of $X_n$. The old base point of $X_n'$ is $\tht_1+1$.

Let
\beq\label{Def-F_n-inN0}
F_n=P_nM_{\infty}(C(X_n))P_n\oplus \bigoplus_{i=2}^{p_n} M_{[n,i]}(\C)\andeqn {\rm let}\\
%\eneq
%This is  similar  to the definition of $F_n$ in \ref{range 0.5} except replacing $M_{[n,1]}(\C)$ by $P_n M_{\infty}(C(X_n))P_n$
%with ${\rm rank}(P_n)=[n,1]$.
%Let
%$$
J_n=\{f\in P_n M_{\infty}(C(X_n))P_n:~ f(\tht_1)=0\}.
\eneq
%$$
Then $J_n$ is an ideal of $F_n$ with
$\hat{F_n}=F_n/J_n=\oplus_{i=1}^{p_n} M_{[n,i]}(\C)$.

Let us denote the spectral of $\hat{F}_n=
\oplus_{i=1}^{p_n}\hat{F}_n^i$ by $\tht_1,\tht_2,...,\tht_{p_n}$.
Then the map $\pi: F_n\to \hat{F}_n$ is given by
$\pi(g,a_2,a_3,...,a_{p_n})=(g(\tht_1),a_2,a_3,...,a_{p_n})$
where $g\in P_n M_{\infty}(C(X_n))P_n$.\\

\end{NN}

\begin{NN}\label{range 0.19}
~~ The map
$$H_n/{\rm Tor}(H_n)~\big(=\Z^{p_n}\big)\lr H_n/G_n~\big(=\Z^{l_n}\big)$$
(induced by the quotient map $H_n\to H_n/G_n$; here we use the fact
that ${\rm Tor}(H_n)={\rm Tor}(G_n)$) can be realized as a difference of two
maps
$$\bb_0, \bb_1:~ \Z^{p_n}\to \Z^{l_n}$$
corresponding to two $l_n\times p_n$ matrices of strictly positive
integer entries $\bb_0=(b_{0,ij}),~\bb_1=(b_{1,ij})$.
%{\bf It looks that $H_n/G_n$ does have an order now. Do we ignore it?} { The map involved here is $\bb_1-\bb_0$ which is not positive, so $H_n/G_n$ still have No order structure, say even you are confused, it seems important to emphasis it has no order there---Gong .}

 Exactly as in \ref{range 0.5} (in which we did  the special
case of torsion free $K_0$ and trivial $K_1$), we can define
\beq\label{150102-1}
\{n,i\}:= \sum_{j=1}^{p_n}b_{0,ij}[n,j]=\sum_{j=1}^{p_n}b_{1,ij}[n,j]~.
\eneq

Let $E_n=\bigoplus_{i=1}^{l_n}M_{\{n,j\}}(\C)$ and let $\bt_0,\bt_1:~
\hat{F}_n\to E_n$ be any homomorphisms
%with
%$$(\bt_0)_{*0},~(\bt_1)_{*0}:~ K_0(\hat{F}_n)\big(\Z^{p_n}\big)
%\lr K_0(E_n)\big(\Z^{l_n}\big)$$ satisfying t
such that $(\bt_0)_{*0}=\bb_0$
and
$(\bt_1)_{*0}=\bb_1$.  Now we define:
$$A_n=\big\{(f,g)\in C([0,1],E_n)\oplus F_n;~ f(0)=\bt_0(\pi(g)),
f(1)=\bt_1(\pi(g))\big\},$$
where {$F_n=P_nM_{\infty}(C(X_n))P_n\oplus \bigoplus_{i=2}^{p_n} M_{[n,i]}(\C)$ as defined in (\ref{Def-F_n-inN0})}
and recall that from \ref{range 0.18}, the map $\pi: F_n\to \hat{F}_n$ is given by
$$\pi(g,a_2,a_3,...,a_{p_n})=(g(\tht_1),a_2,a_3,...,a_{p_n}),$$
where $g\in P_n M_{\infty}(C(X_n))P_n$.

%{\bf As it stands, it is not even well-defined. Note $\bt_0,\bt_1$ are not defined on
%$F_n.$ But I assume that it factors through $\hat{F}.$ We need to be more precise.}{ I think it is very precise, $\pi(g)$ is in ${\hat F}_n$ not $F_n$}

%As same as we already discussed for
Similar to the construction of $C_n$ (see \ref{conditions}), the algebra $A_n$ depends only on $F_n$, $\bb_0,\bb_1: K_0({\hat F}_n)(=K_0(F_n)/\mathrm{Tor}(K_0(F_n)))=\Z^{p_n} \to \Z^{l_n}$.  That is, once $F_n, \bb_0,
\bb_1$ are specified, the construction of the algebra $A_n$ is considered to be done.

Note that from (\ref{150102-1}), $\bt_0,\bt_1$ are
unital homomorphisms and therefore $A_n$ is a unital algebra. Note
that this algebra, in general,  is NOT a direct sum of a homogeneous
algebra and an algebra in ${\cal C}_0.$ {In fact, $A_n\in {\cal D}_2.$}
%In fact the base point of the
%spectra of homogeneous part $F_n^1$ is involved in the definition of
%the boundary of function $f\in C([0,1], E_n)$ above.
Later we will
deal with a nicer special case so that $A_n$ is a direct sum of a  homogeneous \CA\,  and
a \CA\, in ${\cal C}_0.$
%the homogenous part
%and ${\cal C}_0$ part can be separate ? in the definition of $A_n$   \\
\end{NN}

%\pagebreak

\begin{NN}\label{construction} Let $I_n=\{(f,g)\in A_n~|~~ g=0\}$. Then $A_n/I_n =F_n$, and we denote the quotient map $A_n\to A_n/I_n=F_n$ by  $\pi_2$.
Let $ J'_n=\{(f, g)\in A_n~|~~ g\in J_n\subset F_n\}=\pi_2^{-1}(J_n)$. Denote the quotient
algebra $A_n/ J'_n$ by $\bar A_n$ and the quotient map $A_n\to A_n/ J'_n=\bar A_n$ by $\pi_1$. Note that
\begin{equation}\label{construction-1}
\hspace{-0.05in}
{\bar A_n}=A(\hat F_n, E_n, \bt_0,\bt_1)=\big\{(f,g)\in C([0,1],E_n)\oplus \hat F_n;~ f(0)=\bt_0((g))~f(1)=\bt_1(\pi(g))\big\}.\hspace{-0.1in}
\end{equation}

The homomorphism $\phi_{n,n+1}: A_n\to A_{n+1}$ to be constructed should satisfy the conditions $\phi_{n,n+1}(I_n)\subset I_{n+1}$ and $\phi_{n,n+1}( J'_n)\subset J'_{n+1}$, and therefore it induces two homomorphism $\psi_{n,n+1}:A_n/I_n=F_n\to A_{n+1}/I_{n+1}=F_{n+1}$ and $\bar \phi_{n,n+1}: A_n/J'_n=\bar A_n\to A_{n+1}/J'_{n+1}=\bar A_{n+1}$. Let $\bar I_n\triangleq\pi_1(I_n)$. Conversely, if two homomorphisms $\psi_{n,n+1}: F_n\to F_{n+1}$ and $\bar \phi_{n,n+1}: \bar A_n\to \bar A_{n+1}$ satisfy the following conditions:
\begin{enumerate}
\item[(a)] $\psi_{n,n+1}(J_n)\subset J_{n+1}$ and $\bar \phi_{n,n+1} (\bar I_n)\subset \bar I_{n+1}$, and
\item[(b)] the homomorphism $F_n/J_n=\hat F_n \to F_{n+1}/J_{n+1}=\hat F_{n+1}$ induced by $\psi_{n,n+1}$ and the homomorphism $\bar A_n/\bar I_n=\hat F_n \to \bar A_{n+1}/\bar I_{n+1}=\hat F_{n+1}$ induced by $\bar \phi_{n,n+1}$ are the same,
\end{enumerate}
then, there is a unique homomorphism $\phi_{n,n+1}: A_n\to A_{n+1}$ satisfying $\phi_{n,n+1}(I_n)\subset I_{n+1}$, $\phi_{n,n+1}( J'_n)\subset J'_{n+1}$, and it induces the homomorphisms  $\psi_{n,n+1}$ and $\bar \phi_{n,n+1}$. Therefore to construct $\phi_{n,n+1}$, we only needs to construct $\psi_{n,n+1}$ and $\bar \phi_{n,n+1}$ separately with the same restriction map $\psi_{n,n+1}|_{Sp \hat F_n}=\bar \phi_{n,n+1}|_{Sp \hat F_n}: \hat F_n\to \hat F_{n+1}$

Note that $\bar A_n$ is as same as $C_n$ in \ref{range 0.5}--\ref{range 0.12} with $\hat F_n$ and $\bar I_n=\pi_1(I_n)$ in place of $F_n$ and $I_n$ in \ref{range 0.5}--\ref{range 0.12}. Therefore the construction of $\bar \phi_{n,n+1}$ can be carried out as  in \ref{range 0.5}--\ref{range 0.12},  with the map $\hat F_n \to \hat F_{n+1}$ being given by the matrix $\cc$ as in \ref{range 0.16} above---of course, we need to assume that the corresponding maps $\bb_0$, $\bb_1$ in this case (see \ref{range 0.19}) satisfy $\spdd$. So, in what follows, we will focus on the construction of $\psi_{n,n+1}$. But before the construction, we will specify the conditions to be used in the future.

\end{NN}

\begin{NN}\label{condition2} The construction of $A_{n+1}$ and $\phi_{n,n+1}$ will be done by induction. Suppose that, we already have the part of the inductive sequence:
\begin{displaymath}
    \xymatrix{
        A_1 \ar[r]^{\phi_{_{1,2}}} & A_2 \ar[r]^{\phi_{_{2,3}}}&A_3 \ar[r]^{\phi_{_{3,4}}} &\cd \ar[r]^{\phi_{n-1,n}} & A_n,}
\end{displaymath}
satisfies the following conditions: for each $i=1,2,\cd n-1$,
\begin{enumerate}
\item[(a)] $\phi_{i,i+1}(I_i)\subset I_{i+1}$ and $\phi_{i,i+1}(J'_i)\subset J'_{i+1}$ and therefore induces two homomorphisms $\psi_{i,i+1}: F_i\to F_{i+1}$ and $\bar \phi_{i,i+1}: \bar A_i \to \bar A_{i+1}$;
\item[(b)] all the homomorphism $\phi_{i,i+1}$, $\psi_{i,i+1}$ and $\bar \phi_{i,i+1}$ are injective, in particular $\psi_{i,i+1}|_{J_i}: J_i\to J_{i+1}$ is injective;
\item[(c)] the induced map $\bar \phi_{i,i+1}:\bar A_i \to \bar A_{i+1}$ satisfies all the conditions (1)--(6) of \ref{range 0.6} with $i$ in place of $n$, with $\bar A_i$  in place of $C_n$ (and of course naturally with $\bar A_{i+1}$ in place of $C_{n+1}$), and with $G_i/\mathrm{Tor}(G_i)$ (or $G_{i+1}/ \mathrm{Tor}(G_{i+1})$) and $H_i/\mathrm{Tor}(H_i)$ (or $H_{i+1}/\mathrm{Tor}(H_{i+1})$) in place of $G_n$ (or $G_{n+1}$) and $H_n$ (or $H_{n+1}$), respectively;
\item[(d)] the matrices $\bb'_0$ and $\bb'_1$ for each $A_{i+1}$ satisfy the condition $\spdd$ in (\ref{13spdd}) with $i+1$ in place of $n+1$ (this condition will be strengthened below).
\end{enumerate}

Now we need to construct $A_{n+1}$ and late on the homomorphism $\phi_{n,n+1}$.

Choose a finite set $Y\sbs X_n\setminus \{\tht_1\}$ (where $\tht_1$ is the base point of $X_n$) satisfying that for each $i < n,~
\bigcup_{y\in Y}Sp(\phi_{_{i,n}}|_{_y})$ is $\frac1n$-dense in $X_i$.
This can be done since the corresponding map $\psi_{_{i,i+1}}|_{_{J_i}}:~J_i (\subset F_i^1)\to J_{i+1} (\subset F_{i+1}^1)$ is injective for
each $i.$
%and $\pi_1\circ \phi_{_{i,i+1}}=\pi_1\circ \xi_{_{i,i+1}}$
%(see (\ref{150102-p2}) in \ref{range 0.29}).  
Let us denote $t\in(0,1)_j\subset
Sp\big(C([0,1],E_n^j)\big) $ by $t_{(j)}$ to distinguish
%spectral
spectrum
from different direct summands of $C([0,1], E_n)$.
%It is convenient to denote $0\in [0,1]_j$ by
%$0_j$ and $1\in [0,1]_j$ by $1_j$ Note that $0_j$ and $1_j$ do not
%correspond to irreducible representation.  In fact, $0_j$
%corresponds to the direct sum of irreducible
%representation
%representations for the
%set
%$$\left\{ \tht_1^{\sim b_{0,j1}}, \tht_2^{\sim b_{0,j2}},\cd, \tht_n^{\sim b_{0,jp_n}}\right\}$$
%and $1_j$,
%and $1_j$ to the set
%$$\left\{ \tht_1^{\sim b_{1,j1}}, \tht_2^{\sim b_{1,j2}},\cd, \tht_n^{\sim b_{1,jp_n}}\right\}.$$
Let $T\subset Sp(A_n)$ be defined by
$$T=\left\{ (\frac kn)_{(j)};~ j=1,2,\cd,l_n;~ k=1,2,\cd, n-1 \right\}.$$
%We need to modify $\phi_{_{n,n+1}}$ to $\xi_{_{n,n+1}}$ such that
% $$SP\big(\xi_{_{n,n+1}}|_{_{\tht_2'}}\big) \sps T\cup Y.$$
 Let $Y=\{y_1,y_2,... , y_{_{L_1}}\}\subset X_n$ and 
  let $L=l_n\cdot (n-1)+L_1=\#(T\cup Y)$.

Set $A_n=\big\{(f,g)\in C([0,1],E_n)\oplus F_n;~ f(0)=\bt_0(\pi(g)),
f(1)=\bt_1(\pi(g))\big\},$ with $(\bt_0)_{*0}=\bb_0=(b_{0,ij})$
and
$(\bt_1)_{*0}=\bb_1=(b_{1,ij})$.
Put $M=\max\{ b_{0,ij};~ i=1,2,\cd,p_n;~ j=1,2,\cd, l_n\}$.

 First we need that the matrix $\cc=(c_{ij})$ given  by the map $\gm'_{_{n,n+1}}:~ H_n/{\rm Tor}(H_n)=\Z^{p_n}\lr H_{n+1}/{\rm Tor}(H_{n+1})=\Z^{p_{n+1}}$ satisfies
 $$
\spddd_1 \qq\qq\qq\qq\qq\qq\qq c_{ij}> 13 \cdot 2^{2n}\cdot ML\qq\,\,\,
\rforal \,\,i,j.\qq\qq\qq\qq$$
This can be done by increasing  $n+1$ to large enough $m$, and renaming $H_m$ as $H_{n+1}$ and $\gm'_{_{n,m}}$ as $\gm'_{_{n,n+1}}$ (see the end of \ref{range 0.16}).

 Let $A_{n+1}=\big\{(f,g)\in C([0,1],E_{n+1})\oplus F_{n+1}:~ f(0)=\bt'_0(\pi(g)),
f(1)=\bt'_1(\pi(g))\big\}$  be constructed as in \ref{range 0.19}. In this construction, we will choose $\bb_0'=(b_{0,ij}')$ and $\bb_1'=(b_{1,ij}')$ (corresponding to $\bt'_0$ and $\bt'_1$) to  satisfy the
  strengthened condition (stronger than $\spdd$)
$$\spddd_2 ~~~
\left\{ \begin{array}{l}
\td b_{0,il}=\sum_{k=1}^{p_{n+1}} b'_{0,ik}\cdot
c_{kl}> 2^{2n}\left(\sum_{k=1}^{l_{n}}(|d_{ik}|+2)\cdot \{n,k\} \right. \\
\left. \hskip 5cm + (L_1+(n-1)\sum_{k=1}^{l_n}b_{0,kl})\cdot b'_{0,i2}
\right) \andeqn \\
%\\
%$$\spddd_2 \qq
%%\qq\qq\qq
\td b_{1,il}=\sum_{k=1}^{p_{n+1}} b'_{1,ik}\cdot c_{kl}>
2^{2n}\left(\sum_{k=1}^{l_{n}}(|d_{ik}|+2)\cdot \{n,k\}\right. \\
\left. \hskip 5cm +(L_1+(n-1)\sum_{k=1}^{l_n}b_{0,kl})\cdot b'_{1,i2} \right),
\end{array}
\right.
$$
 for
each $l\in \{1,2,..., p_n\}$ and $i\in \{1,2,..., l_{n+1}\}$, where $\td \bb_0=\bb'_0\cdot \cc=(\td b_{0,ij})$ and $\td
\bb_1=\bb'_1\cdot \cc=(\td b_{1,ij})$. To
make the above inequalities hold,  we only need to  increase  the
entries of the third columns of $\bb'_0$ and $\bb'_1$ by adding same
big positive number to make it much larger than the second column of
the matrices (see the end of \ref{range 0.5}).

With the above choice of $A_{n+1}$ satisfying  the conditions $\spddd_1$ and $\spddd_2$, we will begin the construction of $\phi_{n,n+1}$ by constructing $\psi_{n,n+1}$ first.
\end{NN}

\begin{NN}\label{range 0.20}
~~ Recall that $K_0(F_n)=(H_n, u_n)$,
$K_0(F_{n+1})=(H_{n+1},u_{n+1})$ and the map $\td \cc=(\td
c_{ij}):~ H_n\to H_{n+1}$ is as in \ref{range 0.16}.  Assume that
$c_{ij}>13$ for any $ij$ (which is a consequence of $\spddd_1$).  We shall define the unital  homomorphism
$\psi_{_{n,n+1}}:~ F_n\to F_{n+1}$ to satisfy the following
conditions:
%\begin{enumerate}
%\item[(1)] 
 \beq\nonumber
&&\hspace{-1.25in}(1)\hspace{0.1in} (\psi_{_{n,n+1}})_{*0}=\gm_{_{n,n+1}}: K_0(F_n)(=H_n)\lr K_0(F_{n+1})(=H_{n+1})\andeqn\\
\nonumber
&& (\psi_{_{n,n+1}})_{*1}=\chi_{_{n,n+1}}:~ K_1(F_n)(=K_n)\lr K_1(F_{n+1})(=K_{n+1});
 \eneq
%\item[(2)]  $
(2)\hspace{0.2in} $\psi_{_{n,n+1}}(J_n)\subset J_{n+1}$, and the map $\hat{\psi}_{_{n,n+1}}:~{\hat F}_n\to {\hat F}_{n+1} $, induced by $\psi_{_{n,n+1}}$  satisfies $$(\hat{\psi}_{_{n,n+1}})_{*0}=\cc=(c_{ij}):
 K_0({\hat F}_n)(=\Z^{p_n})\to ({\hat F}_{n+1})_{*0}(=\Z^{p_{n+1}}).$$
%\end{enumerate}

     Since a homomorphism from a finite dimensional C*-algebra to another 
     %finite dimensional C*-algebra 
     is completely determined by its K-theory map up to unitary equivalence, the map $\hat{\psi}_{_{n,n+1}}:~{\hat F}_n\to {\hat F}_{n+1} $ is
     completely determined by the condition $(\hat{\psi}_{_{n,n+1}})_{*0}=\cc=(c_{ij})$
     up to unitary equivalence and such a map is necessarily unital (since K-theory map $\cc$ keeps the order unit).

     If $i\geq 2$ and $j\geq 2$, then $F_n^i=\hat F_n^i$ and $F_{n+1}^j=\hat F_{n+1}^j$, and, consequently, $\psi_{n,n+1}^{i,j}$ is $\hat{\psi}_{n,n+1}^{i,j}$.
    Since $\psi_{_{n,n+1}}(J_n)\subset J_{n+1}$ and $J_{n+1}\cap F_{n+1}^j=0$ for any $j\geq 2$, one has that $\psi_{n,n+1}^{1,j}(J_n)=0$.  Therefore each component  $\psi_{_{n,n+1}}^{1,j}: F_n^1\to F_{n+1}^j$ for $j\geq 2$
     (note that $\hat{F}_{n+1}^j=F_{n+1}^j$ for $j\geq 2$) factors through $\hat{F}_n^1=F_n^1/J_n$ and is completely
     determined by $\hat{\psi}_{_{n,n+1}}$.  {Note that}, for any $i$ and any $j\geq 2$, $\psi_{n,n+1}^{i,j}$  is completely
     determined by $\hat{\psi}_{n,n+1}^{i,j}$.

      Therefore
       one only needs to define $\psi_{_{n,n+1}}^{-,1}: F_n\to F_{n+1}^1,$  the map
       from $F_n$ to $F_{n+1},$ then project to the first summand.

  %     { Something is missing}

     First one can find mutually orthogonal projections $Q_1,Q_2, ... , Q_{p_n}$ such that
     $$Q_1+Q_2+\cdots +Q_{p_n}=P_{n+1}$$ and $\gm_{_{n,n+1}}^{i,1}([\e_{F_n^i}])=[Q_i] \in K_0(P_{n+1}M_{\infty}(C(X_{n+1}))P_{n+1})$.
     %{\bf We should avoid $M_{\infty}$ here---this means another notation for a positive integer, can we use $M_{\spadesuit}$ to replace all $M_{\infty}$? }
    Since $c_{ij}>13$ for all $i,j$---that is ${\rm rank}(Q_i)/{\rm rank} (\e_{F_n^i})=c_{1,i}>13$, by \ref{range 0.14} (and \ref{range 0.15}), there are unital homomorphisms $\psi_{_{n,n+1}}^{i,1}:F_n^i\to Q_iF_{n+1}^1Q_i$ which realize the K-theory map $\gm_{_{n,n+1}}^{i,1}
:K_0(F_n^i)\to K_0(F_{n+1}^1)$ and $\chi_{_{n,n+1}}:~
K_1(F_n^1)(=K_n)\to K_1(F_{n+1}^1)(=K_{n+1})$ (note that
$K_1(F_n^i)=0$ for $i\geq 2$) and satisfy Condition (2).
%{\bf {Is it clear that each $Q_i$ satisfies the condition
%in \ref{range 0.14}? It is too far away from here to figure out "how large is $Q_i$. %Perhaps, just say: since....,
%$Q_i$ is ...., therefore \ref{range 0.14} applies...}}.

\end{NN}

\begin{NN}\label{range 0.21}
~~ For the convenience of later construction, we need the projection $Q_i$ and we need the homomorphism
$\psi_{_{n,n+1}}^{i,1}: F_n^i\to Q_iF_{n+1}^1Q_i$ not only satisfies the conditions in \ref{range 0.20}, but also has the simple form described below.
%{\bf This is not a good notation ---if it is what I think it is}
%
%{ I worked very hard to chose notation, if I had a good one, I definitely already used it, another choice is $\psi_{n,n+1}|_{X_{n+1}}$, if you do not have a suggestion, then just keep this one ---Gong}
%which depends on detailed
%description of the map given below.

To {simplify} the
notation, we denote $\psi_{_{n,n+1}}^{-,1}$ by $\psi:
\bigoplus_{i=1}^{p_n}F_n^i \to F_{n+1}^1$ and denoted by $\psi^i
:F_n^i \to Q_iF_{n+1}^1Q_i$ the corresponding partial map.
%{Not particularly good choice
%in a very long proof.----in the sense many notation has been used but defined previously with a structure, and then
%change to something else locally, and then change back}
%
%{ I don't know whether it will be any better if we do not simplify $\psi_{_{n,n+1}}^{~,1}$ to $\psi$, for me I think present way is better---Gong}
Note that
${\rm rank} (Q_i)=c_{1i}[n,i]$. We may require that the projection $Q_i$ (for any $i\geq2$, the case $i=1$ will be discussed later) satisfy the following condition:

(a)~~ $Q_i|_{[\tht_1', \tht_1'+1]}$ is of the form
$$\diag (\0_{c_{11}[n,1]},\0_{c_{12}[n,2]},..., \0_{c_{1~i-1}[n,i-1]},\e_{c_{1i}[n,i]},\0_{c_{1~i+1}[n,i+1]},..., \0_{c_{1p_n}[n,p_n]}), $$
when $P_{n+1}|_{[\tht_1', \tht_1'+1]}$ is identified with
$\e_{[n+1,1]}\in M_{[n+1,1]}(C[\tht_1', \tht_1'+1]) \subset
M_{\infty}(C[\tht_1', \tht_1'+1])$.\linebreak
(In the above, recall that
$[\tht_1', \tht_1'+1]$ is identified with the interval $[0,1]$ in
$X_{n+1}=[0,1]\vee X_{n+1}'$, and $\tht_1'$ is the first element in
$Sp(\hat{F}_{n+1})=\{\tht_1', \tht_2', ... , \tht_{p_{n+1}}'\}.$ Here
we reserve $\{\tht_1, \tht_2,..., \tht_{p_{n}}\}$ for
$Sp(\hat{F}_{n})$. )

Fix $i\geq 2$ (the case of  $i=1$ will be discussed later). Let
{$\{e_{ij}\}$} be the matrix unit of $F_n^i=M_{[n,i]}(\C)$ and let
$q=\psi^1(e_{11})$. It is well known that
$Q_iM_{\infty}(C(X_{n+1}))Q_i$ can be identified with
$qM_{\infty}(C(X_{n+1}))q\otimes M_{[n,i]}(\C)$ so that the
homomorphism $\psi^i$ is given by
\beq\label{0.21*}\psi^i\big( (a_{ij})\big)= q\otimes (a_{ij}) \in qM_{\infty}(C(X_{n+1}))q\otimes M_{[n,i]}(\C).
\eneq
Note that ${\rm rank}(q)=c_{1i}$ as ${\rm rank}(Q_i)=c_{1i}[n,i]$. Denote
$c_{1i}-1$ by $d$.

We can write $q=q_1+q_2+\cd+q_d+p$, where $q_1,q_2,\cd, q_d$ are
mutually equivalent trivial rank $1$ projections and $p$ is a
(possible nontrivial) rank $1$ projection. Under the identification
$Q_i=q\otimes \e_{M_{[n,i]}}$, we denote that $\hat{q}_j=q_j\otimes
\e_{M_{[n,i]}}$ and $\hat{p}=p\otimes \e_{M_{[n,i]}}$.
 From the definition of $\psi^i$ (see (\ref{0.21*}) above), we know that $\psi^i(F_n^i)$
 commutes
 with $\hat{q}_1,\hat{q}_2, ... , \hat{q}_d,$ and $ \hat{p}$.
  We can further require the homomorphisms $\psi^i$ ($i\geq 2$) satisfy

 (b)~~ the above projections $\hat{q}_1,\hat{q}_2,\cd, \hat{q}_d,$ and $ \hat{p}$ for $\psi^i$ can be chosen such that $\hat{q}_1|_{[\tht_1', \tht_1'+1]}$, $\hat{q}_2|_{[\tht_1', \tht_1'+1]}$ ,..., $\hat{q}_d|_{[\tht_1', \tht_1'+1]},$
   and $ \hat{p}|_{[\tht_1', \tht_1'+1]}$ are diagonal
   matrix with $\e_{[n,i]}$ in the correct place when
    $Q_iM_{\infty}(C[\tht_1',\tht_1'+1])Q_i$ identified with
     $M_{c_{i1}[n,i]}(C[\tht_1',\tht_1'+1])$. That is
$$\hat{q}_j= \diag (\underbrace{\0_{[n,i]},
..., \0_{[n,i]}}_{j-1}, \e_{[n,i]},\0_{[n,i]},..., \0_{[n,i]})\quad
 \mbox{and}\quad  \hat{p}=\diag (\underbrace{\0_{[n,i]},..., \0_{[n,i]} }_{d},\e_{[n,i]}).$$

\end{NN}
%{\bf {It would be more helpful to state exact what property of $\psi_{_{n,n+1}}^{-,1}$ we are talking about at this
%stage. It would be even more instructive as where it would be used.  I could not find a sentence in the next proof
%which looks like: Now we use the property mentioned in \ref{range 0.21}.}}

\begin{lem}\label{range 0.22}
 Let $i\geq 2$ and $\psi^i: F_n^i\to
Q_iF_{n+1}^1Q_i$ be as in  $\ref{range 0.21}$ above. In particular, the projection $Q_i$ satisfies Condition (a) and the homomorphism $\psi^i$ satisfies Condition (b) above. Suppose $m\leq d
=c_{1i}-1$. Let $\LD:  Q_iF_{n+1}^1Q_i \to M_m(Q_iF_{n+1}^1Q_i)$ be the 
amplification defined by $\LD(a)=a\otimes {\bf 1}_{m}.$
%\diag(\underbrace{a,...,a}_m)$. 
There
is a projection $R^i\in M_m(Q_iF_{n+1}^1Q_i)$ satisfying the following
conditions:

{\rm (i)} $R^i$ commutes with $\LD(\psi^i(F_n^i))$.

{\rm (ii)} $R^i(\tht_1')= Q_i(\tht_1')\otimes \left(
                                        \begin{array}{cc}
                                          \e_{m-1} & 0 \\
                                          0 & 0 \\
                                        \end{array}
                                      \right)= \diag(\underbrace{Q_i(\tht_1'),..., Q_i(\tht_1')}_{m-1}, 0))\in M_m(F_{n+1}^1|_{\tht_1'}).$
Consequently, $rank (R^i)= c_{i1}(m-1)[n,i]=(d+1)(m-1)[n,i]$.

 Let $\pi: M_m(Q_iF_{n+1}^1Q_i) \to M_{m}(Q_i(\tht_1')\hat{F}_{n+1}^1Q_i(\tht_1')) $ be the map induced by $\pi: F_{n+1}^1 \to \hat{F}_{n+1}^1.$ Then $\pi$ takes $R^iM_m(Q_i F_{n+1}^1 Q_i)R^i$ onto  $M_{m-1}(Q_i(\tht_1')\hat{F}_{n+1}^1Q_i(\tht_1')) \subset M_{m}(Q_i(\tht_1')\hat{F}_{n+1}^1Q_i(\tht_1'))$ (see {\rm (ii)} above).

{\rm (iii)} There is an inclusion unital homomorphism
 $$ \iota: M_{m-1}(Q_i(\tht_1')\hat{F}_{n+1}^1Q_i(\tht_1')) \hookrightarrow R^iM_m(Q_i F_{n+1}^1 Q_i)R^i$$
 such that $\pi\circ \iota=\id|_{M_{m-1}(Q_i(\tht_1')\hat{F}_{n+1}^1Q_i(\tht_1'))}$ and such that $R^i(\LD (\psi^i(F_n^i)))R^i \subset \mathrm{Image}(\iota)$.

\end{lem}

 \begin{proof} In the proof of this lemma, $i\geq 2$ is fixed. So the notation $q, q_1,q_2,..., q_d,p,r$ and $\LD_1$ are for this fixed $i$ and only kept the meaning in this proof (later on, we will use them for the proof of similar lemma for the case $i=1$ in which we will work on the corner $M_m(Q_1F_{n+1}^1Q_1)$).

 The homomorphism $\LD\circ \psi^i: F_n^i= M_{[n,i]}(\C) \to M_m(Q_iF_{n+1}^1Q_i)$ can be regarded as $\LD_1\otimes \id_{[n,i]}$,
 where $\LD_1: \C \to M_m(qF_{n+1}^1q)$
 %change to is
 is the  unital
 homomorphism given by $\LD_1(c)= c\cdot(q\otimes \e_m)$.

  Note that $q=q_1+q_2+\cd +q_d+p$ with $\{q_i\}$ being mutually equivalent rank
  one projections. Furthermore, $p|_{[\tht_1',\tht_1'+1]}$ is also a trivial rank one
  projection. Let $r\in M_m(qF_{n+1}^1q)= qF_{n+1}^1q\otimes M_m(\C)$ be defined
  as below.
  $$r(\tht_1')=q(\tht_1')\otimes \left(
                         \begin{array}{cc}
                           \e_{m-1} & 0 \\
                           0 & 0 \\
                         \end{array}
                       \right)= (q_1(\tht_1')+q_2(\tht_1')+\cd +q_d(\tht_1')+p(\tht_1'))\otimes \left(
                         \begin{array}{cc}
                           \e_{m-1} & 0 \\
                           0 & 0 \\
                         \end{array}
                       \right),$$
$$r(\tht_1'+1)= (q_1(\tht_1'+1)+q_2(\tht_1'+1)+\cd +q_{m-1}(\tht_1'+1))\otimes \e_m +~~~~~~~~~~~~~~$$
$$~~~~~~~~~~~~~+ (q_m(\tht_1'+1)+\cd +q_{d}(\tht_1'+1)) \otimes \left(
                         \begin{array}{cc}
                           \e_{m-1} & 0 \\
                           0 & 0 \\
                         \end{array}
                       \right),$$
and  for $x\in X_{n+1}'\subset X_{n+1}$,
 $$ r(x)= (q_1(x)+q_2(x)+\cd +q_{m-1}(x))\otimes \e_m + (q_m(x)+\cd +q_{d}(x)) \otimes \left(
                         \begin{array}{cc}
                           \e_{m-1} & 0 \\
                           0 & 0 \\
                         \end{array}
                       \right).$$
In the above, between $\tht_1'$ and $\tht_1'+1$, $r(t)$ can be
defined to be any continuous path connecting the projections
$r(\tht_1')$ and $r(\tht_1'+1)$, both of rank
$(d+1)(m-1)=(m-1)m+(d-m+1)(m-1)$. (Note that all $q_i(t)$ and $p(t)$
are constant on $[\tht_1', \tht_1'+1]$.)

Let $R^i=r\otimes \e_{[n,i]}$, under the identification of
$M_m(Q_iF_{n+1}^1Q_i)$ with $M_m(qF_{n+1}^iq)\otimes \e_{[n,i]}$.
Since $\LD_1: \C \to M_m(qF_{n+1}^iq)$ sends $\C$ to the center of
 $M_m(qF_{n+1}^iq)$, we have that $r$ commutes with $\LD_1(\C)$ and consequently,
 $R^i=r\otimes \e_{[n,i]}$ commutes with $\LD(\psi^i(F_n^i))$ as
 $\LD\circ \psi^i= \LD_1\otimes \id_{[n,i]}$. That is, condition (i) holds.
 Condition (ii) follows from definition of $r(\tht_1')$ and
 $R^i(\tht_1')=r(\tht_1')\otimes \e_{[n,i]}$.

Note that $r\in M_m(qF_{n+1}^1q)= M_m(qM_{\infty}(C(X_{n+1}))q)$ is
a trivial projection of rank $(d+1)(m-1)$ and
$r(\tht_1')=q(\tht_1')\otimes \e_{m-1}$. One identifies
$$r(\tht_1')M_m(q(\tht_1')\hat{F}_{n+1}^1q(\tht_1'))r(\tht_1')=M_{m-1}(q(\tht_1')\hat{F}_{n+1}^1q(\tht_1'))\cong
M_{(d+1)(m-1)}(\C).$$
Let $r_{ij}^0, 1\leq i,j\leq (d+1)(m-1)$ be the matrix units
for $M_{(d+1)(m-1)}$.
Since $r$ is a trivial projection, one
can construct
$r_{ij}\in rM_m(qF_{n+1}^1q)r$, $1\leq i,j\leq (d+1)(m-1)$
 with $r_{ij}(\tht_1')=r_{ij}^0$ serving as matrix units for $M_{(d+1)(m-1)}\subset
 rM_m(qF_{n+1}^1q)r \cong M_{(d+1)(m-1)}(C(X_{n+1}))$.
 Here by matrix units, we mean $r_{ij}r_{kl}= \dt_{jk}r_{il}$
 and $r=\sum_{i=1}^{(d+1)(m-1)} r_{ii}$. We can define
 $$\iota_1:  M_{m-1}(q(\tht_1')\hat{F}_{n+1}^1 q(\tht_1'))
 (\cong r(\tht_1') M_{m}(q(\tht_1')\hat{F}_{n+1}^1 q(\tht_1')) r(\tht_1')) \hookrightarrow rM_m(q F_{n+1}^1 q)r$$ by $\iota_1(r_{ij}^0)=r_{ij}$. Finally define $\iota=\iota_1\otimes \id_{[n,i]}$, using the identification $R^i=r\otimes \e_{[n,i]}$ and $Q_i=q\otimes \e_{[n,i]}$.
 %Evidently, (iii) holds.
Then (iii) follows.
\end{proof}

\begin{NN}\label{range 0.23}
~~Now we continue the discussion of \ref{range 0.21} for
%the
homomorphism
 $\psi^1: F_n^1\to Q_1F_{n+1}^1Q_1$.  We know that ${\rm rank} (Q_1)= c_{11}[n,1]$,
 where $[n,1]={\rm rank} (P_n)$ for $F_n^1=P_nM_{\infty}(C(X_n))P_n$.
 Note that $c_{11} > 13$. Denote $d=c_{11} - 13$. We can assume that the projection  $Q_1$ and the homomorphism $\psi^1: F_n^1 \to Q_1 F_{n+1}^1 Q_1$ satisfy the following conditions
\begin{enumerate}
\item[(a)]
 $Q_1=\e_{d[n,i]}\oplus \td{Q}:=Q'\oplus \td{Q} \in M_{\infty}(C(X_{n+1}))$
  such that $\td{Q}|_{[\tht_1',\tht_1'+1]}=\e_{13[n,i]}$ (but in the lower right
  corner of $Q_1|_{[\tht_1',\tht_1'+1]}$);

\item[(b)]
   $\psi^1: F_n^1 \to Q_1 F_{n+1}^1 Q_1$ is decomposed into two parts
   $\psi^1=\psi_1\oplus \psi_2$ with \\
   $\psi_1: F_n^1\to M_{d[n,i]}(C(X_{n+1}))= Q'M_{\infty}(C(X_{n+1}))Q'$ and  $\psi_2:F_n^1\to \td{Q}M_{\infty}(C(X_{n+1}))\td{Q}$ described below:
   \begin{enumerate}
    \item[(b1)]
   the unital homomorphism  $\psi_1: F_n^1\to M_{d[n,i]}(C(X_{n+1}))= Q'M_{\infty}(C(X_{n+1}))Q'$ is
    defined by $\psi_1(f)= \diag(\underbrace{f(\tht_1),f(\tht_1),..., f(\tht_1)}_d)$ as
    a constant function on $X_{n+1}$;

    \item[(b2)] the unital homomorphism
    $\psi_2:F_n^1\to \td{Q}M_{\infty}(C(X_{n+1}))\td{Q}$ is a
    homomorphism satisfying $(\psi_2)_{*0}= \td{c}_{11}-d= c_{11}-d+T_1$
    (where $T_1: H_n^1(=K_0(F_n^1))\to {\rm Tor}(H_{n+1}) \subset H_{n+1}^1$ is as
    in \ref{range  0.16})
    and $(\psi_2)_{*1}=\chi_{n,n+1}: K_1(F_n)(=K_n)\to K_1(F_{n+1})(=K_{n+1})$
    (such $\psi_2$ exists because of  \ref{range 0.14});

    \item[(b3)] furthermore,
     $\psi_2$ is injective (see \ref{range 0.15}) and the definition
    of $\psi_2|_{[\tht_1',\tht_1'+1]}$ is given as below:
    for $t\in [0,\frac12]$,
     $$\psi_2(f)(\tht_1'+t)=\diag(\underbrace{f(\tht_1),f(\tht_1),...,
     f(\tht_1)}_{13})$$
    and for $t\in [\frac12, 1]$,
    $$\psi_2(f)(\tht_1'+t)=\diag(\underbrace{f(\tht_1+ 2t),f(\tht_1+2t),...,
    f(\tht_1+2t)}_{13}).$$
    Here, $f(\tht_1+s)\in P_n(\tht_1+s)M_{\infty}(\C)P_n(\tht_1+s)$ is
    regarded as an $[n,1]\times[n,1]$ matrix for each $s\in [0,1]$ by
    using the fact $P_n|_{[\tht_1, \tht_1+1]}=\e_{[n,1]}$.
    \end{enumerate}
\end{enumerate}

 Let us remark that in the decomposition
 $\psi^1=\psi_1\oplus\psi_2: F_n^1\to Q_1M_{\infty} (C(X_{n+1}))Q_1\subset F_{n+1}^1$,
 the first part $\psi_1: F_n^1\to M_{d[n,i]}(C(X_{n+1}))(=Q'M_{\infty} (C(X_{n+1}))Q'$
  factors through $\hat{F}_n^1=M_{[n,1]}(\C)$, and the
  restriction $\psi^1|_{[\tht_1',\tht_1'+\frac12]}$
 also factors
  through $\hat{F}_n^1$,
  as $$\psi^1(f)(x)=\diag(\underbrace{f(\tht_1),..., f(\tht_1)}_{d+13})\,\,\,{\rm for\,\,\, any}\,\,\,x\in [\tht_1',\tht_1'+\frac12].$$

  \end{NN}

\begin{lem}\label{range 0.24}  Suppose that $Q_1$ and $\psi^1:F_n^1 \to Q_1F_{n+1}^1Q_1$ satisfy Conditions (a) and (b) (including (b1),(b2) and (b3)).
 Suppose that $13m\leq d=c_{11}-13$. Let $\LD:  Q_1F_{n+1}^1Q_1 \to M_m(Q_1F_{n+1}^1Q_1)$
 be amplification  defined by
$\LD(a)=a\otimes {\bf 1}_m.$
%\diag(\underbrace{a,\cd,a}_m)$. 
There is a projection $R^1\in
M_m(Q_1F_{n+1}^1Q_1)$ satisfying the following conditions:

{\rm (i)} $R^1$ commutes with $\LD(\psi^1(F_n^1))$.

{\rm (ii)} $R^1(\tht_1')= Q_1(\tht_1')\otimes \left(
                                        \begin{array}{cc}
                                          \e_{m-1} & 0 \\
                                          0 & 0 \\
                                        \end{array}
                                      \right)= \diag(\underbrace{Q_1(\tht_1'),\cd, Q_1(\tht_1')}_{m-1}, 0))\in M_m(F_{n+1}^1|_{\tht_1'}).$
Consequently, $rank (R^1)= c_{11}(m-1)[n,1]=(d+13)(m-1)[n,1]$.

 Let $\pi: M_m(Q_1F_{n+1}^1Q_1) \to M_{m}(Q_1(\tht_1')\hat{F}_{n+1}^1Q_1(\tht_1')) $ be induced by $\pi: F_{n+1}^1 \to \hat{F}_{n+1}^1$, then $\pi$ takes $R^1M_m(Q_1 F_{n+1}^1 Q_1)R^1$ onto  $M_{m-1}(Q_1(\tht_1')\hat{F}_{n+1}^1Q_1(\tht_1')) \subset M_{m}(Q_1(\tht_1')\hat{F}_{n+1}^1Q_1(\tht_1'))$ (see {\rm (ii)} above).

{\rm (iii)} There is an inclusion unital homomorphism
 $$ \iota: M_{m-1}(Q_1(\tht_1')\hat{F}_{n+1}^1Q_1(\tht_1')) \hookrightarrow R^1M_m(Q_1 F_{n+1}^1 Q_1)R^1$$ such that $\pi\circ \iota=\id|_{M_{m-1}(Q_1(\tht_1')\hat{F}_{n+1}^1Q_1(\tht_1'))}$ and such that $R^1(\LD (\psi^1(F_n^1)))R^1 \subset
 % \mbox
 {\rm Image} (\iota)$.

\end{lem}

The notation $\LD$, $d$, $m$ in the above lemma, and  $q, q_1,q_2,... q_d, p, r $ and $\LD_1$ in the proof below, are also used in Lemma \ref{range 0.22} and
its proof for the case $i\geq 2$ (comparing with $i=1$ here). Since they are used for the same purpose, we choose the same notation.

 \begin{proof} The map
 \begin{displaymath}
\xymatrix{\psi_1: F_n^1\ar[r]^{\pi}&  \hat{F}_n^1 \ar[r] &
M_{d[n,1]}(C(X_{n+1}))=Q'M_{\infty}(C(X_{n+1}))Q'}
\end{displaymath}
(where $Q'=\e_{d[n,1]}$)
%,
can be written as $(\LD_1\otimes
\id_{[n,1]})\circ \pi$, where $\LD_1: \C \to M_d(C(X_{n+1}))$ is
the map sending $c\in \C$ to $c\cdot\e_d$. We write
$\LD_1(1):=q'=q_1+q_2+\cd+q_d$, with each $q_i$ a trivial constant
projection of rank $1$. Here $q'$ is constant subprojection of $Q'$
with $Q'=q'\otimes \e_{[n,1]}$. Consider the map
$\td{\psi}_2:=\psi_2|_{[\tht_1',\tht_1'+\frac12]}: F_n^1 \to
\td{Q}F_{n+1}^1\td{Q}|_{[\tht_1',\tht_1'+\frac12]}$ and $\td{\psi}^1:=\psi^1|_{[\tht_1',\tht_1'+\frac12]}=(\psi_1+\psi_2)|_{[\tht_1',\tht_1'+\frac12]}: F_n^1 \to
{Q_1}F_{n+1}^1{Q_1}|_{[\tht_1',\tht_1'+\frac12]}$. As
%point
pointed out in
\ref{range 0.23}, $\td{\psi}_2$ has the factorization
\begin{displaymath}
\xymatrix{F_n^1\ar[r]^{\pi} &  \hat{F}_n^1 \ar[r] &
M_{13[n,1]}(C[\tht_1',\tht_1'+\frac12]).}
\end{displaymath}
Hence $\td{\psi}^1$ has the factorization
\begin{displaymath}
\xymatrix{F_n^1\ar[r]^{\pi} &  \hat{F}_n^1 \ar[r] &
M_{(d+13)[n,1]}(C[\tht_1',\tht_1'+\frac12]).}
\end{displaymath}
The map $\td{\psi}^1$ can be written as $(\LD_2\otimes \id_{[n,1]})\circ \pi$,
%for $\LD_2: \C \to M_{d+13}(C[\tht_1',\tht_1'+\frac12])$ to be the
where $\LD_2: \C \to M_{d+13}(C[\tht_1',\tht_1'+\frac12])$ is the map defined by
%map
sending $c\in \C$ to $c\cdot\e_{d+13}$. We write
$\LD_2(1):=q=q_1+q_2+\cd+q_d+ p$ with each $q_i$ being the
restriction of $q_i$ appeared in the definition of $\LD_1(1)$
%onto
on $[\tht_1',\tht_1'+\frac12]$, and $p$ is rank $13$ trivial
projection. Here $q$ is a constant projection on
$[\tht_1',\tht_1'+\frac12]$ and
$Q_1|_{[\tht_1',\tht_1'+\frac12]}=q\otimes \e_{[n,1]}$. Let $r\in
M_m(qF_{n+1}^1q)= qF_{n+1}^1q\otimes M_m(\C)$ be defined as below.
  $$r(\tht_1')=q(\tht_1')\otimes \left(
                         \begin{array}{cc}
                           \e_{m-1} & 0 \\
                           0 & 0 \\
                         \end{array}
                       \right)= (q_1(\tht_1')+q_2(\tht_1')+\cd +q_d(\tht_1')+p(\tht_1'))\otimes \left(
                         \begin{array}{cc}
                           \e_{m-1} & 0 \\
                           0 & 0 \\
                         \end{array}
                       \right);$$
for $t\in [\frac12, 1]$,
$$r(\tht_1'+t)= \Big(q_1(\tht_1'+t)+q_2(\tht_1'+t)+\cd +q_{13(m-1)}(\tht_1'+t)\Big)\otimes \e_m + $$
$$+\Big(q_{13(m-1)+1}(\tht_1'+t)+\cd +q_{d}(\tht_1'+t)\Big) \otimes \left(
                         \begin{array}{cc}
                           \e_{m-1} & 0 \\
                           0 & 0 \\
                         \end{array}
                       \right);$$
and for $ x\in X_{n+1}'\subset X_{n+1}$,
 $$ r(x)= (q_1(x)+q_2(x)+\cd +q_{13(m-1)}(x))\otimes \e_m + (q_{13(m-1)+1}(x)+\cd +q_{d}(x)) \otimes \left(
                         \begin{array}{cc}
                           \e_{m-1} & 0 \\
                           0 & 0 \\
                         \end{array}
                       \right).$$
In the above, between $\tht_1'$ and $\tht_1'+\frac12$, $r(t)$ can
 be defined to be any continuous path connecting the projections
 $r(\tht_1')$ and $r(\tht_1'+\frac12)$, both of rank
 $(d+13)(m-1)=13(m-1)m+(d-13(m-1))(m-1)$. (Note that all
 $q_i(x) $ are constant on $x\in X_{n+1}=[\tht_1',\tht_1'+1]\vee X_{n+1}'$
   and $p(t)$ is constant for $t\in [\tht_1', \tht_1'+1]$.)
   Note that for $x \in [\tht_1'+\frac12,\tht_1'+1]\vee X_{n+1}'$,
    $r(x)$ has same form as $r(\tht_1'+\frac12)$ which is constant
    %sub projection
    sub-projection
    of constant projection $q'\otimes \e_m$. We are going to define
    $R^1$
    to be $r\otimes \e_{[n,i]}$
    under certain identification. Note that the
     projection $Q_1$ is identified with $q\otimes \e_{[n,1]}$
     only on interval $[\tht_1', \tht_1'+1]$
     so the definition of $R^1$ will be
     divided into two parts. For the part on $[\tht_1', \tht_1'+\frac12]$, we
     use the identification of $Q_1$ with $q\otimes \e_{[n,1]}$, and for the
      part
      that $x\in [\tht_1'+\frac12,\tht_1'+1]\vee X_{n+1}'$, we use
       the identification of $Q'=\e_{d[n,1]}$ with
       $q'\otimes \e_m$ (of course, we use the fact
       that $r$ is sub-projection of $q'$
       on this part). This is the only difference between the proof of this lemma
       and that of Lemma \ref{range 0.22}. The definition of
       $ \iota: M_{m-1}(Q_1(\tht_1')\hat{F}_{n+1}^1Q_1(\tht_1')) \hookrightarrow R^1M_m(Q_1 F_{n+1}^1 Q_1)R^1$ and
       verification
       that $\iota$ and $R$ satisfy the conditions are exact as same
       as the proof of \ref{range  0.22}, with $(d+1)(m-1)$ replaced by $(d+13)(m-1)$.
 \end{proof}

 Combining \ref{range  0.22} and \ref{range  0.24}, we have the following
 %Theorem
 theorem which is
  used to conclude the algebra $A$ (will be constructed later) satisfies that
   $A\otimes U$ is in ${\cal B}_0.$

% ----{\bf There is a serious problem with file now.-----this above automatically becomes italian and }
  % tracially approximately Elliott-Thomsen algebra---that is $A\in {\cal N}_0$.

\begin{thm}\label{range 0.25}
  Suppose that $1< m\leq {\min} \{\frac{c_{11}-13}{13}, c_{12}-1, c_{13}-1,..., c_{1p_n}-1\}$. Let $\psi: F_n \to F_{n+1}^1$ be the composition
 \begin{displaymath}
\xymatrix{F_n\ar[r]^{\psi_{n,n+1}}& F_{n+1}  \ar[r]^{\pi_1} &
F_{n+1}^1.}
\end{displaymath}
(The map $\pi_1$ is the quotient map to the first block.) Let $\LD:
F_{n+1}^1 \to M_m(F_{n+1}^1)$ be
%amplify map
%the map defined 
the amplification by
$\LD(a)=a\otimes {\bf1}_m.$
%\diag(\underbrace{a,...,a}_m)$. 
There is a projection $R\in
M_m(F_{n+1}^1)=F_{n+1}^1\otimes M_m(\C)$ and there is an inclusion
unital homomorphism
 $ \iota: M_{m-1}(\hat{F}_{n+1}^1) =\hat{F}_{n+1}^1\otimes M_{m-1}(\C) \hookrightarrow RM_m( F_{n+1}^1 )R$, satisfying the following conditions:

{\rm (i)} $R$ commutes with $\LD(\psi(F_n))$.

{\rm (ii)} $R(\tht_1')= \e_{\hat{F}_{n+1}^1}\otimes \left(
                                        \begin{array}{cc}
                                          \e_{m-1} & 0 \\
                                          0 & 0 \\
                                        \end{array}
                                      \right).$
Consequently, the map $\pi: F_{n+1}^1 \to \hat{F}_{n+1}^1$, takes
$RM_m( F_{n+1}^1 )R$ onto  $M_{m-1}(\hat{F}_{n+1}^1)$ .

{\rm (iii)}  $\pi\circ \iota=\id|_{M_{m-1}(\hat{F}_{n+1}^1)}$,  and $R(\LD
(\psi(F_n)))R \subset \mathrm{Image} (\iota)$.

\end{thm}

\begin{proof} Choose $R=\bigoplus_{i=1}^{p_n} R^i\in
M_m(F_{n+1}^1)$, where $R^1\in
M_m(Q_1F_{n+1}^1Q_1)$ is given in Lemma  \ref{range 0.24} and $R^i\in M_m(Q_iF_{n+1}^1Q_i)$ (for $i\geq 2$) are given in Lemma \ref{range 0.22}.
Then the theorem follows.
\end{proof}

\begin{NN}\label{range 0.26}

 Set $F=\varinjlim (F_n, \psi_{n, n+1})$. Since
 $\psi_{n,n+1}(J_n)\subset J_{n+1}$, this procedure also gives
 %up
 an
 inductive limit of quotient algebras $\hat{F}=\lim (\hat{F}_n, \hat{\psi}_{n,n+1})$,
 where $\hat{F}_n=F_n/J_n$. Evidently, $\hat{F}$ is an AF algebra
 with $K_0(\hat{F})= H/{\rm Tor}(H)$.

 \end{NN}

\begin{thm}\label{range 0.27}
If the matrix $\cc=(c_{ij})$ of $\td{\gm}_{n,n+1}: H_n/{\rm Tor}(H_n) \to
H_{n+1}/ {\rm Tor}(H_{n+1})$ satisfies $c_{ij}> 13\cdot 2^{2n}$ for each
$i,j$ (which is true by $\spddd_1$ in \ref{condition2}), then tracial state space
%$TF$
$T(F)$
of $F$ is $T(\hat{F})=\DT$---that
is, the Elliott invariant of $F$ is
$$
((K_0(F),K_0(F)_+, [{\mbox{\large \bf 1}}_F]),K_1(F),T(F), r_F)\cong
((H,H_+,u),K, \Delta,r).$$

\end{thm}

\begin{proof} Note that $\psi_{n,n+1}$ satisfies (1) in \ref{range 0.20}, and consequently,
$$ ((K_0(F),K_0(F)_+, [{\mbox{\large \bf 1}}_F]),K_1(F))\cong
((H,H_+,u),K, ).$$

%{\bf {A hint why $K$-theory does not need to verify should be given. The reader might also be puzzling
%about comparability. I think you sort of did these in previous lemmas and remarks. But it would be good opportunity
%to recall and pointed out which lemma or remark is relevant to this omission at this point.}}
Let $\pi_n: F_n\to \hat F_n$ be the quotient map. Then $\pi_n^{\sharp}: \Aff(T(F_n))=C(X_n, \R)\oplus \R^{p_n-1} \to \Aff (T (\hat F_n)) =\R^{p_n}$
is given by $\pi_n(g,h_2,h_3,...,h_{p_n})=(g(\tht_1), h_2,h_3,..., h_{p_n})$. Define $\GM_n: \Aff (T (\hat F_n)) =\R^{p_n} \to \Aff (T (F_n))=C(X_n, \R)\oplus \R^{p_n-1}$ to be the right inverse of $\pi_n^{\sharp}$ given by $\GM_n(h_1, h_2,h_3,..., h_{p_n})=(g,h_2,h_3,...,h_{p_n})$ with $g$ being constant function $g(x)=h_1$ for all $x\in X_n$. Then with the condition $c_{ij}> 13\cdot 2^{2n}$, we can prove the following claim:

\noindent{\bf Claim:} For any $f\in \Aff (T(F_n))$ with $\|f\|\leq 1$ and $f'=\psi_{n,n+1}^{\sharp} (f)\in \Aff (T(F_{n+1}))$, we have
$$\|\GM_{n+1}\circ\pi_{n+1}^{\sharp}(f')-f'\|<\frac{2}{2^{2n}}.$$

\noindent{Proof of claim:} Write $f=(g, h_2,..., h_{p_n})$ and $f'=(g', h'_2,..., h'_{p_{n+1}})$.
Then $\GM_{n+1}\circ\pi_{n+1}^{\sharp}(f')=(g'', h'_2,..., h'_{p_{n+1}})$ with $g''(x)=g'(\tht'_1)$ for all $x\in X_{n+1}$.
Recall that $\psi_{n,n+1}^{i,1}$ is denoted by $\psi^i: F_n^i \to Q_iF_{n+1}^1Q_i$ and $\psi^1=\psi_1+\psi_2$ with $\psi_1: F_n^1 \to Q'F_{n+1}^1Q' $  and $\psi_2: F_n^1 \to {\td Q}F_{n+1}^1{\td Q} $ as in \ref{range 0.21} and \ref{range 0.23}. Note that
$$\frac{{\rm rank}(Q_i)}{{\rm rank}(P_{n+1})}=\frac{c_{i,1}}{\sum_{j=1}^{p_n}c_{j,1}},~~~~\frac{{\rm rank}(Q')}{{\rm rank}(P_{n+1})}=\frac{c_{1,1}-13}{\sum_{j=1}^{p_n}c_{j,1}},~~~~{\mbox{and}}~~~\frac{{\rm rank}({\td Q})}{{\rm rank}(P_{n+1})}=\frac{13}{\sum_{j=1}^{p_n}c_{j,1}}.$$
Hence
$$g'=\frac{c_{1,1}-13}{\sum_{j=1}^{p_n}c_{j,1}}\psi_1^\#(g)+
\frac{13}{\sum_{j=1}^{p_n}c_{j,1}}\psi_2^\#(g)+
\sum_{i=2}^{p_n}\frac{c_{i,1}}{\sum_{j=1}^{p_n}c_{j,1}}(\psi^i)^\#(h_i).$$
Also from the construction in \ref{range 0.21} and \ref{range 0.23}, we know that $\psi_1^\#(g)$ and $(\psi^i)^\#(h_i)$ ($i\geq 2$) are constant. So we have
$$|g'(x)-g'(\tht'_1)|\leq \frac{2\times 13}{\sum_{j=1}^{p_n}c_{j,1}} < \frac{2}{2^{2n}},$$
and the claim follows.
%{\bf {Do we need the condition \ref{13spdd} here?  Also, we have somewhat different infinite diagrams, do we?
%It looks $A_n$ is replaced by $F_n$ here. Perhaps it is simpler after the reader was informed which diagrams to look----
%I guess what I mean is to display them  here.
%At that point,  you could  tell the reader what these are and why they becomes simpler.  The more complicated (at least part of them) $\Aff$ spaces also involved. You also need to convince the reader that complexity does not matter.  Add all these, \,
%I would say that things are no longer simpler.-----please also note, at this point, one would complain about changing notation
%of $\LD$.}}

The proof of the theorem is then as same as that of \ref{range  0.12} with \ref{range  0.10} replaced by the above claim. Namely, one can prove the following approximately intertwining diagram replacing the diagram  there.

\begin{displaymath}
    \xymatrix{
        \Aff(T(F_1)) \ar@/_/[d]_{\pi_1^\#}\ar[r]^{\psi_{1,2}^\#}  & \Aff(T(F_2)) \ar[r]^{\psi_{2,3}^\#} \ar@/_/[d]_{\pi_2^\#}& \Aff(T(F_3)) \ar[r] \ar@/_/[d]_{\pi_3^\#}& \cd \Aff(T(F ))\\
        \Aff(T({\hat F}_1)) \ar[r]^{\hat{\psi}_{1,2}^\#}\ar@/_/[u]_{\GM_1} & \Aff(T({\hat F}_2)) \ar[r]^{\hat{\psi}_{2,3}^\#}\ar@/_/[u]_{\GM_2}& \Aff(T({\hat F}_3)) \ar[r]\ar@/_/[u]_{\GM_3}& \cd \Aff(T({\hat F}))~, }
\end{displaymath}
where $\pi_i^\#$ is induced by the quotient map $\pi_i: F_i \to {\hat F}_i$.  Exactly the same as the end of the proof of Lemma \ref{range 0.12}, one can get the isomorphism between $((K_0(F),K_0(F)_+, [{\mbox{\large \bf 1}}_F]),K_1(F),T(F), r_F)$ and $((H,H_+,u),K, \Delta,r)$ including the compatibility (note that the quotient map $H_i=K_0(F_i) \to H_i/\mathrm{Tor}(H_i)=K_0(\hat F_i)$ is also induced by quotient map $\pi_i:F_i\to \hat F_i$).
\end{proof}

\begin{NN}\label{range 0.28}~~
 Let $A_n=\Big\{(f,g)\in C([0,1],E_n)\oplus F_n;~ f(0)=\bt_0\pi(g),
f(1)=\bt_1\pi(g)\Big\}$ {be} as in \ref{range 0.19}, where $\bt_0,\bt_1:
\hat{F}_n \to E_n$ are two unital homomorphisms, and $\pi$ is as in the
end of \ref{range 0.18}.

Let $A_{n+1}=\Big\{(f,g)\in
C([0,1],E_{n+1})\oplus F_{n+1};~ f(0)=\bt'_0\pi(g),
f(1)=\bt'_1\pi(g)\Big\}$ be defined with unital homomorphisms
$\bt'_0,\bt'_1: \hat{F}_{n+1} \to E_{n+1}$. Note that $A_{n+1}$ satisfies the conditions $\spddd_1$ and $\spddd_2$ in \ref{condition2}.

Recall that $J_n\subset F_n^1$
 (and $J_{n+1} \subset F_{n=1}^1$, resp.) is the ideal of elements
 vanishing on $\tht_1\in Sp(\hat{F}_n^1)\subset Sp(F_n^1)$
 (and $\tht_1'\in Sp(\hat{F}_{n+1}^1)\subset Sp(F_{n+1}^1)$ respectively).

Let $I_n$ (or $I_{n+1}$), $J'_n$ (or $J'_{n+1}$), $F_n=A_n/I_n$ (or $F_{n+1}=A_{n+1}/I_{n+1}$),  and  $ \bar{A}_n=A_n/J'_n$ (or
 $\bar{A}_{n+1}=A_{n+1}/J'_{n+1}$)  be as in \ref{construction}.

We still use the matrix $\bb_0=(b_{0,ij})$  (or $\bb_1=(b_{1,ij})$ )
to denote $(\bt_0)_{*0}$ (or $(\bt_1)_{*0}$) as before.  Also
$\bb_0'=(b_{0,ij}')=(\bt_0')_{*0}$,
 $ \bb_1'=(b_{1,ij}')=(\bt_1')_{*0}$ as in \ref{range 0.7}.  We still use
 $\cc=(c_{ij})_{p_{n+1}\times p_n}$ to denote the map
 $H_n/{\rm Tor}(H_n)\to H_{n+1}/{\rm Tor}(H_{n+1})$ and
$\dd=(d_{ij})_{l_{n+1}\times l_n}$ to denote the map
 $H_n/G_n \to H_{n+1}/G_{n+1}$.  It follows that, as in \ref{range  0.6} and
 \ref{range  0.7}, if $\bb_0'$ and $\bb_1'$ satisfy the condition $\spd$ (see \ref{13spd-1}) which is weaker than $\spddd_2$,
 then one can use $\hat\psi_{_{n,n+1}}:~ \hat F_n\to \hat F_{n+1}$ to
 define $\bar \phi_{_{n,n+1}}:~ \bar A_n\to \bar A_{n+1}$ (exactly as $C_n \to C_{n+1}$ there) which
 satisfies $\bar \phi_{_{n,n+1}}(I_n)\subset I_{n+1}$.  From the construction, we know that
 $({\bar \phi}_{_{n,n+1}})_{_{*,0}}=
 \gm'_{_{n,n+1}}|_{G_n/\mathrm{Tor}(G_n)}$ and
\beq\label{0.28a}
 \gm'_{_{n,n+1}}|_{G_n/\mathrm{Tor}(G_n)}: K_0(\bar A_n)=G_n/(\mathrm{Tor}(G_n)) \to  K_0(\bar A_{n+1})=G_{n+1}/(\mathrm{Tor}(G_{n+1})),
 \eneq
 where $\gm'_{n,n+1}: H_n/(\mathrm{Tor}(H_n)) \to  H_{n+1}/(\mathrm{Tor}(H_{n+1}))$ is induce by $\gm_{n,n+1}$ as in \ref{range 0.16}.

 Let
 $\pi_1:  A_n\to  A_n/J_n=\bar A_n$  (or $A_{n+1}\to  \bar A_{n+1}$)
 and
$\pi_2:  A_n\to  A_n/I_n= F_n$  (or $A_{n+1}\to  F_{n+1}$) be the
quotient maps as in \ref{construction}.  Then we can combine the above definition of $\bar
\phi_{_{n,n+1}}:~ \bar A_n\to \bar A_{n+1}$ and $\psi_{_{n,n+1}}:~
F_{n}\to F_{n+1}$ to define $\phi_{_{n,n+1}}:~  A_n\to A_{n+1}$, as below.

Let
$f\in A_n$ and $x\in Sp(A_{n+1})$.  If $x\in Sp(\bar A_{n+1})$, then
\beq\label{0.28b}(\phi_{_{n,n+1}}(f))(x)=\bar\phi_{_{n,n+1}}(\pi_1(f)(x));
\eneq
and if $x \in Sp(F_{n+1})$, then
\beq\label{0.28c}
(\phi_{_{n,n+1}}(f))(x)=\psi_{_{n,n+1}}(\pi_2(f)(x)).
\eneq
Note that for $x\in Sp(F_{n+1})\cap Sp(\bar
A_{n+1})=Sp(\hat{F}_{n+1})$,
$$\bar\phi_{_{n,n+1}}(\pi_1(f)(x))=\psi_{_{n,n+1}}(\pi_2(f)(x))=\hat \psi_{_{n,n+1}}(\pi_1(\pi_2(f))),
$$
where $\pi_1(\pi_2(f))=\pi_2(\pi_1(f))\in \hat F_n$.  By $(\ref{0.28a}),(\ref{0.28b}),(\ref{0.28c})$ above and (1) of \ref{range 0.20}, one has
 $$({\phi}_{_{n,n+1}})_{_{*,0}}=
 \gm_{_{n,n+1}}|_{G_n}~~~~~~\mbox{and}~~~~~({\phi}_{_{n,n+1}})_{_{*,1}}=
 \chi_{_{n,n+1}}.$$
 In this way we
define an inductive limit $$A_1\lr A_2\lr A_3\lr \cd \lr A_n\lr \cd
\lr A$$ with $\big(K_0(A),K_0(A)_+, \e_A \big)=(G,G_+,u)$ and
$K_1(A)=K$.
%{\bf {This is of course right. I also think that you have provided ingredients of the proof for the last
%assertion. Perhaps, it is the time to tell people, what is going to be relevant to this assertion. For example, how \ref{range 0.27} is used here. How other remarks used here, $G$ is (if I am correct) a subgroup of $H$ in \ref{range 0.27},
%how an infinite digram play role here.}}

Similar to \ref{range 0.8} (see $(*)$ and $(**)$ there)
%{\bf {I am afraid that I did not find
%$(*)$ and $(**)$ here}},
$$\Aff(T(A_n)) \sbs \bigoplus_{i=1}^{l_n}C([0,1]_i,\R)\oplus C(X_n, \R)\oplus \R^{p_n-1}$$
{consists} of $(f_1,f_2,... ,f_{l_n},g,h_2,... ,h_{p_n})$  which satisfies
the conditions
\beq\nonumber
\qq\qq\qq\qq f_i(0)&=&\frac1{\{n,i\}}\big(b_{_{0,i1}}g(\tht_1)[n,1]+\sum_{j=2}^{p_n} b_{0,ij}h_j\cdot [n,j]\big)\,\,\,\,\,\,\,\,\,\,\andeqn\qq\qq\,\,\,\,(*)\\\nonumber
%$$
%and
%
%$$
\qq\qq\qq\qq f_i(1)&=&\frac1{\{n,i\}}\big(b_{_{1,i1}}g(\tht_1)[n,1]+\sum_{j=2}^{p_n}b_{1,ij}h_j\cdot [n,j]\big),\qq\qq\qq\qq\,(**)
\eneq

For $h=(h_1,h_2,...,h_{p_n})\in \Aff T(\hat F_n)$, let
$\GM_n'(h)(t)= t\cdot \bt_1^{\sharp}(h)+(1-t)\cdot \bt_0^{\sharp}(h)$ which
gives an element $C([0,1],\R^{l_n})=\bigoplus_{i=1}^{l_n}C([0,1]_i,\R)$.
And let
$$
\GM_n:~ \Aff(T(\hat F_{n}))=\R^{p_n}\to \Aff(T(A_{n})) \subset
\bigoplus_{j=1}^{l_n}C([0,1]_j,\R)\oplus C(X_n, \R)\oplus \R^{p_n-1}
$$
be defined by
$$
\GM_n(h_1,h_2,..., h_{p_n})=(\GM_n'(h_1,h_2,..., h_{p_n}), g,h_2,..., h_{p_n} )\in \bigoplus_{j=1}^{l_n}C([0,1]_j,\R)\oplus C(X_n, \R)\oplus \R^{p_n-1},
$$
where $g\in C(X_n,\R)$ is the constant function $g(x)=h_1$.

If we further assume at each step that the K-group
homomorphisms  $\bb_0'=(\bt_0')_{*0}$ and $\bb_1'=(\bt_1')_{*0}$ satisfy
$(\ref{13spdd})$, and assume that $\cc=(c_{ij})$ satisfies $c_{ij}> 13\cdot 2^{2n}$, then by using   \ref{range 0.10} and the  claim in \ref{range 0.27},
we will the following claim:
%{\bf {The situation seems quite different. It is, for example,
%what should happen on $X_{n+1}$ is not very clear.}}
%{\color{Green} (How about calling it a lemma, because it is proved.)We are not going to prove it here,since it is completely same as \ref{range 0.10} }

\noindent{\rm Claim:} For any $f\in \Aff(T (A_n))$ with $\|f\|\leq 1$ and $f'=\phi_{n,n+1}^{\sharp} (f)\in \Aff (T(A_{n+1}))$, we have
$$\|\GM_{n+1}\circ\pi_{n+1}^{\sharp}(f')-f'\|<\frac{4}{2^{2n}},$$
where $\pi_{n+1}^{\sharp}: \Aff(T (A_{n+1})) \to \Aff(T(\hat F_{n+1}))$ is induced by canonical quotient map $\pi_{n+1}: A_{n+1} \to \hat F_{n+1}$.

\noindent Proof of the claim:
For any $n\in\Z_+$, write
\begin{displaymath}
\xymatrix{\GM_n=\GM_n^1\circ\GM_n^2:\Aff(T({\hat F}_n))\ar[r]^{\qq\qq\GM_n^2}& \Aff(T({ F}_n))  \ar[r]^{\GM_n^1} &
\Aff(T({ A}_n))}
\end{displaymath}
with $\GM_n^2: \Aff(T({\hat F}_n))=\R^{p_n} \to \Aff(T({F}_n))=C(X,\R)\oplus \R^{p_{n-1}}$
defined by $$\GM_n^2(h_1,h_2,..., h_{p_n})=(g,h_2,..., h_{p_n}),$$ where $g$ is the constant function $g(x)=h_1$, and with $\GM_n^1: \Aff(T({F}_n))=C(X,\R)\oplus \R^{p_{n-1}} \to \Aff(T({ A}_n))$ defined by $$\GM_n^1(g,h_2,..., h_{p_n})=(\GM_n'(g(\tht_1),h_2,..., h_{p_n}), g,h_2,..., h_{p_n} ).$$ Also write $\pi_n^\#=(\pi_n^2)^\#\circ (\pi_n^1)^\#$, where $\pi_n^1: A_n\to F_n$ and $\pi_n^2: F_n\to {\hat F}_n$ are quotient maps.

For any $f\in \Aff(T({ A}_n))$ with $\|f\|\leq 1$, write $f_1=(\pi_n^1)^\#(f)$, $f'=\phi_{n,n+1}^{\sharp} (f)$ and $f'_1=\psi_{n,n+1}^{\sharp} (f_1)$.
By the condition $c_{ij}>13\cdot2^{2n}$ and  the claim in the proof of Theorem \ref{range 0.27}, we have
\beq\label{0.28d}
\|\GM_{n+1}^2\circ(\pi_{n+1}^2)^{\sharp}(f'_1)-f'_1\|<\frac{2}{2^{2n}}.
\eneq
By condition $\spdd$ (which is weaker than $\spddd_2$), applying Lemma \ref{range 0.10} to ${\bar A}_n \to {\bar A}_{n+1}$ as $C_n \to C_{n+1}$ (note that the definition of $\GM'_n$ is the same as the definition of $\GM'_n$ in \ref{range 0.9}  and $\bar \phi_{n,n+1}: {\bar A}_n \to {\bar A}_{n+1}$ is the same as $\phi_{n,n+1}:C_n \to C_{n+1}$), we have
\begin{equation}\label{0.28e}
\|\left(\GM_{n+1}^1\circ(\pi_{n+1}^1)^{\sharp}(f')-f'\right)|_{Sp({\bar A}_{n+1})}\|<\frac{2}{2^{2n}}.
\end{equation}
Note that $\left(\GM_{n+1}^1\circ(\pi_{n+1}^1)^{\sharp}(f')-f'\right)|_{Sp({ F}_{n+1})}=0$. So we have
\begin{equation}\label{0.28f}
\|\GM_{n+1}^1\circ(\pi_{n+1}^1)^{\sharp}(f')-f'\|<\frac{2}{2^{2n}}.
\end{equation}
Consequently (applying  $(\ref{0.28d}) (\ref{0.28f}$)), we get

\begin{eqnarray*}
&&\hspace{-0.3in}\|\GM_{n+1}\circ\pi_{n+1}^{\sharp}(f')-f'\|
= \|\GM_{n+1}^1\circ\GM_{n+1}^2\circ(\pi_{n+1}^2)^{\sharp}\circ(\pi_{n+1}^1)^{\sharp} (f')-f'\| \\
 & &= \|\GM_{n+1}^1\circ\GM_{n+1}^2\circ(\pi_{n+1}^2)^{\sharp}\circ\psi_{n,n+1}^\#\circ(\pi_{n}^1)^{\sharp} (f)-f'\| \\
 && (\textrm{since $\pi_{n+1}^1\circ\phi_{n,n+1}=\psi_{n,n+1}\circ\pi_n^1$}) \\
& = & \|\GM_{n+1}^1\circ\GM_{n+1}^2\circ(\pi_{n+1}^2)^{\sharp}\circ\psi_{n,n+1}^\#(f_1)-f'\| 
=\|\GM_{n+1}^1\circ\GM_{n+1}^2\circ(\pi_{n+1}^2)^{\sharp}(f'_1)-f'\| \\
& < & \|\GM_{n+1}^1(f'_1)-f'\|+\frac{2}{2^{2n}}\quad\quad(\mbox{by}~(\ref{0.28d})) \\
& = &\|\GM_{n+1}^1\circ\psi_{n,n+1}^\#\circ(\pi_n^1)^\#(f)-f'\|+\frac{2}{2^{2n}}  
=\|\GM_{n+1}^1\circ(\pi_{n+1}^1)^\#\circ\phi_{n,n+1}^\#(f)-f'\|+\frac{2}{2^{2n}} \\
& = &\|\GM_{n+1}^1\circ(\pi_{n+1}^1)^\#(f')-f'\|+\frac{2}{2^{2n}} 
 < \frac{2}{2^{2n}}+\frac{2}{2^{2n}}.
\end{eqnarray*}
%\begin{displaymath}
%\xymatrix{
%\|\GM_{n+1}\circ\pi_{n+1}^{\sharp}(f')-f'\|\qq\qq\qq\qq\qq\qq\qq\qq \\
%=\|\GM_{n+1}^1\circ\GM_{n+1}^2\circ(\pi_{n+1}^2)^{\sharp}\circ(\pi_{n+1}^1)^{\sharp} (f')-f'\|\qq\qq\qq\qq \\
%=\|\GM_{n+1}^1\circ\GM_{n+1}^2\circ(\pi_{n+1}^2)^{\sharp}\circ\psi_{n,n+1}^\#\circ(\pi_{n}^1)^{\sharp} (f)-f'\| \qq\qq\qq\\
%\qq\qq\qq\qq\qq\qq(\mbox{since}~~\pi_{n+1}^1\circ\phi_{n,n+1}=\psi_{n,n+1}\circ\pi_n^1) \\
%=\|\GM_{n+1}^1\circ\GM_{n+1}^2\circ(\pi_{n+1}^2)^{\sharp}\circ\psi_{n,n+1}^\#(f_1)-f'\|\qq\qq\qq\qq \\
%=\|\GM_{n+1}^1\circ\GM_{n+1}^2\circ(\pi_{n+1}^2)^{\sharp}(f'_1)-f'\|\qq\qq\qq\qq\qq\qq \\
%<\|\GM_{n+1}^1(f'_1)-f'\|+\frac{2}{2^{2n}}\qq\qq(\mbox{by}~(\ref{0.28d}))\qq\qq\qq \\
%=\|\GM_{n+1}^1\circ\psi_{n,n+1}^\#\circ(\pi_n^1)^\#(f)-f'\|+\frac{2}{2^{2n}} \qq\qq\qq\qq \\
%=\|\GM_{n+1}^1\circ(\pi_{n+1}^1)^\#\circ\phi_{n,n+1}^\#(f)-f'\|+\frac{2}{2^{2n}} \qq\qq\qq\qq \\
%=\|\GM_{n+1}^1\circ(\pi_{n+1}^1)^\#(f')-f'\|+\frac{2}{2^{2n}}\qq\qq\qq\qq\qq \\
%<\frac{2}{2^{2n}}+\frac{2}{2^{2n}}.\qq\qq\qq\qq\qq\qq\qq\qq\qq
%}
%\end{displaymath}
So the claim is proved.

Use the claim, one can prove the following approximately intertwining diagram:
\begin{displaymath}
    \xymatrix{
        \Aff(T(A_1)) \ar@/_/[d]_{\pi_1^\#}\ar[r]^{\phi_{1,2}^\#}  & \Aff(T(A_2)) \ar[r]^{\phi_{2,3}^\#} \ar@/_/[d]_{\pi_2^\#}& \Aff(T(A_3)) \ar[r] \ar@/_/[d]_{\pi_3^\#}& \cd \Aff(T(F ))\\
        \Aff(T({\hat F}_1)) \ar[r]^{\hat{\psi}_{1,2}^\#}\ar@/_/[u]_{\GM_1} & \Aff(T({\hat F}_2)) \ar[r]^{\hat{\psi}_{2,3}^\#}\ar@/_/[u]_{\GM_2}& \Aff(T({\hat F}_3)) \ar[r]\ar@/_/[u]_{\GM_3}& \cd \Aff(T({\hat F}))~, }
\end{displaymath}
where $\pi_i^\#$ are induced by the quotient map $\pi_i: A_i \to {\hat F}_i$.

With the diagram above and noting that $\Aff(T({\hat F}))=\Aff(\DT)$, exactly the same argument as the proof of \ref{range 0.12} with the above approximately intertwining diagram replacing the diagram  there, we can prove
$$\big((K_0(A),K_0(A)_+, \e_A),\ K_1(A), T(A), r_{\!\!_A} \big) \cong \big((G,G_+,u),~ K, \DT, r \big).$$
% {\bf {It seems that this case is slightly different even from the proof of \ref{range 0.27}, does it not?
%It would be a good idea again, at least, refer to an already display diagram(s).}}
\end{NN}

\begin{NN}\label{range 0.29}
The algebra $A$ in \ref{range 0.28} is not simple, we need to modify
the \hm\,  $\phi_{_{n,n+1}}$ to make the limit algebra simple. Let us
emphasize that every \hm\,  $\phi: A_n\to A_{n+1}$ is completely
determined by $\phi_x=\pi_{_x}\circ \phi$ for each $x\in
Sp(A_{n+1})$, where the map $\pi_x:~ A_{n+1}\to A_{n+1}|_{_x}$ is
the corresponding irreducible representation.

Note that from the definition of $\phi: A_n\to A_{n+1}$ and from the assumption that
$c_{ij}>13$ for each entry of $\cc=(c_{ij})$, we know that for any
$x\in Sp(A_{n+1})$,
\beq\label{150102-p1}
 Sp(\phi_{_{n,n+1}}|_x)\supset Sp(\hat F_n)=(\tht_1,\tht_2,...,\tht_{p_n}).
 \eneq
(See \ref{homrestr} and \ref{ktimes} for notations.  {Also note that each irreducible representation of
$\hat{F}_n$ can be identified with one of $F_n.$}) To see the above is true, one notes that for $x\in Sp(\bar A_{n+1})$, the homomorphism $\phi_{_{n,n+1}}|_x$ {is}  defined to be $\bar \phi_{_{n,n+1}}|_x$, and  in turn, $\bar \phi_{_{n,n+1}}$ is defined in
the proof of \ref{range 0.6} (it was called $\phi_{n,n+1}$ there)
%and
{which} satisfies condition (5)---the same condition as above.
For $x\in Sp(F_{n+1})$, we have $\phi_{_{n,n+1}}|_x=\psi_{_{n,n+1}}|_x $. From the definition of $\psi_{_{n,n+1}}$
(see \ref{range 0.20},\ref{range 0.21} and \ref{range 0.23}), we know that,  if $x\in Sp(F_{n+1}^i)$
(for $i\geq 2$), then $SP(\psi_{_{n,n+1}}|_x )$ contains exactly $c_{ij}$ copies of $\tht_j$; and if $x\in Sp(F_{n+1}^1)$, then $SP(\psi_{_{n,n+1}}|_x )$ contains exactly $c_{1j}$ copies of $\tht_j$ (for $j\geq 2$) and at least $c_{11}-13$ copies of $\tht_1$. Hence the condition also holds for this case.
To make the limit algebra simple, we need to make the set
%$Sp(\phi_{_{n,m}})_x$
$Sp(\phi_{_{n,m}}|_x)$  sufficiently dense
%enough
in $Sp(A_n)$, for any $x\in
Sp(A_m)$, provided $m$ large enough.

Write
$Sp(A_n)=\cup_{j=1}^{l_n}(0,1)_j\cup X_n\cup S_n,$
where $S_n=Sp(\bigoplus_{i=2}^{p_n}M_{[n,i]}(\C)).$
Let $0<d<1/2$ and let $Z\subset Sp(A_n).$

{\it Recall from \ref{density} that, we say
$Z$ is $d$-dense in $Sp(A_n)$ if the following holds:
$Z\cap (0,1)_j$ is $d$-dense in $(0,1)_j$ with the usual metric,
$Z\cap X_n$ is $d$-dense in $X_n$ with a given metric of $X_n$ and
if $\xi\in S_n$ is an isolated point of $Sp(A_n),$ then $\xi\in Y.$ }

Now we will change $\phi_{_{n,n+1}}$ to a map $\xi_{_{n,n+1}}:
A_n\to A_{n+1}$ satisfying:
\begin{enumerate}
\item[(i)]  $\xi_{_{n,n+1}}$ is homotopic equivalent to  $\phi_{_{n,n+1}}$.
\item[(ii)]  $\|\phi_{_{n,n+1}}^{\sharp}-\xi_{_{n,n+1}}^{\sharp}\|\leq \frac1{2^n}$.
\item[(iii)] for any $y\in Sp(A_{n+1})$, if $x\in Sp (A_n)$ satisfies $x\in SP(\phi_{_{n,n+1}}|_y)$, then $x\in SP(\xi_{_{n,n+1}}|_y)$ (we avoid to say
$ SP(\xi_{_{n,n+1}}|_y)\supset SP(\phi_{_{n,n+1}}|_y)$, since multiplicity of the irreducible representation $x\in Sp(A_n)$ in $SP(\xi_{_{n,n+1}}|_y)$ may not be larger or equal to its multiplicity in $SP(\phi_{_{n,n+1}}|_y)$,  when both $ SP(\xi_{_{n,n+1}}|_y)$ and $ SP(\phi_{_{n,n+1}}|_y)$ are regarded as sets with multiplicities), consequently, $\xi_{n,n+1}$ also satisfies that for any
$y\in Sp(A_{n+1})$,
\begin{equation}\label{150102-p2}
 SP(\xi_{_{n,n+1}}|_y)\supset Sp(\hat F_n)=(\tht_1,\tht_2,...,\tht_{p_n}),
 \end{equation}
since $\phi_{n,n+1}$ has the same property.
\item[(iv)]  For any $i\leq n$,
$$\qq\qq\qq
%Sp(\xi_{_{i,n+1}})|_{_{\tht_2'}}
Sp(\xi_{_{i,n+1}}|_{_{\tht_2'}})\ \  \mbox{is  $\frac1n$\,-  dense in  } Sp(A_i),
$$
\end{enumerate}
%\qq\qq\qq\qq (**)$$
where $\tht_2'\in Sp(F^2_{n+1})\sbs
Sp(\hat F_{n+1})=\{\tht_1',\tht_2',...,\tht'_{p_{n+1}}\}$ is the second
point (note that $\tht_1'$ is the base point of
$[\tht_1',\tht_1'+1]\vee X_{n+1}'=X_{n+1}=Sp(F_n^1)$,
and we do not
want to modify this one).
%{\bf  First, are there more typos? $i=n?$} { I do want (iv) holds for all $i\leq n$} {  Second, again, $Sp(\xi_{i,n+1}|_x)$ is not a subset of $Sp(A_i).$
%----perhaps, I know what you meant. But we have changed the meaning of notations.} { Now I changed to $SP(\xi_{i,n+1}|_x)$} {
%Third,  perhaps, I do not know  $Sp(\xi_{i,n+1}|_{\theta'_2})$ means.
%Here is another question, I thought $\theta_1'$ is a point ---at some step, it becomes  an irreducible representation
%of $F_n$ which is perfect fine---but  should $Sp(F_{n+1})$ includes $X_{n+1}?$  In that case, why
%$Sp(F_{n+1})$ is written as $p_{n+1}$ points? May be you nean, the first
%position of $Sp(F_{n+1})$ is just $X_{n+1}$ and others are referred to other direct summands?} { There were  some  typos here, Changed some of $F_n$ to $\td F_{n}$ and $F_{n+1}$ to $\td F_{n+1}$ whose spectrum do not including the point in $X_{n+1}$ other than $\tht_1'$} \\
\begin{enumerate}
\item[(v)] For any $x\in
Sp(F_{n+1})$ satisfying $x\not= \tht_2',$
\begin{equation}\label{150102-p2-1}
%\qq\qq \qq
\phi_{_{n,n+1}}|_{_x} = \xi_{_{n,n+1}}|_{_x}.
\end{equation}
%\qq\qq  (***) $$
In
particular we have,

\item[(v')] $\phi_{n,n+1}|_{Sp(F_{n+1}^1)}=\xi_{n,n+1}|_{Sp(F_{n+1}^1)}$,  or equivalently, for any $x\in X_{n+1}=Sp(F_{n+1}^1),$
\begin{equation}\label{0.29-1}
\phi_{_{n,n+1}}|_{_x}
= \xi_{_{n,n+1}}|_{_x}.
\end{equation}
\end{enumerate}
This property is important for us
%go
to apply
\ref{range 0.25} to prove that the limit algebra $A$ {has} the
property that
 $A\otimes U\in {\cal B}_0$
 for any UHF-algebra $U$.
%\end{enumerate}

Note that, in (iv) we do not need that $Sp(\xi_{_{i,n+1}}|_{_x})$ to be sufficiently dense in $Sp(A_i)$ for \\
 all $x\in Sp(A_{n+1})$, but only need $Sp(\xi_{_{i,n+1}}|_{_{\tht_2'}})$ to be sufficiently dense. Then combined with condition (iii) for $\xi_{n+1,n+2}$, we will have $Sp(\xi_{_{i,n+2}}|_{_x})$ to be sufficiently dense in $Sp(A_i)$ for all $x\in Sp(A_{n+2})$, since
\begin{equation}\label{dense1}
Sp(\xi_{i,n+2}|_x)=\bigcup_{y\in Sp(\xi_{n+1,n+2}|_x)} Sp(\xi_{i,n+1}|_y) \supset Sp(\xi_{i,n+1}|_{\tht_2'}).
\end{equation}
Namely, it follows from  (iii) (for $n+1$ in place of $n$)  and (iv) that for any $i\leq n $
and any  $x\in Sp(A_{n+2})$, the set $Sp(\xi_{i,n+2}|_x)$ is $\frac1n$\,-  dense in $ Sp(A_i).$

%{\bf This seems to contradict the claim
%that you want to change $\phi_{n, n+1}$ so that $Sp(\phi_{n,n+1}|_x)$
%is dense in $Sp(A_n)$ FOR ALL $x\in Sp(A_{n+1})$ ?------Please note, there are few places that I think, by changing the notation, it seems
%that I could guess what you meant. This one really got me.----assume though, again, notation problems.} { I think after add the sentence, it
%should clear your doubt}

 %, that is of tracially approximately an
 %Elliott-Thomsen trivial $K_1$ algebra.\\

\end{NN}

\begin{NN}\label{range 0.30}
Suppose that we
%already
have constructed
%had the definition of
\begin{displaymath}
\xymatrix{ A_1  \ar[r]^{\xi_{_{1,2}}} &A_2  \ar[r]^{\xi_{_{2,3}}} &
\cd \ar[r]^{\xi_{_{n-1,n}}}& A_n   ,}
\end{displaymath}
such that for all $i\leq n-1$, $\xi_{i,i+1}$ satisfy Conditions (i)--(v) (with $i$ in place of $n$) in \ref{range 0.29}. We will
%define
construct the map $\xi_{_{n,n+1}}:~ A_n\to A_{n+1}$.
Let
\beq\nonumber
\vspace{-0.1in}A_n=\big\{(f,g)\in C([0,1],E_n)\oplus F_n;~ f(0)=\bt_0(\pi(g)), f(1)=\bt_1(\pi(g))\big\}\,\,\,{\rm with}\\
%with
%$$
F_n=P_n M_{\infty}(C(X_n))P_n\oplus
%\oplus_{i=2}^{p_n} M_{[n,i]}(\C),
\bigoplus_{i=2}^{p_n} M_{[n,i]}(\C),
\eneq
where ${\rm rank}(P_n)=[n,1]$.

By Condition (iii) of \ref{range 0.29}, applying  to each $j\in\{i,i+1,...,n-1\}$ (in place of $n$), we have that for each $x\in Sp(A_i)$, if $x\in SP(\phi_{_{i,n}}|_{_z})$ for some $z\in Sp(A_n)$, then $x\in SP(\xi_{_{i,n}}|_{_z})$. Recall that the  finite set $Y\subset X_n$ from \ref{condition2} satisfies that for each $i < n$, $
\bigcup_{y\in Y}SP(\phi_{_{i,n}}|_{_y})$ is $\frac1n$-dense in $X_i$. Combine this two facts, we have that

(vi) $\cup_{y\in Y}SP(\xi_{_{i,n}}|_{_y})$ is also $\frac1n$-dense in $X_i$.

%This can be done since the corresponding map $\pi_1\circ
%\phi_{_{i,i+1}}|_{_{F_i^1}}:~F_i^1\to F_{i+1}^1 $ is injective for
%each $i$ and $\pi_1\circ \phi_{_{i,i+1}}=\pi_1\circ \xi_{_{i,i+1}}$
%(see (\ref{150102-p2}) in \ref{range 0.29}).
Also recall that, we  denote $t\in(0,1)_j\sbs
Sp\big(C([0,1],E_n^j)\big) $ by $t_{(j)}$ to distinguish
%spectral
spectrum
from different direct summand of $C([0,1], E_n)$, and we denote  $0\in [0,1]_j$ by
$0_j$ and $1\in [0,1]_j$ by $1_j.$ Note that $0_j$ and $1_j$ do not
correspond to irreducible representations.  In fact, $0_j$
corresponds to the direct sum of irreducible
%representation
representations for the
set
$$\left\{ \tht_1^{\sim b_{0,j1}}, \tht_2^{\sim b_{0,j2}},..., \tht_n^{\sim b_{0,jp_n}}\right\}$$
%and $1_j$,
and $1_j$ corresponds to the set
$$\left\{ \tht_1^{\sim b_{1,j1}}, \tht_2^{\sim b_{1,j2}},..., \tht_n^{\sim b_{1,jp_n}}\right\}.$$
Again recall from \ref{condition2},  $T\sbs Sp(A_n)$ is defined by
$$T=\left\{ (\frac kn)_{(j)};~ j=1,2,... ,l_n;~ k=1,2, ... , n-1 \right\}.$$
 We need to modify $\phi_{_{n,n+1}}$ to $\xi_{_{n,n+1}}$ such that
\begin{equation}\label{0.30-1}
SP\big(\xi_{_{n,n+1}}|_{_{\tht_2'}}\big) \sps T\cup Y.
\end{equation}
 If this is done, then combining with (vi) above, we know that (iv) of \ref{range 0.29} holds.

 Recall that   in conditions $\spddd_1$ and $\spddd_2$,  $L=l_n\cdot (n-1)+L_1=\#(T\cup Y)$ and  $M=\max\{ b_{0,ij};~ i=1,2,...,p_n;~ j=1,2,..., l_n\}$, where we write $Y=\{y_1,y_2,..., y_{_{L_1}}\}\subset X_n$.
 %Assume that
 %matrices $\cc=(c_{ij})$, $\bb_0'=(b_{0,ij}')$ and $\bb_1'=(b_{1,ij}')$ satisfy the
 % strengthened condition
%$$
%\spddd_1 \qq\qq\qq\qq\qq\qq\qq c_{ij}> 13 \cdot 2^{2n}\cdot ML\qq\,\,\,
%\rforal \,\,i,j\qq\qq\qq\qq$$
%$$\spddd_2 ~~~
%\td b_{0,il}=\sum_{k=1}^{p_{n+1}} b'_{0,ik}\cdot
%c_{kl}> 2^{2n}\left(\sum_{k=1}^{l_{n}}(|d_{ik}|+2)\cdot \{n,k\}+
%(L_1+(n-1)\sum_{k=1}^{l_n}b_{0,kl})\cdot b'_{0,i2}
%\right)~\mbox{and}$$
%$$\spddd_2 \qq
%\qq\qq\qq
%\td b_{1,il}=\sum_{k=1}^{p_{n+1}} b'_{1,ik}\cdot c_{kl}>
%2^{2n}\left(\sum_{k=1}^{l_{n}}(|d_{ik}|+2)\cdot \{n,k\}+
%(L_1+(n-1)\sum_{k=1}^{l_n}b_{0,kl})\cdot b'_{1,i2} \right),
%$$
% for
%each $l\in \{1,2,\cd, p_n\}$ and $i\in \{1,2,\cd, l_{n+1}\}$. To
%make the above inequalities hold,  we only need to  increase  the
%entries of the third columns of $\bb'_0$ and $\bb'_1$ by adding same
%big positive number to make it much larger then the second column of
%the matrices (see the end of \ref{range 0.5}).\\

%{\bf Let us emphasis that we assume that, in the original construction of $A_{n+1}$, the matrices $\cc=(c_{ij}),\bb'_0=(b'_{0,ij})$ and $\bb'_1=(b'_{1,ij})$ satisfy the above condition $\spddd_1,\spddd_2$. Under this condition, we will discuss how to construct $\xi_{n,n+1}$ from  $\phi_{n,n+1}$  to satisfy (i)---(iv) in \ref{range 0.29}. So in this procedure, $A_{n+1}$ and $\phi_{n,n+1}$ are kept the same.}

\end{NN}

\begin{NN}\label{range 0.31}~~  We only defined for each $x\in Sp(F_{n+1})$
with $x\not=\tht_2'$, $\xi_{_{n,n+1}}|_{_x} =
\phi_{_{n,n+1}}|_{_x}$. Now we define
 $\xi_{_{n,n+1}}|_{_{\tht_2'}}$. To simplify the notation, write $\xi_{_{n,n+1}}:= \xi$.
Note that
$$SP(\phi_{_{n,n+1}}|_{_{\tht_2'}})=\left\{ \tht_1^{\sim c_{21}}, \tht_2^{\sim c_{22}},..., \tht_n^{\sim c_{2p_n}}\right\}.$$
We use $y_1,y_2,...,y_{_{L_1}}$ to replace $L_1$ copies of $\tht_1$ in $\{\tht_1^{\sim c_{21}}\}$, so $\{\tht_1^{\sim c_{21}}\}$ becomes
 $$\{\tht_1^{\sim (c_{21}-L_1)}, y_1,y_2,... ,y_{_{L_1}} \}.$$  That is, in
 the definition of $\xi_{_{n,n+1}}(f)|_{_{\tht_2'}}$ we use $f(y_i)$ to replace
 one copy of $f(\tht_1)$.  Note that $f(y_i)=P_n(y_i)M_{\infty}(\C)P_n(y_i)$ can
 be identified with $M_{[n,1]}(\C)=\hat F_n^1$.  For late use we also choose a
 path $y_i(s)~(0\leq s\leq 1)$ from $\tht_1$ to $y_i$.  That is,
 $y_i(0)=\tht_1,~ y_i(1)=y_i$, and fix identification of
 $P_n(y_i(s))M_{\infty}(\C)P_n(y_i(s))$ with $M_{[n,1]}(\C)=
 \hat F_n^1=P_n(\tht_1)M_{\infty}(\C)P_n(\tht_1)$---such
 %add an an
 an identification could be chosen to be continuously depending on $s$. (Here,
 we only use the fact that any projection (or vector bundle) over
 % add the
the  interval is trivial to make such identification. Since  the
 projection $P_n$ itself may
 not be trivial, it is possible
 that the paths for different $y_i$ and $y_j$ ($i\not=j$) may intersect at
 $y_i(s_1)=y_j(s_2)$
 and we may use different identification of
 $P_n(y_i(s_1))M_{\infty}(\C)P_n(y_i(s_2))= P_n(y_j(s_2))M_{\infty}(\C)P_n(y_j(s_2))$
 with $M_{[n,1]}(\C)$ for $i$ and $j$.)
 When we talk about $f(y_i(s))$ later, we will considered it to
 be an element in $M_{[n,1]}(\C)$ (rather than in
 $P_n(y_i(s))M_{\infty}(\C)P_n(y_i(s))$).

 Also we use $(\frac kn)_{(j)}\in (0,1)_j$ to replace
 $$0_j=\left\{ \tht_1^{\sim b_{0,j1}}, \tht_2^{\sim b_{0,j2}},..., \tht_n^{\sim b_{0,jp_n}}\right\}.$$

 In summary, we have that
 $SP(\xi|_{_{\tht_2'}}) = \left\{ \tht_1^{\sim a_1}, \tht_2^{\sim a_2},..., \tht_n^{\sim a_{p_n}}\right\} \cup T \cup Y$ with 
 %with
 $$
 \begin{array}{ccl}
   a_1 & = & c_{21}-L_1-\left(\sum_{j=1}^{l_n}b_{0,j1} \right)(n-1) \\
   a_2 & = & c_{22}-\left(\sum_{j=1}^{l_n}b_{0,j2} \right)(n-1) \\
   & \vdots & \\
   a_{p_n} & = & c_{2p_n}-\left(\sum_{j=1}^{l_n}b_{0,jp_n} \right)(n-1) ~.
 \end{array}
$$
From $\spddd_1$ of \ref{range 0.30},
 we know that
\beq\label{150102-e1}
%$$\qq\qq\qq\qq\qq
a_i \geq \frac{2^{2n}-1}{2^{2n}}c_{2i},
%\qq\qq\qq\qq\qq (*)$$
\eneq
and as many as $a_i$ points (counting multiplicity) of $SP(\phi_{_{n,n+1}} |_{_{\tht_2'}})$ remain the
same. The map $\xi|_{_{\tht_2'}}$ is not unique, but is unique up to
unitary equivalence.
 Note that $\phi_{_{n,n+1}} |_{_{\tht_2'}}:~ A_n\to F^2_{n+1}$ is
 homotopic   to $\xi|_{_{\tht_2'}}:~ A_n\to F^2_{n+1}$.  Namely, for each $s\in [0,1]$,
 one can define homomorphism $\Omega(s):A_n\to F^2_{n+1}$ as the definition
 of $\xi|_{_{\tht_2'}}$,
 %add by
 by replacing $y_i$ by $y_i(s)$, and
 $(\frac kn)_{(j)}\in (0,1)_j$ by
 $(s\cdot\frac kn)_{(j)}\in (0,1)_j$.
 It is obvious that $\Omega(0)=\phi_{_{n,n+1}} |_{_{\tht_2'}}$,
 and $\Omega(1)=\xi |_{_{\tht_2'}}$. {Again recall  $\xi$ is simplified notation for $\xi_{n, n+1}$. }\\

 If $A_{n+1}$ is in the case of \ref{range 0.5a}, with
 $p_{n+1}^0 \geq 2$, then, for $i\leq p_{n+1}^0$, $Sp(F_{n+1}^i)$ is a closed and open subset of $Sp(A_{n+1})$ since $F_{n+1}^i$ is separate from $C([0,1], E_{n+1})$ in the definition. On the other hand, $F_{n+1}^2=\hat F_{n+1}^2$ has a single spectrum $\tht_2'$. That is, $\tht_2'$ is an isolate point in
 $Sp(A_{n+1})$.  Then we can define $\xi|_{x}=\phi_{n,n+1}|_{x}$ for
 any $x\in Sp(A_{n+1})\setminus \{\tht_2'\}$ to finish the
 construction of $\xi=\xi_{n,n+1}$ (this case is much
 easier). Also, for this case, we do not need the estimation $\spddd_2$
 (only $\spddd_1$ and the old $\spdd$ (see (\ref{13spdd})) will be enough).

 \end{NN}

\begin{NN}\label{range 0.32}
~~   In \ref{range 0.31}, we have defined the part $\xi$
(restricted to $F_{n+1}$ and also to $\hat F_{n+1}$); that is
\begin{displaymath}
\xymatrix{
 \xi'=\pi\circ \xi:~~ A_n \ar[r]  &A_{n+1}  \ar[r]^{\pi} & \hat F_{n+1}~,}
\end{displaymath}
  where $\pi$ is the quotient map modulo the ideal $I_{n+1}+J'_{n+1} $ (see \ref{construction} for notation of $I_{n}$ and $J'_{n}$).  {Recall that $\xi$ is the simplified notation for  $\xi_{n,n+1}$.} Now we define
$\xi|_{_{[0,1]_j}}  $ for each $[0,1]_j\subset Sp(C[0,1], E^j_{n+1})$.
We already know the definition of
$$\xi|_{_{\hat0_j}}=\pi_j\circ \bt_0'\circ \xi':~ A_n\to E_{n+1}^j\qq\qq \mbox{
and }\qq\qq \xi |_{_{\hat1_j}} = \pi_j\circ \bt_1'\circ \xi':~ A_n\to
E_{n+1}^j,$$
  where $\pi_j:~ E_{n+1}\to E_{n+1}^j$.  Here we use
$\hat0_j,~\hat1_j \in [0,1]_j \sbs Sp(C([0,1],E_{n+1}^j))$, and
reserve $0_j, 1_j$ for the corresponding points in $Sp(C([0,1],E_{n}^j))$.

Now we need to connect $\xi|_{_{\hat0_j}}$ and $\xi|_{_{\hat1_j}}$
to obtain the definition of $\xi|_{_{[0,1]_j}}  $.  Let $\psi^j:~
C([0,1],E_{n}) \to C([0,1],E_{n+1}^j)$ be as defined in the proof of
\ref{range  0.6} (and \ref{range  0.7}).  Let $\GM:~ A_n \to
C([0,1],E_{n})$ be the natural inclusion map.  As in the proof of
\ref{range  0.6}, for the original map $\phi_{_{n,n+1}}$, we have
$$
 SP(\phi_{_{n,n+1}}|_{_{\hat 0_j}})= \left\{SP(\psi^j\circ \GM|_{_{\hat 0_j}} ),~
  \tht_1^{\sim r^j_1}, \tht_2^{\sim r^j_2},..., \tht_{p_n}^{\sim r^j_{p_n}} \right\}~\,\andeqn$$
% and
$$
 SP(\phi_{_{n,n+1}}|_{_{\hat 1_j}})= \left\{SP(\psi^j\circ \GM|_{_{\hat 1_j}} )~
  \tht_1^{\sim r^j_1}, \tht_2^{\sim r^j_2},..., \tht_{p_n}^{\sim r^j_{p_n}} \right\}~,$$
where
 \beq\label{0.32*}  r_l^j=\sum_{k=1}^{p_n+1}b_{0,jk}'c_{kl} -
\left(\sum_{d_{jk}<0}|d_{jk}|b_{1,kl}
+\sum_{d_{jk}>0}|d_{jk}|b_{0,kl}
+\sum_{d_{jk}=0}(b_{0,kl}+b_{1,kl})\right).
\eneq
Note that $SP(\xi|_{_{\hat0_j}})$ and $SP(\xi|_{_{\hat1_j}})$ are
obtained by replacing some subsets
$$\left\{ \tht_1^{\sim b_{0,i1}}, \tht_2^{\sim b_{0,i2}},..., \tht_{p_n}^{\sim b_{0,ip_n}}\right\}=0_i\in Sp(C([0,1],E^i_{n}))$$
in $SP(\phi_{_{n,n+1}}|_{_{\hat 0_j}})$ and  $SP(\phi_{_{n,n+1}}|_{_{\hat 1_j}})$ by $(\frac kn)_{(i)}$ and replacing some of $\tht_1$ each by one of
$y_k~ (k=1,2,..., L_1)$ (see \ref{homrestr} and \ref{ktimes} for notation). We would like to give precise calculation
for the numbers of $\tht_1$ or $0_i$ to be replaced. Recall  that
$${\hat 0}_j=\{(\tht_1')^{\sim b'_{0,j1}}, (\tht_2')^{\sim b'_{0,j2}}, ...,
(\tht_{p_{n+1}}')^{\sim b'_{0,jp_{n+1}}}\},~~~\mbox{and}~~~{\hat
1}_j=\{(\tht_1')^{\sim b'_{1,j1}}, (\tht_2')^{\sim b'_{1,j2}}, ... ,
(\tht_{p_{n+1}}')^{\sim b'_{1,jp_{n+1}}}\}.$$ From the part of $\xi$
already defined, recall that
$\xi|_{\tht_i'}=\phi_{n,n+1}|_{\tht_i'}$ if $i\not=2$. So the set
$SP(\xi|_{{\hat 0}_j})$ is obtained from $SP(\phi_{n,n+1}|_{{\hat
0}_j})$ by doing the following replacement: replace each  group of
$b'_{0,j2}$ copies of $\tht_1$ by same number copies of $y_i$ for
each $i=1,2,..., L_1$(totally  $L_1\cdot b'_{0,j2}$ copies of
$\tht_1$ have been replaced  by some $y_i$'s), and replace each
group of $b'_{0,j2}$ copies of $0_i$ for $i=1,2,...,l_n$, by same
number of copies of $(\frac kn)_{(i)}$ for each $k=1,2,..., n-1$
(totally $(n-1)\cdot b'_{0,j2}$ copies of $0_i=\left\{ \tht_1^{\sim
b_{0,i1}}, \tht_2^{\sim b_{0,i2}},..., \tht_{p_n}^{\sim
b_{0,ip_n}}\right\}$ have been replaced by some $(\frac
kn)_{(i)}$'s). Consequently
$$
 SP(\xi|_{_{\hat 0_j}})= \left\{SP(\psi^j\circ \GM|_{_{\hat 0_j}} ),~
  \tht_1^{\sim s_1}, \tht_2^{\sim s_2},..., \tht_{p_n}^{\sim s_{p_n}} \right\}~
  \cup(T\cup Y)^{\sim b'_{0,j2}}~,$$
where
 $$
 \begin{array}{ccl}
  s_1 & = & r_1^j-(L_1+(n-1)\sum_{i=1}^{l_n}b_{0,i1})\cdot b'_{0,j2} \\
   s_2 & = & r_2^j-((n-1)\sum_{i=1}^{l_n}b_{0,i2})\cdot b'_{0,j2}  \\
   & \vdots & \\
   s_{p_n} & = & r_{p_n}^j-((n-1)\sum_{i=1}^{l_n}b_{0,ip_n})\cdot b'_{0,j2} ~.
 \end{array}
$$
Exactly as in the above argument, we have
$$
 SP(\xi|_{_{\hat 1_j}})= \left\{SP(\psi^j\circ \GM|_{_{\hat 1_j}} ),~
  \tht_1^{\sim t_1}, \tht_2^{\sim t_2},..., \tht_{p_n}^{\sim t_{p_n}} \right\}~
  \cup(T\cup Y)^{\sim b'_{1,j2}}~,$$
where
 $$
 \begin{array}{ccl}
  t_1 & = & r_1^j-(L_1+(n-1)\sum_{i=1}^{l_n}b_{0,i1})\cdot b'_{1,j2} \\
   t_2 & = & r_2^j-((n-1)\sum_{i=1}^{l_n}b_{0,i2})\cdot b'_{1,j2}, \\
   & \vdots & \\
   t_{p_n} & = & r_{p_n}^j-((n-1)\sum_{i=1}^{l_n}b_{0,ip_n})\cdot b'_{1,j2} ~.
 \end{array}
$$
We can pare two sets $ SP(\xi|_{_{\hat 0_j}})$ and $ SP(\xi|_{_{\hat
1_j}})$ now. Of course the part of $ SP(\psi^j\circ \GM|_{_{\hat
0_j}} )$ and $ SP(\psi^j\circ \GM|_{_{\hat 1_j}} )$ already be
paired by the homomorphism $\psi^j$ defined before, we only need to
pair the part
 $\left\{~
  \tht_1^{\sim s_1}, \tht_2^{\sim s_2},..., \tht_{p_n}^{\sim s_{p_n}} \right\}~
  \cup(T\cup Y)^{\sim b'_{0,j2}}~$ of the set $SP(\xi|_{_{\hat 0_j}})$ with the part
 $$\left\{~
  \tht_1^{\sim t_1}, \tht_2^{\sim t_2},..., \tht_{p_n}^{\sim t_{p_n}} \right\}~
  \cup(T\cup Y)^{\sim b'_{1,j2}}~$$ of the set $ SP(\xi|_{_{\hat 1_j}})$.
 If $b'_{1,j2}\geq b'_{0,j2}$ (the other case
$b'_{0,j2}> b'_{1,j2}$ { can be done similarly}), then $s_i\geq t_i$.
First we do the obvious pairing. That is, we can pair $t_i$ copies
of $\tht_i$ from $SP(\xi|_{_{\hat 0_j}})$ with all $t_i$ copies of
$\tht_i$ in $SP(\xi|_{_{\hat 1_j}})$;  and we can also pair, for
each $x\in Y\cup T$, all $b'_{0,j2}$ copies of $x$ in
$SP(\xi|_{_{\hat 0_j}})$ with  some $b'_{0,j2}$ ($\leq b'_{1,j2}$)
copies of $x$ in $ SP(\xi|_{_{\hat 1_j}})$. After this procedure of pairing, what left in $ SP(\xi|_{_{\hat
1_j}})$ is the set $(Y\cup
T)^{\sim (b'_{1,j2}-b'_{0,j2})}$, and what left in $
SP(\xi|_{_{\hat 0_j}})$ is $\{\tht_1^{\sim (s_1-t_1)},\tht_2^{\sim (s_2-t_2)},
\cd, \tht_{p_n}^{\sim (s_{p_n}-t_{p_n})}\}$. And we will pair these parts as below.
 For each $y\in Y$, there are
$b'_{1,j2}-b'_{0,j2}$ copies of $y$ left in $ SP(\xi|_{_{\hat
1_j}})$, and we pair all those copies (for each $y$) with same
number of copies of $\tht_1$ in $SP(\xi|_{_{\hat 0_j}})$---it will
totally cost $L_1\cdot(b'_{1,j2}-b'_{0,j2})$ copies of $\tht_1$; and
for each $(\frac kn)_{(i)}\in T$, there are also
$b'_{1,j2}-b'_{0,j2}$ copies of $(\frac kn)_{(i)}$  left in $
SP(\xi|_{_{\hat 1_j}})$, we pair all those copies (for each $(\frac
kn)_{(i)}$)  with same number of copies of $0_i=\left\{ \tht_1^{\sim
b_{0,i1}}, \tht_2^{\sim b_{0,i2}},..., \tht_{p_n}^{\sim
b_{0,ip_n}}\right\}$---it will totally cost
$$
 \begin{array}{ccl}
 ((n-1)\sum_{i=1}^{l_n}b_{0,i1})\cdot (b'_{1,j2}-b'_{0,j2}) & \mbox{copies of} & \tht_1, \\
 ((n-1)\sum_{i=1}^{l_n}b_{0,i2})\cdot (b'_{1,j2}-b'_{0,j2}) & \mbox{copies of} & \tht_2,  \\
   & \vdots & \\
  ((n-1)\sum_{i=1}^{l_n}b_{0,i{p_n}})\cdot (b'_{1,j2}-b'_{0,j2}) & \mbox{copies of} & \tht_{p_n}.
 \end{array}
$$
Hence, putting
%%%puting
 all them together, totally we will use exactly
\beq
(L_1+(n-1)\sum_{i=1}^{l_n}b_{0,i1})\cdot (b'_{1,j2}-b'_{0,j2})=s_1-t_1~~~ \mbox{copies of} ~~~\tht_1
\andeqn\\
%$$ and
%$$
((n-1)\sum_{i=1}^{l_n}b_{0,ik})\cdot (b'_{1,j2}-b'_{0,j2})=s_k-t_k ~~~\mbox{copies of} ~~~ \tht_k,
\eneq
for each $k\in\{2,3, ..., p_n\}$. That is, the set $(Y\cup
T)^{\sim (b'_{1,j2}-b'_{0,j2})}$ is paired exactly  with
 $$\{\tht_1^{\sim (s_1-t_1)},\tht_2^{\sim (s_2-t_2)},
..., \tht_{p_n}^{\sim (s_{p_n}-t_{p_n})}\}.$$

Then, we can write
$$SP(\xi|_{_{\hat 0_j}})=\left\{SP(\psi^j\circ \GM |_{_{\hat 0_j}}), \Xi_0 \right\}\andeqn 
%and
SP(\xi|_{_{\hat 1_j}})=\left\{SP(\psi^j\circ \GM |_{_{\hat 1_j}}), \Xi_1 \right\},$$
and a  pairing $\Theta$ between $\Xi_0$ and $\Xi_1$ such that each
pair $\sm\in\Tht$ is of the following form: $\sm=\left( x \in\Xi_0,
x \in\Xi_1\right)$ with $x\in Y\cup T\cup \{\tht_1,\tht_2,...
\tht_{p_n}\}$,  $\sm=\left(0_i \in\Xi_0, (\frac
kn)_{(i)}\in\Xi_1\right)$, or $\sm=\left( \tht_1\in\Xi_0,
y_k\in\Xi_1\right)$.  For any $\sm\in\Tht$, define $\psi_\sm:~ A_n
\to C([0,1],E^j_{n})$ by sending $f\in A_n$ to $\psi_\sm(f)$ as:
$$\psi_\sm(f)(t)=\left\{\begin{array}{lll}
                          f(x), & \mbox{if }  \sm=\left( x\in\Xi_0,~ x\in\Xi_1\right) & ({\rm case}~1),\\
                          &\\
                          f((t\frac kn)_{(i)}), & \mbox{if } \sm=\left(0_i \in\Xi_0,~ (\frac kn)_{(i)}\in\Xi_1\right) & ({\rm case}~2),\\
                          &\\
                          f(y_k(t)), & \mbox{if } \sm=\left( \tht_1\in\Xi_0,~
                          y_k\in\Xi_1\right) & ({\rm case}~3).
                        \end{array}
\right.
$$
(Of course, if $b'_{0,j2}>b'_{1,j2}$, then the  case (2) above will
be changed to $\sm=\left( (\frac kn)_{(i)}
\in\Xi_0,~0_i\in\Xi_1\right)$ and the function $\psi_\sm(f)(t)$
should be defined by $\psi_\sm(f)(t)= f(((1-t)\frac kn)_{(i)})$; and
the case (3) above will be changed  to $\sm=\left(y_k
\in\Xi_0,~\tht_1\in\Xi_1\right)$ and the function $\psi_\sm(f)(t)$
should be defined by $\psi_\sm(f)(t)= f(y_k(1-t))$.)

  Finally let $\xi^j(f)=\diag\big(\psi^j(f),
\psi_{\sm_1}(f),\psi_{\sm_2}(f),..., \psi_{\sm_\bullet}(f)  \big)$
where $\{ \sm_1,\sm_2,...,\sm_\bullet\}=\Tht$. Then this $\xi^j$ has
the expected spectrum at $\hat0_j$ and $\hat1_j$.   After conjugating a
unitary, we can get
\beq\label{compatible}\xi^j(f)(\hat 0_j)=\pi_j\circ \bt_0'\circ \xi'(f),~~\mbox{and}~~
\xi^j(f)(\hat 1_j)=\pi_j\circ \bt_1'\circ \xi'(f).
\eneq
Combining this $\xi^j$ together with the previous defined $\xi$ into
$F_{n+1}$, we get the definition of $\xi: A_n\to A_{n+1}$. That is  for $x\in Sp(C_0((0,1),E_n^j))\subset Sp(\bar A_{n+1})$, define $\xi|_x$ to be $\xi^j|_x$; and for
$x\in Sp(F_{n+1})$, define $\xi|_x$ to be as previously defined $\xi|_x$ on this part. The condition \ref{compatible} says these two parts of definition compatible on the their boundary $Sp(\hat F_{n+1})=Sp({\bar A}_{n+1})\cap  Sp(F_{n+1})$.    Notice
that $\psi_\sm$ (part of $\xi$) is homotopic to the constant map,
$$\psi'_\sm(f)(t)=f(0_i)\mbox{~~for case 2 and}~~~ x=(\frac kn)_{(i)}\in T\mbox{~~~in case 1}, \qq \mbox{or}$$
$$\qq \psi'_\sm(f)(t)=f(\tht_1)\mbox{~~for case 3 and}~~~x\in Y\mbox{~~~~in
case 1},$$ which just matches the definition of the corresponding part
of  $\phi_{_{n,n+1}} $ (see the penultimate paragraph  of \ref{range
0.31} also). That is, $\xi_{n,n+1}$ is homotopic to $\phi_{n,n+1}$, which is (i) of \ref{range 0.29}. From $(*)$ above, the definition of $s_l$ and $t_l$, and
condition $\spddd_2$ of \ref{condition2}, we have $$s_l>\frac{2^{2n}-1}{2^{2n}}r_l^j\quad\textrm{and}\quad
 t_l>\frac{2^{2n}-1}{2^{2n}}r_l^j.$$ This implies that, for each
$x\in Sp(A_{n+1})$, $SP(\phi_{_{n,n+1}}|_x) $ and $SP(\xi|_x)$
differ by a fraction at most $\frac1{2^{2n}}$. Thus we have
$$\|\phi_{_{n,n+1}}^{\sharp}-\xi^{\sharp}\|<
\frac1{2^{2n}},$$
which is (ii) of \ref{range 0.29}. Evidently, (iii) and (v) of \ref{range 0.29} follows from the construction.  {Note that we have already made (\ref{0.30-1}) hold. Therefore} (iv) follows, as we mentioned immediately after the equation (\ref{0.30-1}). Therefore our inductive construction is completed. \\

\end{NN}

%{\bf { I noticed that you stated that you wanted to change maps $\phi_n$ so that the inductive
%limit is simple. However, in the process, would it be correct that $A_{n+1}$ is also changed? If so,
%then obvious questions are whether all these lemmas and other requirements are met after one changed $A_{n+1}.$
%These do not seem fully discussed. If $A_{n+1}$ does change, would $\phi_{n, n+1}$ be also changed? or do we still have $A_{n+2}$ to be fit
%for $A_{n+1}?$  It seems to me that a completely different construction which has not so much to do
%with original $A_n$ is constructed---surely it is true, otherwise it would not be simple.  Therefore  I assume
%that $\phi_{n+1}$ is not original one, and (1) and (2)  in \ref{range 0.29} simply are NOT true. You do not really compare with original
%$\phi_{n, n+1}.$ This should make it clearer and the statements like (1) and (2) have to restated.---I already pointed out $\phi_{n,n+1}$ and $A_{n+1}$ never changed
% }}

\begin{NN}\label{range 0.33}
   Let $ B=\lim(A_n, \xi_{_{n,n+1}})$. Recall that $A=\lim(A_n, \phi_{_{n,n+1}})$. Using (i) and (ii), exactly as  the same as the proofs of Theorems \ref{range 0.12} and \ref{range 0.27}, we have that $\big((K_0(B),K_0(B)_+, \e_B),~ K_1(B), TB, r_{\!\!_{B}} \big)$ is isomorphic to $\big((K_0(A),K_0(A)_+, \e_A),~ K_1(A), T(A), r_{\!\!_A} \big)$ (including the compatibility). Hence
\begin{eqnarray*}
\big((K_0(B),K_0(B)_+, \e_B),~ K_1(B), T(B), r_{\!\!_{B}} \big)& \cong & \big((G,G_+,u),~ K, \DT, r \big).
\end{eqnarray*}

On the other hand, we have that, for each $n$, and $m>n+1$,
$Sp(\xi_{_{n,n+1}}|_{_x})$ is $\frac 1{m-2}$ dense in $Sp(A_n)$ for
any $x\in Sp(A_m)$ (see the end of \ref{range 0.29}).  This condition implies that the limit algebra
$B$ is simple (see \ref{simplelimit}).
% {\bf {Let me say this: It is somewhat standard, in
%the construction of inductive limits of homogeneous \CA s that these conditions
%imply the inductive limit is simple.  I do not know  how standard it is to other people. There are actually two
%arguments there. One is actually little tricky since algebraic inductive limit is only dense there. But now
%it involves something else. The density  is more complicated
%than those discussed earlier.}}
Notice that in the definition of $ \xi_{_{n,n+1}}$ we have
$$\xi_{_{n,n+1}}|_{Sp (F^1_{n+1})}= \phi_{_{n,n+1}}|_{Sp (F^1_{n+1})}.$$
Then \ref{range  0.25} implies the following corollary.
\end{NN}

%{\color{Green} (There should be one condition to say that the trace of $R$ is at least $(m-1)/m$. Otherwise, one can choose $R$ to be $0$)Answer: the condition (2), in particular the injectivity of $\iota$ and the choice of $B_{n+1}$ implies that the size of $R$ is exactly $(m-1)/m$ }

\begin{cor}\label{range 0.34}
For any $m>0$ and any $A_i$, there is an integer $n\geq i$ and a
projection $R\in M_m(A_{n+1})$ such that
\begin{enumerate}
\item[(1)] $R$  commutes with $\LD\xi_{_{n,n+1}}(A_n)$, where $\LD:~
A_{n+1}\to M_m(A_{n+1})$ is the amplify map sending $a$ to an
$m\times m$ diagonal matrix:  $\LD(a)=\diag(a,\cd,a)$;

%$$\left(
%                                    \begin{array}{ccc}
%                                      a &  &  \\
%                                       & \ddots &  \\
%                                       &  & a \\
%                                    \end{array}
%                                  \right)_{m\times m~.}
%$$
\item[(2)] Recall  $\bar A_{n+1}=A\left( \hat F_{n+1}, E_{n+1}, \bt_0',\bt_1'
\right)$, where $\bt_0',\bt_1':~ \hat F_{n+1}\to E_{n+1}$ is as in
the definition of $A_{n+1}$ (see (\ref{construction-1}).  Then there is an  injective \hm
$$\iota:~M_{m-1}(\bar A_{n+1}) \lr RM_m(A_{n+1})R$$
such that $R\LD(\xi_{_{n,n+1}}(A_n))R\sbs \iota(M_{m-1}(\bar A_{n+1})).$
\end{enumerate}
\end{cor}

\begin{proof}
We have already defined  the part of $R$ in $M_m(F^1_{n+1})$ with property described in \ref{range 0.25} (the definition is given by combining \ref{range 0.22} and \ref{range 0.24}). Note that this $R$ works for the homomorphism $\xi_{n,n+1}$ in place of $\phi_{n,n+1}$ because $\xi_{_{n,n+1}}|_{Sp (F^1_{n+1})}= \phi_{_{n,n+1}}|_{Sp (F^1_{n+1})}$.  In
particular, $$R(\tht_1') = \e_{_{\hat F_{n+1}^1}}\otimes \left(
                                            \begin{array}{cc}
                                              \e_{m-1} & 0 \\
                                              0 & 0 \\
                                            \end{array}
                                          \right).
                                          $$
We extend the definition of $R$ as follows. For each $x\in
Sp(A_{n+1})\setminus Sp(F^1_{n+1})$, define $$R(x)=
\e_{_{A_{n+1}|_{_x}}}\otimes \left(
                                            \begin{array}{cc}
                                              \e_{m-1} & 0 \\
                                              0 & 0 \\
                                            \end{array}
                                          \right).
                                          $$
One can use the map $\iota$ combining with the identity map $\id:~
M_{m-1}(\bar A_{n+1})\to M_{m-1}(A_{n+1}/J'_{n+1})$ to get the corollary.
(Note that $\bar A_{n+1}=A_{n+1}/J'_{n+1})$.
\end{proof}

%{\bf{ In many Lemmas and Theorems, symbols presumed to refer to previously used symbols, in other words,
%it is difficult to read one lemma without reading all previously stated ones.  This makes reading
%difficult. On the other hand, some time symbols seem
%to change the meaning after some repeated use. Here are few examples: $\LD$ was a map from $\Aff$ spaces. Later it becomes
%maps between algebras.  $F_n$ in this section has at least two quite different forms.
%$A_n$ also may be different. We may try to ease this a bit by using better notations.}}

%{\bf {$1/n$-density:  Beside non-Hausdorff
%%%Hasdorff
% topology used. Note $SP(\phi)$ and $Sp(\phi)$ have different meanings.
%We need not to add any confusion.  It should make clear at the beginning instead stated twice but missed once.}}

%\begin{NN}
\begin{cor}\label{range 0.34a}
Let $B$ be as constructed above. Then $B\otimes U \in {\cal B}_0$  for every UHF-algebra $U$ of infinite dimension.
%{\color{Green} (We don't need $U$ to be of infinite type)} {\red{Yes, we do need it to be infinite dimension}}
\end{cor}

\begin{proof}
In the above corollary, we know that  $\bar A_{n+1} \in {\cal C}_0$ and
therefore  $M_{m-1}(\bar A_{n+1})$ and $ \iota(M_{m-1}(\bar A_{n+1}))$ are in
${\cal C}_0$. Also, for all normalized traces $\tau \in
T(M_m(A_{n+1}))$, we have $\tau(\e-R)=1/m$. 
Fix  an integer $k\ge 1$ and 
any finite subset ${\cal F}\subset B\otimes M_k,$ 
let ${\cal F}_1\subset B$ be a finite subset 
such that $\{(f_{ij})_{k\times k}: f_{ij}\in {\cal F}_1\}\supset {\cal F}.$  Now,  evidently, the
inductive limit algebra $B=\lim(A_n, \xi_{n,m})$ 
% (and $B\otimes M_k$ for any positive integer $k$) 
have the following
property: 
For any finite set ${\cal F}_1\subset B$, $\ep>0$, $\dt>0$, and any
$m>1/\dt$, there is a unital \SCA\,  $C\subset M_m(B)$ which is
in ${\cal C}_0$ such that
\begin{enumerate}
\item[(i)] $\|[\e_C, \diag\{\underbrace{f,...,f}_m\}]\|<\ep/k^2$, for all
$f\in {\cal F}_1$,
\item[(ii)] $\dist(\e_C (\diag\{\underbrace{f,...,f}_m\})\e_C, C)<\ep/k^2$,
for all $f\in {\cal F}_1$, and
\item[(iii)] $\tau(\e-\e_C)=1/m<\dt$ for all $\tau \in T(M_{m}(B))$.
\end{enumerate}
Thus  the above (i), (ii), (iii) hold  by replacing $B$ by $M_k(B),$ 
${\cal F}_1$ by ${\cal F}\subset M_k(B)$ and  $\ep/k^2$ by $\ep.$
%property for  
%$\ep>0,$ $0<\dt<1/m,$ 
%$M_k(B)=B\otimes M_k$ and finite set $F\subset M_k(B)$
%$\ep> 0$ and $\dt<1/m$ 
%follows from the same property for  $\ep/k^2,$ the same $\dt<1/m,$ $B$, $F_1\subset B$ with
%condition that if $(f_{ij})_{k\times k}\in F$, then $f_{ij}\in F_1.$
%, $\frac{\ep}{k^2}$ and same $\dt<1/m$.

Now  $B\otimes U$ can be written as $\lim (B\otimes M_{k_n}, \iota_{n,m})$ with $k_1|k_2|k_3\cd$ and $k_{n+1}/k_n \to \infty$, and $\iota_{n,{n+1}}$ is the amplification
by sending $f\in B\otimes M_{k_n}$ to  $\diag{(f,..., f)}\in B\otimes M_{k_{n+1}},$ where $f$ repeated $k_{n+1}/k_n$
times.
%$\diag\{\underbrace{f,..., f}_{k_{n+1}/k_n}\}\in B\otimes M_{k_{n+1}}$.

To show $B\otimes U\in {\cal  B}_0,$  let ${\cal F} \subset  B\otimes U$
be a finite subset and let $a \in (B\otimes U)_+\setminus  \{0\}.$ There is an integer $m_0>0$ such that $\tau(a)>1/m_0$ for all $\tau\in B\otimes U$.  Without loss of generality, we may assume that
 ${\cal F} \subset B\otimes M_{k_n}$ with $\frac{k_{n+1}}{k_n} >m_0$. Then by 
 what has been proved above
 %the above property 
 for $B\otimes M_{k_n}$ with $m=k_{n+1}/k_n$ (and note that $\iota_{n,n+1}$ is the amplification), there is a unital $C^*$ subalgebra $C\subset B\otimes M_{k_{n+1}}$ with $C\in {\cal C}_0$ such that $\|[\e_C, \iota_{n,n+1}(f) ]\|<\ep$, for all
$f\in {\cal F}$, such that $\dist(\e_C (\iota_{n,n+1}(f))\e_C, C)<\ep$ for all $f\in {\cal F}$, and such that $\tau(\e-\e_C)=1/m<\dt$ for all $\tau \in T(M_{k_{n+1}}(B))$. Then $\iota_{n+1,\infty}(C)$ is the desired subalgebra. (Note that $1-\e_C\lesssim a$ follows from strict comparison property of $B\otimes U$.)
It follows that $B \otimes U\in {\cal B}_0.$
 \end{proof}

%\end{NN}

\begin{thm}\label{range 0.35}
For any simple weakly unperforated Elliott invariant $\big((G,G_+,u),~
K, \DT, r \big)$, there is an unital simple algebra $A\in {\cal
N}_0^{\cal Z}$ which is an inductive limit of $(A_n, \phi_{n,m})$ with $A_n$
described in the end of \ref{range 0.19}, with $\phi_{n,m}$
injective,  such that
$$\big((K_0(A),K_0(A)_+, \e_A),~ K_1(A), T(A), r_{\!\!_A} \big)\cong
\big((G,G_+,u),~ K, \DT, r \big).$$
\end{thm}

\begin{proof} By \ref{range 0.34a}, $A\in {\cal N}_0$. Since $A$ is a unital simple inductive limit of subhomogeneous C*-algebras with no dimension growth. By Corollary 6.5 of \cite{Winter-Z-stable-01}, $A$ is ${\cal Z}$-stable.
\end{proof}

\begin{cor}\label{range 0.35a} For any simple weakly unperforated Elliott invariant $\big((G,G_+,u),~
K, \DT, r \big)$ with $K=0$ and $G$ torsion free,there is an unital ${\cal Z}$-stable simple algebra which is an inductive limit of $(A_n, \phi_{n,m})$ with $A_n$ in ${\cal C}_0$ described in \ref{DfC1}, with $\phi_{n,m}$
injective,  and such that
$$\big((K_0(A),K_0(A)_+, \e_A),~ K_1(A), T(A), r_{\!\!_A} \big)\cong
\big((G,G_+,u),~ 0, \DT, r \big).$$

\end{cor}

\begin{proof} In the construction of $A_n$, just let all the spaces $X_n$ involved to be the space of a single point.

\end{proof}

\section{Models for \CA s in ${\cal A}_0$ with (SP) property}

\begin{NN}\label{range 0.37} For technical  reasons, in the construction of our model algebras,
it is important for us to be able to decompose
$A_n$ into direct sum of two parts: the homogeneous part
which stores the information of $\Inf K_0(A)$ and $K_1(A)$ and the
part of  algebra in ${\cal C}_0$ which stores information of
$K_0(A)/\Inf K_0(A)$, $T(A)$ and paring between them. But this can not
be done in general for the algebras in ${\cal N}_0$ (see \ref{Class0}), since some
algebras in ${\cal N}_0$ has minimal projections (or even the unit
itself is a minimal projection). But we will prove this can be done
if the Elliott invariant satisfies an extra condition (SP) described
below. Note that if $A\in {\cal N}_0$ then the Elliott invariant of $A\otimes M_{\mathfrak{p}}$ has
the condition  (SP) even though the Elliott invariant of $A$
%the Elliott invariant of all  algebras in ${\cal
%N}_0$ satisfy the condition (SP) after tensor with $M_{\mathfrak{p}}$ for
%a supernatural
%%%supernature
% number $\mathfrak{p}$ of infinite type even though
%the algebra
itself may not satisfy the condition.

  Let $\big((G,G_+,u),~
K, \DT, r \big)$ be a weakly unperforated Elliott invariant as
\ref{range 0.1}. We say that it has the (SP) property if for any
real number $s>0$, there is $g\in G_+\setminus \{0\}$ such that
$\tau (g) <s$ for any state $\tau$ on $G$, or equivalently,
$r(\tau)(g)<s$ for any $\tau \in \DT$. In this case, we will prove
that the algebra in \ref{range 0.35} can be chosen to be in class
${\cal B}_0$ (rather than in the larger class ${\cal N}_0=\{A{:} \,
A\otimes M_{\mathfrak{p}}\in {\cal B}_0\}$). Roughly speaking, for
each $A_n$, we will separate the part of homogeneous algebra which
will store all the information of infinitesimal part of $K_0$  and
$K_1$, and it will be in the corner $P_nA_nP_n$ with $P_n$ small compare to
$\e_{A_n}$ in the limit algebra. In fact, the construction of this
case is much easier, since the homogeneous blocks can be separate
out from the part of ${\cal C}_0$---we will first write the group
inclusion $G_n\hookrightarrow H_n$ as in \ref{range 0.5a}.

\end{NN}

\begin{NN}\label{range 0.38}
Let $((G,G_+,u),K,\Delta,r)$ be the one given in \ref{range 0.1} or
\ref{range 0.16}. As in \ref{range  0.3}, let $\rho:~ G\to \Aff\DT$
be dual to the map $r$. Denote by the kernel of the map $\rho$ by
$\Inf(G)$---the infinitesimal part of $G$, that is
$$\Inf(G)=\{g\in G: \rho(g)(\tau)=0,\ \forall \tau\in\DT\}.$$
Let $G^1\subseteq \Aff \DT$ be a countable dense subgroup which is $\Q$ linearly
independent with $\rho(G)$---that is, if $g\in \rho(G)\otimes\Q$ and
$g^1\in G^1\otimes\Q$ satisfy $g+g^1=0$, then both $g$ and $g^1$ are
zero. Note that such $G^1$ exist, since as $\Q$ vector space,
dimension of $\rho(G)\otimes \Q$ is countable, but dimension of $\Aff
\DT$ is uncountable.  Again as in \ref{range 0.3}, let $H=G\oplus
G^1$ with $H_+\setminus \{0\}$ to be the {set} of $(g,f)\in
G\oplus G^1$ with
$$\rho(g)(\tau)+f(\tau)>0  \rforal \tau\in \DT.$$
The scale $u\in G_+$ is regarded as $(u,0)\in G\oplus G^1=H$
as the scale of $H_+$. Since $\rho(u)(\tau)>0,$ it follows  that $u$ is an order unit for $H.$ Since $G^1$ is $\Q$ linearly
independent with $\rho(G)$, we know $\Inf(G)=\Inf(H)$---that is, when
we embed $G$ into $H$, it does not create more elements in the
infinitesimal group. Evidently ${\rm Tor}(G)={\rm Tor}(H)\subset \Inf(G)$. Let
$G'=G/\Inf(G)$ and $ H'=H/\Inf(H)$, then we have the following
diagram:

\begin{displaymath}
    \xymatrix{
       0 \ar[r] & \Inf(G) \ar[r]\ar@{=}[d]&G
         \ar[r]\ar@{_{(}->}[d] &{ G'}\ar[r]\ar@{_{(}->}[d]&0 \\
         0 \ar[r] & \Inf(H) \ar[r]&H
         \ar[r] &{ H'}\ar[r]&0.}
\end{displaymath}

Let ${ G'}^+$ (or $H'^+$) and $u'$  be the image of $G_+$ (or $H_+$)
and $u$ under the quotient map from $G$ to $G'$ (or from $H$ to
$H'$). Then $(G', G'^+, u')$ is weakly unperforated group without
infinitesimal elements.  Note that $G$ and $H$ share same unit $u$, and
therefore $G'$ and $H'$ share same unit $u'$. Since
$r(\tau)|_{\Inf(G)}=0$ for any $\tau \in \DT$,  the map $r:\DT \to
S_u(G)$ induces a map $r': \DT \to S_{u'}(G')$. Hence $((G', G'_+, u'),
\{0\}, \DT, r')$ is a weakly unperforated Elliott Invariant with trivial
$K_1$ group and no infinitesimal elements in the $K_0$-group.
\end{NN}

\begin{NN}\label{range 0.39}
   With the same argument as that of \ref{range 0.5}, we have the following diagram of
   inductive limit{s}:
\begin{displaymath}
    \xymatrix{{ G'}_1 \ar[r]^{\af'_{_{12}}}\ar@{_{(}->}[d]_{\iota_1} & {G'}_2 \ar[r]^{\af'_{_{23}}}\ar@{_{(}->}[d]_{\iota_2}&\cd \ar[r]&{ G'} \ar@{_{(}->}[d]_{\iota}\\
         { H'}_1 \ar[r]^{\gm'_{_{12}}}  & { H'}_2 \ar[r]^{\gm'_{_{23}}}&\cd \ar[r]&{ H'},  }
\end{displaymath}
where each ${ H'}_n$ is a   direct sum of finite copies of ordered
groups $\Z$, $\af_{n, n+1}=\gm_{n, n+1}|_{{G'}_n}$, and
${H'}_n/{G'}_n$ is a free abelian group.

By \ref{range 0.15a}, we can construct an increasing sequence of
finitely generated subgroups
$$\Inf_1\subset \Inf_2\subset \Inf_3\subset\cd\subset \Inf(G),$$
with $\Inf(G) =\bigcup_{i=1}^{\infty} \Inf_n$ and the inductive limit
\begin{displaymath}
\xymatrix{
\Inf_1\oplus H'_1 \ar[r]^-{\gm_{1,2}} &  \Inf_2\oplus H'_2 \ar[r]^-{\gm_{2,3}} & \Inf_3\oplus H'_3 \ar[r]^-{\gm_{3,4}} & \cd \ar[r] & H
}
\end{displaymath}

%$$\Inf_1\oplus H'_1 \stackrel{\gm_{1,2}}\rightarrow
%\Inf_2\oplus H'_2 \stackrel{\gm_{2,3}}\rightarrow \Inf_3\oplus
%H'_3\stackrel{\gm_{3,4}}\rightarrow\cd\rightarrow H.$$

Put $H_n:=\Inf_n\oplus
H'_n$ and $G_n=\Inf_n\oplus G'_n.$ Since  $G'_n$ is a subgroup of
$H'_n$, the group   $G_n$ is also a subgroup of $H_n$. Define $\af_{n,n+1}: G_n\to
G_{n+1}$ by $\af_{n, n+1}=\gm_{n,n+1}|_{G_n}$, which is
compatible with $\af'_{n,n+1}$ in the sense that (ii) of \ref{range
0.15a} holds. Hence we get the following diagram of inductive limits:

\begin{displaymath}
    \xymatrix{{ G}_1 \ar[r]^{\af_{_{12}}}\ar@{_{(}->}[d]_{\iota_1} & {G}_2 \ar[r]^{\af_{_{23}}}\ar@{_{(}->}[d]_{\iota_2}&\cd \ar[r]&{ G} \ar@{_{(}->}[d]_{\iota}\\
         { H}_1 \ar[r]^{\gm_{_{12}}}  & { H}_2 \ar[r]^{\gm_{_{23}}}&\cd \ar[r]&{ H},  }
\end{displaymath}
 with $\af_{n,n+1}(\Inf_n)\subset \Inf_{n+1}$ and $\af_{n,n+1}|_{(\Inf_n)}$
is the inclusion map.

Note that all notations discussed so far in this section will be used for the rest of this section.

\end{NN}

 \begin{lem}\label{range 0.40} Let $(G, G_+, u)=\lim ((G_n,(G_n)_+,u_n), \af_{n,m})$
and $(H,H_+,u)=\lim((H_n,(H_n)_+,u_n), \gm_{n,m})$  be as above. Again
suppose that $\big((G,G_+,u),~ K, \DT, r \big)$ has (SP) property.

 For any $n$ with
$G_n\stackrel{\iota_n}{\hookrightarrow}H_n$,  any positive integer
$L$, and for any $D=\Z^k$ (for $k$ an arbitrary positive integer),
there are $m$ and positive maps $(\kappa, \id): H_n \to D\oplus H_n$,
$(\kappa', \id): G_n \to D\oplus G_n$, $\eta: D\oplus H_n \to H_m$,
and $\eta': D\oplus G_n \to G_m$ such that the following diagram
commutes:
\begin{displaymath}
\xymatrix{
G_n \ar[rr]^{\alpha_{n, m}} \ar[dr]_{({\kappa}', \id)} \ar[ddd]_{\iota_n} &  &  G_m \ar[ddd]^{\iota_m} \\
 &D \oplus G_n  \ar[ur]_{\eta'} \ar@{_{(}->}[d] _{(\id, \iota_n)}&   \\
 &D \oplus H_n \ar[dr]^{\eta}&  \\
H_n\ar[ur]^{({\kappa}, \id)} \ar[rr]_{\gm_{n, m}} &   & H_m, }
\end{displaymath}
and such that the following are true:

{\rm (1)} For any positive element $x\in H_n'^+$, each component of
$\kappa(x)$ in $\Z^k$ is strictly positive and consequently, each
component of $\kappa (u_n)=\kappa'(u_n)$ in $\Z^k$ is strictly
positive.

{\rm (2)} For any $\tau\in \DT$,
$$r(\tau)((\af_{m,\infty}\circ\eta')(\e_D))(=r(\tau)((\gm_{m,\infty}\circ\eta)(\e_D)))<1/L.$$

{\rm (3)} Each component of the map $\eta:D\oplus H'_n
=\Z^k\oplus\Z^{p_n}\to H'_m=\Z^{p_m}$ is  L-large---that is all
entries in the $(k+p_n)\times p_m$ matrix corresponding to the map
are  larger than $L$.

\end{lem}

\begin{proof}
We will use the following  fact several times: the positive cone of
$G_n'$ (and of $H'_n$) is
 finitely generated (note that
 even though $G_n$ and $H_n$ are finitely generated, their positive cone may not be
 finitely generated). For  $H_n$, pick an arbitrary positive nonzero homomorphism $\ld:
H_n\to\mathbb Z$ so that for any nonzero $x\in H_n'^+$, $\ld(x)>0$.
Denote by $\ld'=\ld\circ\iota_n: G_n\to \mathbb Z$. It follows from
positivity that  such map $\ld$ satisfies $\ld(\Inf_n)=0$.

  Since $G$ has the (SP) property,
 there is $p'\in G_+\setminus\{0\}$ such
 that for any $a\in (G_n')_+$ (not use  $(G_n)_+$ here, but we regard it as subset of $(G_n)_+$)
 the element
$$\alpha_{n, \infty}(a)-k\cdot\ld'(a)\cdot p'$$ is positive and for
any $a\in H_n'^+$ (not use $(H_n)_+$ here, but again we regard it as
subset of $(H_n)_+$), the element
$$\gm_{n, \infty}(a)-k\cdot\ld(a)\cdot p'$$ is positive, where
the map $\alpha_{n, \infty}$ and $\gm_{n, \infty}$ are the
homomorphisms from $G_n$ to $G$ and from $H_n$ to $H$ respectively.
Moreover, one may require that
\begin{equation}\label{14star}
r(\tau)(\ld(u_n)\cdot p')<1/2kL  \rforal \tau\in
\DT.
\end{equation}

Since $(G_n')_+$ and $H_n'^+$ are finitely generated, there is an integer $m$
and $p\in G_m^+$ such that 
%$$
\beq\nonumber
\af_{m,\infty}(p)=p',\,\,\,
\alpha_{n, m}(a)-k\cdot\ld'(a)\cdot p \in G_m^+ \rforal a\in (G_n')_+\andeqn\\\nonumber
%and
%$$
\gm_{n, m}(a)-k\cdot\ld(a)\cdot p \in H_m^+  \rforal a\in
H_n'^+.
\eneq
Then define $\tilde{\alpha}_n: G_n \to G_m$ and $\tilde{\gm}_n: H_n
\to H_m$ by
\beq\nonumber
&&\tilde{\alpha}_n: G_n\ni a\mapsto \alpha_{n, m}(a)-k\cdot\ld'(a)\cdot p\in G_m\andeqn\\
&&\tilde{\gm}_n: H_n\ni a\mapsto\gm_{n, m}(a)-k\cdot\ld(a)\cdot p\in H_m.
\eneq
By the choice of $p$, the maps $\tilde{\alpha}_n$ and
$\tilde{\gm}_n$ are positive. (Note that
$\tilde{\alpha}_n=\af_{n,m}$ and $\tilde{\gm}_n=\gm_{n,m}$ on
$\Inf_n$.)

%  $p\in G'_2$, $\tilde{\alpha}_{1}(G_1')_+\subseteq
%(G_2')_+$, and $\tilde{\beta}_{1}(H_1')_+\subseteq (H_2')_+$.

 A direct calculation shows the following diagram commutes:
\begin{displaymath}
\xymatrix{
G_n \ar[rr]^{\alpha_{n, m}} \ar[dr]_{({\kappa}', \id)} \ar[ddd]_{\iota_n} &  &  G_m \ar[ddd]^{\iota_m} \\
 &D \oplus G_n  \ar[ur]_{\eta'} \ar@{_{(}->}[d] _{(\id, \iota_n)}&   \\
 &D \oplus H_n \ar[dr]^{\eta}&  \\
H_n\ar[ur]^{({\kappa}, \id)} \ar[rr]_{\gm_{n, m}} &   & H_m, }
\end{displaymath}
where $D=\Z^k,$
$$\kappa'(a)=(\underbrace{\ld'(a), ..., \ld'(a)}_k)\in D\andeqn
\kappa(a)=(\underbrace{\ld(a), ..., \ld(a)}_k)\in D,$$
$$\eta'((m_1, ..., m_k, g))= (m_1+\cdots+m_k)p+\tilde{\alpha}_n(g)\andeqn$$
$$\eta_1((m_1, ..., m_k, g))= (m_1+\cdots+m_k)p+\tilde{\gm}_n(g).$$
%{\bf (change notation here.)}
 The order 
 %on the groups 
 of $D\oplus G_n$  and $D\oplus H_n$ are the standard
  order on direct sums, i.e., $(a, b)\geq 0$ if and only if $a\geq 0$ and $b\geq 0$.
   Since the maps $\tilde{\alpha}_n$ and $\tilde{\gm}_n$ are positive, the
   maps $\eta'$ and $\eta$ are positive. Condition (1) follows from
   the {construction;} Condition (2) follows from (\ref{14star}), and Condition
   (3) follows from the simplicity of $H$,  if one passes to further stage (choose larger
    $m$).
\end{proof}

\begin{NN}\label{range 0.41}
Write $K$ (the $K_1$ part of invariant) as the union of increasing
sequence of finitely generated abelian subgroups:  $K_1\subset K_2 \subset
K_3\subset\cd\subset K$ with $K=\bigcup_{i=1}^{\infty} K_i$.

For a finitely generated abelian group $G$, we use ${\rm rank}\,
G$ to denote the minimum number of possible generated set of
$G$---that is, $G$ can be written as a direct sum of
${\rm rank}(G)$ cyclic groups (e.g., $\Z$ or $\Z/m\Z$, $m\in \N$).

Let $d_n=\mbox{max}\{2,~1+\mathrm{rank}(\Inf_n)+\mathrm{rank}(K_n)$\}. Apply \ref{range
0.40} with $k=d_n$ (and suitable choice of $L=L_n > 13\cdot 2^n$ to be determined later (\ref{range 0.43-3}))
for each $n$ and pass to subsequence, we can obtain the following
diagram of inductive limits:
\begin{displaymath}
\xymatrix{
\Z^{d_1} \oplus G_1\ar[r]^{\tilde{\alpha}_{1, 2}}\ar@{_{(}->}[d]_{(\id, \iota_1)} & \Z^{d_2} \oplus G_2 \ar@{_{(}->}[d]_{(\id, \iota_2)} \ar[r]^-{\tilde{\alpha}_{2, 3}} & \cdots \ar[r] & G \ar@{_{(}->}[d]_{\iota} \\
\Z^{d_1} \oplus H_1\ar[r]_{\tilde{\gm}_{1, 2}} & \Z^{d_2} \oplus H_2
\ar[r]_-{\tilde{\gm}_{2, 3}} & \cdots \ar[r] & H.}
\end{displaymath}

Consider the ordered group $\Z^{d_n} \oplus G_n= \Z^{d_n} \oplus \Inf_n\oplus
G_n'$ with the scale $(\kappa'(u_n), u_n)$, where
$\kappa'$ is as in \ref{range 0.40}; and let us still denote it by $(G_n, u_n)$. Similarly, consider $
\Z^{d_n} \oplus H_n= \Z^{d_n} \oplus \Inf_n\oplus H_n'$ with the scale
$(\kappa(u_n), u_n)$, where $\kappa$ is
also as in \ref{range 0.40}; and let us still denote it by $(H_n, u_n)$. We will also use $\af_{n,n+1}$ and
$\gm_{n,n+1}$ for $\tilde{\alpha}_{n, n+1}$ and $\tilde{\gm}_{n,
n+1}$ in the above diagram. Let $G_n''=\Z^{d_n}\oplus \Inf_n$, then
with the new notation, we have $G_n =G_n''\oplus G_n'$ and $H_n
=G_n''\oplus H_n'$. Now we have the following diagram
\begin{displaymath}
\xymatrix{
G_1'' \oplus G_1'\ar[r]^{{\alpha}_{1, 2}}\ar@{_{(}->}[d]_{(\id, \iota_1)} &G_2'' \oplus G_2' \ar@{_{(}->}[d]_{(\id, \iota_2)} \ar[r]^-{{\alpha}_{2, 3}} & \cdots \ar[r] & G \ar@{_{(}->}[d]_{\iota} \\
G_1'' \oplus H_1'\ar[r]_{{\gm}_{1, 2}} & G_2''\oplus H_2'
\ar[r]_-{{\gm}_{2, 3}} & \cdots \ar[r] & H.}
\end{displaymath}

The positive cones of $G_n'$ and $H_n'$ are as before.  Write
$G_n''=\bigoplus_{i=1}^{d_n}(G_n'')^i$, with $(G_n'')^i=\Z$ for $i\leq
1+\mathrm{rank}(K_n)$ and $(G_n'')^i=\Z\oplus (\mbox{a cyclic group})$ for
$1+\mathrm{rank}(K_n)<i\leq d_n$, and the direct sum of those cyclic groups is
$\Inf_n$. Here the positive cone of $(G_n'')^i$ is given by the strict
positivity of first coordinate for nonzero positive elements.  And an
element in  $G_n''$ is positive if
%its each
{each of its} components in $(G_n'')^i$
is positive. The order on the groups $G_n''\oplus G_n'$
and $G_n''\oplus H_n'$ are the standard
  order on direct sums, i.e., $(a, b)\geq 0$ if and only if $a\geq 0$ and $b\geq 0$.
  Since
  % the
  each
  %entries
  {entry} of $\gm_{n,n+1}: \Z^{d_n}
\oplus H_n' \to \Z^{d_{n+1}} \oplus H_{n+1}'$ is strictly positive,
the infinitesimal part can be put in any given summand without
affect the order structure of the limit. Let $u_n=(u_n'',u_n')\in
G_n''\oplus G_n'\subset H_n'\oplus H_n''$ be the unit of $G_n$ and
of $ H_n$. Note that
\beq\label{0.41a}
\Inf(G) \subset \bigcup_{n=1}^{\infty}(\Inf_n)\subset \bigcup_{n=1}^{\infty} \alpha_{n,\infty}(G''_n).
\eneq

\end{NN}

\begin{df}\label{AHblock}
A \CA\, is said to be in the class ${\bf H}$ if it is the direct sum
of the algebras of  the form $P(C(X)\otimes M_n)P$, where $X=\{pt\},
[0,1], $ $S^1, $ $S^2,$ $T_{2,k}$ and $T_{3, k}.$
\end{df}

\begin{NN}\label{range 0.43}
For each $n$ as in \ref{range 0.5} applied to $G_n'\subset H_n'$,
one can find  finite dimensional $C^*$-algebras $F_n$ and $E_n$,
unital homomorphisms $\bt_0,~\bt_1: F_n \to E_n$, and the C*-algebra
$$C_n=A(F_n,E_n,\bt_0,\bt_1):= \{(f,a)\in C([0,1], E_n)\oplus F_n;
f(0)=\bt_0(a), f(1)=\bt_1(a)\}$$ such that
$$(K_0(F_n),K_0(F_n)_+, [\e_{F_n}])=(H_n', H_n'^+,u_n'),$$
$$(K_0(C_n), K_0(C_n)_+,\e_{C_n})=(G_n', (G_n)_+,u_n'),\quad K_1(C_n)=\{0\},$$
and furthermore $K_0(C_n)$ is identified with
$$\ker((\bt_1)_{*0}-(\bt_0)_{*0})=\{x\in K_0(F_n); ((\bt_0)_{*0}-(\bt_1)_{*0})(x)=0\in
K_0(E_n)\}.$$

We can also find a unital $C^*$ algebra $B_n \in {\bf H}$ such that
\beq\label{0.43a}
(K_0(B_n), K_0(B_n)_+, \e_{B_n})=(G_n'',G_n''^+,u_n'')
\eneq
 and
$K_1(B_n)=K_n$. More precisely, we have that $B_n=\bigoplus_{i=1}^{d_n}B_n^i$,
with $K_0(B_n^i)=(G_n'')^i$ and $K_1(B_n^i)$ is either  a cyclic group
for the case $2\leq i\leq 1+{\rm rank} (K_n)$ or zero for the other cases.
In particular, the algebra $B_n^1$ can be chosen to be a matrix
algebra over $\C$. And we assume, for at least one block $B_n^2$,
the spectrum is not a single point (note that $d_n\geq2$),
otherwise, we will replace the single point spectrum by interval
$[0,1]$.

Now, we can extend the maps $\bt_0$ and $\bt_1$ to $\bt_0,~\bt_1:
B_n\oplus F_n \to E_n$, by defining them to be zero on $B_n$.

Consider $A_n=B_n\oplus C_n$. Then the C*-algebra $A_n$ can be written as
$$A_n=\{(f,a)\in C([0,1],E_n)\oplus (B_n\oplus F_n); f(0)=\bt_0(a), f(1)=\bt_1(a)\}.$$
%{?} This is very similar to \ref{range 0.19}, of course $B_n$ is like
%the first block $F_n^1$ of $F_n$ in \ref{range 0.19} which is not
%matrix algebras over $\C$ any more.{?}

For each block $B_n^i$, choose a base point $x_{n,i}\in Sp(B_n^i)$.
Let $I_n=C_0\big((0,1),E_n\big)$ be the ideal of $C_n$ as before.
And let $J_n$ be the ideal of $B_n$ consisting of functions vanishing
on all base points $\{x_{n,i}\}_{i=1}^{d_n}$. Applying Proposition
\ref{range 0.14}, and
%completely similar to
exactly as in \ref{range 0.6} and
\ref{range 0.20}, we will have non-simple inductive limit
$A'=\lim(B_n\oplus C_n, \phi_{n,m})$ with
$(\phi_{n,n+1})_{*0}=\af_{n,n+1}$, $K_1(\phi_{n,n+1})$ is the inclusion
from $K_n$ to $K_{n+1}$, ${\rm Ell}(A)\cong ((G, G_+, u), K, \DT,
r),$ and with $\phi_{n,n+1}(I_n)\subset
I_{n+1},~~\phi_{n,n+1}(J_n)\subset J_{n+1}$. Furthermore the map
$\pi\circ\phi_{n,n+1}|_{B_n}$  is injective for projection
$\pi:B_{n+1}\oplus C_{n+1} \to B_{n+1}$, since at least one block of
$B_{n+1}$ is not a single point. (We do not need \ref{range 0.21}
to \ref{range 0.25} and \ref{range 0.34} because the homogeneous
algebra $B_n$ is separate from $C_n\in {\cal C}_0$. It is also true that $B_n$  occupies
relatively smaller space   compare to those occupied by $C_n$ in the
limit algebra---namely, by (2) of \ref{range 0.40} and the choice of $L=L_n\geq 13\cdot 2^n$ in \ref{range 0.41}, one has
$$\frac{\tau(\alpha_{n,\infty}(\e_{B_n}))}{\tau(\alpha_{n,\infty}(\e_{C_n}))}<\frac{1}{13\cdot 2^n-1},$$
for all $\tau$.)

 We need to modify $\phi_{n,n+1}$ to $\psi_{n,n+1}$ to make the
algebra simple. However,  it is much easier than what we did in \ref{range
0.29} to \ref{range 0.33}.
% to modify the construction to make the
%algebra simple.
We only need to modify the partial map from $A_n$ to
$B_{n+1}^1$, the first block of $A_{n+1}=B_{n+1}\oplus C_{n+1}$, and
keep other part of the map to be same---that is, we only need to make
$SP(\psi_{n,n+1}|_{x_{n+1,1}})$ to be dense enough in $Sp(A_n)$. In order to do this, one chooses a finite set $X\subset Sp(A_n)$ dense enough to play the
role of $Y\cup T$ (as in \ref{condition2}) and let $L_n$ in
\ref{range 0.41} satisfies
\beq\label{range 0.43-3}
L_n > 13\cdot
2^n\cdot(\#(X))\cdot(\max\{\mathrm{size}(B_n^i), \mathrm{size}(F_n^i), \mathrm{size}(E_n^i)\}).
\eneq

(Note in \ref{range 0.29} to \ref{range 0.33}, we modify the set
$SP(\phi_{n,n+1}|_{\tht_2'})$, which will force us to change the
definition of the map for the point in $Sp(I_{n+1})$. But now, the
maps $\bt_0$ and $\bt_1$ in the construction of $A_{n+1}$ are zero on
$B_{n+1}$ and $\{x_{n+1,1}\}$ is an isolated point in $Sp(B_{n+1})$ (and
is also isolated in $Sp(A_{n+1})$), so the modification of
$SP(\phi_{n,n+1}|_{x_{n+1,1}})$ will not affect other points; see
the end of  \ref{range 0.31} also.) {\it Let us emphasis that,  in our original construction of $A_{n+1}$, we can assume the $L_n$ involved in the construction satisfying the above condition. Therefore the algebra $A_{n+1}$ will not be changed, when we modify the connecting maps to make the limit algebra simple.} We get the following theorem:
%{\bf I would again ask whether  $A_{n+1}$ is changed in the process.}

\end{NN}

\begin{thm}\label{RangT}
Let $((G, G_+, u), K, \DT, r)$ be  a six-tuple of the following
objects: $(G, G_+, u)$ is a weakly unperforated simple order-unit
group, $K$ is a countable abelian group, $\DT$ is a 
separable 
%compact
Choquet simplex and $r:  \DT\to \mathrm{S}_u(G) $ is surjective
affine map, where $\textrm{S}_u(G)$ the compact convex set of the
states on $(G, G_+, u)$. Assume that $(G, G_+, u)$ has the (SP)
property in the sense that for any real number $s>0$, there is $g\in
G_+\setminus\{0\}$ such that $\tau(g)<s$ for any state $\tau$ on $G$.

Then there is a unital  simple C*-algebra $A\in {\cal B}_0$ which can be written
as $A=\lim_{n\to\infty}(A_i, \psi_{i, i+1})$ with injective
$\psi_{i, i+1}$, where $A_i=B_i\oplus C_i,$ $B_i\in {\bf H},$
$C_i\in {\cal C}_0$ with $K_1(C_i)=\{0\}$ such that
\begin{enumerate}
\item $\lim_{i\to\infty}\sup\{\tau(\psi_{i, \infty}(1_{B_{{{i}}}})): \tau\in T(A)\}=0,$
\item $\ker \rho_A\subset \bigcup_{i=1}^{\infty} [{\psi}_{i, \infty}]_0({\rm ker}\rho_{B_i})$, and
\item ${\rm Ell}(A)\cong ((G, G_+, u), K, \DT, r).$
\end{enumerate}

Moreover, the inductive system $(A_i, \psi_i)$ can be chosen so that $\psi_{i, i+1}=\psi_{i, i+1}^{(0)}\oplus\psi_{i, i+1}^{(1)}$ with $\psi_{i, i+1}^{(0)}: A_i\to A_{i+1}^{(0)}$ and $\psi_{i, i+1}^{(1)}: A_i\to A_{i+1}^{(1)}$ for C*-subalgebras $A_{i+1}^{(0)}$ and $A_{i+1}^{(1)}$ of $A_{i+1}$ with $1_{A_{i+1}^{(0)}}+1_{A_{i+1}^{(0)}}=1_{A_{i+1}}$ % (not necessary with respect to the decomposition $A_i=B_i\oplus C_i$)
such that
%\begin{enumerate}
%\item
%\label{fd-not0}
$A_{i+1}^{(0)}$ is a nonzero finite dimensional C*-algebra,  %$[\psi_{i, i+1}^{(0)}]$ is a point evaluation (in particular, it is homotoplically trivial), and
%\item $\tau(\psi_{i+1, \infty}(1_{A_{i+1}^{(0)}}))\neq 0$ uniformly on $\tau\in T(A)$,
%\item
and $[\psi_{i, \infty}]_1$ is injective.
%\end{enumerate}

\end{thm}

\begin{proof} Condition (2) follows from (\ref{0.41a}) and (\ref{0.43a}).
Now the only item that has not been proved is the assertion that
%We only need to prove
$A\in {\cal B}_0$.  By (1) above and the fact that $C_i \in {\cal C}_0,$ it remains to show that
$A$ has strict comparison for projections. This actually follows from our construction immediately.  {Note, by \ref{2Tg16},  $C_i$ has strictly comparison.  Moreover, dimensions of the underline spaces  do not increase.
A standard argument shows that the inductive limit has strict comparison for projections ( or for positive elements).}
Another way to see this is that since $A$ is a unital simple inductive limit of subhomogeneous C*-algebras with no dimension growth,
it then follows from Corollary 6.5 of \cite{Winter-Z-stable-01} that $A$ is ${\cal Z}$-stable.
Hence $A$ has strict comparison for positive elements.
%{{\bf Condition (2) is not discussed (at least not specified) in the construction of the group.
%Perhaps, in the proof, we should deal with this requirements.  Please note that, even inductive limits
%constructed like $F_n\to F_{n+1}$ could have infinitesimal elements. In other words, we should take more care
%of the issue.}}
\end{proof}

%It follows from \ref{2pg3} and a standard argument, the $C^*$ algebra $A$ has stable rank one and consequently has property of strict comparison {\color{Green} (Stable rank one does not imply comparison; for instance, any AH-algebra with diagonal maps has stable rank one, but some of them have perforation in K-group or Cuntz semigroup. Moreover, it is not straightforward to show that $A$ has stable rank one.)}. Hence it follows from (1) above and $C_i \in {\cal C}_0$ that $A\in {\cal B}_0$.
%\end{proof}

%{\color{Green}
%\begin{proof}
%%We only need to show that $A\in {\cal B}_0$, and we only need to show that $A$ has strict comparison on projections. But this follows from the fact that $A$ has no dimension growth. (I regard it well known that any simple locally ASH algebra with slow dimension growth has strict compariso---very easy to prove for projections, but extra work are needed for Cuntz semigroups.)
%\end{proof}
%}

\begin{rem}
Note that   $A_{i+1}^{(0)}$ can be chosen to be the first block
$B_{i+1}^1$, so we have
$$\lim_{i\to\infty}\tau(\psi_{i+1, \infty}(1_{A_{i+1}^{(0)}}))=0$$
uniformly for $\tau\in T(A)$.
\end{rem}

\begin{rem}\label{range 0.46}
Let $\xi_{n,n+1}$ be the partial map of $\psi_{n,n+1}:B_n \to
B_{n+1}$, and let $\xi_{n,m}: B_n \to B_m$ be the corresponding composition
$\xi_{m-1,m}\circ \xi_{m-2,m-1}\circ \cd \circ \xi_{n,n+1}$. Let
$e_n=\xi_{1,n}(\e_{B_1}).$ Then, from the construction, we know that
the algebra $B=\varinjlim(e_nB_n e_n, \xi_{n,m})$ is simple, as we know
that $SP(\xi_{n,n+1}|_{x_{n,1}})$ is dense enough in $Sp(B_n)$. Note
that the simplicity of $B$ does not follow from simplicity of $A$
itself, since it is not a corner of $A$.

\end{rem}

The following is an analog of Theorem 1.5 of \cite{Lnann}.

\begin{cor}\label{smallmap}
Let $A_1$ be a simple separable C*-algebra in $\mathcal B_1$,
 and let $A=A_1\otimes U$ for {an infinite dimensional}  UHF-algebra $U$. There exists an
  inductive limit algebra $B$ as constructed in Theorem \ref{RangT} such that $A$
  and $B$ have the same  Elliott invariant{.}
  % including}
   %scale ordered K-theory.
   Moreover, the C*-algebra $B$ satisfies the following properties:

Let $G_0$ be a finitely generated subgroup of $K_0(B)$ with decomposition $G_0=G_{00}\oplus G_{01}$, where $G_{00}$  vanishes under all states of $K_0(A)$. Suppose $\mathcal P\subset \underline{K}(B)$ is a finite subset which generates a subgroup $G$ such that $G_0\subset G\cap K_0(B)$.

Then, for any $\epsilon>0$, any finite subset $\mathcal F\subset B$, any $1>r>0$, and any positive integer $K$, there is an
$\mathcal F$-$\epsilon$-multiplicative map $L:B\to B$ such that:
\begin{enumerate}
\item $[L]|_{\mathcal P}$ is well defined.
\item $[L]$ induces the identity maps on the infinitesimal part of  $G\cap K_0(B)$, $G\cap K_1(B)$,
      $G\cap K_0(B,\mathbb Z/k\mathbb Z)$ and $G\cap K_1(B, \mathbb Z/k\mathbb Z)$
      for $k=1, 2, ...$, and $i=0,1$.
\item $\rho_B\circ[L](g)\leq r\rho_B(g)$ for all $g\in
      G\cap K_0(B)$, where $\rho$ is the canonical positive homomorphism from $K_0(A)$ to $\Aff(\mathrm{S}(K_0(A), K_0(A)_+, [1_A]_0))$.
\item For any positive element $g\in G_{01}$, we have $g-[L](g)=Kf$ for some $f\in K_0(B)_+$.
\end{enumerate}
\end{cor}
\begin{proof}
%{ We can write $A=A'\otimes M_{\mathfrak{p}}$ and
%A'=M_{\mathfrak{p'}}$ with both super nature numbers $\mathfrak{p}
% and $\mathfrak{p'}$ of infinite type. Replacing $A_1$ by $A$, we
%can assume that $Ell(A_1)$ has (SP) property}.
{\Wlog, by replacing $A_1$ by $A_1\otimes U,$ we may assume that ${\rm Ell}(A_1)$ has (SP) property.}

 Consider ${\rm Ell}(A_1)$,
%Since $A_1\in \mathcal B_1$, the $(K_0(A_1), K^+_0(A_1), [1_A])$
which satisfies the condition of Theorem \ref{RangT}, and therefore
by the first part of Theorem \ref{RangT}, there is a inductive
system $B_1=\varinjlim(T_i\oplus S_i, \psi_{i, i+1})$ such that
%\begin{enumerate}
%\item

(i) $T_i\in\mathbf H$ and $S_i\in\mathcal C_0$
with $K_1(S_i)=\{0\}$,

%\item
(ii) $\lim \tau(\phi_{i, \infty}(1_{T_i}))= 0$ uniformly on $\tau\in T(B_1)$,

%\item
(iii) $\ker(\rho_{B_1})=\bigcup_{i=1}^\infty (\psi_{i, \infty})_{*0}(\ker(\rho_{T_i}))$, and

%\item
(iv) ${\rm Ell}(B_1)={\rm Ell}(A_1)$.
%\end{enumerate}

Put $B=B_1\otimes U$. Then ${\rm Ell}(A)={\rm Ell}(B)$. Let $\mathcal P\in \underline{K}(B)$ be a finite subset, and let $G$ be the subgroup generated by $\mathcal P$ which contains $G_0$. Then there is a positive integer $M'$ such that $G \cap K_*(B, \mathbb Z/ k\mathbb Z)={\{0\}}$ if $k>M'$. Put $M=M'!$. Then $Mg=0$ for any $g\in G \cap K_*(B, \mathbb Z/ k\mathbb Z)$, $k=1, 2, ...$ .

Let $\ep>0$, $\mathcal F\subseteq B$ and $0< r<1$ be given. Choose a finite subset $\mathcal G\subseteq B$ and ${0<}\epsilon'<\epsilon$ such that $\mathcal F\subseteq \mathcal G$  and for any $\mathcal G$-$\epsilon'$-multiplicative map $L: B\to B$, the map $[L]_{\mathcal P}$ is well defined, and $[L]$ induce a homomorphism on $G$.

%{\bf  Adding another finite subset $\mathcal G$.} For any $\ep>0$, any finite subset $\mathcal F\subseteq B$ and any $0< r<1$,
{By choosing a sufficiently large $i_0,$ we may assume that
$[\psi_{i_0, \infty}](\underline{K}(T_{i_0}\oplus S_{i_0}))\supset G.$
In particular, we may assume that, by (iii) above,
$G\cap {\rm ker}\rho_{B_1}\subset (\psi_{i_0, \infty})_{*0}({\rm ker}\rho_{T_{i_0}}).$ Let $G'\subset \underline{K}(T_{i_0}\oplus S_{i_0})$ be such that $[\psi_{i_0, \infty}]({\cal G}')={\cal G}.$}
{Since $B=B_1\otimes U,$ we may write  that
$U=\varinjlim(M_{m(n)}, \imath_{n, n+1}),$
where $m(n)|m(n+1)$ and $\imath_{n, n+1}: M_{m(n)}\to M_{m(n+1)}$ is defined
by $a\mapsto a\otimes 1_{m(n+1)}.$ }
%By choosing a sufficiently large $i$,
One may assume that for each $f\in\mathcal G$, {there exists $i> i_0$ such that}
\begin{equation}\label{diag-f}
f=
\left(
\begin{array}{ccc}
f_0\oplus f_1 & & \\
 & \ddots & \\
 & & f_0\oplus f_1
\end{array}
\right)\in (T_i'\oplus S_i')\otimes M_m  %\bigcup_k ((T_i\oplus S_i)\otimes {\bf T_{0, k}} \otimes M_n)\cup (T_i\oplus S_i)\otimes {\bf T_{0, k}}\otimes\mathrm{C}(\mathbb T) \otimes M_n
\end{equation}
for some $f_0\in T_i'$, $f_i\in S_i'$, and $m>2MK/r,$ {where
$m=m(i+1)m(i+2)\cdots m(n),$
$T_i'=\psi_{i, \infty}'(T_i\otimes M_{m(i)})$ and $S_i'=\psi_{i, \infty}'(S_i\otimes M_{m(i)})$ and
where $\psi_{i, \infty}=\psi_{i, \infty}\otimes \imath_{i,\infty}.$}
Moreover, one may assume that $\tau(1_{T_i'})<r/2$ for all $\tau\in T(A_1)$. % and if $M[f]=0$ in $K_i(B, \mathbb Z/k\mathbb Z)$ for some $k\neq 0$, so is $M[f_1]$.

%Denote by $\pi_0: (T_i\oplus S_i)\otimes M_n\to T_i\otimes M_n$ and $\pi_1: (T_i\oplus S_i)\otimes M_n\to S_i\otimes M_n$  be the coordinator projection respectively. Define the
Choose a large $n$ so that $m={M_0}+l$ with ${M_0}$ divisible by $KM$ and $0\leq l<KM$. Then define the map $$L: (T_i'\oplus S_i')\otimes M_m \to (T_i'\oplus S_i')\otimes M_m$$ to be
$$
{L((f_{i, j} \oplus g_{i, j})_{m\times m})=(f_{i,j})_{m\times m}\oplus E_l(g_{i,j})_{m\times m}E_l,}
$$
{where $E_l={\rm diag}(\underbrace{1_{S_i'}, 1_{S_i'},...,1_{S_i'}}_l, \underbrace{0,0,...,0}_{M_0}),$}
%$$L((f_{i, j} \oplus g_{i, j})_{i, j})=\mathrm{diag}\{\underbrace{(f_{1, 1}\oplus g_{1, 1}), ..., (f_{l, l}\oplus g_{l, l}}_l), \underbrace{(f_{l+1, l+1}\oplus 0), ..., (f_{n, n}\oplus 0)}_{M_0}\}$$
{which is a \morp\, from $(T_I'\oplus S_i')\otimes M_m$ to $B.$
We then}
 extend $L$ to a completely positive linear map $B\to {(1_B-E_l)B(1_B-E_l)}$.
 Also define $$R: {(T_i'\oplus S_i')\otimes M_m \to T_i'\oplus S_i'}$$ to be
$$R((f_{i, j} \oplus g_{i, j})_{i, j})=g_{1, 1},$$ and extend it to a {\morp}\, $B\to B$,
where ${T_i'\oplus S_i'}$ is regarded as a corner of ${(T_i'\oplus S_i')}\otimes M_m\subseteq B$. Then $L$ and $R$ are $\mathcal G$-$\epsilon'$-multiplicative. Hence $[L]|_\mathcal{P}$ is well defined. Moreover, $$\tau(L(1_A))<\tau(1_{T_1})+\frac{l}{m}<\frac{r}{2}+\frac{MK}{2MK/r}=r\rforal \tau\in T(A).$$

Note that for any $f$ in the form of \eqref{diag-f}, one has
$$f=L(f)+{\overline{R}}(f),$$
{where $\overline{R}(f)$ may be written as}
$$
{\overline{R}(f)=\mathrm{diag}\{\underbrace{0,0,...,0}_l, \underbrace{(0\oplus g_{1,1}), ..., (0\oplus g_{1,1})}_{M_0}\}}.
$$
Hence for any $g\in G$,
$$g=[L](g)+ {M_0}[R](g).$$
%Note that for any $f$ in the form of \eqref{diag-f} , one has
%$$f-L(f)=\mathrm{diag}\{0, ..., 0, \underbrace{f_1, ..., f_1}_{k}\},$$
%with $k$ is divided by $MK$.
%In particular, for any $g\in G$, one has
%$$g-[L](g)=k(g_1)$$
%for some $g_1\in \underline{K}(B)$.
Then, if $g\in\ (G_{0, 1})_+\subseteq (G_0)_+$, one has,
$$g-[L](g)={M_0}[R](g)=K((\frac{M_0}{K})[R](g)).$$
And if $g\in G\cap K_i(B, \mathbb Z/ k\mathbb Z)$ ($i=0,1$), one also has
$$g-[L](g)={M_0}[R](g).
%=\frac{M_0}{M}(M[R](g))
%=\frac{M_0}{M}([R](Mg)).
$$
Since $Mg=0$  and {$M|M_0$}, one has $g-[L](g)=0$.

%By the choice of $M$, one also has that
%$$g-[L](g)=0\rforal g\in G\cap K_0(B, \mathbb Z/ k\mathbb Z),\ k\neq 0.$$ Hence $L$ induces an identity map on $G\cap K_0(B, \mathbb Z/ k\mathbb Z)$, $k\neq 0$.

%Note that {$K_0(S_i)$ is torsion free and} $K_1(S_i)=\{0\}.$  Therefore the restriction of $[R]$ to $G\cap K_1(B, \mathbb Z/n\mathbb Z)$ is zero, and hence $$g-[L](g)=0\rforal g\in G\cap K_1(B, \mathbb Z/ k\mathbb Z).$$
%Hence
%{Therefore} $L$ induces an identity map on $G\cap K_1(B, \mathbb Z/ k\mathbb Z)$.

%By the construction of $B$, there is an inclusion $G\cap \ker(\rho_B)\subset \ker(\rho_{T_i})$. The same argument as above shows that $L$ induces an identity map on $G\cap \ker(\rho_B)$.
{Since $L$ is identity on $\psi_{i, \infty}'(T_i\otimes M_{m(i)})$ and $i>i_0,$ by (iii),
$L$ is identity map on $G\cap \ker \rho_{B}.$  Since $K_1(S_i)=0$ for all $i,$
$L$ induces the identity map on $G\cap K_1(B).$ It follows that}
%Thus,
 $L$ is the desired map.
\end{proof}

Related to the above we have the following decomposition:

{
\begin{prop}\label{lem-compress}
Let $A_1$ be a separable  amenable C*-algebra in ${\cal B}_1$ (or
${\cal B}_0$) and let $A=A_1\otimes U$ for some  infinite dimensional UHF-algebra $U$.
Let $\mathcal G\subseteq A,$
$\mathcal P\subseteq  {\underline{K}}(A)$ be finite subsets, ${\cal
P}_0\subset A\otimes {\cal K}$ be a finite subset of  projections,
and let $\epsilon>0$, $0<r_0<1$ and $M\in\mathbb N$ be arbitrary.
Then there is a projection $p\in A,$ a \SCA\, $B\in {\cal C}$ (or in
${\cal C}_0$) with $p=1_B$ and $\mathcal
G$-$\epsilon$-multiplicative unital \morp s $L_1: A\to  (1-p)A(1-p)$
and $ L_2: A\to B$ such that
\begin{enumerate}
\item
$\|L_1(x)+L_2(x)-x\|<\ep\tforal x\in {\cal G};$
\item $[L_i]|_{\cal P}$  is well defined, $i=1,2;$
\item $[L_1]|_{\cal P}+[\imath\circ L_2]|_{\cal P}=[{\rm id}]|_{\cal P};$
\item $\tau\circ[L_1](g)\leq r_0\tau(g)$ for all $g\in
      {\mathcal P}_0$ and $\tau\in T(A)$;
\item  For any $x\in {\cal P},$ there exists $y\in \underline{K}(B)$ such that
   $x-[L_1](x)=[\imath\circ L_2](x)=M[\imath](y)$
   and,
\item for any $d\in {\cal P}_0,$ there exist positive element  $f\in {K_0}(B)_+$
      such that
      $$d -[L_1](d)=[\imath\circ L_2](d)=M[\imath](f),$$
      where $\imath: B\to A$ is the embedding.
      Moreover, we can require that $1-p\not=0.$
\end{enumerate}
\end{prop}
}
\begin{proof}
Since $A$ is in ${\cal B}_1$ (or ${\cal B}_0$), there is a sequence
of projections $p_n\in A$ and a sequence of \SCA\, $B_n\in {\cal
B}_1$ (${\cal B}_0$) with $1_{B_n}=p_n$ such that
\beq\label{compress-1}
\lim_{n\to\infty}\|(1-p_n)a(1-p_n)+p_nap_n-a\|=0,\\
\lim_{n\to\infty} \dist(p_nap_n, B_n)=0\andeqn\\
\lim_{n\to\infty}\max\{\tau(1-p_n): \tau\in T(A)\}=0 \eneq for all
$a\in A$.
%and some $b_n\in B_n,$ $n=1,2,....$
 Since  each $B_n$ is
amenable, one obtains easily a sequence of unital \morp s $\Psi_n:
A\to B_n$ such that \beq\label{compress-4}
\lim_{n\to\infty}\|p_nap_n-\Psi_n(a)\|=0\rforal a\in A. \eneq In
particular, \beq\label{compress-5}
\lim_{n\to\infty}\|\Psi_n(ab)-\Psi_n(a)\Psi_n(b)\|=0\rforal a, b\in
A. \eneq Let  $j: A\to A\otimes U$ be defined by $j(a)=a\otimes
1_U.$ There is a unital \hm\, $s: A\otimes U\to A$ and  a sequence
of unitaries $u_n\in A\otimes U$ such that \beq\label{compress-6}
\lim_{n\to\infty}\|a-{\rm Ad}\, u_n\circ s\circ j(a)\|=0\rforal a\in
A. \eneq There are  nonzero projection $e_n'\in U$  and $e_n\in U$
such that \beq\label{compress-7} \lim_{n\to\infty}t(e_n)=0\andeqn
1-e_n={\rm diag}(\overbrace{e_n', e_n',...,e_n'}^M),
\eneq
%{ What is $t$ here? Is $t$ here $\tau$ for all $\tau\in TA$?}
{where $t\in T(U)$ is the unique tracial state on $U.$}
Choose
$N\ge 1$ such that
\beq\label{compress-8}
0<t(e_n)<r_0/2\andeqn
\max\{\tau(1-p_n):\tau\in T(A)\}<r_0/2.
\eneq
Define $\Phi_n: A\to
(1-p_n)A(1-p_n)$ by $\Phi_n(a)=(1-p_n)a(1-p_n)$ for all $a\in A.$
Define  $\Phi_n'(a)=\Phi_n(a) \oplus  {\rm Ad}\, u_n\circ s(a\otimes
e_n)$ and $\Psi_n'(a)={\rm Ad}\, u_n\circ s(\Psi(a)\otimes (1-e_n))$
for all $n\ge N.$ Note that $u_n^*s(B_n\otimes (1-e_n))u_n\in {\cal
C}_1$ (or in ${\cal C}_0$). It is then easy to verify that, if we
choose a large $n,$ the maps $L_1=\Phi_n'$ and $L_2=\Psi_n'$ {meet}  the
requirements.
\end{proof}

\section{Positive maps from $K_0$-group of C*-algebras in the class $\mathcal C.$ }

This section contains some technical lemmas about positive \hm s
from $K_0(C)$ for some $C\in {\cal C}.$

\begin{lem}{\rm (Compare 2.8 of \cite{Lnbirr})}\label{multiple-ext} Let
$G\subset \Z^l$ (for some $l>1$) be a subgroup. {{ There is an integer $M>0$ satisfying the following condition:}}
% such that $\Z^l/G$ is {not finite}.
% and $s\ge 1$ be an integer.
Let $1>\sigma_1, \sigma_2>0$ { {be any given numbers.}} There is an integer ${{R}}>0$ {{such that}}:
{{if a set of  $l$ positive numbers $\af_i\in \R_+$ ($i=1,2,\cdots l$) satisfy
%%%satiesfy 
 $\af_i\ge \sigma_1$ for all $i$ and
satisfy }}
\beq\label{mule-0}
\sum_{i=1}^l \af_i m_i\in \Z \tforal (m_1, m_2,...,m_l)\in G,
\eneq
{ {then for any integer $K\geq R$,}} there exist a set of rational positive numbers
$ \bt_i\in \frac{1}{KM}\mathbb Z^+$ ($i=1,2,\cdots l$)
%$\bt_i\in \Z_+/M$
such that
\beq\label{mule-1}
\sum_{i=1}^l|\af_i-\bt_i{{|}}<\sigma_2,\,\,\,i=1,2,...,l\tand {\tilde \phi}|_G=\phi|_G,
\eneq
where $\phi((n_1, n_2,...,n_l))=\sum_{i=1}^l \af_in_i$ and ${\tilde \phi}((n_1,n_2,...,n_l))=\sum_{i=1}^l\bt_in_i$
for all
$(n_1,n_2,...,n_l)\in \Z^l.$
\end{lem}

\begin{proof}
{Denote by $e_j\in \Z^l$ the element having 1 in the $j$-th coordinate and $0$ elsewhere.}
{First we consider the case that $\Z^l/G$ is finite.}
In this case there is an integer $M\ge 1$ such that
$Me_j\in G$ for all $j=1,2,...,l.$  It follows that $\phi(Me_j)\in \phi(G){{\subset \Z}},$ $j=1,2,...,l$.
Hence $\af_j=\phi(e_j)\in \frac{1}{M}\mathbb Z^+$.
%By the assumption, there is $f_j\in {\red{\Z}}$ such that $Mf_j=\phi(Me_j).$
%Define ${\tilde\phi}: \Z^l\to {\red{\Z}}$ by ${\tilde \phi}(e_j)=f_j .$
We choose $\bt_j=\af_j,$ $j=1,2,...,l,$ and $\tilde \phi=\phi$. The lemma  follows{{---that is, for any $\sigma_1, \sigma_2$, we can choose $R=1$.}}

{Now we assume that $\Z^l/G$ is not finite.}

Regard $\Z^l$ as a subset of $\Q^l$ and set
 $H_0$ to be the  vector subspace of $\Q^l$  spanned by elements in $G.$ The assumption
 that $\Z^l/G$ is not finite implies that $H_0$ has dimension  $p<l.$
 Moreover $G\cong \Z^p.$
Let $g_1,g_2,...,g_p\in G$ be free generators of $G.$ View them as elements
in $H_0\subset \Q^l$ and
write
\beq\label{mule-2}
g_i={{\left(g_{i, 1}, g_{i, 2},\cdots, g_{i, l}\right)}},\,\,\,i=1,2,...,p.
\eneq

Define $L: \Q^p\to \Q^l$ to be $L=(f_{i,j})_{l\times p},$ where $f_{i,j}=g_{j,i},$ $i=1,2,...,l$ and
$j=1,2,...,p.$
Then $L^*=(g_{i,j})_{p\times l}.$  We also view $L^*:
\Q^l\to \Q^p.$
Define $T=L^*L: \Q^p\to \Q^p$ which is invertible.  Note that $T=T^*$ and $(T^{-1})^*=T^{-1}.$
Note  also that the matrix representation $(a_{i,j})_{{{l\times p}}}$ of $L\circ T^{-1}$ is {{an}} ${{l\times p}}$ matrix with entries 
%%%entires 
in $\Q.$ There is an integer $M_1\ge 1$ such that $a_{i,j}\in {1\over{M_1}}\Z,$ $i=1,2,...,p$ and
$j=1,2,...,l.$

Let $H_{00}={\rm ker} L^*.$ It has dimension $l-p>0.$ Let $P: \Q^l\to H_{00}$ be an orthogonal projection
which is a $\Q$-linear map. Represent $P$ as a $l\times l$ matrix.
Then its entries are in $\Q.$ There is an integer $M_2\ge 1$ such that
all entries are in
$\frac{1}{M_2}\Z.$
We will use the fact that $L^*=L^*(1-P).$
%$\Z/M_0.$

It is important to note that $M_1$ and $M_2$ depend on
$G$ only and are independent of $\{\af_j: 1\le j\le l\}.$  Let $M=M_1M_2.$

{ {Suppose that $\sigma_1,\sigma_2 \in (0,1)$ are two positive numbers and the numbers $\alpha_i\geq \sigma_1$ ($i=1,2,\cdots,l$) satisfy the condition (\ref{mule-0}). }}

{ {The condition (\ref{mule-0}) is equivalent to that }} $b_i{{:}}=\sum_{j=1}^l \af_jg_{i,j}\in \Z,$ $i=1,2,...,p.$ Put
$b=(b_1,b_2,...,b_p)^T$ and $\af=(\af_1, \af_2,...,\af_l)^T.$   Then
$b=L^*\af.$

%Let $L'=L|_{{\rm ker}L^{\perp}}: {\rm ker}L^{\perp}\to {\rm Im}L$ and
%$L'^{-1}: {\rm Im}L\to {\rm ker}L^{\perp}$ be the inverse.
%Let $v_1, v_2,...,v_p$ be a basis for ${\rm ker}L^{\perp}\subset \Q^l.$
%Put $v_i=(v_{i,1},v_{i,2},...,v_{i,l})\in \Q^l.$ There exists an integer $M_1\ge 1$ such that
%$v_{i,j}\in \frac{1}{M_1}\Z,$
%$v_{i,j}\in \Z/M_1,$
%$i=1,2,...,p$ and $j=1,2,...,l.$ With the basis $\{v_1,v_2,...,v_p\},$ we may write
%$L'^{-1}=(a_{i,j})_{l\times p}$ which has rational entries.
%Denote by $M_2\ge 1$ an integer such that
%$a_{i,j}\in \frac{1}{M_2}\Z.$
%$a_{i,j}\in \Z/M_2.$
%It is important to note that $M_0,$ $M_1$ and $M_2$ depend on
%$G$ only and are independent of $\{\af_j: 1\le j\le l\}.$

If we write
\beq\label{mule-3}
L (T^*)^{-1}b=c=\left(\begin{array}{c} c_1\\ c_2\\ \vdots\\ c_l\end{array}\right),
\eneq
then, since $b\in \Z^p,$ one has that
$c_i\in \frac{1}{M_1}\Z.$
%$c_i\in \Z/M_1M_2.$
Choose an integer ${ {R}}\ge 1$ such that $1/{{R}}<\sigma_1\sigma_2/{{(4l^2)}}$. { {Let $K\geq R$ be any integer.}}
%and put $M=KM_0M_1M_2.$
Note that
\beq\label{mule-3+}
L^*c=L^*L(T^*)^{-1}b=L^*LT^{-1}b=L^*\af.
\eneq
 Thus $\af-c\in{\rm ker} L^*$ as a subspace of $\R^l.$

{{ For the space $\R^l$, we use $\|\cdot\|_1$ and $\|\cdot\|_2$ to denote the $l_1$ and $l_2$ norm on it. Then we have
 $$\|v\|_2 \leq \|v\|_1\leq l \|v\|_2~~~\rforal v\in \R^l.$$}}

 Since $H_{00}$ is dense in the real subspace of
 ${\rm ker} L^*,$
there exists $\xi\in H_{00}$ such that
\beq\label{mule-4}
\|\af-c-\xi\|_2{{\leq \|\af-c-\xi\|_1}}<\sigma_1\sigma_2/4{{l}}.
\eneq
Pick $\eta\in \Q^l$ such that $\xi=P\eta.$
Then there is $\eta_0\in \Q^l$ such that $K\eta_0\in \Z^l$ and
\beq\label{mule-5}
{ {\|\eta_0-\eta\|_2\leq}}\|\eta_0-\eta\|_{{1}}<\sigma_1\sigma_2/2{{l}}.
\eneq
Since $P$ has norm 1 with respect to $l_2$ norm,
\beq\label{mule-6}
\|\af-c-P\eta_0\|_{{1}}{ {\leq l\|\af-c-P\eta_0\|_2 \leq l( \|\af-c-\xi\|_2+ \|P(\eta_0-\eta)\|_2) }}< \sigma_1\sigma_2.
\eneq
Put $\bt=c+P\eta_0=(\bt_1, \bt_2,...,\bt_l)^T.$ Note that $M_2K(P\eta_0)\in \Z^l$ {{, and that $M_1c\in \Z^l$.}}

{ {We have}}  ${ {K}}M\bt\in \Z^l,$ { {and}}
\beq\label{mule-7}
L^*\bt=L^*c=L^*\af\andeqn \|\af-\bt\|_{{1}}<\sigma_1\sigma_2.
\eneq
Moreover, since $\af_i\ge \sigma_1,$
\beq\label{mule-8}
\bt_i>0,\,\,\, i=1,2,...,l.
\eneq
Since $P\eta_0\in H_{00},$  one has that $L ^*\bt =L^*(1-P)\bt=L^*(1-P)c=L^*c=L^*\af=b.$
Define ${{\tilde \phi}}: \Q^l\to \Q$ by
\beq\label{mule-9}
{{\tilde \phi}}(x)=\langle x, \bt \rangle
%=\langle x, LT^{-1} b\rangle=\lange T^{-1}L^*x, b\rangle
\eneq
for all $x\in \Q^l.$
Note $L^*e_i=g_i,$ where
$e_i$ is the element in $\Z^p$ with the $i$-th coordinate being $1$ and $0$ elsewhere.  So
\beq\label{mule-10}
\phi(g_i)&=&\langle Le_i, \af\rangle =\langle e_i, L^*\af\rangle =\langle e_i, L^*\bt\rangle \nonumber \\
&=& \langle Le_i, \bt \rangle = \langle g_i, \bt \rangle ={{\tilde \phi}}(g_i),
%\langle e_i, b\rangle =\langle T^{-1}L^* Le_i, b\rangle\\
  % &=& \langle g_i, L(T^*)^{-1}b\rangle= \langle g_i, LT^{-1}b\rangle=\langle g_i, c\rangle=\phi'(g_i),
   \eneq
$i=1,2,...,p.$
It follows that
${{\tilde \phi}}(g)=\phi(g)$ for all $g\in G.$  { {Hence ${\tilde \phi}|_G=\phi|_{G}.$}}
Note that
${\tilde \phi}(\Z^l)\subset {1\over{{{K}}M}}\Z$, since $\bt_i\in {1\over{{{K}}M}}\Z^{{+}},$ $i=1,2,...,l.$
% Note that
%${\tilde \phi}(\Z^l)\subset {1\over{M_1}}\Z$, since $c_i\in {1\over{M_1}}\Z,$ $i=1,2,...,l.$
\end{proof}

If we do not need to approximate $\{\af_i: 1\le i\le l\},$ then ${{R}}$
can be chosen { {to be $1$, with $M=M_1$  which  only depends}} on $G$ and $l$ (by { {replacing}}  $\bt$ by $c$ in the proof).

From the proof \ref{multiple-ext}, since $L$ and $(T^*)^{-1}$ depends only on $g_1,g_2,...,g_p,$ we have the following:

\begin{lem}\label{ext-norm}
Let $G\subset \Z^l$ be an ordered subgroup with order unit $e$, and let
$g_1,g_2,...,g_p$ ($p\le l$) be a set of free generators of $G.$  For any $\ep>0,$ there exists $\dt>0$ satisfying the
following:
if $\phi: G\to \R$ is a \hm\,
%with $\phi({\red e})=1$ { {I think this restriction should be delete, see the remark after the lemma}}
such that
$$
|\phi(g_i)|<\dt,\,\,\,i=1,2,...,p,
$$
then, there is $\bt=(\bt_1, \bt_2,...,\bt_l)\in \R^l$ with
$|\bt_i|<\ep,$ $i=1,2,...,l,$ such that
$$
\phi(g)=\psi(g)\rforal g\in G,
$$
where $\psi: \Z^l\to \R$ defined by $\psi((m_1,m_2,...,m_l))=\sum_{i=1}^l\bt_im_i$
for all $m_i\in \Z.$
\end{lem}

%{ {The above lemma is tautologically true, and therefore not useful at all. Since one can choose $\dt$ small enough such that $|\phi(g_i)|<\dt$ implies $|\phi(e)|<1$---that is, the condition never be satisfied and therefore the sentence is tautologically true}}
%\begin{cor}\label{Mmultiple-ext}
%Let $G\subset \Z^l$ (for some $l>1$) be a subgroup and let $N\ge 1$ be an integer.
% such that $\Z^l/G$ is not finite.
 %There is an integer $M>0$ satisfying the following:
%For any $\af_i\in \R_+$
%such that
%\beq\label{Mmule-0}
%\sum_{i=1}^l \af_i m_i\in \Z \tforal (m_1, m_2,...,m_l)\in G,
%\eneq
%there exists
%$c_i\in\frac{1}{M} \Z_+$
%$c_i\in \Z_+/M$
%and
%${\tilde \phi}: \Z^l \to \frac{1}{M}\Z$
%${\tilde \phi}: \Z^l \to \Z/M$
%such that
%\beq\label{Mmule-1}
%{\tilde \phi}|_G=\phi_G,
%\eneq
%where $\phi((n_1, n_2,...,n_l))=\sum_{i=1}^l \af_in_i$ and ${\tilde \phi}((n_1,n_2,...,n_l))=\sum_{i=1}^lc_in_i$
%for all
%$(n_1,n_2,...,n_l)\in \Z^l.$
%\end{cor}

\begin{cor}\label{MextC}
Let $G\subset \Z^l$ be an order subgroup. Then, there exists  an integer $M\ge 1$
satisfying the following:
for any positive map $\kappa: G\to \Z^n$ (for any integer $n\ge 1$) with every element in $\kappa(G)$  divisible by $M$, there is  $R_0\ge 1$ such that, for any integer $K\ge R_0,$
there is a positive homomorphism ${\tilde \kappa}: \Z^l\to \Z^n$ such that
${\tilde \kappa}|_G=K\kappa.$
\end{cor}

\begin{proof}
We first prove the case that $n=1.$

Let $S\subset \{1,2,...,l\}$ be a subset and denote by $\Z^{(S)}$ the subset
$$
\Z^{(S)}=\{(m_1, m_2,...,m_l): m_i=0\,\,\,{\rm if}\,\,\, i\not\in S\}.
$$
Let $\Pi_S: \Z^l\to \Z^{{{(}}S{{)}}}$ be
the projection and $G(S)=\Pi_S(G).$

{ {Let $M(S)$ be the  {{integer}} (in place of $M$) as in
 \ref{multiple-ext} associated with $G(S)\subset \Z^{(S)}.$}}
%Let $M_0(S), M_1(S)$ and $M_2(S)$ be the integers (in place of $M_0,M_1, M_2$) in the proof of
% \ref{multiple-ext} associated with $G(S)\subset \Z^{(S)}.$
 %Note that $M_0(S), M_1(S), M_2(S)$ are independent of $\af$ in the proof.
 %Let $K_S\ge 1$ be the constant in the proof  \ref{multiple-ext} corresponding to
% the case  that $\sigma_1=\sigma_2=1/2l.$
 Put $M=\prod_{S\subset \{1,2,...,l\}}{{M(S)}}.$

 Now assume that $\kappa: G\to \Z$ is a positive \hm\, with multiplicity $M${ {---that is every element in $\kappa(G)$ is  divisible by $M$.}}

 By applying   2.8 of \cite{Lnbirr}, we obtain a
positive \hm\, $\bt: \Z^l\to \R$ such that $\bt|_G=\kappa.$
Define $f_i=\bt(e_i),$ where $e_i$ is the element in $\Z^l$ with $1$ at the $i$-th coordinate and
$0$ elsewhere, $i=1,2,...,l.$  Then $f_i\ge 0.$
Choose $S$ such that $f_i>0$ if $i\in S$ and $f_i=0$ if $i\not\in S.$

{{Evidently if $\xi_1,\xi_2\in \Z^l$ satisfy $\Pi_S(\xi_1)=\Pi_S(\xi_2)$, then $\bt(\xi_1)=\bt(\xi_2)$,  and if we  further
assume $\xi_1,\xi_2\in G$ then $\kappa(\xi_1)=\kappa(\xi_2)$. Hence the maps $\bt$ and $\kappa$ induce the maps $\bt': \Z^{(S)} \to R$ and $\kappa': G(S) \to \Z$ such that $\bt=\bt'\circ \Pi_S$ and $\kappa=\kappa'\circ \Pi_S$. In addition, we have $\bt'|_{G(S)}=\kappa'$.}}

{{Let $\sigma_1=2\sigma_2=\frac{\min \{f_i:~ i\in S\}}{2M}$. Apply \ref{multiple-ext} to $\af_i=f_i/M> \sigma_1$ for $i\in S$ and to $G(S)\subset \Z^{(S)},$ we obtain the number $R(\kappa)$ (depending on $\sigma_1$ and $\sigma_2$ and therefore depending on $\kappa$) as in the lemma. For any $K\geq R(\kappa)$, it follows from the lemma there are $\bt_i\in \frac{1}{KM})_+$ (for $i\in S$), such that $\tilde \kappa'|_{G(S)}=\frac{1}{M}\kappa'$, where $\tilde \kappa': Z^{(S)}\to \Q$ is defined by $\tilde \kappa'(\{m_i\}_{i\in S})=\sum_{i\in S}\bt_i m_i$. Evidently $\tilde \kappa=KM(\tilde \kappa'\circ\Pi_S) \Z^l\to \Z$ is as desired for this case.}}

%If this case, we
%Consider $G(S)\subset \Z^{(S)}.$
 %According to the given $\kappa,$
%there is $R_S\ge 1$ such that
%\beq\label{MextC-2}
%R_Sf_i\ge 1/l\rforal t\in S.
%\eneq

%For any integer $N_S\ge R_S,$
%put
%$\af_i=R_Sf_i$
%$\af_i=N_Sf_i$
%for $i\in S.$  Then the proof of \ref{multiple-ext} implies that there exists
%a positive \hm\, ${\tilde \kappa}: \Z^l \to \Z$ such that
%\beq\label{MextC-3}
%{\tilde \kappa}|_G=N_S\kappa.
%\eneq

This prove the case $n=1.$

In general,
 let $s_i: \Z^n\to \Z$ be the projection
to the $i$-th summand, $i=1,2,...,n.$  { {Apply the case $n=1$ to each of the the maps $\kappa_i:=s_i\circ \kappa$ {{(for $i=1,2,\cdots n$)}} to obtain $R(\kappa_i)$. And let $R_0=\max_{i} R(\kappa_i)$.}}  For any $K\ge R_{{0}},$
by what has been proved, we obtain ${\tilde\kappa}_i: \Z^l\to \Z$ such that
\beq\label{MextC-4}
{\tilde \kappa}_i|_G=Ks_i\circ \kappa|_G,\,\,\,i=1,2,....
\eneq
Define ${\tilde \kappa}: \Z^l \to \Z^n$ by
${\tilde \kappa}(z)=({\tilde \kappa}_1(z),{\tilde \kappa}_2(z),...,{\tilde \kappa}_l(z)).$
The lemma follows.
\end{proof}

%{\color{green}  The following is not used. Let's keep it for a while, and let me look the details tomorrow.}

%\begin{cor}[Lemma 3.6 of \cite{Niu-TAS-I}, see also 2.12 of \cite{Lnbirr}]\label{pt-integer}
% Let $C\in\mathcal C$ with $C=A(F_1, F_2, \psi_0, \psi_1).$  Then, there exists  an integer $M\ge 1$
%satisfying the following:
%for any positive map $\kappa: {K_0}(C)\to \Z^n$ (for any integer $n\ge 1$) with every element in $\kappa(K_0(C))$  divisible by $M$,  there exists $R_0\ge 1,$ for any $K\ge R_0,$ there is a homomorphism $\phi: F_1\to F_3$ such that
%$K\kappa=(\phi\circ \pi_e)_{*0},$  where $F_3$ is a direct sum of $n$ copies of full matrix algebras.
%\end{cor}

%\begin{proof}
%Note, by \ref{2Lg13}, $K_0(C)\subset K_0(F_1)=\Z^l$ and the
%embedding is given by $(\pi_e)_{*0}.$ Let  $M_0$ ${\tilde \kappa}$
%be as in \ref{MextC} and ${\tilde
%\kappa}([1_A])=(m_1,m_2,...,m_n)\in \Z^n$ be the positive element.
%Define $F_3=M_{m_1}\oplus M_{m_2}\oplus M_{m_n}.$ It is well known
%that there is a unital \hm\, $\phi: F_1\to F_3$ such that
%$\phi_{*0}={\tilde \kappa}.$   The lemma follows.
%\end{proof}

\begin{lem}[Lemma 3.2 of \cite{EN-Tapprox}]\label{annihilation2}
Let $G={K_0}(S)$, where $S\in \mathcal C$ and let $r: G\to \Z$ be a
strictly positive homomorphism. Then,  for any  order unit $u\in
G_+, $  there exists a natural number $m$ such that if the map
$\theta: G\to G$ is defined by $g\mapsto r(g)u$, then there
exists an integer $m\ge 1$ such that
 the positive homomorphism $\text{id}+m\theta: G\to G$ factors through $\bigoplus_{i=1}^n\Z$ positively for some $n$.
\end{lem}

\begin{proof}
Let $u$ be an order unit of $G$, and define the map $\phi: G\to G$ by $\phi(g)=g+ r(g)u$; that is, $\phi=\id+\theta$.
Define $G_n=G$ and $\phi_n: G\to G$ by $\phi_n(g)=\phi(g)$ for all $g$ and $n.$ Consider the inductive limit
$$
\xymatrix{
G\ar[r]^\phi & G \ar[r]^\phi & \cdots \ar[r] & \varinjlim G.
}
$$
Then the ordered group $\varinjlim G$ has the Riesz decomposition property.
In fact, let $a, b, c\in \varinjlim G_+$ such that
$$
a\le b+c.
$$
Without loss of generality, one may assume that $a\not=b+c.$

We may assume that there are $a', b', c'\in G$ for the $n$-th $G$ such that
$\phi_{n, \infty}(a')=a,$ $\phi_{n, \infty}(b')=b$ and $\phi_{n, \infty}(c')=c,$ {{and furthermore }}
%Therefore, for some large $k\ge 1,$
\beq\label{ann2-n1}
{{a' < b'+c'}}.
\eneq

A straightforward calculation shows that for each $k$,  there is $m(k)\in\mathbb N$ such that
$$
\phi_{n, n+k}(a')=a'+m(k)r(a')u,\quad \phi_{n,n+k}(b')=b'+m(k)r(b')u,\quad \textrm{and}\quad \phi_{n,n+k}(c')=c+m(k)r(c')u.
$$
Moreover, the sequence $(m(k))$ is strictly increasing.
Since $r$ is strictly positive, combing with (\ref{ann2-n1}), we have that
$$
r(a{{'}})< r(b{{'}})+r(c{{'}})\,\,\,{\rm (in} \,\,\, \Z{\rm )}.
$$
There are $l(a{{'}})_i\in \Z_+$ such that
$$
l(a{{'}})_1+l(a{{'}})_2=r(a{{'}}),\,\,l(a{{'}})_1\le r(b{{'}})\andeqn l(a{{'}})_2\le r(c{{'}}).
$$
Without loss of generality, we may assume $d=r(b{{'}})-l(a{{'}})_1>0$ (otherwise
we let $d=r(c{{'}})-l(a{{'}})_2$).
Since $u$ is an order unit, there is $m_1\in \Z_+$ such that
$$
m_1d u>a{{'}}.
$$
Choose $k\ge 1$ such that $m(k)>m_1.$
Let $a_1=a'+m(k)l(a')_1u$ and $a_2=m(k)l(a')_2u.$
Then
$$
a_2{ {=}}m(k)l(a')_2u\le m(k)r(c')u\le c'+m(k)r(c')u=\phi_{n,n+k}(c').
$$
Moreover,
$$
a_1=a'+m(k)l(a')_1u\le m(k)du+m(k)l(a')_1u\le b'+m(k)r(b')u=\phi_{n, n+k}(b').
$$
Note
$$
\phi_{n, n+k}(a')=a_1+a_2\le \phi_{n, n+k}(b')+\phi_{n,n+k}(c').
$$
These  imply that
\beq
a{{=}}\phi_{n+k,\infty}(a_1)+\phi_{n+k{{,\infty}}}(a_2)\le b+c,\\
\phi_{n+k,\infty}(a_1)\le b\andeqn \phi_{n+k,\infty}(a_2)\le c.
\eneq

This implies that the limit group $\varinjlim G$ has the Riesz decomposition property.
Since $G$ is unperforated, so is $\varinjlim G$. It then follows from the Effros-Handelman-Shen Theorem (\cite{EHS}) that  the ordered group $\varinjlim G$ is a dimension group. Therefore, for a sufficiently large $k\in\mathbb N$, the map $\phi^k$ must factor through the ordered group $\bigoplus_n\mathbb Z$ positively for some $n$. Since $\phi^k$ has the desired form $\id+m(k)\theta$, the lemma follows.
\end{proof}

\begin{lem}\label{Cuthm}
Let $(G,G_+,u)$ be  an group
%%%orroup
 with order unit $u$ %which is a  finitely generated  group
such that the positive cone $G_+$ is generated by finitely many positive elements which are smaller than $u.$ Let  $\lambda: G\to K_0(A)$ be an order preserving map such that $\lambda(u)=[1_A]$  and
$\lambda(G_+\setminus \{0\})\subset K_0(A)_+\setminus \{0\},$ where $A\in {\cal B}_1$ ($A\in {\cal B}_0$).
%$A=A_1\otimes U,$ where $U$ is an infinite dimensional UHF-algebra.
Let $a\in K_0(A)_+\setminus \{0\}$ with $a\le [1_A].$ % and for any integer $K>1.$
Let ${\cal P}\subset G_+\setminus \{0\}$ be a finite subset.
Suppose that there exists an integer
$N\ge 1$ such that $N\lambda(x)>[1_A]$ for all $x\in {\cal P}.$

There are two positive homomorphisms $\lambda_0, \lambda_1: G\to
K_0(A)$ and a \SCA\, $S'\subset A$  with $S'\in {\cal C}$ ($S'\in {\cal
C}_0$) satisfying the following: \beq\label{Cuthm-1}
\lambda=\lambda_0+\lambda_1,\,\,\, \lambda_1=\imath_{*0}\circ
\gamma,\,\,\, {{0\leq}}\lambda_0(u)<a\andeqn \gamma(g)>0 \eneq for all
$g\in G_+\setminus \{0\},$ where $\gamma: G\to K_0(S')$ with
$\gamma([{{u}}])=[1_{S'}]$ and where $\imath: S'\to A$ is the
embedding. Moreover, $N\gamma(x)\ge \gamma(u)$ for all $x\in
{\cal P}$. Furthermore, if $A=A_1\otimes U$, where $U$ is an
infinite dimensional UHF-algebra and $A_1\in\mathcal B_1$ (or ${\cal
B}_0$),  then, for any integer $K\ge 1,$ we can require that
$S'=S\otimes M_K$ for some $S\in {\cal C}$ (or ${\cal C}_0$) and
$\gamma$ has multiplicity $K.$
\end{lem}

\begin{proof}
Let $\{g_1, g_2,..., g_m\}\subset G_+$ be the set of generators
of $G_+$ with $g_i < u$. Since $A$ has stable rank one, it is easy to check that
there are projections $q_1, q_2,...,q_m\in A$ such that
$\lambda(g_i)=[q_i],$ $i=1,2,...,m.$ For convenience,  to simplify
notation, \wilog, we may assume that ${\cal P}=\{g_1,g_2,...,g_m\}.$
 Define
$$
Q_i={\rm diag}(\overbrace{q_i, q_i,...,q_i}^{{{N}}}),\,\,\, i=1,2,...,m.
$$
By the assumption, there are $v_i\in M_N(A)$ such that
\beq\label{Cutm-14-n1} v_i^*v_i=1_A\andeqn v_iv_{{i}}^*\le
Q_i,\,\,\,i=1,2,...,m. \eneq
 Since $A\in {\cal B}_1$ (or ${\cal
B}_0$), there exists a sequence of projections $\{p_n\}$ of $A,$ a
sequence of \SCA\, $S_n\in {\cal C}$ (${\cal C}_0$) with
$p_n=1_{S_n}$ and a sequence of unital \morp s $L_n: A\to S_n$ such
that \beq\label{Cuthm-n1}
&&\hspace{-0.2in}\lim_{n\to\infty}\|a-((1-p_n)a(1-p_n)+p_nap_n)\|=0,\,\,\,\lim_{n\to\infty}\|L_n(a)-p_nap_n\|=0\\\label{Cuthm-n1+}
&&\andeqn \lim_{n\to\infty}\sup_{\tau\in T(A)}\{\tau(1-p_n)\}=0.
\eneq

It is also standard to find, for each $i,$  a projection
$e_{i,n}'\in (1-p_n)A(1-p_n),$ a projection $e_{i,n}\in M_{{N}}(S_n)$
and partial isometries $w_{i,n}\in {{M_N(}} S_n{{)}}$ such that
\beq\label{Cuthm-nn1}
&&\lim_{n\to\infty}\|(1-p_n)q_i(1-p_n)-e_{i,n}'\|=0,\\\label{Cuthm-nn2}
&&w_{i,n}^*w_{i,n}=p_n,\,\,\, w_{i,n}w_{i,n}^*\le
e_{i,n},\\\label{Cuthm-nn3}
&&\lim_{n\to\infty}\|{{(}}L_n{{\otimes {\rm id}_{M_N})}}(v_i)-w_{i,n}\|=0\andeqn \lim_{n\to\infty}\|
(L_n\otimes {\rm id}_{M_{{N}}})(Q_i)-e_{i,n}\|=0. \eneq  Let
$\Psi_n: {{A}}\to (1-p_n)A(1-p_n)$ be defined by $\Psi_n(a)=(1-p_n)a(1-p_n)$ for
all $a\in A.$ We will use $[\Psi_n]\circ \lambda$  for $\lambda_0$
and $[L_n\circ \lambda]$ for $\gamma$ for  some large $n.$  The fact
that  $\lambda_0$ and $\gamma$ are \hm s follows from Lemma 7.1 of
\cite{Lin-BDF}. To see that $\lambda_0$ is  positive, we use
(\ref{Cuthm-nn1}) and the fact that $G_+$ is finitely generated. It
follows from (\ref{Cuthm-nn2}) and (\ref{Cuthm-nn3})  that
$N\gamma(x)\ge \gamma(u)$ for all $x\in {\cal P}.$ Since we assume
that the positive cone of $G_+$ is generated by ${\cal P},$ this
also shows that $\gamma(x)>0$ for all $x\in G_+\setminus \{0\}.$
%$\gamma(g)>0$ for all $g\in G_+\setminus \{0\},$  we note that
%an element in  $K_0(S_n)$ is strictly positive if and only if {its image in $\rho_{S_n}(K_0(S_n))$
%is positive.}
%Then, by applying Lemma 9.2 of \cite{LinTAI}, we may assume that, for each $t\in T(S_n),$
%there exists a $\tau'\in T(A)$ such that
%\beq\label{Cuthm-n2}
%\|t(L_n(q_i))-\tau'(q_i)\|<1/n,\,\,\,i=1,2,...,m.
%\eneq
By (\ref{Cuthm-n1+}), we can choose large $n$ so that
 ${{0\leq}}\lambda_0(u)<a.$

It should be noted when $A$ does not have (SP), one can choose $\lambda=\lambda_1${{, and $\lambda_0=0$.}}

If  $A=A_1\otimes U,$ then, \wilog, we may assume that $p_n\in A_1$ for all $n.$
Choose  a sequence of nonzero projections $e_n\in U$ such that $t(1-e_n)=r(n)/K,$
where $t$ is the unique tracial state on $U$ and $r(n)$ are positive rational numbers such that
$\lim_{n\to\infty}t(e_n)=0.$
Thus $S_n\otimes (1-e_n)\subset B_n$ where $B_n\cong S_n\otimes M_K$ and
$p_n\otimes (1-e_n)=1_{B_n}.$ We check that the lemma follows if we replace
$\Psi_n$ by $\Psi_n',$ where $\Psi_n'(a)=(1-p_n)a\otimes 1_U+p_na\otimes e_n$
\end{proof}

\begin{lem} {\rm (see Lemma 3.6 of \cite{EN-Tapprox} or Lemma 2.8 of \cite{Niu-TAS-II})}
\label{decomposition2}
Let $G={K_0}(S)$, where $S\in \mathcal C$. Let $H={K_0}(A)$ for $A=A_1\otimes U,$ where $A_1\in\mathcal B_1$
(or ${\cal B}_0$) and $U$ is  {{an}} infinite dimensional UHF-algebra.
% such that $A$ is tracially approximately divisible. %$A=A_1\otimes U$ for some $A_1\in \mathcal B_1$ and a UHF-algebra $U$.
Let $M_1\ge 1$ be a given integer and $d\in K_0(A)_+\setminus \{0\}.$ Then for any strictly positive homomorphism $\theta: G\to H$ with multiplicity $M_1$, and any integers $M_2\ge 1$ and $K\ge 1$ such that
$K\theta(x)>[1_{A}]$ for all $x\in G_+\setminus\{0\},$ one has a decomposition $\theta=\theta_1+\theta_2$, where $\theta_1$ and $\theta_2$ are positive homomorphisms from $G$ to $H$ such that the following diagrams commute:
\begin{displaymath}
\xymatrix{
G \ar[rr]^{\theta_1} \ar[dr]_{\phi_1} &                      & H & & G \ar[rr]^{\theta_2} \ar[dr]_{\phi_2} &                     & H \\
                                      & G_1 \ar[ur]_{\psi_1} &   & &                                      & G_2 \ar[ur]_{\psi_2} & },
\end{displaymath}
where $\theta_1([1_{{S}}])\le d,$
$G_1\cong\bigoplus_n\Z$ for some natural number $n$ and $G_2=K_0(S')$ for some
%C*subalgebra
\SCA\, $S'$ of $A$ which is in the class $\mathcal C$ (or in ${\cal C}_0$),  $\phi_1, \psi_1$ are
positive \hm s
  and $\psi_2=\imath_{*0},$ where $\imath: S'\to A$  is the embedding.
Moreover, $\phi_1$ has the multiplicity $M_1,$ $\phi_2$ has the multiplicity $M_1M_2,$   $2K\phi_2(x)>\phi_2([1_{{S}}])>0$ for all $x\in G^{+}\setminus\{0\}.$
\end{lem}

\begin{proof}
Let $u=[1_S]$, and let $m$  be as in Lemma \ref{annihilation2}. Suppose that $S={{A}}(F_1, F_2, \psi_0, \psi_1)$ with $F_1=M_{R_1}\oplus M_{R_1}\oplus \cdots \oplus M_{R_l}.$
It is easy to find  a strictly positive  \hm\, $\eta_0: K_0( F_1)\to \Z.$  Define
$r: G\to \Z$ by $r(g)=\eta_0\circ (\pi_e)_{*0}.$
By replacing $S$ with $M_r(S)$  and $A$ by $M_r(A)$ for some integer $r\ge 1, $
\wilog\, we may assume that $S$ has a finite subset of projections
 ${\cal P}=\{p_1, p_2, ..., p_{{l}}\}$  such that
 every projection $q\in S$ is equivalent to one of projections in ${\cal P}$ and
 $\{[p_i]:1\le i\le {{l}}\}$ generates $K_0(S)_+$ (see \ref{FG-Ratn}).
Let
 $$
 \sigma_0=\min\{\rho_A(d)(\tau): \tau\in T(A)\}.
 $$
 Note that since $A$ is simple, one has that $\sigma_0>0.$

Let
$$
\sigma_1=\inf\{\tau(\theta([p]): p\in {\cal P},\,\,\, \tau\in T(A)\}>0.
$$
Since $A=A_1\otimes U,$ $A$  has the (SP) property, there is a projection $f_0\in A_+\setminus \{0\}$ such that
\beq\label{decomp6/7-1}
0<\tau(f_0)< \min\{\sigma_0, \sigma_1\}/8Nr(u)\rforal \tau\in T(A).
\eneq

Since $A=A_1\otimes U,$ we may choose
$f_0$ so that $f_0=M_1\tilde{h}$ for some nonzero ${\tilde h}\in K_0(A)_+.$  Put
$\theta_0': G\to K_0(A)$ by $\theta_0'(g)=r(g){\tilde h}$ for all $g\in G.$ { {And let $\theta'=M_1\theta_0'$. Then
$2\theta'(x) <\theta(x)$ for all $x\in G_+\setminus \{0\}$.}}

Since  $\theta$ has multiplicity $M_{1}$, one has that $\theta(g)-\theta'(g)$ is divisible by $M_1$ for any $g\in G$. By the choice of $\sigma_0,$  one checks that $\theta-\theta'$ is strictly positive.
Moreover,
\beq\label{decom-nn1}
2K\rho_A((\theta(x)-\theta'(x))(\tau) &{{>}}& 2K\rho_A(\theta(x))(\tau)-{{K}}\rho_A(\theta(x))(\tau){{=}} K\rho_A(\theta(x))(\tau)\\
&\ge & \rho_A([1_A])(\tau)\rforal \tau\in T(A).
\eneq

Applying \ref{Cuthm}, one obtains a \SCA\, $S'\in A,$  a \hm\, ${\tilde \theta}_1: G\to K_0(A)$ and
strictly positive \hm\, $\phi_2: G\to K_0(S')$ such that
\begin{equation}\label{decompo6/8-1}
\theta-\theta' = {\tilde \theta}_1+\imath_{*0}\circ \phi_2,
\end{equation}
\begin{equation}\label{decomp6/8-2}
|\tau({\tilde \theta}_1(u))| < \frac{\tau({\tilde h})}{mM_1{{ M_2}}},\quad \tau\in T(A),\\
2K\phi_2(x)>\phi_2([1_S]),\,\,\,\phi_2([1_S])=[1_{S'}]{{,}}
\end{equation}
 where $m$ is from \ref{annihilation2} and $\phi_2$ has multiplicity  $M_1M_2,$
  {{and}} where $\imath: S'\to A$ is the embedding.
Put $$\theta_2=\imath_{*0}\circ \phi_2,\quad \textrm{and}\quad \psi_2=\imath_{*0}.$$

Since $\theta(g)-\theta'(g)$ is divisible by $M_1$ and any element in $\theta_2(G)$ is divisible by $M_1$, one has that any elements in $\tilde{\theta}_1(G)$ is divisible by $M_1$. Therefore, the map $\tilde{\theta}_1$ can be decomposed further as $M_1\theta_1'$, and one has that $\theta-\theta'=M_1\theta_1'+\theta_2$. Therefore, there is a decomposition
$$\theta=\theta'+M_1\theta'_1+\theta_2=M_1\theta'_0+M_1\theta'_1+\theta_2.$$
Put $$\theta_1=M_1(\theta'_0+\theta_1').$$
Then,
$$
\rho_A(\theta_1([1_S]))(\tau)<{{\tau(}}d{{)}}/2\rforal \tau\in T(A).
$$

We then show that $\theta_1$ has the desired factorization property.
For $\theta'_0+\theta'_1$, one has the following farther decomposition: for any $g\in G$,
\begin{eqnarray*}
\theta'_0(g)+\theta'_1(g)&=&r(g){\tilde h}+\theta'_1(g)\\
&=&r(g)({\tilde h}-m\theta'_1(u))+r(g)m\theta'_1(u)+\theta_{{1}}'(g)\\
&=& r(g)({\tilde h}-m\theta'_1(u))+ \theta'_1(mr(g)u)+\theta'_1(g)\\
&=& r(g)(\tilde{h}-m\theta'_1(u))+\theta'_1(mr(g)u+g).
\end{eqnarray*}
By (\ref{decomp6/8-2}), ${\tilde{h}}-m\theta_1'(u))>0.$
By Lemma \ref{annihilation2}, $g\to mr(g)u+g$ factors though $\bigoplus_n\Z$ positively for some $n$. Therefore, the map $M_1(\theta'_0+\theta_1')$ factors though $\bigoplus_{(1+n)M_1}\Z$ positively.
So there are positive homomorphisms $\phi_1: G\to \bigoplus_{(1+n)M_1}\Z$ and $\psi_1: \bigoplus_{(1+n)M_1}\Z.\to K_0(A)$ such that $\theta_1=\psi_1\circ \phi_1$ and $\phi_1$ has multiplicity of $M_1$.
\end{proof}

%Define $\psi: \Z^l\to \Z^n$ by
%$$
%{\tilde \kappa}(x)=(\psi_1(x),\psi_2(x),...,\psi_n(x))\rforal x\in G.
%$$
%It is clear that ${\tilde \kappa}$ is positive and ${\tilde \kappa}|_G=\kappa.$

%\end{proof}

%{\bf We can remove the following too--L}

% In the following lemma $e_j\in \Z^l$ is the element with $1$ at the $j$-th coordinate and zero elsewhere.

%\begin{lem}[cf. Lemma 2.11 of \cite{Lnbirr}]\label{birrL}
%Let $G\subset \Z^l$ (for some $l>1$) be a subgroup which contains an order unit of $\Z^l$ and
%let $T$ be a Choquet simplex.  Let $F$ be dense, divisible and ordered subgroup of $\Aff(T)$ with the strict
%ordering (i.e., $f\ge 0$ if and only if $f=0$ or $f(t)>0$ for all $t\in T$). Suppose that
%$\phi: G\to F$ is a positive \hm\, such that $\phi(g)>0$ for all $g\ge 0$ and $g\not=0.$
%Then there exists a positive \hm\, ${\tilde \phi}: \Z^l\to F$ such that
%${\tilde \phi}|_G=\phi$ and ${\tilde \phi}(e_j)>0$ for $j=1,2,...,l.$
%\end{lem}

%The above  is the same as Lemma 2.11 of \cite{Lnbirr}. Note in the proof, $w_n=z_l^{(1)}c(y_n-l_ny_1).$

%\begin{cor}\label{July22-newlem}
%xxx
%\end{cor}

\section{Existence Theorems for affine maps on tracial state spaces for building blocks}

\begin{lem}\label{discretizeT0}
Let $A$ be a unital separable  \CA\,  with $T(A)\not=\emptyset$ and
%Let $X$ be a compact metric space and
let ${\cal H}\subset A$  be a finite  subset.
Then, for any $\sigma>0,$ there exists an integer $N>0$  and a finite subset
$E\subset \partial_e(T(A))$ satisfying the following:
%For any probability  measure $\mu$ on $X$
For any $\tau\in T(A)$ and any $k\ge N$, there are $\{t_1, t_2,...,t_k\}\subset E$ such that
%$\{\pi_1, \pi_2, ..., \pi_k\}\subseteq {\rm Irr}(A)$ such that
\beq\label{161N-1}
|\tau(h)-\frac{1}{k}\sum_{i=1}^k t_i(h)|<\sigma \tforal h\in\mathcal H.
\eneq
%where $t_i\in E.$
 { {(If $\tau$ is a (possibly unnormalized) trace on $A$ with $\|\tau\|\leq 1$, then there are $\{t_1, t_2,...,t_{k'}\}$ with $k'\leq k$ such that
$$
|\tau(h)-\frac{1}{k}\sum_{i=1}^{k'} t_i(h)|<\sigma \tforal h\in\mathcal H.)
$$}}
  Suppose that $A$ is a subhomogeneous C*-algebras. Then,
there are $\{\pi_1, \pi_2, ..., \pi_k\}\subseteq {\rm Irr}(A)$ such that
\beq\label{161N-2}
|\tau(h)-\frac{1}{k}\left({\rm tr}_1\circ \pi_1(h)+{\rm tr}_2\circ \pi_2(h)+\cdots+{\rm tr}_k\circ \pi_k(h)\right)|<\sigma \tforal h\in\mathcal H,
\eneq
where  { {$\pi_j\in E$ and}} ${\rm tr}_j$ is the canonical trace of $\pi_j(A)$.
Moreover, if, for each ${{l}},$  $\hat{A}_{{l}}$ has no isolated points, then
$\{\pi_1,\pi_2,...,\pi_k\}$ can be required  to consist of distinct points,
where $\hat{A}_{{l}}$ is the set of all irreducible representations of $A$ with rank exactly ${{l}}.$
\end{lem}

%\begin{rem}

{\bf Remark:} Note that in (\ref{161N-2}), $\pi_i$ may not be distinct. But the subset $E$ of   irreducible
representations can be chosen independent of $\tau$ (but dependent of $\sigma$ and
${\cal H}$).
%\end{rem}

\begin{proof}
%The proof is a standard compact argument.
Without loss of generality, one may assume that
$\|f\|\le 1$ for all $f\in {\cal H}.$
Note the tracial state space $T(A)$ is weak *- compact.
Therefore, there are $\tau_1, \tau_2,...,\tau_m\in T(A)$ such that, for any
$\tau\in T(A),$ there is $j\in\{1,2,...,m\}$ such that
\begin{equation}\label{12N1-1}
|\tau(f)-\tau_j(f)|<\sigma/4\tforal f\in {\cal H}.
\end{equation}
Note that  the set of extreme points of $T(A)$ is the set of those tracial
states induced by irreducible representations of $A.$
By the Krein-Milman Theorem, there are $t_1', t_2',...,t_n'\in \partial_e(T(A))$
%there  is a finite subset
%$\{\pi_1', \pi_2',...,\pi_n'\}\subset \hat{A}$
and nonnegative numbers $\{\af_{i,j}\}$  such that
\begin{equation}\label{12N1-2}
|\tau_j(f)-\sum_{i=1}^n \af_{i,j} t_i'(f)|<\sigma/8\andeqn\sum_{i=1}^n\af_{i,j}=1.
\end{equation}
Put $E=\{t_1',t_2',...,t_n'\}.$
Choose $N>32mn/\sigma.$ { {Let $\tau$ be possibly unnormalized trace on $A$ with $0<\tau(1)\leq 1$. Suppose that $j$ is so chosen that $\|\tau(f)/\tau(1)- \tau_j(f)\|<\sigma/4 \tforal f\in {\cal H}$ as in \ref{12N1-1}. }}
Then,  for any $k\ge N,$ there exist positive rational numbers $r_{i,j}$  and positive integers $p_{i,j}$ ($1\le i\le n$ and $1\le j\le m$)  such that,
%for each $j,$
\begin{equation}\label{ratn-001}
{ {\sum_{i=1}^nr_{i,j}\leq 1, ~~{\mbox {or}}}}~~ \sum_{i=1}^nr_{i,j}=1~~ { {{\mbox {if}} ~~\tau(1)=1}}
\end{equation}
\begin{equation}\label{ratn-002}
r_{i,j}=\frac{p_{i,j}}{k}\quad\textrm{and}\quad |{{\tau(1)}}\alpha_{i,j}-r_{i,j}|<\frac{\sigma}{8n},\quad 1\le i\le n,\,\,
%1\le j\le m.
\end{equation}
Note that
\begin{equation}\label{12N1-3}
{ {\sum_{i=1}^np_{i,j}\leq k~~{\mbox {or}}}}~~ \sum_{i=1}^np_{i,j}=k~~{ {{\mbox {if}}~~ \tau(1)=1}}.
\end{equation}

%{\blue {I move the next sentence to two paragraphs later}}
%There exists $\dt>0$ such that
%\begin{equation*}
%|f(x)-f(y)|<\sigma/64k\,\,\,{\rm provided}\,\,\, {\rm dist}(x, y)<\dt \rforal f\in {\cal H}.
%\end{equation*}

%for all $f\in {\cal H}.$

%Now let $\tau\in T(A).$
%Choose $j$ so that (\ref{12N1-1}) holds. One may also assume that (\ref{12N1-2}) holds.
Then, by (\ref{ratn-002}),
\begin{equation*}
|\tau(f)-\sum_{i=1}^n({p_{i,j}\over{k}})t_i'(f)|<\sigma/4+\sigma/8=3\sigma/8\rforal f\in {\cal H}.
\end{equation*}
It is then clear that (\ref{161N-1}) holds by  repeating each $t_i'$ $p_{i,j}$ times.

Now suppose that $A$ is subhomogeneous. Then every
$t_i'$ has the form ${\rm tr}_i\circ \pi_i,$ where $\{\pi_1, \pi_2,..., \pi_n\}\subset {\rm Irr}(A).$
It follows that (\ref{161N-2}) holds.

{ {There exists $\dt>0$ such that for any irreducible representation $x,y\in \hat{A}_l $ with ${\rm dist}(x, y)<\dt$, we have
\begin{equation*}
|f(x)-f(y)|<\sigma/64k
%{\rm provided}\,\,\, {{{\rm dist}(x, y)<\dt}}
{{\rforal f\in {\cal H}}}.
\end{equation*}}}

If $\hat{A}_{{l}}$ has no isolated points, for each $i,$ choose $\pi_{i,j}$ distinct points
 in a neighborhood $O(\pi_i', \dt)$ of $\pi_i'$ (in $\hat{A}_{{l}}$) with diameter
less than $\dt.$
Let $\{\pi_{1,j},\pi_{2,j},..., \pi_{k,j}\}$ be the resulting set of $k$ elements (see (\ref{12N1-3})).
Then, one has
\begin{equation}\label{12n1-5}
|\tau(f)-{1\over{k}}(f(\pi_{1,j})+f(\pi_{2,j})+\cdots+ f(\pi_{k,j}))|<\sigma\rforal f\in {\cal H},
\end{equation}
as desired.
\end{proof}

%\begin{lem}\label{discretizeT}
%Let $\mathcal H$ be a finite subset of $C([0, 1]\times\mathbb T)$, let $\sigma>0$ and $1/2>\dt>0.$ Then there are integer  $N_1\ge 1$
%depending  on $\sigma,$ integer $N_2\ge 1$ depending on $\sigma$ and  $\dt,$
%and a positive integer $N_3\le N_2$ such that for any finite  measure $\mu$ on $[0, 1]\times \mathbb T$
%with $1\ge \|\mu\|\ge \dt$
%and any $k\ge N_1,$  there are $\{x_1, x_2, ..., x_{kN_3}\}\subseteq (0, 1)\times \mathbb T$ such that
%$$|\int h d\mu-\frac{1}{kN_2}(h(x_1)+h(x_2)+\cdots+h(x_{kN_3}))|<\sigma, \quad\forall h\in\mathcal H.$$
%\end{lem}
%Let us consider a unital stably finite C*-algebra $A$ and a matrix algebra $M$. Recall that an order-unit map $\kappa: K_0(A)\to K_0(M)$ is compatible to a tracial state $\tau\in\mathrm{T}(A)$ if $$\mathrm{tr}(\kappa(p))=\tau(p)\rforal p\in K_0(A).$$ In the following, we will show that any tracial state which is almost compatible to a given $K_0$ map can be perturbed to a exactly compatible trace if $A\in\mathcal C_0$ or $A=C\otimes\mathrm{C}(\mathbb T)$ for some $C\in\mathcal C_0$.
%First, one has {the following:}

The following is well known.
%{\color{green} I think we could state it for $C(X)$, and it is well known.}

\begin{lem}\label{oldext}
Let $C=\bigoplus_{i=1}^kC(X_i)\otimes M_{r(i)},$  where each $X_i$ is connected compact metric space. %Let $\Delta: C_+^{q, 1}\setminus\{0\}\to(0, 1)$ be an order preserving map.
Let $\mathcal H\subseteq C$ be a finite subset and let $\sigma>0$. Then there is
%are finite subset $\mathcal H_1\subseteq C^+$ such that for any  $\tau\in T(A)$ satisfying
%$$
%\tau(h) > \Delta(\hat{h}) \rforal h\in \mathcal H_{1},
%$$
an integer $N\ge 1$ satisfying the following:
for any positive homomorphism $\kappa: K_0(C)\to K_0(M_s)=\Z$   with $\kappa([1_{M_{r(i)}}])\ge N$
and any $\tau\in T(C)$ such that
$$
{\rho_C(x)(\tau)}={\rm tr}(\kappa(x))\tforal x\in K_0(C),
$$
where ${\rm tr}$ is the tracial state on $M_s,$
there is a homomorphism $\phi: A\to M_s$ such that
$\phi_{*_0}=\kappa$ and
$$|{\rm tr}\circ\phi(h)-\tau(h)|<\sigma \tforal h\in\mathcal H.$$
\end{lem}

\begin{lem}\label{appoldext}
Let $C=C(\T)\otimes F_1,$ where
$F_1=M_{R(1)}\oplus M_{R(2)}\oplus \cdots \oplus M_{R(l)},$ or $C=F_1.$
%Let $\sigma>0$ be such that $\sigma<1/\sum_{j=1}^l R(j).$
Let $\mathcal H\subseteq C$ be a finite subset, and let $\epsilon>0$.
There is $\dt>0$ satisfying the following:
For any $M_s$, any order-unit map $\kappa: K_0(C)\to K_0(M_s)$
 and any tracial state $\tau\in T(C)$ such that
$$
 |{\rho_{ {M_s}}(\kappa(p))(\mathrm{tr})}-\tau(p)|<\dt
 $$
 for all projections {{$p$}} in $C,$
where $\mathrm{tr}$ is the tracial state  on $M_s$, there is a tracial state $\tilde{\tau}\in T(C)$ such that
$$\mathrm{tr}(\kappa([p]))=\tilde{\tau}(p),$$ and
$$|\tau(h)-\tilde{\tau}(h)|<\epsilon\rforal h\in \mathcal H.$$

\end{lem}

\begin{proof}
We may assume that ${\cal H}$ is in the unit ball of $C.$
Let $\dt=\ep/l.$
We may write that $\tau=\sum_{j=1}^l \af_j \tau_j,$ where $\tau_j$ is a tracial state
on $C(\T)\otimes M_{r(j)},$ $\af_j\in \R_+$ and $\sum_{j=1}^k\af_j=1.$
Let $\bt_j={\rm tr}(\kappa([1_{M_{R(j)}}]),$ $j=1,2...,l.$
Put ${\tilde \tau}=\sum_{j=1}^l \bt_j\tau_j.$
Then
%\beq\label{appoldext-1}
${\rm tr}(\kappa(p))={\tilde \tau}(p)$
%\eneq
for all projections $p\in C;$ and for any $h\in {\cal H},$
$$
|{\tilde \tau}(h)-\tau(h)|\le \sum_{j=1}^l|\bt_j-\af_j|<\ep,
$$
%for all $h\in {\cal H}.$
as desired.
\end{proof}

In the following statement, we keep the notation in \ref{8-N-3} for \CA\, $A=A_m\in {\cal D}_m.$
In particular,
$\Pi_{ {m-1}}=\pi_e^{(m)}: A_m\to A_{m-1},$
$\Pi_{ k}=\pi_e^{(k+1)}\circ \pi_e^{(k+2)}\circ \cdots \circ \pi_e^{(m)}: A_m\to A_k$  and
$\Pi_0=\pi_e^{(1)}\circ \pi_e^{(2)}\circ \cdots \circ \pi_e^{(m)}: A_m\to A_0$ be the quotient map.

\begin{lem}\label{Reduction2015}
Let $A\in {\cal D}_m$ be a unital \CA\, for some $m\ge 0,$
%$\kappa: K_0(A)\to \R$ be a positive \hm\,
%with  $\kappa([1_A])=1$
and
let $\Delta: A^{q, {\bf 1}}_+\setminus \{0\}\to (0,1)$ be an order preserving map.
Let ${\cal H}\subset  A$ be a finite subset, let ${\cal P}\subset M_{m{{'}}}(A)$ be a finite
subset of projections (for some integer $m{{'}}\ge 1$) and let $\ep>0.$
 Then there are $\Delta_0: (A_0)_+^{q,{\bf 1}}\setminus \{0\}\to (0,1),$ a
 finite subset $\mathcal H_1\subseteq A^{\bf{1}}_+\setminus \{0\}$ and a positive integer $K$ such that for any  $\tau\in T(A)$ satisfying
\begin{equation}\label{extnn-0}
\tau(h) > \Delta(\hat{h}) \rforal h\in { {\mathcal H_{1}}}
\end{equation}
 and any positive homomorphism $\kappa: K_0(A)\to {{\Z}}$
with $s=\kappa([1_A]){{\in \N}}$
such that
\begin{equation}\label{extnn-0+}
\rho_A(x)(\tau)=(1/s)(\kappa(x))\tforal x\in K_0(A),
\end{equation}
there are
%integers $K_0\ge 1,$
rational numbers $r_1,r_2,..,r_k\ge 1,$ finite dimensional irreducible representations
$\pi_1, \pi_2,...,\pi_k$ of $A,$
$\pi_{k+1}, \pi_{k+2},...,\pi_{k+l}$ of $A_0$
{{and a  trace $\tau_0$ on $A_0$}}
such that, {{for all $h\in {\cal H}$}},
\beq\label{Reduc-1}
&&\hspace{-0.8in}|\tau(h)-(\sum_{i=1}^k r_i {\rm tr}_i\circ \pi_i(h)+\tau_0\circ \Pi_0(h))|<\ep/2,\,\, |\tau_0\circ \Pi_0(h))-{{\sum_{i=k+1}^{k+l}r_i{\rm tr}_i\pi_i\circ \Pi_0(h)}}|<\ep/2,\\
%\end{equation}
%\tau_0\circ \Pi_0(h))|<\ep,
%\begin{equation}
\label{Reduc-1+}
&&\hspace{-0.4in}\tau(p)=\sum_{i=1}^kr_i{\rm tr}_i \circ \pi_i(p)+\tau_0\circ \Pi_0(p)=\sum_{i=1}^kr_i{\rm tr}_i \circ \pi_i(p)+\sum_{i=k+1}^{k+l}r_i{\rm tr}_i\pi_i\circ \Pi_0(p)
\eneq
%tau_0\circ \Pi_0(p)
%\tforal p\in {\cal P},
%\\\label{Reduc-1++}
%{equation}
for all $p\in{\cal P}$;
  \begin{equation}\label{Reduc-1++}
%\tau_0(\Pi_0(h))
{{\tau_0\circ \Pi_0(h)}}
\ge \Delta_0(\widehat{\Pi_0(h)}) \tforal h\in {\cal H}_1;
\end{equation}
%{\blue{I change  $H_1$ to $H_1'$ above, otherwise, we do not need to introduce $H_1'$---Gong.}}\\
%{\blue {for all $h\in{\cal H}_1$;
and moreover, $Kr_i\in \Z$ for $i=1,2,\cdots, k$, and $sKr_i\in \Z$, for $i=k+1,k+2,\cdots k+l.$

% $K_0r_i\in \Z$ and $K_0k|K.${\blue{It seems in the proof, we did not take care of $K_0k|K$---in particular how to make the number $k$, the number of irreducible representation involved, to be a factor of $K$.  We only prove some thing like $Kr_i\in \Z$. If this is what we need later, then we can eliminate the mention of $K_0$. If we do need $k$ to be factor of $K$ and $(K/k) r_i\in\Z$, we need to make this way, though not very difficult, I will let you guys to do it---Gong.}}
\end{lem}

\begin{proof}
\Wlog, we may assume that $A\in {\cal D}_m'.$
Write $A_m$ and keep the notation in \ref{8-N-3} for $A_m.$  So $X_m=X.$
Put $Y_m=X_m\setminus { Z_m},$ and $Y_k=X_k\setminus {Z_k},$
$k=0, 1, ..., m.$
Let $\pi_e^{(k)}: A_k\to A_{k-1},$ $\lambda_k: A_k\to P_kC(X_k, F_k)P_k$  and
$\Gamma_{ k+1}: A_k\to P_{k+1}C({ Z_{k+1}}, F_{ k+1})P_{k+1}$ be as defined in \ref{8-N-3}.
Let $\Pi_k: A_m\to A_k$ be defined by $\pi_e^{(k+1)}\circ \pi_e^{(k+2)}\cdots \pi_e^{(m)}$
for $k=0,1,2,...,m-1.$  Let $I_k=Q_kC_0(Y_k, F_k)Q_k\subset A_k,$ $k=1,2,...,m${{, where  $Q_k=P_k|_{Y_k}$, as in  \ref{8-N-3} .}} Note that $I_k={\rm ker}\pi_e^{(k)}.$
Note also that we assume  that $P_k(x)\not=0$ for all $x\in X_k,$ $k=1,2,...,m.$

By replacing $A$ by $M_N(A),$ \wilog, we may assume that $K_0(A)$ is generated
by $\{p_1, p_2,...,p_c\},$ where $p_i\in A$ are projections, $i=1,2,...,c.$   Moreover, we may assume
that, \wilog, ${\cal P}= \{p_1,p_2,..., p_c\}$ and ${\cal P}\subset {\cal H}.$  To simplify notation,
we may further assume that $p_1=1_A.$
Write $A_0=\bigoplus_{j=1}^{l(0)}M_{r(0,j)}$ and assume
that $e_j\in M_{r(0,j)}$ is a rank one projection.
%=\bigoplus_{j=1}^{l(0)} P_{0,j}C(X_{0,j}, M_r)P_{0,j}$ as before.
%By choosing even large $N,$ \wilog, we may further assume that $P_{0,j}C(X_{0,j}, M_r)P_{0,j}$
%contains a rank one projection $e_j.$

{{ Recall for any $d>0$, $X_j^d=\{x\in X_j: \dist(x, Z_j)<d\}.$}} Choose $\dt_0>0$ such that there is  continuous map $s_*^{(j,\dt_0)}: \overline{X_j^{\dt_0}}\to { Z_j}$
satisfying
\beq\label{Reduc-2}
s_*^{(j, \dt_0)}(x)=x\rforal x\in { Z_j},\,\,\, {\rm dist}(s_*^{(j,\dt_0)}(x), x)\le \dt_0\rforal x\in X_j^{\dt_0}
\\\label{Reduce-3}
\|\Gamma_{j}(\Pi_{j-1}(h))(s_*^{(k, \dt_0)}(y))-\lambda_j(\Pi_j(h))(y)\|<\ep/2^{{ m}+1} \rforal h\in {\cal H}
\eneq
and $y\in X_j^{\dt_0}=\{y\in X_j: {\rm dist}(y, { Z_j})< \dt_0\},$
$j=1,2,...,{ m}.$
Write
$$
P_kC(X_k, F_k)P_k=\bigoplus_{j=1}^{l(k)} P_{k,j}(C(X_k, M_{s(k,j)})P_{k,j},
$$
where $P_{k,j}\in C(X_k, M_{s(k,j)})$ is a projection of rank $r(k,j)$ at each $x\in X_k,$
%$F_k=M_{r(k,1)}\oplus M_{r(k,2)}\oplus \cdots M_{r(k, l(k))}$
($k=0,1,2,...,m$) and $r(i)=r(0,i),$ $i=1,2,...,l.$  Note that ${{ (K_0(A_k))/\ker(\rho_{A_k})\subset \rho_{A_{k-1}}(K_0(A_{k-1}))/\ker(\rho_{A_{k-1}})}}\subset
\Z^l$ (see { \ref{8-n-K0nccw}}).

%Let $e_j\in M_{r(j)}\subset F_0$ be a rank one projection, $j=1,2,...,l.$
Let $S_{k,i}: P_kC(X_k, F_k)P_k\to P_{k,i}C(X_k, M_{s(k,j)})P_{k,i}$ be the {{quotient map}}, $i=1,2,...,l(k),$ $k=0,1,2,...,m.$
Choose ${{\chi_k}}\in (I_k)_+$ such that  $\|{{\chi_k}}\|\le 1,$ ${{\chi_k(x)}}$ is a scalar multiple of
$\lambda_k(1_{A_k})(x)$ for all $x\in Y_k,$  ${{\chi_k}}(x)=\lambda_k(1_{A_k})(x)$  if ${\rm dist}(x, { Z_k})>\dt_0$ and
${{\chi_k}}(x)=0$ if ${\rm dist}(x, { Z_k})<\dt_0/2$ and for all $x\in X_k,$ $k=1,2,...,m.$

{ Note that} $Y_k\cap X_{ k}^{\dt_0}\not=\emptyset,$
define an element $b_k\in C_0(Y_k)$ such that $b_k(x)=0$ if $x\in { Z_k},$ $b_k(x)>0$ for all $x\in Y_k\cap X_k^{\dt_0}.$

Let $b_{k,i}=S_{k,i}(b_k\cdot Q_k) { {\in {{P_{k,i}}}C_0(Y_k, M_{{{s}}(k,i)}){{P_{k,i}}}, } }$
$i=1,2,...,l(k).$
Let  $h_j\in A_+$ with $\|h_j\|\le 1$ and $\Pi_0(h_j)=e_j$ and
$h_j^{(k)}\in A_+$ with $\|h_j^{(k)}\|\le 1$ such that
$$\Pi_{ k}(h_j^{(k)})=(1-{{\chi_k}})\Pi_{ k}(h_j)=(1-{{\chi_k}})^{1/2}\Pi_{ k}(h_j)(1-{{\chi_k}})^{1/2},$$
 $j=1,2,...,l(0)$ and
$k=0,1,2,...,m-1.$  Let ${\cal H}_0=\{h_j: 1\le j\le l(0)\}.$
For each $b\in A_{ k-1}{ =A_k/I_k}$ define
$$L_k(b)=(1-{{\chi_k}})s^{k,\dt_0}(b)=(1-{{\chi_k}})^{1/2}s^{k,\dt_0}(b)(1-{{\chi_k}})^{1/2}{ \in A_k},$$
for $k=1,2,..., m,$ where
$s^{k, \dt_0}: A_k/I_k\to C(\overline{X_k^{\dt_0}}, F_k)$ defined by
$s^{k, \dt_0}(a)(x)=\Gamma_k(a)(s_*^{k, \dt_0}(x))$ for all $a\in A_k/I_k=A_{k-1}$ and $x\in \overline{X_k^{\dt_0}}.$ Note that $L_k: A_{k-1}\to A_k$ is a \morp.

For each $h\in A_{ {k-1}},$ define $h'_{(k)}=L_m\circ L_{m-1}\circ \cdots\circ L_k({{h}}){ {\in A_m}},$
%=(1-a_k)s^{k,\dt_0}\circ \Pi_k(h)(1-a_k),$
$k=1,2,..., m.$
%Note that $h_{(k)}'\in A_m=A$ for all $k.$
%where
%$s^{k, \dt_0}: A_k/I_k\to C(\overline{X_k^{\dt_0}}, F_k)$ defined by
%$s^{k, \dt_0}(a)(x)=\Gamma_k(\pi_e^{(k)}(a)(s_*^{k, \dt_0}(x))$ for all $a\in A_k$ and $x\in \overline{X_k^{\dt_0}}.$
%where $X_k^{\dt_0}=\{X_k: {\rm dist}(x, \partial(X_k))< \dt_0\}.$
%${{ Note that $h'_{(k)}$ only depends on $ \Pi_{k-1}(h)$ (instead of $h$), we can d}}
Define { $\DT_k: (A_{k-1})^{(q,1)}_+\setminus \{0\} \to (0,1)$ by }
$$
\Delta_k({{\hat h}})=\Delta(\widehat{h_{(k)}'})\rforal h\in { {A_{k-1}^+}},
$$
$k={{1}}, 2, ...,m-1.$ { {Evidently for any $h\in A_{k-2}^+$ ($k\geq 2$) we have}}
\begin{equation}\label{extnn-G}
{ {\Delta_{k-1}({\hat h})=\Delta_{k}({\widehat{L_{k-1}(h)}}).}}
\end{equation}
Put
$$
{\cal H}_1=\{1_A\}\cup {\cal H}\cup {\cal H}_0\cup \cup_{k=1}^m \{h_{(k)}': h\in {\cal H}\cup {\cal H}_0\}\cup
\cup_{k=1}^m\{b_{k,i}':1\le i\le l(k)\},
$$
{ where $b_{k,i}'= L_m\circ L_{m-1}\circ\cdots \circ L_{k+1} (b_{k,i}).$ }
%Let ${\cal H}_1'\supset {\cal H}_1$ be a finite subset.
{ {Put}}
\begin{equation}\label{extnn-2}
\sigma_1=\min\{\Delta(\hat{h}): h\in {\cal H}_1\}/2\quad\mathrm{and} \quad \sigma_2=
\frac{\sigma_1\cdot  \ep\cdot \dt_0}{64m\cdot (\sum_{k=1}^m l(k)+l)}.
\end{equation}
Let $K'\ge 1$ be an integer (in place of {{the product $MR$ there}}) in \ref{multiple-ext} associated
with $G=(\Pi_0)_{*0}(K_0(A))={{(K_0(A))/\ker(\rho_A)}}\subset \Z^{l}$ {{(see (\ref{K0nccw-2}))}}, $\sigma_1$ and $\sigma_2.$  Let $\{d_1, d_2,..., d_m\}$ be the set of ranks of all possible irreducible representations of $A.$
%Put $K_3=\sum_{i=1}^m R_i.$
{{Let $K''=(\prod_{j=1}^md_j).$}}

Let ${{\tilde{\cal H}}}_k=\lambda_k\circ \Pi_{k}({\cal H}),$ $k=1,2,...,m-1.$ Note that
$Y_k\setminus X_k^{\dt_0}=\{{{x}}\in X_k: {\rm dist}(x, {Z_k})\ge \dt_0\}$ is compact.
Let $K_k\ge 1$ be integers (in place of ${{N}}$) in \ref{discretizeT0} for $\sigma_2/2^m$ (in place of ${{\sigma}}$)
and ${{Q_k(C(Y_k\setminus X_k^{\dt_0}, F_k)Q_k}}$ (in place of {{$A$}}){{, and $ {\tilde{\cal H}}_k|_{Y_k\setminus X_k^{\dt_0}}$ (in place of ${\cal H}\subset A)$.}}

Choose an integer $K_0\ge 1$ such that
$1/K_0<\frac{\sigma_2^2\cdot \sigma_1^2}{64 \cdot l(\sum_{k=1}^m l(k)+l)},$  let
${ {{\tilde K}=(K_0^m)(\prod_{k=1}^{m} K_k)}}$ {{and let$K=\tilde KK'K''.$}}

For each $\tau\in {{T(A)}}$, we may write
\begin{equation}\label{extnn-trace}
\tau(f)=\sum_{{ {i}}=1}^{l{ (m)}}\int_{Y_m} a_{m,{ {i}}}{\rm tr}_{{ m},{ {i}}}(f(t))d\mu_j)+ t_{m-1}\circ \pi_e^{(m)}(f)\rforal f\in A,
%\sum_{k=1}^d(\sum_{j=1}^{l(k)}\int_{Y_k} a_{k,j}{\rm tr}_{k,j}(f(t))d\mu_j)+ t_0\circ \Pi_0(f)\rforal f\in A,
\end{equation}
where $\mu_{ {i}}$ is a Borel probability measure on $Y_{{m}},$
${\rm tr}_{m,{ {i}}}$ is the normalized trace on $M_{r(m,{ {i}})}$ {{(here we regard that  $f(t)\in M_{r(m,i)}$ by identifying $ P_{m,i}(t)M_{s(m,i)}P_{m,i}(t)\cong M_{r(m,i)}$), }}$t_{m-1}$ is an unnormalized trace  on $A_{m-1},$ $a_{m,{ {i}}}\in \R_+$ and
$\sum_{{ {i}}}a_{m,{ {i}}}+\|t_{m-1}\|=1.$

We claim that for each $k{{\in\{1,2,\cdots,m+1\}}},$
\begin{eqnarray*}
&& |\tau(h)-(\sum_{s=k}^m({{\sum_{i=1}^{l(s)}}}\sum_{{j}} d(s,i,j) {\rm tr}_{s,{ {i}}}(\LD_s(h)(y_{s,i, j})){)}+t_{k-1}\circ \Pi_{ k-1}(h))|  \le
{(m-k+1)\ep\over{2^{m}}}  \rforal h\in {\cal H},\\
&& \tau(p)  =  \sum_{s=k}^m({{\sum_{i=1}^{l(s)}}}\sum_{{j}} d(s,i,j) {\rm tr}_{s,{ {i}}}(\LD_s(h)(y_{s,i, j})){{)}}+t_{k-1}\circ \Pi_{ k-1}(p) \rforal  p\in {\cal P}, \andeqn\\
&&t_{ k-1}(h)\ge \Delta_k(\hat{h})-{{((m-k+1)\sigma_2)/2^m}} \rforal h\in \Pi_{ k-1}({{\cal H}}_1) { \subset (A_{k-1})_+},
\end{eqnarray*}
%{ {(here again I change  $H_1$ to $H_1'$ above---Gong)}}
where $K_0^{m-k{{+1}}}K_mK_{m-1}\cdots K_k \cdot d(s,i,j)\in \Z_+,$
$\sum_{s,i,j}d(s,i,j)+\|t_k\|=1$ and $y_{s,{ {i,}}j}\in Y_s {{\setminus X_{s-1}^{\dt_0}}},$
${\rm tr}_{s,{ {i}}}$ is the tracial state on $M_{r(s,{ {i}})}$ and
$t_{k-1}$ is a { unnormalized } trace in $(A_{k-1}),$ $s=k, k+1,...,m.$ Here the convention is, when $k=m+1,$
$\tau_m=\tau,$ $\Pi_{m}={\rm id}_A$ and $\sum_{m+1}^m$  is a sum of empty set (so
any such sum is zero). Hence the above holds for $k=m+1.$ { {We will prove the claim by reverse induction.}}
%So the above holds for $k=m+1.$

Suppose that the above holds for $k\le m+1.$
One may write
\beq\label{extnn-15-2}
\hspace{-0.1in}t_{k-1}(a)=\sum_{{ {i}}=1}^{l(k-1)} a_{k-1,{ {i}}}\int_{Y_{k-1}}{\rm tr}_{k-1,{ {i}}} (\lambda_{k-1}(a))d\mu_{k-1,{ {i}}}+t'_{k-2}\circ \Pi_{ k-2}(a)\rforal a\in A_{k-1},
\eneq
%{ {what is $t_{k-1}'$ here, is it $t_{k-1}$ from the induction assumption, or normalized one?}}\\
where $\sum_{{ {i}}=1}^{l(k-1)}a_{k-1,{ {i}}}+\|t'_{k-2}\|=\|t_{k-1}\|,$
%{ {if $t_{k-1}'$ is not normalized one, then you can not have the above}}
${\rm tr}_{k-1,{ {i}}}$ is the tracial state on $M_{r(k-1, { {i}})},$
%\subset F_{k-1},$
 $t_{k-2}'$ is a trace on
  $A_{ k-2}$ and $\mu_{k-1,{ {i}}}$ is a probability Borel measure on $X_{ k-1}.$
 Here we also identify ${\rm tr}_{k-1,{ {i}}}$ with ${\rm tr}_{k-1,{ {i}}}\circ {\tilde \pi}_y\circ S_{k-1,{ {i}}}$ and where ${\tilde \pi}_y: P_{k,{ {i}}}C(X_{k-1}, M_{s(k-1,{ {i}})})P_{k,{ {i}}}\to M_{r(k-1,{ {i}})}$ is the point evaluation at $y\in Y_{k-1}.$
 % where
 %$P_j: F_{k-1}\to M_{r(k-1,j)}.$

%Let $\bt_i=
%\mu_{ k-1,i}(Y_{ k-1}\cap X^{\dt_0}_{ k-1}),$ $i=1,2,...,k(l).$
%Since ${\rm supp} ({\widehat b_{k-1,i}})\subset Y_{ k-1}\cap X^{\dt_0}_{ k-1}$, we have
%\beq\label{extnn-3++}
%\bt_i>t_{k-1}(\widehat{b_{k-1,i}})>\Delta_k(\widehat{b_{k-1,i}})>2\sigma_2.
%\eneq

It follows from \ref{discretizeT0}  that there are
$y_{k-1,i, j}\in Y_{k-1}\setminus X_{k-1}^{\dt_0},$ $j=1,2,...,m(k-1,{ {i}})\le K_{k-1}$ { (for convenient we will denote $y_{k-1,i, j}$ by $t_{i,j}$)} such that
\beq\label{extnn-3}
%\hspace{-0.1in}{\eta(f)\triangleq}
|a_{k-1,i}\hspace{-0.05in}\int_{Y_{k-1}\setminus X_{k-1}^{\dt_0}}{\rm tr}_{k-1,i}(f)d\mu_{k-1,i}(t)-\hspace{-0.1in}
\sum_{j=1}^{m({ k-1,} i)} d(k-1,i,j){\rm tr}_{k-1,i}(f(t_{i,j}))|<{a_{k-1,i}{{\sigma_2^2}}\over{2^{m}}}
\eneq
for all $f\in {{{\cal H}_1}},$ where $K_0K_{k-1}\cdot d(k-1,i,j)\in \Z,$
%and $\sum_jr(k-1,i,j)=a_{k-1,i},$
$i=1,2,...,l(k-1).$
Put
$$
\eta(f)=|a_{k-1,i}\hspace{-0.05in}\int_{Y_{k-1}\setminus X_{k-1}^{\dt_0}}{\rm tr}_{k-1,i}(f)d\mu_{k-1,i}(t)-\hspace{-0.1in}
\sum_{j=1}^{m({ k-1,} i)} d(k-1,i,j){\rm tr}_{k-1,i}(f(t_{i,j}))|.
$$
Then
\beq\label{extnn-3n}
\eta(f)<{a_{k-1,i}\sigma_2^2\over{2^{m}}}\rforal f\in {\cal H}_1.
\eneq
%{ {To make the above equation true, $\af_{k-1,i}$ must be rational numbers with not too large denominators, therefore you must modify $t_{k-1}'$  slightly, and when you do the modification, you need to keep the value of new $t_{k-1}'$ on a projection being same
%as the value of old one on same projection. I guess you can apply Lemma 15.1 to do so.}}

For each $i,$ define $\rho_i
% \rho_i':
: A\to \C$ by
\beq\label{extnn-3+}
&&\hspace{-0.2in}\rho_i(f)=a_{k-1,i}\int_{Y_{k-1}} {\rm tr}_{k-1,i}(\lambda_{k-1}(f))d\mu_{k-1,i}-\sum_{j=1}^{m({ k{{-1}},}i)} d(k-1,i,j){\rm tr}_{k-1,i}(f(t_{i,j}))
%\rho'_i(f)&=& \int_{Y\setminus X^{\dt_0}} {\rm tr}_i(f)d\mu_i-(1/N_1N_2)\sum_{j=1}^{m(i)} {\rm tr}_j(f(t_{i,j}))
\eneq
for all $f\in A.$
Put
$\af_i=\rho_i(1_A)$
and $\bt_i=
\mu_{ k-1,i}(Y_{ k-1}\cap X^{\dt_0}_{ k-1}),$ $i=1,2,...,k(l).$
Since ${\rm supp} ({\widehat b_{k-1,i}})\subset Y_{ k-1}\cap X^{\dt_0}_{ k-1}$, we have
\beq\label{extnn-3++}
\bt_i>t_{k-1}(\widehat{b_{k-1,i}})>\Delta_k(\widehat{b_{k-1,i}})>2\sigma_2.
\eneq
It follows from  (\ref{extnn-3}) and $1_A\in {\cal H}$, that
\beq\label{extnn-3++a}
|a_{k-1,i}\bt_i-\af_i|< a_{k-1,i}\sigma_2^2/2^{m},\, {\rm or}\,\,\, |a_{k-1,i}-{\af_i\over{\bt_i}}|<a_{k-1,i}\sigma_2/2^m
\eneq
and consequently, $\af_i>a_{k-1,i}\bt_i- a_{k-1,i}\sigma_2/2^{m}>0$.
%If $\bt_i\not=0,$  then $\Delta_{ k}(\widehat{b_{k-1,i}})>0.$ Since $\sigma_2<\Delta_k(\widehat{b_{k-1,i}})/2,$ it follows that, $\af_i>0$ when $\bt_i\not=0.$
%If $\bt_i=0,$ define $T_i=0.$ If $\bt_i\not=0,$

 Define
$T_i: A\to \C$ by $$T_i(f)={\af_i\over{\bt_i}}\int_{Y_{k-1}\cap X_{k-1}^{\dt_0}} { {\rm tr}_{k-1,i}}( s^{k-1,\dt_0}\circ { \Pi_{k-{ 2}}}(f)) d\mu_{k-1,i}$$ for all $f\in A.$
{ It follows from (\ref{Reduce-3}) that
\beq\label{extnn6n1}
|T_i(f)-\frac{\af_i}{\bt_i}\int_{Y_{k-1}\cap X_{k-1}^{\dt_0}} { {\rm tr}_{k-1,i}}( { \Pi_{k-1}}(f))d\mu_{k-1,i} |<\frac{\ep\af_i}{2^{m+1}}
\eneq
for each $f\in {\cal H}$. Also,  (\ref{extnn-3n}) implies
 \beq\label{extnn6n1-1}
&&|\rho_i(f)-\frac{\af_i}{\bt_i}\int_{Y_{k-1}\cap X_{k-1}^{\dt_0}} { {\rm tr}_{k-1,i}}( { \Pi_{k-1}}(f))d\mu_{k-1,i} |\\
 &=&\eta(f) +|\left(a_{k-1,i}-\frac{\af_i}{\bt_i}\right) \int_{Y_{k-1}\cap X_{k-1}^{\dt_0}} { {\rm tr}_{k-1,i}}( { \Pi_{k-1}}(f))d\mu_{k-1,i} |\\
& < & a_{k-1,i}\sigma_2/2^{m}+|a_{k-1,i}\bt_i -\af_i|<a_{k-1,i}\sigma_2^2/2^{m-1}
\eneq
for each $f\in H.$

{ Note that $\af_i<a_{k-1,i}$ and $\sigma_2<\ep/64$. We have
\beq\label{extnn6n2}
|T_i(f)-\rho_i(f)|<\frac{a_{k-1,i}\ep}{2^m} ~~~{\mbox{for} }~~~f\in {\cal H}.
\eneq}
It follows from (\ref{extnn-3})  that
\beq\label{extnn-6n}
|a_{k-1,i}\int_{Y_{k-1}\cap X_{k-1}^{\dt_0}} {\rm tr}_{k-1,i}(f)d\mu_i-T_i(f)|{ =|\rho_i(f)-T_i(f)|+{{|}}\eta(f)|<a_{k-1, i}\ep}
%|a_{k-1,i}\bt_i-\af_i|+\af_i\sigma_2/2^{m+1}<a_{k-1, i}\sigma_2
\eneq
for all $f\in {\cal H}.$
Since each projection $\Pi_{k-1}(p)$ is constant on each open set $Y_{k-1}$
% one can define
% $L_{k-1,i,j}={\rm tr}_{k-1,i}(p).$
and since $\Gamma_{k-1}\circ \Pi_{k-{ 2}}(p)={\lambda_{ k-1}(\Pi_{ k-1}(p))|_{Z_{ k-1}}},$ one checks
that
\beq\label{extnn-7}
T_i(p)&=& {\af_i\over{\bt_i}} \int_{Y_{k-1}\cap X_{k-1}^{\dt_0}} {\rm tr}_{k-1,i}(s^{k-1,\dt_0}\circ \Pi_{k-1}(p_j))d\mu_{k-1,i}\\ \label{extnn-7+1}
&=&\af_i\cdot {\rm tr}_{k-1,i}(p)=\rho_i(p)
\eneq
for all $p\in {\cal P}.$
Define
% { {I suppose that $t_{k-1}$ is already given by induction, then construct $t_{k-2}$}}
$$t_{k-1}''(f)=t_{k-2}'(\Pi_{k-{{2}}}(f))+\sum_{i=1}^{l(k-1)}  T_i(f)
%{\af_i\over{\bt_i}} a_{k-1,i}\int_{Y_{k-1}\cap X_{k-1}^{\dt_0}}
 %{\rm tr}_{k-1, i}\circ s^{k-1,\dt_0}(\lambda_{k-1}(f))d\mu_j
 +{ \sum_{i=1}^{l(k-1)}}\sum_{j=1}^{m({{k-1,}} i)} d(k-1,i,j){\rm tr}_{{k-1,i}}(f(t_{i,j}))
 $$
 for all $f\in A$ and {{define}}
 \beq\label{164-n1}
 t_{k-2}(a)=t_{k-2}'(a)+\sum_{i=1}^{l(k-1)}
{{{\af_i\over{\bt_i}} \int}}_{Y_{k-1}\cap X_{k-1}^{\dt_0}}
 {\rm tr}_{k-1, i}\circ s^{k-1,\dt_0}(\Pi_{k-{ 2}}(a))d\mu_{k-1,i}
\eneq
for all $a\in A_{k-2}.$
We estimate that
\beq\label{16-univ-K-10}
|t_{k-1}''(h)-t_{k-1}(h)|{ =\sum_{i=1}^{l(k-1)}|T_i(f)-\rho_i(f)|}<(\sum_{i=1}^{l(k-1)}a_{k-1,i}){ \ep}/2^m\le { \ep}/2^m
\eneq
for all $h\in {\cal H}.$ Therefore
\beq\label{16-univ-K-11}
|\tau(h)-\sum_{s=k-1}^m({{\sum_{i=1}^{l(s)}\sum_{j=1}^{m(s,i)}}} d(s,i,j) {\rm tr}_{s,{{i}}}(\lambda_s(h)(y_{s,i, j}))+t_{k-2}\circ \Pi_{k-{ 2}}(h)|
 \le (m-k+2){\ep}/2^{m}
\eneq
for all $h\in {\cal H}$ and
\beq\label{16-univ-K-12}
\tau(p)=\sum_{s=k-1}^m(\sum_{i,j} d(s,i,j) {\rm tr}_{s,{{i}}}(\lambda_s(p)(y_{s,i, j}))+t_{k-2}\circ \Pi_{k-{ 2}}(p)
\eneq
for all $p\in {\cal P}.$
{{Note that $y_{s,i,j}\in Y_s\setminus X_s^{\dt_0}$ and $(1-\chi_s)(y_{s,i,j})=0.$}
%So we have
%\beq
%\tau_{k-1}(\Pi_{k-1}(h))=\tau_{k-1}((1-\chi_k)
%\eneq
%{{We also note that $\Pi_{k-2}\circ L_{k}\circ \Pi_{k-2}=\Pi_{k-2}.$}}
We {then, by (\ref{extnn-15-2}),(\ref{164-n1}) and using (\ref{extnn-3++a}),}} compute that, for all $ h\in {\cal H}_1,$
\beq\label{164-n2}
\hspace{-0.2in}{ {t}}_{k-1}(L_{k-1}(\Pi_{k-2}(h)))=\sum_{{ {i}}=1}^{l(k-1)} a_{k-1,{ {i}}}\int_{Y_{k-1}}{\rm tr}_{k-1,i} (L_{k-1}(\Pi_{k-2}(h))d\mu_{k-1,{ {i}}}+t'_{k-2}\circ \Pi_{ k-2}(h)\\
=\sum_{{ {i}}=1}^{l(k-1)} a_{k-1,{ {i}}}\int_{Y_{k-1}\cap {{X}}_{k-1}^{\dt_0} }{\rm tr}_{k-1,i} (L_{k-1}(\Pi_{k-2}(h))d\mu_{k-1,{ {i}}}+t'_{k-2}\circ \Pi_{ k-2}(h)\\
\hspace{-0.8in}=\sum_{{ {i}}=1}^{l(k-1)} {\af_i\over{\bt_i}}\int_{Y_{k-1}\cap  {{X}}_{k-1}^{\dt_0} }{\rm tr}_{k-1,i} (L_{k-1}(\Pi_{k-2}(h))d\mu_{k-1,{ {i}}}+t'_{k-2}\circ \Pi_{ k-2}(h)\\
+\sum_{{ {i}}=1}^{l(k-1)} (a_{k-1, i}-{\af\over{\bt_i}})\int_{Y_{k-1}\cap  {{X}}_{k-1}^{\dt_0} }{\rm tr}_{k-1,i} (L_{k-1}(\Pi_{k-2}(h))d\mu_{k-1,{ {i}}}\\
\le \sum_{{ {i}}=1}^{l(k-1)} {\af_i\over{\bt_i}}\int_{Y_{k-1}\cap  {{X}}_{k-1}^{\dt_0} }{\rm tr}_{k-1,i} (s^{k-1,\dt_0}(\Pi_{k-{ 2}}(h))d\mu_{k-1,{ {i}}}\\
+t'_{k-2}\circ \Pi_{ k-2}(h)+\sigma_2/2^m={ {t}}_{k-{{2}}}(\Pi_{k-2}(h)) +\sigma_2/2^m.
\eneq
It follows that {{(using (\ref{extnn-G})}} ,  for all $h\in {\cal H}_1,$
\beq\label{16-univ-K-13}
t_{k-2}\circ \Pi_{k-2}(h)\ge t_{k-1}(L_{k-1}(\Pi_{k-2}(h)))-\sigma_2/2^m\\
 \ge \Delta_{{k}}(\widehat{(L_{k-1}(\Pi_{k-2}(h)))}) -(m-k+2)\sigma_2/2^m\\
 =\Delta_{k-{{1}}}(\widehat{\Pi_{k-2}(h)})-(m-k+2)\sigma_2/2^m.
\eneq
This completes the induction (stop at $k-2=0$) and proves the claim. Then define $\Delta_0=\Delta_1/2.$ { {The inequality  (\ref{Reduc-1+}) and first half of both  inequality (\ref{Reduc-1}) and equation (\ref{Reduc-1+}) follows from the claim for $k=1$. Of course we denote $\Lambda_s(h)(y_{s,i,j})$ by $\pi_i(h)$ for certain irreducible representation $\pi_i$, and denote $d(s,i,j)$ by $r_i$. Furthermore we have $\tilde K r_i\in )_+$, for $ i=1,2,\cdots,k.$ Note so far we only have definition of $r_i$ for first half of (\ref{Reduc-1}).}}

{ {Note that from our construction, there are $h_j\in {\cal H}_1\subset  A_+^{\bf 1}$ with
$\Pi_0(h_j)=e_j,$ $j=1,2,...,l$ .
%Let
%$$
%\sigma_1=\min\{\Delta_0(e_j): 1\le j\le 1\}\andeqn \sigma_2=\sigma_1\sigma/4.
%$$
%}}
%{\red{Note that $\Delta_0(e_j)$ depends on $\ep$ but not $\tau$ or $\kappa.$ Therefore
%$\sigma_1.$}}

%Let $K'\ge 1$ be an integer (in place of {{the product $MR$ there}}) in \ref{multiple-ext} associated
%with $G=(\Pi_0)_{*0}(K_0(A))={{(K_0(A))/\ker(\rho_A)}}\subset \Z^{l}$ {{(see (\ref{K0nccw-2}))}}, $\sigma_1$ and $\sigma_2.$
%Let $\{d_1, d_2,..., d_m\}$ be the set of ranks of all possible irreducible representations of $A.$
%Put $K_3=\sum_{i=1}^m R_i.$
%{{Let $K''=(\prod_{j=1}^md_j)$ and
%let $K=\tilde KK'K''.$}}

 {{ Write $\tau_0=\sum_{j=1}^{l}\af'_jt_j\circ S_j,$ where $t_j$ is a tracial state on $M_{r(0,j)},$
$S_j: A_0\to M_{r(0,j)}$ is the quotient map. Then from (\ref{Reduc-1++}) and $e_j\in \Pi_0({\cal H}_1)$,
we know $\af'_j\geq \sigma_1$. From the first half of (\ref{Reduc-1+}), $\tilde Kr_i\in )_+$ , for $i=1,2,\cdots k$, and (\ref{extnn-0+}),
we have that
$$s\tilde K K''\tau_0\circ \Pi_0(p)=s \tilde K K''\big(\tau(p)-\sum_{i=1}^kr_i{\rm tr}_i \circ \pi_i(p)\big)\in \Z, \tforal p\in {\cal P}.$$
}}

{ {Let ${{\bar \af}}_j =s\tilde K K'' \af'_j$, then from the above, we have for any $(m_1,m_2,\cdots, m_l)\in (\Pi_0)_{*0}K_0(A)\subset K_0(A_0)=\Z^l$, we have $\sum_{j=1}^l{\bar{\af_j}} m_j \in \Z$. Apply (\ref{multiple-ext}), we have ${\bar{\bt}}_j$ ($1\le j\le l$) such that$K'{\bar\bt}_j\in Z^+$,  $\sum_j\|{{\bar\bt}}_j-{{\bar\af}}_j\|<\sigma_2$, and $\sum_{j=1}^l {{\bar\bt}}_j m_j =\sum_{j=1}^l {{\bar\af}}_j m_j $ for all $(m_1,m_2,\cdots, m_l)\in (\Pi_0)_{*0}K_0(A)$. Let $r_{k+j}= \frac1{s\tilde KK''}{\bar \bt}_j$ and $\pi_{k+j}=S_j$ for $j=1,2,\cdots,  l$. The lemma {{now}} follows easily.}}
}}}

\end{proof}

%%%%%%%%%

\begin{lem}\label{ExtTraceMn}
%Let $C\in\mathcal C_0$ with $K_1(C)=\{0\}$. Let $A=C$ or $A=C\otimes\mathrm{C}(\mathbb T)$.
Let $A\in {\cal D}_m$ be  a unital \CA.
%Let $A=C$ for some $C\in\mathcal C$ or
%$A=C\otimes C(\T)$ for some $C\in\mathcal C$.
Let $\Delta: A_+^{q, 1}\setminus\{0\}\to(0, 1)$ be an order preserving map. Let $\mathcal H\subseteq A$ be a finite subset and let $\sigma>0$. Then there are finite subset $\mathcal H_1\subseteq  {{A^{\bf 1}_+}}\setminus \{0\}$ and a positive integer $K$ such that for any  $\tau\in T(A)$ satisfying
\begin{equation}\label{extnn-0-1}
\tau(h) > \Delta(\hat{h}) \rforal h\in \mathcal H_{1}
\end{equation}
 and any positive homomorphism $\kappa: K_0(A)\to K_0(M_s)$
with $s=\kappa([1_A])$
such that
\begin{equation}\label{extnn-0-1+}
\rho_A(x)(\tau)=(1/s)(\kappa(x))\tforal x\in K_0(A),
\end{equation}
there is a unital homomorphism $\phi: A\to M_{sK}$ such that
$\phi_{*0}=K\kappa$ and $$|{\rm tr}'\circ\phi(h)-\tau(h)|<\sigma \tforal h\in\mathcal H,$$
where ${\rm tr}'$ is the tracial state on $M_{sK}.$
\end{lem}

\begin{proof} Lemma follows easily from \ref{Reduction2015} .
%and \ref{multiple-ext}.

Without loss of generality, we may assume that $A\in {\cal D}_m',$ and we may also assume that projections in $A$ generates $K_0(A),$ by replacing $A$ by $M_N(A)$
for some integer $N \ge 1.$ Let ${\cal P}\subset A$ be a finite subset  of projections
such that $\{\rho_A([p]):p\in {\cal P}\}$ generates
 %which generates
${ {\rho_{A}}}(K_0(A))_+.$
Write $A_0=\bigoplus_{j=1}^{l(0)}M_{r(0,j)}$ and $e_j\in M_{r(0,j)}$ is a rank one projection.
%P_{0,j}C(X_{0,j}, F_0)P_{0,j},$
%where $X_{0,j}$ is connected and $P_{0,j}\in C(X_{0,j}, F_0)$ has rank
%$r(0,j).$ By choosing possibly larger
%$m,$ \wilog, we may further assume that $P_{0,j}C(X_{0,j}, F_0))P_{0,j}$ contains a rank one projection
%$e_j.$

%{\red ---that is for any element $x\in K_0(A)_+$,  an element $y$ which is summation of some $K_0$ elements represented by projection in ${\cal P}$, such that $\rho_A(x)=\rho_A(y)$}.

Let $\Delta_0: \Pi_0(A)=A_0: (A_0)_+^{q, \bf 1}\setminus \{0\}\to (0,1),$
${\cal H}_1\subset A_+^{\bf 1}\setminus \{0\}$ and  $K_1$ be the integer $K$ in the statement of
\ref{Reduction2015} for $\sigma/4$ (in place of $\ep$), ${\cal H}$ and $\Delta.$
%We may assume that ${\cal P}$ contains a finite set of projections
%which generates $K_0(A)_+.$
%Write $A_0=\Pi_0(A)=M_{r(1)}\oplus M_{r(2)}\oplus \cdots \oplus M_{r(l)}.$
%Let $e_j\in M_{r(j)}\subset A_0$ be a rank one projection.
%Choose  a finite subset ${\cal H}_1'$ which
%contains ${\cal H}_1$ as well as some $h_j\in A_+^{\bf 1}$ with
%$\Pi_0(h_j)=e_j,$ $j=1,2,...,l$ (in fact, in the proof of \ref{Reduction2015}, we have such $h_j$ in ${\cal H}_1$).
%Let
%$$
%\sigma_1=\min\{\Delta_0(e_j): 1\le j\le 1\}\andeqn \sigma_2=\sigma_1\sigma/4.
%$$
%Let $K_2\ge 1$ be an integer (in place of {{the product $MR$ there}}) in \ref{multiple-ext} associated
%with $G={{(K_0(A))/\ker(\rho_A)}}\subset \Z^{l(0)}$ {{(see (\ref{K0nccw-2}))}}, $\sigma_1$ and $\sigma_2.$
Let $\{d_1, d_2,..., d_m\}$ be the set of ranks of all possible irreducible representations of $A.$
%Put $K_3=\sum_{i=1}^m R_i.$
{{Let $K_2=(\prod_{j=1}^md_j)$ and
let $K=K_1K_2.$}}
% {\red{(\prod_{j=1}^mR_j)}}.$

For given $\kappa$ and $\tau$ as described in the statement,
%there are integers $K_0\ge 1,$
rational numbers $r_1,r_2,..,r_k\ge 1,$ finite dimensional irreducible representations
$\pi_1, \pi_2,...,\pi_k$ of $A,$ {{irreducible representations}} of $A_0$ such that
%and a  trace $\tau_0\in T(A_0)$ such that
\beq\label{Reduc-1-1}
&&|\tau(h)-(\sum_{i=1}^k r_i {\rm tr}_i\circ \pi_i(h)+\sum_{i=k+1}^{k+l}r_{i}{\rm tr}_{i}\pi_i\circ \Pi_0(h)
%\tau_0\circ \Pi_0(h))
|<\sigma/4,\\
&&\tau(p)=\sum_{i=1}^kr_i{\rm tr}_i \circ \pi_i(p)+\sum_{i=k+1}^{k+l}r_{i}{\rm tr}_{i}\pi_i\circ \Pi_0(p)
%\tau_0\circ \Pi_0(p)
%,\quad \forall
\rforal p\in {\cal P},
%\tau_0(\Pi_0(h))
%&&{\red{\sum_{i=k+1}^{k+l}r_{i}{\rm tr}_{i}\pi_i\circ \Pi_0(h)}}\ge \Delta_0({ \widehat{\Pi_0(h)}})
%, \quad \forall
%\rforal h\in {\cal H}_1',
\eneq
and  $sK_1r_i\in \Z.$
%quad \textrm{and}\quad K_0|K_1.$$

Define $r'_i=r_i/R(i)$ and ${\rm Tr}_i=R(i){\rm tr}_i,$  where ${{R(i)}}$ is rank of $\pi_i$.
{{So ${\rm Tr}_i$ is an unnormalized trace on $M_{R(i)}.$}} Then we have
\beq\label{Reduc-1-G}
&&|\tau(h)-(\sum_{i=1}^k r'_i {\rm Tr}_i\circ \pi_i(h)+{{\sum_{i=k+1}^{k+l}r'_{i}{\rm Tr}_{i}\pi_i\circ \Pi_0(h)}}
%\tau_0\circ \Pi_0(h))
|<\sigma/4,\\
&&\tau(p)=\sum_{i=1}^kr'_i{\rm Tr}_i \circ \pi_i(p)+{{\sum_{i=k+1}^{k+l}r'_{i}{\rm Tr}_{i}\pi_i\circ \Pi_0(p)}}
%\tau_0\circ \Pi_0(p)
%,\quad \forall
\rforal p\in {\cal P},
%\tau_0(\Pi_0(h))
%&&{\red{\sum_{i=k+1}^{k+l}r_{i}{\rm tr}_{i}\pi_i\circ \Pi_0(h)}}\ge \Delta_0({ \widehat{\Pi_0(h)}})
%, \quad \forall
%\rforal h\in {\cal H}_1',
\eneq
and $ sK_1K_2r'_i\in \Z.$

Define $T:A\to \C$ by
$T(a)=\sum_{i=1}^k r{{'}}_i {{{\rm Tr}}}_i\circ \pi_i(a)+{{\sum_{i=k+1}^{k+l}  r{{'}}_i {{{\rm Tr}}_i}\circ \pi_i\circ \Pi_0(a)}}$
%\tau_0\circ \Pi_0(a)
for all $a\in A.$ Note that $T\in T(A).$
Therefore
\begin{eqnarray}\label{extnn-9}
(1/s)\circ \kappa(p)=T(p)\tforal p\in K_0(A)\andeqn\\\label{extnn-9+}
|\tau(h)-T(h)|<\sigma/4\tforal h\in {\cal H}.
\end{eqnarray}
Put $m_0=\sum_{i=1}^{k+l} m(i)R(i),$
%+\sum_{i=k+1}^{k+l}m(i)r_{0,i},$
where $R_i$ is rank of $\pi_i$ and $m(i)={{s}}{{K_1K_2}}r{{'}}_i.$  Define
$\phi: A\to M_{m_0}$ by $\Psi(a)=\bigoplus_{i=1}^k {\bar \pi}_i(a)\bigoplus\bigoplus_{i=k+1}^{k+1}{\bar \pi}_i(\Pi_0(a))$
for all $a\in A,$ where ${\bar \pi}_i$ is $m(i)$ copies of $\pi_i.$ Let $\kappa_\phi: K_0(A)\to \Z$ be the map induced  by $\phi.$
Then { {using (\ref{extnn-9}), we have,}} $\kappa_{{\phi}}([1_A])=K_1K_2\kappa([1_A])).$ Therefore $m_0=sK_1K_2=sK,$
$tr'(p)=T(p)$ for all projections $p\in {\cal P},$ where $tr'$ is the normalized trace on $M_{sK}$ and
$tr'(h)=T(h)$ for all $h\in A.$  This completes the proof.

\end{proof}

%{\color{Green}(What is the role of $\sigma_1$ in \ref{appextnn})?}
%{{\bf I think one could do a better job than the following. But I am eager to move beyond that.}} {\color{green} This is okay to me. One also does not need the assumption that $K_1(C)=\{0\}$.}
\begin{lem}\label{appextnn}
%Let $C\in\mathcal C_0$ with $K_1(C)=\{0\}$. Let $A=C$ or $A=C\otimes\mathrm{C}(\mathbb T)$.
Let $A\in {\cal D}_m$ be a unital \CA\, ($m\ge 1$).
Let $\Delta: A_+^{q, {\bf 1}}\setminus\{0\}\to (0,1)$ be an order preserving map.  Let $\mathcal H\subseteq A$ be a finite subset and let $\epsilon>0$.
%Let $\mathcal G$ be a finite generating set for $K_0(A)$ (as an abelian group).
There exist a finite subset $\mathcal H_1\subseteq A_+^{\bf 1}\setminus \{0\}$  and
a finite subset of projections ${\cal P}\subset M_n(A)$ (for some $n\ge 1$),
and
%for any $\sigma_1>0$,
there is $\delta>0$ such that if a tracial state $\tau\in T(A)$ satisfies
$$
\tau(h)>\Delta(\hat{h})\tforal h\in\mathcal H_1
$$
and  any order-unit map $\kappa: K_0(A)\to K_0(M_s)$ satisfying
\begin{equation}\label{extnn-0++}
|\mathrm{tr}(\kappa([p]))-\tau(p)|<\delta
\end{equation}
for all projections $p\in  {\cal P},$
where $\mathrm{tr}$ is the canonical tracial state on $M_s$, there is a tracial state $\tilde{\tau}\in T(A)$ such that
$$\mathrm{tr}(\kappa(p))=\tilde{\tau}(p) \rforal p\in K_0(A)\tand
|\tau(h)-\tilde{\tau}(h)|<\epsilon\rforal h\in \mathcal H.$$
\end{lem}

\begin{proof}
%We may assume that $s>1.$
%To simplify notation, \wilog, by replacing $A$ by $M_n(A)$ if necessary, we may assume
Without loss of generality, we may assume that
%$A\in {\cal D}_m',$ and we may also assume that
projections in $A$ generates $K_0(A),$ by replacing $A$ by $M_N(A)$
for some integer $N \ge 1.$ Let ${\cal P}\subset A$ be a finite subset  of projections
such that $\{\rho_A([p]):p\in {\cal P}\}$ generates
 %which generates
${ {\rho_{A}}}(K_0(A))_+.$
%Let $\sigma_0=\min\{\Delta(\hat{p}): p\in {\cal P}\}.$
Let ${\cal H}_1={\cal H}\cup {\cal P}$ and let
$\sigma_0=\min\{\Delta(\hat{h}): h\in {\cal H}_1\}.$
Let $\dt=\min\{\ep\cdot\sigma_0/128, 1/16\}.$

Now suppose that $\tau$ and $\kappa$ satisfies the assumption for the above mentioned ${\cal H}_1,$
${\cal P}$ and $\dt.$
Recall  $\tau$ may be viewed as an order preserving map from
$K_0(A)$ to $\R.$
%Define $\kappa_1: K_0(A)\to K_0(M_s)=\Z$ by $\kappa_1([p])=s\tau(p)$ for all projection $p\in A.$
%Note that $\kappa_1$  is well-defined, $\kappa_1([1_A])=s$ and order preserving.
Define $\eta=(1-\ep/3)(\kappa/s-\tau): K_0(A)\to \R.$
Let $d_A: K_0(A)\to K_0(A)/{\rm ker}\rho_A$   be a quotient map
and let  $\gamma: K_0(A)/{\rm ker}\rho_A\to \R$ given  such that
$\gamma\circ d_A=(\ep/3s)\kappa+\eta.$

Note that, as in \ref{8-n-K0nccw} ,
\beq\label{166-1}
K_0(A)/{\rm ker}\rho_A\cong (\Pi_0)_{*0}(K_0(A))\subset K_0(A_0).
\eneq
So one may view that $\gamma$ is a \hm\, from  $(\Pi_0)_{*0}(K_0(A))$ to
$\R.$
For each $p\in {\cal P},$ from the assumption (\ref{extnn-0++}), one computes that
\beq\label{166-2}
|\eta([p])|<(1-\ep/3)\dt.
\eneq
Therefore
\beq\label{166-3}
\gamma((\Pi_0)_{*0}([p]))&=&(\ep/3)(\kappa([p])/s)+\eta([p])>
(\ep/3)(\Delta(\hat{p})-\dt)-(1-\ep/3)\dt\\
&\ge & (\ep/3)(1-1/128)\sigma_0-(1-\ep/3)\ep\sigma_0/128\\
&=& \ep\sigma_0[({1\over{3}}-{1\over{3\cdot 128}})-({1\over{128}}-{\ep\over{3\cdot 128}})]>0
\eneq
for all $p\in {\cal P}.$
In other words, $\gamma$ is positive. By applying 2.8 of \cite{Lnbirr}, there
is a positive \hm\, $\gamma_1: K_0(A_0)\to \R$ such
that $\gamma_1\circ (\Pi_0)_{*0}=\gamma.$
It is well-known that there is a (non-normalized) trace  $T_0$ on $A_0$ such
that $\gamma([p])=T_0(p)$ for all projections $p\in A.$

Consider the trace $\tau'=(1-\ep/3)\tau+T_0\circ \Pi_0$ on $A.$
Then, for all projection $p\in A,$
\beq\label{166-4}
\tau'(p)&=&(1-\ep/3)\tau(p)+T_0\circ \Pi_0(p)=(1-\ep/3)\tau(p)+(\ep/3s)\kappa([p])+\eta([p])\\
&=&(1-\ep/3)\tau(p) +(\ep/3s)\kappa([p]) +(1-\ep/3)(\kappa([p])/s-\tau(p))=(1/s)\kappa([p]).
\eneq
Since $(1/s)\kappa([1_A])=1,$ $\tau'\in T(A).$
We also compute that, by  (\ref{166-2})
\beq\label{166-5}
|T_0\circ \Pi_0(1_A)|=|\gamma\circ \rho_A([1_A])|<\ep/2
\eneq
Therefore, we also have
\beq\label{166-6}
|\tau'(h)-\tau(h)|<\ep\rforal h\in {\cal H}.
\eneq

\end{proof}

%{\color{Green} In \ref{APPextnn}, only being approximately compatible to the projection in $A$ is not enough, since a minimal projection of $M_\infty(A)$ might not in $A$; in other words, the K-group might not be generated by projections in $A$. For instance, consider
%$$K_0(A)=\{(x, y, z)\in\Z^3: x+y=2z\}$$
%with $[1_A]=(1, 1, 1)$. Then $A$ is projectionless by itself, but $K_0(A)=\Z^2$ and the other generator is not in $A$. The same problem with \ref{ExtTraceI-D} and maybe others.
%}
\begin{lem}\label{APPextnn}
%Let $C\in\mathcal C_0$ with $K_1(C)=\{0\}$. Let $A=C$ or $A=C\otimes\mathrm{C}(\mathbb T)$.
Let $A\in {\cal D}_m$ be a unital \CA.
%=C$ for some $C\in\mathcal C$ or
%$A=C\otimes C(\T)$ for some $C\in\mathcal C$.
Let $\Delta: A_+^{q, {{\bf1}}}\setminus\{0\}\to(0, 1)$ be an order preserving map. Let $\mathcal H\subseteq A$ be a finite subset and let $\sigma>0$. Then there are  a finite subset $\mathcal H_1\subseteq A_+^{\bf 1}\setminus \{0\},$
%a finite set of projections ${\cal P}\subset M_n(C)$ (for some $n\ge 1$),
$\dt>0,$ a finite subset ${\cal P}\subset K_0(A)$ and a positive integer $K$ such that for any  $\tau\in T(A)$ satisfying
\begin{equation*}
\tau(h) > \Delta(\hat{h}) \rforal h\in \mathcal H_{1}
\end{equation*}
 and any { {positive homomorphism} } $\kappa: K_0(A)\to K_0(M_s)$
with $s=\kappa([1_A])$
such that
\begin{equation*}
|\rho_A(x)(\tau)-(1/s)(\kappa(x))|<\dt
\end{equation*}
for all $x\in {\cal P},$
there is a unital homomorphism $\phi: A\to M_{sK}$ such that
$\phi_{*0}=K\kappa$ and $$|{\rm tr}'\circ\phi(h)-\tau(h)|<\sigma \tforal h\in\mathcal H,$$
where ${\rm tr}'$ is the tracial state on $M_{sK}.$
\end{lem}

\begin{proof}
Note that there is an integer $n\ge 1$ such that projections in $M_n(A)$ generate $K_0(A).$
Therefore this lemma  is a corollary of \ref{ExtTraceMn} and \ref{appextnn}.
\end{proof}

{\begin{rem}\label{Rextnn}
%When $K_0(A)$ is finitely generated, in what follows, we may choose ${\cal P}$ to be a set of generators {\blue{of $K_0(A)}}.

%{The following sentence is not correct, so I delete it and replace it by another sentence and hope it is what we need.}
%Moreover, ${\cal P}$ could be the set of images of all projections of $M_n(A)$ for a suitable integer $n.$
%It should be noted that such ${\cal P}$ is a finite subset of $K_0(A).$

Since $K_0(A)$ is finite generated and $(\Pi_0)_{*0}(K_0(A)_+) \subset \Z_+^l$ is also finitely generated (though $K_0(A)_+$ itself may not be finitely generated), in what follows, we can choose finite subset ${\cal P}$  to be a set of generators of $K_0(A)$ as  {{abelian}} group, and $\{(\Pi_0)_{*0}([p])\}_{p\in{\cal P}}$ generate $(\Pi_0)_{*0}(K_0(A)_+)$ as {{abelian}} semigroup.  In particular, we can choose $n$ such that ${\cal P}\subset M_n(A)$.

\end{rem}}

%{\blue{Let $\phi: A \to B$ be a completely positive map, then the notation $\phi\otimes 1_{M_N(\C)}$ may be one of the following two maps:

%(a)for any $x\in A $ and $a\in M_N(\C)$, $\phi\otimes 1_{M_N(\C)} (x\otimes a) = \phi(x)\otimes a\in B\otimes M_N(\C)$ , which is from $A\otimes M_N(\C)$ to $B\otimes M_N(\C)$;
%
%(b)for any $x\in A $,  $\phi\otimes 1_{M_N(\C)} (x)=\diag\{\underbrace{x,x,\cdots,x}_N\} \in B\otimes M_N(\C)$, which is from $A$ to $B\otimes M_N(\C)$.}}
%{\red{I am not prepared to make announcement (b) in section 2 yet.
  %I will prefer to, whenever you find them, use $\phi(x)\otimes 1_{M_N}.$---does it make sense?}}

%{\blue{This will not cause any confusion, since it will clear what is the domain from the content.}}

\begin{lem}\label{ExtTraceI-D}
%Let $B\in\mathcal C_0$ with $K_1(B)=\{0\}$ and let $C=B$ or
%$C=B\otimes C(\T).$
Let $C\in {\cal D}_m$ be a unital C*-algebra with finitely generated $K_i({ {C}})$ ($i=0,1$).
Let $\Delta: C_+^{q, 1}\setminus\{0\}\to(0, 1)$ be an order preserving map. Let $\mathcal F, \mathcal H\subseteq C$ be  finite subsets, and let $\epsilon>0, \sigma>0$.
%Let $\mathcal H\subseteq C$ be a finite subset and let $\sigma>0$.
Then there are a finite subset $\mathcal H_1\subseteq C_+^{\bf 1}\setminus \{0\}$, $\delta>0$,
a finite subset ${\cal P}\subset K_0(C)$ and a positive integer $K$ such that for any continuous affine map $\gamma: T(C([0, 1]))\to T(C)$ satisfying
$$
\gamma(\tau)(h) > \Delta(\hat{h}) \tforal h\in \mathcal H_{1} \andeqn \rforal \tau\in T(C([0, 1]))
$$
and  any positive homomorphism $\kappa: K_0(C)\to K_0(M_s(C([0, 1])))$ with $\kappa([1_C])=s$
such that
$$
|\rho_A(x)(\gamma(\tau))-(1/s)\tau(\kappa(x))|<\delta  \tforal \tau\in T(C([0, 1]))
$$
for all $x\in {\cal P},$
%projections in $C,$
there is an %homomorphism
$\mathcal F$-$\epsilon$-multiplicative completely positive linear map
$\phi: C\to M_{sK}(C([0, 1]))$ such that
$\phi_0=K\kappa$ and
$$|\tau\circ\phi(h)-\gamma'(\tau)(h)|<\sigma\tforal h\in\mathcal H,$$
where $\gamma': T(M_{sK}(C([0,1])))\to T(C)$ is induced by $\gamma.$ { {Furthermore $\phi_0=\pi_0\circ \phi$ and $\phi_1=\pi_1\circ \phi$ are true homomorphisms.}}
In the case that $C\in\mathcal C$, the map $\phi$ can be chosen to be a homomorphism.
\end{lem}

\begin{proof}
Since any
%C*-algebra
\CA s in $\mathcal C$ are semi-projective,
%(see Theorem 6.22 of \cite{ELP1}),
the second part of the statement follows directly from the first part of the statement. Thus, let us only show the first part of the statement. %Noting that $C$ is finitely generated,
Without loss of generality, one may assume that $\mathcal F\subseteq \mathcal H$.
To simplify notation, without loss of generality, by replacing  $C$ by $M_r(C)$  for some
$r\ge 1,$ we may assume that the  set of projections in $C$ generates $\rho_A(K_0(C)).$

%Since $C$ is weakly semi-projective (see Theorem 6.22 of \cite{ELP1}), there is $\delta'>0$ such that if $L: C\to A$ is a unital $\mathcal H$-$\delta'$-multiplicative map for a unital C*-algebra $A$, there is a homomorphism $\phi: C\to A$ such that $$||\phi(h)-L(h)||<\sigma/2\rforal h\in\mathcal H.$$ One may also assume that $\delta'$ is small enough so that $[L]|_{K_0(C)}$ is well defined and $[L]|_{K_0(A)}=\phi_{*0}$.

Since the K-theory of $C$ is finitely generated, there is $M\in\mathbb N$ such that
$$Mp=0 \rforal p\in {\rm Tor}(K_i(C)), i=0, 1.$$
% n=1, 2,....$$

Let $\mathcal H_{1, 1}\subset C_+^{1}\setminus \{0\}$ (in place of $\mathcal H_{ {1}}$)
%$\mathcal G_{1}\subseteq C_{s.a.}$ (in place of ${ {\mathcal H_2}}$)
and $\sigma_{1}>0$ (in place of ${{\delta}}$) be the finite subsets and the positive constant of { {Theorem \ref{UniqAtoM} }} with respect to $C$ (in the place of ${{A}}$),
% { {we need to change the class ${\cal A}_s$ to ${\cal D}_m$ in \ref{CNewuni1} or whole section 5}}
$\min\{\sigma, \epsilon \}$ (in the place of $\epsilon$) and $\mathcal H$ (in the place of $\mathcal F$), and $\Delta/2$.
{{(We will not need the finite set ${\cal P}$ in Theorem \ref{UniqAtoM}, since $K_0(C)$ is finitely generated and when we apply Theorem \ref{UniqAtoM}, we will require both map induce same $KL$ maps.)}}

Let $\mathcal H_{1, 2}\subseteq C$ (in the place of $\mathcal H_1$)  be a  finite subset, $\delta>0$
%(in place of $\delta$)
be a positive number, ${\cal P}\subset K_0(C)$ be a finite subset
%the finite subset
and $K'$  be an integer required by Lemma \ref{APPextnn} with respect to $C$, $\Delta/2$ (in place of $\Delta$),  $\mathcal H\cup \mathcal H_{1, 1}$
% \cup \mathcal G_1$
(in the place of $\mathcal H$) and
$\min\{\sigma/16, \sigma_1/8,\{\Delta(\hat{h})/4: h\in\mathcal H_{1, 1}\}\}$ (in the place of $\sigma$).
{{W}}e may assume that ${\cal P}$ is the set of projections in $C,$ as mentioned in \ref{Rextnn}.

%Let $M$ be the constant of Lemma \ref{ExtTraceM} with respect to $C$, $\mathcal H\cup \mathcal H_{1, 1} \cup \mathcal G_1$ (in the place of $\mathcal H$) and $\min\{\sigma/16, \Delta(\hat{h})/4, \sigma_1/8;\ h\in\mathcal H_{1, 1}\}$ (in the place of $\sigma$).

%let $\delta>0$ (in the place of $\delta$) be the constant of Lemma \ref{pert-T-K} with respect to $\mathcal H_{1, 2}$ (in the place of $\mathcal H_1$) and $\sigma_{1, 2}$ (in the place of $\sigma_1$).

Put $\mathcal H_1=\mathcal H_{1, 1}\cup\mathcal H_{1, 2}$ and $K=MK'$.
Then, let $\gamma: T(C([0, 1]))\to T(C)$ be a continuous affine map with $$\gamma(\tau)(h) > \Delta(\hat{h})\rforal h\in \mathcal H_{1},$$ and let $\kappa: K_0(C)\to K_0(M_s(C([0, 1])))$
with $\kappa([1_C])=s$
such that
$$
|\rho_\gamma(\tau)(x)-(1/s)\tau(\kappa(x))|<\delta\rforal  x\in {\cal P}
%\mathcal G
\rforal \tau\in T(C([0, 1])).
$$
%for some $s$ divided by $M$.

Since $\gamma$ is continuous, there is a partition
$$0=x_0<x_1<\cdots<x_n=1$$
such that for any $0\leq i\leq n-1$, and any $x\in[x_i, x_{i+1}]$, one has
\begin{equation}\label{eq-par}
|\gamma(\tau_x)(h)-\gamma(\tau_{x_i})(h)|< \min\{\sigma/8, \sigma_1/4\} \rforal h\in \mathcal H_1,
%\cup\mathcal G_1,
\end{equation}
where $\tau_x\in T(M_s(C([0, 1])))$ is the extremal trace which concentrates at $x$.

For any $0\leq i\leq n$, consider the trace $\tilde{\tau}_i=\gamma(\tau_{x_i})\in T(C)$. It is clear that
\beq\nonumber
|\tilde{\tau}_i(x)-\mathrm{tr}(\kappa(x))|<\delta\rforal x\in {\cal P}\andeqn
%\mathcal{G}
%$$
%and
%$$
\tilde{\tau}_i(h)>\Delta(\hat{h})\rforal h\in\mathcal H_{1, 2}.
\eneq
%Then it follows from Corollary \ref{pert-T-K-D} that there is $\tau_i\in T(C)$ such that
%\begin{equation}\label{comp-pair}
%\tau_i(p)=\mathrm{tr}(\kappa(p)),\quad \forall p\in\mathcal{G},
%\end{equation}
%and
%\begin{equation}\label{eq-pert-tr}
%|\tau_i(h)-\tilde{\tau}_i(h)|<\min\{\sigma/16, \Delta(\hat{h})/4, \sigma_1/8;\ h\in\mathcal H_{1, 1}\}\rforal h\in\mathcal H\cup \mathcal H_{1, 1} \cup \mathcal G_1.
%\end{equation}

%Define the map $\gamma_i: T(M_s(\C))\cong\{\mathrm{tr}\}\to T(C)$ by $\gamma_i(\mathrm{tr})=\tau_i$.
%Regard the map $\kappa$ as a map from $K_0(C)$ to $K_0(M_s(\C))$, and by \eqref{comp-pair}, %ne has that $(\kappa, \gamma_i)$ is a compatible pair for $C$ and $M_s(\C)$.

%Since $s$ is divisible by $M$,
By Lemma \ref{APPextnn}, there exists a unital homomorphism $\phi'_i: C\to M_{sK'}(\C)$ such that
$$[\phi_i']_0=K\kappa$$
as we identify $K_0(C([0,1], M_s))$ with $\Z$ and
\begin{equation}\label{pt-exi}
|\mathrm{tr}\circ\phi'_i(h)-\tau_{x_i}(h)|<\min\{\sigma/16, \Delta(\hat{h})/4, \sigma_1/8;\ h\in\mathcal H_{1, 1}\}\rforal h\in\mathcal H\cup \mathcal H_{1, 1}.
% \cup \mathcal G_1.
\end{equation}
%By \eqref{eq-pert-tr}, one has
%\begin{equation}\label{pt-exi}
%|\mathrm{tr}\circ\phi'_i(h)-\gamma(\tau_{x_i})(h)|<\min\{\sigma/8, \Delta(\hat{h})/2, \sigma_1/4;\ h\in\mathcal H_{1, 1}\}\rforal h\in\mathcal H\cup \mathcal H_{1, 1} \cup \mathcal G_1.
%\end{equation}
In particular, by \eqref{pt-exi}, one has that, for any $0\leq i\leq n-1$,
$$
|\mathrm{tr}\circ\phi'_i(h)-\mathrm{tr}\circ\phi'_{i+1}(h)|<\sigma_1\rforal h\in {\cal H}_{1,1}.
%\mathcal G_1.
$$
Note that $\gamma(\tau_{x_i})(h)>\Delta(\hat{h})$ for any $h\in\mathcal H_{1, 1}$ by the assumption. It then also follows from \eqref{pt-exi} that,  for any $0\leq i\leq n$,
$$\mathrm{tr}\circ\phi_i'(h)>\Delta(\hat{h})/2\rforal h\in\mathcal H_{1, 1}.$$
Define {{the amplification $\phi_i''$ as}}
$$\phi''_i:=\phi_i'\otimes 1_{\mathrm{M}_M(\mathbb C)}: C\to M_{sK}(\C){{=M_{M(sK')}(\C)}}.$$
One has that
$$[\phi''_i]=[\phi''_{i+1}]\quad \textrm{in $KL(C, M_{sK})$}.$$
It then follows from {{Theorem \ref{UniqAtoM}}} that there is a unitary $u_1\in M_s(\C)$ such that
$$\|\phi''_{0}(h)-\mathrm{Ad}u_1\circ\phi''_{1}(h)\|<\min\{\sigma, \epsilon\} \rforal h\in\mathcal H.$$
Consider the maps $\mathrm{Ad}u_1\circ\phi''_{1}$ and $\phi''_{2}$.
Applying {{Theorem \ref{UniqAtoM}}} again, one obtains a unitary $u_2\in M_{sK}(\C)$ such that
$$
\|\mathrm{Ad}u_1\circ\phi''_{1}(h)-\mathrm{Ad}u_2\circ\phi''_{2}(h)\|<\min\{\sigma, \epsilon\}\rforal h\in\mathcal H.
$$
Repeat this argument for all $i=1,...,n$, one obtains  unitaries  $u_i\in M_{sK}(\C)$ such that
$$
\|\mathrm{Ad}u_i\circ\phi''_{i}(h)-\mathrm{Ad}u_{i+1}\circ\phi''_{i+1}(h)\|<\min\{\sigma, \epsilon\},\rforal h\in\mathcal H.
$$
Then define $\phi_0=\phi_0''$ and $\phi_i=\mathrm{Ad}u_i\circ\phi''_i$, and one has
\begin{equation}\label{close-ev}
\|\phi_{i}(h)-\phi_{i+1}(h)\|<\min\{\sigma, \epsilon\}\rforal h\in\mathcal H.
\end{equation}

Define the linear map $\phi: C \to M_{sK}([0, 1])$ by
$$\phi(f)(t)=\frac{t-x_i}{x_{i+1}-x_i}\phi_i(f)+\frac{x_{i+1}-t}{x_{i+1}-x_i}\phi_{i+1}(f),\quad\textrm{if $t\in[x_i, x_{i+1}]$}.
$$
Since each $\phi_i$ is a homomorphism, by \eqref{close-ev}, the map $\phi$ is $\mathcal H$-$\epsilon$-multiplicative; in particular, it is $\mathcal F$-$\epsilon$-multiplicative.
%Therefore, there is a homomorphism $\phi: C\to M_{sK}([0, 1])$ such that
%$$\|\phi(h)-L(h)\|<\sigma/2\rforal h\in\mathcal H.$$
It is clear that $\phi_{*0}=K\kappa$. On the other hand, for any $x\in [x_i, x_{i+1}]$ for some $i=1, ..., n-1$, one has that for any $h\in\mathcal H$,
\begin{eqnarray*}
&&|\gamma(\tau_x)(h)-\tau_x\circ \phi(h)| \\
& = &|\gamma(\tau_x)(h)- (\frac{x-x_i}{x_{i+1}-x_i}\mathrm{tr}(\phi_i(f))+\frac{x_{i+1}-x}{x_{i+1}-x_i}\mathrm{tr}(\phi_{i+1}(f)))|\\
&<&|\gamma(\tau_x)(h)- (\frac{x-x_i}{x_{i+1}-x_i}\gamma(\tau_{x_i})(h)+\frac{x_{i+1}-x}{x_{i+1}-x_i}\gamma(\tau_{x_{i+1}})(h))|+\sigma/4 \quad\textrm{{(by \eqref{pt-exi})}}\\
&<&|\gamma(\tau_x)(h)- \gamma(\tau_{x_{i+1}})(h))|+3\sigma/8  \quad \quad \quad \quad \quad\quad \quad  \quad \quad \quad \quad \quad \quad \quad\,\,\,\textrm{{(}by \eqref{eq-par}{)}}\\
&<&\sigma/2,  \quad  \quad \quad \quad \quad \quad\quad\quad \quad  \quad \quad \quad \quad \quad\quad\quad \quad  \quad \quad \quad \quad \quad\quad\quad\,\,\,\,\,\,\,\,\,\,\,\,\,\,\textrm{{(by \eqref{eq-par})}}.
\end{eqnarray*}
Hence for any $h\in\mathcal H$,
$$|\gamma(\tau_x)(h)-\tau_x\circ\phi(h)|<|\gamma(\tau_x)(h)-\tau_x\circ L(h)|+\sigma/2<\sigma,$$ and therefore
$$|\gamma(\tau)(h)-\tau\circ\phi(h)|<\sigma$$
% \,\,\, \rforal h\in\mathcal H$$ and
for any $\tau\in T(M_{sK}(C([0, 1])))$.

{{ Note that the restriction of $\phi$ to the boundaries are $\phi_0$ and $\phi_n$ which are homomorphisms. }}  Thus the map $\phi$ satisfies the statement of the lemma.

\end{proof}

\begin{thm}\label{ExtTraceC-D}
 Let $C\in {\cal D}_m$ be a unital \CA\, with finite generated $K_i(C)$ ($i=0,1$).
 Let $\Delta: C_+^{q, 1}\setminus\{0\}\to (0, 1)$ be an order preserving map. Let $\mathcal F, \mathcal H\subseteq C$ be finite subsets, and let $1>\sigma, \epsilon>0$.
There exist  a finite subset $\mathcal H_1\subseteq C_+^{\bf 1}\setminus \{0\}$, $\delta>0$,
{a finite subset ${\cal P}\subset K_0(C)$}and a positive integer $K$ such that for any continuous affine map $\gamma: T(D)\to T(C)$ satisfying
$$\gamma(\tau)(h) > \Delta(\hat{h})\rforal h\in \mathcal H_1\tforal \tau\in T(D),$$
where $D$ is a C*-algebra in $\mathcal C$,
any positive \hm\, $\kappa: K_0(C)\to K_0(D)$ {with $\kappa([1_C])=s[1_D]$} for some integer $s \ge 1$ satisfying
$$|\rho_C(x)(\gamma(\tau))-{(1/s)}\tau(\kappa(x))|<\delta\tforal \tau\in T(D)$$
and for all $x\in {\cal P},$
there is a $\mathcal F$-$\epsilon$-multiplicative positive linear map $\phi: C\to M_{{sK}}(D)$ such that
$$\phi_{*0}={K}\kappa$$ and
$$|(1/(sK))\tau\circ\phi(h)-\gamma(\tau)(h)|<\sigma\tforal h\in\mathcal H\tand \tau\in T(D).$$

In the case that $C\in\mathcal C$, the map $\phi$ can be chosen to be a homomorphism.
\end{thm}
\begin{proof}
As in the proof of \ref{ExtTraceI-D}, since \CA s in ${\cal C}$ are semi-projective, we will only prove the first part
of the statement.
%Since any
%\CA\,
%C*-algebra
%in $\mathcal C$ is semi-projective, the second part of the statement follows directly from the first part of the statement. Thus, let us only show the first part of the statement. %Noting that $C$ is finitely generated, {\bf finitely
%generated----Do we state/prove this?---We do not really use it, do we?--L}
Without loss of generality, one may {also} assume that $\mathcal F\subseteq \mathcal H$.
By replacing $C$ by $M_m(C)$ for some integer $m\ge 1,$  we may find a finite subset
${\cal P}$  of projections in $C$ { {as in Remark \ref{Rextnn}.}}
%Let ${\cal P}$ be a set of projections in $C$
%such that, for  every projection $q\in C,$ there is a projection
%$p\in {\cal P}$ such that $p$ and $q$ are equivalent.
%Note that ${\cal P}$ is finite.
Without loss of generality, one may also assume that $\mathcal P\subseteq \mathcal H$.

Since the K-group of $C$ is finitely generated (as abelian groups), there is $M\in\mathbb N$ such that
$Mx=0$ for all $x\in \mathrm{Tor}(K_i(C)),$ $i=0,1.$
%$$Mp=0,\quad p\in K_*(C, \mathbb Z/n\mathbb Z), *=0, 1, n=1, 2,,, \ .$$
%Since $C$ is weakly semi-projective, there is $\delta'>0$ such that if $L: C\to A$ is a unital $\mathcal H$-$\delta'$-multiplicative { positive linear}  map for a unital C*-algebra $A$, there is a homomorphism $\phi: C\to A$ such that $$||\phi(h)-L(h)||<\sigma/{4}\tforal h\in\mathcal H.$$ One may also assume that $\delta'$ is small enough so that $[L]|_{K_0(C)}$ is well defined and $[L]|_{K_0(C)}=\phi_{*0}$.

Let $\mathcal H_{1, 1}\subset C_+^1\setminus \{0\}$ (in place of $\mathcal H_{1}$)
%$\mathcal G_1\subseteq C_{s.a.}$ (in place of ${{\mathcal H_2}}$) be finite subsets
and $\sigma_1>0$ (in place of ${{\delta}}$) be  a positive number required by
%the finite subset of
{{Theorem \ref{UniqAtoM}}} with respect to $C$ (in the place of ${{A}}$), $\min\{\sigma/{4}, \epsilon/{2}\}$ (in the place of $\epsilon$), $\mathcal H$ (in the place of $\mathcal F$) and $\Delta$.

%Let ${K}_1$ be the constant of Lemma \ref{ExtTraceM} with respect to $C$, $\mathcal H\cup\mathcal H_{1, 1}\cup\mathcal G_1$ (in the place of $\mathcal H$) and $\frac{1}{2}\min\{\Delta(\hat{h})/2, \sigma_1/4;\ h\in\mathcal H_{1, 1}\}$ (in the place of $\sigma$).

Let $\mathcal H_{1, 2}\subseteq C$ (in place of $\mathcal H_1$) {be a finite subset}, let $\sigma_2$ (in place of $\delta$)
be a positive  number,
%{${\cal P}\subset K_0(C)$ be a finite subset}
{and $K_1$ (in place of $K$) be an integer required by  Lemma \ref{APPextnn}} with respect to $\mathcal H\cup\mathcal H_{1, 1}$
%\cup\mathcal G_1$ (in the place of $\mathcal H$)
and ${\frac{1}{2}\min\{\sigma/16, \sigma_1/4, \min\{\Delta(\hat{h})/2: \ h\in\mathcal H_{1, 1}\}\}}$ (in the place of {$\sigma$}) and $\Delta$. Note we can choose ${\cal P}$ above ---
see \ref{Rextnn}.

Let $\mathcal H_{1, 3}$ (in place of $\mathcal H_1$), $\sigma_3>0$ (in place of $\delta$) and $K_2$
(in place of $K$) be  finite subset and constants required by  Lemma \ref{ExtTraceI-D}  with respect to $C$, $\mathcal H\cup\mathcal H_{1, 1}$
%\cup\mathcal G_1$
(in the place of $\mathcal H$), {$\min\{\sigma/16, \sigma_1/4, \min\{\Delta(\hat{h})/2:\ h\in\mathcal H_{1, 1}\}\}$} (in the place of $\sigma$), $\ep/4$ (in place of $\ep$), ${\cal H}$ (in place of
${\cal F}$)  and $\Delta$ (with the same ${\cal P}$ above).
%$$\min\{\Delta(\hat{h})/2, \sigma_1/4;\ h\in\mathcal H_{1, 1}\},\quad \forall h\in \mathcal H\cup\mathcal H_{1, 1}\cup\mathcal G_1$$

Put $\mathcal H_1=\mathcal H_{1, 1}\cup\mathcal H_{1, 2}\cup\mathcal H_{1, 3}{\cup {\cal P}}$, $\delta=\min\{{\sigma_1/2}, \sigma_2, {1/4}\}$ and $K=MK_1K_2$.
Let
$$
D={{A}}(F_1, F_2, \psi_0, \psi_1)=\{(f,a)\in C([0,1], F_2)\oplus F_1: f(0)=\psi_0(a)\andeqn f(1)=\psi_1(a)\}
$$
be any
%C*-algebra
\CA\, in $\mathcal C$, and let $\gamma: T(D)\to T(C)$ be a given continuous affine map satisfying $$\gamma(\tau)(h) > \Delta(\hat{h})\rforal h\in \mathcal H_{1} \rforal \tau\in T(D).$$
%Consider $M_s(D)$ with $s$ divided by $M$, and l
{L}et $\kappa: K_0(C)\to K_0(M_s(D))$ be any positive map {with $s[1_D]=\kappa([1_C])$} satisfying
$$
|\rho_C(x)(\gamma(\tau))-{(1/s)}\tau(\kappa(x))|<\delta\rforal \tau\in T(D)
$$
{and for all $x\in {\cal P}.$}
Write $C([0, 1], F_2)=I_1\oplus I_2\oplus \cdots \oplus I_{{k}}$ with $I_i=C([0, 1], M_{r_i})$, $i=1, ..., {{k}}$.
%Let $j: D\to J$ be defined by $j((f,a))=f$ for all $(f,a)\in D.$
Note that $\gamma$ induces a continuous affine map $\gamma_i: T(I_i)\to T(C)$ by
%$\gamma^J(\tau)=\gamma(\tau\circ J)$ for all $\tau\in T(J).$
%Therefore   $\gamma$  also induces a continuous affine map
$\gamma_i: T(I_i) \to T(C)$ defined by $\gamma_i(\tau)=\gamma(\tau\circ\pi_i)$ for each $1\leq i\leq {k},$   where $\pi_i$ is the restriction map $D\to I_i$
defined by $(f,a)\to f|_{[0,1]_j}.$  It is clear that
for any $1\leq i\leq {k}$,
one has that
\beq
%\begin{equation}
\label{intv-dense}
\gamma_i(\tau)(h)>\Delta(\hat{h}) \rforal h\in\mathcal H_{1, 3}{\andeqn} \rforal \tau\in T(I_i)
%\end{equation}
%and
%\begin{equation}
\andeqn\\\label{intv-cpt}
|\rho_C(x)(\gamma_i(\tau))-\tau((\pi_i)_{*0}\circ\kappa(x))|<\delta\leq \sigma_3\rforal  \tau\in T(M_s(I_i))
%\end{equation}
\eneq
{and for all $x\in {\cal P}$ and for any $1\le i\le k.$}
Also write $F_1=M_{{R_1}}\oplus\cdots\oplus M_{{R_l}}$ and denote by $\pi'_j: D\to M_{{R_j}}$ the corresponding evaluation of $D$.
%With a slightly abusing notation, one still use $\pi_i$ and $\pi'_j$ to denote the restrictions of $M_s(D)$ to the distinguished intervals and evaluations on the distinguished points.
Since
%$$
\beq\nonumber
&&\gamma(\tau)(h)>\Delta(\hat{h})\rforal h\in\mathcal { {H}}_{1, 2}\,\,\, {\rm and}\rforal \tau\in T(D),
\andeqn\\
%and
%$$
&&|\rho_C(x)(\gamma(\tau))-{(1/s)}\tau(\kappa(x))|<\delta  \rforal \tau\in T(D)
\eneq
{and for all  $x\in {\cal P}, $} one has that, for each $j,$
\beq\nonumber
%$$
\gamma\circ(\pi_j')^*(\mathrm{tr})(h)>\Delta(\hat{h})\rforal h\in\mathcal H_{1, 2}\andeqn\\
%and
%$$
|\rho_C(x)(\gamma\circ(\pi_j')^*(\mathrm{tr'}))-\mathrm{tr}([\pi_j']\circ\kappa(x))|<\delta\leq \sigma_2,
%\rforal \mathcal \tau\in T(M_s(D))
\eneq
where ${\rm tr}$ is the tracial state on $M_{sR_j}$ and ${\rm tr}'$ is the tracial
state on $M_{R_j},$ for all  $x\in K_0(C)$ and
where $\gamma\circ (\pi_j')^*({\rm tr})=\gamma({\rm tr}\circ \pi_j).$
%and for each $j$.
%Applying Lemma \ref{pert-T-K-D}, there is a tracial state $\tilde{\tau}\in T(C)$ such that $([\pi_j']_0\circ\kappa, \tilde{\tau})$ is a compatible pair, and
%\begin{equation}\label{pt-pre-0}
%|\tilde{\tau}(h)-\gamma\circ(\pi_j')^*(h)|< \frac{1}{2}\min\{\Delta(\hat{h})/2, \sigma_1/4;\ h\in\mathcal H_{1, 1}\}\rforal h\in \mathcal H\cup\mathcal H_{1, 1}\cup\mathcal G_1.
%\end{equation}

%Since $s$ is divided by $M$ (hence is divided by $M_1$),
It follows from Lemma \ref{APPextnn} that there is a homomorphism $\phi'_j: C\to M_{{R}_j}\otimes M_{sK_1K_2}$ such that
%\begin{equation}
\beq\label{C-D-nnn1}
&&{(\phi'_j)_{*0}=(\pi_j')_{*0}\circ {K_1K_2}\kappa}\andeqn\\
%\end{equation}
%and
%\begin{equation}\label{pt-pre-1}
%| \mathrm{tr}\circ \phi'_j(h)- \tilde{\tau}(h) | <\frac{1}{2}\min\{\Delta(\hat{h})/2, \sigma_1/4;\ h\in\mathcal H_{1, 1}\},\quad \forall h\in \mathcal H\cup\mathcal H_{1, 1}\cup\mathcal G_1.
%\end{equation}
%Together with \eqref{pt-pre-0}, one has
%\begin{equation}
\label{pt-pre}
&&|\mathrm{tr}\circ \phi'_j(h)- (\gamma \circ (\pi_j')^*)(\mathrm{tr}')(h) | <{\min\{\sigma/16, \sigma_1/4, \min\{\Delta(\hat{h})/2: \ h\in\mathcal H_{1, 1}\}\}}
%\end{equation}
\eneq
for all $h\in \mathcal H\cup\mathcal H_{1, 1},$
%\cup\mathcal G_1,$
where
${\rm tr}$ is the tracial
state on $M_{R_j}\otimes M_{sK}$ and
where ${\rm tr}'$ is the tracial state on $M_{R_j}.$
Denote by $$\phi'=\bigoplus_{j=1}^l \phi'_j: C\to F_1\otimes\mathrm{M}_{sK_1K_2}(\mathbb C).$$
Applying Lemma \ref{ExtTraceI-D} to \eqref{intv-dense} and \eqref{intv-cpt}, one has that, for any $1\leq i\leq {k}$,
there is an ${\cal H}$-$\ep/4$-multiplicative \morp\,
%homomorphism
$\phi_i: C\to I_i\otimes M_{{sK_1K_2}}$ such that ${(\phi_i){*_0}=(\pi_i)_{*0}\circ K_1K_2\kappa}$ and
\begin{equation}\label{int-pre}
|(1/sK_1K_2)\tau \circ \phi_i(h)- ((\gamma\circ(\pi_i)^*)(\tau))(h) |< {\min\{\sigma/16, \sigma/4, \min\{\Delta(\hat{h})/2:\ h\in\mathcal H_{1, 1}\}\}}
\end{equation}
for all $h\in \mathcal H\cup\mathcal H_{1, 1}\cup\mathcal G_1,$ where $\tau\in T({I_i}).$ {{Furthermore, as
{{in the}}  conclusion of Lemma \ref{ExtTraceI-D}, the restrictions of $\phi_i$ to both boundaries are  homomorphisms.}}

%T(M_{sK_1K_2}(D))$
%and $\tau'\in T(D).$

For each $1\leq i\leq {k}$, denote by $\pi_{i, 0}$ and $\pi_{i, 1}$ the evaluations of $I_i\otimes M_s$ at the point $0$ and $1$ respectively. Then one has
\begin{equation}\label{C-D-n1}
{\psi_{0,i}\circ \pi_e=\pi_{i,0}\circ \pi_i.}
\end{equation}
It follows that
\begin{eqnarray}\label{C-Dn-2}
(\psi_{0, i}\circ\phi')_{*0} & = & (\psi_{0,i})_{*0}\circ (\sum_{j=1}^l(\pi_j')_{*0})\circ K_1K_2\kappa\\
&=&(\psi_{0,i})_{*0}\circ (\pi_e)_{*0}\circ K_1K_2\kappa=(\pi_{i,0}\circ \pi_i)_{*0}\circ K_1K_2\kappa\\
&=&(\pi_{i, 0})_{*0}\circ (\phi_i)_{*0}.
\end{eqnarray}
Moreover, note that by \eqref{pt-pre},
\begin{equation}
\mathrm{tr}\circ (\psi_{0, i}\circ\phi')(h)\geq \Delta(\hat{h})/2\rforal h\in \mathcal H_{1, 1},
\end{equation}
and by \eqref{int-pre},
\begin{equation}
\mathrm{tr}\circ (\pi_{i, 0}\circ \phi_i)(h)\geq \Delta(\hat{h})/2\rforal h\in \mathcal H_{1, 1}.
\end{equation}
It also follows from \eqref{pt-pre} and \eqref{int-pre} that
\begin{equation}
|\mathrm{tr}\circ (\psi_{0, i}\circ\phi')(h) - \mathrm{tr}\circ (\pi_{i, 0}\circ \phi_i)(h) | < \sigma_1/2\rforal h\in \mathcal
H_{1,1}.
% G_1.
\end{equation}
%{\bf (More detail calculation will be added later).}
Consider {{amplifications}}
%$$
\beq\nonumber
\phi'_i:&=&\phi_i\otimes 1_{\mathrm{M}_M(\mathbb C)}: C\to I_i\otimes M_{{sK}}\andeqn\\
%$$
%and
%$$
\phi'':&=&\phi'\otimes 1_{\mathrm{M}_M(\mathbb C)}: C\to F_1 \otimes M_{{sK}}.
%$$
\eneq
Then, one has
$$[\psi_{0, i}\circ\phi'']=[(\pi_{i, 0})_{*0}\circ (\phi'_i)]\quad\textrm{in $KL(C, \mathrm{M}_{r_isK})$}.$$
Therefore, by { {Theorem \ref{UniqAtoM}}}, there is a unitary $u_{i, 0}\in M_{r_i} \otimes M_{{sK}}$ such that
$$\|\mathrm{Ad}u_{i, 0}\circ \pi_{i, 0}\circ \phi'_i(f) -  \psi_{0, i}\circ\phi''(f)\|<\min\{\sigma/{4}, \epsilon/{2}\}\rforal f\in\mathcal H.$$
Exactly the same argument shows that
there is a unitary $u_{i, 1}\in M_{r_i} \otimes M_{{sK}}$ such that
$$\|\mathrm{Ad}u_{i, 1}\circ {{\pi}}_{i, 1}\circ \phi_i{{'}}(f) -  \psi_{1, i}\circ\phi''(f)\|<\min\{\sigma/4, \epsilon/2\}\rforal f\in\mathcal H.$$

Choose {two}  path{s} of {unitaries}  ${\{u_{i, 0}(t):t\in [0,1/2]\}\subset} M_{r_i} \otimes M_{{sK}}$
such that $u_{i, 0}(0)=u_{i, 0}$ and {$u_{i, 0}(1/2)=1_{M_{r_i} \otimes M_{sK}},$} {and
$\{u_{i,1}(t):t\in [1/2, 1]\}\subset M_{r_i}\otimes M_{{sK}}$ such that
$u_{i,1}(1/2)=1_{M_{r_i}\otimes M_{sK}}$ and $u_{i,1}(1)=u_{i,1}$ }
{Put $u_i(t)=u_{i,0}(t)$ if $t\in [0,1/2)$ and $u_i(t)=u_{i,1}(t)$ if $t\in [1/2,1].$
Define ${\tilde \phi}_i: C\to I_i\otimes M_{sK}$ by
\begin{equation*}
\pi_t\circ {\tilde \phi_i}={\rm Ad}\, u_i(t)\circ \pi_t\circ \phi'_i,
\end{equation*}
where $\pi_t: I_i\otimes M_{sK}\to M_{r_i}\otimes M_{sK}$ is the point-evaluation at $t\in [0,1].$}

%Then redefine the map $\phi_i$ to be $\mathrm{Ad}u_{i, 0}(t)\circ \phi_i$, and o
{One has that, for each $i,$
%\begin{equation}
\beq\nonumber
&&\| \pi_{i, 0}\circ {\tilde \phi}_i(f) -  \psi_{0, i}\circ\phi''(f)\|<\min\{\sigma/4, \epsilon/2\}
%\rforal f\in\mathcal H
\andeqn\\
%\end{equation}
%and
%\begin{equation}
&&\| \pi_{i, 1}\circ {\tilde \phi_i}(f) -  \psi_{1, i}\circ\phi''(f)\|<\min\{\sigma/4, \epsilon/2\}\rforal f\in\mathcal H.
\eneq
%, i=1,..., k.
%\end{equation}
}
For each $1\leq i\leq {k}$, let $\epsilon_i<1/2$ be a positive number such that
%$$
\beq\nonumber
&&\|{\tilde \phi}_i(f)(t) - \psi_{0, i}\circ\phi''(f)\|<\min\{\sigma/{4}, \epsilon/{2}\} \rforal f\in\mathcal H \rforal t\in [0, \epsilon_i]\andeqn\\
%$$
%and
%$$
&&\|{\tilde \phi}_i(f)(t) -  \psi_{1, i}\circ\phi''(f)\|<\min\{\sigma/{4}, \epsilon/{2}\}\rforal f\in\mathcal H \rforal t\in [1-\epsilon_i, 1].
\eneq
%$$
Define ${\Phi}_i: C\to I_i\otimes M_{{sK}}$ to be
$${\Phi}_i(t)=\left\{
\begin{array}{ll}
\frac{(\epsilon_i-t)}{\epsilon_i}(\psi_{0, i}\circ\phi'') + \frac{t}{\epsilon_i}{{\tilde \phi}}_i(f)(\epsilon_i), & \textrm{if $t\in [0, \epsilon_i]$},\\
{\tilde \phi}_i(f)(t), & \textrm{if $t\in [\epsilon_i, 1-\epsilon_i]$ },\\
\frac{(t-1+\epsilon_i)}{\epsilon_i}(\psi_{1, i}\circ\phi'') + \frac{1-t}{\epsilon_i}{{\tilde\phi}}_i(f)(\epsilon_i), & \textrm{if $t\in [1- \epsilon_i, 1]$}.
\end{array}
 \right.$$
The map ${\Phi}_i$ is not necessarily a homomorphism, but it is $\mathcal H$-$\epsilon$-multiplicative; in particular, it is $\mathcal F$-$\epsilon$-multiplicative. Moreover, it satisfies the relations
%\begin{equation}
\beq\label{C-D-n5}
\pi_{i, 0}\circ {\Phi}_i(f) =  \psi_{0, i}\circ\phi''(f)
%\rforal f\in\mathcal H, i=1, ..., {k},
%\end{equation}
%and
\andeqn
%\begin{equation}\label{C-D-n6}
\pi_{i, 1}\circ {\Phi}_i(f) =  \psi_{1, i}\circ\phi''(f)\rforal f\in\mathcal H, i=1,..., {k}.
\eneq
%\end{equation}
{Define  $\Phi'(f): C\to C([0,1], F_2)\otimes M_{sK}$ by $\pi_{i,t}\circ \Phi'=\Phi_i,$ where
$\pi_{i,t}: C([0,1], F_2)\otimes M_{sK}\to M_{r_i}\otimes M_{sK}$ defined by the point evaluation
at $t\in [0,1]$ (on the $i$-th summand) and define $\Phi'': C\to F_1$ by $\Phi''(f)=\phi'(f)$ for all $f\in C.$  Define
$$\phi(f)=(\Phi'(f),\Phi''(f)).$$
It follows from \eqref{C-D-n5}
%and \eqref{C-D-n6}
that
$\phi$ is $\mathcal F$-$\epsilon$-multiplicative positive linear map from $C$ to $D\otimes M_{sK}$}.
%Since $C$ is semiprojective, by the choice of $\delta'$, there is a homomorphism $\phi: C\to D$ such that
%\beq\label{C-D-n8}
%\|{\Phi(h)-\phi(h)}\|<\sigma/2 \rforal h\in\mathcal H.
%\eneq
It follows from (\ref{C-D-nnn1}) that
\beq\label{C-D-n9}
[\pi_e\circ \phi(p)]=[\phi'(p)]=(\pi_e)_{*0}\circ K\kappa([p]) \rforal p\in {\cal P}.
\eneq
%Since ${\cal P}\subset {\cal H}$ and $\sigma<1/2,$ by (\ref{C-D-n9}),
%\beq\label{C-D-n10}
%[\pi_e\circ \phi(p)]=\pi_{*0}\circ K\kappa([p])\rforal p\in {\cal P}.
%\eneq
Since $(\pi_e)_{*0}: K_0(D)\to \Z^l$ is injective, one has
\begin{equation}\label{C-D-n11}
\phi_{*0}=K\kappa.
\end{equation}
It follows from \eqref{int-pre} and \eqref{pt-pre} that one calculates
that
$$
|(1/sK)\tau\circ \phi(h)-\gamma(\tau)(h)|<\sigma\rforal h\in {\cal H}
$$
and for all $\tau\in T(D).$
%Then a straightforward calculation shows that the homomorphism $\tilde{\phi}$ satisfies the lemma.
\end{proof}

\begin{lem}\label{TtoDelta}
Let $C\in {\cal C}.$ For any $\ep>0$ and any finite subset
${\cal H}\subset C_{s.a.},$ there exists a finite
subset of extremal traces ${\cal T}\subseteq T(C)$  and
a continuous affine map $\lambda: T(C)\to \triangledown,$
where $\triangledown$ is the convex hull of of ${\cal T}$ such that
\beq\label{TtoD-1}
|\lambda(\tau)(h)-\tau(h)|<\ep,\quad h\in\mathcal H,\ \tau\in T(C).
\eneq
\end{lem}

\begin{proof}
Without loss of generality, we may assume that
${\cal H}$ is in the unit ball of $C.$
Write $C=C(F_1, F_2, \psi_0, \psi_1),$ where
$F_1=M_{R_1}\oplus M_{R_2}\oplus\cdots \oplus M_{R_l}$ and
$F_2=M_{r_1}\oplus M_{r_2}\oplus \cdots \oplus M_{r_k}.$
Let $\pi_{e,i}: C\to M_{R_i}$ be the surjective \hm\, defined by the composition
of $\pi_e$ and the projection from $F_1$ onto $M_{R_i},$ and
$\pi_{I,j}: C\to C([0,1],M_{r_j})$ the restriction which may also be
viewed as the restriction of the projection from $C([0,1], F_2)$ to $C([0,1], M_{r_i}).$ Denote by $\pi_t\circ \pi_{I,j}$
the composition of $\pi_{I,j}$ and the point-evaluation at $t\in [0,1].$
There is $\dt>0$ such that, for any $h\in {\cal H},$
\beq\label{TtoD-2}
\|\pi_{I,j}(h)(t)-\pi_{I,j}(h)(t')\|<\ep/16\rforal h\in {\cal H}
\eneq
and $|t-t'|<\dt,$ $t, t'\in [0,1].$

Let $g_1,g_2,..., g_n$ be a partition of unity over interval $[\dt, 1-\dt]$ with respect to
an open cover with order 2
such that each ${\rm supp}(g_i)$ has diameter $<\dt$ and
$g_{s}g_{s'}\not=0$ implies that $|s-s'|\le 1.$
Let $t_s\in {\rm supp}(g_s){{\cap}} [\dt, 1-\dt]$ be a point.  We may assume that
$t_s<t_{s+1},$ We may further choose
$t_1=\dt$ and $t_n=1-\dt$ and assume that $g_1(\dt)=1$ and $g_n(1-\dt)=1$ by
choosing an appropriate open cover of order 2.

Extend $g_s$ to $[0,1]$ by
defining $g_s(t)=0$ if $t\in [0, \dt)\cup (1-\dt,1]$ for $s=2,3,...,n-1$ and
\beq\label{ToD-3}
\hspace{-0.2in}g_1(t)=g_1(\dt)(t/\dt)\,\,\,{\rm for}\,\,\, t\in [0,\dt)\andeqn g_n(t)=g_n(1-\dt)(1-t)/\dt
\,\,\,{\rm for}\,\,\, t\in (1-\dt, 1].
\eneq
Define $g_0=1-\sum_{s=1}^n g_s.$ {{Then $g_0(t)=0 $ for all $t\in[\dt,1-\dt]$.}}
In what follows,
we view $g_s$ as $g_s\cdot {\rm id}_C$ as in the center of $C.$
In particular, $g_0$ is identified
with $(g_0,1_{F_1})$ so that $g_0(0)=\psi_0(1_{F_1})$ and $g_0(1)=\psi_1(1_{F_1}).$
 Let $g_{s,j}=\pi_{I,j}(g_s),$ $s=1,2,...,n,$ $j=1,2,...,k.$
Let $p_i\in F_1$ be the projection corresponding to the summand
$M_{R_i}.$
Choose $d_i\in C([0,1], F_2)$ such that $d_i(t)=\psi_0(p_i)$  for
$t\in [0,\dt]$ and $d_i(t)=\psi_1(p_i)$ for $t\in [1-\dt,1]$ and $0\le d_i(t)\le 1$
for $t\in (\dt, 1-\dt).$ Note that $d_i\in C.$
Moreover,
%\beq\label{DtoT-n1}
$\sum_{i=1}^l g_0d_i=g_0.$
%\eneq
 Without loss of generality,
we may assume that $\{t_s: 1\le s\le n\}$ is a set of distinct points.
Denote by ${\rm tr}_i$ the  tracial state on $M_{R_i}$ and
${\rm tr}_j'$ the tracial state on $M_{r_j},$ $i=1,2,...,l$ and $j=1,2,...,k.$
Let
\beq\label{TtoD-4}
{\cal T}=\{{\rm tr}_i\circ \pi_{e,i}:1\le i\le l\}\cup \bigcup_{s=1}^n\{{\rm tr}_j'\circ \pi_{t_s}\circ \pi_{I,j}: 1\le j\le k\}.
\eneq
Let $\triangledown$ be the convex hull of ${\cal T}.$
Define $\lambda: T(C)\to \triangledown$ by
\beq\label{ToD-5}
\lambda(\tau)(f)=\sum_{j=1}^k \sum_{s=1}^n \tau(g_{s,j}){\rm tr}_j'\circ (\pi_{I,j}(f)(t_s))+
\sum_{i=1}^l \tau(g_0d_i){\rm tr}_i\circ \pi_{e,i}(f){{,}}
\eneq
{{where $g_{s,j}=g_s\cdot 1_{r_j} \in C_0((0,1), M_{r_j})\subset C$,}} for all $f\in C.$
It is clear that $\lambda$ is a continuous affine map.
Note that if $h\in C,$
\beq\label{ToD-6}
\lambda({\rm tr}_j\circ \pi_{e,j})(h)&=&\sum_{i=1}^l {\rm tr}_j\circ \pi_{e,j}(g_0d_i){\rm tr}_i\circ \pi_{e,i}(h)\\
&=&{\rm tr}_j\circ \pi_{e,j}(g_0d_j){\rm tr}_j\circ \pi_{e,j}(h)={\rm tr}_j\circ \pi_{e,j}(h).\\
\eneq
If  $\tau(f)={\rm tr}_j'\circ ( \pi_{I,j}(f)(t))$ with $t\in (\dt, 1-\dt),$ then if $h\in {\cal H},$
\beq\label{ToD-7}
\tau(h)&=& {\rm tr}_j'\circ (\pi_{I,j}(h)(t))=(\sum_{s=1}^n{\rm tr}_j'\circ (\pi_{I,j}(hg_s)(t)))\\
&\approx_{2\ep/16} &\sum_{s=1}^n {\rm tr}_j'\circ \pi_{I,j}(h(t_s)g_s(t))\\
&=&\sum_{s=1}^n g_{s,j}(t){\rm tr}_j'\circ \pi_{I,j}(h(t_s))=\sum_{s=1}^n \tau(g_{s,j}){\rm tr}_j'\circ \pi_{I,j}(h(t_s))\\
&=&\sum_{i=1}^k\sum_{s=1}^n \tau(g_{s,i}){\rm tr}_i'\circ \pi_{I,i}(h(t_s)) +\sum_{i=1}^l 0 \cdot {\rm tr}_i\circ \pi_{e,i}(h)=\lambda(\tau)(h){{,}}
\eneq
where we notice that $\tau(g_{s,i})=0$ if $i\not=j$. If  $\tau$ has the form $\tau(f)={\rm tr}_j'\circ (\pi_{I,j}(f)(t))$ for some fixed $t\in (0, \dt),$
then for $h\in {\cal H}$ with $h=(h_0, h_1),$ where
$h_0\in C([0,1], F_2)$ and $h_1\in F_1$ such that \nolinebreak
$\psi_0(h_1)=h_0(0)=h(0)$ and $\psi_1(h_1)=h_0(1)=h(1),$
\begin{eqnarray*}
\tau(h)&=&{\rm tr}_j'\circ (\pi_{I,j}(h(t)))\nonumber \\
&=&{\rm tr}_j'\circ (\pi_{I,j}(h(t)g_1(t)))+{\rm tr}_j'\circ (\pi_{I,j}(h(t)g_0(t)))\nonumber \\
&\approx_{\ep/8}&{\rm tr}_j'\circ (\pi_{I,j}(h(\dt)g_1(t))+{\rm tr}_j'\circ \pi_{I,j}(h(0)g_0(t)))\nonumber \\
&=& g_1(t) {\rm tr}_j'\circ \pi_{I,j}(h(\dt)) + g_0(t)
{\rm tr}_j'\circ \pi_{I,j}(\psi_0(h_1))\nonumber \\
&=&\tau(g_{1,j}){{{\rm tr}_j'}} \circ \pi_{I,j}(h(t_1))+g_0(t)\sum_{i=1}^l{\rm tr}_j'\circ \pi_{I,j}(\psi_0(h_1p_i))\nonumber\\
&=&\tau(g_{1,j}){\rm tr}_j'\circ \pi_{I,j}(h(t_1))+g_0(t)(\sum_{i=1}^l{\rm tr}_j'\circ \pi_{I,j}(\psi_0(p_i)){\rm tr}_i\circ\pi_{e,i}(h_1))\nonumber\\
&=&\tau(g_{1,j}){\rm tr}_j'\circ \pi_{I,j}(h(t_1))+g_0(t)\sum_{i=1}^l {\rm tr}_j'\circ (\pi_{I,j}(d_i(t))) {\rm tr}_i(\pi_{e,i}(h))\nonumber\\
&=& \tau(g_{1,j}){\rm tr}_j'\circ \pi_{I,j}(h(t_1))+\sum_{i=1}^l {\rm tr}_j'\circ (\pi_{I,j}(g_{{0}}(t)d_i(t))) {\rm tr}_i(\pi_{e,i}(h))\nonumber\\
&=&\tau(g_{1,j}){\rm tr}_j'\circ \pi_{I,j}(h(t_1))+\sum_{i=1}^l \tau(g_0d_i) {\rm tr}_i(\pi_{e,i}(h))\nonumber\\
&=& \lambda(\tau)(h).
\end{eqnarray*}
The same argument as above shows that,
if
$$\tau(f)={\rm tr}_j'\circ (\pi_{I,j}(f)(t)),\quad t\in (1-\dt, 1),$$
then
\begin{equation*}
\tau(h)\approx_{\ep/8} \lambda(\tau)(h)\rforal h\in {\cal H}.
\end{equation*}
It follows that
\begin{equation*}
|\tau(h)-\lambda(\tau)(h)|<\ep/8\tforal h\in {\cal H}
\end{equation*}
and for all extreme points of $\tau\in T(C).$
By the Choquet Theorem, for each $\tau\in T(C),$  there exist a Borel probability measure
$\mu_\tau$ on the extreme points $\partial_e{T(C)}$ of $T(C)$ such that
\begin{equation*}
\tau(f)=\int_{\partial_e{T(C)}}f(t)d\mu_\tau\rforal f\in \Aff(T(C)).
\end{equation*}
Therefore, for
each $h\in {\cal H},$
\begin{equation*}
\tau(h)=\int_{\partial_e(T(C)} h(t)d\mu_{\tau}\approx_{\ep/8} \int_{\partial_e{T(C)} }h(\lambda(t))d\mu_{\tau}\rforal
\tau\in T(C),
\end{equation*}
as desired.
\end{proof}

\begin{lem}\label{cut-trace}
Let $C$ be a unital stably finite C*-algebra, and let $A\in \mathcal B_1$ (or $\mathcal B_0$). Let $\alpha: T(A)\to T(C)$ be a continuous affine map.
\begin{enumerate}
\item\label{cut-trace-a} For any finite subset $\mathcal H\subseteq \Aff(T(C))$,  any $\sigma>0$, there is a C*-subalgebra $D\subseteq A$ and a continuous affine map $\gamma: T(D)\to T(C)$ such that $D\in\mathcal C$ (or $\mathcal C_0,$) and
$$| h(\gamma(\imath(\tau)))- h(\alpha(\tau))|<\sigma\rforal\tau\in T(A) \rforal h\in \mathcal H,$$
where $\imath: T(A)\ni \tau\to \frac{1}{\tau(p)}\tau|_D\in T(D)$, and $p=1_D$.

\item\label{cut-trace-b} If there are a finite subset $\mathcal H_1\subseteq C^+$ and $\sigma_1>0$ such that $$\alpha(\tau)(g)>\sigma_1\rforal g\in\mathcal H_1\rforal \tau\in T(A),$$ the affine map $\gamma$ can be chosen so that $$\gamma(\tau)(g)>\sigma_1\rforal g\in\mathcal H_1$$ for any $\tau\in T(D)$.

\item\label{cut-trace-c} If the positive cone of $K_0(C)$ is generated by a finite subset ${\cal P}$  of projections and there is an order-unit map $\kappa: K_0(C) \to K_0(A)$ which is compatible to $\alpha$, then, for any $\delta>0$, the C*-subalgebra $D$ and $\gamma$ can be chosen so that  there are {positive} homomorphisms $\kappa_0: K_0(C)\to K_0((1-p)A(1-p))$ and $\kappa_1: K_0(C)\to K_0(D)$ such that  {}$\kappa_1$ is strictly positive, $\kappa=\kappa_0+{\imath}\circ \kappa_1,$ {where $\imath: D\to A$ is the embedding,},  and
\begin{equation}\label{june2-nn1}
|\gamma(\tau)(p)-\tau(\kappa_1([p])|<\dt\rforal p\in {\cal P}{\tand\,\,\, \tau\in T(D)}.
\end{equation}

%\and $\kappa_1$ and $\gamma$ are $\delta$-compatible;

\item\label{cut-trace-d} Moreover, if $A\cong A\otimes U$ for some infinite dimensional UHF-algebra, for any given positive integer $K$, the C*-algebra $D$ can be chosen so that  $D=M_K(D_1)$
for some $D_1\in {\cal C}$ ($D_1\in {\cal C}_0$) and $\kappa_1=K\kappa_1',$
where $\kappa_1': K_0(C)\to K_0(D_1)$ is a strictly positive homomorphism. {Furthermore,
$\kappa_0$ can also be chosen to be strictly positive.}
\end{enumerate}

\end{lem}

\begin{proof}

Write $\mathcal H=\{h_1, h_2, ..., h_m\}.$  We may assume that $\|h_i\|\le 1,$ $i=1,2,...,m.$ Choose $f_1, f_2, ..., f_m\in A$ such that that $\tau(f_i)=h_i(\alpha(\tau))$ for all $\tau\in T(A)$ and $\|f_i\|\le 2,$ $i=1,2,...,m$
(see 9.2 of \cite{LinTAI}). Put ${\cal F} =\{1_A, f_1,f_2,...,f_m\}.$

Let  $\dt>0$ and let ${\cal G}_1$ (in place of ${\cal G}$) be a finite subset required by Lemma 9.4 of \cite{LinTAI}  for $A,$ $\sigma/16$ (in place of $\ep$) and  ${\cal F}$.
Let $\sigma_1=\min\{\sigma/16, \dt/16\}.$
We may assume that ${\cal G}_1\supset {\cal F}.$
Put ${\cal G}=\{g,\,gh: g,\, h\in {\cal G}_1\}.$
Since $A\in {\cal B}_1$ {{ (or ${\cal B}_0$)}}  there is a $D\in {\cal C}$ ({{or}}  ${\cal C}_0$),  and $h'\in (1-p)A(1-p)$ and $h''\in D$ with $p=1_D$ such that
\begin{equation}\label{june2-n1}
\|h-(h'+h'')\|<\sigma_1/16,\quad h\in {\cal G}_1\andeqn \tau(1-p)<\sigma_1/2,\quad \tau\in T(A).
\end{equation}
%and
%$$\tau(1-p)<\sigma_1/2,\quad \tau\in T(A).$$

Moreover, since $D$ is amenable, without loss of generality,
we may further assume that there is a unital \morp\,
$L: A\to D$ such that $L(h)=h''$ and $L$ is ${\cal G}$-$\sigma_1/2$-multiplicative.
By the choice of $\dt$ and ${\cal G},$ it follows from Lemma 9.4 of \cite{LinTAI} that, for each $\tau\in  T(D),$ there is $\gamma'(\tau)\in T(A)$ such that
\begin{equation}\label{june2-n2}
|\tau(L(h))-\gamma'(\tau)(h)|<\sigma/16\rforal h\in {\cal F}.
\end{equation}
Applying \ref{TtoDelta},  one obtains $t_1, t_2,...,t_n\in \partial_e{T(D)}$ and
a continuous affine map $\lambda: T(D)\to \Delta$
such that
\begin{equation}\label{june6-n1}
|\tau(f)-\lambda(\tau)(f)|<\sigma_1/16\rforal \tau\in T(D)
\end{equation}
and $f\in {\cal F},$ where $\Delta$ is the convex hull of $\{t_1, t_2,...,t_n\}.$
Define $\lambda_1:  \Delta\to T(A)$ by
\begin{equation}\label{june6-n2}
\lambda_1(t_i)=\gamma'(t_i),\,\,\, i=1,2,...,m.
\end{equation}
Define $\gamma=\af\circ \lambda_1\circ \lambda.$
 Then
\begin{eqnarray*}
h_j(\gamma(\pi(\tau)))&=&h_j(\alpha\circ \lambda_1\circ\lambda(\imath(\tau)))\\
&=& \lambda_1\circ\lambda(\pi(\tau))(f_j)\\
&\approx_{\sigma/16}& \lambda(\imath(\tau))(f_j'')\\
&\approx_{\sigma/16}&\imath(\tau)(f_j'')\\
&\approx_{\sigma/8}& \tau(g_j)= h_j(\alpha(\tau)),
\end{eqnarray*}
and this proves (\ref{cut-trace-a}).
Note that it follows from the construction that $\gamma(\tau)\in\alpha(T(A))$, and hence  (\ref{cut-trace-b}) also holds.  With \ref{Cuthm},  (\ref{cut-trace-c}) and (\ref{cut-trace-d}) follows straightforwardly {except
the "Furthermore" part. }

{To see that, we note that we may choose $D\subset A\otimes 1_U.$ Choose a projection
$e\in U$ such that
$$
0<t_0(e)<\dt_0<\dt-\max\{|\gamma(\tau)(p)-\tau(\kappa_1([p])|: p\in {\cal P} \tand\, \tau\in T(D)\},
$$
where $t_0$ is the unique tracial state of $U.$ We then replace $\kappa_1$ by  $\kappa_2: K_0(A)\to K_0(D_2),$
where $D_2=D\otimes (1-e)$ and $\kappa_2([p])=\kappa_1([p])\otimes [1-e].$
Define $\kappa_3([p])=\kappa_1([p])\otimes [e].$
Then let $\kappa_4: K_0(C)\to K_0((1-p+(1\otimes e))A(1-p+(1\otimes e))$ be defined by
$\kappa_4=\kappa_0+[\imath]\circ \kappa_3,$ where $\imath: D\otimes e\to A\otimes U\cong A$ is the embedding.
We then replace $\kappa_0$ by $\kappa_4.$ Note that, now, $\kappa_4$ is strictly positive.}
\end{proof}

\section{Maps from homogeneous \CA s to \CA s in ${\cal C}.$}

%The following is convent in this section and follows from Lemma 2.15 of \cite{Li-interval} directly.{\bf Checked?
%Could we also remove "easily"---I think I do have reference for this--L}
%
%\begin{thm}\label{uniCMn}
%Let $X$ be a connected compact metric space, and consider $C(X)$. Let $\mathcal F\subseteq C$ be a finite subset, and let $\epsilon>0$ be a constant. There is a finite subset $\mathcal H_1\subseteq C^+$ such that for any $\sigma_1>0$ there is a finite subset $\mathcal H_2\subseteq C$ and $\sigma_2>0$ such that for any  {unital }homomorphisms $\phi, \psi: C\to M_n$ for a matrix algebra $M_n$ satisfying
%\begin{enumerate}
%\item $\phi(h)>\sigma_1$ and $\psi(h)>\sigma_1$ for any $h\in\mathcal H_1$, and
%\item $|{\rm tr}\circ \phi(h)-{\rm tr}\circ \psi(h)|<\sigma_2$ for any $h\in \mathcal H_2$,
%\end{enumerate}
%then there is a unitary $u\in M_n$ such that
%$$\|\phi(f)-u^*\psi(f)u\|<\epsilon\quad \textrm{for any $f\in\mathcal F$}.$$
%\end{thm}

%{{\bf  Original statement has a point-evaluation $\Phi.$ When I tried to verify that part, I realized that existence of $\Phi$ (image in $D$) might be a problem. So I remove it.   Fix a point $x\in X.$ Let $I=C_0(X\setminus \{x\}).$  Let $\imath: I\to C$ be the embedding.  I can show, however, that $[\phi\circ \imath]=0$ as an element in $KK(C,D).$ May I assume that was what meant?}} {\color{green} Yes.}

\begin{lem}\label{ExtTraceH}
Let $X$ be a connected finite CW-complex and let $C=C(X)$. Let $\mathcal H\subseteq C$ be a finite subset, and let $\sigma>0$. There exists a finite subset $\mathcal H_{1, 1}\subseteq C^+$ satisfying the following: for any $\sigma_{1, 1}>0$, there is a finite subset $\mathcal H_{1, 2}\subseteq C^+$ satisfying the following: for any $\sigma_{1, 2}>0$, there is a positive integer $M$ such that for any $D\in\mathcal C$ with the dimension of any irreducible representation of $D$ at least $M$, for any continuous affine map $\gamma: T(D)\to T(C)$ satisfying $$\gamma(\tau)(h) > \sigma_{1, 1}\tforal h\in \mathcal H_{1, 1}\rforal \tau\in T(D),$$ and $$\gamma(\tau)(h) > \sigma_{1, 2} \tforal h\in \mathcal H_{1, 2}\rforal \tau\in T(D),$$ there is a homomorphism $\phi: C\to D$ such that
%$[\phi]_0=[\Phi]_0$ for a point-evaluation map $\Phi$, and
$$|\tau\circ\phi(h)-\gamma(\tau)(h)|<\sigma\rforal h\in\mathcal H.$$
Moreover, if $D\in {\cal C}_0,$ then there is a point-evaluation $\Psi: C\to D$ such that
$[\phi]=[\Psi].$
\end{lem}

\begin{proof}
Without loss of generality, one may assume that every element of $\mathcal H$ has norm at most $1$.

 Let $\eta>0$ such that for any $f\in\mathcal H$ and any $x, x'\in X$ with $d(x, x')<\eta$, one has
$$|f(x)-f(x')|<\sigma/4.$$

 Since $X$ compact, { {one can choose}} a finite subset $\mathcal H_{1, 1}\subseteq C^+$ such that { {for any open ball $O_{\eta/24} \subset X$,  of radius $\frac{\eta}{24}$, there is an nonzero element $h\in  \mathcal H_{1, 1}$ with ${\rm supp}(h) \subset O_{\eta/24}$. We assume that $\|h\|\leq 1 $ for all $h\in \mathcal H_{1, 1}$. Consequently,} }
   if there is $\sigma_{1, 1}>0$ such that $$\tau(h)>\sigma_{1, 1}\rforal h\in\mathcal H_{1, 1},$$ then $$\mu_{\tau}(O_{\eta/24})>\sigma_{1, 1}$$ for any open ball $O_{\eta/24}$ with radius $\eta/24,$ where $\mu_\tau$ is the probability measure
induced by $\tau.$

{ { Fixed $\sigma_{1,1} >0$. }}  Let $\delta$
%(in the place of $\delta$)
and $\mathcal G\subseteq C(X)$ (in the place of $\mathcal G$) be the constant and finite subset of Lemma 6.2 of \cite{LnTAMS12} with respect to $\sigma/2$ (in the place of $\epsilon$), $\mathcal H$ (in the place of $\mathcal F$), and $\sigma_{1, 1}/\eta$ (in the place of $\sigma$).

Let $\mathcal H_{1, 2}\subseteq C(X)$ (in the place of $\mathcal H_1$) be the finite subset of Theorem \ref{uniCMn}
with respect to $\delta$ (in the place of $\epsilon$) and $\mathcal G$ (in the place of $\mathcal F$).

Let $\sigma_{1, 2}>0$. Then let $\mathcal H_2\subseteq C(X)$ (in the place of $\mathcal H_2$) and {$\sigma_2$}
%(in the place of $\sigma_2$)
be the finite subset and positive constant of Theorem \ref{uniCMn} with respect to $\sigma_{1, 2}$ { {and $\mathcal H_{1, 2}$ }}(in the place of $\sigma_1$ { {and $\mathcal H_1$}}).

Let $M$ (in place of $N$) be the constant of Theorem 2.1 of \cite{Li-interval} with respect to $\mathcal H_2\cup \mathcal H_{1, 2}\cup\mathcal H_{1, 1}$ (in the place of $F$) and $\min\{\sigma/4, \sigma_2/4, \sigma_{1, 2}/2, \sigma_{1, 1}/2\}$ (in the place of $\epsilon$).

Let $D=D(F_1, F_2, \psi_0, \psi_1)$ be a \CA\,  in $\mathcal{C}$ with dimensions of the irreducible representations at least $M$. Let $\gamma: T(D)\to T(C)$ be a map satisfying the lemma. Write $C([0, 1], F_2)=I_1\oplus\cdots\oplus I_r$ with $I_i=C([0, 1], M_{r_i})$, $i=1, ..., k$. Then $\gamma$ induces a continuous map $\gamma_i: T(I_i)\to T(C)$ by $\gamma_i(\tau)=\gamma(\tau\circ\pi_i)$, where $\pi_i$ is the restriction map $D\to I_i$. It is then clear that
$$\gamma_i(\tau)(h)>\sigma_{1, 2}\rforal {{~~h\in \mathcal H_{1,2},~~}}\tau\in T(I_i).$$

Also write $F_1=M_{R_1}\oplus\cdots\oplus M_{R_l}$, and denote by $\pi'_j: D\to M_{l_j}$ the corresponding evaluation of $D$.
By Theorem 2.1 of \cite{Li-interval}, for each $1\leq i\leq k$, there is a homomorphism $\phi_i: C(X)\to I_i$ such that
\begin{equation}\label{app-tr-001}
|\tau\circ\phi_i(h)- \gamma_i(\tau)(h) | < \min\{\sigma/4, \sigma_2/4, \sigma_{1, 2}/2, \sigma_{1, 1}/2\}\rforal h\in\mathcal H_2\cup\mathcal H_{1, 2}\cup\mathcal H_{1, 1};
\end{equation}
and for any $j$, there is also a homomorphism $\phi'_j: C(X)\to M_{R_j}$ such that
\begin{equation}\label{app-tr-002}
|\mathrm{tr}\circ\phi'_j(h)- \gamma\circ(\pi_j')^*(\mathrm{tr})(h) | < \min\{\sigma/4, \sigma_2/4, \sigma_{1, 2}/2, \sigma_{1, 1}/2\}\rforal h\in\mathcal H_2\cup\mathcal H_{1, 2}\cup\mathcal H_{1, 1}.
\end{equation}
Denote by $\phi'=\bigoplus_{j}\phi'_j$ and  by $\pi_t: I_i\to M_{r_i}$ the point-evaluation at $t\in [0,1].$
%In particular, i
It follows that
\beq\nonumber
|\mathrm{tr}\circ(\psi_{0, i}\circ \phi')-\mathrm{tr}\circ(\pi_0\circ\phi_i)|\leq \sigma_2/2\rforal h\in\mathcal H_2\andeqn\\
\mathrm{tr}\circ(\psi_{0, i}\circ\phi')(h)\geq \sigma_{1, 2}/2\quad\mathrm{and}\quad \mathrm{tr}\circ(\pi_{0}\circ\phi_i)(h)\geq \sigma_{1, 2}/2\rforal h\in \mathcal H_{1, 2}.
\eneq
By Theorem \ref{uniCMn}, there is a unitary $u_{i, 0}\in M_{r_i}$ such that
$$|| \mathrm{Ad}u_{i, 0}\circ\pi_{ 0}\circ\phi_i(f)-\psi_{0, i}\circ\phi'(f)||<\delta \rforal f\in\mathcal G.$$
Exactly the same argument shows that
there is a unitary $u_{i, 1}\in M_{r_i}$ such that
$$\|\mathrm{Ad}u_{i, 1}\circ \pi_{1}\circ \phi_i(f) -  \psi_{1, i}\circ\phi'(f)\|<\dt\rforal f\in\mathcal H.$$

Choose two  paths of unitaries  $\{u_{i, 0}(t):t\in [0,1/2]\}\subset M_{r_i} $
such that $u_{i, 0}(0)=u_{i, 0}$ and $u_{i, 0}(1/2)=1_{M_{r_i}}$ and
$\{u_{i,1}(t):t\in [1/2, 1]\}\subset M_{r_i}$ such that
$u_{i,1}(1/2)=1_{M_{r_i}}$ and $u_{i,1}(1)=u_{i,1}$
Put $u_i(t)=u_{i,0}(t)$ if $t\in [0,1/2)$ and $u_i(t)=u_{i,1}(t)$ if $t\in [1/2,1].$
Define ${\tilde \phi}_i: C\to I_i$ by
\begin{equation*}
\pi_t\circ {\tilde \phi_i}={\rm Ad}\, u_i(t)\circ \pi_t\circ \phi_i.
\end{equation*}
%where $\pi_t: I_i\to M_{r_i}$ is the point-evaluation at $t\in [0,1].$
%Choose a path of unitaries $u_{i, 0}(t)\in M_{r_i}$ such that $u_{i, 0}(0)=u_{i, 0}$ and $u_{i, 0}(1/2)=1_{M_{r_i}}$, %and redefine the map $\phi_i$ to be $\mathrm{Ad}u_{i, 0}(t)\circ\phi_i$, and
%one has

%$$|| \pi_{i, 0}\circ{\tilde{\phi}}_i(f)-\psi_{0, i}\circ\phi'(f)||<\delta\rforal f\in\mathcal G.$$
%Repeating the same argument for $\pi_{i, 1}$ and all $i=1, ..., r$, one has that
Then
\begin{equation}
\| \pi_{0}\circ {\tilde \phi}_i(f) -  \psi_{0, i}\circ\phi'(f)\|<\delta
\andeqn
\| \pi_{1}\circ {\tilde \phi}_i(f) -  \psi_{1, i}\circ\phi'(f)\|<\delta
\end{equation}
 for all $f\in\mathcal G, i=1,..., k.$

Note that it also follows from \eqref{app-tr-001} and \eqref{app-tr-002} that
$$\mathrm{tr}\circ(\psi_{0, i}\circ\phi')(h)\geq \sigma_{1, 1}/2\quad\mathrm{and}\quad \mathrm{tr}\circ(\pi_{i, 0}\circ{\tilde \phi}_i)(h)\geq \sigma_{1, 1}/2\rforal h\in \mathcal H_{1, 1}/2.$$ Hence
$$\mu_{\tau\circ (\psi_{0, i}\circ\phi')}(O_{\eta/24})\geq\sigma_{1,1}\quad\mathrm{and}\quad \mu_{\tau\circ (\pi_{0}\circ\tilde \phi_i)}(O_{\eta/24})\geq\sigma_{1,1}.$$

Thus, by Lemma 6.2 of \cite{LnTAMS12}, for each $1\leq i\leq k$, there are two unital homomorphisms
$$\Phi_{0, i}, \Phi_{0, i}': C(X)\to C([0, 1], M_{r_i})$$ such that
$$\pi_0\circ\Phi_{0, i}=\psi_{0, i}\circ\phi',\quad \pi_0\circ\Phi_{0, i}'=\pi_{ 0}\circ\tilde\phi_i,$$
$$||\pi_t\circ\Phi_{0, i}(f)-\psi_{0, i}\circ\phi'(f)||<\sigma/2,\quad  ||\pi_t\circ\Phi_{0, i}'(f)-\pi_{0}\circ{\tilde \phi_i} ||<\sigma/2$$
for all $f\in\mathcal H$ and $t\in[0, 1]$, and there is a unitary $w_{i, 0}\in M_{r_i}$ (in the place of $u$) such that
$$\pi_1\circ\Phi_{0, i}=\mathrm{Ad}w_{i, 0}\circ \pi_1\circ \Phi_{0, i}'.$$

The same argument shows that, for each $1\le i\le k,$ there are two unital \hm s
$\Phi_{1,i}, \Phi_{1,i}': C(X)\to C([0,1], M_{r_i})$ such that
$$\pi_1\circ\Phi_{1, i}=\psi_{1, i}\circ\phi',\quad \pi_1\circ\Phi_{1, i}'=\pi_{1}\circ\tilde\phi_i,$$
$$||\pi_t\circ\Phi_{1, i}(f)-\psi_{1, i}\circ\phi'(f)||<\sigma/2,\quad  ||\pi_t\circ\Phi_{1, i}'(f)-\pi_{0}\circ {\tilde \phi}_i ||<\sigma/2$$
for all $f\in\mathcal H$ and $t\in[0, 1]$, and there is a unitary $w_{i, 1}\in M_{r_i}$ (in the place of $u$) such that
$$\pi_0\circ\Phi_{1, i}=\mathrm{Ad}w_{i, 1}\circ \pi_0\circ \Phi_{1, i}'.$$

Choose two continuous paths $\{w_{i,0}(t): t\in [0,1]\},\,\{w_{i,1}(t): t\in [0,1]\}$
in $M_{r_i}$ such that $w_{i,0}(0)=w_{i,0},\,\,\, w_{i,0}(1)=1_{M_{r_i}}$ and
$w_{i,1}(1)=1_{M_{r_i}}$ and $w_{i,1}(0)=w_{i,1}.$

For each $1\leq i\leq k$, by the continuity of $\gamma_i,$ there is  $1>\epsilon_i>0$ such that
$$|\gamma_i(\tau_x)(h)-\gamma_i(\tau_y)(h)|<\sigma/4 \rforal h\in\mathcal H,$$
provided that $|x-y|<\epsilon_i,$ where $\tau_x$ and $\tau_y$ are the extremal trace of $I_i$ concentrated on $x$ and $y$ respectively.

Define the map $\tilde{\tilde{\phi}}_i: C\to I_i$ by
\begin{displaymath}
\pi_t\circ{\tilde{ \tilde{\phi}}}_i=
\left\{
\begin{array}{ll}
\pi_{\frac{3t}{\epsilon_i}}\circ\Phi_{0, i}, & t\in[0, \epsilon_i/3),\\
\mathrm{Ad}(w_{i, 0}(\frac{3t}{\epsilon_i}-1))\circ \pi_1\circ\Phi_{0, i}', & t\in[\epsilon_i/3, 2\epsilon_i/3),\\
\pi_{3-\frac{3t}{\epsilon_i}}\circ\Phi_{0, i}', & t\in[2\epsilon_i/3, \epsilon_i),\\
\pi_{\frac{t-\epsilon_i}{1/2-\epsilon_i}}\circ{\tilde \phi}_i, & t\in[\epsilon_i, 1/2),\\
\pi_{{1 -2\ep_i/3-t\over{1-2\pi_i/3-1/2}}}\circ \Phi_{1,i}' , &   t\in [1/2, 1-2\ep_i/3],\\
\mathrm{Ad}(w_{i, 1}(\frac{(1-\ep_i/3)-t}{\epsilon_i/3}))\circ \pi_0\circ\Phi_{1, i}', & t\in[1-2\epsilon_i/3, 1-\epsilon_i/3],\\
\pi_{\frac{t-1+\ep_i/3}{\epsilon_i/3}}\circ\Phi_{1, i}, & t\in[1-\epsilon_i/3,1].
\end{array}
\right.
\end{displaymath}
Then,
\begin{equation}\label{match-i}
\pi_0\circ\tilde{\tilde{\phi}}_i=\psi_{0, i}\circ\phi' \quad\mathrm{and}\quad \pi_1\circ\tilde{\tilde{\phi}}_i=\psi_{1, i}\circ
\phi'.
\end{equation}
One can also estimate, by the choice of $\ep_i$ and the definition of ${\tilde{\tilde{\phi}}}_i,$ that
\begin{equation}\label{app-tr-003}
|\tau_t\circ \tilde{\tilde{\phi}}_i(h)-\gamma_i(\tau_t)(h)|<\sigma\rforal t\in[0, 1],
\end{equation} where $\tau_t$ is the extremal tracial state of $I_i$ concentrated on $t\in[0, 1].$

Define $\Phi: C(X)\to C([0,1], F_2)$ by
$\Phi(f)=\bigoplus_{i=1}^k {\tilde{\tilde{\phi}}}_i(f)$ for all $f\in C(X).$
Define $\phi: C(X)\to C([0,1],F_2)\oplus F_2$ by $(\Phi(f),\phi'(f)).$
By \eqref{match-i}, $\phi$ is a \hm\, from $C(X)$ to $D.$
By \eqref{app-tr-003} and \eqref{app-tr-002}, one has that
$$|\tau\circ\phi(g)-\gamma(\tau)(h)|<\sigma\rforal h\in\mathcal H$$
and for all $\tau\in T(D),$
as desired.

To see the last part of the lemma, one assumes that  $D\in {\cal C}_0.$
Consider $\pi_e\circ \phi: C\to F_1$ (where $D=D(F_1, F_2, \psi_0, \psi_1)$ as above).
 Since $\pi_e\circ \phi$ has finite dimensional range, it is a point-evaluation.
 We may write $\pi_e\circ \phi(f)=\sum_{i=1}^m f(x_i)p_i$ for all $f\in C(X),$ where
 $\{x_1, x_2,...,x_m\}\subset X$ and $\{p_1, p_2,..., p_m\}\subset F_1$ is a set of mutually orthogonal projections.
%Fix a point $x\in X\setminus \{x_1,x_2,...,x_m\}$ and let
Let $I=\{f\in C(X): f(x_1)=0\}$   and $\imath: I\to C(X)$ be the embedding.
 It follows that $[\pi_e\circ \phi\circ \imath]=0.$ By
 \ref{2Lg13}, $[\pi_e]$ is injective on each $K_i(A)$ and each $K_i(A, \Z/k\Z)$ ($k\ge 2,$ $i=0,1$). Hence $[\phi\circ \imath]=0.$  Choose $\Psi(f)=f(x)\cdot 1_A.$ Since $X$ is connected,
$[\phi]=[\Psi]$ (see the end of the Remark \ref{remark831KL}).
\end{proof}

\begin{cor}\label{ExtTraceH-D}
Let $X$ be a connected finite CW-complex, and put $C=C(X)$. Let $\Delta: C_+^{q,  {\bf 1}}\setminus\{0\}\to (0, 1)$ be an order preserving map. Let $\mathcal H\subseteq C$ be a finite subset and let $\sigma>0$. Then there exists a finite subset $\mathcal H_1\subseteq C_+^{\bf 1}\setminus \{0\}$ and a positive integer $M$ such that  for any $D\in\mathcal C$ with the dimension of any irreducible representation of $D$ at least $M,$ for any continuous affine map $\gamma: T(D)\to T(C)$ satisfying $$\gamma(\tau)(h) > \Delta(\hat{h})\tforal h\in \mathcal H_{1}\tand \tforal \tau\in T(D),$$ there is a homomorphism $\phi: C\to D$ such that
%$[\phi]_0=[\Phi]_0$ for a point-evaluation map $\Phi$, and
$$|\tau\circ\phi(h)-\gamma(\tau)(h)|<\sigma\tforal h\in\mathcal H.$$
\end{cor}
\begin{proof}
Let $\mathcal H_{1, 1}$ be the subset of Lemma \ref{ExtTraceH} with respect to $\mathcal H$ and $\sigma$. Then put
$$\sigma_{1, 1}=\min\{\Delta(\hat{h}):\ h\in\mathcal H_{1, 1}\}.$$
Let $\mathcal H_{1, 2}$ be the finite subset of Lemma \ref{ExtTraceH} with respect to $\sigma_{1, 1}$, and then put
$$\sigma_{1, 2}=\min\{\Delta(\hat{h}):\ h\in\mathcal H_{1, 2}\}.$$
Let $M$ be the positive integer of Lemma \ref{ExtTraceH} with respect to $\sigma_{1, 2}$. Then it follows from Lemma \ref{ExtTraceH} that the finite subset $$\mathcal H_1:=\mathcal H_{1, 1}\cup\mathcal H_{1, 2}$$ and the positive integer $M$ satisfy the statement of the corollary.
\end{proof}

\begin{thm}\label{istTr}
Let  $X$ be a  connected  finite CW complex, and let $A\in {\cal B}_0$ be a unital separable simple C*-algebra.
Suppose that $\gamma: T(A)\to T_{\mathrm f}(C(X))$ is a continuous affine map.
Then, for any $\sigma>0,$ any finite subset ${\cal H}\subset C(X)_{s.a.},$ there exists a
unital \hm\, $h: C(X)\to A$ such that
\beq\label{istTr-1}
[h]=[\Psi]\andeqn
|\tau\circ h(f)-\gamma(\tau)(f)|<\sigma\tforal f\in {\cal H},
\eneq
where $[\Psi]$ is a point-evaluation.
\end{thm}

\begin{proof}
Without loss of generality, one may assume that every element of $\mathcal H$ has norm at most one. Let $\mathcal H_{1, 1}$  be the finite subset of Lemma \ref{ExtTraceH} with respect to $\mathcal H$ (in place of $\mathcal H$), $\sigma/4$ (in place of $\sigma$), and $C$ (in place of $C$). Since $\gamma(T(A))\subseteq T_{\mathrm f}(C(X))$, there is $\sigma_{1, 1}>0$ such that
$$\gamma(\tau)(h)>\sigma_{1, 1}\rforal h\in\mathcal H_{1, 1}\rforal \tau\in T(A).$$

Let $\mathcal H_{1, 2}\subseteq C_+$
%(in the place of $\mathcal H_{1, 2}$)
 be the finite subset of Lemma \ref{ExtTraceH} with respect to $\sigma_{1, 1}$. Again, since $\gamma(T(A))\subseteq T_{\mathrm{f}}(C(X))$, there is $\sigma_{1, 2}>0$ such that
$$\gamma(\tau)(h)>\sigma_{1, 2}\rforal h\in\mathcal H_{1, 2} \rforal \tau\in T(A).$$

Let $M$
% (in place of $M$)
be the constant of Lemma \ref{ExtTraceH} with respect to $\sigma_{1, 2}.$
%(in place of $\sigma_{1, 2}$).
%Since $A_1\in \mathcal B_0$, one has
{ {Note}} that $A\in \mathcal B_0$. By (\ref{cut-trace-a}) and (\ref{cut-trace-b}) of Lemma \ref{cut-trace}, there is a C*-subalgebra $D\in A$ with $D\in\mathcal C_0$, a continuous affine map $\gamma': T(D)\to T(C)$ such that
\begin{equation}\label{istTr-eq1}
|\gamma'(\frac{1}{\tau(p)}\tau|_D)(f)-\gamma(\tau)(f)|<\sigma/4\rforal \tau\in T(A) \rforal f\in\mathcal H,
\end{equation}
where $p=1_D$, $\tau(1-p)<\sigma/(4+\sigma)$,
\beq\label{istTr-eq2}
&&\gamma'(\tau)(h)>\sigma_{1, 1}\rforal \tau\in T(D)\rforal h\in\mathcal H_{1, 1}\andeqn\\
%\end{equation}
%and
%\begin{equation}
\label{istTr-eq3}
&&\gamma'(\tau)(h)>\sigma_{1, 2}\rforal \tau\in T(D)\rforal h\in\mathcal H_{1, 2}.
\eneq
%\end{equation}

Moreover, since $A$ is simple, one may assume that the dimension of any irreducible representation of $D$ is at least $M$ (see \ref{Affon1}). Thus, by \eqref{istTr-eq2} and  \eqref{istTr-eq3}, one applies Lemma \ref{ExtTraceH} to $D$, $C$, and $\gamma'$ (in the place of $\gamma$) to obtain a homomorphism $\phi: C\to D$ such that
%$[\phi]_0=[\Phi]_0$ for a point evaluation map $\Phi$ and
\begin{equation}\label{istTr-eq4}
|\tau\circ\phi(f)-\gamma'(\tau)(f)|<\sigma/4\rforal f\in\mathcal H\rforal \tau\in T(D).
\end{equation}

Moreover, we may assume that $[\phi]=[\Phi_0]$ for some
point evaluation $\Phi_0: C\to D,$ since we assume that $D\in {\cal C}_0.$ Pick a point $x\in X$, and define $h: C\to A$ by $$f \mapsto f(x)(1-p)\oplus\phi(f)\rforal f\in C.$$
%It is clear that $[h]_0$ is a point evaluation map.
For any $f\in\mathcal H$, one has
\begin{eqnarray*}
|\tau\circ h(f)-\gamma(\tau)(f)|&\leq & |\tau\circ \phi(f)-\gamma(\tau)(f)|+\sigma/4\\
&<&|\tau\circ \phi(f)-\gamma'(\frac{1}{\tau(p)}\tau|_D)(f)|+\sigma/2\\
&<&|\tau\circ \phi(f)- \frac{1}{\tau(p)}\tau\circ\phi(f)|+3\sigma/4<\sigma.
%\\
%&<&\sigma.
\end{eqnarray*}
Define $\Phi: C\to A$ by $\Phi(f)=f(x)(1-p)\oplus \Phi_0(f)$ for all $f\in C.$ Then $[h]=[\Phi].$
%To show that there is a point-evaluation $\Phi$ such that $[\Phi]=[h],$
%we consider $\pi_e\circ \phi: C\to F_1,$ where $D=D(F_1, F_2, \psi_0, \psi_1).$  By
% \ref{2Lg13}, $[\pi_e]$ is and embedding.  Let $I=\{f\in C(X): f(x)=0\}$ and $\imath: I\to C(X)$
%be the embedding.
%Since $\pi_e\circ \phi$ has finite dimensional range, it is a point-evaluation.
%It follows that $[\pi_e\circ \phi\circ \imath]=0.$ Hence $[\phi\circ \imath]=0.$
%Therefore $[h\circ \imath]=0.$ Choose $\Phi(f)=f(x)\cdot 1_A.$ Since $X$ is connected,
%$[h]=[\Phi].$ Thus the homomorphism $h$ satisfies the statement.
\end{proof}

\begin{cor}
The statement of Theorem \ref{istTr} holds for $C=PM_m(C(X))P$, where $X$ is a connected finite CW-complex and $P$ is a projection in $M_m(C(X))$.
\end{cor}

%{\bf It seems that in the following lemma, one needs to assume that the trace map $\gamma$ should be compatible to some $K_0$-maps, but not sure if this is automatically. Moreover, since $C$ has stable relation, the map $\phi$ can be chosen to be a homomorphism.}

\section{KK-attainability
%%%attainablity
 of the building blocks}\label{KK-BBlock}

% Recall that
%(\cite{LinTAI})
\begin{df}(9.1 of \cite{LinTAI})
Let ${\cal D}$ be a class of unital \CA s.
A \CA\, $C$ is said to be KK-attainable with respective  to ${\cal D}$ if for any $A\in {\cal D}$ and any $\alpha\in KK(C, A)^{++}$, there exists a sequence of completely positive linear maps  $L_n: C\to A$ such that
\beq
&&\lim_{n\to\infty} \|L_n(ab)-L_n(a)L_n(b)\|=0\rforal a,\, b\in C\andeqn\\
&&{[}L_n{]}=\af.
\eneq

In what follows, we will use ${\cal B}_{u0}$ for the class of those \CA s of the form
$A\otimes U,$ where $A\in {\cal B}_0$ and  $U$ is any UHF-algebras of infinite type.

\end{df}

\begin{thm}{\rm (Theorem 5.9 of \cite{LinTAF2})}\label{kkmaps}
Let $A$ be a separable C*-algebra satisfying UCT and let $B$ be a amenable separable C*-algebra. Assume that $A$ is the closure of an increasing sequence $\{A_n\}$ of residually finite dimensional \SCA s. Then for any $\alpha\in KL(A, B)$, there exist two sequences of completely positive contractions $\phi_n^{(i)}: A\rightarrow{B}\otimes\mathcal K\ (i=1, 2)$ satisfying the following:
\begin{enumerate}
\item $\|{\phi_n^{(i)}(ab)-\phi_n^{(i)}(a)\phi_n^{(i)}(b)}\|\rightarrow0$
      as $n\rightarrow\infty$,
\item for any $n$, the images of $\phi_n^{(2)}$ are contained in a finite
      dimensional sub-C*-algebra of $B\otimes\mathcal K$ and for any finite subset
      $\mathcal P\subset \underline{K}(A)$, $[\phi_n^{(i)}]|_{\mathcal P}$ are well defined for
      sufficiently large $n$,
\item for each finite subset $\mathcal P \subset\underline{K}(A)$, there exists $m>0$
      such that $$[\phi_n^{(1)}]|_\mathcal{P}=\alpha+[\phi_n^{(2)}]|_{\mathcal P} \tforal n\ge m,$$
     % for all $n>m$,
\item for each $n$, we may assume that $\phi_n^{(2)}$ is a homomorphism
      on $A_n$.
\end{enumerate}
\end{thm}

\begin{lem}\label{Nextension}
Let $C={{A}} (F_1, F_2, \phi_0, \phi_1)\in {\cal C}$
%\{(f,a)\in C([0,1], F_2)\oplus F_1: f(0)=\phi_0(a)\andeqn f(1)=\phi_1(a),\, \, a\in F_1\}\in {\cal C}$
and let $A\in {\cal C}$ be another \CA.
Let  $\kappa: K_0(C)\to K_0(A)$ be an order  preserving \hm\,
such that, for any nonzero element $p\in K_0(C)_+,$ there exists an integer $N\ge 1$ such that
$N\kappa(p)>[1_A].$
 Then there is $\sigma>0$ satisfying the following:
For any $\tau\in T(A),$ there exists $t\in T(C)$ such that
\beq\label{Next-1}
t(h)\ge \sigma\int_{[0,1]} {\rm T}(\imath(h(t)))d\mu(t)\tforal h\in C_+\andeqn\\
{\tau(\kappa(p))\over{\tau(\kappa([1_C])}}=t(p)\tforal p\in K_0(C)_+,
\eneq
where $\imath: C\to C([0,1], F_2)\oplus F_1$ is the natural embedding and
${\rm T}(b)=\sum_{i=1}^k{\rm tr}_i(b)$ for all $b\in F_2$ and where ${\rm tr}_i$ is
the normalized tracial state on the $i$-th simple summand of $F_2,$ and
$\mu$ is the Lebesgue measure on $[0,1].$
\end{lem}

\begin{proof}
Write $A={{A(F'_1,F'_2, \phi'_0,\phi'_1)=}}\{(g,c)\in C([0,1], F'_2)\oplus F'_1: g(0)=\phi_0{{'}}(c)\andeqn g(1)=\phi_1{{'}}(c),\, \, c\in F'_1\}$.
Suppose that $F_1'$ has {{$l'$ }} simple summands so that
$K_0(F_1')=\Z^{l'}.$
Let $\pi_i: K_0(F_1')\to \Z$ be the $i$-th projection.
View $\pi_i\circ \kappa$ as a positive \hm\, from $K_0(C)\to \R.$
Since $$N\kappa(x)>[1_A]\rforal x\in K_0(C)_+\setminus\{0\},$$
one has that
\beq\label{18-3n}
(\pi_i\circ\kappa)(x)>0\rforal x\in K_0(C)_+\setminus\{0\}.
\eneq
It follows from \cite{GH-RankFunction} that there are positive homomorphisms $L_i: K_0(C([0,1], F_2)\oplus F_1)
\to \R$ such that $L_i\circ \imath_{*0}=\pi_i\circ \kappa$ and
$
L_i(e_j)=\af_{i,j}{{>}}0,\,\,\,j=1,2,...k+l,\,\,\, i=1,2,...,l',
$
%Note that we do need $ \af_{i,j}>0$ to define $\sigma$ not just $\geq 0$---Gong
where $e_j=\psi_j\circ \imath([1_C])$ and $\psi_j$ is the projection from
$C([0,1], F_2)\oplus F_1$ to the $j$-th summand of $C([0,1],F_2),$ $1\le j\le k,$ and
$\psi_j$ is the projection from $C([0,1], F_2)\oplus F_1$ to
$(j-k)$-th simple summand of $F_1,$ $k+1\le j\le k+l.$  Moreover, by  (\ref{18-3n}), for each $i,$
%\beq\label{18-3-n2}
$\sum_j\af_{i,j}>0.$
%\eneq
Note that
\beq\label{Nextension-n1}
\sum_{j=1}^{k+l}{{\rm rank}\psi_j(p)\over{{\rm rank}\psi_j(1_C)}}\af_{i,j}=\sum_{j=1}^{k+l}{{\rm rank}\psi_j(p)\over{{\rm rank}\psi_j(1_C)}}L_i(e_j)=L_i\circ \imath_{*0}([p])=\pi_i\circ \kappa([p])
\eneq
for any projection $p\in M_m(C),$ where $m\ge 1$ is an integer.
Choose
$$
\sigma={ {\frac{\min\{\af_{i,j}: 1\le j\le k+l,\,\,\, 1\le i\le l'\}}{k\sum_{i,j}\af_{i,j}}}}.
$$
For any $\tau\in T(A),$ by 2.8 of \cite{Lnbirr}, there is $\tau'\in T(F_1')$
such that
$$
\tau'\circ \pi_e(x)=\tau(x)\rforal x\in K_0(A).
$$
Write
$
\tau'=\sum_{i=1}^{l'} \lambda_{i, \tau} {\rm tr}_i,
$
where $0\le \lambda_i\le 1$, $\sum_{i=1}^{l'}\lambda_{i, \tau}=1$ and each ${\rm tr}_i$ is the tracial state of the $i$-th simple summand of $F_1'.$

For each $i,$ define
$$
t_i(s)=({1\over{\sum_{j=1}^{k+l}\af_{i,j}}})(\sum_{j=1}^{k} \af_{i,j}\int_{[0,1]}{ {\rm tr}'_j(s)\over{{\rm tr}'_j(1_C)}} d\mu(t)+\sum_{j=k+1}^{l} \af_{i,j}{{\rm tr}'_j(s)\over{{\rm tr}'_j(1_C)}})
$$
for all $s\in C,$ where each ${\rm tr}'_j$ is the tracial state on the $j$-th simple summand
of $F_2$ for $1\le j\leq k$, and each ${\rm tr}'_j$ is the tracial state on the $j$-th simple summand
of $F_1$ for $k+1\le j\le l$.
In general, if $\tau\in T(A),$
define
$
t_{\tau}=\sum_{i=1}^{l'} \lambda_{i, \tau}t_i.
$
It is straightforward to verify that
$$
t_{\tau}(h)\ge \sigma\int_{[0,1]} {\rm T}(\imath(h))d\mu(t)\rforal h\in C_+.
$$
Moreover, for each $i,$ by (\ref{Nextension-n1}),
\begin{eqnarray*}
t_i(p)&=& ({1\over{\sum_{j}\af_{i,j}}})\left(\sum_{j=1}^{k} \af_{i,j}\int_{[0,1]}{ {\rm tr}'_j(p)\over{{\rm tr}'_j(1_C)}} d\mu(t)+\sum_{j=k+1}^{l} \af_{i,j}{{\rm tr}'_j(p)\over{{\rm tr}'_j(1_C)}}\right)\\
&=& ({1\over{\sum_{j}\af_{i,j}}})\left(\sum_{j=1}^{k} \af_{i,j}{{\rm rank}(\psi_j(p))\over{{\rm rank}(\psi_j(1_C))}}+\sum_{j=k+1}^{l} \af_{i,j}{{\rm rank}\psi_j(p)\over{{\rm rank}\psi_j(1_C))}}\right)\\
&=& ({1\over{\sum_{j}\af_{i,j}}})\left(\sum_{j=1}^{k}L_i(e_i){{\rm rank}\psi_j(p)\over{{\rm rank}\psi_j(1_C)}}
+\sum_{j=k+1}^{l}L_i(e_i){{\rm rank}\psi_j(p)\over{{\rm rank}\psi_j(1_C)}}\right)\\
&=&({1\over{\sum_{j}\af_{i,j}}})(\pi_i\circ \kappa(p))=\pi_i\circ \kappa(p)/\pi_i\circ \kappa([1_C])
\end{eqnarray*}
for all projection $p\in M_m(C).$
This implies that
$$
t_{\tau}(x)=\tau(\kappa(x))/\tau(\kappa([1_C]))\rforal x\in K_0(C)\andeqn
\rforal
\tau\in T(A).
$$
\end{proof}

\begin{prop}\label{pl-lifting}
Let $S\in\mathcal C$ and $N\ge 1.$  There exists an integer $K\ge 1$ satisfying the following:
 For any positive homomorphism  $\kappa: {K_0}(S)\to {K_0}(A)$
 which
 %{\red{preserves}} the order  and
 satisfies { {$\kappa([1_S])\leq [1_A]$ and}} $N\kappa([p])>[1_A]$ for any  $p\in K_0(S)_+\setminus \{0\}, $ where $A\in {\cal C}$,
there exists a  homomorphism $\phi: S\to M_K(A)$ such that $\phi_{*0}=K\kappa.$ {{ If we further assume $\kappa([1_S])=[1_A]$, then $\phi$ can be chosen to be unital.}}

\end{prop}

\begin{proof}
Write
$$
S={ {A(F_1,F_2, \phi_0, \phi_1)=}}\{(f,g): (f,g)\in C([0,1], F_2)\oplus F_1: f(0)=\phi_0(g),\,\,\, f(1)=\phi_1(g){{\}}}.
$$
Denote by $\imath: S\to C([0,1], F_2)\oplus F_1$ the embedding.
To simplify the notation, without loss of generality, by replacing $S$ by $M_r(S)$ for some integer $r\ge 1,$ we may assume
that projections in $S$ generates $K_0(S).$   {{Note that, by the assumption, $\kappa([p])$ associated
with a full projection of $M_N(A).$}}
{{W}}ithout loss of generality, {{by applying \ref{cut-full-pj},}}
we may assume that
that $[1_A]=\kappa([1_S]).$
% cutting down of $A$ by a projection in $A$ is still an algebra in ${\cal C}$.

Let $\sigma>0$ be given by \ref{Nextension} (associated with integer $N$).
%Let $\delta$ be the positive number of Lemma \ref{lbd-intv} with respect to $S$ and $N$.
Define { {$\Delta: S_+^{q,1}\setminus \{0\} \to (0,1)$ by}}
$$\Delta({{\hat h}})=\sigma\int_{[0, 1]}\mathrm{T}(h(t))d\mu(t)$$
for all $h\in S_{{+}},$ where
${\rm T}(c)=\sum_{i=1}^k {\rm tr}_i(c)$ for all $c\in C([0,1], F_2),$ ${\rm tr}_i$
is the normalized tracial state on the $i$-th simple summand of $F_2$ (so we assume
that $F_2$ has {{$k$ }} simple summands), and where
$\mu$ is the Lebesgue measure on $[0,1].$

Let $\mathcal H_1$, $\delta>0$, ${\cal P}$ and $K$  be the finite subset and constants of {{Theorem}} \ref{ExtTraceC-D} with respect to $S$, $\Delta$, an arbitrarily chosen $\mathcal H$ and an arbitrarily chosen $1>\sigma_1>0$
(in place of $\sigma$).

{We may assume that  ${\cal P}\subset S$}
%be
is a finite subset of projections such that
every projection $q\in S$ is equivalent to one of projections in ${\cal P}.$
%To simplify the notation, without loss of generality, we may assume, by the assumption, that
{{Let ${\cal Q}\subset M_{r'}(A)$ (for some $r'\ge 1$ be a finite subset of projections such
that $\kappa({\cal P})$ can be represented by projections in ${\cal Q}.$}}
It follows from \ref{TtoDelta} that there is a finite subset ${\cal T}$ of extreme points of $T({A})$
and there exists a continuous affine map $\gamma': T({A})\to \triangledown$ such that
\beq\label{june7-1}
|\gamma'(\tau)(p)-\tau(p)|<\dt/2\rforal p\in  {{\cal Q} },
\eneq
where $\triangledown$ is the convex hull of ${\cal T}.$

Note that  any \CA\
%C*-algebra
in the class $\mathcal C$ is of type I, it is amenable and in particular it is exact. Therefore, by \cite{Blatrace} and \cite{Haagtrace}, for each $s\in \triangledown,$
 there is tracial state $t_s\in T(S)$ such that
\beq\label{june7-2}
r_S(t_s)(x)=r_A(s)(\kappa(x))\rforal x\in K_0(S),
\eneq
where $r_S: T(S)\to S_{[1_S]}(K_0(S))$ and
$r_A: T(A)\to S_{[1_A]}(K_0(A))$ are the induced maps from
the tracial state spaces to the state space of $K_0,$ respectively.
It follows from \ref{Nextension} that we may choose
$t_s$ such that
\beq\label{june7-3}
t_s(h)\ge \Delta(h)\tforal h\in S_+.
\eneq
For each ${{s}}\in {\cal T},$ define
$\lambda(s)=t_s$ which satisfies (\ref{june7-2}) and (\ref{june7-3}). This extends to a continuous affine map
$\lambda: \triangledown\to T({{S}}).$
Put $\gamma=\lambda\circ \gamma'.$ Then, for any $\tau\in T(S),$
\beq\label{june7-4}
\gamma(\tau)(h)\ge \Delta(h) \rforal h\in {\cal H}_1
\andeqn
\eneq
\beq\label{june7-5}
\hspace{-0.5in}|\gamma(\tau)(q)-\tau(\kappa([q]))|&=&|\lambda(\gamma'(\tau))(q)-\gamma'(\tau)(\kappa([q]))|
+|\gamma'(\tau)(\kappa([q]))-\tau(\kappa([q]))|\\
&=&|\gamma'(\kappa([q]))-\kappa([q])|<\dt/2
\eneq
for all projections $q\in S.$
One then applies \ref{ExtTraceC-D} to obtain a unital \hm\, $\phi: S\to M_K(A)$
such that $[\phi]=K\kappa.$
\end{proof}

\begin{lem}\label{liftingpl-M}
 Let $C\in {\cal C}.$ Then there is $M>0$ satisfying the following: Let $A_1\in\mathcal B_1$  and let
$A=A_1\otimes U$ for some infinite dimensional UHF-algebra
% $U$ of infinite type
and let $\kappa:({K_0}(C), {K_0}^+(C))
%, [1_C])
\to({K_0}(A), {K_0}^+(A))$
%, [1_A])$
 be a strictly positive homomorphism with multiplicity $M$. Then there exists a  homomorphism $\phi: C\to
 {{M_m(A)}}$ {{(for some integer $m\ge 1$)}} such that $\phi_{*0}=\kappa$ and $\phi_{*1}=0.$
 %in particular, $[\phi(1_C)]=\kappa([1_C]).$
 %{\blue {If we replace $\kappa ([1_C])=[1_A]$ by $\kappa ([1_C])=[1_A]$, we still have possible non unital such homomorphism $\phi$.}}

\end{lem}
\begin{proof}
Write $C={{A}}(F_1, F_2, \phi_0, \phi_{{1}}).$
Denote by $M$ the constant of \ref{MextC} for $G=K_0(C)\subset K_0(F_1)=\Z^l.$
Let $\kappa: K_0(C)\to K_0(A)$ be a unital positive homomorphism satisfying the condition of the lemma. Since $K_0(A)$ is finitely generated,  simple and $\kappa$ is strictly positive, there is $N$ such that for any nonzero positive element $x\in K_0(C)_+$, one has that $N\kappa(x)> 2[1_A]$. Let $K$ be the natural number of Proposition \ref{pl-lifting} with respect to $C$ and $N$.

We may also assume that $M_r(C)$ contains { {a set of}} minimal projections such that
every minimal element of $K_0(C)_+\setminus \{0\}$ is represented by {{a}} minimal projection { {from the set}}.

By Lemma \ref{decomposition2}, for any positive map $\kappa$ with multiplicity $M$, one has
$\kappa=\kappa_1+\kappa_2$ and there are  positive \hm s $\lambda_1: K_0(C)\to \Z^n,$
$\gamma_1: \Z^n \to K_0(A),$ $\lambda_2: K_0(C)\to K_0(C')$  such that
$\lambda_1$ has multiplicity $M,$ $\lambda_2$ has multiplicity $MK,$
$\kappa_1=\gamma_1\circ\lambda_1,$ $\kappa_2=\imath_{*0}\circ \lambda_2$ and
$C'\subset A$ is a \SCA\, with $C'\in {\cal C},$
where
$\imath: C'\to A$ is the embedding. Moreover,
\begin{equation}\label{120514-ext-1}
\lambda_2([1_C])=[1_{C'}]\quad\textrm{and}\quad N \lambda_2(x)>\lambda([1_C])>0
\end{equation}
for any $x\in K_0(C)_+\setminus \{0\}$.

Let $R_0$ be as in \ref{MextC} associated with $K_0(C)=G\subset \Z^l$ and $\lambda_1: K_0(C)\to \Z^n$ {{, which has multiplicity $M$.}}
Let $\lambda_1([1_C])=(r_1, r_2,...,r_n),$ where $r_i\in \Z_+,$ $i=1,2,...,n.$
Put $F_3=M_{r_1}\oplus M_{r_2}\oplus \cdots \oplus M_{r_n}.$
Since $A$ has stable rank one,  there is a \hm\, $\psi_0: F_3\to A$ such that
$(\psi_0)_{*0}=\gamma_1.$ Write $U=\lim_{n\to\infty}(M_{R(n)}, h_n),$ where $h_n: M_{R(n)}\to M_{R(n+1)}$ is a unital embedding. Choose $R(n)\ge R_0.$ Consider the unital \hm\, $j_{F_3}: F_3\to F_3\otimes M_{R(n)}$  defined
by $j_{F_3}(a)=a\otimes 1_{M_{R(n)}}$ for all $a\in F_3,$ and consider
the unital \hm\, $\psi_0\otimes h_{n, \infty}: F\otimes M_{R(n)}\to A\otimes U$ defined by
$\psi_0\otimes h_{n, \infty}(a\otimes b)=\psi_0(a)\otimes h_{n, \infty}(b)$ for all $a\in F_3$ and $b\in M_{R(n)}.$
We check, for any projection $p\in F_3,$
\beq\label{141204lift-1}
((\psi_0\otimes h_{n,\infty})\circ j_{F_3})_{*0}([p])=[\psi_0(p)\otimes 1_U]=(\psi_0)_{*0}([p])\in K_0(A).
\eneq
It follows that
\beq\label{141204-lift-2}
(\psi_0\otimes h_{n, \infty})_{*0}\circ (j_{F_3})_{*0}=(\psi_0)_{*0}.
\eneq
Now the map $(j_{F_3})_{*0}\circ {{\lambda}}_1: K_0(C)\to K_0({{F_3}})=\Z^n$ has multiplicity $MR(n).$
Applying \ref{MextC}, we obtain a positive \hm\, $\lambda_1': K_0(F_1)=\Z^l\to K_0(F_3)$
such that $(\lambda_1')_{*0}\circ (\pi_{{e}})_{*0}=(j_{F_3})_{*0}\circ \lambda_1.$
The above can be summarized by the following commutative diagram:
$$
\begin{array}{ccccc}
\Z^l & \stackrel{\lambda_1'}{\rightarrow}  & K_0(F_3)\otimes M_{R(n)} & \stackrel{(\psi_0\otimes h_{n, \infty})_{*0}}{\longrightarrow}& K_0(A\otimes U)\\
\hspace{0.3in}\uparrow_{(\pi_e)_{*0}} && \uparrow_{(j_{F_3})_{*0}} &&   \|\\
K_0(C) & \stackrel{\lambda_1}\longrightarrow  & K_0(F_3) & \stackrel{\gamma_1=(\psi_0)_{*0}}{\longrightarrow} & K_0(A).
\end{array}
$$
 We obtain a \hm\, $h_0: F_1\to F_3\otimes M_{R(n)}$ such that
 $(h_0)_{*0}=\lambda_1'.$
 Define $h_1=h_0\circ \pi_e: C\to F_3\otimes M_{R(n)}$ and
 $h_2=(\psi_0\otimes \psi_{n, \infty})\circ h_1: C\to A\otimes U.
 $
Then, by the %six-term
commutative diagram above,
\beq\label{120412-ext4}
(h_2)_{*0}=\kappa_1.
\eneq

Since $\lambda_2$ has multiplicity $K,$ there exists
$\lambda_2': K_0(C)\to K_0(C')$ such that
$K\lambda_2'=\lambda_2.$ Since $K_0(C')$ is weakly unperforated,
$\lambda_2'$ is positive. Moreover, by  (\ref{120514-ext-1}),
\beq\label{120414-ext4}
KN\lambda_2'(x)>K\lambda_2'([1_C])=\lambda_2([1_C])=[1_{C'}]>0.
\eneq
Since $K_0(C')$ is weakly unperforated, we have
\beq\label{120414-ext5}
N\lambda_2'(x)>\lambda_2'([1_C])>0\rforal x\in K_0(C)_+\setminus \{0\}.
\eneq
There is a projection $e\in M_k(C')$ for some integer $k\ge 1$ such that
$\lambda_2'([1_C])=[e].$  Define $C''=eM_k(C')e.$  By (\ref{120414-ext4}), $e$ is full in $C'.$
In fact $K[e]=[1_{C''}].$ In other words, $M_K(C'')\cong C'.$
By \ref{cut-full-pj}, $C''\in {\cal C}.$  By applying \ref{pl-lifting}, we obtain a  unital \hm\,
$h_2': C\to C''$ such that $(\phi_2')_{*0}=K\lambda_2'=\lambda_2.$
Put $h_2=\imath\circ \lambda_2.$
Note that $[1_A]=\kappa([1_C])=\kappa_1([1_C])+\kappa_2([1_C]),$ by conjugating a unitary, \wilog, we may assume
that $h_1(1_C)+h_2(1_C)=1_A.$
Then it is easy to check that $\phi: C\to A$ defined by $\phi(c)=h_1(c)+h_2(c)$ for all $c\in C$ meets the requirements.
\end{proof}

\begin{lem}\label{smallkkk}{\rm (cf. Lemma 9.8 of \cite{LinTAI})}
Let $A$ be a unital \CA\, and let $B_1$ be a unital separable simple \CA\, in ${\cal B}_0$ and let $B=B_1\otimes U$
for some  UHF-algebra of infinite type and $C\in {\cal C}_0$ be a \SCA\, of $B.$
Let $G\subset \underline{K}(A)$ be a finitely generated subgroup.
Suppose that there exists an ${\cal F}$-$\dt$-multiplicative  \morp\, $\psi: A\to C\subset B$ such that
$[\psi]|_G$ is well defined. Then, for any $\ep>0,$  there exists \SCA\, $C_1\cong C$ of $B$ and an ${\cal F}$-$\dt$-multiplicative \morp\, $L: A\to C_1\subset B$ such that
\beq\label{smallkkk-1}
[L]|_{G\cap K_0(A, \Z/k\Z)}=[\psi]|_{G\cap K_0(A, \Z/k\Z)}\andeqn \tau(1_{C_1})<\ep
\eneq
for all $\tau\in T(B)$ and for all $k\ge 1$ so that $G\cap K_0(A, \Z/k\Z)\not=\{0\},$
where $L$ and $\psi$ are viewed as maps to $B.$ Furthermore,
if $[\psi]|_{G\cap K_0(A)}$ is positive {{then}} so is $[L]|_{G\cap K_0(A)}.$
\end{lem}

\begin{proof}
We use the fact that $U\otimes U\cong U.$ So we may write, without loss of generality,
that $C\subset B_1\otimes U\otimes 1.$
\Wlog, we may assume that $\psi(1_A)=p$ is a projection.
Let $1>\ep>0.$  Suppose that
$$
G\cap K_0(A,\Z/k\Z)=\{0\}\rforal k\ge K
$$
for some integer $K\ge 1.$
Find a projection $e_0\in U$ such that $\tau_0(e_0)<\ep$
and $1_U=e_0+\sum_{i=1}^m p_i,$ where $m=2l K!$ and $1/l<\ep$ and $p_1, p_2,...,p_m$ are mutually orthogonal
and mutually equivalent projections in $U.$
Choose $C_1=C\otimes e_0.$ Then $C_1\cong C.$
Let $\phi: C\to C_1$ be the isomorphism defined by $\phi(c)=c\otimes e_0$ for all $c\in C.$
Put $L=\phi\circ \psi.$ Note that $K_1(C)=K_1(C_1)=\{0\}.$  Both $[L]$ and $[\psi]$ maps
$K_0(A,\Z/k\Z)$ to $K(B)/kK_0(B)$ and factor through $K_0(C, \Z/k\Z).$ It follows that
$$
[L]|_{G\cap K_0(A. \Z/k\Z)}=[\psi]|_{G\cap K_0(A, \Z/k\Z)}.
$$
In case that $[\psi]|_{G\cap K_0(A)}$ is positive, from the definition of $L,$ it is clear
that $[L]|_{G\cap K_0(A)}$ is also positive.
\end{proof}

\begin{thm}\label{preBot2}
Let $C$ and $A$ be unital stably finite \CA s and let  $\af\in KK_{e}(C, A)^{++}.$

{\rm  (i)}\,\, If  $C\in{\bf H}$, or $C\in \mathcal C,$ and $A_1$
%be
is a unital simple C*-algebra in ${\cal B}_{0}$ and
$A=A_1\otimes U$ for some UHF-algebra $U$ of infinite type,
 then there exists a sequence of completely positive linear maps  $L_n: C\to A$ such that
\beq\label{preBot2-1}
\lim_{n\to\infty} \|L_n(ab)-L_n(a)L_n(b)\|=0 \tforal a,\, b\in C\tand\,
{[}L_n{]}=\af;
\eneq

{\rm (ii)}\,\, If $C\in {\cal C}_0$ and $A_1\in {\cal B}_1,$ the above also holds;

{\rm (iii)} If $C=M_n(C(S^2))$ for some integer $n\ge 1,$  $A=A\otimes U$ and  $A_1\in {\cal B}_1,$
then there is a unital \hm\, $h: C\to A$ such that $[h]=\af;$

{\rm (iv)} If $C\in {\bf H}$ with torsion $K_1(C),$ $C\not=M_n(C(S^2)),$ and $A=A_1\otimes Q,$ where $A_1$ is unital and $A$ has stable rank one, 
then there exists a unital \hm\, $h: C\to A$ such that $[h]=\af$;

{\rm (v)} If $C=M_n(C(\T))$ for some integer $n\ge 1,$ then
for any unital \CA\, $A$ with stable rank one, there is a unital \hm\, $h: C\to A$ such that
$[h]=\af.$
\end{thm}

\begin{proof}
Let us first consider (iii). This is a special case of Lemma 2.19 of \cite{Niu-TAS-II}. Let us provide a proof here.
In this case one has that $K_0(C)=\Z\oplus {\rm ker}\rho_C\cong \Z\oplus \Z$ is free and $K_1(C)=\{0\}.$
Write $A=A_1\otimes U,$ where $K_0(U)=D\subset \Q$ is identified with a dense subgroup of $\Q$ and $1_U=1.$
Let $\af_0=\af|_{K_0(C)}.$ Then $\af_0([1_C])=[1_A]$ and $\af_0(x)\in {\rm ker}\rho_A$ for all $x\in {\rm ker}\rho_C.$
Let $\xi\in {\rm ker}\rho_C=\Z$ be a generator and $\af_0(\xi)=\zeta\in {\rm ker}\rho_A.$
Let $B_0$ be a the unital simple AF-algebra with
$$
(K_0(B), K_0(B)_+, [1_{B_0}])=(D\oplus \Z, (D\oplus \Z)_+, (1,0)),
$$
where
$$
(D\oplus \Z)_+=\{(d,m): d>0,m\in \Z\}\cup \{(0,0)\}.
$$
It follows from \cite{EL-Emb}  that there is a unital \hm\, $h_0: C\to B$ such that
$h_{*0}(\xi)=(0, 1).$
There is a positive order-unit preserving \hm\, $\lambda: D\to K_0(A)$ (given by the embedding
$a\to 1_A\otimes a$ from $U\to A_1\otimes U$).
Define a \hm\, $\kappa_0: K_0(B)\to K_0(A)$
by $\kappa_0(r)=\lambda(r)$ for all $r\in D$ and $\kappa_0((0,1))=\zeta.$
Since $A$ has stable rank one,   it is known and easy  to find a unital \hm\, $\phi: B\to A$ such that
\beq\label{12/24/14-1}
\phi_{*0}=\kappa_0.
%((0,1))=\zeta.
\eneq
Define $L=\phi\circ h_0.$ Then, $[L]=\af.$  This proves (iii).

For ({{i}}v), we note that $K_i(A)$ is torsion free and divisible. Write $C=PM_n(C(X))P,$ where $X$ is a connected finite CW complex  and $P\in M_n(C(X))$ is a projection. Note that  $X\not=S^2.$ In this case $K_0(C)=\Z\oplus {\rm Tor}(K_0(C)),$ $K_1(C)=\{0\},$ or $K_0(C)=\Z$ and
$K_1(C)$ is finite.  Suppose that $P$ has rank $r\ge 1.$
Choose $x_0\in X.$ Let $\pi_{x_0}: C\to M_r$ be defined by $\pi_{x_0}(f)=f(x_0)$ for all $f\in C.$
Suppose that $e=(1,0)\in \Z\oplus {\rm Tor}(K_0(C))$ or $e=1\in \Z.$
Choose a projection $p\in A$ such that $[p]=\af_0(e)$ (this is possible since $A$ has stable rank one).
There is a unital \hm\, $h_0: M_r\to A$ such that $h_0(e_{1,1})=p,$ where $e_{1,1}\in M_r$ is a rank one projection.
Define $h: C\to A$ by $h=h_0\circ \pi_{x_0}.$ One verifies that $[h]=\af.$

Now we prove (i) and (ii).
If $C\in\mathbf H$ and $C\neq S^2$, the statement follows from the same argument as that of Lemma 9.9 of \cite{LinTAI}
by replacing Lemma 9.8 of \cite{LinTAI} by \ref{smallkkk} above (and replace $F$ by $C$ and $F_1$ by $C_1$) in
the proof of Lemma 9.9 of \cite{LinTAI}.
%The case $C= C(S^2)$ has been proved.

Assume that $C\in \mathcal C$. By considering ${\rm Ad}\, w\circ L_n|_C$ for suitable
unitary { {$w$}} (in $M_r(A)$),  we may replace $C$ by $M_r(C)$ for some $r\ge 1$
so that $K_0(C)_+$ is generated by minimal projections $\{p_1,p_2,...,p_d\}\subset C$ (see \ref{FG-Ratn}).
Since $A$ is simple and $\af(C_+\setminus \{0\})\subset K_0(A)_+\setminus\{0\},$ there exists an integer
$N\ge 1$ such that
\begin{equation}\label{prebot2141204-n1}
N\af([p])>2[1_A]\rforal [p]\in K_0(C)_+\setminus\{0\}.
\end{equation}

Let $M\ge 1$ (in place of $K$) be the integer given by \ref{pl-lifting} associated with $N$ and $C.$
Since $C$ has a separating family of finite-dimensional representations, by Theorem \ref{kkmaps}, there exist two sequences of completely positive contractions $\phi_n^{(i)}: C\rightarrow{A}\otimes{\cal K}\ (i=0, 1)$ satisfying the following:
%\begin{enumerate}
{\rm (a)} $\|{\phi_n^{(i)}(ab)-\phi_n^{(i)}(a)\phi_n^{(i)}(b)}\|\rightarrow0$, $\rforal a, b\in C$,
      as $n\rightarrow\infty$,

 \noindent
{\rm (b)} for any $n$, $\phi_n^{(1)}$  is a \hm\, with  finite
      dimensional  range
      %\SCA\,  of $A\otimes\cal K$
      and, consequently, for any finite subset
      ${\cal P}\subset \underline{K}(C)$, the map $[\phi_n^{(i)}]|_{\cal P}$ are well defined for
     all sufficiently
     large $n$,

\noindent
%\item
{\rm (c)} for each finite subset ${\cal P}\subset \underline{K}(C)$, there exists $m>0$
      such that $$[\phi_n^{(0)}]|_{\cal P}=\alpha+[\phi_n^{(1)}]|_{\cal P} \tforal n>m,$$

   %   for all $n>m$,
%\item
%\noindent
%{\rm (d)} for each $n$, the map $\phi_n^{(1)}$ is a homomorphism
     % on $C$ with image in a finite dimensional C*-algebra.
      %{\blue{---should we combine this item with item (b)---Gong}----L: We could. But this way we  kept the reference easy}
%\end{enumerate}

Since $C$ is semiprojective and the positive cone of the $K_0$-group is finitely generated, there are homomorphisms $\phi_0$ and $\phi_1$ from $C \to A\otimes\cal K$ such that $$[\phi_0]=\alpha+[\phi_1].$$ Without lose of generality, let us assume that $\phi_0$ and $\phi_1$ are *-homomorphisms from $C$ to $\textrm{M}_r(A)$ for some $r$. Note that $\textrm{M}_r(A)\in \mathcal B_0$ (or in ${\cal B}_1,$ when $C\in {\cal C}_0$).

Since $K_i(C)$ is finitely generated ($i=0,1$), there exists $n_0\ge 1$ such that
every element $\kappa\in KL(C,A)$ is determined by
$\kappa$ on $K_i(C)$ and $K_i(C, \Z/n\Z)$ for $2\le n\le n_0,$ $i=1,2$  (see Corollary 2.11 of \cite{DL}).
Let ${\cal P}\subset \underline{K}(C)$ be a finite subset which generates
$$
\bigoplus_{i=0,1}(K_i(C)\oplus\bigoplus_{2\le n\le n_0}K_i(C, \Z/n\Z)).
$$
Choose $K=n_0!.$
Let $\cal G$ be a finite subset of $\textrm{M}_r(A)$ which contains $\{\phi_0(p_i), \phi_1(p_i); i=1, ..., d\}$.
%Also denote by ${\mathcal P}_0=\{[\phi_0(p_i)], [\phi_1(p_i)], \af([p_i]); i=1, ..., d\}$.
We may assume that ${ {\{[p_i],i=1,2,\cdots d\}}} \subset {\cal P}.$

Let
 $$
 T=\max\{\tau(\phi_0(p_i))+KM\tau(\phi_1(p_i)): 1\le i\le d; \tau\in T(A)\}.
 $$
 Choose $r_0>0$ such that
 \beq\label{prebot2-nn3}
 NTr_0<1/2.
 \eneq
 Let ${\cal Q}=[\phi_0]({\cal P})\cup [\phi_1]({\cal P})\cup\af({\cal P}).$
 Let $1>\ep>0.$
 By Lemma \ref{lem-compress}, for $\ep$ and $r_0$ above,
 there is a non-zero projection $e\in M_r(A),$ a \SCA\, $B\in {\cal C}_0$ {{(or $B\in {\cal C}$ for case (ii))}} with
 $e=1_B,$ $\mathcal{G}$-$\ep$-multiplicative \morp s $L_1 :\mbox{M}_r(A)\to(1-e)\mbox{M}_r(A)(1-e)$
 and $L_2: M_r(A)\to B$ with the following properties:
\begin{enumerate}
\item $\|L_1(a)+L_2(a)-a\|<\ep/2\tforal a\in {\cal G};$
\item $[L_i]|_{\cal Q}$  is well defined, $i=1,2;$
\item $[L_1]|_{\cal Q}+[\imath\circ L_2]|_{\cal Q}=[{\rm id}]|_{\cal Q};$
\item\label{08-10-lem-cond-02} $\tau\circ[L_1](g)\leq r_0\tau(g)$ for all $g\in {\mathcal P}_0$ and $\tau\in T(A)$;
\item\label{08-10-lem-cond-03}  For any $x\in {\cal Q},$ there exists $y\in \underline{K}(B)$ such that
   $x-[L_1](x)=[\imath\circ L_2](x)=KM[\imath](y)$
   and,
%\item\label{08-10-lem-cond-04} for any $d\in {\cal P}_0,$ there exist positive element  $f\in {K_0}(B)_+$
%      such that
   %   $$d -[L_1](d)=[\imath\circ L_2](d)=M[\imath](f),$$
%\item $\tau\circ[L](g)\leq r_0\tau(g)$ for all $g\in
   %   {\mathcal P}$ and $\tau\in T(A)$;
\item\label{cond3} There exist positive elements $\{f_i\}\subset{K_0}(B)_+$
      such that for  $i=1,...,n,$
      $$\af([p_i])-[L_1](\af([p_i])=[\imath\circ L_2](\alpha([p_i]))=KM\imath_{*0}(f_i).$$
      \end{enumerate}
      where $\imath: B\to A$ is the embedding.
      %Note that, with small $r_0,$ we may assume that $f_i\in K_0(A)_+\setminus\{0\}.$
%\item the image of $\mathrm{Id}-L$ is in a sub-C*-algebra $S\subset \textrm{M}_r(A)$ with $S\in\mathcal C_0$. (Note that by Condition \eqref{cond3}, any irreducible representation of $S$ must be divisible by M.)

By (\ref{08-10-lem-cond-03}), since $K=n_0!,$
\beq\label{prebot2-nn4}
[\imath\circ L_2\circ \phi_0]|_{K_i(C, \Z/n\Z)\cap {\cal P}}=[\imath\circ L_2\circ \phi_1]|_{K_i(C, \Z/n\Z)\cup {\cal P}}=0,\,\,\,i=0,1,
\eneq
and $n=1,2,...,n_0.$  It follows that
\beq\label{prebot2-n1}
%&&[L_1\circ \phi_0]|_{K_1(C)\cap {\cal P}}=[\phi_0]|_{K_1(C)\cup{\cal P}},\,\,\, [L_1\circ \phi_1]|_{K_1(C)}=[\phi_1]|_{K_1(C)},\\
&&{[}L_1\circ\phi_0{]}|_{K_i(C, \Z/n\Z)\cap {\cal P}}=[\phi_0]|_{K_i(C, \Z/n\Z)\cap {\cal P}}\quad\textrm{and}\quad\\
&& [L_1\circ\phi_1]|_{K_i(C, \Z/n\Z)\cap {\cal P}}=[\phi_1]|_{K_i(C, \Z/n\Z)\cap {\cal P}},
\eneq
$i=0,1$ and $n=1, 2, ...,n_0.$
Furthermore, { {for the case $B\in {\cal C}_0$ we have}} $K_1(B)=0,$ { {and consequently}}
\beq\label{prebot2-nn6}
[\imath\circ L_2]|_{K_1(C)\cap {\cal P}}=0.
\eneq

It follows that
\beq\label{prebot2-nn5}
[L_1\circ \phi_0]|_{K_1(C)\cap {\cal P}}=[\phi_0]|_{K_1(C)\cap{\cal P}},\,\,\, [L_1\circ \phi_1]|_{K_1(C)\cap {\cal P}}=[\phi_1]|_{K_1(C)\cap {\cal P}}.\\
\eneq
In the second case when  we assume that $C\in {\cal C}_0$ and $A_1\in {\cal B}_1,$ then $K_1(C)=0.$  Therefore
(\ref{prebot2-nn5}) above also holds.

Denote by $\Psi:=\phi_0\oplus\bigoplus_{KM-1}\phi_1$. One then has
\begin{eqnarray*}
[L_1\circ\Psi]_{K_i(C,\ \mathbb Z/n\mathbb Z)\cap {\cal P}} & = & [L_1\circ\phi_0]_{K_i(C,\ \mathbb Z/n\mathbb Z)\cap {\cal P}}+ (KM-1)[L_1\circ\phi_1]_{K_i(C,\ \mathbb Z/n\mathbb Z)\cap {\cal P}} \\
& = & [L_1\circ\phi_0]_{K_i(C,\ \mathbb Z/n\mathbb Z)\cap {\cal P}}-[L_1\circ\phi_1]_{K_i(C,\ \mathbb Z/n\mathbb Z)\cap {\cal P}}\\
&=&[\phi_0]_{K_i(C,\ \mathbb Z/n\mathbb Z)\cap {\cal P}}-[\phi_1]_{K_i(C,\ \mathbb Z/n\mathbb Z)\cap {\cal P}}\\
&=&\alpha|_{K_i(C,\ \mathbb Z/n\mathbb Z)\cap {\cal P}},
\end{eqnarray*}
where $i=0, 1$, $n=1, 2, ...,n_0!$

%Since $A$ is simple, there exists $\delta>0$ such that $\tau(\alpha([p_i]))>\delta$ for any $\tau\in T(A)$. One then chooses $r_0$ sufficiently small such that $\tau\circ[L]\circ[\Psi]([p_i])<\delta/2$ for all $\tau\in T(A)$, and hence
By  (\ref{prebot2141204-n1}), (\ref{08-10-lem-cond-02}) and (\ref{prebot2-nn3}),
\beq\label{prebot1204-n2}
N(\tau(\alpha([p_i])-[L_1\circ\Psi]([p_i]))\ge 2-Nr_0T\ge {{3/2}}  \rforal \tau\in T(A).
\eneq
Since the strict order on ${K_0}(A)$ is determined by traces, one has that $\alpha([p_i])-[L_1\circ\Psi]([p_i])>[1_A].$

Moreover, one also has
$$\begin{array}{ll}
  & \alpha([p_i])-[L_1\circ\Psi]([p_i]) \\
  =&\alpha([p_i])-([L_1\circ\alpha]([p_i])+KM[L_1\circ\phi_1]([p_i])) \\
  =&(\alpha([p_i])-[L_1\circ\alpha]([p_i]))-KM[L_1\circ\phi_1]([p_i]) \\
  =& KM(\imath_{*0}(f_j)-[L_1]\circ[\phi_1]([p_i])) \\
  =& KMf'_j, \quad \mbox{where}\ f'_j=f_j-[L_1]\circ[\phi_1]([p_i]).
\end{array}$$
Note that $f_j'\in K_0(A)_+\setminus \{0\},$ $j=1,2,...,d.$
Define a \hm\, $\bt: K_0(C)\to K_0(A)$ by
$\bt([p_j])=Mf_j',$ $j=1,2,...,d.$  Since the $K_0(C)_+$  is generated by $[p_1], [p_2],...,[p_d],$
$\bt(K_0(C)_+\setminus\{0\})\subset K_0(A)_+\setminus\{0\}.$
%Consider the map $$\beta: \frac{1}{K}(\alpha|_{K_0}-[L\circ\Psi]_0).$$
Since $\bt$ has  multiplicity $M$, By the choice of $M$ and by Lemma \ref{liftingpl-M}, there  exists a *-homomorphism $h: C\to \textrm{M}_r(A)$ (may not be unital) such that
$$
h_{*0}=\beta\quad\mathrm{and}\quad h_{*1}=0.
$$

Consider the map $\phi':=L_{{1}}\circ\Psi\oplus (\bigoplus_{i=1}^K h): C \to A\otimes{\cal K}$, one has that
$$[\phi']|_{K_0(C)\cap {\cal P}}=[L_{{1}}\circ\Psi]|_{K_0(C)\cap {\cal P}}+K\beta=\alpha|_{K_0(C)\cap {\cal P}}.$$ It is clear that
$$K[h]|_{K_i(C,\ \mathbb Z/n\mathbb Z)\cap {\cal P}}=0,\quad i=0, 1,\ n=1, 2,...,n_0$$
and therefore
$$[\phi']_{K_i(C,\ \mathbb Z/n\mathbb Z)\cap {\cal P}}=[L_{{1}}\circ\Psi]|_{K_i(C,\ \mathbb Z/n\mathbb Z)\cap {\cal P}}=\alpha|_{K_i(C,\ \mathbb Z/n\mathbb Z)\cap {\cal P}},$$
where $i=0, 1$, $n=1, 2, ..., n_0$ .
We also have  that
$$[\phi']_{K_1(C)\cap {\cal P}}=[L_{{1}}\circ\Psi]|_{K_1(C)\cap {\cal P}}=\alpha|_{K_1(C)\cap {\cal P}}.$$
Therefore
$$[\phi']|_{\cal P}=\alpha|_{\cal P}.$$

Since $[\phi'(1_C)]=[1_A]$ and $A$ has stable rank one, there is a unitary $u$ in a matrix algebra of $A$ such that the map $\phi=\mathrm{ad}(u)\circ\phi'$ satisfies $\phi(1_C)=1_A$, as desired.

Case (v) is standard and is well known.
\end{proof}

\begin{cor}\label{est-m2a}
Any C*-algebra $A$ of Theorem \ref{RangT} is KK-attainable with respective to ${\cal B}_{u0}$.
\end{cor}

\begin{proof}
Note that $A$ is an inductive limit of C*-algebras  which are finite direct sums of \CA s in $\mathbf H$ and $\mathcal C_0$. Since KK-attainability passes to inductive limits, by Theorem \ref{preBot2}, $A$ is KK-attainable with respective to ${\cal B}_{u0}.$
\end{proof}

 %{ Compare the following with \ref{liftingpl-M}: {\bf The original statement follows from \ref{liftingpl-M}. The current one is new .}}
\begin{cor}\label{C0ext}
Let $C\in {\cal C},$   let
$A\in {\cal B}_{u0}$ and $\af: KK(C,A)^{++}$ be  such that
$\af([1_C])=[p]$ for some projection $p\in A$ and $\alpha$ is strictly positive. Then
there exists a \hm\, $\phi: C\to A$ such that
$\phi_{*0}=\af.$
\end{cor}

\begin{proof}
This is a special case of Theorem \ref{preBot2} since $C$ is semiprojective.
\end{proof}

\begin{cor}\label{Cistbk}
%Let $A_1$ be a unital separable simple C*-algebra in $\mathcal B_1$, and let $A=A_1\otimes U$ for some UHF-algebra $U$.
Let $A\in\mathcal B_{u0}$.
Then there exists a unital simple C*-algebra $B_1=\varinjlim(C_n, \phi_n),$ where
each $C_n$ is in ${\cal C}_0,$ and a UHF algebra $U$ of infinite type   such that { {for}} $B=B_1\otimes {{U}}$ , {{we have}}
$$
(K_0(B), K_0(B)_+, [1_{B}], T(B), r_{B})
=(\rho_A(K_0(A)), \rho_A(K_0(A)_+), \rho_A([1_A]), T(A), r_A).
$$
Moreover, for each $n,$ there is a unital \hm\, $h_n: C_n\otimes U\to A$ such that
\begin{equation}\label{Cistbk-1}
\rho_A\circ (h_n)_{*0}=(\phi_{n, \infty}\otimes\mathrm{id}_U)_{*0}.
\end{equation}
%there is a unital monomorphism $h: B_1\otimes U\to A$ such that
%$h$ carries the above the identity.
\end{cor}
\begin{proof}
Consider the tuple
$$(\rho_A(K_0(A)), \rho_A(K_0(A)_+), \rho_A([1_A]), T(A), r_A).$$
Since $A\cong A\otimes U_1,$ { {for a UHF algebra $U_1$ of infinite type, }} it has the (SP) property (see \cite{BKR-ADiv}), and therefore the ordered group $(K_0(A), K_0(A)_+, [1_A])$ has the (SP) property in the sense of Theorem \ref{RangT}; that is, for any positive real number $0<s<1$, there is $g\in K_0(A)_+$ such that $\tau(g)<s$ for any $\tau\in T(A)$. Then it is clear that the scaled ordered group $(\rho_A(K_0(A)), (\rho_A(K_0(A)_+), \rho_A(1_A))$ also has the (SP) property in the sense of Theorem \ref{RangT}. Therefore, by Theorem \ref{RangT}, there is a simple unital C*-algebra $B_1=\varinjlim(C_n, \phi_n)$, where each $C_n\in\mathcal C_0$  such that
$$(K_0(B_1), K_0(B_1)_+, [1_{B_1}], T(B_1), r_{B_1})\cong(\rho_A(K_0(A)), \rho_A(K_0(A)_+), \rho_A(1_A), T(A), r_A).$$

{ {Let $U=U_1$ and}}  $B=B_1\otimes U$.
%Since $U$ is of  infinite type, one has that
{ {Recall that $A\otimes U=A$.}}
{{T}}herefore $$(K_0(B), K_0(B)_+, [1_{B}], T(B), r_{B})\cong(\rho_A(K_0(A)), \rho_A(K_0(A)_+), \rho_A(1_A), T(A), r_A).$$ Clearly, $B$ has the inductive limit decomposition $$B=\varinjlim(C_n\otimes U, \phi_n\otimes \mathrm{id}_U).$$

For each $n$, consider the positive homomorphism $(\phi_{n, \infty})_{*0}: K_0(C_n)\to K_0(B_1)\cong\rho_A(K_0(A))$. Since $K_0(C_n)$ is torsion free, $K_0(C_n)_+$ is finitely generated, and the strict order on the projections of $A$ is determined by traces, there is a positive homomorphism $\kappa_n: K_0(C_n) \to K_0(A)$ such that
$$\rho_A\circ\kappa_n=(\phi_{n, \infty})_{*0} \andeqn \kappa([1_{C_n}]_0)=[1_A]_0.$$

By Corollary \ref{C0ext}, there is a unital homomorphism $h'_n: C_n \to A$ such that $(h'_n)_{*0}=\kappa_n.$ It is clear that $h_n:=h'_n\otimes \mathrm{id}_U$ satisfies the desired condition.
\end{proof}

\begin{lem}\label{ExtTrace}
{Let $C\in {\cal C}$.
Let  $\sigma>0$  and let
${\cal H}\subset C_{s.a.}$ be any finite subset.
Let $A\in\mathcal B_{u0}$. Then for any $\kappa\in KL_e(C,A)^{++}$
%strictly positive $\kappa\in \Hom(K_*(C), K_*(A))$
and any continuous affine map $\gamma: T(A)\to T_{\rm f}(C)$ which is compatible to $\kappa$, there is a unital homomorphism $\phi: C\to A$ such that
$$
[\phi]_*=\kappa\quad\textrm{and}\quad |\tau\circ \phi(h)-\gamma(\tau)(h)|<\sigma\tforal h\in {\cal H}.
$$
Moreover,
the above also holds if $C\in {\cal C}_0$ and $A\in B_{u1}.$}
\end{lem}

\begin{proof}
%Since $U$ is an UHF-algebra, without loss of generality, one may assume that $A=B$. That is, $B\cong B\otimes U$.
Without loss of generality, one may assume that every element of $\mathcal H$ has norm at most one.
Let $\kappa$ and  $\gamma$ be given.  Define $\Delta: C_+^{q, 1}\setminus\{0\}\to (0, 1)$ by
$$\Delta(\hat{h})=\inf\{\gamma(\tau)(h)/2:\ \tau\in T({{A}})\}.$$

Let $\mathcal H_1\subseteq C^+$, $\delta$, and $K$ be the finite subset, positive constant and the positive integer of \ref{ExtTraceC-D} with respect to $C$, $\Delta$, $\mathcal H$ and $\sigma/8$ (in place of $\sigma$).
Let ${\cal P}\subset K_0(C)_+$ be a finite subset which generates $K_0(C)_+.$

By Lemma \ref{cut-trace}, there is a \SCA\, $D\subseteq  A$ with $D\in\mathcal C$ and with
$1_D=p\in A,$ a continuous affine map $\gamma': T(D)\to T(C)$ such that
\begin{equation}\label{istTrC-eq1}
|\gamma'(\frac{1}{\tau(p)}\tau)(f)-\gamma(\tau)(f)|<\sigma/8\rforal \tau\in T(A) \rforal f\in\mathcal H,
\end{equation}
where $p=1_D'$, $\tau(1-p)<\sigma/(8+\sigma)$,
\begin{equation}\label{istTrC-eq2}
\gamma'(\tau)(h)>\Delta(\hat{h})\rforal \tau\in T(D) \rforal h\in\mathcal H_{1}.
\end{equation}
Denote by $\imath: D\to pAp$  the embedding.
Moreover, by \ref{cut-trace}, there are positive homomorphisms $\kappa_{0,0}: K_0(C)\to K_0((1-p){{A}}(1-p))$ and $\kappa_{1,0}: K_0(C)\to K_0(D)$ such that  $\kappa_{1,0}$ is strictly positive and has the multiplicity $K$,
\beq\nonumber
\kappa|_{K_0(A)}=\kappa_{0,0}+\imath_{*0}\circ \kappa_{1,0} \andeqn
|\gamma'(\tau)(p)-\tau(\kappa_{{1,0}}(p))|<\delta,\quad p\in{\mathcal P},\ \tau\in T(D).
\eneq
Since we assume that $A\otimes U\cong A,$ by the last part of \ref{cut-trace}, we may assume
that $\kappa_{0,0}$ is also strictly positive.
%Note that, since $D\in\mathcal C_0$, one has $\kappa_1|_{K_1(C)}=0$.
%Let $\kappa_1\in KL
Therefore, by Lemma \ref{ExtTraceC-D} (for suitable modified projection $p$  {{which still satisfies}} (\ref{istTrC-eq1}) and ({{\ref{istTrC-eq2}}})), there is a homomorphism $\phi_1: C\to D\otimes U\subset A\otimes U$ such that
$(\phi_1)_{*0}=\kappa_{1,0}$ and $$|\tau\circ\phi_1(h)-\gamma'(\tau)(h)|<\sigma/4\rforal h\in\mathcal H.$$
Since $A$ is simple, $K_i((1-p)A(1-p))=K_i(A),$ $i=0,1.$
Let $\kappa_0=\kappa-[\imath\circ \phi_1]\in KL(C,A)=KL(C, (1-p)A(1-p)).$
%By the UCT, there is $\kappa_0\in KL(C,A)$ such that
Then
$\kappa_0|_{K_0(C)}=\kappa_{0,0}.$
%and $\kappa_0|_{K_1(C)}=(\kappa-[\phi_1])|_{K_1(C)}.$
Since $\kappa_0|_{K_0(C)}
=\kappa_{0.0},$ it is strictly positive. Note that $(1-p)A(1-p)\otimes U\cong (1-p)A(1-p).$ Therefore,
by Theorem \ref{preBot2}, since $C$ is semi-projective, there is a homomorphism $\phi_0: C\to (1-p)A(1-p)$ such that $[\phi_0]=\kappa_0$. Note that, this holds for both the case that $A\in {\cal B}_{u0}$ and
the case that $C\in {\cal C}_0$ and $A\in {\cal B}_{u1}.$ Consider the homomorphism $$\phi:=\phi_0\oplus\imath\circ \phi_1: C\to (1-p)A(1-p)\oplus D\subseteq A.$$ One has that $[\phi]=\kappa$ and, for all $h\in {\cal H}$,
\begin{eqnarray*}
|\tau\circ \phi(h)-\gamma(\tau)(h)|&\leq & |\tau\circ \phi_1(h)-\gamma(\tau)(h)|+\sigma/4\\
&<&|\tau\circ \phi_1(h)-\gamma'(\frac{1}{\tau(p)}\tau|_D)(f)|+\sigma/2\\
&<&|\tau\circ \phi_1(h)- \frac{1}{\tau(p)}\tau\circ\phi_1(h)|+3\sigma/4<\sigma\\
%&<&\sigma,\quad \forall h\in {\cal H},
\end{eqnarray*}
as desired.%for any $h\in\mathcal H$.
\end{proof}
It turns out that the KK-attainability implies the following existence theorem.
\begin{prop}\label{add-tr}
Let $A\in \mathcal B_0$, and assume that $A$ is KK-attainable with respective to ${\cal B}_{u0}$. Then for any $B\in \mathcal B_{u0}$, any $\alpha\in KL^{++}(A, B)$, and any affine continuous map $\gamma: T(B)\to T(A)$ which is compatible to $\alpha$, there is a sequence of completely positive linear maps $L_n: A\to B$ such that
 $$ \lim_{n\to\infty} ||L_n(ab)-L_n(a)L_n(b)||  =  0
 %,\quad \forall
 \tforal a, b \in A, $$
  $$ [L_n]  =  \alpha \andeqn $$
$$ \lim_{n\to\infty}\sup\{|\tau\circ L_n(f)-\gamma(\tau)(f)|: \tau\in T(A)\} = 0 \tforal f\in A.$$
\end{prop}

\begin{proof}
The proof is the same as that of Proposition 9.7 of \cite{LinTAI}. Instead of using Lemma 9.6 of \cite{LinTAI}, one uses Lemma \ref{ExtTrace}.
\end{proof}

\begin{rem}\label{RUHF-TAD}
The condition that $A=A\otimes U$ ($A\in {\cal B}_1$, or $A\in {\cal B}_0$) in the statements in Section \ref{KK-BBlock} can be easily eased to the condition
that $A\in {\cal B}_1$ (or ${\cal B}_0$) such that $A$ is tracially approximately divisible.  As a consequence, with almost no additional efforts, one can also ease the the same condition for $B.$ Since we eventually do not need to assume $A\cong A\otimes U,$ to shorten the length of this article, we conveniently use the stronger assumption.
\end{rem}

\section{The class ${\cal N}_1$}

Let $A$ be a unital \CA\, such that $A\otimes Q\in {\cal B}_1.$
In this section,  we will show that  $A\otimes Q\in {\cal B}_0.$ Note that this is proved without assuming 
$A\otimes Q$ is nuclear.
However, it implies that ${\cal N}_1={\cal N}_0.$ 
If we assume that $A\otimes Q$ has finite nuclear dimension, using \ref{ExtTrace} and a characterization 
of of $TAS$ by Winter, a much shorter proof of \ref{B0=B1}could easily given here. 
%Note that this is proved without assuming 
%$A\otimes Q$ is nuclear.
%let us show that any unital simple C*-algebra which can be tracially approximated by subhomogeneous C*-algebras with dimensions of spectra at most one are in fact can be rationally tracially approximated by Elliott-Thomsen algebras with trivial $K_1$-group; that is, its tensor product with an UHF algebra is in the class $\mathcal B_1$. In particular, this shows that any simple inductive limits of subhomogeneous C*-algebras with dimensions of spectra at most one are classifiable.

\begin{lem}\label{K1inj}
Let $A\in {\cal B}_1$ such that $A\cong A\otimes Q.$
Then the following holds:
For any $\ep>0,$  any two non-zero mutually orthogonal elements $a_1, a_2\in A_+$ and any finite subset ${\cal F}\subset A,$
there exists a projection $q\in A$ and a \SCA\, $C_1\in {\cal C}$ with $1_{C_1}=q$ such that
\begin{enumerate}
\item $\|[x, q]\|<\ep/16$, $x\in\mathcal F$,
\item $pxp\in_{\ep/16} C_1$, $x\in\mathcal F$, and
\item $1-q\lesssim a_1.$
\end{enumerate}
%Suppose that $K_1(C_1)=\Z^m\oplus G_0$ such that
%$j_{*1}(G_0)=\{0\}$ and $j_{*1}|_{\Z^m}$ is injective, where

Suppose $\Delta: (C_1)_{{+}}^{q,{\bf 1}}\setminus \{0\}\to (0,1)$ is an order preserving map
such that
$$
\tau(c)\ge \Delta(\hat{c})\rforal {{~~\tau \in T(A)~~ \mbox{and}~~}} c\in (C_1)^{\bf 1}_+\setminus \{0\}.
$$
%{\blue {(
%Such 
{\rm(} By \ref{Ldet} such $\Delta$  always exists.{\rm )}
%%%exiats
% since $A$ is simple.
%)}}---{\red{\bf Zhuang: Could you help
%to find a reference---I did somewhere---could be in this paper}}
Suppose also that
${\cal H}\subset (C_1)_+\setminus \{0\}$
% is a finite subset
%that $j: C_1\to A$ is the embedding, and suppose
%and that
and ${\cal F}_1\subset C_1$ are finite subsets. Then, there exists  another projection $p\in A$
and a \SCA\, $C_2\in {\cal C}$ with $p=1_{C_2}$ such that $p\le q,$ and a  unital \hm\,
$H: C_1\to C_2$ such that

\begin{enumerate}\setcounter{enumi}{3}
\item $\|[x, p]\|<\ep/16$, $x\in {\cal F}$,
%\item
\item  $\|H(y)-{{p}}y{{p}}\|<\ep/16$, $y\in {\cal F}_1$, and
%\item

\item $1-p\lesssim a_1+a_2$.
\end{enumerate}
%
%\beq\label{K1inj-2}
%\|[x, p]\|<\ep/16\tforal x\in {\cal F},\\
%\|H(y)-y\|<\ep/16\tforal y\in {\cal F}_1,\\
%1-p\lesssim a_1+a_2,
%\eneq
Moreover $K_1(C_1)=\Z^m\oplus G_0$ such that
$ H_{*1}(G_0)=\{0\},$ and $H_{*1}|_{\Z^m}$ and $(j\circ H)_{*1}|_{\Z^m}$ are both  injective,
where $j: C_2\to A$ is the embedding ($m$ could be zero, in this case, $G_0=K_1(C_1)$).
Furthermore, we may assume that
\begin{equation*}
\tau(j\circ H(c))\ge  3\Delta(\hat{c})/4\tforal c\in {\cal H}\tforal \tau\in T(A).
\end{equation*}

\end{lem}
\begin{proof}
%The case that $A=\overline{\cup_{n=1}^{\infty} C_n},$ where $C_n\subset C_{n+1}$ is a sequence of \SCA s
%in ${\cal C}_1$ follows from the definition easily. In fact, in this case, we may choose $p=q=1_A.$

%We now assume that this is not the case.

Since $A\in {\cal B}_1,$
%For any $\ep>0,$  $a\in A_+\setminus \{0\}$ and any finite subset ${\cal F}\subset A,$
there exists a projection $q\in A$ and a \SCA\, $C_1\in {\cal C}$ with $1_{C_1}=q$ such that
%\begin{enumerate}
%\item

(a) $\|[x, q]\|<\ep/16$, $x\in\mathcal F$,
%\item

(b) $pxp\in_{\ep/16} C_1$, $x\in\mathcal F$, and
%\item

(c) $1-q\lesssim a_1.$
%\end{enumerate}
%\beq\label{K1inj-1}
%\|[x, q]\|<\ep/16\tand qxq\in_{\ep/16} C_1 \tforal x\in {\cal F},\\\tand
%1-q\lesssim a_1
%\eneq

There are two non-zero mutually orthogonal elements $a_2'$ and $a_3\in \overline{a_2Aa_2}.$
Note that $A\cong A\otimes Q.$ Therefore $K_1(A)$ is torsion free. Denote by $j: C_1\to qAq$ the embedding.
Since $K_1(C_1)$ is finitely generated,
we may write $K_1(C_1)=G_1\oplus G_0,$ where $G_1\cong \Z^{m_1},$ $j_{*1}|_{G_1}$ is injective and $j_{{*1}}|_{G_0}=0.$
Define
$$
\sigma={{\min}}\{\Delta(\hat{h})/{{16}}: h\in {\cal H}\}>0.
$$
Choose an element $a_2''\in \overline{a_2'Aa_2'}$ such that
$d_\tau(a_2'')<\sigma$ for all $\tau\in T(A)$ {{, where recall $d_{\tau} (a) =\lim_{n\to \infty} \tau (a^{\frac1n})$.}}

Suppose that $G_0$ is generated by $[v_1], [v_2],...,[v_l],$ where $v_i\in U(C_1)$ (note $C_1$ has stable rank one).
Then, we may write $v_k=\prod_{s=1}^{l(k)} \exp (i h_{s,k}),$ { {(since $j_{*1}([v_k])=0$ in$K_1(A)$)}} where $h_{s,k}\in (qAq)_{s.a},$ $s=1,2,...,l(k),$ $k=1,2,...,l.$

Let ${\cal F}_1$ be a finite subset of $C_1$ which also has the following property:
if $x\in {\cal F},$ there is $y\in {\cal F}_1$ such that $\|{{p}}x{{p}}-y\|<\ep/16.$

Let ${\cal F}_2$ be a finite subset  of $qAq$ to be determined which at least contains ${\cal F}\cup {\cal F}_1\cup{\cal H}$ and
$h_{s,k}, \exp(i h_{s,k}),$ $s=1,2,...,l(k),$ $k=1,2,...,l.$

Let $0<\dt<\min\{\ep/64, \sigma/4\}$ to be determined.   Since $qAq\in {\cal B}_1,$
one obtains a non -zero projection $q_1\in qAq$ and a \SCA\, $C_2{{\subset}} qAq$  such that
%\begin{enumerate}
%\item

(d) $\|[x,\, q_1]\|<\dt$, $x\in {\cal F}_2$,
%\item

(e) $q_1xq_1\in_{\dt} C_2$, $x\in {\cal F}_2$,
%\item

(f) $1-q_1\lesssim a_2{{''}}.$
%\end{enumerate}
%\beq\label{K1inj-2}
%\|[x,\, q_1]\|<\dt\rforal x\in {\cal F}_2,\\
%q_1xq_1\in_{\dt} C_2\rforal x\in {\cal F}_2\andeqn\\
%1-q_1\lesssim a_2'.
%\eneq
%We can assume that $\|q_1xq_1\|\ge \|x\|/2$ for all
%$x\in {\cal F}_1.$

With sufficiently small $\dt$ and large ${\cal F}_2,$ using the semi-projectivity of $C_1,$
we obtain unital \hm s $h_1: C_1\to C_2$
% and $h_2: C_1\to (q-q_1)A(q-q_1)$
such that
\beq\label{K1inj-3}
\|h_1(a)-q_1aq_1\|<\min\{\ep/16, \sigma\}\rforal a\in {\cal F}_1\cup {\cal H}.
\eneq
%We may assume that both are  injective{\bf----I do not get a chance to check that--L} , by choosing sufficiently large
%if ${\cal F}_2$  and sufficiently small $\dt.$
One computes that
\beq\label{K1inj-4}
\tau(j\circ h_1(c))\ge {{7}}\Delta(\hat{c})/{{8}} \rforal c\in {\cal H},
\eneq
where we also use $j$ for the embedding from $C_2$ into $qAq.$
Note that, when ${\cal F}_2$ is large enough, $(h_1)_{*1}(G_0)=\{0\}.$

We may write {{$K_1(C_1)=G_2\oplus G_{2,0}\oplus G_{2,0,0},$ where $G_2\cong \Z^{m_2}$ with $m_2\le m_1,$
$G_2$ is a subgroup of $G_1,$ $G_{2,0,0}\supset G_0,$  $(h_1)_{*1}(G_{2,0,0})=\{0\},$
$(j\circ h_1)_{*1}|_{G_{2,0}}=0$ and $(j\circ h_1)_{*1}|_{G_2}$ is injective.
Here we use the fact that $K_1(A)$ is torsion free.}}
If $G_{2,0}=\{0\}${{,}} we are done.
Otherwise $m_2<m_1.$ We also note that $(j\circ h_1)_{*1}(G_{2,0}\oplus G_{2,0,0})=\{0\}.$

{{We will  repeat the process}} to construct $h_2$, and consider  $h_2\circ h_1.$
Then we may write $K_1(C_1)=G_3\oplus G_{3,0}\oplus G_{3,0,0}$ with
$G_3\cong \Z^{m_3}$ ($m_3\le m_2$) $G_3\subset G_2,$ $G_{3,0,0}\supset G_{2,0}\oplus G_{2,0,0},$
$(j\circ h_2\circ h_1)_{*1}(G_{3,0})=\{0\}$ and $(j\circ h_2\circ h_1)_{*1}|_{G_3}$ is injective.
Again, if $G_{3,0}=\{0\},$ we are done (choose $H=h_2\circ h_1$).  Otherwise, $m_3<m_2<m_1.$
We continue this process. Since $G_1$ is a commutative noetherian ring, this process stops at  a finite stage.
This proves the lemma.

%We may write $K_1(C_1)=G_2\oplus G_{2,0},$ where $G_2\cong \Z^{m_2}$ with $m_2\le m_1,$
%$G_2$ is a subgroup of $G_1$ and $G_{2,0}\supset G_0,$ and $(h_1)_{*1}(G_{2,0})=\{0\}.$ Moreover,
%there is a subgroup $H_1\subset G_1$ such that $G_1/H_1\cong G_2$ and $G_{2,0}=H_1\oplus G_0\supset G_0.$
%If $(j\circ h_1)_{*1}|_{G_2}$ is injective, we are done.
%
%Otherwise, we continue this process. Since $G_1$ is a commutative noetherian ring, this process stops at  a finite stage.
%This proves the lemma.
\end{proof}

\begin{thm}\label{B0=B1}
Let $A_1$ be a unital separable amenable \CA\,  such that $A_1\otimes Q\in {\cal B}_1.$  Then $A_1\otimes Q\in {\cal B}_0.$
\end{thm}

\begin{proof}
Let $A=A_1\otimes Q.$ Suppose that $A\in {\cal B}_1.$
Let $\ep>0,$  let $a\in A_+\setminus\{0\}$ and let ${\cal F}\subset A.$
Since $A$ has property (SP), we obtain three non-zero and mutually orthogonal projections $e_0, e_1, e_2\in \overline{aAa}.$
There exists a projection $q_1\in A$ and a \SCA\, $C_1\in {\cal C}$ with $1_{C_1}=q_1$ such that
\beq\label{B0=B1-n1}
\|[x,\, q_1]\|<\ep/16\tand q_1xq_1\in_{\ep/16} C_1 \tforal x\in {\cal F},\\\tand
1-q_1\lesssim e_0.
\eneq
Let ${\cal F}_1\subset C_1$ be a finite subset
such that, for any $x\in {\cal F},$ there is $y\in {\cal F}_1$ such that
$\|q_1xq_1-y\|<\ep/16.$
For each $h\in (C_1)_+\setminus\{0\},$ define
$$
\Delta(\hat{h})=(1/2)\inf\{\tau(h): \tau\in T(A)\}.
$$

Then $\Delta: (C_1)_+^{q, {\bf 1}}\setminus \{0\}\to (0,1)$ preserves the order.
Let ${\cal H}_1\subset (C_1)^{\bf 1}_+\setminus \{0\}$ be a finite subset, $\gamma_1, \gamma_2>0,$  $\dt>0,$ ${\cal G}\subset C_1$ be a finite subset,
${\cal P}\subset \underline{K}(C_1)$ be a finite subset, ${\cal H}_2\subset (C_1)_{s.a.}$ be a finite subset and
let ${\cal U}\subset J_c^{(1)}(U(C_1)/U_0(C_1))$  be a finite subset required  by \ref{UniCtoA} and \ref{RemUniCtoA} for $C=C_1,$
for $\ep/16$ (in place of $\ep$)  and ${\cal F}_1$ (in place of ${\cal F}$) and $\Delta/2$ (in place of $\Delta$).

%Suppose that $K_1(C_1)=\Z^m\oplus G_{00}$ such that
%$j_{*1}(G_{00})=\{0\}$ and $j_{*1}|_{\Z^m}$ is injective.
%and suppose that ${\cal F}_1\subset C_1$ is a finite subset.
By \ref{K1inj}, there exists  another projection $q_2\in A$
and a \SCA\, $C_2\in {\cal C}$ with $q_2=1_{C_2}$ such that $q_2\le  q_1,$ and a unital \hm\,
$H: C_1\to C_2$ such that
\beq\label{K1inj-2}
\|[x, q_2]\|<\ep/16\tforal x\in {\cal F},\\\label{K1inj-2+}
\|H(y)-q_2yq_2\|<\ep/16\tforal y\in {\cal F}_1,\\
\tau(j\circ H(c))\ge 3 \Delta(\hat{c})/4\tforal c\in {\cal H}\andeqn\\
1-q_2\lesssim e_0+e_1.
\eneq
Moreover, we may write $K_1(C_1)=\Z^m\oplus G_{{00}},$ where
$H_{*1}(G_{00})=\{0\},$ $H_{*1}|_{\Z^m}$ and $( j\circ H)_{*1}|_{\Z^m}$ are injective, where
$j: C_2\to A$ is the embedding.
Let $A_2=q_2Aq_2$ and denote by $j_1: C_2\to A_2$ the embedding.

By \ref{RangT}, there exists a unital  simple \CA\, $B\cong B\otimes Q$  such that
$B=\lim_{n\to\infty} (B_n, \imath_n)$ such that each $B_n=B_{n,0}\oplus B_{n,1}$ with $B_{n,0}\in {\bf H}$ and
$B_{n,1}\in {\cal C}_0,$  $\imath_n$ is injective,
\beq\label{B0=B1-2}
\lim_{n\to\infty}\max\{\tau(1_{B_{n,0}}):\tau\in T(B)\}=0\andeqn
{\rm Ell}(B)={\rm Ell}(A_2).
\eneq

We may assume that ${\cal U}={\cal U}_1\cup {\cal U}_0,$
such that $\pi({\cal U}_1)$ generates $\Z^m,$ and $\pi({\cal U}_0)\subset G_0,$
where $\pi: U(C_1)/CU(C_1)\to K_1(C_1)$ is the quotient map. Here
$J_c^{(1)}: K_1(C_1)\to U(C_1)/CU(C_1)$ is a fixed splitting map defined in \ref{Dcu}.
Suppose that ${\bar v_1},{\bar v_2},...,{\bar v_m}$ form a set of free generators for $J_c^{(1)}(\Z^m).$ \Wlog, we may assume
that ${\cal U}_1=\{{\bar v_1},{\bar v_2},...,{\bar v_m}\}.$

Put
$$
\gamma_3=\min\{{{(}}\Delta/2{{)}}(\hat{h}): h\in{\cal H}_1\}.
$$
Note $H^{\ddag}({\cal U}_0)\subset U_0(C_2)/CU(C_2).$
Choose a finite subset ${\cal H}_3\subset (C_2)_{s.a.}$  and $\sigma>0$ which have the following property:
for any two unital \hm s $h_1, h_2: C_2\to D$ (for any unital \CA\, $D$ of stable rank one),
if
\beq\label{B0=B1-n3}
|\tau\circ h_1(g)-\tau\circ h_2(g)|<\sigma\tforal g\in {\cal H}_3,
\eneq
then
\beq\label{B0=B1-n4}
{\rm dist}(h_1^{\ddag}({\bar v}), h_2^{\ddag}({\bar v}))<\gamma_2/8
\eneq
for all ${\bar v}\in H^{\ddag}( {\cal U}_0)\subset U_0(C_2)/CU(C_2).$
%where we identify ${\bar v}$ with $\overline{v\oplus (q_2-q_1)}\subset U_0(C_2)/CU(C_2).$
\Wlog, we may assume that $\|h\|\le 1$ for all $h\in {\cal H}_1\cup {\cal H}_2\cup {\cal H}_3.$

Let $\kappa: {\rm Ell}(A_2)\to {\rm Ell}(B)$  be the above identification. So
$\kappa\circ [j_2]\in KK_e(C_2, B)^{++}.$
It follows from \ref{ExtTrace} that there exists a unital \hm\, $\phi: C_2\to B$
such that
\beq\label{B0=B1-4}
[\phi]=\kappa\circ [j_2]\andeqn\\\label{B0=B1-4+}
 |\tau(\phi(h))-\gamma(\tau)(h)|<\min\{\gamma_1, \gamma_2,\gamma_3,\sigma\}/8
\eneq
for all $h\in  H({\cal H}_1)\cup H({\cal H}_2)\cup {\cal H}_3,$ where $\gamma: T(B)\to T_f(C_1)$ is induced
by $\kappa$ and the embedding $j_1.$
In particular, $\phi_{*1}$ is injective on $H_{*1}(\Z^m).$

%Let $p_0=\phi(q_1).$ Since $B$ is simple, there are
%$x_1,x_2,...,x_l\in B$ such that
%\beq\label{B0=B1-5-1}
%\sum_{i=1}^lx_i^*p_0x_i=1_B.
%\eneq

Since $C_2$ is semi-projective, \wilog, we may assume that $\phi(C_2)\subset B_1.$
%Moreover, by considering (\ref{B0=B1-5-1}) and by choosing $B_n$ for sufficiently large $n,$
%\wilog, we may assume that $p_0$ is {\it full} in $B_1.$
We may also assume that $B_1=B_{1,0}\oplus B_{1,1}$
and
\beq\label{B0=B1-5}
\tau(1_{B_{1,0}})<\min\{\tau(\kappa([e_2]))/4, \gamma_1/8, \gamma_2/8, \gamma_3/8, \sigma/8\}\tforal \tau\in T(B).
\eneq
Furthermore, identifying ${\rm Ell}(A_2)$ and $ {\rm Ell}(B)$,  we have  $(\imath_{1, \infty} \circ \phi)_{*1}=\phi_{*1}$, which is is injective on $H_{*1}(\Z^m)$. Consequently one has
% by replacing $B_1$ by $B_n$ for sufficiently large $n,$ we mayassume, \wilog,
that $(\imath_{1, \infty})_{*1}$ is injective on $\phi_{*1}(H_{*1}(\Z^m)).$

Let $G_1=H^{\ddag}(J_c^{(1)}(\Z^m))\subset U(C_2)/CU(C_2)$ and
$G_0=H^{\ddag}(J_c^{(1)}(G_{00})).$ Note, by the construction, $G_0\subset U_0(C_2)/CU(C_2).$
Since $\kappa|_{K_1(C)}$ is an isomorphism and
$(\imath_{1, \infty})_{*1}$ is injective on $\phi_{*1}(H_{*1}(\Z^m)),$ $\phi^{\ddag}|_{G_1}$ is injective
Let ${\cal H}_4=P(\phi({\cal H}_1\cup {\cal H}_2\cup {\cal H}_3)),$ where $P: B_1\to B_{1,1}$  is the projection.

Let $e_3{{\in}} A$ be a projection such that $[e_3]=\kappa^{-1}([\imath_{1, \infty})(1_{B_{1,1}})]).$
It follows from the second part of \ref{ExtTrace}
%\ref{CorCinjj}
that there is a unital monomorphism
$\psi_1: B_{1,1}\to e_3Ae_3$ such that
\beq\label{B0=B1-6}
\hspace{-0.6in}&&[\psi_1]=\kappa^{-1}\circ [(\imath_{1,\infty})|_{B_{1,1}}]\andeqn\\\label{B0=B1-6+}
\hspace{-0.6in}&&|\tau\circ \psi_1(g)-\gamma'(\tau)(\imath_{1,\infty}(g))|<\min\{\tau(([e_3]))/4, \gamma_1/8, \gamma_2/8, \gamma_3/8, \sigma/8\}\rforal g\in {\cal H}_4
\eneq
for all $\tau\in T(A),$ where $\gamma': T(A)\to T_f(B_{1,1})$ induced by $\kappa^{-1}$ and $\imath_{1, \infty}.$

Write $B_{1,0}=B_{1,0,1}\oplus B_{1,0,2},$ where $B_{1,0,1}$ is a finite direct sum of circle algebras and
$K_1(B_{1,0,2})$ is finite.
Since $B\cong B\otimes Q,$  we may assume that $(\imath_{1,\infty})_{*1}|_{K_1(B_{1,0.2})}=\{0\}.$

Since $A_2\cong A_2\otimes Q,$ by part  (iii) and (iv) of \ref{preBot2},
there exists a unital \hm\, $\psi_2: B_{1,0,2}\to e_4A_2e_4$
such that $[\psi_{{2}}]=\kappa^{-1}\circ [(\imath_{1,\infty})_{B_{1,0}}],$ where $e_4$ is a projection
orthogonal to $e_2$ and $[e_4]=\kappa^{-1}\circ [\imath_{1,\infty}](1_{B_{1,0,2}}).$ As in the proof of \ref{preBot2}, $(\psi_2)_{*1}=0.$
Let $\psi_3: B_{1,1}\oplus B_{1,0,2}\to (e_3+e_4)A_2(e_3+e_4)$ by $\psi_3={{\psi}}_1\oplus \psi_2.$

Let $P_1: B_1\to B_{1,0,1}.$ Then, since $(\imath_{1, \infty})_{*1}|_{K_1(B_{1,0,2})}=\{0\},$ $(P_1)_{*1}|_{\phi_{*1}(H_{*1}(\Z^m))}$
is injective. Also $P_1^{\ddag}$ is injective on $\phi^{\ddag}\circ H^{\ddag}(J_c^{(1)}(\Z^m)).$
Put $G_1'=P_1^{\ddag}\circ \phi^{\ddag}\circ H^{\ddag}(J_c^{(1)}(\Z^m)) { {\subset U(B_{1,0,1})/CU(B_{1,0,1})}}.$ Then $G_1'\cong \Z^m.$

It follows from \ref{preBot2} that there is a unital \hm\, $\psi_4': B_{1,0,1}\to (1-e_3-e_4)A_2(1-e_3-e_4)$
such that $[\psi_4']=\kappa^{-1}\circ [\imath_{1, \infty}|_{B_{1,0,1}}].$
Let
\beq\label{B0=B1-15}
z_i= P_1^{\ddag}\circ \phi^{\ddag}\circ H^{\ddag}({\bar v}_i)\andeqn
\xi_i=\psi_3^{\ddag}\circ (1_{B_1}-P_1)^{\ddag}\circ \phi^{\ddag}\circ H^{\ddag}({\bar v}_i),
\eneq
$i=1,2,...,m.$
It should be noted that, since $(\psi_1)_{*1}=(\psi_2)_{*1}=0,$ $\xi_i\in U_0(A_2)/CU(A_2),$ $i=1,2,...,m.$
Moreover, since
\beq\label{B0=B1-16}
(\psi_4')_{*1}\circ (P_1)_{*1}\circ \phi_{*1}(x)=j_{*1}(x)\tforal x\in \Z^m\subset K_1(C_1),
\eneq
\beq\label{B0=B1-17}
\pi({\bar v}_i)\pi((\psi_4')^{\ddag}(z_i))^{-1}=0\,\,\, {\rm in}\,\,\, K_1(A).
\eneq

Define a \hm\, $\lambda: G_1'\to U(A_2)/CU(A_2)$ by
\beq\label{B0=B1-18}
\lambda(z_i)={\bar v}_i((\psi_4')^{\ddag}(z_i)\xi_i)^{-1},\,\,\,i=1,2,...,m.
\eneq
Note that $\lambda(z_i)\in U_0(A)/CU(A),$  $i=1,2,...,m.$
By \ref{UCUdiv},  $U_0(A)/CU(A)$ is divisible. There exists  a \hm\, $\bar{\lambda}: U(B_{1,0,1})/CU(B_{1,0,1})
\to U(A_2)/CU(A_2)$
such that ${\bar \lambda}|_{G_1'}=\lambda.$
Define a \hm\, $\lambda_1:  U(B_{1,0,1})/CU(B_{1,0,1})\to U(A_2)/CU(A_2)$ by
$\lambda_1(x)=(\psi_4')^{\ddag}(x){\bar \lambda}(x)$ for all $x\in U(B_{1,0,1})/CU(B_{1,0,1}).$
By \ref{UCUiso}, the  \hm\,
$$
U((1-e_3-e_4)A(1-e_3-e_4))/CU((1-e_3-e_4)A_2(1-e_3-e_4))\to
U(A_2)/CU(A_2)
$$
is an isomorphism.  Since $B_{1,0,1}$ is a circle algebra, one easily obtains a unital \hm\,
$$
\psi_4: B_{1,0,1}\to (1-e_3-e_4)A_2(1-e_3-e_4)
$$
such that
\beq\label{B0=B1-19}
[\psi_4]=[\psi_4']=\kappa^{-1}\circ [\imath_{1, \infty}|_{B_{1,0,1}}]\andeqn \psi_4^{\ddag}=\lambda_1.
\eneq
Define $\psi: B_1\to A_2$ by $\psi=\psi_3\oplus \psi_4.$
Then, we have
\beq\label{B0=B1-20}
[\psi]=\kappa^{-1}\circ [\imath_{1, \infty}].
\eneq
Since $[\imath_{1, \infty}\circ \phi]=\kappa\circ [j],$
we compute  that
\beq\label{B0=B1-21}
[\psi\circ \phi\circ H]=\kappa^{-1}\circ [\imath_{1, \infty}]\circ [\phi\circ H]=[j\circ H]\,\,\,{{\mbox{in}~~}}KK(C_1,A).
\eneq
Put $C_1'=H(C_1)$ and $\psi'=\psi\circ \phi|_{C_1'}.$ Then $\psi': C_1'\to A_2$ is a monomorphism.
%By replacing $\psi'$ by ${\rm Ad} \, w\circ \psi'$ for some suitable unitary $w,$ we may assume
%that $\psi'(1_{C_1'})=q_2.$
We have
\beq\label{B0=B1-22}
[\psi']=[j|_{C_1'}].
\eneq	
By (\ref{B0=B1-4+}), (\ref{B0=B1-5}) and (\ref{B0=B1-6+}),
\beq\label{B0=B1-23}
|\tau\circ \psi\circ \phi(g)-\tau(g)|<\min\{\sigma, \gamma_1, \gamma_2, \gamma_3,\sigma\}\rforal g\in H({\cal H}_1)\cup  H({\cal H}_2)\cup {\cal H}_3.
\eneq
In particular,
\beq\label{B0=B1-24}
|\tau\circ \psi'(g)-\tau(j(g))|<\min\{ \gamma_1, \gamma_2, \gamma_4\}\rforal g\in H({\cal H}_1)\cup H({\cal H}_2).
\eneq
We then compute that
\beq\label{B0=B1-24+}
\tau\circ \psi'(h)\ge \Delta(\hat{h})/2\rforal h\in {\cal H}_1.
\eneq

By  (\ref{B0=B1-18}) and the definition of $\lambda_1,$ we have
\beq\label{B0=B1-25}
(\psi')^{\ddag}({\bar v}_i)={\bar v}_i,\,\,\,i=1,2,...,m.
\eneq
By the choice of ${\cal H}_3$ and $\sigma,$
we also have
\beq\label{B0=B1-26}
{\rm dist}((\psi')^{\ddag}({\bar v}), (j\circ H)^{\ddag}({\bar v}))<\sigma_2\rforal {\bar v}\in {\cal U}_0.
\eneq
It follows from \ref{UniCtoA} and \ref{RemUniCtoA} that there is a unitary $U\in q_2Aq_2$ such that
\beq\label{B0=B1-27}
\|{\rm Ad}\, U\circ \psi'(x)-x\|<\ep/16\rforal x\in {\cal F}_1.
\eneq
%Now let $p'=\p1_{B_{1,1}}.$ It is a projection. It is a full projection of $B_{1,1}.$
%Put $C_0=p'B_{1,1}p'.$ Then, by \ref{cut-full-pj}, $C_0\in {\cal C}_0.$
Let $C_3={\rm Ad}\, U\circ \psi(B_{1,1})$ and $p=1_{C_3}.$
%Then $C_3\in {\cal C}_0$
and  $p={\rm Ad}\, U\circ \psi(1_{B_{1,1}}).$
%In particular, $p\le q_2=\psi(p').$
For any $y\in {\cal F}_1,$ $\phi(y)1_{B_{1,1}}=1_{B_{1,1,}}\phi(y).$
Therefore ${\rm Ad}\, U\circ \psi'(y)p=p{\rm Ad}\, U\circ \psi\circ \phi(y)=p{\rm Ad}\, U\circ \psi'(y).$ Hence,
\beq\label{B0=B1-28}
\|py-yp\|\le \|py-p{\rm Ad}\,U\circ \psi'(y)\|+\|p{\rm Ad}\,U\circ \psi'(y)-yp\|<\ep/8+\ep/8=\ep/4
\eneq
for all $y\in {\cal F}_1.$
Thus, combining with  (\ref{B0=B1-n1})
 and (\ref{K1inj-2}),
%and (\ref{K1inj-2+}),
\beq\label{B0=B1-29}
\|px-xp\|=\|pq_2x-xq_2p\|<2\ep/16+\|pq_2xq_2-q_2xq_2p\|<\ep\rforal x\in {\cal F}.
\eneq
Let $x\in {\cal F}.$ Choose $y\in {\cal F}_1$ such that $\|q_2xq_2-q_2yq_2\|<\ep/16.$ Then, by (\ref{B0=B1-27}),
\beq\label{B0=B1-30}
\|pxp-p{\rm Ad}\,U\circ \psi'(y)p\| &\le & \|pxp-pq_2xq_2p\|\\
&&+\|pq_2xq_2p-pyp\|+\|pyp-p{\rm Ad}\, U\circ \psi'(y)p\|\\
&&<\ep/16+\ep/16=\ep/8.
\eneq	
However, $p{\rm Ad}\, U\circ \psi'(y)p={\rm Ad}\, U\circ(\psi(\phi(q_2)\phi(y)\phi(q_2))\in C_3$ for all $y\in {\cal F}_1.$
Therefore
\beq\label{B0=B1-31}
pxp\in_{\ep} C_3.
\eneq
We then estimate that
\beq\label{B0=B1-32}
[1-p]\le [1-q_2]+[\psi(1_{B_{1,0}})]\le [e_0\oplus e_1\oplus e_2]\le [a].
\eneq
Since $B_{1,1}\in {\cal C}_0,$ by applying \ref{subapprox},  $A\in {\cal B}_0.$
\end{proof}

%{\color{Green} In \ref{B0=B1tU}, one does not need $U$ to have infinite type. In fact, any UHF algebra works.}
\begin{cor}\label{B0=B1tU}
If $A$ is a unital separable amenable simple \CA\, such that $A\otimes Q\in {\cal B}_1,$  then, for any infinite dimensional UHF-algebra $U,$ $A\otimes U\in {\cal B}_0.$
\end{cor}

\begin{proof}
It follows from \ref{B0=B1} that $A\otimes Q\in {\cal B}_0.$ Then, by \ref{subapprox} and by \cite{LS}, $A\otimes U\in {\cal B}_0$ for
every UHF-algebra of infinite type.
\end{proof}

\section{KK-attainability of the C*-algebras in $\mathcal B_0$}

In the following, let us show an existence theorem for the maps from an algebra in the class $\mathcal B_0$ to a C*-algebra in Theorem \ref{RangT}. The procedure is similar to that of Section 2 of \cite{Lnduke}, and roughly, we will construct a map factors through the C*-subalgebras (in $\mathcal C_0$) of the given C*-algebra in $\mathcal B_0$, and also require this map to carry the given KL-element. But since positive cone of the $K_0$-groups of a C*-algebra in $\mathcal C_0$ in general is not free, extra work has to be done to take care of this issue.

\begin{NN}\label{201} Let us proceed as that of Section 2 of \cite{Lnduke}. Let $A\in\mathcal B_0$, and assume that $A$ has the (SP) property. By Lemma \ref{MF}, the C*-algebra $A$ can be embedded as a C*-subalgebra of $\prod M_{n_k}/\bigoplus M_{n_k}$ for some $(n_k)$, and therefore $A$ is MF in the sense of Blackadar and Kirchberg (Theorem 3.2.2 of \cite{BK-inductive}). Since $A$ is assumed to be amenable, by Theorem 5.2.2 of \cite{BK-inductive}, the C*-algebra $A$ is strong NF, and hence, by Proposition 6.1.6 of \cite{BK-inductive}, there is an increasing family of RFD sub-C*-algebras $\{A_n\}$ such that their union is dense in $A$.

Let $\{x_1, x_2,..., x_n,...\}$ be a dense sequence of elements in the unit ball of $A$. Let $\mathcal P_0\subseteq \mathrm{M}_\infty(A)$ be a finite subset of projections. We assume that $x_1$ and $\mathcal P_0\subseteq \mathrm{M}_\infty(A_1)$. Consider a finite subset $\mathcal F_0\subset A_1\subset A$ with $\{1, x_1\}\subseteq \mathcal F_0$, $\delta_0>0$, and a homomorphism $h_0$ from $A_1$ to a finite-dimensional C*-algebra $F_0$ which is non-zero on $\mathcal F_0$. Since $A$ is assumed to have the (SP) property,
% {and is quasidiaogonal, {\color{Green} (one does not need the quasidiagnality in Lemma 2.1, it was used to show $A$ is strong NF, and hence $\{A_n\}$ exists)}}
by Lemma 2.1 of \cite{Niu}, there is a non-zero homomorphism $h': F_0\to A$ such that
\begin{enumerate}
\item $||e_0x-xe_0||<\delta_0/256$\ \ and
\item $|| h'\circ h_0(x)-e_0xe_0 ||<\delta_0/256$
%\end{enumerate}
for all $x\in\mathcal F_0$, where $e_0=h'(1)$.

%{\bf Which reference?  This looks familiar and a combination of TAC and a result of BK---L. ----Just (SP) property and some functional calculus.}

Since $F_0$ has finite dimension, it follows from Arveson's Extension Theorem that the homomorphism $h_0: A_1 \to F_0$ can be extended to a
\morp\,  from $A$ to $F_0$, and let us still denote it by $h_0$.

   Put $H=h'\circ h_0: A\to A.$  Note that $e_0=H(1)$. Since the hereditary C*-subalgebra $(1-e_0)A(1-e_0)$ is in the class $\mathcal B_0$ again, there is a projection $q_1'\leq 1-e_0$ and a C*-subalgebra $S_1'\in \mathcal C_0$ with $1_{S_1'}=q_1'$ such that
%\begin{enumerate}
\item $|| q'_1x-xq'_1 ||<\delta_0/256$ for any $x\in (1-e_0){\cal F}_0(1-e_0),$
\item $\textrm{dist}(q_1'xq'_1, S_1')<\delta_0/256$ for any $x\in (1-e_0)\mathcal F_0(1-e_0)$, and
\item ${\tau(1-e_0-q_1')} < 1/16$ for any tracial state $\tau$ on $A$.
%\end{enumerate}
Put $q_1=q_1'+e_0$ and $S_1=S'_1\oplus h'(F_0)$. One has
%\begin{enumerate}
\item $|| q_1x-xq_1 ||<\delta_0/{64}$ for any $x\in \{ab: a,\, b\in{ \cal F}_0\},$
\item $\textrm{dist}(q_1xq_1, S_1)<\delta_0/64$ for any $x\in\ \{ab: a,\, b\in {\cal  F}_0\}$, and
\item $\tau(1-q_1)=\tau(1-q'_1-e_0)  <{ 1/16}$ for any tracial state $\tau$ on $A$.
%\end{enumerate}

{Let ${\cal F}_0'\subset S_1$ be a finite subset such that ${\rm dist}(q_1yq_1, {\cal F}_0')<\dt_0/16$
for all $y\in \{ab: a, b\in {\cal F}_{ 0}\}.$ }
Let $\mathcal G_1$ be a finite generating set of $S_1$ which is in the unit ball.

Since $S_1$ is amenable, there is a \morp\,  $L_0': q_1Aq_1\to S_1$ such that
\beq\label{June13-1}
\|L_0'(s)-s\|<\dt_0/256\tforal s\in   {\cal G}_1\cup {\cal F}_0'.
\eneq
Set $L_0(a)=L_0'(q_1aq_1)$ for any $a\in A$. Then $L_0$ is a completely positive contraction from $A$ to $S_1$
such that
\beq\label{June13-2}
\|L_0(s)-s\|<\dt_0/256\rforal s\in {\cal G}_1\cup{\cal F}_0'.
\eneq
%Note that  the map $sen so L'_0$ can be chosen
We estimate that $L_0$ is $\{x_1\}$-$\delta_0/16$-multiplicative.
%{by (6)}.
{Since $S_1$ is semi-projective, there is $\delta_1'>0$  and a finite subset ${\cal G}_1'\supset {\cal G}_1$ in the unit ball} such that for any $\mathcal G_1'$-$\delta_1'$-multiplicative \morp\,  $L$ from $S_1$ to {any} C*-algebra, there is a homomorphism $h$ from $S_1$ to that C*-algebra such that
$$
|| L(a)-h(a) ||<\delta_0/2 \quad\textrm{for any}\ a\in \{L_0(x_1)\}\cup \mathcal G_1.
$$
Set $\mathcal F_1=\{x_2\}\cup\mathcal F_0\cup {\mathcal G_1'}$.

%Since $S_1$ is semi-projective, there is $\delta_1'>0$ such that for any $\mathcal G_1$-$\delta_1'$-multiplicative \morp\,  $L$ from $S_1$ to a C*-algebra, there is a homomorphism $h$ from $S_1$ to that C*-algebra such that
%$$
%|| L(a)-h(a) ||<\delta_0/2 \quad\textrm{for any}\ a\in \{L_0(x_1)\}\cup \mathcal G_1.
%$$
Set $\delta_1=\min\{\delta_1', \delta_0/256\}$.  Also there is a projection $q_2\in A$ and a C*-subalgebra  $S_2\in\mathcal C_0$ with $1_{S_2}=q_2$ such that
%\begin{enumerate}
\item  $|| q_2x-xq_2 ||<{\delta_1/16}$ for any $x\in {\{ab: a,b\in {\cal F}_1\}},$
\item  $\textrm{dist}(q_2xq_2, S_2)<{\delta_1/16}$ for any $x\in {\{ab: a,b\in {\cal F}_1\}}$, and
\item   $\tau(1-q_2)\leq 1/32$ for any tracial state $\tau$ on $A$.
%\end{enumerate}
%Moreover, it follows from the semi-projectivity of $S_1$ that there is a homomorphism
%$h_{1, 2}: S_1\to S_2$
%such that
%$$\|L_2(a) - h_2(a)\|<\delta_0/16,\quad a\in \{L_1(x_1)\}\cup\mathcal G_1.$$

{Let ${\cal F}_1'\subset S_2$ be a finite subset such that ${\rm dist}(q_2yq_2, {\cal F}_1')<\dt_1/16$
for all $y\in \{ab: a, b\in {\cal F}_1\}.$ }
Let ${\cal G}_2$ be a finite generating set of $S_2$ which is in the unit ball.

Therefore there  exists a completely positive contraction $L_1': q_2Aq_2\to S_2$
such that
\begin{equation*}
\|L_1'(a)-a\|<\dt_1/16 \rforal a\in {\cal G}_2{\cup {\cal F}_1'}.
\end{equation*}
 Define  $L_1(a)=L_1'(q_2aq_2)$ for all $a\in A.$ It is a completely positive contraction from $A$ to $S_2$
 %with the restriction to $S_2$ being identity,
 such that
 \begin{equation*}
 \|L_1(a)-a\|<\dt_1/16\rforal a\in {\cal G}_2{\cup {\cal F}_1'},
 \end{equation*}
 %Note that $L_1'$ can be chosen so that
 {We estimate, from (9), the choice of ${\cal F}_1'$  and the above inequality, that} $L_1$ is $\mathcal F_1$-$\delta_1/2$-multiplicative.  Therefore, there is a homomorphism $h_{1}: S_1\to S_2$ such that $$|| L_1(a)-h_{1}(a) ||<\delta_0/2 \quad\textrm{for any}\ a\in \{L_0(x_1)\}\cup \mathcal G_1.$$
 { Since $S_2$ is semi-projective, there is  $\delta_2'>0$  and a finite subset ${\cal G}_2'\supset {\cal G}_2$ in the unit ball} such that for any $\mathcal G_2'$-$\delta_2'$-multiplicative \morp\,  $L$ from $S_2$ to {any} C*-algebra, there is a homomorphism $h'$ from $S_2$ to that C*-algebra such that
$$
\| L(a)-h'(a) \|<\delta_1/2,\quad a\in \{L_1(x_1), L_1(x_2)\}\cup \mathcal G_2.
$$
Put
${\cal F}_2={\{x_3\}\cup {\cal F}_1\cup {\cal F}_1'\cup {\cal G}_1'\cup {\cal G}_2'}.$

%  Since $S_2$ is semi-projective, there is $\delta_2'>0$ such that for any $\mathcal G_2$-$\delta_2'$-multiplicative \morp\,  $L$ from $S_2$ to a C*-algebra, there is a homomorphism $h'$ from $S_2$ to that C*-algebra such that
%$$
%\| L(a)-h'(a) \|<\delta_1,\quad a\in \{L_1(x_1), L_1(x_2)\}\cup \mathcal G_2.
%$$
Set $\delta_2=\min\{\delta_2', \delta_1/256\}$. Also there is a projection $q_3\in A$ and a C*-subalgebra  $S_3\in\mathcal C_0$ with $1_{S_3}=q_3$ such that
%\begin{enumerate}
\item $|| q_3x-xq_3 ||<\delta_2/16$ for any $x\in  \{ab: a,b\in {\cal F}_2\},$
\item $\textrm{dist}(q_3xq_3, S_3)<\delta_2/16$ for any $x\in \{ab: a,\,b\in {\cal F}_2\}$, and
\item $\tau(1-q_2)\leq 1/64$ for any tracial state $\tau$ on $A$.
%\end{enumerate}

{Let ${\cal F}_2'\subset S_3$ be a finite subset such that ${\rm dist}(q_3yq_3, {\cal F}_2')<\dt_3/16$
for all $y\in \{ab: a, b\in {\cal F}_2\}.$ }
Let ${\cal G}_3$ be a finite generating set of $S_3$ which is in the unit ball.

There then exists a completely positive contraction $L_2': q_3Aq_3\to S_3$
such that
\begin{equation*}
\|L_2'(a)-a\|<\dt_2/16\rforal a\in {\cal G}_3\cup {\cal F}_2'.
\end{equation*}
 Define  $L_2(a)=L_2'(q_3aq_3)$ for all $a\in A.$ It is a completely positive contraction from $A$ to $S_2$
 %with the restriction to $S_2$ being identity,
 such that
 \begin{equation*}
 \|L_2(a)-a\|<\dt_2/16\rforal a\in  {\cal G}_3\cup {\cal F}_2'.
 \end{equation*}
% Note that $L_2'$ can be chosen
 As estimated above, $L_2$ is $\mathcal F_2$-$\delta_2/2$-multiplicative. Since $S_3$ is semi-projective, there is  $\delta_3'>0$  and a finite subset ${\cal G}_3'\supset {\cal G}_3$ in the unit ball such that for any $\mathcal G_3'$-$\delta_3'$-multiplicative \morp\,  $L$ from $S_3$ to a C*-algebra, there is a homomorphism $h'$ from $S_3$ to that C*-algebra such that
$$
\| L(a)-h'(a) \|<\delta_2/2,\quad a\in \{L_2(x_1), L_2(x_2), L_2(x_3)\}\cup \mathcal G_3.
$$
Put
${\cal F}_3={\{x_4\}\cup {\cal F}_2\cup {\cal F}_2'\cup {\cal G}_2'\cup {\cal G}_3'}.$

{On the other hand,} there is a homomorphism %are  homomorphisms $h_{1,2}: S_1\to S_3$
 and $h_{2}: S_2\to S_3$ such that
 $$\|L_2(a)-h_{2}(a)\|<\delta_1,\quad a\in \{L_1(x_1), L_1(x_2)\}\cup \mathcal G_2.$$

 %Note that if define $h_{1, 3}:=h_{2, 3} \circ h_{1, 2}$, then, for any $a\in \{L_0(x_1)\}\cup\mathcal G_1$, one has that
 %
% \begin{eqnarray*}
 %\| L_2 (a) - h_{1, 3}(a) \| & = & \| L_2 (a) - h_{2, 3} (h_{1, 2}(a)) \| \\
% & = & \| L_2 (a) - h_{2, 3} (L_1(a)) \| + \delta_0/2 \\
% & = & \| L_2 (a) - L_2 (L_1(a)) \| + \delta_0/2 +\delta_1 \\
% & = &
% \end{eqnarray*}

% \beq\label{June15-1}
% || L_2(a)-h_{1,3}(a) ||<\delta_0/256\rforal  a\in \{L_0(x_1)\}\cup \mathcal G_1\andeqn\\
% \|L_2(a)-h_{2,3}(a)\|<\dt_1
% \rforal a\in \{L_1(x_1), L_1(x_2), L_1(x_3)\}\cup \mathcal G_2.
% \eneq

Repeating this construction, one obtains a sequence of finite subsets $\mathcal F_0, \mathcal F_1, ...,$ with dense union in the unit ball of $A$,  a decreasing sequence of positive numbers $\{\delta_n\}$ {with $\sum_{n=1}^{\infty}\dt_n<\infty,$} a sequence of projections $\{q_n\}\subset A$, a sequence of C*-subalgebras $S_n\in\mathcal C_0$ with $1_{S_n}=q_n$
with finite generating subsets ${\cal G}_n\subset S_n,$ and a sequence of homomorphisms $h_{n}:S_n\to S_{n+1}$
{and a sequence of  ${\cal F}_n$-$\dt_n/2$-multiplicative \morp s $L_n: A\to S_{n+1}$} such that
%\begin{enumerate}
\item $|| q_nx-xq_n ||<\delta_{n-1}/16$ for all $x\in\mathcal F_n$.
\item dist$(q_nx_iq_n,S_n)<\delta_{n-1}/16$, $i=1, ... , n$.
\item $\tau(1-q_n)<1/2^{n+1}$ for all tracial states $\tau$ on $A$.
\item $\mathcal G_n\subset\mathcal F_{n+1}$, where $\mathcal G_n$ is a finite set of generators for $S_n$.
\item $|| L_{n+1}(a)-h_{n+1}(a)||<\dt_n$ for
      all $a\in L_n(\mathcal F_n) \cup \mathcal G_n$. %$i=1,2,...,n,$
      %where $L_{n+1}:A\to S_{n}$ is a
      %contractive positive linear map
      %with
      and
      $\|L_{n}(a)-a\|<\dt_n$ for all $a\in {\cal G}_{n+1}.$
%and where $h_{i,n+1}=h_{n+1}\circ h_n\circ\cdots \circ h_i,$ $i=1,2,...,n.$
\end{enumerate}

\end{NN}

\begin{NN}\label{June13-constr}
{\rm
Let $\Psi_n: A\to(1-q_{{n+1}})A(1-q_{{n+1}})$ denote the cut-down map sending $a$ to $(1-q_{{n+1}})a(1-q_{{n+1}})$, and let $J_n:A\to A$ denote the map sending $a$ to $\Psi_n(a)\oplus L_n(a)$. Note that $\Psi_n$ and $J_n$ are $\mathcal F_n$-$\delta_n/2$-multiplicative. Set $J_{m,n}=J_{n-1}\circ\cdots\circ J_m$ and
$h_{m,n}=h_{n-1}\circ\cdots\circ h_{m}: S_{m}\to S_{n}$. Note that $J_{m,n}$ is $\mathcal F_m$-$\delta_m$-multiplicative. We also use $L_n,
\Psi_n, J_n, J_{m,n}, h_m, $ and $h_{m,n}$ for their extensions on a matrix algebra over $A$.

%{{Write ${\tilde \lambda}_{k,k+j}: S_k\to S_{k+j}$ for the map
%$L_{k+j}\circ \Psi_{k+j-1}|_{S_k}$ ($j\ge 1$).   Define, for $n>j,$
%${\tilde h}_{k+j,k+n}: S_k\to S_{k+n}$ by ${\tilde h}_{k+j, k+n}=h_{k+j+1, k+n-1}\circ {\tilde \lambda}_{k,k+j}.$}}
}
\end{NN}

Using the same argument as that of Lemma 2.7 of \cite{Lnduke}, one has the following lemma.

\begin{lem}[Lemma 2.7 of \cite{Lnduke}]\label{traces}
Let $\mathcal P\subset \text{M}_k(A)$ be a finite set of projections. Assume that $\mathcal F_{\bf 1}$ is sufficiently large and $\delta_0$ is
sufficiently small such that $[L_{n}\circ J_{1,n}]|_{\mathcal P}$ and $[L_{n}\circ J_{1,n}]|_{G_0}$ are well defined, where $G_0$ is the subgroup generated by $\mathcal P$. Then
$$\lim_{n\to\infty}\sup_{\tau\in T(A)}|\tau([\iota_{n+1}\circ L_{n}\circ J_{1,n}]([p]))-\tau([p])|=0$$
for any $p\in\mathcal P.$

Furthermore, for any projection ${{p \in {\cal P}}}$  and $k\ge 1,$ we have
$$
{{|\tau(h_{k, k+n+1}\circ [L_{k-1}]([p]))-\tau(h_{k, k+n}\circ [L_{k-1}]([p])|<(1/2)^{n+k}}}
%\lim_{n\to\infty}\sup_{\tau\in T(A)}|\tau(\imath_n\circ {{L_n}}(q))-\tau(q)|=0
$$
for all $\tau\in T(A)$ and
$$
{{\lim_{n\to\infty} \tau(h_{k,k+n}\circ [L_{k-1}]([p])\ge (1-\sum_{i=1}^n 1/2^{i+k})\tau([L_{k-1}]([p]))>0}}
$$
for all $p\in {\cal P}$ and $\tau\in T(A).$
\end{lem}

\begin{rem}\label{Kstate}
Since $A$ is stably finite and assumed to be amenable, therefore exact, any positive state of $K_0(A)$ is the restriction of a tracial state of
$A$ (\cite{Blatrace} and \cite{Haagtrace}).  Thus, the lemma above still holds if one replaces the trace $\tau$ by any positive state $\tau_0$ on $K_0(A)$.
\end{rem}

\begin{NN}\label{Constr-1508}
Fix a finite subset ${\cal P}$ of projections of $M_r(A)$ (for some $r\ge 1$)
and an integer $N\ge 1$ such that $[L_{N+i}]|_{\cal P},$
$[J_{N+i}]|_{\cal P}$ and $[\Psi_{N+i}]|_{\cal P}$ are all well defined.
Keep notation in \ref{June13-constr}. Then, on ${\cal P},$
\beq\nonumber
[L_{N+1}\circ J_N]&=&[L_{N+1}\circ L_N]\oplus [L_{N+1}\circ \Psi_N]\\
&=& [h_{N+1}\circ L_N]\oplus [L_{N+1}\circ  \Psi_N], \andeqn
\eneq
\beq\nonumber
\hspace{-0.6in}[L_{N+2}\circ J_{N, N+2}]&=&[L_{N+2}\circ L_{N+1}\circ J_N]\oplus [L_{N+2}\circ \Psi_{N+1}\circ J_N]\\
&=&[L_{N+2}\circ L_{N+1}\circ L_N]\oplus [L_{N+2}\circ L_{N+1}\circ \Psi_N]\\
&&\oplus [L_{N+2}\circ  \Psi_{N+1}\circ J_N]\\\hspace{-0.1in}&=&[h_{N+1, N+3}]\circ [L_N]\oplus [L_{N+2}\circ L_{N+1}\circ \Psi_N]\oplus
[L_{N+2}\circ  \Psi_{N+1}\circ J_N].\\
\eneq
Moreover,  on ${\cal P},$
\beq\nonumber
[L_{N+n}\circ J_{N, N+n-1}] &=& [h_{N+1, N+n+1}]\circ [L_N]\oplus [L_{N+n}\circ \Psi_{N+n-1}\circ J_{N, N+n-1}]\\
&&\oplus [L_{N+n}\circ L_{N+n-1}\circ \Psi_{N+n-2}\circ J_{N, N+n-2}]\\
&&\oplus  [L_{N+n}\circ L_{N+n-1}\circ L_{N+n-2}\circ \Psi_{N+n-3}\circ J_{N, N+n-3}]\\
&&\oplus \cdots { \oplus[L_{N+n}\circ L_{N+n-1}\circ\cdots \circ L_{N+2}\circ \Psi_{N+1}\circ J_N}]\\
&& \oplus [L_{N+n}\circ L_{N+n-1}\circ\cdots \circ L_{N+1}\circ \Psi_N].
\eneq

Set $\psi^N_{N}=L_N,$ $\psi^N_{N+1}=L_{N+1}\circ \Psi_N,$
$\psi^N_{N+2}=L_{N+2}\circ \Psi_{N+1},..., \psi^N_{N+n}=L_{N+n}\circ \Psi_{N+n-1},$
$n=1,2,....$

\end{NN}

\begin{NN}\label{Constr-Sn}
%\rm
For each $S_n$, since the abelian group $K_0(S_n)$ is finitely generated and torsion free, there is a set of free generators $\{e^{n}_1, e^{n}_2, ..., e^{n}_{l_n}\}\subseteq K_0(S_n)$. By Theorem \ref{FG-Ratn}, the positive cone of the $K_0(S_n)$ is finitely generated; denote a set of generators by $\{s^n_1, s^n_2, ..., s^n_{r_n}\}\subseteq K_0^+(S_n)$. Then there is an
$r_n\times l_n$ integer-valued matrix $R'_n$ such that
$$\vec{r}_n=R'_n\vec{e}_n,$$ where $\vec{r}_n=(s^n_1, s^n_2, ..., s^n_{r_n})^{\mathrm T}$ and $\vec{e}_n=(e^{n}_1, e^{n}_2, ..., e^{n}_{l_n})^{\mathrm T}$. In particular, for any ordered group $H$, and any elements $h_1, h_2, ..., h_{l_n}\in H$, the map $e^n_i\mapsto h_i$, $i=1, ..., l_n$ induces an abelian-group homomorphism $\phi: K_0(S_n)$ to $H$, and the map $\phi$ is positive (or strictly positive) if and only if
$$R'_n\vec{h}\in H^{r_n}_+ \quad \textrm{(or $R_n'\vec{h}\in {(H_+\setminus \{0\})^{r_n}}$ ),}$$
where $\vec{h}=(h_1, h_2, ..., h_{l_n})^{{T}} \in H^{l_n}$.
Moreover, for each $e^n_i$, write it as $e^n_i=(e^n_i)_+-(e^n_i)_-$ for $(e^n_i)_+, (e^n_i)_-\in K_0(S_n)_+$ and fix this decomposition.
Define  a $r_n\times 2l_n$ matrix
\begin{displaymath}
R_n=
R'_n
\left(
\begin{array}{ccccccc}
1 & -1 & 0 & 0 &\cdots & 0 & 0 \\
0 & 0 & 1 & -1 & \cdots & 0 & 0 \\
\vdots & \vdots & \vdots & \vdots & \ddots & \vdots & \vdots \\
0 & 0 & 0 & 0 & \cdots & 1 & -1
\end{array}
\right).
\end{displaymath}
Then one has $$\vec{r}_n=R_n\vec{e}_{n, \pm},$$ where $\vec{e}_{n, \pm}=((e^{n}_{1})_+, (e^{n}_{1})_-, ..., (e^{n}_{l_n})_+, (e^{n}_{l_n})_-)^{\mathrm T}$. Hence, for any ordered group $H$, and any elements $h_{1, +}, h_{1, -}, ..., h_{l_n, +}, h_{l_n, -}\in H$, the map $e^n_i\mapsto (h_{i, +}-h_{i, -})$, $i=1, ..., l_n$ induces a positive (or strictly positive) homomorphism if and only if
$$R_n\vec{h}_{{\pm}}\in H^{r_n}_+\quad \textrm{ (or $R_n\vec{h}_{{\pm}}\in (H_+\setminus\{0\})^{r_n}$) },$$
where $\vec{h}_{\pm}=(h_{1, +}, h_{1, -}, ..., h_{l_n, +}, h_{l_n, -})^{{T}}   \in H^{l_n}$.

Since $\{e^{n}_1, e^{n}_2, ..., e^{n}_{l_n}\}$ is a set of generators of $K_0(S_n)$, for any projection $p$ in a matrix algebra of $S_n$, there are integers $m^n_1({{[p]}}), ..., m^n_{l_n}({{[p]}})$ such that for any homomorphism $\tau: K_0(S_n)\to\mathbb R$, one has
$$
\tau([p])=\langle\vec{m}_n(p), \tau(\vec{e}_n)\rangle=\sum_{i=1}^{l_n} m_i^{{n}}({{[p]}})\tau(e^n_i)=\sum_{i=1}^{l_n} m^n_i({{[p]}})\tau((e^n_i)_+)-m^n_i({{[p]}})\tau((e^n_i)_-),$$
where $\vec{m}_n{{([p])}}=(m^n_1({{[p]}}), ..., m^n_{l_n}({{[p]}}))^{{T}}$ and $\vec{e_n}=(e^{n}_1, e^{n}_2, ..., e^{n}_{l_n}).$

%Keep the notation in \ref{June13-constr}, define
%$\psi_{k,k}=L_{k-1}: A\to S_{k},$
%$\psi_{k,n}={L_{n-1}\circ \Psi_{n-1}},$ $n>k,$ $k=2,3,...\ .$  %({?})
For each $p\in M_m(A),$ for some integer $m\ge1,$ denote by $[\psi_{{k}, k+j}(p)]$ an element
in $K_0(S_{{k+j+1}})$ associated with $\psi_{{k}, k+j}(p).$  Let $\imath_n: S_n\to A$ be the imbedding.
Denote by
$${\overline{(\imath_n)_{*0}}}: \vec{e_n}\mapsto (((\imath_n)_{*0}(e_1^n), (\imath_n)_{*0}(e_2^n),...,(\imath_n)_{*0}(e_{l_n}^n)).$$

\end{NN}
Then, by Lemma \ref{traces} and Remark \ref{Kstate}, one has the following lemma.

\begin{lem}\label{rho-uniform}
With the notion same as the above, for any $p\in\mathcal P_0,$  for each fixed $k$, one has that
$$
\tau(p)=
\lim_{n\to\infty}\sum_{j=1}^n(\sum_{i=1}^{{{l_{k+j+1}}}} m_i^{{k+j+1}}([{{\psi^k_{k+j}}}(p)])\tau((\imath_{k+n}\circ h_{k+j,k+n})_{*0}(e_i^{{k+j+1}}))) $$
$$~~~~~-m_i^{k+j+1}({{\psi^k_{k+j}}}(p)])\tau((\imath_{k+n}\circ h_{k+j,k+n})_{*0}(e_i^{{k+j+1}})_-))
$$
uniformly on $\mathrm{S}(K_0(A)).$  Moreover,  $\rho_A\circ (\imath_n)_{*0}\circ h_{k+j,k+n}(e_{i,\pm}^{k+j})$
converge{{s}} to a strictly positive element in $\mathrm{Aff}(\mathrm{S}(K_0(A)))$ as $n\to\infty$ uniformly.
\end{lem}

\begin{proof}
We first compute that, if $j>1,$
\beq\label{rho-comput-1508}
%&&(\sum_{i=1}^{l_n} m^n_i([\psi_{k,k+j}]([p]))\tau((\imath_n\circ {\tilde h_{k+j,n}})_{*0}(e^k_i)_+)-m^n_i([\psi_{k,k+j}]([p]))\tau((\imath_n\circ h_{k,n})_{*0}(e^k_i)_-))\\
&&\sum_{i=1}^{l_{k+j+1}} m^{k+j+1}_i([\psi^k_{k+j}]([p]))\tau((\imath_{k+n}\circ h_{k+j,k+n})_{*0}(e^{k+j+1}_i))\\
&=& \tau([L_{k+n}\circ\cdots \circ L_{k+j+1}\circ L_{k+j}\circ \Psi_{k+j-1}]([p]))
\eneq
and, if $j=1,$
\beq\label{rho-compute-1508-1}
&&\sum_{i=1}^{l_{k+1}} m^{k+1}_i([\psi^k_{k+1}]([p]))\tau((\imath_{k+n}\circ h_{k+1,k+n})_{*0}(e^{k+1}_i))\\
&=& \tau([L_{k+n}\circ\cdots \circ L_{k+j+1}\circ L_k]([p])).
\eneq
Thus  (see \ref{Constr-1508})
\beq\label{comput-1508-2}
&&\sum_{j=1}^n (\sum_{i=1}^{l_{k+j}} m^{k+j}_i([\psi^k_{k+j}(p)])\tau((\imath_n\circ h_{k+j,k+n})_{*0}(e^{k+j}_i)))\\
&=& \tau([L_{k+n}\circ\cdots \circ L_{k+j+1}\circ L_k]([p]))\\
&&+\sum_{j=2}^n\tau([L_{k+n}\circ\cdots \circ L_{k+j+1}\circ L_{k+j}\circ \Psi_{k+j-1}]([p]))\\
&=& \tau([L_{k+n}\circ\cdots \circ L_{k+j+1}\circ L_k]([p]))\\
&&+\sum_{j=2}^n\tau([L_{k+n}\circ\cdots \circ L_{k+j+1}\circ L_{k+j}\circ \Psi_{k+j-1}\circ J_{k,k+j-1}]([p]))\\
&=&\tau([L_{k+n}\circ J_{k, n+k-1}]([p])).
\eneq
Thus the first part of the lemma then follows from \ref{traces}. The second part also follows.

\end{proof}

One then has the following
\begin{cor}\label{rho}
Let $\mathcal P$ be a finite subset of projections in a matrix algebra over $A$,  let $G_0$ be the subgroup of $K_0(A)$ generated by $\mathcal P$ {{ and let $k\ge 1$ be an integer}}.  Denote by $\tilde{\rho}: G_{{0}}\to\Pi\mathbb Z$ the map defined by
\begin{eqnarray}
&&\hspace{-0.7in}[p] \mapsto \nonumber\\
&&\hspace{-0.7in} ({m}^{k}_1(q_0), -m^{k}_1(q_0), {m}^{k}_2(q_0), -m^{k}_2(q_0),\cdots m^{k}_{l_k}(q_0), -m^k_{l_k}(q_0), \nonumber\\
&&m^{k+1}_1(q_1), -m^{k+1}_1(q_1),m^{k+1}_2(q_1), -m^{k+1}_2(q_1), ...,m^{k+1}_{l_{k+1}}(q_1), -m^{k+1}_{l_{k+1}}(q_1),\nonumber\\
&&m^{k+2}_1(q_2), -m^{k+2}_1(q_2),m^{k+2}_2(q_2), -m^{k+2}_2(q_2), ...,m^{k+2}_{l_{k+2}}(q_2), -m^{k+2}_{l_{k+2}}(q_2), \cdots), \label{rho-1}
\end{eqnarray}
where $q_i=[\psi^{{{k}}}_{k+i}(p)],$ $i=0,1,2,....$
If $\tilde{\rho}(g)=0$, then $\tau(g)=0$ for any trace over $A$.
\end{cor}

By the definition of the map $\tilde{\rho}$ and $H=h'\circ h_0: A\to F_0\to  A$, using the same argument as that of Lemma 2.12 of \cite{Lnduke}, one has the following lemma.
\begin{lem}\label{ker}
Let $\mathcal P$ be a finite subset of projections in $M_k(A_1)\subseteq M_k(A)$. Then there is a finite subset
$\mathcal F_1\subset A_1$ and $\delta_0>0$ such that if the above construction starts with $\mathcal F_1$ and $\delta_0$, then
$$\ker\tilde{\rho}\subset\ker[H]\ \ \mbox{and}\ \ \ker\tilde{\rho}\subset\ker[h_0].$$
\end{lem}

The $K_0$-part of the existence theorem will almost factor through the map $\tilde{\rho}$, and this lemma will help us to handle the elements of $K_0(A)$ which vanish under $\tilde\rho$. Moreover, to get a such $K_0$-homomorphism, one also needs to find a copy of the generating set of the positive cone of $K_0(S)$ inside  the image of  ${{\tilde \rho}}$
as an ordered group
%the codomain's ordered group
for certain algebra $S\in \mathcal C_0$. In order to do so, one {needs} the following technical lemma, which is  {essentially} Lemma 3.4 of \cite{Lnduke}.

%{\bf The following Lemma should be rewritted.}
\begin{lem}[Lemma 3.4 of \cite{Lnduke}]\label{solveeq}
Let $S$ be a compact convex set, and let $\Aff(S)$ be the affine continuous functions on $S$. Let $\mathbb{D}$ be a dense ordered subgroup of $\Aff(S)$, and let $G$ be an ordered group with the strict order determined by a surjective homomorphism $\rho:G\to\mathbb{D}$. Let $\{x_{ij}\}_{1\leq i\leq r, 1\leq j<\infty}$ be an $r\times\infty$ matrix having rank $r$ and with $x_{ij}\in\mathbb Z$ for each $i,j$.  Let $g_j^{(n)}\in G$ such that $\rho(g_j^{(n)})=a_j^{(n)}$, where $\{a_j^{(n)}\}$ is a sequence of positive elements in $\mathbb{D}$ such that $a_j^{(n)}\to a_j (>0)$ uniformly on $S$ as $n\to\infty$.
% For each $n$, there is an $s(n)$ such that

{Further suppose that there is a sequence of integers $s(n)$ satisfying the following condition:}

{Let $\widetilde{v_n}=(g_j^{(n)})_{s(n)\times1}$ be the part of $(g_j^{(n)})_{1\leq j<\infty}$
%from the above
and let
$$\widetilde{y_n}=(x_{ij})_{r\times s(n)}\widetilde{v_n}$$
Denote by $y_n=\rho^{(r)}(\widetilde{y_n}).$
%---that is, by writing $\widetilde{y_n}=(\widetilde{b}_j^{(n)})$, we have $y_n=(b_j^{(n)})$ with $b_j^{(n)}=\rho(\tilde{b}_j^{(n)})$.
Then there exists $z=(z_j)_{r\times1}$ such that $y_n\to z$ on $S$ uniformly.}
%where $\widetilde{v_n}=(g_j^{(n)})_{s(n)\times1}$ and $\widetilde{y_n}=(\widetilde{b}_j^{(n)})\in G^r$. Set $b_j^{(n)}=\rho(\tilde{b}_j^{(n)})$ and $y_n=(b_j^{(n)})$.
%Suppose that $y_n\to z$ on $S$ uniformly for some $z=(z_j)_{r\times1}$.

{With the above condition,  there exist $\delta>0$ and a positive integer $K>0$ satisfying the following:}

{For some sufficiently large $n$, if $M$ is a positive integer,  and if $\tilde{z}'\in ({1\over{K^3}}G)^r$ (i.e.,  there is $\tilde{z}''\in G^r$ such that $K^3\tilde{z}''=\tilde{z}'$) satisfies $|| z-M{z}' ||<\delta$, where $z'=(z'_1,z'_2,\cdots, z'_r)$ with  $z'_j=\rho(\tilde{z}'_j)$
%if $\tilde{z}'=(\tilde{z}'_1,...,\tilde{z}'_r)$, and $M$ is a positive integer
 then there is a $\tilde{u}=(\tilde{c}_j)_{s(n)\times1}\in G^{s(n)}_+$ such that
  $$(x_{ij})_{r\times s(n)}\tilde{u}=\tilde{z}'.$$}

 {Moreover, if each $s(n)$ can be written as $s(n)=\sum_{k=1}^n l_k$, where $l_k$ are positive integers,  and for each $k$,   $R_k$ is a $r_k\times l_k$ matrix with entries in $\Z$ so that }
$$
{\bar R}_n={\rm diag}(R_1, R_2,...,R_n)
$$
satisfies
\beq\label{June14-1}
{\bar R}_n{\bar g}_n>0
\eneq
as an element in { $G^{\sum_{k=1}^n r_k},$} where
%$s(n)=\sum_{k=1}^nl_k$ and
${\bar g}_n=(g_1^n,g_2^n,...,g_{s(n)}^n)^T,$
$n=1,2,...,$,  then we may choose $\tilde u$ that
\beq\label{JUne14-2}
{\bar R}_n{\tilde u}>0.
\eneq
%Moreover, with $(g_i^{(n)})$ the elements as those of Lemma \ref{rho-uniform}, one can choose $$\tilde{u}=(c_{1, +}^{(n)}, c_{1, -}^{(n)}, ..., c_{s(n), +}^{(n)}, c_{s(n), -}^{(n)})$$ such that for each $n$, if $\{e^n_i\}$ are the generators of $K_0(S_n)$ fixed previously, then the map $e_i^n\mapsto (c_{i, +}^{(n)}-c_{i, -}^{(n)})$ induces a positive homomorphism from the $K_0(S)$ to $G$.
\end{lem}
\begin{proof}
The proof is repeating the argument of Lemma 3.4 of \cite{Lnduke}, and we need to show that if  (\ref{June14-1}) holds, then $u=(\tilde{c}_j)_{s(n)\times1}$ can be chosen to make (\ref{JUne14-2}) hold.

%Assume that $\{g_{j}^{(n)}\}_{j=1}^{s(n)}$ are positive elements of $K_0(S)$ which forms a basis of $K_0(S)$ as an abelian group. Since the positive cone of $K_0(S)$ is finitely generated, there are integer $k(n)$ and $k(n)\times s(n)$ integer matrix $A_n$ such that any entry of $A_n(g_{1}^{(n)},..., g_{2s(n)}^{(n)})^T$ is strictly positive. Hence for any ordered group $H$ and strict positive element $h_1, h_2, ..., h_{s(n)}\in H$, if any entry of $A_n(h_{1},..., h_{s(n)})^T$ is strictly positive, then the map $g^{(n)}_j\mapsto h_j$ extend to a positive map $K_0(S)\to H$.

%Note that if $(g_i^{(n)})$ are as those of Lemma \ref{rho-uniform},
By (\ref{June14-1}), any entry of $R_k(\frac{a_{{s(k-1)}+1}^{(n)}}{M}, ..., \frac{a_{{s(k)}}^{(n)}}{M})^T$ is strictly positive. Since
$a_j^{(n)}\to a_j$ and the strict order on $G$ is determined by $\rho,$
%and the positive cone of $K_0(S)$ is finitely generated,
there exists an integer $L_k\ge 1$ and $\delta'_k>0$ such that if $n\ge L_k$ and
$$||\rho(f_{j}^{(k)})-\frac{\rho(g_{j}^{(n)})}{M}||<\delta'_k/2\andeqn \|a_j/M-a_j^{(n)}/M\|<\dt_k'/4,\quad j=s(k-1)+1,..., s(k)$$
for some $f_{j}^{(k)}\in G_+$, then any entry of $R_k(f_{{s(k-1)}+1}^{(k)}, ..., f_{{s(k)}}^{(k)})^T$ is positive.
%and hence $e_i^n\mapsto (\tilde{d}_{j}^{(n)}-\tilde{d}_{j, -}^{(n)})$ induced a positive homomorphism.
%Fix $\delta'_k$, and
We may assume that $L_{k+1}> L_k,$ $k=1,2,.....$
One may assume that $\{\delta'_k\}$ is decreasing.

Without loss go generality, let us assume that $(x_{ij})_{r\times r}$ has rank $r$. Set
$A'_n=(x_{ij})_{r\times s(n)}.$ Then there is an invertible matrix $B'\in M_{r}(\mathbb Q)$ such that $B'A'_n=C_n,$ where $C_n=(I_r, D_n')$ for some $r\times (s(n)-r)$ matrix $D_n'$. Moreover, there is an integer $K$ such that all entries of $KB'$ and $K(B')^{-1}$ are integers. {Evidently, $KD_n'$ is also a matrix with integer entries.}

Since $\mathbb D$ is dense in $\mathrm{Aff}(S)$, there are $\xi_n\in G^{s(n)}$ such that $\xi_n=(\tilde{d}_j^{(n)})$ and,
for all $n\ge L_n,$
\begin{equation}\label{equ5001}
||K^3\rho(\tilde{d}_j^{(n)})-\frac{a_j^{(n)}}{M}||<\min\{\frac{1}{s(n)^2\cdot 2^n}, \delta_n'/4\},\,\,\, j=1,2,....,s(n).
\end{equation}

Let $\tilde{w}_n=K^3\xi_n$, and let $\rho(\tilde{d}_j^{(n)})=d_j^{(n)}$. Then
$$K^3 d_j^{(n)}\to K^3 d_j=\frac{a_j}{M}>0$$
uniformly on $S$.
Set $\tilde{y}_n'=A_n'\tilde{w}_n=K^3 A_n'\xi_n$, and set $y'_n=\rho^{(s(n))}(\tilde{y}_n')$, where $\rho^{(m)} G^m \to \mathbb D^m$ is the nature map induced by $\rho: G\to \mathbb D$.  Then one has $y_n'\to z/M$ uniformly on $S$.

We then has
$$C_n\tilde{w}_n=K^3C_n\xi_n=K^3B'A_n'\xi_n=B'\tilde{y}_n'$$
and
$$I_r\tilde{v}_n'=B'\tilde{y}_n'-D_n\tilde{w}_n,$$
where $\tilde{v}_n'=(K^3\tilde{d}_1^{(n)}, ..., K^3\tilde{d}_r^{(n)} )^{{T}}$ and
$D_n=(0, D_n')$.

Since $d_j^{(n)}\to d_j>0$ uniformly on $S$, there is $N_1>0$ such that
$$d_j^{(n)}\geq \inf\{\frac{d_j}{2}(\tau):\ \tau\in S\}>0$$
for all $n\geq N_1$ and $j=1, 2, ..., r$.
Let $k_r$ be a natural number such that $ r < s(k_r)$, and set
$$0 <\epsilon < \min\{\delta'_{k_r}/8, \inf\{d_j(\tau): \tau\in S\}/32K^4:\ j=1, 2, ..., r\}.$$

There is $N_2$ such that if $n\geq N_2$, then
$$||B'y_n-B'z ||_\infty<\epsilon/4;$$
and there is $N_3>0$ such that
$$||B'(y_n')-B'(z/M)||_\infty<\epsilon/4,\quad n\geq N_3.$$

There is $\delta>0$ depending only on $B'$ such that if $||z-Mz'||<\delta$, one has
$$||B'y'_n-B'z'||_\infty<\epsilon/2,\quad n\geq N_3.$$

Let $\tilde{z}'$ and $\tilde{z}''$ be described in the lemma. Set $N=\max\{N_1, N_2, N_3\}$. {Since both $KB'$ and $KD_n$ are matrix over $\mathbb Z$, $u'=KB'\tilde{z}''-KD_n\xi_n\in G^{r}$.} Let
\begin{equation}\label{equ4999}
B'\tilde{z}'-D_n\tilde{w}_n=K^3B'\tilde{z}''-K^3D_n\xi_n=K^2u',
\end{equation}
where $u'=(\tilde{c}_1, \tilde{c_2}, ..., \tilde{c}_r)\in G^{r}$ and $n\ge N.$ Set $u''=K^2u'$, and let $$\rho^{(r)}(u')=(c_1, c_2, ..., c_r)\in\mathbb D^r.$$ One may write
$$I_r u''=B'\tilde{z}'-D_{{n}}\tilde{w}_n,$$ and one has
$$||\rho^{(r)}(u'')-\rho^{(r)}(\tilde{v}_n')||_\infty=||\rho^{(r)}(B'\tilde{z}'-D_n\tilde{w}_n)-\rho^{(r)}(B'\tilde{y}_n'-D_n\tilde{w}_n)||_\infty<\epsilon.$$

Set $$\tilde{u}=(K^2\tilde{c}_1,..., K^2\tilde{c}_r, K^3\tilde{d}_{r+1}^{({n})},..., K^3\tilde{d}_{s(n)}^{({n})})^T.$$
%and the same argument as that of Lemma 3.4 of \cite{Lnduke} shows that
{ We calculate that
$$B'A'_n\tilde{u}=(I_r, D_n')\tilde{u}=u''+(0,D_n')\tilde{u}=B'\tilde{z}'-D_n\tilde{w}_n+(0,D_n')\tilde{u}=B'\tilde{z}'.$$
Consequently,}
$A'_n\tilde{u}=\tilde{z}',$
{since $B'$ is invertible.}
%$$\bar{u}=(\tilde{c}_1,..., \tilde{c}_r, K\tilde{d}_{r+1}^{(n)},..., K\tilde{d}_{s(n)}^{(n)}).$$
%
%
%With such $\xi_n$ and $\epsilon$, one obtains the solution (with the notation of \cite{Lnduke})
%$$\tilde{u}=(K^2\tilde{c}_1,..., K^2\tilde{c}_r, K^3\tilde{d}_{r+1}^{(n)},..., K^3\tilde{d}_{s(n)}^{(n)}).$$

In order to show that
%the map $e_i^n\mapsto (c_{i, +}^{(n)}-c_{i, -}^{(n)})$ induces a positive homomorphism,
%\beq\label{Jun14-4}
${\bar R}_n {\tilde u}>0,$
%\eneq
it suffices to
%one only to
verify that
\begin{equation}\label{equ5002}
|| K^2\rho(\tilde{c}_j)-\frac{\rho(g_j^{(L_{k_r})})}{M} ||<\delta'_{k_r}/2, \quad j=1,..., r
\end{equation}
and, if $n\ge L_n,$
\begin{equation}\label{equ5003}
|| K^3\rho(\tilde{d}_j^{(n)})-\frac{\rho(g_j^{(L_n)})}{M} ||<\delta'_{n}/2, \quad j=r+1,..., s(n).
\end{equation}

Equation \eqref{equ5003} follows from 
%Equation 
\eqref{equ5001} directly. Let us verify 
%Equation 
\eqref{equ5002}. By \eqref{equ4999}, one has
\begin{eqnarray*}
|| K^2\rho(u')-\rho(g')||&=&|| B'\rho(\tilde{z}')-D_n\rho(\tilde{w})_n-\rho(g')||\\
&\leq&|| B'y'_n-D_n\rho(\tilde{w})_n-\rho(g')||+\epsilon/2\\
&\leq&|| B'\rho(A'_n\tilde{w}_n)-D_n\rho(\tilde{w})_n-\rho(g')||+\delta'_{k_r}/16\\
&=&|| \rho((B'A_n'-D_n)(\tilde{w}_n)-g') ||+\delta'_{k_r}/16\\
&=&|| \rho((C_n-D_n)(\tilde{w}_n)-g') ||+\delta'_{k_r}/16\\
&=&|| K^3\rho((\tilde{d_j^{(n)}})_{r\times 1})-\rho(g')  ||+\delta'_{k_r}/16\le
\delta'_{k_r}/2,
\end{eqnarray*}
where $g'=(g_1^{(L_k)}/M, g_2^{(L_{k_r})}/M,...,g_r^{(L_{k_r})}/M),$ and $\rho^{r}: G^r\to \mathbb D^r$ is denoted by $\rho$.  This verifies Equation \eqref{equ5002}. Thus, $\tilde{u}$ is the desired solution.
\end{proof}

\begin{df}
A unital stably finite C*-algebra $A$ is said to have {the $K_0$-}density property if the image $\rho(K_0(A))$ is dense in $\mathrm{Aff}(\mathrm{S}_{[1]}(K_0(A))$, where $\mathrm{S}_{[1]}{(K_0(A))}$ is the convex set of the states of $K_0(A),$ i.e.,
the convex set of all
%convex of the
positive homomorphisms $r: K_0(A)\to \mathbb R$ satisfying
%$r(K_0^+(A))\subseteq\mathbb{R}$ and
{$r([1])=1$.}
\end{df}

\begin{rem}
By Corollary 7.9 of \cite{Goodearl}, the linear space spanned by $\rho(K_0(A))$ is always dense in $\mathrm{Aff}(\mathrm{S}_{[1]}(K_0(A))$. Therefore, the unital stably finite C*-algebra in the form of $A\otimes U$ for a UHF-algebra $U$ always has the {$K_0$-}density property. Moreover, any unital stably finite C*-algebra $A$ which is tracially approximately divisible has the {$K_0$-}density property.
\end{rem}

\begin{rem}
Not all C*-algebras in $\mathcal B_1$ with the (SP) property satisfy the density property. The following is one example:
Consider
\beq\nonumber
G=\{(a,b)\in {\mathbb Q}\oplus {\mathbb Q}:  2a-b \in \mathbb Z\}\andeqn\\
G_+\setminus \{0\}=({\mathbb Q}_+\setminus\{0\}\oplus \Q_+\setminus\{0\} )\cap G
\eneq
and $1=(1,1)\in G_+$ as unit.

%For each rational number $r$, consider the function $g_r: [0, 1]\to \mathbb R$ defined by
%$$g_r(x)=rx+r.$$  Consider the group $G\subseteq \{f: [0, 1]\to\mathbb R\}$ generated by
%$$\{g_r:\ r\in\mathbb Q\}\cup\{1\}$$
%with the positive cone defined by
%$$\{f:\ \textrm{$f(x)>0$ for all $x\in[0, 1]$}\}\cup\{0\}.$$

Then $(G, G_+, 1)$ is a weakly unperforated rational Riesz simple ordered group (but not a  Riesz group---see \cite{LNjfa}).
Evidently $G$ has (SP) property; but the image of $G$ is not dense in $\mathrm{Aff}(\mathrm{S}_{[1]}(G)=\mathbb R\oplus \mathbb R$, as $(1/2,1/2)\in \mathbb R\oplus \mathbb R$ is not in the closure of the image of $G$. We leave the details to the readers.

%Also note that for each $x\in[0, 1]$, the evaluation $f\mapsto f(x)$ is a state of $G$. Denote it by $\tau_x$.

%The ordered group $G$ has the (SP) property, as for any positive $N\in\mathbb N$, there is $r\in\mathbb Q^+$ such that $Ng_r< 1$, and hence $\tau(g_r)<1/N$ for any state $\tau$ of $G$. By Theorem \ref{RangT}, there is a simple C*-algebra $A\in\mathcal B_1$ such that $K_0(A)=(G, G^+, 1)$.

%One asserts that there does not exists $f\in G$ such that
%\begin{equation}\label{approx-1/2}
%|\tau(f)-1/2|<1/32\rforal \tau\in \mathrm S(G),
%\end{equation}
%and hence the image of $G$ cannot be dense in $\mathrm{Aff}(\mathrm{S}(G))$.
%If such $f$ exists, there are $m_1, ..., m_k, n,\in\mathbb Z$ and $r_1, r_2, ..., r_k\in\mathbb Q$ such that
%\begin{eqnarray*}
%f&=&m_1g_{r_1}+\cdots+m_kg_{r_k}+n\\
%&=&(m_1r_1+\cdots+m_kr_k)x+ (m_1r_1+\cdots+m_kr_k) +n.
%\end{eqnarray*}
%Applying \eqref{approx-1/2} to $\tau_x$, one has
%\begin{equation}\label{n-est-01001}
%|(m_1r_1+\cdots+m_kr_k)x- ((m_1r_1+\cdots+m_kr_k) +n+1/2)|<1/32\rforal x\in[0, 1].
%\end{equation}
%In particular, it implies the function $x\mapsto (m_1r_1+\cdots+m_kr_k)x$ is almost constant up to $1/32$, and hence
%\begin{equation}\label{n-est-01002}
%|m_1r_1+\cdots+m_kr_k|<1/16.
%%\end{equation}
%Then \eqref{n-est-01001} and \eqref{n-est-01002} implies
%$$(m_1r_1+\cdots+m_kr_k) +n+1/2\in (-3/32, 3/32),$$
%and by \eqref{n-est-01002} again, one has
%$$n+1/2\in(-5/32, 5/32),$$
%which is absurd. So, the C*-algebra $A$ does not have the density property.
\end{rem}

\begin{prop}\label{est-a2m}
Let $A\in\mathcal B_{0}$ satisfy the density property and let $B_1$ be an inductive limit C*-algebra in Theorem \ref{RangT} such that $$(K_0(A),{K_0(A)_+},[1_A],K_1(A))\cong(K_0(B),{K_0(B)_+},[1_B],K_1(B)),$$
where $B=B_1\otimes U$ for a UHF-algebra $U$ of infinite type.
%{\color{Green} (I'm not sure if we need the assumption $B=B_1\otimes U$ for a UHF algebra $U$ of infinite type in the proof. Implicitly, this would imply that $A$ is UHF-stable with some UHF algebra, which is not true in general (for instance, consider $A$ to be the irrational rotation algebra). Of course, this will not affect the classification theorem.)} {\red{L: I think it is sufficient to have infinite dimensional UHF-algebra $U.$ }}
 Let $\alpha\in KL(A, B)$ be an element which implements the isomorphism above. Then, for any $\mathcal P\subset P(A)$, there is a sequence of completely positive linear maps $L_n: A\to B$ such that
$$\lim_{n\to\infty} ||L_n(ab)-L_n(a)L_n(b)||= 0\rforal a, b\in A$$ and $[L_n]|_{\mathcal P}=\alpha|_{\mathcal P}$ as $n\to\infty$.
\end{prop}

\begin{proof}
By Lemma \ref{MF}, $A$ is the closure of an increasing union of RFD \SCA s $\{A_n\}$. We may assume
${\mathcal P}\subset \underline{K}(A_1)$. Let $G=G({\cal P})$ be the subgroup generated by ${\cal P}$ and let ${\mathcal P}_0\subset{\mathcal P}$ be such that ${\mathcal P}_0$ generate $G\cap{K_0}(A).$
%where $G({\mathcal P})$ is the group generated by ${\mathcal P}$.
Write ${\mathcal P}_0=\{p_1,..., p_l\},$ where $p_1, ..., p_l$ are
projections in a matrix algebra over $A$.  Let $G_0$ be the group generated by ${\cal P}_0.$ Let $\mathcal F_1$ be a finite subset of $A_1$ and let $\delta_0>0$ be such that
any $\mathcal F_1$-$\delta_0$-multiplicative linear map, the map $[L]|_{\mathcal P}$ is well defined. Moreover, one requires that $\mathcal F_1$ and
$\delta_0$ satisfy Lemma \ref{ker}. Let $k_0$ be an integer such that
$G({\mathcal P})\cap\mbox{K}_i(A,\mathbb Z/k\mathbb Z)=\{0\}$ for any $k\geq k_0$, $i=0,1.$

By Theorem \ref{kkmaps}, there are two $\mathcal F_1$-$\delta_0/2$ multiplicative \morp s
$\Phi_0, \Phi_1$ from $A$ to $B\otimes\mathcal {\mathcal K}$ such that
$$[\Phi_0]|_{\mathcal P}=\alpha|_{\mathcal P}+[\Phi_1]|_{\mathcal P}$$
and the image of $\Phi_1$ is in a finite dimensional \SCA.
Moreover, we may also  assume that $\Phi_1$ is a homomorphism when it is restricted on $A_1$, and the image is a finite-dimensional C*-algebra.
With $\Phi_1$ in the role of $h_0$, we can proceed with the construction as described at the beginning of this
section. We will keep the same notation.

Consider the map $\tilde\rho:G({\mathcal P})\cap\mbox{K}_0(A)\to l^\infty(\mathbb Z)$ defined in Corollary \ref{rho}. The
linear span of $\{\tilde\rho(p_1), ..., \tilde\rho(p_l)\}$ over $\mathbb Q$ will have finite rank, say $r$. So,
we may assume that $\{\tilde\rho(p_1), ..., \tilde\rho(p_r)\}$ are linearly independent and the $\mathbb Q$-linear span of them
give us the whole subspace. Therefore, there is an integer $M$ such that for any $g\in\tilde\rho(G_0)$, the element $Mg$
is in the subgroup generated by $\{\tilde\rho(p_1), ..., \tilde\rho(p_r)\}$. Let
$x_{ij}=(\tilde\rho(p_i))_j$, and $z_i=\rho_B(\af([p_i])\in\mathbb{D}$, where $\mathbb{D}=\rho_B(K_0(B))$ in
$\text{Aff}(S_{[1]}(K_0(B))$. Since $A$ is assumed to have the density property, so is $B$. Therefore the image $\mathbb{D}$ is a dense subgroup of $\text{Aff}(\mathrm{S}_{[1]}({K_0}(B)))$.

Let  $\{S_j\}$ be the sequence of \SCA s in ${\cal C}_0$ in the construction at the beginning of this section.
Fix $k\ge 1.$
Let $e_i^{{{k+j}}}, e_{i,\pm}^{{{k+j}}}\in K_0(S_{k+j}), $ $i=1,2,..., l_{k+j}$  and
$R_{{{k+j}}}$ be {$r_{k+j}\times 2l_{k+j}$ matrix} as described in \ref{Constr-Sn}.
{Let $s(j)=\sum_{i=1}^j2l_{k+i},$ $j=1,2,....$} Put
\beq\label{June14n2-1}
\af([ \imath_n\circ {\tilde h_{k+j, n}}(e_{i,+}^{{k+j}})])=g^{(n)}_{{s(j-1)+2i-1}},\,\,\,
%2l_{k-1}+2i-1}^{(n)},\,\,\,
\af([\imath_n\circ {\tilde h_{k+j,N_0+n}}(e_{i,-}^{{k+j}})])=g^{(n)}_{s(j-1)+2i},
 \eneq
 $i=1,2,...,l_k,$
and $a_j^{(n)}=\rho_B(g_j^{(n)}),$ $j=1,2,..., s(n)=\sum_{j=1}^n l_{k+j},$ $n=1,2,....$ Note that $a_j^{(n)}\in\mathbb{D}^+\backslash\{0\}.$ It follows from Lemma \ref{traces} that  $\lim_{n\to\infty}a_{{s(j-1)+2i}}^{(n)}=a_{s(j-1)+2i}=\rho_B(\af(g_{s(j-1)+2i}^{(n)}))>0$  and  $\lim_{n\to\infty}a_{s(j-1)+2i-1}^{(n)}=a_{s(j-1)+2i-1}=\rho_B(\af(g_{s(j-1)+2i-1}^{(n)}))>0$ uniformly.
 Moreover, by \ref{rho-uniform},
$\sum_{j=1}^nx_{ij}a_j^{(n)}\to z_i$ uniformly.  Furthermore, by \ref{Constr-Sn},
${\bar R}_n {\bar g}_n>0,$ ${\bar g}_n=(g_1^n, g_2^n,...,g_{{s(n)}}^n)^T.$

So, Lemma \ref{solveeq} applies. Fix $K$ and $\delta$ obtained from Lemma \ref{solveeq}.
%Also note that if $g\in\ker\tilde\rho$, then $\tau(g)=0$ for all traces $\tau$. By Lemma \ref{ker}, one has that
%$g\in\ker[H]$ and $g\in\ker[\Phi_1]$, and hence $g$ is in the kernel of any direct sum of $[\Phi_1]$.

Let $\Psi:=\Phi_0\oplus(\overbrace{\Phi_1\oplus\cdots\oplus\Phi_1}^{MK^3(k_0+1)!-1})$. Since $\Phi_1$
factors through a finite-dimensional C*-algebra, it is zero when restricted to $K_1(A)\cap G$  and
${K_1}(A,\mathbb Z/k\mathbb Z)\cap G$ for
$2\le k\le k_0.$ Moreover,  the map $(\overbrace{\Phi_1\oplus\cdots\oplus\Phi_1}^{MK^3(k_0+1)!})$ vanishes on ${K_0}(A,\mathbb Z/k\mathbb Z)$ for $2\le k\le k_0.$ Therefore we have $$[\Psi]|_{K_1(A)\cap G}=\alpha|_{K_1(A)\cap G},\quad [\Psi]|_{K_1(A,\mathbb Z/k\mathbb T)\cap G}=\alpha|_{K_1(A,\mathbb Z/k\mathbb Z)\cap G}$$ and $[\Psi]|_{K_0(A,\mathbb Z/k\mathbb Z)\cap
G}=\alpha|_{K_0(A,\mathbb Z/k\mathbb Z)\cap G}$. We may assume $\Psi(1_A)$ is a projection in M$_r(B)$ for some integer $r$.
% Then the map $\Psi$ induces the desired{---what is the desired} maps on all the invariants except on the ${K_0}$ group. In the
%following, we will compress the map $\Psi$ to a small corner of $B\otimes{\mathcal K}$, and then find a suitable map $h$ to fill the $K_0$-part which is missing under $\Psi$.

We may also assume  that there exist projections $\{p'_1, ..., p'_l\}$ in  $B\otimes{\mathcal K}$ which are sufficiently close
to $\{\Psi(p_1), ..., \Psi(p_l)\}$ respectively, so that $[p'_i]=[\Psi(p_i)]$. Note that $B\in\mathcal B_0$, and hence the strict order on the projections of $B$ is determined by traces. Thus there is a projection $q'_i\le p'_i$ such that $[q'_i]=MK^3(k_0+1)![\Phi_1(p_i)]$. Set $e'_i=p'_i-q'_i$, and let
${\mathcal P}_1=\Psi({\mathcal P})\cup\Phi_1({\mathcal P})\cup\{p'_i,q'_i,e'_i;i=1,...,l\}$. Denote by $G_1$ the group generated
by ${\mathcal P}_1$. Recall that $G_0=G({\mathcal P})\cap{K_0}(A)$, and decompose it as $G_{00}\oplus G_{01}$, where
$G_{00}$ is the infinitesimal part of $G_0$. Fix this decomposition and denote by $\{d_1, ..., d_t\}$ the positive elements which generate $G_{01}$.

Applying Corollary \ref{smallmap} to M$_r(B)$ with any finite subset $\mathcal{G}$, any $\epsilon>0$ and any
$0<r_0<\delta<1$, one has a $\mathcal{G}$-$\epsilon$-multiplicative map $L:\mbox{M}_r(B)\to\mbox{M}_r(B)$ with the following properties:
\begin{enumerate}
\item $[L]|_{{\mathcal P}_1}$ and $[L]|_{G_1}$ are well defined;
\item $[L]$ induces the identity maps on the infinitesimal part
      of $G_1\cap$K$_0(B)$, $G_1\cap\mbox{K}_1(B)$,
      $G_1\cap\mbox{K}_0(B,\mathbb Z/k\mathbb Z)$ and $G_1\cap\mbox{K}_1(B, \mathbb Z/k\mathbb Z)$
      for the $k$ with $G_1\cap\mbox{K}_i(B, \mathbb Z/k\mathbb Z)\neq\{0\}$,
      $i=0,1$;
\item $\tau\circ[L](g)\leq r_0\tau(g)$ for all $g\in
      G_1\cap\mbox{K}_0(B)$ and $\tau\in T(B)$;
\item There exist positive elements $\{f_i\}\subset{K_0}(B)_+$
      such that for  $i=1,...,t,$
      $$\alpha(d_i)-[L](\alpha(d_i))=MK^3(k_0+1)!f_i.$$
\end{enumerate}

 Using the compactness of $T(B)$ and the strict comparison for positive elements for  $B$, the positive number $r_0$ can be chosen sufficiently small such that
 $\tau\circ[L]\circ[\Psi]([p_i])<\delta/2$ for all $\tau\in T(B)$,
 and $\alpha([p_i])-[L\circ\Psi]([p_i])>0,$  $i=1,2,...,l.$

Let $[p_i]=\sum_{j=1}^tm_jd_j+s_j,$ where $m_j\in\Z$ and $s_j\in G_{00}.$
Note, by (2) above, $(\af-[L]\circ \af)({s_i})=0.$  Then we have
$$\begin{array}{ll}
  & \alpha([p_i])-[L\circ\Psi]([p_i]) \\
  =&\alpha([p_i])-([L\circ\alpha]([p_i])+MK^3(k_0+1)![L\circ\Phi_1]([p_i])) \\
  =&(\alpha(\sum m_jd_j)-[L\circ\alpha](\sum m_jd_j))-MK^3(k_0+1)![L\circ\Phi_1]([p_i]) \\
  =& MK^3(k_0+1)!(\sum m_jf_j-[L]\circ[\Phi_1]([p_i]))=MK^3(k_0+1)!f_i',
 \end{array}
 $$
 {where $f_i'=\sum m_jf_j-[L]\circ[\Phi_1]([p_i]),$}
for $i=1,2,...,l.$
 Define $\bt([p_i])=K^3(k_0+1)!f_i',$ $i=1,2,...,l.$
%\beq\label{June16-10}
%\bt|_{G_{00}}=0,\,\,\, \bt(d_i)=K^3(k_0+1)!(f_i-[L\circ \Phi_1](d_i)),\,\,\,
%i=1,2,...,t.
%\eneq
%Let $x=\sum_{i=1}^tm_id_i+s_x,$ where $m_i\in\Z$ and $s_x\in G_{00}.$
%Note, by (2) above, $(\af-[L]\circ \af)(s_x)=0.$  Then we have
%$$\begin{array}{ll}
 % & \alpha(x)-[L\circ\Psi](x) \\
 % =&\alpha(x)-([L\circ\alpha](x)+MK^3(k_0+1)![L\circ\Phi_1](x)) \\
 % =&(\alpha(\sum m_id_i)-[L\circ\alpha](\sum m_id_i))-MK^3(k_0+1)![L\circ\Phi_1](x) \\
  %=& MK^3(k_0+1)!(\sum m_if_i-[L]\circ[\Phi_1](x)) \\
  %=& M\bt(x).
  %%, \quad \mbox{where}\quad f'_j=\sum
 % m_j^{(i)}f_j-[L]\circ[\Phi_1]([p_i]).
%\end{array}$$
%Define $\bt: G_0\to K_0(B)$ by $\bt|_{G_{00}}=0,$ $\bt(d_i)=K^3(k_0+1)!f_i,$
%$i=1,2,...,t.$

%Note that  $\alpha([p_i])-[L\circ\Psi])([p_i])>0,$ $ i=1,2,...,r.$
%and ${K_0}(B)$ is
%weakly unperforated, one has that $f'_j>0$.

%Then define
%$\beta:G\cap{K_0}(A)\to{K_0}(B)$ by
%$$\beta=(\alpha-[L\circ\Psi])/M.$$ Note that $\beta([p_i])=K^3(k_0+1)!f'_i$.

Let us now construct a map $h': A \to B.$
%which carries the map $\beta$.
It will be constructed by
factoring through the $K_0$-group of some \CA\, in the class $\mathcal C_0$ in the construction given at the beginning of this section.
Let $\tilde{z_i}'=\beta([p_i])$, and $z'_i=\rho_B(\tilde{z_i}')\in\Aff(S_{[1]}(K_0(B)))$.
%Recall that we identify
%the ${K_0}(A)$ and ${K_0}(B)$ and $\alpha$ is the identity.
Then we have:
$$\begin{array}{lll}
    ||{Mz'-z}||_\infty & = & \max_i\{||\rho_B(\alpha([p_i])-[L\circ\Psi]([p_i]))-\rho(\af([p_i]))||\}\\
     & = & \max_i\{\sup_{\tau\in T(B)}\{\tau\circ[L]\circ[\Psi]([p_i])\}
     \le \delta/2,
  \end{array}$$
  where $z=(z_1,z_2,...,z_r)$ and $z'=(z_1',z_2',...,z_r').$
By Lemma \ref{solveeq},  for sufficiently large $n,$ one obtains
$\tilde{u}=(u_1,u_2,...,u_{s(n)})\in K_0(B)^{{s(n)}}_+$
%=\{u_1, u_2, ..., u_{r_1+\cdots+r_{s(n)}}\}=\{\{u^{(1)}_1, ...,u^{(1)}_{r_1}\}, ..., \{u^{(s(n))}_1, ...,u^{(s(n))}_{r_{s(n)}}\}\},$$ where $u^{(j)}_i$'s are positive elements of ${K_0}(B)$
such that

\begin{equation}\label{June-16-n1}
\sum x_{ij}u_j=\tilde{z}'_i.
\end{equation}
More importantly,
\beq\label{June16-3}
{\bar R}_n{\tilde u}>0.
\eneq
%Moreover, the map $g^{(k)}_i \mapsto u^{(k)}_i$, $1\leq i\leq r_k$ induces a positive map from ${K_0}(S_k)\to K_0(B)$, where $\{g^{(k)}_1, ..., g^{(k)}_{r_k}\}$ is the generator of $K_0(S_k)$ which is fixed previously. Therefore, with $D=S_1\oplus\cdots\oplus S_{s(n)}$,
It follows from \ref{Constr-Sn} that the maps $$e_i^{k+j}\mapsto (u_{{s(j-1)}+2i-1}-u_{s(j-1)+2i}),\quad 1\leq j\leq n, 1\leq i\leq l_{k+j}$$  defines
%an order preserving
strictly positive homomorphism, $\kappa_0^{(k+j)}$ from $K_0(S_{k+j})$ to $K_0(B)$ which defines a strictly positive homomorphism from ${K_0}(D)$ to ${K_0}(B)$, where $D=S_{{{k+1}}}\oplus\cdots\oplus S_{{k+n}}$. Since $B\in\mathcal B_{u0}$, by Corollary \ref{C0ext},  there is a
 \hm\, $h':D\to \mbox{M}_m(B)$ for some large $m$
 %(note that we use the semi-projectivity of $D$ twice!)
 such that $h'_{*0}|_{S_k}=\kappa_0^{(k)}.$
By \eqref{June-16-n1}, one has, keeping the notation in the construction at the beginning of this section,
$$h'_{*0}([\psi^{{{k}}}_{k+1}(p_i)], [{{\psi^k_{k+2}}}(p_i)], ...,[{{\psi^k_{n+k}}}(p_i)])=\beta([p_i]),\quad i=1,...,r.$$
Now, define $h'': A\to  {{D \to}} \mathrm{M}_{{{m}}}(B)$ by $$h''=h'\circ ({{{\psi^k_{k+1}}}}\oplus {{\psi^k_{k+2}}}{\oplus}\cdots\oplus {{\psi^k_{k+n}}}).$$
%where $\phi_{1, k}: {?}A\to S_k$, $k=1, ..., n$, are the maps constructed in \ref{Constr-Sn}.
%regard $h'$ as a map from $A$ to M$_k(B)$ by consider it with the map $\psi_1\oplus\cdots\oplus\psi_{s(n)}:A\to D$.
Then $h''$ is $\mathcal F$-$\delta$-multiplicative.

For any $x\in\ker\tilde{\rho}$, by Lemma \ref{ker}, $x\in\ker \rho_B\circ\alpha\cap\ker [H]$ and $x\in\ker[h_0]=\ker[\Phi_1]$. Therefore, we have $[\Phi_1](x)=0$ and $[\Psi](x)=\alpha(x).$ Note that $\alpha(x)$ also vanishes under any state of $({K_0}(B), {K_0}^+(B))$, we have $[L]\circ\alpha(x)=\alpha(x)$. So, we get $$\alpha(x)-[L\circ\Psi](x)=0.$$ Therefore $\alpha-[L\circ\Psi]|_{{\rm ker}{\tilde \rho}}=0.$
%gives us a homomorphism on $G_0$.
Therefore we may view $\af-[L\circ \Psi]$ as a \hm\, {from} ${\tilde \rho}(G_0).$
Since $Mg$ is in the subgroup generated by ${\tilde \rho}([p_1]),...,{\tilde \rho}([p_r])$ for any $g\in {\tilde \rho}(G_0).$
{Recall} that
\beq\label{June16-11}
(\af-[L\circ \Psi])([p_i])=M\bt([p_i]),\,\,\,i=1,2,...,r.
%,...,l.
\eneq

Set $h$ to be $M$ copies of $h''$. The map $h$ is $\mathcal F$-$\delta$-multiplicative, and $$[h]([p_i])=\alpha([p_i])-[L]\circ[\Psi]([p_i])\quad i=1,...,{r}.$$ Note that $[h]$ also has the multiplicity $MK^3(k_0+1)!$, and $D\in\mathcal C_0$ with trivial ${K_1}$ groups. One can concludes that $h$ induces zero map on $G\cap{K_1}(A)$, $G\cap{K_1}(A,\mathbb Z/k\mathbb Z)$ and $G\cap{K_1}(A,\mathbb Z/k\mathbb Z)$ for $k\leq k_0$. Therefore, we have $$[h]|_{\mathcal P}=\alpha|_{\mathcal P}-[L]\circ[\Psi]|_{\mathcal P}.$$

Set $L_1=(L\circ\Psi)\oplus h$. It is $\mathcal F$-$\delta$ multiplicative and $$[L_1]|_{\mathcal P}=[h]|_{\mathcal P}+[L]\circ[\Psi]|_{\mathcal P}=\alpha|_{\mathcal P}.$$
We may assume $L_1(1_A)=1_B$ by taking a conjugation with a partial isometry. Then $L_1$ is an $\mathcal F$-$\delta$-multiplicative map from $A$ to $B$, and $[L_1]|_{\mathcal P}=\alpha|_{\mathcal P}$.
\end{proof}

\begin{cor}\label{kk-attain}
Let $A$ be a amenable C*-algebra in the class $\mathcal B_1$, and assume that $A$ satisfies the {$K_0$-}density property and the UCT. Then $A$ is KK-attainable with respect to ${\cal B}_{u0}.$
\end{cor}
\begin{proof}
Let $C$ be any C*-algebra in $\mathcal B_{u0}$, and let $\alpha\in KL^{++}(A, C)$. We may write that $C=C_1\otimes U$ for some $C_1\in {\cal B}_0$ and for some UHF-algebra of infinite type.
By Theorem \ref{RangT}, there is a C*-algebra $B$ which is an inductive limit of C*-algebras in the class $\mathcal C_0$ together with homogeneous \CA s in the class $\mathbf H$
%{\color{Green} \bf ----Need to be revised  according to new Guihua's file---Let's do it later}
such that
$$(K_0(A), {K_0(A)_+}, [1_A]_0, K_1(A)) \cong (K_0(B), {K_0(B)_+}, [1_B]_0, K_1(B)).$$
Since $A$ satisfies the UCT, there is an invertible $\beta\in KL^{++}(A, B)$ such that $\beta$ carries the isomorphism of K-theories of $A$ and $B$. Applying Proposition \ref{est-a2m} to $\beta$ and applying Corollary \ref{est-m2a} to $\alpha\circ\beta^{-1}$, one has the desired conclusion.
\end{proof}

\begin{thm}\label{MEST}
Let $A\in\mathcal B_1$ be {{a}}  amenable  {{\CA}} satisfying the {$K_0$-}density property and the UCT, and let $B\in {\mathcal B}_{u0}$. Then for any $\alpha\in KL^{++}(A, B)$, and any $\gamma: T(B)\to T(A)$ which is compatible to $\alpha$, there is a sequence of completely positive linear maps $L_n: A\to B$ such that
$$\lim_{n\to\infty} ||L_n(ab)-L_n(a)L_n(b)||= 0\rforal a, b\in A$$
$$[L_n]=\alpha\quad\mathrm{and}\quad \lim_{n\to\infty}\sup_{\tau\in T(B)}|\tau\circ L_n(f)-\gamma(\tau)(f)|= 0\rforal f\in A.$$
\end{thm}

\begin{proof}
It follows from Corollary \ref{kk-attain} and Proposition \ref{add-tr} directly.
\end{proof}

\section{The isomorphism theorem}

\begin{df}\label{bdbk-k1}
Put $C'=P\text{M}_n(\text{C}(X))P$, where $\overbrace{ \T\sqcup \cdots \sqcup \T}^s\sqcup Y$
for some connected finite CW-complex $Y$ with torsion $K_1$-group (no restriction on $K_0(C(Y))$) and dimension no more than $3$, and $P$ is a projection in $\text{M}_n(\text{C}(X))$ with rank $r\geq6$. Then $K_1(P\text{M}_n(\text{C}(X))P)=\mathrm{Tor}(K_1(C'))\oplus G_1$ for some torsion free group $G_1\cong\mathbb Z^s$.  Let $D'=\bigoplus_{i=1}^s \mathrm{M}_r(\mathrm{C}(\mathbb T))$. Let  $\Pi'_i: PM_n(C(X))P\to E_i\, (E_i=M_r(C(\T)))$ {be} defined
by $\Pi_i'(f)(x)=f|_{\T_i},$ where $\T_i$ is the $i$-th circle, for all $f\in C',$ $i=1,2,...,s.$ Define
$\Pi': P\text{M}_n(\text{C}(X))P\to D'$ by $\Pi'(f)=(\Pi_1'(f),\Pi_2'(f),...,\Pi_s'(f))$ for all $f\in C'.$
We have that $K_1(D')\cong G_1.$

Denote by $C$ a finite direct sum of the C*-algebras of the form $C'$ above, matrix algebras {and}  \CA s in $\mathcal C_0$ (with trivial $K_1$-group).  Denote by $D$ the direct sum of $D'$ corresponding to  the C*-algebras in the form $C'$.
In other words, $C=D\oplus C_0,$ where $C_0$ is a direct sum of \CA s in ${\cal C}_0$ and
those with the form $PM_n(C(Y))P,$ where $Y$ is a connected finite CW complex
with $K_1(C(Y))$ is a finite abelian group.
Then, one has that $$\text{U}(C)/\text{CU}(C)\cong\text{U}_0(C)/\text{CU}(C)\oplus K_1(D) \oplus \mathrm{Tor}(K_1(C)).$$
Here we identify $K_1(D)\oplus {\rm Tor}(K_1(C))$ with a subgroup of $U(C)/CU(C).$
Denote by $\pi_0, \pi_1, \pi_2$  the projection maps from $\text{U}(C)/\text{CU}(C)$ to each component according to the decomposition above. Define $P': C\to PM_n(C(X))P$ {to}  be the projection and $\Pi=\Pi'\circ P'.$

We will frequently refer to the above notation later in  this section.
\end{df}

As in \cite{LinTAI}, we have the following lemmas to control the maps from  $U(C)/CU(C)$ in the approximate intertwining argument in the proof of \ref{IST0}. The proofs are the repetition of the corresponding arguments in \cite{LinTAI}.

\begin{lem}[See Lemma 7.2 of \cite{LinTAI}]\label{LinTAI-72}
Let $C$ be as above,  let $\mathcal U\subset \mathrm{U}(C)$ be a finite subset, and let $F$ be the group generated by $\mathcal U$. Suppose that $G$ is a subgroup of $\mathrm{U}(C)/\mathrm{CU}(C)$ which contains $\overline{F},$ the image of $F$ in $U(C)/CU(C),$  also contains $\pi_1( \mathrm{U}(C)/\mathrm{CU}(C))$, and $\pi_2(\mathrm{U}(C)/\mathrm{CU}(C))$. Suppose that the composition map $\gamma: \overline{F}\to \mathrm{U}(D)/\mathrm{CU}(D)\to \mathrm{U}(D)/\mathrm{U}_0(D)$ is injective{{---that is, if $x,y \in \overline{F}$ and $x\not=y$, then $[x]\not=[y]$ in $U(D)/U_0(D)$.}} Let $B$ be a unital C*-algebra and $\Lambda: G\to \mathrm{U}(B)/\mathrm{CU}(B)$ be a homomorphism such that $\Lambda(G\cap(\mathrm{U}_0(C)/\mathrm{CU}(C)))\subset \mathrm{U}_0(B)/\mathrm{CU}(B)$.
Let $\theta: \pi_2(\mathrm{U}(C)/\mathrm{CU}(C))\to \mathrm{U}(B)/\mathrm{CU}(B)$ be defined by
$\theta(g)=\Lambda|_{\pi_2(\mathrm{U}(C)/\mathrm{CU}(C))}(g^{-1})$ for any $g\in\pi_2(\mathrm{U}(C)/\mathrm{CU}(C))$.
Then there is a homomorphism $\beta: \mathrm{U}(D)/\mathrm{CU}(D)\to \mathrm{U}(B)/\mathrm{CU}(B)$ with $$\beta(\mathrm{U}_0(D)/\mathrm{CU}(D))\subset \mathrm{U}_0(B)/\mathrm{CU}(B),$$ and  such that $$\beta\circ\Pi^{\ddagger}\circ\pi_1(\bar{w})=\Lambda(\bar{w})(\theta\circ\pi_2(\bar{w}))\tforal w\in F.$$ 
%for any $w\in F.$
If furthermore $B\cong B_1\otimes U$ for a unital C*-algebra $B_1\in\mathcal B_0$ and a UHF-algebra $U$, and
$\Lambda(G)\subset \mathrm{U}_0(B)/\mathrm{CU}(B)$, then $\beta\circ\Pi^{\ddagger}\circ(\pi_1)|_{\bar{F}}=\Lambda|_{\bar{F}}$.
\end{lem}
The above may be summarized by the following commutative diagram:
\begin{displaymath}
\xymatrix{
{\bar F} \ar[d]^{\pi_1}  \ar[r]^{\mathrm{inclusion}} & G \ar[d]^{\Lambda+\theta\circ\pi_2} \\
\pi_1(\bar F) \ar[d]^{\Pi^\ddag} &  U(B)/CU(B) \\
U(D)/CU(D) \ar[ur]_{\beta}&
}
\end{displaymath}
%
%$$
%{\begin{array}{cccc}
%\hspace{-0.1in}{\bar F} &\stackrel{\gamma}{\longrightarrow}&  G \\
%\downarrow_{\pi_1} && & \stackrel{\Lambda+\theta\circ \pi_2}{\searrow}\\
%\pi_1({\bar F}) && & \hspace{0.1in} U(B)/CU(B)\\
%\downarrow_{\Pi^{\ddag}} &&\\
%U(D)/CU(D) && {\vspace{-0.1in}\nearrow_{\bt}}\\
%\end{array}
%}
%$$

\begin{proof}
The proof is exactly the same as that of Lemma 7.2 of \cite{LinTAI}.

Let $\kappa_1: U(D)/CU(D)\to K_1(D)\subset U(C)/CU(C)$ be  the quotient map. Let $\eta: \pi_1(U(C)/CU(C))\to K_1(D)$ be the map defined by $\eta=\kappa_1\circ\Pi^{\ddagger}|_{\pi_1(U(C)/CU(C))}$. Note that $\eta$ is an isomorphism{{----as we regard both $\pi_1(U(C)/CU(C))$ and $K_1(D)$ as the same subgroup of $U(C)/CU(C)$, this $\eta$ is the identity map.}} Since $\gamma$ is injective and $\gamma(\overline{F})$ is free, we conclude that $\kappa_1\circ\Pi^{\ddagger}\circ\pi_1$ is also injective on $\overline{F}$. Since $U_0(C)/CU(C)$ is divisible (6.6 of \cite{LinTAI}), there is a homomorphism $\lambda: K_1(D)\to U_0(C)/CU(C)$ such that
$$\lambda|_{\kappa_1\circ\Pi^\ddagger\circ\pi_1(\overline{F})}=\pi_0\circ((\kappa_1\circ\Pi^{\ddagger}\circ\pi_1)|_{\overline{F}})^{-1}{{,}}$$
{{ where $((\kappa_1\circ\Pi^{\ddagger}\circ\pi_1)|_{\overline{F}})^{-1}: ~~\eta\circ \pi_1(\overline{F})\to \overline{F}$ is the inverse map of the injective map $(\kappa_1\circ\Pi^{\ddagger}\circ\pi_1)|_{\overline{F}}$.}} This could be viewed as the following commutative diagram:
\begin{displaymath}
\xymatrix{
& K_1(D) \ar[dr]^\lambda& \\
(\eta\circ \pi_1)({\bar F}) \ar[rr]^{\pi_0\circ(\eta\circ\pi_1)^{-1}} \ar[ur] & & U_0(C)/CU(C). \\
}
\end{displaymath}
%$$
%\begin{array}{cccc}
%& K_1(D)  \\
%\nearrow &&  \stackrel{\lambda}{\searrow} \\
%\hspace{-0.1in}(\eta\circ \pi_1)({\bar F}) &\stackrel{\pi_0\circ (\eta\circ \pi_1)^{-1}} {\longrightarrow}& \hspace{0.2in}U_0(C)/CU(C)\\
%%\hspace{0.4in}\searrow_{(\eta\circ \pi_1)^{-1}}&& \nearrow_{\pi_0}\\
%%& {\bar F}\\
%\end{array}
%$$

Define $$\beta=\Lambda((\eta^{-1}\circ\kappa_1)\oplus(\lambda\circ\kappa_1)).$$
Then, for any $\overline{w}\in\overline{F}$,
$$\beta(\Pi^{\ddagger}\circ\pi_1(\overline{w}))=\Lambda(\eta^{-1}(\kappa_1\circ\Pi^{\ddagger}(\pi_1(\overline{w})))\oplus \lambda\circ\kappa_1(\Pi^\ddagger(\pi_1(\overline{w}))))=\Lambda(\pi_1(\overline{w})\oplus\pi_0(\overline{w})).$$
Define
$\theta: \pi_2(U(C)/CU(C))\to U(B)/CU(B)$
by
$$\theta(x)=\Lambda(x^{-1})\rforal x\in\pi_2(U(C)/CU(C)).$$
Then
$$\beta(\Pi^{\ddagger}(\pi_1(\overline{w})))=\Lambda(\overline{w})\theta(\pi_2(\overline{w}))\rforal w\in F.$$

For the second part of the statement, one assume that $\Lambda(G)\subseteq U_0(B)/CU(B)$. Then $\Lambda(\pi_2(U(C)/CU(C)))$ is a torsion subgroup of $U_0(B)/CU(B)$. But $\mathrm{U}_0(B)/\mathrm{CU}(B)$ is torsion free by Lemma \ref{UCUdiv}, and hence $\theta=0$.
\end{proof}

\begin{lem}[See Lemma 7.3 of \cite{LinTAI}]\label{exp-length3}
Let $B\in\mathcal B_1$ be a separable simple C*-algebra, and let $C$ be as above. Let $\mathcal{U}\subset \text{U}(B)$ be a finite subset, and  let $F$ be the subgroup generated by $\mathcal{U}$ such that $\kappa_1(\bar{F})$ is free, where $\kappa_1: \text{U}(B)/\text{CU}(B)\to K_1(B)$ is the quotient map. Suppose that $\alpha :K_1(C)\to K_1(B)$ is an injective homomorphism and $L: \bar{F}\to \text{U}(C)/\text{CU}(C)$ is an injective homomorphism with ${L(\bar{F}\cap \text{U}_0(B)/\text{CU}(B))\subset \text{U}_0(C)/\text{CU}(C) }$ such that $\pi_1\circ L$ is one-to-one and $$\alpha\circ\kappa'_1\circ L(g)=\kappa_1(g)\quad\mbox{for all}\quad g\in\bar{F},$$ where $\kappa'_1: \text{U}(C)/\text{CU}(C)\to K_1(C)$ is the quotient map. Then there exists a homomorphism $\beta: \text{U}(C)/\text{CU}(C)\to \text{U}(B)/\text{CU}(B)$ with $\beta(\text{U}_0(C)/\text{CU}(C))\subset \text{U}_0(B)/\text{CU}(B)$ such that
$$\beta\circ L(f)=f\rforal f\in\bar{F}.$$
\end{lem}
\begin{proof}
The proof is exactly the same as that of Lemma 7.3 of \cite{LinTAI}.

Let $G$ be the preimage of $\alpha\circ\kappa_1'(U(C)/CU(C))$ under $\kappa_1$. So we have the short exact sequence
$$0\to U_0(B)/CU(B) \to G \to \alpha\circ\kappa_1'(U(C)/CU(C)) \to 0.$$
Since $U_0(B)/CU(B)$ is divisible, there is an injective homomorphism
$$\gamma: \alpha\circ\kappa_1'(U(C)/CU(C)) \to G$$
such that $\kappa_1\circ\gamma(g)=g$ for any $g\in\alpha\circ\kappa_1'(U(C)/CU(C))$. Since $\alpha\circ\kappa_1'\circ L(f)=\kappa_1(f)$ for any $f\in\overline{F}$, we have $\overline{F}\subseteq G$.  Moreover, note that
$$(\gamma\circ\alpha\circ\kappa_1'\circ L(f))^{-1} f\in U_0(B)/CU(B)\rforal f\in \overline{F}.$$
Define $\psi: L(\overline{F})\to U_0(B)/CU(B)$ by
$$
\psi(x)=\gamma\circ\alpha\circ\kappa_1'(x^{-1})L^{-1}(x)
%\psi(x)=\gamma\circ\alpha\circ\kappa_1'\circ L(((L)^{-1}(x))^{-1})L^{-1}(s)
$$
for $x\in L(\overline{F})$. Since $U_0(B)/CU(B)$ is divisible, there is a homomorphism $\tilde{\psi}: U(C)/CU(C)\to U_0(B)/CU(B)$ such that $\tilde{\psi}|_{L(\overline{F})}=\psi$. Now define
$$\beta(x)=\gamma\circ\alpha\circ\kappa_1'(x)\tilde{\psi}(x).$$
Hence $\beta(L(f))=f$ for $f\in \overline{F}$.
\end{proof}

%{\bf We will move this to section 11?----- I did not do this because I do not want mixed up the file. If you are now in position do it, please move it.}
\begin{lem}\label{FDCU}
Let $A$ be a unital separable \CA\, such that
$\{\rho_A([p]): p\in A \,\,\, {\rm projections}\}$ is dense in
real linear span of $\{\R\rho_A([p]): p\in A\,\,\, {\rm projections}\}.$
Then, for any finite dimensional \SCA\, $B\subset A,$
$U(B)\subset CU(A).$
\end{lem}
\begin{proof}
Let $u\in U(B).$  Since $B$ has finite dimensional,  $u=\exp(i h)$ for some $h\in B_{s.a.}.$ We may write
$h=\sum_{i=1}^n \lambda_i p_i,$ where $\lambda_i\in \R$ and $\{p_1, p_2,...,p_n\}$ is a set of mutually orthogonal
projections.  By the assumption and applying Proposition 3.6 of \cite{GLX-ER}, one has
$u\in CU(A).$
\end{proof}

\begin{lem}[See Lemma 7.4 of \cite{LinTAI}]\label{exp-length2}
Let $B\cong A\otimes U$, where $A\in\mathcal B_0$ and $U$ is an infinite dimensional  UHF-algebra. Let
$C=C_1\oplus C_2,$ where $C_1=PM_n(C(X))P$ with $X$ be as defined in \ref{bdbk-k1}, $C_2\in {\cal C}_0.$  Let $F$ be a group generated by a finite subset $\mathcal{U}\subset \text{U}(C)$ such that $(\pi_1)|_{\bar{F}}$ is one-to-one. Let $G$ be a subgroup containing $\bar{F}$, $\pi_1(\text{U}(C)/\text{CU}(C))$ and $\pi_2(\text{U}(C)/\text{CU}(C))$. Suppose that $\alpha: \text{U}(C)/\text{CU}(C)\to \text{U}(B)/\text{CU}(B)$ is a homomorphism  such that $\alpha(\text{U}_0(C)/\text{CU}(C))\subset \text{U}_0(B)/\text{CU}(B)$. Then, for any $\epsilon>0$, there { are $\delta>0$ and finite subset $\mathcal G\subseteq C$} satisfying the following: if $\phi=\phi_0\oplus\phi_1: C\to B$ {{ (by such decomposition, we mean there is a projection $e_0$ such that $\phi_0: C \to e_0Be_0$ and $\phi_1: C\to (1_B-e_0)B(1_B-e_0)$)}} is a $\mathcal G$-${\bf \delta}$-multiplicative completely positive linear contraction such that
\begin{enumerate}
\item  $\phi_0$ maps identity of each summand of $C$ to a projection,
%$\phi_0$ is a \hm,
%and  $\phi_1$ are $\mathcal G$-${\bf \delta}$-multiplicative,
\item $\mathcal G$ is sufficiently large and ${\bf \delta}$ is sufficiently small depending only on $F$ and $C$ (such that $\phi^\ddagger$ is well defined on a subgroup of $\text{U}(C)/\text{CU}(C)$ containing all of $\bar{F}$, $\pi_0(\bar{F})$, $\pi_1(\text{U}(C)/\text{CU}(C))$, and  $\pi_2(\text{U}(C)/\text{CU}(C))$),
\item $\phi_0$ is homotopic to a homomorphism with finite dimensional image,
$[\phi_0]|_{K_0(C)}$ is well defined
and $[\phi]|_{K_1(C)}=\alpha_*$, where $\af_*: K_1(C)\to K_1(B)$ is induced map,
\item $\tau(\phi_0(1_C))<\delta$ for all $\tau\in T(B)$ (assume $e_0=\phi_0(1_C)$),
\end{enumerate}
then there is a homomorphism $\Phi: C\to e_0Be_0$ such that
%\begin{enumerate}
%\item

{\rm (i)} $\Phi|_{C_1}$ is {homotopic to a \hm\,  with finite dimensional image}  and $(\Phi)_{*0}=[\phi_0]|_{K_0(C)}$ and

{\rm (ii)}  $\alpha(\bar{w})^{-1}(\Phi\oplus \phi_1)^\ddagger(\bar{w})=\bar{g_w}$ where $g_w\in \text{U}_0(B)$ and $\mbox{cel}(g_w)<\epsilon$ for any $w\in\mathcal{U}$.
%\end{enumerate}
\end{lem}
\begin{proof}

The argument is exactly the same as that of Lemma 7.4 of \cite{LinTAI} {{(see also \cite{EGL-AH} and \cite{NT})}} . Since the source algebra and target algebra
in this lemma are different from the ones in 7.4 of \cite{LinTAI}, we will repeat some of the argument here.
We will retain the notations in \ref{bdbk-k1}.  By \ref{LinTAI-72}, there are \hm s $\bt_1, \bt_2: U(D)/CU(D)\to U(B)/CU(B)$
with $\bt_i(U_0(D)/CU(D))\subset U_0(B)/CU(B)$ ($i=1,2$) and \hm s
$$
\theta_1,\theta_2: \pi_2(U(C)/CU(C))\to U(B)/CU(B)
$$
such that
\beq\label{74-e1}
\bt_1\circ \Pi^{\ddag}(\pi_1({\bar w}))=\af({\bar w})\theta_1(\pi_2({\bar w}))\andeqn
\bt_2\circ \Pi^{\ddag}(\pi_1({\bar w}))=\phi^{\ddag}({\bar w}^*)\theta_2(\pi_2({\bar w}))
\eneq
for all ${\bar w}\in {\bar F}.$
Moreover, $\theta_{{1}}(g)=\af(g^{{-1}})$ and $\theta_2(g)={{\phi}}^{\ddag}(g)$ for all $g\in \pi_2({\bar F}).$
Since ${{\phi}}_0$ is homotopic
%%%homtopic
 to a \hm\, with finite dimensional range, {{using condition (3),}}
we compute that
\beq\label{74-e2}
\theta_1(g)\theta_2(g)\in U_0(B)/CU(B)\rforal g\in \pi_2({\bar F}).
\eneq
Since $\pi_2(U(C)/CU(C))$ is torsion free and $U_0(B)/CU(B)$ is torsion free (see \ref{UCUdiv}), we conclude that
\beq\label{72-e3}
\theta_1(g)\theta_2(g)=1\rforal g\in \pi_2({\bar F}).
\eneq
To simplify notation, \wilog, we may write
$C=D\oplus C_0$ as in \ref{bdbk-k1}.  To keep the same notation as in the proof
of 7.4 of \cite{LinTAI}, we also write $D=C^{(1)}$ and $C_0=\oplus_{j=2}^{l_1} C^{(j)},$ where each
$C^{(j)}$ is in ${\cal C}_0$ and is {{minimal}} (cannot be written as finite direct sum of more than one
copies of \CA s in ${\cal C}_0$), or a \CA\, of the form $PM_n(C(Y))P$ with connected spectrum
$Y.$ Note $K_1(C^{(j)})=\{0\}$ for all $j\ge 2.$

We then proceed the proof of 7.4 of \cite{LinTAI} and construct $\Phi_1$ exactly the same way as
in the proof of 7.4 of \cite{LinTAI}.
% {{In the construction of $\Phi_1$, do we need the existence theorem such as 18.8 or 20.16. I think it is better to spell out the construction instead of refer to \cite{LinTai}.}}
We will keep the notation used in the proof of
7.4 of \cite{LinTAI}.
We will again use the fact, by Corollary \ref{Unotrosion}, the group $\mathrm{U}_0(B)/\mathrm{CU}(B)$ is torsion free, and by Theorem \ref{UCUiso}, the map $$U(eBe)/CU(eBe)\to U(B)/CU(B)$$ is an isomorphism, where $e\in B$ is a nonzero projection.
By the assumption, $\phi_0|_{C_0}$ is homotopy to $\phi_{00}: C_0\to \phi_0(1_C)B\phi_0(1_C)$
such that $B_0=\phi_{00}(C_0)$ is a finite dimensional \SCA\, in $(e_0-E_1)B(e_0-E_1)$ with
$1_{B_0}=e_0-E_1,$ where $E_1\le e_0$ is a non-zero projection such that $[E_1]=(\phi_0)_{*0}([1_D]).$
Let $\Phi_2: C\to B_0$ be defined by $\Phi_2(f, g)=\phi_{00}(g)$ for all $f\in  D$ and $g\in C_0.$
It is important to note (as the main difference from this lemma and that of 7.4 of \cite{LinTAI}) that, with the assumption,
by \ref{FDCU}, for any $w\in U(C),$ $\Phi_2(w)\in CU(B).$ In other words, the general case could be reduced to the case
that $C=D.$
The rest of the proof is exactly the same as that of 7.4 of \cite{LinTAI}.

%
%By Lemma \cite{LinTAI-72}, there are homomorphisms $\beta_1, \beta_2: U(D)/CU(D)\to U(B)/CU(B)$ with
%$$\beta_i(U_0(D))\subseteq U_0(B)/CU(B),\quad i=1, 2$$
%and homomorphisms
%$\theta_1, \theta_2: \pi_2(U(C)/CU(C))\to U(B)/CU(B)$
%such that
%$$\beta\circ\Pi^{\ddager}(\pi(\bar{w}))=\alpha(\bar{w})\theta_1(\pi(\bar{w}))\quad\textrm{and}\quad \beta\circ\Pi^{\ddager}(\pi(\bar{w}))=\phi_1^{\ddager}(\bar{w}^*)\theta_2(\pi(\bar{w}))$$
%for all $\bar{w}\in\overline{F}$.
%Moreover, $\theta_1(g)=\alpha(g^{-1})$ and $\theta_2(g)=\phi_1^{\ddager}(g)$ if $g\in\pi_2(\overline{F})$. Since $\phi_0$ is homotopic to a map with finite dimensional image, one has that $\theta_1(g)\theta_2(g)\in U_0(B)/CU(B)$ for all $g\in \pi_2(\overline{F})$. Since $\pi_2(U(C)/CU(C))$ is torsion and $U_0(B)/CU(B)$ is torsion free, one concludes that $\theta_1(g)\theta_2(g)=\bar{1}$ for all $g\in \pi_2(\overline{F})$.  Write $C=\bigoplus_{l=1}^{l_1+1}$
%
\end{proof}

\begin{lem}[See Lemma 7.5 of \cite{LinTAI}]\label{exp-length1}
Let $B\cong A\otimes U$, where $A \in\mathcal B_1$ and $U$ an {{infinite dimensional}}  UHF-algebra. Let $\mathcal{U}\subset \text{U}(B)$ be a finite subset and $F$ be the subgroup generated by $\mathcal{U}$ such that $\kappa_1(\bar{F})$ is free, where $\kappa_1: \text{U}(B)/\text{CU}(B)\to K_1(B)$ is the quotient map. Let $C$ be as above and let $\phi: C\to B$ be a homomorphism such that $(\phi)_{*1}$ is one-to-one. Suppose that $j, L: \bar{F}\to \text{U}(C)/\text{CU}(C)$ are two injective homomorphisms with $j(\overline{F\cap U_0(B)})$, $L(\overline{F\cap U_0(B)})\subset \text{U}_0(C)/\text{CU}(C)$ such that $\kappa_1\circ\phi^\ddagger\circ L=\kappa_1\circ\phi^\ddagger\circ j=\kappa_1|_{\bar{F}}$, and they are one-to-one.

Then, for any $\epsilon>0$, there exists $\delta>0$ such that if {{$\phi$ can be decomposed as}}
%{\red {It seems readers may think that $\delta$ depends on $\phi$ and the condition (1) of $\phi$ below depends on $\delta$,so it is a cyclic depending. I guess when we apply the lemma, our $\phi$ can be decomposed as such for any $\delta$, so I add above to avoid the confusion--Gong}}
$\phi=\phi_0\oplus\phi_1: C\to B$, where $\phi_0$ and $\phi_1$ are homomorphisms satisfying the following:
\begin{enumerate}
\item $\tau(\phi_0(1_C))<\delta$ for all $\tau\in T(B)$ and
\item $\phi_0$ is {homotopic to a \hm\, with finite dimensional image},
\end{enumerate}
then there is a homomorphism $\psi: C\to e_0Be_0$ ($e_0=\phi_0(1_C)$) such that
\begin{enumerate}\setcounter{enumi}{3}
\item  $[\psi]=[\phi_0]$ in $KL(C, B)$ and
\item  $(\phi^\ddagger\circ j(\bar{w}))^{-1}(\psi\oplus\phi_1)^\ddagger(L(\bar{w}))=\bar{g_w}$ where $g_w\in \text{U}_0(B)$ and $\mbox{cel}(g_w)<\epsilon$ for any $w\in\mathcal{U}$.
\end{enumerate}
\end{lem}
\begin{proof}
The proof is the same as that of Lemma 7.5 of \cite{LinTAI}. Note that instead of using Lemma 7.4 of \cite{LinTAI}, one uses Lemma \ref{exp-length2}.
\end{proof}

\begin{rem}
Roles of Lemma \ref{exp-length2} and Lemma \ref{exp-length1} played   in the proof of the following  isomorphism theorem (\ref{IST0})  are the same as those in the proof of Theorem 10.4 in \cite{LinTAI}.
\end{rem}

{The following statement is well known. For the reader's convenience, we include a proof.
\begin{lem}\label{indlim-inv}
Let $(A_n, \phi_{n, n+1})$ be a unital inductive sequence of separable C*-algebras, and denote by $A=\varinjlim A_n$. Assume that $A$ is amenable. Let $\mathcal F\subseteq A$ be a finite subset, and let $\epsilon>0$. Then there is  an integer $m\ge 1$ and a unital completely positive linear map $\Psi:A\to A_m$ such that
$$
\|\phi_{m, \infty}\circ\Psi(f)-f\|<\epsilon
\rforal f\in\mathcal F.
$$
\end{lem}
\begin{proof}

%Without loss of generality, one may assume that $\mathcal F\subseteq \phi_{n, \infty}(A_n)$.
Regard $A$ as the C*-subalgebra of $\prod A_n/\bigoplus A_n$ generated by the equivalence classes of the sequences $(x_1, x_2, ..., x_n, ...)$ satisfying {{that there is $N$ with}}
$x_{n+1}=\phi_n(x_n),\quad n={{N, N+1}}, ...\ .$
%Note that $$\phi_{n, \infty}(A_n)=\{\overline{(\underbrace{0, ..., 0}_{n-1}, a, \phi_{n, n+1}(a), \phi_{n, n+2}(a), ...)}: a\in A_n\}.$$
Since $A$ is amenable, by the Choi-Effros lifting theorem, there is a unital completely positive linear map $\Phi: A\to \prod A_n$ such that
$\pi\circ\Phi=\id_A,$
where $\pi$ is the quotient map. In particular, this implies that
\begin{equation}\label{coincidelim}
\lim_{k\to\infty}\|\pi_k\circ\Phi(a)-a_k\|=0,
\end{equation}
if $a=\overline{(a_1, a_2, ..., a_k, ...)}\in A$.

Write $\mathcal F=\{f_1, f_2, ..., f_l\}$, and for each $f_i$, fix a representative
$$f_i=\overline{(f_{i, 1}, f_{i, 2}, ..., f_{i, k}, ...)}.$$ In particular
\begin{equation}\label{convlim}
\lim_{k\to\infty} \phi_{k, \infty}(f_{i, k})=f_i.
\end{equation}
Then, for each $f_i$, one has
\begin{eqnarray*}
&&\limsup_{k\to\infty}\|\phi_{k, \infty}\circ\pi_k\circ\Phi(f_i)-f_i\| \\
& = & \limsup_{k\to\infty} \|\phi_{k, \infty}\circ\pi_k\circ\Phi(f_i)-\phi_{k, \infty}(f_{i, k})\| \quad(\textrm{by \eqref{convlim}}) \\
&\leq&\limsup_{k\to\infty}\|\pi_k\circ\Phi(f_i)-f_{i, k}\|
= 0\quad \textrm{(by \eqref{coincidelim}}).
\end{eqnarray*}
There then {{exists}} $m\in \N$ such that
$$\|\phi_{m, \infty}\circ\pi_m\circ\Phi(f_i)-f_i\|<\epsilon,\quad 1\leq i\leq l.$$
Thus, the unital completely positive linear map $$\Psi:=\pi_m\circ\Phi$$
satisfies the lemma.
\end{proof}
}

\begin{thm}\label{IST0}
Let $A_1\in \mathcal B_{{0}}$ be a unital separable simple C*-algebra satisfying the UCT, and denote by $A=A_1\otimes U$ for a UHF-algebra {$U$} of infinite type. Let $C$ be a C*-algebra in Theorem \ref{RangT}.
%{{\rm (} and
%\ref{ReRangT}).}{\color{Green} (I can't find it.)}
If ${\rm Ell}(C)\cong {\rm Ell}(A)$, then there is an isomorphism $\phi: C\to A$ which carries the identification of ${\rm Ell}(C) \cong {\rm Ell}(A)$.

Moreover, if there is a homomorphism $\Gamma: {\rm Ell}(C)\to {\rm Ell}(A)$, then there is a *-homomorphism $\phi: C\to A$ such that $\phi$ induces $\Gamma$.
\end{thm}
\begin{proof}
We only prove  the first part of the statement. The second part can be proved in a similar way (one only has to do one-sided intertwining arguments in this case).

Let $\alpha\in KL(C, A)$ {with $\alpha^{-1}\in KL(A, C)$} and $\gamma: T(A)\to T(C)$ be given by the isomorphism
${\rm Ell}(C)\cong {\rm Ell}(A).$

Assume that $C=\varinjlim (C_n, \iota_n)$
be as in \ref{RangT}.
%and {\ref{ReRangT}}.
Let $\mathcal G_1\subseteq \mathcal G_2\subseteq \cdots \subseteq C$ and $\mathcal F_1\subseteq \mathcal F_2\subseteq\cdots \subseteq A$ be increasing sequences of finite subsets with dense union. Let $\ep_1>\ep_2> \cdots >0$ be a decreasing sequence of positive numbers with finite sum.

We will repeatedly apply Theorem \ref{MUN1}. Let $\delta_c^{(1)}>0$ (in place of $\delta$), $\mathcal G^{(1)}_c\subseteq C$ (in place of $\mathcal G$), $\sigma^{(1)}_{c, 1}, \sigma^{(1)}_{c, 2}>0$ (in place of $\sigma_1$ and $\sigma_2$ respectively), $\mathcal P^{(1)}_c\subseteq\underline{K}(C)$ (in place of $\mathcal P$), $\mathcal U^{(1)}_c\subseteq U(C)$ (in place of $\mathcal U$)
and $\mathcal H^{(1)}_c\subseteq C_{s.a}$ (in place of $\mathcal H_2$) be as required by \ref{MUN1} for $C$ (in the place of $A$ with $X$ being a point), $\ep_1$ (in the place of $\ep$), $\mathcal G_1$ (in the place of $\mathcal F$).  Note, in this case,
$X$ is a point, by Remark \ref{ReMUN1}, we do not introduce the map $\Delta$ and ${\cal H}_1.$

As in the remark of \ref{ReMUN1}, with sufficiently large $n\ge 1,$ we may assume that ${\cal U}_c^{(1)}$ is in the image of $U(C_n)$ for some large  $n\ge 1.$ Moreover, as \ref{ReMUN1}, we let ${\cal U}_c^{(1)}\subset U(C)$ (instead in
$U(M_2(C))$).

Denote by $F_c\subseteq U(C)$ the subgroup generated by $\mathcal U^{(1)}_c$. Write $\overline{F_c}=(\overline{F_{c}})_0 \oplus \mathrm{Tor}(\overline{F_c})$ according to the decomposition described in \ref{bdbk-k1}, where $(\overline{F_{c}})_0$ is torsion free.
Without loss of generality (by choosing smaller $\sigma_{c,2}^{(1)}$), one may assume that
$$\mathcal U^{(1)}_c=\mathcal U^{(1)}_{c, 0}\sqcup\mathcal U^{(1)}_{c, 1}$$ where $\overline{\mathcal U^{(1)}_{c, 0}}$ generates $(\overline{F_{c}})_0$ and $\overline{\mathcal U^{(1)}_{c, 1}}$ generates $\mathrm{Tor}(\overline{F_c})${{---namely, we can choose $\mathcal U^{(1)}_{c, 0}$ and $\mathcal U^{(1)}_{c, 1}$ so that $\mathcal U^{(1)}_c\subset\mathcal U^{(1)}_{c, 0}\cdot\mathcal U^{(1)}_{c, 1}$; then by choosing smaller $\sigma_{c,2}^{(1)}$ one can replace $\mathcal U^{(1)}_c$ by  $\mathcal U^{(1)}_{c, 0}\sqcup\mathcal U^{(1)}_{c, 1} .$}} Note that for each $u\in \mathcal U^{(1)}_{c, 1}$, one has that $u^k \in CU(C)$, where $k$ is the order of $\overline{u}$.

By Theorem \ref{MEST}, there is a $\mathcal G^{(1)}_c$-$\delta^{(1)}_c$-multiplicative map $L_1: C\to A$ such that
\beq\label{eq-kk-001}
&&[L_1]|_{\mathcal P^{(1)}_c}=\alpha|_{\mathcal{P}^{(1)}_c}\andeqn\\\label{eq-tr-001}
&& |\tau\circ L_1(f)-\gamma(\tau)(f)|<\sigma^{(1)}_{c, 1}/3\rforal f\in\mathcal H^{(1)}_c\rforal \tau\in T(A).
 \eneq
 Without loss of generality, one may assume that $L_1^\ddagger$ is well defined and injective on $(\overline{F_c})_0$.  Moreover, if $k$ is the order of $\overline{u}, $  one may also assume that
 $$\mathrm{dist}(L_1(u^k), CU(A))<\sigma_{c, 2}^{(1)}/2\rforal u\in\mathcal U_{c, 2}^{(1)}.$$
 %where $k$ is the order of $\overline{u}$.

To apply  Theorem \ref{MUN1}  second time, let $\delta^{(1)}_a>0$ (in the place of $\delta$), $\mathcal G^{(1)}_a\subseteq {{C}}$ (in the place of $\mathcal G$), $\sigma^{(1)}_{a, 1}, \sigma^{(1)}_{a, 2}>0$ (in the place of $\sigma_1$ and $\sigma_2$ respectively), $\mathcal P^{(1)}_a\subseteq\underline{K}(A)$ (in the place of $\mathcal P$), $\mathcal U^{(1)}_a\subseteq U(A)$ (in the place of $\mathcal U$) and $\mathcal H^{(1)}_a\subseteq A_{s.a}$ (in the place of $\mathcal H_2$) be as required by Theorem \ref{MUN1} for
$A$ (in the place of $A$ with $X$ being a point), $\ep_1$ (in the place of $\ep$), and $\mathcal F_1$ (in the place of $\mathcal F$).  Again, note that $X$ is a point, in the application of Theorem \ref{MUN1}.

Denote by $F_a\subseteq U(A)$ the subgroup generated by $\mathcal U_a^{(1)}$. Since ${{\mathcal U_a^{(1)}}}$ is finite, we can write $\overline{F_a}=(\overline{F_{a}})_0 \oplus \mathrm{Tor}(\overline{F_a})$, where $(\overline{F_{a}})_0$ is torsion free. Fix this decomposition. Without loss of generality (by choosing smaller $\sigma_{{{a}},2}^{(1)}$),
one may assume that
$$\mathcal U^{(1)}_a=\mathcal U^{(1)}_{a, 0}\sqcup\mathcal U^{(1)}_{a, 1}$$ where ${ {\mathcal U}}^{(1)}_{a, 0}$ generates $(\overline{F_{a}})_0$ and $\mathcal U^{(1)}_{a, 1}$ generates $\mathrm{Tor}(\overline{F_{{a}}})$.
{ {By enlarging $\mathcal P^{(1)}_a$, we can assume $\mathcal P^{(1)}_a\supset \kappa_{1, A}((\overline{F_{a}})_0)$(in $ K_1(A)$), where $\kappa_{1, A}: U(A)/CU(A)\to K_1(A)$
is the quotient map.}}
Note that for each $u\in \mathcal U^{(1)}_{a, 1}$, one has that $u^k \in CU(A)$, where $k$ is the order of $\overline{u}$.

By Theorem \ref{MEST} and amenability of $C,$  there are a finite subset ${\cal G}'\supset \mathcal G^{(1)}_a,$ a positive number $\dt'<\dt_a^{(1)},$ a sufficiently large integer $n\ge 1$ and there is a $\mathcal G' $-$\delta'$-multiplicative map $\Phi'_1: A\to C_n$ such that, {{with $\sigma_{a,c,1}^{(1)}=\min\{\sigma^{(1)}_{a, 1}, \sigma^{(1)}_{c, 1}\}/3,$}}
\beq\label{eq-kk-002}
&&\hspace{-0.2in}[\imath_n\circ \Phi'_1]_{\mathcal P^{(1)}_a\cup [L_1](\mathcal P^{(1)}_c)}=\alpha^{-1}|_{\mathcal P^{(1)}_a \cup [L_1](\mathcal P^{(1)}_c)}\andeqn\\
%\end{equation}
%and
%\begin{equation}
\label{eq-tr-002}
&&\hspace{-0.2in}|\tau\circ \imath_n\circ \Phi'_1(f)-\gamma^{-1}(\tau)(f)|<\sigma_{a,c,1}^{(1)}
%{\min\{\sigma^{(1)}_{a, 1}, \sigma^{(1)}_{c, 1}\}\over{3}}
\rforal f\in\mathcal H^{(1)}_a\cup L_1(\mathcal H^{(1)}_c)\,\,\,{\rm and}\,\,\tau\in T(C).
\eneq
%\end{equation}
Moreover, one  may assume that $\Phi'_1\circ L_1$ is $\mathcal G^{(1)}_c$-$\delta^{(1)}_c$-multiplicative, $(\Phi_1')^\ddagger$ is defined and injective on $(\overline{F_a})_0$, and $(\Phi'_1\circ L_1)^\ddagger$ is well defined and injective on $(\overline{F_c})_0$.
Furthermore,  one may also assume that
 \begin{equation}\label{almost-cd-01}
 \mathrm{dist}(\imath_n\circ \Phi_1'(u^k), CU(C))<\sigma_{a, 2}^{(1)}/2\rforal u\in \mathcal U_{a, 1}^{(1)},
 \end{equation}
 where $k$ is the order of $\overline{u}$, and
 \begin{equation}\label{almost-cd-02}
 \mathrm{dist}((\imath_n\circ \Phi_1'\circ L_1)(v^{k'}), CU(A))<\sigma_{c, 2}^{(1)}\rforal v\in \mathcal U_{c, 1}^{(1)},
 \end{equation}
 where $k'$ is the order of $\overline{v}$.
It then follows from \eqref{eq-kk-001} and \eqref{eq-kk-002} that
\begin{equation}\label{comp-001}
[\imath_n\circ \Phi'_1\circ L_1]|_{\mathcal P^{(1)}_c}=[\id]|_{\mathcal P^{(1)}_c};
\end{equation}
and it follows from \eqref{eq-tr-001} and \eqref{eq-tr-002}  that
\begin{equation}\label{comp-002}
|\tau\circ\imath_n\circ \Phi'_1\circ L_1(f)-\tau(f)|<2\sigma^{(1)}_{c, 1}/3\rforal f\in\mathcal H^{(1)}_{c} \rforal \tau\in T(C).
\end{equation}

Recall that $(\overline{F_c})_0\subseteq U(C)/CU(C)$ is the subgroup generated by $\overline{\mathcal U^{(1)}_{c, 0}}$.  Since we have assumed that $\overline{\mathcal U^{(1)}_{c, 0}}$ is in the image of $U(C_n)/CU(C_n),$
there is an injective homomorphism
$j: (\overline{F_c})_0 \to U(C_n)/CU(C_n)$
such that
\begin{equation}\label{june-29-n-001}
\iota_{n}^\ddagger\circ j=\id|_{\overline{(F_c)_0}}.
\end{equation}
Moreover, by (\ref{comp-001}),  $\kappa_{1,C}\circ\iota_n^\ddagger\circ(\Phi_1'\circ L_1)^\ddagger|_{\overline{(F_c)_0}}=\kappa_{1,A}\circ\iota_n^\ddagger\circ j=\kappa_{1, C}|_{\overline{(F_c)_0}},$ where
$\kappa_{1,C}: U(C)/CU(C)\to K_1(C)$ be the quotient map.

Let $\delta$ be the constant of Lemma \ref{exp-length1} with respect to $C_n$ (in place of $C$), $C$ (in place of $B$), $\sigma_{c, 2}$ (in place of $\ep$), $\iota_n$ (in place of $\phi$), $j$ and $(\Phi_1'\circ L_1)^\ddagger |_{\overline{(F_c)_0}}$ (in place of $L$). By the construction of $C$, one has a decomposition ${ \iota_n=\iota_n^{(0)}\oplus\iota_n^{(1)}}$ such that
\begin{enumerate}
\item $\tau(\iota_n^{(0)}(1_{C_n}))<\min\{\delta, \sigma^{(1)}_{c, 1}/3\}$ for all $\tau\in T(C)$, and
\item $\iota_n^{(0)}$ has finite dimensional range.
%\end{enumerate}
Then, by Lemma \ref{exp-length1}, there is a homomorphism $h: C_n\to e_0Ce_0$, where $e_0=\iota_n^{(0)}(1_{C_n})$, such that
%\begin{enumerate}
\item\label{purt-sm-cor} $[h]=[\iota_n^{(0)}]$ in $KL(C_n, C)$, and
\item for each $u\in\mathcal U^{(1)}_{c, 0}$, one has that
\begin{equation}\label{eq-ak-001}
(\iota_n^\ddagger\circ j(\overline{u}))^{-1}(h\oplus\iota_n^{(1)})^\ddagger((\Phi'_1\circ L_1)^\ddagger(\overline{u}))=\overline{g_u}\end{equation}
for some $g_u\in U_0(C)$ with $\textrm{cel}(g_u)<\sigma^{(1)}_{c, 2}$.
\end{enumerate}

Define  $\Phi_1=(h\oplus \iota_n^{({{1}})})\circ\Phi'_1$. By \eqref{eq-kk-002} and \eqref{purt-sm-cor}, one has
\begin{equation}\label{June20-n-002}
[\Phi_1]_{\mathcal P^{(1)}_a\cup [L_1](\mathcal P^{(1)}_c)}=\alpha^{-1}|_{\mathcal P^{(1)}_a \cup [L_1](\mathcal P^{(1)}_c)}
\end{equation}

Note that $\Phi_1$ is still $\mathcal G_a^{(1)}$-$\delta_a^{(1)}$-multiplicative, and hence \eqref{almost-cd-01} and \eqref{almost-cd-02} still hold with $\Phi_1'$ replaced by $\Phi_1$. That is
\begin{equation}\label{almost-cd-001}
 \mathrm{dist}(\Phi_1(u^k), CU(C))<\sigma_{a, 2}^{(1)}/2\rforal u\in \mathcal U_{a, 1}^{(1)},
 \end{equation}
 where $k$ is the order of $\overline{u}$, and
 \begin{equation}\label{almost-cd-002}
 \mathrm{dist}((\Phi_1\circ L_1)(v^{k'}), CU(A))<\sigma_{c, 2}^{(1)}\rforal v\in \mathcal U_{c, 1}^{(1)},
 \end{equation}
 where $k'$ is the order of $\overline{v}$.
By \eqref{comp-001} and \eqref{comp-002}, one has
%\begin{equation}
\beq\label{uniq-kk-001}
[\Phi_1\circ L_1]|_{\mathcal P^{(1)}_c}=[\id]|_{\mathcal P^{(1)}_c}\andeqn\\
%\end{equation}
%and
%\begin{equation}
\label{uniq-tr-001}
|\tau\circ\Phi_1\circ L_1(f)-\tau(f)|<\sigma^{(1)}_{c, 1}\rforal f\in\mathcal H^{(1)}_{c} \rforal \tau\in T(C).
\eneq
%\end{equation}
Moreover, for any $u\in \mathcal U^{(1)}_{c, 0}$, one has (by \eqref{june-29-n-001} and \eqref{eq-ak-001})
\begin{eqnarray}\label{alg-k1-c-tf}
(\Phi_1 \circ L_1)^\ddagger(\overline{u})&=&(\iota_n^\ddagger\circ j(\overline{u}))\cdot\overline{g_u}=\overline{u}\cdot\overline{g_u}\approx_{\sigma^{(1)}_{c, 2}}\overline{u}.
\end{eqnarray}
Let $u\in \mathcal{U}_{c, 1}^{(1)}$ with order $k$. By \eqref{almost-cd-002}, there is a self-adjoint element $b\in C$ with $||b||<\sigma_{c, 2}^{(1)}$ such that
\begin{equation*}
(u^*)^k(\Phi_1 \circ L_1)(u^k) \exp(2\pi i b)\in CU(C),
\end{equation*}
{{(where, we notice that  $ (u^*)^k\in CU(C)$)}} and hence
\begin{equation*}
((u^*)(\Phi_1 \circ L_1)(u)\exp(2\pi i b/k))^k\in CU(C).
\end{equation*}
Note that
$$(u^*)(\Phi_1 \circ L_1)(u)\exp(2\pi i b/k)\in U_0(C)$$
and $U_0(C)/CU(C)$ is torsion free (Corollary \ref{Unotrosion}). One has that
$$(u^*)(\Phi_1 \circ L_1)(u)\exp(2\pi i b/k)\in CU(C).$$
In particular, this implies that
\begin{equation}\label{alg-k1-c-tor}
\mathrm{dist}((\Phi_1 \circ L_1)^{\ddagger}({\bar u}), {\bar u})<\sigma_{c, 2}^{(1)}/k \rforal  {\bar u}\in \mathcal{U}_{c, 1}^{(1)}
\end{equation}
Together with \eqref{alg-k1-c-tf}, one has that
\begin{equation}\label{alg-k1-c}
\mathrm{dist}((\Phi_1 \circ L_1)^{\ddagger}({\bar u}), {\bar u})<\sigma_{c, 2}^{(1)} \rforal u\in \mathcal U_{c}^{(1)}.
\end{equation}

Therefore, by \eqref{uniq-kk-001}, \eqref{uniq-tr-001}, and \eqref{alg-k1-c}, applying Theorem \ref{MUN1},
one obtains a unitary $U_1$ such that
$$||U^*(\Phi_1\circ L_1 (f))U-f||<\ep_1 \rforal f\in\mathcal G_1.$$
By replacing $\Phi_1$  by $\mathrm{Ad}(U)\circ \Phi_1$, \wilog,  one may assume that
 $$||\Phi_1\circ L_1 (f)-f||<\ep_1 \rforal f\in\mathcal G_1.$$ In other words, one has the following diagram:
\begin{displaymath}
\xymatrix{
C \ar[r]^{\id} \ar[d]_{L_1} & C\\
A \ar[ur]_{\Phi_1}
}
\end{displaymath}
which is approximately commutative on the subset ${\cal G}_1$ within $\ep_1.$

We will continue to  apply Theorem \ref{MUN1}.
Let $\delta^{(2)}_c>0$ (in the place of $\delta$), $\mathcal G^{(2)}_c\subseteq C$ (in the place of $\mathcal G$), $\sigma^{(2)}_{c, 1}, \sigma^{(2)}_{c, 2}>0$ (in the place of $\sigma_1$ and $\sigma_2$ respectively), $\mathcal P^{(2)}_c\subseteq\underline{K}(C)$ (in the place of $\mathcal P$), $\mathcal U^{(2)}_c\subseteq U(C)$ (in the place of $\mathcal U$) and $\mathcal H^{(2)}_c\subseteq C_{s.a}$ (in the place of $\mathcal H_2$) be as required by Theorem \ref{MUN1}
for  $C$ (in the place of $A$ with $X$ being a point), $\ep_2$ (in the place of $\ep$), and $\mathcal G_2$ (in the place of $\mathcal F$).

Denote by $F^{(2)}_c\subseteq U(C)$ the subgroup generated by $\mathcal U^{(2)}_c$. Since ${\cal U}^{(2)}$ is finite,  we can write $\overline{F^{(2)}_c}=(\overline{F^{(2)}_{c}})_0 \oplus \mathrm{Tor}(\overline{F^{(2)}_c})$, where $(\overline{F^{(2)}_{c}})_0$ is torsion free. Fix this decomposition. Without loss of generality (by choosing smaller $\sigma_{c,2}^{(2)}$), one may assume that
$$\mathcal U^{(2)}_c=\mathcal U^{(2)}_{c, 0}\sqcup\mathcal U^{(2)}_{c, 1}$$ where $\overline{\mathcal U^{(2)}_{c, 0}}$ generates $(\overline{F^{(2)}_{c}})_0$ and $\overline{\mathcal U^{(2)}_{c, 1}}$ generates $\mathrm{Tor}(\overline{F^{(2)}_c})$.
Note that, for each $u\in \mathcal U^{(2)}_{c, 1}$, one has that $u^k \in CU(C)$, where $k$ is the order of $\overline{u}$.

There are a finite subset ${\cal G}_0\subset C$ and a positive number $\dt_0>0$  such that, for any two
${\cal G}_0$-$\dt_0$-multiplicative \morp s $L_1'', L_2'': C\to A,$\\ {{if $\|L_1''(c)-L_2''(c)\|<\dt_0\rforal c\in {\cal G}_0$, then}}
%\begin{equation*}%\label{June20-n1}
\beq\nonumber
&&[L_1'']|_{{\cal P}_c^{(2)}}=[L_2'']|_{{\cal P}_c^{(2)}}\andeqn\\\nonumber
%\end{equation*}
%and
%\begin{equation*}%\label{June20-n1-1}
&&\hspace{-0.1in}|\tau\circ L_1''(h)-\tau\circ L_2''(h)|<\min\{\sigma^{(2)}_{c, 1}/3, \sigma^{(1)}_{a, 1}/3 \}/2 \rforal h\in \mathcal H^{(2)}_c\cup\Phi_1(\mathcal H^{(1)}_a) \rforal \tau\in T(A){{.}}
\eneq
%\end{equation*}
%provided that
%\begin{equation*}%\label{June24-1}
%\|L_1''(c)-L_2''(c)\|<\dt_0\rforal c\in {\cal G}_0.
%\end{equation*}

Note, by Lemma \ref{indlim-inv}, for any finite subset ${\cal G}''\subset C$ and any $\dt''>0,$ there exists a
large $m$ and a unital \morp\, $L_{0,2}: C\to C_m$ such that
\begin{equation}\label{June20-2}
\|\imath_m\circ L_{0,2}(g)-g\|<\dt''\rforal g\in {\cal G}''.
\end{equation}
Let $\kappa_{1,C_m}: U(C_m)/CU(C_m)\to K_1(C_m)$ and $\kappa_{1, A}: U(A)/CU(A)\to K_1(A)$
be the quotient maps, respectively.
We may assume that, with sufficiently large ${\cal G}''$ and sufficiently small $\dt'',$  $(L_{0,2}\circ \Phi_1)^\ddagger$ is defined and injective on $\overline{(F_a)_0}$, and moreover,
$$\kappa_{1, C_m}\circ(L_{0, 2}\circ\Phi_1)^{\ddagger}(g)=[L_{0, 2}\circ\Phi_1](\kappa_{1, A}(g)), \quad g\in\overline{(F_a)_0},$$
and by \eqref{June20-n-002}, for any $g\in\overline{(F_a)_0}$ {{(note that $\mathcal P^{(1)}_a\supset \kappa_{1, A}((\overline{F_{a}})_0)$(in $K_1(A)$)),}}
\begin{eqnarray*}
\alpha \circ [\imath_m] \circ [L_{0, 2}\circ\Phi_1](\kappa_{1,A}(g)) & = & \alpha \circ [\imath_m \circ L_{0, 2}]\circ[\Phi_1](\kappa_{1,A}(g))\\
&=& \alpha \circ [\Phi_1](\kappa_{1,A}(g))
= \kappa_{1, A}(g).
\end{eqnarray*}
Hence,
$$\alpha\circ[\iota_m]\circ \kappa_{1, C_m}\circ(L_{0, 2}\circ\Phi_1)^{\ddagger}(g)=\alpha\circ[\iota_m]\circ[L_{0, 2}\circ\Phi_1](\kappa_{1, A}(g))=\kappa_{1, A}(g) \rforal g\in\overline{(F_a)_0}.$$
It follows from Lemma \ref{exp-length3} that there is a homomorphism $\beta: U(C_m)/CU(C_m)\to U(A)/CU(A)$ with $\beta(U_0(C_m)/CU(C_m))\subseteq U_0(A)/CU(A)$ such that
\begin{equation}\label{eq-ak-lift}
\beta\circ (L_{0,2}\circ \Phi_1)^\ddagger(f)=f\rforal f\in\overline{(F_a)_0}.
\end{equation}

By Theorem \ref{MEST}, there is a
${\cal G}''$-$\dt''$-multiplicative map $L'_2: C\to A$ such that
\beq
%\begin{equation*}%\label{eq-kk-003}
&&[L'_2]|_{\mathcal P^{(2)}_c}=\alpha|_{\mathcal{P}^{(2)}_c}\andeqn\\
%\end{equation*}
 %and
% \begin{equation*}%\label{eq-tr-003}
&& |\tau\circ L'_2(f)-\gamma(\tau)(f)|<\min\{\sigma^{(2)}_{c, 1}/3, \sigma^{(1)}_{a, 1}/3 \}/2\rforal f\in\mathcal H^{(2)}_c\cup\Phi_1(\mathcal H^{(1)}_a)
% \end{equation*}
\eneq
 and for all  $\tau\in T(A).$
 We may choose that
 \begin{equation*}%\label{eq-June24-2}
 {\cal G''}\supset {\cal G}_0\cup \mathcal G_c^{(2)}\andeqn \dt''<\min\{\dt_0, \delta_c^{(2)}\}.
 \end{equation*}
 Define $L_2'': C\to A$ by
 $L_2''(c)=L_2'\circ \imath_m\circ L_{0,2}.$ Then, by (\ref{June20-2}), since $L_2'$ is contractive,
 \begin{equation*}%\label{June24-4}
 \|L_2''(c)-L_2'(c)\|<\dt_0\rforal c\in {\cal G}_0.
 \end{equation*}
 It follows that
 \begin{eqnarray*}%\label{June24-5}
&& [L_2'']|_{\mathcal P^{(2)}_c}=[L'_2]|_{\mathcal P^{(2)}_c}=\alpha|_{\mathcal{P}^{(2)}_c}\andeqn\\
 && |\tau\circ L''_2(f)-\gamma(\tau)(f)|<\min\{\sigma^{(2)}_{c, 1}/3, \sigma^{(1)}_{a, 1}/3 \}\rforal f\in\mathcal H^{(2)}_c\cup\Phi_1(\mathcal H^{(1)}_a)
 \end{eqnarray*}
 and for all $\tau\in T(A).$

By choosing ${\cal G}''$ sufficiently large and $\delta''$ sufficiently small, one may assume that $L_2''\circ\Phi_1$ is $\mathcal G_a^{(1)}$-$\delta_a^{(1)}$-multiplicative, and
\begin{equation*}%\label{almost-cd-003}
 \mathrm{dist}((L''_2\circ\Phi_1)(u^k), CU(A))<\sigma_{a, 2}^{(1)}\rforal u\in \mathcal U_{a, 1}^{(1)},
 \end{equation*}
 where $k$ is the order of $\overline{u}$ (see \eqref{almost-cd-002}).

Moreover, by the construction of $C$ (see \ref{RangT}), one may assume that $L_2'\circ\iota_m=h_0\oplus \iota_m^{(1)}$, where $h_0, \iota_m^{(1)}: C_m\to A$ are $\mathcal G'$-$\eta'$-multiplicative for a sufficiently large $\mathcal G'$ and sufficiently small $\eta'$ (only depends on $C_m$ and $(\Phi_1)^\ddagger(\overline{(F_a})_0)$) so that $(L_2\circ\iota_m)^\ddagger$ is defined on a subgroup of $U(C_m)/CU(C_m)$ containing $(\Phi_1)^\ddagger(\overline{(F_a)_0})$, $\pi_0((\Phi_1)^\ddagger(\overline{(F_a)_0}))$, $\pi_1(U(C_m)/CU(C_m))$ and $\pi_2(U(C_m)/CU(C_m))$. Moreover, since every finite dimensional
\CA\, is semi-projective and since $L_2'$ is chosen after $C_m$ is chosen, we may assume that  the map $h_0$ is
a \hm\, and has finite dimensional range, and $\tau(h_0(1_{C_m}))<\min\{\delta', \sigma_{a, 1}/3\}$ for any $\tau\in T(A)$, where $\delta'$ is the constant (in place of $\delta$) of Lemma \ref{exp-length2} with respect to $\sigma_{a, 2}$ (in place of $\ep$).

Then, by Lemma \ref{exp-length2}, there is is a homomorphism $\psi_0: C_m\to e_0'Ae_0'$, where $e_0'=\psi_0(1_{C_m})$, such that
%\begin{enumerate}
%\item

(i) $\psi_0$ is homotopically trivial, and $[\psi_0]_0=[h_0]_0$, and

%\item
(ii) for any $u\in \mathcal U_{a, 0}^{(1)}$, one has
\begin{equation}\label{eq-ak-002}
\beta(\Phi_1^\ddagger(\overline{u}))^{-1}(\psi_0\oplus\iota_m^{(1)})^\ddagger(\Phi_1^\ddagger(\overline{u}))=\overline{g_u}
\end{equation}
 for some $g_u\in U_0(A)$ with $\mathrm{cel}(g_u)<\sigma_{a, 2}$.
%\end{enumerate}

Define $L_2=(\psi_0\oplus\iota_m^{1})\circ L_{0,2}: C_m\to A$. Then, for any $u\in\mathcal U_a$, by \eqref{eq-ak-002} and \eqref{eq-ak-lift},  one then has
\begin{equation}\label{alg-k1-a-tf}
(L_2\circ\Phi_1)^\ddagger(u)=\beta(\Phi_1^\ddagger(\overline{u}))\cdot\overline{g_u}=\overline{u}\cdot \overline{g_u}\approx_{\sigma_{a, 2}}\overline{u} \rforal u\in\mathcal U_{a, 0}^{(1)}.
\end{equation}
Moreover, it is also clear that
\begin{equation}\label{uniq-kk-002}
[L_2\circ \Phi_1]|_{\mathcal P^{(1)}_a}=[\id]|_{\mathcal P^{(1)}_a}\andeqn
\end{equation}
%and
\begin{equation}\label{uniq-tr-002}
|\tau\circ L_2\circ \Phi_1(f)-\tau(f)|<\sigma^{(1)}_{a, 1} \rforal f\in\mathcal H^{(1)}_{a} \rforal \tau\in T(A).
\end{equation}
Note that $L_2$ is still $\mathcal G_c^{(2)}$-$\delta_c^{(2)}$-multiplicative. One then has that for arbitrary $u\in \mathcal U_{a, 1}^{(1)}$ (with $\overline{u}$ has order $k$),
\begin{equation*}%\label{almost-cd-004}
 \mathrm{dist}((L_2\circ\Phi_1)(u^k), CU(A))<\sigma_{a, 2}^{(1)}.
 \end{equation*}
Therefore, there is a self-adjoint element $h\in A$ with $||h||<\sigma_{a, 2}^{(1)}$ such that
\begin{equation*}
(u^*)^k(L_2\circ \Phi_1 )(u^k) \exp(2\pi i h)\in CU(A),
\end{equation*}
and hence
\begin{equation*}
((u^*)(L_2 \circ \Phi_1)(u)\exp(2\pi i h/k))^k\in CU(A).
\end{equation*}
Note that
$$(u^*)(L_2\circ \Phi_1 )(u)\exp(2\pi i h/k)\in U_0(A)$$
and $U_0(A)/CU(A)$ is torsion free (Corollary \ref{Unotrosion}). On has that
$$(u^*)(L_2\circ \Phi_1)(u)\exp(2\pi i h/k)\in CU(A).$$
In particular, this implies that
\begin{equation*}%\label{alg-k1-a-tor}
\mathrm{dist}((L_2\circ \Phi_1)^{\ddagger}(u), u)< \sigma_{a, 2}^{(1)}.
\end{equation*}
Together with \eqref{alg-k1-a-tf}, one has that
\begin{equation}\label{alg-k1-a}
\mathrm{dist}((L_2\circ \Phi_1)^{\ddagger}(u), u)<\sigma_{a, 2}^{(1)} \rforal u\in \mathcal U_{a}^{(1)}.
\end{equation}

Then, applying Theorem \ref{MUN1} with \eqref{uniq-kk-002}, \eqref{uniq-tr-002} and \eqref{alg-k1-a}, one obtains a unitary $W\in A$ such that
$$||W^*(L_2 \circ \Phi_1(f))W-f||<\ep_1 \rforal f\in\mathcal F_1.$$
Redefine $L_2$ to be $\mathrm{Ad}(W)\circ L_2$, and one has
$$||L_2 \circ \Phi_1(f)-f||<\ep_1 \rforal f\in\mathcal F_1.$$ That is, one has the following diagram
\begin{displaymath}
\xymatrix{
C \ar[r]^{\id} \ar[d]_{L_1} & C \ar[d]^{L_2}\\
A \ar[ur]_{\Phi_1}\ar[r]_\id & A
}
\end{displaymath}
with the upper triangle approximately commutes on $\mathcal G_1$ up to $\ep_1$ and the lower triangle approximately commutes on $\mathcal F_1$ up to $\ep_1$.
Since $\dt_c^{(2)},$ ${\cal G}_c^{(2)},$
$\sigma_{c,1}^{(2)},$ $\sigma_{c,2}^{(2)},$  ${\cal P}_c^{(2)}$ and ${\cal H}_c^{(2)}$ have been chosen and
embedded into the construction of $L_2,$ the construction can continue.

By repeating  this argument, one obtains the following approximate intertwining diagram
\begin{displaymath}
\xymatrix{
C \ar[r]^{\id} \ar[d]_{L_1} & C \ar[d]^{L_2}  \ar[r]^{\id} & C \ar[d]^{L_3}  \ar[r]^{\id}  & C \ar[d]^{L_4} \ar[r] & \cdots\\
A \ar[ur]_{\Phi_1}\ar[r]_\id & A \ar[ur]_{\Phi_2} \ar[r]_\id & A \ar[ur]_{\Phi_3} \ar[r]_\id & A \ar[r] \ar[ur]& \cdots ,
}
\end{displaymath}
where
$$||\Phi_{n} \circ L_n(g)-g||<\ep_{n}\rforal g\in\mathcal G_n\andeqn$$
%and
$$||L_{n+1} \circ \Phi_n(f)-f||<\ep_{n}\rforal f\in\mathcal F_n.$$
By the choices of ${\cal G}_n,$ ${\cal F}_n$ and the fact that $\sum_{n=1}^{\infty}\ep_n<\infty,$
the standard Elliott approximate intertwining argument applies, which shows that $A\cong B$, as desired.  This proves the first part of the proof.

%To see that
{{The}} proof of the second part is basically the same but simpler since we only need to
have a one-sided  approximate intertwining.  In particular, we do not need to construct $\Phi_1.$
Thus, once $L_1$ is constructed, we can go to construct $L_2.$ We can first construct $L_2'$ as above right
after (\ref{eq-ak-lift}).  It is important that we do not need to assume that $L_1^{\ddag}$ is injective on
$\overline{F_c}_0,$ since we will apply \ref{exp-length2} but not \ref{exp-length1}.
\end{proof}

\begin{thm}\label{IST1}
Let $A_1, B_1\in {\cal B}_0$ be two unital separable simple \CA s which satisfying the UCT.
Let $A=A_1\otimes U$ and $B=B_1\otimes U,$ where $U$ is a UHF-algebra of infinite type.
Suppose that ${\rm Ell}(A)={\rm Ell}(B).$ Then there exists an isomorphism $\phi: A\to B$
which carries the identification of ${\rm Ell}(A)={\rm Ell}(B).$
\end{thm}

\begin{proof}
By Theorem \ref{RangT}, there is  a C*-algebra $C$ {{ as in \ref{RangT}}}, such that ${\rm Ell}(C)\cong {\rm Ell}(A)\cong {\rm Ell}(B)$. By the first part of Proposition \ref{IST0}, one has that $C\cong A$ and $C\cong B$. In particular, $A\cong B$.
\end{proof}

\begin{cor}\label{Cist}
Let $A$ and $B$ be as in \ref{IST1}. Suppose that there is a \hm\, $\Gamma$
from ${\rm Ell}(A)$ to ${\rm Ell}(B)$ such that $\Gamma([1_A])=[1_B].$ Then there exists a
unital \hm\, $\phi: A\to B$ such that $\phi$ induces $\Gamma.$
\end{cor}

\begin{proof}
By Theorem \ref{IST1}, one may assume that $A$ is a unital C*-algebra described in \ref{RangT}. Then the Corollary follows from the second part of Theorem \ref{IST0}.
\end{proof}

\section{More existence theorems}
%{\color{Green} In the following lemma, where do we need that $U$ is of infinite type?}{{---I think, in the original proof, it not only uses that $U$ is of infinite type, it actually  uses that $U$ is self absorb $U\otimes U= U$. But in general, $U$ is of infinite type does not imply self absorb. In this section, there are many places involving using this stronger assumption. I hope that I correct all of them---Gong}}
\begin{lem}\label{preBot1}
Let $X$ be a finite CW complex, $C=PM_k(C(X))P$ and let $A_1\in {\cal B}_0$ be a unital simple \CA. Assume that $A=A_1\otimes U$ for {{an infinite dimensional}}  UHF-algebra $U.$
%of infinite type.
Let $\af\in KK_e(C, A)^{++}.$
Then there exists a unital  monomorphism $\phi: C\to A$ such that $[\phi]=\af.$ Moreover
we may write $\phi=\phi_n'\oplus \phi_n'',$
where $\phi_n': C\to (1-p_n)A(1-p_n)$ is a unital monomorphism,  $\phi_n'': C\to p_nAp_n$ is a unital homomorphism with $[\phi_n'']=[\Phi]$ in $KK(C, p_nAp_n)$ for some point evaluation map $\Phi$ and
$$\lim_{n\to\infty} \max\{\tau(1-p_n):\tau\in T(A)\}=0\rforal \tau\in T(A),$$
where $p_n\in A$ is a sequence of projections.
\end{lem}

\begin{proof}
To simplify the matter, we may assume that $X$ is connected. It is also easy to check that the general case
can be reduced to the case that $C=C(X).$ {{It is routine to prove that  $\alpha\in KK(C,A)^{++}$ implies that $\alpha ({\rm ker}\rho_C) \subset {\rm ker}\rho_A$.}}

Since $K_i(C)$ is finitely generated, $i=0,1,$ $KK(C,A)=KL(C,A).$
Let $\af\in KL_e(C,A)^{++}$ which we will identify  with an element in ${\rm Hom}_{\Lambda}(\underline{K}(C),
\underline{K}(A))$ by a result 
%of Dadarlat
%%%Dadalart
% and Loring 
in (\cite{DL}).

{ {Write  $A=\lim_{n\to\infty}(A_1\otimes M_{r_n}, \imath_{n,n+1}),$  where $r_n|r_{n+1},$
$r_{n+1}=m_nr_n$ and
$\imath_{n,n+1}(a)=a\otimes 1_{M_{m_n}},$ $n=1,2,....$ Since $K_*(C)$ is finitely generated and consequently, $\underline{K}(C)$ is finitely generated modulo Bockstein $\Lambda$ operations, there is an element  $\alpha_1\in KK(C, A_1\otimes M_{r_n})$ such that $\alpha=\alpha_1\times [  \imath_n],$ where $ [  \imath_n] \in KK(A_1\otimes M_{r_n}, A)$ is induced by the inclusion $ \imath_n: A_1\otimes M_{r_n} \to A$.  By increasing $n$, we can assume that $\alpha_1 ({\rm ker}\rho_C) \subset {\rm ker}\rho_{A_1\otimes M_{r_n}}$ and further that $\alpha_1\in KK_e(C, A_1\otimes M_{r_n})^{++}$. Replacing  $A_1$ by $A_1\otimes M_{r_n}$, without lose of generality we can assume that $\alpha=\alpha_1 \times [\imath]$, where $\alpha_1\in KK_e(C, A_1)^{++}$ and $ \imath: A_1 \to A$ is the inclusion.}}

It induces an element ${{\tilde \af}}_1\in KL(C\otimes U, A\otimes U).$
Let $K_0(U)=\D,$ a dense subgroup of $\Q$.
% with the property that $\D\cdot \D=\D.$
%The above using stronger assumption that $U$ is self absorb $U\otimes U =U$---Gong
Note that $K_i(C\otimes U)=K_i(C)\otimes \D,$ $i=0,1,$ by the Kunneth formula.

We verify that ${{\tilde \af}}_1(K_0(C\otimes U)_+\setminus \{0\})\subset K_0(A\otimes U)_+\setminus\{0\}.$
Consider $x=\sum_{i=1}^m x_i\otimes d_i\in K_0(C\otimes U)_+\setminus \{0\}$ with $x_i\in K_0(C)$ and $d_i\in \D,$ $i=1,2,...,m,$   There  is a projection $p\in M_r(C)$ for some $r\ge 1$ such that $[p]=x.$
 Let $t\in T(C),$ then
\beq\label{prebot1-n1}
\sum_{i=1}^m t(x_i)d_i>0.
\eneq
It should be noted that, since $C=C(X)$ and $X$ is connected, $t(x_i)\in \Z$ and $t(x_i)=t'(x_i)$ for all $t, t'\in T({{C}}).$
Since $\af_{ {1}}([1_C])=[1_{A_{{1}}}],$ $\tau\circ \af_{{1}}(x_i)=t(x_i)$ for any $\tau\in T(A_{{1}})$ and $t\in T(C).$
By (\ref{prebot1-n1}),
\beq\label{prebot1-n2}
\tau({{\tilde \af}}_1(x))=\sum_{i=1}^m \tau\circ\af_{{1}}(x_i)d_i=\sum_{i=1}^m t(x_i)d_i>0
\eneq
for all $\tau\in T(A_{{1}}).$  This shows that ${{\tilde \af}}_1$ is strictly positive.
For any \CA\, $A',$ in this proof, we will use $j_{A'}: A'\to A'\otimes U$ for the \hm\, $j_{A'}(a)=a\otimes 1_U$
for all $a\in A'.$
%Let $j_C: C\to C\otimes U$ be defined by $j(c)=c\otimes 1_U$ and let $j_A: A\to A\otimes U.$
%There exists a unital \hm\, $s: A\otimes U\to A$ such that
%$s\circ j_A$ is approximately inner.
Evidently,
\beq\label{prebot1-n3}
\alpha= {\tilde \af}_1\circ j_C=j_{A_1}\circ \af_1.
\eneq

%We obtain
%\beq\label{prebot1-n3}
%[s]\circ \af_1\circ j_C=\af.
%\eneq

Write  $U=\lim_{n\to\infty}(M_{r_n}, \imath_n),$  where $r_n|r_{n+1},$
$r_{n+1}=m_nr_n$ and
$\imath_n(a)=a\otimes 1_{M_{m_n}},$ $n=1,2,....$   We may assume that $r_1=1.$  { {The UHF algebra $U$ corresponds to the {super natural number} $\Pi_{i=1}^{\infty} m_i$. Evidently, we can choose $m_i$ carefully so that we can write the {super natural} number $\Pi_{i=1}^{\infty} m_i$ in another way $\Pi_{i=1}^{\infty} m_i=\Pi_{i=1}^{\infty} l_i$
with $l_i|m_i$ and $\lim_{i\to \infty} \frac{m_i}{l_i}=\infty$.}}

 Let $\{x_n\}$ be a sequence of points in $X$ such that $\{x_k,x_{k+1},...,x_n,...\}$ is dense in $X$ for each $k$
and  each point in $\{x_n\}$ repeats infinitely many times.
Let ${ {B}}=\lim_{n\to\infty} ({{B_n=}}M_{r_n}(C), \psi_n),$ where
$$
\psi_n(f)={\rm diag}({ {\underbrace{f,f,\cdots f}_{l_n}}}, f(x_1),f(x_2),...,f(x_{m_n-{{l_n}}}))\tforal f\in M_{r_n}(C),
$$
$n=1,2,....$  Note that $\psi_n$ is injective.
Denote $e_n={\rm diag}(1_{M_{r_n{{\cdot l_n}}}},0,...,0)\in M_{r_{n+1}}(C),$ $n=1,2,...$

It is standard that ${ {B}}$ has tracial rank zero and { {$K_*(B)=K_*(C\otimes U)$.}}
%$K_0(B_1)=\D\oplus {\rm ker}\rho_C$
%and $K_1(B_1)=K_1(C).$
%Put $B=B_1\otimes U.$
Note that
 $B$ is a unital  simple AH-algebra with no dimension growth, with real rank zero and with a unique tracial state.
 %Define ${\red{h}}: C\to {\red {B}}$ by ${\red=\psi_{1, \infty}.$
 %Define $h: C\otimes U\to B$ by
% $h(c\otimes a)=h'(c)\otimes a$ for all $c\in C$ and $a\in U.$  Since $\psi_n$ is injective, $h'$ is a monomorphism.
 %We have
 %\beq\label{prebot1-n4}
 %h\circ j_C=j_{B_1}\circ h'.
 %\eneq
 Note that
$$
(K_0(B), K_0(B)_+, [1_B], K_1(B))=(K_0(C\otimes U), K_0(C\otimes U)_+, [1_{C\otimes U}], K_1(C\otimes U)).
$$
Thus we obtain a $KK$-equivalence ${{\kappa}}\in KL_e(B, C\otimes U)^{++}$ (by the UCT).  { {It is standard to construct a unital homomorphism $h: C\to B$ such that $KK(h)={\kappa}^{-1}\circ j_C$ (see \cite{Li-K-theory}). In particular if $N$ is lager enough, we can choose homomorphism $h': C\to B_N$ such that $h=\imath_N\circ h'$, where $\imath_N: B_N\to B$ is the inclusion. }}

% such that
%\beq\label{prebot1-n5}
%\kappa_0\circ [h]\circ [j_C]=[j_C].
%\eneq
We also have that ${{{\tilde\af}_1\circ \kappa}}\in KL_e(B, {{A}})^{++}, $ {{where recall $A=A_1\otimes U.$}}
We also note that $B$ has a unique tracial state. Let $\gamma: T(A)\to T(B)$ by
$\gamma(\tau)=t_0$ where $t_0\in T(B)$ is the unique tracial state. It follows that
${{{\tilde\af}_1\circ \kappa}}$ and $\gamma$ is compatible.  By the second part of Theorem \ref{IST0},
there is a unital \hm\, $H: B\to {{A}}$ such that
$[H]={{{\tilde\af}_1\circ \kappa.}}$
Define $\phi: C\to A$ by ${{ \phi= H\circ h=H\circ \imath_N\circ h'.}}$ Then,  $\phi$ is injective and by (\ref{prebot1-n3}) and {{${{[h]}}={\kappa}^{-1}\circ j_C$ , we have}}
$[\phi]=\af.$

To show the last part,  define $q_n=\phi_{n+1, \infty}{{(e_n)}}\in B,$  $n={{ N+1, N+2,\cdots.}}$
Define $p_n={ {1- H(q_n)}},$ $n={{ N+1, N+2,\cdots.}}$ One checks that
\beq\label{prebot1-n6}
\lim_{n\to\infty}\max\{\tau(1-p_n):\tau\in T(A)\}={{\lim \frac{l_n}{m_n}=}}0
\eneq
{{Note that for $n>N$, $q_n$ commutes with image of $h$ and the homomorphism $(1-q_n)h(1-q_n): C\to (1-q_n)B(1-q_n)$ is defined by point evaluation. }} Define $\phi_n': C\to {{ (1-p_n)A(1-p_n)}}$ by $\phi_n'(f)= H(q_n)H\circ {{ h(f)H(q_n)}}.$
%Note that ${\tilde \phi}_n(f)=(1-e_n)\psi_{1,n+1}(f)$ is a point-evaluation map.
 Define
{ {$\phi_n''(f)= p_nH\circ  h(f)p_n,$  which is }}a point-evaluation map. The lemma follows.

\end{proof}

%{\color{Green} I think it might be better to give a reference or the outline for the following lemma: It is clear with real rank zero, but just with the (SP) property, it is not clear why the non-trivial K-theory is in a small corner.} {\red {I tried to add a proof for the following, but seems it is not so easy, the trouble is I can not find a single homomorphism which can be decomposed as in the lemma for sequence $n$. The best I  can do is that for each estimation of trace, we can find a homomorphism as in lemma to satisfies that the non point evaluation part is small according to the trace estimation. But we can not find a same one. I suggested to remove it if not needed in the future---Gong.}}

We also have the following:
\begin{lem}\label{Circlef-1}
Let $C=M_k(C(\T))$ and $A$ be a unital  infinite dimensional simple \CA\, with stable rank one and with the (SP) property.  Then the conclusion of \ref{preBot1}
also holds.
\end{lem}

\begin{proof}
The existence of $\af$ implies that $A$ contains mutually equivalent and mutually orthogonal projections
$e_1, e_2,...,e_k$ such that $\sum_{i=1}^k e_i=1_A.$ Since $e_iAe_i$ are unital infinite dimensional simple \CA\, with stable rank one and with
(SP), the general case can be reduced to the case that $k=1.$
Fix $1>\dt>0.$ Choose a non-zero projection $p\in A$ such that $\tau(p)<\dt$ for all $\tau\in T(A).$
Note that $K_1(pAp)=K_1(A),$ since $A$ is simple.  Let $\af_1: K_1(C(\T))\to K_1(pAp)$ be the \hm\, given by $\af.$
Let $z\in C(\T)$ be the standard unitary generator. Let $x=\af_1([z])\in K_1(pAp).$ Since
$pAp$ has stable rank one, there is a unitary $u\in pAp$ such that $[u]=x$ in $K_1(pAp)=K_1(A).$
Define $\phi': C(\T)\to pAp$ by $\phi'(f)=f(u)$ for all $f\in C(\T).$ Define $\phi'': C(\T)\to (1-p)A(1-p)$ by
$\phi''(f)=f(1)(1-p)$ for all $f\in C(\T)$ (where $f(1)$ is a point-evaluation at $1$ on the unit circle).
Define $\phi=\phi'\oplus \phi': C(\T)\to A.$  The lemma follows.

\end{proof}

%{\red{Since the new statement of Lemma \ref{ExtTraceH} does not say the homomorphism is  defined by point evaluation, the proof is not working. I make significant change in the proof, please read it---Gong}} {\red{I changed the proof almost back to original proof with a modified \ref{ExtTraceH} which mentioned point evaluation--Gong}}

%{\color{Green} It seems that $U$ only has to be a UHF algebra of infinite dimensional (it always can be written as tensor product of other two UHF algebras of infinite dimensional). It seems that in many other places, the UHF algebra only has to be infinite dimensional. There are also a few places where $U$ is indeed need to be of infinite type; trying to recall where it is.} {\red{Since in section 21, we use U to be infinite type, and we use the result from section 21, so we can not simply change infinite type to infinite diemensional---Gong----Fine L}}

\begin{cor}\label{istBotC}
Let  $X$ be a  connected  finite CW complex, $C=PM_m(C(X))P,$
where $P\in M_m(C(X))$ is a projection, let $A_1\in {\cal B}_0$ be a unital separable simple \CA\, which satisfies
the UCT and let $A=A_1\otimes U,$ where $U$ is a UHF-algebra of infinite type.
Suppose that $\af\in KK(C, A)^{++}$ and $\gamma: T(A)\to T_f(C(X))$ is a continuous affine map.
Then
there exists a
sequence of \morp  s $h_n: C\to A$ such that
\begin{enumerate}
\item $\lim_{n\to\infty}\|h_n(ab)-h_n(a)h_n(b)\|=0,$ for any $a,b\in C$,
\item for each $h_n$, the map $[h_n]$ is well defined and $[h_n]=\alpha$, and
\item $\lim_{n\to\infty}\max\{|\tau\circ h_n(f)-\gamma(\tau)(f)|: \tau\in T(A)\}=0$  for any $f\in C$.
\end{enumerate}
\end{cor}

\begin{proof}
By Theorem \ref{IST1}, one may assume that $A$ is a unital C*-algebra described in \ref{RangT}.  It follows from Lemma \ref{preBot1} that there is a unital \hm\,  $h_n: C\to A$ such that
%\begin{equation*}
%\lim_{n\to\infty}\|h_n(ab)-h_n(a)h_n(b)\|=0\tforal a,b\in C
%\end{equation*}and
$[h_n]=\alpha$. Moreover, $$h_n=h_n'\oplus h_n'',$$ where $h''_n: C\to p_nAp_n$ is a homomorphism with $[h''_n]=[\Phi^{{'}}]$ in $KK(C, p_nAp_n)$ for some point evaluation map $\Phi^{{'}}$, where $p_n$ is a projection in $A$ with $\tau(1-p_n)$ converge to $0$ uniformly as $n\to\infty$. We will modify the map $ h_n=h_n'\oplus h_n''$ to get our \hm.

Assert that for any finite subset $\mathcal H\subseteq C_{s.a}$, and $\epsilon>0$, and any sufficiently large $n$, there is a unital homomorphism $\tilde{h}_n: C\to p_nAp_n$ such that  $[\tilde{h}_n]=[\Phi]$ in $KK(C,p_n Ap_n)$ for some point-evaluation $\Phi$, and
$$|\tau\circ \tilde{h}_n(f)-\gamma(\tau)(f)|<\epsilon\rforal \tau\in T(A)$$ for all $f\in\mathcal H$. The corollary then  follows by replacing the map $h_n''$ by the map $\tilde{h}_n${{----of course, we use $\lim_{n\to \infty}\tau(1-p_n)=0$.}}

Let $\mathcal H_{1, 1}$ (in place of $\mathcal H_{1, 1}$) be the finite subset of Lemma \ref{ExtTraceH} with
 respect to $\mathcal H$ (in place of $\mathcal H$), $\epsilon/{8}$ (in the place of $\sigma$), and $C$
 (in the place of $C$). Since $\gamma(T(A))\subseteq T_f(C(X))$, there is $\sigma_{1, 1}>0$ such that $$\gamma(\tau)(h)>\sigma_{1, 1}\rforal h\in\mathcal H_{1, 1}\rforal \tau\in T(A).$$

Let $\mathcal H_{1, 2}\subseteq C^+$ (in the place of $\mathcal H_{1, 2}$) be the finite subset of Lemma \ref{ExtTraceH} with respect to $\sigma_{1, 1}$. Since $\gamma(T(A))\subseteq T_f(C(X))$, there is $\sigma_{1, 2}>0$ such that $$\gamma(\tau)(h)>\sigma_{1, 2}\rforal h\in\mathcal H_{1, 2} \rforal \tau\in T(A).$$

Let $M$ be the constant of Lemma \ref{ExtTraceH} with respect to $\sigma_{1, 2}$. Using a same argument as that of Lemma \ref{cut-trace}, for sufficiently large $n$, there is a C*-subalgebra $D\subseteq p_nAp_n\subseteq A$ such that $D\in \mathcal C_0$, a continuous affine map $\gamma': T(D)\to T(C)$ such that
$$|\gamma'(\frac{1}{\tau(p)}\tau|_D)(f)-\gamma(\tau)(f)|<\epsilon/4 \rforal \tau\in T(A)\rforal f\in\mathcal H,$$
where $p=1_D$, $\tau(1-p)<\epsilon/(4+\epsilon)$, $\rforal \tau\in T(A)$,
\beq
&&\gamma'(\tau)(h)>\sigma_{1, 1}\rforal \tau\in T(D) \rforal h\in\mathcal H_{1, 1}\andeqn\\
&&\gamma'(\tau)(h)>\sigma_{1, 2}\rforal \tau\in T(D)\rforal h\in\mathcal H_{1, 2}.
\eneq
Since $A$ is simple and not of elementary, one may assume that the dimensions of the irreducible representations of $D$ are at least $M$. Thus, by Lemma \ref{ExtTraceH}, there is a homomorphism $\phi: C\to D$ such that
$[\phi]=[\Phi]$ in $KK(C, D)$ for a point evaluation map $\Phi$, and that
%{\red{(in the statement of \ref{ExtTraceH}, we did not mention that $[\phi]=[\Phi]$for a point evaluation map $\Phi$, may be we should make such a remark---Gong---I found out it is in old version but was removed for the new version }},and
 $$|\tau\circ\phi(f)-\gamma'(\tau)(f)|<\epsilon/{{4}}\rforal f\in \mathcal H\rforal \tau\in T(D).$$

Pike a point $x\in X$, and define ${\tilde h}: C\to p_nAp_n$ by
$$f\mapsto f(x)(p_n-p)\oplus\phi(f)\rforal f\in C.$$ Then a calculation as in the proof of Theorem \ref{istTr} shows that the homomorphism $h_n'\oplus {\tilde h}$ satisfies the assertion.

\end{proof}

%{\color{Green} This Lemma is never used.}
%\begin{lem}\label{Circlef-2}
%{Let $C=M_n(C(\T)),$  let $A_1\in {\cal B}_1$ and let $A=A_1\otimes U$ for some UHF-algebra of infinite type.
%Let $\af\in KL_e(C, A)^{++},$ $\gamma: T(A)\to T_f(C)$ and $\lambda: U(C)/CU(C)\to U(A)/CU(A)$ such that
%$(\af, \gamma, \lambda)$ is compatible. Then, for any $\ep>0$ and any finite subset ${\cal H}\subset C_{s.a.},$ there
%exists a unital \hm\, $h: C\to A$ such that
%$[h]=\af,$ $h^{\ddag}=\lambda$ and
%\beq\label{Circlef-3}
%|\tau\circ h(f)-\gamma(\tau)(f)|<\ep\tforal f\in {\cal H}.
%\eneq
%}
%\end{lem}

%\begin{proof}

%\end{proof}

\begin{cor}\label{istBotC+}
Let   $C\in {\bf H}$ {{(where ${\bf H}$ is defined in \ref{AHblock})}}
and let $A_1\in {\cal B}_0$ be a unital separable simple \CA\, which satisfies
the UCT and let $A=A_1\otimes U$ { for some UHF-algebra $U$ of infinite type.}
Suppose that $\af\in KK_e(C, A)^{++},$ $\lambda: U(C)/CU(C)\to U(A)/CU(A)$ is
a continuous \hm\, and $\gamma: T(A)\to T_f({{C}})$ is a continuous affine map such that
$\af, \lambda, $ and $\gamma$ are compatible.
Then there exists a
sequence of  unital \morp  s $h_n: C\to A$ such that
\begin{enumerate}
\item $\lim_{n\to\infty}\|h_n(ab)-h_n(a)h_n(b)\|=0$ for any $a,b\in C$,
\item for each $h_n$, the map $[h_n]$ is well defined and $[h_n]=\alpha$,
\item $\lim_{n\to\infty}\max\{|\tau\circ h_n(f)-\gamma(\tau)(f) |: \tau\in T(A)\}=0$ for all $f\in C,$  and,
\item $\lim_{n\to\infty}{\rm dist}(h_n^{\ddag}({\bar u}), \lambda({\bar u}))=0$
for any $u\in U(C).$
\end{enumerate}
\end{cor}

\begin{proof}
Let $\epsilon>0$. Let $\mathcal U$ be a finite generating set of $J_c(K_1(C))$, {{where $J_c(K_1(C))$ is as in \ref{Dcu}.}} Let $\delta>0$ and $\mathcal G$ be the constant and finite subset of Lemma \ref{exp-length2} with respect to $\mathcal U$, $\epsilon$ and $\lambda$ (in the place of $\alpha$). Without loss of generality, one may assume that $\delta<\epsilon$.

Let $\mathcal F$ be a finite subset such that $\mathcal F\supset \mathcal G$. Let $\mathcal H\subseteq C$ be a finite subset of self-adjoint elements with norm at most one. By Corollary \ref{istBotC}, there is a positive completely linear map $h': C\to A$ such that $h$ is $\mathcal F$-$\delta$-multiplicative, $[h']$ is well-defined and $[h']=\alpha$, and
\begin{equation}\label{tr-eqn-09-27}
|\tau(h'(f))-\gamma(\tau)(f)|<\epsilon,\quad\tau\in\mathrm{T}(A), \ f\in\mathcal H.
\end{equation}

By Theorem \ref{IST0}, the C*-algebra $A$ is isomorphic to one of the model algebras constructed in Theorem \ref{RangT}, and therefore there is an inductive limit decomposition $A=\varinjlim(A_i, \phi_i)$, where $A_i$ and $\phi_i$ are described as in Theorem \ref{RangT}. Without loss of generality, one may assume that $h'(C)\subseteq A_i$. Therefore, by Theorem \ref{RangT}, the map $\phi_{1, \infty}\circ h$ has a decomposition $$\phi_{1, \infty}\circ h'=\psi_0\oplus\psi_1$$ such that $\psi_0, \psi_1$ satisfy the (1)-(4) Lemma \ref{exp-length2} with the above $\delta$.

It then follows from Lemma \ref{exp-length2} that there is a homomorphism $\Phi: C\to e_0Ae_0$, where $e_0=\psi_0(1_C)$, such that
%\begin{enumerate}
%\item

(i) $\Phi$ is homotopic to a homomorphism with finite dimensional range and
\begin{equation}\label{same-k0}
[\Phi]_{*0}=[\psi_0]\andeqn
\end{equation}
%\item

(ii) for each $w\in \mathcal U$, there is $g_w\in \text{U}_0(B)$ with $\mbox{cel}(g_w)<\epsilon$ such that
        \begin{equation}\label{k1-eqn-09-27-01}
         \lambda(\bar{w})^{-1}(\Phi\oplus \psi_1)^\ddagger(\bar{w})=\bar{g_w}
         \end{equation}
%\end{enumerate}

Consider the map $h:=\Phi\oplus\psi_1.$ Then $h$ is $\mathcal F$-$\epsilon$-multiplicative. By \eqref{same-k0}, one has
$$[h]=[\psi_0]\oplus[\psi_1]=[h']=\alpha.$$
By \eqref{tr-eqn-09-27} and Condition (4) of Lemma \ref{exp-length2}, one has that, for all $f\in {\cal H},$
$$|\tau(h(f))-\gamma(\tau)(f)|\leq |\tau(h'(f))-\gamma(\tau)(f)|+\delta <\epsilon+\delta<2\epsilon.
%,\quad f\in\mathcal H.
$$
It follows from \eqref{k1-eqn-09-27-01} that, for all $u\in {\cal U},$ 
$$\mathrm{dist}(\overline{h(u)}, \lambda(\overline{u}))<\epsilon.$$
Since $\mathcal F$, $\mathcal H$, and $\epsilon$ are arbitrary, this shows the Corollary.
\end{proof}

\begin{cor}\label{istBotC++}
Let $C\in {\bf H}$
and let $A_1\in {\cal B}_0$ be a unital separable simple \CA\, which satisfies
the UCT and let $A=A_1\otimes U$ { for some UHF-algebra $U$ of infinite type.}
Suppose that $\af\in KL_e(C, A)^{++},$ $\lambda: U(C)/CU(C)\to U(A)/CU(A)$ is
a continuous \hm\, and $\gamma: T(A)\to T_f({{C}})$ is a continuous affine map such that
$\af, \lambda, $ and $\gamma$ are compatible.
Then
there exists a  unital \hm\, $h: C\to A$
 such that
\begin{enumerate}
 \item $[h]=\af, $
\item $\tau\circ h(f)=\gamma(\tau)(f) $ for any $f\in C,$ and,
\item $h_n^{\ddag}=\lambda.$
\end{enumerate}
\end{cor}

\begin{proof}
Let us construct a sequence of homomorphisms $h_n: C\to A$ which satisfies (1)--(4) of Corollary \ref{istBotC+}, and moreover, the sequence $\{h_n(f)\}$ is Cauchy for any $f\in\mathcal C$. Then the limit map $h=\lim_{n\to\infty} h_n$ is the desired homomorphism.

To construct such a sequence of homomorphisms, it is enough to construct a sequence of homomorphisms satisfying (1)-(4) of Corollary \ref{istBotC+} such that $\{h_n(f)\}$ is Cauchy for any
$f\in C.$

Let $\{{\cal F}_n\}$ be an increasing sequence of the unit ball of $C$ such that its union is dense in  the unit ball of $C.$
Define $\Delta(a)=\min\{\gamma(\tau)(a): \tau\in\mathrm{T}(A)\}$. Since $\gamma$ is continuous and $\mathrm{T}(A)$ is compact, the map $\Delta$ is an order preserving map from $C_+^{1, q}\setminus\{0\}$ to $(0, 1)$.
Let $\mathcal G(n), \mathcal H_1(n), \mathcal H_2(n) \subseteq C$,
$\mathcal U(n)\subseteq U_\infty(C)$, $\mathcal P(n)\subseteq\underline{K}(C)$, $\gamma_1(n)$, $\gamma_2(n)$, $\delta(n)$ be the finite subset and constants of Theorem \ref{UniCtoA} with respect to ${\cal F}_n$, $1/2^{n+1}$ and $\Delta/2$. Without loss of generality, one may assume that $\delta(n)$ decrease to $0$ if $n\to\infty$  and $\mathcal P(n)\subset \mathcal P(n+1),$ $n=1,2,...,$ and
$\cup_{n=1}^{\infty} {\cal P}(n)=
\underline{K}(C).$

Let $\mathcal G_1\subseteq\mathcal G_2\subseteq\cdots$ be an increasing sequence of finite subset of $C$ such that $\bigcup\mathcal  G_n$ is dense in $C$, and let $\mathcal U_1\subseteq\mathcal U_2\subseteq\cdots$ be an increasing sequence of finite subset of $U(C)$ such that $\bigcup\mathcal  U_n$ is dense in $U(C)$. One may assume that $\mathcal G_{n}\supseteq \mathcal G(n)\cup \mathcal G(n-1)$, $\mathcal G_{n}\supseteq {\cal H}_1(n)\cup {\cal H}_1(n+1)\cup \mathcal H_2(n)\cup \mathcal H_2(n-1) $, and $\mathcal U_n\supseteq \mathcal U(n)\cup \mathcal U(n-1)$.

By Corollary  \ref{istBotC+}, there is a $\mathcal G_1$-$\delta(1)$-multiplicative map $h'_1: C\to A$ such that
\begin{enumerate}\setcounter{enumi}{3}
\item the map $[h'_1]$ is well defined and $[h_1]=\alpha$,
\item $|\tau\circ h_n(f)-\gamma(\tau)(f)|<\min\{\gamma_1(1), \frac{1}{2}\Delta(f): f\in \mathcal H_{1}\}$ for any $f\in \mathcal G_1$  and,
\item  ${\rm dist}(h_n^{\ddag}({\bar u}), \lambda({\bar u}))<\gamma_2(1)$
for any $u\in \mathcal U_n.$
\end{enumerate}

Define $h_1=h_1'$. Assume that $h_1, h_2, ..., h_n: C\to A$ are constructed such that
\begin{enumerate}\setcounter{enumi}{6}
\item $h_i$ is $\mathcal G_i$-$\delta(i)$-multiplicative, $i=1, ..., n$,
\item the map $[h_i]$ is well defined and $[h_i]=\alpha$, $i=1, ..., n$,
\item  $|\tau\circ h_i(f)-\gamma(\tau)(f)|<\min\{\frac{1}{2}\gamma_1(i), \frac{1}{2}\Delta(f): f\in \mathcal H_1({i})\}$ for any $f\in \mathcal G_i$, $i=1, ..., n$,
\item ${\rm dist}(h_i^{\ddag}({\bar u}), \lambda({\bar u}))<\frac{1}{2}\gamma_2(i)$ for any $u\in \mathcal U_i,$ $i=1, ..., n$, and
\item $\|h_{i-1}(g)-h_{i}(g)\|<\frac{1}{2^{i-1}}$ for all $g\in {\cal G}_{i-1},$ $i=2, 3, ..., n$.
\end{enumerate}

Let us construct $h_{n+1}: C\to A$ such that
\begin{enumerate}\setcounter{enumi}{11}
\item $h_{n+1}$ is $\mathcal G_{n+1}$-$\delta(n+1)$-multiplicative,
\item the map $[h_{n+1}]$ is well defined and $[h_{n+1}]=\alpha$,
\item $|\tau\circ h_{n+1}(f)-\gamma(\tau)(f)|<\min\{\frac{1}{2}\gamma_1(n+1), \frac{1}{2}\Delta(f): f\in \mathcal H_1({n+1})\}$ for any $f\in \mathcal G_{n+1}$,
\item  ${\rm dist}(h_{n+1}^{\ddag}({\bar u}), \lambda({\bar u}))<\frac{1}{2}\gamma_2(n+1)$ for any $u\in \mathcal U_,$ $i=1, ..., n$, and
\item $\|h_{n}(g)-h_{n+1}(g)\|<\frac{1}{2^{n}}$ for all $g\in {\cal F}_n$.
\end{enumerate}

Then the statement follows.

By Corollary \ref{istBotC+}, there is $\mathcal G(n+1)$-$\delta(n+1)$-multiplicative map $h'_{n+1}: C\to A$ such that
%\begin{enumerate}
%\item

 $h'_{n+1}$ is $\mathcal G_{n+1}$-$\delta(n+1)$-multiplicative,
%\item
 the map $[h'_{n+1}]$ is well defined and $[h_{n+1}']=\alpha$,
%\item

 \begin{equation}\label{eqn-n-09-27-02}
        |\tau\circ h'_{n+1}(f)-\gamma(\tau)(f)|<\min\{\frac{1}{2}\gamma_1(n+1), \frac{1}{2}\Delta(f): f\in \mathcal H_2({n+1})\}
        \end{equation}
        for any $f\in \mathcal G_{n+1}$
and
$$
{\rm dist}((h'_{n+1})^{\ddag}({\bar u}), \lambda({\bar u}))<\frac{1}{2}\gamma_2(n+1)
$$ for any $u\in \mathcal U_,$ $i=1, ..., n$.
%\end{enumerate}
In particular, this implies that
$$[h'_{n+1}]|_{\mathcal P}=[h_n]|_{\mathcal P},$$
and for any $f\in \mathcal H_2(n)$ (note that $\mathcal H_2(n)\subseteq \mathcal G_n$)
\begin{eqnarray*}
|\tau\circ h_n(f)-\tau\circ h'_{n+1}(f)| & < & \gamma_1(n)+|\gamma(\tau)(f)-\tau\circ h'_{n+1}(f)| \\
&< & \gamma_1(n)/2 + \gamma_1(n+1)/2
<\gamma_1(n).
\end{eqnarray*}

Also by \eqref{eqn-n-09-27-02}, for any $f\in\mathcal H_1(n)$, one has that
$$\tau(h'_{n+1}(f))\geq \gamma(\tau)(f)-\frac{1}{2}\Delta(f)> \frac{1}{2}\Delta(f).$$
By the inductive hypothesis, a same argument shows that
$$\tau(h_{n}(f))\geq \gamma(\tau)(f)-\frac{1}{2}\Delta(f)> \frac{1}{2}\Delta(f)\rforal f\in\mathcal H_1(n).$$

For any $u\in\mathcal U(n)$, one has
\begin{eqnarray*}
\mathrm{dist}(\overline{h'_{n+1}(u)}, \overline{h_{n}(u)}) & < & \frac{1}{2}\gamma_2(n+1)+\mathrm{dist}(\gamma(\overline{u}),\overline{h_{n}(u)})\\
&< & \frac{1}{2}\gamma_2(n+1) + \frac{1}{2}\gamma_2(n) <\gamma_1(n).
\end{eqnarray*}

Note that both $h'_{n+1}$ and $h_n$ are $\mathcal G(n)$-$\delta(n)$-multiplicative,  by Theorem \ref{UniCtoA}, there is a unitary $W\in A$ such that
$$\|W^*h_{n+1}'(g)W-h_n(g)\|<1/2^n \rforal g\in {\cal F}_n.$$
Then the map $h_{n+1}:={{AdW\circ h'_{n+1} }}$ satisfies the desired conditions, and the statement is proved.
\end{proof}

\begin{lem}\label{CtimesText00}
Let $C\in\mathcal C_0.$  Let $\ep>0,$ ${\cal F}\subset  C$ be
 any finite subset.
Suppose that $B$ is a unital separable simple \CA\, in ${\mathcal B}_0,$  $A=B\otimes U$ for some UHF-algebra of infinite
type, $\af\in KK_e(C\otimes C(\T),A)^{++}$.
Then there is a unital $\ep$-${\cal F}$-multiplicative \morp\, $\phi: C\otimes C(\T)\to A$ such that
\beq
[\phi]=\af.
\eneq
\end{lem}
\begin{proof}
Denote by $\alpha_0$ and $\alpha_1$ the induced maps induced by $\alpha$ on $K_0$-groups and $K_1$-groups.

By Theorem \ref{kkmaps}, there exist an $\mathcal F$-$\epsilon$-multiplicative map $\phi_1: C\otimes C(\T)\to A\otimes\cal K$ and a homomorphism $\phi_2: C\otimes C(\T)\to A\otimes\cal K$ with finite dimensional range such that
$$[\phi_{1}]=\alpha+[\phi_{2}]\ \mathrm{in}\ KK(C, A).$$ In particular, one has that $(\phi_1)_{*1}=\alpha_1$.
Without loss of generality, one may assume that $\phi_1, \phi_2: C \to M_r(A)$ for some integer $r$.

Since $M_r(A)\in \mathcal B_0$, for any finite subset $\mathcal G\subseteq M_r(A)$ and any $\epsilon'>0$, there are  $\mathcal G$-$\epsilon'$-multiplicative map $L_1: M_r(A)\to (1-p)M_r(A)(1-p)$ and $L_2: M_r(A)\to S_0\subseteq  pM_r(A)p$ for a C*-subalgebra $S_0\in \mathcal C_0$ with $1_{S_0}=p$ such that
\begin{enumerate}
\item $||a-L_1(a)\oplus L_2(a)||<\epsilon'$ for any $a\in\mathcal G$, and
\item $\tau((1-p))<\epsilon'$ for any $\tau\in T(M_r(A)).$
\end{enumerate}
Since $K_1(S_0)=\{0\}$, by choosing $\mathcal G$ sufficiently large and $\epsilon'$ sufficiently small, one may assume that $L_1\circ\phi_1$ is {{$\cal F$-$\ep$}}-multiplicative, and
$$[L_1\circ \phi_1]|_{K_1(C\otimes C(\T))}=(\phi_1)_{*1}=\alpha_1.$$ Moreover, since the positive cone of $K_0(C\otimes C(\T))$ is finitely generated, by choosing $\epsilon'$ even smaller, one may assume that the map
$$\kappa:=\alpha_0-[L_1\circ\phi_1]|_{K_0(C\otimes C(\T)}: K_0(C\otimes C(\T))\to K_0(A)$$ is positive.

Pick a point $x_0\in \T$, and consider the evaluation map $$\pi: C\otimes C(\T) \in f\otimes g \mapsto f\cdot g(x_0)\in C.$$ Then $\pi_{*0}: K_0(C\otimes C(\T))\to K_0(C)$ is an order isomorphism.

Pick a projection $q\in A$ with $[q]=\kappa([1])$. Since $qAq\in\mathcal B_0$, by Theorem \ref{preBot2}, there is a unital homomorphism $h: C\to qAq$ such that $$[h]_0=\kappa\circ \pi_{*0}^{-1}\quad\textrm{on $K_0(C)$},$$ and hence one has $$(h\circ\pi)_{*0}=\kappa,\quad \textrm{on $K_0(C\otimes C(\T))$.}$$
Put $\phi=(L_1\circ\phi_1)\oplus (h\circ\pi): C\otimes C(\T)\to A$, and it is clear that
$$
\phi_{*0}=[L_1\circ\phi_1]|_{K_0(C\otimes C(\T))}+\kappa=[L_1\circ\phi_1]|_{K_0(C\otimes C(\T))} + \alpha_0-[L_1\circ\phi_1]|_{K_0(C\otimes C(\T))}=\alpha_0
$$
%and
$$ \andeqn[\phi]_1=[L_1\circ\phi_1]|_{K_1(C\otimes C(\T))}=\alpha_1.$$
Since $K_*(C\otimes C(\T))$ is finitely generated and torsion free, one has that $[\phi]=\alpha$ in $KK(C\otimes C(\T), A)$.
\end{proof}

\begin{lem}\label{CtimesText0}
Let $C\in\mathcal C_0.$
%with $K_1(C)=\{0\}.$ %be as in \ref{ExtTrace}.
Let $\ep>0,$ ${\cal F}\subset  C\otimes\mathrm{C}(\T)$ be
 any finite subset, $\sigma>0,$
${\cal H}\subset (C\otimes C(\T))_{s.a.}$ be a finite subset.
Suppose that $A$ is a unital 
%separable simple 
\CA\, in ${\cal B}_0,$  $B=A\otimes U$ for some UHF-algebra {$U$} of infinite
type, $\af\in KK_e(C\otimes C(\T), B)^{++}$, and $\gamma: T(B)\to T_{\rm f}(C\otimes C(\mathbb T))$ is a continuous affine
map such that $\af$ and $\gamma$ are compatible.
There is a unital {{$\cal F$-$\ep$}}-multiplicative \morp\, $\phi: C\otimes C(\T)\to B$ such that
\begin{enumerate}
\item $[\phi]=\alpha$ and
\item  $|\tau\circ \phi(h)-\gamma(\tau)(h)|<\sigma$ for any $h\in {\cal H}$.
\end{enumerate}
Moreover, if $A\in {\cal B}_1,$ $\bt\in KK_e(C, A)^{++},$ $\gamma': T(A)\to T_f(C)$ is a continuous affine map
which is compatible with $\bt$  and ${\cal H}'\subset C_{s.a.}$ is a finite subset,
then there is also a unital \hm\, $\psi: C\to A$ such that
\beq\label{Ct-2}
[\psi]=\bt \tand |\tau\circ \psi(h)-\gamma'(\tau)(h)|<\sigma \rforal f\in {\cal H}'.
\eneq
%for all $h\in {\cal H}'.$
\end{lem}

\begin{proof}
Since $K_*(C\otimes C(\T))$ is finitely generated and torsion free, by the UCT, the element $\alpha\in KK(C\otimes C(\T), A)$ is determined by the induced maps $\alpha_0\in \mathrm{Hom}(K_0(C\otimes C(\T)), K_0(A))$ and $\alpha_1\in \textrm{Hom}(K_1(C\otimes C(\T)), K_1(A))$.
{\Wlog\, we may assume that projections in $M_2(C\otimes C(\T))$ generates $K_0(C\otimes C(\T)).$}

To simplify the proof, without loss of generality, we may assume
that $\|h\|\le 1$ for all $h\in {\cal H}.$ Fix a finite generating set $\mathcal G$ of $K_0(C\otimes C(\T))$. Since $\gamma(\tau)\in T_{\mathrm{f}}(C\otimes\mathrm{C}(\mathbb T))$ for all $\tau\in T(B)$ and $\tau(B)$ is compact, one is able to define $\Delta: {{(C\otimes C(\T))}}_+^{q, 1}\setminus\{0\} \to (0, 1)$
by
$$\Delta(\hat{h})=\frac{1}{2}\inf\{\gamma(\tau)(h):\ {{\tau}}\in T(B)\}.$$
Fix a finite generating set $\mathcal G$ of $K_0(C\otimes C(\T))$. Let $\mathcal H_1\subseteq C\otimes C(\mathbb T)$, $\delta>0$, and $K\in\mathbb N$ be the finite subset and the constants of Lemma \ref{ExtTraceC-D} with respect to $\mathcal F$, $\mathcal H$,  $\epsilon$, $\sigma/4$ (in the place of $\sigma$), and $\Delta$.

Since $A\in \mathcal B_0$ and  $U$ is {{of infinite type,}}
%self-absorbing,
for any finite subset $\mathcal G'\subseteq B$ and any $\epsilon'>0$, there are  $\mathcal G'$-$\epsilon'$-multiplicative maps $L_1: B\to (1-p)B(1-p)$ and $L_2: B\to D\otimes 1_{M_K}\subset D\otimes M_K\subseteq  pBp$ for a C*-subalgebra $D\in \mathcal C_0$ with $1_{D\otimes M_K}=p$ such that
\begin{enumerate}\setcounter{enumi}{2}
\item $||a-L_1(a)\oplus L_2(a)||<\epsilon'$ for any $a\in\mathcal G'$, and
\item $\tau((1-p))<\min\{\epsilon', \sigma/4\}$ for any $\tau\in T(B).$
\end{enumerate}
Put $S=D\otimes M_K.$
By choosing $\mathcal G'$ large enough and $\epsilon'$ small enough, one may assume
$[L_1]$ and $[L_2]$ are well-defined on $\alpha(\underline{K}(C\otimes C(\T)))$, and
\begin{equation}\label{July20-1}
\alpha=[L_1]\circ\alpha+[j]\circ [L_2]\circ\alpha,
\end{equation}
where $j: S\to A$ is the embedding. Note that since $K_1(S)=\{0\}$, one has that
$$\alpha_1=[L_1]\circ\alpha|_{K_1(C\otimes C(\T))}.$$
Define $\kappa'=[L_2]\circ\alpha|_{K_0(C\otimes C(\T))}$, which  is a homomorphism  from $K_0(C\otimes C(\T))$ to $K_0(D)$
which mapping  $[1_{C\otimes C(\T)}]$ to $[1_D].$ Let $\{e_{i,j}: 1\le i,j\le K\}$ be a system of matrix units for $M_K.$
View  $e_{i,j}\in D\otimes M_K.$ Then $e_{i,j}$ commutes with the image of $L_2.$
Define $L_2': B\to D\otimes e_{1,1}$ by $L_2'(a)=L_2(a)\otimes e_{1,1}$ for all $a\in M_r(A).$

Put $\kappa=[L_2']\circ \af|_{K_0(C\otimes C(\T))}.$ Put $D'=D\otimes e_{1,1}.$

Moreover, by choosing $\mathcal G'$ larger and $\epsilon'$ smaller,  if necessarily, there is an continuous affine map $\gamma': T(D')\to T(C\otimes C(\T))$ such that, for all $\tau\in T(A),$
\begin{enumerate}\setcounter{enumi}{4}
\item $|\gamma'(\frac{1}{\tau(e_{1,1})}\tau|_{D'})(f)-\gamma(\tau)(f)|<\sigma/4$ for any $f\in\mathcal H$,
\item $\gamma'(\tau)(h)>\Delta(\hat{h})$ for any $h\in\mathcal H_{1}$, and
\item $|\gamma'(\frac{1}{\tau(e_{1,1})}\tau|_{D'})(p)-\tau(\kappa([p]))|<\dt$ for all projections $p\in M_2(C\otimes C(\T)).$
\end{enumerate}
Then it follows from Lemma \ref{ExtTraceC-D} that there is an ${\cal F}$-$\epsilon$-multiplicative map $\phi_2: C\otimes C(\T)\to M_K(D)=S$ such that
\begin{equation*}%\label{July20-2}
(\phi_2)_{*0}=K\kappa=\kappa'
\end{equation*}
and
\begin{equation*}
|(1/K)t\circ  \phi_2(h)-\gamma'(t)(h)|<\sigma/4,\quad h\in\mathcal H,\ t\in T(D').
\end{equation*}
%for any $h\in\mathcal H$ and any $t\in T(D')$.

On the other hand, since $(1-p)A(1-p)\in \mathcal B_0$, by Lemma \ref{CtimesText00}, there is a unital  $\mathcal F$-$\epsilon$-multiplicative map $\phi_1: C\otimes C(\T)\to (1-p)A(1-p)$ such that
$$[\phi']=[L_1]\circ\alpha\quad\mathrm{in}\ KK(C\otimes C(\T), A).$$
Define $\phi=\phi_1\oplus j\circ \phi_2: C\otimes C(\T)\to (1-p)A(1-p)\oplus S\subseteq A$. Then one has
$$\phi_{*0}=(\phi_1)_{*0}+(j\circ \phi_2)_{*0}=([L_1]\circ\alpha)|_{K_0(C\otimes C(\T))}+([j\circ L_2]\circ\alpha)|_{K_0(C\otimes C(\T))}=\alpha_0$$ and
$$\phi_{*1}=(\phi_1)_{*1}+(j\circ \phi_2)_{*1}=([L_1]\circ\alpha|_{K_1(C\otimes C(\T))}=\alpha_1.$$ Hence $[\phi]=\alpha$ in $KK(C\otimes C(\T))$.

For any $h\in \mathcal H$ and any $\tau\in T(A)$, one has (note that
we have assumed $\|h\|\le 1$ for all $h\in {\cal H}$),
\begin{eqnarray*}
&&|\tau\circ\phi(h)-\gamma(\tau)(h)|\\
&<& |\tau\circ\phi(h)-\tau\circ j\circ \phi_2(h)|+|\tau\circ j\circ \phi_2(h)-\gamma(\tau)(h)|\\
&<& \sigma/4+|\tau\circ j\circ \phi_2(h)-\gamma'({1\over{\tau(e_{1,1})}}\tau|_{D'})(h)|+
|\gamma'({1\over{\tau(e_{1,1})}}\tau|_{D'})(h)-\gamma(\tau)(h)|\\
&<& |\tau\circ\phi(h)-\gamma'(\frac{1}{\tau(p)}\tau|_S)(h)|+ |\gamma'(\frac{1}{\tau(p)}\tau|_S)(h)-  \gamma(\tau)(h)|\\
&<& \sigma/4+\sigma/4+\sigma/4<\sigma.
\end{eqnarray*}

Hence the map $\phi$ satisfies the lemma.

To see the last part of the lemma holds, we note that, when $C\otimes C(\T)$ is replaced by $C$ and $A$ is assumed
to be in ${\cal B}_1,$ the only difference is that we may not use \ref{CtimesText00}.  But then we can appeal to \ref{preBot2}
to obtain $\phi_1.$   The semi-projectivity of $C$ allows us actually obtain a unital homomorphism.
\end{proof}

%{\color{Green} One should be a little bit careful that the statement of \ref{CtimesText0} is for the algebras tensoring with circle.} {\bf More specific?}

\begin{cor}\label{CorCinjj}
Let $C\in\mathcal C_0.$
%Let $\ep>0,$ ${\cal F}\subset  C$ be
%a finite subset, $1>\sigma_1>0,$  $1>\sigma_2>0,$  $\overline{{\cal U}}\subset J_c(K_1(C\otimes C(\T)))\subset U(C\otimes C(\T))/CU(C\otimes C(\T))$ be any finite subset (see \ref{Dcu}) and
%${\cal H}\subset (C\otimes C(\T))_{s.a.}$ be a finite subset.
Suppose that $A$ is a unital separable simple \CA\, in $\mathcal B_{{0}}$,  $B=A\otimes U$ for some UHF-algebra of infinite
type, $\af\in KK_{e}(C, B)^{++}$ and
%$\lambda: J_c(K_1(C)) \to U(B)/CU(B)$ is a homomorphism, and
 $\gamma: T(B)\to T_{\rm f}(C))$ is a continuous affine
map. Suppose that $(\af, \lambda, \gamma)$ is a compatible triple. Then there is a unital
\hm\, $\phi: C\to B$ such that
$$
%\begin{enumerate}
[\phi]=\af \andeqn \phi_T=\gamma.
$$
%\item  $\phi^{\ddag}=\lambda$ and
%${\rm dist}(\phi^{\ddag}(x), \lambda(x))<\sigma_1$, for any $x\in \overline{{\cal U}}$, and
%\item   $\phi_T=\gamma.$
%$|\tau\circ \phi(h)-\gamma(\tau)(h)|<\sigma_2$, for any $h\in {\cal H}$.
%\end{enumerate}
In particular, $\phi$ is a monomorphism.
\end{cor}

\begin{proof}
The proof is exactly the same argument employed in \ref{istBotC++} by using  the second part of \ref{CtimesText0}.
The reason $\phi$ is a monomorphism because $\gamma(\tau)$ is faithful for each $\tau\in T(A).$
\end{proof}

\begin{lem}\label{L85}
Let $B\in {\cal B}_0$ which satisfies the UCT, $A_1\in {\cal B}_0,$ let $C=B\otimes U_1$ and $A=A_1\otimes U_2,$ where $U_1$ and $U_2$ are unital infinite dimensional UHF-algebras. Suppose that $\kappa\in KL_e(C,A)^{++},$
$\gamma: T(A)\to T(C)$ is a continuous affine map and $\af: U(C)/CU(C)\to
U(A)/CU(A)$ is a continuous \hm\, for which $\gamma,\, \af$ and $\kappa$ are compatible. Then there exists a unital monomorphism $\phi: C\to A$ such that
\begin{enumerate}
\item $[\phi]=\kappa$ in $KL_e(C,A)^{++}$,
\item $\phi_T=\gamma$ and $\phi^{\ddag}=\af.$
\end{enumerate}
%\beq\label{L85-1}
%[\phi]=\kappa \,\,\,{\rm in}\,\,\, KL_e(C,A)^{++},\\
%\phi_T=\gamma\andeqn \phi^{\ddag}=\af.
%\eneq
\end{lem}
\begin{proof}
The proof follows a same line as that of Lemma 8.5 of \cite{Lnclasn}. By the classification theorem, one can write
$$C=\varinjlim(C_n, \phi_{n, n+1})$$
where $C_n$ is a direct sum of C*-algebras in $\mathcal C_0$ or {{in}} $\mathbf H$. Let $\kappa_n=\kappa\circ[\phi_{n, \infty}],$ $\alpha_n=\alpha\circ\phi_{n, \infty}^\ddag$, and $\gamma_n=(\phi_{n, \infty})_T\circ\gamma$. By Corollary \ref{istBotC++} or Corollary \ref{CorCinjj}, there are unital monomorphism
$$\psi_n: C_n\to A$$ such that
$$[\phi_n]=\alpha_n, \quad \psi_n^{\ddag}=\alpha_n,\quad\textrm{and}\quad (\psi_n)_T=\gamma_n.$$
In particular, the sequence of monomorphisms  $\psi_n$ satisfies
$$[\psi_{n+1}\circ\phi_{n, n+1}]=[\psi_n],\quad \psi_{n+1}^\ddag\circ\phi_{n, n+1}=\psi_{n}^\ddag,
\quad\textrm{and}\quad
(\phi_{n+1}\circ\phi_{n, n+1})=(\psi_n)_T.$$

Let $\mathcal F_n\subseteq C_n$ be a finite subset such that $\phi_{n, n+1}(\mathcal F_n)\subseteq \mathcal F_{n+1}$ and $\bigcup\phi_{n, \infty}(\mathcal F_n)$ is dense in $C$. Applying Theorem \ref{UniCtoA} with $\Delta(h)=\inf\{\gamma(\tau)(\phi_{n, \infty}(h)): \tau\in T(A)\}$, $h\in C_n^+$, there is a sequence of unitaries $u_n\in A$ such that
$$\mathrm{Ad} u_{n+1}\circ\psi_{n+1}\circ\phi_{n, n+1}\approx_{1/2^n}\mathrm{Ad} u_n\circ \psi_{n}\quad \textrm{on $\mathcal F_n$}.$$

The the maps $\{\mathrm{Ad} u_n\circ \psi_{n}: n=1, 2, ...\}$ converge to a unital homomorphism $\phi: C\to A$ which satisfies the lemma.
\end{proof}

\begin{lem}\label{alg-cut}
Let $C$ be a unital C*-algebra.
Let $p\in C$ be a full projection. Then, for any $u\in U_0(C)$, there is a unitary $v\in pCp$ such that $$\overline{u}=\overline{v\oplus (1-p)}\quad\textrm{in $U_0(C)/CU(C)$}.$$
If, furthermore, $C$ is separable and has stable rank one, then, for any $u\in U(C)$, there is a unitary $v\in pCp$ such that $$\overline{u}=\overline{v\oplus (1-p)}\quad\textrm{in $U(C)/CU(C)$}.$$
\end{lem}

\begin{proof}
It is sufficient to prove the first part of the statement. It is essentially contained in the proof of  4.5 and 4.6 of \cite{GLX-ER}.
As in the proof of 4.5 of \cite{GLX-ER}, for any $b\in C_{s.a.},$ there is $c\in pCp$ such that
$b-c\in C_0,$ where
$C_0$ is the closed subspace of $A_{s.a.}$ consisting of elements
of the form $x-y,$ where $x=\sum_{n=1}^{\infty}c_n^*c_n$ and $y=\sum_{n=1}^{\infty}c_nc_n^*$
(converge in norm) for  some sequence $\{c_n\}$ in $C.$

Now let $u=\prod_{k=1}^n\exp(ib_k)$ for some $b_k\in C_{s.a.},$ $k=1,2,...,n.$
Then there are $c_k\in pCp$ such that $b_k-c_k\in  C_0,$ $k=1,2,...,n.$
Put $v=p(\prod_{k=1}^n \exp(ic_k))p.$ Then $v\in U_0(pCp)$ and
$v+(1-p)=\prod_{k=1}^n\exp(i c_k).$   By 3.1 of \cite{Thomsen-rims},
$u^*(v+(1-p))\in CU(C).$
\end{proof}

\begin{lem}\label{CtimesText}
Let $C\in\mathcal C_0.$
Let $\ep>0,$ ${\cal F}\subset  C$ be
a finite subset, $1>\sigma_1>0,$  $1>\sigma_2>0,$  $\overline{{\cal U}}\subset J_c(K_1(C\otimes C(\T)))\subset U(C\otimes C(\T))/CU(C\otimes C(\T))$ be any finite subset (see \ref{Dcu}) and
${\cal H}\subset (C\otimes C(\T))_{s.a.}$ be a finite subset.
Suppose that $A$ is a unital separable simple \CA\, in $\mathcal B_0$,  $B=A\otimes U$ for some UHF-algebra {{$U$ of infinite}}
type, $\af\in KK_{e}(C\otimes C(\T), B)^{++},$
$\lambda: J_c(K_1(C\otimes C(\T))) \to U(B)/CU(B)$ is a homomorphism, and $\gamma: T(B)\to T_{\rm f}(C\otimes C(\mathbb T))$ is a continuous affine
map. Suppose that $(\af, \lambda, \gamma)$ is a compatible triple. Then there is a unital ${\cal F}$-$\ep$-multiplicative \morp\, $\phi: C\otimes C(\T)\to B$ such that
\begin{enumerate}
\item $[\phi]=\alpha$,
\item ${\rm dist}(\phi^{\ddag}(x), \lambda(x))<\sigma_1$, for any $x\in \overline{{\cal U}}$, and
\item  $|\tau\circ \phi(h)-\gamma(\tau)(h)|<\sigma_2$, for any $h\in {\cal H}$.
\end{enumerate}

%Moreover, theorem also holds, if $C\otimes C(\T)$ is replaced by $C$ and $B\in {\cal B}_1.$
\end{lem}

\begin{proof} {{Note that $\underline{K}(C\otimes C(\T))$ is finitely generated modulo Bockstein operations and $K_0(C\otimes C(\T))_+$ is a finitely generated semi group.  Using the inductive limit $B=\lim_{n\to\infty}(A\otimes M_{r_n}, \imath_{n,n+1}),$  one can find, for $n$ large enough,  $\alpha_n\in KK_e(C\otimes C(\T), A\otimes M_{r_n})^{++}$  such that $\alpha=\alpha_n\times [  \imath_n]$ where $ [  \imath_n] \in KK(A\otimes M_{r_n}, B)$ is induced by the inclusion $ \imath_n: A\otimes M_{r_n} \to B$. Replacing  $A$ by $A\otimes M_{r_n}$, without lose of generality we can assume that $\alpha=\alpha_1 \times [\imath]$, where $\alpha_1\in KK_e(C\otimes C(\T), A)^{++}$ and $ \imath: A \to A\otimes U=B$ is the inclusion.  Note that  $J_c(K_1(C\otimes C(\T)))$, same argument as above, we know that if the integer $n$ above is large enough, then there is a  $\lambda_n: J_c(K_1(C\otimes C(\T))) \to U(A\otimes M_{r_n})/CU(A\otimes M_{r_n})$ such that $|\imath_n^{\ddag}\circ \lambda_n(u)-\lambda(u)|$ arbitrarily small (e.g smaller than {$\frac{\sigma_1}4$}) for all $u\in  \overline{{\cal U}}$. When we replace $A$ by $ A\otimes M_{r_n}$, we can assume $\lambda=\imath^{\ddag}\circ\lambda_n$, with $\lambda_1: J_c(K_1(C\otimes C(\T))) \to U(A)/CU(A)$ and $ \imath^{\ddag}: U(A)/CU(A)\to U(B)/CU(B)$ being induced by the inclusion map.  Furthermore, we can assume $\lambda_1$ being compatible with $\af_1$.}
}

Without loss of generality, we may assume that $\|h\|\le 1$ for all $h\in {\cal H}.$ Let $p_i, q_i\in M_k(C)$ be projections such that
$\{[p_1]-[q_1], ..., [p_d]-[q_d]\}$ forms a set of  standard generators of $K_0(C)$ (as an abelian group)
for some integer $k\ge 1$.
By choosing a specific $J_c,$
without loss of generality, one may assume that
$$
\overline{\mathcal U}=\{(({\bf 1}_k-p_i)+p_i\otimes z)(({\bf 1}_k-q_i)+q_i\otimes z^*): 1\le i\le d\},
$$
where $z\in C(\T)$ is the identity function on the unit circle.
Put $u_i'=({\bf 1}_k-p_i)+p_i\otimes z)(({\bf 1}_k-q_i)+q_i\otimes z^*).$ Hence $\{[u_1'], ..., [u_d']\}$ is a set of standard generators of $K_1(C\otimes C(\T))\cong K_0(C)\cong \mathbb Z^d.$ Then $\lambda$ is a homomorphism from $\mathbb{Z}^d$ to $U(B)/CU(B)$.

Let $\pi_e: C\to F_1=\bigoplus_{i=1}^l M_{n_i}$ be defined in \ref{DfC1}. By \ref{2Lg13},  the map $(\pi_e)_{*0}$ induces an embedding of $K_0(C)$ to $\mathbb{Z}^l$, and the map $(\pi_e\otimes\id)_{*1}$ induces an embedding of $K_1(C\otimes C(\T))\cong \mathbb{Z}^d$ to $K_1(\bigoplus_{i=1}^lM_{n_i}\otimes C(\T)) \cong \mathbb{Z}^l$. Choose $J_c({{K_1(}}\bigoplus_{i=1}^lM_{n_i}\otimes C(\T){{)}})$ to be the subgroup generated by $\{\overline{e_i\otimes z_i\oplus (1-e_i)}; i=1, ..., l\}$, where $e_i$ is a rank one projection of $M_{n_i}$ and $z_i$ is the standard unitary generator of  $i$-th copy of $C(\T).$  Note that the image of $J_c(K_1(C\otimes C(\T)))$ under $\pi_e$ is contained in $J_c(K_1(\bigoplus_{i=1}^lM_{n_i}\otimes C(\T)))$. Denote by $w_j=e_j\otimes z_j\oplus (1-e_j)$, $1\leq j\leq l$.

We write $B=B_0\otimes U_{{2}},$ and {{ $B_0=A\otimes U_{1},$ with $U=U_1\otimes  U_2$, both $U_1$ and $U_2$ being UHF algebra of infinite type. Denoted by $\imath_1:A \to  B_0$, $\imath_2:B_0 \to  B$ and $\imath=\imath_2\circ \imath_1:A \to  B$ the inclusion maps. Recall $\af=\af_1\times [\imath] \in KK(C,B)$.}}

%where $B_0=A\otimes U\cong B.$---we do not have the condition $B$ is self absorbing---Gong

By applying Lemma \ref{CtimesText0},
one obtains a unital ${\cal F}'$-$\ep'$-multiplicative \morp\, $\psi: C\otimes\mathrm{C}(\T) \to B_0$ such that
\begin{equation}\label{July22-a1}
[\psi]=\alpha_{{1}}{{\times [\imath_1]}}\andeqn
\end{equation}
\begin{equation}\label{July22-a2}
|\tau\circ \psi(h)-\gamma(\tau)(h)|<\min\{\sigma_1, \sigma_2\}/3 \rforal  h\in {\cal H}\andeqn \rforal \tau\in\mathrm{T}(B_0),
\end{equation}
where $\ep/2>\ep'>0$ and ${\cal F}_1\supset {\cal F}.$ { {(Note that $T(B_0)=T(A)=T(B)$, and the map $\gamma: T(B)\to T_{\rm f}(C\otimes C(\mathbb T))$ can be regarded as a map with domain $T(B_0)$.}}
We may assume that $\ep'$ is sufficiently   small and ${\cal F}_1$ is sufficiently  large so that
not only (\ref{July22-a1}) and (\ref{July22-a2}) make sense
but also  that $\psi^{\ddag}$ is well defined on $\bar {\mathcal U}$ and induces a
homomorphism from $J_c(K_1(C\otimes C(\T)))$ to $U(B_{{0}})/CU(B_{{0}}).$

Let $M$ be the integer in \ref{MextC}.

For any $\ep''>0$ and any finite subset ${\cal F}''\subset B_0,$ since $B_0$ has the Popa condition and has (SP) property,
there exist a non-zero projection $e\in B_0$ and a unital ${\cal F}''$-$\ep''$-multiplicative \morp\, $L_0: B_0\to F\subset eB_0e,$
where $F$ is a finite  dimensional and $1_{{F}}=e$ and a  unital ${\cal F}''$-$\ep''$-multiplicative \morp\, $L_1: B_0\to
(1-e)B_0(1-e)$ such that
\beq\label{July22-b1}
\|b-\imath\circ L_0(b)\oplus L_1(b)\|<\ep''\rforal b\in {\cal F}'',\\
\|L_0(b)\| \ge \|b\|/2\rforal b\in {\cal F}'' \andeqn\\ \label{Dec1407-n1}
\tau(e)<\min\{\sigma_1/2, \sigma_2/2\}\rforal \tau\in T(B_0),
\eneq
where $\imath: F\to eB_0e$ is the embedding and $L_1(b)=(1-p)b(1-p)$ for all $b\in B_0.$

Since the positive cone of $K_0(C\otimes C(\T))$ is finitely generated, with sufficiently small
$\ep''$ and sufficiently large ${\cal F}''$, one may assume that
$[L_0\circ \psi]|_{K_0(C\otimes C(\T))}$ is positive. Moreover, one may assume that $(L_0\circ \psi)^{\ddag}$ and
$(L_1\circ \psi)^{\ddag}$ are  well defined and induce homomorphisms from $J_c(K_1(C\otimes C(\T)))$ to $U(B_{{0}})/CU(B_{{0}})$. One may also assume that $[L_1\circ \psi]$ is well defined.  Moreover, we may assume that $L_i\circ \psi$ is ${\cal F}$-$\ep$-multiplicative for $i=0,1.$

There is a projection $E_c\in U_{{2}}$ such that $E_c$ is a direct sum of $M$ copies of some non-zero projections $E_{c,0}\in U_{{2}}.$
Put $E=1_{{U_2}}-E_c.$

Define $\phi_0: C\otimes C(\T)\to F\otimes EU_{{2}}E\to eB_0e\otimes EU_{{2}}E$ by $\phi_0(c)= L_0\circ \psi(c)\otimes E{{(\in B)}}$ for all $c\in C\otimes C(\T)$ and
define $\phi_1': C\to F\otimes E_cU_{{2}}E_c$ by $\phi_1'(c)=L_0\circ \psi(c)\otimes E_c$ for all $c\in C.$
Note that
$\phi_0$ is also { {${\cal F}$-$\ep$-}}multiplicative and $\phi_0^{\ddag}$ is also well defined as
$(L_0\circ \psi)^{\ddag}$ is.  Moreover $[\phi_1']$ is well defined.  Define
$$L_2={{\imath_2 \circ}}L_1\circ \psi+\phi_0{{: C\otimes C(\T) \to \big((1-e)B_0(1-e)\otimes 1_{U_2}\big)\oplus \big(eB_0e\otimes EU_2E\big)~~(\subset B)}}.$$

Denote by
$$\lambda_0=\lambda-L_2^{\ddag}=\lambda-\phi_0^{\ddag}-({{\imath_2 \circ}}L_1\circ \psi)^{\ddag}: J_c(K_1(C\otimes C(\T)))\to U(B)/CU(B).$$
{ { Note that $L_0$ factors through the finite dimensional algebra $F$ and therefore $[L_0]=0$ on {{$K_1(B_0)$}}, consequently $[\phi_0]|_{K_1(C\otimes C(\T))}=0$  and $[L_1\circ \psi]=[\imath_1]\circ[\alpha_1]$ on ${K_1(C\otimes C(\T))}$. Hence $[\imath_2\circ L_1\circ \psi]=\alpha$ on $K_1$. Furthermore $\alpha$ is compatible with $\lambda$. We know that the image of $\lambda_0$ is in $U_0(B)/CU(B).$}}
%I add the above, since in general $U(B)/CU(B) is not divisible, only U_0(B)/CU(B) is, and we need to use it below--Gong

Note that, by \ref{UCUdiv},
the group $U_{{0}}(B)/CU(B)$ is divisible. It is an injective abelian group.
Therefore there  is a homomorphism ${{\tilde \lambda}}: J_c(\bigoplus_{i=1}^lM_{n_i}\otimes C(\T))\to U_{{0}}(B)/CU(B)$
such that
\begin{equation}\label{July22-a4}
{{\tilde \lambda}}\circ (\pi_e)^{\ddag}=\lambda_0-L_2^{\ddag}.
\end{equation}

Let $\bt=[{ {L_0\circ \psi}}]|_{K_0(C)}{{:K_0(C) \to}}K_0(F)=\Z^n.$  Let $R_0\ge 1$ be the integer given by \ref{MextC}.
There is a unital \SCA\, $M_{MK}\subset E_cU_{{2}}E_c$ such that $K\ge R_0$ { {and such that $E_cU_{2}E_c$ can be written as $M_{MK}\otimes U_3$}}.  It follows from \ref{MextC}
that there is a positive \hm\, $\bt_1: K_0(F_1)\to K_0(F)$ such that
$\bt_1\circ (\pi_e)_{*0}=MK\bt.$
Let $h: F_1\to F\otimes M_{MK}$ be the unital \hm\, such that $h_{*0}=\beta_1.$
Put $\phi_1''=h\circ \pi_e{{: C\to F\otimes M_{MK}}}$, and then one has that $(\phi_1'')_{*0}=MK\bt.$ Let $J: M_{MK}\to E_cU_{{2}}E_c$ be the embedding.
One verifies that
\beq\label{120714-18nn2}
(\imath_{{F}}\otimes J)_{*0}\circ (\phi_1'')_{*0}=(\imath_{{F}}\otimes J)_{*0}\circ MK\bt={{\tilde\imath}}_{*0}\circ (\phi_1')_{*0}{{,}}
\eneq
{ { where $\imath_{F}:F\to eB_0e$ and ${\tilde \imath}: F\otimes E_cU_2E_c \to eB_0e\otimes E_cU_2E_c $
be the unital embedding.
%%%embeddings.
}}

Choose a unitary  $y_i\in { {(\imath_F\otimes J\circ h)}}(e_j)B{ {(\imath_F\otimes J\circ h)}}(e_j)$ such that
$${\bar y}_j={{\tilde \lambda}}(w_j),\,\,\,j=1,2,...,l{{,}}$$
{ {where { we} recall that $w_j=e_j\otimes z_j\oplus (1-e_j)\in F_1\otimes C(\T)=\oplus_j M_{n_j}\otimes C(\T)$ is the standard generator of $K_1(M_{n_j}\otimes C(\T)$. Let $1_j$ be the unit of $M_{n_j}\subset F_1$, then $1_j=\underbrace{e_j\oplus e_j\oplus \cdots\oplus e_j}_{n_j}$.}}

Define ${ {\tilde y}}_j={\rm diag}(\overbrace{y_j,y_j,...,y_j}^{{n_j}}){{\in  (\imath_F\otimes J\circ h)(1_j)B (\imath_F\otimes J\circ h)(1_j)}},$ $j=1,2,...,l.$ { {Then ${\tilde y}_j$ commutes with $\imath_F\otimes J (F_1)$.}}

Define ${\tilde  \phi}_1: F_1\otimes C(\T)\to { {\big((\imath_F\otimes J)\circ \phi_1''\big)}}(1_C)B{ {\big((\imath_F\otimes J)\circ \phi_1''\big)}}(1_C)$ by ${\tilde \phi}_1(c_j\otimes f)={ {\big((\imath_F\otimes J)\circ \phi_1''\big)}}(c_j)f({{\tilde y}}_j)$ for
all $c_j\in M_{n_j}$ and $f\in C(\T).$ Define $\phi_1={\tilde \phi}_1\circ (\pi_e\otimes {\rm id}_{C(\T)}).$
Then, by identifying $K_0(C\otimes C(\T))$ with $K_0(C),$ one has
\begin{equation}\label{July23-4}
(\phi_1)_{*0}={{\tilde \imath}}_{*0}\circ (\phi_1')_{*0}\andeqn (\phi_1)^{\ddag}={{\tilde \lambda}}.
\end{equation}
Define $\phi=\phi_0\oplus \phi_1\oplus {{\imath_2 \circ}} L_1\circ \psi.$
We verify that, by \eqref{July22-a2} and \eqref{Dec1407-n1},
\begin{equation*}
|\tau\circ \phi(h)-\gamma(\tau)(h)|<\sigma_2/3+\sigma_2/3=2\sigma_2/3,\quad h\in {\cal H}.
\end{equation*}

It is ready to verify that
\beq\label{120714-18-nn3}
\phi_{*0}=\af|_{K_0(C\otimes C(\T))}\andeqn \phi^{\ddag}=\lambda.
\eneq
Thus, since $\lambda$ is compatible with $\af,$
\beq\label{120714-18-nn4}
\phi_{*1}=\af|_{K_1(C\otimes C(T))}.
\eneq
Since  $K_{*i}(C\otimes C(\T))\cong K_0(C)$ is free and finitely generated,
one concludes that
\begin{equation*}
[\phi]=\af.
\end{equation*}

\end{proof}

\begin{cor}\label{CorCinj}
Let $C\in\mathcal C_0$ and $C_1=C\otimes C(\T).$
%Let $\ep>0,$ ${\cal F}\subset  C$ be
%a finite subset, $1>\sigma_1>0,$  $1>\sigma_2>0,$  $\overline{{\cal U}}\subset J_c(K_1(C\otimes C(\T)))\subset U(C\otimes C(\T))/CU(C\otimes C(\T))$ be any finite subset (see \ref{Dcu}) and
%${\cal H}\subset (C\otimes C(\T))_{s.a.}$ be a finite subset.
Suppose that $A$ is a unital separable simple \CA\, in $\mathcal B_0$,  $B=A\otimes U$ for some UHF-algebra of infinite
type, $\af\in KK_{e}(C_1, B)^{++},$
$\lambda: J_c(K_1(C)) \to U(B)/CU(B)$ is a homomorphism, and
 $\gamma: T(B)\to T_{\rm f}(C_1))$ is a continuous affine
map. Suppose that $(\af, \lambda, \gamma)$ is a compatible triple. Then there is a unital
\hm\, $\phi: C_1\to B$ such that
%\begin{enumerate}
$$ [\phi]=\af,
 \phi^{\ddag}=\lambda\andeqn
%${\rm dist}(\phi^{\ddag}(x), \lambda(x))<\sigma_1$, for any $x\in \overline{{\cal U}}$, and
\phi_T=\gamma.$$
%$|\tau\circ \phi(h)-\gamma(\tau)(h)|<\sigma_2$, for any $h\in {\cal H}$.
%\end{enumerate}
In particular, $\phi$ is a monomorphism.
\end{cor}

\begin{proof}
The proof is exactly the same argument employed in \ref{istBotC++} by using   \ref{CtimesText}.
\end{proof}

\begin{cor}\label{CCorCinj}
Let $C\in\mathcal C_0$ and let $C_1=C$ or $C_1=C\otimes C(\T).$
%Let $\ep>0,$ ${\cal F}\subset  C$ be
%a finite subset, $1>\sigma_1>0,$  $1>\sigma_2>0,$  $\overline{{\cal U}}\subset J_c(K_1(C\otimes C(\T)))\subset U(C\otimes C(\T))/CU(C\otimes C(\T))$ be any finite subset (see \ref{Dcu}) and
%${\cal H}\subset (C\otimes C(\T))_{s.a.}$ be a finite subset.
Suppose that $A$ is a unital separable simple \CA\, in $\mathcal B_0$,  $B=A\otimes U$ for some UHF-algebra of infinite
type, $\af\in KK_{e}(C_1, B)^{++}$
%$\lambda: J_c(K_1(C\otimes C(\T))) \to U(B)/CU(B)$ is a homomorphism,
and $\gamma: T(B)\to T_{\rm f}(C_1)$ is a continuous affine
map. Suppose that $(\af, \gamma)$ is  compatible. Then there is a unital
\hm\, $\phi: C_1\to B$ such that
%\begin{enumerate}
$$
[\phi]=\alpha, \,\,\,
%\andeqn \phi_T=\gamma.
%\item
\phi^{\ddag}=\lambda \,\tand
%${\rm dist}(\phi^{\ddag}(x), \lambda(x))<\sigma_1$, for any $x\in \overline{{\cal U}}$, and
%\item
\phi_T=\gamma.$$
%$|\tau\circ \phi(h)-\gamma(\tau)(h)|<\sigma_2$, for any $h\in {\cal H}$.
%\end{enumerate}
In particular, $\phi$ is a monomorphism.
\end{cor}

\begin{proof}
To apply \ref{CorCinj}, one needs $\lambda.$
Note that $J_c(K_1(C_1)$ is isomorphic to $K_1(C_1)$ which is  finitely generated.
Let $J_c^{(1)}: K_1(B)\to U(B)/CU(B)$ be the splitting map defined in \ref{Dcu}.
Define $\lambda=J_c^{(1)}\circ \af|_{K_1(C_1)}.$ Then
$(\af, \lambda, \gamma)$ is compatible. This corollary then follows from the previous one.
\end{proof}

\begin{thm}\label{Next2014/09}
Let $X$ be a finite CW complex and let $C=PM_n(C(X))P,$ where
$n\ge 1$ is an integer and $P\in M_n(C(X))$ is a projection.
Let $A_1\in {\cal B_{{0}}}$ and let $A=A_1\otimes U$ for UHF-algebra {{$U$}} of infinite type.
Suppose $\af\in KL_e(C, A)^{++},$
$\lambda: U_{\infty}(C)/CU_{\infty}(C)\to  U(A)/CU(A)$ be a continuous \hm\, and $\gamma: T(A)\to T_f(C)$ be a
continuous affine map which are compatible. Then there exists a unital \hm\,
$h: C\to A$ such that
\beq\label{Next1409-1}
[h]=\af,\,\,\, h^{\ddag}=\lambda \andeqn
h_T=\gamma.
\eneq
\end{thm}

\begin{proof}
The proof is similar to that of 6.6 of \cite{Lnclasn}.
To simplify the notation, without loss of generality, we may assume
that $X$ is connected. Furthermore, a standard argument shows that the general case case
can be reduced to the case that $C=C(X).$
We assume that $U(M_N(C))/U_0(M_N(C))=K_1(C).$
Therefore, in this case, $$U(M_N(C))/CU(M_N(C))=U_{\infty}(C)/CU_{\infty}(C).$$
Write $K_1(C)=G_1\oplus Tor(K_1(C)),$  where $G_1$ is the free part of $K_1(C).$
Fix a point $\xi\in X$ let $C_0=C_0(X\setminus \{\xi\}).$  Note that $C_0$ is an ideal of $C$ and
$C/C_0{{\cong \C}}$. Write
\beq\label{Next1409-2}
K_0(C)=\Z\cdot [1_C]\oplus K_0(C_0).
\eneq
Let $B\in {\cal B}_0$ be a unital separable simple \CA\, constructed in \ref{smallmap}
such that
\beq\label{Next1409-3}
(K_0(B), K_0(B)_+, [1_B], T(B), r_B)=(K_0({{A}}), K_0({{A}}), [1_{{A}}], T({{A}}), r_{{A}})
\eneq
and $K_1(B)=G_1{{\oplus}} Tor(K_1({{A}})).$
%Let $B_1=B\otimes U$ and $\imath: B\to B\otimes 1_U\subset B_1$ be the embedding.

Put
\beq\label{Next1510-1}
\Delta(\hat{g})=\inf\{\gamma(\tau)(g): \tau\in T(A)\}.
\eneq
For each $g\in C_+\setminus\{0\},$ since $\gamma(\tau)\in T_f(C),$ $\gamma$ is continuous and $T(A)$ is compact,
$\Delta(\hat{g})>0.$

Let $\ep>0,$ ${\cal F}\subset C$ be a finite subset,  let $1>\sigma_1, \sigma_2>0,$
${\cal H}\subset C_{s.a.}$ be a finite
subset,  ${\cal U}\subset U(M_N(C))/CU(M_N(C))$ be a finite subset.
We assume that ${\cal U}={\cal U}_0\cup {\cal U}_1,$ where
${\cal U}_0\subset U_0(M_N(C))/CU(M_N(C))$ and ${\cal U}_1\subset  J_c(K_1(C))\subset U(M_N(C))/CU(M_N(C)).$

For each $u\in {\cal U}_0,$ write $u=\prod_{j=1}^{n(u)}\exp(\sqrt{-1} a_i(u)),$ where
$a_i(u)\in M_N(C)_{s.a.}.$  Write
\beq\label{Next1409-4}
a_i(u)=(a_i^{(k,j)}(u))_{N\times N},\,\,\,i=1,2,...,n(u).
\eneq
Write
\beq\label{Next1409-5}
c_{i,k,j}(u)={a_i^{(k,j)}(u)+(a_i^{(k,j)})^*\over{2}}\andeqn d_{i,k,j}={a_i^{(k,j)}(u)-(a_i^{(k,j)})^*\over{2i}}.
\eneq
Put
\beq\label{Next1409-6}
M=\max\{\|c\|, \|c_{i,k,j}(u)\|,\|d_{i,k,j}(u)\|: c\in {\cal H}, u\in {\cal U}_0\}.
\eneq
Choose a non-zero projection $e\in {{B}}$ such that
\begin{equation*}
\tau(e)<{\min\{\sigma_1,\sigma_2\}\over{16N^2(M+1)\max\{n(u): u\in {\cal U}_0\}}}\rforal \tau\in T({{B}}).
\end{equation*}
Let $B_2=(1-e){{B}}(1-e).$
%Note that $T(B_1)=T(B).$

In what follows we will use the identification (\ref{Next1409-3}).
Define $\kappa_0\in Hom(K_0(C). K_0(B_2))$ as follows.
Define $\kappa_0(m[1_C])=m[1-e]$ for $m\in \Z$ and
$\kappa_0|_{K_0(C_0)}=\af|_{K_0(C_0)}.$ { {Note that $K_1(B)=G_1\oplus Tor(K_1(A))$ and that $\af$ { induces} a map $[\af]|_{Tor(K_1(C))}: Tor(K_1(C)) \to Tor(K_1(A))$. Using the given decomposition $K_1(C)=G_1\oplus
Tor(K_1(C)$, we can define $\kappa_1: K_1(C) \to K_1(B)$ by $\kappa_1|_{G_1}={\rm id}$ and $\kappa_1|_{Tor(K_1(C))}=[\af]|_{Tor(K_1(C))}$.}}
%define $\kappa_1\in Hom(K_1(C), K_1(B_1))$ as follows
%$\kappa_1|_{K_1(C)}=\imath_{*1}\circ {\rm id}_{K_1(C)}.$
%Let $\kappa_0'=\af|_{K_0(C)}$ and let $\kappa_1'={\rm id}_{K_1(C)},$ where
%we identify $K_1(B)$ with $K_1(C).$

 By the Universal Coefficient Theorem, there is
$\kappa\in KL(C, B_2)$ which gives the above two \hm s {{$\kappa_0, \kappa_1$}}.
Note that $\kappa\in KL_e(C,B_2)^{++},$  since $K_0(C_0)= {{\rm ker}\rho_C(K_0(C))}$.
Choose
$$
{\cal H}_1={\cal H}\cup\{c_{i,k,j}(u), d_{i,k,j}(u): u\in {\cal U}_0\}.
$$
Every tracial state $\tau'$ of $B_2$ has the form
$\tau'(b)=\tau(b)/\tau(1-e)$ for all $b\in B_2$ { { for some $\tau\in T(B)$}}.  Let $\gamma': T(B_2)\to
T({{C}})$ be defined { {as follows. For $\tau'\in T(B_2)$ as above, define $\gamma'(\tau')(f)=\gamma(\tau)(f)$ for $f\in C$.}}
%by $\gamma'(t)(b)=\gamma(t)(b)/\tau(1-e)$ for all $b\in B_2$ and$t\in T(C).$

It follows from \ref{istBotC} that
there exists a sequence of unital \morp s $h_n: C\to B_2$ such that
\begin{eqnarray*}
\lim_{n\to\infty}\|h_n(ab)-h_n(a{{)h_n(}}b)\|=0\rforal a,\, b\in  C,\\
{[}h_n{]}=\kappa \,\,\,{\text (-} K_*(C)\,\,\text{is\,\,\, finitely \,\, generated\,\,)}\andeqn\\
\lim_{n\to\infty}\max\{ |\tau\circ h_n(c)-\gamma'(\tau)(c)|: \tau\in T(B_2)\}=0{{ \rforal c\in C}}.
\end{eqnarray*}
Here we may assume that $[h_n]$ is well defined for all $n$ and
\beq\label{Next1409-8}
|\tau\circ h_n(c)-\gamma(\tau)(c)|<{\min\{\sigma_1, \sigma_2\}\over{8N^2}},\,\,\, n=1,2,....
\eneq
for all $c\in {\cal H}_1$ and for all $\tau\in T(B_2).$
Choose $\theta\in KL({{B}}, {{A}})$ such that
it gives the identification of (\ref{Next1409-3}), and, $\theta|_{G_1}=\af|_{G_1}$
and  $\theta|_{Tor(K_1({{A}}))}={\rm id}_{Tor(K_1({{A}}))}.$ {{ Let $e'\in A $ be a projection such that $[e']\in K_0(A)$ corresponds $[e]\in K_0(B)$ under the identification \ref{Next1409-3}.}}
Let $\bt=\af-\theta\circ \kappa.$ Then
\beq\label{Next1409-9}
\bt([1_C])={ {[}}e{{']}},\,\,\,\bt_{K_0(C_0)}=0\andeqn \bt_{K_1(C)}=0.
\eneq
Then $\bt\in KL_e(C, {{e'Ae'}}).$ It follows \ref{istBotC} that there exists a sequence of unital \morp s
$\phi_{0,n}: C\to {{e'Ae'}}$ such that
\beq\label{Next1409-10}
\lim_{n\to\infty}\|\phi_{0,n}(ab)-\phi_{0,n}(a)\phi_{0,n}(b)\|=0\andeqn [\phi_{0,n}]=\bt.
\eneq
Note that, for each $u\in U(M_N(C))$ with ${\bar u}\in {\cal U}_0,$
\beq\label{Next1409-11}
D_{C}(u)=\overline{\sum_{i=1}^{n(u)}\widehat{a_j(u)}},
\eneq
where $\widehat{c}(\tau)=\tau(c)$ for all $c\in C_{s.a.}$ and $\tau\in T(C).$
Since $\kappa$ and $\lambda$ are compatible,  we compute, for ${\bar u}\in {\cal U}_0,$
\beq\label{Next1409-12}
{\rm dist}((h_n)^{\ddag}({\bar u}), \lambda({\bar u}))<\sigma_2/8.
\eneq
Fix a pair of large integers $n,m,$
define $\chi_{n,m}: J_c(G_1)\to \Aff(T({{A}}))/\overline{\rho_{{{A}}}(K_0({{A}}))}$ by
\beq\label{Next1409-13}
\lambda|_{J_c(G_1)}-(h_n)^{\ddag}|_{J_c(G_1)}-(\phi_{0,m}^{\ddag}|_{J_c(G_1)}.
\eneq
Viewing $J_c(G_1)$ as subgroup of $J_c(K_1({{B}})){{=J_c(K_1(B_2))}},$ define
$\chi_{n,m}$ on $Tor(K_1(B_2))$  to be zero, we obtain a \hm\,
$\chi_{n,m}: J_c(K_1(B_2))\to \Aff({{T(A)}})/\overline{\rho_{{{A}}}(K_0({{A}}))}.$
It follows from \ref{L85} that there is a unital \hm\, $\psi: B_2\to (1-e{{'}}){{A}}(1-e{{'}})$ such that
\beq\label{Next1409-14}
[\psi]=\theta,\,\, \psi_T={\rm id}_{T({{A}})}\andeqn\\
\psi^{\ddag}|_{J_c(K_1(B_2))}=\chi_{n,m}|_{J_c(K_1(B_2))}+J_c\circ \theta|_{K_1(B_2)},
\eneq
where we identify $K_1(B_2)$ with $K_1({{B}})$
By (\ref{Next1409-14}),
\beq\label{Next1409-15}
\psi^{\ddag}|_{\Aff(T(B_2))/\overline{\rho_{B_2}(K_0(B_2))}}={\rm id}.
\eneq
Define $L(c)=\phi_{0,m}(c)\oplus \psi\circ h_n(c)$ for all $c\in C.$
It follows, by choosing sufficiently large $m$ and $n,$
one has that $L$ is $\ep$-${\cal F}$-multiplicative,
\beq\label{Next1409-16}
[L]=\af,\\
\max\{|\tau\circ \psi(f)-\gamma(\tau)(f)|: \tau\in T({{A}})\}<\sigma_1\rforal f\in {\cal H}\andeqn\\
{\rm dist}(L^{\ddag}({\bar u}), \lambda({\bar u}))<\sigma_2.
\eneq
This implies that that there is a sequence of \morp s $\psi_n: C\to {{A}}$ such that
\beq\label{Next1409-17}
\lim_{n\to\infty}\|\psi_n(ab)-\psi_n(a)\psi_n(b)\|=0\rforal a,b\in C,\\
{[}\psi_n{]}=\af,\\
\lim_{n\to\infty}\max\{|\tau\circ \psi_n(c)-\gamma(\tau)(c)|: \tau\in T(A_1)\}=0\rforal c\in C_{s.a.}\andeqn\\
\lim_{n\to\infty}{\rm dist}(\psi_n^{\ddag}({\bar u}), \lambda({\bar u}))=0\rforal u\in U(M_N(C))/CU(M_N(C)).
\eneq
Finally, by applying \ref{UniCtoA}, as in the proof of \ref{istBotC++},
using $\Delta/2$ above, we obtain a unital \hm\,
$h: C\to {{A}}$ such that
\beq\label{Next1409-18}
{[}h{]}=\af,\,\,\, h_T=\gamma\andeqn h^{\ddag}=\lambda
\eneq
as desired.

\end{proof}

\begin{thm}\label{BB-exi}
Let  $C\in \mathcal C_0$
and let
 $G= K_0(C)$. Write
 ${G}=\Z^k$ with
 $\Z^k$ generated by
$$
\{x_{1}=[p_1]-[q_1], x_2=[p_2]-[q_2], ..., x_{k}=[p_{k}]-[q_{k}]\},
$$
where $p_i, q_i\in M_n(C)$ (for some integer $n\ge 1$) are projections,
$i=1,...,k$.

Let $A$ be a simple C*-algebra in $\mathcal B_0$, and let $B=A\otimes U$ for a UHF algebra $U$ of infinite type. Suppose that $\phi: C\to B$ is a monomorphism. Then, for any finite subsets  $\mathcal F\subseteq C$ and $\mathcal P\subseteq \underline{K}(C)$, any $\ep>0$ and $\sigma>0$, any homomorphism $$\Gamma: \Z^k \to U_0({{B}})/CU({{B}}),$$ there is a unitary $w\in B$ such that
\begin{enumerate}
\item $\|[\phi(f), w]\|<\ep$, for any $f\in\mathcal F$,
\item $\mathrm{Bott}(\phi, w)|_{\mathcal P}=0$, and
\item $\mathrm{dist}(
\overline{\langle(({\mathbf 1}_n-\phi(p_i))+\phi(p_i)\tilde{w})(({\mathbf 1}_n-\phi(q_i))+\phi(q_i)\tilde{w}^*)\rangle},
 \Gamma(x_i)))<\sigma$, for any $1\leq i\leq k,$
{ where $\tilde{w}=\mathrm{diag}(\overbrace{w, ..., w}^n)$.}
\end{enumerate}
%\beq
%\|[\phi(f), w]\|<\ep\quad\forall f\in\mathcal F,\,\
%\mathrm{Bott}(\phi, w)|_{\mathcal P}=0,\andeqn
%\eneq
%$$\mathrm{dist}(
%\overline{\langle(({\mathbf 1}_n-\phi(p_i))+\phi(p_i)\tilde{w})(({\mathbf 1}_n-\phi(q_i))+\phi(q_i)\tilde{w}^*\rangle},
% \Gamma(x_i)))<\sigma, \quad\forall 1\leq i\leq k,$$
%{ where $\tilde{w}=\mathrm{diag}(\overbrace{w, ..., w}^n)$.}
\end{thm}

\begin{proof}{{ Write $B=\lim_{{n\to \infty}}(A\otimes M_{r_n},\imath_{n,n+1})$. Using the fact that $C$ is
%generated by a stable relation
{ semi-projective} (see \cite{ELP1}), one can construct a sequence of homomorphisms $\phi_n: C\to A\otimes M_{r_n}$ such that $\imath_n\circ\phi_n(c) \to \phi(c)$ for all $c\in C$. Without loss of generality, we can assume $\phi=\imath\circ \phi_1$ for a homomorphism $\phi_1 C\to A$ (replacing $A$ by $A\otimes M_{r_n}$), where $\imath: A \to A\otimes U=B$ is the standard inclusion.}}

Without loss of generality, one may assume that $||f||\leq 1$ for any $f\in\mathcal F$.

For any nonzero positive element $h\in C$ with norm at most $1$,  define $$\Delta(h)=\inf\{\tau(\phi(h));\ \tau\in T(B)\}.$$ Since $B$ is simple, one has that $\Delta(h)\in (0, 1)$.

Let $\mathcal H_1\subseteq A_+\setminus\{0\}$, $\mathcal G\subseteq A$, $\delta>0$, $\mathcal P\subseteq\underline{K}(A)$, $\mathcal H_2\subseteq A_{s.a.}$ and $\gamma_1>0$ be the finite subsets and constants of Corollary \ref{UniCtoA} with respect to $C$ (in the place of $C$), $\mathcal F$, $\epsilon/2$ and $\Delta/2$ (since $K_1(C)=\{0\}$, one does not need $\mathcal U$ and $\gamma_2$).

Note that $B={{A}}\otimes U.$ Pick a unitary $z\in U$
with ${\rm sp}(u)=\T$ and consider the map $\phi': C\otimes C(\mathbb T)\to B\otimes U$ by $$a\otimes f \mapsto \phi(a)\otimes f(z).$$
{{(Recall that $\phi(a)=\phi_1(a)\otimes 1_U$.)}}
Denote by $$\gamma=(\phi')_*: T(B)\to T_{\mathrm{f}}(C\otimes C(\mathbb{T})).$$ Also define $$\alpha:=[\phi'] \in KK(C\otimes C(\mathbb T), B).$$

Note that $K_1(C\otimes C(\mathbb T))=K_0(C)=\mathbb Z^k$. Identifying $J_c(K_1(C\otimes C(\mathbb T)))$ with $\mathbb Z^k$, and define a map $\lambda: J_c(K_1(U(C\otimes C(\mathbb T))))\to U_{{0}}(B)/CU(B)$
%{\color{Green} (should  $J_c(U(C\otimes C(\mathbb T)))$ be $J_c(K_1(C\otimes C(\mathbb T)))$?) }
by $\lambda(a)=\Gamma(a)$ for any $a\in\mathbb Z^k$.

Denote by
$$\mathcal U=\{(1_n-p_i+p_i\tilde{z'})(1_n-q_i+q_i\tilde{z'}^*);\ i=1, ..., k\}\subseteq J_c(U(C\otimes C(\mathbb T))),$$ where $z'$ is the standard generator of $C(\mathbb T)$, and denote by
$$\delta=\min\{\Delta(h)/4;\ h\in\mathcal H_1\}.$$
Applying Lemma \ref{CtimesText}, one obtains a $\mathcal F$-$\epsilon/4$-multiplicative map $\Phi: C\otimes C(\mathbb T)\to B$ such that
\begin{equation*}
[\Phi]=\af,\,\,\,{\rm dist}(\Phi^{\ddag}(x), \lambda(x))<\sigma,\quad x\in \overline{{\cal U}}
\end{equation*}
and
\begin{equation}\label{eq-012-02}
|\tau\circ \Phi(h\otimes 1)-\gamma(\tau)(h\otimes 1)|<\min\{\gamma_1,\delta \},\quad h\in {\mathcal H_1\cup\mathcal H_2}.
\end{equation}

Let $\psi$ denote the restriction of $\Phi$ to $C\otimes 1$. Then one has
$$[\psi]|_{\mathcal P}=[\phi]|_{\mathcal P}.$$
By \eqref{eq-012-02}, one has that for any $h\in\mathcal H_1$,
$$\tau(\psi(h))>\gamma(\tau)(h)-\delta=\tau(\phi'(h\otimes 1))-\delta=\tau(\phi(h))-\delta>\Delta(h)/2,$$
and it is also clear that $$\tau(\phi(h))>\Delta(h)/2\rforal h\in \mathcal H_1.$$
Moreover, for any $h\in\mathcal H_2$, one has
\begin{eqnarray*}
|\tau\circ\psi(h)-\tau\circ\phi(h)|&=& |\tau\circ\Phi(h\otimes 1)-\tau\circ\phi'(h\otimes 1)|\\
&=&  |\tau\circ\Phi(h\otimes 1)-\gamma(\tau)(h\otimes 1)|\\
&\leq & \gamma_1.
\end{eqnarray*}
Therefore, by Corollary \ref{UniCtoA}, there is a unitary $W\in B$ such that
$$||W^*\psi(f)W-\phi(f)||<\epsilon/2\rforal f\in\mathcal F.$$ The the element $$w=W^*\Phi(1\otimes z')W$$ is the desired unitary.
\end{proof}

\begin{thm}\label{BB-exi+}
Let  $C$ be a unital  \CA\, which is a finite direct sum of \CA s in $\mathcal C_0$  and \CA s with
the form $PM_n(C(X))P,$  where $X$ is a finite CW complex, and let
 $G= K_0(C)$. Write
 ${G}=\Z^k\bigoplus \mathrm{Tor}(G)$ with
 $\Z^k$ generated by
$$
\{x_{1}=[p_1]-[q_1], x_2=[p_2]-[q_2], ..., x_{k}=[p_{k}]-[q_{k}]\},
$$
where $p_i, q_i\in M_n(C)$ (for some integer $n\ge 1$) are projections,
$i=1,...,k$.

Let $A$ be a simple C*-algebra in $\mathcal B_0$, and let $B=A\otimes U$ for a UHF algebra $U$ of infinite type. Suppose that $\phi: C\to B$ is a monomorphism. Then, for any finite subsets  $\mathcal F\subseteq C$ and $\mathcal P\subseteq \underline{K}(C)$, any $\ep>0$ and $\sigma>0$, any homomorphism $$\Gamma: \Z^k \to U_0(M_n(B))/CU(M_n(B)),$$ there is a unitary $w\in B$ such that
\begin{enumerate}
\item $\|[\phi(f), w]\|<\ep$, for any $f\in\mathcal F$,
\item $\mathrm{Bott}(\phi, w)|_{\mathcal P}=0$, and
\item $\mathrm{dist}(
\overline{\langle(({\mathbf 1}_n-\phi(p_i))+\phi(p_i)\tilde{w})(({\mathbf 1}_n-\phi(q_i))+\phi(q_i)\tilde{w}^*)\rangle},
 \Gamma(x_i)))<\sigma$, for any $1\leq i\leq k,$
{ where $\tilde{w}=\mathrm{diag}(\overbrace{w, ..., w}^n)$.}
\end{enumerate}
%\beq
%\|[\phi(f), w]\|<\ep\quad\forall f\in\mathcal F,\,\
%\mathrm{Bott}(\phi, w)|_{\mathcal P}=0,\andeqn
%\eneq
%$$\mathrm{dist}(
%\overline{\langle(({\mathbf 1}_n-\phi(p_i))+\phi(p_i)\tilde{w})(({\mathbf 1}_n-\phi(q_i))+\phi(q_i)\tilde{w}^*\rangle},
% \Gamma(x_i)))<\sigma, \quad\forall 1\leq i\leq k,$$
%{ where $\tilde{w}=\mathrm{diag}(\overbrace{w, ..., w}^n)$.}
\end{thm}

\begin{proof}
By \ref{BB-exi}, it suffices to prove the case that $C=PM_n(C(X))P,$   where
$X$ is a finite CW complex, $n\ge 1$ is an integer and $P\in M_n(C(X))$ is a projection.
Proof follows the same lines of the proof as that of \ref{BB-exi} but one  will apply \ref{Next2014/09}
instead of \ref{CtimesText}.
\end{proof}

\section{A pair of almost commuting unitaries}

%{\color{Green} Should we keep the following paragraph or delete it?}

%{\bf One could obtain an existence theorem for Bott element in much more generality.
%Note that, in the following as well as the actual use,  the constants $n$ and $\dt$ appeared  in  \ref{orderK0}, \ref{Vpair1}
%and \ref{Vpair2} need to be  universal.
%Therefore we construct with bare hand first from $C(\T^2).$ }

\begin{lem}\label{orderK0}
Let $C\in {\cal C}.$ There exists a constant $M_C>0$ satisfying the following:
For any $\ep>0,$ any $x\in K_0(C)$ and any $n\ge M_C/\ep,$ if
\beq\label{orderK0-0}
|\rho_C(x)(\tau)| <\ep\tforal \tau\in T(C\otimes M_n),
\eneq
then, there are mutually inequivalent
%{\color{Green} (nonequivalent?)}
and mutually orthogonal minimal projections $p_1, p_2,...,p_{k_1}$ and $q_1,q_2,...,q_{k_2}$ in
$C\otimes M_n$  and positive integers $l_1, l_2,...,l_{k_1}, m_1, m_2,...,m_{k_2}$ such that
\beq\label{orddeK0-1}
x=[\sum_{i=1}^{k_1}l_ip_i]-[\sum_{j=1}^{k_2}m_jq_{{j}}]\tand\\
\tau(\sum_{i=1}^{k_1}l_ip_i)<4\ep\andeqn \tau(\sum_{j=1}^{k_2}m_jq_{{j}})<4\ep
\eneq
for all $\tau\in T(C\otimes M_n).$
%where $l_i$ and $m_j$ are positive integers, $i=1,2,...,k_1$ and
%$j=1,2,..., k_2.$

\end{lem}

\begin{proof}
Let $C=C(F_1, F_2, \phi_1, \phi_2)$ and $F_1=\bigoplus_{i=1}^l M_{r(i)}.$
%let $\{1\}, \{2\}, ..., \{l\}$ be as defined in \ref{K0D1}.
By {{\ref{FG-Ratn}}},
there is an integer $N(C)>0$ such that
every projection in $C\otimes {\cal K}$ is equivalent to  a finite direct sum of   projections {{from a set of finitely many mutually inequivalent minimal projections (some of them may repeat in the direct sum)}} in $M_{N(C)}(C).$ { {By enlarge the number $N(C)$, we may assume this set of mutually inequivalent projections are sitting in $M_{N(C}(C)$, orthogonally.}}
We also assume that, as in \ref{DfC1}, $C$ is minimal.
Let
$$
M=N(C)+2(r(1)\cdot r(2)\cdots r(l))
$$
Suppose that $n\ge M/\ep.$
With the canonical embedding of $K_0(C)$ into $K_0(F_1)\cong\mathbb Z^l$, write
\beq\label{odderK0-2}
x=
\begin{pmatrix}
%\begin{array}{c}
x_1\\
x_2\\
           \vdots \\
           x_l
           \end{pmatrix}\in\mathbb Z^l
         %\end{array}
\eneq
By (\ref{orderK0-0}),
for any irreducible representation $\pi$ of $C$ and any tracial state $t$ on $M_n(\pi(C)),$
\beq\label{orderK0-3}
|t\circ \pi(x)|<\ep.
\eneq
It follows that
\beq\label{orderK0-4}
|x_s|/r(s)n<\ep,\,\,\, s=1,2,...,l.
\eneq
Let
\beq\label{orderK0-5}
T=\max\{|x_s|/r(s): 1\le s\le l\}.
\eneq
Define
\beq\label{orderK0-6}
y=x+T\begin{pmatrix}
                   r(1)\\
                   r(2)\\
           \vdots \\
           r(l)
         \end{pmatrix} \andeqn z=T\begin{pmatrix}
                   r(1)\\
                   r(2)\\
           \vdots \\
           r(l)
         \end{pmatrix}.
         \eneq
It is clear that $z\in K_0(C)_+$ (see \ref{2Lg13}). It follows that $y\in K_0(C).$  One also computes that
$y\in K_0(C)_+.$ It follows that there are projections $p, \, q\in M_L(C)$ for some integer $L\ge 1$ such that $[p]=y$ and $[q]=z.$
Moreover, $x=[p]-[q].$
One also computes that
\beq\label{orderK0-9}
\tau(q)<T/n<\ep\tforal \tau\in T(C\otimes M_n).
\eneq
One also has
\beq\label{orderK0-10}
\tau(p)<2\ep \tforal \tau\in T(C\otimes M_n).
\eneq

There are {{two sets of}} mutually orthogonal minimal projections ${{\{}}p_1, p_2,...,p_{k_1}{{\}}}$ and $
{{\{}}q_1,q_2,...,q_{k_2}{{\}}}$ in $C\otimes M_n$ (since $n>N({{C}})$) such that
\beq\label{orderK0-11}
[p]=\sum_{i=1}^{k_1}l_i[p_i]\andeqn [q]=\sum_{j=1}^{k_2}m_j[q_j].
\eneq
Therefore
\beq\label{orderK0-12}
x=\sum_{i=1}^{k_1}l_i[p_i]-\sum_{j=1}^{k_2}m_j[q_j].
\eneq
%Note that, we may further assume that $p_1, p_2,...,p_{k_1}, q_1, q_2, ...,q_{k_2}$ are mutually inequivalent, by allowing the error
%to be no more that $4\ep$ instead of $2\ep.$

\end{proof}

\begin{lem}\label{Vpair1}
Let $C\in \cal C.$  There is an integer $M_C>0$ satisfying the following: For any $\ep>0$
%and for any unital \CA\,
and for any
$x\in K_0(C)$ with
$$
|\tau(\rho_C(x))|<\ep/24\pi
$$
for all  $\tau\in T(C\otimes M_n),$ where $n\ge 2M_C\pi/\ep,$
there exists a pair of unitaries $u$ and $v\in C\otimes M_n$ such that
\beq\label{Vpair1-1}
\|uv-vu\|<\ep\andeqn \tau({\rm bott}_1(u,\, v))=\tau(x).
\eneq
\end{lem}

\begin{proof}
To simplify the proof, without loss of generality, we may assume that $C$ is minimal.
By applying \ref{orderK0}, there are mutually orthogonal and mutually inequivalent  minimal projections $p_1, p_2,...,p_{k_1}, q_1, q_2,...,q_{k_2}\in C\otimes M_n$ such that
$$
\sum_{i=1}^{k_1}l_i[p_i]-\sum_{j=1}^{k_2}m_j[q_j]=x,
$$
where $l_1,l_2,...,l_{k_1},$ $m_1, m_2,...,m_{k_2}$ are positive integers.
Moreover,
%Since $p_1, p_2,..., p_{k_1},q_1, q_2,...,q_{k_2}$ are mutually in-equivalent and minimal,
\beq\label{Vpair1-2}
\sum_{i=1}^{k_1}l_i\tau(p_i)<\ep/6\pi\andeqn \sum_{j=1}^{k_2}m_j\tau(q_j)<\ep/6\pi
\eneq
for all $\tau\in T(C\otimes M_n).$
%{\bf This requires a lemma}.
Choose $N\le n$ such that $N=[2\pi/\ep]+1.$
By (\ref{Vpair1-2}),
\beq\label{Vpair1-3}
\sum_{i=1}^{k_1}Nl_i\tau(p_i)+\sum_{j=1}^{k_2}Nm_j\tau(q_j)<1/2\tforal \tau\in T(C\otimes M_n).
\eneq
It follows  that there are mutually orthogonal  projections
$d_{i,k}, d_{j,k}'{{\in C\otimes M_n}}, $ $k=1,2,...,N,$ $i=1,2,...,k_1,$ and $j=1,2,...,k_2$ such that
\beq\label{Vpair1-4}
[d_{i,k}]=l_i[p_i]\andeqn [d_{j,k}']=m_j[q_j], \,\,\,i=1,2,...,k_1, \, j=1,2,...,k_2
\eneq
and $k=1,2,...,N.$
Let $D_i=\sum_{k=1}^Nd_{i,k}\andeqn D_j'=\sum_{k=1}^N d_{j,k}',$ $i=1,2,...,k_1$ and $j=1,2,...,k_2.$
There is a partial isometry $s_{i,k},\, s_{j,k}'\in C\otimes M_n$ such that
\beq\label{July28-2}
s_{i,k}^*d_{i,k}s_{i,k}=d_{i,k+1},\,,\,\, (s_{j,k}')^*d_{j,k}'s_{j,k}'=d_{j,k{{+1}}}'\,\,\, k=1,2,...,N-1,\\
s_{i, N}^*d_{i, N}s_{i,N}=d_{i,1}\andeqn (s_{j,N}')^* d_{j,N}'s_{j,N}'=d_{j,1}',
\eneq
 $i=1,2,...,k_1$ and  $j=1,2,...,k_2.$
Thus we obtain unitary $u_i\in D_i(C\otimes M_n)D_i$ and $u_j'=D_j'(C\otimes M_n)D_j'$ such that
\beq\label{Vpair1-5}
u_i^*d_{i,k}u_i=d_{i,k+1},\,\,\,u_i^*d_{i,N}u_i=d_{i,1},\,\,\,
(u_j')^*d_{j,k}'u_j=d_{j,k+1}'\andeqn (u_j')^*d_{j,N}'u_j'=d_{j,1}',
\eneq
$i=1,2,...,k_1,$ $j=1,2,...,k_2.$
Define
$$
v_i=\sum_{k=1}^N e^{\sqrt{-1}(2k\pi/N)}d_{i,k}\andeqn v_j'=\sum_{k=1}^N e^{\sqrt{-1}(2k\pi/N)}d_{j,k}'.
$$
We compute that
\beq\label{Vpair1-6}
\|u_iv_i-v_iu_i\|<\ep\andeqn \|u_j'v_j'-v_j'u_j'\|<\ep,\\
{1\over{2\pi \sqrt{-1}}}\tau(\log v_iu_iv_i^*u_i^*)=l_i\tau(p_i)\andeqn\\
{1\over{2\pi \sqrt{-1}}}\tau(\log v_j'u_j'(v_j')^*(u_j')^*)=m_j\tau(q_j),
\eneq
for $\tau\in T(C\otimes M_n),$ $i=1,2,...,k_1$ and $j=1,2,...,k_2.$
Now define
\beq\label{Vpair1-7}
u=\sum_{i=1}^{k_1}u_i+\sum_{j=1}^{k_2}u_j'+(1_{C\otimes M_n}-\sum_{i=1}^{k_2}D_i-\sum_{j=1}^{k_2}D_j')\andeqn\\
v=\sum_{i=1}^{k_1}v_i+\sum_{j=1}^{k_2}(v_j')^*+(1_{C\otimes M_n}-\sum_{i=1}^{k_2}D_i-\sum_{j=1}^{k_2}D_j').
\eneq
We then compute that
\beq\label{Vpair1-8}
\hspace{-0.4in}\tau({\rm bott}_1(u, v))=\sum_{i=1}^{k_1}{1\over{2\pi \sqrt{-1}}}\tau(\log (v_iu_iv_i^*u_i^*))-\sum_{j=1}^{k_2}{1\over{2\pi \sqrt{-1}}}\tau(\log v_j'u_j'(v_j')^*(u_j')^*)\\
=\sum_{i=1}^{k_1} l_i\tau(p_i) -\sum_{j=1}^{k_2}m_j\tau(q_j)=\tau(x)
\eneq
for all $\tau\in T(C\otimes M_n).$
\end{proof}

\begin{lem}\label{Vpair2}
Let $\ep>0.$ There exists $\sigma>0$ satisfying the following:
Let $A=A_1\otimes U,$ where $U$ is a UHF-algebra of infinite type  and $A_1\in {\cal B}_0$, let $u\in U(A)$ be a unitary with $sp(u)=\T,$ and let
$x\in K_0(A)$ with
 $|\tau(\rho_A(x))|<\sigma$ for all $\tau\in T(A)$ and $y\in K_1(A) .$ Then there exists a unitary
$v\in U(A)$ such that
\beq\label{Vpair2-1}
\|uv-vu\|<\ep,\,\,\,{\rm bott}_1(u,v)=x\andeqn [v]=y.
\eneq
\end{lem}

\begin{proof}
Let $\phi_0: C(\T)\to A$ be the unital monomorphism
defined by $\phi_0(f)=f(u)$ for all $f\in C(\T).$
Let $\Delta_0: A_+^{q,{\bf 1}}\setminus \{0\}\to (0,1)$ be defined by
$\Delta_0(\hat{f})=\inf \{\tau(f): \tau\in T(A)\}.$
Let $\ep>0$ be given.
Choose $0<\ep_1<\ep$ such that
$$
{\rm bott}_1(z_1,z_2)={\rm bott}_1(z_1', z_2')
$$
if $\|z_1-z_1'\|<\ep_1$ and $\|z_2-z_2'\|<\ep_1$ for any two pairs of unitaries
$z_1, z_2, $ and $z_1', z_2'$ which also have the property that $\|z_1z_2-z_2z_1\|<\ep_1$ and
$\|z_1'z_2'-z_2'z_1'\|<\ep_1.$

Let ${\cal H}_1\subset C(\T)_+^{\bf 1}\setminus \{0\}$
be a finite subset, $\gamma_1>0,$ $\gamma_{{2}}>0$  and ${\cal H}_2\subset C(\T)_{s.a.}$
be a finite subset  as required by \ref{Unitaryuni} (for
$\ep_1/4$ and $\Delta_0/2$).

Let
$$\dt_1=\min\{\gamma_1/16, \gamma_2/16, \min\{\Delta(\hat{f}): f\in {\cal H}_1\}/4\}.$$
 Let
$\dt=\min\{\dt_1/16, (\dt_1/16)(\ep_1/32\pi)\}.$

Let $e\in 1\otimes U \subseteq A $ be a non-zero projection such that
$\tau(e)<\dt_1$ for all $\tau\in T(A).$
Let $B=eAe$ (then $B\cong A\otimes U'$ for some UHF-algebra $U'$).
It follows from \ref{Cistbk} that there is a unital simple \CA\, $C'=\lim_{n\to\infty}(C_n, \psi_n),$ where $C_n\in {\cal C}_0$ and $C=C'\otimes U$ such that
$$
(K_0(C), K_0(C)_+, [1_C], T(C), r_C)=(\rho_A(K_0(A)), (\rho_A(K_0(A)))_+,
[e], T(A), r_A).
$$
Moreover, we may assume that all $\psi_n$ are unital.

Now suppose that $x\in K_0(A)$ with $|\tau(\rho_A(x))|<\dt$
for all $\tau\in T(A)$ and suppose that $y\in K_1(A).$   Let $z=\rho_A(x){{\in K_0(C)}}.$ We identify $z$ with the element in $K_0(C)$ in the above identification.
We claim that,
there is $n_0\ge 1$ such that there is $x'\in K_0(C_{n_0}\otimes U)$ such that ${{z}}= (\psi_{n_0, \infty})_{*0}(x'){{\in K_0(C)}}$
and $|t(\rho_{C_{n_0}\otimes U})(x')|<\dt$ for all $t\in T(C_{n_0}\otimes U).$

Otherwise, there is an increasing sequence $n_k,$ $x_k\in K_0(C_{n_k}\otimes U)$ such that
\beq\label{Vpair2-2}
(\psi_{n_k,\infty})_{*0}(x_k)={{z\in K_0(C)}}\andeqn |t_k(\rho_{C_{n_k}\otimes U})(x_k)|\ge \dt
\eneq
for some $t_k\in T(C_{n_k}\otimes U),$ $k=1,2,....$
Let $L_k: C\to C_{n_k}\otimes U$ such that
$$
\lim_{n\to\infty}\|\psi_{n,\infty}\circ L_n(c)-c\|=0
$$
for all $c\in \psi_{k,\infty}( C_{n_k}\otimes U),$ $k=1,2,....$ It follows that
the limit points of $t_k\circ L_k$ is a tracial state of $C.$  Let $t_0$ be one of such limit.
Then, by (\ref{Vpair2-2}),
$$
t_0(\rho_C(z))\ge \dt.
$$
This proves the claim.

Write $U=\lim_{n\to\infty}(M_{r(m)}, \imath_m),$ where $\imath_m: M_{r(m)}\to M_{r(m+1)}$ is a unital
embedding. By repeating the above argument, we obtain $m_0\ge 1$  and $y'\in K_0(C_{n_0}\otimes M_{r(m_0)})=
K_0(C_{n_0})$  such that
$(\imath_{m_0, \infty})_{*0}(y')=x'$ and $|t(\rho_{C_{n_0}}(y'))|<\dt$ for all $t\in T(C_{n_0}\otimes M_{r(m_0)}).$
Let $M_{C_{N_0}}$ be as the constant in \ref{Vpair1}. Choose $r(m_1)\ge \max\{{{48}}M_{C_{n_0}}/\dt, r(m_0)\}$
and let $y''=(\imath_{m_0, m_1})_{*0}(y').$ Then, we compute that
$$
|t(\rho_{C_{n_0}}(y''))|<\dt\rforal t\in T(C_{n_0}\otimes M_{r(m_1)}).
$$

It follows from \ref{Vpair1} that there exists a pair of unitaries
$u_1', v_1'\in C_{n_0}\otimes M_{r(m_1)}$ such that
\beq\label{Vpair2-3-}
\|u_1'-v_1'\|<\ep_1/4\andeqn {\rm bott}_1(u_1', v_1')=y''.
\eneq
Put $u_1=\imath_{m_1, \infty}(u_1')$ and $u_2=\imath_{m_1, \infty}(v_1').$
Then (\ref{Vpair2-3-}) implies that
\beq\label{Vpair2-3}
\|u_1-v_1\|<\ep_1/4\andeqn {\rm bott}_1(u_1, v_1)=x'.
\eneq

Let $h_0: C_{n_0}\otimes U\to eAe$ be a unital \hm\, given by \ref{Cistbk}
such that
\beq\label{Vpair2-4}
\rho_A\circ (h_0)_{*0}=(\psi_{n_0, \infty})_{*0}.
\eneq
It follows that
\beq\label{Vpair2-5}
\rho_A((h_0)_{*0}(x')-x)=0.
\eneq
Let $u_2=h_0(u_1)$ and $v_2=h_0(v_1).$
We have that
\beq\label{Vpair2-6}
\rho_A({\rm bott}_1(u_2, v_2)-x)=0.
\eneq
Choose another non-zero projection $e_1\in A$ such that $e_1e=ee_1=0$ and $\tau(e_1)<\dt_1/16$
for all $\tau\in T(A).$
It follows from \ref{preBot1} that there is a  sequence of unital \morp s $L_n: C(\T^2)\to e_1Ae_1$ such that
\beq\label{Vpair2-7}
\lim_{n\to\infty}\|L_n(fg)-L_n(f)L_n(g)\|=0\andeqn [L_n](b)=x-{\rm bott}_1(u_2, v_2){{,}}
\eneq
{ {where $b$ is the Bott element in $K_0(C(T^2))$.}} (In fact, we can also apply \ref{Next2014/09} here.)
Thus we obtain a pair of unitaries $u_3, v_3\in e_1Ae_1$ such that
\beq\label{Vpair2-8}
\|u_3v_3-v_3u_3\|<\ep_1/4\andeqn {\rm bott}_1(u_3, v_3)=x-{\rm bott}_1(u_2, v_2).
\eneq
Let $e_2, e_3\in (1-e-e_1)A(1-e-e_1)$ be a pair of non-zero mutually orthogonal projections such that $\tau(e_2)<\dt_1/32$
and $\tau(e_3)<\dt_1/32$
for all $\tau\in T(A).$
It follows from \ref{istTr}{{, applied to $X=\T$,}}  that there is a unitary $u_4\in (1-e-e_1-e_2-e_3)A(1-e-e_1-e_2-e_3)$ such that
\beq\label{Vpair2-9}
|\tau\circ (f(u_4))-\tau\circ f(u)\|<\dt_1/4\tforal f\in {\cal H}_2\cup {\cal H}_1\andeqn \tforal   \tau\in T(A).
\eneq
Let $w=u_2+u_3+u_4+(1-e-e_1-e_2-e_3).$
It follows from Theorem 3.10 of \cite{GLX-ER}  that there exists $u_5\in U(e_2Ae_2)$ such that
\beq\label{Vpair2-10}
\overline{u_5}=\bar{u}{\bar w^*}\in U(A)/CU(A).
\eneq
%Suppose that $u_5\in e_2Ae_2$ such that ${\bar u_5}=w_1.$
Since $A$ is simple and has stable rank one, there exists
a unitary
%unitray 
$v_4\in e_3Ae_3$ such that$[v_4]=y-[v_2+v_3+({{e_2+}}e_3)]\in K_1(A).$
Now define
$$
u_6=u_2+u_3+u_4+u_5+e_3\andeqn v_6=v_2+v_3+(1-e-e_1-{{e_2-}}e_3)+{{e_2+}}v_4.
$$
Then
\beq\label{Vpair2-11}
\|u_6v_6-v_6u_6\|<\ep_1/2,\,\,\,{\rm bott}_1(u_6, v_6)=x\andeqn
[v_6]=y.
\eneq
Moreover,
\beq\label{Vpair2-12}
\tau(f(u_6))\ge \Delta(\hat{f})/2\tforal f\in {\cal H}_1\\
|\tau(f(u)-\tau(f(u_6))|<\gamma_1\andeqn
\bar{u_6}=\bar{u}.
\eneq
It follows from \ref{Unitaryuni} that there exists a unitary $W\in A$ such that
\beq\label{Vpair2-13}
\|W^*u_6W-u\|<\ep_1/2.
\eneq
Now let $v=W^* v_6W.$
We compute that
\beq\label{Vpair2-14}
\|uv-vu\|<\ep,\,\,\,{\rm bott}_1(u, v)={\rm bott}_1(u_6, v_6)=x
\andeqn [v]=y.
\eneq

\end{proof}

\section{More existence theorems for Bott elements}

Using \ref{Vpair2},
%\ref{Prebot1} {\bf
\ref{preBot1}, \ref{Cist}, \ref{ExtTrace} and \ref{MUN1} we can show the following:

\begin{lem}\label{Extbot1}
Let  $A=A_1\otimes U_1,$ where $A_1$ is as in {\rm \ref{RangT}} and $B=B_1\otimes U_2,$ where $B_1\in {\cal B}_0$ %{\color{Green} (it seems to me should be $\mathcal B_0$, since Lemma \ref{Vpair2} uses $\mathcal B_0$)},
and $U_1, U_2$ are two UHF-algebras of infinite type.
For any $\ep>0,$ any finite subset ${\cal F}\subset A,$ any finite subset ${\cal P}\subset \underline{K}(A),$ there exist $\dt>0$ and a finite subset ${\cal Q}\subset K_1(A)$ satisfying the following:
Let $\phi: A\to B$ be a unital \hm\, and $\af\in KL(A\otimes C(\T), B)$ such that
\beq\label{Extbot1-1}
|\tau\circ \rho_B(\af(\boldsymbol{\bt}(x)))|<\dt \tforal x\in {\cal Q}\tand \tforal \tau\in T(B){{,}}
\eneq
there exists a unitary $u\in B$ such that
\beq\label{Extbot1-2}
\|[\phi(x), \, u]\|&<&\ep \tforal x\in {\cal F}\andeqn\\
{\rm Bott}(\phi, u)|_{\cal P}&=&\af(\boldsymbol{\bt})|_{\cal P}{{.}}
\eneq
\end{lem}

\begin{proof}
Let $\ep_1>0$ and let ${\cal F}_1\subset A$ be a finite subset satisfying the following:
If $$L, L': A\otimes C(\T)\to B$$ are
two {{unital}} ${\cal F}_1'$-$\ep_1$-multiplicative \morp s
such that
\beq\label{Extbot1-4}
\|L(f)-L'(f)\|<\ep_1\tforal f\in {\cal F}_1',
\eneq
where
$$
{\cal F}'_1=\{a\otimes g: a\in {\cal F}_1\andeqn g\in \{z, z^*,1_{C(\T)}\}\},
$$
then
\begin{equation}\label{Extbot1-3}
[L]|_{\boldsymbol{\bt}({\cal P})}=[L']|_{\boldsymbol{\bt}({\cal P})}.
\end{equation}

%{\color{Green} (I forgot if we decide to use wedge product or disjoint union)}
Let $B_{1,n}=M_{m(1,n)}(C(\T))\oplus M_{m(2,n)}(C(\T))\oplus\cdots \oplus M_{m(k_1(1),n)}(C(\T)),$
$B_{2,n}=PM_{r_1(n)}(C(X_n))P,$ where
$X_n$ is a finite disjoint union of
$S^2, T_{0,k}$ and $T_{1,k}$ (for various $k\ge 1$). Let $B_{3,n}$
be a finite direct sum of \CA s in ${\cal C}_0$ (with trivial $K_1$ and
${\rm ker}\rho_{B_{3,n}}=\{0\}$), $n=1,2,....$
Put $C_n=B_{1,n}\oplus B_{2,n}\oplus B_{3,n},$ $n=1,2,....$
We may write that $A=\lim_{n\to\infty}(C_n,\imath_n)$ as in \ref{RangT}. So we also assume that
$\imath_n$ are injective,
\beq\label{Extbot1-5}
&&{\rm ker}\rho_A\subset (\imath_{n, \infty})_{*0}({\rm ker}\rho_{C_n})\andeqn\\
&&\lim_{n\to\infty}\sup\{\tau(1_{B_{1,n}}\oplus 1_{B_{2,n}}): \tau\in T(B)\}=0
\eneq

Let $\ep_2=\min\{\ep_1/4, \ep/4\}$ and
let ${\cal F}_2={\cal F}_1\cup {\cal F}.$
% We will apply \ref{MUN1}
%for $A$ (with $X$ being a single point). Note in this case $\Delta$ and ${\cal H}_1$  will not be needed. Let $\dt_1>0$ (in place of $\dt$), ${\cal G}_1\subset A$ be a finite subset (in place of ${\cal G},$ $\sigma_1, \sigma_2>0$ ${\cal P}_1\subset \underline{K}(A)$ (in place of ${\cal P}$) be a finite subset,
%${\cal U}\subset U(A)$ be a finite subset and ${\cal H}_2\subset A_{a.s.}$ be as required by \ref{MUN1} for $\ep_2$ and ${\cal F}_2.$

Let ${\cal P}_{1,1}\subset \underline{K}(B_{1,n_1}),$
${\cal P}_{2,1}\subset \underline{K}(B_{2,n_1})$ and ${\cal P}_{3,1}\subset \underline{K}(B_{3,n_1})$ be finite subsets
such that
$$
{\cal P}\subset [\imath_{n_1, \infty}]({\cal P}_{1,1})\cup [\imath_{n_1, \infty}]({\cal P}_{2,1})\cup
[\imath_{n_1, \infty}])({\cal P}_{3,1})
$$
for some $n_1\ge 1.$
Let  ${\cal Q}'$ be a finite set of generators of $K_1(C_n)$ and
let ${\cal Q}=[\imath_{n_1,\infty}]({\cal Q}').$

Without loss of generality, we may assume that ${\cal F}_1\cup {\cal F}\subset
\imath_{n_1,\infty}(C_{n_1}).$ Let
${\cal F}_{1,1}\subset B_{1, n_1},$ ${\cal F}_{2,1}\subset B_{2, n_2}$ and
${\cal F}_{3,1}\subset B_{3, n_1}$ be finite subsets such that
%We also assume that (possibly with larger $n_1$)
\beq\label{Extbot1-6}
{\cal F}_1\cup {\cal F}\subset \imath_{n_1, \infty}({\cal F}_{1,1}\cup{\cal F}_{2,1}\cup {\cal F}_{3,1}).
%\tau(1_{B_{1, n_1}}+1_{B_{2, n_2}})<\min\{\sigma_1/4, \sigma_2\}\tforal \tau\in T(A).
\eneq
Let $e_1=\imath_{n_1, \infty}(1_{B_{1, n_1}}),$ $e_2=\imath_{n_1, \infty}(1_{B_{2, n_1}})$ and $e_3=1-e_1-e_2.$
Let $\Delta_1: (B_{2,n_1})_+^{q, {\bf 1}}\setminus \{0\}\to (0,1)$ be defined by
$$
\Delta_1(\hat{h})=(1/2)\inf\{\tau(\phi(h)): \tau\in T(A)\}\tforal h\in (B_{2, n_1})_+^{{\bf 1}}\setminus \{0\}.
$$
Let $\Delta_2: B_{{{3}}, n_1}^{q, {\bf 1}}\setminus \{0\}\to (0,1)$ be defined by
$$
\Delta_2(\hat{h})=(1/2)\inf\{\tau(\phi(h)): \tau\in T(A)\}\tforal h\in (B_{3, n_1})_+^{{\bf 1}}\setminus \{0\}.
$$
Note that $B_{2, n_1}$ has the form $C$ in \ref{UniCtoA}. So we will apply \ref{UniCtoA}.
Let ${\cal H}_{2,1}\subset (B_{2, n_1}^{\bf 1})_+\setminus \{0\}$ (in place of ${\cal H}_1$) be a finite subset, $\gamma_{2,1}>0$
(in place of $\gamma_1$),
%$\gamma_{2,2}>0$ (in place of $\gamma_2$),
$\dt_{2,1}>0$ (in place of $\dt$),
${\cal G}_{2,1}\subset B_{2, n_1}$ (in place of ${\cal G}$) be a finite subset, ${\cal P}_{2,2}\subset \underline{K}(B_{2, n_1})$
(in place of ${\cal P}$), ${\cal H}_{2,2}\subset (B_{2, n_1})_{s.a.}$ (in place of ${\cal H}_2$) be a finite subset
%${\cal U}_{2,1}\subset J_c(K_1(C))$ (in place of ${\cal U}$) be a finite subset for which
%$[{\cal U}]\subset {\cal P}_{2,2}$
 {{as}} required by \ref{UniCtoA} for  $\ep_2/16,$ ${\cal F}_{2,1}$ and
$\Delta_1$ (see also the Remark  \ref{RemUniCtoA} since
$K_1(B_{2,n_1})$ is torsion or zero).

Now let $\sigma>0$ be required by  \ref{Vpair2} for $\ep_2/4$ (in place of $\ep$).
Let $\dt=\sigma\cdot \inf\{\tau(e_1): \tau\in T(A)\}.$
It follows from \ref{Vpair2} that { {if $|\tau\circ\rho_B(\af(\bt(x)))|<\dt$ for all $x\in [\imath_{n_1, \infty}])({\cal P}_{1,1})$ then}} there is a unitary $v_1\in e_1Be_1$ such that
\beq\label{Extbot1-7}
{\rm Bott}(\phi\circ \imath_{n_1, \infty}, v_1)|_{{\cal P}_{1,1}}={{\af\circ{\bt}\circ[\imath_{n_1, \infty}]|_{({\cal P}_{1,1})}}}.
\eneq
 Note that $K_1(B_{2, n_1})$
 is a finite group.
Therefore
\beq\label{Extbot1-8}
\af(\boldsymbol{\bt}([\imath_{n_1, \infty}])(K_1(B_{2, n_1}))\subset {\rm ker}\rho_B
%\andeqn
%\\\label{Extbot1-8+}
%\af(\boldsymbol{\bt}([\imath_{n_1, \infty}]))({\rm ker}\rho_{B_{2, n_1}})\subset {\rm ker}\rho_B.
\eneq
Define $\kappa_1\in KK(B_{2, n_1}\otimes C(\T){{, A}})$ by
$\kappa_1|_{\underline{K}(B_{2, n_1})}=[\phi\circ \imath_{n_1, \infty}|_{B_{2, n_1}}]$ and
$\kappa_1|_{\boldsymbol{\bt}(\underline{K}(B_{2, n_1})})=\af|_{ \boldsymbol{\bt}(\underline{K}(B_{2, n_1})}.$
Since $\imath_{n_1, \infty}$ is injective, by (\ref{Extbot1-8}),
%and (\ref{Extbot1-8+}),
$\kappa_1\in KK_e(B_{2, n_1}\otimes C(\T), e_2Be_2)^{++}.$

Let
$$
\sigma_0=\min\{\gamma_{ 2,1}/2,  \min\{\Delta_1(\hat{h}): h\in {\cal H}_{2,1}\}\cdot \inf\{\tau(e_2): \tau\in T(A)\}.
$$
Define $\gamma_0: T(e_2Ae_2)\to T_f(B_{2,n_1}{{\otimes C(\T)}})$ by
$\gamma_0(\tau)(f\otimes 1_{C(\T)})=\tau\circ \phi\circ \imath_{n_1, \infty}(f)$ for all $f\in B_{2,n_1}$ and
$\gamma_0(1\otimes g)=\int_{\T}g(t)dt$ for all $g\in C(\T).$
It follows from \ref{Next2014/09}, {{applied to the space $X_{n_1}\times \T$,}} that there
is a unital monomorphism $\Phi: B_{2,n_1}\otimes C(\T)\to e_2Ae_2$
such that $[\Phi]=\kappa_1$ and $\Phi_T=\gamma_0.$  Put $L_2=\Phi|_{B_{2,n_1}}${{(by identifying $B_{2,n_1}$ with ${ B_{2,n_1}\otimes 1_{C(\T)}}$)}}  and $v_2'=\Phi(1\otimes z),$
where $z\in C(\T)$ is the identity function on the unit circle.
Then $L_2$ is a unital monomorphism  from $ B_{2, n_1}$ to $ e_2Ae_2.$ We also have the following:
\beq\label{Extbot-9}
&&[L_2]=[\phi\circ \imath_{n_1, \infty}],\,\,\,\|[L_2(f),\, v_2']\|=0\\
&&{\rm Bott}(L_2,\, v_2')|_{{\cal P}_{2,2}}=\af(\boldsymbol{\bt}( [\imath_{n_1, \infty}]))|_{{\cal P}_{2,2}}\andeqn\\\label{Extbot-9++}
&&|\tau\circ L_2(f)-\tau\circ \phi\circ \imath_{n_1, \infty}(f)|=0 \tforal f\in {\cal H}_{2,1}\cup {\cal H}_{2,2}
\eneq
and for all $\tau\in T(e_2Ae_2).$
It follows from (\ref{Extbot-9++}) that
\beq\label{Extbot-10}
\tau(L_2(f))\ge \Delta_1(\hat{f})\cdot \tau(e_2).
\tforal f\in {\cal H}_{2,1}\andeqn \tau\in  T(A).
\eneq
By \ref{UniCtoA}(see also \ref{RemUniCtoA}),  there exists a unitary $w\in e_2Ae_2$ such that
\beq\label{Extbot-11}
\|{\rm Ad}\, w\circ L_2(f)-\phi\circ \imath_{n_1, \infty} (f)\|<\ep_2/16\tforal  f\in {\cal F}_{2,1}.
\eneq
Define $v_2=w^*v_2'w.$
Then
\beq\label{Extbot-12}
\|[\phi\circ \imath_{n_1,\infty}(f),\,v_2]\|<\ep_2/8\andeqn
{\rm Bott}(\phi\circ \imath_{n_1, \infty},\, v_2)|_{{\cal P}_{2,1}}=\af(\boldsymbol{\bt}([\imath_{n_2, \infty}])|_{{\cal P}_{2,1}}
\eneq

Note that $B_{3, n_1}$ has the form $C$ in \ref{UniCtoA}.
  Let ${\cal H}_{3,1}\subset (B_{3, n_1})_+^{\bf 1}\setminus \{0\}$
(in place of ${\cal H}_1$) be a finite subset, $\gamma_{3,1}>0$ (in place of $\gamma_1$),
%$\gamma_{3,2}>0$ (in place of $\gamma_2$),
$\dt_{3,1}>0$ (in place of $\dt$), ${\cal G}_{3,1}\subset B_{3, N_1})$ (in place of ${\cal G}$) be a finite subset,  ${\cal P}_{3,2}\subset \underline{K}(B_{3,n_1})$ (in place of ${\cal P}$) be a finite subset and ${\cal H}_{3,2}\subset
(B_{3, n_1})_{s.a.}$ (note that $K_1(B_{3, n_1})=\{0\}$) be  {{as}} required by \ref{UniCtoA} for
$\ep_2/16,$ ${\cal F}_{3,1}$ and $\Delta_2$ (see also Remark \ref{RemUniCtoA}).

Let
$$
\sigma_1=(\gamma_{3,1}/2)\min\{\tau(e_3): \tau\in T(A)\}\cdot \min\{ \Delta_1(\hat{f}): f\in {\cal H}_{3,1}\}.
$$

Note that ${\rm ker}\rho_{B_{3, n_1}}=\{0\}$ and $K_1(B_{3, n_1})=\{0\}.$
Therefore ${\rm ker}\rho_{B_{3, n_1}\otimes C(\T)}={\rm ker}\rho_{B_{3, n_1}}=\{0\}.$
Define
$\kappa_2\in KK(B_{3, n_1}\otimes C(\T))$ as follows
$$
\kappa_2|_{\underline{K}(B_{3, n_1})}=[\phi\circ \imath_{n_1, \infty}]|_{B_{3, n_1}}\andeqn
\kappa_2|_{\boldsymbol{\bt}(\underline{K}(B_{3, n_1})}=\af(\boldsymbol{\bt}(\imath_{n_1, \infty})|_{\underline{K}(B_{3, n_1})}.
$$
Thus $\kappa_2\in KK_e(B_{3,n_1}\otimes C(\T), e_3Ae_3)^{++}.$
It follows from \ref{CtimesText00} that there is  a unital ${\cal G}_{3,1}$-$\min\{\ep_2/16,\dt_{3,1}/2\}$-multiplicative
\morp\, $L_3: B_{3, n_1}\to e_3Ae_3$ and a unitary $v_3'\in e_3Ae_3$ such that
\beq\label{Extbot-13}
&&\hspace{-1.4in}[L_3]=[\phi],\,\,\, \|[L_3(f),\, v_3']\|<\ep_2/16\tforal f\in {\cal G}_{3,1},\\
&&\hspace{-1.4in}{\rm Bott}(L_3,\, v_3')|_{{\cal P}_{3,1}}=\kappa_2|_{\boldsymbol{\bt}({{\cal P}_{3,2}})}\andeqn\\\label{Extbot-13++}
&&\hspace{-1.4in}|\tau\circ L_3(f)-\tau\circ \phi\circ \imath_{n_1, \infty}(f)|<\sigma_1\tforal f\in {\cal H}_{3,1}\cup {\cal H}_{3,2}\\\andeqn
\tforal \tau\in T(e_3AE_3).
\eneq
It follows that  (\ref{Extbot-13++}) that
\beq\label{Extbot-14}
\tau(L_3(f))\ge \Delta_1(\hat{f})\tau(e_3)\tforal f\in {\cal H}_{3,1}\andeqn \tforal \tau\in T(A).
\eneq
It follows from \ref{UniCtoA} and its remark that there exists a unitary $w_1\in e_3Ae_3$ such that
\beq\label{Extbot-15}
\|{\rm Ad}\, w_1\circ L_2(f)-\phi\circ \imath_{n_1, \infty}(f)\|<\ep_2/16\tforal f\in {\cal F}_{3,1}.
\eneq
Define $v_3=w_1^*v_3'w_1.$ Then
\beq\label{Extbot-16}
\|[\phi\circ \imath_{n_1, \infty}(f), \, v_3]\|<\ep_2/8\andeqn {\rm Bott}(\phi\circ \imath_{n_1, \infty},\, v_3)|_{{\cal P}_{3,1}}
={\rm Bott}(L_3,\, v_3')|_{{\cal P}_{3,1}}.
\eneq
Let $v=v_1+v_2+v_3.$ Then
\beq\label{Extbot-17}
\|[\phi(f), \, v]\|<\ep\tforal f\in {\cal F}.
\eneq
Moreover, we compute that
\beq\label{Extbot-18}
{\rm Bott}(\phi, \, v)|_{\cal P}=\af|_{\boldsymbol{\bt}({\cal P})}.
\eneq
\end{proof}

We actually prove the following:
\begin{lem}\label{Extbot2}
Let  $A=A_1\otimes U_1,$ where $A_1$ {{is}}  as in {\rm  \ref{RangT}} and $B=B_1\otimes U_2,$ where $B_1\in {\cal B}_{{0}}$ {{is}}  unital simple \CA\, and where $U_1, U_2$ are two UHF-algebras of infinite type. Write
$A=\lim_{n\to\infty} (C_n, \imath_n)$ as described in {\rm \ref{RangT}}.
For any $\ep>0,$ any finite subset ${\cal F}\subset A,$ any finite subset ${\cal P}\subset \underline{K}(A),$ there exists an integer $n\ge 1$
such that ${\cal P}\subset [\imath_{n,\infty}](\underline{K}(C_n)) $
% unital \SCA\, $A_0\subset A_1$ with $1_{A_0}=1_A$ and with finitely generated $K_i(A_0)$
%($i=0,1$)
and  there is a finite subset
${\cal Q}\subset K_1(C_n)$ which generates $K_1(C_n)$
%with $[\imath]({\cal Q})\supset{\cal P},$
%where $\imath: A_0\to A$ is the embedding,
and there exists $\dt>0$ satisfying the following:
Let $\phi: A\to B$ be a unital \hm\, and let $\af\in KK(C_n\otimes C(\T), B)$ such that
\begin{equation*}
|\tau\circ \rho_B([\imath_{n, \infty}]\circ \af(\boldsymbol{\bt}(x)))|<\dt \tforal x\in {\cal Q}\tand \tforal \tau\in T(B),
\end{equation*}
%where ${\cal Q}_1$ is a finite set of generators of $K_1(A_n)$
%and
%$\imath: A_0\to A$ is the embedding,
there exists a unitary $u\in B$ such that
\begin{equation*}
\|[\phi(x), \, u]\|<\ep \tforal x\in {\cal F}
\andeqn
{\rm Bott}(\phi\circ [\imath_{n,\infty}], u)=\af(\boldsymbol{\bt})|_{\cal P}.
\end{equation*}
\end{lem}

%{\color{Green}{
\begin{rem}\label{RemB2} Note that, in the above statement, if an integer $n$ works, any integer $m\ge n$ also
works.
In the terminology of Definition 3.6 of \cite{L-N}, the above also implies that
$B$ has property (B1) and (B2) associated with $C.$
\end{rem}
%}}

\begin{cor}\label{L86}
Let $B\in {\cal B}_{{0}}$ which satisfies the UCT, $A_1\in {\cal B}_{{0}},$ let $C=B\otimes U_1$ and $A=A_1\otimes U_2,$ where $U_1$ and $U_2$ are unital infinite dimensional UHF-algebras. Suppose that $\kappa\in KK_e(C,A)^{++},$
$\gamma: T(A)\to T(C)$ is a continuous affine map and $\af: U(C)/CU(C)\to
U(A)/CU(A)$ is a continuous \hm\, for which $\gamma,\, \af$ and $\kappa$ are compatible. Then there exists a unital monomorphism $h: C\to A$ such that
\begin{enumerate}
\item $[h]=\kappa$ in $KK_e(C,A)^{++}$,
\item $h_T=\gamma$ and $h^{\ddag}=\af.$
\end{enumerate}
\end{cor}

\begin{proof}
The proof follows the same line as that of Theorem 8.6 of \cite{Lnclasn}. Denote by $\overline{\kappa}\in KL(C, A)$ be the image of $\kappa$. It follows from Lemma \ref{L85} that there is a unital monomorphism $\phi: C\to A$ such that
$$[\phi]=\overline{\kappa}, \quad \phi^\ddag=\alpha,\quad\textrm{and}\quad (\phi)_T=\gamma.$$
Note that it follows from the UCT that (as an element  in $KK(C,A)$)
$$\kappa-[\phi]\in\textrm{Pext}(K_*(C), K_{*+1}(A)).$$
By Lemma \ref{Extbot2}, the C*-algebra $A$ has Property (B1) and Property (B2) associated with C in the sense of \cite{L-N}. Since $A$ contains a copy of $U_2,$  it is infinite dimensional, simple and antiliminal. It follows
from a result in \cite{AS} that $A$ contains an element $b$ with $sp(b)=[0,1].$ Moreover,
$A$ is approximately divisible. By Theorem 3.17 of \cite{L-N}, there is a unital monomorphism $\psi_0: A\to A$ which is approximate inner such that
$$[\psi_0\circ\phi]-[\phi]=\kappa-[\phi]\quad\textrm{in $KK(C, A)$}.$$ Then the map
$$h:=\psi_0\circ\phi$$
satisfy the corollary.
\end{proof}

\begin{lem}\label{Extbot3}
Let  $A=A_1\otimes U_1,$ where $A_1$ {{is}} as in \ref{RangT} and $B=B_1\otimes U_2,$ where $B_1\in {\cal B}_{{0}}$ {{is}} unital simple \CA\, and where $U_1, U_2$ are two UHF-algebras of infinite type. Let
$A=\lim_{n\to\infty} (C_n, \imath_n)$ be as described in \ref{RangT}.
For any $\ep>0,$ any $\sigma>0,$ any finite subset ${\cal F}\subset A,$ any finite subset ${\cal P}\subset \underline{K}(A),$
and any projections $p_1, p_2,...,p_k,q_1,q_2,...,q_k\in A$
such that $\{x_1, x_2,...,x_k\}$
generates a free subgroup  $G$ of $K_0(A),$  where $x_i=[p_i]-[q_i],$ $i=1,2,...,k,$
there exists an integer $n\ge 1$
such that ${\cal P}\subset [\imath_{n,\infty}](\underline{K}(C_n)) $
% unital \SCA\, $A_0\subset A_1$ with $1_{A_0}=1_A$ and with finitely generated $K_i(A_0)$
%($i=0,1$)
and  there is a finite subset
${\cal Q}\subset K_1(C_n)$ which generates $K_1(C_n)$
%with $[\imath]({\cal Q})\supset{\cal P},$
%where $\imath: A_0\to A$ is the embedding,
and there exists $\dt>0$ satisfying the following:
Let $\phi: A\to B$ be a unital \hm,  let
$\Gamma: G\to U_0(B)/CU(B)$ be  a \hm\, and  let $\af\in KK(C_n\otimes C(\T), B)$
 such that
\begin{equation*}
|\tau\circ \rho_B([\imath_{n, \infty}]\circ \af(\boldsymbol{\bt}(x)))|<\dt \tforal x\in {\cal Q}\tand \tforal \tau\in T(B),
\end{equation*}
%where ${\cal Q}_1$ is a finite set of generators of $K_1(A_n)$
%and
%$\imath: A_0\to A$ is the embedding,
there exists a unitary $u\in B$ such that
\begin{equation*}
\|[\phi(x), \, u]\| < \ep \tforal x\in {\cal F},\,\,\,
%\end{equation*}
%\begin{equation*}
{\rm Bott}(\phi\circ [\imath_{n,\infty}], u) = \af(\boldsymbol{\bt})|_{\cal P},
\end{equation*}
and
\begin{equation*}
{\rm dist}(\overline{\langle ((1-\phi(p_i))+\phi(p_i)u)((1-\phi(q_i))+\phi(q_i)u^*)\rangle}, \Gamma(x_i))<\sigma,\,i=1,2,...,k.
\end{equation*}
\end{lem}

\begin{proof}
This follows from \ref{Extbot2} and \ref{BB-exi+}.
In fact,  for any $0<\ep_1<\ep/2$ and finite subset ${\cal F}_1\supset {\cal F},$  by applying \ref{Extbot2}, there exists $\dt,$ $n\ge 1,$  ${\cal Q}\subset K_1(C_n)$ and $\dt$ described above, and a unitary $u_1\in U_0(B)$  such that
\begin{equation*}
\|[\phi(x), \, u_1]\|<\ep_1\tforal x\in {\cal F}_1
\end{equation*}
and
\begin{equation*}
{\rm Bott}(\phi\circ \imath_{n,\infty}, u_1)  = \af(\boldsymbol{\bt})|_{\cal P}.
\end{equation*}

By choosing smaller $\ep_1$ and larger ${\cal F}_1,$ {{if}} necessary, we may assume that
$$\overline{\langle ((1-\phi(p_i))+\phi(p_i)u_1)((1-\phi(q_i))+\phi(q_i)u_1^*)\rangle}\in U_0(B)/CU(B)$$
is well defined for all $1\le i\le k.$
Define {{a}}  map $\Gamma_1: G\to U_0(B)/CU(B)$ by
\begin{equation}\label{Extbot2-4}
\Gamma_1(x_i)=\overline{\langle ((1-\phi(p_i))+\phi(p_i)u_1)(1-\phi(q_i))+\phi(q_i)u_1^*)\rangle},\quad i=1, 2, ..., k.
\end{equation}
%$i=1,2,...,k.$

%
%there exists
%a \hm\, $\Gamma_1: G\to U_0(B)/CU(B)$ such that
%\beq\label{Extbot2-4}
%{\rm dist}(\overline{\langle ((1-\phi(p_i))+\phi(p_i)u_1)((1-\phi(q_i))+\phi(q_i)u_1^*)\rangle}, \Gamma_1(x_i))&<&\sigma/2,
%\eneq
%$i=1,2,...,k.$
%{\bf \red{ Why can we just define that  $\Gamma_1(x_i)=\overline{\langle ((1-\phi(p_i))+\phi(p_i)u_1)(1-\phi(q_i))+\phi(q_i)u_1^*)\rangle}$.}}{\color{Green}(I think we can, as long as the the corresponding elements are well defined.)}
By choosing a  large $n,$ \wilog,    we may assume that there are projections $p_1', p_2',...,p_k',$
$q_1',q_2,'...,q_k'\in C_n$ such that $\imath_{n, \infty}(p_i')=p_i$ and
$\imath_{n, \infty}(q_i')=q_i,$ $i=1,2,...,k.$ Moreover, we may assume that ${\cal F}_1\subset \imath_{n, \infty}(C_n).$

Let $\Gamma_2: G\to U_0(B)/CU(B)$ {{be defined}} by
$\Gamma_2(x_i)=\Gamma_1(x_i)^*\Gamma(x_i),$ $i=1,2,...,k.$
It follows  \ref{BB-exi+} that is a unitary $v\in U_0(B)$ such that
\beq\label{Extbot2-5}
\|[\phi(x), \, v]\|<\ep/2\rforal x\in {\cal F},\\
{\rm Bott}(\phi\circ \imath_{n,\infty}, v)=0\andeqn\\
{\rm dist}(\overline{\langle ((1-\phi(p_i))+\phi(p_i)v)((1-\phi(q_i))+\phi(q_i)v^*)\rangle}, \Gamma_2(x_i))&<&\sigma,
%\sigma/2,
\eneq
$i=1,2,...,k.$ Define $u=u_1v,$
\beq
X_i&=&\overline{\langle ((1-\phi(p_i))+\phi(p_i)u_1)((1-\phi(q_i))+\phi(q_i)u_1^*)\rangle}\andeqn\\
Y_i&=&\overline{\langle ((1-\phi(p_i))+\phi(p_i)v)((1-\phi(q_i))+\phi(q_i)v^*)\rangle},
\eneq
 $i=1,2,...,k.$
We then compute that
\begin{equation}\label{Extbot2-6}
\|[\phi(x), \, u]\|<\ep_1+\ep/2<\ep\tforal x\in {\cal F},
\end{equation}
\begin{equation}
{\rm Bott}(\phi\circ \imath_{n,\infty}, u)={\rm Bott}(\phi\circ \imath_{n,\infty}, u_1)=\af(\boldsymbol{\bt})|_{\cal P}
\end{equation}
and
\begin{eqnarray}
&&{\rm dist}(\overline{\langle ((1-\phi(p_i))+\phi(p_i)u)((1-\phi(q_i))+\phi(q_i)u^*)\rangle}, \Gamma(x_i))\\
&\le & {\rm dist}(X_iY_i, \Gamma_1(x_i)Y_i)+{\rm dist}(\Gamma_1(x_i)Y_i, \Gamma(x_i))\\
&=&{\rm dist}(X_i, \Gamma_1(x_i))+{\rm dist}(Y_i, \Gamma_2(x_i))< \sigma,
\end{eqnarray}
for $i=1,2,...,k$.
\end{proof}

\section{Another Basic Homotopy Lemma}

\begin{lem}\label{densesp}
Let $A$  be a unital C*-algebra and let $U$ be an infinite dimensional  UHF-algebra. Then there is a unitary $w\in U$ such that for
any unitary $u\in A$, one has
\begin{equation}
\tau(f(u\otimes w))=\tau(f(1_A\otimes w))=\int_{\T}f dm,\quad f\in\mathbb C(\T),\ \tau\in T(A\otimes U)
\end{equation}
where $m$ is the normalized Lebesgue
%%%Lesbegue
 measure on $\T.$
Furthermore, for any $a\in A$  and $\tau\in T(A\otimes U)$, $\tau (a\otimes w^j)=0$ if $j\not=0$.
%(Note that $w\in CU(U)$.)
\end{lem}

\begin{proof}
Denote by $\tau_U$ the unique trace of $U$. Then any trace $\tau\in T(A\otimes U)$ is a product trace, i.e.,
$$\tau(a\otimes b)=\tau(a\otimes 1)\otimes\tau_U(b),\quad a\in A, b\in U.$$

Pick a unitary $w\in U$ such that the spectral measure of $w$ is the Lebesgue measure (a Haar unitary). Such a unitary always exists (it can be constructed directly; or, one can consider a strictly ergodic Cantor system $(\Omega, \sigma)$ such that $K_0(\mathrm{C}(\Omega)\rtimes_\sigma\Z)\cong K_0(U)$, and note that the canonical unitary in $\mathrm{C}(\Omega)\rtimes_\sigma\Z$ is a Haar unitary. Then by embedding $\mathrm{C}(\Omega)\rtimes_\sigma\Z$ into $U$, one obtains a Haar unitary in $U$).
%{\bf Zhuang:  If you find the reference, then we keep this and delete next one since that was all I did }
Then one has
$$
\tau_U(w^n)=
\left\{
\begin{array}{ll}
1, & \textrm{if $n=0$},\\
0, & \textrm{otherwise}.
\end{array}
\right.
$$
Hence, for any $\tau\in T(A\otimes U)$, one has
$$\tau((u\otimes w)^n)=\tau(u^n\otimes w^n)=\tau(u^n\otimes 1)\tau_U(w^n)=
\left\{
\begin{array}{ll}
1, & \textrm{if $n=0$},\\
0, & \textrm{otherwise};
\end{array}
\right.
$$
and therefore
$$\tau(P(u\otimes w))=\tau(P(1\otimes w))=\int_{\T}P(z)dm$$
for any polynomial $P$.
Since polynomials are dense in $\mathrm{C}(\T)$, one has
$$\tau(f(u\otimes w))=\tau(f(1\otimes w))=\int_{\T}fdm,\quad f\in \mathrm{C}(\T),$$
as desired.
\end{proof}

\begin{lem}\label{Defuv}
Let $A$ be a unital separable amenable \CA\, and let
$L: A\otimes C(\T)\to B$ be a unital completely positive linear map, where $B$ is another
unital amenable \CA. Suppose that $C$ is a unital \CA\, and $u\in C$ is a
unitary. Then, there is a  unique unital completely positive linear map $\Phi: A\otimes C(\T)\to B\otimes C$
such that
$$\Phi|_{A\otimes 1_{C(\T)}}=\imath\circ L|_{A\otimes 1_{C(\T)}}\quad \textrm{and}\quad \Phi(a\otimes z^j)=L(a\otimes z^j)\otimes u^j$$
for any $a\in A$ and any  integer $j,$ where $\imath: B\to B\otimes C$ is the standard inclusion.
%Then $\Phi$ is a unital \morp.

Furthermore,
if $\dt>0$ and ${\cal G}\subset A\otimes C(\T)$ is a finite subset, there is a $\dt_1>0$ and finite set  ${\cal G}_1\subset A\otimes C(\T)$  (which does not  depend on $L$) such that  if
$L$ is {{${\cal G}_1$-$\dt_1$-}}multiplicative, then
$\Phi$ is  { {${\cal G}$-$\dt$-}}multiplicative.

%and if
%$L$ is $\dt$-${\cal G}$-multiplicative, then
%$\Phi$ is also $\dt$-${\cal G}$-multiplicative.

\end{lem}

\begin{proof}
%%If $x=\sum_{i,j} a_i\otimes z^j,$ then $\Phi(x)=\sum_{i,j}L(a_i\otimes z^j)\otimes u^j.$ where $a_i\in A.$
%We first show that $\Phi$ is unital and positive.
%Let $a_{-n}, a_{-n+1},..., a_{-1}, a_0, a_1,a_2,...,a_n\in A.$ We will show that
%\beq\label{Defuv-2}
%\Phi(\sum_{-n\le i\le n} a_i\otimes z^i)^*(\sum_{-n\le i\le n}a_i\otimes z^i))\ge 0.
%\eneq
%%We compute that
%%\beq\label{Defuv-3}
%%&&\Phi(\sum_{i=1}^n a_i\otimes z^i)^*(\sum_{i=1}a_i\otimes z^i))=
%%\Phi(\sum_{i,j}a_i^*a_j\otimes z^{-i+j})\\
%%&=& \sum_{i,j}L(a_i^*a_j\otimes z^{-i+j})\otimes u^{-i+j}.
%%\eneq
%Let $C_0$ be the unital \SCA\, of $C$ generated by $u.$ Then
%$C_0$ is a commutative \CA.
%Define $L_1: A\otimes C(\T)\otimes C_0\to B\otimes C_0$ by
%$L_1=L\otimes {\id}_{C_0}.$ Since $L$ is positive, so is $L_1.$
%Note that
%\beq\label{Defuv-4}
%(\sum_{-n\le i\le n} (a_i\otimes z^i)\otimes u^i)^*(\sum_{-n\le i\le n}(a_i\otimes z^i)\otimes u^i)\\
%=\sum_{i,j}(a_i^*a_j\otimes z^{-i+j})\otimes u^{-i+j}\in (A\otimes C)_+
%\eneq
%Therefore
%\beq\label{Defuv-5}
%&&\Phi(\sum_{-n\le i\le n} a_i\otimes z^i)^*(\sum_{-n\le i\le n}a_i\otimes z^i))=
%\Phi(\sum_{i,j}a_i^*a_j\otimes z^{-i+j})\\
%&=& \sum_{i,j}L(a_i^*a_j\otimes z^{-i+j})\otimes u^{-i+j}\\
%&=& L_1(\sum_{i,j}(a_i^*a_j\otimes z^{-i+j})\otimes u^{-i+j})\ge 0.
%\eneq
%It follows that $\Phi$ on $A\otimes_{alg} C(\T)$ is unital and positive.
%Thus $\Phi$ is uniquely extended to $A\otimes C(\T)$ which is unital and
%positive.
%
%Since $L$ is positive on $M_n(A\otimes C(\T))=M_n(A)\otimes C(\T),$ the same argument holds by replacing $A$ by $M_n(A).$ It follows that $\Phi$ is \morp.
Denote by $C_0$ the unital C*-subalgebra of $C$ generated by $u$, then the tensor product map
$$L \otimes\id_{C_0}: A\otimes C(\T)\otimes C_0 \to B\otimes C_0$$
is unital and completely positive (see, for example, Theorem 3.5.3 of \cite{BO-Book}). Define the homomorphism $\psi: C(\T) \to C(\T)\otimes C_0$ by
$$\psi(z)=z\otimes u.$$ By Theorem 3.5.3 of \cite{BO-Book} again, the tensor product map
$$\id_A\otimes\psi: A\otimes C(\T) \to A\otimes C(\T)\otimes C_0$$
is unital and completely positive. Then the map
$$\Phi:=(L\otimes\id_{C_0})\circ(\id_A\otimes \psi)$$
satisfies the first part of the lemma.

Let us consider the second part of the lemma. Let $\dt>0$ and ${\cal G}\subset A\otimes C(\T)$ be a finite subset.
Suppose that $L$ is $\dt$-${\cal G}$-multiplicative. To simplify notation, without loss of generality, we may assume that
elements in ${\cal G}$ has the form $\sum_{-n\le i\le n} a_i\otimes z^i.$
{{Let $N=\max\{  n: \sum_{-n\le i\le n} a_i\otimes z^i\in {\cal G}\},$
let $\dt_1=\dt/2N^2$ and let
${\cal G}_1\supset \{a_i\otimes z^i: -n\le i\le n: \sum_{-n\le i\le n} a_i\otimes z^i\in {\cal G}\}.$ }}

Then
\beq\label{Defuv-6}
%\hspace{-0.5in}
&&\hspace{-0.5in}\Phi((\sum_{-n\le i\le n} a_i\otimes z^i)(\sum_{-n\le i\le n} b_i\otimes z^i))\\
&=&\sum_{i,j} \Phi(a_ib_j\otimes z^{i+j})\\
&=& \sum_{i,j}L(a_ib_j\otimes z^{i+j})\otimes u^{i+j}\\
%&\approx_{\dt {\color{Green} (2n+1)\delta?}}& \sum_{i,j} L(a_i\otimes z^i)L(b_j\otimes z^j)\otimes u^{i+j}\\
&\approx_{\dt}& (\sum_{-n\le i\le n}L(a_i\otimes z^i)\otimes u^i)(\sum_{-n\le i\le n}L(b_i\otimes z^i)\otimes u^i)\\
&=&\Phi(\sum_{-n\le i\le n} a_i\otimes z^i)\Phi(\sum_{-n\le i\le n} b_i\otimes z^i),
\eneq
if $\sum_{-n\le i\le n}a_i\otimes z^i, \sum_{-n\le i\le n}b_i\otimes z^i\in {\cal G}.$
It follows that $\Phi$ is {{${\cal G}$-$\dt$}}-multiplicative.

Let $P(\T)=\{\sum_{i=-n}^n c_iz^i, c_i\in \C\}$ be the algebra of all Laurent  polynomials.  The uniqueness follows from that $A\otimes C(\T)$ is the closure of the algebraic tensor product $A\otimes_{alg} P(\T)$.
\end{proof}

%{\color{Green} I did corresponding change on \ref{densesp3}.}

The following follows immediately from \ref{Defuv} and \ref{densesp}.

\begin{cor}\label{densesp3}
Let $C$ be a unital \CA\, and let $U$ be an infinite dimensional UHF-algebra.
%\, in ${\cal A}^b.$
%Suppose that $L: C\otimes C(\T)\to A$ is a unital \morp.
{{For any $\dt>0$ and any finite subset ${\cal G}\subset
C\otimes C(\T),$ there exist $\dt_1>0$ and a finite subset
 ${\cal G}_1\subset C\otimes C(\T)$ satisfying the following:}}
For any $1>\sigma_1,\,\sigma_2>0,$ any finite subset ${\cal H}_1\subset C(\T)_+\setminus \{0\}$ and any finite subset ${\cal H}_2\subset  (C\otimes C(\T))_{s.a.},$   %and there exists an integer $N\ge 1$
% such that, for
 %any integer $n\ge N,$
 any unital
 ${\cal G}_1$-$\dt_1$-multiplicative \morp\, $L: C\otimes C(\T)\to A,$
 where $A$ is another unital \CA,
 there exists unitary $w\in U$ satisfying the following:
 \beq\label{densesp3-1}
 |\tau(L_1(f))-\tau(L_2(f))|<\sigma_1\tforal f\in {\cal H}_2,\ \tau\in T(B),\ \textrm{and}\ \label{densesp3-2} \\
 \tau(g(1_A\otimes w))\ge \sigma_2(\int g dm)\tforal g\in {\cal H}_1,\ \tau\in T(B),
 \eneq
where $B=A\otimes U$ and $m$ is the normalized Lebesgue
%%%Lesbegue
 measure
 on $\T,$ and $L_1, L_2: C\otimes C(\T)\to A\otimes U$ are $\mathcal G$-$\dt$-multiplicative \morp s as $\Phi$ given by  {\rm  \ref{Defuv}} such that
 $L_i(c\otimes 1_{C(\T)})=L(c\otimes 1_{C(\T)})\otimes 1_U$ ($i=1,2$),  $L_1(c\otimes z^j)=L(c\otimes z^j)\otimes w^j,$
 and $L_2(c\otimes z^j)=L(c)\otimes w^j$ for all $c\in C$ and all
 integer {{$j.$}}
 %given by  {\rm  \ref{Defuv}}.
 %Moreover, one can choose $v$ and $w$ so that ${\rm det}(v)={\rm det}(w)=1.$
%
\end{cor}

%{\color{Green} From the proof of \ref{BHfull}, it seems that $\mathcal H_2$ above should be a subset of $C\otimes C(\T)$. Otherwise, it seems that $f$ only should be applied to $1\otimes z$ (in this case, one has the exact statement rather than a perturbation result).}

\begin{lem}\label{BHfull}
Let $A=A_1\otimes U_1,$
where $A_1\in {\cal B}_0$ which satisfies the UCT and $U_1$ is a UHF-algebra
of infinite type.
%and  $B=B_1\otimes U_2,$
 %where $A_1, B_1\in {\cal B}_1$  which satisfies the UCT, and $U_1$ and $U_2$ are  UHF-algebras of infinite type.
% Suppose that $\phi: A\to B$ is a unital \hm.
% Let $\Delta: C(\T)_+^{\bf 1}\setminus \{0\}\to (0,1)$ be an order preserving map.
For any $1>\ep>0$ and any finite subset ${\cal F}\subset A,$ there exist  $\dt>0,$ $\sigma>0$, a finite subset
${\cal G}\subset A,$ a finite subset $\{p_1,p_2,...,p_k, q_1,q_2,...,q_k\}$ of projections of $A$ such that
$\{[p_1]-[q_1],[p_2]-[q_2],...,[p_k]-[q_k]\}$ generates a free subgroup $G_u$ of $K_0(A),$
%and $m[1_A]\in G_u$  for some $m\ge 1,$
and a finite subset ${\cal P}\subset \underline{K}(A),$
%a finite subset ${\cal V}\subset K_1(A),$
%and a finite subset
%${\cal H}\subset C(\T)_+^{\bf 1}\setminus \{0\}$
satisfying the following:

Let $B=B_1\otimes U_2,$
where $B_1\in {\cal B}_0$  which satisfies the UCT and  $U_2$ are  UHF-algebras of infinite type.
Suppose that $\phi: A\to B$ is a unital \hm.

For  unitary $u\in U(B)$ such that
\beq\label{BHfull-1}
%\tau(f(u))\ge \Delta(\hat{f})\tforal f\in {\cal H},\\
&&\|[\phi(x),\, u]\|<\dt\tforal x\in {\cal G},\\\label{BHfull-1n}
&&{\rm Bott}(\phi,\, u)|_{\cal P}=0{{,}}\\\label{BHfull-1n1}
&&{\rm dist}(\overline{\langle ((1-\phi(p_i))+\phi(p_i)u)(1-\phi(q_i))+\phi(q_i)u^*)\rangle}, {\bar 1})<\sigma\andeqn\\\label{BHfull-1n2}
&&\mathrm{dist}(\bar{u}, \bar{1})<\sigma,
%{\rm dist}({\rm Bu}(\boldsymbol{\bt}(y)),{\bar{1}})<\eta
%\tforal y\in  {\cal V},
\eneq
there exists a continuous path of unitaries $\{u(t): t\in [0,1]\}\subset U(B)$ such
that
\beq\label{BHTL-3}
&&u(0)=u,\,\,\, u(1)=1_B{{,}}\\
&&{\rm dist}(u(t), CU(A))<\ep {\tforal} t\in [0,1],\\
&&\|[\phi(a),\, u(t)]\|<\ep\tforal a\in {\cal F}\tand for\,\, all\,\, t\in [0,1]\\
&&\tand {\rm length}(\{u(t)\})\le 2\pi+\ep.
\eneq
\end{lem}

\begin{proof}
Without loss of generality, one only has to prove the statement with assumption that $u\in CU(B)$.

In what follows we will use the fact that every C*-algebra in ${\cal B}_0$ has stable rank one.
% {\color{Green} (But we tensor everything with a UHF algebra here.)}
Define
$$
\Delta(f)=(1/2)\int f dm \tforal f\in C(\T)_+^{\bf 1}\setminus \{0\},
$$
where $m$ is the normalized Lebesgue measure on the unit circle $\T.$
Let $A_2=A\otimes C(\T).$ Let ${\cal F}_1=\{x\otimes f: x\in {\cal F}, f=1, z, z^*\}.$
To simplify notation, {\wilog}, we may assume that ${\cal F}$ is a subset of the unit ball of $A.$
Let $1>\dt_1>0$ (in place of $\dt$), ${\cal G}_1\subset A_2$ be a finite subset
(in place of ${\cal G}$), $1/4>\sigma_1>0, \,1/4>\sigma_2>0,$ ${\cal P}\subset \underline{K}(A_2)$
be a finite subset, ${\cal H}_1\subset C(\T)_+^{\bf 1}\setminus \{0\}$ be a finite subset, ${\cal H}_2\subset (A_2)_{s.a.}$
be a finite subset and ${\cal U}\subset U(M_2(A_2))/CU(M_2(A))$ be a finite subset as required by \ref{MUN1} for
$\ep/4$ (in place of $\ep$), ${\cal F}_1$ (in place of ${\cal F}$),  $\Delta$ and $A_2$ (in place of $A$).

We may assume, without loss of generality, that
$$
{\cal G}_1=\{a\otimes f: a\in {\cal G}_2\andeqn  f=1,z,z^*\},
$$
where ${\cal G}_2\subset A$ is a finite subset, and
$
{\cal P}={\cal P}_1\cup \boldsymbol{\bt}({\cal P}_2),
$
where ${\cal P}_1, {\cal P}_2\subset \underline{K}(A)$ are finite subsets.

We assume that $(2\dt_1, {\cal P}, {\cal G}_1)$ is a $KL$-triple for $A_2$, $(2\dt_1, {\cal P}_1, {\cal G}_2)$ is a $KL$-triple for $A$, and $1_{A_1} \otimes \mathcal H_1\subseteq \mathcal H_2$.

%We may also assume that
%$$
%{\cal H}_2=\{h_1\otimes h_2: h_1\in {\cal H}_3\andeqn h_2\in {\cal H}_4\},
%$$
%where ${\cal H}_3\subset A_{s.a.}$ and ${\cal H}_4\subset C(\T)_{s.a.}$ are finite subsets and ${\cal H}_1\subset {\cal H}_4.$

We may also choose $\sigma_1$ and $\sigma_2$ such that
\beq\label{Bufull-9}
\max\{\sigma_1, \sigma_2\}<(1/4)\inf\{\Delta(f): f\in {\cal H}_1\}.
\eneq

Let $\dt_2$ (in place of $\dt_1$) and a finite subset  ${\cal G}_3$ (in place of ${\cal G}_1$)
be as required by  \ref{densesp3} for $A$ (in place of $C$), $\dt_1/4$ (in place of $\dt$) and
${\cal G}_1$ (in place of ${\cal G}$).
By choosing even smaller $\dt_2,$
\wilog, we may assume that
${\cal G}_3=\{a\otimes f: g\in {\cal G}_2'\andeqn f=1,z,z^*\}$ with a large finite subset
${\cal G}_2'\supset {\cal G}_2.$
%We also assume $\dt_1$ is smaller than $\dt$ and ${\cal G}_1$ is larger
%than ${\cal G}$ that required by \ref{densesp3} for $A$ (in place of $C$),
%$\sigma_1$ %$\sigma_1/8$
%(in place of $\sigma_1$), 1/2 %$1-\sigma_2/16$
%(in place of $\sigma_2$), ${\cal H}_1$ and $\mathcal H_2$ %${\cal H}_2\cup {\cal H}_1$
%mentioned above.
We may assume that $\dt_2<\dt_1.$
%and ${\cal G}_3\supset {\cal G}_1.$

We may further assume
that,
\beq\label{BUfull-8}
{\cal U}={\cal U}_1\cup \{\overline{1\otimes z}\}\cup  {\cal U}_2,
\eneq
where ${\cal U}_1=\{\overline{a\otimes 1}: a\in {\cal U}_1'\subset U(A)\}$  and
${\cal U}_1'$ is a finite subset,
${\cal U}_2\subset U(A_2)/CU(A_2)$ is a finite subset
whose elements represent a finite subset of $\boldsymbol{\bt}(K_0(A)).$
So we may assume that ${\cal U}_2\in J_c(\boldsymbol{\bt}(K_0(A))).$
As in \ref{ReMUN1},  we may assume that the subgroup of $J_c(\boldsymbol{\bt}(K_0(A)))$
generated by ${\cal U}_2$ is free.  Let ${\cal U}_2'$ be a finite subset of unitaries such
that $\{{\bar x}: x\in {\cal U}_2'\}={\cal U}_2.$ We may also
assume that unitaries in ${\cal U}_2'$ has the form
\beq\label{Bufull-10-1}
((1-p_i)+p_i\otimes z)(1-q_i)+q_i\otimes z^*),\,\,\, i=1,2,...,k.
\eneq

We further assume that $p_i\otimes z\in {\cal G}_1,$ $i=1,2,...,k.$
Choose $\dt_3>0$ and a finite subset ${\cal G}_4'\subset A$ {{ (and denote ${\cal G}_4:=\{g\otimes f: g\in {\cal G}_4', f=1,z,z^*\}.$)}} such that, for any two
unital ${\cal G}_4$-$\dt_3$-multiplicative \morp s $\Psi_1, \Psi_2: A\otimes C(\T)\to C$ (any unital \CA\, $C$),
any ${\cal G}_4'$-$\dt_3$-multiplicative \morp s $\Psi_0: A \to C$ and unitary $V\in C$ ($1\le i\le k$), { {if
\beq\label{Bf-10+1}
&&\|\Psi_0(g)-\Psi_1(g\otimes 1)\|<\dt_3\rforal g\in {\cal G}_4'\\\label{Bf-10+2}
 &&\|\Psi_1(z)-V\|<\dt_3\andeqn
\|\Psi_1(g)-\Psi_2(g)\|<\dt_3\rforal g\in {\cal G}_4,
\eneq then}}
\beq\label{Bufull-10-1+}
&&\langle (1-\Psi_0(p_i)+\Psi_0(p_i)V)(1-\Psi_0(q_i)+\Psi_0(q_i) V^*\rangle\\
&&\hspace{0.3in}\approx_{{\sigma_2\over{2^{10}}}}\langle \Psi_1(((1-p_i)+p_i\otimes z){{(}}(1-q_i)+q_i\otimes z^*)\rangle,\\
&&\|\langle \Psi_1(x)\rangle -\langle \Psi_2(x)\rangle\|<{{\sigma_2/2^{10}}}\rforal x\in {\cal U}_2',\\
&&\Psi_1(((1-p_i)+p_i\otimes z)(1-q_i)+q_i\otimes z^*))\\\label{Bf-1510}
&&\hspace{0.3in}\approx_{{\sigma_2\over{2^{10}}}}
\Psi_1({{((1-p_i)+p_i\otimes z)}})\Psi_1((1-q_i)+q_i\otimes z^*)),
%\\\label{Bf-15107}
\eneq
{ {furthermore for $d_i^{(1)}=p_i,$ $d_i^{(2)}=q_i,$, there are projections ${\bar d}_i^{(j)}\in C$ and unitaries  ${\bar z}^{(j)}_i\in {\bar d}_i^{(j)}C {\bar d}_i^{(j)}$ such that}}
\beq\label{Bufull-10-1++}
&&\Psi_1(((1-d_i^{(j)})+d_i^{(j)}\otimes z))\approx_{{\sigma_2\over{2^{12}}}}(1-{\bar d}_i^{(j)})+{\bar z}^{(j)}_i\andeqn\\\label{Bf-15107+1-1}
&&{\bar d}_i^{(j)}\approx_{\sigma_2\over{2^{12}}} \Psi_1(d_i^{(j)}),\,\,{\bar z}_i^{(1)}\approx_{\sigma_2\over{2^{12}}} \Psi_1(p_i\otimes z),\andeqn
{\bar z}_i^{(2)}\approx_{\sigma_2\over{2^{12}}}  \Psi_1(q_i\otimes z^*),
\eneq
where
%$d_i^{(1)}=p_i,$ $d_i^{(2)}=q_i,$
%${\bar p}_i,$ ${\bar q}_i$ are projections and ${\bar z}^{(j)}_i\in {\bar d}_i^{(j)}C {\bar d}_i^{(j)}$ is a
%unitary,
$1\le i\le k,$ $j=1,2${{.}}
%provided that
%\beq\label{Bf-10+1}
%&&\|\Psi_0(g)-\Psi_1(g\otimes 1)\|<\dt_3\rforal g\in {\cal G}_4'\\\label{Bf-10+2}
 %&&\|\Psi_1(z)-V\|<\dt_3\andeqn
%\|\Psi_1(g)-\Psi_2(g)\|<\dt_3\rforal g\in {\cal G}_4,
%\eneq
%where ${\cal G}_4=\{g\otimes f: g\in {\cal G}_4', f=1,z,z^*\}.$ \Wlog, we may assume that ${\cal G}_2'\subset {\cal G}_4'.$
Choose $\sigma>0$ so it is smaller than $\min\{\sigma_1/16, \ep/16, \sigma_2/16, \dt_2/16, \dt_3/16\}.$

Choose $\dt_5>0$ and a finite subset ${\cal G}_5\subset A$ satisfying the following:
there is a unital  ${{{\cal G}_4}}$-$\sigma/8$-multiplicative \morp\,
$L: A\otimes C(\T)\to B'$ such that
\beq\label{Bufull-10}
&&\|L(a\otimes 1)-\phi'(a)\|<\sigma/8
%\min\{\dt_2/8, \ep/16\}
\tforal a\in {{{\cal G}_4'}}\andeqn
%\\\label{Bufull-10+}
\|L(1\otimes z)-u'\|<\sigma/8
%\min\{\sigma/8, \ep/16\}
\eneq
for any unital \hm\, $\phi': A\to B'$ and any unital $u'\in B'$ so that
$$
\|\phi'(g)u'-u'\phi'(g)\|<\dt_5\tforal  g\in {\cal G}_5.
$$
Let $\dt=\min\{\dt_5/4, \sigma\}$ and ${\cal G}={\cal G}_5\cup {\cal G}_4'\cup {\cal G}_2'.$
%By choosing smaller $\dt_1$ and larger ${\cal G}_1$ we may assume
%that $(2\dt_1, {\cal G}_1, {\cal P}_1)$ is a $KL$-triple
%for $A.$

Now suppose that $\phi: A\to B$ is a unital \hm\, and $u\in CU(B)$
which satisfy the assumption (\ref{BHfull-1}) to (\ref{BHfull-1n1}) for the above mentioned
$\dt,$ $\sigma,$ ${\cal G},$  ${\cal P},$ $p_i,$ and $q_i.$
There is an isomorphism $s: U_2\otimes U_2\to U_2.$ Moreover,
$s\circ\imath$ is approximately unitarily equivalent to the identity map on $U_2,$
where   ${\imath}: U_2\to U_2\otimes U_2$ defined by $\imath(a)=a\otimes 1$ (for all $a\in U_2$).
To simplify notation, without loss of generality, we may assume that
$\phi(A)\subset B\otimes 1\subset B\otimes U_2.$
Suppose that $u\in U(B)\otimes 1_{U_2}$ is a unitary which satisfies the assumption. {{As mentioned at the beginning,  we may assume that $u\in CU(B)\otimes 1_{U_2}.$
\Wlog, we may further assume that $u=\prod_{j=1}^{m_1}c_jd_jc_j^*d_j^*,$
where $c_j, d_j\in U(B)\otimes 1_{U_2},$ $1\le j\le m.$
Let ${\cal  F}_1= \{c_j, d_j: 1\le j\le m_1\}.$ }}

Let $L: A\otimes C(\T)\to B$ be a unital ${\cal G}_4$-$\dt_2/8$-multiplicative \morp\,
such that
\beq\label{BUfull11}
\|L(a\otimes 1)-\phi(a)\|<\sigma/8
%\min\{\dt_2/8, \ep/16\}
\tforal a\in {{{\cal G}_4'}}\andeqn
\|L(1\otimes z)-u\|<\sigma/8.
%\min\{\dt_2/8, \ep/16\}.
\eneq
Since ${\rm Bott}(\phi, u)|_{\cal P}=0$, we can also assume that
that
\beq\label{Bufull11+}
[L]|_{{\cal P}_1}=[\phi]|_{{\cal P}_1}\andeqn
[L]|_{\boldsymbol{\bt}({\cal P}_2)}=0.
\eneq

%For any $v\in U_2,$ there is a continuous path of unitaries
%$\{v(t): t\in [0,1]\}\subset 1\otimes U_2$ such that
%$v(0)=1_{U_2}$ and $v(1)=v.$
Since $B$ is in ${\cal B}_0,$ there is a projection $p\in B$ and a unital \SCA\,
$C\in {\cal C}_0$ with $1_C=p$ satisfying the following:
\beq\label{BUfull-12}
&&\|L(g)-[(1-p)L(g)(1-p)+L_1(g)]\|<{{\sigma^2/32(m_1+1)}}
%\min\{{\dt_2^2\over{{{32(m+1)}}}}, {\ep^2\over{{\red{32(m+1)}}}}\}
\tforal g\in {{ {\cal G}_4}}\\\label{BUfull-12+}
&&{{\andeqn \|(1-p)x-x(1-p)\|<\sigma^2/32(m_1+1)
%\min\{{\dt_2^2\over{32(m+1)}}, {\ep^2\over{32(m+1)}}\}
\rforal x\in {\cal F}_1,}}
\eneq
where $L_1: A{{\otimes C(\T)}}\to C$ is a unital ${{{\cal G}_4}}$-$\min\{\dt_2/8, \ep/8\}$-multiplicative
\morp,
\beq\label{BUfull-13}
\tau(1-p)<\min\{\sigma_1/16, \sigma_2/16\}\tforal \tau\in T(B)
\eneq
and,  using (\ref{BHfull-1n1}), (\ref{BHfull-1n2}), (\ref{BUfull11}) and (\ref{Bufull-10-1+}) to (\ref{Bf-10+2}),
%(and assuming that $\dt_1$ and $\sigma$
%are small enough comparing with $\sigma_2$),
we have that
\beq\label{Bufull-14}
&&{\rm dist}(L_2^{\ddag}(x), {\bar 1})<\sigma_2/4\tforal x\in \{1\otimes {\bar z}\}\cup {\cal U}_2\andeqn \\
&&{\rm dist}(L_2^{\ddag}(x), \overline{\phi(x')\otimes 1_{C(\T)}})<\sigma_2/4
\tforal  x\in {\cal U}_1,
\eneq
where $\overline{x'\otimes 1_{C(\T)}}=x,$
$L_2(a)=(1-p)L(a)(1-p)+L_1(a)$ for all $a\in A\otimes C(\T).$
%and ${\cal G}_2$ is a finite subset of $A$
%containing ${\cal G}_1.$
Note that  we also have
\beq\label{Bufull-15}
\|\phi(g)-L_2(g\otimes 1)\|<\sigma/2
%\min\{\dt_2/2, \ep/8\}
\rforal g\in {{{\cal G}_4'}}\andeqn
%\label{Bufull-16}
[L_2|_A]|_{{\cal P}_1}=[\phi]|_{{\cal P}_1}.
\eneq
%Moreover, we may assume that
%\beq\label{Bufull-17}
%{\rm dist}(L_1^{\ddag}(
%\eneq
By  (\ref{BUfull-12+}) and the choice of ${\cal F}_1$ there is a unitary $v_0\in CU(C)$ and a unitary\\ $v_{00}\in CU((1-p)B(1-p))$
 such that
\beq\label{Bufull-17}
&&\|L_1(1\otimes z)-v_0\|<\min\{\dt_2/2, \ep/8\}\andeqn \\\label{Bufull-17+}
&&\|(1-p)L(1\otimes z)(1-p)-v_{00}\|<\min\{\dt_2/2, \ep/8\}.
\eneq
By applying \ref{densesp3}, we obtain  a unitary $w\in U(U_2)=U_0(U_2)=CU(U_2)$ such that
%%two unitaries $v_2, w\in CU(U_2)$ such that
%and a continuous path of unitaries
%$\{v(t): t\in [1/4, 1/2]\}$ such that
\beq\label{Bufull-18}
%v(1/4)=1_{U_2},\,\,\, v(1/2)=v_0\otimes v_2\oplus (1-p)\otimes v_2 \andeqn\\\label{Bufull-18--}
&&|t(L_3(g)))-t(\Phi(g))|<\sigma_1,\quad g\in\mathcal H_2, % \min\{\sigma_1/4, \sigma_2/4\}\tforal g\in {\cal H}_2%\cup {\cal H}_1
\andeqn\\ \label{Bufull-18-+}
&&t(g(1\otimes w))\ge %(1-\sigma_2/16)
\frac{1}{2}\int_{\T} g dm\tforal g\in {\cal H}_1
\eneq
for all $t\in T(B\otimes U_2),$
where
$L_3: A\otimes C(\T)\to B\otimes U_2$ is  a unital ${\cal G}_1$-$\dt_1/4$-completely positive linear map defined by
\beq\label{Bufull-18-}
L_3(a\otimes 1)=L_2(a\otimes 1)\otimes 1_{U_2} \andeqn
L_3(a\otimes z^j)=L_2(a\otimes z^j)\otimes (w)^j
\eneq
for all $a\in A$ and all integers $j$ as given by \ref{Defuv},
%$v_2(1-p)=(1-p)v_2=(1-p),$
$w(1-p)=(1-p)w=(1-p)$, as we consider both $w$, $1-p$ as elements in $B\otimes U_2$ as $1_B\otimes w$ and $(1-p)\otimes 1_{U_2}$, respectively, and
$\Phi: A\otimes C(\T)\to B\otimes U_2$  is defined by
$\Phi(a\otimes 1)=\phi(a)\subset B\otimes 1$ for all $a\in A$ and
$\Phi(1\otimes f)=f(w)$ for all $f\in C(\T).$ Moreover,
%by the virtue of
%(\ref{BUfull-13}), we may assume that $(1-p)w=w(1-p)=1-p$ and
$\Phi(1\otimes f)=
f(\lambda)((1-p)\otimes 1_{U_2})+f(p w)$
for all $f\in C(\T)$  and for some $\lambda\in \T.$
%By \ref{2Lg9}, there is a continuous path of unitaries  in $CU(C)$
%connecting $p$ to $v_0$ with length at most $4\pi.$
Note that $CU(U_2)=U(U_2).$
%It is also known (by working in matrices)  that there is a continuous
%path of unitaries in $CU(U_2)$ connecting $1_{U_2}$ to $v_2$ with length
%no more than $\pi+\ep/256.$
One obtains a continuous path of unitaries
$\{v(t): t\in [1/4, 1/2]\}\subset CU(U_2)$ such that
\beq\label{Bufull-1409}
v(1/4)=1_{U_2},\,\,\, v(1/2)=w\andeqn
{\rm length}(\{v(t): t\in [1/4, 1/2]\})\le \pi+\ep/256.
\eneq
Note that $\phi(a)\Phi(1\otimes z)=\Phi(1\otimes z)\phi(a)$
for all $a\in A.$ So, in particular, $\Phi$ is a unital \hm\, and
\beq\label{Bf-moved}
[\Phi]|_{\boldsymbol{\bt}(\underline{K}(A))}=0.
\eneq
%Note  that $L_3$ is $\dt_1/2$-${\cal G}_1$-multiplicative.
Define a unital \morp\,
$L_t: A_2\to C([2,3], B\otimes U_2)$ by
$$
L_t(f\otimes 1)=L_2(f\otimes 1)\andeqn L_t(a\otimes z^j)=L_2(a\otimes z^j)\otimes (v((t-2)/4+1/4))^j
$$
for all $a\in A$ and integers $j$ and $t\in [2,3]$.
%as in \ref{Defuv}.
Moreover, $L_t(1\otimes z)=(v_0\oplus v_{00})\otimes v((t-2)/4+1/4)),$ and,  since $v(s)\in CU(U_2),$
$L_t(1\otimes z)\in CU(B\otimes U_2)$ for all $t\in [2, 3].$
Note, as in the proof
of \ref{Defuv},  that $L_t$ are ${\cal G}_1$-$\dt_1/4$-multiplicative.
Note at $t=2,$ $L_t=L_2$ and at $t=3,$ $L_t=L_3.$
It follows that
\beq\label{Bufull-18+}
&&[L_3]|_{{\cal P}_1}=[L_2]|_{{\cal P}_1}=[\phi]|_{{\cal P}_1},\,\,\,
[L_3]|_{{\boldsymbol{\bt}}({\cal P}_2)}=0\andeqn\\\label{Bufull-18++}
&&L_3^{\ddag}(x)=L_2^{\ddag}(x)\tforal x\in {\cal U}_1.
\eneq
If $v=(e\otimes z)+(1-e)$ for some projection $e\in A,$ then
\beq\label{Bufull-19-1}
L_3(v)=L_2(e\otimes z)\otimes w+L_2((1-e)).
\eneq
Since $w\in CU(U_2),$
%ying \ref{UCUiso},
one computes  from (\ref{Bf-1510}) %and (\ref{Bf-15107})
%{\color{Green} (I erase {Bf-15107}, since it is the same as \eqref{Bf-1510} now)}
 that
%and (\ref{Bufull-10-1})
that, with $x=((1-p_i)+p_i\otimes z)((1-q_i)+q_i\otimes z^*),$
\beq\label{Bf-15107+1}
&&\hspace{-0.4in}{{\overline{\langle L_3(x)\rangle}\approx_{\sigma_2/2^{10}} \overline{({\bar z}_i^{(1)}\otimes w+(1-{\bar p}_i))({\bar z}_i^{(2)}\otimes w+(1-{\bar q}_i))}}}\\
&&\hspace{-0.1in}{{=\overline{({\bar z}_i^{(1)}+(1-{\bar p}_i))({\bar p}_i\otimes w+(1-{\bar p}_i)\otimes 1_{U_2})({\bar z}_i^{(2)}+(1-{\bar q}_i))
({\bar q}_i\otimes w+(1-{\bar q}_i))}}}\\
&&{{=\overline{({\bar z}_i^{(1)}+(1-{\bar p}_i))({\bar z}_i^{(2)}+(1-{\bar q}_i))}={\overline{\langle L_2(x)\rangle}}}},
\eneq
where ${\bar p}_i, {\bar q}_i, {\bar z}_i^{(1)}, {\bar z}_i^{(2)}$ are as above (see the lines below (\ref{Bf-1510})), replacing $\Psi_1$ by $L_2.$  It follows that
\beq\label{Bufull-19-2}
{\rm dist}(L_3^{\ddag}(x), {\bar 1})<\sigma_2/2 \tforal x\in \{\overline{1\otimes z}\}\cup {\cal U}_2.
\eneq
Note that, since $w\in CU(U_2)$ and $\phi(q)\in B\otimes 1_{U_2},$
\beq\label{Bufull-19+3}
\Phi(q\otimes z+(1-q)\otimes 1)=
\phi(q)\otimes w+\phi(1-q)\in CU(B\otimes U_2)
\eneq
for any projection $q\in A.$
It follows  that
\beq\label{Bufull-19}
%&&[\Phi]|_{\boldsymbol{\bt}(\underline{K}(A))}=0,\\
&&\Phi^{\ddag}(x)\in CU(B\otimes U_2)\tforal x\in \{\overline{1\otimes z}\}\cup {\cal U}_2.
\eneq
Therefore (see also (\ref{Bf-moved}))
\beq\label{Bufull-20}
[L_3]|_{{\cal P}}=[\Phi]|_{\cal P}\andeqn
{\rm dist}(\Phi^{\ddag}(x), L_3^{\ddag}(x))<\sigma_2\tforal x\in {\cal U}.
\eneq
It follows from \eqref{Bufull-18-+} that
\begin{equation}\label{Bufull-21}
\tau(\Phi(f))\ge \Delta(f),\quad f\in {\cal H}_1,\ \tau\in T(B\otimes U_2),
\end{equation}
and it follows from \eqref{Bufull-18} that
\begin{equation}
|\tau(\Phi(f))-\tau(L_3(f))|<\sigma_1,\quad f\in\mathcal H_2,\ \tau\in T(B\otimes U_2).
\end{equation}
%Moreover, by (\ref{Bufull-18--}),
%\beq\label{Bufull-22}
%\|
By applying \ref{MUN1}, we obtain a unitary $w_1\in B\otimes U_2$ such that
\beq\label{Bufull-22}
\|w_1^*\Phi(f)w_1-L_3(f)\|<\ep/4\tforal f\in {\cal F}_1.
\eneq
Since $w\in U_2,$ there is a continuous path of unitaries $\{w(t):t\in [3/4, 1]\}\subset CU(U_2)$ such that
\beq\label{Bufull-23}
\hspace{-0.2in}w(3/4)=\Phi(1\otimes z)=w,\,\,\, w(1)=1_{U_2}\andeqn {\rm length}(\{w(t):t\in [3/4,1]\})\le \pi+\ep/256.
\eneq
Note that
\begin{equation}\label{Bufull-24}
\Phi(a)w(t)=w(t)\Phi(a)\tforal a\in A\andeqn t\in [3/4,1].
\end{equation}
It follows from (\ref{Bufull-22}) that there exists
a continuous path of unitaries $\{u(t): t\in [1/2, 3/4]\}\subset B\otimes U_2$ such that
\beq\label{Bufull-25}
u(1/2)=(v_{00}+(v_0))\otimes w,\,\,\,
u(3/4)=w_1^*\Phi(1\otimes z)w_1\andeqn\\\label{Bufull-25+}
\|u(t)-u(1/2)\|<\ep/4\tforal t\in [1/2, 3/4].
\eneq
It follows from (\ref{Bufull-10})
%{(\color{Green} change from Bufull-10+ to Bufull-10)}
and (\ref{Bufull-17+}) that
there exists a continuous path of unitaries
$\{u(t): t\in [0, 1/4]\}\subset B$ such that
\beq\label{Bufull-26}
u(0)=u,\,\,\,u(1/4)=v_{00}+v_0\andeqn\\\label{Bufull-26+}
\|u(t)-u\|<\ep/4\tforal t\in [0, 1/4].
\eneq
Now define
\beq\label{Bufull-27}
u(t)=w_1^*w(t)w_1\tforal t\in [3/4,1].
\eneq
Then $\{u(t): t\in [0,1]\}\subset B\otimes U_2$ is a continuous path of unitaries such that $u(0)=u$ and $u(1)=1.$ Moreover, by (\ref{Bufull-22}),
(\ref{Bufull-23}),  (\ref{Bufull-26+}), (\ref{Bufull-1409}),  (\ref{Bufull-23}) and (\ref{Bufull-25+}),
\beq\label{Bufull-28}
\|\phi(f)u(t)-u(t)\phi(f)\|<\ep\tforal f\in {\cal F}\andeqn {\rm length}(\{u(t)\})\le 2\pi+\ep.
\eneq

\end{proof}

\begin{rem}\label{RBHfull}

Let $A$ be a unital simple separable amenable \CA\, with stable rank one.
%{\color{Green} (It seems that one uses that $A$ is stably finite rather than stable rank one. One also does not use simplicity of $A$.)}
Let $G_0\subset K_0(A)$ be a finitely generated subgroup containing $[1_A].$
Let $G_r=\rho_A(G_0).$ Then $\rho_A([1_A])\not=0$ {and} $G_r$ is a finitely generated free
group. Then we may write $G_0=G_0\cap {\rm ker}\rho_A\oplus G_r',$
where $\rho_A(G_r')=G_r$ and $G_r'\cong G_r.$
Note that $G_0\cap {\rm ker}\rho_A$ is {{a}} finitely generated subgroup.
We may write $G_0\cap {\rm ker}\rho_A=G_{00}\oplus G_{01},$ where
$G_{00}$ is a torsion group and $G_{01}$ is free.
Note that $G_{01}\oplus G_r'$ is  free.
Therefore  $G_0=Tor(G_0)\oplus F,$ where $F$ is a finitely generated free subgroup.
Note that there is an integer $m\ge 1$ such that $m[1_A]\in F.$
Let $z\in C(\T)$ be the standard unitary generator. Consider $A\otimes C(\T).$
Then $\boldsymbol{\bt}(G_0)\subset \boldsymbol{\bt}(K_0(A))$ is a subgroup
of $K_1(A\otimes C(\T)).$  Moreover $\boldsymbol{\bt}([1_A])$ may be identified with
$[1\otimes z].$

If we choose ${\cal U}_2$ in the above proof
% {\color{Green} (I cannot find $\mathcal U_2$)}
which generates $\boldsymbol{\bt}(F),$
then, for any $\sigma_1>0,$ we may assume that
\beq\label{RBH-1}
{\rm dist}({\overline{ u^m}}, {\bar 1})<\sigma_1/m
\eneq
provided that (\ref{BHfull-1n1}) holds for a sufficiently small $\sigma.$
Since we assume that $u\in U_0(B)$ as ${\cal P}$ may be large enough in (\ref{BHfull-1n}),
by (\ref{RBH-1}),
%\beq\label{RBH-2}
${\rm dist}({\bar u}, {\bar 1})<\sigma_1.$
%\eneq
This implies that (with sufficiently small $\sigma$) the condition (\ref{BHfull-1n2}) is redundant
and therefore can be  {{omitted.}}
% removed.

\end{rem}

\section{Stably results}

\begin{lem}\label{stablehomtp}
Let $C$ be a unital separable \CA\,  which
%whose irreducible representations have bounded dimensions {\color{Green}(it seems that one only has
%to assume that $C$
is residually finite dimensional and satisfies the UCT.  For any $\ep>0,$ any finite subset ${\cal F}\subset C,$ any finite subset ${\cal P}\subset \underline{K}(C),$ any unital \hm\, $h: C\to A,$ where $A$ is any unital \CA, and any ${\kappa}\in Hom_{\Lambda}(\underline{K}({SC}), \underline{K}(A)),$ there exists an integer $N\ge 1,$ a unital \hm\,
$h_0: C\to M_N(\C)\subset M_N(A)$  and  a unitary $u\in U(M_{N+1}(A))$ such that
\beq\label{sthomp-1}
\|H(c),\, u]\|<\ep\tforal c\in {\cal F}\andeqn {\rm Bott}(H,\, u)|_{\cal P}=\kappa,
\eneq
where $H(c)={\rm daig}(h(c), h_0(c))$ for all $c\in C.$

\end{lem}

\begin{proof}
Define $S=\{z, 1_{C(\T)}\},$ where $z$ is identity function on the unit circle. Define $x\in {\rm Hom}_{\Lambda}(\underline{K}(C\otimes C(\T)), \underline{K}(A))$ as follows:
\beq\label{sthomp-2}
x|_{\underline{K}(C)}=[h]\andeqn x|_{\boldsymbol{\bt}(\underline{K}(C))}=\kappa.
\eneq
Fix a finite subset ${\cal P}_1\subset \boldsymbol{\bt}(\underline{K}(C)).$
Choose $\ep_1>0$ and a finite subset ${\cal F}_1\subset C$ satisfying the following:
\beq\label{sthomp-3}
[L']|_{{\cal P}_1}=[L'']|_{{\cal P}_1}
\eneq
for any pair of {{(${\cal F}_1\otimes S)-\ep_1$}}-multiplicative \morp s $L',L'':C\otimes C(\T)\to B$ (for any unital \CA\, $B$), provided  that
\beq\label{sthomtp-4}
L'\approx_{\ep_1} L''\,\,\,{\rm on} \,\, {\cal F}_1\otimes S.
\eneq

Let a positive number $\ep>0,$  a finite subset ${\cal F}$  and a finite subset ${\cal P}\subset \underline{K}(C)$ be given.
We may assume, without loss of generality,  that
\beq\label{sthomtp-4+1}
{\rm Bott}(H',\, u')|_{\cal P}={\rm Bott}(H',\, u'')|_{\cal P}
\eneq
provided $\|u'-u''\|<\ep$ for any unital \hm\, {{$H'$}} from $C.$
Put
$\ep_2=\min\{\ep/2, \ep_1/2\}$ and ${\cal F}_2={\cal F}\cup {\cal F}_1.$

Let $\dt>0,$ ${\cal G}\subset C$ be a finite {{subset}} and ${\cal P}_0\subset
\underline{K}(C)$ (in place of ${\cal P}$) be as required by \ref{Newstableuniq}
%{\red{warning: I modified  \ref{Newstableuniq} }}
for $\ep_2/2$ (in place of $\ep$) and ${\cal F}_2$ (in place of ${\cal F}$).
Without loss of generality, we may assume that ${\cal F}_2$ and ${\cal G}$ are in the unit ball of $C$ and $\dt<\min\{1/2, \ep_2/16\}.$
Fix another finite subset ${\cal P}_2\subset \underline{K}(C)$ and defined
${\cal P}_3={\cal P}_0\cup {\boldsymbol{\bt}}({\cal P}_2)$ (as a subset of
$\underline{K}(C\otimes C(\T))$). We may assume that ${\cal P}_1\subset
{\boldsymbol{\bt}}({\cal P}_2).$

It follows from {\ref{kkmaps}} that there are integer $N_1\ge 1,$ a unital \hm\, $h_1:C\otimes C(\T)\to M_{N_1}(\C)\subset M_{N_1}(A)$ and a {{$({\cal G}\otimes S)$-$\dt/2$}}-multiplicative \morp\, $L: C\otimes C(\T)\to M_{N_1+1}(A)$ such that
\beq\label{sthomtp-5}
[L]|_{{\cal P}_3}=(x+[h_1])|_{{\cal P}_3}.
\eneq
We may assume that there is a unitary $v_0\in M_{N_1+1}(A)$ such that
\beq\label{sthomtp-5+}
\|L(1\otimes z)-v_0\|<\ep_2/2.
\eneq

Define $H_1: C\to M_{N_1+1}(A)$ by
\beq\label{sthomtp-6}
H_1(c)=h(c)\oplus h_1(c\otimes 1)\tforal c\in C.
\eneq
Define $L_1: C\to M_{N_1+1}(A)$ by $L_1(c)=L(c\otimes 1)$ for all $c\in C.$ Note that
\beq\label{sthomtp-7}
[L_1]|_{{\cal P}_0}=[H_1]|_{{\cal P}_0}.
\eneq
It follows from \ref{Newstableuniq} that there exists an integer $N_2\ge 1,$ a unital \hm\, $h_2: C\to M_{N_2(N_1+1)}(\C){{\subset M_{N_2(N_1+1)}(A)}}$ and a unitary $W\in M_{(N_2+1)(1+N_1)}(A)$ such that
\beq\label{sthomtp-8}
W^*(L_1(c)\oplus h_2(c))W\approx_{\ep/4} H_1(c)\oplus h_2(c)\rforal c\in {\cal F}_2.
\eneq
Put $N=N_2(N_1+1)+N_1.$ Now define $h_0: C\to M_N(\C)$ and $H: C\to M_{N+1}(A)$ by
\beq\label{sthomtp-9}
h_0(c)=h_1(c\otimes 1)\oplus h_2(c)\andeqn H(c)=h(c)\oplus h_0(c)
\eneq
for all $c\in C.$
Define
%\beq\label{sthomtp-10}
$u=W^*(v_0\oplus 1_{M_{N_2(N_1+1)}})W.$
%\eneq
Then, by ({{\ref{sthomtp-8}}}), { {and $L_1$ being $({\cal G}\otimes S)$-$\dt/2$- multiplicative, we have }}
\beq\label{sthomtp-11}
&&\hspace{-0.2in}\|[H(c),\, u]\| \le \|(H(c)-{\rm Ad}\, W\circ (L_1(c)\oplus h_2(c))) u]\|\\
&&+\|{\rm Ad}\, W\circ (L_1(c)\oplus h_2(c)),\, u]\|+
\|u(H(c)-{\rm Ad}\, W\circ (L_1(c)\oplus h_2(c)))\|\\
&&<\ep/4+\dt/2+\ep/4<\ep\rforal c\in {\cal F}_2.
\eneq
Define $L_2: C\to M_{N+1}(A)$ by $L_2(c)=L_1(c)\oplus h_2(c)$ for all $c\in C.$ Then, we compute that
\beq\label{sthomtp-12}
{\rm Bott}(H,\, u)|_{\cal P}&=& {\rm Bott}({\rm Ad}\, W\circ L_2,\,u)|_{\cal P}=
 {\rm Bott}(L_2,\, v_0\oplus 1_{M_{N_2(N_1+1)}})|_{\cal P}\\
&=& {\rm Bott}(L_1,\, v_0)|_{\cal P}+{\rm Bott}(h_2,\, 1_{M_{N_2(N_1+1)}})|_{\cal P}\\
&=& [L]|_{{\boldsymbol{\bt}({\cal P})}}+0=
(x+[h])|_{{\boldsymbol{\bt}}({\cal P})}=\kappa|_{\cal P}.
\eneq

\end{proof}

\begin{thm}\label{STHOM}
Let $C$ be a unital amenable separable \CA\,  which  is residually finite dimensional and satisfies the UCT.  For any $\ep>0$ and any finite subset ${\cal
F}\subset C,$ there is $\dt>0,$ a finite subset ${\cal G}\subset
C,$ a finite subset ${\cal P}\subset \underline{K}(C)$
 satisfying the following:

Suppose that $A$ is a unital \CA, suppose $h: C\to A$ is a unital
\hm\, and suppose that $u\in U(A)$ is a unitary such that
\beq\label{Shomp1}
\|[h(a), u]\|<\dt\tforal a\in {\cal G}\tand
{\rm{Bott}}(h,u)|_{\cal P}=0.
\eneq
Then there exists an integer $N\ge 1$ and a continuous path of
unitaries $\{U(t): t\in [0,1]\}$ in $M_{N+1}(A)$ such that
\beq\label{Shomp2}
U(0)=u',\,\,\, U(1)=1_{M_{N+1}(A)}\andeqn \|[h'(a),
U(t)]\|<\ep\tforal a\in {\cal F},
\eneq
where
$$
u'={\rm diag}(u, H_0(1\otimes z))
$$
and $h'(f)=h(f)\oplus H_0(f\otimes 1)$ for $f\in C,$ where $H_0:
C\otimes C(\T)\to M_N(\C)$ {{($\subset M_N(A)$)}} is a unital \hm\, (with finite dimensional range)
and $z\in C(\T)$ is the identity function on the unit circle.

Moreover,
\beq\label{Shomp2+}
\rm{Length}(\{U(t)\})\le \pi+\ep.
\eneq

\end{thm}

\begin{proof}
Let $\ep>0$ and ${\cal F}\subset C$  be given. Without loss of
generality, we may assume that ${\cal F}$ is in the unit ball of
$C.$

Let $\dt_1>0,$  ${\cal G}_1\subset C\otimes C(\T),$ ${\cal
P}_1\subset \underline{K}(C\otimes C(\T))$  be required by
\ref{Newstableuniq} for $\ep/4$ and ${\cal F}\otimes S.$  Without loss of
generality, we may assume that ${\cal G}_1={\cal G}_1'\otimes S,$
where ${\cal G}_1'$ is in the unit ball of $C$ and
$S=\{1_{C(\T)}, z\}\subset C(\T).$ Moreover, without loss of
generality, we may assume that ${\cal P}_1={\cal P}_2\cup{\cal
P}_3,$ where ${\cal P}_2\subset \underline{K}(C)$ and ${\cal
P}_3\subset {\boldsymbol{\bt}}(\underline{K}(C)).$ Let ${\cal P}= {\cal P}_2 \cup \beta^{-1}({\cal P}_3) \subset \underline{K}(C)$. Furthermore, we
may assume that any $\dt_1$-${\cal G}_1$-multiplicative \morp\,
$L'$ from $C\otimes C(\T)$ to a unital \CA\, well defines $[L']|_{{\cal P}_1}.$

Let $\dt_2>0$ and ${\cal G}_2\subset C$ be a finite subset
required by 2.8 of \cite{LnHomtp} for $\dt_1/2$ and ${\cal G}_1'$
above.

Let $\dt=\min\{\dt_2/2, \dt_1/2, \ep/2\}$ and ${\cal G}={\cal
F}\cup {\cal G}_2.$

Suppose that $h$ and $u$ satisfy the assumption with above $\dt,$
${\cal G}$ and ${\cal P}.$ Thus, by 2.8 of \cite{LnHomtp}, there is
$\dt_1/2$-${\cal G}_1$-multiplicative \morp\, $L: C\otimes
C(\T)\to A$ such that
\beq\label{Shomp4}
&&\|L(f\otimes 1)-h(f)\|<\dt_1/2\rforal f\in {\cal G}_1'\\
&& \|L(1\otimes z)-u\|<\dt_1/2.
\eneq

Define $y\in Hom_{\Lambda}(\underline{K}(C\otimes C(\T)),
\underline{K}(A))$ as follows:
$$
y|_{\underline{K}(C)}=[h]|_{\underline{K}(C)}\andeqn
y|_{\boldsymbol{\bt}(\underline{K}(C))}=0.
$$
{It follows from ${\rm{Bott}}(h,u)|_{\cal P}=0$ that $[L]|_{\bt({\cal P})}=0$.}

Then
\beq\label{Shomp4+1}
[L]|_{{\cal P}_{ 1}}=y|_{{\cal P}_{ 1}}.
\eneq

Define $H: C\otimes C(\T)\to A$ by
$$
H(c\otimes g)=h(c)\cdot g(1)\cdot 1_A
$$
for all $c\in C$ and $g\in C(\T),$ where $\T$ is identified with
the unit circle (and $1\in \T$).

It follows that
\beq\label{Shomp3}
[H]|_{{\cal P}_{ 1}}= y|_{{\cal P}_{ 1}}=[L]|_{{\cal P}_{ 1}}.
\eneq

It follows from \ref{Newstableuniq} that there is an integer $N\ge 1,$ a
unital \hm\, $H_0: C\otimes C(\T)\to M_N(\C)$ {{($\subset M_N(A)$)}} with finite
dimensional range and a unitary $W\in U(M_{1+N}(A))$ such that
\beq\label{Shomp5}
W^*(H(c)\oplus H_{0}(c))W\approx_{\ep/4} L(c)\oplus H_{0}(c)\rforal c\in {\cal F}\otimes S.
\eneq
%for all $c\in {\cal F}\otimes S.$

Since $H_0$ has finite dimensional range {{and since $H_0(1\otimes z)$ is in the center of range($H_0$) $\subset M_N(\C)$, }} it is easy to construct
a continuous path $\{V'{ (t)}: t\in [0,1]\}$ in a finite dimensional
\SCA\, of $M_N(\C)$  such that
\beq\label{Shomp6}
&&V'(0)=H_0(1\otimes z),\,\,\,
V'(1)=1_{M_{N}(A)}\andeqn\\
&&H_0(c\otimes 1)V'(t)=V'(t)H_0(c\otimes 1)
\eneq
for all $c\in C$ and $t\in [0,1].$ Moreover,
\beq\label{Shomp6+}
\text{Length}(\{V'(t)\})\le \pi.
\eneq

Now define $U(1/4+3t/4)=W^*{\rm diag}(1, V'(t))W$ for $t\in [0,1]$
and
$$
u'=u\oplus  H_0(1_A\otimes z)\andeqn h'(c)=h(c)\oplus H_0(c\otimes
1)
$$
for $c\in C$ for $t\in [0,1].$  Then
\beq\label{Shomp7}
\|u'-U(1/4)\|<\ep/4\andeqn \|[U(t),\, h'(a)]\|<\ep/4
\eneq
for all $a\in {\cal F}$ and $t\in [1/4,1].$ The theorem follows by
connecting $U(1/4)$ with $u'$ with a short path as follows: There
is a self-adjoint element $a\in M_{1+N}(A)$ with $\|a\|\le {\ep
\pi\over{8}}$ such that
\beq\label{Shomp-8}
\exp(i a)=u'U(1/4)^*
\eneq
Then the path of unitaries $U(t)=\exp(i (1-4t) a)U(1/4)$ for $t\in [0,1/4)$ satisfy the above.
\end{proof}

\begin{lem}\label{VuV}
Let $C$ be a unital separable \CA\,  whose irreducible representations have bounded dimensions and $B$ be a unital \CA\, with $T(B)\not=\emptyset.$
Suppose  $\phi_1, \, \phi_2:C\to B$ are two unital monomorphisms
such that
$$
[\phi_1]=[\phi_2]\,\,\,{\rm in} \,\,\, KK(C,B),\\
%(\phi_1)_T=(\phi_2)_T\tand\,{\overline{R}_{\phi_1, \phi_2}}=0.
$$
Let $\theta: \underline{K}(C)\to \underline{K}(M_{\phi_1, \phi_2})$ be the splitting map defined in  {\rm(\ref{Aug2-2})}.
%such that $[\pi_0]\circ \theta=[\phi_1].$

For any $1/2>\ep>0,$ any finite subset ${\cal F}\subset C$ and any
finite subset ${\cal P}\subset \underline{K}(C),$ there are
integers $N_1\ge 1,$  an $\ep/2$-${\cal F}$-multiplicative \morp\,
$L: C\to M_{1+N_1}(M_{\phi_1,\phi_2}),$  a unital \hm\,
$h_0: C\to M_{N_1}(\C),$
%unital \hm\, $h_0'': C\to M_{N_1+N_2-1}(\C),$
and a continuous path of unitaries $\{V(t): t\in [0,1-d]\}$ of
$M_{1+N_1}(B)$ for some $1/2>d>0,$ such that $[L]|_{\cal P}$ is
well defined, $V(0)=1_{M_{1+N_1}(B)},$
\beq\label{VuV1}
[L]|_{\cal P}&=&(\theta+[h_0])|_{\cal P},\\
%\eneq
%\beq
\label{VuV2}
\pi_t\circ L&\approx_{\ep}&{\rm ad}\, V(t)\circ
(\phi_1\oplus h_0)\,\,\,\text{on}\,\,\,{\cal F} \tforal t\in (0,1-d],\\
%\eneq
%for all $t\in (0,1-d],$
%\beq
\label{VuV3}
\pi_t\circ L&\approx_{\ep}&{\rm ad}\, V(1-d)\circ
(\phi_1\oplus h_0)\,\,\,\text{on}\,\,\,{\cal F} \tforal t\in (1-d, 1]\andeqn\\
%\eneq
%for all $t\in (1-d,1),$  and
%\beq
\label{VuV33}
\pi_1\circ L&\approx_{\ep}&\phi_2\oplus h_0\,\,\,\text{on}\,\,\,{\cal F},
\eneq
where $\pi_t: M_{\phi_1,\phi_2}\to B$ is the point-evaluation at
$t\in (0,1).$

\end{lem}

\begin{proof}
Let $\ep>0$ and let ${\cal F}\subset C$ be a finite subset.
Let $\dt_1>0,$ ${\cal G}_1\subset C$ be a finite subset and ${\cal
P}\subset \underline{K}(C)$ be a finite subset required by
\ref{STHOM} for $\ep/4$ and ${\cal F}$ above.
In particular, we assume that $\dt_1<\dt_{\cal P}$ (see \ref{Dbeta}).
We may further assume that $\dt_1$ is sufficiently small such that
\beq\label{1510-n1}
{\rm Bott}(\Phi, U_1U_2U_3)|_{\cal P}=\sum_{i=1}^3{\rm Bott}(\Phi, U_i)|_{\cal P},
\eneq
provided that
$\|[\Phi, U_i]\|<\dt_1,\,\,\,i=1,2,3.$

Let $\ep_1=\min\{\dt_1/2, \ep/4\}$ and ${\cal F}_1={\cal F}\cup
{\cal G}_1.$ We may assume that ${\cal F}_1$ is in the unit ball
of $C.$ We may also assume that $[L']|_{\cal P}$ is well defined
for any $\ep_1$-${\cal F}_1$-multiplicative \morp\, from $C$ to
any unital C*-algebra.

Let $\dt_2>0$ and ${\cal G}\subset C$ be a finite subset and
${\cal P}_1\subset \underline{K}(C)$ be finite subset required by
\ref{Newstableuniq}  for $\ep_1/2$ and ${\cal F}_1.$ We may assume that
$\dt_2<\dt_1/2,$ ${\cal G}\supset {\cal F}_1$ and ${\cal
P}_1\supset {\cal P}.$ We also assume that ${\cal G}$ is in the
unit ball of $C.$

It follows from \ref{kkmaps} that there exists an integer $K_1\ge
1,$ a unital \hm\, $h_0': C\to M_{K_1}(\C)$ and a
$\dt_2/2$-${\cal G}$-multiplicative \morp\, $L_1:  C\to
M_{K_1+1}(M_{\phi_1, \phi_2})$ such that
\beq\label{VuV-1}
[L_1]|_{{\cal P}_1}=(\theta+[h_0'])|_{{\cal P}_1}.
\eneq

%Note that the point-evaluation at $0,$ $\pi_0$ from  $M_{\phi_1, \phi_2}$ has an image in
%$\phi_0(C).$ We also have, by viewing $\pi_0$ as a \hm\, from
%$M_{\phi_1,\phi_2}$ to $B,$
%\beq\label{VuV-2}
%[\pi_0]\circ \theta\circ [\psi]=[\psi]\,\,\,\text{in}\,\,\, KK(C,B).
%\eneq
%
%Moreover, by viewing $\pi_0$ as a map from $M_{\phi_1,\phi_2}$ to
%$B,$
%\beq\label{VuV-3}
%[\pi_0]\circ \theta\circ
%[\psi]=[\phi_1\circ\psi]\,\,\,\text{in}\,\,\, KK(C, B).
%\eneq
Note that  $[\pi_0]\circ \theta=[\phi_1]$ and $[\pi_1]\circ \theta =[\phi_2]$ and, for each $t\in (0,1),$
\beq\label{VuV-4}
[\pi_t]\circ \theta=[\phi_1].
\eneq
By \ref{Newstableuniq}, we obtain an integer $K_0,$ a unitary $V\in
U(M_{1+K_1+K_0}((C)))$ and a unital \hm\, $h: C\to
M_{K_0}(\C)$ such that
\beq\label{VuV-5--1}
{\rm ad}\, V\circ (\pi_e\circ L_1\oplus h)\approx_{\ep_1/2}
(\id\oplus h_0'\oplus h)\,\,\,\text{on}\,\,\,{\cal F}_1,
\eneq
where $\pi_e: M_{\phi_1,\phi_2} \to C$ is the canonical projection.

(Here and below, we will identify a homomorphism mapped to $M_k(\C)$ as a homomorphism mapped to $M_k(A)$ for any unital $C^*$ algebra $A$, without introducing a new notation.)

%By \ref{Newstableuniq}, we obtain an integer $K_0,$ a unitary $V_{00}\in
%U(M_{1+K_1+K_0}(\phi_1(C)))$ and a unital \hm\, $h: C\to
%M_{K_0}(\C)$ such that
%\beq\label{VuV-5--1}
%{\rm ad}\, V_{00}\circ (\pi_0\circ L_1\oplus h)\approx_{\ep_1/2}
%(\phi_1\oplus h_0'\oplus h)\,\,\,\text{on}\,\,\,{\cal F}_1.
%\eneq
Write $V_{00}=\phi_1(V)$ and
$V_{00}'=\phi_2(V).$ The assumption that
$[\phi_1]=[\phi_2]$ implies that $[V_{00}]=[V_{00}']$ in $K_1(B).$
By adding another $h$ in (\ref{VuV-5--1}), replacing $K_0$ by $2K_0$  and replacing $V$
by $V\oplus 1_{M_{K_0}},$ if necessary, we may assume that
$V_{00}$ and $V_{00}'$ are in the same component of
$U(M_{1+K_1+K_0}(B)).$

One obtains a continuous path of unitaries $\{Z(t): t\in [0,1]\}$
in $M_{1+K_1+K_0}(B)$ such that
\beq\label{VuV-5--2}
Z(0)=V_{00}\andeqn Z(1)=V_{00}'.
\eneq
It follows that $Z\in M_{1+K_1+K_0}(M_{\phi_1,\phi_2}).$ By
replacing $L_1$ by ${\rm ad}\, Z\circ (L_1\oplus h)$ and using a
new $h_0',$ we may assume that
\beq\label{VuV-5--3}
\pi_0\circ L_1\approx_{\ep_1/2}\phi_1\oplus
h_0'\,\,\,\text{on}\,\,\,{\cal F}_1\andeqn
%\eneq
%and that
%\beq
%\label{VuV-5--4}
\pi_1\circ L_1\approx_{\ep_1/2} \phi_2\oplus
h_0'\,\,\,\text{on}\,\,\, {\cal F}_1,\,\,\,i=1,2.
\eneq

Define $\lambda: C\to M_{1+K_1+K_0}(C)$ by $\lambda(c)={\rm diag}(c, h_0'(c)),$ where
we also identify $M_{K_0+K_1}(\C)$ with the scalar matrices in $M_{K_0+K_1}(C).$
In particular, since $\phi_i$ is  unital, $\phi_i\otimes {\rm id}_{M_{K_1+K_0}}$ is identity on  $M_{K_0+K_1}(\C),$
$i=1,2.$ {{Consequently $(\phi_i\otimes {\rm id}_{M_{K_0+K_1}})\circ h_0'=h_0'$.}}
Therefore, in this way, one may write
\begin{equation*}
\phi_i(c)\oplus h_0'(c)=(\phi_i\otimes {\rm id}_{M_{K_0+K_1{{+1}}}})\circ \lambda(c)\rforal c\in C.
\end{equation*}
%It follows from 10.1 of \cite{Lnamj} that we may also assume that
%\beq\label{VuV-5--4}
%\pi_1\circ L_1\approx_{\ep_1/2} \phi_2\oplus
%h_0'\,\,\,\text{on}\,\,\, {\cal F}_1,\,\,\,i=1,2.
%\eneq

There is a partition
%\beq\label{VuV-5}
$0=t_0<t_1<\cdots <t_n=1$
%\eneq
such that
\beq\label{VuV-6}
\pi_{t_{i}}\circ L_1\approx_{\dt_2/8}\pi_t\circ
L_1\,\,\,\text{on}\,\,\,{\cal G} \rforal t_i\le t\le t_{i+1},\,\,\,i=1,2,...,n-1.
\eneq
%for all $t_i\le t\le t_{i+1},$ $i=1,2,...,n-1.$
By applying \ref{Newstableuniq}  again, we
obtain an integer $K_2\ge 1, $
%K_4\ge N(K_1)$ with $K_1+K_4=1+K_3,$
a unital \hm\, $h_{00}: C\to M_{K_2}(\C),$
%and a unital \hm\, $h_{00}': C\to M_{K_4}(\C)$
 and  a unitary $V_{t_i}\in M_{1+K_0+K_1+K_2}(B)$ such that
\beq\label{VuV-7}
{\rm ad}\, V_{t_i}\circ (\phi_1\oplus h_0'\oplus
h_{00})\approx_{\ep_1/2} (\pi_{t_i}\circ L_1 \oplus
h_{00})\,\,\,\text{on}\,\,\, {\cal F}_1.
\eneq

Note that, by (\ref{VuV-6}), (\ref{VuV-7}) and (\ref{VuV-5--3}),
$$
\|[\phi_1\oplus h_0'\oplus h_{00}(a), V_{t_i}V_{t_{i+1}}^*]\|<\dt_2/4+\ep_1\rforal
a\in {\cal F}_1.
$$

Denote by  $\eta_{-1}=0$ and
$$
\eta_k=\sum_{i=0}^k {\rm Bott}(\phi_1\oplus h_0'\oplus h_{00}, V_{t_i}V_{t_{i+1}}^*)|_{\cal P},\,\,k=0,1,...,n-1.
$$

Now we will construct, for each $i,$ unital \hm\, $F_i: C\to M_{J_i}(\C)\subset M_{J_i}(B)$ and
a unitary $W_i\in M_{1+K_0+K_1+K_2+\sum_{k=1}^iJ_i}(B)$ such that
\beq\label{150110-L2}
\|[H_i,\, W_i]\|<\dt_2/4\andeqn {\rm Bott}(H_i,\, W_i)=\eta_{i-1},
\eneq
where $H_i=\phi_{{1}}\oplus h_0'\oplus h_{00}\oplus\oplus_{k=1}^i F_i,$ $i=1,2,...,n-1.$

Let $W_0=1_{M_{1+K_0+K_1+K_2}}.$
It follows from \ref{stablehomtp} that there is an integer $J_1\ge 1,$
a unital \hm\, $F_1: C\to M_{J_1}(\C)$ and a unitary $W_0\in
U(M_{1+K_0+K_1+K_2+J_1}(B))$ such that
\beq\label{VuV-8}
\|[H_1(a), \, W_1]\|<\dt_2/4\rforal a\in {\cal F}_1\andeqn
\text{Bott}(H_1, W_1)=\eta_0,
\eneq
where $H_1=\phi_1\oplus h_0'\oplus h_{00}\oplus F_1.$

Assume that, we have construct required $F_i$ and $W_i$ for $i=0,1,...,k<n-1.$
It follows from \ref{stablehomtp} that there is an integer $J_{k+1}\ge 1,$
a unital \hm\, $F_{k+1}: C\to M_{J_{k+1}}(\C)$ and a unitary $W_{k+1}\in
U(M_{1+K_0+K_1+K_2+\sum_{i=1}^{k+1}J_i}(B))$ such that
\beq\label{VuV-9}
\|[H_{k+1}(a), \, W_{k+1}]\|<\dt_2/{ 4}\rforal a\in {\cal F}_1\andeqn
\text{Bott}(H_{k+1}, W_{k+1})=\eta_{k}
\eneq
where $H_{k+1}=\phi_1\oplus h_0'\oplus h_{00}\oplus \oplus_{i=1}^{k+1}F_i.$

Now define $F_{00}=h_{00}\oplus \oplus_{i=0}^{n-1}F_i$ and define
$K_3=1+K_0+K_1+K_2+\sum_{i=1}^{n-1}J_i.$
Define
$$
v_{t_k}= {\rm diag}(W_k{\rm diag}(V_{t_k}, {\rm id}_{1_{M_{\sum_{i=1}^kJ_i}}}), 1_{M_{\sum_{i=k+1}^{n-1} J_i}}),
$$
$k=1,1,...,n-1$ and
$v_{t_0}=1_{M_{1+K_0+K_1+K_2+\sum_{i=1}^{n-1}J_i}}.$
Then
\beq\label{VuV-14}
&&{\rm ad}\, v_{t_i}\circ (\phi_1\oplus h_0'\oplus
F_{00})\approx_{\dt_2+\ep_1} \pi_{t_i}\circ (L_1\oplus
F_{00})\,\,\,\text{on}\,\,\,{\cal F}_1,\\
&&\|[\phi_1\oplus  h_0'\oplus F_{00}(a),\,v_{t_i}v_{t_{i+1}}^*]\|<\dt_2/2+2\ep_1\rforal a\in {\cal F}_1\andeqn\\
%\eneq
%\beq
\label{VuV-15}
%&&\|[\phi_1\oplus  h_0'\oplus F_{00}(a),\,V_iV_{i+1}^*]\|<\dt_2+2\ep_1\rforal a\in {\cal F}_1\andeqn
%&&\text{Bott}(\phi_1\oplus h_0'\oplus F_{00}, V_1)=0\andeqn
&&\text{Bott}(\phi_1\oplus h_0'\oplus F_{00}, v_{t_i}v_{t_{i+1}}^*)\\
&=&\text{Bott}(\phi_1', W_i')+
\text{Bott}(\phi_1', V_{t_i}'(V_{t_{i+1}}')^*) +\text{Bott}(\phi_1', (W_{i+1}')^*)\\
&=& \eta_{i-1}+\text{Bott}(\phi_1', V_{t_i}V_{t_{i+1}}^*)-\eta_{i}=0,
\eneq
where $\phi'_{{1}}=\phi_1\oplus h_0'\oplus F_{00},$   $W_i'={\rm diag}(W_i, 1_{M_{\sum_{j=i+1}^{n-1}J_i}})$
and $V_{t_i}'={\rm diag}(V_{t_i}, 1_{M_{\sum_{i=1}^{n-1}J_i}}),$
$i=0,1,2,...,n-2.$

It follows from \ref{STHOM} that there is an integer $N_1\ge 1,$
another unital \hm\, $F_0': C\to M_{N_1}(\C)$ and a continuous path of
unitaries $\{w_i(t): t\in [t_{i-1}, t_i]\}$ such that
\beq\label{VuV-17}
&&\hspace{-0.4in}
%w_1(0)=v_v_{t_1}, \,\,\,w_1(t_1)=1,\,\,\,
w_i(t_{i-1})=v_{i-1}'(v_{i}')^*, w_i(t_i)=1, \,\,\, i=1,2,...,n-1\andeqn\\
&&\|[\phi_1\oplus h_0'\oplus F_{00}\oplus F_{0}'(a), \,w_i(t)]\|<\ep/2\rforal a\in {\cal F},
\eneq
$i=1,2,...,n-1,$
where $v_i'={\rm diag}(v_i, 1_{M_{N_1}}(B)),$ $i=1,2,...,n-1.$
Define $V(t)=w_i(t)v_i'$ for $t\in [t_{i-1}, t_i],$ $i=1,2,...,n-1.$
Then $V(t)\in C([0,t_{n-1}], M_{N_1}(B)).$ Moreover,
\beq\label{VuV-18}
{\rm ad}\, V(t)\circ (\phi_1\oplus h_0'\oplus F_{00}\oplus
F_0')\approx_{\ep/2} \pi_t\circ L_1\oplus F_{00}\oplus F_0'
\,\,\,\text{on}\,\,\,{\cal F}.
\eneq

Define  $h_0=h_0'\oplus F_{00}\oplus F_0',$ $L=L_1\oplus
F_{00}+F_0'$ and $d=1-t_{n-1}.$ Then, by (\ref{VuV-18}),
(\ref{VuV2}) and (\ref{VuV3}) hold. From (\ref{VuV-5--3}),
(\ref{VuV33}) also holds.

\end{proof}

\section{Asymptotically unitary equivalence}

\begin{lem}\label{MUN2}
Let $C_1$ and {$A_1$} be two unital separable simple \CA s in {${\cal B}_1,$} let
$U_1$ and $U_2$ be two UHF-algebras of infinite type and
let $C=C_1\otimes U_1$ and $A=A_1\otimes U_2.$  Suppose that
%$C_1$ satisfies the UCT and
$\phi_1, \phi_2: C\to A$ are two unital monomorphisms.
Suppose also that
\beq\label{MUN2-1}
[\phi_1]=[\phi_2]\,\,\,{\rm in}\,\,\, KL(C,A),\\
(\phi_1)_T=(\phi_2)_T\andeqn \phi_1^{\ddag}=\phi_2^{\ddag}.
\eneq
Then $\phi_1$ and $\phi_2$ are approximately unitarily equivalent.
\end{lem}

\begin{proof}
This follows  immediately from \ref{MUN1}. Note that both $A$ and $C$ are 
%unital simple C*-algebras  
in $\mathcal B_1.$
\end{proof}

\begin{lem}\label{botgroup}
Let $B$ be a unital \CA\, and let $u_1, u_2,...,u_n\subset U(B)$ be unitaries. Suppose that $v_1, v_2,...,v_m\subset U(B)$ are also unitaries
such that $[v_j]\subset G,$ where $G$ is the subgroup of $K_1(B)$ generated by $[u_1], [u_2],...,[u_n].$ There exists  $\dt>0$ satisfying the following:
For any unital \CA\, $A$ and any unital \hm s $\phi_1, \phi_2: B\to A${{,}}
if there is a unitary $w\in U(B)$ such that
\beq\label{botgroup-1}
\| w^*\phi_1(u_i)w-\phi_2(u_i)\|<\dt,
\eneq
then there exists a group \hm\, $\af: G\to \Aff(T(A))$ such that
\beq\label{botgroup-2}
{1\over{2\pi i}}\tau(\log(\phi_2(v_j)w^*\phi_1(v_j^{ *})w)&=&\af([v_j]){{(\tau)}},
\eneq
{{any $\tau\in T(A),$}} $i=1,2,...,n$ and $j=1,2,...,m.$
\end{lem}

\begin{proof}
The proof is essentially contained in the proof of 6.1, 6.2 and 6.3 of \cite{Lnind}.
\end{proof}

\begin{lem}\label{inv71}
Let $C_1$ be a unital simple \CA\,  in {\rm \ref{RangT}}, let $A_1$ be a unital separable simple \CA\, in {${\cal B}_0,$} and let $U_1$ and $U_2$ be two UHF-algebras of infinite type. Let $C=C_1\otimes U_1$ and $A=A_1\otimes U_2.$ Suppose that $\phi_1, \phi_2: C\to A$ are two unital monomorphisms. Suppose also that
\beq\label{71-1}
&&[\phi_1]=[\phi_2]\,\,\,{\rm in}\,\,\, KL(C, A),\\\label{71-2}
&&\phi_1^{\ddag}=\phi_2^{\ddag}, \,\,\, (\phi_1)_T=(\phi_2)_T\tand\\\label{71-3}
&&R_{\phi, \psi}(K_1(M_{\phi_1,\phi_2})){ \subset} \rho_A(K_0(A)).
\eneq
Then, for any increasing sequence of finite subsets $\{{\cal F}_n\}$ of $C$ whose union is dense in $C,$ any increasing sequence of finite
subsets ${\cal P}_n$ of $K_1(C)$ with
$\cup_{n=1}^{\infty} {\cal P}_n=K_1(C)$ and any decreasing sequence of positive numbers $\{\dt_n\}$ with $\sum_{n=1}^{\infty} \dt_n<\infty,$ there exists a sequence of unitaries $\{u_n\}$ in $U(A)$ such that
\beq\label{71-4}
{\rm Ad}u_n\circ \phi_1\approx_{\dt_n} \phi_2\,\,\,{\rm on}\,\,\,{\cal F}_n\andeqn\\
\rho_A({\rm bott}_1(\phi_2, u_n^*u_{n+1})(x))=0\rforal x\in {\cal P}_n
\eneq
and for all sufficiently large $n.$

\end{lem}

\begin{proof}
Note that $A\cong A\otimes U_2.$ Moreover, {there is a unital \hm\, $s: A\otimes U_2\to A$ such that
$s\circ \imath$} is approximately unitarily equivalent to the identity map on $A,$ where
 ${\imath}: A\to A\otimes U_2$ defined by $a\to a\otimes 1_{U_2}$ {for all $a\in A$}.
 Therefore we
may assume that $\phi_1(C),\, \phi_2(C)\subset A\otimes 1_{U_2}.$
It follows from \ref{MUN2} that there exists a sequence of unitaries
$\{v_n\}\subset A$ such that
$$
\lim_{n\to\infty} {\rm Ad}\, v_n\circ \phi_1(c)=\phi_2(c)\tforal c\in C.
$$

%Without loss of generality, we may assume
%that ${\cal P}_n=\{z_1, z_2,...,z_n\}.$
We may assume that ${\cal F}_n$ are in the unit ball and $\cup_{n=1}^{\infty} {\cal F}_n$ is dense in the unit ball of $C.$

Put $\ep_n'=\min\{1/2^{n+1}, \dt_n/2\}.$ Let
$C_n\subset C$ be a unital \SCA\, (in place of ${ C_n}$) which
has finitely generated $K_i(C_n)$ ($i=0,1$),
%let ${\cal P}_n\subset \underline{K}(C_n)$ be a finite subset
and let ${\cal Q}_n$ be a finite set of generators of $K_1(C_n),$
%with finite subset ${\cal Q}_n\subset K_1(C_n)$ such that
%$[\imath_n]{\cal Q}_n={\cal P}_n,$ where
 %$\imath_n: C_n\to A$ is the embedding
 let $\dt_n'>0$ (in place of $\dt$) be as in \ref{Extbot2} for $C$ (in place of $A$), $\ep_n'$ (in
place of $\ep$), ${\cal F}_n$ (in place of ${\cal F}$) and $[\imath_n]({\cal Q}_{n-1})$
(in place of ${\cal P}$), where $\imath_n: C_n\to C$ is the embedding.
Note that, we assume that
\beq\label{71-9}
[\imath_{n+1}]({\cal Q}_{n+1})\supset {\cal P}_{n+1}\cup [\imath_n]({\cal Q}_n).
\eneq

%Suppose that the subgroup generated by ${\cal P}_n$
%is $G_n$ which can be written
Write $K_1(C_n)=G_{n,f}\oplus {\rm Tor}(K_1(C_n)),$ {where $G_{n,f}$ is a finitely generated free { abelian} group}.
Let $z_{1,n}, z_{2,n},...,z_{f(n),n}$ be the free generators of $G_{n,f}$ and $z'_{1,n},z'_{2,n},...,z_{t(n),n}'$ be generators of ${\rm Tor}(K_1(C_n)).$
We may assume that
$$
{\cal Q}_n=\{z_{1,n}, z_{2,n},...,z_{f(n),n},z'_{1,n},z'_{2,n},...,z_{t(n),n}'\}.
$$

Let $1/2>\ep_n''>0$ so that ${\text{bott}}_1(h', u')|_{K_1(C_n)}$ is well defined group \hm, ${\text{bott}}_1(h', u')|_{{\cal Q}_n}$ is well defined and
$(\text{bott}_1(h',u')|_{K_1(C_n)})|_{{\cal Q}_n}={\text{bott}}_1(h', u')|_{{\cal Q}_n}$ for any unital \hm\, $h': C\to A$ and any unitary $u'\in A$ for which
\beq\label{71-10}
\|[h'(c), u']\|<\ep_n''\tforal c\in {\cal G}_n
\eneq
for some finite subset ${\cal G}_n\subset C$ which contains ${\cal F}_n.$

Let $w_{1,n},w_{2,n},...,w_{f(n),n}, w_{1,n}',w_{2,n}',...,w_{t(n),n}'\in C$ be unitaries (note that $C$ has stable rank one)
such that $[w_{i,n}]=(\imath_n)_{*1}(z_{i,n})$ and
$[w_{j,n}']=(\imath_n)_{*1}(z_{j,n}'),$ $i=1,2,...,f(n),$
$j=1,2,...,t(n)$ and  $n=1,2,....$
To simplify notation, without loss of generality, we may assume that $w_{i,n}\in {\cal G}_n,$ $n=1,2,....$

Let $\dt_1''=1/2$ and, for $n\ge 2,$ let $\dt_n''>0$ (in place of $\dt$) be as in \ref{botgroup}  associated
with $w_{1,n},w_{2,n},...,w_{f(n),n}, w_{1,n}',w_{2,n}',...,w_{t(n),n}'$
(in place of $u_1, u_2,..,u_n$) and
$$
\{w_{1,{n-1}}, w_{2,n-1},...,w_{f(n-1), n-1}, w_{1,n-1}',w_{2,n-1}',...,w_{t(n-1), n-1}'\}
$$
(in place of $v_1,v_2,...,v_m$).

Put $\ep_n=\min\{\ep_n''/2, \ep_n'/2,\dt'_n, \dt_n''/2\}.$ We may assume that
\beq\label{71-11}
{\rm Ad}\, v_n\circ \phi_1\approx_{\ep_n} \phi_2\,\,\,{\rm on\,\,\,} {\cal G}_n,\,\,\, n=1,2,....
\eneq
Thus ${\rm bott}_1(\phi_2\circ \imath_n, v_n^*v_{n+1})$ is well defined.
Since $\Aff(T(A))$ is torsion free,
\beq\label{71-12-1}
\tau\big({\rm bott}_1(\phi_2\circ \imath_n, v_n^*v_{n+1})|_{{\rm Tor}(K_1(C_n))}\big)=0.
\eneq

We have
\beq\label{71-12}
\|\phi_2(w_{j,n}){\rm Ad}\, v_n(\phi_1(w_{j,n})^*)-1\|<(1/4)\sin(2\pi \ep_n)<\ep_n,\,\,\, n=1,2,....
\eneq
%$n=1,2,....$
%It follows from \ref{botgroup} that there exists a group
%\hm\, $\gamma_n: G_n\to \tau
Define
\beq\label{71-13}
h_{j,n}={1\over{2\pi i}}\log(\phi_2(w_{j,n}){\rm Ad}\, v_n(\phi_1(w_{j,n})^*)),\,\,\,j=1,2,...,f(n), n=1,2,...
\eneq
Then, for any $\tau\in T(A),$
%\begin{equation*}
$|\tau(h_{j,n})|<\ep_n<\dt_n',\,\,\, j=1,2,...,f(n), n=1,2,....$
%\end{equation*}
Since $\Aff(T(A))$ is torsion free, it follows from
\ref{botgroup} that
\beq\label{71-13+}
\tau({1\over{2\pi i}}\log(\phi_2(w_{j,n}'){\rm Ad}\, v_n(\phi_1(w_{j,n}^{'*}))))=0,
\eneq
$j=1,2,...,t(n)$ and $n=1,2,....$
By the assumption that
$R_{\phi_1, \phi_2}({ K_1}(M_{\phi_1,\phi_2})){ \subset} \rho_A(K_0(A)),$ by
the Exel's formula (see \cite{HL}) and by Lemma 3.5 of \cite{Linajm}, we conclude that
\begin{equation*}
\widehat{h_{j,n}}(\tau)=\tau(h_{j,n})\in R_{\phi_1, \phi_2}({ K_1}(M_{\phi_1,\phi_2})){ \subset} \rho_A(K_0(A)).
\end{equation*}
Now define $\af_n': K_1(C_n)\to \rho_A(K_0(A))$ by
\beq\label{71-15}
\af_n'(z_{j,n})(\tau)=\widehat{h_{j,n}}(\tau)=\tau(h_{j,n}),\,\,\, j=1,2,...,f(n)\andeqn
\af_n'(z_{j,n}')=0,\,\,\,j=1,2,...,t(n),
\eneq
$n=1,2,....$
 Since $\af_n'(K_1(C_n))$ is free,
%Since $\Aff(T(A))$ is divisible {\color{Green} (it is not clear to me how to use the divisibility of $\Aff(T(A))$; seems to me, what used here is the projectivity of the image of $K_1(C_n)$---it is isomorphic to $\mathbb Z^k$ for some $k$)},
it follows that there is a \hm\, $\af_n^{(1)}: K_1(C_n)\to K_0(A)$ such that
\beq\label{71-16}
{{(}}\rho_A\circ \af_n^{(1)}(z_{j,n}))(\tau)=\tau(h_{j,n}),\,\,\, j=1,2,...,f(n),\,\,{{\tau\in T(A)}} \,\,\andeqn\\
\af_{n}^{(1)}(z_{j,n}')=0,\,\,\, j=1,2,...,t(n).
\eneq
Define
$\af_n^{(0)}: K_0(C_n)\to K_1(A)$ by $\af_n^{(0)}=0.$ By the UCT,
there is $\kappa_n\in KL({ S}C_n, A)$ such that $\kappa_n|_{K_i(C_n)}=\af_n^{(i)},$ $i=0,1,$
where $SC_n$ is the suspension of $C_n$ (here, we also {identify} $K_i(C_n)$ with $K_{i+1}(SC_n)$).

By the UCT again, there is $\af_n\in KL(C_n\otimes C(\T), A)$ such that
$\af_n\circ \boldsymbol{\bt}|_{\underline{K}(C_n)}=\kappa_n.$ In particular,
$\af_n\circ \boldsymbol{\bt}|_{K_1(C_n)}=\af_n^{(1)}.$ It follows from \ref{Extbot2} that there exists a unitary $U_n\in U_0(A)$ such that
\beq\label{71-17}
\|[\phi_{ 2}(c),\, U_n]\|<\ep_n''\rforal c\in {\cal F}_n\andeqn\\\label{71-17+}
\rho_A({\text{bott}_1}(\phi_2,\, U_n)(z_{j,n}))=-\rho_A\circ \af_n(z_{j,n}),
\eneq
$j=1,2,...,f(n).$ We also have
\beq\label{71-18}
\rho_A({\text{bott}_1}(\phi_2,\, U_n)(z_{j,n}'))=0,\,\,\, j=1,2,...,t(n).
\eneq
By the Exel trace formula {{(see \cite{HL})}}, {{(\ref{71-16}) and (\ref{71-17+}),}}
 we have
\beq\label{71-19}
\tau(h_{j,n})&=&-\rho_A({\rm bott}_1(\phi_2, U_n)(z_{j,n})(\tau)
%\\\label{71-20}
    =-\tau({1\over{2\pi i}}\log(U_n\phi_2(w_{j,n})U_n^*\phi_2(w_{j,n}^*)))
\eneq
for all $\tau\in T(A),$ $j=1,2,...,f(n).$
Define $u_n=v_nU_n,$ $n=1,2,....$ By 6.1 of \cite{Lnind}, (\ref{71-19}) and
(\ref{71-17+}), we compute that
\beq\label{71-21}
&&\hspace{-0.2in}\tau({1\over{2\pi i}}\log(\phi_2(w_{j,n}){\rm Ad} u_n(\phi_{ 1}(w_{j,n}^*)))))\\
&=&\tau({1\over{2\pi i}}\log(U_n\phi_2(w_{j,n})U_n^*v_n^*\phi_1(w_{j,n}^*)v_n)))\\
&=&\tau({1\over{2\pi i}}\log(U_n\phi_2(w_{j,n})U_n^*\phi_2(w_{j,n}^*)\phi_2(w_{j,n})v_n^*
\phi_1(w_{j{{,n}}}^*)v_n)))\\
&=&\tau({1\over{2\pi i}}\log(U_n\phi_2(w_{j,n})U_n^*\phi_2(w_{j,n}^*))))
+\tau({1\over{2\pi i}}\log(\phi_2(w_{j,n})v_n^*\phi_1(w_{j,n}^*)v_n)))\\
%\\%\label{71-22}
&=&\rho_A({\rm bott}_1(\phi_2, U_n)(z_{j,n}))(\tau)+\tau(h_{j,n})=0
\eneq
for all $\tau\in T(A),$ $j=1,2,...,f(n)$ and $n=1,2,....$
By (\ref{71-13+}) and (\ref{71-18}),
\begin{equation}\label{71-22}
\tau({1\over{2\pi i}}\log(\phi_2(w_{j,n}'){\rm Ad} u_n(\phi_{1}((w_{j,n}')^*))))=0,
\end{equation}
$j=1,2,...,t(n)$ and $n=1,2,....$
Let
\beq\label{71-23}
b_{j,n}&=&{1\over{2\pi i}}\log(u_n\phi_2(w_{j,n})u_n^*\phi_1(w_{j,n}^*)),\\
b_{j,n}'&=&{1\over{2\pi i}}\log(\phi_2(w_{j,n})u_n^*u_{n+1}\phi_2(w_{j,n}^*)u_{n+1}^*u_n)\andeqn\\
b_{j,n+1}''&=&{1\over{2\pi i}}\log(u_{n+1}\phi_2(w_{j,n})u_{n+1}^*\phi_1(w_{j,n}^*)).
\eneq
$j=1,2,...,f(n)$ and $n=1,2,....$ We have, by (\ref{71-22}),
\beq\label{71-24}
\tau(b_{j,n})&=&\tau({1\over{2\pi i}}\log(u_n\phi_2(w_{j,n})u_n^*\phi_1(w_{j,n})))\\
%&=&\tau({1\over{2\pi i}}\log(u_n^*u_n\phi_2(w_{j,n})u_n^*\phi_1(w_{j,n}^*)u_n))\\
&=&\tau({1\over{2\pi i}}\log(\phi_2(w_{j,n}^*)u_n^*\phi_1(w_{j,n}^*)u_n))=0
\eneq
for all $\tau\in T(A),$ $j=1,2,..., f(n)$ and $n=1,2,....$ Note that
$\tau(b_{j, n+1})=0$ for all $\tau\in T(A),$ $j=1,2,...,f(n+1).$
It follows from \ref{botgroup} {and (\ref{71-9})} that
\begin{equation*}
\tau(b_{j,n+1}'')=0\tforal \tau\in T(A),\,\,\,j=1,2,...,f(n),\,\,\,n=1,2,....
\end{equation*}
Note  also
that
\begin{equation*}
u_ne^{2\pi i b_{j,n}'}u_n^*=e^{2\pi i b_{j,n}}\cdot e^{-2\pi i b_{j,n+1}''},\,\,\,j=1,2,...,f(n).
\end{equation*}
Thus, by 6.1 of \cite{Lnind}, we compute that
\beq\label{71-26}
\tau(b_{j, n}')=\tau(b_{j,n})-\tau(b_{j,n+1}'')=0\tforal \tau\in T(A).
\eneq
By the Exel formula {{(see \cite{HL})}} and (\ref{71-26}),
\beq\label{71-27}
&&\hspace{-0.3in}\rho_A({\rm bott}_1(\phi_2, u_n^*u_{n+1}))(w_{j,n}^*)(\tau)
=\tau({1\over{2\pi i}}\log(u_n^*u_{n+1}\phi_2(w_{j,n})u_{n+1}^*u_n\phi_2(w_{j,n}^*)))\\
&=&\tau({1\over{2\pi i}}\log(\phi_2(w_{j,n})u_{n}^*u_{n+1}\phi_2(w_{j,n}^*)u_{n+1}^*u_n))=0
\eneq
for all $\tau\in T(A)$ and $j=1,2,...,f(n).$
Thus
\beq\label{71-28}
\rho_A({\rm bott}_1(\phi_2,u_n^*u_{n+1})(w_{j,n})(\tau)=0\tforal \tau\in T(A),
\eneq
$j=1,2,...,f(n)$ and $n=1,2,....$ We also have
\beq\label{71-29}
\rho_A({\rm bott}_1(\phi_2, u_n^*u_{n+1})(w_{j,n}')(\tau)=0\tforal \tau\in T(A),
\eneq
$j=1,2,...,f(n)$ and $n=1,2,....$ By \ref{botgroup}, we have that
\beq\label{71-30}
\rho_A({\rm bott}_1(\phi_2, u_n^*u_{n+1})(z)=0\tforal z\in {\cal P}_n,
\eneq
$n=1,2,....$
\end{proof}

\begin{thm}\label{Tm72}
Let $C_1$ be a unital simple \CA\, as in {\rm \ref{RangT}},  let $A_1$ be a unital separable simple \CA\, in ${{\cal B}_0},$ let $C=C_1\otimes U_1$ and let $A=A_1\otimes U_2,$ where $U_1$ and $U_2$ are UHF-algebras of infinite type. Suppose that $\phi_1,\, \phi_2: C\to A$ are two unital monomorphisms.
Then they are asymptotically unitarily equivalent if and only if
\beq\label{Tm72-1}
[\phi_1]=[\phi_2]\,\,\,{\rm in}\,\,\, KK(C,A),\\
\phi^{\ddag}=\psi^{\ddag},\,\,\,(\phi_1)_T=(\phi_2)_T \andeqn {\overline{R_{\phi_1,\phi_2}}}=0.
\eneq
\end{thm}

\begin{proof}
We will prove the ``if " part only. The ``only if" part follows
from 4.3 of \cite{Lnclasn}.
{Note $C=C_1\otimes U_1$ can be also regarded as a \CA\, in \ref{RangT}.}
Let $C=\lim_{n\to\infty}(C_n, \imath_n)$ be as in \ref{RangT}, where $\imath_n: C_n\to C_{n+1}$ is injective \hm. Let ${\cal F}_n\subset C$ be
an increasing sequence of subsets of $C$ such that
${\cup_{n=1}^{\infty}{\cal F}_n}$ is dense in $C.$
Put
$$
M_{\phi_1, \phi_2}=\{(f, {  c})\in C([0,1], A){ \oplus C}: f(0)=\phi_1(c)\andeqn f(1)=\phi_2(c)\}.
$$
Since $C$ satisfies the UCT, the assumption that $[\phi_1]=[\phi_2]$ in
$KK(C,A)$ implies that the following exact sequence splits:
\beq\label{72-1}
0\to \underline{K}(SA)\to \underline{K}(M_{\phi_1, \phi_2}) \rightleftharpoons_{\theta}^{\pi_e} \underline{K}(C)\to 0
\eneq
for some $\theta\in Hom(\underline{K}(C), \underline{K}(A)),$ where
$\pi_e: M_{\phi_1,\phi_2}\to C$ is { the projection  to $C$ defined in \ref{Mappingtorus}}.  Furthermore, since $\tau\circ \phi_1=\tau\circ \phi_2$ for all $\tau\in T(A)$ and ${\overline{R_{\phi_1,\phi_2}}}=0,$ we may also assume that
\beq\label{72-2}
R_{\phi_1,\phi_2}(\theta(x))=0\rforal x\in K_1(C).
\eneq
By \cite{DL}, we have that
\beq\label{72-3}
\lim_{n\to\infty}(\underline{K}(C_n),[\imath_n])=\underline{K}(C).
\eneq
Since $K_i(C_n)$ is finitely generated, there exists $K(n)\ge 1$ such that
\beq\label{72-3+}
Hom_{\Lambda}(F_{K(n)}\underline{K}(C_n),\, F_{K(n)}\underline{K}(A))=Hom_{\Lambda}(\underline{K}(C_n), \underline{K}(A))
\eneq
(see also \cite{DL} for the notation $F_m$ there).
%{\red{ can we recall where the notation $F_{K(n)}$ appeared before.}}

Let $\dt_n'>0$ (in place of $\dt$), $\sigma_n'>0$ (in place of $\sigma$), ${\cal G}_n'\subset C$ (in place of ${\cal G}$), \linebreak ${\{p_{1,n}', p_{2,n}',...,p'_{I(n),n)}, q_{1,n}', q_{2,n}',...,q_{I(n),n}'\}}$ (in place of
$\{p_1, p_2,...,p_k, q_1, q_2,...,q_k\}),$ ${\cal P}_n'\subset  {\underline{K}(C)}$
(in place of ${\cal P}$) corresponding to $1/2^{n+2}$ (in place of $\ep$) and ${\cal F}_n$ (in place of ${\cal F}$) as required by \ref{BHfull} (see also \ref{RBHfull}).
Note that, by the choice as in \ref{BHfull}, we assume that $G_{u,n}',$ the subgroup generated by
${\{[p'_{i,n}]-[q'_{i,n}]: 1\le i\le I(n)\}}$ is free.

Without loss of generality, we may assume that
${\cal G}_n'\subset \imath_{n, \infty}({\cal G}_n)$ and ${\cal P}_n'\subset [\imath_{n, \infty}]({\cal P}_n)$ for some finite subset ${\cal G}_n\subset C_n,$ and for some finite subset ${\cal P}_n\subset \underline{K}(C_n),$ we may assume that $p_{i,n}'=\imath_{n, \infty}(p_{i,n})$ and $q_{i,n}'=\imath_{n,\infty}(q_{i,n})$
for some projections $p_{i,n}, q_{i,n}\in C_n,$ $i=1,2,...,I(n).$
{\Wlog\, we may assume that the subgroup $G_{n,u}$  generated by $\{[p_{i,n}]-[q_{i,n}]: 1\le i\le I(n)\}$
is free and $p_{i,n}, q_{i,n}\in {\cal G}_n,$ $n=1,2,...,I(n).$}

 We may assume that ${\cal P}_n$ contains a set of generators
of $F_{K(n)}\underline{K}(C_n),$ ${\cal F}_n\subset {\cal G}_n'$ and
$\dt_n'<1/2^{n+3}.$
We may also assume that
${\rm Bott}(h', u')|_{{\cal P}_n}$ is well defined whenever $\|[h'(a),\, u']\|<\dt_n'$ for all $a\in {\cal G}_n'$ and for any unital \hm\, $h'$ from $C_n$ and a unitary $u'$ in the target algebra. Note that ${\rm Bott}(h',\,u')|_{{\cal P}_n}$ defines ${\rm Bott}(h'\, u').$
%Moreover,
%we may also assume that there exists a \hm\, $\Lambda: G_{u,n}\to U(A)/CU(A)$ such that
%\beq\label{72-n201410}
%\hspace{-0.2in}{\rm dist}(\overline{\langle ((h'(1-p_{i,n})+h'(p_{i,n})u')(h'(1-q_{i,n})+h'(q_{i,n})(u')^*)\rangle}, %\Lambda(x_{i,n}))<\sigma_n'/32,
%\eneq
%{\red{did we choose $\{x_{i,n}\}$ be the set of free generators of $G_{u,n}$ and can we just define
%$$\Lambda(x_{i,n}=\overline{\langle ((h'(1-p_{i,n})+h'(p_{i,n})u')(h'(1-q_{i,n})+h'(q_{i,n})(u')^*)\rangle}.$$}}
% where $x_{i,n}=[p_{i,n}]-[q_{i,n}],$ $i=1,2,...,I(n).$
We further assume that
\beq\label{72-4}
{\rm Bott}(h,\, u)|_{{\cal P}_n}={\rm Bott}(h',u)|_{{\cal P}_n}
\eneq
provided that $h\approx_{\dt_n'}h'$ on ${\cal G}_n'.$
%{Moreover, if $h'=h''\circ \imath_{n, \infty}$ for some \hm\, from $C,$ we may assume
%that $\Lambda$ factors through $G_{u,n}'.$ }
We may also assume
that $\dt_n'$ is smaller than {$\dt/16$} for the $\dt$ defined in 2.15 of \cite{Lnclasn} for $C_n$ (in place of $A$) and ${\cal P}_n$ (in place of ${\cal P}$).
{Let $k(n)\ge n$ (in place of $n$), $\eta_n'>0$ (in place of $\dt$) and ${\cal Q}_{k(n)}\subset K_1(C_{k(n)})$ be required by \ref{Extbot3}
for
$\dt_{k(n)}'/4$ (in place of $\ep$),
$\imath_{n, \infty}({\cal G}_{k(n)})$ (in place of ${\cal F}$),  ${\cal P}_{k(n)}$ (in place of ${\cal P}$) and
%Let $k(n)''\ge n$ (in place of $n$) and $\eta_n''>0$ (in place of $\dt$) be required by \ref{Extbot3} for $\imath_{n, \infty}({\cal G}_{k(n)''})$ (in place of ${\cal F}$),
$\{p_{i,n}, q_{i,n},: i=1,2,...,k(n)\}$ (in place of $\{p_i, q_i: i=1,2,...,k\}$),
%${\cal P}_{k(n)}'$ (in place of ${\cal P}$),
%{ $\dt_{k(n)}'/4$ (in place of $\ep$)
{and $\sigma_{k(n)}'/16$ (in place of $\sigma$)}.
%By changing the notation, we may assume that $k(n)=k(n)'=k(n)''.$
%Put $\eta_n=\min\{\eta_n', \eta_n''\}.$
 We may assume that ${\cal Q}_{k(n)}$ form a generator set of $K_1(C_{k(n)}).$
Since ${\cal P}$ generates $F_{K(n)}\underline{K}(C_{k(n+1)}),$ we may assume that ${\cal Q}_n\subset {\cal P}_{k(n)}.$}

For $C_n,$ since $K_i(C_n)$ ($i=0,1$) is finitely generated, by
(\ref{72-3+}), we may further assume that
$[\imath_{k(n), \infty}]$ is injective on $[\imath_{n, k(n)}](\underline{K}(C_n)),$ $n=1,2,....$ By passing to a subsequence, to simplify notation, we may also assume that
$k(n)=n+1.$ Let $\dt_n=\min\{\eta_n, \sigma_n',\dt_n'/2\}.$ By \ref{inv71}, there are unitaries $v_n\in U(A)$ such that
\beq\label{72-5}
&&{\rm Ad}\, v_n\circ \phi_1\approx_{\dt_{n+1}/4} \phi_2\,\,\,{\rm on}\,\,\,
\imath_{n, \infty}({\cal G}_{n+1}),\\
&&\rho_A({\rm bott}_1(\phi_2, v_n^*v_{n+1}))(x)=0
\,\,\,\rforal x\in (\imath_{n, \infty})_{*1}(K_1(C_{n+1}))\andeqn\\
&&\|[\phi_2(c),\, v_n^*v_{n+1}]\|<\dt_{n+1}/2\,\,\rforal a\in \imath_{n, \infty}({\cal G}_{n+1})
\eneq
(note that $K_1(C_{n+1})$ is finitely generated).
Note that, by (\ref{72-4}),
we may also assume that
\beq\label{72-7}
{\rm Bott}(\phi_1,\, v_{n+1}v_n^*)|_{[\imath_{n, \infty}]({\cal P}_n)}&=&{\rm Bott}(v_n^*\phi_1v_n, \, v_n^*v_{n+1})|_{[\imath_{n,\infty}]({\cal P}_n)}\\
&=&{\rm Bott}(\phi_2,\, v_n^*v_{n+1}])|_{[\imath_{n, \infty}]({\cal P}_n)}.
\eneq
In particular,
\beq\label{72-8}
{\rm bott}_1(v_n^*\phi_1v_n,\, v_n^*v_{n+1})(x)={\rm bott}_1(\phi_2,\, v_n^*v_{n+1})(x)
\eneq
for all $x\in (\imath_{n, \infty})_{*1}(K_1(C_{n+1})).$

By applying 10.4 and 10.5 of \cite{Linajm}, without loss of generality,
we may assume that the pair $\phi_1$ and $v_n$ defines
an element $\gamma_n\in {\rm Hom}_{\Lambda}(\underline{K}(C_{n+1}),\underline{K}(M_{\phi_1, \phi_2}))$ and $[\pi_0]\circ \gamma_n=[{\rm id}_{C_{n+1}}]$ (see 10.4 and 10.5 of \cite{Linajm} for the definition of $\gamma_n$). Moreover, by 10.4 and 10.5 of \cite{Linajm}, we may assume, without loss of {generality}, that
\beq\label{72-9}
|\tau(\log(\phi_2\circ \imath_{n, \infty}(z_j^*){\tilde v_n}\phi_1\circ \imath_{n, \infty}(z_j){\tilde v_n}))|<\dt_{n+1},
\eneq
$j=1,2,...,r(n),$ where $\{z_1,z_2,...,z_{r(n)}\} \subset U(M_k(C_{n+1}))$
which forms a set of generators of $K_1(C_{n+1})$ and where
${\tilde v_n}={\rm diag}(\overbrace{v_n, v_n,...,v_n}^k).$ {We may assume that
$z_j\in {\cal Q}_n\subset {\cal P}_n,$ $j=1,2,...,r(n).$}

Let $H_n=[\imath_{n+1}](\underline{K}(C_{n+1})).$ Since $\bigcup_{n=1}[\imath_{n+1, \infty}](\underline{K}(C_n))=\underline{K}(C)$
and $[\pi_0]\circ \gamma_n=[{{\rm id}}_{C_{n+1}}],$ we conclude that
\beq\label{72-10}
\underline{K}(M_{\phi_1, \phi_2})=\underline{K}(SA)+\cup_{n=1}^{\infty}\gamma_n(H_n).
\eneq
Thus, by passing to a subsequence, we may further assume that
\beq\label{72-11}
\gamma_{n+1}(H_n)\subset \underline{K}(SA)+\gamma_{n+2}(H_{n+1}),\,\,\,n=1,2,....
\eneq

By identifying $H_n$ with $\gamma_{n+1}(H_n),$ we may write $j_n: \underline{K}(SA)\oplus H_n\to \underline{K}(SA)\oplus H_{n+1}$ {for the inclusion in  (\ref{72-11}).} By
(\ref{72-10}), the inductive limit is $\underline{K}(M_{\phi_1,\phi_2}).$
From the definition of $\gamma_n,$ we note that $\gamma_n-\gamma_{n+1}\circ [\imath_{n+1}]$ maps $\underline{K}(C_{n+1})$ into $\underline{K}(SA).$
By 10.6 of \cite{Linajm},
$$
\Gamma({\rm Bott}(\phi_1,\, v_nv_{n+1}^*))|_{H_n}=(\gamma_{n+1}-\gamma_{n+2}\circ [\imath_{n+2}])|_{H_n}
$$
%{\red{should we recall what is $\Gamma$ here}}
(see 10.4, 10.5 and 10.6  of \cite{Linajm} for the definition of $\Gamma$) gives a \hm\, $\xi_n: H_n\to \underline{K}(SA).$ Put $\zeta_n=\gamma_{n+1}|_{H_n}.$ Then
\beq\label{72-13}
j_n(x,y)=(x+\xi_n(y),[\imath_{n+2}](y))
\eneq
for all $(x,y)\in \underline{K}(SA)\oplus H_n.$ Thus we obtain the following diagram:
\beq\label{72-14}
\begin{array}{ccccccc}
0 \to  & \underline{K}(SA)  &\to & \underline{K}(SA)\oplus H_n &\to
& H_n &\to 0\\\nonumber
 &\| & &\hspace{0.4in}\| \hspace{0.15in}\swarrow_{\xi_n} \hspace{0.05in}\downarrow_{[\imath_{n+2,\infty}]} &&
 \hspace{0.2in}\downarrow_{[\imath_{n+2,\infty}]} &\\\label{72-15}
 0 \to  & \underline{K}(SA)  &\to & \underline{K}(SA)\oplus H_{n+1} &\to & H_{n+1} &\to 0\\\nonumber
 &\| & &\hspace{0.4in}\| \hspace{0.1in}\swarrow_{\xi_{n+1}}\downarrow_{[\imath_{n+3,\infty}]} &&
\hspace{0.2in} \downarrow_{[\imath_{n+3,\infty}]} &\\
 0 \to  & \underline{K}(SA)  &\to & \underline{K}(SA)\oplus H_{n+2} &\to & H_{n+2} &\to 0\\
 \end{array}
\eneq
By the assumption that ${\bar R}_{\phi_1, \phi_2}=0,$ map $\theta$ also gives the following decomposition:
\beq\label{72-16}
{\rm ker}R_{\phi_1,\phi_2}={\rm ker}\rho_A\oplus K_1(C).
\eneq
Define $\theta_n=\theta\circ [\imath_{n+2, \infty}]$ and $\kappa_n=\zeta_n-\theta_n.$
Note that
\beq\label{72-17}
\theta_n=\theta_{n+1}\circ [\imath_{n+2}].
\eneq
We also have that
\beq\label{72-18}
\zeta_n-\zeta_{n+1}\circ [\imath_{n+2}]=\xi_n.
\eneq
Since $[\pi_0]\circ (\zeta_n-\theta_n)|_{H_n}=0,$ $\kappa_n$ maps
$H_n$ into $\underline{K}(SA).$ It follows that
\beq\label{72-19}
\kappa_n-\kappa_{n+1}\circ [\imath_{n+2}] &=&
\zeta_n-\theta_n-\zeta_{n+1}\circ [ \imath_{n+2}]+\theta_{n+1}\circ
[\imath_{n+2}]\\{\label{72-19+}}
&=&\zeta_n-\zeta_{n+1}\circ[\imath_{n+2}]=\xi_n{{.}}
\eneq
It follows from \ref{VuV} that there are integers
$N_1\ge 1,$  a ${\dt_{n+1}\over{4}}$-$\imath_{n+1}({\cal
G}_{n+1})$-multiplicative \morp\, $L_n: \imath_{n,
\infty}(C_{n+1})\to M_{1+N_1}(M_{\phi_1,\phi_2}),$ a unital \hm\,
$h_0: \imath_{n+1, \infty}(C_{n+1})\to M_{N_1}(\C),$
%unital \hm\, $h_0'': C\to M_{N_1+N_2-1}(\C),$
and a continuous path of unitaries $\{V_n(t): t\in [0,3/4]\}$ of
$M_{1+N_1}(A)$ such that $[L_n]|_{{\cal P}_{n+1}'}$ is well
defined, $V_n(0)=1_{M_{1+N_1}(A)},$
\begin{equation*}
[L_n\circ \imath_{n,\infty}]|_{{\cal P}_n}=(\theta\circ
[\imath_{n+1,\infty}]+[h_0\circ \imath_{n+1,\infty}])|_{{\cal P}_n},
\end{equation*}
and
\begin{equation*}
\pi_t\circ L_n\circ
\imath_{n+1,\infty}\approx_{\dt_{n+1}/4} {\rm ad}\, V_n(t)\circ
((\phi_1\circ \imath_{n+1,\infty})\oplus (h_0\circ
\psi_{n+1,\infty}))
\end{equation*}
on $\imath_{n+1,\infty}({\cal G}_{n+1})$ for all $t\in (0,3/4],$
\begin{equation*}
\pi_t\circ L_n\circ
\imath_{n+1,\infty}\approx_{\dt_{n+1}/4} {\rm ad}\, V_n(3/4)\circ
((\phi_1\circ \imath_{n+1,\infty})\oplus (h_0\circ
\imath_{n+1,\infty}))
\end{equation*}
on $\imath_{n+1,\infty}({\cal G}_{n+1})$ for all $t\in (3/4,1),$ and
\begin{equation*}
\pi_1\circ L_n\circ
\imath_{n+1,\infty}\approx_{\dt_{n+1}/4}\phi_2\circ
\imath_{n+1,\infty}\oplus h_0\circ \imath_{n+1,\infty}
\end{equation*}
on $\imath_{n+1,\infty}({\cal G}_{n+1}),$ where $\pi_t:
M_{\phi_1,\phi_2}\to A$ is the point-evaluation at $t\in (0,1).$

Note that $R_{\phi_1,\phi_2}(\theta(x))=0$ for all $x\in
\imath_{n+1,\infty}(K_1(C_{n+1})).$  As computed in 10.4 of
\cite{Linajm},
\beq\label{72-24}
\tau(\log((\phi_2(x)\oplus h_0(x)^*V_n(3/4)^*(\phi_1 (x)\oplus
h_0(x))V_n(3/4)))=0
\eneq
for $x=\imath_{n+1, \infty}(y),$ where $y$ is in a set of generators
of $K_1(C_{n+1})$ and for all $\tau\in T(A).$

Define $W_n'={\rm daig}(v_n,1)\in M_{1+N_1}(A).$ Then
$\text{Bott}((\phi_1\oplus h_0)\circ \imath_{n+1, \infty},\,
W_n'(V_n(3/4)^*)$ defines a \hm\, ${\tilde \kappa}_n\in
{\rm Hom}_{\Lambda}(\underline{K}(C_{n+1}),\underline{K}(SA)).$ By
(\ref{72-9})
\beq\label{72-25}
|\tau(\log((\phi_2\oplus h_0)\circ \imath_{n+1,\infty}(z_j)^*(W_n')^*(\phi_1\oplus h_0)\circ \imath_{n+1,\infty}(z_j)W_n'))|<\dt_{n+1},
\eneq
$j=1,2,...,r(n).$
Put ${\tilde V}_n={\rm diag}(V_n(3/4),1).$

Let
\beq\label{72-26}
&&\hspace{-0.7in}b_{j,n}={1\over{2\pi i}}\log({\tilde V}_n^*(\phi_1\oplus h_0)\imath_{n+1,\infty}(z_j){\tilde V}_n (\phi_2\oplus h_0)\circ \imath_{n+1, \infty}(z_j)^*),\\
&&\hspace{-0.7in}b_{j,n}'={1\over{2\pi i}}\log((\phi_1\oplus h_0)\circ \imath_{n+1,\infty}(z_j){\tilde V}_n(W_n')^*(\phi_1\oplus h_0)\circ \imath_{n+1, \infty}(z_j)^*W_n'{\tilde V}_n^*)\andeqn\\
&&\hspace{-0.7in}b_{j,n}''={1\over{2\pi i}}\log((\phi_2\oplus h_0)\imath_{n+1,\infty}(z_j)(W_n')^* (\phi_1\oplus h_0)\circ \imath_{n+1, \infty}(z_j)^*W_n'),
\eneq
$j=1,2,...,r(n).$
By (\ref{72-24}) and (\ref{72-25}),
\beq\label{72-27}
\tau(b_{j,n})=0 \andeqn |\tau(b_{j,n}'')|<\dt_{n+1}
\eneq
for all $\tau\in T(A).$
Note that
\beq\label{72-28}
{\tilde V}_n^*e^{2\pi i b_{j,n}'}{\tilde V}_n=
e^{2\pi i b_{j,n}} e^{2\pi i b_{j,n}''}{{.}}
\eneq
Then, by 6.1 of \cite{Lnind} and  by (\ref{72-27})
\beq\label{72-29}
\tau(b_{j,n}')&=&\tau(b_{j,n})-\tau(b_{j,n}'')
=\tau(b_{j,n}'')
\eneq
for all $\tau\in T(A).$
It follows from this and (\ref{72-7}) that
\beq\label{72-30}
|\rho_A({\tilde \kappa}_n(z_j))(\tau)| <\dt_{n+1},\,\,\,j=1,2,....
\eneq
for all $\tau\in T(A).$
It follows from {\ref{Extbot3}} that there is a unitary $w_n'\in U(A)$
such that
\beq\label{72-31}
&&\|[\phi_1(a), w_n']\|<\dt_{n+1}'/4\rforal a\in \imath_{n+1,\infty}({\cal G}_{n+1})\andeqn\\
&&\text{Bott}(\phi_1\circ \imath_{n+1,\infty},\, w_n')=-{\tilde
\kappa}_n\circ[\imath_{n+1}].
\eneq
By (\ref{72-4}),
\beq\label{72-32}
\text{Bott}(\phi_2\circ \imath_{n+1,\infty},\,v_n^*w_n'v_n)|_{{\cal
P}_n}=-{\tilde \kappa}_n\circ [\imath_{n+1}]|_{{\cal P}_n}.
\eneq
Put $w_n=v_n^*w_n'v_n.$ It follows from 10.6 of \cite{Linajm} that
\beq\label{72-33}
\Gamma(\text{Bott}(\phi_1\circ \imath_{n+1, \infty}, w_n'))&=&-\kappa_n\circ [\imath_{n+1}]
\andeqn\\
 \Gamma(\text{Bott}(\phi_1\circ \imath_{n+2, \infty},
w_{n+1}'))&=&-\kappa_{n+1}\circ [\imath_{n+2}].
\eneq
We also have
\beq\label{72-34}
\Gamma(\text{Bott}(\phi_1\circ \imath_{n+1,\infty},
v_nv_{n+1}^*))|_{H_n}=\zeta_n-\zeta_{n+1}\circ [\imath_{n+2}]=\xi_n.
\eneq
But, by (\ref{72-19}) {and (\ref{72-19+})},
\beq\label{72-35}
(-\kappa_n +\xi_n+\kappa_{n+1}\circ [{\imath_{n+2}}])=0.
\eneq
By 10.6 of \cite{Linajm}, $\Gamma({\rm Bott}(.,.))=0$ if and only
if ${\rm Bott}(.,.)=0.$ Thus, by (\ref{72-32}), (\ref{72-33}) and (\ref{72-34}),
\beq\label{72-36}
{\hspace{-0.5in}} -\text{Bott}(\phi_1\circ \imath_{n+1, \infty},\,w_n')
+\text{Bott}(\phi_1\circ \imath_{n+1, \infty}, \,v_nv_{n+1}^*)
+\text{Bott}(\phi_1\circ\imath_{n+1,\infty}, {w_{n+1}'}) =0.
\eneq

Define $u_n=x_nv_nw_n^*,$ $n=1,2,....$ Then, by (\ref{72-5}) and
(\ref{72-31}),
\beq\label{72-37}
{\rm ad}\, u_n\circ \phi_1\approx_{\dt_n'/2} \phi_2\rforal a\in
\imath_{n+1, \infty}({\cal G}_{n+1}).
\eneq

From (\ref{72-7}), (\ref{72-4}) and (\ref{72-36}), we compute
that
\beq\label{72-38}
&&\hspace{-0.6in}\text{Bott}(\phi_2\circ \imath_{n+1, \infty},u_n^*u_{n+1})
= \text{Bott}(\phi_2\circ \imath_{n+1, \infty}, w_nv_n^*v_{n+1}w_{n+1}^*)\\
&=&\hspace{-0.1in} \text{Bott}(\phi_2\circ \imath_{n+1, \infty}, w_n)+\text{Bott}(\phi_2\circ \imath_{n+1, \infty}, v_n^*v_{n+1})
%&&\hspace{1.8in}
+\text{Bott}(\phi_2\circ \imath_{n+1, \infty}, w_{n+1}^*)\\
&=&\text{Bott}(\phi_1\circ \imath_{n+1, \infty},w_n')+\text{Bott}(\phi_1\circ \imath_{n+1, \infty}, (w_{n+1}')^*)\\
&=&-[-\text{Bott}(\phi_1\circ \imath_{n+1, \infty}, w_n')+\text{Bott}(\phi_1\circ \imath_{n+1, \infty}, v_nv_{n+1}^*)\\
&&\hspace{2.1in}+\text{Bott}(\phi_1\circ \imath_{n+1, \infty}, {w_{n+1}'})]=0.
%&=&0.
\eneq
Let $x_{i,n}=[p_{i,n}]-[q_{i,n}], \,\, 1\le i\le I(n).$ Note that we assume that
$G_{u,n}$ is a free group generated by  $\{x_{i,n}: 1\le i\le I(n)\}.$ \Wlog, we may assume that these are free generators.
Define, for each $n\ge 1,$   a \hm\,
%By (\ref{72-n201410}),
%let
$\Lambda_n: G_{u,n}\to U(A)/CU(A)$   by
$$\Lambda_n(x_{i,n})=(\overline{\langle ((1-e_{i,n})+e_{i,n}u_n)((1-e'_{i,n})+e'_{i,n}u_n^*)\rangle},$$
%be a \hm\, {\red {why $\Lambda_n$ not $\Lambda$ as in (\ref{72-n201410}),}} such that
%\beq\label{72-n201410-1}
%\hspace{-0.2in}{{\rm dist}(\overline{\langle ((1-e_{i,n})+e_{i,n}u_n)((1-e'_{i,n})+e'_{i,n}u_n^*)\rangle}, \Lambda_{\red{n}}%(x_{i,n}))}<\sigma_n'/16,
%\eneq
where $e_{i,n}=\phi_2\circ \imath_{n+1,\infty}(p_{i,n}),$
$e'_{i,n}=\phi_2\circ \imath_{n+1,\infty}(q_{i,n}),$
%and $x_{i,n}=[p_{i,n}]-[q_{i,n}],$
$i=1,2,..., I(n).$
In what follows, we will construct unitaries $s_1, s_2, ..., {s_n,...}$ in $A$ such that
%\begin{equation}
\beq\label{n-72-37}
||[\phi_2\circ\iota_{n+1, \infty}(f), s_n]||<\delta'_{n+1}/4\rforal f\in\mathcal G_{n+1},\\
%\end{equation}
%\begin{equation}
\label{n-72-38}
\mathrm{Bott}(\phi_2\circ\iota_{n+1, \infty}, s_n)|_{\mathcal P_n}=0\andeqn\\
%\end{equation}
%and
%\beq
\label{n-72-39}
\hspace{-0.2in}{\rm dist}(\overline{\langle ((1-e_{i,n})+e_{i,n}s_n)((1-e'_{i,n})+e'_{i,n}s_n^*)\rangle},
\Lambda_n(-x_{i,n}))<\sigma_n'/16,
\eneq
%{\red{again $\Lambda_n$ or $\Lambda$}} {As mentioned earlier, we may assume that $\Lambda$ factors through $G_{u,n}'.$ }
%where $e_{i,n}=\phi_2\circ \imath_{n+1,\infty}(p_{i,n})$ and
%$e_{i,n}=\phi_2\circ \imath_{n+1,\infty}(q_{i,n}),$ $i=1,2,...,I(n).$
%\mathrm{dist}(\overline{s^*_ns_{n+1}, \overline{u^*_nu_{n+1}})<\sigma'_{n+1}.

Let $s_1=1$, and assume that $s_2, s_3, ..., s_n$ are already constructed. Let us construct $s_{n+1}$. Note that by \eqref{72-38}, the $K_1$ class of the unitary $u_n^*u_{n+1}$ is trivial. In particular, the $K_1$ class of $s_nu_n^*u_{n+1}$ is trivial.  {Since $\Lambda$ factors through $G_{u,n}',$ applying Theorem \ref{Extbot3}} to $\phi_2\circ\iota_{n+2, \infty}$, one obtains  a unitary $s_{n+1}\in B$ such that
%\begin{equation}
\beq
||[\phi_2\circ\iota_{n+2, \infty}(f), s_n]||<\delta'_{n+1}/4\rforal f\in\mathcal G_{n+2},\\
%\end{equation}
%\eneq
%\begin{equation}
%\beq
\mathrm{Bott}(\phi_2\circ\iota_{n+2, \infty}, s_{n+1})|_{\mathcal P_n}=0\andeqn
\eneq
%\end{equation}
%and
%\begin{equation}
\beq\label{n-72-39+}
\hspace{-0.2in}{\rm dist}(\overline{\langle ((1-e_{i,n+1})+e_{i,n+1}s_{n+1})((1-e'_{i,n+1})+e'_{i,n+1}s_{n+1}^*)\rangle},
\Lambda_{{{n+1}}}(-x_{i,n+1}))<\sigma_n'/16,
\eneq
%\end{equation}
%{\red{again $\Lambda_{n+1}$ or $\Lambda$}}
$i=1,2,...,I(n+1).$
Then  $s_1, s_2, ..., s_{n+1}$ satisfies \eqref{n-72-37}, \eqref{n-72-38} and \eqref{n-72-39}.

Put $\widetilde{u_n}=u_ns_n^*$. Then by \eqref{72-37} and \eqref{n-72-37}, one has
\beq\label{n-72-40}
{\rm ad}\, \widetilde{u_n}\circ \phi_1\approx_{\dt_n'} \phi_2\rforal a\in
\imath_{n+1, \infty}({\cal G}_{n+1}).
\eneq
By \eqref{72-38} and \eqref{n-72-38}, one has
\begin{equation}\label{n-72-41}
\mathrm{Bott}(\phi_2\circ\iota_{n+1, \infty}, (\widetilde{u_{n}})^* \widetilde{u_{n+1}})|_{\mathcal P_n}=0.
\end{equation}
Note that
\beq\label{72-n2014103}
\overline{\langle (1-e_{i,n})+e_{i,n}s_{n+1}u_{n+1}\rangle \langle(1-e_{i,n}')+e_{i,n+1}u_{n+1}^*s_{n+1}^*\rangle}
=\overline{c_1c_2c_3c_4}=\overline{c_1c_4c_2c_3},
\eneq
where
\beq
c_1=\langle(1-e_{i,n+1})+e_{i,n+1}s_{n+1}\rangle,\,\, c_2=\langle (1-e_{i,n+1})+e_{i,n+1}u_{n+1}\rangle,\\
c_3=\langle (1-e'_{i,n+1})+{e_{i,n+1}'}u_{n+1}^*\rangle,\,\, c_4=\langle (1-e'_{i,n+1})+e'_{i,n+1}{{s_{n+1}^*}}\rangle.
\eneq
Therefore,  by \eqref{n-72-39+},
%and by \eqref{72-n201410-1},
one has
\beq\label{n-72-42}
&&\hspace{-0.2in}\mathrm{dist}(\overline{\langle ((1-e_{i,n+1})+e_{i,n+1}{\widetilde{u_{n+1}}})((1-e_{i,n+1}')+e_{i,n+1}'{\widetilde{u_{n+1}}^*)}\rangle}), {\bar 1})\\
&&<\sigma_{n+1}'/16+{\rm dist}(\Lambda(-x_{i,n+1})\Lambda(x_{i,n+1}),{\bar 1})=\sigma_{n+1}'/16,
\eneq
$i=1,2,...,I(n).$
%&&<{\rm dist}(\overline{((1-e_{i,n+1})+e_{i,n+1}s_{n+1}((1-e_{i,n+1}')+e_{i,n+1}'s_{n+1}^*)), {\bar 1})\\
%&&+{\rm dist}(\overline{((1-e_{i,n+1})+e_{i,n+1}s_{n+1}((1-e_{i,n+1}')+e_{i,n+1}'s_{n+1}^*)), {\bar 1})
%\sigma_{n+1}.
%\end{equation}
Therefore, by \ref{BHfull} (and Remark \ref{RBHfull}), there exists a piece-wise smooth and
continuous path of unitaries $\{z_n(t): t\in [0,1]\}$ of $A$ such
that
\beq\label{72-39}
&&z_n(0)=1,\,\,\, z_n(1)=(\widetilde{u_n})^*\widetilde{u_{n+1}}\andeqn\\\label{72-40}
&&\|[\phi_2(a),\, z_n(t)]\|<1/2^{n+2}\rforal a\in {\cal F}_n\andeqn
t\in [0,1].
\eneq
Define
$$
u(t+n-1)=\widetilde{u_n}z_{n+1}(t)\,\,\,t\in (0,1].
$$
Note that $u(n)=\widetilde{u_{n+1}}$ for all integer $n$ and $\{u(t):t\in [0,
\infty)\}$ is a continuous path of unitaries in $A.$ One estimates
that, by (\ref{72-37}) and (\ref{72-40}),
\beq\label{72-41}
{\rm ad}\, u(t+n-1)\circ \phi_1 \approx_{\dt_n'} {\rm
ad}\,z_{n+1}(t)\circ \phi_2
 \approx_{1/2^{n+2}} \phi_2
\,\,\,\,\,\,\,\text{on}\,\,\,\,{\cal
F}_n
\eneq
 for all $t\in (0,1).$
It then follows that
\beq\label{72-42}
\lim_{t\to\infty}u^*(t)\phi_1(a) u(t)=\phi_2(a)\rforal a\in C.
\eneq

\end{proof}

\section{Rotation maps and strong asymptotic equivalence}

\begin{lem}\label{L91}
{Let {{$A$}} be a unital separable simple \CA\, of stable rank one.}
%Let $C\in {\cal B}_1$ be a unital separable simple \CA\, and let $A=C\otimes U,$ where $U$ is an infinite dimensional UHF-algebra.
Suppose that
$u\in CU(A).$ Then, for any piecewise smooth  and continuous path
$\{u(t): t\in [0,1]\}\subset U(A)$ with $u(0)=u$ and $u(1)=1_A,$
\beq\label{L91-1}
D_A(\{u(t)\})\in \overline{\rho_A(K_0(A))}.\,\,\,\,\,\,\hspace{0.4in}{\rm (} recall\,\, \ref{Dcu}\,\,for\,\,  D_A {\rm)}
\eneq
%{\rm (recall \ref{Dcu} for  $D_A$).}
%{where $D_A(\{u(t)\})(\tau)=({1\over{2\pi i}})\int_0^1 \tau({du(t)\over{dt}}u(t)^* )dt$ for all $\tau\in T(A).$}
\end{lem}

\begin{proof}
It follows from \ref{UCUiso} that the map
$j: u\mapsto {\rm diag}(u,1,..., 1)$ from $U(A)$ to $U(M_n(A))$ induces
an isomorphism from $U(A)/CU(A)$ to $U(M_n(A))/CU(M_n(A)).$ Then
the lemma follows from 3.1 and 3.2 of \cite{Thomsen-rims}.
\end{proof}

\begin{lem}\label{L92}
{Let {{$A$}}  be a unital separable simple \CA\, of stable rank one.}
%and let
%$A=C\otimes U,$ where $U$ is an infinite dimensional UHF-algebra.
%Let $C\in {\cal B}_1$ be a unital separable simple \CA\, and let
%$A=C\otimes U,$ where $U$ is an infinite dimensional UHF-algebra.
{Suppose that $B$ is a unital separable \CA\,}
% such that $U(B)/U_0(B)=K_1(B)$
and suppose that $\phi,\,\psi: B \to A$ are two unital monomorphisms such that
\beq\label{L92-1}
[\phi]=[\psi]\,\,\,{\rm in}\,\,\,KK(B,A)\\\label{L92-2}
\phi_T=\psi_T\andeqn \phi^{\ddag}=\psi^{\ddag}.
\eneq
Then
\beq\label{L92-3}
R_{\phi, \psi}\in {\rm Hom}(K_1(B), \overline{\rho_A(K_0(A))}).
\eneq

\end{lem}

\begin{proof}
Let $z\in K_1(B)$ be represented by a unitary $u\in {U(M_m(B))}$ {for some integer $m$}.
Then, by (\ref{L92-2}),
$$
{(\phi\otimes {\rm id}_{M_m})(u)(\psi\otimes {\rm id}_{M_m})(u)^*\in CU(M_m(A)).}
$$
Suppose that $\{u(t): t\in [0,1]\}$ is a piecewise smooth and continuous path in ${M_m(U(A))}$ such that
$u(0)={(\phi\otimes {\rm id}_{M_m})}(u)$ and $u(1)=(\psi\otimes{\rm id}_{M_m})(u).$ Put $w(t)={(\psi\otimes {\rm id}_{M_m})(u)^*u(t)}.$  Then $w(0)={(\psi\otimes {\rm id}_{M_m})(u)^*(\phi\otimes {\rm id}_{M_{m}})(u)}\in CU(A)$ and $w(1)=1_A.$ Thus
\beq\label{L92-4}
R_{\phi, \psi}(z)(\tau)&=&{1\over{2\pi i}}\int_0^1 \tau({du(t)\over{dt}}u^*(t))dt
= {1\over{2\pi i}}\int_0^1 \tau(\psi(u)^*{du(t)\over{dt}}u^*(t)\psi(u))dt\\
&=& {1\over{2\pi i}}\int_0^1 \tau({dw(t)\over{dt}}w^*(t))dt
\eneq
for all $\tau\in T(A).$ By \ref{L91},
\beq\label{L92-5}
R_{\phi, \psi}(z)\in \overline{\rho_A(K_0(A))}.
\eneq
It follows that
\beq\label{L92-6}
R_{\phi, \psi}\in {\rm Hom}(K_1(B), \overline{\rho_A(K_0(A))}).
\eneq

\end{proof}

\begin{thm}\label{T94}
Let $C_1, C_2\in {\cal B}_0$ be unital separable simple \CA s, $A=C_1\otimes U_1,$ $B=C_2\otimes U_2,$ where $U_1$ and $U_2$ are  UHF-algebras of infinite type.  Suppose that $B$ is a  {unital \SCA\,} of $A,$ and denote by $\imath$ the embedding. For any $\lambda\in  Hom(K_{{1}}(B), \overline{\rho_A(K_0(A))}),$ there exists $\phi\in {\overline{{\rm Inn}}}(B,A)${{(see \ref{dInn})}} such that there are \hm s $\theta_i: K_i(B)\to K_i(M_{\imath, \phi})$ with $(\pi_0)_{*i}\circ \theta_i={\rm id}_{K_i({{B}})},$ $i=0,1,$ and
the rotation map $R_{\imath, \phi}: K_1({{M_{\imath, \phi}}})\to {\rm Aff}(T(A))$ given by
\beq\label{T94-1}
R_{\imath, \phi}(x)=\rho_A(x-\theta_1(\pi_0)_{*1}(x))+\lambda\circ (\pi_0)_{*1}(x))
\eneq
for all $x\in K_1(M_{\imath, \phi}).$ In other words,
\beq\label{T94-2}
[\phi]=[\imath]\,\,\,{\rm in}\,\,\, KK(B,A)
\eneq
and the rotation map $R_{\imath, \phi}: K_1(M_{\imath, \phi})\to
{\rm Aff}(T(A))$ is given by
\beq\label{T94-3}
R_{\imath, \phi}(a,b)=\rho_A(a)+\lambda(b)
\eneq
for some identification of $K_1(M_{\imath, \phi})$ with $K_0(A)\oplus K_1(B).$
\end{thm}

\begin{proof}
The proof is exactly the same as that of Theorem 4.2 of \cite{L-N}.
In 4.2 of \cite{L-N}, it is assumed that $\rho_A(K_0(A))$ is dense in
${\rm Aff}(T(A)).$ However, it is that $\lambda(K_1(B))\subset \overline{\rho_A(K_0(A))}$ is used which we assume here. In 4.2 of \cite{L-N}, it is also assume that $A$ has the property (B1) and (B2) associated with $B$ (defined in {3.6} of \cite{L-N}). But this follows from
\ref{Extbot1} (see also \ref{Extbot2}).

\end{proof}

\begin{df}\label{d95}
Let $A$ be a unital \CA\, and let $C$ be a unital separable
\CA. Denote by ${\rm{Mon}}_{asu}^e(C,A)$  the set of all asymptotically unitary equivalence classes of unital monomorphisms from $C$ into $A.$ Denote by $\small{{\boldsymbol{K}}}:
% the map from
\textrm{Mon}_{asu}^e(C, A)\to {KK}_e(C,A)^{++}$  the map defined by
$$
\phi\mapsto [\phi]\tforal \phi\in \mathrm{Mon}_{asu}^e(C,A).
$$
Let $\kappa\in {KK}_e(C,A)^{++}.$ Denote by $\langle \kappa
\rangle$ the classes of $\phi\in \mathrm{Mon}_{asu}^e(C,A)$ such
that $\small{{\boldsymbol{ K}}}(\phi)=\kappa.$

Denote by $KKUT_e(A,B)^{++}$ the set of triples $(\kappa,
\af,\gamma)$ for which $\kappa\in KK_e(A,B)^{++},$ $\af:
U(A)/CU(A)\to U(B)/CU(B)$ is a \hm\, and $\gamma: T(B)\to T(A)$ is an
affine continuous map and both $\af$ and $\gamma$ are compatible
with $\kappa.$ Denote by $\boldsymbol{\mathfrak{K}}$ the map from
$\textrm{Mon}_{asu}^e(C, A)$ into ${KKUT}(C,A)^{++}$ defined by
$$
\phi\mapsto ([\phi],\phi^{\ddag}, \phi_T)\tforal \phi\in
\mathrm{Mon}_{asu}^e(C,A).
$$
Denote by $\langle \kappa, \af,\gamma \rangle $ the subset of
$\phi\in \mathrm{Mon}_{asu}^e(C,A)$ such that
$\boldsymbol{\mathfrak{K}}(\phi)=(\kappa,\,\af,\,\gamma).$

\end{df}

\begin{thm}\label{T96}
Let $C$ and $A$ be two  unital separable  amenable \CA s. Suppose that $\phi_1, \phi_2, \phi_3: {C\to A}$
are three unital monomorphisms for which
\beq\label{Mul-1}
[\phi_1]=[\phi_2]=[\phi_3]\,\,\,{\rm in}\,\,\, KK({C,A}))\andeqn
(\phi_1)_T=(\phi_2)_T=(\phi_3)_T.
\eneq
Then
\beq\label{Mul-2}
\overline{R}_{\phi_1, \phi_2}+\overline{R}_{\phi_2,
\phi_3}=\overline{R}_{\phi_1,\phi_3}.
\eneq
\end{thm}

\begin{proof}
The proof is exactly the same as that of Theorem 9.6 of \cite{Lnclasn}.
\end{proof}

\begin{lem}\label{L97}
Let $A$ {{and}} $B$ be two unital separable amenable \CA s.
Suppose that $\phi_1, \phi_2: A\to B$ are  two
unital monomorphisms such that
$$
[\phi_1]=[\phi_2]\,\,\,{\rm in}\,\,\, KK(A,B)\andeqn
(\phi_1)_T=(\phi_2)_T.
$$
Suppose that $(\phi_2)_T: T(B)\to T(A)$ is an affine
homeomorphism. Suppose also that there is $\af\in Aut(B)$ such that
%$\boldsymbol{\mathfrak{K}}(\af)={\boldsymbol{\mathfrak{K}}}({\rm
%id}_B).$
$$[\af]=[{\rm id}_B]\,\,\,{\rm in}\,\,\, KK(B,B)\andeqn \af_T={\rm id}_T.
$$
Then
\beq\label{Group1-0}
\overline{R}_{\phi_1,\af\circ \phi_2} =\overline{R}_{{\rm
id}_B,\af}\circ (\phi_2)_{*1}+\overline{R}_{\phi_1,\phi_2}
% \eta_{\phi_1}(\af\circ
%\phi_2)=\eta_{{\rm id}_B}(\af)\circ
%(\phi_2)_{*1}+\eta_{\phi_1}(\phi_2)
\eneq
in ${\rm Hom}(K_1(A), \overline{\rho_B(K_0(B))})/{\cal R}_0.$

\end{lem}

\begin{proof}
By \ref{T96}, we compute that
\beq\nonumber
\overline{R}_{\phi_1, \af\circ \phi_2}&=&\overline{R}_{\phi_1, \phi_2}+\overline{R}_{\phi_2,\af\circ \phi_2}=
 \overline{R}_{\phi_1, \phi_2}+\overline{R}_{{\rm id}_B, \af}\circ (\phi_2)_{*1}.
\eneq
\end{proof}

\begin{thm}\label{MT2}
Let $B\in {\cal N}$ be a unital   simple \CA\, in ${\cal B}_0,$ let $C=B\otimes U_1,$ where $U_1$ is a UHF-algebra
of infinite type{,} { let $A_1$ is a unital separable amenable simple \CA\, in ${\cal B}_0$ and let $A=A_1\otimes U_2,$
where $U_2$ is another UHF-algebra of infinite type. }
% and let $A$ be a unital separable simple \CA\, in ${\cal B}_0$.
 Then the map $\boldsymbol{\mathfrak{K}}:
\mathrm{Mon}_{asu}^e(C,A)\to {KKUT}(C,A)^{++}$ is surjective.
Moreover, for each $(\kappa, \af, \gamma)\in {KKUT}(C,A)^{++}$,
there exists a bijection
$$
\eta: \langle \kappa,\af,\gamma \rangle \to \mathrm{Hom}({K}_1(C),
{\overline{\rho_A(K_0(A))}})/{\cal R}_0.
$$
%where $\langle{\cal R}_0\rangle$ denotes the subgroup generated by $\mathcal R_0$.
\end{thm}

\begin{proof}
It follows from \ref{L86} that $\boldsymbol{\mathfrak{K}}$ is
surjective.

Fix a triple $(\kappa,\af, \gamma)\in {KKT}(C,A)^{++}$ and choose
a unital monomorphism $\phi: C\to A$ such that $[\phi]=\kappa$,
$\phi^{\ddag}=\af$ and  $\phi_\mathrm{T}=\gamma.$
If $\phi_1: C\to A$ is another unital monomorphism such that
$\boldsymbol{\mathfrak{K}}(\phi_1)=\boldsymbol{\mathfrak{K}}({{\phi}}),$ then by \ref{L92},
\beq\label{T98-1}
{\overline{R}}_{\phi, \phi_1}\in {\rm Hom}(K_1(C), \overline{\rho_A(K_0(A))})/{\cal R}_0.
\eneq

Let $\lambda\in \mathrm{Hom}({K}_1(C), \overline{\rho_A(K_0(A))})$
be a \hm.  It follows from \ref{T94} that there is a unital monomorphism
$\psi\in \overline{{{\rm Inn}}}(\phi(C), A)$ with $[\psi\circ \phi]=[\phi]$ in $KK(C,A)$ such that there exists a homomorphism
$\theta: {K}_1(C)\to K_1(M_{\phi, \psi\circ \phi})$ with
$(\pi_0)_{*1}\circ \theta={\rm{id}}_{{K}_1(C)}$ for which
$R_{\phi, \psi\circ \phi}\circ \theta=\lambda.$ Let
$\beta=\psi\circ \phi.$ Then $R_{\phi, \beta}\circ
\theta=\lambda.$ Note also since $\psi\in
{\overline{\rm{Inn}}}(\phi(C), A),$ $\beta^{\ddag}=\phi^{\ddag}$
and $\beta_\mathrm{T}=\phi_\mathrm{T}.$ In particular, $\boldsymbol{\mathfrak{K}}(\beta)={\boldsymbol{\mathfrak{K}}}(\phi).$

Thus{{, for each unital monomorphism $\phi$,}} we obtain a well-defined  and surjective map $$\eta_{{\phi}}: \langle
[\phi], \phi^{\ddag}, \phi_T\rangle \to \mathrm{Hom}({K}_1(A),
\overline{\rho_A(K_0(A))})/{\cal R}_0.$$
To see it is one to one,  consider two monomorphisms
%let
%
$\phi_1, \phi_2: C\to A$
% be  two unital monomorphisms 
in $ \langle
[\phi], \phi^{\ddag},\phi_T\rangle$  such that
$$
\overline{R}_{\phi, \phi_1}=\overline{R}_{\phi, \phi_2}.
$$
Then, by \ref{L92},
\beq
\overline{R}_{\phi_1,\phi_2}=\overline{R}_{\phi_1,\phi}+\overline{R}_{\phi,
\phi_2}
=-\overline{R}_{\phi, \phi_1}+\overline{R}_{\phi, \phi_2}=0.
\eneq
It follows from \ref{Tm72} that $\phi_1$ and $\phi_2$ are
asymptotically unitarily equivalent. The map $\eta_{{\phi}}$ is the desired bijection $\eta$ as $\langle
[\phi], \phi^{\ddag}, \phi_T\rangle=\langle \kappa,\af,\gamma \rangle$.
\end{proof}

\begin{df}
{\rm  Denote by $KKUT_e^{-1}(A,A)^{++}$ the subgroup of those
elements ${{(}} \kappa, \af, \gamma{{)} }\in KKUT_e(A,A)^{++}$
for which $\kappa|_{K_i(A)}$ is an isomorphism $(i=0,1$), $\af$ is an
isomorphism and $\gamma$ is {an}  affine homeomorphism.  {{Recall from the proof of \ref{MT2}
$\eta_{{\rm id}_A}: {\langle [{\rm id}_A], {\rm
id}_A^{\ddag}, ({\rm id}_A)_T\rangle} \to \mathrm{Hom}({K}_1(A),
\overline{\rho_A(K_0(A))})/{\cal R}_0$ is a bijection.}}

Denote by $\langle {\rm id}_A\rangle $ the class of those
automorphisms $\psi$ which are asymptotically unitarily equivalent
to ${\rm id}_A${{---this subset of ${\rm Aut}(A)$ gives up a single element in $\mathrm{Mon}_{asu}^e(A,A)$ which should not be confused with the subset  ${\langle [{\rm id}_A], {\rm
id}_A^{\ddag}, ({\rm id}_A)_T\rangle}\subset\mathrm{Mon}_{asu}^e(A,A)$.}}
Note that, if $\psi\in \langle {\rm id}_A\rangle
,$ then $\psi$ is {\it asymptotically inner}, i.e., there exists a
continuous path of unitaries $\{u(t): t\in [0,\infty)\}\subset A$
such that
$$
\psi(a)=\lim_{t\to\infty}u(t)^*au(t)\tforal a\in A.
$$

{{ Note that $\langle {\rm id}_A\rangle $ is a normal subgroup of ${\rm Aut}(A)$.}}

}
\end{df}

\begin{cor}\label{C910}
Let $A_1\in {\cal N}\cap {\cal B}_0$ be a unital simple  \CA\, and let
$A=A_1\otimes U$ for some  UHF-algebra $U$ of infinite type. Then one
% group ${\rm Aut}(A)/\langle {\rm id}_A\rangle$
has the following short exact sequence:
\beq\label{Inn1}
\hspace{-0.4in}0 \to {\rm Hom}(K_1(A), \overline{\rho_A(K_0(A))})/{\cal
R}_0\stackrel{\eta_{{\rm id}_A}^{-1}}{\to}{\rm Aut}(A)/\langle
{\rm id}_A\rangle \stackrel{{\boldsymbol{\mathfrak{K}}}}{\to}
KKUT_e^{-1}(A,A)^{++}\to 0.
\eneq

In particular, if $\phi, \psi\in Aut(A)$ such that
$$
{\boldsymbol{\mathfrak{K}}}(\phi)={\boldsymbol{\mathfrak{K}}}(\psi)={\boldsymbol{\mathfrak{K}}}({\rm
id}_A),
$$
Then
$$
\eta_{{\rm id}_A}(\phi\circ \psi)=\eta_{{\rm
id}_A}(\phi)+\eta_{{\rm id}_A}(\psi).
$$

\end{cor}

\begin{proof}
It follows from \ref{L86} that,  for any $\langle \kappa, \af,
\gamma\rangle,$ there is a unital monomorphism $h: A\to A$ such
that ${\boldsymbol{\mathfrak{K}}}(h)=\langle \kappa, \af,
\gamma\rangle.$ The fact that $\kappa\in KK_e^{-1}(A,A)^{++}$
implies that there is $\kappa_1\in KK_e^{-1}(A,A)^{++}$ such that
$$
\kappa\times \kappa_1=\kappa_1\times \kappa=[{\rm id}_A].
$$
By \ref{L86}, choose $h_1: A\to A$ such that
$$
{\boldsymbol{\mathfrak{K}}}(h)=\langle \kappa_1, \af^{-1},
\gamma^{-1}\rangle.
$$

It follows from \ref{MUN2} that $h_1\circ h$
and $h\circ h_1$ are approximately unitarily equivalent. Applying
a standard approximate intertwining argument of G. A. Elliott, one
obtains two isomorphisms $\phi$ and $\phi^{-1}$ such that there is
a sequence of unitaries $\{u_n\}$ in $A$ such that
$$
\phi(a)=\lim_{n\to\infty}{\rm Ad}\, u_{2n+1}\circ h(a)\andeqn
\phi^{-1}(a)=\lim_{n\to\infty}{\rm Ad}\, u_{2n}\circ h_1(a)
$$
for all $a\in A.$ Thus $[\phi]=[h]$ in $KL(A,A)$ and
$\phi^{\ddag}=h^{\ddag}$ and $\phi_T=h_T.$ Then, as in the proof
of \ref{L86}, there is $\psi_0\in {\overline{{\rm Inn}}}(A,A)$ such that
$[\psi_0\circ \phi]=[{\rm id}_A]$ in $KK(A,A)$ as well as
$(\psi_0\circ \phi)^{\ddag}=h^{\ddag}$ and $(\psi_0\circ
\phi)_T=h_T.$  So we have $\psi_0\circ \phi\in Aut(A,A)$ such that
${\boldsymbol{\mathfrak{K}}}(\psi_0\circ \phi)=\langle \kappa, \af,
\gamma\rangle.$ {This implies that ${\mathfrak{K}}$ is surjective.}

Now let $\lambda\in {\rm Hom}(K_1(C), \overline{\text{Aff}(T(A))})/{\cal
R}_0.$  The proof \ref{MT2} says that there is $\psi_{00} \in
\overline{{\rm Inn}}(A,A)$ (in place of $\psi$) such that
${\boldsymbol{\mathfrak{K}}}(\psi_{00}\circ {\rm
id}_A)={\boldsymbol{\mathfrak{K}}}({\rm id}_A)$ and
$$
\overline{R}_{{\rm id}_A, \psi_{00}}=\lambda.
$$
Note that $\psi_{00}$ is again an automorphism.
The last part of the lemma then follows from \ref{L97}.
\end{proof}

%\section{Strongly asymptotic unitary equivalence}

\begin{df}{\rm (Definition 10.2 of \cite{Linajm} and see also \cite{LnTAMS12})}\label{Dsu}
Let $A$ be a unital \CA\, and $B$ be another \CA. Recall
(\cite{LnTAMS12}) that
$$
H_1(K_0(A), K_1(B))=\{x\in K_1(B): \phi([1_A])=x\, \, \text{for\,some }\,\phi\in
{\rm Hom}(K_0(A), K_1(B))\}.
$$
\end{df}

\begin{prop}{\rm (Proposition 12.3 of \cite{Linajm})}\label{Sup}
Let $A$  be a unital separable \CA\, and let $B$ be a unital \CA.
Suppose that $\phi: A\to B$ is a unital \hm\, and $u\in U(B)$ is a
unitary. Suppose that there is a continuous path of unitaries
$\{u(t): t\in [0,\infty)\}\subset B$ such that
\beq\label{sup-1}
u(0)=1_B\andeqn \lim_{t\to\infty}{\rm ad}\, u(t)\circ \phi(a)={\rm
ad}\, u\circ \phi(a)
\eneq
for all $a\in A.$  Then
$$
[u]\in H_1(K_0(A), K_1(B)).
$$

\end{prop}

%{\bf The proof of the following follows exactly the same lines as those of
%Theorem 10.5 of \cite{Lnclasn}.}

\begin{lem}\label{fix-unitary}
Let $C=C'\otimes U$ for some $C'=\varinjlim(C_n, \psi_n)$ and a UHF algebra $U$ of infinite type, where each $C_n$ is a direct sum of C*-algebras in $\mathcal C_0$ and $\mathbf H$. Assume that $\psi_n$ is unital and injective. Let $A\in\mathcal B_1$. Let $\phi_1, \phi_2: C\to A$ be two monomorphisms such that there is an increasing sequence of finite subsets $\mathcal F_n\subseteq C$ with dense union, an increasing sequence of finite subsets $\mathcal P_n\subseteq K_1(C)$ with union to be $K_1(C)$, a sequence of positive numbers $(\delta_n)$ with $\sum\delta_n<1$ and a sequence of unitaries $\{u_n\}\subseteq A$, such that
$$\mathrm{Ad}u_n\circ\phi_1\approx_{\delta_n} \phi_2\quad\textrm{on}\ \mathcal F_n\tand\,\,
%and
%$$
\rho_A(\mathrm{bott}_1(\phi_2, u_n^*u_{n+1}))=0\tforal  x\in\mathcal P_n.$$
Suppose that $H_1(K_0(C), K_1(A))=K_1(A).$
Then there exists a sequence of unitaries $v_n\in U_0(A)$ such that
\beq\label{strong-1}
\mathrm{Ad}v_n\circ\phi_1\approx_{\delta_n} \phi_2\quad\textrm{on}\ \mathcal F_n\andeqn\\
\rho_A(\mathrm{bott}_1(\phi_2, v_n^*v_{n+1}))=0,\quad x\in\mathcal P_n.
\eneq
\end{lem}

\begin{proof}
Let $x_n=[u_n]\in K_1(A)$. Since $H_1(K_0(C), K_1(A))=K_1(A)$, there is a homomorphism $$\kappa_{n, 0}: K_0(C)\to K_1(A)$$
such that $\kappa_{n, 0}([1_C])=-x_n$. Since $C$ satisfies the Universal Coefficient Theorem, there is $\kappa_n\in KL(C\otimes C(\T), A)$ such that
$$
(\kappa_n)|_{{\boldsymbol{\bt}}(K_0(C))}=\kappa_{n, 0}\quad\mathrm{and}\quad (\kappa_n)|_{{\boldsymbol{\bt}}(K_1(C))}=0.
$$
\Wlog, we may assume that
%${\cal F}_n\subset \psi_{n, \infty}(C_n)$ and
$[1_{C}]\in {\cal P}_n,$  $n=1,2,.....$
For each $\delta_n$, choose a positive number $\eta_n<\delta_n,$
such that
$$\mathrm{Ad}u_{n}\circ\phi_1\approx_{\eta_n}\phi_2\quad \textrm{on}\ \mathcal F_n. $$
By \ref{Extbot1}, there is a unitary $w_n\in U(A)$ such that
\beq\nonumber
\|[\phi_2(a), w_n]\|<(\delta_n-\eta_n)/2\quad \rforal a\in\mathcal F_n\andeqn
%and
%$$
\mathrm{Bott}(\phi_2, w_n)|_{{\cal P}_n}=\kappa_n|_{{\boldsymbol{\bt}}({\cal P}_n)}.
\eneq
%$$
Put $v_n=u_nw_n$, $n=1, 2, ...$. Then
$$\mathrm{Ad} v_n\circ\phi_1\approx_{\delta_n} \phi_2\quad\textrm{on}\ \mathcal F_n,\,\,\,
\rho_A(\mathrm{bott}_1(\phi_2, v_n^*v_{n+1}))|_{{\cal P}_n}=0$$
and, since $[1_C]\in {\cal P}_n,$
$$[v_n]=[u_n]-x_n=0,$$
as desired.
\end{proof}

\begin{thm}\label{T105}
Let $B\in {\cal B}_1$ be a unital separable simple \CA\, which satisfies the UCT,
let $A_1\in {\cal B}_1$ be a unital separable
simple \CA\, and let $C=B\otimes U_1$ and $A=A_1\otimes U_2,$ where $U_1$ and $U_2$ are two unital infinite dimensional UHF-algebras. Suppose that
$H_1(K_0(C),K_1(A))=K_1(A)$ and suppose that $\phi_1, \phi_2: C\to
A$ are two unital monomorphisms which are asymptotically unitarily
equivalent. Then there exists a continuous path of unitaries
$\{u(t): t\in [0, \infty)\}\subset A$ such that
$$
u(0)=1\andeqn \lim_{t\to\infty}{\rm ad}u(t)\circ
\phi_1(a)=\phi_2(a)\rforal a\in C.
$$
\end{thm}
\begin{proof}
By 4.3 of \cite{Lnclasn}, one has
\begin{eqnarray*}
[\phi_1]=[\phi_2]\quad {\rm in}\,\,\, KK(C,A),\\
\phi^{\ddag}=\psi^{\ddag},\,\,\,(\phi_1)_T=(\phi_2)_T \andeqn {\overline{R_{\phi_1,\phi_2}}}=0.
\end{eqnarray*}
Then by Lemma \ref{fix-unitary}, one may assume that $v_n\in U_0(A)$, $n=1, 2, ...$ in the proof \ref{Tm72}. It follows that $\xi_n([1_C])=0$, $n=1, 2, ...$, and therefore $\kappa_n([1_C])=0$. This implies that $\gamma_n\circ\beta([1_C])=0$. Hence $w_n\in U_0(A)$, and $u_n\in U_0(A)$. Therefore, the continuous path of unitaries $\{u(t)\}$ constructed in \ref{Tm72} is in $U_0(A)$, and then one may require that $u(0)=1_A$.
\end{proof}

\section{The classification theorem}

\begin{lem}\label{L112}
Let $A_1\in {\cal B}_0$ be a unital separable simple \CA, let $A=A_1\otimes U$ for some infinite dimensional UHF-algebra and
let $\mathfrak{p}$ be a supernatural number of infinite type. Then
the \hm\, $\imath: a\mapsto a\otimes 1$ induces an isomorphism
from $U_0(A)/CU(A)$ to $U_0(A\otimes M_{\mathfrak{p}})/CU(A\otimes
M_{\mathfrak{p}}).$
\end{lem}

\begin{proof}
There are  sequences of positive integers $ \{m(n)\}$ and
$\{k(n)\}$ such that $A\otimes
M_{\mathfrak{p}}=\lim_{n\to\infty}(A\otimes M_{m(n)}, \imath_n),$
where
$$
\imath_n: M_{m(n)}(A)\to M_{m(n+1)}(A)
$$
is defined by $\imath(a)={\rm diag}(\overbrace{a,a,...,a}^{k(n)})$
for all $a\in M_{m(n)}(A),$ $n=1,2,....$ Note, $M_{m(n)}(A)=M_{m(n)}(A_1)\otimes U$ and $M_{m(n)}(A_1)\in {\cal B}_0.$ Let $$j_n:
U(M_{m(n)}(A))/CU(M_{m(n)}(A)))\to
U(M_{m(n+1)}(A))/CU(M_{m(n+1)}(A))
$$ be defined by
$$
j_n({\bar u})=\overline{{\rm
diag}(u,\underbrace{1,1,...,1}_{k(n)-1})}\tforal u\in
U((M_{m(n)}(A)).
$$
It follows from Theorem \ref{UCUiso} that $j_n$ is
an isomorphism.
 By \ref{Unotrosion}, $U_0(M_{m(n)}(A))/CU(M_{m(n)}(A))$
is divisible. For each $n$ and $i,$ there is a unitary $U_i\in
M_{m(n+1)}(A)$ such that
$$
U_i^*E_{1,1}U_i=E_{i,i},\,\,\,i=2,3,...,k(n),
$$
where $E_{i,i}=\sum_{j=(i-1)m(n)+1}^{im(n)}e_{j,j}$  and
$\{e_{i,j}\}$ is a matrix unit for $M_{m(n+1)}.$ Then
$$
\imath_n(u)=u'U_2^*u'U_2\cdots U_{k(n)}^*u'U_{k(n)},
$$
where $u'={\rm diag}(u,\overbrace{1,1,...,1}),$ for all $u\in
M_{m(n)}(A).$ Thus
$$
\imath_n^{\ddag}({\bar u})=k(n)j_n(\bar u).
$$
It follows that
$\imath_n^{\ddag}|_{U_0(M_{m(n)}(A))/CU(M_{m(n)}(A))}$ is
injective, since $U_0(M_{m(n+1)}(A))/CU(M_{m(n+1)}(A))$ is torsion
free (see \ref{UCUdiv}). For each $z\in
U_0(M_{m(n+1)}(A)/CU(M_{m(n+1)}),$ there is a unitary $v\in
M_{m(n+1)}(A)$ such that
$$
j_n({\bar v})=z,
$$
since $j_n$ is an isomorphism. By the divisibility of
$U_0(M_{m(n)}(A)/CU(M_{m(n)}),$ there is $u\in M_{m(n)}(A)$ such
that
$$
\overline{u^{k(n)}}=\overline{u}^{k(n)}=\overline{v}.
$$
As above,
$$
\imath_n^{\ddag}({\bar u})=k(n)j_n(\bar v)=z.
$$
So $\imath_n^{\ddag}|_{U_0(M_{m(n)}(A))/CU(M_{m(n)}(A))}$ is a
surjective.  It follows that $\imath_{n,
\infty}^{\ddag}|_{U_0(M_{m(n)}(A))/CU(M_{m(n)}(A))}$ is also an
isomorphism. One then concludes that
$\imath^{\ddag}|_{U_0(A)/CU(A)}$ is an isomorphism.
\end{proof}

\begin{lem}\label{L113}
Let $A_1$ and $B_1$ be two unital separable simple \CA s  in ${\cal B}_0,$ let $A=A_1\otimes U_1$ and let $B=B_1\otimes U_2,$ where $U_1$ and $U_2$ are two infinite dimensional UHF-algebras.
Let $\phi: A\to B$ be an isomorphism and let $\bt:
B\otimes M_{\mathfrak{p}}\to B\otimes M_{\mathfrak{p}}$ be an
automorphism such that $\bt_{*1}={\rm id}_{K_1(B\otimes
M_{\mathfrak{p}})}$ for some supernatural number $\mathfrak{p}$ of infinite type.
Then
$$\psi^{\ddag}(U(A)/CU(A))=
(\phi_0)^{\ddag}(U(A)/CU(A))=U(B)/CU(B),
$$
where $\phi_0=\imath\circ \phi,$ $\psi=\bt\circ \imath\circ \phi$
and where $\imath: B\to B\otimes M_{\mathfrak{p}}$ is defined by
$\imath(b)=b\otimes 1$ for all $b\in B.$ Moreover there is an
isomorphism $\mu: U(B)/CU(B)\to U(B)/CU(B)$ with
$\mu(U_0(B)/CU(B))\subset U_0(B)/CU(B)$ such that
$$
\imath^{\ddag}\circ \mu\circ \phi^{\ddag}=\psi^{\ddag}\andeqn
q_1\circ \mu=q_1,
$$
where $q_1: U(B)/CU(B)\to K_1(B)$ is the quotient map.
\end{lem}

\begin{proof}
The proof is exactly the same as that of Lemma 11.3 of \cite{Lnclasn}.
\end{proof}

\begin{lem}\label{L114}
Let $A_1$ and $B_1$ be two unital   simple amenable \CA s  in
${\cal N}\cap {\cal B}_0,$ let $A=A_1\otimes U_1$ and let
$B=B_1\otimes U_2,$ where $U_1$ and $U_2$ are
UHF-algebras {{of infinite type}}.  Suppose that
$\phi_1, \phi_2: A\to B$ are two isomorphisms such that
$[\phi_1]=[\phi_2]$ in $KK(A,B).$ Then there exists an
automorphism $\bt: B\to B$ such that $[\bt]=[{\rm id}_B]$ in
$KK(B,B)$ and $\bt\circ \phi_2$ is asymptotically unitarily
equivalent to $\phi_1.$ Moreover, if $H_1(K_0(A), K_1(B))=K_1(B),$
they are strongly asymptotically unitarily equivalent.
\end{lem}

\begin{proof}
It follows from \ref{MT2} that there is an automorphism $\bt_1:
B\to B$ satisfying the following:
\beq
[\bt_1]=[{\rm id}_B]\,\,\,{\rm in}\,\,\,KK(B,B),\\
\bt_1^{\ddag}=\phi_1^{\ddag}\circ (\phi_2^{-1})^{\ddag}\andeqn
(\bt_1)_T=(\phi_1)_T\circ (\phi_2)_T^{-1}.
\eneq
By \ref{C910}, there is automorphism $\bt_2\in Aut(B)$ such that
\beq
[\bt_2]=[{\rm id}_B]\,\,\, {\rm in}\,\,\, KK(B,B),\\
\bt_2^{\ddag}={\rm id}_B^{\ddag},\,\,\,(\bt_2)_T=({\rm
id}_B)_T\andeqn\\
 \overline{R}_{{\rm id}_B,\bt_2}=-\overline{R}_{\phi_1,\bt_1\circ \phi_2}\circ
(\phi_2)_{*1}^{-1}.
\eneq
Put $\bt=\bt_2\circ \bt_1.$
 It follows that
\beq
[\bt\circ \phi_2]=[\phi_1]\,\,\,{\rm in}\,\,\, KK(A,B),
(\bt\circ \phi_2)^{\ddag}=\phi_1^{\ddag}\andeqn (\bt\circ
\phi_2)_T=(\phi_1)_T.
\eneq
Moreover, by \ref{L97},
\beq
\overline{R}_{\phi_1,\bt\circ \phi_2}&=&\overline{R}_{{\rm
id}_B,\bt_2}\circ
(\phi_2)_{*1}+\overline{R}_{\phi_1,\bt_1\circ\phi_2}\\
&=&(-\overline{R}_{\phi_1,\bt_1\circ \phi_2}\circ
(\phi_2)_{*1}^{-1})\circ (\phi_2)_{*1}
+\overline{R}_{\phi_1,\bt_1\circ\phi_2}=0.
\eneq
It follows from \ref{MT2} that $\bt\circ \phi_2$ and $\phi_1$ are
asymptotically unitarily equivalent.

In the case that $H_1(K_0(A), K_1(B))=K_1(B),$ it follows from
\ref{T105} that $\bt\circ \phi_2$ and $\phi_1$ are strongly
asymptotically unitarily equivalent.
\end{proof}

\begin{lem}\label{L115}
Let $A_1$ and let $B_1$ be two  unital simple amenable \CA s in ${\cal
N}\cap {\cal B}_0$ and let $A=A\otimes U_1$ and
 $B=B_1\otimes U_2$ for some  UHF-algebras $U_1$ and $U_2$ of infinite type.  Let $\phi: A\to B$ be an
isomorphism. Suppose that $\bt\in Aut(B\otimes M_{\mathfrak{p}})$
for which
$$[\bt]=[{\rm id}_{B\otimes
M_{\mathfrak{p}}}]\,\,\,{\rm in}\,\,\,KK(B\otimes
M_{\mathfrak{p}},B\otimes M_{\mathfrak{p}}) \andeqn \bt_T=({\rm
id}_{B\otimes M_{\mathfrak{p}}})_T
$$
for some supernatural number $\mathfrak{p}$ of infinite type.

Then there exists an automorphism $\af\in Aut(B)$ with
$[\af]=[{\rm id}_{B}]$ in $KK(B,B)$ such that $ \imath\circ
\af\circ \phi $ and $\bt\circ \imath\circ \phi$ are asymptotically
unitarily equivalent, where $\imath: B\to B\otimes
M_{\mathfrak{p}}$ is defined by $\imath(b)=b\otimes 1$ for all
$b\in B.$
\end{lem}

\begin{proof}
It follows from \ref{L113} that there is an isomorphism $\mu:
U(B)/CU(B)\to U(B)/CU(B)$ such that
$$
\imath^{\ddag}\circ \mu\circ \phi^{\ddag}=(\bt\circ \imath\circ
\phi)^{\ddag}.
$$
Note that $\imath_T: T(B\otimes M_{\mathfrak{p}})\to T(B)$ is an
affine homeomorphism.

It follows from \ref{MT2} that there is an automorphism $\af: B\to
B$ such that
\beq\label{l1-1}
&&[\af]=[{\rm id}_B]\,\,\,{\rm in}\,\,\,KK(B,
B),\\
&&\af^{\ddag}=\mu,\,\,\, \af_T=(\bt\circ \imath\circ \phi)_T\circ
((\imath\circ
\phi)_T)^{-1}=({\rm id}_{B\otimes M_{\mathfrak{p}}})_T\andeqn\\
&&\overline{R}_{{\rm
id}_B,\af}(x)(\tau)=-\overline{R}_{\bt\circ\imath\circ \phi,\,
\imath\circ \phi}(\phi_{*1}^{-1}(x))(\imath_T(\tau))\tforal x\in
K_1(A)
\eneq
and for all $\tau\in T(B).$

Denote by $\psi=\imath\circ \af\circ \phi.$ Then we have, by
\ref{L97},
\beq
&&[\psi]= [\imath\circ \phi]=[\bt\circ\imath\circ \phi]\,\,\,{\rm
in}\,\,\, KK(A, B\otimes M_{\mathfrak{p}})\\
&&\psi^{\ddag}=\imath^{\ddag}\circ\mu\circ \phi^{\ddag}=(\bt\circ
\imath\circ \phi)^{\ddag},\\
&&\psi_T=(\imath\circ \af\circ \phi)_T=(\imath\circ
\phi)_T=(\bt\circ \imath\circ \phi)_T.
\eneq
Moreover, for any $x\in K_1(A)$ and $\tau\in T(B\otimes
M_{\mathfrak{p}}),$
\beq
\overline{R}_{\bt\circ \imath\circ \phi, \psi}(x)(\tau)&=&
\overline{R}_{\bt\circ \imath\circ\phi,
\imath\circ\phi}(x)(\tau)+\overline{R}_{\imath, \imath\circ \af}\circ
\phi_{*1}(x)(\tau)\\
&=&\overline{R}_{\bt\circ \imath\circ\phi,
\imath\circ\phi}(x)(\tau)+\overline{R}_{{\rm id}_B, \imath\circ \af}\circ
\phi_{*1}(x)(\imath_T^{-1}(\tau))\\
&=&\overline{R}_{\bt\circ \imath\circ\phi,
\imath\circ\phi}(x)(\tau)-\overline{R}_{\bt\circ \imath\circ\phi,
\imath\circ\phi}(\phi_{*1}^{-1})(\phi_{*1}(x))(\tau)=0
\eneq
It follows from \ref{MT2}  that $\imath\circ \af\circ \phi$ and
$\bt\circ \imath\circ \phi$ are asymptotically unitarily
equivalent.
\end{proof}

%{\bf It appears that $Ell(A)$ is first time defined here. We need to move earlier.}

%\begin{NN}
%{\rm Let $A$ be a unital separable simple \CA. By the Elliott
%invariant we mean the 6-tuple:
%$$
%Ell(A)=(K_0(A), K_0(A)_+, [1_A], K_1(A), T(A),\rho_A).
%$$
%
%Let $B$ be another unital separable simple \CA. We say that there
%is an isomorphism
%$$
%\Gamma: Ell(A)\to Ell(B),
%$$
%if there is an  order isomorphism $\kappa_0: (K_0(A),
%K_0(A)_+,[1_A])\to (K_0(B), K_0(B)_+, [1_B]),$  there is an
%isomorphism $\kappa_1: K_1(A)\to K_1(B)$ and there is an affine
%homeomorphism $\mu: T(B)\to T(A)$ for which
%$$
%\mu(\tau)([p])=\tau(\kappa_0([p]))
%$$
%for all projection $p\in M_k(A)$ (for all $k\ge 1$) and all
%tracial states $\tau\in T(B).$

%}
%\end{NN}

  \begin{thm}\label{CMT1}
Let $A$ and $B$ be two unital separable  simple
\CA s in ${\cal N}.$  Suppose that there is an isomorphism
$$
\Gamma: {\rm Ell}(A)\to {\rm Ell}(B).
$$
Suppose also that, for some pair of relatively prime supernatural
numbers $\mathfrak{p}$ and $\mathfrak{q}$ of infinite type such
that $M_{\mathfrak{p}}\otimes M_{\mathfrak{q}}\cong Q,$
$A\otimes M_{\mathfrak{p}}\in {\cal B}_0,B\otimes
M_{\mathfrak{p}}\in {\cal B}_0$, $A\otimes M_{\mathfrak{q}}\in {\cal B}_0$ and
$B\otimes M_{\mathfrak{q}}\in {\cal B}_0$ Then,
$$
A\otimes {\cal Z}\cong B\otimes {\cal Z}.
$$
\end{thm}

\begin{proof}
The proof is almost identical to that of 11.7 of \cite{Lnclasn} with a few necessary modifications.
Note that $\Gamma$ induces an isomorphism
$$
\Gamma_{\mathfrak{p}}: {\rm Ell}(A\otimes M_{\mathfrak{p}})\to
{\rm Ell}(B\otimes M_{\mathfrak{p}}).
$$
 Since $A\otimes M_{\mathfrak{p}} \in {\cal B}_0$ and
$B\otimes M_{\mathfrak{p}}\in {\cal B}_0$, by Theorem \ref{IST1},  there is an isomorphism $\phi_{\mathfrak{p}}:
A\otimes M_{\mathfrak{p}}\to B\otimes M_{\mathfrak{p}}.$ Moreover,
(by the proof of \ref{IST1}),
$\phi_{\mathfrak{p}}$ carries $\Gamma_{\mathfrak{p}}.$  For
exactly the same reason, $\Gamma$ induces an isomorphism
$$
\Gamma_{\mathfrak{q}}:{\rm Ell}(A\otimes M_{\mathfrak{q}})\to
{\rm Ell}(B\otimes M_{\mathfrak{q}})
$$
and there is an isomorphism $\psi_{\mathfrak{q}}: A\otimes
M_{\mathfrak{q}}\to B\otimes M_{\mathfrak{q}}$ which induces
$\Gamma_{\mathfrak{q}}.$

 Put $\phi=\phi_{\mathfrak{p}}\otimes {\rm
id}_{M_{\mathfrak{q}}}: A\otimes Q\to B\otimes Q$ and
$\psi=\psi_{\mathfrak{q}}\otimes {\rm
id}_{{M_{\mathfrak{p}}}}: A\otimes Q\to B\otimes Q.$
Note that
$$
(\phi)_{*i}=(\psi)_{*i}\,\,{\rm (} i=0,1 {\rm )} \andeqn
\phi_T=\psi_T
$$
(they are induced by $\Gamma$). Note that $\phi_T$ and $\psi_T$
are affine homeomorphisms. Since $K_{*i}(B\otimes Q)$ is
divisible, we in fact have $[\phi]=[\psi]$ (in $KK(A\otimes Q,
B\otimes Q)$). It follows from \ref{L114} that there is an
automorphism $\bt: B\otimes Q\to B\otimes Q$ such that
$$
[\bt]=[{\rm id}_{B\otimes Q}]\,\,\,KK(B\otimes Q, B\otimes Q)
$$
{and} such that $\phi$ and $\bt\circ \psi$ are asymptotically unitarily
equivalent.  Since $K_1(B\otimes Q)$ is divisible,
$H_1(K_0(A\otimes Q), K_1(B\otimes Q))=K_1(B\otimes Q).$ It
follows that $\phi$ and $\bt\circ \psi$ are strongly
asymptotically unitarily equivalent. Note also in this case
$$
\bt_T=({\rm id}_{B\otimes Q})_T.
$$
Let $\imath: B\otimes M_{\mathfrak{q}}\to B\otimes Q$ defined by
$\imath(b)=b\otimes 1$ for $b\in B.$ We consider the pair
$\bt\circ \imath\circ {\psi}_{\mathfrak{q}}$ and $\imath \circ
{\psi}_{\mathfrak{q}}.$ By applying \ref{L115}, there exists an
automorphism $\af: B\otimes M_{\mathfrak{q}}\to B\otimes
M_{\mathfrak{q}}$ such that $\imath\circ \af\circ
\psi_{\mathfrak{q}}$ and $\bt\circ \imath\circ
\psi_{\mathfrak{q}}$ are asymptotically unitarily equivalent (in
$B\otimes Q$). So they are strongly asymptotically unitarily
equivalent. Moreover,
$$
[\af]=[{\rm id}_{B\otimes M_{{\mathfrak{q}}}}]\,\,\,{\rm in}\,\,\,
KK(B\otimes M_{\mathfrak{q}},B\otimes M_{\mathfrak{q}}).
$$

We will show that $\bt\circ \psi$ and $\af\circ
{\psi}_{\mathfrak{q}}\otimes {\rm id}_{M_{\mathfrak{p}}}$ are
strongly asymptotically unitarily equivalent. Define
$\bt_1=\bt\circ \imath\circ\psi_{\mathfrak{q}}\otimes {\rm
id}_{M_{\mathfrak{p}}}: B\otimes Q\otimes M_{\mathfrak{p}}\to
B\otimes Q\otimes M_{\mathfrak{p}}.$ Let $j: Q\to Q\otimes
M_{\mathfrak{p}}$ be defined by $j(b)=b\otimes 1.$ There is an
isomorphism $s: M_{\mathfrak{p}}\to M_{\mathfrak{p}}\otimes
M_{\mathfrak{p}}$ with $({\rm id}_{M_{\mathfrak{q}}}\otimes
s)_{*0}=j_{*0}.$ In this case $[{\rm id}_{M_{\mathfrak{q}}}\otimes s]=[j].$ Since
$K_1(M_{\mathfrak{p}})=0{,}$ {by} \ref{Tm72}, ${\rm
id}_{M_{\mathfrak{q}}}\otimes s$ is strongly asymptotically
%uniatrily
{{unitarily}} equivalent to $j.$ It follows that $\af\circ
\psi_{\mathfrak{q}}\otimes {\rm id}_{M_{\mathfrak{p}}}$ and
$\bt\circ \imath\circ \psi_{\mathfrak{q}}\otimes {\rm
id}_{M_{\mathfrak{p}}}$ are strongly asymptotically unitarily
equivalent. Consider the \SCA\, $C=\bt\circ \psi(1\otimes
M_{\mathfrak{p}})\otimes M_{\mathfrak{p}}\subset B\otimes Q\otimes
M_{\mathfrak{p}}.$ In $C,$  $\bt\circ \phi|_{1\otimes
M_{\mathfrak{p}}}$ and $j_0$ are strongly asymptotically unitarily
equivalent, where $j_0: M_{\mathfrak{p}}\to C$ is defined by $j_0(a)=1\otimes
a$ for all $a\in M_{\mathfrak{p}}.$ There exists a continuous path
of unitaries $\{v(t): t\in [0,\infty)\}\subset C$ such that
\beq\label{CM1-1}
\lim_{t\to\infty}{\rm ad}\, v(t) \circ\bt\circ \phi(1\otimes
a)=1\otimes a\tforal a\in M_{\mathfrak{p}}.
\eneq
It follows that $\bt\circ \psi$ and $\bt_1$ are strongly
asymptotically unitarily equivalent. Therefore $\bt\circ \psi$ and
$\af\circ \psi_{\mathfrak{q}}\otimes {\rm id}_{M_{\mathfrak{p}}}$
are strongly asymptotically unitarily equivalent. Finally, we
conclude that $\af\circ \psi_{\mathfrak{q}}\otimes {\rm
id}_{\mathfrak{p}}$ and $\phi$ are strongly asymptotically
unitarily equivalent. Note that $\af\circ \psi_{\mathfrak{q}}$ is
an isomorphism which induces $\Gamma_{\mathfrak{q}}.$

Let $\{u(t): t\in [0,1)\}$ be a continuous path of unitaries in
$B\otimes Q$ with $u(0)=1_{B\otimes Q}$ such that
$$
\lim_{t\to\infty}{\rm ad}\, u(t)\circ \phi(a)=\af\circ
\psi_{\mathfrak{q}}\otimes {\rm id}_{{M_{\mathfrak{p}}}}(a)\tforal
a\in A\otimes Q.
$$
One then obtains a unitary suspended isomorphism which lifts
$\Gamma$ along $Z_{p,q}$ (see \cite{Winter-Z}). It follows from Theorem
7.1 of \cite{Winter-Z} that $A\otimes {\cal Z}$ and $B\otimes {\cal Z}$
are isomorphic.
\end{proof}

\begin{df}\label{Class}
{ Denote by ${\cal N}_0$ the class of those unital  simple
\CA s $A$  in ${\cal N}$ for which $A\otimes
M_{\mathfrak{p}}\in {\cal N}\cap {\cal B}_0$ for any supernatural number
${\mathfrak{p}}$ of infinite type.

Of course ${\cal N}_0$  contains all unital simple amenable \CA s
in ${\cal B}_0$  which satisfy the UCT.
It contains all unital simple inductive limits of \CA s in ${\cal C}_0.$
}

\end{df}

\begin{cor}\label{CM1}
Let $A$ and $B$ be two \CA s in ${\cal N}_0.$ Then $A\otimes {\cal
Z}\cong B\otimes {\cal Z}$ if and only if ${\rm Ell}(A\otimes {\cal
Z})\cong {\rm Ell}(B\otimes {\cal Z}).$
\end{cor}

\begin{proof}
This follows from \ref{CMT1} immediately.
\end{proof}

%\section{The invariant range}
%Given  a countable weakly unperforated ordered  group $G_0$ with order unit $u$  (rationally Riesz?---No), a countable abelian group $G_1,$ a separable Choquet simplex $\Delta,$ a surjective affine continuous map $r: \Delta\to S_{[u]}(G_0)$
%where $S_{[u]}(G_0)$ is the state space of the ordered group $G_0$ with the order unit $u,$  such that $r(\tau)(g)=\tau(g)$ for all $\tau\in \Delta$ and $g\in G_0.$

%There is a unital simple C*-algebra $A$ in ${\cal N}_1$ such that
%$$
%Ell(A)=(G_0, (G_0)_+, u, G_1, \Delta, r).
%$$

% {\bf  Do not forget to prove, in \ref{RangT}, $A$ constructed actually in ${\cal B}_0!$-----Add a proof, please--L, already done--Gong}

%{\color{Green} In fact, \ref{MFTh} and \ref{MFTC} need the classification theorem, which will be proved in Part II. So, how about to move them into Part II.} {\bf {Let us move them---soon as the second part was revised.}}

\begin{thm}\label{MFTh}
Let $A$ and $B$ be two unital separable simple amenable ${\cal Z}$-stable \CA s which satisfy the UCT.
Suppose that $gTR(A\otimes Q)\le 1$ and $gTR(B\otimes Q)\le 1.$ Then
$A\cong B$ if and only if
\begin{equation*}
{\rm Ell}(A)\cong {\rm Ell}(B).
\end{equation*}
\end{thm}

\begin{proof}
It follows from \ref{B0=B1tU} that $A\otimes U, B\otimes U\in {\cal B}_0$ for any UHF-algebra $U$ of infinite type.
Therefore the theorem follows immediately from \ref{CM1}.
\end{proof}

\begin{cor}\label{MFTC}
Let $A$ and $B$ be two unital separable simple \CA s which satisfy the UCT.
Suppose that $gTR(A)\le 1$ and $gTR(B)\le 1.$ Then
$A\cong B$ if and only if
\begin{equation*}
{\rm Ell}(A)\cong {\rm Ell}(B).
\end{equation*}
\end{cor}

\begin{cor}\label{CM2}
Let $A$ and $B$ be two unital simple \CA s in ${\cal B}_1 \cap {\cal N}.$
Then $A\cong B$ if and only if
$$
{\rm Ell}(A)\cong {\rm Ell}(B).
$$
\end{cor}

\begin{proof}
It follows from \ref{Zstable} that $A\otimes {\cal Z}\cong A$ and
$B\otimes {\cal Z}\cong B.$ Thus the corollary follows from \ref{MFTh}.
%\ref{CM1}.
\end{proof}

{\bf Added in Proof}:  As this paper  was revised in 2015, it has been shown (see \cite{EGLN})  that all unital separable simple
\CA s $A$ with finite decomposition rank  which satisfy the UCT have $gTR(A\otimes U)\le 1$ for all
UHF-algebras of infinite type.  Therefore, by Theorem \ref{MFTh}, they are classified by the Elliott invariant.

{\small
\bibliographystyle{amsplain}
%\bibliography{operator_algebras}

\providecommand{\bysame}{\leavevmode\hbox to3em{\hrulefill}\thinspace}
\providecommand{\MR}{\relax\ifhmode\unskip\space\fi MR }
\providecommand{\MRhref}[2]{%
  \href{http://www.ams.org/mathscinet-getitem?mr=#1}{#2}
}
\providecommand{\href}[2]{#2}

}
ghgong@gmail.com,
hlin@uoregon.edu,
zniu@uwyo.edu.

\end{document}